\documentclass{report}%
\usepackage{amssymb}
\usepackage{amsfonts}
\usepackage{svsing2e}
\usepackage{amsmath}
\usepackage[a4paper,nomarginpar,verbose]{geometry}
\usepackage{graphicx}%
\newtheorem{theorem}{Theorem}

\newtheorem{conclusion}[theorem]{Conclusion}

\newtheorem{conjecture}[theorem]{Conjecture}
\newtheorem{corollary}[theorem]{Corollary}

\newtheorem{definition}[theorem]{Definition}
\newtheorem{example}[theorem]{Example}

\newtheorem{lemma}[theorem]{Lemma}
\newtheorem{notation}[theorem]{Notation}

\newtheorem{proposition}[theorem]{Proposition}
\newtheorem{remark}[theorem]{Remark}

\newtheorem{summary}[theorem]{Summary}
\newenvironment{proof}[1][Proof]{\noindent\textbf{#1.} }{\ \rule{0.5em}{0.5em}}
\makeatletter
\def\@crosshairs{\vbox to0pt{}}
\makeatother
\begin{document}
%
\begin{titlepage}%
%

\title
{A weight function theory of positive order basis function interpolants and smoothers}%
%

\author{\\Phillip Y. Williams\\\\
5/31 Moulden Court, Canberra, ACT 2617, Australia.\\\\
\texttt{phil.y.williams@bigpond.com.au}}%
%

\date{\today}%
%

\maketitle
%

\thanks
{Thanks to my Masters degree supervisors Dr Markus Hegland and Dr Steve Roberts
of the Centre for Mathematics and its Applications (CMA) at the Australian National University.}%
%

\end{titlepage}%
%

\begin{abstract}%

In this document I develop a weight function theory of positive order basis
function interpolants and smoothers.

In \textbf{Chapter 1} the basis functions and data spaces are defined directly
using weight functions. The data spaces are used to formulate the variational
problems which define the interpolants and smoothers discussed in later
chapters. The theory is illustrated using some standard examples of radial
basis functions and two classes of weight functions I will call the tensor
product extended B-splines and the central difference weight functions.

\textbf{Chapter 2} The goal of Chapter 2 is to derive `modified'
inverse-Fourier transform formulas for the basis functions and the data
functions (native spaces) and to use these formulas to obtain bounds for the
rates of increase of these functions and their derivatives near infinity.

\textbf{Chapter 3} shows how to prove functions are basis functions without
using the awkward space of test functions \$S\_\{0,n\}\$ which are infinitely
smooth functions of rapid decrease with several zero-valued derivatives at the
origin. Worked examples include several classes of well-known radial basis functions.

In \textbf{Chapter 4} we prove the existence and uniqueness of a solution to
the minimal seminorm interpolation problem. We then derive orders for the
pointwise convergence of the interpolant to its data function as the density
of the data increases.

In \textbf{Chapter 5} a well-known non-parametric variational smoothing
problem will be studied with special interest in the order of pointwise
convergence of the smoother to its data function. This smoothing problem is
the minimal norm interpolation problem stabilized by a smoothing coefficient.

In \textbf{Chapter 6} a non-parametric, scalable, variational smoothing
problem will be studied, again with special interest in its order of pointwise
convergence to its data function.

In \textbf{Chapter 7} we study the bounded linear functionals on the data spaces.%

\end{abstract}%
%

\tableofcontents

\section{Introduction}

\pagenumbering{arabic}In this document I develop a weight function theory of
\textbf{positive order} basis function interpolants and smoothers. Note that
the Appendix of this document contains a list of basic notation, definitions
and properties also used in this document.

This document had its genesis in the development of a \textit{scalable
algorithm} for \textit{Data Mining} applications. Data Mining is the
extraction of complex information from large databases, often having tens of
millions of records. Scalability means that the time of execution is linearly
dependent on the number of records processed and this is necessary for the
algorithms to have practical execution times. One approach is to develop
additive regression models and these require the approximation of large
numbers of data points by surfaces. Here one is concerned with approximating
data by surfaces of the form $y=f\left(  x\right)  $, where $x\in
\mathbb{R}^{d}$, $y\in\mathbb{R}$ and $d$ is any dimension. Smoothing
algorithms are one way of approximating surfaces and in particular we have
decided to use a class of non-parametric smoothers called basis function
smoothers, which solve a variational smoothing problem over a semi-Hilbert
space of continuous \textit{data functions} (\textit{often called native
spaces e.g. Schaback \cite{Schaback99}, Wendland \cite{Wendland05} etc.}) and
express the solution in terms of a single \textit{basis function}. Indeed,
several years ago I succeeded in developing a scalable smoothing
algorithm(unpublished) of the basis function type. This algorithm was derived
by approximating a minimal smoothing problem on a regular rectangular grid. In
this document we develop some theoretical tools to construct and analyze it.
This theoretical approach applies in any dimension but the smoothing algorithm
is only practical up to three dimensions. In higher dimensions the matrices
are too large to put into computer memory.

I started my Masters degree (supervised by Dr. Markus Hegland and Dr. Steve
Roberts at the ANU, Canberra, Australia) searching for a scalable basis
function smoothing algorithm and had the good fortune to devise such an
algorithm by approximating, on a regular grid, the convolution in Definition
38 of the space $J_{G}$ in Dyn's review article \cite{Dyn89}. For the case of
zero order basis functions this algorithm is analyzed in Chapter 4 of the
document Williams \cite{WilliamsZeroOrdSmthV4}, and the positive order basis
function case is studied in Chapter \ref{Ch_Approx_smth} of this document.

Dr. Hegland was particularly interested in using the tensor product
hat(triangular) function as a basis function. At about this time we had a
visit by the late Professor Will Light who showed me his paper
\cite{LightWayneX98Weight} which defined basis functions in terms of weight
functions using the Fourier transform. Light and Wayne's weight function
properties, given in Definition \ref{Def_LW_wt_fn_prop} below, were designed
for positive order basis functions which excluded the tensor product hat
function. They were designed for the well-known `classical' radial weight
functions. I therefore developed a version of his theory designed to generate
zero order basis functions, including both tensor product and radial types.
This theory was developed in Chapter 1 of Williams
\cite{WilliamsZeroOrdSmthV4} and requires that the basis functions have
Fourier transforms which can take zero values outside the origin since this is
a property of hat functions.

??? Need to add more discussion about the native space approach i.e. relate my
work to other's.\textbf{\medskip}

\textbf{Chapter by chapter:\medskip}

\textbf{Chapter }\ref{Ch_weight_fn_exten} \textbf{Extensions of Light's class
of weight functions} The goal of this chapter is to extend the theoretical
work of Light and Wayne in \cite{LightWayneX98Weight} to allow classes of
weight functions defined using integrals and which allow tensor product weight
functions. Two \textit{weight function} definitions are introduced, one
designed for the standard radial weight functions and one designed for tensor
product functions. These definitions allow for basis functions which have
Fourier transforms which are zero on a closed set of measure zero and are
analogous to the zero order definitions introduced in Williams
\cite{WilliamsZeroOrdSmthV4}. These definitions involve the positive integer
order parameter $\theta$ and the smoothness parameter $\kappa$ which can take
any non-negative real value.

The classes of tensor product weight functions are introduced which I call the
\textit{extended B-splines} and the \textit{central difference weight
functions}.

Both the semi-inner product \textit{data spaces} $\left(  X_{w}^{\theta
},\left\vert \cdot\right\vert _{w,\theta}\right)  $ and the \textit{basis
functions} $G$ are defined in terms of the weight function. We then prove some
completeness and smoothness properties of the data space as well as some
smoothness and positive definiteness properties of the basis function. In fact
the data space consists of continuous functions and the basis functions are
continuous, and the data space is used to define the variational interpolation
and smoothing problems of later chapters and the basis functions are used to
express their (unique) solutions.

The \textit{extended B-spline} basis functions are shown to be the
convolutions of hat (triangle) functions.\textbf{\medskip}

\textbf{Chapter }\ref{Ch_MoreBasisTheory} \textbf{More basis function and
semi-Hilbert data space theory} The goal of this chapter is to derive
`modified' inverse-Fourier transform formulas for basis functions and the data
functions (often native spaces) and to use these formulas to obtain bounds for
the rates of increase near infinity of these functions and their derivatives.
This will be done by proving a general inverse-Fourier transform formula for a
subspace of the distributions and then applying it to both the basis functions
and the data functions. Some other basis function properties are derived from
their weight function properties.\textbf{\medskip}

\textbf{Chapter }\ref{Ch_basis_no_So2n} \textbf{Calculating the weight
function from the basis function without using }$S_{\emptyset,2n}$ Shows how
to prove functions are basis functions without using the awkward test
functions $S_{\emptyset,2\theta}=\left\{  \phi\in S:D^{\alpha}\phi\left(
0\right)  =0,\text{ }\left\vert \alpha\right\vert <2\theta\right\}  $ where
$S$ is the $C^{\infty}$ space of rapidly decreasing functions given in
Definition \ref{Def_Distributions}. Worked examples include several classes of
well-known radial basis functions, the choice here following Dyn \cite{Dyn89}:
the thin-plate splines, the shifted thin-plate splines, the multiquadric and
inverse multiquadric functions and the Gaussian. In the last section I will
illustrate the method using a non-radial example: the fundamental solutions of
homogeneous elliptic differential operators.\textbf{\medskip}

\textbf{Chapter }\ref{Ch_Interpol} \textbf{The basis function interpolant and
its convergence to the data function} Studies the \textit{minimum seminorm
interpolation problem}. Topics include unisolvency and Lagrange polynomial
interpolation and related operators and matrices, the existence and uniqueness
of the basis function solution, a matrix equation for the solution and the
pointwise convergence of the interpolant to its data function on a bounded region.

We use a \textit{data function} $f_{d}\in X_{w}^{\theta}$ and
\textit{unisolvent independent data} $X=\left\{  x^{\left(  k\right)
}\right\}  _{k=1}^{N}$ contained in a bounded \textit{data region} $\Omega$ to
generate $N$ data points $\left(  x^{\left(  k\right)  },f_{d}\left(
x^{\left(  k\right)  }\right)  \right)  $ and the interpolation problem
requires us to minimize $\left\vert \cdot\right\vert _{w,\theta}$ over all
$f\in X_{w}^{\theta}$ which interpolate these data points. The seminorm
$\left\vert \cdot\right\vert _{w,\theta}$ is converted into a \textit{Light
norm} $\left\Vert \cdot\right\Vert _{w,\theta}$ and the problem is
reformulated as the \textit{minimum norm interpolation problem}. Geometric
Hilbert space theory using orthogonal projections shows the solution, denoted
$\mathcal{I}_{X}f_{d}$, is unique. Indeed, this solution is identical to that
of the seminorm problem and it lies in the finite dimensional basis function
space $W_{G,X}$ (Definition \ref{Def_Wgx}). It is called a \textit{basis
function solution}.

We are now interested in what happens to the interpolant as the data $X$
becomes denser in $\Omega$ when we use the \textit{cavity measure} of data
density%
\[
h_{X}=\sup_{x\in\Omega}\operatorname*{dist}\left(  x,X\right)  .
\]

Theorem \ref{Thm_interpol_converg} establishes an error estimate of order
$\eta=\min\left\{  \theta,\frac{1}{2}\left\lfloor 2\kappa\right\rfloor
\right\}  $ for the interpolant in the sense that there are positive constants
$k_{I}$ and $h$ such that%
\begin{equation}
\left\vert f_{d}\left(  x\right)  -\mathcal{I}_{X}f_{d}\left(  x\right)
\right\vert \leq\left\vert f_{d}\right\vert _{w,\theta}k_{I}\left(
h_{X}\right)  ^{\eta},\quad x\in\overline{\Omega},\text{ }f_{d}\in
X_{w}^{\theta},\label{a01}%
\end{equation}

when $h_{X}<h$.

Note that in the case of the thin-plate, the shifted thin-plate and the
Gaussian radial basis functions the Light and Wayne weight function definition
of \cite{LightWayneX98Weight} yields the same orders of convergence as
obtained above. However, slightly improved orders of convergence are obtained
later: see Corollary \ref{Cor_Thm_int_rx(x)_bnd_better_order} and examples
where an improvement of 1/2 is obtained for the shifted thin-plate
splines.\textbf{\medskip}

\textbf{Chapter }\ref{Ch_ExactSmth} \textbf{The Exact smoother and its
convergence to the data function} Studies the well-known \textit{Exact
smoother problem }(our terminology) which stabilizes the interpolant by adding
a smoothing coefficient to the seminorm functional. Topics include the
existence and uniqueness of the basis function solution i.e. the \textit{Exact
smoother}, a matrix equation for the solution and the convergence of the
solution to its data function.

Using the same data as the interpolation problem the Exact smoother problem
requires us to minimize the functional
\begin{equation}
J_{e}\left[  f\right]  =\rho\left\vert f\right\vert _{w,\theta}^{2}+\frac
{1}{N}\sum_{i=1}^{N}\left\vert f\left(  x^{(i)}\right)  -f_{d}\left(
x^{(i)}\right)  \right\vert ^{2},\label{a02}%
\end{equation}

for $f\in X_{w}^{\theta}$. \textbf{Note} that this smoothing functional is
different from that used in Narcowich, Ward and Wendland
\cite{NarcWardWend2004} which is%
\[
\rho\left\Vert f\right\Vert _{w,0}^{2}+\sum_{i=1}^{N}\left\vert f(x^{(i)}%
)-y_{i}\right\vert ^{2}.
\]

To obtain their error estimates you simply replace $\rho N$ by $\rho$ in our
error formulas.

Our proofs will be carried out within a Hilbert space framework by formulating
(Definition \ref{Def_ex_Hilbert_smoothing}) the smoothing functional in terms
of an inner product $\left(  u,v\right)  _{V}=\rho\left\langle u_{1}%
,v_{1}\right\rangle _{w,\theta}+\frac{1}{N}\left(  \widetilde{u}%
_{2},\widetilde{v}_{2}\right)  _{\mathbb{C}^{N}}$ on the product space
$V=X_{w}^{\theta}\otimes\mathbb{C}^{N}$ and the continuous operator
$\mathcal{L}_{X}:X_{w}^{\theta}\rightarrow V$ defined by $\mathcal{L}%
_{X}f=\left(  f,\widetilde{\mathcal{E}}_{X}f\right)  $ where $\widetilde
{\mathcal{E}}_{X}f=\left(  f\left(  x^{\left(  k\right)  }\right)  \right)  $.
We can now write $J_{e}\left[  f\right]  =\left\Vert \mathcal{L}_{X}f-\left(
0,\widetilde{\mathcal{E}}_{X}f\right)  \right\Vert _{V}^{2}$ and using the
technique of orthogonal projection it is shown in Theorem
\ref{Thm_smooth_Exact} that a unique solution to the smoothing problem exists
given by $s_{e}=\frac{1}{N}\left(  \mathcal{L}_{X}^{\ast}\mathcal{L}%
_{X}\right)  ^{-1}\widetilde{\mathcal{E}}_{X}^{\ast}\widetilde{\mathcal{E}%
}_{X}f_{d}$. Like the interpolant, the Exact smoother lies in the finite
dimensional basis function space $W_{G,X}$ and a matrix equation \ref{h68} is
derived for the coefficients of the basis functions.

To derive our estimates for the Exact smoother error $\left\vert s_{e}\left(
x\right)  -f_{d}\left(  x\right)  \right\vert $ on a bounded data region we
make use of the operators $\mathcal{L}_{X}$ and $\mathcal{L}_{X}^{\ast}$. Now
$\left\Vert \mathcal{L}_{X}f\right\Vert _{V}$ is an equivalent norm to
$\left\Vert f\right\Vert _{w,\theta}$ which implies that $X_{w}^{\theta}$ is
also a reproducing kernel Hilbert space under the inner product $\left(
\mathcal{L}_{X}f,\mathcal{L}_{X}f\right)  _{V}$. It then follows that there
exists a unique $R_{V,x}\in V$ such that $f\left(  x\right)  =\left(
\mathcal{L}_{X}f,R_{V,x}\right)  _{V}$ and $R_{x}=\mathcal{L}_{X}^{\ast
}R_{V,x}$. Thus, using the geometric properties of the orthogonal projection
\begin{align*}
\left\vert s_{e}\left(  x\right)  -f_{d}\left(  x\right)  \right\vert
=\left\vert \left(  \mathcal{L}_{X}\left(  s_{e}-f_{d}\right)  ,R_{V,x}%
\right)  _{V}\right\vert  &  \leq\left\Vert \mathcal{L}_{X}\left(  s_{e}%
-f_{d}\right)  \right\Vert _{V}\left\Vert R_{V,x}\right\Vert _{V}\\
&  \leq\left\Vert \mathcal{L}_{X}f_{d}-\left(  0,\widetilde{\mathcal{E}}%
_{X}f_{d}\right)  \right\Vert _{V}\left\Vert R_{V,x}\right\Vert _{V}\\
&  =\left\vert f_{d}\right\vert _{w,\theta}\sqrt{\rho}\left\Vert
R_{V,x}\right\Vert _{V}.
\end{align*}

We then estimate $\left\Vert R_{V,x}\right\Vert _{V}$ and show (eq'n
\ref{h104}) that there exist positive constants $k_{e},k_{e}^{\prime}$ and
$h^{\prime}$ such that%
\[
\left\vert s_{e}\left(  x\right)  -f_{d}\left(  x\right)  \right\vert
\leq\left\vert f_{d}\right\vert _{w,\theta}k_{e}\left(  k_{e}^{\prime}\left(
h_{X}\right)  ^{\eta}+\sqrt{N\rho}\right)  ,\quad x\in\overline{\Omega},\text{
}f_{d}\in X_{w}^{\theta},
\]

when $h_{X}<h^{\prime}$. Here $X$ has $N_{X}$ points and $k_{e}k_{e}^{\prime
}=k_{I}$. The latter equation means that when $\rho=0$ the smoother estimate
coincides with the interpolant error estimate. As with the interpolant
estimates, slightly higher orders of convergence can be obtained - see
Corollary \ref{Cor_Thm_int_rx(x)_bnd_better_order} and
examples.\textbf{\medskip}

\textbf{Chapter }\ref{Ch_Approx_smth} \textbf{The Approximate smoother and its
convergence to the Exact smoother and it's data function} Studies the
\textit{scalable\ Approximate smoother} (our suggested name). Topics include
the derivation of the \textit{Approximate smoother problem} and the existence,
uniqueness of its solution i.e. the Approximate smoother, a matrix equation
for the smoother, its scalability and the pointwise convergence of the
smoother to the Exact smoother and to its data function.

Using data generated by a data function $f_{d}\in X_{w}^{\theta}$ and a
\textbf{fixed} set $X^{\prime}\subset\mathbb{R}^{d}$ the Approximate smoother
problem requires us to minimize the Exact smoother functional \ref{a02} over
the finite dimensional space $W_{G,X^{\prime}}$. Using Hilbert space
techniques and then matrix techniques the smoother $s_{a}$ will be shown to
exist, be unique and satisfy the regular matrix equation \ref{s16} i.e. if
$s_{a}\left(  x\right)  =\sum\limits_{i=1}^{N^{\prime}}\alpha_{i}^{\prime
}G\left(  x-x_{i}^{\prime}\right)  +\sum\limits_{i=1}^{M}\beta_{i}^{\prime
}p_{i}\left(  x\right)  $ then%
\begin{equation}
\left(
\begin{array}
[c]{lll}%
\left(  2\pi\right)  ^{\frac{d}{2}}N\rho G_{X^{\prime},X^{\prime}%
}+G_{X^{\prime},X}G_{X,X^{\prime}} & G_{X^{\prime},X}P_{X} & P_{X^{\prime}}\\
P_{X}^{T}G_{X,X^{\prime}} & P_{X}^{T}P_{X} & O_{M}\\
P_{X^{\prime}}^{T} & O_{M} & O_{M}%
\end{array}
\right)  \left(
\begin{array}
[c]{c}%
\alpha^{\prime}\\
\beta^{\prime}\\
\gamma^{\prime}%
\end{array}
\right)  =\left(
\begin{array}
[c]{l}%
G_{X^{\prime},X}\\
P_{X}^{T}\\
O_{M,N}%
\end{array}
\right)  y,\label{a03}%
\end{equation}

where: $G_{X^{\prime},X^{\prime}}=G\left(  x_{i}^{\prime}-x_{j}^{\prime
}\right)  $, $G_{X^{\prime},X}=G\left(  x_{i}^{\prime}-x^{\left(  j\right)
}\right)  $, $G_{X^{\prime},X}=G\left(  x^{\left(  i\right)  }-x_{j}^{\prime
}\right)  $, $P_{X}=\left(  p_{j}\left(  x^{\left(  i\right)  }\right)
\right)  $, $P_{X^{\prime}}=\left(  p_{j}\left(  x_{i}^{\prime}\right)
\right)  $ and $\left\{  p_{i}\right\}  $ is a basis for the polynomials
$P_{\theta-1}$ of degree less than $\theta$.

The \textit{Approximate smoother matrix is} $N+N^{\prime}$ square and this
dependence on the number of data points $N$ implies \textit{computational
scalability} i.e. the computational effort required to construct and solve the
matrix equation depends linearly on the number of data points $N$.

For a bounded data region $\Omega$ we first derive some uniform, pointwise
convergence results which do not involve orders of convergence e.g. in
Corollary \ref{Cor_Thm_Jsd[sig(Zk)]toJsd(sig(Z))_2} it is shown that as the
points in $X^{\prime}$ get closer to those in $X$ the Approximate smoother
converges uniformly to the Exact smoother on $\overline{\Omega}$.

Orders of pointwise convergence are derived for the convergence of the
Approximate smoother to the data function e.g. estimate \ref{s46} says that
for some constants $B,C,r>0$,%
\begin{align*}
\left\vert s_{e}\left(  x\right)  -s_{a}\left(  x\right)  \right\vert \leq &
\left\vert f_{d}\right\vert _{w,\theta}\left(  1+\left(  1+K_{\Omega,\theta
}^{\prime}\right)  \sqrt{c_{G,\eta,\sigma}}\frac{\left(  c_{\Omega,\theta
}h_{X^{\prime}}\right)  ^{\eta_{G}}}{\sqrt{\rho}}\right)  \times\\
&  \times\left(  1+K_{\Omega,\theta}^{\prime}\right)  \left(  \sqrt
{c_{G,\eta,\sigma}}\left(  c_{\Omega,\theta}h_{X}\right)  ^{\eta_{G}}%
+\sqrt{N_{X}\rho}\right)  ,\quad x\in\overline{\Omega},\text{ }f_{d}\in
X_{w}^{\theta},
\end{align*}

when $h_{X^{\prime}},h_{X}<r$. Here $h_{X}=\sup\limits_{\omega\in\Omega
}\operatorname*{dist}\left(  \omega,X\right)  $ and $h_{X^{\prime}}%
=\sup\limits_{\omega\in\Omega}\operatorname*{dist}\left(  \omega,X^{\prime
}\right)  $ measure the density of the point sets $X$ and $X^{\prime}$.

These estimates are then added to the Exact smoother error formulas of Chapter
\ref{Ch_ExactSmth} to obtain error estimates for the Approximate smoother e.g.
estimate \ref{s25}.\textbf{\medskip}

\textbf{Chapter }\ref{Ch_bnd_lin_fnal_Xwth} \textbf{The bounded linear
functionals on the data space }$X_{w}^{\theta}$ In this chapter we study the
spaces $X_{1/w}^{-\theta}$ where $\theta\geq1$. Note that we use the notation
$X_{1/w}^{-\theta}$ instead of $X_{w}^{-\theta}$ so we can distinguish between
the space $X_{1/w}^{0}$ and the space $X_{w}^{0}$. The spaces $X_{1/w}%
^{-\theta}$ are thus not the negative order versions of $X_{w}^{\theta}$ but
are actually isometrically isomorphic to the spaces of \textbf{bounded linear
functionals} on $X_{w}^{\theta}$, denoted $\left(  X_{w}^{\theta}\right)
^{\prime}$. The properties of the spaces $X_{1/w}^{-\theta}$ are established
by constructing inverse isometric mappings between $X_{1/w}^{-\theta}$ and
$X_{w}^{\theta}$.

\begin{notation}
\label{Not_distribution}Regarding notation for distributions, in
\cite{LightWayneX98Weight} Light and Wayne used an operator notation which I
feel significantly reduced the readability of an already difficult paper,
especially replacing $\left(  ix\right)  ^{\gamma}$ by $V_{\gamma}\left(
x\right)  $ and replacing the exponential $e^{ixy}$ by $e_{x}\left(  y\right)
$. For me this turned old mathematical friends into strangers. So I have
settled on the notion of an \textbf{action} \textbf{variable} e.g. in the
equations%
\[
\left[  \xi^{\alpha}f,\phi\right]  =\left[  f,\xi^{\alpha}\phi\right]
,\text{\quad}\left[  e^{i\xi}f,\phi\right]  =\left[  f,e^{i\xi}\phi\right]  .
\]
involving a distribution $f$ and a test function $\phi$, $\xi$ will be called
the action variable. I will also use the dot notation e.g.%
\[
\left[  \left\vert \cdot\right\vert ^{2}f,\phi\right]  =\left[  f,\left\vert
\cdot\right\vert ^{2}\phi\right]  ,\text{\quad}\left[  f\left(  \cdot
-x\right)  ,\phi\right]  =\left[  f,\phi\left(  \cdot+x\right)  \right]  .
\]

My aim is to retain as much as possible of \textbf{the form} of `conventional'
mathematical expressions.
\end{notation}

\chapter{Extensions of Light's class of weight
functions\label{Ch_weight_fn_exten}}

\section{Introduction\label{Sect_Introd}}

The goal of this chapter is to extend the theoretical work of Light and Wayne
in \cite{LightWayneX98Weight} to allow classes of weight functions analogous
to the zero order weight functions developed in Section 1.2 of Williams
\cite{WilliamsZeroOrdSmthV4} i.e. weight functions which can generate positive
order tensor product basis functions and the radial basis functions which have
Fourier transforms with zeros outside the origin. A class of weight functions
which I call the extended B-splines of positive order is used to illustrate
the weight function properties and the basis function theory. Besides the
basis function theory I have shown Light and Wayne's semi-Hilbert function
spaces (often called native spaces) are still valid for the expanded weight
function class but I give a different proof of completeness that uses single
mappings in a manner analogous to Sobolev space theory.

The material of this chapter can be summarized as:

\begin{enumerate}
\item Define the positive order weight function properties and present
equivalent properties and relationships between them. Introduce the
\textit{extended B-spline weight functions }which are tensor product functions
with bounded support.

\item Define the semi-inner product space of continuous \textit{data functions
(native spaces)} and prove its completeness and smoothness properties.

\item Define the basis distributions of positive order and prove continuity
and positive definiteness properties. Derive convolution formulas for the
extended B-spline basis functions.
\end{enumerate}

The theory of this document lays the foundations for the study of the basis
function interpolation and smoothing problems introduced in later chapters.

\section{An extended class of positive order weight functions
\label{Sect_gen_wt_funcs}}

In this section we introduce our extended class of positive order weight functions.

\subsection{Zero order weight functions and basis
functions\label{SbSect_ZeroOrdWeightBasis}}

In the Chapter 1 of \cite{WilliamsZeroOrdSmthV4} Williams developed a theory
of zero order weight functions and basis functions. Here we shall summarize
the relevant parts of this theory which only needs simple $L^{1}$ Fourier
theory e.g. Section 2.2 of Petersen \cite{Petersen83}. However, here we will
allow non-integer values of the smoothness parameter $\kappa$ to match the
positive order weight function definitions below.

\begin{definition}
\label{Def_wt_fn zero}\textbf{Zero order weight functions with smoothness
parameter} $\kappa$

A weight function's properties are defined with reference to a\textbf{\ }set
$\mathcal{A}\subset\mathbb{R}^{d}$ which is closed and has measure zero. A
weight function $w$ is a mapping $w:\mathbb{R}^{d}\rightarrow\mathbb{R}$ which
has the properties:

\begin{description}
\item[W01] There exists a closed set $\mathcal{A}$ with measure zero such that
$w$ is continuous and positive outside $\mathcal{A}$ i.e. $w\in C^{\left(
0\right)  }\left(  \mathbb{R}^{d}\setminus\mathcal{A}\right)  $ and $w>0$ on
$\mathbb{R}^{d}\setminus\mathcal{A}$.

\item[W02] For some (possibly non-integer) $\kappa\in\mathbb{R}^{1}$,
$\kappa\geq0$,%
\[
\int\dfrac{\left\vert x\right\vert ^{2s}}{w\left(  x\right)  }dx<\infty
,\quad0\leq s\leq\kappa.
\]

\item[W03] For some $\kappa\in\mathbb{R}^{d}$, $\kappa\geq0$,%
\begin{equation}
\int\dfrac{x^{2\lambda}}{w\left(  x\right)  }dx<\infty,\quad0\leq\lambda
\leq\kappa.\label{p62}%
\end{equation}

\end{description}
\end{definition}

The \textit{basis function of order zero} is defined by%
\begin{equation}
G=\left(  \frac{1}{w}\right)  ^{\vee}.\label{p63}%
\end{equation}

It follows directly from Corollary 2.12 of Petersen \cite{Petersen83} that
\begin{equation}
G\in C_{B}^{\left(  \left\lfloor 2\kappa\right\rfloor \right)  },\label{p73}%
\end{equation}

and so $\kappa$ can be called the weight function \textit{smoothness
parameter}.

\subsection{Motivation for the extended weight function class properties
\label{SbSect_motiv_weight_fn}}

I will now state the Light and Wayne weight function properties extracted from
Section 3 of Light and Wayne \cite{LightWayneX98Weight}. To this list
\textbf{I have added property B3.5} - \textbf{not} \textbf{defined by Light
and Wayne}. The space $X$ of functions mentioned in this definition were used
to define the variational interpolant. The basis function is used to construct
the interpolant.

\begin{definition}
\label{Def_LW_wt_fn_prop}\textbf{Light and Wayne's weight function properties
}\cite{LightWayneX98Weight}

A \textbf{weight function} $w$ is a mapping $w:\mathbb{R}^{d}\rightarrow
\mathbb{R}$ which has properties A3.1 and A3.2 below, as well as combinations
of the other properties. Properties A3.1 and A3.2 are used to define the
semi-inner product space $X$ of distributions.

\begin{description}
\item[A3.1] $w\in C^{\left(  0\right)  }(\mathbb{R}^{d}\setminus0)$.

\item[A3.2] $w(x)>0$ on $\mathbb{R}^{d}\setminus0$.\medskip

Properties A3.3 and A3.4 are used to prove the completeness of\ the $X$ space
and to allow the definition of a basis distribution.

\item[A3.3] $1/w\in L_{loc}^{1}$.

\item[A3.4] There exist $\mu$, $R>0$ and $C_{R}>0$ such that
\[
\dfrac{1}{w(x)}\leq C_{R}|x|^{-2\mu},\quad|x|\geq R.
\]

\item[B3.5] We say that $w$ has property B3.5 for parameters $\mu$ and
\textbf{order} $\theta$ if $w$ satisfies property A3.4 for some $\mu$
satisfying $\mu+\theta>d/2$.\smallskip

\textbf{I have constructed Property B3.5 }from Light and Wayne\textbf{\ }%
\cite{LightWayneX98Weight} and it allows the space $X$ to be embedded in the
continuous functions $C^{\left(  k\right)  }$ where $k$ is the largest integer
such that $k<\mu+\theta-d/2$ (Theorem 2.18), and the basis distributions of
order $\theta$ to be continuous functions $C^{\left(  j\right)  }$ where $j $
is the largest integer such that $j<2\left(  \mu+\theta\right)  -d$ (Theorem
3.14). Property B3.5 allows these weight function properties to be more easily
compared with the `extended' properties given below.\medskip
\end{description}
\end{definition}

Light and Wayne's weight function properties were designed to generate the
positive order radial basis functions and the aim of this section is to extend
Light and Wayne's weight function properties so that they can generate
positive order tensor product basis functions and also to formulate the weight
function properties in terms of integrals.

In Subsection 1.2.9 of Williams \cite{WilliamsZeroOrdSmthV4} a class of weight
functions, which I called the extended B-spline weight functions, was used to
illustrate the weight function properties and the basis function theory. In
this document a class of weight functions which I will call the
\textit{extended B-spline weight functions of positive order} is used to
illustrate the weight function properties and the basis function theory.

In more detail, the extended weight function properties given below in
Definition \ref{Def_extend_wt_fn} were formulated by taking the following
considerations into account:\medskip

\textbf{1}. I was interested in properties expressed in terms of integrals.
Also, these properties should be scalable.\medskip

\textbf{2}. When a function $w_{1}$ defined on $\mathbb{R}^{1}$ has the weight
function properties every tensor product $\prod\limits_{i=1}^{d}w_{1}\left(
x_{i}\right)  $ should also have these properties but perhaps for different
parameter values. Now suppose $w$ is a tensor product function of the
1-dimensional function $w_{1}$ i.e.%
\[
w\left(  x\right)  =\prod\limits_{i=1}^{d}w_{1}\left(  x_{i}\right)  ,\quad
x\in\mathbb{R}^{d},
\]

and suppose that $w_{1}$ has Light and Wayne properties A3.1 and A3.2. Now we
want $w$ to be a weight function i.e. satisfy A3.1 and A3.2. In general this
is not possible since if $w_{1}$ is discontinuous at zero then $w$ is
discontinuous on the closed set $\bigcup\limits_{i=1}^{d}\left\{
x:x_{i}=0\right\}  $. This motivates extended weight function property W1 of
Definition \ref{Def_extend_wt_fn} i.e. there exists a closed set
$\mathcal{A}\subset\mathbb{R}^{d}$ of measure zero such that $w\in C^{\left(
0\right)  }\left(  \mathbb{R}^{d}\setminus\mathcal{A}\right)  $ and $w>0$ on
$\mathbb{R}^{d}\setminus\mathcal{A}$.\medskip

\textbf{3}. I want to generalize the zero order extended B-spline weight
functions studied in Section 1.2 of Williams \cite{WilliamsZeroOrdSmthV4} to
the positive order case. This will mean allowing basis functions with Fourier
transforms which have zeros outside the origin i.e. weight function with poles
(discontinuities) outside the origin.\textbf{\ }Indeed, the one-dimensional
hat function is defined by
\begin{equation}
\Lambda\left(  x\right)  =\left\{
\begin{array}
[c]{ll}%
1-\left\vert x\right\vert , & \left\vert x\right\vert \leq1,\\
0, & \left\vert x\right\vert >1,
\end{array}
\right.  \text{\quad}x\in\mathbb{R}^{1},\label{p04}%
\end{equation}

and in higher dimensions it is defined as the tensor product
\begin{equation}
\Lambda\left(  x\right)  =\prod_{i=1}^{d}\Lambda\left(  x_{i}\right)
,\text{\quad}x\in\mathbb{R}^{d}.\label{p06}%
\end{equation}

It is well known that
\begin{equation}
\widehat{\Lambda}\left(  \xi\right)  =\prod\limits_{i=1}^{d}\widehat{\Lambda
}\left(  \xi_{i}\right)  =\left(  2\pi\right)  ^{-d/2}\prod\limits_{i=1}%
^{d}\left(  \dfrac{\sin\left(  \xi_{i}/2\right)  }{\xi_{i}/2}\right)
^{2},\text{\quad}\xi\in\mathbb{R}^{d}.\label{p07}%
\end{equation}

The zero order extended B-spline weight functions were created by generalizing
the right side of \ref{p07} to the two-parameter class of functions
\begin{equation}
\frac{1}{w\left(  \xi\right)  }=\prod\limits_{i=1}^{d}\dfrac{\sin^{2n}\xi_{i}%
}{\xi_{i}^{2k}},\quad n,k\geq1,\label{p30}%
\end{equation}

and then restricting the choice of parameters $l$ and $n$. The positive order
weight functions will involve a different choice of these parameters. Here
$1/w$ has zeros outside the origin on a closed set of measure zero.\medskip

\textbf{4}. Property A3.4 is $1/w\in L_{loc}^{1}$ and so is already defined by
an integral and this becomes extended weight function property W2.1 of
Definition \ref{Def_extend_wt_fn}. Property A4.4 states: there exists $\mu
\in\mathbb{R}^{1}$ and $R,C_{R}>0$ such that $\dfrac{1}{w(x)}\leq
C_{R}|x|^{-2\mu}$ for $|x|\geq R$. Now Lemma \ref{Thm_basis_smth_W3.2_r3_pos}
requires that $1/w\in S^{\prime}$ i.e. $1/w$ be a tempered distribution
(Appendix \ref{Sect_tempered_distrib}) and from part 2 of Appendix
\ref{SbSect_property_S'} we see that A3.4 and A3.4 imply that $1/w$ is a
regular, tempered distribution. Further it is clear that extended property
W2.2 i.e. $\int\limits_{\left\vert \cdot\right\vert \geq r_{2}}\frac
{1}{w\left\vert \cdot\right\vert ^{2\sigma}}<\infty$ for some $\sigma>0$ and
$r_{2}>0$, combined with property W2.1 also implies that $1/w\in S^{\prime}$.
Happily, property W2.2 is a good substitute for property A3.4.\medskip

\textbf{5}. Some clarification of the relationship between property B3.5 and
the new properties W3.1 and W3.2 will be provided by Theorem
\ref{Thm_L&W_implies_extend} in Subsection
\ref{SbSect_extend_propert_LightWayne}.

\subsection{The extended class of positive order weight functions
\label{SbSect_extend_wt_fn_propert}}

The weight function class of Light and Wayne will be extended as follows:

\begin{definition}
\label{Def_extend_wt_fn}\textbf{The weight function properties}

The properties are defined with reference to a set $\mathcal{A}\subset
\mathbb{R}^{d}$ which is a closed set of measure zero. The weight function is
a mapping $w:\mathbb{R}^{d}\rightarrow\mathbb{R}$ which has at least property W1:

\begin{description}
\item[W1] $w\in C^{\left(  0\right)  }\left(  \mathbb{R}^{d}\setminus
\mathcal{A}\right)  $ and $w>0$ on $\mathbb{R}^{d}\setminus\mathcal{A}$.

Property W1 will be used to define the semi-inner product distribution space
$X_{w}^{\theta}$ for positive integer \textbf{order }$\theta$. This property
is identical to the zero order weight function property W01 (Definition
\ref{Def_wt_fn zero} above).\medskip

\item[W2] Property W2 is satisfied if the following two sub-properties are
satisfied:\medskip

\item[W2.1] $1/w\in L_{loc}^{1}$.

\item[W2.2] $\int\limits_{\left\vert \cdot\right\vert \geq r_{2}}\frac
{1}{w\left\vert \cdot\right\vert ^{2\sigma}}<\infty$ for some $\sigma>0$ and
some $r_{2}>0$.

Property W2 is used to prove completeness and $C^{\infty}$ density results
regarding the $X_{w}^{\theta}$ spaces, as well as allowing the definition of
the basis distributions. Condition W2 is satisfied iff the function $\lambda$
defined by:%
\begin{equation}
\lambda\left(  x\right)  =\left\{
\begin{array}
[c]{ll}%
0, & \left\vert x\right\vert \leq r_{2},\\
\sigma, & r_{2}<\left\vert x\right\vert ,
\end{array}
\right. \label{p70}%
\end{equation}

satisfies $\int\dfrac{dx}{w\left(  x\right)  \left\vert x\right\vert
^{2\lambda\left(  x\right)  }}<\infty$. The function $\lambda$ is useful for
the concise expression of integral inequalities involving weight functions
e.g. Lemma \ref{Lem_functnal_phi_sq}.\bigskip

Now to introduce properties W3.1 and W3.2. These will allow the $X_{w}%
^{\theta}$ spaces to be continuously embedded in the continuous functions and
the basis distributions to be continuous functions. Suppose $\theta\geq1$ is a
positive integer:\medskip

\item[W3.1] $w$ has property W3.1 for order $\theta$ and some parameter
\fbox{$\kappa\in\mathbb{R}^{d}$}, $\kappa\geq0$ if there exists a multi-index
$\alpha$ such that $\left\vert \alpha\right\vert =\theta$ and%
\begin{equation}
\int\dfrac{x^{2\lambda}}{w\left(  x\right)  x^{2\alpha}}dx<\infty,\quad
0\leq\lambda\leq\kappa.\label{p36}%
\end{equation}

Here $x^{2\lambda}:=\left(  x_{1}^{2}\right)  ^{\lambda_{1}}\times\cdots
\times\left(  x_{d}^{2}\right)  ^{\lambda_{d}}=\left\vert x_{1}\right\vert
^{2\lambda_{1}}\times\cdots\times\left\vert x_{d}\right\vert ^{2\lambda_{d}%
}=\left(  \left\vert x_{k}\right\vert \right)  ^{2\lambda}=x_{+}^{2\lambda}$.

\item[W3.1*] $w$ has property W3.1* for order $\theta$ and $\kappa$ if it has
property W3.1 for $\theta$ and $\kappa$ and\textbf{\ all} $\alpha$ such that
$\left\vert \alpha\right\vert =\theta$.

\item[W3.2] $w$ has property W3.2 for order $\theta$ and some \fbox{$\kappa
\in\mathbb{R}^{1}$}, $\kappa\geq0$ if there exists some $r_{3}\geq0$ such
that
\[
\int\limits_{\left\vert x\right\vert \geq r_{3}}\dfrac{\left\vert x\right\vert
^{2\kappa}}{w\left(  x\right)  \left\vert x\right\vert ^{2\theta}}dx<\infty.
\]

\item[W3.3] $w$ has property W3.3 for order $\theta$ and some \fbox{$\kappa
\in\mathbb{R}^{d}$}, $\kappa\geq0$ if
\[
\int\dfrac{x^{2\lambda}}{w\left(  x\right)  \left\vert x\right\vert ^{2\theta
}}dx<\infty,\quad0\leq\lambda\leq\kappa.
\]

\item[W3] $w$ has property W3 if it has either property W3.1 or W3.2 or W3.3.
\end{description}
\end{definition}

\begin{remark}
\label{Rem_Def_extend_wt_fn_1}\ 

\begin{enumerate}
\item This definition permits the weight function to have
\textbf{discontinuities outside the origin} on a closed set $\mathcal{A}$ of
measure zero. This allows the Fourier transform of the \textbf{basis function
to have zeros outside the origin}. For example, the extended B-spline tensor
product weight functions introduced in Subsection
\ref{SbSect_ext_splin_wt_fn_W3.1*} will generate basis functions which have
zeros outside the origin on a set of measure zero.

\item The extended properties have been framed with a desire to use
\textbf{integrals to formulate properties} which generalize Light and Wayne's properties.

\item Like the Light and Wayne properties A3.3 and A3.4, the extended
properties W2.1 and W2.2 ensure that $1/w$\textbf{\ is a regular tempered
distribution} (or generalized function of slow growth) as per part 2 of
Appendix \ref{SbSect_property_S'}.

\item Property W3.1 is designed for use with \textbf{tensor product weight
functions} such as the tensor product extended B-spline weight functions
introduced in Subsection \ref{SbSect_ext_splin_wt_fn_W3.1*}.

With regard to point 2 of Subsection \ref{SbSect_motiv_weight_fn} above, when
a weight function $w_{1}$ defined on $\mathbb{R}^{1}$ has any of the weight
function properties W2.1, W2.2, W3.1, W3.1* etc. then every homogeneous tensor
product $\prod\limits_{i=1}^{d}w_{1}\left(  x_{i}\right)  $ also has the same
properties but perhaps for different parameter values.

\item Property W3.2 is designed to \textbf{handle radial weight functions}.

\item In one dimension properties W3.1 and W3.1* are identical but in higher
dimensions property W3.1* may be stronger than property W3.1. However,
property W3.1* is easier to handle analytically but W3.1 will give ?? better results?.??

\item The expressions used to define properties W3.1 and W3.2 may look a
little strange because we have not combined the \textbf{exponents from the
denominator and the numerator} but these expressions allow direct comparison
of the two definitions and point to important similarities between them. Also,
if a weight function has property W3 for order $\theta$ then we will only
define the function spaces $X_{w}^{\theta}$ and the basis functions for this order.

\item I turns out that $X_{w}^{\theta}\subset C_{B}^{\left(  \left\lfloor
\kappa\right\rfloor \right)  }$ and $G\in C_{B}^{\left(  \left\lfloor
2\kappa\right\rfloor \right)  }$ so we shall sometimes refer to $\kappa$ as
the \textbf{smoothness parameter}.

\item Observe that all the \textbf{weight function properties are scalable}
w.r.t. both the dependent and independent variables.

\item If we set $\theta=0$ in the definition of the extended properties we
obtain a set of \textbf{zero order weight function} properties: property W1
for the positive order weight functions is identical to the zero order weight
function property W01 of Definition \ref{Def_wt_fn zero}. Property W3.3
becomes equivalent to zero order property W03. Properties W3.2 and W2 are
equivalent to zero order property W02. When $\theta=0$ positive order property
W2 stays the same and ensures that $1/w$ is a tempered distribution. It allows
the definition of a \textbf{zero order basis distribution} which corresponds
to the positive order distributions of Definition \ref{Def_basis_distrib}.

\item \textbf{Property W3.3} is used for constructing positive order basis
functions from zero order basis functions. The basis functions have similar
differentiability to those associated with W3.1. See Theorem
\ref{Thm_pos_wt_from_zero_wt}.
\end{enumerate}
\end{remark}

\subsection{Relationships between the extended properties and Light and
Wayne's weight function properties\label{SbSect_extend_propert_LightWayne}}

The next theorem justifies the use of the terminology \textit{extended weight
function properties}.

\begin{theorem}
\label{Thm_L&W_implies_extend}Suppose for parameters $\mu\in\mathbb{R}^{1}$
and $R>0$ a weight function $w$ satisfies properties A3.1 to A3.4 of Light and
Wayne's Definition \ref{Def_LW_wt_fn_prop}. Then:

\begin{enumerate}
\item $w$ satisfies the extended weight function properties W1 and W2 of
Definition \ref{Def_extend_wt_fn} with $\mathcal{A}=\left\{  0\right\}  $,
when $\sigma>d-2\mu$, $\sigma\geq0$ and $r_{2}=R$.

\item Suppose $w$ also satisfies property B3.5 for $\mu$ and order $\theta$
i.e. $\mu+\theta>d/2$. Then the weight function satisfies property W3.2 of
Definition \ref{Def_extend_wt_fn} for order $\theta$ and all $\kappa$ such
that $0\leq\kappa<\mu+\theta-d/2$.
\end{enumerate}
\end{theorem}

\begin{proof}
\textbf{Part 1} Clearly $\mathcal{A}=\left\{  0\right\}  $ is a closed set of
measure zero. It is also clear that $w$ has properties W1 and W2. Now to prove
property W2.2. Property A3.4 requires that for some $R>0$
\[
\dfrac{1}{w(x)}\leq C_{R}|x|^{-2\mu},\quad when\text{ }|x|\geq R.
\]

Regarding the definition of property W2.2, choose $\sigma>d-2\mu$, $\sigma
\geq0$ and $r=R$. Then
\[
\int\limits_{\left\vert x\right\vert \geq r}\dfrac{dx}{\left\vert x\right\vert
^{\sigma}w\left(  x\right)  }\leq C_{R}\int\limits_{\left\vert x\right\vert
\geq r}\dfrac{dx}{\left\vert x\right\vert ^{\sigma+2\mu}}<\infty,
\]

since $\sigma+2\mu>d$. Thus $w$ has property W2.2.\medskip

\textbf{Part 2} Choose any $r^{\prime}>0$. Then, since $w$ has property A3.4
\[
\int\limits_{\left\vert x\right\vert \geq r^{\prime}}\dfrac{\left\vert
x\right\vert ^{2\kappa}dx}{\left\vert x\right\vert ^{2\theta}w\left(
x\right)  }=\int\limits_{\left\vert x\right\vert \geq r^{\prime}}\dfrac
{dx}{\left\vert x\right\vert ^{2\left(  \theta-\kappa\right)  }w\left(
x\right)  }\leq C_{R}\int\limits_{\left\vert x\right\vert \geq r^{\prime}%
}\frac{\left\vert x\right\vert ^{-2\mu}dx}{\left\vert x\right\vert ^{2\left(
\theta-\kappa\right)  }}=C_{R}\int\limits_{\left\vert x\right\vert \geq
r^{\prime}}\frac{dx}{\left\vert x\right\vert ^{2\left(  \mu+\theta
-\kappa\right)  }}.
\]

But the condition $0\leq\kappa<\mu+\theta-d/2$ implies $\mu+\theta-\kappa>d/2$
and so the last integral exists.
\end{proof}

\subsection{Equivalent criteria for property W3.2}

\begin{theorem}
\label{Thm_equiv_W3.2}The following criteria are \textbf{equivalent} to weight
function property W3.2\textbf{\ }for order $\theta$ and $\kappa$:

\begin{enumerate}
\item $\int\limits_{\left\vert x\right\vert \geq r_{3}}\dfrac{\left\vert
x\right\vert ^{2\left(  \kappa-\left\lfloor \kappa\right\rfloor \right)
}x^{2\beta}}{\left\vert x\right\vert ^{2\theta}w\left(  x\right)  }%
dx<\infty,\quad\left\vert \beta\right\vert =\left\lfloor \kappa\right\rfloor
.$

\item $\int\limits_{\left\vert x\right\vert \geq r_{3}}\dfrac{x_{k}^{2\left(
\kappa-\left\lfloor \kappa\right\rfloor \right)  }x^{2\beta}}{\left\vert
x\right\vert ^{2\theta}w\left(  x\right)  }dx<\infty,\quad\left\{
\begin{array}
[c]{l}%
\left\vert \beta\right\vert =\left\lfloor \kappa\right\rfloor ,\\
k=1,\ldots,d.
\end{array}
\right.  $
\end{enumerate}
\end{theorem}

\begin{proof}
Part 1 is proved by using the identity%
\[
\left\vert x\right\vert ^{2\kappa}=\left\vert x\right\vert ^{2\left(
\kappa-\left\lfloor \kappa\right\rfloor \right)  }\left\vert x\right\vert
^{2\left\lfloor \kappa\right\rfloor }=\left\vert x\right\vert ^{2\left(
\kappa-\left\lfloor \kappa\right\rfloor \right)  }\sum\limits_{\left\vert
\beta\right\vert =\left\lfloor \kappa\right\rfloor }\frac{\left\lfloor
\kappa\right\rfloor !}{\beta!}x^{2\beta}=\sum\limits_{\left\vert
\beta\right\vert =\left\lfloor \kappa\right\rfloor }\frac{\left\lfloor
\kappa\right\rfloor !}{\beta!}\left\vert x\right\vert ^{2\left(
\kappa-\left\lfloor \kappa\right\rfloor \right)  }x^{2\beta},
\]

and part 2 is proved by applying the inequalities,

$\frac{1}{d}\sum\limits_{k=1}^{d}x_{k}^{2p}\leq\left\vert x\right\vert
^{2p}\leq\sum\limits_{k=1}^{d}x_{k}^{2p},$\quad$0\leq p\leq1$, with
$p=\kappa-\left\lfloor \kappa\right\rfloor $ to the criterion of part 1.
\end{proof}

\subsection{Tensor product weight functions and weight function property W3.1}

These results are \textbf{much simpler} than results for property W3.1*.

\begin{theorem}
\label{Thm_ten_prod_wt_fn_W3.1} Suppose $w=%
{\textstyle\bigotimes\limits_{k=1}^{n}}
w_{k}$ is a tensor product of weight functions with $w_{k}\in\mathbb{R}%
^{d_{k}}$ and $w\in\mathbb{R}^{d}$. Then:

\begin{enumerate}
\item Each $w_{k}\in W2.1$ iff $w\in W2.1$.\medskip

Regarding property W3.1, $w\in W3.1$ iff each $w_{k}\in W3.1$. In fact:

\item[2.1] Suppose each $w_{k}\in W3.1$ for $\alpha^{\left(  k\right)  }$,
$\theta_{k}$ and $\kappa^{\left(  k\right)  }$.

Then $w\in W3.1$ for $\alpha=\left(  \alpha^{\left(  1\right)  }%
,\alpha^{\left(  2\right)  },\ldots,\alpha^{\left(  n\right)  }\right)  $,
$\left\vert \alpha\right\vert =\theta=\theta_{1}+\ldots+\theta_{n}$ and
$\kappa=\left(  \kappa^{\left(  1\right)  },\kappa^{\left(  2\right)  }%
,\ldots,\kappa^{\left(  n\right)  }\right)  $.

\item[2.2] Suppose $w\in W3.1$ for $\alpha$, $\theta$ and $\kappa$. Write
$\alpha=\left(  \alpha^{\left(  1\right)  },\alpha^{\left(  2\right)  }%
,\ldots,\alpha^{\left(  n\right)  }\right)  $ and $\theta_{k}=\left\vert
\alpha_{k}\right\vert $ where $\alpha^{\left(  k\right)  }\in\mathbb{R}%
^{d_{k}}$. Write $\kappa=\left(  \kappa^{\left(  1\right)  },\kappa^{\left(
2\right)  },\ldots,\kappa^{\left(  n\right)  }\right)  $ where $\kappa
^{\left(  k\right)  }\in\mathbb{R}^{d_{k}}$.

Then each $w_{k}\in W3.1$ for $\alpha^{\left(  k\right)  }$, $\theta_{k}$ and
$\kappa^{\left(  k\right)  }$.
\end{enumerate}
\end{theorem}

\begin{proof}
\textbf{Part 1} Suppose $K\subset\mathbb{R}^{d}$ is compact. Set
$r=\max\limits_{\xi\in K}\left\vert \xi\right\vert $ and let $\overline
{R}\left(  a,b\right)  $ denote the closed (compact) rectangle which has
left-most point $a$ and right-most point $b$. Then%
\[
\int\limits_{K}\frac{1}{w}\leq\int\limits_{\overline{R}\left(  -r\mathbf{1}%
_{d},r\mathbf{1}_{d}\right)  }\frac{1}{w}=\prod\limits_{k=1}^{n}%
\int\limits_{R\left(  -r\mathbf{1}_{d_{k}},r\mathbf{1}_{d_{k}}\right)  }%
\frac{1}{w_{k}}<\infty.
\]

Conversely, suppose $w\in W2.1$ i.e. $\int_{K}\frac{1}{w}<\infty$ for any
compact $K\subset\mathbb{R}^{d}$. Now choose $K_{k}\subset\mathbb{R}^{d_{k}}$
to be any compact sets. Then $K=%
{\textstyle\bigotimes\limits_{k=1}^{n}}
K_{k}$ is a compact set and $\prod\limits_{k=1}^{n}\int_{K_{k}}\frac{1}{w_{k}%
}=\int_{K}\frac{1}{w}<\infty$.\medskip

\textbf{Part 2.1} There exist multi-indexes $\alpha^{\left(  k\right)  }$ such
that $\left\vert \alpha^{\left(  k\right)  }\right\vert =\theta_{k}$ and
\[
\int\limits_{\mathbb{R}^{d_{k}}}\frac{\xi^{2\tau^{\left(  k\right)  }}d\xi
}{w_{k}\left(  \xi\right)  \xi^{2\alpha^{\left(  k\right)  }}}<\infty
,\text{\quad}\left\vert \alpha^{\left(  k\right)  }\right\vert =\theta
_{k},\text{ }0\leq\tau^{\left(  k\right)  }\leq\kappa^{\left(  k\right)  }.
\]

Set $\alpha=\left(  \alpha^{\left(  1\right)  },\alpha^{\left(  2\right)
},\ldots,\alpha^{\left(  n\right)  }\right)  \in\mathbb{R}^{d}$ and
$\tau=\left(  \tau^{\left(  1\right)  },\tau^{\left(  2\right)  },\ldots
,\tau^{\left(  n\right)  }\right)  \in\mathbb{R}^{d}$. Then $\left\vert
\alpha\right\vert =\theta$ and $0\leq\tau\leq\kappa$ and%
\[
\int\limits_{\mathbb{R}^{d}}\frac{\xi^{2\tau}d\xi}{w\left(  \xi\right)
\xi^{2\alpha}}=\prod\limits_{k=1}^{n}\int\limits_{\mathbb{R}^{d_{k}}}\frac
{\xi^{2\tau^{\left(  k\right)  }}d\xi}{w_{k}\left(  \xi\right)  \xi
^{2\alpha^{\left(  k\right)  }}}<\infty.
\]
\medskip

\textbf{Part 2.2} ?? \textbf{FINISH}! ??
\end{proof}

\subsection{Relationships between the extended weight function properties}

The next result demonstrates some relationships between the weight function properties.

\begin{theorem}
\label{Thm_weight_property_relat}Suppose $w$ has weight function property W1. Then:

\begin{enumerate}
\item Property W3.1 for $\kappa$ and order $\theta$ implies property W3.2 for
$\underline{\kappa}=\min\kappa$, order $\theta$ and any $r_{3}>0$.

\item Property W3.1 implies property W2 for any $\sigma\geq\theta$ and any
$r_{2}>0$.

\item Property W3.2 implies property W2.2 for $\sigma=\theta$ and any
$r_{2}\geq r_{3}$.

\item Property W3.3 for $\kappa$ and order $\theta$ implies property W2.2 and
property W3.2 for $\underline{\kappa}=\min\kappa$, order $\theta$ and any
$r_{3}>0$.
\end{enumerate}
\end{theorem}

\begin{proof}
\textbf{Part 1} Choosing any $r_{3}>0$ and $0\leq t\leq\min\kappa$ we get%
\[
\int\limits_{\left\vert x\right\vert \geq r_{3}}\dfrac{\left\vert x\right\vert
^{2t}dx}{w\left(  x\right)  \left\vert x\right\vert ^{2\theta}}=\int%
\limits_{\left\vert x\right\vert \geq r_{3}}\dfrac{\left\vert x\right\vert
^{2t}dx}{w\left(  x\right)  \sum\limits_{\left\vert \mu\right\vert =\theta
}\frac{\theta!}{\mu!}x^{2\mu}}\leq\int\limits_{\left\vert x\right\vert \geq
r_{3}}\dfrac{\left\vert x\right\vert ^{2t}dx}{w\left(  x\right)  \frac
{\theta!}{\alpha!}x^{2\alpha}}=\frac{\alpha!}{\theta!}\int\limits_{\left\vert
x\right\vert \geq r_{3}}\dfrac{\left\vert x\right\vert ^{2t}dx}{w\left(
x\right)  x^{2\alpha}}.
\]

But there exist constants $a_{t}$, $b_{t}>0$, independent of $x$, such that%
\begin{equation}
a_{t}\sum_{i=1}^{d}x_{i}^{2t}\leq\left\vert x\right\vert ^{2t}\leq b_{t}%
\sum_{i=1}^{d}x_{i}^{2t},\text{\quad}x\in\mathbb{R}^{d},\label{a40}%
\end{equation}

so that%
\begin{equation}
\int\limits_{\left\vert x\right\vert \geq r_{3}}\dfrac{\left\vert x\right\vert
^{2t}dx}{w\left(  x\right)  \left\vert x\right\vert ^{2\theta}}\leq
\frac{\alpha!}{\theta!}\int\limits_{\left\vert x\right\vert \geq r_{3}}%
\dfrac{b_{t}\sum_{i=1}^{d}x_{i}^{2t}}{w\left(  x\right)  x^{2\alpha}}%
dx=\frac{\alpha!}{\theta!}b_{t}\sum_{i=1}^{d}\int\limits_{\left\vert
x\right\vert \geq r_{3}}\dfrac{x_{i}^{2t}dx}{w\left(  x\right)  x^{2\alpha}%
},\label{a00}%
\end{equation}

and the existence of each of these integrals is guaranteed by property
W3.1.\medskip

\textbf{Part 2} We prove properties W2.1 and W2.2. Suppose $K$ is compact.
Then for some $r>0$, $K\subset B\left(  0;r\right)  $ and%
\[
\int\limits_{\left\vert x\right\vert \leq r}\frac{dx}{w\left(  x\right)
}=\int\limits_{\left\vert x\right\vert \leq r}\frac{x^{2\alpha}dx}{w\left(
x\right)  x^{2\alpha}}\leq\int\limits_{\left\vert x\right\vert \leq r}%
\frac{\left\vert x\right\vert ^{2\theta}dx}{w\left(  x\right)  x^{2\alpha}%
}\leq r^{2\theta}\int\limits_{\left\vert x\right\vert \leq r}\frac
{dx}{w\left(  x\right)  x^{2\alpha}}.
\]

The last integral exists since $w$ has property W3.1. Thus $\dfrac{1}{w}\in
L_{loc}^{1}$.

We next show that W2.2 is true for $\sigma=2\theta$ and any $r_{2}>0$. In
fact, the identity $\left\vert x\right\vert ^{2\theta}=\sum\limits_{\left\vert
\beta\right\vert =\theta}\frac{\theta!}{\beta!}x^{2\beta}$ implies $\left\vert
x\right\vert ^{2\theta}\geq\frac{\theta!}{a!}x^{2\alpha}$, and so
\[
\int\limits_{\left\vert x\right\vert \geq r_{2}}\dfrac{dx}{w\left(  x\right)
\left\vert x\right\vert ^{2\theta}}\leq\frac{\alpha!}{\theta!}\int%
\limits_{\left\vert x\right\vert \geq r_{2}}\dfrac{dx}{w\left(  x\right)
x^{2\alpha}}.
\]

The last integral then exists by property W3.1.\medskip

\textbf{Part 3 }By inspection of the integrals defining property W3.2.\medskip

\textbf{Part 4} $w$ has Property W3.3 for order $\theta$ and some $\kappa
\in\mathbb{R}_{\oplus}^{d}$ if
\[
\int\dfrac{x^{2\lambda}dx}{w\left(  x\right)  \left\vert x\right\vert
^{2\theta}}<\infty,\quad0\leq\lambda\leq\kappa.
\]

The proof of part 1 now implies this part.
\end{proof}

Other results:

\begin{theorem}
\label{Thm_property_wt_fn_W3}\ 

\begin{enumerate}
\item Suppose a weight function $w$ has property W3.1 for order $\theta$ and
parameter $\kappa$. Then%
\[
\int\dfrac{x^{2\lambda}dx}{w\left(  x\right)  \left\vert x\right\vert
^{2\theta}}<\infty,\quad0\leq\lambda\leq\kappa,
\]

and%
\[
\int\dfrac{\left\vert x\right\vert ^{2s}dx}{w\left(  x\right)  \left\vert
x\right\vert ^{2\theta}}<\int\dfrac{\left\vert x\right\vert ^{2s}dx}{w\left(
x\right)  x^{2\alpha}}<\infty,\quad s\leq\underline{\kappa},\text{ }\left\vert
\alpha\right\vert =\theta,
\]

where $\underline{\kappa}:=\min\kappa$.

\item Suppose a weight function $w$ has property W3.2 for order $\theta$ and
parameter $\kappa$. Then%
\[
\int\limits_{\left\vert x\right\vert \geq r_{3}}\dfrac{\left\vert x\right\vert
^{2t}dx}{w\left(  x\right)  \left\vert x\right\vert ^{2\theta}}<\infty
,\quad0\leq t\leq\kappa.
\]

\end{enumerate}
\end{theorem}

\begin{proof}
\textbf{Part 1} For the first set of inequalities use the identity \ref{p08}
i.e. $\left\vert x\right\vert ^{2k}=\sum\limits_{\left\vert \alpha\right\vert
=k}\frac{k!}{\alpha!}x^{2\alpha}$.

Regarding the second set of inequalities: set $t=\underline{\kappa
}-\left\lfloor \underline{\kappa}\right\rfloor $. Then%
\[
\int\dfrac{\left\vert x\right\vert ^{2s}dx}{w\left(  x\right)  x^{2\alpha}%
}=\int_{\left\vert \cdot\right\vert \leq1}\dfrac{\left\vert \cdot\right\vert
^{2s}}{wx^{2\alpha}}+\int_{\left\vert \cdot\right\vert \geq1}\dfrac{\left\vert
\cdot\right\vert ^{2s}}{wx^{2\alpha}}\leq\int_{\left\vert \cdot\right\vert
\leq1}\dfrac{1}{wx^{2\alpha}}+\int_{\left\vert \cdot\right\vert \geq1}%
\dfrac{\left\vert \cdot\right\vert ^{2\underline{\kappa}}}{wx^{2\alpha}}.
\]

Now
\begin{align*}
\int_{\left\vert \cdot\right\vert \geq1}\dfrac{\left\vert \cdot\right\vert
^{2\underline{\kappa}}}{wx^{2\alpha}}  & =\int_{\left\vert \cdot\right\vert
\geq1}\dfrac{\left\vert \cdot\right\vert ^{2\left\lfloor \underline{\kappa
}\right\rfloor }\left\vert \cdot\right\vert ^{2t}}{wx^{2\alpha}}\leq
\int_{\left\vert \cdot\right\vert \geq1}\dfrac{\left\vert \cdot\right\vert
^{2\left\lfloor \underline{\kappa}\right\rfloor }\left(  \left\vert
x_{1}\right\vert ^{2t}+\ldots+\left\vert x_{d}\right\vert ^{2t}\right)
}{wx^{2\alpha}}\\
& =\sum\limits_{\left\vert \beta\right\vert =\left\lfloor \underline{\kappa
}\right\rfloor }\frac{\left\lfloor \underline{\kappa}\right\rfloor !}{\beta
!}\int_{\left\vert \cdot\right\vert \geq1}\dfrac{x^{2\beta}\left(  \left\vert
x_{1}\right\vert ^{2t}+\ldots+\left\vert x_{d}\right\vert ^{2t}\right)
}{wx^{2\alpha}}\\
& =\sum\limits_{\left\vert \beta\right\vert =\left\lfloor \underline{\kappa
}\right\rfloor }\frac{\left\lfloor \underline{\kappa}\right\rfloor !}{\beta
!}\sum\limits_{k=1}^{d}\int_{\left\vert \cdot\right\vert \geq1}\dfrac
{\xi^{2\beta}\left\vert \xi_{k}\right\vert ^{2t}}{wx^{2\alpha}}\\
& =\sum\limits_{\left\vert \beta\right\vert =\left\lfloor \underline{\kappa
}\right\rfloor }\frac{\left\lfloor \underline{\kappa}\right\rfloor !}{\beta
!}\sum\limits_{k=1}^{d}\int_{\left\vert \cdot\right\vert \geq1}\dfrac
{\xi^{2\left(  \beta+t\mathbf{e}_{k}\right)  }}{wx^{2\alpha}}\\
& \leq\sum\limits_{\left\vert \beta\right\vert =\left\lfloor \underline
{\kappa}\right\rfloor }\frac{\left\lfloor \underline{\kappa}\right\rfloor
!}{\beta!}\sum\limits_{k=1}^{d}\int\dfrac{\xi^{2\left(  \beta+t\mathbf{e}%
_{k}\right)  }}{wx^{2\alpha}},
\end{align*}

but $\beta+t\mathbf{e}_{k}\leq\left\lfloor \underline{\kappa}\right\rfloor
\mathbf{1}+t\mathbf{e}_{k}\leq\underline{\kappa}\mathbf{1}\leq\kappa$ so
$\int\dfrac{\xi^{2\left(  \beta+t\mathbf{e}_{k}\right)  }}{w\left\vert
\cdot\right\vert ^{2\theta}}<\infty$ for all $k$.\smallskip

\textbf{Part 2}%
\[
\int\limits_{\left\vert x\right\vert \geq r_{3}}\dfrac{\left\vert x\right\vert
^{2t}dx}{w\left(  x\right)  \left\vert x\right\vert ^{2\theta}}=\int%
\limits_{\left\vert x\right\vert \geq r_{3}}\frac{1}{\left\vert x\right\vert
^{2\left(  \kappa-t\right)  }}\dfrac{\left\vert x\right\vert ^{2\kappa}%
dx}{w\left(  x\right)  \left\vert x\right\vert ^{2\theta}}\leq\frac{1}%
{r_{3}^{2\left(  \kappa-t\right)  }}\int\limits_{\left\vert x\right\vert \geq
r_{3}}\dfrac{\left\vert x\right\vert ^{2\kappa}dx}{w\left(  x\right)
\left\vert x\right\vert ^{2\theta}}<\infty.
\]

\end{proof}

\subsection{Tensor product central difference weight
functions\label{Sect_wt_fn_central_diff}}

The tensor product central difference weight functions were introduced as zero
order weight functions in Chapter 5 of \cite{WilliamsZeroOrdSmthV4}. They are
closely related to the tensor product extended B-spline weight functions which
were studied in the same document and are studied below. Here we exhibit
necessary and sufficient conditions under which the tensor product of $d$
identical 1-dimensional central difference weight functions satisfies
condition W3.1.

\begin{definition}
\label{Def_central_diff_wt_func}\textbf{Central difference weight functions}

Suppose that $q\in L^{1}\left(  \mathbb{R}^{1}\right)  $, $q\neq0$, $q\left(
\xi\right)  \geq0$ and $l,n\geq0$ are integers. The \textbf{univariate}
central difference weight function is defined by
\begin{equation}
w\left(  \xi\right)  =\frac{\xi^{2n}}{\Delta_{2l}\widehat{q}\left(
\xi\right)  },\quad\xi\in\mathbb{R}^{1},\label{a962}%
\end{equation}

where $\Delta_{2l}$ is central difference operator%
\begin{equation}
\Delta_{2l}f\left(  \xi\right)  =\sum_{k=-l}^{l}\left(  -1\right)  ^{k}%
\tbinom{2l}{k+l}f\left(  -k\xi\right)  ,\text{\quad}l=1,2,3,\ldots,\text{ }%
\xi\in\mathbb{R}^{1}.\label{a914}%
\end{equation}

In $d$ dimensions the central difference weight function with parameters $n$
and $l$ is the tensor product of the univariate weight function.
\end{definition}

The following results are quoted from the unpublished zero order weight
function document Williams \cite{WilliamsZeroOrdSmthV4}. This next result
justifies the central weight function definition.

\begin{theorem}
\label{Thm_wt_func_dim1_2}Suppose $w$ is the function on $\mathbb{R}^{1}$
introduced in Definition \ref{Def_central_diff_wt_func}. Then $w$ is an even
function satisfying weight function property W1 for $\mathcal{A}=\left\{
0\right\}  $.
\end{theorem}

\begin{proof}
$w$ is an even function since it is clear from equation \ref{a914} that
$\Delta_{2l}\widehat{q}$ is even. Now $\Delta_{2l,\xi}e^{-i\xi t}=2^{2l}%
\sin^{2l}\left(  \xi t/2\right)  $ so that
\begin{align}
\frac{1}{w\left(  \xi\right)  }=\frac{\Delta_{2l}\widehat{q}\left(
\xi\right)  }{\xi^{2n}} &  =\tfrac{1}{\sqrt{2\pi}}\frac{1}{\xi^{2n}}%
\int\left(  \Delta_{2l,\xi}e^{-i\xi t}\right)  q\left(  t\right)
dt\nonumber\\
&  =\tfrac{2^{2l}}{\sqrt{2\pi}}\frac{1}{\xi^{2n}}\int\sin^{2l}\left(  \xi
t/2\right)  q\left(  t\right)  dt.\label{a943}%
\end{align}

Since $q\in L^{1}$ implies $\widehat{q}\in C_{B}^{\left(  0\right)  }$, from
the definition of $w$ we have $w\in C^{\left(  0\right)  }\left(
\mathbb{R}^{1}\setminus0\right)  $ and $w>0$ on $\mathbb{R}^{1}\setminus0$.
\end{proof}

The central difference weight function is related to the extended B-spline
weight function as follows: suppose $w_{\Lambda}$ is the extended B-spline
weight function with parameters $n,l$. Suppose $w_{c}$ is a central difference
weight function with parameters $n,l,q\left(  \cdot\right)  $. Then \ref{a943}
can be written as%
\begin{equation}
\frac{1}{w\left(  s\right)  }=\tfrac{2^{2\left(  l-n\right)  }}{\sqrt{2\pi}%
}\int_{\mathbb{R}^{1}}\frac{t^{2n}q\left(  t\right)  }{w_{\Lambda}\left(
st/2\right)  }dt.\label{a968}%
\end{equation}

\subsection{Central difference weight functions with property
W3.1\label{SbSect_CentDiffWtW3.1}}

\begin{lemma}
\label{Lem_ten_prod_wt_fn_is_W3.1}Suppose $w$ is the (homogeneous) tensor
product of $d$ identical 1-dimensional functions $w_{1}\in W1$.

Then $w\in W3.1$ for $\theta=\left\vert \alpha\right\vert \geq1$ and parameter
$\kappa=\kappa_{1}\mathbf{1}$ iff for some $r>0$,%
\begin{equation}
\int_{\left\vert s\right\vert \leq r}\frac{ds}{s^{2\overline{\alpha}}%
w_{1}\left(  s\right)  }<\infty,\text{ }\int_{\left\vert s\right\vert \geq
r}\dfrac{s^{2\kappa_{1}}ds}{s^{2\underline{\alpha}}w_{1}\left(  s\right)
}<\infty,\label{p02}%
\end{equation}

where we have used the multi-index notation%
\begin{equation}
\overline{\alpha}:=\max_{k}\alpha_{k},\quad\underline{\alpha}:=\min_{k}%
\alpha_{k}.\label{p18}%
\end{equation}

\end{lemma}

\begin{proof}
Property W3.1 requires that $\int_{\mathbb{R}^{d}}\frac{\xi^{2\tau}d\xi}%
{\xi^{2\alpha}w\left(  \xi\right)  }<\infty$ for $\mathbf{0}\leq\tau\leq
\kappa_{1}\mathbf{1}$. This is true iff
\[
\int_{\mathbb{R}^{1}}\frac{s^{2t}ds}{s^{2\alpha_{k}}w_{1}\left(  s\right)
}<\infty,\quad0\leq t\leq\kappa_{1},\text{ }1\leq k\leq d,
\]

which is true iff for some $r>0$,%
\[
\int_{\left\vert s\right\vert \leq r}\frac{s^{2t}ds}{s^{2\alpha_{k}}%
w_{1}\left(  s\right)  }<\infty\text{ and }\int_{\left\vert s\right\vert \geq
r}\frac{s^{2t}ds}{s^{2\alpha_{k}}w_{1}\left(  s\right)  }<\infty,\quad\left\{
\begin{array}
[c]{l}%
0\leq t\leq\kappa_{1},\\
1\leq k\leq d,
\end{array}
\right.
\]

which is true iff
\[
\int_{\left\vert s\right\vert \leq r}\frac{ds}{s^{2\alpha_{k}}w_{1}\left(
s\right)  }<\infty\text{ and }\int_{\left\vert s\right\vert \geq r}%
\frac{s^{2\kappa_{1}}ds}{s^{2\alpha_{k}}w_{1}\left(  s\right)  }<\infty
,\quad1\leq k\leq d,
\]

which is true iff%
\[
\int_{\left\vert s\right\vert \leq r}\frac{ds}{s^{2\overline{\alpha}}%
w_{1}\left(  s\right)  }<\infty\text{ and }\int_{\left\vert s\right\vert \geq
r}\dfrac{s^{2\kappa_{1}}ds}{s^{2\underline{\alpha}}w_{1}\left(  s\right)
}<\infty.
\]

\end{proof}

\begin{lemma}
\label{Lem_low_bnd_integ_sinx_dev_x}If $n+\lambda>1/2$ and $x>0$, then
\begin{equation}%
\begin{array}
[c]{ll}%
\frac{1}{2\left(  n+\lambda\right)  -1}\frac{2^{-l}}{\left(  x+\frac{5\pi}%
{4}\right)  ^{2\left(  n+\lambda\right)  -1}} & \leq\int\limits_{x}^{\infty
}\frac{\sin^{2l}u}{u^{2\left(  n+\lambda\right)  }}du\leq\frac{1}{2\left(
n+\lambda\right)  -1}\frac{1}{x^{2\left(  n+\lambda\right)  -1}}.
\end{array}
\label{a099}%
\end{equation}

If further $x\geq\rho>0$ then%
\begin{equation}%
\begin{array}
[c]{l}%
\frac{2^{-l}}{2\left(  n+\lambda\right)  -1}\frac{1}{\left(  \left(
1+\frac{5\pi}{4\rho}\right)  x\right)  ^{2\left(  n+\lambda\right)  -1}}%
\leq\int_{x}^{\infty}\frac{\sin^{2l}u}{u^{2\left(  n+\lambda\right)  }}du.
\end{array}
\label{a098}%
\end{equation}

Finally, if $n+\lambda\neq l+1/2$ and $0<x<\pi/2$,%
\begin{equation}
\left.
\begin{array}
[c]{ll}%
\left(  \frac{2}{\pi}\right)  ^{2}\frac{\left(  \frac{\pi}{2}\right)
^{2\left(  l-n-\lambda\right)  +1}-x^{2\left(  l-n-\lambda\right)  +1}%
}{2\left(  l-n-\lambda\right)  +1} & \leq\int\limits_{x}^{\frac{\pi}{2}}%
\frac{\sin^{2l}u}{u^{2\left(  n+\lambda\right)  }}du\\
& \leq\frac{\left(  \frac{\pi}{2}\right)  ^{2\left(  l-n-\lambda\right)
+1}-x^{2\left(  l-n-\lambda\right)  +1}}{2\left(  l-n-\lambda\right)  +1}%
\end{array}
\right\} \label{a100}%
\end{equation}

\end{lemma}

\begin{proof}
\textbf{Upper bound of }\ref{a099}:%
\begin{align*}
\int_{x}^{\infty}\frac{\sin^{2l}u}{u^{2\left(  n+\lambda\right)  }}du\leq
\int_{x}^{\infty}\frac{1}{u^{2\left(  n+\lambda\right)  }}du  & =\left[
\frac{1}{1-2\left(  n+\lambda\right)  }\frac{1}{u^{2\left(  n+\lambda\right)
-1}}\right]  _{x}^{\infty}\\
& =\frac{1}{2\left(  n+\lambda\right)  -1}\frac{1}{x^{2\left(  n+\lambda
\right)  -1}}.
\end{align*}

\textbf{Lower bound of} \ref{a099}:%
\begin{align*}
\int_{x}^{\infty}\frac{\sin^{2l}u}{u^{2\left(  n+\lambda\right)  }}du  &
>\sum_{k=\left\lceil x/\pi\right\rceil }^{\infty}\int_{\pi k}^{\pi\left(
k+1\right)  }\frac{\sin^{2l}u}{u^{2\left(  n+\lambda\right)  }}du\\
& >\sum_{k=\left\lceil x/\pi\right\rceil }^{\infty}\int_{\pi k+\pi/4}%
^{\pi\left(  k+1\right)  -\pi/4}\frac{\sin^{2l}u}{u^{2\left(  n+\lambda
\right)  }}du\\
& >\sin^{2l}\frac{\pi}{4}\sum_{k=\left\lceil x/\pi\right\rceil }^{\infty}%
\int_{\pi k+\pi/4}^{\pi\left(  k+1\right)  -\pi/4}\frac{du}{u^{2\left(
n+\lambda\right)  }}\\
& =\frac{2^{-l}}{2\left(  n+\lambda\right)  -1}\sum_{k=\left\lceil
x/\pi\right\rceil }^{\infty}\left[  \frac{-1}{u^{2\left(  n+\lambda\right)
-1}}\right]  _{\pi k+\pi/4}^{\pi\left(  k+1\right)  -\pi/4}\\
& =\frac{2^{-l}}{2\left(  n+\lambda\right)  -1}\frac{1}{\left(  \pi\left\lceil
\frac{x}{\pi}\right\rceil +\frac{\pi}{4}\right)  ^{2\left(  n+\lambda\right)
-1}}\\
& >\frac{2^{-l}}{2\left(  n+\lambda\right)  -1}\frac{1}{\left(  \pi\left(
\frac{x}{\pi}+1\right)  +\frac{\pi}{4}\right)  ^{2\left(  n+\lambda\right)
-1}}\\
& =\frac{2^{-l}}{2\left(  n+\lambda\right)  -1}\frac{1}{\left(  x+\frac{5\pi
}{4}\right)  ^{2\left(  n+\lambda\right)  -1}}.
\end{align*}

\textbf{Proof of }\ref{a100}:%
\begin{align*}
\int_{x}^{\frac{\pi}{2}}\frac{\sin^{2l}u}{u^{2\left(  n+\lambda\right)  }}du
& \leq\int_{x}^{\frac{\pi}{2}}\frac{du}{u^{2\left(  n+\lambda\right)  }}\\
& =\frac{1}{1-2\left(  n+\lambda\right)  }\left[  \frac{1}{u^{2\left(
n+\lambda\right)  -1}}\right]  _{x}^{\frac{\pi}{2}}\\
& =\frac{1}{1-2\left(  n+\lambda\right)  }\left[  \frac{1}{u^{2\left(
n+\lambda\right)  -1}}\right]  _{x}^{\frac{\pi}{2}}\\
& =\frac{1}{1-2\left(  n+\lambda\right)  }\left(  \frac{1}{\left(
\pi/2\right)  ^{2\left(  n+\lambda\right)  -1}}-\frac{1}{x^{2\left(
n+\lambda\right)  -1}}\right) \\
& =\frac{1}{2\left(  n+\lambda\right)  -1}\left(  \frac{1}{x^{2\left(
n+\lambda\right)  -1}}-\frac{1}{\left(  \pi/2\right)  ^{2\left(
n+\lambda\right)  -1}}\right)  ,
\end{align*}

and%
\[
\int_{x}^{\frac{\pi}{2}}\frac{\sin^{2l}u}{u^{2\left(  n+\lambda\right)  }%
}du=\int_{x}^{\frac{\pi}{2}}\frac{u^{2l}}{u^{2\left(  n+\lambda\right)  }%
}\frac{\sin^{2l}u}{u^{2l}}du=\int_{x}^{\frac{\pi}{2}}u^{2\left(
l-n-\lambda\right)  }\frac{\sin^{2l}u}{u^{2l}}du.
\]

But%
\[
\left(  \frac{2}{\pi}\right)  ^{2}\leq\frac{\sin^{2}u}{u^{2}}\leq
1,\quad\left\vert u\right\vert \leq\frac{\pi}{2},
\]

so%
\begin{align*}
\left(  \frac{2}{\pi}\right)  ^{2}\int_{x}^{\frac{\pi}{2}}u^{2\left(
l-n-\lambda\right)  }du  & \leq\int_{x}^{\frac{\pi}{2}}\frac{\sin^{2l}%
u}{u^{2\left(  n+\lambda\right)  }}du\\
& \leq\int_{x}^{\frac{\pi}{2}}u^{2\left(  l-n-\lambda\right)  }du,
\end{align*}

and%
\begin{align*}
\left(  \frac{2}{\pi}\right)  ^{2}\frac{\left(  \frac{\pi}{2}\right)
^{2\left(  l-n-\lambda\right)  +1}-x^{2\left(  l-n-\lambda\right)  +1}%
}{2\left(  l-n-\lambda\right)  +1}  & \leq\int\limits_{x}^{\frac{\pi}{2}}%
\frac{\sin^{2l}u}{u^{2\left(  n+\lambda\right)  }}du\\
& \leq\frac{\left(  \frac{\pi}{2}\right)  ^{2\left(  l-n-\lambda\right)
+1}-x^{2\left(  l-n-\lambda\right)  +1}}{2\left(  l-n-\lambda\right)  +1}.
\end{align*}

\end{proof}

\begin{theorem}
\label{Thm_centdiff_wtfn_in_W3.1}Suppose $w$ is the tensor product of $d$
identical 1-dimensional \textbf{central difference} weight functions $w_{1}$. Then:

\begin{enumerate}
\item The \textbf{first condition} of \ref{p02} holds iff for some $R>0$,
\begin{equation}
\overline{\alpha}\leq l-n\quad and\quad\int_{\left\vert t\right\vert \geq
R}\left\vert t\right\vert ^{2\left(  n+\overline{\alpha}\right)  -1}q\left(
t\right)  dt<\infty,\label{a101}%
\end{equation}

\item The \textbf{second condition} of \ref{p02} holds iff
\begin{equation}
\underline{\alpha}>\kappa_{1}-n+\frac{1}{2}.\label{a107}%
\end{equation}

In other words, $w\in W3.1$ for $\theta=\left\vert \alpha\right\vert \geq1 $
and parameter $\kappa=\kappa_{1}\mathbf{1}$ iff \ref{a101} and \ref{a107} hold.
\end{enumerate}
\end{theorem}

\begin{proof}
\fbox{\textbf{Part 1}} Suppose the \textbf{first condition} of\textbf{\ }%
\ref{p02} holds. Then from \ref{a968},%
\begin{align*}
\int_{\left\vert s\right\vert \leq r}\frac{ds}{s^{2\overline{\alpha}}%
w_{1}\left(  s\right)  }  & =\int_{\left\vert s\right\vert \leq r}\frac
{1}{s^{2\overline{\alpha}}}\tfrac{2^{2\left(  l-n\right)  }}{\sqrt{2\pi}}%
\int_{\mathbb{R}^{1}}\frac{t^{2n}q\left(  t\right)  }{w_{\Lambda}\left(
st/2\right)  }dt\,ds\\
& =\tfrac{2^{2\left(  l-n\right)  }}{\sqrt{2\pi}}\int_{\mathbb{R}^{1}}%
t^{2n}q\left(  t\right)  \left(  \int_{\left\vert s\right\vert \leq r}%
\frac{ds}{s^{2\overline{\alpha}}w_{\Lambda}\left(  st/2\right)  }\right)  dt\\
& =\tfrac{2^{2\left(  l-n\right)  +1}}{\sqrt{2\pi}}\int_{\mathbb{R}^{1}}%
t^{2n}q\left(  t\right)  \left(  \int_{0}^{r}\frac{ds}{s^{2\overline{\alpha}%
}w_{\Lambda}\left(  st/2\right)  }\right)  dt\\
& =\tfrac{2^{2\left(  l-n\right)  +1}}{\sqrt{2\pi}}\int_{\mathbb{R}^{1}}%
t^{2n}q\left(  t\right)  \left(  \int_{0}^{r}\frac{ds}{s^{2\overline{\alpha}%
}w_{\Lambda}\left(  s\left\vert t\right\vert /2\right)  }\right)  dt\\
& =\tfrac{2^{2\left(  l-n\right)  +1}}{\sqrt{2\pi}}\int_{\mathbb{R}^{1}}%
t^{2n}q\left(  t\right)  \left(  \int_{0}^{r}\frac{ds}{s^{2\overline{\alpha}%
}w_{\Lambda}\left(  s\left\vert t\right\vert /2\right)  }\right)  dt.
\end{align*}

Change of variables: $u=\frac{\left\vert t\right\vert }{2}s$, $du=\frac
{\left\vert t\right\vert }{2}ds$ yields%
\begin{equation}
\int\limits_{\left\vert s\right\vert \leq r}\frac{ds}{s^{2\overline{\alpha}%
}w_{1}\left(  s\right)  }=\tfrac{2^{2\left(  l-n-\overline{\alpha}+1\right)
}}{\sqrt{2\pi}}\int\limits_{\mathbb{R}^{1}}\left\vert t\right\vert ^{2\left(
n+\overline{\alpha}\right)  -1}q\left(  t\right)  \int\limits_{0}^{r\left\vert
t\right\vert /2}\frac{du}{u^{2\overline{\alpha}}w_{\Lambda}\left(  u\right)
}dt.\label{a105}%
\end{equation}

Since%
\[
\int_{0}^{r\left\vert t\right\vert /2}\frac{du}{u^{2\overline{\alpha}%
}w_{\Lambda}\left(  u\right)  }=\int_{0}^{r\left\vert t\right\vert /2}%
\frac{\sin^{2l}u}{u^{2\left(  \overline{\alpha}+n\right)  }}du,
\]

this integral exists iff
\begin{equation}
2\left(  n+\overline{\alpha}\right)  -2l<1\Leftrightarrow2\left(
n+\overline{\alpha}\right)  -2l\leq0\Leftrightarrow n+\overline{\alpha}\leq
l,\label{a102}%
\end{equation}

and, since $n\geq1$ implies $2\left(  n+\overline{\alpha}\right)  \geq2$, we
further have%
\[
\int_{0}^{\infty}\frac{du}{u^{2\overline{\alpha}}w_{\Lambda}\left(  u\right)
}=\int_{0}^{\infty}\frac{\sin^{2l}v}{v^{2\left(  n+\overline{\alpha}\right)
}}dv<\infty.
\]

Thus, if $n+\overline{\alpha}\leq l$ then \ref{a105} implies%
\begin{align*}
&  \int_{\left\vert s\right\vert \leq r}\frac{ds}{s^{2\overline{\alpha}}%
w_{1}\left(  s\right)  }\\
&  <\tfrac{2^{2\left(  l-n-\overline{\alpha}+1\right)  }}{\sqrt{2\pi}}%
\int_{\mathbb{R}^{1}}\left\vert t\right\vert ^{2\left(  n+\overline{\alpha
}\right)  -1}q\left(  t\right)  \int_{0}^{\infty}\frac{du}{u^{2\overline
{\alpha}}w_{\Lambda}\left(  u\right)  }dt\\
&  =\tfrac{2^{2\left(  l-n-\overline{\alpha}+1\right)  }}{\sqrt{2\pi}}\int%
_{0}^{\infty}\frac{du}{u^{2\overline{\alpha}}w_{\Lambda}\left(  u\right)
}\int_{\mathbb{R}^{1}}\left\vert \cdot\right\vert ^{2\left(  n+\overline
{\alpha}\right)  -1}q\\
&  =\tfrac{2^{2\left(  l-n-\overline{\alpha}+1\right)  }}{\sqrt{2\pi}}%
\int\limits_{0}^{\infty}\frac{du}{u^{2\overline{\alpha}}w_{\Lambda}\left(
u\right)  }\left(  \int\limits_{\left\vert \cdot\right\vert \leq R}\left\vert
\cdot\right\vert ^{2\left(  n+\overline{\alpha}\right)  -1}q+\int%
\limits_{\left\vert \cdot\right\vert \geq R}\left\vert \cdot\right\vert
^{2\left(  n+\overline{\alpha}\right)  -1}q\right) \\
&  \leq\tfrac{2^{2\left(  l-n-\overline{\alpha}+1\right)  }}{\sqrt{2\pi}}%
\int\limits_{0}^{\infty}\frac{du}{u^{2\overline{\alpha}}w_{\Lambda}\left(
u\right)  }\left(  R^{2\left(  n+\overline{\alpha}\right)  -1}\int%
\limits_{\left\vert \cdot\right\vert \leq R}q+\int\limits_{\left\vert
\cdot\right\vert \geq R}\left\vert \cdot\right\vert ^{2\left(  n+\overline
{\alpha}\right)  -1}q\right) \\
&  <\infty.
\end{align*}

On the other hand, if $\int_{\left\vert s\right\vert \leq r}\frac
{ds}{s^{2\overline{\alpha}}w_{c}\left(  s\right)  }<\infty$ then \ref{a105}
implies $\int_{0}^{r\left\vert t\right\vert /2}\frac{du}{u^{2\overline{\alpha
}}w_{\Lambda}\left(  u\right)  }<\infty$ for almost all $t$ and so
$n+\overline{\alpha}\leq l$.

Thus the first condition of \ref{p02} holds iff $n+\overline{\alpha}\leq l$
and $\int_{\left\vert t\right\vert \geq R}\left\vert t\right\vert ^{2\left(
n+\overline{\alpha}\right)  -1}q\left(  t\right)  dt$.\medskip

\fbox{\textbf{Part 2}} Using \ref{a968},
\begin{align*}
\int_{\left\vert s\right\vert \geq r}\frac{ds}{s^{2\overline{\alpha}}%
w_{1}\left(  s\right)  }  & =\int_{\left\vert s\right\vert \geq r}%
\dfrac{s^{2\kappa_{1}}}{s^{2\underline{\alpha}}}\tfrac{2^{2\left(  l-n\right)
}}{\sqrt{2\pi}}\left(  \int_{\mathbb{R}^{1}}\frac{t^{2n}q\left(  t\right)
}{w_{\Lambda}\left(  st/2\right)  }dt\right)  ds\\
& =\tfrac{2^{2\left(  l-n\right)  }}{\sqrt{2\pi}}\int_{\mathbb{R}^{1}}%
t^{2n}q\left(  t\right)  \left(  \int_{\left\vert s\right\vert \geq r}%
\frac{s^{2\kappa_{1}}ds}{s^{2\underline{\alpha}}w_{\Lambda}\left(
st/2\right)  }\right)  dt\\
& =\tfrac{2^{2\left(  l-n\right)  +1}}{\sqrt{2\pi}}\int_{\mathbb{R}^{1}}%
t^{2n}q\left(  t\right)  \left(  \int_{r}^{\infty}\frac{s^{2\kappa_{1}}%
ds}{s^{2\underline{\alpha}}w_{\Lambda}\left(  st/2\right)  }\right)  dt\\
& =\tfrac{2^{2\left(  l-n\right)  +1}}{\sqrt{2\pi}}\int_{\mathbb{R}^{1}}%
t^{2n}q\left(  t\right)  \left(  \int_{r}^{\infty}\frac{s^{2\kappa_{1}}%
ds}{s^{2\underline{\alpha}}w_{\Lambda}\left(  s\left\vert t\right\vert
/2\right)  }\right)  dt\\
& =\tfrac{2^{2\left(  l-n\right)  +1}}{\sqrt{2\pi}}\int_{\mathbb{R}^{1}}%
t^{2n}q\left(  t\right)  \left(  \int_{r}^{\infty}\frac{s^{2\kappa_{1}}%
ds}{s^{2\underline{\alpha}}w_{\Lambda}\left(  s\left\vert t\right\vert
/2\right)  }\right)  dt.
\end{align*}

The change of variables: $u=\frac{\left\vert t\right\vert }{2}s$,
$du=\frac{\left\vert t\right\vert }{2}ds$ yields%
\begin{align*}
\int_{r}^{\infty}\frac{s^{2\kappa_{1}}ds}{s^{2\underline{\alpha}}w_{\Lambda
}\left(  s\left\vert t\right\vert /2\right)  }  & =\int_{2r/\left\vert
t\right\vert }^{\infty}\frac{\left(  \frac{2u}{\left\vert t\right\vert
}\right)  ^{2\kappa_{1}}\frac{2}{\left\vert t\right\vert }du}{\left(
\frac{2u}{\left\vert t\right\vert }\right)  ^{2\underline{\alpha}}w_{\Lambda
}\left(  u\right)  }\\
& =2^{2\left(  \kappa_{1}-\underline{\alpha}\right)  +1}\left\vert
t\right\vert ^{2\left(  \underline{\alpha}-\kappa_{1}\right)  -1}%
\int_{2\left\vert t\right\vert r}^{\infty}\frac{u^{2\kappa_{1}}du}%
{u^{2\underline{\alpha}}w_{\Lambda}\left(  u\right)  },
\end{align*}

so that%
\begin{align}
&  \int_{\left\vert s\right\vert \geq r}\frac{ds}{s^{2\overline{\alpha}}%
w_{1}\left(  s\right)  }\nonumber\\
&  =\tfrac{2^{2\left(  l-n\right)  +1}}{\sqrt{2\pi}}\int\limits_{\mathbb{R}%
^{1}}t^{2n}q\left(  t\right)  \left(  2^{2\left(  \kappa_{1}-\underline
{\alpha}\right)  +1}\left\vert t\right\vert ^{2\left(  \underline{\alpha
}-\kappa_{1}\right)  -1}\int\limits_{2\left\vert t\right\vert r}^{\infty}%
\frac{u^{2\kappa_{1}}du}{u^{2\underline{\alpha}}w_{\Lambda}\left(  u\right)
}\right)  dt\nonumber\\
&  =\tfrac{2^{2\left(  l-n+\kappa_{1}-\underline{\alpha}+1\right)  }}%
{\sqrt{2\pi}}\int\limits_{\mathbb{R}^{1}}\left\vert t\right\vert ^{2\left(
n+\underline{\alpha}-\kappa_{1}\right)  -1}q\left(  t\right)  \int%
\limits_{2\left\vert t\right\vert r}^{\infty}\frac{u^{2\kappa_{1}}%
du}{u^{2\underline{\alpha}}w_{\Lambda}\left(  u\right)  }dt.\label{a106}%
\end{align}

Since%
\begin{equation}
\int_{2\left\vert t\right\vert r}^{\infty}\frac{u^{2\kappa_{1}}du}%
{u^{2\underline{\alpha}}w_{\Lambda}\left(  u\right)  }=\int_{2\left\vert
t\right\vert r}^{\infty}\frac{u^{2\kappa_{1}}\sin^{2l}u}{u^{2\underline
{\alpha}}u^{2n}}du=\int_{2\left\vert t\right\vert r}^{\infty}\frac{\sin
^{2l}u\text{ }du}{u^{2\left(  n+\underline{\alpha}-\kappa_{1}\right)  }%
},\label{a108}%
\end{equation}

this integral exists for any $\left\vert t\right\vert >0$ iff
\[
2\left(  n+\underline{\alpha}-\kappa_{1}\right)  >1\Leftrightarrow
n+\underline{\alpha}>\kappa_{1}+\frac{1}{2},
\]

so by assumption this integral exists. Now writing%
\[
\int\limits_{\left\vert s\right\vert \geq r}\frac{ds}{s^{2\overline{\alpha}%
}w_{1}\left(  s\right)  }=\tfrac{2^{2\left(  l-n+\kappa_{1}-\underline{\alpha
}+1\right)  }}{\sqrt{2\pi}}\int\limits_{\mathbb{R}^{1}}\left\vert t\right\vert
^{2\left(  n+\underline{\alpha}-\kappa_{1}\right)  -1}q\left(  t\right)
\int\limits_{2\left\vert t\right\vert r}^{\infty}\frac{\sin^{2l}u\text{ }%
du}{u^{2\left(  n+\underline{\alpha}-\kappa_{1}\right)  }}dt,
\]

an application of inequality \ref{a098} of Lemma
\ref{Lem_low_bnd_integ_sinx_dev_x} gives%
\begin{align*}
&  \int_{\left\vert s\right\vert \geq r}\frac{ds}{s^{2\overline{\alpha}}%
w_{1}\left(  s\right)  }\\
&  \leq\tfrac{2^{2\left(  l-n+\kappa_{1}-\underline{\alpha}+1\right)  }}%
{\sqrt{2\pi}}\int\limits_{\mathbb{R}^{1}}\left\vert t\right\vert ^{2\left(
n+\underline{\alpha}-\kappa_{1}\right)  -1}q\left(  t\right)  \tfrac
{1}{2\left(  n+\underline{\alpha}-\kappa_{1}\right)  -1}\frac{dt}{\left(
2\left\vert t\right\vert r\right)  ^{2\left(  n+\underline{\alpha}-\kappa
_{1}\right)  -1}}\\
&  =\tfrac{1}{\sqrt{2\pi}}\tfrac{1}{2\left(  n+\underline{\alpha}-\kappa
_{1}\right)  -1}\frac{2^{2\left(  l-n+\kappa_{1}-\underline{\alpha}\right)
+1}}{\left(  2r\right)  ^{2\left(  n+\underline{\alpha}-\kappa_{1}\right)
-1}}\int\limits_{\mathbb{R}^{1}}q\\
&  =\tfrac{1}{\sqrt{2\pi}}\tfrac{1}{2\left(  n+\underline{\alpha}-\kappa
_{1}\right)  -1}\frac{4^{l}}{\left(  4r\right)  ^{2\left(  n+\underline
{\alpha}-\kappa_{1}\right)  -1}}\int\limits_{\mathbb{R}^{1}}q\\
&  <\infty.
\end{align*}

Now suppose $\int_{\left\vert s\right\vert \leq r}\frac{ds}{s^{2\overline
{\alpha}}w_{c}\left(  s\right)  }<\infty$. Then from \ref{a106} and \ref{a108}
we have%
\begin{align*}
\int_{\left\vert s\right\vert \leq r} &  \frac{ds}{s^{2\overline{\alpha}}%
w_{1}\left(  s\right)  }\\
&  =\lim_{\varepsilon\rightarrow0^{+}}\tfrac{2^{2\left(  l-n+\kappa
_{1}-\underline{\alpha}+1\right)  }}{\sqrt{2\pi}}\int_{\left\vert t\right\vert
\geq\varepsilon}\left\vert t\right\vert ^{2\left(  n+\underline{\alpha}%
-\kappa_{1}\right)  -1}q\left(  t\right)  \int_{2\left\vert t\right\vert
r}^{\infty}\frac{\sin^{2l}u\text{ }du}{u^{2\left(  n+\underline{\alpha}%
-\kappa_{1}\right)  }}dt,
\end{align*}

which implies $\int_{2\left\vert t\right\vert r}^{\infty}\frac{\sin^{2l}%
u}{u^{2\left(  n+\underline{\alpha}-\kappa_{1}\right)  }}du<\infty$ for almost
all $t$ and so \ref{a107} must hold.
\end{proof}

\begin{corollary}
\label{Cor_Thm_centdiff_wtfn_in_W3.1}Suppose $w_{c}$ is the (homogeneous)
tensor product of $d$ identical 1-dimensional \textbf{central difference}
weight functions $w_{1}$ for which $q$ \textbf{has bounded support}. Then
$w_{c}\in W3.1$ for $\theta=\left\vert \alpha\right\vert \geq1$ and parameter
$\kappa=\kappa_{1}\mathbf{1}$ iff%
\[
\kappa_{1}-n+1/2<\underline{\alpha}\leq\overline{\alpha}\leq l-n.
\]

In fact the statement of Theorem \ref{Thm_ExtNatSplin_wt_fn_W3.1} holds if we
replace \textbf{extended B-spline weight function} by \textbf{central
difference weight function}.
\end{corollary}

\begin{remark}
See Theorem \ref{Thm_ExtNatSplin_wt_fn_W3.1} and Remark
\ref{Rem_Thm_ExtNatSplin_wt_fn_W3.1} for information regarding choosing
$\alpha$ in order to maximize the smoothness $\kappa_{1}$.
\end{remark}

\subsection{Extended B-spline weight functions with property
W3.1*\label{SbSect_ext_splin_wt_fn_W3.1*}}

In this subsection we will introduce the tensor product, extended B-spline
weight functions. These are two-parameter, tensor product weight functions and
first we will derive necessary and sufficient conditions for a weight function
to have property W3.1*\textbf{.}

\begin{lemma}
\label{Lem_ten_prod_wt_fn_is_W3.1*}Now suppose $w$ is the tensor product of
1-dimensional functions $\left\{  w_{k}\right\}  $ each with property W1. Then:

\begin{enumerate}
\item $w$ has property W3.1* for order $\theta$ and $\kappa$ iff for all
$r,R>0$ and all $k$,%
\begin{equation}
\int_{\left\vert s\right\vert \leq r}\dfrac{ds}{s^{2\theta}w_{k}\left(
s\right)  }<\infty,\text{\quad}\int_{\left\vert s\right\vert \geq R}%
\dfrac{s^{2\kappa_{k}}ds}{w_{k}\left(  s\right)  }<\infty.\label{p69}%
\end{equation}

\item If $\kappa=\kappa_{1}\mathbf{1}$ and the $w_{k}$ are identical the
conditions \ref{p69} become%
\begin{equation}
\int_{\left\vert s\right\vert \leq r}\dfrac{ds}{s^{2\theta}w_{1}\left(
s\right)  }<\infty,\text{\quad}\int_{\left\vert s\right\vert \geq R}%
\dfrac{s^{2\kappa_{1}}ds}{w_{1}\left(  s\right)  }<\infty.\label{p71}%
\end{equation}

\end{enumerate}
\end{lemma}

\begin{proof}
\textbf{Part 1} From the definition of property W3.1*, $w$ has property W3.1*
iff for each $k$,
\[
\int\dfrac{s^{2\beta_{k}}ds}{s^{2\alpha_{k}}w_{k}\left(  s\right)  }%
<\infty,\quad0\leq\beta_{k}\leq\kappa_{k},\text{ all }\left\vert
\alpha\right\vert =\theta.
\]

Near zero: $\int_{\left\vert s\right\vert \leq r}\dfrac{s^{2\beta_{k}}%
ds}{s^{2\alpha_{k}}w_{k}\left(  s\right)  }<\infty$ iff $\int_{\left\vert
s\right\vert \leq r}\dfrac{ds}{s^{2\theta}w_{k}\left(  s\right)  }<\infty$.

Then near $\infty$: $\int_{\left\vert s\right\vert \geq R}\dfrac{s^{2\beta
_{k}}}{s^{2\alpha_{k}}w_{k}\left(  s\right)  }ds<\infty$ for $0\leq\beta
_{k}\leq\kappa_{k}$ iff $\int_{\left\vert s\right\vert \geq R}\dfrac
{s^{2\kappa_{k}}}{w_{k}\left(  s\right)  }ds<\infty$.\smallskip

\textbf{Part 2} Set $k=1$ in \ref{p69}.
\end{proof}

\begin{theorem}
\label{Thm_ExtNatSplin_wt_fn_W3.1*}\textbf{Tensor product, extended B-spline
weight functions} Suppose $d\geq1$ and $\theta\geq1$ is an integer. For given
integers $l,n\geq0$ the tensor product, extended B-spline weight function $w$
is defined by
\begin{equation}
w\left(  x\right)  :=\prod_{i=1}^{d}\frac{x_{i}^{2n}}{\sin^{2l}x_{i}%
}.\label{p32}%
\end{equation}

Then:\smallskip

\frame{When $d=1$} there exists a closed set $\mathcal{A}_{1}$ of measure zero
such that $w$ is a function with properties W1 and W2. It also has property
\textbf{W3.1}* for order $\theta$ and scalar $\kappa$ iff $n$ and $l$ satisfy%
\begin{equation}
\kappa+1/2<n+\theta,\quad\theta\leq l-n.\label{p16}%
\end{equation}

Consequently%
\begin{equation}
\left\lfloor \kappa\right\rfloor \leq n+\theta-1,\quad\left\lfloor
2\kappa\right\rfloor \leq2\left(  n+\theta-1\right)  .\label{p05}%
\end{equation}

\frame{When $d\geq2$} there exists a closed set $\mathcal{A}_{d}$ of measure
zero such that $w$ is a function with properties W1 and W2. It also has
property W3.1* for order $\theta$ and scalar $\kappa$ iff $n$ and $l$ satisfy
\begin{equation}
\kappa+1/2<n,\quad\theta\leq l-n.\label{p17}%
\end{equation}

Consequently%
\begin{equation}
\left\lfloor \kappa\right\rfloor \leq n-1,\quad\left\lfloor 2\kappa
\right\rfloor \leq2\left(  n-1\right)  .\label{p15}%
\end{equation}

\end{theorem}

\begin{proof}
Clearly for each parameter pair $l,n$ there exists a closed set $\mathcal{A}%
_{d}$ of measure zero such that property W1 is satisfied. Further, by part 2
of Theorem \ref{Thm_weight_property_relat} $w$ has property W2 and so we need
only prove property W3.1* which we do so by using the conditions \ref{p69} of
part 1 Lemma \ref{Lem_ten_prod_wt_fn_is_W3.1*} with $\frac{1}{w_{k}\left(
s\right)  }=\frac{\sin^{2l}s}{s^{2n}}$. Now for given $k$ the integral%
\[
\int_{\left\vert s\right\vert \leq r}\dfrac{ds}{s^{2\theta}w_{k}\left(
s\right)  }=\int_{\left\vert s\right\vert \leq r}\frac{\sin^{2l}s}{s^{2n}%
}\dfrac{ds}{s^{2\theta}}=\int_{\left\vert s\right\vert \leq r}\frac{\sin
^{2l}s}{s^{2l}}\frac{ds}{s^{2\left(  n+\theta-l\right)  }},
\]

exists%
\[%
\begin{array}
[c]{ll}%
iff & 2\left(  n+\theta-l\right)  <1,\\
iff & n+\theta-l<1/2,\\
iff & n+\theta-l\leq0,\\
iff & \theta\leq l-n.
\end{array}
\]

Also the integral%
\[
\int_{\left\vert s\right\vert \geq r}\dfrac{ds}{s^{2\left(  \underline{\theta
}-\kappa_{k}\right)  }w_{k}\left(  s\right)  }=\int_{\left\vert s\right\vert
\geq r}\frac{\sin^{2l}s}{s^{2n}}\frac{ds}{s^{2\left(  \underline{\theta
}-\kappa_{k}\right)  }}=\int_{\left\vert s\right\vert \geq r}\frac{\sin
^{2l}s\text{\thinspace}ds}{s^{2\left(  n+\underline{\theta}-\kappa_{k}\right)
}},
\]

exists%
\[%
\begin{array}
[c]{ll}%
iff & 2\left(  n+\underline{\theta}-\kappa_{k}\right)  >1,\\
iff & n+\underline{\theta}-\kappa_{k}>1/2,\\
iff & n+\theta-\kappa>1/2\text{ }when\text{ }d=1,\text{ }and\text{ }%
n-\kappa>1/2\text{ }when\text{ }d\geq2\\
iff & \kappa+1/2<n+\theta\text{ }when\text{ }d=1,\text{ }and\text{ }%
\kappa+1/2<n\text{ }when\text{ }d\geq2\\
iff & \max\kappa+1/2<n+\theta\text{ }when\text{ }d=1,\text{ }and\text{ }%
\max\kappa+1/2<n\text{ }when\text{ }d\geq2.
\end{array}
\]

\end{proof}

\subsection{Extended B-spline weight functions with property
W3.1\label{SbSect_ext_splin_wt_fn_W3.1}}

We now apply Lemma \ref{Lem_ten_prod_wt_fn_is_W3.1} to the extended B-spline
weight functions to obtain necessary and sufficient conditions that they have
property W3.1.

\begin{theorem}
\label{Thm_ExtNatSplin_wt_fn_W3.1}Suppose $w$ is the \textbf{extended B-spline
weight function} \ref{p32} with parameters $l$ and $n$. Then:

\begin{enumerate}
\item $w\in W3.1$ for $\theta=\left\vert \alpha\right\vert $ and
$\kappa=\kappa_{1}\mathbf{1}$ iff $n$ and $l$ satisfy
\begin{equation}
\kappa_{1}-n+1/2<\underline{\alpha}\leq\overline{\alpha}\leq l-n.\label{a103}%
\end{equation}

so that
\begin{equation}
\left\lfloor \kappa_{1}\right\rfloor \leq n+\underline{\alpha}-1,\quad
\left\lfloor 2\kappa_{1}\right\rfloor \leq2n+2\underline{\alpha}%
-2.\label{a104}%
\end{equation}

When we choose $n+\underline{\alpha}-1\leq\kappa_{1}<n+\underline{\alpha}-1/2$
the largest values of $\left\lfloor \kappa_{1}\right\rfloor $ and
$\left\lfloor 2\kappa_{1}\right\rfloor $ are attained.

\item There exists $\alpha$ such that $w\in W3.1$ for $\theta=\left\vert
\alpha\right\vert $ and $\kappa=\kappa_{1}\mathbf{1}$ iff
\begin{equation}%
\begin{array}
[c]{ll}%
\kappa_{1}-n+1/2<\theta/d\leq l-n, & if\text{ }\operatorname{rem}\left(
\theta,d\right)  =0,\\
\kappa_{1}-n+1/2<\left\lfloor \theta/d\right\rfloor \leq l-n-1, & if\text{
}1\leq\operatorname{rem}\left(  \theta,d\right)  \leq d-1.
\end{array}
\label{p03}%
\end{equation}

If \ref{p03} holds then we can choose $\alpha=\beta$ or some permutation of
$\beta$, where $\beta$ is given by \ref{p13} in the proof. Thus
\begin{equation}
\left\lfloor \kappa_{1}\right\rfloor \leq n-1+\left\lfloor \theta
/d\right\rfloor ,\quad\left\lfloor 2\kappa_{1}\right\rfloor \leq2\left(
n-1+\left\lfloor \theta/d\right\rfloor \right)  .\label{p21}%
\end{equation}

When we choose $n+\theta/d-1\leq\kappa_{1}<n-1/2+\theta/d$ or $n+\left\lfloor
\theta/d\right\rfloor -1\leq\kappa_{1}<n-1/2+\left\lfloor \theta
/d\right\rfloor $, whichever is applicable, the largest values of
$\left\lfloor \kappa_{1}\right\rfloor $ and $\left\lfloor 2\kappa
_{1}\right\rfloor $ are attained.
\end{enumerate}
\end{theorem}

\begin{proof}
\textbf{Part 1} We apply Lemma \ref{Lem_ten_prod_wt_fn_is_W3.1}. From
\ref{p32}, $w\left(  x\right)  =\prod\limits_{i=1}^{d}\frac{x_{i}^{2n}}%
{\sin^{2l}x_{i}}$ so here $\frac{1}{w_{1}\left(  s\right)  }=\frac{\sin^{2l}%
s}{s^{2n}}$. Regarding the conditions of \ref{p02}:%
\[
\int_{\left\vert s\right\vert \leq r}\frac{ds}{s^{2\overline{\alpha}}%
w_{1}\left(  s\right)  }=\int_{\left\vert s\right\vert \leq r}\frac{\sin
^{2l}s}{s^{2n}}\dfrac{ds}{s^{2\overline{\alpha}}}=\int_{\left\vert
s\right\vert \leq r}\frac{\sin^{2l}s}{s^{2l}}\dfrac{ds}{s^{2\left(
n-l+\overline{\alpha}\right)  }},
\]

which exists iff $2\left(  n-l+\overline{\alpha}\right)  <1$ iff
$\overline{\alpha}<l-n+1/2$ iff $\overline{\alpha}\leq l-n$. Also%
\[
\int_{\left\vert s\right\vert \geq r}\dfrac{s^{2\kappa_{1}}ds}{s^{2\underline
{\alpha}}w_{1}\left(  s\right)  }=\int_{\left\vert s\right\vert \geq r}%
\dfrac{s^{2\kappa_{1}}ds}{s^{2\underline{\alpha}}w_{1}\left(  s\right)  }%
=\int_{\left\vert s\right\vert \geq r}\dfrac{\sin^{2l}sds}{s^{2\left(
n+\underline{\alpha}-\kappa_{1}\right)  }},
\]

which exists iff $2\left(  n+\underline{\alpha}-\kappa_{1}\right)  >1$ iff
$\underline{\alpha}>\kappa_{1}-n+1/2$. These two inequalities give \ref{a103}
directly and thus \ref{a104}.\medskip

\textbf{Part 2} \textbf{Suppose }\ref{p03}\textbf{\ holds}. Define
$\beta=\left(  \beta_{i}\right)  $\ by%
\begin{equation}
\beta_{i}=\left\{
\begin{array}
[c]{ll}%
1+\left\lfloor \theta/d\right\rfloor , & 1\leq i\leq\operatorname{rem}\left(
\theta,d\right)  ,\\
\left\lfloor \theta/d\right\rfloor , & \operatorname{rem}\left(
\theta,d\right)  <i\leq d.
\end{array}
\right. \label{p13}%
\end{equation}
\medskip

\fbox{If $d=1$} then $\left\lfloor \theta/d\right\rfloor =\theta$ and
$\operatorname{rem}\left(  \theta,d\right)  =0$, and so
\begin{align*}
\beta_{1}  & =\left\{
\begin{array}
[c]{ll}%
1+\left\lfloor \theta\right\rfloor , & 1\leq i\leq\operatorname{rem}\left(
\theta,1\right)  ,\\
\left\lfloor \theta\right\rfloor , & \operatorname{rem}\left(  \theta
,d\right)  <i\leq1,
\end{array}
\right. \\
& =\left\{
\begin{array}
[c]{ll}%
1+\theta, & 1\leq i\leq0,\\
\theta, & 1\leq i\leq1,
\end{array}
\right. \\
& =\theta.
\end{align*}

so $\left\vert \beta\right\vert =\theta$ and $\underline{\beta}=\overline
{\beta}=\theta=\left\lfloor \theta/d\right\rfloor $.\medskip

\fbox{If $d\geq2$} there are two cases:\medskip

\underline{\textbf{Case:} $\operatorname{rem}\left(  \theta,d\right)  =0$}
Then $\left\lfloor \theta/d\right\rfloor =\theta/d$ and so\smallskip%
\[
\beta_{i}=\left\{
\begin{array}
[c]{ll}%
1+\theta/d, & 1\leq i\leq0,\\
\theta/d, & 1\leq i\leq d,
\end{array}
\right.
\]

$\beta=\left(  \theta/d\right)  \mathbf{1}$, $\left\vert \beta\right\vert
=\theta$ and $\underline{\beta}=\overline{\beta}=\left\lfloor \theta
/d\right\rfloor $.\medskip

\underline{\textbf{Case:} $1\leq\operatorname{rem}\left(  \theta,d\right)
\leq d-1$} Here $\theta=\left\lfloor \theta/d\right\rfloor
d+\operatorname{rem}\left(  \theta,d\right)  $ so that \ref{p13}
implies\smallskip\
\begin{align*}
\left\vert \beta\right\vert  & =\left(  1+\left\lfloor \theta/d\right\rfloor
\right)  \operatorname{rem}\left(  \theta,d\right)  +\left\lfloor
\theta/d\right\rfloor \left(  d-\left(  1+\operatorname{rem}\left(
\theta,d\right)  \right)  +1\right) \\
& =\left(  1+\left\lfloor \theta/d\right\rfloor \right)  \operatorname{rem}%
\left(  \theta,d\right)  +\left\lfloor \theta/d\right\rfloor \left(
d-\operatorname{rem}\left(  \theta,d\right)  \right) \\
& =\operatorname{rem}\left(  \theta,d\right)  +\left\lfloor \theta
/d\right\rfloor d\\
& =\theta,
\end{align*}

and $\underline{\beta}=\left\lfloor \theta/d\right\rfloor $ and $\overline
{\beta}=1+\left\lfloor \theta/d\right\rfloor $.

We have thus shown that if $\beta$ is given by \ref{p13} then for all $d\geq
1$, $\left\vert \beta\right\vert =\theta$ and
\begin{equation}
\left\{
\begin{array}
[c]{ll}%
\underline{\beta}=\overline{\beta}=\left\lfloor \theta/d\right\rfloor , &
if\text{ }\operatorname{rem}\left(  \theta,d\right)  =0,\\
\underline{\beta}=\overline{\beta}-1=\left\lfloor \theta/d\right\rfloor , &
if\text{ }1\leq\operatorname{rem}\left(  \theta,d\right)  \leq d-1.
\end{array}
\right. \label{p00}%
\end{equation}

Because \ref{p03} holds, $\beta$ satisfies \ref{p00} and $\left\vert
\beta\right\vert =\theta$ which means that \ref{a103} holds for $a=\beta$ and
part 1 implies that $w\in W3.1$ for $\theta$ and $\kappa_{1}\mathbf{1}$.

\textbf{To prove the converse} assume that \ref{a103} holds for some $\alpha$
such that $\left\vert \alpha\right\vert =\theta$. We will show that the
$\alpha=\beta$ solves the problems $\max\limits_{\left\vert \alpha\right\vert
=\theta}\underline{\alpha}$ and $\min\limits_{\left\vert \alpha\right\vert
=\theta}\overline{\alpha}$ because this would mean that $\underline{\alpha
}\leq\underline{\beta}$ and $\overline{\beta}\leq\overline{\alpha}$ and hence
that \ref{a103} and \ref{p00}\ imply \ref{p03}.

If $\left\lfloor \theta/d\right\rfloor =\theta/d$ then these problems are
clearly solved by \ref{p13} i.e. by $\beta=\left(  \theta/d\right)
\mathbf{1}$. Now suppose $\left\lfloor \theta/d\right\rfloor <\theta/d$ i.e.
$\operatorname{rem}\left(  \theta,d\right)  >0$. If $\underline{\alpha}%
\geq1+\left\lfloor \theta/d\right\rfloor $ then $\left\vert \alpha\right\vert
\geq d+d\left\lfloor \theta/d\right\rfloor =d+\theta-\operatorname{rem}\left(
\theta,d\right)  >\theta$, and if $\overline{\alpha}\leq\left\lfloor
\theta/d\right\rfloor $ then $\left\vert \alpha\right\vert <d\left\lfloor
\theta/d\right\rfloor =\theta-\operatorname{rem}\left(  \theta,d\right)
<\theta$. Thus $\underline{\alpha}\leq\left\lfloor \theta/d\right\rfloor $ and
$\overline{\alpha}\geq1+\left\lfloor \theta/d\right\rfloor $ and so \ref{p13}
solves both $\max\limits_{\left\vert \alpha\right\vert =\theta}\underline
{\alpha} $ and $\min\limits_{\left\vert \alpha\right\vert =\theta}%
\overline{\alpha} $.
\end{proof}

Concerning the last two theorems:

\begin{remark}
\label{Rem_Thm_ExtNatSplin_wt_fn_W3.1}Observe that when $\kappa=\kappa
_{1}\mathbf{1}$, $\max\left\lfloor \kappa_{1}\right\rfloor $ and
$\max\left\lfloor 2\kappa_{1}\right\rfloor $ are sometimes larger when W3.1 is
used instead of W3.1*.

Specifically, under W3.1$^{\text{*}}$,%
\[%
\begin{array}
[c]{lll}%
\max\left\lfloor \kappa_{1}\right\rfloor =n-1+\theta, & \max\left\lfloor
2\kappa_{1}\right\rfloor =2\left(  n-1\right)  +2\theta, & d=1,\\
\max\left\lfloor \kappa_{1}\right\rfloor =n-1, & \max\left\lfloor 2\kappa
_{1}\right\rfloor =2\left(  n-1\right)  , & d\geq2,
\end{array}
\]

and under W3.1,%
\[%
\begin{array}
[c]{lll}%
\max\left\lfloor \kappa_{1}\right\rfloor =n-1+\theta, & \max\left\lfloor
2\kappa_{1}\right\rfloor =2\left(  n-1\right)  +2\theta, & d=1,\\
\max\left\lfloor \kappa_{1}\right\rfloor =n-1+\left\lfloor \theta
/d\right\rfloor , & \max\left\lfloor 2\kappa_{1}\right\rfloor =2\left(
n-1\right)  +2\left\lfloor \theta/d\right\rfloor , & d\geq2.
\end{array}
\]

Thus for $d=1$, $\max\left\lfloor \kappa_{1}\right\rfloor $ and $\max
\left\lfloor 2\kappa_{1}\right\rfloor $ are the same for both W3.1$^{\text{*}%
}$ and W3.1, and when $\theta\geq d$ they are larger for property W3.1.

This has smoothness consequences later in Corollary
\ref{Cor_Thm_Xwth_W3_smooth}.
\end{remark}

In Section \ref{SbSect_basis_Bsplin_W3.1} we will derive expressions for the
basis functions of positive order generated by the weight functions of the
last theorem.

\subsection{Constructing positive order weight functions from zero order
weight functions.}

Zero order weight functions were introduced in Definition \ref{Def_wt_fn zero}%
, and the positive order weight functions in Definition \ref{Def_extend_wt_fn}.

The next result shows that all the weight functions introduced in the zero
order document Williams \cite{WilliamsZeroOrdSmthV4} can be used to construct
weight functions of any positive order. Property W3.3 was motivated by part 2
of this theorem.

\begin{theorem}
\label{Thm_pos_wt_from_zero_wt}\textbf{Constructing positive order weight
functions from any zero order weight function}

\begin{enumerate}
\item Suppose the \textbf{zero order} weight function $w_{0}$ has property W02
of Definition \ref{Def_wt_fn zero} for parameter $\kappa\in\mathbb{R}^{1}$.
Then $w=\frac{w_{0}}{\left\vert \cdot\right\vert ^{2\theta}}$ is a weight
function with property W2 and also property \textbf{W3.2} for any positive
order $\theta$ and parameter $\kappa$.

\item Suppose the \textbf{zero order} weight function $w_{0}$ has property W03
of Definition \ref{Def_wt_fn zero} for parameter $\kappa\in\mathbb{R}^{d}$.
Then $w=\frac{w_{0}}{\left\vert \cdot\right\vert ^{2\theta}}$ is a weight
function with property W2, and also has property \textbf{W3.3} for $r_{3}=0$,
parameter $\kappa$ and any positive order $\theta$.
\end{enumerate}
\end{theorem}

\begin{proof}
\textbf{Part 1} Property W02 is: $\int\dfrac{\left\vert \cdot\right\vert
^{2s}}{w_{0}}<\infty$ for $0\leq s\leq\kappa$. Thus $1/w_{0}\in L_{loc}^{1}$
and so $\frac{1}{w}=\frac{\left\vert \cdot\right\vert ^{2\theta}}{w_{0}}\in
L_{loc}^{1}$ which is property W2.1. Also%
\[
\int_{\left\vert \cdot\right\vert \geq r_{3}}\frac{\left\vert \cdot\right\vert
^{2\kappa}}{w\left\vert \cdot\right\vert ^{2\theta}}=\int_{\left\vert
\cdot\right\vert \geq r_{3}}\frac{\left\vert \cdot\right\vert ^{2\kappa}%
}{w_{0}}<\infty,
\]

so property W3.2 holds and thus by Theorem \ref{Thm_weight_property_relat},
property W2.2 holds.\medskip

\textbf{Part 2} Property W03 is: $\int\dfrac{x^{2\lambda}}{w_{0}}<\infty$ for
$0\leq\lambda\leq\kappa$. Thus $1/w_{0}\in L_{loc}^{1}$ and so $\frac{1}%
{w}=\frac{\left\vert \cdot\right\vert ^{2\theta}}{w_{0}}\in L_{loc}^{1}$ which
is property W2.1. Also, if $0\leq\lambda\leq\kappa$,%
\[
\int\frac{x^{2\lambda}}{w\left(  x\right)  \left\vert x\right\vert ^{2\theta}%
}dx=\int\frac{x^{2\lambda}}{w_{0}}<\infty,
\]

so property W3.2 holds and thus property W2.2.
\end{proof}

\section{The spaces $S_{\emptyset,n}$ and related spaces
\label{Sect_So,n_S'o,n}}

We now define some key spaces and exhibit some of their properties.

\begin{definition}
\label{Def_So,n}\textbf{The spaces} $S_{\emptyset,n}$ \textbf{and}
$C_{\emptyset,n}^{\infty}$
\begin{equation}
S_{\emptyset,0}=S,\qquad S_{\emptyset,n}=\left\{  \phi\in S:D^{\alpha}%
\phi\left(  0\right)  =0,\text{\ }\left\vert \alpha\right\vert <n\right\}
,\;n=1,2,3,\ldots,\label{p60}%
\end{equation}

and we endow $S_{\emptyset,n}$ with the subspace topology induced by the space
$S$. $S$ is the space of $C^{\infty}$ functions of rapid decrease used as test
functions for the tempered distributions. $S$ is endowed with the countable
seminorm topology described in Appendix \ref{Def_Distributions}.%
\[
C_{\emptyset,0}^{\infty}=C^{\infty},\qquad C_{\emptyset,n}^{\infty}=\left\{
\phi\in C^{\infty}:D^{\beta}\phi\left(  0\right)  =0,\text{\ }|\beta
|<n\right\}  ,\;n=1,2,3,\ldots,
\]

so the space $C_{\emptyset,n}^{\infty}$ retains the constraints of
$S_{\emptyset,n}$ near the origin.
\end{definition}

The next result gives a simple upper bound for functions in $S_{\emptyset,n} $
near the origin. This follows directly from the estimate \ref{p20} (Appendix
\ref{Sect_taylor_expansion}) of the integral remainder term of the Taylor
series expansion.%
\begin{equation}
\left\vert u(x)\right\vert \leq\left(  \sum\limits_{\left\vert \alpha
\right\vert =n}\left\Vert D^{\alpha}u\right\Vert _{\infty}\right)  \left\vert
x\right\vert ^{n},\quad u\in S_{\emptyset,n},\text{ }x\in\mathbb{R}%
^{d}.\label{p33}%
\end{equation}

Noting that $C_{BP}^{\infty}$ is the space of $C^{\infty}$ functions for which
each derivative is bounded by a polynomial we will need the following useful results:

\begin{theorem}
\label{Thm_product_of_Co,k_funcs}\ 

\begin{enumerate}
\item If $f,g\in C_{BP}^{\infty}$ and $\psi\in S$ then $f\psi\in S$ and $fg\in
C_{BP}^{\infty}$.

\item $S_{\emptyset,k}=S\cap C_{\emptyset,k}^{\infty}$.

\item For $k,l\geq0$, $\phi\in C_{\emptyset,k}^{\infty}$ and $\psi\in
C_{\emptyset,l}^{\infty}$ implies that $\phi\psi\in C_{\emptyset,k+l}^{\infty
}$.

\item If $\left\vert \alpha\right\vert =k$ then $x^{\alpha}\in C_{\emptyset
,k}^{\infty}$.

\item If $k$ is a non-negative integer then $\left\vert \cdot\right\vert
^{2k}\in C_{\emptyset,2k}^{\infty}$.
\end{enumerate}
\end{theorem}

In Appendix \ref{Def_Some_basic_spaces} we defined $P_{n}$ to be the space of
polynomials of degree at most $n$ with complex coefficients, and $P$ to be the
space of all polynomials with complex coefficients.

\begin{notation}
\label{Not_Fourier_P}\textbf{Fourier transforms of polynomial spaces}%
\[
\widehat{P}_{n}=\left\{  \widehat{p}:p\in P_{n}\right\}  ,\qquad\widehat
{P}=\left\{  \widehat{p}:p\in P\right\}  .
\]

\end{notation}

The next theorem proves two relationships between the tempered distributions
$\widehat{P_{n-1}}$ and the functions $S_{\emptyset,n}$.

\begin{theorem}
\label{Thm_So,n_and_Pnhat}Suppose that $n$ is a non-negative integer and $u\in
S^{\prime}$. Then using Notation \ref{Not_Fourier_P}:

\begin{enumerate}
\item $u\in\widehat{P}_{n-1}\ $iff$\ \left[  u,\phi\right]  =0$ for all
$\phi\in S_{\emptyset,n}$.

\item $u\in\widehat{P}_{n-1}$\ iff\ $\phi u=0$ for all $\phi\in S_{\emptyset
,n}$.
\end{enumerate}
\end{theorem}

\begin{proof}
\textbf{Part 1} Suppose $\left[  u,\phi\right]  =0$ for all $\phi\in
S_{\emptyset,n}$. This implies that the $\operatorname*{supp}u\subset\left\{
0\right\}  $ and by a well known theorem in distribution theory, $u\in
\widehat{P}$. Thus $u=\widehat{p}$ for some polynomial $p$. Suppose $\deg p>n
$. For each $\left\vert \beta\right\vert >n$ choose $\phi_{\beta}\in
S_{\emptyset,n}$ such that $\phi_{\beta}\left(  x\right)  =x^{\beta}/\beta!$
in a neighborhood of zero. Then we have $\left(  D^{\alpha}\phi_{\beta
}\right)  \left(  0\right)  =\delta_{\alpha,\beta}$ and if the coefficients of
$p$ are $b_{\alpha}$,

$0=\left[  \widehat{p},\phi_{\beta}\right]  =\left[  p(-D)\left(  \phi_{\beta
}\right)  \right]  (0)=\left(  -1\right)  ^{\left\vert \beta\right\vert
}b_{\beta}$, and thus $\deg p<n$ and $u\in\widehat{P}_{n-1}$.\medskip

\textit{Conversely}, suppose $u\in\widehat{P}_{n-1}$. Then there exists $p\in
P_{n-1}$ such that $u=\widehat{p}=\left(  2\pi\right)  ^{^{\frac{d}{2}}%
}p(iD)\delta$. Hence, if $\phi\in S_{\emptyset,n}$
\begin{align*}
\left[  u,\phi\right]  =\left[  \left(  2\pi\right)  ^{\frac{d}{2}}%
p(iD)\delta,\phi\right]  =\left(  2\pi\right)  ^{^{\frac{d}{2}}}\left[
p(iD)\delta,\phi\right]  =\left(  2\pi\right)  ^{^{\frac{d}{2}}}\left[
p(iD)\delta,\phi\right]   &  =\left(  2\pi\right)  ^{^{\frac{d}{2}}}\left[
\delta,p(-iD)\phi\right] \\
&  =\left(  2\pi\right)  ^{^{\frac{d}{2}}}\left[  p(-iD)\phi\right]  \left(
0\right) \\
&  =0.
\end{align*}
\medskip

\textbf{Part 2} Suppose that $\phi u=0$ for all $\phi\in S_{\emptyset,n}$.
Then for all $\psi\in S$
\[
\left[  \phi u,\psi\right]  =\left[  u,\phi\psi\right]  =\left[  \psi
u,\phi\right]  =0,
\]

Thus $\left[  \psi u,\phi\right]  =0$ for all $\phi\in S_{\emptyset,n}$ and so
by part 1, $\psi u\in\widehat{P}_{n-1}$. Choose a point $x_{0}\in
\mathbb{R}^{d}$ and an open ball containing $x_{0}$ but not the origin. Choose
$\psi\in S$ so that $\psi=1$ on this ball. Then $u=\psi u\in\widehat{P}_{n-1}$
in this neighborhood and so $u=0$ in this neighborhood. Thus $u=0$ on
$\mathbb{R}^{d}\setminus0$. Choose an open ball containing the origin and
choose $\psi\in S$ so that $\psi=1$ on this ball. Then $u=\psi u\in\widehat
{P}_{n-1}$ in this neighborhood and so $u\in\widehat{P}_{n-1}$ in this
neighborhood, say $u=\widehat{p}$ where $p\in P_{n-1}$. But then $\widehat
{p}=0$ on $\mathbb{R}^{d}\setminus0$. Hence $u=\widehat{p}$ on $\mathbb{R}%
^{d}$.\medskip

Conversely, suppose $u\in\widehat{P}_{n-1}$. Then, if $\phi\in S_{\emptyset
,n}$ and $\psi\in S$ it follows that $\phi\psi\in S_{\emptyset,n}$. Thus by
part 1 of this theorem $\left[  \phi u,\psi\right]  =\left[  u,\phi
\psi\right]  =0$.
\end{proof}

\subsection{The functionals $S_{\emptyset,n}^{\prime}$\label{SbSect_S'on}}

The following lemma will enable us to define the space of continuous
functionals $S_{\emptyset,n}^{\prime}$ which in turn will allow us to define
the basis distributions in Section \ref{Sect_basis_distrib} and the operator
$\mathcal{J}$ of Section \ref{SbSect_I_J}.

\begin{lemma}
\label{Lem_convex_tls_extend_2}\ 

\begin{enumerate}
\item The space of rapidly decreasing functions $S$ is a locally convex
topological space when endowed with the countable seminorm topology defined in
Appendix \ref{Def_Distributions}.

\item A linear functional $f$ defined on a subspace $\mathbb{M}$ of the space
$S$ is continuous iff there exists an integer $n\geq0$ and a constant $C$ such
that%
\[
\left\vert \left[  f,\psi\right]  \right\vert \leq C\left\Vert \sum
\limits_{\left\vert \alpha\right\vert \leq m}\left(  1+\left\vert
\cdot\right\vert \right)  ^{m}\left\vert D^{\alpha}\psi\right\vert \right\Vert
_{\infty},\quad\phi\in\mathbb{M}.
\]

We write $f\in\mathbb{M}^{\prime}$.

\item Any continuous linear functional $f$ on a subspace $\mathbb{M}$ of a
locally convex topological vector space $\mathbb{T}$ can be extended
non-uniquely to a continuous linear functional $f^{e}$ on $\mathbb{T}$
($f^{e}\in\mathbb{T}^{\prime}$).

The set of extensions is $f^{e}+\mathbb{M}^{\bot}$ where $\mathbb{M}^{\bot} $
is the set of annihilators of $\mathbb{M}$ i.e. the members of $\mathbb{T}%
^{\prime}$ which are zero on each member of $\mathbb{M}$.
\end{enumerate}
\end{lemma}

\begin{proof}
This lemma can be proved, for example, by using the results and definitions of
Chapter V., Volume I of Reed and Simon \cite{ReedSimon72}.
\end{proof}

Recall Definition \ref{Def_So,n} of the space $S_{\emptyset,n}$.

\begin{definition}
\label{Def_functnl_on_Son}\textbf{The spaces} $S_{\emptyset,n}^{\prime},$
$n=1,2,3,\ldots$.

$S_{\emptyset,n}^{\prime}$ is the space of continuous linear functionals on
$S_{\emptyset,n}$, where $S_{\emptyset,n}$ has the subspace topology induced
by $S$.
\end{definition}

The next theorem gives the criterion we shall use to prove that a linear
functional on $S_{\emptyset,n}$ is continuous i.e. is a member of
$S_{\emptyset,n}^{\prime}$. It also supplies the key extension result needed
to define the basis distributions.

\begin{theorem}
\label{Thm_prop_functnl_on_Son}\ 

\begin{enumerate}
\item A linear functional $f$ on $S_{\emptyset,n}$ is member of $S_{\emptyset
,n}^{\prime}$ iff for some integer $m\geq0$ and some constant $C>0$
\[
\left\vert \left[  f,\psi\right]  \right\vert \leq C\left\Vert \sum
\limits_{\left\vert \alpha\right\vert \leq m}\left(  1+\left\vert
\cdot\right\vert \right)  ^{m}\left\vert D^{\alpha}\psi\right\vert \right\Vert
_{\infty},\quad\phi\in S_{\emptyset,n}.
\]

\item If $f\in S_{\emptyset,n}^{\prime}$ then there exists $f^{e}\in
S^{\prime}$ such that $f^{e}=f$ on $S_{\emptyset,n}$. The set of extensions is
$f^{e}+\widehat{P_{n-1}}$. We sometimes say that $f$ can be extended to $S$ as
a member $f^{e}$ of $S^{\prime}$.
\end{enumerate}
\end{theorem}

\begin{proof}
\textbf{Part 1} follows from part 2 of Lemma \ref{Lem_convex_tls_extend_2}
with $\mathbb{M}=S_{\emptyset,n}$.\medskip

\textbf{Part 2} We use part 3 of Lemma \ref{Lem_convex_tls_extend_2} when
$\mathbb{M}=S_{\emptyset,n}$ and $\mathbb{T}=S$. In this case the the set of
extensions is $f^{e}+S_{\emptyset,n}^{\bot}$ and $g\in S_{\emptyset,n}^{\bot}$
iff $g\in S^{\prime}$ and $\left[  g,\phi\right]  =0$ for all $\phi\in
S_{\emptyset,n}$ i.e. $g\in\widehat{P_{n-1}}$ by part 1 of Theorem
\ref{Thm_So,n_and_Pnhat}.
\end{proof}

\section{The data (native) spaces $X_{w}^{\theta},$ $\theta=1,2,3,..$%
\label{Sect_Xwo_Xwth}}

\subsection{Introduction}

In Section 2 of \cite{LightWayneX98Weight} Light and Wayne introduced the
reproducing kernel semi-Hilbert spaces $X$ (native spaces, see Schaback
\cite{Schaback99}, Wendland \textit{\cite{Wendland05}} etc.) to formulate the
minimal seminorm interpolation problem. The space $X$ was defined using a
positive weight function $w$ and had positive order $\theta$. In fact, the
spaces $X$ can be described as Beppo-Levi spaces \cite{DenyLions54}
generalized using a positive weight function. In Section \ref{Sect_Xwo_Xwth}
we extended the class of weight functions to include functions which were
positive and continuous except on a set of measure zero. In this section we
will use the same definition for our function spaces as Light and Wayne but I
will use the notation $X_{w}^{\theta}$ for the $X$ spaces of positive order
$\theta$ generated by the weight function $w$. Only weight function property
W1 is required to define $X_{w}^{\theta}$. I will also introduce an
alternative definition that I have found quite useful.

By analogy with Sobolev space theory mappings will be constructed between
$X_{w}^{\theta}$ and $L^{2}$ in order to prove various properties of
$X_{w}^{\theta}$. To prove that $X_{w}^{\theta}$ is a semi-Hilbert space we
construct the mappings $\mathcal{I}:X_{w}^{\theta}\rightarrow L^{2}$ and
$\mathcal{J}:L^{2}\rightarrow X_{w}^{\theta}$ which are adjoints, inverses and
isometric isomorphisms in the seminorm sense. Here weight function property W2
is also required. Finally we prove some smoothness results for $X_{w}^{\theta
}$ when the weight function has properties W2 and W3.

\begin{definition}
\label{Def_Xwth}\textbf{The semi-inner product data (native) spaces}
$X_{w}^{\theta}$ \textbf{of order} $\theta=1,2,3,\ldots$.

Suppose $w$ is a weight function i.e. it has property W1. We now define
$X_{w}^{\theta}$ by%
\begin{equation}
X_{w}^{\theta}=\left\{  f\in S^{\prime}:\widehat{D^{\alpha}f}\in L_{loc}%
^{1}\left(  \mathbb{R}^{d}\right)  ,\text{ }\int w\left\vert \widehat
{D^{\alpha}f}\right\vert ^{2}<\infty\text{ }for\text{ }all\text{ }\left\vert
\alpha\right\vert =\theta\right\}  ,\label{p37}%
\end{equation}

and endow $X_{w}^{\theta}$ with a semi-inner product and seminorm (part 3
Theorem \ref{Thm_properties_Xwm})
\begin{equation}
\left\langle f,g\right\rangle _{w,\theta}=\sum\limits_{\left\vert
\alpha\right\vert =\theta}\frac{\theta!}{\alpha!}\int w\widehat{D^{\alpha}%
f}\,\overline{\widehat{D^{\alpha}g}},\text{\quad}\left\vert f\right\vert
_{w,\theta}=\sqrt{\left\langle f,f\right\rangle _{w,\theta}}.\label{p61}%
\end{equation}

\end{definition}

To prove that $X_{w}^{\theta}$ is non-empty we need the following lemma:

\begin{lemma}
\label{Lem_fe=1_on_set}Let $\mathcal{F}$ be a set of points in $\mathbb{R}%
^{d}$. Then for any $\eta>0$ there exists a function $f_{\eta}$ such that:

\begin{enumerate}
\item $f_{\eta}\in C^{\infty}$,

\item $0\leq f_{\eta}\left(  x\right)  \leq1$,

\item $f_{\eta}\left(  x\right)  =1$ when $x\in\mathcal{F}_{\eta}$,

\item $f_{\eta}\left(  x\right)  =0$ when $x\notin\mathcal{F}_{3\eta}$.
\end{enumerate}

Here $\mathcal{F}_{\eta}$ and $\mathcal{F}_{3\eta}$ denote $\eta-$ and
$3\eta-$neighborhoods of the set $\mathcal{F},$ as defined in Appendix
\ref{Def_topol_notation}.
\end{lemma}

\begin{proof}
Choose $\eta>0$. There exists a mollifier $\omega\in C_{0}^{\infty}$
satisfying
\[
\operatorname*{supp}\omega\subseteq B\left(  0;1\right)  ;\quad\int%
\omega=1;\quad\omega\geq0.
\]

Now define the scaled function $\omega_{\eta}\left(  x\right)  =\eta
^{-d}\omega\left(  x/\eta\right)  $ and let $\chi_{2\eta}$ be the
characteristic function of the set $\mathcal{F}_{2\eta}$. Then it can be shown
that the function $f_{\eta}\left(  x\right)  =\int\chi_{2\eta}\left(
y\right)  \omega_{\eta}\left(  x-y\right)  dy$ has the required properties.
\end{proof}

\begin{theorem}
\label{Thm_Xwth_non_empty}Suppose the weight function $w$ has property W1
w.r.t. the set $\mathcal{A}$. Then
\begin{equation}
\left\{  \left(  u\left\vert \cdot\right\vert ^{-\theta}\right)  ^{\vee}:u\in
C_{0}^{\infty}\text{ }and\text{ }\operatorname*{supp}u\subset\mathbb{R}%
^{d}\setminus\left(  \mathcal{A}\cup\left\{  0\right\}  \right)  \right\}
\subset X_{w}^{\theta},\label{p80}%
\end{equation}

and the set on the left is non-empty.
\end{theorem}

\begin{proof}
If $\Omega$ is an open set we define $C_{0}^{\infty}\left(  \Omega\right)
=\left\{  \phi\in C_{0}^{\infty}:\operatorname*{supp}\phi\subset
\Omega\right\}  $. A weight set is a closed set of measure zero so
$\mathcal{B}=\mathcal{A}\cup\left\{  0\right\}  $ is a closed set of measure
zero. Then there exists $\eta>0$ such that the open set $\mathbb{R}%
^{d}\setminus\mathcal{B}_{\eta}$ is non-empty where $\mathcal{B}_{\eta}$ is
the neighborhood set of $\mathcal{B}$ defined by $\eta$.

We now use Lemma \ref{Lem_fe=1_on_set} with $\mathcal{F}=\mathcal{B}$ to show
the set on the left of \ref{p80} is non-empty. Fix $\eta>0$. Then there exists
a function $f_{\eta}$ with properties 1 to 4 of Lemma \ref{Lem_fe=1_on_set}.
Therefore, $v\in C_{0}^{\infty}$ implies $\left(  1-f_{\eta}\right)  v\in
C_{0}^{\infty}$ and $\left(  1-f_{\eta}\right)  v=0$ on $\mathcal{B}_{\eta}$
i.e. $\left(  1-f_{\eta}\right)  v\subset C_{0}^{\infty}\left(  \mathbb{R}%
^{d}\setminus\mathcal{B}\right)  $. We conclude that the set on the left of
\ref{p80} is non-empty.

To prove the inclusion \ref{p80} holds let $f=\left(  \frac{u}{\left\vert
\cdot\right\vert ^{\theta}}\right)  ^{\vee}$ where $u\in C_{0}^{\infty}\left(
\mathbb{R}^{d}\setminus\mathcal{B}\right)  $. Because $0\in\mathcal{B}$ it
follows that $\frac{1}{\left\vert \cdot\right\vert ^{\theta}}\in C^{\infty
}\left(  \mathbb{R}^{d}\setminus\mathcal{B}\right)  $ and $\frac{u}{\left\vert
\cdot\right\vert ^{\theta}}\in C_{0}^{\infty}\left(  \mathbb{R}^{d}%
\setminus\mathcal{B}\right)  $. Thus $f\in S^{\prime}$, $\xi^{\alpha}%
\widehat{f}=\xi^{\alpha}\frac{u}{\left\vert \cdot\right\vert ^{\theta}}\in
L_{loc}^{1}$ and $\int w\left\vert \cdot\right\vert ^{2\theta}\left\vert
\widehat{f}\right\vert ^{2}=\int\limits_{\operatorname*{supp}u}w\left\vert
u\right\vert ^{2}$. But $w\in C^{\left(  0\right)  }\left(  \mathbb{R}%
^{d}\setminus\mathcal{A}\right)  $ so $w\in C^{\left(  0\right)  }\left(
\mathbb{R}^{d}\setminus\mathcal{B}\right)  $ and so, $w\left\vert u\right\vert
^{2}=wu\overline{u}$ is continuous on $\operatorname*{supp}u$. This means that
$\int\limits_{\operatorname*{supp}u}w\left\vert u\right\vert ^{2}<\infty$ and
so $\int w\left\vert \cdot\right\vert ^{2\theta}\left\vert \widehat
{f}\right\vert ^{2}<\infty$. We conclude that $f\in X_{w}^{\theta}$.
\end{proof}

\subsection{Some properties of the data spaces $X_{w}^{\theta}$}

The space $X_{w}^{\theta}$ is obviously contained in the space%
\begin{equation}
X^{\theta}:=\left\{  f\in S^{\prime}:\widehat{D^{\alpha}f}\in L_{loc}%
^{1}\text{ }whenever\text{ }\left\vert \alpha\right\vert =\theta\right\}
.\label{p66}%
\end{equation}

We will now prove some basic properties of this space of tempered distributions.

\begin{lemma}
\label{Lem_properties_Xwm_distrib}The distributions in the space $X^{\theta} $
of \ref{p66} have the properties:

\begin{enumerate}
\item $\widehat{f}\left(  \xi\right)  =\dfrac{\left(  -i\right)  ^{\theta}%
}{\left\vert \xi\right\vert ^{2\theta}}\sum\limits_{\left\vert \alpha
\right\vert =\theta}\dfrac{\theta!}{\alpha!}\xi^{\alpha}\widehat{D^{\alpha}%
f}\left(  \xi\right)  $,\quad$\xi\neq0$.

\item $\widehat{f}\in L_{loc}^{1}\left(  \mathbb{R}^{d}\setminus0\right)  $.

\item $\left\vert \widehat{f}\left(  \xi\right)  \right\vert \leq\dfrac
{1}{\left\vert \xi\right\vert ^{\theta}}\sum\limits_{\left\vert \alpha
\right\vert =\theta}\dfrac{\theta!}{\alpha!}\left\vert \widehat{D^{\alpha}%
f}\left(  \xi\right)  \right\vert $,\quad$\xi\neq0$.

\item Define the function $f_{F}$ a.e. on $\mathbb{R}^{d}$ by\textbf{:}
$f_{F}=\widehat{f}$ on $\mathbb{R}^{d}\setminus0$.

Then $\left\vert \cdot\right\vert ^{\theta}f_{F}\in L_{loc}^{1}$ and
$\xi^{\alpha}f_{F}\in L_{loc}^{1}$ when $\left\vert \alpha\right\vert =\theta$.

\item Also $\xi^{\alpha}f_{F}=\xi^{\alpha}\widehat{f}$ a.e. when $\left\vert
\alpha\right\vert =\theta$.
\end{enumerate}
\end{lemma}

\begin{proof}
\textbf{Part 1} is proved by using the second identity of part 5 of Summary
\ref{Sum_Multi_index} to write

$\left\vert \xi\right\vert ^{2\theta}\widehat{f}\left(  \xi\right)  =\left(
\sum\limits_{\left\vert \alpha\right\vert =\theta}\dfrac{\theta!}{\alpha!}%
\xi^{2\alpha}\right)  \widehat{f}\left(  \xi\right)  $. \textbf{Parts 2 to 4}
have easy proofs and these have been omitted. However, some readers may find
the definition of the the function $f_{F}$ in part 4 confusing. The function
$f_{F}$ is \textbf{not} a function on $\mathbb{R}^{d}$ in the distribution
sense i.e. it is not an $L_{loc}^{1}\left(  \mathbb{R}^{d}\right)  $ function.
It is a function in the elementary sense which happens to be a member of
$L_{loc}^{1}\left(  \mathbb{R}^{d}\setminus0\right)  $ and being a member of
$L_{loc}^{1}\left(  \mathbb{R}^{d}\setminus0\right)  $ defines a function a.e.
on $\mathbb{R}^{d}\setminus0$ which defines a function a.e. on $\mathbb{R}%
^{d}$.\smallskip

\textbf{Part 5} True since $\xi^{\alpha}f_{F}\in L_{loc}^{1}$, $\xi^{\alpha
}\widehat{f}\in L_{loc}^{1}$, $\xi^{\alpha}f_{F}=\xi^{\alpha}\widehat{f}$ on
$\mathbb{R}^{d}\setminus0$ so $\xi^{\alpha}f_{F}=\xi^{\alpha}\widehat{f}$ a.e.
\end{proof}

The next theorem derives some properties of $X_{w}^{\theta}$ and of the
function $f_{F}$ defined in part 4 of the previous lemma. The proof of the
formula \ref{p10} for $\left\vert f\right\vert _{w,\theta}$ reveals the
motivation for using the coefficients $\frac{\theta!}{\alpha!}$ to define the
seminorm $\left\vert f\right\vert _{w,\theta}$. This formula will prove very
useful because of its algebraic simplicity. Part 2 of the next theorem
provides a way of verifying that a function is in $X_{w}^{\theta}$ by making
use of this formula.

\begin{theorem}
\label{Thm_properties_Xwm}Suppose $w$ is a weight function i.e. $w$ has
property W1. If $f\in X_{w}^{\theta}$ for some integer $\theta\geq1$ then
$\widehat{f}\in L_{loc}^{1}\left(  \mathbb{R}^{d}\setminus0\right)  $ and we
can a.e. define the function $f_{F}:\mathbb{R}^{d}\rightarrow\mathbb{C}$ by
$f_{F}=\widehat{f}$ on $\mathbb{R}^{d}\setminus0$. Further:

\begin{enumerate}
\item The seminorm \ref{p61} satisfies
\begin{equation}
\int w\left\vert \cdot\right\vert ^{2\theta}\left\vert f_{F}\right\vert
^{2}=\left\vert f\right\vert _{w,\theta}^{2}.\label{p10}%
\end{equation}

\item An \textbf{alternative definition} of $X_{w}^{\theta}$ is:
\begin{equation}
X_{w}^{\theta}=\left\{  f\in S^{\prime}:\xi^{\alpha}\widehat{f}\in L_{loc}%
^{1}\text{ }if\text{ }\left\vert \alpha\right\vert =\theta;\text{\thinspace
}\int w\left\vert \cdot\right\vert ^{2\theta}\left\vert f_{F}\right\vert
^{2}<\infty\right\}  .\label{p40}%
\end{equation}

This definition actually makes sense because by part 2 of Corollary
\ref{Lem_properties_Xwm_distrib} the condition $\xi^{\alpha}\widehat{f}\in
L_{loc}^{1}$ for all $\left\vert \alpha\right\vert =\theta$ implies that
$\widehat{f}\in L_{loc}^{1}\left(  \mathbb{R}^{d}\setminus0\right)  $, and so
$f_{F}$ can be defined a.e. on $\mathbb{R}^{d}$.

\item The functional $\left\vert \cdot\right\vert _{w,\theta}$ is a seminorm.
In fact $\operatorname*{null}\left\vert \cdot\right\vert _{w,\theta}%
=P_{\theta-1}$ and we also have $X_{w}^{\theta}\cap P=P_{\theta-1}$.

\item Suppose $X^{\theta}$ is the space defined by \ref{p66}. Then%
\[
X_{w}^{\theta}=\left\{  f\in X^{\theta}:\int w\left\vert \cdot\right\vert
^{2\theta}\left\vert f_{F}\right\vert ^{2}<\infty\right\}  ,
\]

where the space $X_{w}^{\theta}$ is independent of the weight function.
\end{enumerate}
\end{theorem}

\begin{proof}
\textbf{Part 1} Since $f_{F}$ is the function defined in part 4 of Lemma
\ref{Lem_properties_Xwm_distrib} we have from part 5 of the same lemma that
$\xi^{\alpha}f_{F}\in L_{loc}^{1}$, $\xi^{\alpha}\widehat{f}\in L_{loc}^{1}$
and $\xi^{\alpha}f_{F}=\xi^{\alpha}\widehat{f}$ a.e. Hence
\begin{align*}
\int w\left\vert \cdot\right\vert ^{2\theta}\left\vert f_{F}\right\vert
^{2}=\int w\left(  \sum\limits_{\left\vert \alpha\right\vert =\theta}%
\frac{\theta!}{\alpha!}\xi^{2\alpha}\right)  \left\vert f_{F}\right\vert
^{2}=\sum\limits_{\left\vert \alpha\right\vert =\theta}\frac{\theta!}{\alpha
!}\int w\left\vert \xi^{\alpha}f_{F}\right\vert ^{2} &  =\sum
\limits_{\left\vert \alpha\right\vert =\theta}\frac{\theta!}{\alpha!}\int
w\left\vert \xi^{\alpha}\widehat{f}\right\vert ^{2}\\
&  =\sum\limits_{\left\vert \alpha\right\vert =\theta}\frac{\theta!}{\alpha
!}\int w\left\vert \widehat{D^{\alpha}f}\right\vert ^{2}\\
&  =\left\vert f\right\vert _{w,\theta}^{2}<\infty,
\end{align*}

since $f\in X_{w}^{\theta}$.\medskip

\textbf{Part 2} Now let $\mathbb{U}$ be the space defined by the right hand
side of equation \ref{p40}.\smallskip

\qquad\textbf{Step 1} Prove that $X_{w}^{\theta}\subset\mathbb{U}$:
By\textbf{\ }part 1 of\ this theorem, if $f\in X_{w}^{\theta}$ then $\int
w\left\vert \cdot\right\vert ^{2\theta}\left\vert f_{F}\right\vert
^{2}=\left\vert f\right\vert _{w,\theta}^{2}<\infty.$ Also, $f\in
X_{w}^{\theta}$ implies $\widehat{f}\in L_{loc}^{1}\left(  \mathbb{R}%
^{d}\setminus0\right)  $ and the other requirements for $f$ to be in
$\mathbb{U}$ follows from Lemma \ref{Lem_properties_Xwm_distrib}.

\qquad\textbf{Step 2} Prove that $\mathbb{U}\subset X_{w}^{\theta}$: Assume
$f\in\mathbb{U}$. We need to show that, $\sum\limits_{\left\vert a\right\vert
=\theta}\dfrac{\theta!}{\alpha!}\int w\left\vert \widehat{D^{\alpha}%
f}\right\vert ^{2}<\infty$ when $\left\vert \alpha\right\vert =$ $\theta$. But
from the proof of part 1, $\xi^{\alpha}\widehat{f}=\xi^{\alpha}f_{F}\in
L_{loc}^{1}$ and so
\begin{align*}
\sum\limits_{\left\vert a\right\vert =\theta}\frac{\theta!}{\alpha!}\int
w\left\vert \widehat{D^{\alpha}f}\right\vert ^{2}=\sum\limits_{\left\vert
a\right\vert =\theta}\frac{\theta!}{\alpha!}\int w\left\vert \xi^{\alpha
}\widehat{f}\right\vert ^{2}=\sum\limits_{\left\vert a\right\vert =\theta
}\frac{\theta!}{\alpha!}\int w\left\vert \xi^{\alpha}f_{F}\right\vert ^{2} &
=\sum\limits_{\left\vert a\right\vert =\theta}\frac{\theta!}{\alpha!}\int
w\xi^{2\alpha}\left\vert f_{F}\right\vert ^{2}\\
&  =\int w\left\vert \cdot\right\vert ^{2\theta}\left\vert f_{F}\right\vert
^{2}<\infty.
\end{align*}
\medskip

\textbf{Part 3} Suppose that $\left\vert f\right\vert _{w,\theta}=0$. Then
$f\in S^{\prime}$, $\widehat{D^{\alpha}f}\in L_{loc}^{1}$ for $\left\vert
\alpha\right\vert =$ $\theta$, and so
\[
\left\vert f\right\vert _{w,\theta}^{2}=\sum_{\left\vert \alpha\right\vert
=\theta}\frac{\theta!}{\alpha!}\int w\left\vert \widehat{D^{\alpha}%
f}\right\vert ^{2}=0.
\]

Now $w>0$ a.e. implies that when $\left\vert \alpha\right\vert =$ $\theta$,
$\widehat{D^{\alpha}f}=0$ a.e. and hence $D^{\alpha}f=0$ a.e. But by a well
known result from distribution theory e.g. corollary to Theorem VI, p.60
Schwartz \cite{Schwartz66}, we can conclude that $f\in P_{\theta-1}$. Finally,
it is clear that $f\in P_{\theta-1}$ implies $\left\vert f\right\vert
_{w,\theta}=0$. Thus we have proved $\operatorname*{null}\left\vert
\cdot\right\vert _{w,\theta}=P_{\theta-1}$.

Clearly $P_{\theta-1}\subset X_{w}^{\theta}\cap P$. To prove the converse,
suppose $f\in X_{w}^{\theta}\cap P$ and $f\notin P_{\theta-1}$. Then
$\left\vert \alpha\right\vert =$ $\theta$ implies $D^{\alpha}f\in
P\setminus\mathbf{0}$ and $\widehat{D^{\alpha}f}\notin L_{loc}^{1}$, but $f\in
X_{w}^{\theta}$ which implies $\widehat{D^{\alpha}f}\in L_{loc}^{1}$, a
contradiction. Thus $X_{w}^{\theta}\cap P=P_{\theta-1}$.\medskip

\textbf{Part 4} A direct consequence of \ref{p40}.
\end{proof}

The function $f_{F}$, introduced in the previous theorem, is of central
importance to the theory of this document. We now derive some important
properties of this function. To do this we need the following lemma:

\begin{lemma}
\label{Lem_functnal_phi_sq}Suppose the weight function $w$ has property W2.
Then for any integer $\theta\geq1$ there exists a constant $c_{r_{2},\theta}$
independent of $\phi\in S_{\emptyset,\theta}$ such that%
\[
\left(  \int\frac{\left\vert \phi\right\vert ^{2}}{w\left\vert \cdot
\right\vert ^{2\theta}}\right)  ^{1/2}\leq c_{r_{2},\theta}\left(  \int%
\frac{1}{w\left\vert \cdot\right\vert ^{2\lambda\left(  \cdot\right)  }%
}\right)  ^{1/2}\sum\limits_{\left\vert \alpha\right\vert \leq n}\left\Vert
\left(  1+\left\vert \cdot\right\vert \right)  ^{n}D^{\alpha}\phi\right\Vert
_{\infty},
\]

where $n=\operatorname{ceil}\left\{  \theta,\sigma\right\}
=\operatorname{ceil}\max\left\{  \theta,\sigma\right\}  $. Here $\lambda$ is
the function introduced in the definition of weight function property W2.
\end{lemma}

\begin{proof}
Suppose $r_{2}$ is the parameter in the definition of weight function property
W2. Then for $\phi\in S_{\emptyset,\theta}$ we write%
\begin{align}
\int\frac{\left\vert \phi\right\vert ^{2}}{w\left\vert \cdot\right\vert
^{2\theta}}\leq\int\limits_{\left\vert \cdot\right\vert \leq r_{2}}%
\frac{\left\vert \phi\right\vert ^{2}}{w\left\vert \cdot\right\vert ^{2\theta
}}+\int\limits_{\left\vert \cdot\right\vert \geq r_{2}}\frac{\left\vert
\phi\right\vert ^{2}}{w\left\vert \cdot\right\vert ^{2\theta}}= &
\int\limits_{\left\vert \cdot\right\vert \leq r_{2}}\frac{\left\vert
\phi\right\vert ^{2}}{\left\vert \cdot\right\vert ^{2\theta}}\frac{1}{w}%
+\int\limits_{\left\vert \cdot\right\vert \geq r_{2}}\frac{\left\vert
\cdot\right\vert ^{2\sigma}\left\vert \phi\right\vert ^{2}}{\left\vert
\cdot\right\vert ^{2\theta}}\frac{1}{w\left\vert \cdot\right\vert ^{2\sigma}%
}\nonumber\\
&  \leq\left\Vert \frac{\phi^{2}}{\left\vert \cdot\right\vert ^{2\theta}%
}\right\Vert _{\infty;\leq r_{2}}\int\limits_{\left\vert \cdot\right\vert \leq
r_{2}}\frac{1}{w}+\left\Vert \frac{\left\vert \cdot\right\vert ^{2\sigma}%
\phi^{2}}{\left\vert \cdot\right\vert ^{2\theta}}\right\Vert _{\infty
;\left\vert \cdot\right\vert \geq r_{2}}\int\limits_{\left\vert \cdot
\right\vert \geq r_{2}}\frac{1}{w\left\vert \cdot\right\vert ^{2\sigma}%
}\nonumber\\
&  =\left\Vert \frac{\phi}{\left\vert \cdot\right\vert ^{\theta}}\right\Vert
_{\infty;\leq r_{2}}^{2}\int\limits_{\left\vert \cdot\right\vert \leq r_{2}%
}\frac{1}{w}+\left\Vert \frac{\left\vert \cdot\right\vert ^{\sigma}\phi
}{\left\vert \cdot\right\vert ^{\theta}}\right\Vert _{\infty;\geq r_{2}}%
^{2}\int\limits_{\left\vert \cdot\right\vert \geq r_{2}}\frac{1}{w\left\vert
\cdot\right\vert ^{2\sigma}}\nonumber\\
&  \leq\left(  \max\left\{  \left\Vert \frac{\phi}{\left\vert \cdot\right\vert
^{\theta}}\right\Vert _{\infty;\leq r_{2}}^{2},\left\Vert \frac{\left\vert
\cdot\right\vert ^{\sigma}\phi}{\left\vert \cdot\right\vert ^{\theta}%
}\right\Vert _{\infty;\geq r_{2}}^{2}\right\}  \right)  ^{2}\int\frac
{1}{w\left\vert \cdot\right\vert ^{2\lambda\left(  \cdot\right)  }%
}.\label{p52}%
\end{align}

which exists by weight function property W2. Since $n=\operatorname*{ceil}%
\left\{  \theta,\sigma\right\}  $ we can apply the inequality \ref{p33} (see
remark below) to get
\[
\left\Vert \frac{\phi}{\left\vert \cdot\right\vert ^{\theta}}\right\Vert
_{\infty;\leq r_{2}}\leq\sum\limits_{\left\vert \alpha\right\vert =\theta
}\left\Vert D^{\alpha}\phi\right\Vert _{\infty}\leq\sum\limits_{\left\vert
\alpha\right\vert =\theta}\left\Vert \left(  1+\left\vert \cdot\right\vert
\right)  ^{n}D^{\alpha}\phi\right\Vert _{\infty}\leq\sum\limits_{\left\vert
\alpha\right\vert \leq n}\left\Vert \left(  1+\left\vert \cdot\right\vert
\right)  ^{n}D^{\alpha}\phi\right\Vert _{\infty}.
\]

and, since $\phi\in S$
\[
\left\Vert \frac{\left\vert \cdot\right\vert ^{\sigma}\phi}{\left\vert
\cdot\right\vert ^{\theta}}\right\Vert _{\infty;\geq r_{2}}\leq r_{2}%
^{-\theta}\left\Vert \left\vert \cdot\right\vert ^{\sigma}\phi\right\Vert
_{\infty}^{2}\leq r_{2}^{-\theta}\left\Vert \left(  1+\left\vert
\cdot\right\vert \right)  ^{n}\phi\right\Vert _{\infty}\leq r_{2}^{-\theta
}\sum\limits_{\left\vert \alpha\right\vert \leq n}\left\Vert \left(
1+\left\vert \cdot\right\vert \right)  ^{n}D^{\alpha}\phi\right\Vert _{\infty
}.
\]

Substituting these inequalities into the right side of \ref{p52} gives the
estimate of this lemma where%
\begin{equation}
c_{r_{2},\theta}=\max\left\{  1,r_{2}^{-\theta}\right\}  .\label{p74}%
\end{equation}

Finally, by part 1 of Theorem \ref{Thm_prop_functnl_on_Son}, $\left(
\int\frac{\left\vert \phi\right\vert ^{2}}{w\left\vert \cdot\right\vert
^{2\theta}}\right)  ^{1/2}\in S_{\emptyset,\theta}^{\prime}$.
\end{proof}

\begin{theorem}
\label{Thm_property_(g)_F}Suppose $f\in X_{w}^{\theta}$. We define the
function $f_{F}$ by: $f_{F}=\widehat{f}$ on $\mathbb{R}^{d}\setminus0$. Now
suppose the weight function $w$ also has property W2. Then $f_{F}$ has the
following properties:

\begin{enumerate}
\item $f_{F}=0$ iff $f\in P_{\theta-1}$.

\item $f_{F}\in S_{\emptyset,\theta}^{\prime}$ with action $\int f_{F}\phi$,
$\phi\in S_{\emptyset,\theta}$. Also $f_{F}=\widehat{f}$ on $S_{\emptyset
,\theta}$.

\item For any compact set $K$
\[
\int_{K}\left\vert \cdot\right\vert ^{\theta}f_{F}\leq\left(  \int_{K}\frac
{1}{w}\right)  ^{1/2}\left\vert f\right\vert _{w,\theta}<\infty.
\]

\item If $r_{2}$ is the parameter in the definition of weight function
property W2.2%
\[
\int_{\left\vert \cdot\right\vert \geq r_{2}}\frac{\left\vert \cdot\right\vert
^{\theta}\left\vert f_{F}\right\vert }{\left(  1+\left\vert \cdot\right\vert
\right)  ^{\sigma}}\leq\left(  \int_{\left\vert \cdot\right\vert \geq r_{2}%
}\frac{1}{w\left\vert \cdot\right\vert ^{2\sigma}}\right)  ^{1/2}\left\vert
f\right\vert _{w,\theta}<\infty.
\]

\item $\left\vert \cdot\right\vert ^{\theta}f_{F}$ is a regular tempered
distribution in the sense of Appendix \ref{SbSect_property_S'}. Further, if
$\left\vert \alpha\right\vert =\theta$ then $\xi^{\alpha}f_{F}$ is also a
regular tempered distribution.
\end{enumerate}
\end{theorem}

\begin{proof}
\textbf{Part 1} By parts 1 and 3 of Theorem \ref{Thm_properties_Xwm},
$f_{F}=0$ iff $\left\vert f\right\vert _{w,\theta}=0$ iff $f\in P_{\theta-1}
$.\medskip

\textbf{Part 2} First we show that $\int f_{F}\phi$ exists for $\phi\in
S_{\emptyset,\theta}$. By part 1 of Theorem \ref{Thm_properties_Xwm},
$\sqrt{w}\left\vert \cdot\right\vert ^{\theta}f_{F}\in L^{2}$ and $\left\vert
f\right\vert _{w,\theta}^{2}=\int w\left\vert \cdot\right\vert ^{2\theta
}\left\vert f_{F}\right\vert ^{2}$, so for $\phi\in S_{\emptyset,\theta}$
\[
\left\vert \int f_{F}\phi\right\vert \leq\int\left\vert f_{F}\phi\right\vert
=\int\sqrt{w}\left\vert \cdot\right\vert ^{\theta}\left\vert f_{F}\right\vert
\frac{\left\vert \phi\right\vert }{\sqrt{w}\left\vert \cdot\right\vert
^{\theta}}\leq\left\vert f\right\vert _{w,\theta}\left(  \int\frac{\left\vert
\phi\right\vert ^{2}}{w\left\vert \cdot\right\vert ^{2\theta}}\right)  ^{1/2},
\]

having used the Cauchy-Schwartz inequality. By Lemma
\ref{Lem_functnal_phi_sq}
\[
\left(  \int\frac{\left\vert \phi\right\vert ^{2}}{w\left\vert \cdot
\right\vert ^{2\theta}}\right)  ^{1/2}\leq c_{r_{2},\theta}\left(  \int%
\frac{1}{w\left\vert \cdot\right\vert ^{2\lambda\left(  \cdot\right)  }%
}\right)  ^{1/2}\sum\limits_{\left\vert \alpha\right\vert \leq n}\left\Vert
\left(  1+\left\vert \cdot\right\vert \right)  ^{n}D^{\alpha}\phi\right\Vert
_{\infty},
\]

where $n=\operatorname*{ceil}\left\{  \theta,\sigma\right\}  $ and $\lambda$
is the function defined in weight function property W2. By Theorem
\ref{Thm_prop_functnl_on_Son} this inequality implies $f_{F}\in S_{\emptyset
,\theta}^{\prime}$.\smallskip

Now to show that $f_{F}=\widehat{f}$ on $S_{\emptyset,\theta}$. Since
$f_{F}\in S_{\emptyset,\theta}^{\prime}$ Theorem \ref{Thm_prop_functnl_on_Son}
implies $f_{F}$ has a non-unique extension $f_{F}^{e}$ to $S^{\prime}$. Define
$g=\left(  f_{F}^{e}\right)  ^{\vee}\in S^{\prime}$. Next we want to show that
$g\in X_{w}^{\theta}$. In fact if $\left\vert \alpha\right\vert =\theta$ and
$\phi\in S$ then, since $\xi^{\alpha}\phi\in S_{\emptyset,\theta}$, $\left[
\xi^{\alpha}\widehat{g},\phi\right]  =\left[  \xi^{\alpha}f_{F}^{e}%
,\phi\right]  =\left[  f_{F}^{e},\xi^{\alpha}\phi\right]  =\left[  f_{F}%
,\xi^{\alpha}\phi\right]  =\left[  \xi^{\alpha}f_{F},\phi\right]  $ and so
$\xi^{\alpha}\widehat{g}=\xi^{\alpha}f_{F}$. But $\xi^{\alpha}f_{F}\in
L_{loc}^{1}$ by part 4 of Lemma \ref{Lem_properties_Xwm_distrib} which means
that $\xi^{\alpha}\widehat{g}\in L_{loc}^{1}$ and hence by part 2 of Lemma
\ref{Lem_properties_Xwm_distrib} that $\widehat{g}\in L_{loc}^{1}\left(
\mathbb{R}^{d}\setminus0\right)  $. We can now define the function $g_{F}$ by
$g_{F}=\widehat{g}$ on $\mathbb{R}^{d}\setminus0$ and $g_{F}\left(  0\right)
=0$. Now if $\phi\in C_{0}^{\infty}\left(  \mathbb{R}^{d}\setminus0\right)  $
then $\phi\in S_{\emptyset,\theta}$ and it follows that $\left[  g_{F}%
,\phi\right]  =\left[  f_{F}^{e},\phi\right]  =\left[  f_{F},\phi\right]  $
i.e. $g_{F}=f_{F}$ a.e. on $\mathbb{R}^{d}\setminus0$ and hence $g_{F}=f_{F}$
a.e. on $\mathbb{R}^{d}$. Thus $\left\vert g\right\vert _{w,\theta}=\left\vert
f\right\vert _{w,\theta}$ and so $g\in X_{w}^{\theta}$. Further $\left(
g-f\right)  _{F}=0$ and therefore $g-f\in P_{\theta-1}$ by part 1. Finally, by
Theorem \ref{Thm_So,n_and_Pnhat}, $u\in\widehat{P}_{\theta-1}\ $iff$\ \left[
u,\phi\right]  =0$ for all $\phi\in S_{\emptyset,\theta}$. Hence $\widehat
{g}=\widehat{f}$ on $S_{\emptyset,\theta}$ so that for $\phi\in S_{\emptyset
,\theta}$, $\left[  \widehat{f},\phi\right]  =\left[  \widehat{g},\phi\right]
=\left[  f_{F}^{e},\phi\right]  =\left[  f_{F},\phi\right]  $ so that
$\widehat{f}=f_{F}$ on $S_{\emptyset,\theta}$, as required.\medskip

\textbf{Part 3} If $K$ is compact then by using the Cauchy-Schwartz
inequality
\[
\int\limits_{K}\left\vert \cdot\right\vert ^{\theta}\left\vert f_{F}%
\right\vert =\int\limits_{K}\frac{1}{\sqrt{w}}\sqrt{w}\left\vert
\cdot\right\vert ^{\theta}\left\vert f_{F}\right\vert \leq\left(
\int\limits_{K}\frac{1}{w}\right)  ^{1/2}\left(  \int\limits_{K}w\left\vert
\cdot\right\vert ^{2\theta}\left\vert f_{F}\right\vert ^{2}\right)
^{1/2}=\left(  \int\limits_{K}\frac{1}{w}\right)  ^{1/2}\left\vert
f\right\vert _{w,\theta}<\infty,
\]

since weight function property W2 implies that $w\in L_{loc}^{1}$.\medskip

\textbf{Part 4} Using the Cauchy-Schwartz inequality
\begin{align*}
\int_{\left\vert \cdot\right\vert \geq r_{2}}\frac{\left\vert \cdot\right\vert
^{\theta}\left\vert f_{F}\right\vert }{\left(  1+\left\vert \cdot\right\vert
\right)  ^{\sigma}}\leq\int_{\left\vert \cdot\right\vert \geq r_{2}}%
\frac{\left\vert \cdot\right\vert ^{\theta}\left\vert f_{F}\right\vert
}{\left\vert \cdot\right\vert ^{\sigma}} &  =\int_{\left\vert \cdot\right\vert
\geq r_{2}}\sqrt{w}\left\vert \cdot\right\vert ^{\theta}\left\vert
f_{F}\right\vert \frac{1}{\sqrt{w}\left\vert \cdot\right\vert ^{\sigma}}\\
&  \leq\left(  \int_{\left\vert \cdot\right\vert \geq r_{2}}w\left\vert
\cdot\right\vert ^{2\theta}\left\vert f_{F}\right\vert ^{2}\right)
^{1/2}\left(  \int_{\left\vert \cdot\right\vert \geq r_{2}}\frac
{1}{w\left\vert \cdot\right\vert ^{2\sigma}}\right)  ^{1/2}\\
&  \leq\left\vert f\right\vert _{w,\theta}\left(  \int_{\left\vert
\cdot\right\vert \geq r_{2}}\frac{1}{w\left\vert \cdot\right\vert ^{2\sigma}%
}\right)  ^{1/2}\\
&  <\infty,
\end{align*}

since $w$ has property W3.2.\medskip

\textbf{Part 5} That $\left\vert \cdot\right\vert ^{\theta}f_{F}$ is a regular
tempered distribution follows immediately from parts 3 and 4. If $\left\vert
\alpha\right\vert =\theta$ then $\left\vert \xi^{\alpha}f_{F}\right\vert
\leq\left\vert \cdot\right\vert ^{\theta}\left\vert f_{F}\right\vert $ and so
$\xi^{\alpha}f_{F}$ is also a regular tempered distribution.
\end{proof}

\begin{corollary}
\label{Cor_1_Thm_property_(g)_F}Suppose the weight function $w$ has property
W2, and suppose $f\in S^{\prime}$, $\widehat{f}\in L_{loc}^{1}\left(
\mathbb{R}^{d}\setminus0\right)  $ and $\int w\left\vert \cdot\right\vert
^{2\theta}\left\vert f_{F}\right\vert ^{2}<\infty$, where $f_{F}$ was defined
in Theorem \ref{Thm_property_(g)_F}. Then:

\begin{enumerate}
\item $f_{F}=0$ iff $f\in P$.

\item $f_{F}\in S_{\emptyset,\theta}^{\prime}$ with action $\int f_{F}\phi$,
$\phi\in S_{\emptyset,\theta}$.

\item $f_{F}$ has properties 3 and 4 of Theorem \ref{Thm_property_(g)_F}.

\item $f_{F}$ has property 5 of Theorem \ref{Thm_property_(g)_F} i.e.
$\left\vert \cdot\right\vert ^{\theta}f_{F}$ is a regular tempered
distribution and $\xi^{\alpha}f_{F}$ is also a regular tempered distribution
when $\left\vert \alpha\right\vert =\theta$.
\end{enumerate}
\end{corollary}

\begin{proof}
\textbf{Part 1} Since $S^{\prime}\subset\mathcal{D}^{\prime}$, $f_{F}=0$
implies $\widehat{f}\in\mathcal{D}^{\prime}$ and $\operatorname*{supp}%
\widehat{f}=\left\{  0\right\}  $ which implies $f\in P$. Conversely, $f\in P$
implies $f\in S^{\prime}$ and $\operatorname*{supp}\widehat{f}=\left\{
0\right\}  $ i.e. $f_{F}=0$.\medskip

\textbf{Parts 2 to 4} follow from an examination of the proof of Theorem
\ref{Thm_property_(g)_F}.
\end{proof}

The last corollary allows another definition of the space $X_{w}^{\theta}$
when the weight function has property W2.

\begin{corollary}
\label{Cor_2_Thm_property_(g)_F}Suppose the weight function $w$ has property
W2. Then as sets%
\begin{align}
X_{w}^{\theta}  & =\left\{  f\in S^{\prime}:\widehat{f}\in L_{loc}^{1}\left(
\mathbb{R}^{d}\setminus0\right)  ,\text{ }\int w\left\vert \cdot\right\vert
^{2\theta}\left\vert f_{F}\right\vert ^{2}<\infty,\text{ }and\text{
}\left\vert \alpha\right\vert =\theta\text{ }implies\text{\ }\xi^{\alpha
}\widehat{f}=\xi^{\alpha}f_{F}\text{ }on\text{ }S\right\}  ,\nonumber\\
& \label{p67}%
\end{align}

where $f_{F}:\mathbb{R}^{d}\rightarrow\mathbb{C}$ is the function defined by
$f_{F}=\widehat{f}$ on $\mathbb{R}^{d}\setminus0$.

This definition makes sense since by part 4 of Corollary
\ref{Cor_1_Thm_property_(g)_F} the first two constraints imply that $f_{F}\in
S_{\emptyset,\theta}^{\prime}$ and when $\left\vert \alpha\right\vert =\theta
$, $\xi^{\alpha}f_{F}$ is a regular tempered distribution in the sense of
Appendix \ref{SbSect_property_S'} with action $\int\xi^{\alpha}f_{F}\phi$,
$\phi\in S$.
\end{corollary}

\begin{proof}
Definition \ref{p40} of $X_{w}^{\theta}$ is%
\[
\left\{  f\in S^{\prime}:\xi^{\alpha}\widehat{f}\in L_{loc}^{1}\text{
}if\text{ }\left\vert \alpha\right\vert =\theta;\text{\thinspace}\int
w\left\vert \cdot\right\vert ^{2\theta}\left\vert f_{F}\right\vert ^{2}%
<\infty\right\}  .
\]

Suppose $f$ is a member of the the right side of \ref{p67}. Since $w$ has
property W2, part 4 of Corollary \ref{Cor_1_Thm_property_(g)_F} implies that
when $\left\vert \alpha\right\vert =\theta$, $\xi^{\alpha}f_{F}$ is a regular
tempered distribution in the sense of Appendix \ref{SbSect_property_S'} and so
$\xi^{\alpha}f_{F}\in L_{loc}^{1}$. But $\xi^{\alpha}\widehat{f}=\xi^{\alpha
}f_{F}$ as distributions so $\xi^{\alpha}\widehat{f}\in L_{loc}^{1}$ and thus
$f\in X_{w}^{\theta}$.

On the other hand, if $f\in X_{w}^{\theta}$ then by part 2 of Theorem
\ref{Thm_properties_Xwm}, $\widehat{f}\in L_{loc}^{1}\left(  \mathbb{R}%
^{d}\setminus0\right)  $. Further, by part 2 of Theorem
\ref{Thm_property_(g)_F}, $f_{F}\in S_{\emptyset,\theta}^{\prime}$ and
$f_{F}=\widehat{f}$ on $S_{\emptyset,\theta}$. But Theorem
\ref{Thm_product_of_Co,k_funcs} implies $\xi^{\alpha}\psi\in S_{\emptyset
,\theta}$ when $\left\vert \alpha\right\vert =\theta$ so that $\xi^{\alpha
}\widehat{f}=\xi^{\alpha}f_{F}$ on $S$, and we have shown that $f\in
X_{w}^{\theta}$ implies $f$ is a member of the right side of \ref{p67}.
\end{proof}

\begin{remark}
\textbf{1}. ?? \textbf{A word of caution} ?? about definition \ref{p67} of
$X_{w}^{\theta}$. It is clear that $f\in S^{\prime}$ and $\widehat{f}\in
L_{loc}^{1}\left(  \mathbb{R}^{d}\setminus0\right)  $ imply that $f_{F}$
exists. Further, the conditions $f\in S^{\prime}$, $\widehat{f}\in L_{loc}%
^{1}\left(  \mathbb{R}^{d}\setminus0\right)  $ and $\int w\left\vert
\cdot\right\vert ^{2\theta}\left\vert f_{F}\right\vert ^{2}<\infty$
automatically imply that $f_{F}\in S_{\emptyset,\theta}^{\prime}$ and
$\xi^{\alpha}f_{F}\in S^{\prime}$.

However, if we prove that $f\in S^{\prime}$ and $\widehat{f}\in L_{loc}%
^{1}\left(  \mathbb{R}^{d}\setminus0\right)  $ so that $f_{F}$ is defined, and
then try to prove that $\left\vert \alpha\right\vert =\theta$ implies
$\xi^{\alpha}\widehat{f}=\xi^{\alpha}f_{F}$ on $S$, we will first have to show
that $\xi^{\alpha}f_{F}\in S^{\prime}$. There is a question of order
here.\medskip

\textbf{2}. If $w\in W2$ we define%
\begin{equation}
\widetilde{X}^{\theta}:=\left\{  f\in S^{\prime}:\widehat{f}\in L_{loc}%
^{1}\left(  \mathbb{R}^{d}\setminus0\right)  ;\text{ }f_{F}\in S_{\emptyset
,\theta}^{\prime};\text{ }\left\vert \alpha\right\vert =\theta\text{
}implies\text{\ }\xi^{\alpha}\widehat{f}=\xi^{\alpha}f_{F}\text{ }on\text{
}S\right\}  ,\label{p023}%
\end{equation}

then%
\[
X_{w}^{\theta}=\left\{  f\in\widetilde{X}^{\theta}:\int w\left\vert
\cdot\right\vert ^{2\theta}\left\vert f_{F}\right\vert ^{2}<\infty\right\}  .
\]

This follows directly from Corollary \ref{Cor_2_Thm_property_(g)_F}.\medskip

\textbf{3}. The definition of $\widetilde{X}^{\theta}$ in part 2 prompts me to
ask the following question: Suppose $g\in S_{\emptyset,\theta}^{\prime} $ and
$\xi^{\alpha}g=0$ $\forall\left\vert \alpha\right\vert =\theta$. \textbf{What
can we conclude about} $g$? \textbf{NOTE}: See the proof of part 2 of Theorem
\ref{Thm_property_(g)_F}.

By part 2 of Theorem \ref{Thm_prop_functnl_on_Son} there exists a (non-unique)
extension $g^{e}$ of $g$ to $S^{\prime}$ and all the extensions are given by
the Fourier transform formula $g^{e}+\widehat{P}_{\theta-1}$. Thus by Theorem
\ref{Thm_product_of_Co,k_funcs}, $\left[  g^{e},\xi^{\alpha}\phi\right]
=\left[  g,\xi^{\alpha}\phi\right]  =0$ $\forall\phi\in S$ and hence
$\xi^{\alpha}g^{e}=0$ $\forall\left\vert \alpha\right\vert =\theta$ i.e.
$D^{\alpha}\widehat{g^{e}}=0$ $\forall\left\vert \alpha\right\vert =\theta$
and by a well-known distribution result this is true iff $\widehat{g^{e}}\in
P_{\theta-1}$ iff $g^{e}\in\widehat{P}_{\theta-1}$. But by part 1 of Theorem
\ref{Thm_So,n_and_Pnhat}, $\widehat{P}_{\theta-1}\subset S_{\emptyset,\theta
}^{\prime}$ and so \textbf{we can conclude that}%
\begin{equation}
\left\{  g\in S_{\emptyset,\theta}^{\prime}:\xi^{\alpha}g=0\text{ }%
\forall\left\vert \alpha\right\vert =\theta\right\}  =\widehat{P}_{\theta
-1}.\label{p022}%
\end{equation}

We can now conclude that $f\in\widetilde{X}^{\theta}$ implies $f_{F}%
-\widehat{f}\in\widehat{P}_{\theta-1}$ and indeed we have%
\begin{equation}
\widetilde{X}^{\theta}=\left\{  f\in S^{\prime}:\widehat{f}\in L_{loc}%
^{1}\left(  \mathbb{R}^{d}\setminus0\right)  ;\text{ }f_{F}\in S_{\emptyset
,\theta}^{\prime};\text{ }f_{F}-\widehat{f}\in\widehat{P}_{\theta-1}\right\}
.\label{p024}%
\end{equation}

Also from part 1 of Theorem \ref{Thm_So,n_and_Pnhat},%
\[
f_{F}-\widehat{f}\in\widehat{P}_{\theta-1}\text{ }iff\text{ }f_{F}=\widehat
{f}\text{ }on\text{ }S_{\emptyset,\theta},
\]

and consequently%
\begin{equation}
\widetilde{X}^{\theta}=\left\{  f\in S^{\prime}:\widehat{f}\in L_{loc}%
^{1}\left(  \mathbb{R}^{d}\setminus0\right)  ;\text{ }f_{F}\in S_{\emptyset
,\theta}^{\prime};\text{ }f_{F}=\widehat{f}\text{ }on\text{ }S_{\emptyset
,\theta}\right\}  .\label{p025}%
\end{equation}

\end{remark}

The results of the previous remark can be summarized as:

\begin{theorem}
\label{Thm_Xwm_using_squigXm}If $w\in W2$ then%
\[
X_{w}^{\theta}=\left\{  f\in\widetilde{X}^{\theta}:\int w\left\vert
\cdot\right\vert ^{2\theta}\left\vert f_{F}\right\vert ^{2}<\infty\right\}
,\quad f_{F}=\widehat{f}\text{ }on\text{ }\mathbb{R}^{d}\setminus0,
\]

where $\widetilde{X}^{\theta}$ is given by any of \ref{p023}, \ref{p024} and
\ref{p025} each of which is independent of the weight function $w$.
\end{theorem}

??

\begin{remark}
\textbf{1}. Note that if $f\in\widetilde{X}^{\theta}$ then $g=\left(
f_{F}\right)  ^{\vee}$ implies $\widehat{g}\in S_{\emptyset,\theta}^{\prime}$
and $\left(  \widehat{g}\right)  ^{e}\in L_{loc}^{1}\left(  \mathbb{R}%
^{d}\setminus0\right)  \cap S^{\prime}$. We thus define%
\[
\widetilde{Y}^{\theta}:=\left\{  h\in S_{\emptyset,\theta}^{\prime}:r_{0}%
h^{e}\in L_{loc}^{1}\left(  \mathbb{R}^{d}\setminus0\right)  ,\text{ }%
h=r_{0}h^{e}\text{ }on\text{ }S_{\emptyset,\theta}\right\}  ,
\]

where
\[
r_{0}g:=restriction\text{ }of\text{ }g\in S^{\prime}\text{ }to\text{
}\mathbb{R}^{d}\setminus0\text{ }as\text{ }a\text{ }distribution.
\]

If $r_{\emptyset}:S^{\prime}\rightarrow S_{\emptyset,\theta}^{\prime}$ is the
restriction operator then is $r_{\emptyset}\left(  \left(  \widetilde
{X}^{\theta}\right)  ^{\wedge}\right)  \subset\widetilde{Y}^{\theta}$? Is this
mapping onto, 1-1? How to endow $\widetilde{X}^{\theta}$ and $\widetilde
{Y}^{\theta}$ with topologies? Show that%
\begin{equation}
\left(  \widetilde{X}^{\theta}\right)  ^{\wedge}=\left\{  g\in S^{\prime
}:r_{0}g\in L_{loc}^{1}\left(  \mathbb{R}^{d}\setminus0\right)  ;\text{ }%
r_{0}g=g\text{ }on\text{ }S_{\emptyset,\theta}\right\}  .\label{p020}%
\end{equation}

Note that $f_{F}=r_{0}\widehat{f}$. Now suppose $g\in\left(  \widetilde
{X}^{\theta}\right)  ^{\wedge}$. Clearly $r_{\emptyset}g\in S_{\emptyset
,\theta}^{\prime}$. Also $r_{\emptyset}g=r_{0}g\in L_{loc}^{1}\left(
\mathbb{R}^{d}\setminus0\right)  $ so $r_{\emptyset}\left(  \left(
\widetilde{X}^{\theta}\right)  ^{\wedge}\right)  \subset\widetilde{Y}^{\theta
}$. Now prove equality. If $h\in\widetilde{Y}^{\theta}$ then $h\in
S_{\emptyset,\theta}^{\prime}$ and $r_{0}h^{e}\in L_{loc}^{1}\left(
\mathbb{R}^{d}\setminus0\right)  $ and $h=r_{0}h^{e}$ $on$ $S_{\emptyset
,\theta}$. By inspection $h^{e}\in\left(  \widetilde{X}^{\theta}\right)
^{\wedge}$. Thus%
\[
r_{\emptyset}\left(  \left(  \widetilde{X}^{\theta}\right)  ^{\wedge}\right)
=\widetilde{Y}^{\theta}.
\]

Is $r_{\emptyset}$ 1-1? From Theorem \ref{Thm_So,n_and_Pnhat}, $r_{\emptyset
}g=0$ implies $g\in\widehat{P}_{\theta-1}$. Thus%
\[
\ker r_{\emptyset}=\widehat{P}_{\theta-1}.
\]
\medskip

\textbf{2}. ?? RELATE $\widetilde{X}^{\theta}$ and $X_{w}^{\theta}$ to
$\left(  S_{\emptyset,\theta}^{\prime}\right)  ^{\vee}=\left(  \widehat
{S}_{\emptyset,\theta}\right)  ^{\prime}$. We have a mapping $\Phi
:f\rightarrow\left(  f_{F}\right)  ^{\vee}$ and $f=\left(  f_{F}\right)
^{\vee}$ on ?? $\widehat{S}_{\emptyset,\theta}$.??
\end{remark}

\subsection{The operators $\mathcal{I}:X_{w}^{\theta}\rightarrow L^{2}$ and
$\mathcal{J}:L^{2}\rightarrow X_{w}^{\theta}$\label{SbSect_I_J}}

In this section, by analogy with the study of Sobolev spaces, we will define
inverse, isometric operators (in the seminorm sense) between $X_{w}^{\theta} $
and $L^{2}$ for $\theta\geq1$. In the next section these operators will be
used to prove the completeness of $X_{w}^{\theta}$, again in the seminorm
sense, without referring to any other space such as $X_{w}^{0}$. If the weight
function also has property W2 these mapping will turn out to be, in the
seminorm sense, adjoints, inverses and isometric isomorphisms. It is the
properties of these operators that count and the reader can avoid their
definitions and the associated lemmas and just study the properties given in
Theorems \ref{Thm_properties_op_Itheta}, \ref{Thm_properties_op_Jtheta} and
\ref{Thm_Itheta_Jtheta}.

\begin{definition}
\label{Def_Itheta}\textbf{The operator} $\mathcal{I}:X_{w}^{\theta}\rightarrow
L^{2}$, $\theta=1,2,3,\ldots$.

Suppose that the function $w$ has weight function property W1. Using the
definition of $X_{w}^{\theta}$ given in part 2 of Theorem
\ref{Thm_properties_Xwm} we define the linear operator $\mathcal{I}%
:X_{w}^{\theta}\rightarrow L^{2}$ by
\begin{equation}
\mathcal{I}f=\left(  \sqrt{w}\left\vert \cdot\right\vert ^{\theta}%
f_{F}\right)  ^{\vee}.\label{p26}%
\end{equation}

where the function $f_{F}$ is defined by: $f_{F}=\widehat{f}$ on
$\mathbb{R}^{d}\setminus0$ and $f_{F}\left(  0\right)  =\left\{  0\right\}  $.
\end{definition}

\begin{theorem}
\label{Thm_properties_op_Itheta}Properties of\textbf{\ }$\mathcal{I}%
$\textbf{:}

\begin{enumerate}
\item $\mathcal{I}$ is an isometric mapping from $X_{w}^{\theta}$ to $L^{2}$
in the seminorm sense.

\item $\operatorname*{null}\mathcal{I}=P_{\theta-1}$.
\end{enumerate}
\end{theorem}

\begin{proof}
\textbf{Property 1} From part 1 of Theorem \ref{Thm_properties_Xwm} we have

$\left\Vert \mathcal{I}f\right\Vert _{2}=\left\Vert \sqrt{w}\left\vert
\cdot\right\vert ^{\theta}f_{F}\right\Vert _{2}=\left\vert f\right\vert
_{w,\theta}$.\smallskip

\textbf{Property 2} Now $\mathcal{I}f=0$ iff$\sqrt{w}\left\vert \cdot
\right\vert ^{\theta}f_{F}=0$ iff $f_{F}=0$. But $f_{F}=0$ iff $f\in
P_{\theta-1}$ by part 1 of Theorem \ref{Thm_property_(g)_F}.
\end{proof}

The next step is to construct an inverse of $\mathcal{I}$, which we will
denote by $\mathcal{J}$ and this will require the theory of the spaces
$S_{\emptyset,\theta}$ and $S_{\emptyset,\theta}^{\prime}$ discussed above in
Section \ref{Sect_So,n_S'o,n}. We will be looking for an operator which makes
rigorous the formal operator $\left(  \frac{\widehat{g}}{\sqrt{w}\left\vert
\cdot\right\vert ^{\theta}}\right)  ^{\vee}$, $g\in L^{2}$. To justify the
inverse-Fourier transform we show that $\frac{\widehat{g}}{\sqrt{w}\left\vert
\cdot\right\vert ^{\theta}}\in S_{\emptyset,\theta}^{\prime}$ and then extend
$\frac{\widehat{g}}{\sqrt{w}\left\vert \cdot\right\vert ^{\theta}}$ to $S$ as
a member of $S^{\prime}$. This operator from $L^{2}$ to $S^{\prime}$ will
become our inverse operator $\mathcal{J}:L^{2}\rightarrow X_{w}^{\theta}$. We
start by showing $\left(  \frac{\widehat{g}}{\sqrt{w}\left\vert \cdot
\right\vert ^{\theta}}\right)  ^{\vee}\in S_{\emptyset,\theta}^{\prime}$ when
$g\in L^{2}$. Choose $\phi\in S_{\emptyset,\theta}$ and apply the
Cauchy-Schwartz inequality to obtain%
\[
\left\vert \int\frac{\widehat{g}}{\sqrt{w}\left\vert \cdot\right\vert
^{\theta}}\phi\right\vert \leq\left(  \int\left\vert g\right\vert ^{2}\right)
^{1/2}\left(  \int\frac{\left\vert \phi\right\vert ^{2}}{w\left\vert
\cdot\right\vert ^{2\theta}}\right)  ^{1/2}=\left\Vert g\right\Vert
_{2}\left(  \int\frac{\left\vert \phi\right\vert ^{2}}{w\left\vert
\cdot\right\vert ^{2\theta}}\right)  ^{1/2}.
\]

Now observe that by Lemma \ref{Lem_functnal_phi_sq}, $\left(  \int%
\frac{\left\vert \phi\right\vert ^{2}}{w\left\vert \cdot\right\vert ^{2\theta
}}\right)  ^{1/2}\in S_{\emptyset,\theta}^{\prime}$ and so $\frac{\widehat{g}%
}{\sqrt{w}\left\vert \cdot\right\vert ^{\theta}}\in S_{\emptyset,\theta
}^{\prime}$. Theorem \ref{Thm_prop_functnl_on_Son} now allows us to extend the
functional $\frac{\widehat{g}}{\sqrt{w}\left\vert \cdot\right\vert ^{\theta}}$
to $S$ as a member of $S^{\prime}$.\ We can now define the operator
$\mathcal{J}$.

\begin{definition}
\label{Def_operator_Ttheta}\textbf{The operator} $\mathcal{J}:L^{2}\rightarrow
S^{\prime}$, $\theta\geq1$. Suppose that the weight function $w$ has property
W2 and $g\in L^{2}$. Then $\frac{\widehat{g}}{\sqrt{w}\left\vert
\cdot\right\vert ^{\theta}}\in S_{\emptyset,\theta}^{\prime}$. Hence by
Theorem \ref{Thm_prop_functnl_on_Son}, this functional can be extended (up to
a member of $\widehat{P_{\theta-1}}$) to $S$ as an element of $S^{\prime}$,
say $f$. Now define the class of mappings $\mathcal{J}:L^{2}\rightarrow
S^{\prime}$ by
\[
\mathcal{J}g=\overset{\vee}{f}.
\]

Note that $\mathcal{J}$ is not linear but the next theorem shows it is linear
modulo a polynomial in $P_{\theta-1}$.
\end{definition}

\begin{theorem}
\label{Thm_properties_op_Jtheta}\textbf{Properties of }$\mathcal{J}$:

\begin{enumerate}
\item $\left(  \mathcal{J}g\right)  _{F}=\frac{\widehat{g}}{\sqrt{w}\left\vert
\cdot\right\vert ^{\theta}}$,\quad$g\in L^{2}$.

\item $\mathcal{J}:L^{2}\rightarrow X_{w}^{\theta}$ and is an isometry in the
seminorm sense.

\item $\mathcal{J}$ is linear modulo a polynomial in $P_{\theta-1}$ i.e.%
\[
\mathcal{J}\left(  \lambda_{1}g_{1}+\lambda_{2}g_{2}\right)  -\lambda
_{1}\mathcal{J}g_{1}-\lambda_{2}\mathcal{J}g_{2}\in P_{\theta-1}.
\]

\item $\mathcal{J}g\in P_{\theta-1}$ iff $g=0$.
\end{enumerate}
\end{theorem}

\begin{proof}
\textbf{Properties 1 \& 2} To prove $\mathcal{J}g\in X_{w}^{\theta}$ we use
the definition of $X_{w}^{\theta}$ given in Corollary
\ref{Cor_2_Thm_property_(g)_F} i.e. $\mathcal{J}g\in S^{\prime}$,
$\widehat{\mathcal{J}g}\in L_{loc}^{1}\left(  \mathbb{R}^{d}\setminus0\right)
$, and for $\left(  \mathcal{J}g\right)  _{F}$ defined by: $\left(
\mathcal{J}g\right)  _{F}=\widehat{\mathcal{J}g}$ on $\mathbb{R}^{d}%
\setminus0$ and $\left(  \mathcal{J}g\right)  _{F}\left(  0\right)  =\left\{
0\right\}  $, it is required that $\int w\left\vert \cdot\right\vert
^{2\theta}\left\vert \left(  \mathcal{J}g\right)  _{F}\right\vert ^{2}<\infty$
and $\xi^{\alpha}\widehat{\mathcal{J}g}=\xi^{\alpha}\left(  \mathcal{J}%
g\right)  _{F}$ on $S$ when $\left\vert \alpha\right\vert =\theta$.\medskip

From the definition of $\mathcal{J}g$ we have $\widehat{\mathcal{J}g}%
=\frac{\widehat{g}}{\sqrt{w}\left\vert \cdot\right\vert ^{\theta}}$ on
$\mathbb{R}^{d}\setminus0$ and if $K\subset\mathbb{R}^{d}\setminus0$ is
compact
\[
\int_{K}\frac{\left\vert \widehat{g}\right\vert }{\sqrt{w}\left\vert
\cdot\right\vert ^{\theta}}\leq\left(  \int_{K}\left\vert \widehat
{g}\right\vert ^{2}\right)  ^{1/2}\left(  \int_{K}\frac{1}{w\left\vert
\cdot\right\vert ^{2\theta}}\right)  ^{1/2}\leq\left\Vert g\right\Vert
_{2}\left(  \int_{K}\frac{1}{w\left\vert \cdot\right\vert ^{2\theta}}\right)
^{1/2}.
\]

Since $0\in0$ we have $\operatorname*{dist}\left(  0;K\right)  >0$ and so%
\[
\int_{K}\frac{1}{w\left\vert \cdot\right\vert ^{2\theta}}\leq\frac{1}{\left(
\operatorname*{dist}\left(  0;K\right)  \right)  ^{2\theta}}\int_{K}\frac
{1}{w}<\infty,
\]

because $1/w\in L_{loc}^{1}$ by weight function property W2.1. We now have
$\left(  \mathcal{J}g\right)  _{F}=\frac{\widehat{g}}{\sqrt{w}\left\vert
\cdot\right\vert ^{\theta}}$ which proves property 1. Next%
\[
\left\vert \mathcal{J}g\right\vert _{w,\theta}^{2}=\int w\left\vert
\cdot\right\vert ^{2\theta}\left\vert \left(  \mathcal{J}g\right)
_{F}\right\vert ^{2}=\int w\left\vert \cdot\right\vert ^{2\theta}\left\vert
\frac{\widehat{g}}{\sqrt{w}\left\vert \cdot\right\vert ^{\theta}}\right\vert
^{2}=\left\Vert g\right\Vert _{2}^{2}<\infty,
\]

so that $\mathcal{J}$ is isometric. Now if $\left\vert \alpha\right\vert
=\theta$ and $\psi\in S$ then $\xi^{\alpha}\psi\in S_{\emptyset,\theta}$ by
Theorem \ref{Thm_product_of_Co,k_funcs}, and since the operator $\mathcal{J}$
was defined by extending $\frac{\widehat{g}}{\sqrt{w}\left\vert \cdot
\right\vert ^{\theta}}\in S_{\emptyset,\theta}^{\prime}$ from $S_{\emptyset
,\theta}$ to $S$ as a member of $S^{\prime}$ it follows that
\[
\left[  \xi^{\alpha}\widehat{\mathcal{J}g},\psi\right]  =\left[
\widehat{\mathcal{J}g},\xi^{\alpha}\psi\right]  =\left[  \frac{\widehat{g}%
}{\sqrt{w}\left\vert \cdot\right\vert ^{\theta}},\xi^{\alpha}\psi\right]
=\left[  \xi^{\alpha}\left(  \mathcal{J}g\right)  _{F},\psi\right]  ,
\]

so that $\xi^{\alpha}\widehat{\mathcal{J}g}=\xi^{\alpha}\left(  \mathcal{J}%
g\right)  _{F}$ on $S$, confirming that $\mathcal{J}g\in X_{w}^{\theta}%
$.\medskip

\textbf{Property 3} If $\phi\in S_{\emptyset,\theta}$ then
\begin{align*}
\left[  \left(  \mathcal{J}\left(  \lambda_{1}g_{1}+\lambda_{2}g_{2}\right)
-\lambda_{1}\mathcal{J}g_{1}-\lambda_{2}\mathcal{J}g_{2}\right)  ^{\wedge
},\phi\right]   & =\left[  \frac{\widehat{\lambda_{1}g_{1}+\lambda_{2}g_{2}}%
}{\sqrt{w}\left\vert \cdot\right\vert ^{\theta}}-\frac{\widehat{\lambda
_{1}g_{1}}}{\sqrt{w}\left\vert \cdot\right\vert ^{\theta}}-\frac
{\widehat{\lambda_{2}g_{2}}}{\sqrt{w}\left\vert \cdot\right\vert ^{\theta}%
},\phi\right] \\
& =0.
\end{align*}

Part 2 of Theorem \ref{Thm_So,n_and_Pnhat} gives the required result.\medskip

\textbf{Property 4} The details of this proof are very similar to those of the
proof of part 2 of Theorem \ref{Thm_properties_op_Itheta}. Suppose
$\mathcal{J}g\in P_{\theta-1}$. Then $\left(  \mathcal{J}g\right)  _{F}=0$,
$\left\vert \mathcal{J}g\right\vert _{w,\theta}=0$ and thus $\left\Vert
g\right\Vert _{2}=0$ which implies $g=0$. The argument is easily reversible.
\end{proof}

Having proved some properties of $\mathcal{I}$ and $\mathcal{J}$ we now study
how they interact:

\begin{theorem}
\label{Thm_Itheta_Jtheta}Suppose the weight function $w$ has property W2. Then
for $\theta\geq1$ the operators $\mathcal{I}:X_{w}^{\theta}\rightarrow L^{2}$
and $\mathcal{J}:L^{2}\rightarrow X_{w}^{\theta}$ interact as follows:

\begin{enumerate}
\item $\left(  \mathcal{JI}f\right)  _{F}=f_{F}$ when $f\in X_{w}^{\theta}$.
Also, $\mathcal{JI}:X_{w}^{\theta}\rightarrow X_{w}^{\theta}$ is an isometry
in the seminorm sense.

\item $\mathcal{JI}f-f\in P_{\theta-1}$.

\item For all choices of $\mathcal{J}$, $\mathcal{IJ}=I$ on $L^{2}$.

\item $\mathcal{I}$ and $\mathcal{J}$ are adjoints in the sense that
$\left\langle \mathcal{J}g,f\right\rangle _{w,\theta}=\left(  g,\mathcal{I}%
f\right)  _{2}$.
\end{enumerate}
\end{theorem}

\begin{proof}
\textbf{Part 1} From part 1 of Theorem \ref{Thm_properties_op_Itheta}
$\widehat{\mathcal{I}f}=\sqrt{w}\left\vert \cdot\right\vert ^{\theta}f_{F}$
and from Theorem \ref{Thm_properties_op_Jtheta} $\left(  \mathcal{J}g\right)
_{F}=\frac{\widehat{g}}{\sqrt{w}\left\vert \cdot\right\vert ^{\theta}}$. Hence
$\left(  \mathcal{JI}f\right)  _{F}=\frac{\widehat{\mathcal{I}f}}{\sqrt
{w}\left\vert \cdot\right\vert ^{\theta}}=f_{F}$. Further, $\left\vert
\mathcal{JI}f\right\vert _{w,\theta}^{2}=\int w\left\vert \cdot\right\vert
^{2\theta}\left\vert \left(  \mathcal{JI}f\right)  _{F}\right\vert ^{2}=\int
w\left\vert \cdot\right\vert ^{2\theta}\left\vert f_{F}\right\vert
^{2}=\left\vert f\right\vert _{w,\theta}^{2}.$\medskip

\textbf{Part 2} From part 1, $\mathcal{JI}f-f\in X_{w}^{\theta}$ and $\left(
\mathcal{JI}f-f\right)  _{F}=0$. But by part 1 of Theorem
\ref{Thm_property_(g)_F} this implies that $\mathcal{JI}f-f\in P_{\theta-1}%
$.\medskip

\textbf{Part 3} Suppose $g\in L^{2}$. Then $\mathcal{J}g\in X_{w}^{\theta}$
and $\left(  \mathcal{J}g\right)  _{F}=\frac{\widehat{g}}{\sqrt{w}\left\vert
\cdot\right\vert ^{\theta}}$. Thus

$\mathcal{IJ}f=\left(  \sqrt{w}\left\vert \cdot\right\vert ^{\theta}\left(
\mathcal{J}g\right)  _{F}\right)  ^{\vee}=g$.\medskip

\textbf{Part 4}
\[
\left\langle \mathcal{J}g,f\right\rangle _{w,\theta}=\int w\left\vert
\cdot\right\vert ^{2\theta}\left(  \mathcal{J}g\right)  _{F}\overline{f_{F}%
}=\int\sqrt{w}\left\vert \cdot\right\vert ^{\theta}\widehat{g}\,\overline
{f_{F}}=\int\widehat{g}\overline{\sqrt{w}\left\vert \cdot\right\vert ^{\theta
}f_{F}}=\left(  g,\mathcal{I}f\right)  _{2}.
\]

\end{proof}

\subsection{The completeness of $X_{w}^{\theta}$\label{SbSect_Xwth_complete}}

In this section, by analogy with Sobolev space theory, we use the operators
$\mathcal{I}:X_{w}^{\theta}\rightarrow L^{2}$ and $\mathcal{J}:L^{2}%
\rightarrow X_{w}^{\theta}$ where $\theta\geq1$, studied in the previous
subsection to prove that when $w$ has weight function properties W1 and W2 the
semi-inner product space $X_{w}^{\theta}$ is a semi-Hilbert space i.e. it is
complete in the seminorm sense.

Light and Wayne \cite{LightWayneX98Weight} do not define the operator
$\mathcal{J}:L^{2}\rightarrow X_{w}^{\theta}$. They only define the operators
$\mathcal{I}$ and $\mathcal{J}$ between $X_{w}^{0}$ and $L^{2}$ and their
proof of the completeness of $X_{w}^{\theta}$ uses multiple operators and
spaces which correspond to each multi-index $\alpha$ such that $\left\vert
\alpha\right\vert =\theta$. In fact, in Definition 2.11 they define the
semi-inner product spaces $\left(  Y_{\alpha},\left\vert f\right\vert
_{\alpha}\right)  _{\left\vert \alpha\right\vert =\theta}$ given by%
\[
Y_{\alpha}=\left\{  f\in S^{\prime}:\widehat{D^{\alpha}f}\in L_{loc}^{1}\text{
}and\text{ }\int w\left\vert \widehat{D^{\alpha}f}\right\vert ^{2}%
<\infty\right\}  ,\quad\left\vert f\right\vert _{\alpha}=\sqrt{\int
w\left\vert \widehat{D^{\alpha}f}\right\vert ^{2}},
\]

and show that $D^{\alpha}:Y_{\alpha}\rightarrow X_{w}^{0}$ is onto and an
isometric isomorphism, which implies each $Y_{\alpha}$ is complete. Then in
Definition 2.14 they define their positive order semi-inner product space
$\left(  X,\left\vert f\right\vert \right)  $ by
\[
X=\bigcap\limits_{\left\vert \alpha\right\vert =\theta}Y_{\alpha}%
,\quad\left\vert f\right\vert ^{2}=\sum_{\left\vert \alpha\right\vert =\theta
}\frac{1}{\alpha!}\int w\left\vert \widehat{D^{\alpha}f}\right\vert ^{2},
\]
and argue in Theorem 2.15 that the completeness of each $Y_{\alpha}$ implies
that $Y$ is complete.

So there are two ways to prove completeness and although the definition of
$\mathcal{J}$ may be difficult, the operators $\mathcal{I}$ and $\mathcal{J}$
have nice properties and it is really only these properties that are important.

\begin{definition}
\label{Def_seminorm_complete}\textbf{Completeness in the seminorm sense}

Suppose that $\mathbb{U}$ is a linear space equipped with a seminorm
$\left\vert \cdot\right\vert $. We will refer to $\mathbb{U}$ as being
complete in the seminorm sense, if to each Cauchy sequence $\left\{
u_{j}\right\}  \subset\mathbb{U}$ there corresponds an element $u\in
\mathbb{U}$ such that $\left\vert u-u_{j}\right\vert \rightarrow0$ as
$j\rightarrow\infty$. Note that $u$ is no longer uniquely defined by the
sequence $\left\{  u_{j}\right\}  $.
\end{definition}

\begin{theorem}
\label{Thm_Xw_complete}Suppose the weight function $w$ only has property W2.
Then for $\theta\geq1$, $X_{w}^{\theta}$ is complete in the seminorm sense of
Definition \ref{Def_seminorm_complete}.
\end{theorem}

\begin{proof}
First note that the operator $\mathcal{J}$ is not linear. Here we use the
results of Section \ref{SbSect_I_J} concerning the operators $\mathcal{I}$ and
$\mathcal{J}$. Suppose $\left\{  f_{k}\right\}  $ is Cauchy in $X_{w}^{\theta
}$. Then $\left\{  \mathcal{I}f_{k}\right\}  $ is Cauchy in $L^{2}$ since
$\mathcal{I}$ is an isometry. Since $L^{2}$ is complete, $\mathcal{I}%
f_{k}\rightarrow g$ for some $g\in L^{2}$. But $\mathcal{J}\left(
\mathcal{I}f_{k}-g\right)  \in X_{w}^{\theta}$ and since $\mathcal{J}$ is an
isometry, $\left\vert \mathcal{J}\left(  \mathcal{I}f_{k}-g\right)
\right\vert _{w,\theta}=\left\vert \mathcal{I}f_{k}-g\right\vert
_{2}\rightarrow0$.

However, by part 2 of Theorem \ref{Thm_Itheta_Jtheta}, $\mathcal{JI}f-f\in
P_{\theta-1}$. so that
\[
\left\vert \mathcal{J}\left(  \mathcal{I}f_{k}-g\right)  \right\vert
_{w,\theta}^{2}=\left\vert \mathcal{JI}f_{k}-\mathcal{J}g\right\vert
_{w,\theta}^{2}=\left\vert f_{k}-\mathcal{J}g\right\vert _{w,\theta}^{2},
\]

and so $f_{k}\rightarrow\mathcal{J}g$ as $k\rightarrow\infty$.
\end{proof}

The next corollary relates weight function property W3.1 to the completeness
of $X_{w}^{\theta}$. Here weight function property W3.1 replaces property W2.

\begin{corollary}
\label{Cor_thm_Xw_complete_W3.1}Suppose the weight function $w$ only has
property W3.1 for order $\theta$ and parameter $\kappa$. Then $X_{w}^{\theta}$
is complete in the seminorm sense of Definition \ref{Def_seminorm_complete}.
\end{corollary}

\begin{proof}
By part 2 of Theorem \ref{Thm_weight_property_relat}, property W3.1 implies
property W2. Thus the conditions of Theorem \ref{Thm_Xw_complete} are
satisfied and $X_{w}^{\theta}$ is complete.
\end{proof}

\subsection{The operators $\mathcal{T}:X_{w}^{\theta}\rightarrow X_{w}^{0}$
and $\mathcal{U}:X_{w}^{0}\rightarrow X_{w}^{\theta}$}

Following my thesis I will construct the operators $\mathcal{T}:X_{w}^{\theta
}\rightarrow X_{w}^{0}$ and $\mathcal{U}:X_{w}^{0}\rightarrow X_{w}^{\theta}$.
These mapping will turn out to be, in the seminorm sense, adjoints, inverses
and isometric isomorphisms. These mappings are indicated on Figure
\ref{Fig_ops_I_J_U_T_Trho_U_rho}.

These operators are not actually used in the sequel but I found them
interesting. We start with several lemmas.

\begin{lemma}
\label{Lem_prep_for_op_T} Suppose that the weight function $w$ has properties
W1 and W2. If $f\in X_{w}^{\theta}$, $\widehat{f}\in L_{loc}^{1}\left(
\mathbb{R}^{d}\setminus0\right)  $ and we can define the function $f_{F}$ on
$\mathbb{R}^{d}$ by $f_{F}=\widehat{f}$ on $\mathbb{R}^{d}\setminus0$ and
$f_{F}\left(  0\right)  =0$. Then,

\begin{enumerate}
\item $\left\vert \cdot\right\vert ^{\theta}f_{F}\in S^{\prime}\cap
\mathbb{L}_{loc}^{1}$ i.e. a regular tempered distribution.

\item The mapping $f\rightarrow\left\vert \cdot\right\vert ^{\theta}f_{F}$ is
a continuous operator from $X_{w}^{\theta}$ to $S^{\prime}$. In fact%
\[
\int\left\vert \cdot\right\vert ^{\theta}\left\vert f_{F}\phi\right\vert \leq
c_{r_{2}}\left\vert f\right\vert _{w,\theta}\left(  \int\frac{1}{w\left\vert
\cdot\right\vert ^{2\lambda\left(  \cdot\right)  }}\right)  ^{1/2}\left\Vert
\left(  1+\left\vert \cdot\right\vert \right)  ^{\left\lceil \sigma
\right\rceil }\phi\right\Vert _{\infty},
\]

where $c_{r_{1}}=2\sqrt{2}\max\left\{  1,r_{1}^{2}\right\}  $ and the function
$\lambda$ is given in the definition of property W2.
\end{enumerate}
\end{lemma}

\begin{proof}
\textbf{Part 1}. By part 4 of Lemma \ref{Lem_properties_Xwm_distrib},
$\left\vert \cdot\right\vert ^{\theta}f_{F}\in L_{loc}^{1}$ so that
$\left\vert \cdot\right\vert ^{\theta}f_{F}$ is a regular distribution. But by
part 4 of Theorem \ref{Thm_property_(g)_F}, $%
{\displaystyle\int_{\left\vert \cdot\right\vert \geq r_{2}}}
\dfrac{\left\vert \cdot\right\vert ^{\theta}\left\vert f_{F}\right\vert
}{\left(  1+\left\vert \cdot\right\vert \right)  ^{\sigma}}<\infty$ where
$\sigma$ and $r_{2}$ are in the definition of weight function property W2.2.
Thus $\left\vert \cdot\right\vert ^{\theta}f_{F}$ is a regular tempered
distribution. See, for example, Section 2.3 of Vladimirov \cite{Vladimirov}%
.\medskip

\textbf{Part 2}. To prove the continuity of the mapping $f\rightarrow
\left\vert \cdot\right\vert ^{\theta}f_{F}$ we will show that for all $\phi\in
S$, $\left\vert \left[  \left\vert \cdot\right\vert ^{\theta}f_{F}%
,\phi\right]  \right\vert \leq\left\vert f\right\vert _{w,\theta}\left\vert
\phi\right\vert _{S}$ where $\left\vert \cdot\right\vert _{S}$ is a finite
linear combination of the countable set of seminorms, independent of $\phi$,
which define the topology of $S$. Here we use the directed set of seminorms
$\left\{  p_{n,\alpha}:\left\vert \alpha\right\vert \leq n,\text{
}n=0,1,2,\ldots\right\}  $ where $p_{n,\alpha}\left(  \psi\right)  =\left\Vert
\left(  1+\left\vert \cdot\right\vert \right)  ^{n}D^{\alpha}\psi\right\Vert
_{\infty}$. This set was introduced in part 1 of Theorem
\ref{Thm_prop_functnl_on_Son}. In fact%
\begin{align*}
\int\left\vert \cdot\right\vert ^{\theta}\left\vert f_{F}\phi\right\vert  &
\leq\int_{\left\vert \cdot\right\vert \leq r_{2}}\left(  \left\vert
\cdot\right\vert ^{\theta}\left\vert f_{F}\right\vert \right)  \left\vert
\phi\right\vert +\int_{\left\vert \cdot\right\vert \geq r_{2}}\frac{\left\vert
\cdot\right\vert ^{\theta}\left\vert f_{F}\right\vert }{\left\vert
\cdot\right\vert ^{\sigma}}\left(  1+\left\vert \cdot\right\vert \right)
^{\sigma}\left\vert \phi\right\vert \\
& \leq\left\Vert \phi\right\Vert _{\infty}\int_{\left\vert \cdot\right\vert
\leq r_{2}}\left\vert \cdot\right\vert ^{\theta}\left\vert f_{F}\right\vert
+\left\Vert \left(  1+\left\vert \cdot\right\vert \right)  ^{\left\lceil
\sigma\right\rceil }\phi\right\Vert _{\infty}\int_{\left\vert \cdot\right\vert
\geq r_{2}}\frac{\left\vert \cdot\right\vert ^{\theta}\left\vert
f_{F}\right\vert }{\left\vert \cdot\right\vert ^{\sigma}}\\
& \leq\left(  \int_{\left\vert \cdot\right\vert \leq r_{2}}\left\vert
\cdot\right\vert ^{\theta}\left\vert f_{F}\right\vert +\int_{\left\vert
\cdot\right\vert \geq r_{2}}\frac{\left\vert \cdot\right\vert ^{\theta
}\left\vert f_{F}\right\vert }{\left\vert \cdot\right\vert ^{\sigma}}\right)
\left\Vert \left(  1+\left\vert \cdot\right\vert \right)  ^{\left\lceil
\sigma\right\rceil }\phi\right\Vert _{\infty},
\end{align*}

and using the inequalities of parts 3 and 4 of Theorem
\ref{Thm_property_(g)_F}.
\begin{align*}
\int\left\vert \cdot\right\vert ^{\theta}\left\vert f_{F}\phi\right\vert  &
\leq\left(  \left(  \int\limits_{\left\vert \cdot\right\vert \leq r_{2}}%
\frac{1}{w}\right)  ^{1/2}\left\vert f\right\vert _{w,\theta}+\left(
\int\limits_{\left\vert \cdot\right\vert \geq r_{2}}\frac{1}{w\left\vert
\cdot\right\vert ^{2\sigma}}\right)  ^{1/2}\left\vert f\right\vert _{w,\theta
}\right)  \left\Vert \left(  1+\left\vert \cdot\right\vert \right)
^{\left\lceil \sigma\right\rceil }\phi\right\Vert _{\infty}\\
& \leq\sqrt{2}\left\vert f\right\vert _{w,\theta}\left(  \int_{\left\vert
\cdot\right\vert \leq r_{2}}\frac{1}{w}+\int_{\left\vert \cdot\right\vert \geq
r_{2}}\frac{1}{w\left\vert \cdot\right\vert ^{2\sigma}}\right)  ^{1/2}%
\left\Vert \left(  1+\left\vert \cdot\right\vert \right)  ^{\left\lceil
\sigma\right\rceil }\phi\right\Vert _{\infty}\\
& \leq\sqrt{2}\left\vert f\right\vert _{w,\theta}\left(  \int%
\limits_{\left\vert \cdot\right\vert \leq r_{2}}\frac{1}{w}+\int%
\limits_{\left\vert \cdot\right\vert \geq r_{2}}\frac{1}{w\left\vert
\cdot\right\vert ^{2\sigma}}\right)  ^{1/2}\left\Vert \left(  1+\left\vert
\cdot\right\vert \right)  ^{\left\lceil \sigma\right\rceil }\phi\right\Vert
_{\infty}\\
& \leq\sqrt{2}\left\vert f\right\vert _{w,\theta}\left(  \int\frac
{1}{w\left\vert \cdot\right\vert ^{2\lambda\left(  \cdot\right)  }}\right)
^{1/2}\left\Vert \left(  1+\left\vert \cdot\right\vert \right)  ^{\left\lceil
\sigma\right\rceil }\phi\right\Vert _{\infty},
\end{align*}

where the last line follows from the definition of the function $\lambda$
(weight function property W2).
\end{proof}

\begin{remark}
In the above lemma we use the function $f_{F}$ instead of the tempered
distribution $\widehat{f}\in S^{\prime}$ because $\left\vert \cdot\right\vert
^{\theta}\notin C_{BP}^{\infty}$ and so $\left\vert \cdot\right\vert ^{\theta
}\widehat{f}\notin S^{\prime}$ on $\mathbb{R}^{d}$. On the face of it, the
product looks rather confusing.
\end{remark}

The next lemma is required for the definition of the operator $\mathcal{U}%
:X_{w}^{0}\rightarrow X_{w}^{\theta}$.

\begin{lemma}
\label{Lem_prep_for_op_U}Suppose $w$ has property W1 and W2. The function
$\left\vert \cdot\right\vert ^{-\theta}\widehat{g}$ defines a continuous
linear functional on $S_{\emptyset,\theta}$ by the action
\[
\left[  \left\vert \cdot\right\vert ^{-\theta}\widehat{g},\phi\right]
=\int\left\vert \cdot\right\vert ^{-\theta}\widehat{g}\phi,\text{\quad}\phi\in
S_{\emptyset,\theta}.
\]

\end{lemma}

\begin{proof}
By definition, $g\in X_{w}^{0}$ implies $\widehat{g}\in L_{loc}^{1}$ and hence
$\left\vert \cdot\right\vert ^{-\theta}\widehat{g}\in L_{loc}^{1}\left(
\mathbb{R}^{d}\setminus0\right)  $. We then have for $\phi\in S_{\emptyset
,\theta}$
\[
\int\left\vert \cdot\right\vert ^{-\theta}\left\vert \widehat{g}%
\phi\right\vert =\int\sqrt{w}\left\vert \widehat{g}\right\vert \frac
{\left\vert \phi\right\vert }{\sqrt{w}\left\vert \cdot\right\vert ^{\theta}%
}\leq\left\Vert g\right\Vert _{w,0}\left(  \int\frac{\left\vert \phi
\right\vert ^{2}}{w\left\vert \cdot\right\vert ^{2\theta}}\right)  ^{1/2}.
\]

Now the estimate of $\left(  \int\frac{\left\vert \phi\right\vert ^{2}%
}{w\left\vert \cdot\right\vert ^{2\theta}}\right)  ^{1/2}$ derived in Lemma
\ref{Lem_functnal_phi_sq} yields this lemma.
\end{proof}

We now define the mappings $\mathcal{T}$ and $\mathcal{U}$, which we discussed
before the previous lemma. The mapping $\mathcal{T}$ will be shown to be an
isometry, in the seminorm sense, from $X_{w}^{\theta}$ onto $X_{w}^{0}$. This
is the analogue of the isometric isomorphic Sobolev space mapping $\left(
\left(  1+\left\vert \cdot\right\vert \right)  ^{\theta/2}\widehat{f}\right)
^{\vee}$, between the Hilbert spaces $H^{\theta}$ and $H^{0}=L^{2}$. Such a
mapping allows certain properties of $H^{\theta}$ to be proved if the same
properties have already been established for the simpler space $L^{2}$.
Completeness is the property we are concerned with here.

\begin{definition}
\label{Def_mapping_T} \textbf{The mapping }$\mathcal{T}:X_{w}^{\theta
}\rightarrow S^{\prime}$

Suppose that the weight function $w$ has properties W1 and W2. If $f\in
X_{w}^{\theta}$ then $\widehat{f}\in L_{loc}^{1}\left(  \mathbb{R}%
^{d}\setminus0\right)  $ and we define the function $f_{F}$ by: $f_{F}%
=\widehat{f}$ on $\mathbb{R}^{d}\setminus0$ and $f_{F}\left(  0\right)  =0$.
From part 1 of Lemma \ref{Lem_prep_for_op_T}, $\left\vert \cdot\right\vert
^{\theta}f_{F}\in S^{\prime}$ and so we can define the linear operator
$\mathcal{T}:X_{w}^{\theta}\rightarrow S^{\prime}$ by
\[
\mathcal{T}f=\left(  \left\vert \cdot\right\vert ^{\theta}f_{F}\right)
^{\vee}.
\]

\end{definition}

In the definition above we have used the function $f_{F}$ instead of the
tempered distribution $\widehat{f}\in S^{\prime}$ because $\left\vert
\cdot\right\vert ^{\theta}\notin C_{BP}^{\infty}$ and so $\left\vert
\cdot\right\vert ^{\theta}\widehat{f}\notin S^{\prime}$ on $\mathbb{R}^{d}$.
On the face of it, the product looks rather confusing.

\begin{theorem}
$\mathcal{T}:X_{w}^{\theta}\rightarrow X_{w}^{0}$ and is isometric in the
seminorm sense.
\end{theorem}

\begin{proof}
Suppose $f\in X_{w}^{\theta}$. From the definition of the operator
$\mathcal{T}$, $\mathcal{T}f=\left(  \left\vert \cdot\right\vert ^{\theta
}f_{F}\right)  ^{\vee}\in S^{\prime}$ where the function $f_{F}$ is defined by
$f_{F}=\widehat{f}$ on $\mathbb{R}^{d}\setminus0$ and $f_{F}\left(  0\right)
=0$. From part 4 of Lemma \ref{Lem_properties_Xwm_distrib}, $\widehat
{\mathcal{T}f}=\left\vert \cdot\right\vert ^{\theta}f_{F}\in L_{loc}^{1}$ and
so
\[
\left\Vert \mathcal{T}f\right\Vert _{w,0}^{2}=\int w\left\vert \widehat
{\mathcal{T}f}\right\vert ^{2}=\int w\left\vert \cdot\right\vert ^{2\theta
}\left\vert f_{F}\right\vert ^{2}=\left\vert f\right\vert _{w,\theta}^{2},
\]

where the last equality is part 1 of Theorem \ref{Thm_properties_Xwm}.
\end{proof}

The next step is to define a continuous mapping $\mathcal{U}:X_{w}%
^{0}\rightarrow X_{w}^{\theta}$ which will turn to be the inverse and adjoint
of $\mathcal{T}$.

\begin{definition}
\label{Def_mapping_U} \textbf{The mapping }$\mathcal{U}:X_{w}^{0}\rightarrow
S^{\prime}$

Suppose that the weight function $w$ has properties W1 and W2 and $g\in
X_{w}^{0}$. Since $\left\vert \cdot\right\vert ^{-\theta}\widehat{g}$
satisfies the inequality of the proof of Lemma \ref{Lem_prep_for_op_U}, part 1
of Theorem \ref{Thm_prop_functnl_on_Son} implies that $\left\vert
\cdot\right\vert ^{-\theta}\widehat{g}$ defines a continuous linear functional
on $S_{\emptyset,\theta}$. Hence, by part 2 of Theorem
\ref{Thm_prop_functnl_on_Son}, this functional can be extended (non-uniquely)
to $S$ as an element of $S^{\prime}$, say $f$. Now define the mapping
$\mathcal{U}:X_{w}^{0}\rightarrow S^{\prime}$ by
\[
\mathcal{U}g=\overset{\vee}{f}.
\]

In general, $\mathcal{U}$ is \textbf{non-linear}. However, the next lemma
shows that $\mathcal{U}$ is linear modulo a polynomial in $P_{\theta}$.
\end{definition}

The operator $\mathcal{U}$ has the following properties.

\begin{lemma}
\label{Lem_U_hat} Suppose the weight function $w$ has properties W1 and W2.

\begin{enumerate}
\item Suppose $g\in X_{w}^{0}$. Then $\widehat{\mathcal{U}g}=\left\vert
\cdot\right\vert ^{-\theta}\widehat{g}$ on $S_{\emptyset,\theta}$, and
$\widehat{\mathcal{U}g}\in L_{loc}^{1}\left(  \mathbb{R}^{d}\setminus0\right)
$.

\item $g\in X_{w}^{0}$ implies $\widehat{D^{\alpha}\mathcal{U}g}\in
L_{loc}^{1}$ when $\left\vert \alpha\right\vert =$ $\theta$.

\item $\mathcal{U}$ is linear modulo a polynomial in $P_{\theta}$ i.e. for
scalars $\lambda_{1}$ and $\lambda_{2}$%
\[
\mathcal{U}\left(  \lambda_{1}g_{1}+\lambda_{2}g_{2}\right)  -\lambda
_{1}\mathcal{U}g_{1}-\lambda_{2}\mathcal{U}g_{2}\in P_{\theta}.
\]

\item $\mathcal{U}:X_{w}^{0}\rightarrow X_{w}^{\theta}$ and is isometric in
the seminorm sense.
\end{enumerate}
\end{lemma}

\begin{proof}
\textbf{Part 1}. That $\widehat{\mathcal{U}g}=\left\vert \cdot\right\vert
^{-\theta}\widehat{g}$ on $S_{\emptyset,\theta}$ follows directly from the
definition of $\mathcal{U}$.

By definition of $X_{w}^{0}$, $\widehat{g}\in L_{loc}^{1}$, and so $\left\vert
\cdot\right\vert ^{-\theta}\widehat{g}\in L_{loc}^{1}\left(  \mathbb{R}%
^{d}\setminus0\right)  $.

Hence $\widehat{\mathcal{U}g}=\left\vert \cdot\right\vert ^{-\theta}%
\widehat{g}$ on $\mathbb{R}^{d}\setminus0$, and $\widehat{\mathcal{U}g}\in
L_{loc}^{1}\left(  \mathbb{R}^{d}\setminus0\right)  $.\medskip

\textbf{Part 2}. Suppose $g\in X_{w}^{0}$. We have $\widehat{D^{\alpha
}\mathcal{U}g}\in L_{loc}^{1}\left(  \mathbb{R}^{d}\setminus0\right)  $ and so
if $K$ is compact
\begin{align*}
\int_{K}\left\vert \widehat{D^{\alpha}\mathcal{U}g}\right\vert  & =\int%
_{K}\frac{\left\vert \left(  -i\xi\right)  ^{\alpha}\widehat{g}\left(
\xi\right)  \right\vert }{\left\vert \xi\right\vert ^{\theta}}d\xi\\
& \leq\int_{K}\left\vert \widehat{g}\right\vert \\
& =\int_{K}\frac{1}{\sqrt{w}}\sqrt{w}\left\vert \widehat{g}\right\vert \\
& \leq\left(  \int_{K}\frac{1}{w}\right)  ^{1/2}\left\Vert g\right\Vert
_{w,0}\\
& <\infty,
\end{align*}

since weight function property W2 implies $1/w\in L_{loc}^{1}$.\medskip

\textbf{Part 3}. If $\phi\in S_{\emptyset,\theta}$ then by part 1%
\begin{align*}
& \left[  \left(  \mathcal{U}\left(  \lambda_{1}g_{1}+\lambda_{2}g_{2}\right)
-\lambda_{1}\mathcal{U}g_{1}-\lambda_{2}\mathcal{U}g_{2}\right)  ^{\wedge
},\phi\right] \\
& =\left[  \widehat{\mathcal{U}\left(  \lambda_{1}g_{1}+\lambda_{2}%
g_{2}\right)  }-\lambda_{1}\widehat{\mathcal{U}g_{1}}-\lambda_{2}%
\widehat{\mathcal{U}g_{2}},\phi\right] \\
& =\left[  \left\vert \cdot\right\vert ^{-\theta}\widehat{\lambda_{1}%
g_{1}+\lambda_{2}g_{2}}-\lambda_{1}\left\vert \cdot\right\vert ^{-\theta
}\widehat{g_{1}}-\lambda_{2}\left\vert \cdot\right\vert ^{-\theta}%
\widehat{g_{2}},\phi\right] \\
& =0.
\end{align*}

Part 1 of Theorem \ref{Thm_So,n_and_Pnhat} gives the required result.\medskip

\textbf{Part 4}. Suppose $f\in X_{w}^{\theta}$. Then by part 1 of Lemma
\ref{Lem_U_hat}, $\widehat{\mathcal{T}f}=\left\vert \cdot\right\vert ^{\theta
}f_{F}\in S^{\prime}\cap L_{loc}^{1}$ and so
\[
\left\Vert \mathcal{T}f\right\Vert _{w,0}^{2}=\int w\left\vert \widehat
{\mathcal{T}f}\right\vert ^{2}=\int w\left\vert \cdot\right\vert ^{2\theta
}\left\vert f_{F}\right\vert ^{2}=\left\vert f\right\vert _{w,\theta}^{2},
\]

where the last equality is part 1 of Theorem \ref{Thm_properties_Xwm}. Thus,
by part 3 of Theorem \ref{Thm_properties_Xwm}, $\operatorname*{null}%
\mathcal{T}=\operatorname*{null}\left\vert \cdot\right\vert _{w,\theta
}=P_{\theta}$.\medskip
\end{proof}

The next theorem shows the mappings $\mathcal{T}$ and $\mathcal{U}$ are
inverses and adjoints.

\begin{theorem}
\label{Thm_op_T_properties} Suppose that the weight function $w$ has
properties W1 and W2. Then the operators $\mathcal{T}$ and $\mathcal{U}$ are
related as follows:
\end{theorem}

\begin{enumerate}
\item $\mathcal{TU}=I$ on $X_{w}^{0}$.

\item $\mathcal{UT}f-f\in P_{\theta}$ for $f\in X_{w}^{\theta}$.

\item $\mathcal{U}$ and $\mathcal{T}$ are adjoints i.e. $\left\langle
\mathcal{U}g,f\right\rangle _{w,\theta}=\left(  g,\mathcal{T}f\right)  _{w,0}$.
\end{enumerate}

\begin{proof}
\textbf{Part 1}. Suppose $g\in X_{w}^{0}$. Then by definition of $\mathcal{T}
$, $\widehat{\mathcal{TU}g}=\left\vert \cdot\right\vert ^{\theta}f_{F}\in
L_{loc}^{1}$ where $f_{F}=\widehat{\mathcal{U}g}$ on $\mathbb{R}^{d}%
\setminus0$ and $f_{F}\left(  0\right)  =0$. Now by part 2 of Lemma
\ref{Lem_U_hat}, $\widehat{\mathcal{U}g}=\left\vert \cdot\right\vert
^{-\theta}\widehat{g}$ on $\mathbb{R}^{d}\setminus0$, so $f_{F}=\left\vert
\cdot\right\vert ^{-\theta}\widehat{g}$ on $\mathbb{R}^{d}\setminus0$.

Thus on $\mathbb{R}^{d}\setminus0$, $\widehat{\mathcal{TU}g}=\left\vert
\cdot\right\vert ^{\theta}f_{F}=\widehat{g}$, and since $g\in X_{w}^{0}$
implies $\widehat{g}\in L_{loc}^{1}$, it follows that $\widehat{\mathcal{TU}%
g}=\widehat{g}$ in the distribution sense. We conclude that $\mathcal{TU}%
g=g$.\medskip

\textbf{Part 2}. Suppose $f\in X_{w}^{\theta}$. By part 2 of Lemma
\ref{Lem_U_hat}, $\widehat{\mathcal{UT}f}=\left\vert \cdot\right\vert
^{-\theta}\widehat{\mathcal{T}f}$, on $S_{\emptyset,\theta}$. But by
definition of $\mathcal{T}$, $\widehat{\mathcal{T}f}=\left\vert \cdot
\right\vert ^{\theta}\widehat{f}$ on $\mathbb{R}^{d}\setminus0.$ Thus, if
$\phi\in S_{\emptyset,\theta}$%
\[
\left[  \widehat{\mathcal{UT}f},\phi\right]  =\int\left\vert \cdot\right\vert
^{-\theta}\widehat{\mathcal{T}f}\phi=\int\widehat{f}\phi,
\]

and so $\widehat{\mathcal{UT}f}=\widehat{f}$ on $S_{\emptyset,\theta}$. An
application of part 1 of Theorem \ref{Thm_So,n_and_Pnhat} gives the
result.\medskip

\textbf{Part 3}. From part 1, $\widehat{\mathcal{U}g}=\left\vert
\cdot\right\vert ^{-\theta}\widehat{g}$ on $\mathbb{R}^{d}\setminus0$. Also,
$\widehat{f}=f_{F}$ on $\mathbb{R}^{d}\setminus0$ and $\widehat{\mathcal{T}%
f}=\left\vert \cdot\right\vert ^{\theta}f_{F}$ on $\mathbb{R}^{d}\setminus0$.
Hence
\[
\left\langle \mathcal{U}g,f\right\rangle _{w,\theta}=\int w\left\vert
\cdot\right\vert ^{2\theta}\widehat{\mathcal{U}g}\overline{f_{F}}=\int
w\left\vert \cdot\right\vert ^{2\theta}\left\vert \cdot\right\vert ^{-\theta
}\widehat{g}\overline{f_{F}}=\int w\widehat{g}\left(  \overline{\left\vert
\cdot\right\vert ^{\theta}f_{F}}\right)  =\int w\widehat{g}\overline
{\widehat{\mathcal{T}f}}=\left(  g,\mathcal{T}f\right)  _{w,0}.
\]

\end{proof}

\subsubsection{The operators $\mathcal{T}_{\rho}:X_{w}^{\theta}\rightarrow
X_{w}^{0}$ and $\mathcal{U}_{\rho}:X_{w}^{0}\rightarrow X_{w}^{\theta}$}

Inspired by Light and Wayne's Lemma 2.16, we will introduce the simple adjoint
operators $\mathcal{U}_{\rho}$ and $\mathcal{T}_{\rho}$, and use them to prove
the density of various $C^{\infty}$ function subspaces in $X_{w}^{\theta}$.
These mappings are indicated on Figure \ref{Fig_ops_I_J_U_T_Trho_U_rho}.
Unlike the operators $\mathcal{U}$ and $\mathcal{T}$, which were used to prove
the completeness of $X_{w}^{\theta}$ and assumed properties W1 and W2, these
operators have simple definitions for arbitrary $\theta\geq1$, assuming only
the properties W1. However, the isometry property, which is a key property of
$\mathcal{U}$ and $\mathcal{T}$ , is no longer satisfied by $\mathcal{U}%
_{\rho}$ and $\mathcal{T}_{\rho}$.

From part 3 Theorems \ref{Thm_op_Up_properties} and \ref{Thm_op_Tt_properties}
below, it is easy to see that in a sense these operators can be used to
approximate $\mathcal{U}$ and $\mathcal{T}$ away from the origin, in the
Fourier transform space.

\begin{definition}
\label{Def_op_Up} \textbf{The operator} $\mathcal{U}_{\rho}:X_{w}%
^{0}\rightarrow S^{\prime}$.

Suppose the weight function $w$ has property W1. Suppose $\rho\in
C_{BP}^{\infty}$, $\operatorname*{supp}\rho\subset\mathbb{R}^{d}\setminus0$
and $\rho$ is \textbf{real-valued}.

Then for $\theta\geq1$, $\dfrac{\rho}{\left\vert \cdot\right\vert ^{\theta}%
}\in C_{BP}^{\infty}$ and we can define the linear operator $\mathcal{U}%
_{\rho}:X_{w}^{0}\rightarrow S^{\prime}$ by
\[
\mathcal{U}_{\rho}g=\left(  \frac{\rho}{\left\vert \cdot\right\vert ^{\theta}%
}\widehat{g}\right)  ^{\vee},\text{\qquad}g\in X_{w}^{0}.
\]

\end{definition}

From the point of view of proving that $X_{w}^{\theta}\cap\widehat
{C_{0}^{\infty}}$ is dense in $X_{w}^{\theta}$, the mapping $\mathcal{U}%
_{\rho}$ has the following nice properties.

\begin{theorem}
\label{Thm_op_Up_properties} Suppose that $\rho\in C_{BP}^{\infty}\cap
C_{B}^{\left(  0\right)  }$. Then the operator $\mathcal{U}_{\rho}$ has the
following properties\textbf{:}

\begin{enumerate}
\item $\mathcal{U}_{\rho}:X_{w}^{0}\rightarrow X_{w}^{\theta}$ is continuous
and
\[
\left\vert \mathcal{U}_{\rho}g\right\vert _{w,\theta}\leq\left\Vert
\rho\right\Vert _{\infty}\left\Vert g\right\Vert _{w,0},\text{\quad}f\in
X_{w}^{0}.
\]

\item $\mathcal{U}_{\rho}:X_{w}^{0}\cap\widehat{C_{0}^{\infty}}\rightarrow
X_{w}^{\theta}\cap\widehat{C_{0}^{\infty}}$.

\item If $w$ also has properties W2 then, $\mathcal{U}_{\rho}g=\left(
\rho\widehat{\mathcal{U}g}\right)  ^{\vee}$ when $f\in X_{w}^{0}$.
\end{enumerate}
\end{theorem}

\begin{proof}
\textbf{Part 1}. Clearly $\widehat{\mathcal{U}_{\rho}g}\in L_{loc}^{1}$
because $\widehat{g}\in L_{loc}^{1}$. Hence $\xi^{\alpha}\widehat
{\mathcal{U}_{\rho}g}\in L_{loc}^{1}$ when $\left\vert \alpha\right\vert
=\theta$. Also
\begin{align}
\int w\left\vert \cdot\right\vert ^{2\theta}\left\vert \widehat{\mathcal{U}%
_{\rho}g}\right\vert ^{2}=\int w\left\vert \cdot\right\vert ^{2\theta
}\left\vert \frac{\rho}{\left\vert \cdot\right\vert ^{\theta}}\widehat
{g}\right\vert ^{2}  & \leq\left\vert \rho\right\vert _{\infty}^{2}\int
w\left\vert \widehat{g}\right\vert ^{2}\nonumber\\
& =\left\Vert \rho\right\Vert _{\infty}^{2}\left\Vert g\right\Vert _{w,0}%
^{2}.\label{p1.044}%
\end{align}

Thus, by part 2 of Theorem \ref{Thm_properties_Xwm}, $\mathcal{U}_{\rho}g\in
X_{w}^{\theta}$ and so by part 1 of the same theorem,

$\int w\left|  \cdot\right|  ^{2\theta}\left|  \widehat{\mathcal{U}_{\rho}%
g}\right|  ^{2}=\left|  \mathcal{U}_{\rho}g\right|  _{w,\theta}$ and we have
the desired inequality.\smallskip

\textbf{Part 2}. Suppose $g\in X_{w}^{0}\cap\widehat{C_{0}^{\infty}}$ i.e.
$g\in X_{w}^{0}$ and $\widehat{g}\in C_{0}^{\infty}$. From part 1,
$\mathcal{U}_{\rho}g\in X_{w}^{\theta}$. Clearly, since $\dfrac{\rho
}{\left\vert \cdot\right\vert ^{\theta}}\in C_{BP}^{\infty}$ it follows that
$\widehat{\mathcal{U}_{\rho}g}=\dfrac{\rho}{\left\vert \cdot\right\vert
^{\theta}}\widehat{g}\in C_{0}^{\infty}$ and $\mathcal{U}_{\rho}g\in
X_{w}^{\theta}\cap\widehat{C_{0}^{\infty}}$.\smallskip

\textbf{Part 3}. Choose $\phi\in S$. By part 2 of Lemma \ref{Lem_U_hat},
$\widehat{\mathcal{U}g}=\left\vert \cdot\right\vert ^{-\theta}\widehat{g}$ on
$S_{\emptyset,\theta}$. Hence, since $\rho\phi\in S_{\emptyset,\theta}$%
\begin{align*}
\left[  \rho\widehat{\mathcal{U}g},\phi\right]  =\left[  \widehat
{\mathcal{U}g},\rho\phi\right]  =\left[  \left\vert \cdot\right\vert
^{-\theta}\widehat{g},\rho\phi\right]  =\int\left\vert \cdot\right\vert
^{-\theta}\widehat{g}\rho\phi & =\int\frac{\rho}{\left\vert \cdot\right\vert
^{\theta}}\widehat{g}\phi\\
& =\left[  \widehat{\mathcal{U}_{\rho}g},\phi\right]  .
\end{align*}

Hence $\mathcal{U}_{\rho}g=\left(  \rho\widehat{\mathcal{U}g}\right)  ^{\vee}$.
\end{proof}

\begin{definition}
\label{Def_op_Tt} \textbf{The operator} $\mathcal{T}_{\rho}:X_{w}^{\theta
}\rightarrow S^{\prime}$

Suppose the weight function $w$ has properties W1 for weight set $\mathcal{A}
$ i.e. $w$ is continuous and positive on $\mathbb{R}^{d}\setminus\mathcal{A}
$. Suppose $\rho\in C_{BP}^{\infty}$, $\operatorname*{supp}\rho\subset
\mathbb{R}^{d}\setminus0$ and $\rho$ is \textbf{real-valued}.

Then for $\theta\geq0$, $\rho\left\vert \cdot\right\vert ^{\theta}\in
C_{BP}^{\infty}$ and we can define the linear operator $\mathcal{T}_{\rho
}:X_{w}^{\theta}\rightarrow S^{\prime}$ by
\[
\mathcal{T}_{\rho}f=\left(  \rho\left\vert \cdot\right\vert ^{\theta}%
\widehat{f}\right)  ^{\vee},\text{\qquad}f\in X_{w}^{\theta}.
\]

\end{definition}

The mapping $\mathcal{T}_{\rho}$ has the following properties.

\begin{theorem}
\label{Thm_op_Tt_properties} Suppose that $\rho\in C_{BP}^{\infty}\cap
C_{B}^{\left(  0\right)  }$. Then the operator $\mathcal{T}_{\rho}$ has the
following properties \textbf{:}

\begin{enumerate}
\item $\mathcal{T}_{\rho}:X_{w}^{\theta}\rightarrow X_{w}^{0}$ is continuous
and
\[
\left\Vert \mathcal{T}_{\rho}f\right\Vert _{w,0}\leq\left\Vert \rho\right\Vert
_{\infty}\left\vert f\right\vert _{w,\theta},\text{\quad}f\in X_{w}^{\theta}.
\]

\item $\mathcal{T}_{\rho}:X_{w}^{\theta}\cap\widehat{C_{0}^{\infty}%
}\rightarrow X_{w}^{0}\cap\widehat{C_{0}^{\infty}}$.

\item If $w$ also has property W2 then, $\mathcal{T}_{\rho}f=\left(
\rho\widehat{\mathcal{T}f}\right)  ^{\vee}$ when $f\in X_{w}^{\theta}$.
\end{enumerate}
\end{theorem}

\begin{proof}
\textbf{Part 1}. Clearly $\widehat{\mathcal{T}_{\rho}f}\in L_{loc}^{1}$
because $\widehat{f}\in L_{loc}^{1}$. Hence $\xi^{\alpha}\widehat
{\mathcal{T}_{\rho}f}\in L_{loc}^{1}$ when $\left\vert \alpha\right\vert =$
$\theta$. Also
\[
\left\Vert \mathcal{T}_{\rho}f\right\Vert _{w,0}^{2}=\int w\left\vert
\widehat{\mathcal{T}_{\rho}f}\right\vert ^{2}=\int w\left\vert \rho\left\vert
\cdot\right\vert ^{\theta}\widehat{f}\right\vert ^{2}\leq\left\Vert
\rho\right\Vert _{\infty}^{2}\int w\left\vert \cdot\right\vert ^{2\theta
}\left\vert \widehat{f}\right\vert ^{2}=\left\Vert \rho\right\Vert _{\infty
}^{2}\left\vert f\right\vert _{w,\theta}^{2}.
\]
\smallskip

\textbf{Part 2}. Suppose $f\in X_{w}^{\theta}\cap\widehat{C_{0}^{\infty}}$
i.e. $f\in X_{w}^{\theta}$ and $\widehat{f}\in C_{0}^{\infty}$. From part 1,
$\mathcal{T}_{\rho}f\in X_{w}^{0}$. Clearly, since $\rho\left\vert
\cdot\right\vert ^{\theta}\in C_{BP}^{\infty}$ it follows that $\widehat
{\mathcal{T}_{\rho}f}=\rho\left\vert \cdot\right\vert ^{\theta}\widehat{f}\in
C_{0}^{\infty}$ and $\mathcal{T}_{\rho}f\in X_{w}^{\theta}\cap\widehat
{C_{0}^{\infty}}$.\smallskip

\textbf{Part 3}. Since $w$ has properties W1 and W2 we can define
$\mathcal{T}$. From the definition of $\mathcal{T}$, $\widehat{\mathcal{T}%
f}=\left\vert \cdot\right\vert ^{\theta}\widehat{f}$ on $\mathbb{R}%
^{d}\setminus0$, and hence $\operatorname*{supp}\left(  \widehat{\mathcal{T}%
f}-\left\vert \cdot\right\vert ^{\theta}\widehat{f}\right)  \subset\left\{
0\right\}  $. But $\operatorname*{supp}\rho\subset\mathbb{R}^{d}\setminus0$,
so $\rho\left(  \left\vert \cdot\right\vert ^{\theta}\widehat{f}%
-\widehat{\mathcal{T}f}\right)  =0$ i.e. $\rho\left\vert \cdot\right\vert
^{\theta}\widehat{f}=\rho\widehat{\mathcal{T}f}$. But $\widehat{\mathcal{T}%
_{\rho}f}=\rho\left\vert \cdot\right\vert ^{\theta}\widehat{f}$ so
$\widehat{\mathcal{T}_{\rho}f}=\rho\widehat{\mathcal{T}f}$, as required.
\end{proof}

The mappings $\mathcal{T}_{\rho}$ and $\mathcal{U}_{\rho}$ are related as follows.

\begin{theorem}
\label{Thm_ops_Tt_Up_properties}. The mappings $\mathcal{T}_{\rho}$ and
$\mathcal{U}_{\rho}$ are related as follows:

\begin{enumerate}
\item $\widehat{\mathcal{T}_{\rho}\mathcal{U}_{\rho}g}=\rho^{2}\widehat{g}
$,\quad$g\in X_{w}^{0}$.

\item $\widehat{\mathcal{U}_{\rho}\mathcal{T}_{\rho}f}=\rho^{2}\widehat{f}
$,\quad$f\in X_{w}^{\theta}$.

\item $\mathcal{T}_{\rho}$ and $\mathcal{U}_{\rho}$ are adjoints i.e. $\left(
\mathcal{T}_{\rho}f,g\right)  _{w,0}=\left\langle f,\mathcal{U}_{\rho
}g\right\rangle _{w,\theta}$,\quad$f\in X_{w}^{\theta}$, $g\in X_{w}^{0}$.

\item If $f\in X_{w}^{\theta}$ and $g\in X_{w}^{0}$ then
\[
\left\vert f-\mathcal{U}_{\rho}g\right\vert _{w,\theta}^{2}+\int\left(
1-\rho^{2}\right)  w\left\vert \widehat{g}\right\vert ^{2}=\left\Vert
\mathcal{T}_{\rho}f-g\right\Vert _{w,0}^{2}+\int\left(  1-\rho^{2}\right)
w\left\vert \cdot\right\vert ^{2\theta}\left\vert \widehat{f}\right\vert ^{2}.
\]

\end{enumerate}
\end{theorem}

\begin{proof}
\textbf{Part 1}. $\widehat{\mathcal{T}_{\rho}\mathcal{U}_{\rho}g}%
=\rho\left\vert \cdot\right\vert ^{\theta}\widehat{\mathcal{U}_{\rho}g}%
=\rho\left\vert \cdot\right\vert ^{\theta}\frac{\rho}{\left\vert
\cdot\right\vert ^{\theta}}\widehat{g}=\rho^{2}\widehat{g}$.

\textbf{Part 2}. $\widehat{\mathcal{U}_{\rho}\mathcal{T}_{\rho}f}=\frac{\rho
}{\left|  \cdot\right|  ^{\theta}}\widehat{\mathcal{T}_{\rho}f}=\frac{\rho
}{\left|  \cdot\right|  ^{\theta}}\rho\left|  \cdot\right|  ^{\theta}%
\widehat{f}=\rho^{2}\widehat{f}$.

\textbf{Part 3}. Since $\rho$ is real-valued
\[
\left\langle f,\mathcal{U}_{\rho}g\right\rangle _{w,\theta}=\int w\left|
\cdot\right|  ^{2\theta}\widehat{f}\overline{\widehat{\mathcal{U}_{\rho}g}%
}=\int w\left|  \cdot\right|  ^{2\theta}\widehat{f}\frac{\rho}{\left|
\cdot\right|  ^{\theta}}\overline{\widehat{g}}=\int w\left(  \rho\left|
\cdot\right|  ^{\theta}\widehat{f}\right)  \overline{\widehat{g}}=\int
w\widehat{\mathcal{T}_{\rho}f}\overline{\widehat{g}}=\left(  \mathcal{T}%
_{\rho}f,g\right)  _{w,0}.
\]

\textbf{Part 4}. Using part 3 of this theorem
\begin{align*}
\left\vert f-\mathcal{U}_{\rho}g\right\vert _{w,\theta}^{2}-\left\Vert
\mathcal{T}_{\rho}f-g\right\Vert _{w,0}^{2}  & =\left\langle f-\mathcal{U}%
_{\rho}g,f-\mathcal{U}_{\rho}g\right\rangle _{w,\theta}-\left(  \mathcal{T}%
_{\rho}f-g,\mathcal{T}_{\rho}f-g\right)  _{w,0}\\
& =\left\langle f,f\right\rangle _{w,\theta}-\left\langle f,\mathcal{U}_{\rho
}g\right\rangle _{w,\theta}-\left\langle \mathcal{U}_{\rho}g,f\right\rangle
_{w,\theta}+\left\langle \mathcal{U}_{\rho}g,\mathcal{U}_{\rho}g\right\rangle
_{w,\theta}-\\
& -\left(  \mathcal{T}_{\rho}f,\mathcal{T}_{\rho}f\right)  _{w,0}+\left(
\mathcal{T}_{\rho}f,g\right)  _{w,0}+\left(  g,\mathcal{T}_{\rho}f\right)
_{w,0}-\left(  g,g\right)  _{w,0}\\
& =\left\langle f,f\right\rangle _{w,\theta}+\left\langle \mathcal{U}_{\rho
}g,\mathcal{U}_{\rho}g\right\rangle _{w,\theta}-\left(  \mathcal{T}_{\rho
}f,\mathcal{T}_{\rho}f\right)  _{w,0}-\left(  g,g\right)  _{w,0}\\
& =\left\langle f,f\right\rangle _{w,\theta}+\left(  \mathcal{T}_{\rho
}\mathcal{U}_{\rho}g,g\right)  _{w,0}-\left(  f,\mathcal{U}_{\rho}%
\mathcal{T}_{\rho}f\right)  _{w,\theta}-\left(  g,g\right)  _{w,0}\\
& =\left\langle f,f-\mathcal{U}_{\rho}\mathcal{T}_{\rho}f\right\rangle
_{w,\theta}+\left(  \mathcal{T}_{\rho}\mathcal{U}_{\rho}g-g,g\right)  _{w,0}\\
& =\int w\left\vert \cdot\right\vert ^{2\theta}\widehat{f}\overline
{\widehat{\left(  f-\mathcal{U}_{\rho}\mathcal{T}_{\rho}f\right)  }}+\int
w\widehat{\left(  \mathcal{T}_{\rho}\mathcal{U}_{\rho}g-g\right)  }%
\overline{\widehat{g}}\\
& =\int\left(  1-\rho^{2}\right)  w\left\vert \cdot\right\vert ^{2\theta
}\left\vert \widehat{f}\right\vert ^{2}-\int\left(  1-\rho^{2}\right)
w\left\vert \widehat{g}\right\vert ^{2},
\end{align*}

where the last step used parts 1 and 2 of this theorem.
\end{proof}

\subsubsection{Dense $C^{\infty}$ subspaces of $X_{w}^{\theta}$}

A feature here is the use of the simple mappings $\mathcal{U}_{\rho}$ and
$\mathcal{T}_{\rho}$ of Definitions \ref{Def_op_Tt} and \ref{Def_op_Up} and
the corresponding results for $X_{w}^{0}$ to prove the density of various
$C^{\infty}$ function subspaces in $X_{w}^{\theta}$. Unlike the operators
$\mathcal{U}$ and $\mathcal{T}$ these operators are easily defined for
arbitrary $\theta\geq0$.

The next theorem corresponds to Theorem 2.17 of Light and Wayne
\cite{LightWayne95ErrEst}. Because we have used the operators $\mathcal{U}%
_{\rho} $ and $\mathcal{T}_{\rho}$ instead of $\mathcal{U}$ and $\mathcal{T}$
we will only need to assume weight function properties W1, which correspond to
properties A3.1 and A3.2 of Light and Wayne \cite{LightWayneX98Weight}.

\begin{theorem}
\label{Thm_Co_inf_dense_Xwth}Suppose the weight function $w$ has properties W1
and W2.

Then $X_{w}^{\theta}\cap\widehat{C_{0}^{\infty}}$ and $X_{w}^{\theta}\cap S$
are dense in $X_{w}^{\theta}$ for $\theta\geq1$.
\end{theorem}

\begin{proof}
Choose $\varepsilon>0$ and $f\in X_{w}^{\theta}$. Then $\mathcal{T}_{\rho}f\in
X_{w}^{0}$ and, since $X_{w}^{0}\cap\widehat{C_{0}^{\infty}}$ is dense in
$X_{w}^{0}$ (Subsection 1.3.4 of the zero order document Williams
\cite{WilliamsZeroOrdSmthV4}), there exists $g_{\rho}\in X_{w}^{0}\cap
\widehat{C_{0}^{\infty}}$ such that $\left\Vert \mathcal{T}_{\rho}f-g_{\rho
}\right\Vert _{w,0}\leq\varepsilon/4$. Now $\mathcal{U}_{\rho}g_{\rho}\in
X_{w}^{\theta}\cap\widehat{C_{0}^{\infty}}$, and it will be shown that we can
choose $\rho$ such that $\left\vert f-\mathcal{U}_{\rho}g_{\rho}\right\vert
_{w,\theta}\leq\varepsilon$.

From part 4 of Theorem \ref{Thm_ops_Tt_Up_properties}
\begin{equation}
\left\vert f-\mathcal{U}_{\rho}g_{\rho}\right\vert _{w,\theta}\leq\left\Vert
\mathcal{T}_{\rho}f-g_{\rho}\right\Vert _{w,0}+\left(  \int\left(  1-\rho
^{2}\right)  w\left\vert \cdot\right\vert ^{2\theta}\left\vert \widehat
{f}\right\vert ^{2}\right)  ^{1/2}.\label{p1.054}%
\end{equation}

At the start of this proof $g_{\rho}$ was chosen to satisfy $\left\Vert
\mathcal{T}_{\rho}f-g_{\rho}\right\Vert _{w,0}\leq\varepsilon/2$, so now
\[
\left\vert f-\mathcal{U}_{\rho}g_{\rho}\right\vert _{w,\theta}\leq
\varepsilon/2+2\left(  \int\left(  1-\rho^{2}\right)  w\left\vert
\cdot\right\vert ^{2\theta}\left\vert \widehat{f}\right\vert ^{2}\right)
^{1/2}.
\]

Because $\int w\left\vert \cdot\right\vert ^{2\theta}\left\vert \widehat
{f}\right\vert ^{2}$ exists, $\eta>0$ can be selected such that $\int%
\limits_{\left\vert \cdot\right\vert \leq\eta}w\left\vert \cdot\right\vert
^{2\theta}\left\vert \widehat{f}\right\vert ^{2}\leq\left(  \varepsilon
/2\right)  ^{2}$. Then, if $\rho$ satisfies $\rho\left(  x\right)  =1$ when
$\left\vert x\right\vert \geq\eta$
\[
\int\left(  1-\rho^{2}\right)  w\left\vert \cdot\right\vert ^{2\theta
}\left\vert \widehat{f}\right\vert ^{2}=\int\limits_{\left\vert \cdot
\right\vert \leq\eta}w\left\vert \cdot\right\vert ^{2\theta}\left\vert
\widehat{f}\right\vert ^{2}\leq\left(  \varepsilon/2\right)  ^{2},
\]

and $\left\vert f-\mathcal{U}_{\rho}g_{\rho}\right\vert _{w,\theta}%
\leq\varepsilon/2+\varepsilon/2=\varepsilon$ as required.
\end{proof}

??? Next we prove that $\left(  C_{0}^{\infty}\left(  \mathbb{R}^{d}%
\setminus\mathcal{A}\right)  \right)  ^{\vee}$ is dense in $X_{w}^{\theta}$
for $\theta\geq1$. Use the above operators and the corresponding result for
$X_{w}^{0}$ proven in Subsection 1.3.4 of the zero order document Williams
\cite{WilliamsZeroOrdSmthV4}.

\subsection{The space $X_{w:K}^{\theta}$}

\begin{definition}
\label{Def_XowK}If $K\subset\mathbb{R}^{d}$ is a closed set then define%
\[
X_{w:K}^{\theta}:=\left(  X_{w}^{\theta}\right)  _{K}:=\left\{  f\in
X_{w}^{\theta}:\operatorname*{supp}f\subseteq K\right\}  ,
\]

which we endow with the subspace norm.

This notation is analogous to the Sobolev space notation of Petersen
\cite{Petersen83} following Exercise 3:4.10 on p243.
\end{definition}

??? Is $X_{w:K}^{\theta}$ a semi-Hilbert space? An attempted proof based on
the (zero order) proof for $X_{w:K}^{0}$ and which used the Light norm
\ref{p917} failed.

\subsection{Summary figure for mappings, completeness and density results}

Assuming weight properties W1 and W2 hold, Figure
\ref{Fig_ops_I_J_U_T_Trho_U_rho} and subsequent notes summarizes many of the
results of the previous sections of this chapter concerning the spaces
$X_{w}^{0}$ and $X_{w}^{\theta}$. One can think of this as a ladder of
mappings going from simple to more complex spaces and from lower to higher
order spaces. It is an analogy with the Sobolev spaces of non-negative order.%
\begin{equation}%
\begin{array}
[c]{ccccccc}
& \underset{(Hilbert)}{L^{2}} & \underset{\mathcal{J}=\mathcal{I}%
^{-1}}{\overset{\mathcal{I}}{%
{\tt\setlength{\unitlength}{0.5pt}
\begin{picture}(160,20)
\thinlines\put(0,0){\vector(1,0){160}}
\put(160,15){\vector(-1,0){160}}
\end{picture}}%
}} & \underset{(Hilbert)}{X_{w}^{0}} & \underset{\mathcal{U},\text{
}\mathcal{U}_{\rho}}{\overset{\mathcal{T},\text{ }\mathcal{T}_{\rho}}{%
{\tt\setlength{\unitlength}{0.5pt}
\begin{picture}(170,20)
\thinlines\put(0,0){\vector(1,0){170}}
\put(170,15){\vector(-1,0){170}}
\end{picture}}%
}} & \underset{(semi-Hilbert)}{X_{w}^{\theta}} & \\
& \underset{dense\text{ }embed}{%
{\tt\setlength{\unitlength}{0.5pt}
\begin{picture}(5,70)
\thinlines\put(0,0){\vector(0,1){50}}
\end{picture}}%
} &  & \underset{dense\text{ }embed}{%
{\tt\setlength{\unitlength}{0.5pt}
\begin{picture}(5,70)
\thinlines\put(0,0){\vector(0,1){50}}
\end{picture}}%
} &  & \underset{dense\text{ }embed}{%
{\tt\setlength{\unitlength}{0.5pt}
\begin{picture}(5,70)
\thinlines\put(0,0){\vector(0,1){50}}
\end{picture}}%
} & \\
&
{\tt\setlength{\unitlength}{0.5pt}
\begin{picture}(5,70)
\thinlines\put(0,0){\line(0,1){50}}
\end{picture}}%
&  &
{\tt\setlength{\unitlength}{0.5pt}
\begin{picture}(5,70)
\thinlines\put(0,0){\line(0,1){50}}
\end{picture}}%
&  &
{\tt\setlength{\unitlength}{0.5pt}
\begin{picture}(5,70)
\thinlines\put(0,0){\line(0,1){50}}
\end{picture}}%
& \\
& L^{2}\cap\sqrt{w}\widehat{C_{0}^{\infty}} & \underset{\mathcal{I}}{%
{\tt\setlength{\unitlength}{0.5pt}
\begin{picture}(160,20)
\thinlines\put(160,10){\vector(-1,0){160}}
\end{picture}}%
} & X_{w}^{0}\cap\widehat{C_{0}^{\infty}} & \underset{\mathcal{U}_{\rho
}}{\overset{\mathcal{T}_{\rho}}{%
{\tt\setlength{\unitlength}{0.5pt}
\begin{picture}(170,20)
\thinlines\put(0,0){\vector(1,0){170}}
\put(170,15){\vector(-1,0){170}}
\end{picture}}%
}} & X_{w}^{\theta}\cap\widehat{C_{0}^{\infty}} & \\
&  &  &  &  &  &
\end{array}
\label{Fig_ops_I_J_U_T_Trho_U_rho}%
\end{equation}

Notes on Figure \ref{Fig_ops_I_J_U_T_Trho_U_rho}:

\begin{enumerate}
\item The definition of operators $\mathcal{I}:X_{w}^{\theta}\rightarrow
L^{2}$ and $\mathcal{J}:L^{2}\rightarrow X_{w}^{\theta}$ requires weight
function properties W1 and W2.

\item In the seminorm sense, $\mathcal{I}$ and $\mathcal{J}$ are isometric
isomorphisms, adjoints and inverses.

\item The definition of operators $\mathcal{T}:X_{w}^{0}\rightarrow
X_{w}^{\theta}$ and $\mathcal{U}:X_{w}^{\theta}\rightarrow X_{w}^{0}$ requires
weight function properties W1 and W2.

\item In the seminorm sense, $\mathcal{T}$ and $\mathcal{U}$ are isometric
isomorphisms, adjoints and inverses

\item The definition of operators $\mathcal{T}_{\rho}:X_{w}^{0}\rightarrow
X_{w}^{\theta}$ and $\mathcal{U}_{\rho}:X_{w}^{\theta}\rightarrow X_{w}^{0}$
only requires weight function property W1.

\item $\mathcal{T}_{\rho}$ and $\mathcal{U}_{\rho}$ are adjoints.

\item ??? Replace $C_{0}^{\infty}$ by $S$?
\end{enumerate}

\subsection{The smoothness of functions in $X_{w}^{\theta}$%
\label{SbSect_Xw,th_smooth}}

This theorem provides information about the $L^{1}$ behavior of the Fourier
transform of functions in $X_{w}^{\theta}$, $\theta\geq1$. This information
will be used to prove smoothness properties for $X_{w}^{\theta}$ when the
weight function has property W3.1 or properties W2.1 and W3.2 or W3.3.

\begin{theorem}
\label{Thm_Xwm_is_L1loc(Rd_0)}If a weight function $w$ has property W2 then
given $f\in X_{w}^{\theta}$ we can define a.e. the function $f_{F}%
:\mathbb{R}^{d}\rightarrow\mathbb{C}$ by $f_{F}=\widehat{f}$ on $\mathbb{R}%
^{d}\setminus0$. Now we have the following results:

\begin{enumerate}
\item If $w$ has property \textbf{W3.1} for order $\theta$ and $\kappa$, then
$\xi^{\beta}f_{F}\in L^{1}$ when $\beta\leq\kappa$.

\item Suppose $w$ has property W2.1 and also property \textbf{W3.2} for order
$\theta$ and $\kappa$. Choose $\rho\in C_{0}^{\infty}$ such that, $0\leq
\rho\leq1$, $\rho\left(  x\right)  =1$ when $\left\vert x\right\vert \leq
r_{3}$. Then, $\left(  1-\rho\right)  \xi^{\beta}f_{F}\in L^{1}$ when
$\left\vert \beta\right\vert \leq\kappa$.

\item Suppose $w$ has property W2.1 and also property \textbf{W3.3} for order
$\theta$ and $\kappa$, then $\xi^{\beta}f_{F}\in L^{1}$ when $\beta\leq\kappa$.
\end{enumerate}
\end{theorem}

\begin{proof}
\textbf{Part 1} First note that by part 2 of Theorem
\ref{Thm_weight_property_relat} property W3.1 implies property W2 and so
$f_{F}$ is defined. Now by the definition of property W3.1, for all $\alpha$
satisfying $\left\vert \alpha\right\vert =\theta$, $%
{\displaystyle\int}
\dfrac{\xi^{2\lambda}}{w\,\xi^{2\alpha}}<\infty$ for $0\leq\lambda\leq\kappa$.
So if $\beta\leq\kappa$, by using the Cauchy-Schwartz inequality, we obtain
\[
\int\left\vert \xi^{\beta}f_{F}\right\vert =\int\frac{\left\vert \xi^{\beta
}\right\vert }{\sqrt{w}\left\vert \xi^{\alpha}\right\vert }\sqrt{w}\left\vert
\xi^{\alpha}\right\vert \left\vert f_{F}\right\vert \leq\int\frac{\left\vert
\xi^{\beta}\right\vert }{\sqrt{w}\left\vert \xi^{\alpha}\right\vert }\sqrt
{w}\left\vert \cdot\right\vert ^{\theta}\left\vert f_{F}\right\vert =\left(
\int\frac{\xi^{2\beta}}{w\text{ }\xi^{2\alpha}}\right)  ^{\frac{1}{2}%
}\left\vert f\right\vert _{w,\theta},
\]

where the last step used part 1 of Theorem \ref{Thm_properties_Xwm}. The term
on the last line is finite since $w$ has property W3.1.\medskip

\textbf{Part 2} First note that, by part 3 of Theorem
\ref{Thm_weight_property_relat}, property W3.2 implies property W2.2. Hence
$w$ has property W2 and so $f_{F}$ is defined.

Since $\widehat{f}\in L_{loc}^{1}\left(  \mathbb{R}^{d}\setminus0\right)  $,
$\left(  1-\rho\right)  f_{F}\in L_{loc}^{1}$ and $\left(  1-\rho\right)
\xi^{\beta}f_{F}\in L_{loc}^{1}$. Further, since $\rho\left(  \xi\right)  =1$
when $\left\vert \xi\right\vert \leq r_{3}>0$, we can use the Cauchy-Schwartz
inequality to write
\begin{align*}
\int\left\vert \left(  1-\rho\left(  \xi\right)  \right)  \xi^{\beta}%
f_{F}\left(  \xi\right)  d\xi\right\vert \leq\int\limits_{\left\vert
\cdot\right\vert \geq r_{3}}\left\vert \cdot\right\vert ^{\left\vert
\beta\right\vert }\left\vert f_{F}\right\vert  &  =\int\limits_{\left\vert
\cdot\right\vert \geq r_{3}}\frac{\left\vert \cdot\right\vert ^{\kappa}%
}{\left\vert \cdot\right\vert ^{\left(  \kappa-\left\vert \beta\right\vert
\right)  }\sqrt{w}\left\vert \cdot\right\vert ^{\theta}}\sqrt{w}\left\vert
\cdot\right\vert ^{\theta}\left\vert f_{F}\right\vert \\
&  \leq\frac{1}{\left(  r_{3}\right)  ^{\left(  \kappa-\left\vert
\beta\right\vert \right)  }}\int\limits_{\left\vert \cdot\right\vert \geq
r_{3}}\frac{\left\vert \cdot\right\vert ^{\kappa}}{\sqrt{w}\left\vert
\cdot\right\vert ^{\theta}}\sqrt{w}\left\vert \cdot\right\vert ^{\theta
}\left\vert f_{F}\right\vert \\
&  \leq\frac{1}{\left(  r_{3}\right)  ^{\left(  \kappa-\left\vert
\beta\right\vert \right)  }}\left(  \int\limits_{\left\vert \cdot\right\vert
\geq r_{3}}\frac{\left\vert \cdot\right\vert ^{2\kappa}}{w\left\vert
\cdot\right\vert ^{2\theta}}\right)  ^{\frac{1}{2}}\left(  \int%
\limits_{\left\vert \cdot\right\vert \geq r_{3}}w\left\vert \cdot\right\vert
^{2\theta}\left\vert f_{F}\right\vert ^{2}\right)  ^{\frac{1}{2}}\\
&  \leq\frac{1}{\left(  r_{3}\right)  ^{\left(  \kappa-\left\vert
\beta\right\vert \right)  }}\left(  \int\limits_{\left\vert \cdot\right\vert
\geq r_{3}}\frac{\left\vert \cdot\right\vert ^{2\kappa}}{w\left\vert
\cdot\right\vert ^{2\theta}}\right)  ^{\frac{1}{2}}\left\vert f\right\vert
_{w,\theta}.
\end{align*}

Since $w$ has property W3.2 for order $\theta$ and $\kappa$, the last integral
is finite.\medskip

\textbf{Part 3} First note that by part 4 of Theorem
\ref{Thm_weight_property_relat} property W3.3 implies property W2.2. Thus
$w\in W2$ and $f_{F}$ is defined. Now by the definition of property W3.3, $%
{\displaystyle\int}
\dfrac{\xi^{2\lambda}}{w\left\vert \cdot\right\vert ^{2\theta}}<\infty$ for
$0\leq\lambda\leq\kappa$. So if $\beta\leq\kappa$, by using the
Cauchy-Schwartz inequality, we obtain
\[
\int\left\vert \xi^{\beta}f_{F}\right\vert =\int\frac{\left\vert \xi^{\beta
}\right\vert }{\sqrt{w}\left\vert \cdot\right\vert ^{\theta}}\sqrt
{w}\left\vert \cdot\right\vert ^{\theta}\left\vert f_{F}\right\vert
\leq\left(  \int\frac{\xi^{2\beta}}{w\left\vert \cdot\right\vert ^{2\theta}%
}\right)  ^{\frac{1}{2}}\left\vert f\right\vert _{w,\theta},
\]

where the last step used part 1 of Theorem \ref{Thm_properties_Xwm}. The term
on the last line is finite since $w$ has property W3.3.
\end{proof}

\begin{definition}
\label{Def_Ca_smth_space}\textbf{Some} \textbf{spaces of continuous
functions:} For each multi-index $\alpha$ define:%
\begin{align*}
C^{\left(  \alpha\right)  }  & :=\left\{  u\in C^{\left(  0\right)  }%
:D^{\beta}u\in C^{\left(  0\right)  }\text{ }for\text{ }\beta\leq
\alpha\right\}  ,\\
C_{B}^{\left(  \alpha\right)  }  & :=\left\{  u\in C_{B}^{\left(  0\right)
}:D^{\beta}u\in C_{B}^{\left(  0\right)  }\text{ }for\text{ }\beta\leq
\alpha\right\}  ,
\end{align*}

and note that%
\[
C^{\left(  \alpha\right)  }\subset C^{\left(  \underline{\alpha}\right)  },
\]

where $\underline{\alpha}=\min\limits_{i}\alpha_{i}$.
\end{definition}

This theorem corresponds to Light and Wayne's \cite{LightWayneX98Weight}
Theorem 2.18 and represents our main smoothness result for $X_{w}^{\theta}$ spaces.

\begin{theorem}
\label{Thm_Xwth_W3_smooth}\textbf{Smoothness of }$X_{w}^{\theta}$ Suppose the
weight function $w$ has properties W2.1 and W3 for order $\theta$ and $\kappa
$, as described in Definition \ref{Def_extend_wt_fn}. Then, if $\left\lfloor
\kappa\right\rfloor :=\left(  \left\lfloor \kappa_{k}\right\rfloor \right)
$,
\[
X_{w}^{\theta}\subset C_{BP}^{\infty}+C_{B}^{\left(  \left\lfloor
\kappa\right\rfloor \right)  }\subset C_{BP}^{\left(  \left\lfloor
\kappa\right\rfloor \right)  },
\]

and also%
\[
X_{w}^{\theta}\subset C_{BP}^{\infty}+C_{B}^{\left(  \left\lfloor
\underline{\kappa}\right\rfloor \right)  }\subset C_{BP}^{\left(  \left\lfloor
\underline{\kappa}\right\rfloor \right)  }.
\]

\end{theorem}

\begin{proof}
\fbox{\textbf{Case 1} $w$ has property W3.1} First observe that property
W3.1$^{\text{*}}$ implies property W3.1. By part 3 of Theorem
\ref{Thm_weight_property_relat} property W3.1 implies property W2.

Now suppose $\rho_{0}\in C_{0}^{\infty}$, $0\leq\rho_{0}\leq1$, and $\rho
_{0}=1$ on the ball $B\left(  0;1\right)  $. By part 1 of Theorem
\ref{Thm_Xwm_is_L1loc(Rd_0)}, $f\in X_{w}^{\theta}$ implies $\xi^{\beta}%
f_{F}\in L^{1}$ when $\beta\leq\kappa$ and so $\xi^{\beta}\left(  1-\rho
_{0}\right)  \widehat{f}=\xi^{\beta}\left(  1-\rho_{0}\right)  f_{F}\in L^{1}$
when $\beta\leq\left\lfloor \kappa\right\rfloor $. Hence $D^{\beta}\left(
\left(  1-\rho_{0}\right)  \widehat{f}\right)  ^{\vee}\in C_{B}^{\left(
0\right)  }$ for $\beta\leq\left\lfloor \kappa\right\rfloor $ and so $\left(
\left(  1-\rho_{0}\right)  \widehat{f}\right)  ^{\vee}\in C_{B}^{\left(
\left\lfloor \kappa\right\rfloor \right)  }$.

We now write $\widehat{f}=\rho_{0}\widehat{f}+\left(  1-\rho_{0}\right)
\widehat{f}$. Taking the inverse Fourier transform of this equation yields
\[
f=\left(  \rho_{0}\widehat{f}\right)  ^{\vee}+\left(  \left(  1-\rho
_{0}\right)  \widehat{f}\right)  ^{\vee}.
\]

The first term on the right is the inverse Fourier transform of a tempered
distribution with compact support and so by Appendix \ref{SbSect_property_E'},
it is a $C_{BP}^{\infty}$ function. Thus $f\in C_{BP}^{\infty}+C_{B}^{\left(
\left\lfloor \kappa\right\rfloor \right)  }\subset C_{BP}^{\left(
\left\lfloor \kappa\right\rfloor \right)  }$ and, from Definition
\ref{Def_Ca_smth_space}, $f\in C_{BP}^{\infty}+C_{B}^{\left(  \left\lfloor
\underline{\kappa}\right\rfloor \right)  }\subset C_{BP}^{\left(  \left\lfloor
\underline{\kappa}\right\rfloor \right)  }$.\medskip

\fbox{\textbf{Case 2} $w$ has properties W2.1 and W3.2} By part 3 of Theorem
\ref{Thm_weight_property_relat} properties W2.1 and W3.2 imply property W2.
Now suppose $\rho_{0}\in C_{0}^{\infty}$, $0\leq\rho_{0}\leq1$, and $\rho
_{0}=1$ on the ball $B\left(  0;1\right)  $. By part 2 of Theorem
\ref{Thm_Xwm_is_L1loc(Rd_0)}, $f\in X_{w}^{\theta}$ implies $\xi^{\beta
}\left(  1-\rho_{0}\right)  \widehat{f}=\xi^{\beta}\left(  1-\rho_{0}\right)
f_{F}\in L^{1}$ when $\left\vert \beta\right\vert \leq\left\lfloor
\kappa\right\rfloor $. Hence $D^{\beta}\left(  \left(  1-\rho_{0}\right)
\widehat{f}\right)  ^{\vee}\in C_{B}^{\left(  0\right)  }$ for $\left\vert
\beta\right\vert \leq\left\lfloor \kappa\right\rfloor $ and so $\left(
\left(  1-\rho_{0}\right)  \widehat{f}\right)  ^{\vee}\in C_{B}^{\left(
\left\lfloor \kappa\right\rfloor \right)  }$.

We now write $\widehat{f}=\rho_{0}\widehat{f}+\left(  1-\rho_{0}\right)
\widehat{f}$. Taking the inverse Fourier transform of this equation yields
\[
f=\left(  \rho_{0}\widehat{f}\right)  ^{\vee}+\left(  \left(  1-\rho
_{0}\right)  \widehat{f}\right)  ^{\vee}.
\]

The first term on the right is the inverse Fourier transform of a tempered
distribution with compact support and so by Appendix \ref{SbSect_property_E'},
it is a $C_{BP}^{\infty}$ function. Thus $f\in C_{BP}^{\infty}+C_{B}^{\left(
\left\lfloor \kappa\right\rfloor \right)  }\subset C_{BP}^{\left(
\left\lfloor \kappa\right\rfloor \right)  }$ and, from Definition
\ref{Def_Ca_smth_space}, $f\in C_{BP}^{\infty}+C_{B}^{\left(  \left\lfloor
\underline{\kappa}\right\rfloor \right)  }\subset C_{BP}^{\left(  \left\lfloor
\underline{\kappa}\right\rfloor \right)  }$.\medskip

\fbox{\textbf{Case 3} $w$ has properties W2.1 and W3.3} Since $w\in W2.1$ by
part 4 of Theorem \ref{Thm_weight_property_relat}, $w\in W2$. Now suppose
$\rho_{0}\in C_{0}^{\infty}$, $0\leq\rho_{0}\leq1$, and $\rho_{0}=1$ on the
ball $B\left(  0;1\right)  $. By part 3 of Theorem
\ref{Thm_Xwm_is_L1loc(Rd_0)}, $f\in X_{w}^{\theta}$ implies $\xi^{\beta}%
f_{F}\in L^{1}$ when $\beta\leq\kappa$ and so $\xi^{\beta}\left(  1-\rho
_{0}\right)  \widehat{f}=\xi^{\beta}\left(  1-\rho_{0}\right)  f_{F}\in L^{1}$
when $\beta\leq\left\lfloor \kappa\right\rfloor $. Hence $D^{\beta}\left(
\left(  1-\rho_{0}\right)  \widehat{f}\right)  ^{\vee}\in C_{B}^{\left(
0\right)  }$ for $\beta\leq\left\lfloor \kappa\right\rfloor $ and so $\left(
\left(  1-\rho_{0}\right)  \widehat{f}\right)  ^{\vee}\in C_{B}^{\left(
\left\lfloor \kappa\right\rfloor \right)  }$.

We now write $\widehat{f}=\rho_{0}\widehat{f}+\left(  1-\rho_{0}\right)
\widehat{f}$. Taking the inverse Fourier transform of this equation yields
\[
f=\left(  \rho_{0}\widehat{f}\right)  ^{\vee}+\left(  \left(  1-\rho
_{0}\right)  \widehat{f}\right)  ^{\vee}.
\]

The first term on the right is the inverse Fourier transform of a tempered
distribution with compact support and so by Appendix \ref{SbSect_property_E'},
it is a $C_{BP}^{\infty}$ function. Thus $f\in C_{BP}^{\infty}+C_{B}^{\left(
\left\lfloor \kappa\right\rfloor \right)  }\subset C_{BP}^{\left(
\left\lfloor \kappa\right\rfloor \right)  }$ and, from Definition
\ref{Def_Ca_smth_space}, $f\in C_{BP}^{\infty}+C_{B}^{\left(  \left\lfloor
\underline{\kappa}\right\rfloor \right)  }\subset C_{BP}^{\left(  \left\lfloor
\underline{\kappa}\right\rfloor \right)  }$.
\end{proof}

\begin{remark}
We could also choose $\rho_{0}\in S_{1;\theta}$ and use part 2 of Definition
\ref{Def_Convolution}.
\end{remark}

The next corollary shows that in the case of the \textbf{extended B-splines},
by \textbf{enlarging the space of weight functions from W3.1* to W3.1} we can
prove that the \textbf{data space} $X_{w}^{\theta}$\textbf{\ has more
differentiability for dimension }$d>1$\textbf{\ and the order }$\theta\geq d $.

\begin{corollary}
\label{Cor_Thm_Xwth_W3_smooth}Suppose $w$ is the \textbf{extended B-spline}
weight function \ref{p32} with parameters $l$ and $n$. Then:

\begin{enumerate}
\item If $w$ has property W3.1 for order $\theta$ and $\kappa=\kappa
_{1}\mathbf{1}$ it follows that $X_{w}^{\theta}\subset C_{BP}^{\left(  \left(
n-1+\left\lfloor \theta/d\right\rfloor \right)  \mathbf{1}\right)  }\subset
C_{BP}^{\left(  n-1+\left\lfloor \theta/d\right\rfloor \right)  }$.

\item If $w$ has property W3.1* for order $\theta$ and $\kappa=\kappa
_{1}\mathbf{1}$ it follows that $X_{w}^{\theta}\subset C_{BP}^{\left(
n-1+\theta\right)  }$ when $d=1$ and $X_{w}^{\theta}\subset C_{BP}^{\left(
\left(  n-1\right)  \mathbf{1}\right)  }\subset C_{BP}^{\left(  n-1\right)  }$
when $d>1$.
\end{enumerate}
\end{corollary}

\begin{proof}
\textbf{Part 1} From \ref{p21} of Theorem \ref{Thm_ExtNatSplin_wt_fn_W3.1},
$\max\left\lfloor \kappa_{1}\right\rfloor =n-1+\left\lfloor \theta
/d\right\rfloor $ so that Theorem \ref{Thm_Xwth_W3_smooth} implies
$X_{w}^{\theta}\subset C_{BP}^{\left(  \left(  n-1+\left\lfloor \theta
/d\right\rfloor \right)  \mathbf{1}\right)  }\subset C_{BP}^{\left(
n-1+\left\lfloor \theta/d\right\rfloor \right)  }$.\smallskip

\textbf{Part 2} Apply Theorem \ref{Thm_Xwth_W3_smooth} to \ref{p15} and
\ref{p05}.
\end{proof}

\section{The function $\frac{1}{w\left\vert \cdot\right\vert ^{2\theta}}$}

In this section we will prove some properties of the function $\frac
{1}{w\left\vert \cdot\right\vert ^{2\theta}}$, where $w$ is a weight function
and $\theta$ is a positive integer, which prepare for the definition of the
basis distribution and basis function in the next section. The results in this
section study the functional $\int\frac{\phi}{w\left\vert \cdot\right\vert
^{2\theta}}$, $\phi\in S_{\emptyset,2\theta}$, and show how it can be extended
to a member of the tempered distributions $S^{\prime}$. Note that the function
$\frac{1}{w\left\vert \cdot\right\vert ^{2\theta}}$ was introduced in
Subsection \ref{SbSect_I_J} where it was used to define the operator
$\mathcal{J}:L^{2}\rightarrow X_{w}^{\theta}$.

We now prove some properties of the function $\frac{1}{w\left\vert
\cdot\right\vert ^{2\theta}}$ where $w$ is a weight function with property W2.
The next result shows that $\frac{1}{w\left\vert \cdot\right\vert ^{2\theta}%
}\in S_{\emptyset,2\theta}^{\prime}$ and will allow basis distributions to be
defined in the next section. Compare this result with a property of $\frac
{1}{w\left\vert \cdot\right\vert ^{2\theta}}$that was proved in Lemma
\ref{Lem_functnal_phi_sq}.

\begin{theorem}
\label{Thm_1/w|.|^2u_in_So,2u'}Suppose the weight function $w$ also has
property W2. Then the functional $\int\frac{\phi}{w\left\vert \cdot\right\vert
^{2\theta}}$ defined for $\phi\in S_{\emptyset,2\theta}$ is a member of
$S_{\emptyset,2\theta}^{\prime}$. In fact there exists a constant
$c_{r_{2},\theta}$, independent of $\phi\in S_{\emptyset,2\theta}$, such that%
\begin{equation}
\int\frac{\left\vert \phi\right\vert }{w\left\vert \cdot\right\vert ^{2\theta
}}\leq\left(  c_{r_{2},\theta}\right)  ^{2}\left(  \int\frac{1}{w\left\vert
\cdot\right\vert ^{2\lambda\left(  \cdot\right)  }}\right)  \sum
\limits_{\left\vert \alpha\right\vert \leq2n}\left\Vert \left(  1+\left\vert
\cdot\right\vert \right)  ^{2n}D^{\alpha}\phi\right\Vert _{\infty},\label{p11}%
\end{equation}

where $n=\operatorname{ceil}\left\{  \theta,\sigma\right\}  $. Here $\lambda$
is the function \ref{p70} introduced in the definition of weight function
property W2, and $c_{r_{2},\theta}=\max\left\{  1,r_{2}^{-\theta}\right\}  $
is the constant in the estimate of Lemma \ref{Lem_functnal_phi_sq}.
\end{theorem}

\begin{proof}
Suppose $r_{2}$ is the parameter in the definition of weight function property
W2. Then for $\phi\in S_{\emptyset,2\theta}$ we write%
\begin{align}
\int\frac{\left\vert \phi\right\vert }{w\left\vert \cdot\right\vert ^{2\theta
}}  & =\int\limits_{\left\vert \cdot\right\vert \leq r_{2}}\frac{\left\vert
\phi\right\vert }{w\left\vert \cdot\right\vert ^{2\theta}}+\int%
\limits_{\left\vert \cdot\right\vert \geq r_{2}}\frac{\left\vert
\phi\right\vert }{w\left\vert \cdot\right\vert ^{2\theta}}\nonumber\\
& =\int\limits_{\left\vert \cdot\right\vert \leq r_{2}}\frac{\left\vert
\phi\right\vert }{\left\vert \cdot\right\vert ^{2\theta}}\frac{1}{w}%
+\int\limits_{\left\vert \cdot\right\vert \geq r_{2}}\frac{\left\vert
\cdot\right\vert ^{2\sigma}\left\vert \phi\right\vert }{\left\vert
\cdot\right\vert ^{2\theta}}\frac{1}{w\left\vert \cdot\right\vert ^{2\sigma}%
}\nonumber\\
& \leq\left\Vert \frac{\phi}{\left\vert \cdot\right\vert ^{2\theta}%
}\right\Vert _{\infty;\leq r_{2}}\int\limits_{\left\vert \cdot\right\vert \leq
r_{2}}\frac{1}{w}+\left\Vert \frac{\left\vert \cdot\right\vert ^{2\sigma}\phi
}{\left\vert \cdot\right\vert ^{2\theta}}\right\Vert _{\infty;\geq r_{2}}%
\int\limits_{\left\vert \cdot\right\vert \geq r_{2}}\frac{1}{w\left\vert
\cdot\right\vert ^{2\sigma}}\nonumber\\
& \leq\max\left\{  \left\Vert \frac{\phi}{\left\vert \cdot\right\vert
^{2\theta}}\right\Vert _{\infty;\leq r_{2}},\left\Vert \frac{\left\vert
\cdot\right\vert ^{2\sigma}\phi}{\left\vert \cdot\right\vert ^{2\theta}%
}\right\Vert _{\infty;\geq r_{2}}\right\}  \int\frac{1}{w\left\vert
\cdot\right\vert ^{2\lambda\left(  \cdot\right)  }}\label{p53}%
\end{align}

where the integrals exist by weight function property W2. Since
$n=\operatorname*{ceil}\left\{  \theta,\sigma\right\}  $ we can apply
inequality \ref{p33} to get
\[
\left\Vert \frac{\phi}{\left\vert \cdot\right\vert ^{2\theta}}\right\Vert
_{\infty;\leq r_{2}}\leq\left\Vert \frac{\phi}{\left\vert \cdot\right\vert
^{2\theta}}\right\Vert _{\infty}\leq\sum\limits_{\left\vert \alpha\right\vert
=2\theta}\left\Vert D^{\alpha}\phi\right\Vert _{\infty}\leq\sum
\limits_{\left\vert \alpha\right\vert =2\theta}\left\Vert \left(  1+\left\vert
\cdot\right\vert \right)  ^{2n}D^{\alpha}\phi\right\Vert _{\infty}\leq
\sum\limits_{\left\vert \alpha\right\vert \leq2n}\left\Vert \left(
1+\left\vert \cdot\right\vert \right)  ^{2n}D^{\alpha}\phi\right\Vert
_{\infty},
\]

and, since $\phi\in S$,
\[
\left\Vert \frac{\left\vert \cdot\right\vert ^{2\sigma}\phi}{\left\vert
\cdot\right\vert ^{2\theta}}\right\Vert _{\infty;\geq r_{2}}\leq
r_{2}^{-2\theta}\left\Vert \left\vert \cdot\right\vert ^{2\sigma}%
\phi\right\Vert _{\infty;\geq r_{2}}\leq r_{2}^{-2\theta}\left\Vert \left(
1+\left\vert \cdot\right\vert \right)  ^{2n}\phi\right\Vert _{\infty}\leq
r_{2}^{-2\theta}\sum\limits_{\left\vert \alpha\right\vert \leq2n}\left\Vert
\left(  1+\left\vert \cdot\right\vert \right)  ^{2n}D^{\alpha}\phi\right\Vert
_{\infty}.
\]

Substituting these inequalities into the right side of \ref{p53} gives the
estimate \ref{p11} of this lemma. Finally, part 1 of Theorem
\ref{Thm_prop_functnl_on_Son} implies that the functional $\int\frac{\phi
}{w\left\vert \cdot\right\vert ^{2\theta}}$ defined on $S_{\emptyset,2\theta}$
is a member of $S_{\emptyset,2\theta}^{\prime}$.
\end{proof}

?? \textbf{NOTE}: Theorem \ref{Thm_Tay_rem_zeros} is \textbf{much} more
suitable for radial functions.?? \textbf{BUT} in general $\phi\in
S_{\emptyset,2\theta}$ and is not radial.

The last theorem now allows us to extend $\int\frac{\phi}{w\left\vert
\cdot\right\vert ^{2\theta}}$ to $S$ as a member of $S^{\prime}$.

\begin{corollary}
\label{Cor_Thm_1/w|.|^2u_in_So,2u'}Suppose the weight function $w$ has
property W2. Then the functional $\int\frac{\phi}{w\left\vert \cdot\right\vert
^{2\theta}}$ of Theorem \ref{Thm_1/w|.|^2u_in_So,2u'} is a member of
$S_{\emptyset,2\theta}^{\prime}$ and can be extended from $S_{\emptyset
,2\theta}$ to $S$ as a member of $S^{\prime}$, say $\chi^{e}$. The set of
extensions is $\chi^{e}+\widehat{P}_{2\theta}$.
\end{corollary}

\begin{proof}
A direct consequence of Theorem \ref{Thm_1/w|.|^2u_in_So,2u'} and Theorem
\ref{Thm_prop_functnl_on_Son}.
\end{proof}

\section{Basis distributions and basis functions\label{Sect_basis_distrib}}

It is now time to define the (in general complex-valued) basis distributions
that are generated by a weight function with weight property W2. However, it
will be shown that if the weight function also has property W3.1 or properties
W2.1 and W3.2 or W3.3 then the basis distribution is a continuous function
that will be used in the later chapters to construct the basis function
interpolants and smoothers studied in those documents. Following Light and
Wayne \cite{LightWayneX98Weight} we will define basis distributions
\textit{directly} i.e. without reference to the variational problems which
define the basis function interpolants and smoothers to be studied in later documents.

We will also calculate the basis functions which correspond to the (tensor
product) extended B-spline weight functions.

\subsection{Definition of a basis distribution}

We first define a tempered basis distribution of positive order by only
assuming weight function properties W1 and W2. Later property W3 will be
applied and this ensures the basis distributions are continuous basis functions.

\begin{definition}
\label{Def_basis_distrib}\textbf{Basis distributions and basis functions}

Suppose the weight function $w$ also has property W2. Then by Corollary
\ref{Cor_Thm_1/w|.|^2u_in_So,2u'}, $\chi=\dfrac{1}{w\left\vert \cdot
\right\vert ^{2\theta}}\in S_{\emptyset,2\theta}^{\prime}$ for each positive
integer $\theta$, and can be extended to $S$ non-uniquely as a member of
$S^{\prime}$ which we denote by $\chi^{e}$.

If $G\in S^{\prime}$ and satisfies $\widehat{G}=\chi^{e}$, we say $G$ is a
\textbf{(tempered) basis distribution} of order $\theta$ generated by $w$.

Thus a basis distribution of order $\theta$ is any tempered distribution $G$
which satisfies
\begin{equation}
\left[  \widehat{G},\phi\right]  =\int\frac{\phi}{w\left\vert \cdot\right\vert
^{2\theta}}\text{ }for\text{ }all\text{ }\phi\in S_{\emptyset,2\theta
}.\label{p01}%
\end{equation}

For the purposes of this series of documents we will call a basis distribution
a \textbf{basis function} if it is continuous. We \textbf{note} that usually
we would call a tempered basis distribution a basis function if it were a
regular tempered distribution in the sense of Appendix
\ref{SbSect_property_S'}.
\end{definition}

From their definition basis distributions of a given order are not unique. In fact:

\begin{theorem}
\label{Thm_basis_fn_set}Suppose the weight function $w$ also has property W2.
Suppose $G$ is a basis distribution of order $\theta\geq1$ generated by $w$.

Then the set of basis distributions of order $\theta$ is $G+P_{2\theta-1}$.
\end{theorem}

\begin{proof}
This is a direct consequence of the definition of a basis distribution and
Corollary \ref{Cor_Thm_1/w|.|^2u_in_So,2u'}.
\end{proof}

The next theorem will require the following lemma which we give without proof.

\begin{lemma}
\label{Lem_lin_comb_exp(i)_equals_0}Suppose $\left\{  x^{\left(  k\right)
}\right\}  _{k=1}^{n}$ is a set of distinct points in $\mathbb{R}^{d}$ and
$\left\{  \lambda_{k}\right\}  _{k=1}^{n}$ is a set of complex numbers. Then
the function $f\left(  \xi\right)  =\sum\limits_{j=1}^{n}\lambda
_{k}e^{-ix^{\left(  k\right)  }\xi}$ has the properties:

\begin{enumerate}
\item $f\left(  \xi\right)  =0$ a.e. implies $\lambda_{k}=0$ for all $k$.

\item The null space of $f$ is a closed set of measure zero.
\end{enumerate}
\end{lemma}

\begin{theorem}
\label{Thm_G(*-Xk)_lin_indep}Suppose the points $\left\{  x^{\left(  k\right)
}\right\}  _{k=1}^{N}$ are distinct. Then the translated basis distributions

$\left\{  G\left(  \cdot-x^{\left(  k\right)  }\right)  \right\}  _{k=1}^{N}$
are linearly independent w.r.t. the complex scalars.
\end{theorem}

\begin{proof}
Suppose $\sum\limits_{k=1}^{n}\lambda_{k}G\left(  \cdot-x^{\left(  k\right)
}\right)  =0$ and not all $\lambda_{k}\neq0$. Then the definition of $G$
implies $\widehat{G}=\frac{1}{w\left(  \xi\right)  \left\vert \xi\right\vert
^{2\theta}}$ on $\mathbb{R}^{d}\setminus0$ and so on taking the Fourier
transform%
\[
0=\sum\limits_{k=1}^{n}\lambda_{k}e^{-ix^{\left(  k\right)  }\xi}\widehat
{G}=\left(  \sum\limits_{k=1}^{n}\lambda_{k}e^{-ix^{\left(  k\right)  }\xi
}\right)  \widehat{G}=\frac{\sum\nolimits_{k=1}^{n}\lambda_{k}e^{-ix^{\left(
k\right)  }\xi}}{w\left(  \xi\right)  \left\vert \xi\right\vert ^{2\theta}},
\]

a.e. on $\mathbb{R}^{d}\setminus0$. Hence, since $w\left(  \xi\right)  >0$
a.e., it follows that $\sum\limits_{k=1}^{n}\lambda_{k}e^{-i\xi x^{\left(
k\right)  }}=0$ a.e. and on applying Lemma \ref{Lem_lin_comb_exp(i)_equals_0}
we conclude that $\lambda_{k}=0$ for all $k$.
\end{proof}

\subsection{The smoothness of basis distributions; continuous basis functions
\label{SbSect_basis_smth}}

In this subsection we look at the smoothness and growth of basis distributions
when the weight function has either property W3.1 or properties W2.1 and W3.2
of Definition \ref{Def_extend_wt_fn}. In the case of property W3.1 we will
also derive a simple inverse Fourier transform formula for the basis function
and its (bounded) derivatives which only uses $L^{1}$ Fourier transform
theory. When the weight function has property W3.2 much more effort is
required to obtain a `modified' inverse Fourier transform formula and in
Chapter \ref{Ch_MoreBasisTheory} will deal with this.

However, the next theorem will deal with the W3.1 case and it will require a
lemma, a standard $L^{1}$ inverse Fourier transform result due to Laurent
Schwartz. See, for example, Theorem 4.2, p.150 of Malliavin \cite{Malliavin95}.

\begin{lemma}
\label{Lem_L1_Fourier_contin}If $f\in S^{\prime}$ and $\widehat{f}\in L^{1}$,
then $f\in C_{B}^{\left(  0\right)  }$ and
\[
f\left(  x\right)  =\left(  2\pi\right)  ^{-d/2}\int e^{ix\xi}\widehat
{f}\left(  \xi\right)  d\xi.
\]

\end{lemma}

\begin{theorem}
\label{Thm_basis_smth_W3.1}Suppose the weight function $w$ has property W3.1
for order $\theta$ and smoothness parameter $\kappa$. Then $\frac
{1}{w\left\vert \cdot\right\vert ^{2\theta}}\in L^{1}$, $G=\left(  \frac
{1}{w\left\vert \cdot\right\vert ^{2\theta}}\right)  ^{\vee}$ is a basis
function of order $\theta$ generated by $w$, and $G\in C_{B}^{\left(
\left\lfloor 2\kappa\right\rfloor \right)  }$. Further, the set of basis
functions is $G+P_{2\theta-1}$ and%
\begin{equation}
D^{\gamma}G\left(  x\right)  =\left(  2\pi\right)  ^{-\frac{d}{2}}\int
e^{ix\xi}\frac{\left(  i\xi\right)  ^{\gamma}}{w\left(  \xi\right)  \left\vert
\xi\right\vert ^{2\theta}}d\xi,\quad\gamma\leq\left\lfloor 2\kappa
\right\rfloor .\label{p41}%
\end{equation}

\end{theorem}

\begin{proof}
By part 2 of Theorem \ref{Thm_weight_property_relat}, $w$ has property W2 and
so the basis distributions are defined. By Theorem \ref{Thm_property_wt_fn_W3}%
, $\frac{1}{w\left\vert \cdot\right\vert ^{2\theta}}\in L^{1}\subset
S^{\prime}$. Hence $\left(  \frac{1}{w\left\vert \cdot\right\vert ^{2\theta}%
}\right)  ^{\vee}\in C_{B}^{\left(  0\right)  }\subset S^{\prime}$. Set
$G=\left(  \frac{1}{w\left\vert \cdot\right\vert ^{2\theta}}\right)  ^{\vee}$.
Then $G\in S^{\prime}$ and $\widehat{G}=\frac{1}{w\left\vert \cdot\right\vert
^{2\theta}}$ a.e., implying $G\in C_{B}^{\left(  0\right)  }$ and so $G$
satisfies%
\[
\left[  \widehat{G},\phi\right]  =\int\frac{\phi}{w\left\vert \cdot\right\vert
^{2\theta}},\quad\phi\in S_{\emptyset,2\theta},
\]

and by Definition \ref{Def_basis_distrib} $G$ is a basis function. By Theorem
\ref{Thm_basis_fn_set} the set of basis functions is $G+P_{2\theta-1}$.

Since $w$ has property W3.1, part 1 of Theorem \ref{Thm_property_wt_fn_W3}
implies that for order $\theta$ and $\kappa$, $\frac{x^{2\gamma}}{w\left(
x\right)  \left\vert x\right\vert ^{2\theta}}\in L^{1}$ when $0\leq\gamma
\leq\kappa$. Hence%
\[
\widehat{D^{\gamma}G}\left(  \xi\right)  =\left(  i\xi\right)  ^{\gamma
}\widehat{G}\left(  \xi\right)  =\frac{\left(  i\xi\right)  ^{\gamma}%
}{w\left(  \xi\right)  \left\vert \xi\right\vert ^{2\theta}}\in L^{1}%
,\quad\gamma\leq\left\lfloor 2\kappa\right\rfloor ,
\]

and applying Lemma \ref{Lem_L1_Fourier_contin} with $f=D^{\gamma}G\in
S^{\prime}$ we obtain \ref{p41} and $D^{\gamma}G\in C_{B}^{\left(  0\right)
}$ for $\gamma\leq\left\lfloor 2\kappa\right\rfloor $ i.e. $G\in
C_{B}^{\left(  \left\lfloor 2\kappa\right\rfloor \right)  }$.
\end{proof}

Thus when the weight function has property W3.1 for order $\theta$ and
$\kappa$, the basis distribution of order $\theta$ must have differentiability
of at least $\left\lfloor 2\kappa\right\rfloor $ and each derivative is
bounded. This is at least twice the minimum smoothness of the functions in
$X_{w}^{\theta}$, which have differentiability of at least $\left\lfloor
\kappa\right\rfloor $. Continuous basis functions will be used to construct
the basis function interpolants and smoothers which are discussed in Chapters
\ref{Ch_Interpol}, \ref{Ch_ExactSmth} and \ref{Ch_Approx_smth}.

We will now consider the question of smoothness when a weight function has
property W3.2 and then property W3.3. The next theorem shows that for property
W3.2 a basis distribution of order $\theta$ is a $C_{BP}^{\left(  \left\lfloor
2\kappa\right\rfloor \right)  }$ basis function. This is at least twice the
minimum smoothness of the functions in $X_{w}^{\theta}$, which were shown to
have differentiability of at least $\left\lfloor \kappa\right\rfloor $.

\begin{theorem}
\label{Thm_basis_smth_W3.2_r3_pos}\textbf{Property W3.2 and basis functions}

\begin{enumerate}
\item If a weight function $w$ has property W2 then $1/w\in S^{\prime}\cap
L_{loc}^{1}$.\medskip

Now suppose the weight function $w$ has properties W2.1 and \textbf{W3.2} for
order $\theta$ and $\kappa$. Then:\smallskip

\item Suppose $G$ is a basis distribution of order $\theta$ generated by $w$.
Then $\left\vert \cdot\right\vert ^{2\theta}\widehat{G}=\frac{1}{w}$ as
tempered distributions.

\item Any $f\in S^{\prime}$ such that $\left\vert \cdot\right\vert ^{2\theta
}\widehat{f}=\frac{1}{w}$ satisfies $f\in C_{BP}^{\left(  \left\lfloor
2\kappa\right\rfloor \right)  }$.

\item $G\in C_{BP}^{\left(  \left\lfloor 2\kappa\right\rfloor \right)  }$ and
we say that $G$ is a basis function of order $\theta$.

\item The set of all basis functions is $G+P_{2\theta-1}$.
\end{enumerate}
\end{theorem}

\begin{proof}
\textbf{Part 1} Suppose weight function $w$ has property W2. Noting that
property W2 is defined by the conditions
\[
W2.1:1/w\in L_{loc}^{1},\quad W2.2:\int_{\left\vert \cdot\right\vert \geq
r_{2}}\frac{1}{w\left\vert \cdot\right\vert ^{2\sigma}}<\infty\text{
}for\text{ }some\text{ }r_{2}>0,
\]

we calculate that%
\[
\left\vert \int\frac{\phi}{w}\right\vert \leq2\max\left\{  \int_{\left\vert
\cdot\right\vert \leq R_{2}}\frac{1}{w},\int_{\left\vert \cdot\right\vert \leq
R_{2}}\frac{1}{w\left\vert \cdot\right\vert ^{2\sigma}}\right\}  \left\Vert
\left(  1+\left\vert \cdot\right\vert \right)  ^{2\left\lceil \sigma
\right\rceil }\phi\right\Vert _{\infty},\quad\phi\in S,
\]

where $\left\lceil \sigma\right\rceil $ denotes the \textit{ceiling} of
$\sigma$. Hence $1/w\in S^{\prime}$, since $\left\Vert \left(  1+\left\vert
\cdot\right\vert \right)  ^{2\left\lceil \sigma\right\rceil }\phi\right\Vert
_{\infty}$ is one of the seminorms used in Definition \ref{Def_Distributions}
to specify the topology of $S^{\prime}$.\medskip

\textbf{Part 2} \ By part 4 of Theorem \ref{Thm_weight_property_relat}, $w$
has property W2 and so the basis distribution is defined and by part 1,
$1/w\in S^{\prime}$. If $\phi\in S$ then by Theorem
\ref{Thm_product_of_Co,k_funcs}, $\left\vert \cdot\right\vert ^{2\theta}%
\phi\in S_{\emptyset,2\theta}$ and so by \ref{p01} of the basis distribution
definition
\[
\left[  \left\vert \cdot\right\vert ^{2\theta}\widehat{G},\phi\right]
=\left[  \widehat{G},\left\vert \cdot\right\vert ^{2\theta}\phi\right]
=\left[  \frac{1}{\left\vert \cdot\right\vert ^{2\theta}w},\left\vert
\cdot\right\vert ^{2\theta}\phi\right]  =\int\frac{1}{\left\vert
\cdot\right\vert ^{2\theta}w}\left\vert \cdot\right\vert ^{2\theta}\phi
=\int\frac{1}{w}\phi=\left[  \frac{1}{w},\phi\right]  .
\]
\medskip

\textbf{Part 3} In part 2 it was shown that $1/w\in S^{\prime}$ and that there
exists $f\in S^{\prime}$ such that $\left\vert \cdot\right\vert ^{2\theta
}\widehat{f}=\frac{1}{w}$. The rest of this part is based on the fact that
$g\in L^{1}$ implies $\widehat{g}\in C_{B}^{\left(  0\right)  }$. Since
$r_{3}>0$ there exists a function $\psi\in C_{0}^{\infty}$ such that
$0\leq\psi\leq1$, $\psi=1$ on $B\left(  0;r_{3}\right)  $ and $\psi=0$ outside
$B\left(  0;2r_{3}\right)  $. Hence $\dfrac{1-\psi}{\left\vert \cdot
\right\vert ^{2\theta}}\in L_{loc}^{1}\cap C_{B}^{\infty}$ and since
$\left\vert \cdot\right\vert ^{2\theta}\widehat{f}=\frac{1}{w}\in S^{\prime}$
we have
\begin{equation}
\dfrac{1-\psi}{\left\vert \cdot\right\vert ^{2\theta}}\left\vert
\cdot\right\vert ^{2\theta}\widehat{f}=\left(  1-\psi\right)  \widehat
{f}=\frac{1-\psi}{w\left\vert \cdot\right\vert ^{2\theta}}\in S^{\prime}\cap
L_{loc}^{1}.\label{p180}%
\end{equation}

The next step is to prove that $x^{\alpha}\left(  1-\psi\right)  \widehat
{f}\in L^{1}$ when $\left\vert \alpha\right\vert \leq\left\lfloor
2\kappa\right\rfloor $. But
\begin{equation}
x^{\alpha}\left(  1-\psi\right)  \widehat{f}=x^{\alpha}\dfrac{1-\psi
}{w\left\vert \cdot\right\vert ^{2\theta}},\label{p182}%
\end{equation}

and if $\left\vert \alpha\right\vert \leq\left\lfloor 2\kappa\right\rfloor $,%
\[
\int\left\vert x^{\alpha}\frac{1-\psi}{w\left\vert \cdot\right\vert ^{2\theta
}}\right\vert dx\leq\int\left\vert \cdot\right\vert ^{\left\vert
\alpha\right\vert }\frac{1-\psi}{w\left\vert \cdot\right\vert ^{2\theta}}%
=\int\limits_{\left\vert \cdot\right\vert \geq r_{3}}\left\vert \cdot
\right\vert ^{\left\vert \alpha\right\vert }\frac{1-\psi}{w\left\vert
\cdot\right\vert ^{2\theta}}\leq\int\limits_{\left\vert \cdot\right\vert \geq
r_{3}}\frac{\left\vert \cdot\right\vert ^{\left\vert \alpha\right\vert }%
}{w\left\vert \cdot\right\vert ^{2\theta}},
\]

and
\[
\int\limits_{\left\vert \cdot\right\vert \geq r_{3}}\frac{\left\vert
\cdot\right\vert ^{\left\vert \alpha\right\vert }}{w\left\vert \cdot
\right\vert ^{2\theta}}=\int\limits_{\left\vert \cdot\right\vert \geq r_{3}%
}\frac{\left\vert \cdot\right\vert ^{2\kappa}}{\left\vert \cdot\right\vert
^{2\kappa-\left\vert \alpha\right\vert }w\left\vert \cdot\right\vert
^{2\theta}}\leq\frac{1}{\left(  r_{3}\right)  ^{2\kappa-\left\vert
\alpha\right\vert }}\int\limits_{\left\vert \cdot\right\vert \geq r_{3}}%
\frac{\left\vert \cdot\right\vert ^{2\kappa}}{w\left\vert \cdot\right\vert
^{2\theta}}<\infty,
\]

by property W3.2. Hence $x^{\alpha}\dfrac{1-\psi}{w\left\vert \cdot\right\vert
^{2\theta}}\in L^{1}$ when $\left\vert \alpha\right\vert \leq\left\lfloor
2\kappa\right\rfloor $, and so $D^{\alpha}\left(  \dfrac{1-\psi}{w\left\vert
\cdot\right\vert ^{2\theta}}\right)  ^{\vee}$

$\in C_{B}^{\left(  0\right)  }$ when $\left\vert \alpha\right\vert
\leq\left\lfloor 2\kappa\right\rfloor $ i.e.
\begin{equation}
\left(  \dfrac{1-\psi}{w\left\vert \cdot\right\vert ^{2\theta}}\right)
^{\vee}\in C_{B}^{\left(  \left\lfloor 2\kappa\right\rfloor \right)
}.\label{p121}%
\end{equation}

Also, from equation \ref{p180}
\[
\left(  \dfrac{1-\psi}{w\left\vert \cdot\right\vert ^{2\theta}}\right)
^{\vee}=\left(  \left(  1-\psi\right)  \widehat{f}\right)  ^{\vee}=f-\left(
\psi\widehat{f}\right)  ^{\vee},
\]

or on rearranging
\begin{equation}
f=\left(  \dfrac{1-\psi}{w\left\vert \cdot\right\vert ^{2\theta}}\right)
^{\vee}+\left(  \psi\widehat{f}\right)  ^{\vee}.\label{p181}%
\end{equation}

We already know from \ref{p121} that the first term on the right is a
$C_{B}^{\left(  \left\lfloor 2\kappa\right\rfloor \right)  }$ function and
since $\psi$ has bounded support $\psi\widehat{f}$ is a distribution with
compact support and by Appendix \ref{SbSect_property_E'} the inverse Fourier
transform is a $C_{BP}^{\infty}$ function. Thus $f\in C_{B}^{\left(
\left\lfloor 2\kappa\right\rfloor \right)  }+C_{BP}^{\infty}\subset
C_{BP}^{\left(  \left\lfloor 2\kappa\right\rfloor \right)  }$.\medskip

\textbf{Part 4} Follows directly from parts 2 and 3.\medskip

\textbf{Part 5} Follows directly from Theorem \ref{Thm_basis_fn_set}.
\end{proof}

This result is similar to that of Theorem \ref{Thm_basis_smth_W3.1} for weight
functions with property W3.1.

\begin{theorem}
\label{Thm_basis_smth_W3.3}\textbf{Property W3.3 and basis functions} Suppose
the weight function $w$ has property W2.1 and also property W3.3 for order
$\theta$ and smoothness parameter $\kappa$.

Then $\frac{1}{w\left\vert \cdot\right\vert ^{2\theta}}\in L^{1}$, $G=\left(
\frac{1}{w\left\vert \cdot\right\vert ^{2\theta}}\right)  ^{\vee}$ is a basis
function of order $\theta$ generated by $w$, and $G\in C_{B}^{\left(
\left\lfloor 2\kappa\right\rfloor \right)  }\subset C_{B}^{\left(
\left\lfloor 2\underline{\kappa}\right\rfloor \right)  }$. Further, the set of
basis functions is $G+P_{2\theta-1}$ and%
\[
D^{\gamma}G\left(  x\right)  =\left(  2\pi\right)  ^{-\frac{d}{2}}\int
e^{ix\xi}\frac{\left(  i\xi\right)  ^{\gamma}d\xi}{w\left(  \xi\right)
\left\vert \xi\right\vert ^{2\theta}},\quad\gamma\leq\left\lfloor
2\kappa\right\rfloor .
\]

\end{theorem}

\begin{proof}
By part 4 of Theorem \ref{Thm_weight_property_relat}, $w$ has property W2 and
so the basis distributions are defined. By Theorem \ref{Thm_property_wt_fn_W3}%
, $\frac{1}{w\left\vert \cdot\right\vert ^{2\theta}}\in L^{1}\subset
S^{\prime}$. Hence $\left(  \frac{1}{w\left\vert \cdot\right\vert ^{2\theta}%
}\right)  ^{\vee}\in C_{B}^{\left(  0\right)  }\subset S^{\prime}$. Set
$G=\left(  \frac{1}{w\left\vert \cdot\right\vert ^{2\theta}}\right)  ^{\vee}$.
Then $G\in S^{\prime}$ and $\widehat{G}=\frac{1}{w\left\vert \cdot\right\vert
^{2\theta}}$ a.e., implying $G\in C_{B}^{\left(  0\right)  }$ and so $G$
satisfies%
\[
\left[  \widehat{G},\phi\right]  =\int\frac{\phi}{w\left\vert \cdot\right\vert
^{2\theta}},\quad\phi\in S_{\emptyset,2\theta},
\]

and by Definition \ref{Def_basis_distrib} $G$ is a basis function. By Theorem
\ref{Thm_basis_fn_set} the set of basis functions is $G+P_{2\theta-1}$.

Since $w$ has property W3.1, part 1 of Theorem \ref{Thm_property_wt_fn_W3}
implies that for order $\theta$ and $\kappa$, $\frac{x^{2\gamma}}{w\left(
x\right)  \left\vert x\right\vert ^{2\theta}}\in L^{1}$ when $0\leq\gamma
\leq\kappa$. Hence%
\[
\widehat{D^{\gamma}G}\left(  \xi\right)  =\left(  i\xi\right)  ^{\gamma
}\widehat{G}\left(  \xi\right)  =\frac{\left(  i\xi\right)  ^{\gamma}%
}{w\left(  \xi\right)  \left\vert \xi\right\vert ^{2\theta}}\in L^{1}%
,\quad\gamma\leq\left\lfloor 2\kappa\right\rfloor ,
\]

and applying Lemma \ref{Lem_L1_Fourier_contin} with $f=D^{\gamma}G\in
S^{\prime}$ we obtain \ref{p41} and $D^{\gamma}G\in C_{B}^{\left(  0\right)
}$ for $\gamma\leq\left\lfloor 2\kappa\right\rfloor $ i.e. $G\in
C_{B}^{\left(  \left\lfloor 2\kappa\right\rfloor \right)  }$.
\end{proof}

We summarize our basis function smoothness results in:

\begin{summary}
\label{Sum_basis_func_smooth}Suppose a weight function $w$ has properties W2
and W3 for order $\theta$ and smoothness parameter $\kappa$. Suppose $G$ is a
basis distribution of order $\theta$ generated by $w$. Then:

\begin{enumerate}
\item If $w$ has property \textbf{W3.1} then $G\in C_{B}^{\left(  \left\lfloor
2\kappa\right\rfloor \right)  }$ (Theorem \ref{Thm_basis_smth_W3.1}).

\item If $w$ has properties W2.1 and \textbf{W3.2} then $G\in C_{BP}^{\left(
\left\lfloor 2\kappa\right\rfloor \right)  }$ (part 4 Theorem
\ref{Thm_basis_smth_W3.2_r3_pos}).

\item If $w$ has properties W2.1 and \textbf{W3.3} then $G\in C_{B}^{\left(
\left\lfloor 2\kappa\right\rfloor \right)  }$ (Theorem
\ref{Thm_basis_smth_W3.3}).
\end{enumerate}
\end{summary}

\begin{remark}
Note that if a weight function has property W2 then basis distributions are
defined for all positive orders. However, if the weight function also has
property W3 we have only proved that the basis distribution is
continuous\textbf{\ if the order of the basis distribution matches the order
of the weight function}.

Continuous basis distributions can be used to express the solution of the
interpolation and smoothing problems discussed later.
\end{remark}

\subsection{Summary diagram\label{SbSect_summary_diag}}

Figure \ref{Fig_prop_basis} illustrates the relationships between the weight
function properties and the basis distribution properties proved in the last section.

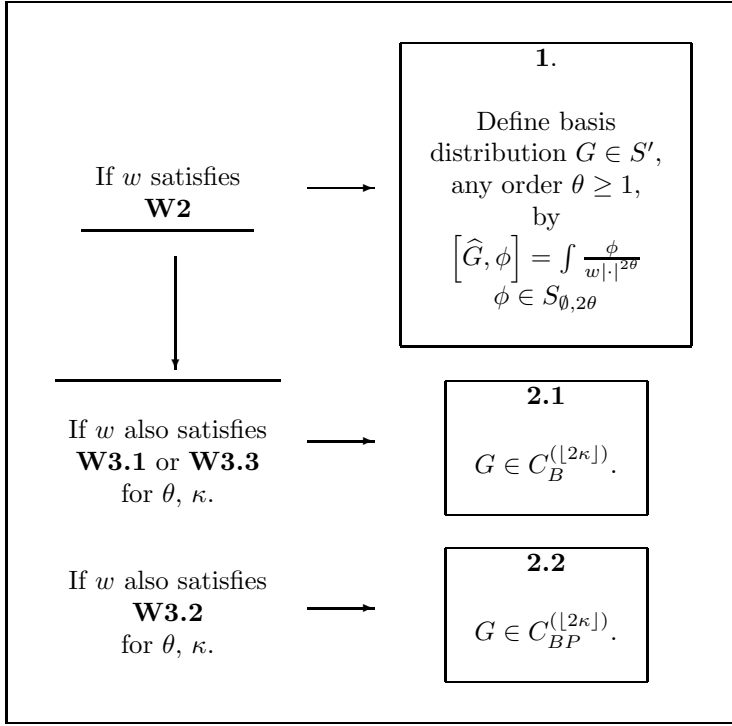
\fbox{$%
\begin{array}
[c]{ccccc}
&  &  &  & \\
& \underline{%
\begin{array}
[c]{c}%
\text{If }w\text{ satisfies}\\
\text{\textbf{W2}}%
\end{array}
} &
{\tt\setlength{\unitlength}{0.5pt}
\begin{picture}(50,20)
\thinlines\put(0,10){\vector(1,0){50}}
\end{picture}}%
&
\begin{tabular}
[c]{|l|}\hline
$%
\begin{array}
[c]{c}%
\text{\textbf{1}.}\\
\\
\text{Define basis}\\
\text{distribution }G\in S^{\prime}\text{,}\\
\text{any order }\theta\geq1\text{,}\\
\text{by}\\
\left[  \widehat{G},\phi\right]  =\int\frac{\phi}{w\left\vert \cdot\right\vert
^{2\theta}}\\
\phi\in S_{\emptyset,2\theta}\\
\quad
\end{array}
$\\\hline
\end{tabular}
& \\
& \text{%
{\tt\setlength{\unitlength}{0.5pt}
\begin{picture}(1,1)
\thinlines\put(1,85){\vector(0,-1){85}}
\end{picture}}%
} &  &  & \\
& \overline{%
\begin{array}
[c]{c}%
\\
\text{If }w\text{ also satisfies}\\
\text{\textbf{W3.1} or \textbf{W3.3}}\\
\text{for }\theta\text{, }\kappa\text{.}\\
\end{array}
} &
{\tt\setlength{\unitlength}{0.5pt}
\begin{picture}(50,20)
\thinlines\put(0,10){\vector(1,0){50}}
\end{picture}}%
&
\begin{tabular}
[c]{|l|}\hline
$%
\begin{array}
[c]{c}%
\text{\textbf{2.1}}\\
\\
G\in C_{B}^{\left(  \left\lfloor 2\kappa\right\rfloor \right)  }\text{.}\\
\text{\quad}%
\end{array}
$\\\hline
\end{tabular}
& \\
&  &  &  & \\
&
\begin{array}
[c]{c}%
\text{If }w\text{ also satisfies}\\
\text{\textbf{W3.2}}\\
\text{for }\theta\text{, }\kappa\text{.}\\
\end{array}
&
{\tt\setlength{\unitlength}{0.5pt}
\begin{picture}(50,20)
\thinlines\put(0,10){\vector(1,0){50}}
\end{picture}}%
&
\begin{tabular}
[c]{|l|}\hline
$%
\begin{array}
[c]{c}%
\text{\textbf{2.2}}\\
\\
G\in C_{BP}^{\left(  \left\lfloor 2\kappa\right\rfloor \right)  }\text{.}\\
\text{\quad}%
\end{array}
$\\\hline
\end{tabular}
& \\
&  &  &  &
\end{array}
$}\label{Fig_prop_basis}

\begin{description}
\item[Box 1] A tempered basis distribution $G$ of order $\theta\geq1$ is
defined for weight functions $w$ which satisfy property W2. When $\theta$ is
positive this coincides with Light and Wayne's definition.\smallskip

\item[Box 2.1] Suppose $w$ also has the property W3.1 or W3.3 for order
$\theta\geq1$ and smoothness parameter $\kappa\in\mathbb{R}^{d}$. Then by
Theorem \ref{Thm_basis_smth_W3.1} or Theorem \ref{Thm_basis_smth_W3.3} the
basis distribution $G$ of order $\theta$ is a basis function in $C_{B}%
^{\left(  \left\lfloor 2\mathbf{\kappa}\right\rfloor \right)  }$.

\item[Box 2.2] Suppose $w$ also has the property W3.2 for order $\theta\geq1$
and scalar smoothness parameter $\kappa\geq0$. Then by Theorem
\ref{Thm_basis_smth_W3.2_r3_pos} the basis distribution $G$ of order $\theta$
is a basis function in $C_{BP}^{\left(  \left\lfloor 2\kappa\right\rfloor
\right)  }$. There are no growth estimates near infinity but see Theorem
\ref{Thm_deriv_Grho} below.
\end{description}

\section{The spaces $S_{\emptyset,\alpha}$ and $C_{\emptyset,\alpha}^{\infty}%
$\label{Sect_So,a}}

We now define some key spaces and exhibit some of their properties.

\begin{definition}
\label{Def_So,a}\textbf{The spaces} $S_{\emptyset,\alpha}$ \textbf{and}
$C_{\emptyset,\alpha}^{\infty}$\textbf{\ where }$\alpha\geq\mathbf{0}%
$\textbf{\ is a multi-index}
\begin{equation}
S_{\emptyset,0}=S,\qquad S_{\emptyset,\alpha}=\left\{  \phi\in S:D^{\beta}%
\phi\left(  0\right)  =0,\text{\ }\beta.<\alpha\right\}  ,\label{p46}%
\end{equation}

and we endow $S_{\emptyset,\alpha}$ with the subspace topology induced by the
space $S$. Here $S$ is the Schwartz space of $C^{\infty}$ functions of rapid
decrease used as test functions for the tempered distributions. $S$ is endowed
with the countable seminorm topology described in Appendix
\ref{Def_Distributions}.%
\begin{equation}
C_{\emptyset,0}^{\infty}=C^{\infty},\qquad C_{\emptyset,\alpha}^{\infty
}=\left\{  \phi\in C^{\infty}:D^{\beta}\phi\left(  0\right)  =0,\text{\ }%
\beta.<\alpha\right\}  ,\label{p47}%
\end{equation}

so the space $C_{\emptyset,\alpha}^{\infty}$ retains the constraints of
$S_{\emptyset,\alpha}$ near the origin.
\end{definition}

We will need the following useful results:

\begin{theorem}
\label{Thm_product_of_Co,a_funcs}\ 

\begin{enumerate}
\item $S_{\emptyset,\alpha}=S\cap C_{\emptyset,\alpha}^{\infty}$.

\item For $\alpha,\beta\geq0$, $\phi\in C_{\emptyset,\alpha}^{\infty}$ and
$\psi\in C_{\emptyset,\beta}^{\infty}$ implies that $\phi\psi\in
C_{\emptyset,\alpha+\beta}^{\infty}$.

\item Regarding monomials: $x^{\alpha}\in C_{\emptyset,\alpha}^{\infty}$.

\item $S_{\emptyset,n\mathbf{1}}\subset S_{\emptyset,n}$ and $C_{\emptyset
,n\mathbf{1}}^{\infty}\subset C_{\emptyset,n}^{\infty}$.
\end{enumerate}
\end{theorem}

\subsection{The extended B-spline basis functions when $w\in W3.1$%
\label{SbSect_basis_Bsplin_W3.1}}

Suppose $w$ is an extended B-spline weight function \ref{p32} which has
property W3.1 for order $\theta$ and smoothness parameter $\kappa$. Then by
Theorem \ref{Thm_basis_smth_W3.1}, $\left(  \frac{1}{w\left\vert
\cdot\right\vert ^{2\theta}}\right)  ^{\vee}$ is a basis function of order
$\theta$. However the formula $\left(  \frac{1}{w\left\vert \cdot\right\vert
^{2\theta}}\right)  ^{\vee}$ does not in general lend itself to the
calculation of a `closed form' suitable for numeric calculations involving
basis functions e.g. basis function smoothers, so we will derive an
alternative convolution form for the basis functions in Theorem
\ref{Thm_ExtNatSplin_convol}\textbf{.} \textbf{In fact, these basis functions
turn out to be derivatives of the zero order extended B-spline basis function
convolved with a thin-plate spline basis function}. To show this we will need
three lemmas. The first lemma gives some properties of the zero order extended
B-splines proved in Subsection 1.4.4 of Williams \cite{WilliamsZeroOrdSmthV4}.

\begin{lemma}
\label{Lem_basis_ord0_Bsplin_W3}Suppose $w$ is the \textbf{extended B-spline
weight function} \ref{p32} and suppose $1\leq n\leq l$. Then $1/w\in L^{1}$
and the basis function $G_{0}=\left(  \frac{1}{w}\right)  ^{\vee}$ of
\textbf{order zero} generated by $w$ is the tensor product $G_{0}\left(
x\right)  =%
{\textstyle\prod\limits_{k=1}^{d}}
H\left(  x_{k}\right)  $ where
\[
H\left(  t\right)  =\left(  -1\right)  ^{l-n}\tfrac{\left(  2\pi\right)
^{l/2}}{2^{2\left(  l-n\right)  +1}}\left(  D^{2\left(  l-n\right)  }\left(
\ast\Lambda\right)  ^{l}\right)  \left(  \tfrac{t}{2}\right)  ,\text{\quad
}t\in\mathbb{R}^{1},
\]

and $\left(  \ast\Lambda\right)  ^{l}$ denotes the convolution of $l$
1-dimensional hat functions. Here $\operatorname*{supp}H\subset\left[
-l,l\right]  $.

Further, if $n<l$ we have
\[
D^{2\left(  l-n\right)  }\left(  \ast\Lambda\right)  ^{l}=\left(  -1\right)
^{l-n}\left(  \ast\Lambda\right)  ^{l-n}\ast\sum_{k=-\left(  l-n\right)
}^{l-n}\left(  -1\right)  ^{\left\vert k\right\vert }\tbinom{l-n}{\left\vert
k\right\vert }\Lambda\left(  \cdot-k\right)  ,
\]

$H\in C_{0}^{\left(  2n-2\right)  }\left(  \mathbb{R}^{1}\right)  $,
$D^{2n-1}H$ is a bounded, piecewise constant function and $D^{2n}H$ is the sum
of a finite number of translated delta functions.

Further $G_{0}\in C_{0}^{\left(  \left(  2n-2\right)  \mathbf{1}\right)
}\left(  \mathbb{R}^{d}\right)  $, the derivatives $\left\{  D^{\alpha}%
G_{0}:\left(  2n-2\right)  \mathbf{1}<\alpha\leq\left(  2n-1\right)
\mathbf{1}\right\}  $ are bounded with bounded support, and $D^{2n\mathbf{1}%
}G_{0}$ is the finite sum of translated delta functions.
\end{lemma}

\begin{lemma}
\label{Lem_ExtNatSplin_1/w_in_Co,2d(l-n)}If $w$ is the extended B-spline
weight function \ref{p32} and $1\leq n<l$ then $\frac{1}{w}\in C_{\emptyset
,2\left(  l-n\right)  \mathbf{1}}^{\infty}\cap C_{BP}^{\infty}$.
\end{lemma}

\begin{proof}
The basis function $G_{0}$ generated by $w$ is a distribution with bounded
support so by Appendix \ref{SbSect_property_E'}, $\widehat{G_{0}}=\frac{1}%
{w}\in C_{BP}^{\infty}$. Further, if $\mathbf{1}=\left(  1,\ldots,1\right)
\in\mathbb{R}^{d}$ we write%
\begin{align*}
\frac{1}{w\left(  x\right)  }=\left(  2\pi\right)  ^{\frac{d}{2}}\prod
_{i=1}^{d}\frac{\sin^{2l}x_{i}}{x_{i}^{2n}}=\left(  2\pi\right)  ^{\frac{d}%
{2}}\prod_{i=1}^{d}x_{i}^{2\left(  l-n\right)  }\frac{\sin^{2l}x_{i}}%
{x_{i}^{2l}} &  =\left(  2\pi\right)  ^{\frac{d}{2}}x^{2\left(  l-n\right)
\mathbf{1}}\prod_{i=1}^{d}\frac{\sin^{2l}x_{i}}{x_{i}^{2l}}\\
&  =\left(  2\pi\right)  ^{\frac{d}{2}}x^{2\left(  l-n\right)  \mathbf{1}%
}\left(  \prod_{i=1}^{d}\frac{\sin^{2}x_{i}}{x_{i}^{2}}\right)  ^{l}\\
&  =\left(  2\pi\right)  ^{\frac{d}{2}}x^{2\left(  l-n\right)  \mathbf{1}%
}\left(  \left(  2\pi\right)  ^{\frac{d}{2}}\widehat{\Lambda}\left(
2x\right)  \right)  ^{l},
\end{align*}

where the last step used property \ref{p07} of the tensor product hat function
defined by \ref{p04} and \ref{p06}. But by \ref{p04}, $\Lambda$ has bounded
support so by Appendix \ref{SbSect_property_E'}, $\widehat{\Lambda}\in
C_{BP}^{\infty}$. Theorem \ref{Thm_product_of_Co,a_funcs} can now be applied
to prove that $x^{2\left(  l-n\right)  \mathbf{1}}\in C_{\emptyset,2\left(
l-n\right)  \mathbf{1}}^{\infty}\cap C_{BP}^{\infty}$ and consequently that
$\frac{1}{w}\in C_{\emptyset,2\left(  l-n\right)  \mathbf{1}}^{\infty}\cap
C_{BP}^{\infty}$.
\end{proof}

\begin{lemma}
\label{Lem_thin_plate_splin}\textbf{Thin-plate spline properties} Suppose
$T_{m}$\ is the thin-plate spline function with integer order $m>d/2$ i.e.
from \ref{m19},%
\begin{equation}
T_{m}\left(  x\right)  :=\left\{
\begin{array}
[c]{ll}%
\left(  -1\right)  ^{m-\frac{d}{2}+1}\left\vert x\right\vert ^{2m-d}%
\log\left\vert x\right\vert , & d\text{ }is\text{ }even,\\
\left(  -1\right)  ^{m-\frac{1}{2}\left(  d-1\right)  }\left\vert x\right\vert
^{2m-d}, & d\text{ }is\text{ }odd,
\end{array}
\right. \label{p24}%
\end{equation}

for all $x\in\mathbb{R}^{d}$. Then:

\begin{enumerate}
\item If $c_{m}=e\left(  m-d/2\right)  $ where the function $e$ is given by
equation \ref{m25} or \ref{m29}\ we have%
\begin{equation}
\widehat{T_{m}}=\frac{c_{m}}{\left\vert \cdot\right\vert ^{2m}}\text{
}on\text{ }S_{\emptyset,2m},\label{p25}%
\end{equation}

\item and $T_{m}\in C_{BP}^{\left(  2m-d-1\right)  }\cap C^{\infty}\left(
\mathbb{R}^{d}\setminus0\right)  $.

\item When $d$ \textbf{is odd}
\[
\left\vert \left(  \widehat{a}D\right)  ^{n}T_{m}\left(  x\right)  \right\vert
\leq k_{n,2m-d}\left\vert x\right\vert ^{2m-d-n},
\]

where $k_{n,2m-d}=k_{n,2\left(  m-\frac{d+1}{2}\right)  +1}$ is given by the
\textbf{conjectures} \ref{p12}. Also, using the Pochhammer symbol,%
\[
\left\vert D\right\vert ^{2n}T_{m}\left(  x\right)  =\left\{
\begin{array}
[c]{ll}%
2^{2n}\left(  -\left(  m-\frac{d-1}{2}\right)  \right)  _{n}\left(  -\left(
m-\frac{1}{2}\right)  \right)  _{n}\left\vert x\right\vert ^{2m-d-2n}, & n\leq
m-\frac{d-1}{2},\\
0, & n>m-\frac{d-1}{2},
\end{array}
\right.
\]

and%
\[
\left(  \widehat{x}D\right)  ^{n}T_{m}\left(  x\right)  =\left\{
\begin{array}
[c]{ll}%
\left(  -1\right)  ^{n+m-\frac{d-1}{2}}\left(  -\left(  2m-d\right)  \right)
_{n}, & n\leq2m-d,\\
0, & n>2m-d,
\end{array}
\right\}  \left\vert x\right\vert ^{2m-d-n}.
\]

\item When $d$ \textbf{is even}, using the Pochhammer symbol,%
\[
\frac{1}{n!}\left\vert \left(  \widehat{a}D\right)  ^{n}T_{m}\left(  x\right)
\right\vert \leq\left(  \frac{k_{n,2m-d}}{n!}\left\vert \log\left\vert
x\right\vert \right\vert +\sum\limits_{j=0}^{n-1}\left(  2-\frac{1}{j}\right)
\frac{k_{j,2m-d}}{j!}\right)  \left\vert x\right\vert ^{2m-d-n},\text{ }%
n\geq1,
\]

where the $k_{n-j,2m-d}$ are given by the \textbf{conjectures} \ref{p14}. Also%
\[
\left\vert D\right\vert ^{2n}T_{m}\left(  x\right)  =\left(  -1\right)
^{m-\frac{d}{2}+1}2^{2n}\left(  -m+\frac{d}{2}\right)  _{n}\left(
-m+1\right)  _{n}\left\vert x\right\vert ^{2m-d-2n},\text{ }n\geq1,
\]

and for $n\geq0$,%
\[
\frac{1}{n!}\left(  \widehat{\cdot}D\right)  ^{n}T_{m}\left(  x\right)
=\left\{
\begin{array}
[c]{ll}%
\left(  \tbinom{2m-d}{n}\log\left\vert x\right\vert +\sum\limits_{j=1}%
^{n}\frac{\left(  -1\right)  ^{j+1}}{j}\tbinom{2m-d}{n-j}\right)  , & 1\leq
n\leq2m-d,\\
\left(  -1\right)  ^{n-1}B\left(  2m-d+1,n-2m+d\right)  , & n>2m-d.
\end{array}
\right\}  \left\vert x\right\vert ^{2m-d-n}.
\]

\end{enumerate}
\end{lemma}

\begin{proof}
\textbf{Part 1} From \ref{m19} and Theorem \ref{Thm_surf_spline_ext_basis}%
\ \textbf{below} we have: the thin-plate spline (or surface spline) function
$H $ defined by
\[
H\left(  x\right)  =\left\{
\begin{array}
[c]{ll}%
\left(  -1\right)  ^{s+1}\left\vert x\right\vert ^{2s}\log\left\vert
x\right\vert , & s=1,2,3,\ldots,\\
\left(  -1\right)  ^{\left\lceil s\right\rceil }\left\vert x\right\vert
^{2s}, & s>0,\text{ }s\neq1,2,3\ldots,
\end{array}
\right.
\]

satisfies%
\[
\widehat{H}=\frac{e\left(  s\right)  }{\left\vert \cdot\right\vert ^{2s+d}%
}\text{ }on\text{ }\mathbb{R}^{d}\setminus0.
\]

Further, the function%
\[
w\left(  \xi\right)  =\frac{1}{e\left(  s\right)  }\left\vert \xi\right\vert
^{-2\theta+2s+d},\quad\xi\neq0.,
\]

is a weight function with order $\theta$ and smoothness parameter $\kappa$ iff
$0\leq\kappa<s<\theta$. In this case $H$ is a basis function of order $\theta$
generated by $w$.

Now if $s=m-d/2$ and $\theta=m$ then there exists $\kappa$ such that
$0\leq\kappa<m-d/2<m$ and we have $H=T_{m}$ and
\[
\widehat{T_{m}}=\frac{c_{m}}{\left\vert \cdot\right\vert ^{2m}}\text{
}on\text{ }S_{\emptyset,2m}.
\]

\textbf{Part 2} From part 4 Theorem \ref{Thm_surf_spline_ext_basis},
\begin{align*}
\left\lfloor 2\kappa\right\rfloor  & \leq\left\{
\begin{array}
[c]{ll}%
2s-1, & s=1,2,3,\ldots,\\
2\left\lfloor s\right\rfloor , & n<s\leq n+1/2,\text{ }n=0,1,2,\ldots,\\
2\left\lfloor s\right\rfloor +1, & n+1/2<s<n+1,\text{ }n=0,1,2,\ldots.
\end{array}
\right. \\
& =\left\{
\begin{array}
[c]{ll}%
2\left(  m-d/2\right)  -1, & d\text{ }is\text{ }even,\\
2\left\lfloor m-d/2\right\rfloor , & d\text{ }is\text{ }odd,
\end{array}
\right. \\
& =\left\{
\begin{array}
[c]{ll}%
2m-d-1, & d\text{ }is\text{ }even,\\
2m-d-1, & d\text{ }is\text{ }odd,
\end{array}
\right. \\
& =2m-d-1,
\end{align*}

so by part 4 Theorem \ref{Thm_basis_smth_W3.2_r3_pos}, $T_{m}\in
C_{BP}^{\left(  2m-d-1\right)  }\cap C^{\infty}\left(  \mathbb{R}^{d}%
\setminus0\right)  $.\medskip

\textbf{Part 3} $d$\textbf{\ is odd} The bound on $\left(  \widehat
{a}D\right)  ^{n}T_{m}\left(  x\right)  $ follows directly from the
\textbf{conjectures} of Corollary \ref{Cor_Thm_(aD)^m_|x|^n}:%
\[
\left\vert \left(  \widehat{a}D\right)  ^{n}T_{m}\left(  x\right)  \right\vert
=\left\vert \left(  \widehat{a}D\right)  ^{n}\left(  -1\right)  ^{m-\frac
{1}{2}\left(  d-1\right)  }\left\vert x\right\vert ^{2m-d}\right\vert
=\left\vert \left(  \widehat{a}D\right)  ^{n}\left\vert x\right\vert
^{2m-d}\right\vert \leq k_{n,2m-d}\left\vert x\right\vert ^{2m-d-n},
\]

where $k_{n,2m-d}$ is given by \ref{p12}.

From part 9 of Lemma \ref{Lem_iterLaplacian_rad}, for integer $n,k\geq0$,%
\[
\left\vert D\right\vert ^{2n}\left(  \left\vert x\right\vert ^{2k-1}\right)
=2^{2n}\left(  -k+\frac{1}{2}\right)  _{n}\left(  -k-\frac{d}{2}+\frac{3}%
{2}\right)  _{n}\left\vert x\right\vert ^{2k-1-2n},
\]

so that%
\begin{align*}
\left\vert D\right\vert ^{2n} &  T_{m}\left(  x\right) \\
&  =\left\vert D\right\vert ^{2n}\left\vert x\right\vert ^{2m-d}\\
&  =\left\vert D\right\vert ^{2n}\left\vert x\right\vert ^{2m-d}\\
&  =\left\{
\begin{array}
[c]{ll}%
2^{2n}\left(  -\left(  m-\frac{d-1}{2}\right)  \right)  _{n}\left(  -\left(
m-\frac{d-1}{2}\right)  -\frac{d}{2}+1\right)  _{n}\left\vert x\right\vert
^{2m-d-2n}, & n\leq m-\frac{d-1}{2},\\
0, & n>m-\frac{d-1}{2},
\end{array}
\right. \\
&  =\left\{
\begin{array}
[c]{ll}%
2^{2n}\left(  -\left(  m-\frac{d-1}{2}\right)  \right)  _{n}\left(  -\left(
m-\frac{1}{2}\right)  \right)  _{n}\left\vert x\right\vert ^{2m-d-2n}, & n\leq
m-\frac{d-1}{2},\\
0, & n>m-\frac{d-1}{2}.
\end{array}
\right.
\end{align*}

From part 1 of Lemma \ref{Lem_deriv_rad_funcs}, when $d$ is odd,
\begin{align*}
\left(  \widehat{x}D\right)  ^{n}T_{m}\left(  x\right)   &  =\left(
\widehat{x}D\right)  ^{n}\left(  -1\right)  ^{m-\frac{d-1}{2}}\left\vert
x\right\vert ^{2m-d}\\
&  =\left(  -1\right)  ^{m-\frac{d-1}{2}}\left(  D^{n}s^{2m-d}\right)  \left(
\left\vert x\right\vert \right) \\
&  =\left\{
\begin{array}
[c]{ll}%
\left(  -1\right)  ^{n+m-\frac{d-1}{2}}\left(  -\left(  2m-d\right)  \right)
_{n}, & n\leq2m-d,\\
0, & n>2m-d,
\end{array}
\right\}  \left\vert x\right\vert ^{2m-d-n}.
\end{align*}
\medskip

\textbf{Part 4} $d$\textbf{\ is even} The bound on $\left(  \widehat
{a}D\right)  ^{2n}T_{m}\left(  x\right)  $ follows directly from the
\textbf{conjectures} of Corollary \ref{Cor_Thm_(aD)^m_|x|^n}.

From part 8 of Lemma \ref{Lem_iterLaplacian_rad}, when $2n\leq2m-d$,%
\begin{align*}
\left\vert D\right\vert ^{2n}T_{m}\left(  x\right)   & =\left(  -1\right)
^{m-\frac{d}{2}+1}\left\vert D\right\vert ^{2n}\left(  \left\vert x\right\vert
^{2m-d}\log\left\vert x\right\vert \right) \\
& =\left(  -1\right)  ^{m-\frac{d}{2}+1}2^{2n}\left(  -k\right)  _{n}\left(
-k-\frac{d}{2}+1\right)  _{n}\left\vert x\right\vert ^{2k-2n}\\
& =\left(  -1\right)  ^{m-\frac{d}{2}+1}2^{2n}\left(  -\frac{2k}{2}\right)
_{n}\left(  -\frac{1}{2}\left(  2k+d\right)  +1\right)  _{n}\left\vert
x\right\vert ^{2k-2n}\\
& =\left(  -1\right)  ^{m-\frac{d}{2}+1}2^{2n}\left(  -\frac{2m-d}{2}\right)
_{n}\left(  -m+1\right)  _{n}\left\vert x\right\vert ^{2m-d-2n}\\
& =\left(  -1\right)  ^{m-\frac{d}{2}+1}2^{2n}\left(  -\frac{2m-d}{2}\right)
_{n}\left(  -m+1\right)  _{n}\left\vert x\right\vert ^{2m-d-2n}\\
& =\left(  -1\right)  ^{m-\frac{d}{2}+1}2^{2n}\left(  -m+\frac{d}{2}\right)
_{n}\left(  -m+1\right)  _{n}\left\vert x\right\vert ^{2m-d-2n}.
\end{align*}
\smallskip

Finally, the formula for $\frac{1}{n!}\left(  \widehat{x}D\right)  ^{n}%
T_{m}\left(  x\right)  $ follows from an application of part 1 of Lemma
\ref{Lem_deriv_rad_funcs} to part 7 of the same lemma to give:%
\begin{align*}
D^{n}\left(  s^{k}\log s\right)   & =\left\{
\begin{array}
[c]{ll}%
n!\left(  \tbinom{k}{n}\log s+\sum\limits_{j=1}^{n}\frac{\left(  -1\right)
^{j+1}}{j}\tbinom{k}{n-j}\right)  s^{k-n}, & 1\leq n\leq k,\\
\left(  -1\right)  ^{n-k-1}k!\left(  n-k-1\right)  !s^{k-n}, & n>k,
\end{array}
\right. \\
\frac{1}{n!}D^{n}\left(  s^{k}\log s\right)   & =\left\{
\begin{array}
[c]{ll}%
\left(  \tbinom{k}{n}\log s+\sum\limits_{j=1}^{n}\frac{\left(  -1\right)
^{j+1}}{j}\tbinom{k}{n-j}\right)  s^{k-n}, & 1\leq n\leq k,\\
\left(  -1\right)  ^{n-k-1}B\left(  k+1,n-k\right)  s^{k-n}, & n>k,
\end{array}
\right.
\end{align*}

so that when $k=2m-d$,%
\begin{align*}
\frac{1}{n!}\left(  \widehat{x}D\right)  ^{n}T_{m}\left(  x\right)   &
=\frac{1}{n!}\left(  \widehat{x}D\right)  ^{n}\left(  \left\vert x\right\vert
^{2m-d}\log\left\vert x\right\vert \right) \\
& =\frac{1}{n!}\left(  D^{n}\left(  s^{2m-d}\log s\right)  \right)  \left(
s=\left\vert x\right\vert \right) \\
& =\left\{
\begin{array}
[c]{ll}%
\tbinom{2m-d}{n}\log\left\vert x\right\vert +\sum\limits_{j=1}^{n}%
\frac{\left(  -1\right)  ^{j+1}}{j}\tbinom{2m-d}{n-j}, & 1\leq n\leq2m-d,\\
\left(  -1\right)  ^{n-1}B\left(  2m-d+1,n-2m+d\right)  , & n>2m-d.
\end{array}
\right\}  \left\vert x\right\vert ^{2m-d-n}%
\end{align*}

\end{proof}

\begin{theorem}
\label{Thm_ExtNatSplin_convol}Suppose $w$ is an \textbf{extended B-spline}
weight function \ref{p32} which has property W3.1 for order $\theta$ and
smoothness $\kappa$. Then $G_{0}=\left(  \frac{1}{w}\right)  ^{\vee}$ is the
\textbf{zero order} extended B-spline basis function and $G_{\theta}=\left(
\frac{1}{w\left\vert \cdot\right\vert ^{2\theta}}\right)  ^{\vee}$ is an order
$\theta$ basis function. In fact, the basis function $G_{\theta}$ has the
form
\begin{equation}
G_{\theta}=\left\{
\begin{array}
[c]{ll}%
\frac{1}{e\left(  \theta-d/2\right)  }G_{0}\ast T_{\theta}, & 2\theta>d,\\
\frac{\left(  -1\right)  ^{\frac{d+1}{2}-\theta}}{e\left(  1/2\right)
}\left(  \left\vert D\right\vert ^{2}\right)  ^{\frac{d+1}{2}-\theta}\left(
G_{0}\ast T_{\frac{d+1}{2}}\right)  , & 2\theta<d,\text{ }d\text{ }odd,\text{
}d\geq3,\\
\frac{\left(  -1\right)  ^{\frac{d}{2}+1-\theta}}{e\left(  1\right)  }\left(
\left\vert D\right\vert ^{2}\right)  ^{\frac{d}{2}+1-\theta}\left(  G_{0}\ast
T_{\frac{d}{2}+1}\right)  , & 2\theta\leq d,\text{ }d\text{ }even,
\end{array}
\right. \label{p23}%
\end{equation}

where $T_{m}$\ is the thin-plate spline basis function of integer order
$m>d/2$ given in Lemma \ref{Lem_thin_plate_splin}. Here $T_{\frac{d+1}{2}%
}=-\left\vert \cdot\right\vert $ when $d$ is odd and $T_{\frac{d}{2}%
+1}=\left\vert \cdot\right\vert ^{2}\log\left\vert \cdot\right\vert $ when $d$
is even. The values of $e\left(  \cdot\right)  $ are given in Remark
\ref{Rem_Lem_e(s)_pos}.

The convolutions of \ref{p23} can be written%
\begin{equation}
\left(  G_{0}\ast T_{m}\right)  \left(  x\right)  =\left(  2\pi\right)
^{-d/2}\int G_{0}\left(  y\right)  T_{m}\left(  x-y\right)  dy,\quad
2m>d,\label{p19}%
\end{equation}

where $G_{0}\ast T_{m}$ is the regular tempered distribution defined in
Appendix \ref{SbSect_property_S'}.
\end{theorem}

\begin{proof}
Since $w$ has property W3.1 for parameters $\theta$ and $\kappa$ it follows
from Lemma \ref{Lem_basis_ord0_Bsplin_W3} that $\frac{1}{w}\in L^{1}$ and that
$G=\left(  \frac{1}{w}\right)  ^{\vee}$ is the unique zero order basis
function. Also, from Theorem \ref{Thm_basis_smth_W3.1} we have that $\left(
\frac{1}{w\left\vert \cdot\right\vert ^{2\theta}}\right)  ^{\vee}$ is a basis
function of order $\theta$. Now suppose that $\phi\in S$. By Lemma
\ref{Lem_ExtNatSplin_1/w_in_Co,2d(l-n)}, $\frac{1}{w}\in C_{\emptyset
,2\theta\mathbf{1}}^{\infty}\cap C_{BP}^{\infty}$, and so by part 4 Theorem
\ref{Thm_product_of_Co,a_funcs},
\begin{equation}
\frac{1}{w}\phi\in S_{\emptyset,2\theta\mathbf{1}}\subset S_{\emptyset
,2\theta}.\label{p22}%
\end{equation}
\medskip

\fbox{\textbf{Case} $2\theta>d$} By Lemma \ref{Lem_thin_plate_splin},
$\frac{1}{c_{\theta}}\widehat{T_{\theta}}=\frac{1}{e\left(  \theta-d/2\right)
}\widehat{T_{\theta}}=\frac{1}{\left\vert \cdot\right\vert ^{2\theta}}$ on
$S_{\emptyset,2\theta}$, and thus for $\phi\in S$,%
\begin{align*}
\left[  \widehat{G_{\theta}},\phi\right]  =\left[  \frac{1}{w\left\vert
\cdot\right\vert ^{2\theta}},\phi\right]  =\left[  \frac{1}{\left\vert
\cdot\right\vert ^{2\theta}},\frac{1}{w}\phi\right]   &  =\left[  \frac
{1}{e\left(  \theta-d/2\right)  }\widehat{T_{\theta}},\frac{1}{w}\phi\right]
\\
&  =\left[  \frac{1}{e\left(  \theta-d/2\right)  }\frac{1}{w}\widehat
{T_{\theta}},\phi\right] \\
&  =\left[  \frac{1}{e\left(  \theta-d/2\right)  }\widehat{G_{0}}%
\widehat{T_{\theta}},\phi\right] \\
&  =\left[  \frac{1}{e\left(  \theta-d/2\right)  }\widehat{G_{0}\ast
T_{\theta}},\phi\right]  ,
\end{align*}

so that $G_{\theta}=\frac{1}{e\left(  \theta-d/2\right)  }G_{0}\ast T_{\theta
}$.\medskip

\fbox{\textbf{Case} $2\theta<d$ and $d$ is odd} By Lemma
\ref{Lem_thin_plate_splin}, $\widehat{T_{\frac{d+1}{2}}}=c_{\frac{d+1}{2}%
}\frac{1}{\left\vert \cdot\right\vert ^{2\left(  \frac{d+1}{2}\right)  }%
}=\frac{e\left(  1/2\right)  }{\left\vert \cdot\right\vert ^{2\left(
\frac{d+1}{2}\right)  }}$ on $S_{\emptyset,d+1}$, and thus \ref{p22} implies
for $\phi\in S$,
\begin{align*}
\left[  \widehat{G_{\theta}},\phi\right]  =\left[  \frac{1}{w\left\vert
\cdot\right\vert ^{2\theta}},\phi\right]  =\left[  \frac{1}{\left\vert
\cdot\right\vert ^{2\theta}},\frac{\phi}{w}\right]   &  =\left[
\frac{\left\vert \cdot\right\vert ^{d+1-2\theta}}{\left\vert \cdot\right\vert
^{d+1}},\frac{\phi}{w}\right] \\
&  =\left[  \frac{1}{\left\vert \cdot\right\vert ^{2\left(  \frac{d+1}%
{2}\right)  }},\frac{\left\vert \cdot\right\vert ^{d+1-2\theta}}{w}\phi\right]
\\
&  =\left[  \frac{1}{e\left(  1/2\right)  }\widehat{T_{\frac{d+1}{2}}}%
,\frac{\left\vert \cdot\right\vert ^{d+1-2\theta}}{w}\phi\right] \\
&  =\left[  \frac{1}{e\left(  1/2\right)  }\frac{\left\vert \cdot\right\vert
^{d+1-2\theta}}{w}\widehat{T_{\frac{d+1}{2}}},\phi\right] \\
&  =\left[  \frac{1}{e\left(  1/2\right)  }\left\vert \cdot\right\vert
^{d+1-2\theta}\widehat{G_{0}}\text{ }\widehat{T_{\frac{d+1}{2}}},\phi\right]
\\
&  =\left[  \frac{1}{e\left(  1/2\right)  }\left\vert \cdot\right\vert
^{d+1-2\theta}\left(  G_{0}\ast T_{\frac{d+1}{2}}\right)  ^{\wedge}%
,\phi\right] \\
&  =\left[  \frac{\left(  -1\right)  ^{\frac{d+1}{2}-\theta}}{e\left(
1/2\right)  }\left(  \left\vert D\right\vert ^{d+1-2\theta}G_{0}\ast
T_{\frac{d+1}{2}}\right)  ^{\wedge},\phi\right]  ,
\end{align*}

so that
\[
G_{\theta}=\frac{\left(  -1\right)  ^{\frac{d+1}{2}-\theta}}{e\left(
1/2\right)  }\left\vert D\right\vert ^{d+1-2\theta}G_{0}\ast T_{\frac{d+1}{2}%
}.
\]

From \ref{p24}, $T_{\frac{d+1}{2}}=\left(  -1\right)  ^{\frac{d+1}{2}%
-\frac{d-1}{2}}\left\vert \cdot\right\vert ^{2\frac{d+1}{2}-d}=-\left\vert
\cdot\right\vert $.\medskip

\fbox{\textbf{Case} $2\theta\leq d$ and $d$ even} In a very similar manner to
the previous case we can obtain%
\[
G_{\theta}=\frac{\left(  -1\right)  ^{\frac{d}{2}+1-\theta}}{e\left(
1\right)  }\left\vert D\right\vert ^{d+2-2\theta}G_{0}\ast T_{\frac{d}{2}+1}.
\]

From \ref{p24}, $T_{\frac{d}{2}+1}=\left(  -1\right)  ^{\frac{d}{2}%
+1-\frac{d-2}{2}}\left\vert \cdot\right\vert ^{2\left(  \frac{d}{2}+1\right)
-d}\log\left\vert \cdot\right\vert =\left\vert \cdot\right\vert ^{2}%
\log\left\vert \cdot\right\vert $.

The last three cases combine to prove \ref{p23}.

Since Theorem \ref{Lem_basis_ord0_Bsplin_W3} implies $G_{0}\in C_{0}^{\left(
\left(  2n-2\right)  \mathbf{1}\right)  }\left(  \mathbb{R}^{d}\right)  $ and
Lemma \ref{Lem_thin_plate_splin} shows that $T_{m}\in C_{BP}^{\left(
2m-d-1\right)  }$, \ref{p19} follows from the convolution formulas of part 5
of Appendix \ref{Def_Convolution}.
\end{proof}

\begin{remark}
In \ref{p23} we considered the case \fbox{$2\theta<d$, $d$ odd, $d\geq3$.}
From \ref{p24} of Lemma \ref{Lem_thin_plate_splin}:
\[
T_{m}\left(  x\right)  =\left(  -1\right)  ^{m-\frac{1}{2}\left(  d-1\right)
}\left\vert x\right\vert ^{2m-d},\quad m>d/2,
\]

so that%
\[
T_{\frac{d+1}{2}}\left(  x\right)  =-\left\vert x\right\vert .
\]

From part 9 of Lemma \ref{Lem_iterLaplacian_rad},%
\[
\left\vert D\right\vert ^{2n}\left\vert x\right\vert ^{2k-1}=2^{2n}\left(
-k+\frac{1}{2}\right)  _{n}\left(  -k-\frac{d}{2}+\frac{3}{2}\right)
_{n}\left\vert x\right\vert ^{2k-1-2n},\quad n,k\geq1,
\]

so%
\[
\left\vert D\right\vert ^{2n}T_{\frac{d+1}{2}}=-\left\vert D\right\vert
^{2n}\left\vert x\right\vert =-2^{2n}\left(  -\frac{1}{2}\right)  _{n}\left(
-\frac{d-1}{2}\right)  _{n}\left\vert x\right\vert ^{1-2n},\quad n\geq1.
\]

Hence%
\begin{align*}
\left\vert D\right\vert ^{d+1-2\theta}T_{\frac{d+1}{2}}  & =\left\vert
D\right\vert ^{2\left(  \frac{d+1}{2}-\theta\right)  }T_{\frac{d+1}{2}}\\
& =-2^{d+1-2\theta}\left(  -\frac{1}{2}\right)  _{\frac{d+1}{2}-\theta}\left(
-\frac{d-1}{2}\right)  _{\frac{d+1}{2}-\theta}\left\vert x\right\vert
^{2\theta-d}\\
& \in L_{loc}^{1},
\end{align*}

since $2\theta-d>-d$. Further%
\begin{align*}
G_{\theta} &  =\frac{\left(  -1\right)  ^{\frac{d+1}{2}-\theta}}{e\left(
1/2\right)  }\left\vert D\right\vert ^{d+1-2\theta}\left(  G_{0}\ast
T_{\frac{d+1}{2}}\right) \\
&  =\frac{\left(  -1\right)  ^{\frac{d+1}{2}-\theta}}{e\left(  1/2\right)
}G_{0}\ast\left\vert D\right\vert ^{d+1-2\theta}T_{\frac{d+1}{2}}\\
&  =\frac{\left(  -1\right)  ^{\frac{d+1}{2}-\theta}}{e\left(  1/2\right)
}\left(  -1\right)  2^{d+1-2\theta}\left(  -\frac{1}{2}\right)  _{\frac
{d+1}{2}-\theta}\left(  -\frac{d-1}{2}\right)  _{\frac{d+1}{2}-\theta
}\left\vert x\right\vert ^{2\theta-d}\\
&  =\left(  \frac{\left(  -1\right)  ^{\frac{d-1}{2}-\theta}}{e\left(
1/2\right)  }2^{d+1-2\theta}\left(  -\frac{1}{2}\right)  _{\frac{d+1}%
{2}-\theta}\left(  -\frac{d-1}{2}\right)  _{\frac{d+1}{2}-\theta}\right)
G_{0}\ast\frac{1}{\left\vert \cdot\right\vert ^{d-2\theta}}.
\end{align*}

Since $D^{2n\mathbf{1}}G_{0}$ is the finite sum of translated delta functions
we have%
\begin{align*}
D^{2n\mathbf{1}}G_{\theta}  & =\left(  \frac{\left(  -1\right)  ^{\frac
{d-1}{2}-\theta}}{e\left(  1/2\right)  }2^{d+1-2\theta}\left(  -\frac{1}%
{2}\right)  _{\frac{d+1}{2}-\theta}\left(  -\frac{d-1}{2}\right)  _{\frac
{d+1}{2}-\theta}\right)  \left(  D^{2n\mathbf{1}}G_{0}\right)  \ast\frac
{1}{\left\vert \cdot\right\vert ^{d-2\theta}}\\
& \in L_{loc}^{1}.
\end{align*}

From Remark \ref{Rem_Lem_e(s)_pos}, since $d$ is odd,%
\[
e\left(  1/2\right)  =\pi^{\frac{d-1}{2}}2^{d}\left(  \frac{d-1}{2}\right)  !
\]

\end{remark}

\begin{remark}
?? Can this formula be used to characterize the data space locally?
\end{remark}

\subsection{Convolution formulas for basis functions generated by weight
functions with property W3.1\label{SbSect_convol_W3.1}}

In this subsection we will develop convolution formulas designed to facilitate
the expression of positive order basis functions in a closed form which can be
used for numerical work. Here the weight function has property W3.1.

\begin{lemma}
\label{Lem_W3.1_to_W2_zero}Suppose the weight function $w$ has property W3.1
for some $\left\vert \alpha\right\vert =\theta$ and $\kappa$. Then:

\begin{enumerate}
\item $w_{\alpha}\left(  x\right)  =w\left(  x\right)  x^{2\alpha}$ is a
\textbf{zero order weight function} with property W02 for parameter $\kappa$.

\item $G_{\alpha}=\left(  1/w_{\alpha}\right)  ^{\vee}\in C_{B}^{\left(
\left\lfloor 2\kappa\right\rfloor \right)  }$ is the unique zero order basis
function generated by $w_{\alpha}$.

\item $G_{0}=\left(  1/w\right)  ^{\vee}\in S^{\prime}$ and we call this the
zero order distribution generated by $w$.

\item $G_{\theta}=\left(  \frac{1}{w\left\vert \cdot\right\vert ^{2\theta}%
}\right)  ^{\vee}\in C_{BP}^{\left(  \left\lfloor 2\kappa\right\rfloor
\right)  }$ is a basis function of order $\theta$ generated by $w$.
\end{enumerate}
\end{lemma}

\begin{proof}
The properties W01 and W02 of a zero order weight function are given in
Subsection \ref{SbSect_ZeroOrdWeightBasis}. Observe that the properties W01
and W1 for the zero and positive order cases are identical.\medskip

\textbf{Part 1} Since $w$ has property W1 it has property W01 as a zero order
weight function and must be continuous and positive outside some closed set
$\mathcal{A}$ with measure zero. Hence $w_{\alpha}$ is continuous and positive
outside the closed set $\mathcal{A}\cup\bigcup\limits_{k=1}^{d}\left\{
x:x_{k}=0\right\}  $. Since this set also has measure zero, $w_{\alpha}$ also
has property W01 i.e. it is a zero order weight function.

Property W3.1 is $\int\dfrac{x^{2\lambda}}{w\left(  x\right)  x^{2\alpha}%
}dx<\infty$ when $0\leq\lambda\leq\kappa$ and for some $\left\vert
\alpha\right\vert =\theta$, and thus $w_{\alpha}$ has the zero order weight
function property W03 for parameter $\kappa$.\medskip

\textbf{Part 2} From Subsection \ref{SbSect_ZeroOrdWeightBasis} the zero order
basis function is $\left(  1/w_{\alpha}\right)  ^{\vee}$ and it is a function
in $C_{B}^{\left(  \left\lfloor 2\kappa\right\rfloor \right)  }$.\medskip

\textbf{Part 3} By part 1 of Theorem \ref{Thm_basis_smth_W3.2_r3_pos}, $1/w\in
S^{\prime}$.\medskip

\textbf{Part 4} By Theorem \ref{Thm_basis_smth_W3.1}, $G_{\theta}=\left(
\frac{1}{w\left\vert \cdot\right\vert ^{2\theta}}\right)  ^{\vee}\in
C_{BP}^{\left(  \left\lfloor 2\kappa\right\rfloor \right)  }$ is a basis
function of order $\theta$ generated by $w$.
\end{proof}

We begin with the case where $G_{\alpha}$\ has bounded support, which
simplifies the convolution.

\begin{theorem}
Suppose the weight function $w$ has property W3.1 for some $\left\vert
\alpha\right\vert =\theta$ and parameter $\kappa$. Assume the weight function
$w_{\alpha}$, the basis functions $G_{\alpha}$, $G_{\theta}$ and the basis
distribution $G_{0}$ are as defined above in Lemma \ref{Lem_W3.1_to_W2_zero}.
Then if $G_{\alpha}$ has bounded support the convolution formulas
\begin{equation}
G_{\theta}=\left\{
\begin{array}
[c]{ll}%
\frac{\left(  -1\right)  ^{\theta}}{e\left(  \theta-d/2\right)  }D^{2\alpha
}G_{\alpha}\ast T_{\theta}, & 2\theta>d,\\
\frac{\left(  -1\right)  ^{\frac{d+1}{2}}}{e\left(  1/2\right)  }\left\vert
D\right\vert ^{d+1-2\theta}D^{2\alpha}G_{\alpha}\ast T_{\frac{d+1}{2}}, &
2\theta<d,\text{ }d\text{ }odd,\\
\frac{\left(  -1\right)  ^{\frac{d+2}{2}}}{e\left(  1\right)  }\left\vert
D\right\vert ^{d+2-2\theta}D^{2\alpha}G_{\alpha}\ast T_{\frac{d}{2}+1}, &
2\theta\leq d,\text{ }d\text{ }even,
\end{array}
\right. \label{p38}%
\end{equation}

hold. Here $T_{m}$\ is the thin-plate spline basis function of integer order
$m>d/2$ defined in Lemma \ref{Lem_thin_plate_splin}. The values of the
function $e$ are given in Remark \ref{Rem_Lem_e(s)_pos}.

The convolutions of \ref{p38} can be written%
\begin{equation}
\left(  G_{\alpha}\ast T_{m}\right)  \left(  x\right)  =\left(  2\pi\right)
^{-d/2}\int G_{\alpha}\left(  y\right)  T_{m}\left(  x-y\right)  dy,\quad
m>d/2,\label{p56}%
\end{equation}

with $G_{\alpha}\ast T_{m}$ being a regular tempered distribution (Appendix
\ref{SbSect_property_S'}).

Finally
\begin{equation}
\left(  -1\right)  ^{\theta}D^{2\alpha}G_{\alpha}=G_{0}.\label{p57}%
\end{equation}

\end{theorem}

\begin{proof}
Since $G_{\alpha}$ has bounded support it follows from Appendix
\ref{SbSect_property_E'} and Lemma \ref{Lem_W3.1_to_W2_zero} that
\begin{equation}
\widehat{G_{\alpha}}=\frac{1}{w_{\alpha}}\in C_{BP}^{\infty},\label{p27}%
\end{equation}

and by Theorem \ref{Thm_product_of_Co,k_funcs},%
\begin{equation}
\frac{\xi^{2\alpha}}{w_{\alpha}\left(  \xi\right)  }\phi\left(  \xi\right)
\in S_{\emptyset,2\theta},\quad\phi\in S.\label{p55}%
\end{equation}

\fbox{\textbf{Case} $2\theta>d$} Thus if $\xi$ is the \textit{action} variable%
\[
\left[  \widehat{G_{\theta}},\phi\right]  =\left[  \frac{1}{\left\vert
\cdot\right\vert ^{2\theta}w},\phi\right]  =\left[  \frac{1}{\left\vert
\cdot\right\vert ^{2\theta}w_{\alpha}},\xi^{2\alpha}\phi\right]  =\left[
\frac{1}{\left\vert \cdot\right\vert ^{2\theta}},\frac{\xi^{2\alpha}%
}{w_{\alpha}}\phi\right]  .
\]

But Lemma \ref{Lem_thin_plate_splin} tells us that $T_{\theta}\in
C_{BP}^{\left(  2\theta-d-1\right)  }$ and $\widehat{T_{\theta}}%
=\frac{c_{\theta}}{\left\vert \cdot\right\vert ^{2\theta}}$ on $S_{\emptyset
,2\theta}$. Hence%
\[
\left[  \frac{1}{\left\vert \cdot\right\vert ^{2\theta}},\frac{\xi^{2\alpha}%
}{w_{\alpha}}\phi\right]  =\left[  \frac{1}{c_{\theta}}\widehat{T_{\theta}%
},\frac{\xi^{2\alpha}}{w_{\alpha}}\phi\right]  =\left[  \frac{1}{c_{\theta}%
}\frac{\xi^{2\alpha}}{w_{\alpha}}\widehat{T_{\theta}},\phi\right]  =\left[
\frac{1}{c_{\theta}}\xi^{2\alpha}\widehat{G_{\alpha}}\text{ }\widehat
{T_{\theta}},\phi\right]  ,
\]

and, since $G_{\alpha}$ has bounded support, by part 4 of Appendix
\ref{Def_Convolution}%
\[
\left[  \frac{1}{c_{\theta}}\xi^{2\alpha}\widehat{G_{\alpha}}\text{ }%
\widehat{T_{\theta}},\phi\right]  =\left[  \frac{1}{c_{\theta}}\xi^{2\alpha
}\widehat{G_{\alpha}\ast T_{\theta}},\phi\right]  =\left[  \frac{\left(
-1\right)  ^{\theta}}{c_{\theta}}\left(  D^{2\alpha}G_{\alpha}\ast T_{\theta
}\right)  ^{\wedge},\phi\right]  ,
\]

so that $\widehat{G_{\theta}}=\frac{\left(  -1\right)  ^{\theta}}{c_{\theta}%
}\left(  D^{2\alpha}G_{\alpha}\ast T_{\theta}\right)  ^{\wedge}$ and
$G_{\theta}=\frac{\left(  -1\right)  ^{\theta}}{c_{\theta}}D^{2\alpha
}G_{\alpha}\ast T_{\theta}$.

Finally, from part 1 of Lemma \ref{Lem_thin_plate_splin}, $c_{\theta}=e\left(
\theta-d/2\right)  $.\medskip

\fbox{\textbf{Case} $2\theta<d$ and $d$ is odd} If $\xi$ is the distribution
action variable then since $\frac{1}{w_{\alpha}}\in C_{BP}^{\infty}$
\begin{align*}
\left[  \widehat{G_{\theta}},\phi\right]  =\left[  \frac{1}{\left\vert
\cdot\right\vert ^{2\theta}w},\phi\right]  =\left[  \frac{1}{\left\vert
\cdot\right\vert ^{2\left(  \frac{d+1}{2}\right)  }w},\left\vert
\cdot\right\vert ^{2\left(  \frac{d+1}{2}-\theta\right)  }\phi\right]   &
=\left[  \frac{1}{\left\vert \cdot\right\vert ^{2\left(  \frac{d+1}{2}\right)
}w_{\alpha}},\xi^{2\alpha}\left\vert \cdot\right\vert ^{2\left(  \frac{d+1}%
{2}-\theta\right)  }\phi\right] \\
&  =\left[  \frac{1}{\left\vert \cdot\right\vert ^{2\left(  \frac{d+1}%
{2}\right)  }},\frac{\xi^{2\alpha}}{w_{\alpha}}\left\vert \cdot\right\vert
^{2\left(  \frac{d+1}{2}-\theta\right)  }\phi\right]  .
\end{align*}

By Theorem \ref{Thm_product_of_Co,k_funcs}, $\xi^{2\alpha}\in S_{\emptyset
,2\theta}$, $\left\vert \cdot\right\vert ^{2\left(  \frac{d+1}{2}%
-\theta\right)  }\in S_{\emptyset,d+1-2\theta}$, $\xi^{2\alpha}\left\vert
\cdot\right\vert ^{2\left(  \frac{d+1}{2}-\theta\right)  }\in S_{\emptyset
,d+1}$ and $\xi^{2\alpha}\left\vert \cdot\right\vert ^{2\left(  \frac{d+1}%
{2}-\theta\right)  }\phi\in S_{\emptyset,d+1}$. Thus
\begin{equation}
\frac{\xi^{2\alpha}}{w_{\alpha}}\left\vert \cdot\right\vert ^{2\left(
\frac{d+1}{2}-\theta\right)  }\phi\in S_{\emptyset,d+1}.\label{p48}%
\end{equation}

Also by part 1 of Lemma \ref{Lem_thin_plate_splin},
\begin{equation}
\widehat{T_{\frac{d+1}{2}}}=c_{\frac{d+1}{2}}\frac{1}{\left\vert
\cdot\right\vert ^{2\left(  \frac{d+1}{2}\right)  }}=e\left(  1/2\right)
\frac{1}{\left\vert \cdot\right\vert ^{2\left(  \frac{d+1}{2}\right)  }}\text{
}on\text{ }S_{\emptyset,d+1},\label{p49}%
\end{equation}

so that%
\begin{align*}
\left[  \widehat{G_{\theta}},\phi\right]  =\left[  \frac{1}{\left\vert
\cdot\right\vert ^{2\left(  \frac{d+1}{2}\right)  }},\frac{\xi^{2\alpha}%
}{w_{\alpha}}\left\vert \cdot\right\vert ^{2\left(  \frac{d+1}{2}%
-\theta\right)  }\phi\right]   &  =\left[  \frac{1}{e\left(  1/2\right)
}\widehat{T_{\frac{d+1}{2}}},\frac{\xi^{2\alpha}}{w_{\alpha}}\left\vert
\cdot\right\vert ^{2\left(  \frac{d+1}{2}-\theta\right)  }\phi\right] \\
&  =\left[  \frac{1}{e\left(  1/2\right)  }\xi^{2\alpha}\left\vert
\cdot\right\vert ^{2\left(  \frac{d+1}{2}-\theta\right)  }\frac{1}{w_{\alpha}%
}\widehat{T_{\frac{d+1}{2}}},\phi\right] \\
&  =\left[  \frac{1}{e\left(  1/2\right)  }\xi^{2\alpha}\left\vert
\cdot\right\vert ^{d+1-2\theta}\widehat{G_{\alpha}}\text{ }\widehat
{T_{\frac{d+1}{2}}},\phi\right] \\
&  =\left[  \frac{1}{e\left(  1/2\right)  }\xi^{2\alpha}\left\vert
\cdot\right\vert ^{d+1-2\theta}\left(  G_{\alpha}\ast T_{\frac{d+1}{2}%
}\right)  ^{\wedge},\phi\right] \\
&  =\left[  \frac{\left(  -1\right)  ^{\frac{d+1}{2}}}{e\left(  1/2\right)
}\left(  \left\vert D\right\vert ^{d+1-2\theta}D^{2\alpha}G_{\alpha}\ast
T_{\frac{d+1}{2}}\right)  ^{\wedge},\phi\right]  ,
\end{align*}

and hence%
\[
G_{\theta}=\frac{\left(  -1\right)  ^{\frac{d+1}{2}}}{e\left(  1/2\right)
}\left\vert D\right\vert ^{d+1-2\theta}D^{2\alpha}G_{\alpha}\ast T_{\frac
{d+1}{2}}.
\]
\medskip

\fbox{\textbf{Case} $2\theta\leq d$ and $d$ is even} This case can be proved
in a very similar manner to the previous case.

The last three cases combine to prove \ref{p38}. Equation \ref{p56} follows
directly from part 5 of Appendix \ref{Def_Convolution}.\medskip

Since $\widehat{G_{\alpha}}=\frac{1}{w_{\alpha}}\in L^{1}$ we have
$x^{2\alpha}\widehat{G_{\alpha}}=\frac{1}{w}=\widehat{G_{0}}$, and hence
$\left(  -1\right)  ^{\theta}D^{2\alpha}G_{\alpha}=G_{0}\in S^{\prime}%
\cap\mathcal{E}^{\prime}$.
\end{proof}

\begin{remark}
The reciprocal of an extended B-spline weight function is $L^{1}$ so the basis
function $G_{0}$ is continuous. Substituting \ref{p57} i.e. $G_{0}=\left(
-1\right)  ^{\theta}D^{2\alpha}G_{\alpha}$ into \ref{p38} we obtain the
equations \ref{p23} for the extended B-spline basis functions.
\end{remark}

We now will prove a more general result which \textbf{does not assume that
}$G_{\alpha}$\textbf{\ has bounded support}. Instead we assume that
$G_{\alpha}\in L^{1}$ and that $\int G_{\alpha}\left(  y\right)  T_{m}\left(
x-y\right)  dy$ defines a regular tempered distribution. To prove our theorem
we will use the following result in which we use the function notation - a
notational necessity - for the action of distributions:

\begin{lemma}
\label{Lem_Thm_W3.1_calc_basis}(Theorem 2.7.5 Vladimirov \cite{Vladimirov})
Suppose $f,g\in\mathcal{D}^{\prime}$ and $g$ has bounded support
($g\in\mathcal{E}^{\prime}$). Then the convolution $f\ast g$ exists and%
\begin{equation}
\left[  f\ast g,\phi\right]  =\left[  f\left(  x\right)  g\left(  y\right)
,\eta\left(  y\right)  \phi\left(  x+y\right)  \right]  ,\quad\phi\in
C_{0}^{\infty},\label{p58}%
\end{equation}

where $\eta$ is any test function equal to one in a neighborhood of
$\operatorname*{supp}g$.

\textbf{Note:} the distribution on the right side of \ref{p58} is called a
\textbf{tensor or direct product} and in the literature it is variously
denoted $f\otimes g$, $f\times g$ or $f\cdot g$ in `operator notation' and
when `function notation with action variables' is used we can write $f\left(
x\right)  g\left(  y\right)  $ or $f\left(  x\right)  \cdot g\left(  y\right)
$, for example. The tensor product is defined for $f,g\in\mathcal{D}^{\prime
}\left(  \mathbb{R}^{d}\right)  $ by%
\begin{equation}
\left[  f\otimes g,\phi\right]  =\left[  f\left(  x\right)  ,\left[  g\left(
y\right)  ,\phi\left(  x,y\right)  \right]  \right]  ,\quad\phi\in
C_{0}^{\infty}\left(  \mathbb{R}^{d}\times\mathbb{R}^{d}\right)  ,\label{p59}%
\end{equation}

and has the properties you would expect e.g. commutivity, associativity,
continuity. See for example Section 2.7 Vladimirov \cite{Vladimirov}.
\end{lemma}

\begin{theorem}
\label{Thm_W3.1_calc_basis}Suppose the weight function $w$ has property W3.1
for $\left\vert \alpha\right\vert =\theta$ and $\kappa$. Then inequality
\ref{p36} holds for the multi-index $\alpha$. Let $G_{\alpha}$ be the basis
function of order zero generated by the zero order weight function $w_{\alpha
}\left(  \xi\right)  =\xi^{\alpha}w\left(  \xi\right)  $ and assume
$G_{\alpha}\in L^{1}$. Suppose $T_{m}$\ is the thin-plate spline basis
function of integer order $m>d/2$ given in Lemma \ref{Lem_thin_plate_splin}.

We shall show that if $G_{\theta}$ is the basis function $\left(  \frac
{1}{\left\vert \cdot\right\vert ^{2\theta}w}\right)  ^{\vee}$ of order
$\theta$ generated by $w$, and $H_{\theta}$ is defined by the convolutions
\begin{equation}
H_{\theta}=\left\{
\begin{array}
[c]{ll}%
\frac{\left(  -1\right)  ^{\theta}}{e\left(  \theta-d/2\right)  }D^{2\alpha
}G_{\alpha}\ast T_{\theta}, & \theta>d/2,\\
\frac{\left(  -1\right)  ^{\frac{d+1}{2}}}{e\left(  1/2\right)  }\left\vert
D\right\vert ^{d+1-2\theta}D^{2\alpha}G_{\alpha}\ast T_{\frac{d+1}{2}}, &
\theta<d/2,\text{ }d\text{ }odd,\\
\frac{\left(  -1\right)  ^{\frac{d+2}{2}}}{e\left(  1\right)  }\left\vert
D\right\vert ^{d+2-2\theta}D^{2\alpha}G_{\alpha}\ast T_{\frac{d+2}{2}}, &
\theta\leq d/2,\text{ }d\text{ }even,
\end{array}
\right. \label{p39}%
\end{equation}

where
\begin{equation}
\left(  G_{\alpha}\ast T_{m}\right)  \left(  x\right)  =\left(  2\pi\right)
^{-d/2}\int G_{\alpha}\left(  y\right)  T_{m}\left(  x-y\right)  dy,\quad
m>d/2,\label{p42}%
\end{equation}

then $G_{\theta}=H_{\theta}$ as tempered distributions whenever the
convolution integral defines a regular tempered distribution. The values of
the function $e\left(  \cdot\right)  $ are given in \ref{Rem_Lem_e(s)_pos}.
\end{theorem}

\begin{proof}
In this proof we will refer to Lemma \ref{Lem_W3.1_to_W2_zero} for the
properties of the zero order weight and basis functions used in this theorem.
Since $G_{\alpha}$ no longer has bounded support we will use a partition of
unity $\left\{  \psi_{\varepsilon},1-\psi_{\varepsilon}\right\}  $: choose
$\psi\in C_{0}^{\infty}$ so that $\operatorname*{supp}\psi\subset B\left(
0;2\right)  $, $0\leq\psi\leq1$ and $\psi=1$ on $B\left(  0;1\right)  $. Now
define $\psi_{\varepsilon}\left(  \xi\right)  =\psi\left(  \varepsilon
\xi\right)  $ for $\varepsilon>0$.\medskip

\fbox{\textbf{Case} $2\theta>d$} Since $\psi_{\varepsilon}G_{\alpha}$ has
bounded support, $\widehat{\psi_{\varepsilon}G_{\alpha}}\in C_{BP}^{\infty}$
and if $\xi$ is the \textit{action} variable%
\begin{align*}
\left[  \widehat{G_{\theta}},\phi\right]  =\left[  \frac{1}{\left\vert
\cdot\right\vert ^{2\theta}w},\phi\right]   &  =\left[  \frac{1}{\left\vert
\cdot\right\vert ^{2\theta}w_{\alpha}},\xi^{2\alpha}\phi\right] \\
&  =\left[  \frac{\widehat{G_{\alpha}}}{\left\vert \cdot\right\vert ^{2\theta
}},\xi^{2\alpha}\phi\right] \\
&  =\left[  \frac{\widehat{\psi_{\varepsilon}G_{\alpha}}}{\left\vert
\cdot\right\vert ^{2\theta}},\xi^{2\alpha}\phi\right]  +\left[  \frac{\left(
\left(  1-\psi_{\varepsilon}\right)  G_{\alpha}\right)  ^{\wedge}}{\left\vert
\cdot\right\vert ^{2\theta}},\xi^{2\alpha}\phi\right] \\
&  =\left[  \frac{1}{\left\vert \cdot\right\vert ^{2\theta}},\widehat
{\psi_{\varepsilon}G_{\alpha}}\xi^{2\alpha}\phi\right]  +\left[  \frac
{\xi^{2\alpha}\left(  \left(  1-\psi_{\varepsilon}\right)  G_{\alpha}\right)
^{\wedge}}{\left\vert \cdot\right\vert ^{2\theta}},\phi\right]  .
\end{align*}

But Lemma \ref{Lem_thin_plate_splin} tells us that $\widehat{T_{\theta}}%
=\frac{c_{\theta}}{\left\vert \cdot\right\vert ^{2\theta}}=\frac{e\left(
\theta-d/2\right)  }{\left\vert \cdot\right\vert ^{2\theta}}$ on
$S_{\emptyset,2\theta}$, and since $\widehat{\psi_{\varepsilon}G_{\alpha}}%
\xi^{2\alpha}\phi\in S_{\emptyset,2\theta}$ when $\phi\in S$
\begin{align*}
\left[  \frac{1}{\left\vert \cdot\right\vert ^{2\theta}},\widehat
{\psi_{\varepsilon}G_{\alpha}}\xi^{2\alpha}\phi\right]  =\left[  \frac
{1}{e\left(  \theta-d/2\right)  }\widehat{T_{\theta}},\widehat{\psi
_{\varepsilon}G_{\alpha}}\xi^{2\alpha}\phi\right]   &  =\left[  \frac
{1}{e\left(  \theta-d/2\right)  }\xi^{2\alpha}\widehat{\psi_{\varepsilon
}G_{\alpha}}\text{ }\widehat{T_{\theta}},\phi\right] \\
&  =\left[  \frac{\left(  -1\right)  ^{\theta}}{e\left(  \theta-d/2\right)
}\left(  D^{2\alpha}\left(  \psi_{\varepsilon}G_{\alpha}\right)  \ast
T_{\theta}\right)  ^{\wedge},\phi\right]  ,
\end{align*}

so that%
\begin{equation}
\widehat{G_{\theta}}=\frac{\left(  -1\right)  ^{\theta}}{e\left(
\theta-d/2\right)  }\left(  D^{2\alpha}\left(  \psi_{\varepsilon}G_{\alpha
}\right)  \ast T_{\theta}\right)  ^{\wedge}+\frac{\xi^{2\alpha}\left(  \left(
1-\psi_{\varepsilon}\right)  G_{\alpha}\right)  ^{\wedge}}{\left\vert
\cdot\right\vert ^{2\theta}}.\label{p09}%
\end{equation}

We now want to show that the second term on the right converges to zero as a
tempered distribution. In fact%
\begin{align*}
\left\vert \left[  \xi^{2\alpha}\frac{\left(  \left(  1-\psi_{\varepsilon
}\right)  G_{\alpha}\right)  ^{\wedge}}{\left\vert \cdot\right\vert ^{2\theta
}},\phi\right]  \right\vert =\left\vert \left[  \xi^{2\alpha}\frac{\left(
\left(  1-\psi_{\varepsilon}\right)  G_{\alpha}\right)  ^{\wedge}}{\left\vert
\cdot\right\vert ^{2\theta}},\phi\right]  \right\vert  &  =\left\vert \int%
\xi^{2\alpha}\frac{\left(  \left(  1-\psi_{\varepsilon}\right)  G_{\alpha
}\right)  ^{\wedge}}{\left\vert \cdot\right\vert ^{2\theta}}\phi\right\vert \\
&  \leq\int\left\vert \cdot\right\vert ^{2\theta}\frac{\left\vert \left(
\left(  1-\psi_{\varepsilon}\right)  G_{\alpha}\right)  ^{\wedge}\right\vert
}{\left\vert \cdot\right\vert ^{2\theta}}\left\vert \phi\right\vert \\
&  =\int\left\vert \left(  \left(  1-\psi_{\varepsilon}\right)  G_{\alpha
}\right)  ^{\wedge}\right\vert \left\vert \phi\right\vert \\
&  \leq\left(  \int\left\vert \left(  \left(  1-\psi_{\varepsilon}\right)
G_{\alpha}\right)  ^{\wedge}\right\vert \right)  \left\Vert \phi\right\Vert
_{\infty},
\end{align*}

but since $G_{\alpha}\in L^{1}$ we have $\left\vert \left(  \left(
1-\psi_{\varepsilon}\right)  G_{\alpha}\right)  ^{\wedge}\right\vert \leq
\int_{\left\vert \cdot\right\vert \geq\frac{1}{\varepsilon}}\left\vert
G_{\alpha}\right\vert \rightarrow0$ as $\varepsilon\rightarrow0^{+}$. Thus as
a tempered distributions
\begin{equation}
\xi^{2\alpha}\frac{\left(  \left(  1-\psi_{\varepsilon}\right)  G_{\alpha
}\right)  ^{\wedge}}{\left\vert \cdot\right\vert ^{2\theta}}\text{ }as\text{
}\varepsilon\rightarrow0^{+}\label{p51}%
\end{equation}

and so \ref{p09} implies that $\frac{\left(  -1\right)  ^{\theta}}{e\left(
\theta-d/2\right)  }\left(  D^{2\alpha}\left(  \psi_{\varepsilon}G_{\alpha
}\right)  \ast T_{\theta}\right)  ^{\wedge}$ converges to $\widehat{G_{\theta
}}$ in the sense of tempered distribution i.e.
\begin{equation}
\frac{\left(  -1\right)  ^{\theta}}{e\left(  \theta-d/2\right)  }D^{2\alpha
}\left(  \psi_{\varepsilon}G_{\alpha}\right)  \ast T_{\theta}\rightarrow
G_{\theta}\text{ }as\text{ }\varepsilon\rightarrow0^{+}.\label{p35}%
\end{equation}

Since $\psi_{\varepsilon}G_{\alpha}$ has bounded support we can use equation
\ref{p58} of Lemma \ref{Lem_Thm_W3.1_calc_basis} to express the convolution in
terms of the direct product. Thus%
\[
\left(  2\pi\right)  ^{\frac{d}{2}}\left[  \left(  \psi_{\varepsilon}%
G_{\alpha}\right)  \ast T_{\theta},\phi\right]  =\left[  \left(
\psi_{\varepsilon}G_{\alpha}\right)  \left(  x\right)  T_{\theta}\left(
y\right)  ,\eta_{\varepsilon}\left(  x\right)  \phi\left(  x+y\right)
\right]  .
\]

for any $\eta_{\varepsilon}\in C_{0}^{\infty}$ equal to one in a neighborhood
of $\operatorname*{supp}\left(  \psi_{\varepsilon}G_{\alpha}\right)  $. Now%
\begin{align*}
\left[  \left(  \psi_{\varepsilon}G_{\alpha}\right)  \left(  x\right)
T_{\theta}\left(  y\right)  ,\eta_{\varepsilon}\left(  x\right)  \phi\left(
x+y\right)  \right]   & =\int\int\psi_{\varepsilon}\left(  x\right)
G_{\alpha}\left(  x\right)  T_{\theta}\left(  y\right)  \phi\left(
x+y\right)  dydx\\
& =\int\int\psi_{\varepsilon}\left(  x^{\prime}-y^{\prime}\right)  G_{\alpha
}\left(  x^{\prime}-y^{\prime}\right)  T_{\theta}\left(  y^{\prime}\right)
\phi\left(  x^{\prime}\right)  dy^{\prime}dx^{\prime}\\
& =\int\int\psi_{\varepsilon}\left(  x^{\prime}-y^{\prime}\right)  G_{\alpha
}\left(  x^{\prime}-y^{\prime}\right)  T_{\theta}\left(  y^{\prime}\right)
dy^{\prime}\phi\left(  x^{\prime}\right)  dx^{\prime}\\
& =\int\int G_{\alpha}\left(  x^{\prime}-y^{\prime}\right)  T_{\theta}\left(
y^{\prime}\right)  dy^{\prime}\phi\left(  x^{\prime}\right)  dx^{\prime}-\\
& -\int\limits_{\mathbb{R}^{d}}\underset{\left\vert x^{\prime}-y^{\prime
}\right\vert \geq1/\varepsilon}{\int}\left(  1-\psi_{\varepsilon}\left(
x^{\prime}-y^{\prime}\right)  \right)  G_{\alpha}\left(  x^{\prime}-y^{\prime
}\right)  T_{\theta}\left(  y^{\prime}\right)  dy^{\prime}\phi\left(
x^{\prime}\right)  dx^{\prime}.
\end{align*}

Since we have assumed that $\int G_{\alpha}\left(  x^{\prime}-y^{\prime
}\right)  T_{\theta}\left(  y^{\prime}\right)  dy^{\prime}$ defines a regular
tempered distribution the first integral of the last equation exists and hence
the second integral also exists and so converges to zero as $\varepsilon
\rightarrow0$ since the region of integration lies outside the ball $B\left(
0;\frac{1}{\sqrt{d}\varepsilon}\right)  $. Thus $\left(  \psi_{\varepsilon
}G_{\alpha}\right)  \ast T_{\theta}$ converges to

$\left(  2\pi\right)  ^{-\frac{d}{2}}\int G_{\alpha}\left(  x^{\prime
}-y^{\prime}\right)  T_{\theta}\left(  y^{\prime}\right)  dy^{\prime}$ as
$\varepsilon\rightarrow0$. But by \ref{p35}, $H_{\theta}=\frac{\left(
-1\right)  ^{\theta}}{e\left(  \theta-d/2\right)  }D^{2\alpha}\left(
\psi_{\varepsilon}G_{\alpha}\right)  \ast T_{\theta}\rightarrow G_{\theta}$ as
$\varepsilon\rightarrow0^{+}$ and we can conclude that

$G_{\theta}=\left(  2\pi\right)  ^{-\frac{d}{2}}\frac{\left(  -1\right)
^{\theta}}{e\left(  \theta-d/2\right)  }D^{2\alpha}\int G_{\alpha}\left(
x^{\prime}-y^{\prime}\right)  T_{\theta}\left(  y^{\prime}\right)  dy^{\prime
}$.\medskip

\fbox{\textbf{Case} $2\theta<d$ and $d$ is odd} If $\xi$ is the
\textit{action} variable%
\begin{align}
\left[  \widehat{G_{\theta}},\phi\right]  =\left[  \frac{1}{\left\vert
\cdot\right\vert ^{2\theta}w},\phi\right]   &  =\left[  \frac{1}{\left\vert
\cdot\right\vert ^{2\left(  \frac{d+1}{2}\right)  }w},\left\vert
\cdot\right\vert ^{d+1-2\theta}\phi\right] \nonumber\\
&  =\left[  \frac{1}{\left\vert \cdot\right\vert ^{2\left(  \frac{d+1}%
{2}\right)  }w_{\alpha}},\xi^{2\alpha}\left\vert \cdot\right\vert
^{d+1-2\theta}\phi\right] \nonumber\\
&  =\left[  \frac{\widehat{G_{\alpha}}}{\left\vert \cdot\right\vert ^{2\left(
\frac{d+1}{2}\right)  }},\xi^{2\alpha}\left\vert \cdot\right\vert
^{d+1-2\theta}\phi\right] \nonumber\\
&  =\left[  \frac{\widehat{\psi_{\varepsilon}G_{\alpha}}}{\left\vert
\cdot\right\vert ^{2\left(  \frac{d+1}{2}\right)  }},\xi^{2\alpha}\left\vert
\cdot\right\vert ^{d+1-2\theta}\phi\right]  +\left[  \frac{\left(  \left(
1-\psi_{\varepsilon}\right)  G_{\alpha}\right)  ^{\wedge}}{\left\vert
\cdot\right\vert ^{2\left(  \frac{d+1}{2}\right)  }},\xi^{2\alpha}\left\vert
\cdot\right\vert ^{d+1-2\theta}\phi\right] \nonumber\\
&  =\left[  \frac{\widehat{\psi_{\varepsilon}G_{\alpha}}}{\left\vert
\cdot\right\vert ^{2\left(  \frac{d+1}{2}\right)  }},\xi^{2\alpha}\left\vert
\cdot\right\vert ^{d+1-2\theta}\phi\right]  +\left[  \frac{\xi^{2\alpha}%
}{\left\vert \cdot\right\vert ^{2\theta}}\left(  \left(  1-\psi_{\varepsilon
}\right)  G_{\alpha}\right)  ^{\wedge},\phi\right]  ,\label{p50}%
\end{align}

so that by \ref{p48} and \ref{p49}%
\begin{align*}
\left[  \frac{\widehat{\psi_{\varepsilon}G_{\alpha}}}{\left\vert
\cdot\right\vert ^{2\left(  \frac{d+1}{2}\right)  }},\xi^{2\alpha}\left\vert
\cdot\right\vert ^{d+1-2\theta}\phi\right]   & =\left[  \frac{1}{\left\vert
\cdot\right\vert ^{2\left(  \frac{d+1}{2}\right)  }},\widehat{\psi
_{\varepsilon}G_{\alpha}}\xi^{2\alpha}\left\vert \cdot\right\vert
^{d+1-2\theta}\phi\right] \\
& =\left[  \frac{1}{e\left(  1/2\right)  }\widehat{T_{\frac{d+1}{2}}}%
,\widehat{\psi_{\varepsilon}G_{\alpha}}\xi^{2\alpha}\left\vert \cdot
\right\vert ^{d+1-2\theta}\phi\right] \\
& =\left[  \frac{1}{e\left(  1/2\right)  }\xi^{2\alpha}\left\vert
\cdot\right\vert ^{d+1-2\theta}\widehat{\psi_{\varepsilon}G_{\alpha}}\text{
}\widehat{T_{\frac{d+1}{2}}},\phi\right] \\
& =\left[  \frac{1}{e\left(  1/2\right)  }\xi^{2\alpha}\left\vert
\cdot\right\vert ^{d+1-2\theta}\left(  \psi_{\varepsilon}G_{\alpha}\ast
T_{\frac{d+1}{2}}\right)  ^{\wedge},\phi\right] \\
& =\left[  \frac{\left(  -1\right)  ^{\frac{d+1}{2}}}{e\left(  1/2\right)
}\left(  D^{2\alpha}\left\vert D\right\vert ^{d+1-2\theta}\left(
\psi_{\varepsilon}G_{\alpha}\right)  \ast T_{\frac{d+1}{2}}\right)  ^{\wedge
},\phi\right]  .
\end{align*}

and \ref{p50} now becomes%
\[
\left[  \widehat{G_{\theta}},\phi\right]  =\left[  \frac{\left(  -1\right)
^{\frac{d+1}{2}}}{e\left(  1/2\right)  }\left(  D^{2\alpha}\left\vert
D\right\vert ^{d+1-2\theta}\left(  \psi_{\varepsilon}G_{\alpha}\right)  \ast
T_{\frac{d+1}{2}}\right)  ^{\wedge},\phi\right]  +\left[  \frac{\xi^{2\alpha}%
}{\left\vert \cdot\right\vert ^{2\theta}}\left(  \left(  1-\psi_{\varepsilon
}\right)  G_{\alpha}\right)  ^{\wedge},\phi\right]  ,
\]

and thus%
\[
\widehat{G_{\theta}}=\frac{\left(  -1\right)  ^{\frac{d+1}{2}}}{e\left(
1/2\right)  }\left(  D^{2\alpha}\left\vert D\right\vert ^{d+1-2\theta}\left(
\psi_{\varepsilon}G_{\alpha}\right)  \ast T_{\frac{d+1}{2}}\right)  ^{\wedge
}+\frac{\xi^{2\alpha}}{\left\vert \cdot\right\vert ^{2\theta}}\left(  \left(
1-\psi_{\varepsilon}\right)  G_{\alpha}\right)  ^{\wedge}.
\]

Comparing this equation with equation \ref{p09} of the previous case we see
that the subsequent calculations are also valid here and we can conclude that
\ref{p51} holds i.e. $\xi^{2\alpha}\frac{\widehat{\left(  1-\psi_{\varepsilon
}\right)  G_{\alpha}}}{\left\vert \cdot\right\vert ^{2\theta}}$ converges to
zero as a tempered distribution. Thus as tempered distributions%
\[
\frac{\left(  -1\right)  ^{\frac{d+1}{2}}}{e\left(  1/2\right)  }D^{2\alpha
}\left\vert D\right\vert ^{d+1-2\theta}\left(  \psi_{\varepsilon}G_{\alpha
}\right)  \ast T_{\frac{d+1}{2}}\rightarrow G_{\theta}\text{ }as\text{
}\varepsilon\rightarrow0^{+}.
\]

Finally, the calculations of the previous part following \ref{p35} are valid
for this case when $\theta=\left(  d+1\right)  /2$ and these complete the
proof of this case.\medskip

\fbox{\textbf{Case} $2\theta<d$ and $d$ is even} This case can be proved in a
very similar manner to the previous case.
\end{proof}

\begin{remark}
The condition $G_{\alpha}\in L^{1}$ of the last result can be used to deduce
information about $w_{\alpha}$ and $w$. In fact, Corollary 3.7 of Petersen
\cite{Petersen83} implies that
\begin{equation}
\frac{1}{w_{\alpha}\left(  \xi\right)  }=\left(  2\pi\right)  ^{-\frac{d}{2}%
}\int e^{i\xi x}G_{\alpha}\left(  x\right)  dx\label{p65}%
\end{equation}

holds a.e. and that $\frac{1}{w_{\alpha}}$ can be modified on a set of measure
zero so that it is a continuous, bounded function which converges to zero at
infinity and is such that \ref{p65} holds everywhere. Thus
\begin{equation}
\frac{1}{w_{\alpha}\left(  \xi\right)  }\leq\left(  2\pi\right)  ^{-\frac
{d}{2}}\left\Vert G_{\alpha}\right\Vert _{1}\label{p64}%
\end{equation}

and so $\lim\limits_{\left\vert \xi\right\vert \rightarrow\infty}w_{\alpha
}\left(  \xi\right)  =\infty$, $w_{\alpha}\in C^{\left(  0\right)  }\left(
\mathbb{R}^{d}\right)  $ and $w_{\alpha}$ is always positive which implies
that $w_{\alpha}$ has property W1 for the set $\mathcal{A}=\left\{
{}\right\}  $.

Further, $\frac{1}{w\left(  \xi\right)  }=\frac{\xi^{2\alpha}}{w_{\alpha
}\left(  \xi\right)  }$ so that by \ref{p64}, $w\left(  \xi\right)  \geq
\frac{\left(  2\pi\right)  ^{\frac{d}{2}}}{\left\Vert G_{\alpha}\right\Vert
_{1}}\frac{1}{\xi^{\alpha}}$ which implies $w$ is always positive and can only
have discontinuities on the axes $\xi_{k}=0$. Indeed, $w$ is discontinuous on
$\xi_{k}=0$ iff $\alpha_{k}=0$.
\end{remark}

The last theorem required that the convolution integral of \ref{p42} be a
regular tempered distribution and that $G_{\alpha}\in L^{1}$. We complete this
subsection by providing a single condition \ref{p44} for which those
requirements are satisfied.

\begin{corollary}
\label{Cor_Thm_W3.1_calc_basis}Suppose $G_{\alpha}$ and $T_{m}$ are as given
in Theorem \ref{Thm_W3.1_calc_basis} and suppose $G_{\alpha}$ also satisfies
the inequality%
\begin{equation}
\left\vert G_{\alpha}\left(  y\right)  \right\vert \leq c\left(  1+\left\vert
y\right\vert \right)  ^{-s},\label{p44}%
\end{equation}

for some $s>2m$ and constant $c$. Then $G_{\alpha}\in L^{1}$ and%
\begin{equation}
\left(  2\pi\right)  ^{-d/2}\int G_{\alpha}\left(  y\right)  T_{m}\left(
x-y\right)  dy,\quad m>d/2,\label{p43}%
\end{equation}

is absolutely convergent for all $x$ and defines a continuous function of
polynomial increase. It is also a regular tempered distribution (as defined in
Appendix \ref{SbSect_property_S'}).
\end{corollary}

\begin{proof}
Since $G_{\alpha}$ satisfies inequality \ref{p44} and $m>d/2$ it follows that
$s>d$ and so $G_{\alpha}\in L^{1}$. Further
\begin{align}
\left\vert \int G_{\alpha}\left(  y\right)  T_{m}\left(  x-y\right)
dy\right\vert  &  =\left\vert \int G_{\alpha}\left(  x-y\right)  T_{m}\left(
y\right)  dy\right\vert \nonumber\\
&  \leq\int\left\vert G_{\alpha}\left(  x-y\right)  T_{m}\left(  y\right)
\right\vert dy\nonumber\\
&  \leq c\int\frac{\left\vert T_{m}\left(  y\right)  \right\vert }{\left(
1+\left\vert x-y\right\vert \right)  ^{s}}dy\nonumber\\
&  \leq c\int\frac{\left(  1+\left\vert x\right\vert \right)  ^{s}}{\left(
1+\left\vert y\right\vert \right)  ^{s}}\left\vert T_{m}\left(  y\right)
\right\vert dy\nonumber\\
&  =c\left(  1+\left\vert x\right\vert \right)  ^{s}\int\frac{\left\vert
T_{m}\left(  y\right)  \right\vert }{\left(  1+\left\vert y\right\vert
\right)  ^{s}}dy\nonumber\\
&  =c\left(  1+\left\vert x\right\vert \right)  ^{s}\int\nolimits_{\left\vert
y\right\vert \leq1}\frac{\left\vert T_{m}\left(  y\right)  \right\vert
}{\left(  1+\left\vert y\right\vert \right)  ^{s}}dy+c\left(  1+\left\vert
x\right\vert \right)  ^{s}\int_{\left\vert y\right\vert \geq1}\frac{\left\vert
T_{m}\left(  y\right)  \right\vert }{\left(  1+\left\vert y\right\vert
\right)  ^{s}}dy.\label{p45}%
\end{align}

Since $T_{m}$ is continuous the first integral of the last term exists. From
\ref{p24}%
\[
T_{m}\left(  y\right)  =\left\{
\begin{array}
[c]{ll}%
\left(  -1\right)  ^{m-\left(  d-2\right)  /2}\left\vert y\right\vert
^{2m-d}\log\left\vert y\right\vert , & d\text{ }even,\\
\left(  -1\right)  ^{m-\left(  d-1\right)  /2}\left\vert y\right\vert
^{2m-d}, & d\text{ }odd.
\end{array}
\right.
\]

Thus if $d$ is odd then%
\[
\int_{\left\vert y\right\vert \geq1}\frac{\left\vert T_{m}\left(  y\right)
\right\vert }{\left(  1+\left\vert y\right\vert \right)  ^{s}}dy\leq
\int\limits_{\left\vert y\right\vert \geq1}\frac{\left\vert y\right\vert
^{2m-d}}{\left(  1+\left\vert y\right\vert \right)  ^{s}}dy\leq\int%
\limits_{\left\vert y\right\vert \geq1}\left\vert y\right\vert ^{\left(
2m-s\right)  -d}dy,
\]

which exists since $s>2m$. If $d$ is even and $s=2m+2\varepsilon$ then%
\begin{align*}
\int\limits_{\left\vert y\right\vert \geq1}\left(  1+\left\vert y\right\vert
\right)  ^{-s}\left\vert T_{m}\left(  y\right)  \right\vert dy\leq
\int\limits_{\left\vert y\right\vert \geq1}\left(  1+\left\vert y\right\vert
\right)  ^{-s}\left\vert y\right\vert ^{2m-d}\log\left\vert y\right\vert dy &
\leq\int\limits_{\left\vert y\right\vert \geq1}\left\vert y\right\vert
^{2m-s-d}\log\left\vert y\right\vert dy\\
&  =\int\limits_{\left\vert y\right\vert \geq1}\left\vert y\right\vert
^{-d-2\varepsilon}\log\left\vert y\right\vert dy\\
&  =\int\limits_{\left\vert y\right\vert \geq1}\frac{1}{\left\vert
y\right\vert ^{d+\varepsilon}}\frac{\log\left\vert y\right\vert }{\left\vert
y\right\vert ^{\varepsilon}}dy\\
&  \leq\sup_{r\geq1}\frac{\log r}{r^{\varepsilon}}\int\limits_{\left\vert
y\right\vert \geq1}\frac{1}{\left\vert y\right\vert ^{d+\varepsilon}}\\
&  <\infty.
\end{align*}

We now know that if $s>2m$ then for a constant $c^{\prime}$ independent of $x$%
\begin{equation}
\int\left\vert G_{\alpha}\left(  y\right)  T_{m}\left(  x-y\right)
\right\vert dy\leq c^{\prime}\left(  1+\left\vert x\right\vert \right)
^{s},\quad x\in\mathbb{R}^{d},\label{p54}%
\end{equation}

and so \ref{p43} is absolutely convergent and with polynomial growth at infinity.

To prove that $\int G_{\alpha}\left(  y\right)  T_{m}\left(  x-y\right)  dy$
is a continuous function of $x$ we use the Lebesgue-dominated convergence
theorem. First note that $\int G_{\alpha}\left(  y\right)  T_{m}\left(
x-y\right)  dy=\int G_{\alpha}\left(  x-y\right)  T_{m}\left(  y\right)  dy$
and that in the proof of \ref{p54} it was shown that
\[
\left\vert G_{\alpha}\left(  x-y\right)  T_{m}\left(  y\right)  \right\vert
\leq c\left(  1+\left\vert x\right\vert \right)  ^{s}\frac{\left\vert
T_{m}\left(  y\right)  \right\vert }{\left(  1+\left\vert y\right\vert
\right)  ^{s}},
\]

and that $\frac{T_{m}\left(  y\right)  }{\left(  1+\left\vert y\right\vert
\right)  ^{s}}\in L^{1}$. Thus if $x_{k}\rightarrow x$ is any sequence in the
ball $B\left(  x;1\right)  $ then%
\[
\left\vert G_{\alpha}\left(  x_{k}-y\right)  T_{m}\left(  y\right)
\right\vert \leq c\left(  1+\left\vert x_{k}-x\right\vert +\left\vert
x\right\vert \right)  ^{s}\frac{\left\vert T_{m}\left(  y\right)  \right\vert
}{\left(  1+\left\vert y\right\vert \right)  ^{s}}\leq c\left(  2+\left\vert
x\right\vert \right)  ^{s}\frac{\left\vert T_{m}\left(  y\right)  \right\vert
}{\left(  1+\left\vert y\right\vert \right)  ^{s}},
\]

and by the Lebesgue-dominated convergence theorem%
\[
\lim_{x_{k}\rightarrow x}\int G_{\alpha}\left(  y\right)  T_{m}\left(
x_{k}-y\right)  dy=\lim_{x_{k}\rightarrow x}\int G_{\alpha}\left(
x_{k}-y\right)  T_{m}\left(  y\right)  dy=\int G_{\alpha}\left(  x-y\right)
T_{m}\left(  y\right)  dy,
\]

and so proving continuity. A function $f\in L_{loc}^{1}$ for which
$\int\left\vert f\left(  x\right)  \right\vert \left(  1+\left\vert
x\right\vert \right)  ^{-\lambda}dx$ exists for some $\lambda\geq0$ is called
a regular tempered distribution. Clearly a continuous function of polynomial
growth at infinity is a regular tempered distribution.
\end{proof}

\subsection{Positive definite and conditionally positive definite basis
distributions}

In his work, Duchon did not mention the concept of conditional positive
definite basis functions, but used as his starting point a weight function.
Subsequent to Duchon's work various other classes of functions were used as
basis functions in numerical work and these were studied from the point of
view of (conditional) positive definiteness. This second major strand to the
development of the theory of basis function interpolation has involved the use
of the theory surrounding conditionally positive definite functions to develop
the correct setting for a variational approach. Work on this approach includes
papers by Schoenberg \cite{Schoenberg38}, Micchelli \cite{Micchelli86}, Madych
and Nelson \cite{MadychNelson88}, \cite{MadychNelson90}, Wu and Schaback
\cite{WuZongSchab93} and in a series of papers by Schaback which the reader
can find in the survey \cite{Schaback95}.

In their work Light and Wayne \cite{LightWayneX98Weight} returned to the
approach of Duchon and the concept of (conditional) positive definiteness
again plays a peripheral\textbf{\ }role. In the work of these authors it is
only after their basis distributions have been defined, and their properties
elucidated, that it is shown that they are positive definite or conditionally
positive definite. They showed in Theorem 4.3 of \cite{LightWayneX98Weight}
that the (integer) order of the basis distribution is the same as the order of
the conditional positive definiteness.

Following Chapter 4 of Light and Wayne \cite{LightWayneX98Weight}, the
objective of this section is to show that the basis distributions $G$ of order
$\theta$ are conditionally positive definite distributions of order $\theta$.
The definition of a conditionally positive definite tempered distribution
involves the concept of a homogeneous polynomial over the complex numbers
$\mathbb{C}$.

\begin{definition}
\label{Def_homog_poly_over_C}\textbf{Homogeneous polynomial over }$\mathbb{C}
$

A homogeneous polynomial $p$, of degree $\theta$ over $\mathbb{C}$, has the
form
\[
p\left(  x\right)  =\sum_{\left\vert \beta\right\vert =\theta}a_{\beta
}x^{\beta},\quad x\in\mathbb{R}^{d},\text{ }a_{\beta}\in\mathbb{C}.
\]

\end{definition}

Following Gelfand and Vilenkin \cite{GelfandVilenkin64} and Light and Wayne
\cite{LightWayneX98Weight} we define a conditionally strictly positive
definite tempered distribution as follows:

\begin{definition}
\label{Def_csposdef_S'}\textbf{Conditionally strictly positive definite
tempered distributions}

A distribution $F\in S^{\prime}$ is said to be conditionally strictly positive
definite of order $\theta\geq1$ if the inequality $\left[  p\bar{p}\widehat
{F},\psi\bar{\psi}\right]  >0$, holds for all $\psi\in S$, $\psi\neq0$ and all
homogeneous polynomials $p\neq0$, of order $\theta$.
\end{definition}

\begin{theorem}
\label{Thm_cpd_basis_fn}Assume $w$ is a weight function with property W2. Now
suppose $G\in S^{\prime}$ is a basis distribution of order $\theta$ generated
by the weight function $w$.

Then we show that $G$ is a conditionally positive definite tempered
distribution of order $\theta$.
\end{theorem}

\begin{proof}
Suppose $p=\sum\limits_{\left\vert \beta\right\vert <\theta}a_{\beta}x^{\beta
}$ is a homogeneous polynomial of order $\theta$ and $\psi\in S$, $\psi\neq0$.

Then $\left[  p\bar{p}\widehat{G},\psi\bar{\psi}\right]  =\left[  \widehat
{G},p\bar{p}\psi\bar{\psi}\right]  $ and
\[
p\bar{p}\psi\bar{\psi}=\left(  \sum_{\left\vert \alpha\right\vert =\theta
}a_{\alpha}x^{\alpha}\right)  \overline{\left(  \sum_{\left\vert
\beta\right\vert =\theta}a_{\beta}x^{\beta}\right)  }\psi\bar{\psi}%
=\sum_{\left\vert \alpha\right\vert =\theta}\sum_{\left\vert \beta\right\vert
=\theta}a_{\alpha}\overline{a_{\beta}}x^{\alpha+\beta}\psi\bar{\psi}.
\]

Now $\left\vert \alpha+\beta\right\vert =2\theta$ implies $x^{\alpha+\beta
}\psi\bar{\psi}\in S_{\emptyset,2\theta}$ and so $p\bar{p}\psi\bar{\psi}\in
S_{\emptyset,2\theta}$.

But $G$ has order $\theta$ so, $\left[  \widehat{G},\phi\right]  =\int%
\frac{\phi}{w\left\vert \cdot\right\vert ^{2\theta}}$ for all $\phi\in
S_{\emptyset,2\theta}$. Thus for all $\psi\in S$, $p\bar{p}\psi\bar{\psi}\in
S_{\emptyset,2\theta}$ and
\[
\left[  p\bar{p}\widehat{G},\psi\bar{\psi}\right]  =\left[  \widehat{G}%
,p\bar{p}\psi\bar{\psi}\right]  =\int\frac{\left(  p\bar{p}\psi\bar{\psi
}\right)  }{w\left\vert \cdot\right\vert ^{2\theta}}=\int\frac{\left\vert
p\psi\right\vert ^{2}}{w\left\vert \cdot\right\vert ^{2\theta}}.
\]
\medskip

Clearly $\int\dfrac{\left\vert p\psi\right\vert ^{2}}{w\left\vert
\cdot\right\vert ^{2\theta}}=0$ implies that $p\psi=0$ a.e. Does this imply
$p=0 $ or $\psi=0$? Certainly, in one dimension $p$ will have a finite number
of zeros, and so $\psi=0$. But what about in higher dimensions?

This result will be proved by induction on the number of dimensions. Suppose
the result is true for $d$ dimensions, and that $p\psi=0$ in $d+1$ dimensions
for some $p\neq0$. Thus $p\left(  x^{\prime},x_{d}\right)  \psi\left(
x^{\prime},x_{d}\right)  =0$ for all $x=\left(  x^{\prime},x_{d}\right)  $.
For fixed $x_{d}$, $p\left(  x^{\prime},x_{d}\right)  $ is a polynomial in
$d-1$ dimensions and $\psi\left(  x^{\prime},x_{d}\right)  \in S$ in $d-1$
dimensions. Thus, for each $x_{d}$, $\psi\left(  x^{\prime},x_{d}\right)  =0$
for all $x^{\prime}$ and so $\psi=0$.
\end{proof}

\chapter{More basis function and semi-Hilbert data space
theory\label{Ch_MoreBasisTheory}}

\section{Introduction}

This chapter is based on the theory developed in Chapter
\ref{Ch_weight_fn_exten} which extended the work of Light and Wayne
\cite{LightWayneX98Weight} The results given here are also closely related to
the work of Madych and Nelson \cite{MadychNelson90} as is indicated in the
remarks to several of our theorems. However, Madych and Nelson deal with the
\textit{reciprocal} of our weight function so that \textit{their}
$X_{w}^{\theta}$ space is formally $X_{1/w}^{-\theta}=\left\{  u\in S^{\prime
}:\int\frac{\left\vert \widehat{u}\right\vert ^{2}}{w\left\vert \cdot
\right\vert ^{2\theta}}<\infty\right\}  $ so that in order to calculate their
data space they need to construct a mapping from $X_{1/w}^{-\theta}$ to
$X_{w}^{\theta}$ Section 3. In a future document I will discuss the space
$X_{1/w}^{-\theta}$ and construct isometric isomorphisms from $X_{w}^{\theta}$
to $X_{1/w}^{-\theta}$ which show that $X_{1/w}^{-\theta}$ is isomorphic to
the space of bounded linear functionals on $X_{w}^{\theta}$ i.e. $\left(
X_{w}^{\theta}\right)  ^{\prime}=X_{1/w}^{-\theta}$.

In Chapter \ref{Ch_weight_fn_exten} smoothness and growth properties were
derived for basis functions and smoothness, growth and completeness properties
were exhibited for the data spaces $X_{w}^{\theta}$. It was shown that every
weight function generates a $C^{\left(  \left\lfloor 2\kappa\right\rfloor
\right)  }$ basis function where $\left\lfloor \cdot\right\rfloor $ denotes
the floor function. However, the rate of growth of basis functions near
infinity was only determined for a subclass of the weight functions. This
subclass generates bounded basis functions with bounded derivatives. As for
the functions in $X_{w}^{\theta}$, it was shown that $X_{w}^{\theta}\subset
C_{BP}^{\left(  \left\lfloor \kappa\right\rfloor \right)  }$ for all weight
functions but no bounds on growth rates near infinity were obtained.

The goal of this chapter is to derive `modified' inverse-Fourier transform
formulas for basis functions and the data functions $X_{w}^{\theta}$ and to
use these formulas to obtain bounds for the rates of increase of these
functions and their derivatives near infinity. This will be done by proving a
general inverse-Fourier transform formula for a subspace of the distributions
and then applying it to both the basis functions and the data functions. From
these formulas we will be able to show that all basis functions are either
bounded or have a rate of increase of at most $\left\vert \cdot\right\vert
^{\left\lfloor 2\theta\right\rfloor }$ near infinity and that all data
functions have a rate of increase of at most $\left\vert \cdot\right\vert
^{\left\lfloor \kappa\right\rfloor }$ near infinity.

We will also show that:

\begin{enumerate}
\item There always exists a conjugate-even (complex-valued) basis function.

\item If the weight function is even then there exists an even, real valued
basis function.

\item If the weight function is radial then there exists a radial basis function.
\end{enumerate}

The key operators used in this document are the projections $\mathcal{P}%
_{\emptyset,n}$ and $\mathcal{Q}_{\emptyset,n}$ which are given by
\[
\left(  \mathcal{P}_{\emptyset,n}u\right)  \left(  x\right)  =\rho\left(
x\right)  \sum_{\left\vert \alpha\right\vert <n}\frac{x^{\alpha}}{\alpha
!}D^{\alpha}u\left(  0\right)  ,\quad\mathcal{Q}_{\emptyset,n}=I-\mathcal{P}%
_{\emptyset,n},
\]

where $\rho\in S,$ $\rho\left(  0\right)  =1,$ $D^{\beta}\rho\left(  0\right)
=0$ for\ $1\leq|\beta|<n$. Here $S$ is the space of $C^{\infty}$ test
functions for the tempered distributions - see Appendix
\ref{Sect_tempered_distrib}. Thus they are based on the Taylor series
expansion about the origin.\medskip

\textbf{Section by section in brief:\smallskip}

In \textbf{Section }\ref{Sect_Pon_Qon} the null spaces and ranges of
$\mathcal{P}_{\emptyset,n}$ and $\mathcal{Q}_{\emptyset,n}$ are
derived.\smallskip

In \textbf{Section }\ref{Sect_adjoint_Pon_Qon} the distribution adjoints of
$\mathcal{P}_{\emptyset,n}^{\ast}$ and $\mathcal{Q}_{\emptyset,n}^{\ast}$ are
calculated. These are projections and we calculate their null spaces and
ranges.\smallskip

In \textbf{Section} \ref{Sect_Taylor_g_n} We derive a tempered distribution
Taylor series expansion for several spaces of continuous functions.

In \textbf{Section }\ref{Sect_Qon(exp(i(x,z)))} we prove some upper bounds for
the function $\left\vert D_{x}^{\gamma}\mathcal{Q}_{\emptyset,n,\xi}\left(
e^{i\left(  x,\xi\right)  }\right)  \right\vert $.\smallskip

In \textbf{Section }\ref{Sect_invers_Fourier} an inverse Fourier transform
result is derived for a class of distributions which will be applied in the
next two sections.\smallskip

In \textbf{Section }\ref{Sect_explicit_basis} we prove an inverse Fourier
transform formula for the basis functions. This is formula is used to derive
the rates of increase of the basis functions near infinity and to derive the
properties of a basis function given various properties of the weight
function.\smallskip

In \textbf{Section }\ref{Sect_another_inv_Four_X} an inverse Fourier transform
theorem is proved for functions in $X_{w}^{\theta}$ and then used to estimate
their rates of increase near infinity.\smallskip

In \textbf{Section} \ref{Sect_data_fn_Taylor_W3.2} Taylor series expansions of
(data) functions in $X_{w}^{\theta}$ when $w\in W3.2$.\smallskip

In \textbf{Section} \ref{Sect_data_fn_Taylor_W3.1} Taylor series expansions of
(data) functions in $X_{w}^{\theta}$ when $w\in W3.1$.

\section{The operators $\mathcal{P}_{\emptyset,n}$ and $\mathcal{Q}%
_{\emptyset,n}=I-\mathcal{P}_{\emptyset,n}$\label{Sect_Pon_Qon}}

The space $S_{1,n}$ is required to define the operators $\mathcal{P}%
_{\emptyset,n}$ and $\mathcal{Q}_{\emptyset,n}$ which are central to this document:

\begin{definition}
\label{Def_S1,n}\textbf{The spaces} $S_{1,n}$
\[
S_{1,n}=\left\{
\begin{array}
[c]{ll}%
S, & n=0,\\
\left\{  \phi\in S:\phi\left(  0\right)  =1;\text{ }D^{\beta}\phi\left(
0\right)  =0,\text{\quad}1\leq|\beta|<n\right\}  , & n=1,2,3,\ldots
\end{array}
\right.
\]

\end{definition}

\begin{definition}
\label{Def_Projections_Po_Qo}\textbf{The space }$\rho P_{n}$\textbf{\ and the
operators }$\mathcal{P}_{\emptyset,n}$\textbf{,\ }$\mathcal{Q}_{\emptyset,n}$.

Suppose $\rho\in S_{1,n}$. Then for integer $n\geq0$ we define:

\begin{enumerate}
\item The space $\rho P_{n}=\left\{  \rho p:p\in P_{n}\right\}  $.

\item The mappings $\mathcal{P}_{\emptyset,n}:C^{\left(  n\right)
}\rightarrow C^{\left(  n\right)  }$ and $\mathcal{Q}_{\emptyset
,n}=I-\mathcal{P}_{\emptyset,n}:C^{\left(  n\right)  }\rightarrow C^{\left(
n\right)  }$ are given by
\begin{align*}
\left(  \mathcal{P}_{\emptyset,n}u\right)  \left(  x\right)   &  =\rho\left(
x\right)  \sum_{\left\vert \alpha\right\vert <n}\frac{x^{\alpha}}{\alpha
!}D^{\alpha}u\left(  0\right)  =\rho\left(  x\right)  \sum_{k=0}^{n-1}\left(
\left(  xD_{y}\right)  ^{k}u\right)  \left(  y=0\right)  .\\
\left(  \mathcal{Q}_{\emptyset,n}u\right)  \left(  x\right)   &  =u\left(
x\right)  -\rho\left(  x\right)  \sum_{\left\vert \alpha\right\vert <n}%
\frac{x^{\alpha}}{\alpha!}D^{\alpha}u\left(  0\right)  .
\end{align*}

\end{enumerate}
\end{definition}

The next theorem proves some important properties of the operators
$\mathcal{P}_{\emptyset,n}$\textbf{\ }and\textbf{\ }$\mathcal{Q}_{\emptyset
,n}$ as well as demonstrating relationships between $S$ and $S_{\emptyset,n}$.

\begin{theorem}
\label{Thm_S=So+pP2m-1}Suppose $\rho\in S_{1,n}$. Then:

\begin{enumerate}
\item $\mathcal{P}_{\emptyset,n}$ and $\mathcal{Q}_{\emptyset,n}$ are
continuous linear mappings from $S$ into $S$.

\item $\mathcal{P}_{\emptyset,n}$ and $\mathcal{Q}_{\emptyset,n}$ are
projections into $S$ which satisfy:

$%
\begin{array}
[c]{lll}%
\mathcal{P}_{\emptyset,n}:S\rightarrow\rho P_{n} & is\text{ }onto\text{ }and &
\operatorname*{null}\mathcal{P}_{\emptyset,n}=S_{\emptyset,n},\\
\mathcal{Q}_{\emptyset,n}:S\rightarrow S_{\emptyset,n} & is\text{ }onto\text{
}and & \operatorname*{null}\mathcal{Q}_{\emptyset,n}=\rho P_{n}.
\end{array}
$

\item $S=S_{\emptyset,n}\oplus\rho P_{n}$.

\item $S=\widehat{S_{\emptyset,n}}\oplus\widehat{\rho P_{n}}=\widehat
{S_{\emptyset,n}}\oplus\left\{  p\left(  D\right)  \widehat{\rho}:p\in
P_{n}\right\}  $.
\end{enumerate}
\end{theorem}

\begin{proof}
\textbf{Part 1} It needs to be shown that if $\phi_{k}\rightarrow0$ in $S$
then $\mathcal{P}_{\emptyset,n}\phi_{k}\rightarrow0$ in $S$

i.e. $x^{\alpha}D^{\beta}\left(  \mathcal{P}_{\emptyset,n}\phi_{k}\right)
\left(  x\right)  \rightarrow0$ for all $x,\alpha,\beta$. Since
\[
\mathcal{P}_{\emptyset,n}\phi_{k}\left(  x\right)  =\rho\left(  x\right)
\sum_{\left\vert \alpha\right\vert <n}\frac{x^{\alpha}}{\alpha!}\left(
D^{\alpha}\phi_{k}\right)  \left(  0\right)  =\rho\left(  x\right)
\sum_{\left\vert \alpha\right\vert <n}\frac{\left[  D^{\alpha}\delta,\phi
_{k}\right]  }{\alpha!}x^{\alpha},
\]

we have using Leibniz's rule \ref{p68} in the Appendix.%
\begin{align*}
D^{\beta}\left(  \mathcal{P}_{\emptyset,n}\phi_{k}\right)  \left(  x\right)
& =D^{\beta}\left(  \rho\left(  x\right)  \sum_{\left\vert \alpha\right\vert
<n}\frac{\left[  D^{\alpha}\delta,\phi_{k}\right]  }{\alpha!}x^{\alpha}\right)
\\
& =\sum_{\gamma\leq\beta}\binom{\beta}{\gamma}\left(  D^{\beta-\gamma}%
\rho\right)  \left(  x\right)  D^{\gamma}\left(  \sum_{\left\vert
\alpha\right\vert <n}\frac{\left[  D^{\alpha}\delta,\phi_{k}\right]  }%
{\alpha!}x^{\alpha}\right) \\
& =\sum_{\gamma\leq\beta}\binom{\beta}{\gamma}\left(  D^{\beta-\gamma}%
\rho\right)  \left(  x\right)  \left(  \sum_{\left\vert \alpha\right\vert
<n}\frac{\left[  D^{\alpha}\delta,\phi_{k}\right]  }{\alpha!}D^{\gamma
}x^{\alpha}\right) \\
& =\sum_{\gamma\leq\beta}\binom{\beta}{\gamma}\left(  D^{\beta-\gamma}%
\rho\right)  \left(  x\right)  \left(  \sum_{\substack{\left\vert
\alpha\right\vert <n \\\alpha\leq\gamma}}\frac{\left[  D^{\alpha}\delta
,\phi_{k}\right]  }{\alpha!}D^{\gamma}x^{\alpha}\right)  .
\end{align*}

If $\phi_{k}\rightarrow0$ in $S$ then $\left[  D^{\alpha}\delta,\phi
_{k}\right]  \rightarrow0$ since $D^{\alpha}\delta\in S^{\prime}$. Hence
$D^{\beta}\left(  \mathcal{P}_{\emptyset,n}\phi_{k}\right)  \left(  x\right)
\rightarrow0$ for all $x$ and consequently $x^{\alpha}D^{\beta}\left(
\mathcal{P}_{\emptyset,n}\phi_{k}\right)  \left(  x\right)  \rightarrow0$ for
all $x$, $\alpha$, $\beta$.\smallskip

\textbf{Parts 2, 3}%
\begin{align*}
\mathcal{P}_{\emptyset,n}\mathcal{P}_{\emptyset,n}\phi=\rho\left(  x\right)
\sum_{\left\vert \beta\right\vert <n}\frac{x^{\beta}}{\beta!}\left(  D^{\beta
}\mathcal{P}_{\emptyset,n}\phi\right)  \left(  0\right)   &  =\rho\left(
x\right)  \sum_{\left\vert \beta\right\vert <n}\frac{x^{\beta}}{\beta
!}D^{\beta}\left(  \rho\left(  x\right)  \sum_{\left\vert \alpha\right\vert
<n}\frac{x^{\alpha}}{\alpha!}D^{\alpha}\phi\left(  0\right)  \right)  \left(
0\right) \\
&  =\rho\left(  x\right)  \sum_{\left\vert \beta\right\vert <n}\frac{x^{\beta
}}{\beta!}\rho\left(  x\right)  D^{\beta}\left(  \sum_{\left\vert
\alpha\right\vert <n}\frac{x^{\alpha}}{\alpha!}D^{\alpha}\phi\left(  0\right)
\right)  \left(  0\right)  ,
\end{align*}

since $\left(  D^{\alpha}\rho\right)  \left(  0\right)  =0$ for $1\leq
\left\vert \alpha\right\vert <n$. Further, $\rho\left(  0\right)  =1$ implies%
\begin{align*}
\mathcal{P}_{\emptyset,n}\mathcal{P}_{\emptyset,n}\phi=\rho\left(  x\right)
\sum_{\left\vert \beta\right\vert <n}\frac{x^{\beta}}{\beta!}D^{\beta}\left(
\sum_{\left\vert \alpha\right\vert <n}\frac{x^{\alpha}}{\alpha!}D^{\alpha}%
\phi\left(  0\right)  \right)  \left(  0\right)   &  =\rho\left(  x\right)
\sum_{\left\vert \beta\right\vert <n}\frac{x^{\beta}}{\beta!}\left(
\frac{D^{\beta}\left(  x^{\beta}\right)  }{\beta!}D^{\beta}\phi\left(
0\right)  \right)  \left(  0\right) \\
&  =\rho\left(  x\right)  \sum_{\left\vert \beta\right\vert <n}\frac{x^{\beta
}}{\beta!}D^{\beta}\phi\left(  0\right) \\
&  =\mathcal{P}_{\emptyset,n}\phi\left(  x\right)  ,
\end{align*}

and so $\mathcal{P}_{\emptyset,n}$ is a projection. Hence $\mathcal{Q}%
_{\emptyset,n}$ is also a projection, $S=\operatorname*{null}\mathcal{Q}%
_{\emptyset,n}\oplus\operatorname*{range}\mathcal{P}_{\emptyset,n}$ and
clearly $\operatorname*{range}\mathcal{P}_{\emptyset,n}=\rho P_{n}$.

Now to show that $\operatorname*{null}\mathcal{P}_{\emptyset,n}=S_{\emptyset
,n}$. Since $\rho\left(  0\right)  =1$ there exists $r>0$ such that
$\rho\left(  x\right)  >0$ for $\left\vert x\right\vert <r$. Hence if
$\mathcal{P}_{\emptyset,n}\phi=0$
\[
\rho\left(  x\right)  \sum_{\left\vert \alpha\right\vert <n}\frac{x^{\alpha}%
}{\alpha!}D^{\alpha}\phi\left(  0\right)  =\phi\left(  x\right)  ,\text{\quad
}\left\vert x\right\vert <r,
\]

and so $\sum\limits_{\left\vert \alpha\right\vert <n}\dfrac{x^{\alpha}}%
{\alpha!}D^{\alpha}\phi\left(  0\right)  =0$ for $\left\vert x\right\vert <r$.
Thus $D^{\alpha}\phi\left(  0\right)  =0$ if $\left\vert \alpha\right\vert
<n$, and so $\phi\in S_{\emptyset,n}$. Conversely, $\phi\in S_{\emptyset,n}$
clearly implies $\mathcal{P}_{\emptyset,n}\phi=0$ and since $\mathcal{P}%
_{\emptyset,n}=I-\mathcal{Q}_{\emptyset,n}$ we have the results
\begin{align*}
\operatorname*{range}\mathcal{Q}_{\emptyset,n}  & =\operatorname*{null}%
\mathcal{P}_{\emptyset,n}=S_{\emptyset,n},\\
\operatorname*{null}\mathcal{Q}_{\emptyset,n}  & =\operatorname*{range}%
\mathcal{P}_{\emptyset,n}=\rho P_{n}.
\end{align*}
\smallskip

\textbf{Part 4} From part 3, $S=S_{\emptyset,n}\oplus\rho P_{n}$. Hence
\[
S=\widehat{S}=\widehat{S_{\emptyset,n}}\oplus\widehat{\rho P_{n}}%
=\widehat{S_{\emptyset,n}}\oplus\left\{  p\left(  D\right)  \widehat{\rho
}:p\in P_{n}\right\}  .
\]

\end{proof}

\section{The tempered distribution adjoints of $\mathcal{P}_{\emptyset,n}$ and
$\mathcal{Q}_{\emptyset,n}$\label{Sect_adjoint_Pon_Qon}}

The projections $\mathcal{P}_{\emptyset,n}$ and $\mathcal{Q}_{\emptyset,n}$
were studied in the previous section and here we will derive their tempered
distribution adjoints $\mathcal{P}_{\emptyset,n}^{\ast}$ and $\mathcal{Q}%
_{\emptyset,n}^{\ast}$. These are projections into $S^{\prime}$ and several
properties are proved which relate to their ranges and null spaces.

\begin{definition}
\label{Def_P*o_Q*o}\textbf{The distribution adjoints }$\mathcal{P}%
_{\emptyset,n}^{\ast}$\textbf{\ and }$\mathcal{Q}_{\emptyset,n}^{\ast}$

The adjoints of $\mathcal{P}_{\emptyset,n},\mathcal{Q}_{\emptyset
,n}:S\rightarrow S$ are denoted by $\mathcal{P}_{\emptyset,n}^{\ast
},\mathcal{Q}_{\emptyset,n}^{\ast}$ and are defined for $u\in S^{\prime}$ and
$\varphi\in S$ by
\begin{align*}
\left[  \mathcal{Q}_{\emptyset,n}^{\ast}u,\varphi\right]   & =\left[
u,\mathcal{Q}_{\emptyset,n}\varphi\right]  ,\\
\left[  \mathcal{P}_{\emptyset,n}^{\ast}u,\varphi\right]   & =\left[
u,\mathcal{P}_{\emptyset,n}\varphi\right]  ,
\end{align*}

so that $\mathcal{P}_{\emptyset,n}^{\ast}:S^{\prime}\rightarrow S^{\prime}$
and $\mathcal{Q}_{\emptyset,n}^{\ast}:S^{\prime}\rightarrow S^{\prime} $ are continuous.
\end{definition}

\begin{theorem}
\label{Thm_P*_and_Q*}The operators $\mathcal{P}_{\emptyset,n}^{\ast}$ and
$\mathcal{Q}_{\emptyset,n}^{\ast}$ have the following properties:

\begin{enumerate}
\item $\mathcal{P}_{\emptyset,n}^{\ast}$ and $\mathcal{Q}_{\emptyset,n}^{\ast
}$ are projections and $\mathcal{P}_{\emptyset,n}^{\ast}+\mathcal{Q}%
_{\emptyset,n}^{\ast}=I$.

\item If $u\in S^{\prime}$ then $\mathcal{P}_{\emptyset,n}^{\ast}%
u=p_{n-1;u}\left(  -iD\right)  \delta$ where
\begin{equation}
p_{n-1;u}\left(  \xi\right)  =\sum_{\left\vert \alpha\right\vert <n}%
\dfrac{b_{\alpha,u}}{\alpha!}\xi^{\alpha}\in P_{n-1},\quad b_{\alpha
,u}=\left[  u,\left(  -ix\right)  ^{\alpha}\rho\right]  .\label{a91}%
\end{equation}

\item $\operatorname*{null}\mathcal{Q}_{\emptyset,n}^{\ast}%
=\operatorname*{range}\mathcal{P}_{\emptyset,n}^{\ast}=\left\{  p\left(
D\right)  \delta:p\in P_{n-1}\right\}  =\widehat{P}_{n-1}$.

\item $\operatorname*{range}\mathcal{Q}_{\emptyset,n}^{\ast}%
=\operatorname*{null}\mathcal{P}_{\emptyset,n}^{\ast}=\left\{  u\in S^{\prime
}:\left[  u,x^{\alpha}\rho\right]  =0\text{ }when\text{ }\left\vert
\alpha\right\vert <n\right\}  $.

\item $S^{\prime}=\widehat{P}_{n-1}\oplus\left\{  u\in S^{\prime}:\left[
u,x^{\alpha}\rho\right]  =0\text{ }when\text{ }\left\vert \alpha\right\vert
<n\right\}  $.

\item $\varphi u=\varphi\mathcal{Q}_{\emptyset,n}^{\ast}u$, where $u\in
S^{\prime}$ and $\varphi\in C_{\emptyset,n}^{\infty}\cap C_{BP}^{\left(
0\right)  }$.
\end{enumerate}
\end{theorem}

\begin{proof}
\textbf{Part 1} is true since $\mathcal{P}_{\emptyset,n}$ and $\mathcal{Q}%
_{\emptyset,n}$ are projections.

\textbf{Part 2} If $u\in S^{\prime}$ and $\varphi\in S$ then
\begin{align*}
\left[  u,\mathcal{P}_{\emptyset,n}\varphi\right]  =\left[  u,\rho
\sum_{\left\vert \alpha\right\vert <n}\frac{x^{\alpha}}{\alpha!}D^{\alpha
}\varphi\left(  0\right)  \right]   &  =\sum_{\left\vert \alpha\right\vert
<n}\frac{\left[  u,x^{\alpha}\rho\right]  }{\alpha!}D^{\alpha}\varphi\left(
0\right) \\
&  =\sum_{\left\vert \alpha\right\vert <n}\frac{\left[  u,x^{\alpha}%
\rho\right]  }{\alpha!}\left[  \delta,D^{\alpha}\varphi\right] \\
&  =\sum_{\left\vert \alpha\right\vert <n}\frac{\left[  u,x^{\alpha}%
\rho\right]  }{\alpha!}\left(  -1\right)  ^{\left\vert \alpha\right\vert
}\left[  D^{\alpha}\delta,\varphi\right] \\
&  =\left[  \left(  \sum_{\left\vert \alpha\right\vert <n}\frac{\left[
u,x^{\alpha}\rho\right]  }{\alpha!}\left(  -1\right)  ^{\left\vert
\alpha\right\vert }D^{\alpha}\right)  \delta,\varphi\right] \\
&  =\left[  \left(  \sum_{\left\vert \alpha\right\vert <n}\frac{\left[
u,\left(  -ix\right)  ^{\alpha}\rho\right]  }{\alpha!}\left(  -iD\right)
^{\alpha}\right)  \delta,\varphi\right] \\
&  =\left[  \mathcal{P}_{\emptyset,n}^{\ast}u,\varphi\right]  ,
\end{align*}

and thus
\[
\mathcal{P}_{\emptyset,n}^{\ast}u=\left(  \sum_{\left\vert \alpha\right\vert
<n}\frac{b_{\alpha,u}}{\alpha!}\left(  -iD\right)  ^{\alpha}\right)  \delta\in
S^{\prime}.
\]
\medskip

\textbf{Part 3} From the formulas of part 2 it is clear that $\mathcal{P}%
_{\emptyset,n}^{\ast}:S^{\prime}\rightarrow\widehat{P}_{n-1}$ so that
$\operatorname*{range}\mathcal{P}_{\emptyset,n}^{\ast}\subset\widehat{P}%
_{n-1}$. We now show that $\operatorname*{null}\mathcal{Q}_{\emptyset,n}%
^{\ast}=\widehat{P}_{n-1}$. If $u\in\widehat{P}_{n-1}$ then $u=p(D)\delta$ for
some $p\in P_{n-1}$ and, since $\varphi\in S$ implies $\mathcal{Q}%
_{\emptyset,n}\varphi\in S_{\emptyset,n}$,%
\begin{align*}
\left[  \mathcal{Q}_{\emptyset,n}^{\ast}u,\varphi\right]  =\left[
u,\mathcal{Q}_{\emptyset,n}\varphi\right]  =\left[  p\left(  D\right)
\delta,\mathcal{Q}_{\emptyset,n}\varphi\right]   &  =\left[  \delta,p\left(
-D\right)  \left(  \mathcal{Q}_{\emptyset,n}\varphi\right)  \right] \\
&  =\left[  p\left(  -D\right)  \left(  \mathcal{Q}_{\emptyset,n}%
\varphi\right)  \right]  \left(  0\right) \\
&  =0,
\end{align*}

by the definition of $S_{\emptyset,n}$. This means $\mathcal{Q}_{\emptyset
,n}^{\ast}u=0$ and hence that $\widehat{P}_{n-1}\subset\operatorname*{null}%
\mathcal{Q}_{\emptyset,n}^{\ast}$.\medskip

\textbf{Part 4} From part 2, $\mathcal{P}_{\emptyset,n}^{\ast}u=p_{u}\left(
D\right)  \delta$ so $\mathcal{P}_{\emptyset,n}^{\ast}u=0$ iff $\widehat
{\mathcal{P}_{\emptyset,n}^{\ast}u}=0$ iff $p_{u}=0$ iff $\left[  u,x^{\alpha
}\rho\right]  =0$ when $\left\vert \alpha\right\vert <n$. Here $x$ was an
\textit{action} variable - see Notation \ref{Not_distribution}.\medskip

\textbf{Part 5} Part 1 implies that $S^{\prime}=\operatorname*{null}%
\mathcal{Q}_{\emptyset,n}^{\ast}\oplus\operatorname*{range}\mathcal{Q}%
_{\emptyset,n}^{\ast}$.\medskip

\textbf{Part 6} Suppose $u\in S^{^{\prime}}$ and $\varphi\in C_{\emptyset
,n}^{\infty}\cap C_{BP}^{\left(  0\right)  }$. Then $\psi\in S$ implies
$\varphi\psi\in S_{\emptyset,n}$ and so%
\[
\left[  \varphi u,\psi\right]  =\left[  u,\varphi\psi\right]  =\left[
u,\mathcal{Q}_{\emptyset,n}\left(  \varphi\psi\right)  \right]  =\left[
\mathcal{Q}_{\emptyset,n}^{\ast}u,\varphi\psi\right]  =\left[  \varphi
\mathcal{Q}_{\emptyset,n}^{\ast}u,\psi\right]  .
\]

\end{proof}

This theorem is used to estimate the polynomial term in \ref{a93}.

\begin{theorem}
\label{Thm_Adel_poly_S1n}Suppose $a,b\in\mathbb{R}^{d}$ are constants. Then:

\begin{enumerate}
\item $aD_{x}\frac{\left(  bx\right)  ^{k}}{k!}=ab\frac{\left(  bx\right)
^{k-1}}{\left(  k-1\right)  !}$;

\item $\left(  aD_{x}\right)  ^{j}\frac{\left(  bx\right)  ^{k}}{k!}=\left(
ab\right)  ^{j}\frac{\left(  bx\right)  ^{k-j}}{\left(  k-j\right)  !}$;

\item Suppose $p_{n-1;u}$ is the polynomial described in \ref{a91}. Then
$\left(  aD\right)  ^{j}p_{n-1;u}=p_{n-j-1;v_{j}}$ where $v_{j}\left(
x\right)  =\left(  -iax\right)  ^{j}u\left(  x\right)  $.
\end{enumerate}
\end{theorem}

\begin{proof}
\textbf{Parts 1 and 2} Part 2 follows directly from part 1 and part 1 is an
simple calculation.\medskip

\textbf{Part 3} By definition%
\begin{align*}
p_{n-1;u}\left(  \xi\right)   & =\sum_{\left\vert \alpha\right\vert <n}%
\dfrac{\left[  u,\left(  -ix\right)  ^{\alpha}\rho\right]  }{\alpha!}%
\xi^{\alpha}=\left[  u_{x},\sum_{\left\vert \alpha\right\vert <n}\frac{\left(
-i\xi\right)  ^{\alpha}x^{\alpha}}{\alpha!}\rho\left(  x\right)  \right]  =\\
& =\left[  u_{x},\sum_{k=0}^{n-1}\sum_{\left\vert \alpha\right\vert =k}%
\frac{\left(  -i\xi\right)  ^{\alpha}x^{\alpha}}{\alpha!}\rho\left(  x\right)
\right]  =\sum_{k=0}^{n-1}\left[  u_{x},\frac{\left(  -i\xi x\right)  ^{k}%
}{k!}\rho\left(  x\right)  \right]  ,
\end{align*}

and thus for $j=1$, part 2 gives,%
\begin{align*}
aDp_{n-1;u}\left(  \xi\right)   & =\sum_{k=0}^{n-1}\left[  u\left(  x\right)
,aD_{\xi}\frac{\left(  -i\xi x\right)  ^{k}}{k!}\rho\left(  x\right)  \right]
=\sum_{k=1}^{n-1}\left[  u\left(  x\right)  ,-iax\frac{\left(  -i\xi x\right)
^{k-1}}{\left(  k-1\right)  !}\rho\left(  x\right)  \right]  =\\
& =\sum_{k=1}^{n-1}\left[  v_{1}\left(  x\right)  ,\frac{\left(  -i\xi
x\right)  ^{k-1}}{\left(  k-1\right)  !}\rho\left(  x\right)  \right]
=\sum_{k=0}^{n-2}\left[  v_{1}\left(  x\right)  ,\frac{\left(  -i\xi x\right)
^{k}}{k!}\rho\left(  x\right)  \right]  =\\
& =p_{n-2;v_{1}}\left(  \xi\right)  ,
\end{align*}

and several more applications of $aD$ proves this part.
\end{proof}

\section{Taylor series expansions and the function $g_{n}$%
\label{Sect_Taylor_g_n}}

We start by deriving the tempered distribution Taylor series expansion
\ref{a1.55} with Fourier transform remainder \ref{a50.5} by means of the
1-dimensional Taylor series expansion with integral remainder applied to
$e^{ia\xi}$. Indeed, for $x\in\mathbb{R}^{1}$,%
\begin{align*}
e^{x}=\sum_{k\leq n}\frac{x^{n}}{n!}+\left(  \mathcal{R}_{n+1}e^{x}\right)
\left(  0,x\right)   &  =\sum_{k\leq n}\frac{x^{k}}{k!}+\frac{x^{n+1}}{n!}%
\int_{0}^{1}\left(  1-t\right)  ^{n}\left(  D^{n+1}e^{s}\right)  \left(
tx\right)  dt\\
&  =\sum_{k\leq n}\frac{x^{k}}{k!}+\frac{x^{n+1}}{n!}\int_{0}^{1}\left(
1-t\right)  ^{n}e^{tx}dt,
\end{align*}

so that if we define $g_{n}$ by%
\begin{equation}
g_{n}\left(  t\right)  =\left\{
\begin{array}
[c]{ll}%
0, & t<0,\\
\left(  1-t\right)  ^{n}, & 0\leq t\leq1,\\
0, & t>1,
\end{array}
\right.  \text{\quad}n=0,1,2,\ldots,\label{a117}%
\end{equation}

it follows that%
\begin{align*}
e^{ix}=\sum_{k\leq n}\frac{x^{k}}{k!}+\frac{x^{n+1}}{n!}\int_{0}^{1}\left(
1-t\right)  ^{n}e^{itx}dt  & =\sum_{k\leq n}\frac{x^{k}}{k!}+\frac{x^{n+1}%
}{n!}\int_{0}^{1}g_{n}\left(  t\right)  e^{itx}dt\\
& =\sum_{k\leq n}\frac{x^{k}}{k!}+\frac{\sqrt{2\pi}}{n!}x^{n+1}\overline
{\widehat{g_{n}}}\left(  x\right) \\
& =\sum_{k\leq n}\frac{x^{k}}{k!}+\frac{\sqrt{2\pi}}{n!}x^{n+1}\overset{\vee
}{g_{n}}\left(  x\right)  ,
\end{align*}

and%
\[
\overset{\vee}{g_{n}}\left(  x\right)  =\frac{n!}{\sqrt{2\pi}}\frac
{e^{ix}-\sum\limits_{k\leq n}\frac{x^{k}}{k!}}{x^{n+1}}.
\]

Also, for $a,\xi\in\mathbb{R}^{d}$,%
\begin{align}
e^{ia\xi}  & =\sum_{k\leq n}\frac{\left(  ia\xi\right)  ^{k}}{k!}%
+\frac{\left(  ia\xi\right)  ^{n+1}}{n!}\int_{0}^{1}\left(  1-t\right)
^{n}e^{ita\xi}dt\nonumber\\
& =\sum_{k\leq n}\frac{\left(  ia\xi\right)  ^{k}}{k!}+\frac{\sqrt{2\pi}}%
{n!}\left(  ia\xi\right)  ^{n+1}\overline{\widehat{g_{n}}}\left(  a\xi\right)
.\label{a114}%
\end{align}

Clearly%
\[
g_{n}\left(  t\right)  =\left(  1-t\right)  ^{n-k}g_{k}\left(  t\right)  .
\]

More properties of $g_{n}$ and $\widehat{g_{n}}$ are:

\begin{lemma}
\label{Lem_gm_properties_2}The function $g_{n}$ given by \ref{a117} has the
following properties:

\begin{enumerate}
\item $\left\Vert g_{n}\right\Vert _{1}=\frac{1}{n+1}$ and $\left\Vert
g_{n}\right\Vert _{2}=\frac{1}{\sqrt{2n+1}}$.

\item $\widehat{g_{n}}\in C_{B}^{\infty}$ and%
\begin{equation}
\left\vert \widehat{g_{n}}\left(  t\right)  \right\vert \leq\frac{1}%
{\sqrt{2\pi}}\frac{1}{n+1},\quad n=0,1,2,\ldots\label{a50.3}%
\end{equation}

\item
\[
Dg_{n}=\left\{
\begin{array}
[c]{ll}%
\delta-\delta\left(  \cdot-1\right)  , & n=0,\\
-ng_{n-1}+\delta, & n=1,2,3\ldots,
\end{array}
\right.
\]

and%
\begin{equation}
\widehat{Dg_{n}}=\left\{
\begin{array}
[c]{ll}%
\frac{2i}{\sqrt{2\pi}}e^{-it/2}\sin\frac{t}{2}, & n=0,\\
-n\widehat{g_{n-1}}+\frac{1}{\sqrt{2\pi}}, & n=1,2,3\ldots.
\end{array}
\right. \label{a2.04}%
\end{equation}

\item If we define%
\begin{equation}
c_{n}=2+\frac{1}{n+1},\quad n=0,1,2,\ldots,\label{a2.08}%
\end{equation}

then
\begin{equation}
\left\vert \widehat{g_{n}}\left(  t\right)  \right\vert \leq\frac{c_{n}}%
{\sqrt{2\pi}}\frac{1}{1+\left\vert t\right\vert },\text{\quad}t\in
\mathbb{R}^{1},\text{ }n=0,1,2,\ldots.\label{a2.07}%
\end{equation}

\item For each $a\in\mathbb{R}^{d}$ such that $a.\neq\mathbf{0}$ we have
$\widehat{g_{n}}\left(  a\xi\right)  \in C_{B}^{\infty}\left(  \mathbb{R}%
^{d}\right)  $.

\item Define the \textbf{tilda operator} $\symbol{126}$ by $\widetilde
{u}\left(  x\right)  =\overline{u\left(  -x\right)  }$. Then%
\[
\left\Vert g_{n}\ast\widetilde{Dg_{n}}\right\Vert _{1}\leq\left\{
\begin{array}
[c]{ll}%
\frac{2}{\sqrt{2\pi}}, & n=0,\\
\left(  1+\frac{1}{\sqrt{2\pi}}\right)  \frac{1}{n+1}, & n\geq1.
\end{array}
\right.
\]

\item
\[
\left\vert \widehat{g_{n}}\left(  t\right)  \right\vert \leq\left\{
\begin{array}
[c]{ll}%
\frac{1}{\sqrt{2\pi}}\frac{1}{n+1}, & \left\vert t\right\vert \leq2\left(
n+1\right)  ,\\
\frac{1}{\sqrt{2\pi}}\left(  2+\frac{1}{n+1}\right)  \frac{1}{1+\left\vert
t\right\vert }, & \left\vert t\right\vert \geq2\left(  n+1\right)  .
\end{array}
\right.
\]

\end{enumerate}
\end{lemma}

\begin{proof}
\textbf{Part 1} A simple calculation.\smallskip

\textbf{Part 2} Since $g_{n}$ is a distribution with bounded support,
$\widehat{g_{n}}\in C_{BP}^{\infty}$ and the estimates \ref{a50.3} follow
directly from the formula, $\widehat{g_{n}}\left(  t\right)  =\frac{1}%
{\sqrt{2\pi}}\int_{0}^{1}e^{-ist}\left(  1-s\right)  ^{n}ds$.\smallskip

\textbf{Part 3} If $n=0$ then $Dg_{0}=\left\{  Dg_{0}\right\}  +\delta
-\delta\left(  \cdot-1\right)  =0+\delta-\delta\left(  \cdot-1\right)  $.

Hence $\widehat{Dg_{0}}=\widehat{\delta}-\widehat{\delta\left(  \cdot
-1\right)  }=\frac{1}{\sqrt{2\pi}}\left(  1-e^{-it}\right)  =\frac{2i}%
{\sqrt{2\pi}}e^{-it/2}\left(  \frac{e^{it/2}-e^{-it/2}}{2i}\right)  =\frac
{2i}{\sqrt{2\pi}}e^{-it/2}\sin\frac{t}{2}$.

If $n\geq1$ then $Dg_{n}=\left\{  Dg_{n}\right\}  +\delta=-n\left(
1-t\right)  ^{n-1}g_{0}+\delta=-ng_{n-1}+\delta$.\smallskip

\textbf{Part 4} Suppose $n\geq1$. Then from part 3 and then part 2:%
\begin{align*}
\left\vert t\right\vert \left\vert \widehat{g_{n}\left(  t\right)
}\right\vert =\left\vert \widehat{Dg_{n}}\left(  t\right)  \right\vert  &
\leq\left\{
\begin{array}
[c]{ll}%
\frac{2}{\sqrt{2\pi}}, & n=0,\\
n\left\vert \widehat{g_{n-1}}\left(  t\right)  \right\vert +\frac{1}%
{\sqrt{2\pi}}, & n\geq1,
\end{array}
\right. \\
& \leq\frac{2}{\sqrt{2\pi}},
\end{align*}

so that%
\[
\left(  1+\left\vert t\right\vert \right)  \left\vert \widehat{g_{n}\left(
t\right)  }\right\vert \leq\left\vert \widehat{g_{n}\left(  t\right)
}\right\vert +\frac{2}{\sqrt{2\pi}}\leq\frac{1}{\sqrt{2\pi}}\left(  2+\frac
{1}{n+1}\right)  =\frac{c_{n}}{\sqrt{2\pi}}.
\]

\textbf{Part 5} True since $\widehat{g_{n}}\in C_{B}^{\infty}\left(
\mathbb{R}^{1}\right)  $.\smallskip

\textbf{Part 6} From part 3,%
\[
Dg_{n}=\left\{
\begin{array}
[c]{ll}%
\delta-\delta\left(  \cdot-1\right)  , & n=0,\\
-ng_{n-1}+\delta, & n=1,2,3\ldots.
\end{array}
\right.
\]

Hence%
\begin{align*}
\left\Vert g_{0}\ast\widetilde{Dg_{0}}\right\Vert _{1}=\left\Vert g_{0}%
\ast\left(  \delta-\delta\left(  -\cdot-1\right)  \right)  \right\Vert _{1}  &
=\frac{1}{\sqrt{2\pi}}\left\Vert g_{0}-g_{0}\left(  \cdot+1\right)
\right\Vert _{1}\\
& =\frac{1}{\sqrt{2\pi}}\left(  \left\Vert g_{0}\right\Vert _{1}+\left\Vert
g_{0}\left(  \cdot+1\right)  \right\Vert _{1}\right) \\
& =\frac{1}{\sqrt{2\pi}}\left(  \left\Vert g_{0}\right\Vert _{1}+\left\Vert
g_{0}\left(  \cdot+1\right)  \right\Vert _{1}\right) \\
& =\frac{2}{\sqrt{2\pi}}\left\Vert g_{0}\right\Vert _{1}\\
& \leq\frac{2}{\sqrt{2\pi}},
\end{align*}

and if $n\geq1$,%
\begin{align*}
\left\Vert g_{n}\ast\widetilde{Dg_{n}}\right\Vert _{1}  & =\left\Vert
g_{n}\ast\left(  -n\widetilde{g}_{n-1}+\delta\right)  \right\Vert _{1}\\
& =\left\Vert -ng_{n}\ast\widetilde{g}_{n-1}+\frac{1}{\sqrt{2\pi}}%
g_{n}\right\Vert _{1}\\
& \leq n\left\Vert g_{n}\ast\widetilde{g}_{n-1}\right\Vert _{1}+\frac{1}%
{\sqrt{2\pi}}\left\Vert g_{n}\right\Vert _{1}\\
& \leq n\left\Vert g_{n}\right\Vert _{1}\left\Vert \widetilde{g}%
_{n-1}\right\Vert _{1}+\frac{1}{\sqrt{2\pi}}\left\Vert \widetilde{g}%
_{n}\right\Vert _{1}\\
& =n\left\Vert g_{n}\right\Vert _{1}\left\Vert g_{n-1}\right\Vert _{1}%
+\frac{1}{\sqrt{2\pi}}\left\Vert g_{n}\right\Vert _{1}\\
& \leq n\frac{1}{n\left(  n+1\right)  }+\frac{1}{\sqrt{2\pi}}\frac{1}{n+1}\\
& =\left(  1+\frac{1}{\sqrt{2\pi}}\right)  \frac{1}{n+1}.
\end{align*}
\smallskip

\textbf{Part 7} Regarding parts 2 and 4, $\frac{1}{\sqrt{2\pi}}\frac{1}%
{n+1}=\frac{c_{n}}{\sqrt{2\pi}}\frac{1}{1+\left\vert t\right\vert }=\frac
{1}{\sqrt{2\pi}}\left(  2+\frac{1}{n+1}\right)  \frac{1}{1+\left\vert
t\right\vert }$ when $\left\vert t\right\vert =2n+2$.
\end{proof}

\begin{remark}
From parts 1 and 4,%
\[
\left\vert \widehat{g_{n}}\left(  t\right)  \right\vert \leq\frac{1}%
{\sqrt{2\pi}}\frac{1}{n+1}\min\left\{  1,\frac{2n+3}{1+\left\vert t\right\vert
}\right\}  ,\quad t\in\mathbb{R}^{1},
\]

and the two bounds intersect when $t=2\left(  n+1\right)  $.
\end{remark}

Thus if $f\in S^{\prime}$ and $\xi$ is the \textbf{action} variable when
$a.\neq\mathbf{0}$ we can write \ref{a114} as%
\begin{equation}
\left(  e^{ia\xi}-\sum_{k\leq n}\frac{\left(  ia\xi\right)  ^{k}}{k!}\right)
\widehat{f}=\frac{\sqrt{2\pi}}{n!}\left(  ia\xi\right)  ^{n+1}\overline
{\widehat{g_{n}}}\left(  a\xi\right)  \widehat{f},\label{a2.11}%
\end{equation}

or on using the first of the multi-index identities \ref{p08},
\begin{equation}
\frac{\left(  a\xi\right)  ^{k}}{k!}=\sum_{\left\vert \beta\right\vert
=k}\frac{a^{\beta}\xi^{\beta}}{\beta!},\quad\frac{d^{k}}{k!}=\sum_{\left\vert
\beta\right\vert =k}\frac{1}{\beta!},\label{a2241}%
\end{equation}

we get%
\[
\left(  f\left(  \cdot+a\right)  -\sum_{k\leq n}\frac{\left(  aD\right)  ^{k}%
}{k!}f\right)  ^{\wedge}=\left(  f\left(  \cdot+a\right)  -\sum_{\left\vert
\beta\right\vert \leq n}\frac{a^{\beta}}{\beta!}D^{\beta}f\right)  ^{\wedge
}=\frac{\sqrt{2\pi}}{n!}\left(  ia\xi\right)  ^{n+1}\overline{\widehat{g_{n}}%
}\left(  a\xi\right)  \widehat{f},
\]

which implies the tempered distribution Taylor series expansion
\begin{equation}
f\left(  \cdot+a\right)  -\sum_{\left\vert \beta\right\vert \leq n}%
\frac{a^{\beta}}{\beta!}D^{\beta}f=f\left(  \cdot+a\right)  -\sum_{k\leq
n}\frac{\left(  aD\right)  ^{k}}{k!}f=\left(  \mathcal{R}_{n+1}f\right)
\left(  \cdot,a\right)  ,\quad f\in S^{\prime},\label{a1.55}%
\end{equation}

where%
\begin{equation}
\left(  \mathcal{R}_{n+1}f\right)  \left(  \cdot,a\right)  =\frac{\sqrt{2\pi}%
}{n!}\left(  \left(  ia\xi\right)  ^{n+1}\overline{\widehat{g_{n}}}\left(
a\xi\right)  \widehat{f}\right)  ^{\vee},\quad f\in S^{\prime},\text{
}n=0,1,2,\ldots\label{a50.5}%
\end{equation}

Suppose $f\in C_{B}^{\left(  n\right)  }$ and $\left(  aD\right)  ^{n+1}f\in
L^{1}$. For clarity set $u=\left(  aD\right)  ^{n+1}f$ where $D=\left(
D_{k}\right)  _{k=1}^{d}$, so that
\begin{align*}
\left(  ia\xi\right)  ^{n+1}\overline{\widehat{g_{n}}}\left(  a\xi\right)
\widehat{f}\left(  \xi\right)  =\overline{\widehat{g_{n}}}\left(  a\xi\right)
\widehat{u}\left(  \xi\right)   &  =\frac{1}{\sqrt{2\pi}}\int_{0}^{1}%
e^{isa\xi}g_{n}\left(  s\right)  \widehat{u}\left(  \xi\right)  ds\\
&  =\frac{1}{\sqrt{2\pi}}\int_{0}^{1}g_{n}\left(  s\right)  F_{x}\left[
u\left(  x+sa\right)  \right]  \left(  \xi\right)  ds\\
&  =\frac{1}{\sqrt{2\pi}}F_{x}\left[  \int_{0}^{1}g_{n}\left(  s\right)
u\left(  x+sa\right)  ds\right]  \left(  \xi\right)  ,
\end{align*}

since the Fourier transform can be transposed with the integral because $u\in
L^{1}$ allows an application of Fubini's theorem. Thus%
\[
\left(  \mathcal{R}_{n+1}f\right)  \left(  x,a\right)  =\frac{1}{n!}\int%
_{0}^{1}g_{n}\left(  s\right)  \left(  \left(  aD\right)  ^{n+1}f\right)
\left(  x+sa\right)  ds,\quad f\in C_{B}^{\left(  n\right)  },\text{ }\left(
aD\right)  ^{n+1}f\in L^{1},
\]

If, in addition, $\left(  \widehat{a}D\right)  ^{n+1}f\in L^{\infty}\left(
\left[  x,x+a\right]  \right)  $ where $\widehat{a}=a/\left\vert a\right\vert
$, then
\begin{equation}
\left\vert \left(  \mathcal{R}_{n+1}f\right)  \left(  x,a\right)  \right\vert
\leq\frac{1}{n!}\left\Vert \left(  aD\right)  ^{n+1}f\right\Vert _{\infty}%
\int_{0}^{1}g_{n}\left(  s\right)  ds=\frac{\left\vert a\right\vert ^{n+1}%
}{\left(  n+1\right)  !}\left\Vert \left(  \widehat{a}D\right)  ^{n+1}%
f\right\Vert _{\infty,\left[  x,x+a\right]  }.\label{p929}%
\end{equation}

Thus we have proved:

\begin{lemma}
\label{Lem_Taylor_estim_L1_Linf}Suppose $f\in C_{B}^{\left(  n\right)
}\left(  \mathbb{R}^{d}\right)  $ and in the distribution sense $\left(
aD\right)  ^{n+1}f\in L^{1}\left(  \mathbb{R}^{d}\right)  $ when $\left\vert
a\right\vert =1$.

Then for $0\leq k\leq n$,%
\begin{equation}
\left(  \mathcal{R}_{k+1}f\right)  \left(  x,a\right)  =\frac{\left\vert
a\right\vert ^{k+1}}{k!}\int_{0}^{1}g_{k}\left(  s\right)  \left(  \left(
\widehat{a}D\right)  ^{k+1}f\right)  \left(  x+sa\right)  ds,\quad
a,x\in\mathbb{R}^{d}.\label{a111}%
\end{equation}

If, in addition, $\left(  \widehat{a}D\right)  ^{n+1}f\in L^{\infty}\left(
\left[  x,x+a\right]  \right)  $ then we have the estimate \ref{p929}.
\end{lemma}

\section{Bounds for $\mathcal{Q}_{\emptyset,n,\xi}\left(  e^{i\left(
x,\xi\right)  }\right)  $\label{Sect_Qon(exp(i(x,z)))}}

In Section \ref{Sect_Pon_Qon} the operator $\mathcal{Q}_{\emptyset,n}$ was
introduced and now we will derive some upper bounds for the function
$\left\vert \mathcal{Q}_{\emptyset,n,\xi}\left(  e^{i\left(  x,\xi\right)
}\right)  \right\vert $, where $\mathcal{Q}_{\emptyset,n,\xi}$ acts on the
variable $\xi$ in $e^{i\left(  x,\xi\right)  }$. These estimates will be used
in Section \ref{Sect_explicit_basis} to prove an explicit formula for a basis
function in terms of its weight function, and used again in Section
\ref{Sect_another_inv_Four_X} to prove an inverse Fourier transform for a
member of $X_{w}^{\theta}$.

\begin{theorem}
\label{Thm_bound_on_g(e,x)}The function $\mathcal{Q}_{\emptyset,n,\xi}\left(
e^{i\left(  x,\xi\right)  }\right)  $ has the following properties:

\begin{enumerate}
\item There exist constants $\left\{  C_{m}=C_{m}\left(  \rho\right)
\right\}  _{m\leq n}$, defined by \ref{a76}, independent of $x$ and $\xi$,
such that%
\begin{equation}
\left\vert D_{x}^{\gamma}\mathcal{Q}_{\emptyset,n,\xi}\left(  e^{i\left(
x,\xi\right)  }\right)  \right\vert \leq\left\{
\begin{array}
[c]{ll}%
C_{n-\left\vert \gamma\right\vert }\left\vert \xi\right\vert ^{\left\vert
\gamma\right\vert }\left(  1+\left\vert x\right\vert \right)  ^{n-\left\vert
\gamma\right\vert -1}, & \left\vert \gamma\right\vert <n,\\
\left\vert \xi\right\vert ^{\left\vert \gamma\right\vert }, & \left\vert
\gamma\right\vert \geq n.
\end{array}
\right. \label{a27}%
\end{equation}

This is a useful estimate for \textbf{large }$x$\textbf{\ and large }$\xi
$.\smallskip

\item Given $r>0$ there exist constants $\left\{  C_{n,m,r}=C_{n,m,r}\left(
\rho\right)  \right\}  _{m\leq n}$, defined by \ref{a97}, such that for all $x
$ and $\left\vert \xi\right\vert \leq r$,
\begin{equation}
\left\vert D_{x}^{\gamma}\mathcal{Q}_{\emptyset,n,\xi}\left(  e^{i\left(
x,\xi\right)  }\right)  \right\vert \leq\left\{
\begin{array}
[c]{ll}%
C_{n,\left\vert \gamma\right\vert ,r}\left\vert \xi\right\vert ^{n}\left(
1+\left\vert x\right\vert \right)  ^{n-\left\vert \gamma\right\vert }, &
\left\vert \gamma\right\vert <n,\\
\left\vert \xi\right\vert ^{\left\vert \gamma\right\vert }, & \left\vert
\gamma\right\vert \geq n.
\end{array}
\right. \label{a25}%
\end{equation}

This estimate is useful \textbf{near} $\xi=\mathbf{0}$.
\end{enumerate}
\end{theorem}

\begin{proof}
We first derive a formula for $D_{x}^{\gamma}\mathcal{Q}_{\emptyset,n,\xi
}\left(  e^{i\left(  x,\xi\right)  }\right)  $. For some $\rho\in S_{1,n}$
\begin{align}
\mathcal{Q}_{\emptyset,n,\xi}\left(  e^{i\left(  x,\xi\right)  }\right)
=e^{i\left(  x,\xi\right)  }-\rho(\xi)\sum_{\left\vert \alpha\right\vert
<n}\frac{i^{\left\vert \alpha\right\vert }x^{\alpha}\xi^{\alpha}}{\alpha!} &
=e^{i\left(  x,\xi\right)  }-\rho(\xi)\sum_{k<n}\sum_{\left\vert
\alpha\right\vert =k}\frac{i^{\left\vert \alpha\right\vert }x^{\alpha}%
\xi^{\alpha}}{\alpha!}\nonumber\\
&  =e^{i\left(  x,\xi\right)  }-\rho(\xi)\sum_{k<n}\frac{\left(
ix,\xi\right)  ^{k}}{k!},\label{a74}%
\end{align}

so that
\begin{equation}
D_{x}^{\gamma}\mathcal{Q}_{\emptyset,n,\xi}\left(  e^{i\left(  x,\xi\right)
}\right)  =\left\{
\begin{array}
[c]{ll}%
\left(  i\xi\right)  ^{\gamma}\mathcal{Q}_{\emptyset,n-\left\vert
\gamma\right\vert ,\xi}\left(  e^{i\left(  x,\xi\right)  }\right)  , &
\left\vert \gamma\right\vert <n,\\
\left(  i\xi\right)  ^{\gamma}, & \left\vert \gamma\right\vert \geq n.
\end{array}
\right. \label{a21}%
\end{equation}
\medskip

\textbf{Part 1} If $\left\vert \gamma\right\vert \geq n$ then the inequality
is a very simple consequence of the formula%
\[
D_{x}^{\gamma}\mathcal{Q}_{\emptyset,n,\xi}\left(  e^{i\left(  x,\xi\right)
}\right)  =\left(  i\xi\right)  ^{\gamma}e^{i\left(  x,\xi\right)  }.
\]

On the other hand, if $\left\vert \gamma\right\vert <n$ then from equations
\ref{a21} and then \ref{a74},
\begin{align*}
\left\vert D_{x}^{\gamma}\mathcal{Q}_{\emptyset,n,\xi}\left(  e^{i\left(
x,\xi\right)  }\right)  \right\vert  & =\left\vert \left(  i\xi\right)
^{\gamma}\mathcal{Q}_{\emptyset,n-\left\vert \gamma\right\vert ,\xi}\left(
e^{i\left(  x,\xi\right)  }\right)  \right\vert \\
& \leq\left\vert \xi\right\vert ^{\left\vert \gamma\right\vert }\left(
1+\sum_{k<n-\left\vert \gamma\right\vert }\frac{\left\vert \rho(\xi
)\right\vert \left\vert \xi\right\vert ^{k}\left\vert x\right\vert ^{k}}%
{k!}\right) \\
& \leq\left\vert \xi\right\vert ^{\left\vert \gamma\right\vert }\left(
1+\sum_{k<n-\left\vert \gamma\right\vert }\frac{\left\vert \rho(\xi
)\right\vert \left\vert \xi\right\vert ^{k}\left(  1+\left\vert x\right\vert
\right)  ^{k}}{k!}\right) \\
& \leq\left\vert \xi\right\vert ^{\left\vert \gamma\right\vert }\left(
1+\sum_{k<n-\left\vert \gamma\right\vert }\frac{\left\vert \rho(\xi
)\right\vert \left\vert \xi\right\vert ^{k}}{k!}\right)  \left(  1+\left\vert
x\right\vert \right)  ^{n-\left\vert \gamma\right\vert -1}\\
& \leq\max_{j<n-\left\vert \gamma\right\vert }\left\Vert \left\vert
\cdot\right\vert ^{j}\rho\right\Vert _{\infty}\left\vert \xi\right\vert
^{\left\vert \gamma\right\vert }\left(  1+\sum_{k<n-\left\vert \gamma
\right\vert }\frac{1}{k!}\right)  \left(  1+\left\vert x\right\vert \right)
^{n-\left\vert \gamma\right\vert -1}\\
& \leq\max_{j<n-\left\vert \gamma\right\vert }\left\Vert \left\vert
\cdot\right\vert ^{j}\rho\right\Vert _{\infty}\left\vert \xi\right\vert
^{\left\vert \gamma\right\vert }\left(  1+\left\vert x\right\vert \right)
^{n-\left\vert \gamma\right\vert -1}\\
& =C_{n-\left\vert \gamma\right\vert }\left\vert \xi\right\vert ^{\left\vert
\gamma\right\vert }\left(  1+\left\vert x\right\vert \right)  ^{n-\left\vert
\gamma\right\vert -1},
\end{align*}

where
\begin{equation}
C_{m}=\max\limits_{k<m}\left\Vert \left\vert \cdot\right\vert ^{k}%
\rho\right\Vert _{\infty},\label{a76}%
\end{equation}

as required.\medskip

\textbf{Part 2} If $\left\vert \gamma\right\vert \geq n$ the inequality
follows from part 1. There remains the case $\left\vert \gamma\right\vert <n
$. Using the Taylor series expansion about zero (Appendix
\ref{Sect_taylor_expansion}) we define the "remainder" function $\mu_{m}$ by
$e^{it}=\sum\limits_{k<m}\frac{\left(  it\right)  ^{k}}{k!}+\left(  it\right)
^{m}\mu_{m}(t)$ (note \ref{a114}, \ref{a117} imply $\mu_{m}=\frac{\sqrt{2\pi}%
}{m!}\overline{\widehat{g_{m}}}$) and note that $\left\Vert \mu_{m}\right\Vert
_{\infty}\leq\frac{1}{m!}$. The following calculations%
\begin{align*}
D_{x}^{\gamma}\mathcal{Q}_{\emptyset,n,\xi}\left(  e^{ix\xi}\right)   &
=\left(  i\xi\right)  ^{\gamma}\mathcal{Q}_{\emptyset,n-\left\vert
\gamma\right\vert ,\xi}\left(  e^{ix\xi}\right) \\
&  =\left(  i\xi\right)  ^{\gamma}\left(  e^{ix\xi}-\rho(\xi)\sum
_{k<n-\left\vert \gamma\right\vert }\frac{\left(  ix\xi\right)  ^{k}}%
{k!}\right) \\
&  =\left(  i\xi\right)  ^{\gamma}\left(  e^{ix\xi}-\rho(\xi)\left(  e^{ix\xi
}-\left(  ix\xi\right)  ^{n-\left\vert \gamma\right\vert }\mu_{n-\left\vert
\gamma\right\vert }(ix\xi)\right)  \right) \\
&  =\left(  i\xi\right)  ^{\gamma}\left(  e^{ix\xi}\left(  1-\rho(\xi)\right)
+\rho(\xi)\left(  ix\xi\right)  ^{n-\left\vert \gamma\right\vert }%
\mu_{n-\left\vert \gamma\right\vert }(ix\xi)\right)  ,
\end{align*}

now imply that when $\left\vert \xi\right\vert \leq r$:
\begin{align*}
\left\vert D_{x}^{\gamma}\mathcal{Q}_{\emptyset,n,\xi}\left(  e^{ix\xi
}\right)  \right\vert  &  \leq\left\vert \xi\right\vert ^{\left\vert
\gamma\right\vert }\left(  \left\vert 1-\rho(\xi)\right\vert +\left\vert
\rho(\xi)\right\vert \left\vert \xi\right\vert ^{n-\left\vert \gamma
\right\vert }\frac{\left\vert x\right\vert ^{n-\left\vert \gamma\right\vert }%
}{\left(  n-\left\vert \gamma\right\vert \right)  !}\right) \\
&  =\left\vert \xi\right\vert ^{\left\vert \gamma\right\vert }\left(
\left\vert 1-\rho(\xi)\right\vert +\left\vert \rho(\xi)\right\vert \left\vert
\xi\right\vert ^{n-\left\vert \gamma\right\vert }\frac{\left\vert x\right\vert
^{n-\left\vert \gamma\right\vert }}{\left(  n-\left\vert \gamma\right\vert
\right)  !}\right) \\
&  \leq\left(  \left\vert \xi\right\vert ^{\left\vert \gamma\right\vert }%
\frac{\left\vert 1-\rho(\xi)\right\vert }{\left\vert \xi\right\vert ^{n}%
}+\left\vert \rho(\xi)\right\vert \frac{\left\vert x\right\vert ^{n-\left\vert
\gamma\right\vert }}{\left(  n-\left\vert \gamma\right\vert \right)
!}\right)  \left\vert \xi\right\vert ^{n}\\
&  \leq\left(  r^{\left\vert \gamma\right\vert }\frac{\left\vert 1-\rho
(\xi)\right\vert }{\left\vert \xi\right\vert ^{n}}+\left\Vert \rho\right\Vert
_{\infty;\leq r}\left\vert x\right\vert ^{n-\left\vert \gamma\right\vert
}\right)  \left\vert \xi\right\vert ^{n}.
\end{align*}

But $\rho\in S_{1,n}$ implies $1-\rho\in C_{\emptyset,n}^{\infty}$ and by
Theorem \ref{Thm_Tay_rem_zeros} in Appendix \ref{Sect_taylor_expansion}, given
$r>0$,
\[
\left\vert 1-\rho(\xi)\right\vert \leq k_{\rho;n,r}\left\vert \xi\right\vert
^{n},\quad\left\vert \xi\right\vert \leq r.
\]

when%
\[
k_{\rho;n;r}=\frac{1}{n!}\left\Vert \left(  \widehat{\cdot}D\right)  ^{n}%
\rho\right\Vert _{\infty;\leq r},\quad n=0,1,2,\ldots.
\]

Set%
\begin{equation}
C_{n,m,r}=\max\left\{  r^{m}k_{\rho;n,r},\left\Vert \rho\right\Vert
_{\infty;\leq r}\right\}  ,\quad m\leq n.\label{a97}%
\end{equation}

Then when $x\in\mathbb{R}^{d}$ and $\left\vert \xi\right\vert \leq r$,
\[
\left\vert D_{x}^{\gamma}\mathcal{Q}_{\emptyset,n,\xi}\left(  e^{ix\xi
}\right)  \right\vert \leq C_{n,\left\vert \gamma\right\vert ,r}\left\vert
\xi\right\vert ^{n}\left(  1+\left\vert x\right\vert ^{n-\left\vert
\gamma\right\vert }\right)  \leq C_{n,\left\vert \gamma\right\vert
,r}\left\vert \xi\right\vert ^{n}\left(  1+\left\vert x\right\vert \right)
^{n-\left\vert \gamma\right\vert }.
\]

\end{proof}

We modify the proofs of Theorem \ref{Thm_bound_on_g(e,x)} to obtain the
following more delicate estimates:

\begin{theorem}
\label{Thm_bound_on_g(e,x)_2}The function $\mathcal{Q}_{\emptyset,n,\xi
}\left(  e^{ix\xi}\right)  $ has the following properties:

\begin{enumerate}
\item Given $r\geq0$ for all $x$ and $\left\vert \xi\right\vert \geq r$,%
\begin{align}
\left\vert D_{x}^{\gamma}\mathcal{Q}_{\emptyset,n,\xi}\left(  e^{ix\xi
}\right)  \right\vert  & \leq\left\{
\begin{array}
[c]{ll}%
\overline{C}_{n-\left\vert \gamma\right\vert ,r}^{\left(  \rho\right)
}\left(  1+\sum\limits_{k<n-\left\vert \gamma\right\vert }\frac{\left\vert
x\widehat{\xi}\right\vert ^{k}}{k!}\right)  \left\vert \xi^{\gamma}\right\vert
, & \left\vert \gamma\right\vert <n,\\
\left\vert \xi^{\gamma}\right\vert , & \left\vert \gamma\right\vert \geq n,
\end{array}
\right.  ,\label{a27.2}\\
& \smallskip\nonumber\\
\left\vert D_{x}^{\gamma}\mathcal{Q}_{\emptyset,n,\xi}\left(  e^{ix\xi
}\right)  \right\vert  & \leq\left\{
\begin{array}
[c]{ll}%
\overline{C}_{n-\left\vert \gamma\right\vert ,r}^{\left(  \rho\right)
}\left(  1+\sum\limits_{\left\vert \alpha\right\vert <n-\left\vert
\gamma\right\vert }\frac{\left\vert x^{\alpha}\right\vert }{\alpha!}\right)
\left\vert \xi^{\gamma}\right\vert , & \left\vert \gamma\right\vert <n,\\
\left\vert \xi^{\gamma}\right\vert , & \left\vert \gamma\right\vert \geq n,
\end{array}
\right. \label{a27.1}%
\end{align}

where $\overline{C}_{n-\left\vert \gamma\right\vert ,r}^{\left(  \rho\right)
}$ is defined by \ref{a2.03}. This is a useful estimate for \textbf{large} $x$
and \textbf{large }$\xi$.\smallskip

\item Given $r>0$ for all $x$ and $\left\vert \xi\right\vert \leq r$,%
\begin{align}
\left\vert D_{x}^{\gamma}\mathcal{Q}_{\emptyset,n,\xi}\left(  e^{ix\xi
}\right)  \right\vert  & \leq\left\{
\begin{array}
[c]{ll}%
\underline{C}_{n,r}^{\left(  \rho\right)  }\left\vert \xi\right\vert
^{n}\left(  \left\vert \xi^{\gamma}\right\vert +\frac{\left\vert x\widehat
{\xi}\right\vert ^{n-\left\vert \gamma\right\vert }}{\left(  n-\left\vert
\gamma\right\vert \right)  !}\right)  , & \left\vert \gamma\right\vert <n,\\
\left\vert \xi^{\gamma}\right\vert , & \left\vert \gamma\right\vert \geq n.
\end{array}
\right. \label{a702}\\
& \smallskip\nonumber\\
\left\vert D_{x}^{\gamma}\mathcal{Q}_{\emptyset,n,\xi}\left(  e^{ix\xi
}\right)  \right\vert  & \leq\left\{
\begin{array}
[c]{ll}%
\underline{C}_{n,r}^{\left(  \rho\right)  }\left\vert \xi\right\vert
^{n}\left(  \left\vert \xi^{\gamma}\right\vert +\sum\limits_{\left\vert
\alpha\right\vert =n-\left\vert \gamma\right\vert }\frac{\left\vert x^{\alpha
}\right\vert }{\alpha!}\right)  , & \left\vert \gamma\right\vert <n,\\
\left\vert \xi^{\gamma}\right\vert , & \left\vert \gamma\right\vert \geq n.
\end{array}
\right. \label{a77}%
\end{align}

where $\underline{C}_{n,r}^{\left(  \rho\right)  }$ is given by \ref{a75}.
\end{enumerate}
\end{theorem}

\begin{proof}
We first derive a formula for $D_{x}^{\gamma}\mathcal{Q}_{\emptyset,n,\xi
}\left(  e^{ix\xi}\right)  $. For some $\rho\in S_{1,n}$,
\begin{align}
\mathcal{Q}_{\emptyset,n,\xi}\left(  e^{ix\xi}\right)  =e^{ix\xi}-\rho
(\xi)\sum_{\left\vert \alpha\right\vert <n}\frac{i^{\left\vert \alpha
\right\vert }x^{\alpha}\xi^{\alpha}}{\alpha!} &  =e^{ix\xi}-\rho(\xi
)\sum_{k<n}\sum_{\left\vert \alpha\right\vert =k}\frac{i^{\left\vert
\alpha\right\vert }x^{\alpha}\xi^{\alpha}}{\alpha!}\nonumber\\
&  =e^{i\left(  x,\xi\right)  }-\rho(\xi)\sum_{k<n}\frac{\left(  ix\xi\right)
^{k}}{k!},\label{a741}%
\end{align}

so that
\begin{equation}
D_{x}^{\gamma}\mathcal{Q}_{\emptyset,n,\xi}\left(  e^{ix\xi}\right)  =\left\{
\begin{array}
[c]{ll}%
\left(  i\xi\right)  ^{\gamma}\mathcal{Q}_{\emptyset,n-\left\vert
\gamma\right\vert ,\xi}\left(  e^{ix\xi}\right)  , & \left\vert \gamma
\right\vert <n,\\
\left(  i\xi\right)  ^{\gamma}, & \left\vert \gamma\right\vert \geq n.
\end{array}
\right. \label{a211}%
\end{equation}
\medskip

\textbf{Part 1} If $\left\vert \gamma\right\vert \geq n$ then clearly
$\left\vert D_{x}^{\gamma}\mathcal{Q}_{\emptyset,n,\xi}\left(  e^{ix\xi
}\right)  \right\vert \leq\left\vert \xi^{\gamma}\right\vert $. On the other
hand, if $\left\vert \gamma\right\vert <n$ then from equations \ref{a211} and
then \ref{a741},
\begin{align*}
\left\vert D_{x}^{\gamma}\mathcal{Q}_{\emptyset,n,\xi}\left(  e^{ix\xi
}\right)  \right\vert  & =\left\vert \left(  i\xi\right)  ^{\gamma}%
\mathcal{Q}_{\emptyset,n-\left\vert \gamma\right\vert ,\xi}\left(  e^{ix\xi
}\right)  \right\vert \\
& \leq\left\vert \xi^{\gamma}\right\vert \left(  1+\sum_{k<n-\left\vert
\gamma\right\vert }\left\vert \rho(\xi)\right\vert \frac{\left\vert
x\xi\right\vert ^{k}}{k!}\right) \\
& \leq\left\vert \xi^{\gamma}\right\vert \left(  1+\sum_{k<n-\left\vert
\gamma\right\vert }\frac{1}{k!}\left\vert \xi\right\vert ^{k}\left\vert
\rho(\xi)\right\vert \left\vert x\widehat{\xi}\right\vert ^{k}\right) \\
& \leq\left\vert \xi^{\gamma}\right\vert \left(  1+\sum_{k<n-\left\vert
\gamma\right\vert }\frac{1}{k!}\left\Vert \left\vert \cdot\right\vert ^{k}%
\rho\right\Vert _{\infty;\geq r}\left\vert x\widehat{\xi}\right\vert
^{k}\right) \\
& \leq\max\left\{  1,\max_{j<n-\left\vert \gamma\right\vert }\left\Vert
\left\vert \cdot\right\vert ^{j}\rho\right\Vert _{\infty;\geq r}\right\}
\left\vert \xi^{\gamma}\right\vert \left(  1+\sum_{k<n-\left\vert
\gamma\right\vert }\frac{1}{k!}\left\vert x\widehat{\xi}\right\vert
^{k}\right) \\
& \leq\overline{C}_{n-\left\vert \gamma\right\vert ,r}^{\left(  \rho\right)
}\left\vert \xi^{\gamma}\right\vert \left(  1+\sum_{k<n-\left\vert
\gamma\right\vert }\frac{1}{k!}\left\vert x\widehat{\xi}\right\vert
^{k}\right)  ,
\end{align*}

where%
\begin{equation}
\overline{C}_{m,r}^{\left(  \rho\right)  }=\max\left\{  1,\max_{j<m}\left\Vert
\left\vert \cdot\right\vert ^{j}\rho\right\Vert _{\infty;\geq r}\right\}
,\quad m\geq1.\label{a2.03}%
\end{equation}

But%
\[
\sum_{k<n-\left\vert \gamma\right\vert }\frac{1}{k!}\left\vert x\widehat{\xi
}\right\vert ^{k}=\sum_{k<n-\left\vert \gamma\right\vert }\left\vert
\sum\limits_{\left\vert \alpha\right\vert =k}\frac{1}{\alpha!}x^{\alpha
}\widehat{\xi}^{\alpha}\right\vert \leq\sum_{k<n-\left\vert \gamma\right\vert
}\sum\limits_{\left\vert \alpha\right\vert =k}\frac{\left\vert x^{\alpha
}\right\vert }{\alpha!}=\sum\limits_{\left\vert \alpha\right\vert
<n-\left\vert \gamma\right\vert }\frac{\left\vert x^{\alpha}\right\vert
}{\alpha!}.
\]
\smallskip

\textbf{Part 2} If $\left\vert \gamma\right\vert \geq n$ the inequality
follows from part 1.

There remains the case $\left\vert \gamma\right\vert <n$. Using the Taylor
series expansion about zero (Appendix \ref{Sect_taylor_expansion}) we define
the "remainder" function $\mu_{m}$ by $e^{it}=\sum\limits_{k<m}\frac{\left(
it\right)  ^{k}}{k!}+\left(  it\right)  ^{m}\mu_{m}(t)$ (note \ref{a114},
\ref{a117} imply $\mu_{m}=\frac{\sqrt{2\pi}}{m!}\overline{\widehat{g_{m}}}$)
and note that $\left\Vert \mu_{m}\right\Vert _{\infty}\leq\frac{1}{m!}$. The
following calculations%
\begin{align*}
D_{x}^{\gamma}\mathcal{Q}_{\emptyset,n,\xi}\left(  e^{ix\xi}\right)  =\left(
i\xi\right)  ^{\gamma}\mathcal{Q}_{\emptyset,n-\left\vert \gamma\right\vert
,\xi}\left(  e^{ix\xi}\right)   &  =\left(  i\xi\right)  ^{\gamma}\left(
e^{ix\xi}-\rho(\xi)\sum_{k<n-\left\vert \gamma\right\vert }\frac{\left(
ix\xi\right)  ^{k}}{k!}\right) \\
&  =\left(  i\xi\right)  ^{\gamma}\left(  e^{ix\xi}-\rho(\xi)\left(  e^{ix\xi
}-\left(  ix\xi\right)  ^{n-\left\vert \gamma\right\vert }\mu_{n-\left\vert
\gamma\right\vert }(x\xi)\right)  \right) \\
&  =\left(  i\xi\right)  ^{\gamma}\left(  e^{ix\xi}\left(  1-\rho(\xi)\right)
+\rho(\xi)\left(  ix\xi\right)  ^{n-\left\vert \gamma\right\vert }%
\mu_{n-\left\vert \gamma\right\vert }(x\xi)\right)  ,
\end{align*}

now imply that when $\left\vert \xi\right\vert \leq r$,
\begin{align*}
\left\vert D_{x}^{\gamma}\mathcal{Q}_{\emptyset,n,\xi}\left(  e^{ix\xi
}\right)  \right\vert  &  \leq\left\vert \xi^{\gamma}\right\vert \left(
\left\vert 1-\rho(\xi)\right\vert +\left\vert \rho(\xi)\right\vert \left\vert
\xi\right\vert ^{n-\left\vert \gamma\right\vert }\frac{\left\vert
x\widehat{\xi}\right\vert ^{n-\left\vert \gamma\right\vert }}{\left(
n-\left\vert \gamma\right\vert \right)  !}\right) \\
&  =\left\vert \xi^{\gamma}\right\vert \left\vert 1-\rho(\xi)\right\vert
+\left\vert \rho(\xi)\right\vert \left\vert \xi^{\gamma}\right\vert \left\vert
\xi\right\vert ^{n-\left\vert \gamma\right\vert }\frac{\left\vert
x\widehat{\xi}\right\vert ^{n-\left\vert \gamma\right\vert }}{\left(
n-\left\vert \gamma\right\vert \right)  !}\\
&  \leq\left\vert \xi^{\gamma}\right\vert \left\vert 1-\rho(\xi)\right\vert
+\left\vert \rho(\xi)\right\vert \left\vert \xi\right\vert ^{\left\vert
\gamma\right\vert }\left\vert \xi\right\vert ^{n-\left\vert \gamma\right\vert
}\frac{\left\vert x\widehat{\xi}\right\vert ^{n-\left\vert \gamma\right\vert
}}{\left(  n-\left\vert \gamma\right\vert \right)  !}\\
&  =\left\vert \xi^{\gamma}\right\vert \left\vert 1-\rho(\xi)\right\vert
+\left\vert \rho(\xi)\right\vert \left\vert \xi\right\vert ^{n}\frac
{\left\vert x\widehat{\xi}\right\vert ^{n-\left\vert \gamma\right\vert }%
}{\left(  n-\left\vert \gamma\right\vert \right)  !}\\
&  \leq\left\vert \xi\right\vert ^{n}\left(  \left\vert \xi^{\gamma
}\right\vert \frac{\left\vert 1-\rho(\xi)\right\vert }{\left\vert
\xi\right\vert ^{n}}+\left\vert \rho(\xi)\right\vert \frac{\left\vert
x\widehat{\xi}\right\vert ^{n-\left\vert \gamma\right\vert }}{\left(
n-\left\vert \gamma\right\vert \right)  !}\right)  .
\end{align*}

But $\rho\in S_{1,n}$ implies $1-\rho\in C_{\emptyset,n}^{\infty}$ and by
Theorem \ref{Thm_Tay_rem_zeros} in Appendix \ref{Sect_taylor_expansion},%
\[
\left\vert 1-\rho(\xi)\right\vert \leq\frac{1}{n!}\underline{C}_{n,r}^{\left(
\rho\right)  }\left\vert \xi\right\vert ^{n},\quad\left\vert \xi\right\vert
\leq r,
\]

where
\begin{equation}
\underline{C}_{n,r}^{\left(  \rho\right)  }=\frac{1}{n!}\left\Vert \left(
\widehat{\cdot}D\right)  ^{n}\rho\right\Vert _{\infty;\leq r},\quad
n=0,1,2,\ldots.\label{a75}%
\end{equation}

Thus when $\left\vert \xi\right\vert \leq r$,
\begin{align*}
\left\vert D_{x}^{\gamma}\mathcal{Q}_{\emptyset,n,\xi}\left(  e^{i\left(
x,\xi\right)  }\right)  \right\vert  & \leq\left\vert \xi\right\vert
^{n}\left(  \underline{C}_{n,r}^{\left(  \rho\right)  }\left\vert \xi^{\gamma
}\right\vert +\left\Vert \rho\right\Vert _{\infty;\leq r}\frac{\left\vert
x\widehat{\xi}\right\vert ^{n-\left\vert \gamma\right\vert }}{\left(
n-\left\vert \gamma\right\vert \right)  !}\right) \\
& \leq\underline{C}_{n,r}^{\left(  \rho\right)  }\left\vert \xi\right\vert
^{n}\left(  \left\vert \xi^{\gamma}\right\vert +\frac{\left\vert x\widehat
{\xi}\right\vert ^{n-\left\vert \gamma\right\vert }}{\left(  n-\left\vert
\gamma\right\vert \right)  !}\right)  ,
\end{align*}

Finally, from the identity \ref{p08} we have%
\[
\frac{\left\vert x\widehat{\xi}\right\vert ^{n-\left\vert \gamma\right\vert }%
}{\left(  n-\left\vert \gamma\right\vert \right)  !}=\left\vert \sum
\limits_{\left\vert \alpha\right\vert =n-\left\vert \gamma\right\vert }%
\frac{x^{\alpha}\widehat{\xi}^{\alpha}}{\alpha!}\right\vert \leq
\sum\limits_{\left\vert \alpha\right\vert =n-\left\vert \gamma\right\vert
}\frac{\left\vert x^{\alpha}\widehat{\xi}^{\alpha}\right\vert }{\alpha!}%
\leq\sum\limits_{\left\vert \alpha\right\vert =n-\left\vert \gamma\right\vert
}\frac{\left\vert x^{\alpha}\right\vert }{\alpha!},
\]

so that%
\[
\left\vert D_{x}^{\gamma}\mathcal{Q}_{\emptyset,n,\xi}\left(  e^{i\left(
x,\xi\right)  }\right)  \right\vert \leq\underline{C}_{n,r}^{\left(
\rho\right)  }\left\vert \xi\right\vert ^{n}\left(  \left\vert \xi^{\gamma
}\right\vert +\sum\limits_{\left\vert \alpha\right\vert =n-\left\vert
\gamma\right\vert }\frac{\left\vert x^{\alpha}\right\vert }{\alpha!}\right)
,\quad\left\vert \xi\right\vert \leq r.
\]

\end{proof}

\begin{remark}
\label{Rem_Thm_bound_on_g(e,x)_2}\ 

\begin{enumerate}
\item We use forms involving sums of $\left\vert x^{\alpha}\right\vert $
instead of $\left\vert x\right\vert $ because of their simple interaction with
the Fourier transform i.e. $\left\vert x^{\alpha}\right\vert \left\vert
\widehat{u}\right\vert =\left\vert x^{\alpha}\widehat{u}\right\vert
=\left\vert \widehat{D^{\alpha}u}\right\vert $.

\item The Stirling formula
\[
k!=\left(  2\pi\right)  ^{1/2}k^{1/2}\left(  \frac{k}{e}\right)  ^{k}%
\exp\left(  \frac{1}{12k}-\frac{1}{360k^{3}}+\ldots\right)  >\left(
2\pi\right)  ^{1/2}k^{1/2}\left(  \frac{k}{e}\right)  ^{k}>\left(  \frac{k}%
{e}\right)  ^{k},
\]

implies that
\begin{equation}
1+\frac{\left\vert s\right\vert ^{k}}{k!}<1+\left(  \frac{e}{k}\left\vert
s\right\vert \right)  ^{k}<\left(  1+\frac{e}{k}\left\vert s\right\vert
\right)  ^{k}.\label{a021}%
\end{equation}

\end{enumerate}
\end{remark}

\section{A useful inverse Fourier transform theorem\label{Sect_invers_Fourier}%
}

The next theorem will be used in Section \ref{Sect_explicit_basis} to derive
an inverse Fourier transform formula for a basis function and in Section
\ref{Sect_another_inv_Four_X} to derive an inverse Fourier transform formula
for any member of $X_{w}^{\theta}$. This theorem will use the following lemma
which shows how, when $\phi\in S$, the expressions $\mathcal{P}_{\emptyset
,n}\phi$ and $\mathcal{Q}_{\emptyset,n}\phi$ can be expressed in terms of
$\mathcal{P}_{\emptyset,n,x}\left(  e^{i\left(  x,\xi\right)  }\right)  $ and
$\mathcal{Q}_{\emptyset,n,x}\left(  e^{i\left(  x,\xi\right)  }\right)  $.

\begin{lemma}
\label{Lem_Pon(f)_in_terms_of_Pon(exp)}If $\psi\in S$ then:

\begin{enumerate}
\item $\left(  \mathcal{P}_{\emptyset,n}\overset{\vee}{\psi}\right)
(\xi)=\left(  2\pi\right)  ^{-\frac{d}{2}}\int\mathcal{P}_{\emptyset,n,\xi
}\left(  e^{i\left(  x,\xi\right)  }\right)  \psi\left(  x\right)  dx$.

\item $\left(  \mathcal{Q}_{\emptyset,n}\overset{\vee}{\psi}\right)
(\xi)=\left(  2\pi\right)  ^{-\frac{d}{2}}\int\mathcal{Q}_{\emptyset,n,\xi
}\left(  e^{i\left(  x,\xi\right)  }\right)  \psi\left(  x\right)  dx$.
\end{enumerate}
\end{lemma}

\begin{proof}
If $\psi\in S$ then%
\begin{align*}
\left(  \mathcal{P}_{\emptyset,n}\overset{\vee}{\psi}\right)  (\xi
)=\rho\left(  \xi\right)  \sum_{\left\vert \alpha\right\vert <n}\frac
{\xi^{\alpha}D^{\alpha}\overset{\vee}{\psi}(0)}{\alpha!} &  =\rho\left(
\xi\right)  \sum_{\left\vert \alpha\right\vert <n}\frac{\xi^{\alpha}}{\alpha
!}\left(  \left(  ix\right)  ^{\alpha}\psi\right)  ^{\vee}(0)\\
&  =\left(  2\pi\right)  ^{-\frac{d}{2}}\rho\left(  \xi\right)  \sum
_{\left\vert \alpha\right\vert <n}\frac{\xi^{\alpha}}{\alpha!}\int\left(
ix\right)  ^{\alpha}\psi\left(  x\right)  dx\\
&  =\left(  2\pi\right)  ^{-\frac{d}{2}}\int\left(  \sum_{\left\vert
\alpha\right\vert <n}\rho\left(  \xi\right)  \frac{\left(  ix\right)
^{\alpha}\xi^{\alpha}}{\alpha!}\right)  \psi\left(  x\right)  dx\\
&  =\left(  2\pi\right)  ^{-\frac{d}{2}}\int\mathcal{P}_{\emptyset,n,\xi
}\left(  e^{i\left(  x,\xi\right)  }\right)  \psi\left(  x\right)  dx,
\end{align*}

and hence%
\begin{align*}
\left(  \mathcal{Q}_{\emptyset,n}\overset{\vee}{\psi}\right)  (\xi)  &
=\overset{\vee}{\psi}\left(  \xi\right)  -\left(  \mathcal{P}_{\emptyset
,n}\overset{\vee}{\psi}\right)  (\xi)\\
& =\left(  2\pi\right)  ^{-\frac{d}{2}}\int e^{i\xi x}\psi\left(  x\right)
dx-\left(  2\pi\right)  ^{-d/2}\int\mathcal{P}_{\emptyset,n,\xi}\left(
e^{i\left(  x,\xi\right)  }\right)  \psi\left(  x\right)  dx\\
& =\left(  2\pi\right)  ^{-\frac{d}{2}}\int\mathcal{Q}_{\emptyset,n,\xi
}\left(  e^{i\left(  x,\xi\right)  }\right)  \psi\left(  x\right)  dx\\
& =\left(  2\pi\right)  ^{-\frac{d}{2}}\int\mathcal{Q}_{\emptyset,n,\xi
}\left(  e^{i\left(  x,\xi\right)  }\right)  \psi\left(  x\right)  dx.
\end{align*}

\end{proof}

To prove the next theorem we will require another lemma concerning
differentiation under the integral sign.

\begin{lemma}
\label{Lem_diff_under_integral_2}(Prop 7.8.4 of Malliavin \cite{Malliavin95}%
)Suppose $f:\mathbb{R}^{m+n}\rightarrow\mathbb{C}$ and we write $f\left(
\xi,x\right)  $ where $\xi\in\mathbb{R}^{m}$ and $x\in\mathbb{R}^{n}$. Further
suppose that:

\begin{enumerate}
\item For each $\xi$, $f\left(  \xi,\cdot\right)  \in C^{\left(  k\right)
}\left(  \mathbb{R}^{n}\right)  $.

\item For each $x$, $\int\left\vert D_{\xi}^{\alpha}f\left(  \xi,x\right)
d\xi\right\vert <\infty$ for $\left\vert \alpha\right\vert \leq k$.
\end{enumerate}

Then we have
\[
D_{x}^{\alpha}\int f\left(  \xi,x\right)  d\xi=\int D_{x}^{\alpha}f\left(
\xi,x\right)  d\xi,\;when\text{ }\left\vert \alpha\right\vert \leq k,
\]

and $\int f\left(  \xi,\cdot\right)  d\xi\in C^{\left(  k\right)  }\left(
\mathbb{R}^{n}\right)  $.
\end{lemma}

\begin{theorem}
\label{Thm_integ(|DQ(exp)fF|)}Suppose $f\in S^{\prime}$ and $\widehat{f}\in
L_{loc}^{1}\left(  \mathbb{R}^{d}\setminus0\right)  $. Define the function
$f_{F}:\mathbb{R}^{d}\rightarrow\mathbb{C}$ a.e. by $f_{F}=\widehat{f}$ on
$\mathbb{R}^{d}\setminus0$.

Then:

\begin{enumerate}
\item If the action $\int f_{F}\phi$ on $\phi\in S_{\emptyset,n}$ defines a
member of $S_{\emptyset,n}^{\prime}$, and if $f_{F}=\widehat{f}$ on
$S_{\emptyset,n}$, it follows that for all multi-indexes $\gamma$%
\begin{equation}
\left[  \widehat{D^{\gamma}f},\psi\right]  =\int\left(  i\xi\right)  ^{\gamma
}\left(  \mathcal{Q}_{\emptyset,n}\psi\right)  \left(  \xi\right)
f_{F}\left(  \xi\right)  d\xi+\left(  2\pi\right)  ^{-\frac{d}{2}}\left[
\widehat{p_{\widehat{D^{\gamma}f}}},\psi\right]  ,\quad\psi\in S,\label{a19}%
\end{equation}

where for $u\in S^{\prime}$, $p_{u}\in P_{n-1}$ is defined by
\begin{equation}
p_{u}\left(  x\right)  :=\sum\limits_{\left\vert \alpha\right\vert <n}%
\dfrac{b_{\alpha,u}}{\alpha!}x^{\alpha},\quad b_{\alpha,u}=\left[  u,\left(
-i\xi\right)  ^{\alpha}\rho\right]  .\label{a20}%
\end{equation}

Further, $\widehat{D^{\gamma}f}=\left(  i\xi\right)  ^{\gamma}f_{F}$ on
$S_{\emptyset,n}$, and $\left(  i\xi\right)  ^{\gamma}f_{F}\in L_{loc}%
^{1}\left(  \mathbb{R}^{d}\setminus0\right)  $.\medskip

\item For all $\gamma$,%
\begin{equation}
D^{\gamma}f=\left(  \mathcal{Q}_{\emptyset,n}^{\ast}\widehat{D^{\gamma}%
f}\right)  ^{\vee}+\left(  2\pi\right)  ^{-\frac{d}{2}}p_{\widehat{D^{\gamma
}f}}=D^{\gamma}\left(  \mathcal{Q}_{\emptyset,n}^{\ast}\widehat{f}\right)
^{\vee}+\left(  2\pi\right)  ^{-\frac{d}{2}}D^{\gamma}p_{\widehat{f}}.\quad
f\in S^{\prime}\text{.}\label{a89}%
\end{equation}

\item Now also assume that for a given multi-index $\gamma$ there exist
constants $s_{\gamma}\geq0$ and $k_{\gamma}>0$ such that%
\begin{equation}
\int\left\vert D_{x}^{\gamma}\mathcal{Q}_{\emptyset,n,\xi}\left(  e^{i\left(
x,\xi\right)  }\right)  f_{F}\left(  \xi\right)  \right\vert d\xi\leq
k_{\gamma}\left(  1+\left\vert x\right\vert \right)  ^{s_{\gamma}},\quad
x\in\mathbb{R}^{d}.\label{a37}%
\end{equation}

Then $D^{\gamma}f\in C_{BP}^{\left(  0\right)  }$,%
\begin{equation}
D^{\gamma}f\left(  x\right)  =\left(  2\pi\right)  ^{-\frac{d}{2}}\int
D_{x}^{\gamma}\mathcal{Q}_{\emptyset,n,\xi}\left(  e^{i\left(  x,\xi\right)
}\right)  f_{F}\left(  \xi\right)  d\xi+\left(  2\pi\right)  ^{-\frac{d}{2}%
}D^{\gamma}p_{\widehat{f}}\left(  x\right)  ,\label{a24}%
\end{equation}

and%
\begin{align}
\int D_{x}^{\gamma}\mathcal{Q}_{\emptyset,n,\xi}\left(  e^{i\left(
x,\xi\right)  }\right)  f_{F}\left(  \xi\right)  d\xi & =\left\{
\begin{array}
[c]{ll}%
\int\left(  i\xi\right)  ^{\gamma}\mathcal{Q}_{\emptyset,n-\left\vert
\gamma\right\vert ,\xi}\left(  e^{i\left(  x,\xi\right)  }\right)
f_{F}\left(  \xi\right)  d\xi, & \left\vert \gamma\right\vert <n,\\
\int e^{i\left(  x,\xi\right)  }\left(  i\xi\right)  ^{\gamma}f_{F}\left(
\xi\right)  d\xi, & \left\vert \gamma\right\vert \geq n,
\end{array}
\right. \label{a26}\\
& =\left\{
\begin{array}
[c]{ll}%
\int\mathcal{Q}_{\emptyset,n-\left\vert \gamma\right\vert ,\xi}\left(
e^{i\left(  x,\xi\right)  }\right)  \left(  D^{\gamma}f\right)  _{F}\left(
\xi\right)  d\xi, & \left\vert \gamma\right\vert <n,\\
\int e^{i\left(  x,\xi\right)  }\left(  D^{\gamma}f\right)  _{F}\left(
\xi\right)  d\xi, & \left\vert \gamma\right\vert \geq n,
\end{array}
\right. \label{a36}%
\end{align}

where the function $\left(  D^{\gamma}f\right)  _{F}:\mathbb{R}^{d}%
\rightarrow\mathbb{C}$ is a.e. defined by $\left(  D^{\gamma}f\right)
_{F}=\widehat{D^{\gamma}f}$ on $\mathbb{R}^{d}\setminus0$.

\item Suppose, in addition to the conditions imposed in part 1, we assume that
there are constants $k_{\gamma}^{\prime}>0$ and $s_{\gamma}^{\prime}\geq0$
such that%
\[
\int\left\vert \left(  i\xi\right)  ^{\gamma}\mathcal{Q}_{\emptyset,n,\xi
}\left(  e^{i\left(  x,\xi\right)  }\right)  f_{F}\left(  \xi\right)
\right\vert d\xi\leq k_{\gamma}^{\prime}\left(  1+\left\vert x\right\vert
\right)  ^{s_{\gamma}^{\prime}},\quad x\in\mathbb{R}^{d}.
\]

Then%
\begin{align}
D^{\gamma}f\left(  x\right)   & =\left(  2\pi\right)  ^{-\frac{d}{2}}%
\int\left(  i\xi\right)  ^{\gamma}\mathcal{Q}_{\emptyset,n,\xi}\left(
e^{i\left(  x,\xi\right)  }\right)  f_{F}\left(  \xi\right)  d\xi+\left(
2\pi\right)  ^{-\frac{d}{2}}p_{\widehat{D^{\gamma}f}}\left(  x\right)
\label{a90}\\
& =\left(  2\pi\right)  ^{-\frac{d}{2}}\int\mathcal{Q}_{\emptyset,n,\xi
}\left(  e^{i\left(  x,\xi\right)  }\right)  \left(  D^{\gamma}f\right)
_{F}\left(  \xi\right)  d\xi+\left(  2\pi\right)  ^{-\frac{d}{2}}%
p_{\widehat{D^{\gamma}f}}\left(  x\right)  .\nonumber
\end{align}

\end{enumerate}
\end{theorem}

\begin{proof}
\textbf{Part 1} We first prove equation \ref{a19}. From part 2 of Theorem
\ref{Thm_S=So+pP2m-1} we know that $\mathcal{Q}_{\emptyset,n}:S\rightarrow
S_{\emptyset,n}$ and since the current theorem assumes that $f_{F}=\widehat
{f}$ on $S_{\emptyset,n}$ and that $\int f_{F}\phi$ on $\phi\in S_{\emptyset
,n}$ defines a member of $S_{\emptyset,n}^{\prime}$, we have for $\psi\in S$
and the \textit{action} variable $\xi$
\[
\left[  \widehat{D^{\gamma}f},\psi\right]  =\left[  \left(  i\xi\right)
^{\gamma}\widehat{f},\psi\right]  =\left[  \left(  i\xi\right)  ^{\gamma
}\widehat{f},\mathcal{Q}_{\emptyset,n}\psi\right]  +\left[  \left(
i\xi\right)  ^{\gamma}\widehat{f},\mathcal{P}_{\emptyset,n}\psi\right]  .
\]

But%
\begin{align*}
\left[  \left(  i\xi\right)  ^{\gamma}\widehat{f},\mathcal{Q}_{\emptyset
,n}\psi\right]  =\left[  \widehat{f},\left(  i\xi\right)  ^{\gamma}%
\mathcal{Q}_{\emptyset,n}\psi\right]   &  =\left[  f_{F},\left(  i\xi\right)
^{\gamma}\mathcal{Q}_{\emptyset,n}\psi\right] \\
&  =\int\left(  i\xi\right)  ^{\gamma}\mathcal{Q}_{\emptyset,n}\psi\left(
\xi\right)  f_{F}\left(  \xi\right)  d\xi,
\end{align*}

and%
\[
\left[  \left(  i\xi\right)  ^{\gamma}\widehat{f},\mathcal{P}_{\emptyset
,n}\psi\right]  =\left[  \mathcal{P}_{\emptyset,n}^{\ast}\left(  \left(
i\xi\right)  ^{\gamma}\widehat{f}\right)  ,\psi\right]  =\left[
\mathcal{P}_{\emptyset,n}^{\ast}\widehat{D^{\gamma}f},\psi\right]  ,
\]

so that%
\begin{equation}
\left[  \widehat{D^{\gamma}f},\psi\right]  =\int\left(  i\xi\right)  ^{\gamma
}\mathcal{Q}_{\emptyset,n}\psi\left(  \xi\right)  f_{F}\left(  \xi\right)
d\xi+\left[  \mathcal{P}_{\emptyset,n}^{\ast}\widehat{D^{\gamma}f}%
,\psi\right]  .\label{a17}%
\end{equation}

The next step is to simplify the second term on the right. From part 1 of
Theorem \ref{Thm_P*_and_Q*}%
\[
\mathcal{P}_{\emptyset,n}^{\ast}\widehat{f}=\left(  \sum\limits_{\left\vert
\alpha\right\vert <n}\dfrac{b_{\alpha,\widehat{f}}}{\alpha!}\left(
-iD\right)  ^{\alpha}\right)  \delta,\quad b_{\alpha,\widehat{f}}=\left[
\widehat{f},\left(  -i\xi\right)  ^{\alpha}\rho\right]  .
\]

Hence%
\[
\left(  \mathcal{P}_{\emptyset,n}^{\ast}\widehat{f}\right)  ^{\vee}\left(
x\right)  =\left(  2\pi\right)  ^{-\frac{d}{2}}\sum\limits_{\left\vert
\alpha\right\vert <n}\dfrac{b_{\alpha,\widehat{f}}}{\alpha!}x^{\alpha}=\left(
2\pi\right)  ^{-\frac{d}{2}}p_{\widehat{f}}\left(  x\right)  ,
\]

so that%
\[
\left[  \mathcal{P}_{\emptyset,n}^{\ast}\widehat{f},\psi\right]  =\left[
\left(  \mathcal{P}_{\emptyset,n}^{\ast}\widehat{f}\right)  ^{\vee}%
,\widehat{\psi}\right]  =\left(  2\pi\right)  ^{-\frac{d}{2}}\left[
p_{\widehat{f}},\widehat{\psi}\right]  =\left(  2\pi\right)  ^{-\frac{d}{2}%
}\left[  \widehat{p_{\widehat{f}}},\psi\right]  ,
\]

and%
\begin{equation}
\mathcal{P}_{\emptyset,n}^{\ast}\widehat{D^{\gamma}f}=\left(  2\pi\right)
^{-\frac{d}{2}}\widehat{p_{\widehat{D^{\gamma}f}}},\label{a22}%
\end{equation}

so that \ref{a17} becomes \ref{a19}. It remains to show that $\widehat
{D^{\gamma}f}=\left(  i\xi\right)  ^{\gamma}f_{F}$ on $S_{\emptyset,n}$. In
fact, if $\phi\in S_{\emptyset,n}$ then, since $f_{F}=\widehat{f}$ on
$S_{\emptyset,n}$, $\left[  \widehat{D^{\gamma}f},\phi\right]  =\left[
\left(  i\xi\right)  ^{\gamma}\widehat{f},\phi\right]  =\left[  \widehat
{f},\left(  i\xi\right)  ^{\gamma}\phi\right]  =\left[  f_{F},\left(
i\xi\right)  ^{\gamma}\phi\right]  =\left[  \left(  i\xi\right)  ^{\gamma
}f_{F},\phi\right]  $ as required. Also, we are given that $f_{F}=\widehat{f}
$ on $\mathbb{R}^{d}\setminus0$ and $\widehat{f}\in L_{loc}^{1}\left(
\mathbb{R}^{d}\setminus0\right)  $ so $f_{F}\in L_{loc}^{1}\left(
\mathbb{R}^{d}\setminus0\right)  $ and $\left(  i\xi\right)  ^{\gamma}f_{F}\in
L_{loc}^{1}\left(  \mathbb{R}^{d}\setminus0\right)  $. \medskip

\textbf{Part 2} From \ref{a22},%
\begin{align*}
D^{\gamma}f=\left(  \widehat{D^{\gamma}f}\right)  ^{\vee}=\left(
\mathcal{Q}_{\emptyset,n}^{\ast}\widehat{D^{\gamma}f}\right)  ^{\vee}+\left(
\mathcal{P}_{\emptyset,n}^{\ast}\widehat{D^{\gamma}f}\right)  ^{\vee}  &
=\left(  \mathcal{Q}_{\emptyset,n}^{\ast}\widehat{D^{\gamma}f}\right)  ^{\vee
}+\left(  \left(  2\pi\right)  ^{-\frac{d}{2}}\widehat{p_{\widehat{D^{\gamma
}f}}}\right)  ^{\vee}\\
& =\left(  \mathcal{Q}_{\emptyset,n}^{\ast}\widehat{D^{\gamma}f}\right)
^{\vee}+\left(  2\pi\right)  ^{-\frac{d}{2}}p_{\widehat{D^{\gamma}f}}.
\end{align*}

Thus $f=\left(  \mathcal{Q}_{\emptyset,n}^{\ast}\widehat{f}\right)  ^{\vee
}+\left(  2\pi\right)  ^{-\frac{d}{2}}p_{\widehat{f}}$ and hence
\[
D^{\gamma}f=D^{\gamma}\left(  \mathcal{Q}_{\emptyset,n}^{\ast}\widehat
{f}\right)  ^{\vee}+\left(  2\pi\right)  ^{-\frac{d}{2}}D^{\gamma}%
p_{\widehat{f}}.
\]
\medskip

\textbf{Part 3} The next step is to prove equation \ref{a24}. If $\psi\in S$
then since $\mathcal{Q}_{\emptyset,n}+\mathcal{P}_{\emptyset,n}=I$
\begin{align*}
\left[  D^{\gamma}f,\psi\right]  =\left(  -1\right)  ^{\left\vert
\gamma\right\vert }\left[  f,D^{\gamma}\psi\right]  =\left(  -1\right)
^{\left\vert \gamma\right\vert }\left[  \widehat{f},\left(  D^{\gamma}%
\psi\right)  ^{\vee}\right]   &  =\left(  -1\right)  ^{\left\vert
\gamma\right\vert }\left[  \widehat{f},\mathcal{Q}_{\emptyset,n}\left(
\left(  D^{\gamma}\psi\right)  ^{\vee}\right)  \right]  +\\
&  \qquad+\left(  -1\right)  ^{\left\vert \gamma\right\vert }\left[
\widehat{f},\mathcal{P}_{\emptyset,n}\left(  \left(  D^{\gamma}\psi\right)
^{\vee}\right)  \right]  .
\end{align*}

By part 2 of Theorem \ref{Thm_S=So+pP2m-1} $\mathcal{Q}_{\emptyset
,n}:S\rightarrow S_{\emptyset,n}$ and since this theorem assumes that
$f_{F}=\widehat{f}$ on $S_{\emptyset,n}$ and that $\int f_{F}\phi$ on $\phi\in
S_{\emptyset,n}$ defines a member of $S_{\emptyset,n}^{\prime}$, we have
\[
\left[  \widehat{f},\mathcal{Q}_{\emptyset,n}\left(  \left(  D^{\gamma}%
\psi\right)  ^{\vee}\right)  \right]  =\left[  f_{F},\mathcal{Q}_{\emptyset
,n}\left(  \left(  D^{\gamma}\psi\right)  ^{\vee}\right)  \right]
=\int\mathcal{Q}_{\emptyset,n}\left(  \left(  D^{\gamma}\psi\right)  ^{\vee
}\right)  \left(  \xi\right)  f_{F}\left(  \xi\right)  d\xi.
\]

Also by \ref{a22} with $\gamma=0$%
\[
\left[  \widehat{f},\mathcal{P}_{\emptyset,n}\left(  \left(  D^{\gamma}%
\psi\right)  ^{\vee}\right)  \right]  =\left[  \left(  \mathcal{P}%
_{\emptyset,n}^{\ast}\widehat{f}\right)  ^{\vee},D^{\gamma}\psi\right]
=\left(  2\pi\right)  ^{-\frac{d}{2}}\left[  p_{\widehat{f}},D^{\gamma}%
\psi\right]  =\left(  2\pi\right)  ^{-\frac{d}{2}}\left[  \left(  -D\right)
^{\gamma}p_{\widehat{f}},\psi\right]  ,
\]

so that now%
\[
\left[  D^{\gamma}f,\psi\right]  =\left(  -1\right)  ^{\left\vert
\gamma\right\vert }\int\mathcal{Q}_{\emptyset,n}\left(  \left(  D^{\gamma}%
\psi\right)  ^{\vee}\right)  \left(  \xi\right)  f_{F}\left(  \xi\right)
d\xi+\left(  2\pi\right)  ^{-\frac{d}{2}}\left[  D^{\gamma}p_{\widehat{f}%
},\psi\right]  .
\]

From Lemma \ref{Lem_Pon(f)_in_terms_of_Pon(exp)}, for action variable $x$,
\begin{align*}
\mathcal{Q}_{\emptyset,n}\left(  \left(  D^{\gamma}\psi\right)  ^{\vee
}\right)  (\xi)  & =\left(  2\pi\right)  ^{-\frac{d}{2}}\int\mathcal{Q}%
_{\emptyset,n,\xi}\left(  e^{i\left(  x,\xi\right)  }\right)  D^{\gamma}%
\psi\left(  x\right)  dx\\
& =\left(  2\pi\right)  ^{-\frac{d}{2}}\left[  \mathcal{Q}_{\emptyset,n,\xi
}\left(  e^{i\left(  \cdot,\xi\right)  }\right)  ,D^{\gamma}\psi\right] \\
& =\left(  2\pi\right)  ^{-\frac{d}{2}}\left[  \left(  -D\right)  ^{\gamma
}\mathcal{Q}_{\emptyset,n,\xi}\left(  e^{i\left(  \cdot,\xi\right)  }\right)
,\psi\right] \\
& =\left(  2\pi\right)  ^{-\frac{d}{2}}\int\left(  -D\right)  _{x}^{\gamma
}\left(  \mathcal{Q}_{\emptyset,n,\xi}\left(  e^{i\left(  x,\xi\right)
}\right)  \right)  \psi\left(  x\right)  dx,
\end{align*}

so that%
\begin{align*}
\left[  D^{\gamma}f,\psi\right]   & =\left(  2\pi\right)  ^{-\frac{d}{2}}%
\int\int D_{x}^{\gamma}\left(  \mathcal{Q}_{\emptyset,\xi}\left(  e^{i\left(
x,\xi\right)  }\right)  \right)  \psi\left(  x\right)  dx\text{ }f_{F}\left(
\xi\right)  d\xi+\left(  2\pi\right)  ^{-\frac{d}{2}}\left[  D^{\gamma
}p_{\widehat{f}},\psi\right] \\
& =\left(  2\pi\right)  ^{-\frac{d}{2}}\int\int D_{x}^{\gamma}\left(
\mathcal{Q}_{\emptyset,\xi}\left(  e^{i\left(  x,\xi\right)  }\right)
\right)  f_{F}\left(  \xi\right)  \psi\left(  x\right)  dxd\xi+\left(
2\pi\right)  ^{-\frac{d}{2}}\left[  D^{\gamma}p_{\widehat{f}},\psi\right]  .
\end{align*}

But from the assumptions of this theorem we know the double integral is
absolutely integrable. Thus%
\[
\left[  D^{\gamma}f,\psi\right]  =\left(  2\pi\right)  ^{-\frac{d}{2}}%
\int\left(  \int D_{x}^{\gamma}\mathcal{Q}_{\emptyset,\xi}\left(  e^{i\left(
x,\xi\right)  }\right)  f_{F}\left(  \xi\right)  d\xi\right)  \psi\left(
x\right)  dx+\left(  2\pi\right)  ^{-\frac{d}{2}}\left[  D^{\gamma}%
p_{\widehat{f}},\psi\right]  ,
\]

and assumption \ref{a37} implies $\int D_{x}^{\gamma}\mathcal{Q}%
_{\emptyset,n,\xi}\left(  e^{i\left(  x,\xi\right)  }\right)  f_{F}\left(
\xi\right)  d\xi\in L_{loc}^{1}$, so we can conclude that%
\[
D^{\gamma}f\left(  x\right)  =\left(  2\pi\right)  ^{-\frac{d}{2}}\int
D_{x}^{\gamma}\mathcal{Q}_{\emptyset,n,\xi}\left(  e^{i\left(  x,\xi\right)
}\right)  f_{F}\left(  \xi\right)  d\xi+\left(  2\pi\right)  ^{-\frac{d}{2}%
}D^{\gamma}p_{\widehat{f}}\left(  x\right)  ,
\]

which proves \ref{a24}. In addition,

$D_{x}^{\gamma}\mathcal{Q}_{\emptyset,n,\xi}\left(  e^{i\left(  x,\xi\right)
}\right)  f_{F}\left(  \xi\right)  \in C^{\left(  0\right)  }$ for each $\xi$,
and by inequality \ref{a37},

$D_{x}^{\gamma}\mathcal{Q}_{\emptyset,n,\xi}\left(  e^{i\left(  x,\xi\right)
}\right)  f_{F}\left(  \xi\right)  \in L^{1}$ for each $x$. Thus by Lemma
\ref{Lem_diff_under_integral}, $D^{\gamma}f\in C^{\left(  0\right)  }$,

and the upper bound \ref{a37} applied to the right side of \ref{a24} implies
$D^{\gamma}f\in C_{BP}^{\left(  0\right)  }$.

Further, from \ref{a21}
\[
D_{x}^{\gamma}\mathcal{Q}_{\emptyset,n,\xi}\left(  e^{i\left(  x,\xi\right)
}\right)  =\left\{
\begin{array}
[c]{ll}%
\left(  i\xi\right)  ^{\gamma}\mathcal{Q}_{\emptyset,n-\left\vert
\gamma\right\vert ,\xi}\left(  e^{i\left(  x,\xi\right)  }\right)  , &
\left\vert \gamma\right\vert <n,\\
\left(  i\xi\right)  ^{\gamma}, & \left\vert \gamma\right\vert \geq n,
\end{array}
\right.
\]

which proves \ref{a26}.

In part 1$\ $it was shown that $\widehat{D^{\gamma}f}=\left(  i\xi\right)
^{\gamma}f_{F}$ on $S_{\emptyset,n}$ and $\left(  i\xi\right)  ^{\gamma}%
f_{F}\in L_{loc}^{1}\left(  \mathbb{R}^{d}\setminus0\right)  $. Thus
$\widehat{D^{\gamma}f}\in L_{loc}^{1}\left(  \mathbb{R}^{d}\setminus0\right)
$ and it is meaningful to define the function $\left(  D^{\gamma}f\right)
_{F}$. In fact $\left(  D^{\gamma}f\right)  _{F}=\left(  i\xi\right)
^{\gamma}f_{F}$ a.e. on $\mathbb{R}^{d}\setminus0$ and hence $\left(
D^{\gamma}f\right)  _{F}=\left(  i\xi\right)  ^{\gamma}f_{F}$ a.e. which means
the right side of \ref{a26} implies the right side of \ref{a36}.\medskip

\textbf{Part 4} From part 2, if $\psi\in S$ then%
\[
\left[  D^{\gamma}f,\psi\right]  =\left[  \left(  \mathcal{Q}_{\emptyset
,n}^{\ast}\widehat{D^{\gamma}f}\right)  ^{\vee},\psi\right]  +\left[  \left(
2\pi\right)  ^{-\frac{d}{2}}p_{\widehat{D^{\gamma}f}},\psi\right]  .
\]

But by Theorem \ref{Thm_S=So+pP2m-1}, $\mathcal{Q}_{\emptyset,n}:S\rightarrow
S_{\emptyset,n}$, and since $\widehat{f}=f_{F}$ on $S_{\emptyset,n}$ we
continue as:%
\begin{align*}
&  \left[  \left(  \mathcal{Q}_{\emptyset,n}^{\ast}\widehat{D^{\gamma}%
f}\right)  ^{\vee},\psi\right]  =\left[  \mathcal{Q}_{\emptyset,n}^{\ast
}\widehat{D^{\gamma}f},\overset{\vee}{\psi}\right]  =\left[  \widehat
{D^{\gamma}f},\mathcal{Q}_{\emptyset,n}\overset{\vee}{\psi}\right]  =\\
&  =\left[  \left(  i\xi\right)  ^{\gamma}\widehat{f},\mathcal{Q}%
_{\emptyset,n}\overset{\vee}{\psi}\right]  =\left[  \widehat{f},\left(
i\xi\right)  ^{\gamma}\mathcal{Q}_{\emptyset,n}\overset{\vee}{\psi}\right]
=\left[  \widehat{f},\left(  i\xi\right)  ^{\gamma}\mathcal{Q}_{\emptyset
,n}\overset{\vee}{\psi}\right]  =\\
&  =\left[  f_{F},\left(  i\xi\right)  ^{\gamma}\mathcal{Q}_{\emptyset
,n}\overset{\vee}{\psi}\right]  =\int\left(  i\xi\right)  ^{\gamma}\left(
\mathcal{Q}_{\emptyset,n}\overset{\vee}{\psi}\right)  \left(  \xi\right)
f_{F}\left(  \xi\right)  d\xi.
\end{align*}

and substituting part 2 of Lemma \ref{Lem_Pon(f)_in_terms_of_Pon(exp)} and
using \ref{a37} to change the order of integration by means of Fubini's
theorem we get
\begin{align*}
D^{\gamma}f\left(  x\right)   & =\left(  2\pi\right)  ^{-\frac{d}{2}}%
\int\left(  i\xi\right)  ^{\gamma}\mathcal{Q}_{\emptyset,n,\xi}\left(
e^{i\left(  x,\xi\right)  }\right)  f_{F}\left(  \xi\right)  d\xi+\left(
2\pi\right)  ^{-\frac{d}{2}}p_{\widehat{D^{\gamma}f}}.\left(  x\right) \\
& =\int\left(  \mathcal{Q}_{\emptyset,n}\overset{\vee}{\psi}\right)  \left(
\xi\right)  \left(  D^{\gamma}f\right)  _{F}\left(  \xi\right)  d\xi+\left(
2\pi\right)  ^{-\frac{d}{2}}p_{\widehat{D^{\gamma}f}}.\left(  x\right)  .
\end{align*}

\end{proof}

\begin{remark}
\label{Rem_Thm_integ(|DQ(exp)fF|)}\ 
\end{remark}

\begin{enumerate}
\item Equation \ref{a19} can be regarded as a modified inverse Fourier
transform for a subspace of the tempered distributions. Equation \ref{a24} can
be regarded as a `modified' inverse Fourier transform for a subspace of
$C_{BP}^{\left(  0\right)  }$. Instead of using the Fourier transform as the
argument its $L_{loc}^{1}$ restriction to $\mathbb{R}^{d}\setminus0$ is used,
and instead of the exponential $e^{i\left(  x,\xi\right)  }$, the exponential
has a Taylor series-like term subtracted.

\item In part 3 the polynomial term depends on the order of the derivative. In
part 4 the polynomial term is in general always of constant degree $n-1$.
\end{enumerate}

The next corollary will be directly applicable to both basis functions and to
functions belonging to $X_{w}^{\theta}$.

\begin{corollary}
\label{Cor_Thm_integ(|DQ(exp)fF|)}Suppose $m$ and $n$ are positive integers. Further:

\begin{enumerate}
\item Suppose $f\in S^{\prime}$ and $\widehat{f}\in L_{loc}^{1}\left(
\mathbb{R}^{d}\setminus0\right)  $. Define a.e. the function $f_{F}%
:\mathbb{R}^{d}\rightarrow\mathbb{C}$ by $f_{F}=\widehat{f}$ on $\mathbb{R}%
^{d}\setminus0$.

\item Assume the action $\int f_{F}\phi$ on $\phi\in S_{\emptyset,n}$ defines
a member of $S_{\emptyset,n}^{\prime}$, and that $f_{F}=\widehat{f} $ on
$S_{\emptyset,n}$.

\item Assume that if $\left\vert \gamma\right\vert \leq m$ there exist
$s_{\gamma}\geq0$ and $k_{\gamma}>0$ independent of $\xi$ such that%
\begin{equation}
\int\left\vert D_{x}^{\gamma}\mathcal{Q}_{\emptyset,n,\xi}\left(  e^{i\left(
x,\xi\right)  }\right)  f_{F}\left(  \xi\right)  \right\vert d\xi\leq
k_{\gamma}\left(  1+\left\vert x\right\vert \right)  ^{s_{\gamma}}.\label{a34}%
\end{equation}

\end{enumerate}

Let $f_{\rho}=\left(  \mathcal{Q}_{\emptyset,n}^{\ast}\widehat{f}\right)
^{\vee}\in S^{\prime}$ so that%
\begin{equation}
f=f_{\rho}+\left(  2\pi\right)  ^{-\frac{d}{2}}p_{\widehat{f}},\label{a31}%
\end{equation}

where $p_{\widehat{f}}\in P_{n-1}$ is defined by \ref{a20} and the subscript
$\rho\in S_{1,n}$ is the function used to define the operator $\mathcal{P}%
_{\emptyset,n}$.

Then $f_{\rho}\in C_{BP}^{\left(  m\right)  }\cap S^{\prime}$ and if
$\left\vert \gamma\right\vert \leq m$ then%
\begin{align}
D^{\gamma}f_{\rho}\left(  x\right)   & =\left(  2\pi\right)  ^{-\frac{d}{2}%
}\int D_{x}^{\gamma}\mathcal{Q}_{\emptyset,n,\xi}\left(  e^{i\left(
x,\xi\right)  }\right)  f_{F}\left(  \xi\right)  d\xi\label{a29}\\
& =\left\{
\begin{array}
[c]{ll}%
\left(  2\pi\right)  ^{-\frac{d}{2}}\int\mathcal{Q}_{\emptyset,n-\left\vert
\gamma\right\vert ,\xi}\left(  e^{i\left(  x,\xi\right)  }\right)  \left(
i\xi\right)  ^{\gamma}f_{F}\left(  \xi\right)  d\xi, & \left\vert
\gamma\right\vert <n,\\
\left(  2\pi\right)  ^{-\frac{d}{2}}\int e^{i\left(  x,\xi\right)  }\left(
i\xi\right)  ^{\gamma}f_{F}\left(  \xi\right)  d\xi, & \left\vert
\gamma\right\vert \geq n,
\end{array}
\right. \label{a30}\\
& =\left\{
\begin{array}
[c]{ll}%
\left(  2\pi\right)  ^{-\frac{d}{2}}\int\mathcal{Q}_{\emptyset,n-\left\vert
\gamma\right\vert ,\xi}\left(  e^{i\left(  x,\xi\right)  }\right)  \left(
D^{\gamma}f\right)  _{F}\left(  \xi\right)  d\xi, & \left\vert \gamma
\right\vert <n,\\
\left(  2\pi\right)  ^{-\frac{d}{2}}\int e^{i\left(  x,\xi\right)  }\left(
D^{\gamma}f\right)  _{F}\left(  \xi\right)  d\xi, & \left\vert \gamma
\right\vert \geq n,
\end{array}
\right. \label{a38}%
\end{align}

and%
\begin{equation}
\left\vert D^{\gamma}f_{\rho}\left(  x\right)  \right\vert \leq k_{\gamma
}\left(  1+\left\vert x\right\vert \right)  ^{s_{\gamma}},\quad\left\vert
\gamma\right\vert \leq m.\label{a32}%
\end{equation}

\end{corollary}

\begin{proof}
Equation \ref{a24} of Theorem \ref{Thm_integ(|DQ(exp)fF|)} can be written%
\[
D^{\gamma}\left(  f\left(  x\right)  -\left(  2\pi\right)  ^{-\frac{d}{2}%
}p_{\widehat{f}}\left(  x\right)  \right)  =\left(  2\pi\right)  ^{-\frac
{d}{2}}\int D_{x}^{\gamma}\mathcal{Q}_{\emptyset,n,\xi}\left(  e^{i\left(
x,\xi\right)  }\right)  f_{F}\left(  \xi\right)  d\xi,
\]

so the definition of $f_{\rho}$ implies equations \ref{a29} to \ref{a36}.
Theorem \ref{Thm_integ(|DQ(exp)fF|)} implies $D^{\gamma}f\in C^{\left(
0\right)  }$ for $\left\vert \gamma\right\vert <m$ i.e. $f\in C^{\left(
m\right)  }$. Further, the upper bound \ref{a34} applied to the right side of
equation \ref{a29} implies the upper bound \ref{a32} and hence that $f_{\rho
}\in C_{BP}^{\left(  m\right)  }$.
\end{proof}

\section{An inverse Fourier transform for basis
functions\label{Sect_explicit_basis}}

From Theorems \ref{Thm_basis_smth_W3.1} and \ref{Thm_basis_smth_W3.3} we see
that when a weight function has properties W3.1 or W3.3 the basis functions
have bounded derivatives. However, for property W3.2 we only know that $G\in
C_{BP}^{\left(  \left\lfloor 2\kappa\right\rfloor \right)  }$ i.e. $G$ has
polynomial growth at infinity but with no upper bound on the rate of growth
near infinity. This situation will be rectified in Theorem
\ref{Thm_deriv_Grho} where we will use the upper bounds for $\left\vert
D_{x}^{\gamma}\mathcal{Q}_{\emptyset,2\theta,\xi}\left(  e^{i\left(
x,\xi\right)  }\right)  \right\vert $ derived in the next lemma to give an
upper bound for the basis function growth near infinity and to prove some
inverse Fourier transform formulas for the basis functions and their
derivatives. We then explore some other consequences of these results.

But first we need the lemma which will supply the estimates \ref{a37} required
by Corollary \ref{Cor_Thm_integ(|DQ(exp)fF|)}.

?? \textbf{COMPARE\ WITH SECTION }\ref{Sect_Tay_expan_basis_W3.2}%
\textbf{\ BELOW} ???

See Theorem \ref{Thm_bnd_on_integ(|DQ(exp)|/w|x|)_radial} below which is the
radial version of this lemma.

\begin{lemma}
\label{Lem_bound_on_integ(|DQ(exp)|/w|x|)}Suppose the weight function $w$ has
properties W2.1 and W3.2 for order $\theta$ and parameter $\kappa$. Then if
$\left\vert \gamma\right\vert \leq\left\lfloor 2\kappa\right\rfloor $ there
exist constants $\left\{  \upsilon_{k,\gamma}^{\left(  w\right)  }\right\}
_{k=0}^{2\theta-\left\vert \gamma\right\vert }$, independent of $x$, such
that
\begin{equation}
\int\dfrac{\left\vert D_{x}^{\gamma}\mathcal{Q}_{\emptyset,2\theta,\xi}\left(
e^{i\left(  x,\xi\right)  }\right)  \right\vert }{w\left(  \xi\right)
\left\vert \xi\right\vert ^{2\theta}}d\xi\leq\left\{
\begin{array}
[c]{ll}%
\max\left\{  \underline{C}_{2\theta,r_{3}}^{\left(  \rho\right)  }%
,\overline{C}_{2\theta-\left\vert \gamma\right\vert ,r_{3}}^{\left(
\rho\right)  }\right\}  \sum\limits_{k\leq2\theta-\left\vert \gamma\right\vert
}\upsilon_{k,\gamma}^{\left(  w\right)  }\frac{\left\vert x\right\vert ^{k}%
}{k!}, & \left\vert \gamma\right\vert <2\theta,\\
\int\frac{\left\vert \xi^{\gamma}\right\vert }{w\left\vert \cdot\right\vert
^{2\theta}}, & \left\vert \gamma\right\vert \geq2\theta.
\end{array}
\right. \label{a33}%
\end{equation}

Here the $\upsilon_{k,\gamma}^{\left(  w\right)  }$ are given by \ref{a09} and
only depend on the weight function $w$ and on the parameters $\theta
,\kappa,r_{3}$ which define weight function property W3.2.

Finally, inequality \ref{a33} implies that for all multi-indexes $\gamma$,
there exist constants $C_{w,\gamma}>0$ such that%
\[
\int\dfrac{\left\vert D_{x}^{\gamma}\mathcal{Q}_{\emptyset,2\theta,\xi}\left(
e^{i\left(  x,\xi\right)  }\right)  \right\vert }{w\left(  \xi\right)
\left\vert \xi\right\vert ^{2\theta}}d\xi\leq\left\{
\begin{array}
[c]{ll}%
C_{w,\gamma}\left(  1+\left\vert x\right\vert \right)  ^{2\theta-\left\vert
\gamma\right\vert }, & \left\vert \gamma\right\vert <2\theta,\\
C_{w,\gamma}, & \left\vert \gamma\right\vert \geq2\theta.
\end{array}
\right.  ,
\]

where $C_{w,\gamma}=\int\frac{\left\vert \xi^{\gamma}\right\vert }{w\left\vert
\cdot\right\vert ^{2\theta}}$ when $\left\vert \gamma\right\vert \geq2\theta$.
\end{lemma}

\begin{proof}
There are two cases to be considered: $\left\vert \gamma\right\vert <2\theta$
and $\left\vert \gamma\right\vert \geq2\theta$. In both cases we will split
the range of integration into the two concentric regions defined by the sphere
$S\left(  0;r_{3}\right)  $.

We use the estimate of part 2 of Theorem \ref{Thm_bound_on_g(e,x)_2} for
$\left\vert D_{x}^{\gamma}\mathcal{Q}_{\emptyset,2\theta,\xi}\left(
e^{i\left(  x,\xi\right)  }\right)  \right\vert $ in the closed ball
$\overline{B}\left(  0;r_{3}\right)  $ and the estimate of part 1 of Theorem
\ref{Thm_bound_on_g(e,x)_2} outside the closed ball.\medskip

\fbox{\textbf{Case 1} $\left\vert \gamma\right\vert <2\theta$ and $\left\vert
\gamma\right\vert \leq\left\lfloor 2\kappa\right\rfloor $} Thus%
\begin{align*}
\int\limits_{\left\vert \cdot\right\vert \leq r_{3}}\frac{\left\vert
D_{x}^{\gamma}\mathcal{Q}_{\emptyset,2\theta,\xi}\left(  e^{ix\xi}\right)
\right\vert }{w\left(  \xi\right)  \left\vert \xi\right\vert ^{2\theta}}d\xi &
\leq\int\limits_{\left\vert \cdot\right\vert \leq r_{3}}\frac{\underline
{C}_{2\theta,r_{3}}^{\left(  \rho\right)  }\left\vert \cdot\right\vert
^{2\theta}\left(  \left\vert \xi^{\gamma}\right\vert +\frac{\left\vert
x\widehat{\xi}\right\vert ^{2\theta-\left\vert \gamma\right\vert }}{\left(
2\theta-\left\vert \gamma\right\vert \right)  !}\right)  }{w\left\vert
\cdot\right\vert ^{2\theta}}\\
& =\underline{C}_{2\theta,r_{3}}^{\left(  \rho\right)  }\left(  \int%
\limits_{\left\vert \cdot\right\vert \leq r_{3}}\frac{\left\vert \xi^{\gamma
}\right\vert }{w}+\frac{1}{\left(  2\theta-\left\vert \gamma\right\vert
\right)  !}\int\limits_{\left\vert \cdot\right\vert \leq r_{3}}\frac
{\left\vert x\widehat{\xi}\right\vert ^{2\theta-\left\vert \gamma\right\vert
}}{w}\right) \\
& =\underline{C}_{2\theta,r_{3}}^{\left(  \rho\right)  }\left(  \int%
\limits_{\left\vert \cdot\right\vert \leq r_{3}}\frac{\left\vert \xi^{\gamma
}\right\vert }{w}+\frac{\left\vert x\right\vert ^{2\theta-\left\vert
\gamma\right\vert }}{\left(  2\theta-\left\vert \gamma\right\vert \right)
!}\int\limits_{\left\vert \cdot\right\vert \leq r_{3}}\frac{\left\vert
\widehat{x}\widehat{\xi}\right\vert ^{2\theta-\left\vert \gamma\right\vert }%
}{w}\right)  .
\end{align*}

and the integrals exist since property W2.1 is $1/w\in L_{loc}^{1}$. Further%
\begin{align*}
\int\limits_{\left\vert \cdot\right\vert \geq r_{3}}\frac{\left\vert
D_{x}^{\gamma}\mathcal{Q}_{\emptyset,2\theta,\xi}\left(  e^{i\left(
x,\xi\right)  }\right)  \right\vert }{w\left(  \xi\right)  \left\vert
\xi\right\vert ^{2\theta}}d\xi & \leq\int\limits_{\left\vert \cdot\right\vert
\geq r_{3}}\frac{\overline{C}_{2\theta-\left\vert \gamma\right\vert ,r_{3}%
}^{\left(  \rho\right)  }\left(  1+\sum\limits_{k<2\theta-\left\vert
\gamma\right\vert }\frac{\left\vert x\widehat{\xi}\right\vert ^{k}}%
{k!}\right)  \left\vert \xi^{\gamma}\right\vert }{w\left\vert \cdot\right\vert
^{2\theta}}\\
& \leq\overline{C}_{2\theta-\left\vert \gamma\right\vert ,r_{3}}^{\left(
\rho\right)  }\int\limits_{\left\vert \cdot\right\vert \geq r_{3}}%
\frac{\left(  1+\sum\limits_{k<2\theta-\left\vert \gamma\right\vert }%
\frac{\left\vert x\widehat{\xi}\right\vert ^{k}}{k!}\right)  \left\vert
\xi^{\gamma}\right\vert }{w\left\vert \cdot\right\vert ^{2\theta}}\\
& =\overline{C}_{2\theta-\left\vert \gamma\right\vert ,r_{3}}^{\left(
\rho\right)  }\int\limits_{\left\vert \cdot\right\vert \geq r_{3}}%
\frac{\left\vert \xi^{\gamma}\right\vert }{w\left\vert \cdot\right\vert
^{2\theta}}+\overline{C}_{2\theta-\left\vert \gamma\right\vert ,r_{3}%
}^{\left(  \rho\right)  }\sum\limits_{k<2\theta-\left\vert \gamma\right\vert
}\frac{1}{k!}\int\limits_{\left\vert \cdot\right\vert \geq r_{3}}%
\frac{\left\vert x\widehat{\xi}\right\vert ^{k}\left\vert \xi^{\gamma
}\right\vert }{w\left\vert \cdot\right\vert ^{2\theta}}\\
& =\overline{C}_{2\theta-\left\vert \gamma\right\vert ,r_{3}}^{\left(
\rho\right)  }\int\limits_{\left\vert \cdot\right\vert \geq r_{3}}%
\frac{\left\vert \xi^{\gamma}\right\vert }{w\left\vert \cdot\right\vert
^{2\theta}}+\overline{C}_{2\theta-\left\vert \gamma\right\vert ,r_{3}%
}^{\left(  \rho\right)  }\sum\limits_{k<2\theta-\left\vert \gamma\right\vert
}\frac{\left\vert x\right\vert ^{k}}{k!}\int\limits_{\left\vert \cdot
\right\vert \geq r_{3}}\frac{\left\vert \widehat{x}\widehat{\xi}\right\vert
^{k}\left\vert \xi^{\gamma}\right\vert }{w\left\vert \cdot\right\vert
^{2\theta}},
\end{align*}

And by Theorem \ref{Thm_property_wt_fn_W3} the last integral exists. Adding
the two estimates we get%
\begin{align*}
& \int\frac{\left\vert D_{x}^{\gamma}\mathcal{Q}_{\emptyset,2\theta,\xi
}\left(  e^{ix\xi}\right)  \right\vert }{w\left(  \xi\right)  \left\vert
\xi\right\vert ^{2\theta}}d\xi\\
& \leq\underline{C}_{2\theta,r_{3}}^{\left(  \rho\right)  }\left(
\int\limits_{\left\vert \cdot\right\vert \leq r_{3}}\frac{\left\vert
\xi^{\gamma}\right\vert }{w}+\frac{\left\vert x\right\vert ^{2\theta
-\left\vert \gamma\right\vert }}{\left(  2\theta-\left\vert \gamma\right\vert
\right)  !}\int\limits_{\left\vert \cdot\right\vert \leq r_{3}}\frac
{\left\vert \widehat{x}\widehat{\xi}\right\vert ^{2\theta-\left\vert
\gamma\right\vert }}{w}\right)  +\\
& +\overline{C}_{2\theta-\left\vert \gamma\right\vert ,r_{3}}^{\left(
\rho\right)  }\int\limits_{\left\vert \cdot\right\vert \geq r_{3}}%
\frac{\left\vert \xi^{\gamma}\right\vert }{w\left\vert \cdot\right\vert
^{2\theta}}+\overline{C}_{2\theta-\left\vert \gamma\right\vert ,r_{3}%
}^{\left(  \rho\right)  }\sum\limits_{k<2\theta-\left\vert \gamma\right\vert
}\frac{\left\vert x\right\vert ^{k}}{k!}\int\limits_{\left\vert \cdot
\right\vert \geq r_{3}}\frac{\left\vert \widehat{x}\widehat{\xi}\right\vert
^{k}\left\vert \xi^{\gamma}\right\vert }{w\left\vert \cdot\right\vert
^{2\theta}}\\
& =\left(  \underline{C}_{2\theta,r_{3}}^{\left(  \rho\right)  }%
\int\limits_{\left\vert \cdot\right\vert \leq r_{3}}\frac{\left\vert
\xi^{\gamma}\right\vert }{w}+2\overline{C}_{2\theta-\left\vert \gamma
\right\vert ,r_{3}}^{\left(  \rho\right)  }\int\limits_{\left\vert
\cdot\right\vert \geq r_{3}}\frac{\left\vert \xi^{\gamma}\right\vert
}{w\left\vert \cdot\right\vert ^{2\theta}}\right)  +\sum\limits_{k<2\theta
-\left\vert \gamma\right\vert }\left(  \overline{C}_{2\theta-\left\vert
\gamma\right\vert ,r_{3}}^{\left(  \rho\right)  }\int\limits_{\left\vert
\cdot\right\vert \geq r_{3}}\frac{\left\vert \widehat{x}\widehat{\xi
}\right\vert ^{k}\left\vert \xi^{\gamma}\right\vert }{w\left\vert
\cdot\right\vert ^{2\theta}}\right)  \frac{\left\vert x\right\vert ^{k}}%
{k!}+\\
& \qquad+\left(  \underline{C}_{2\theta,r_{3}}^{\left(  \rho\right)  }%
\int\limits_{\left\vert \cdot\right\vert \leq r_{3}}\frac{\left\vert
\widehat{x}\widehat{\xi}\right\vert ^{2\theta-\left\vert \gamma\right\vert }%
}{w}\right)  \frac{\left\vert x\right\vert ^{2\theta-\left\vert \gamma
\right\vert }}{\left(  2\theta-\left\vert \gamma\right\vert \right)  !}\\
& \leq\max\left\{  \underline{C}_{2\theta,r_{3}}^{\left(  \rho\right)
},\overline{C}_{2\theta-\left\vert \gamma\right\vert ,r_{3}}^{\left(
\rho\right)  }\right\}  \times\\
& \times\left(
\begin{array}
[c]{c}%
\left(  \int\limits_{\left\vert \cdot\right\vert \leq r_{3}}\frac{\left\vert
\xi^{\gamma}\right\vert }{w}+2\int\limits_{\left\vert \cdot\right\vert \geq
r_{3}}\frac{\left\vert \xi^{\gamma}\right\vert }{w\left\vert \cdot\right\vert
^{2\theta}}\right)  +\sum\limits_{k<2\theta-\left\vert \gamma\right\vert
}\left(  \int\limits_{\left\vert \cdot\right\vert \geq r_{3}}\frac{\left\vert
\widehat{x}\widehat{\xi}\right\vert ^{k}\left\vert \xi^{\gamma}\right\vert
}{w\left\vert \cdot\right\vert ^{2\theta}}\right)  \frac{\left\vert
x\right\vert ^{k}}{k!}+\\
+\left(  \int\limits_{\left\vert \cdot\right\vert \leq r_{3}}\frac{\left\vert
\widehat{x}\widehat{\xi}\right\vert ^{2\theta-\left\vert \gamma\right\vert }%
}{w}\right)  \frac{\left\vert x\right\vert ^{2\theta-\left\vert \gamma
\right\vert }}{\left(  2\theta-\left\vert \gamma\right\vert \right)  !}%
\end{array}
\right)  ,
\end{align*}

so that%
\begin{equation}
\int\frac{\left\vert D_{x}^{\gamma}\mathcal{Q}_{\emptyset,2\theta,\xi}\left(
e^{ix\xi}\right)  \right\vert }{w\left(  \xi\right)  \left\vert \xi\right\vert
^{2\theta}}d\xi\leq\max\left\{  \underline{C}_{2\theta,r_{3}}^{\left(
\rho\right)  },\overline{C}_{2\theta-\left\vert \gamma\right\vert ,r_{3}%
}^{\left(  \rho\right)  }\right\}  \sum\limits_{k\leq2\theta-\left\vert
\gamma\right\vert }\upsilon_{k,\gamma}^{\left(  w\right)  }\frac{\left\vert
x\right\vert ^{k}}{k!},\quad\left\vert \gamma\right\vert <2\theta,\text{
}\left\vert \gamma\right\vert \leq\left\lfloor 2\kappa\right\rfloor
\label{a08}%
\end{equation}

where%
\begin{equation}
\upsilon_{k,\gamma}^{\left(  w\right)  }=\left\{
\begin{array}
[c]{ll}%
\int\limits_{\left\vert \cdot\right\vert \leq r_{3}}\frac{\left\vert
\xi^{\gamma}\right\vert }{w}+2\int\limits_{\left\vert \cdot\right\vert \geq
r_{3}}\frac{\left\vert \xi^{\gamma}\right\vert }{w\left\vert \cdot\right\vert
^{2\theta}}, & k=0,\\
\int\limits_{\left\vert \cdot\right\vert \geq r_{3}}\frac{\left\vert
\widehat{x}\widehat{\xi}\right\vert ^{k}\left\vert \xi^{\gamma}\right\vert
}{w\left\vert \cdot\right\vert ^{2\theta}}, & 0<k<2\theta-\left\vert
\gamma\right\vert ,\\
\int\limits_{\left\vert \cdot\right\vert \leq r_{3}}\frac{\left\vert
\widehat{x}\widehat{\xi}\right\vert ^{2\theta-\left\vert \gamma\right\vert }%
}{w}, & k=2\theta-\left\vert \gamma\right\vert .
\end{array}
\right.  .\label{a09}%
\end{equation}
\medskip

\fbox{\textbf{Case 2} $\left\vert \gamma\right\vert \geq2\theta$ and
$\left\vert \gamma\right\vert \leq\left\lfloor 2\kappa\right\rfloor $} From
Theorem \ref{Thm_bound_on_g(e,x)_2},%
\[
\int\frac{\left\vert D_{x}^{\gamma}\mathcal{Q}_{\emptyset,\xi}\left(
e^{i\left(  x,\xi\right)  }\right)  \right\vert }{w\left(  \xi\right)
\left\vert \xi\right\vert ^{2\theta}}\leq\int\frac{\left\vert \xi^{\gamma
}\right\vert }{w\left\vert \cdot\right\vert ^{2\theta}}.
\]

\end{proof}

If we only assume the weight function has property W2 then part 1 of Theorem
\ref{Thm_integ(|DQ(exp)fF|)} allows us to convert the definition \ref{p01} of
a basis distribution $G\in S^{\prime}$ i.e. $\left[  \widehat{G},\phi\right]
=\int\frac{\phi}{w\left\vert \cdot\right\vert ^{2\theta}}$ for all $\phi\in
S_{\emptyset,2\theta}$, into the explicit formulas \ref{a23} of the next theorem.

\begin{theorem}
\label{Thm_basis_distrib_properties}Suppose the weight function $w$ has
property W2 and let $G$ be a basis \textbf{distribution} of order $\theta
\geq1$ generated by $w$. Then for any multi-index $\gamma$%
\begin{equation}
\left[  \widehat{D^{\gamma}G},\phi\right]  =\int\left(  i\xi\right)  ^{\gamma
}\dfrac{\mathcal{Q}_{\emptyset,2\theta}\phi\left(  \xi\right)  }{w\left(
\xi\right)  \left\vert \xi\right\vert ^{2\theta}}d\xi+\left(  2\pi\right)
^{-\frac{d}{2}}\left[  \left(  p_{\widehat{D^{\gamma}G}}\right)  ^{\wedge
},\phi\right]  ,\quad\phi\in S,\label{a23}%
\end{equation}

where for $u\in S^{\prime}$%
\begin{equation}
p_{u}\left(  x\right)  =\sum_{\left\vert \alpha\right\vert <2\theta}%
\frac{b_{u,\alpha}}{\alpha!}x^{\alpha},\quad b_{u,\alpha}=\left[  u,\left(
-i\xi\right)  ^{\alpha}\rho\right]  ,\label{a45}%
\end{equation}

and $\xi$ is the action variable.
\end{theorem}

\begin{proof}
From weight function property W2.1, $\frac{1}{w}\in L_{loc}^{1}$ and so by
Definition \ref{Def_basis_distrib} we have $G\in S^{\prime}$, $\widehat
{G}=\frac{1}{w\left\vert \cdot\right\vert ^{2\theta}}$ on $S_{\emptyset
,2\theta}$ and $\frac{1}{w\left\vert \cdot\right\vert ^{2\theta}}\in
L_{loc}^{1}\left(  \mathbb{R}^{d}\setminus0\right)  \cap S_{\emptyset,2\theta
}^{\prime}$. Now set $G_{F}=\widehat{G}$ on $\mathbb{R}^{d}\setminus0$ so that
$G_{F}=\frac{1}{w\left\vert \cdot\right\vert ^{2\theta}}\in L_{loc}^{1}\left(
\mathbb{R}^{d}\setminus0\right)  \cap S_{\emptyset,2\theta}^{\prime}$.
Equation \ref{a23} now follows from equation \ref{a19} of Theorem
\ref{Thm_integ(|DQ(exp)fF|)} with $f=G$.
\end{proof}

\begin{remark}
This result is closely related to Theorem 2.1 of Madych and Nelson
\cite{MadychNelson90}. Indeed, in the comments following Theorem 2.1 Madych
and Nelson illustrate Theorem 2.1 by choosing $d\mu\left(  \xi\right)
=w\left(  \xi\right)  d\xi$ where $w$ corresponds to $1/w$ in this document
i.e. to the \textbf{reciprocal} of our weight function.
\end{remark}

If we further assume that the weight function has property W3.2 then we can
apply the bounds derived in the previous lemma to Corollary
\ref{Cor_Thm_integ(|DQ(exp)fF|)} and show that every basis function $G$ is the
sum of a polynomial in $P_{2\theta-1}$ and a special basis function $G_{\rho}$
which depends only on the choice of a function $\rho\in S_{1,n}$. In turn
$G_{\rho}$ satisfies the modified inverse-Fourier transform equations
\ref{a47} and the growth estimates \ref{a46}.

\begin{theorem}
\label{Thm_deriv_Grho}Suppose the weight function $w$ has property W2.1 and
property \textbf{W3.2} for order $\theta$ and parameter $\kappa$. Further
suppose that $G$ is a basis distribution of order $\theta$ and the operator
$\mathcal{P}_{\emptyset,\theta}$ is defined using the function $\rho\in
S_{1,2\theta}$.

Then $G\in C_{BP}^{\left(  \left\lfloor 2\kappa\right\rfloor \right)  }$ and
we define the non-unique basis distribution $G_{\rho}$ by
\begin{equation}
G_{\rho}:=G-\left(  2\pi\right)  ^{-\frac{d}{2}}p_{\widehat{G}},\label{a44}%
\end{equation}

where $p_{\widehat{G}}\in P_{2\theta-1}$ is defined using \ref{a45}. For all
$\left\vert \gamma\right\vert \leq\left\lfloor 2\kappa\right\rfloor $:
\begin{align}
D^{\gamma}G_{\rho}\left(  x\right)   & =\left(  2\pi\right)  ^{-\frac{d}{2}%
}\int\frac{D_{x}^{\gamma}\mathcal{Q}_{\emptyset,2\theta,\xi}\left(
e^{i\left(  x,\xi\right)  }\right)  }{w\left(  \xi\right)  \left\vert
\xi\right\vert ^{2\theta}}d\xi\label{a43}\\
& =\left\{
\begin{array}
[c]{ll}%
\left(  2\pi\right)  ^{-\frac{d}{2}}\int\mathcal{Q}_{\emptyset,2\theta
-\left\vert \gamma\right\vert ,\xi}\left(  e^{i\left(  x,\xi\right)  }\right)
\frac{\left(  i\xi\right)  ^{\gamma}}{w\left(  \xi\right)  \left\vert
\xi\right\vert ^{2\theta}}d\xi, & \left\vert \gamma\right\vert <2\theta,\\
\left(  2\pi\right)  ^{-\frac{d}{2}}\int e^{i\left(  x,\xi\right)  }%
\frac{\left(  i\xi\right)  ^{\gamma}}{w\left(  \xi\right)  \left\vert
\xi\right\vert ^{2\theta}}d\xi, & \left\vert \gamma\right\vert \geq2\theta,
\end{array}
\right. \label{a431}\\
& =\left\{
\begin{array}
[c]{ll}%
\left(  2\pi\right)  ^{-\frac{d}{2}}\int\mathcal{Q}_{\emptyset,n-\left\vert
\gamma\right\vert ,\xi}\left(  e^{i\left(  x,\xi\right)  }\right)  \left(
D^{\gamma}G_{\rho}\right)  _{F}\left(  \xi\right)  d\xi, & \left\vert
\gamma\right\vert <2\theta,\\
\left(  2\pi\right)  ^{-\frac{d}{2}}\int e^{i\left(  x,\xi\right)  }\left(
D^{\gamma}G_{\rho}\right)  _{F}\left(  \xi\right)  d\xi, & \left\vert
\gamma\right\vert \geq2\theta,
\end{array}
\right. \label{a47}%
\end{align}

and $D^{\gamma}G_{\rho}$ satisfies the growth estimate%
\begin{equation}
\left\vert D^{\gamma}G_{\rho}\left(  x\right)  \right\vert \leq\left\{
\begin{array}
[c]{ll}%
C_{w,\gamma}\left(  1+\left\vert x\right\vert \right)  ^{2\theta-\left\vert
\gamma\right\vert }, & \left\vert \gamma\right\vert <2\theta,\\
C_{w,\gamma}, & \left\vert \gamma\right\vert \geq2\theta,
\end{array}
\right. \label{a46}%
\end{equation}

where the constant $C_{w}$ is given in Lemma
\ref{Lem_bound_on_integ(|DQ(exp)|/w|x|)}.
\end{theorem}

\begin{proof}
Since $p_{\widehat{G}}\in P_{\theta-1}$ the basis distribution Definition
\ref{Def_basis_distrib} implies $G_{\rho}$ is also a basis distribution, and
by part 4 of Theorem \ref{Thm_basis_smth_W3.2_r3_pos} it lies in
$C_{BP}^{\left(  \left\lfloor 2\kappa\right\rfloor \right)  }$.

We now want to apply Corollary \ref{Cor_Thm_integ(|DQ(exp)fF|)} with
$f=G$.\ However, inspection of the proof of Theorem
\ref{Thm_basis_distrib_properties} and the inequality \ref{a46} shows that
when $m=\left\lfloor 2\kappa\right\rfloor $, $n=2\theta$, $f_{F}=\frac
{1}{w\left\vert \cdot\right\vert ^{2\theta}}$, $k_{\gamma}=C_{w}$ and
$s_{\gamma}=2\theta-\left\vert \gamma\right\vert $ for $\left\vert
\gamma\right\vert \leq\left\lfloor 2\kappa\right\rfloor $, all the conditions
of Corollary \ref{Cor_Thm_integ(|DQ(exp)fF|)} are satisfied and we have our result.
\end{proof}

The special basis function $G_{\rho}$ defined by \ref{a44} of the last theorem
will enable us to relate various weight function properties to the
corresponding basis function properties e.g. the weight function is radial
implies the function $G_{\rho}$ is radial. But first we need a lemma.

\begin{lemma}
\label{Lem_radial_func}This lemma will employ the technique of Section 4.1,
Stein and Weiss \cite{SteinWeiss71} which defines radial functions in terms of
orthogonal transformations.

Stein and Weiss observe that a function $f$ is radial if and only if $f\left(
\mathcal{O}x\right)  =f\left(  x\right)  $ for any linear, orthogonal
transformation $\mathcal{O}:\mathbb{R}^{d}\rightarrow\mathbb{R}^{d}$ and any
$x\in\mathbb{R}^{d}$.

Note that an orthogonal transformation $\mathcal{O}$ satisfies $\mathcal{O}%
^{T}=\mathcal{O}^{-1}$ where $\left(  \mathcal{O}x,y\right)  =\left(
x,\mathcal{O}^{T}y\right)  $ for the Euclidean inner product. An orthogonal
transformation has a Jacobian of one.
\end{lemma}

\ 

\begin{theorem}
\label{Thm_Grho}Suppose the weight function $w$ has properties W2.1 and W3.2
for order $\theta$ and $\kappa$. Then the special basis function $G_{\rho} $
of order $\theta$ introduced in Theorem \ref{Thm_deriv_Grho} has the following properties:

\begin{enumerate}
\item If $\rho_{1}$ and $\rho_{2}$ are in $S_{1,2\theta}$ then $\rho_{2}%
-\rho_{1}\in S_{\emptyset,2\theta}$ and
\[
G_{\rho_{1}}\left(  x\right)  -G_{\rho_{2}}\left(  x\right)  =\left(
2\pi\right)  ^{-d/2}\sum_{\left\vert \alpha\right\vert <2\theta}\left(
\int\frac{\rho_{2}(\xi)-\rho_{1}(\xi)}{w\left(  \xi\right)  \left\vert
\xi\right\vert ^{2\theta}}\left(  i\xi\right)  ^{\alpha}d\xi\right)
\frac{x^{\alpha}}{\alpha!}\in P_{2\theta-1}.
\]

\item $G_{\rho}\left(  -x\right)  =\overline{G_{\rho}\left(  x\right)  }$ i.e.
$G_{\rho}$ is conjugate-even.

\item If $w$ and $\rho$ are even functions then $G_{\rho}$ is a real valued
even function.

\item If $w$ and $\rho$ are radial functions then $G_{\rho}$ is also a radial function.

\item If $w$ is homogeneous of order $s$ then%
\begin{align*}
G_{\rho}\left(  tx\right)   & =t^{s+2\theta-d}G_{\rho\left(  t^{-1}\xi\right)
}\left(  x\right) \\
& =t^{s+2\theta-d}G_{\rho}\left(  x\right)  +\frac{t^{s+2\theta-d}}{\left(
2\pi\right)  ^{d/2}}\sum_{\left\vert \alpha\right\vert <2\theta}\left(
\int\frac{\rho(\xi)-\rho(t^{-1}\xi)}{w\left(  \xi\right)  \left\vert
\xi\right\vert ^{2\theta}}\left(  i\xi\right)  ^{\alpha}d\xi\right)
\frac{x^{\alpha}}{\alpha!}\\
& =t^{s+2\theta-d}G_{\rho}\left(  x\right)  +\frac{1}{\left(  2\pi\right)
^{d/2}}\sum_{\left\vert \alpha\right\vert <2\theta}\left(  \int\frac
{\rho(t\eta)-\rho(\eta)}{w\left(  \eta\right)  \left\vert \eta\right\vert
^{2\theta}}\left(  i\eta\right)  ^{\alpha}d\eta\right)  \frac{\left(
tx\right)  ^{\alpha}}{\alpha!}.
\end{align*}

\end{enumerate}
\end{theorem}

\begin{proof}
\textbf{Part 1} By part 5 of Theorem \ref{Thm_bound_on_g(e,x)} the integrands
defining $G_{\rho_{1}}$ and $G_{\rho_{2}}$ are absolutely integrable. Hence
\begin{align*}
\left(  2\pi\right)  ^{d/2}\left(  G_{\rho_{1}}\left(  x\right)  -G_{\rho_{2}%
}\left(  x\right)  \right)   &  =\int\left(  e^{i\left(  x,\xi\right)  }%
-\rho_{1}(\xi)\sum_{\left\vert \alpha\right\vert <2\theta}\frac{\left(
ix\right)  ^{\alpha}\xi^{\alpha}}{\alpha!}\right)  \frac{d\xi}{w\left(
\xi\right)  \left\vert \xi\right\vert ^{2\theta}}-\\
&  \qquad-\int\left(  e^{i\left(  x,\xi\right)  }-\rho_{2}(\xi)\sum
_{\left\vert \alpha\right\vert <2\theta}\frac{x^{\alpha}\left(  i\xi\right)
^{\alpha}}{\alpha!}\right)  \frac{d\xi}{w\left(  \xi\right)  \left\vert
\xi\right\vert ^{2\theta}}\\
&  =\int\left(  \left(  \rho_{2}(\xi)-\rho_{1}(\xi)\right)  \sum_{\left\vert
\alpha\right\vert <2\theta}\frac{x^{\alpha}\left(  i\xi\right)  ^{\alpha}%
}{\alpha!}\right)  \frac{d\xi}{w\left(  \xi\right)  \left\vert \xi\right\vert
^{2\theta}}\\
&  =\int\sum_{\left\vert \alpha\right\vert <2\theta}\frac{\rho_{2}(\xi
)-\rho_{1}(\xi)}{w\left(  \xi\right)  \left\vert \xi\right\vert ^{2\theta}%
}\frac{x^{\alpha}\left(  i\xi\right)  ^{\alpha}}{\alpha!}d\xi\\
&  =\sum_{\left\vert \alpha\right\vert <2\theta}\left(  \int\frac{\rho_{2}%
(\xi)-\rho_{1}(\xi)}{w\left(  \xi\right)  \left\vert \xi\right\vert ^{2\theta
}}\left(  i\xi\right)  ^{\alpha}d\xi\right)  \frac{x^{\alpha}}{\alpha!},
\end{align*}

if the integrals on the last line all exist. However, by parts 3 and 4 of
Theorem \ref{Thm_product_of_Co,k_funcs}, $\rho_{2}(\xi)-\rho_{1}(\xi)\in
S_{\emptyset,2\theta}$ implies $\left(  \rho_{2}(\xi)-\rho_{1}(\xi)\right)
\left(  i\xi\right)  ^{\alpha}\in S_{\emptyset,2\theta}$. But Theorem
\ref{Thm_1/w|.|^2u_in_So,2u'} showed that if $\phi\in S_{\emptyset,2\theta}$
then $\int\frac{\left\vert \phi(\xi)\right\vert }{w\left(  \xi\right)
\left\vert \xi\right\vert ^{2\theta}}<\infty$ and so $\frac{\rho_{2}(\xi
)-\rho_{1}(\xi)}{w\left(  \xi\right)  \left\vert \xi\right\vert ^{2\theta}%
}\left(  i\xi\right)  ^{\alpha}\in L^{1}$.\medskip

\textbf{Part 2} From \ref{a43} and Definition \ref{Def_S1,n},%
\[
\overline{G_{\rho}\left(  x\right)  }=\left(  2\pi\right)  ^{-\frac{d}{2}}%
\int\left(  e^{-i\left(  x,\xi\right)  }-\rho(\xi)\sum_{\left\vert
\alpha\right\vert <2\theta}\frac{\left(  -ix\right)  ^{\alpha}\xi^{\alpha}%
}{\alpha!}\right)  \frac{d\xi}{w\left(  \xi\right)  \left\vert \xi\right\vert
^{2\theta}}=G_{\rho}\left(  -x\right)  .
\]

\textbf{Part 3} Suppose $w$ and $\rho$ are even functions. Then $\rho
(-x)=\rho(x)$ and
\begin{align*}
\left(  2\pi\right)  ^{\frac{d}{2}}G_{\rho}\left(  -x\right)   & =\int\left(
e^{-i\left(  x,\xi\right)  }-\rho(\xi)\sum_{\left\vert \alpha\right\vert
<2\theta}\frac{\left(  -ix\right)  ^{\alpha}\xi^{\alpha}}{\alpha!}\right)
\frac{d\xi}{w\left(  \xi\right)  \left\vert \xi\right\vert ^{2\theta}}\\
& =\int\left(  e^{-i\left(  x,\xi\right)  }-\rho(-\xi)\sum_{\left\vert
\alpha\right\vert <2\theta}\frac{\left(  ix\right)  ^{\alpha}\left(
-\xi\right)  ^{\alpha}}{\alpha!}\right)  \frac{d\xi}{w\left(  -\xi\right)
\left\vert \xi\right\vert ^{2\theta}}\\
& =\int\left(  e^{i\left(  x,\xi\right)  }-\rho(\xi)\sum_{\left\vert
\alpha\right\vert <2\theta}\frac{\left(  ix\right)  ^{\alpha}\xi^{\alpha}%
}{\alpha!}\right)  \frac{d\xi}{w\left(  \xi\right)  \left\vert \xi\right\vert
^{2\theta}}\\
& =\left(  2\pi\right)  ^{d/2}G_{\rho}\left(  x\right)  .
\end{align*}

Part 2 implies $G_{\rho}$ is real.\medskip

\textbf{Part 4} First note that $S_{1,2\theta}$ contains radial functions e.g.
the standard example of a distribution test function, the
`cap-shaped'\ function%
\[
\left\{
\begin{array}
[c]{ll}%
e\exp\left(  -\frac{1}{1-\left\vert x\right\vert ^{2}}\right)  , & \left\vert
x\right\vert \leq1,\\
0, & \left\vert x\right\vert >1,
\end{array}
\right.
\]

lies in
\[
S_{1,\infty}:=\left\{  \phi\in S:\phi\left(  0\right)  =1,\text{ }D^{\alpha
}\phi\left(  0\right)  =0\text{ }for\text{ }all\text{ }\alpha\neq0\right\}  .
\]

We now use the definition of a radial function given in Lemma
\ref{Lem_radial_func} to prove that $G_{\rho}$ is radial if $\rho$ is radial.
Now
\[
G_{\rho}\left(  x\right)  =\left(  2\pi\right)  ^{-\frac{d}{2}}\int\left(
e^{i\left(  x,\xi\right)  }-\rho(\xi)\sum_{k=0}^{2\theta}\frac{i^{k}\left(
x,\xi\right)  ^{k}}{k!}\right)  \frac{d\xi}{w\left(  \xi\right)  \left\vert
\xi\right\vert ^{2\theta}}.
\]

Thus
\begin{align*}
G_{\rho}\left(  \mathcal{O}x\right)   &  =\left(  2\pi\right)  ^{-\frac{d}{2}%
}\int\left(  e^{i\left(  \mathcal{O}x,\xi\right)  }-\rho(\xi)\sum
_{k=0}^{2\theta}\frac{i^{k}\left(  \mathcal{O}x,\xi\right)  ^{k}}{k!}\right)
\frac{d\xi}{w\left(  \xi\right)  \left\vert \xi\right\vert ^{2\theta}}\\
&  =\left(  2\pi\right)  ^{-\frac{d}{2}}\int\left(  e^{i\left\langle
x,\mathcal{O}^{T}\xi\right\rangle }-\rho(\xi)\sum_{k=0}^{2\theta}\frac
{i^{k}\left(  x,\mathcal{O}^{T}\xi\right)  ^{k}}{k!}\right)  \frac{d\xi
}{w\left(  \xi\right)  \left\vert \xi\right\vert ^{2\theta}},
\end{align*}

and then making the change of variables $\eta=\mathcal{O}^{T}\xi$ with
$d\eta=d\xi$ and $\xi=\mathcal{O}\eta$
\begin{align*}
G_{\rho}\left(  \mathcal{O}x\right)   &  =\left(  2\pi\right)  ^{-\frac{d}{2}%
}\int\left(  e^{i\left\langle x,\eta\right\rangle }-\rho(\mathcal{O}\eta
)\sum_{k=0}^{2\theta}\frac{i^{k}\left(  x,\eta\right)  ^{k}}{k!}\right)
\frac{d\eta}{w\left(  \mathcal{O}\eta\right)  \left\vert \mathcal{O}%
\eta\right\vert ^{2\theta}}\\
&  =\left(  2\pi\right)  ^{-\frac{d}{2}}\int\left(  e^{i\left\langle
x,\eta\right\rangle }-\rho(\eta)\sum_{k=0}^{2\theta}\frac{i^{k}\left(
x,\eta\right)  ^{k}}{k!}\right)  \frac{d\eta}{w\left(  \mathcal{O}\eta\right)
\left\vert \mathcal{O}\eta\right\vert ^{2\theta}}\\
&  =G_{\rho}\left(  x\right)  ,
\end{align*}

since $\rho$ and $w$ are radial. Hence $G_{\rho}$ is also radial.\medskip

\textbf{Part 5} If $t\in\mathbb{R}_{+}^{1}$ then%
\begin{align*}
G_{\rho}\left(  tx\right)   & =\left(  2\pi\right)  ^{-\frac{d}{2}}\int\left(
e^{i\left(  tx,\xi\right)  }-\rho(\xi)\sum_{\left\vert \alpha\right\vert
<2\theta}\frac{\left(  itx\right)  ^{\alpha}\xi^{\alpha}}{\alpha!}\right)
\frac{d\xi}{w\left(  \xi\right)  \left\vert \xi\right\vert ^{2\theta}}\\
& =\left(  2\pi\right)  ^{-\frac{d}{2}}\int\left(  e^{i\left(  x,t\xi\right)
}-\rho(\xi)\sum_{\left\vert \alpha\right\vert <2\theta}\frac{\left(
ix\right)  ^{\alpha}\left(  t\xi\right)  ^{\alpha}}{\alpha!}\right)
\frac{d\xi}{w\left(  \xi\right)  \left\vert \xi\right\vert ^{2\theta}}\\
& =\left(  2\pi\right)  ^{-\frac{d}{2}}\int\left(  e^{i\left(  x,\eta\right)
}-\rho(t^{-1}\eta)\sum_{\left\vert \alpha\right\vert <2\theta}\frac{\left(
ix\right)  ^{\alpha}\eta^{\alpha}}{\alpha!}\right)  \frac{t^{-d}d\eta
}{w\left(  t^{-1}\eta\right)  \left\vert t^{-1}\eta\right\vert ^{2\theta}}\\
& =t^{s+2\theta-d}\left(  2\pi\right)  ^{-\frac{d}{2}}\int\left(  e^{i\left(
x,\eta\right)  }-\rho(t^{-1}\eta)\sum_{\left\vert \alpha\right\vert <2\theta
}\frac{\left(  ix\right)  ^{\alpha}\eta^{\alpha}}{\alpha!}\right)  \frac{d\xi
}{w\left(  \eta\right)  \left\vert \eta\right\vert ^{2\theta}}\\
& =t^{s+2\theta-d}G_{\rho\left(  t^{-1}\cdot\right)  }\left(  x\right)  .
\end{align*}

Further, using part 1,%
\begin{align*}
G_{\rho\left(  t^{-1}\cdot\right)  }\left(  x\right)   & =G_{\rho}\left(
x\right)  +\left(  G_{\rho\left(  t^{-1}\cdot\right)  }\left(  x\right)
-G_{\rho}\left(  x\right)  \right) \\
& =G_{\rho}\left(  x\right)  +\left(  2\pi\right)  ^{-d/2}\sum_{\left\vert
\alpha\right\vert <2\theta}\left(  \int\frac{\rho(\xi)-\rho(t^{-1}\xi
)}{w\left(  \xi\right)  \left\vert \xi\right\vert ^{2\theta}}\left(
i\xi\right)  ^{\alpha}d\xi\right)  \frac{x^{\alpha}}{\alpha!},
\end{align*}

so that%
\[
G_{\rho}\left(  tx\right)  =t^{s+2\theta-d}G_{\rho}\left(  x\right)
+\frac{t^{s+2\theta-d}}{\left(  2\pi\right)  ^{d/2}}\sum_{\left\vert
\alpha\right\vert <2\theta}\left(  \int\frac{\rho(\xi)-\rho(t^{-1}\xi
)}{w\left(  \xi\right)  \left\vert \xi\right\vert ^{2\theta}}\left(
i\xi\right)  ^{\alpha}d\xi\right)  \frac{x^{\alpha}}{\alpha!}.
\]

But the change of variables $\eta=t^{-1}\xi$, $d\eta=t^{-d}d\xi$ yields%
\begin{align*}
\int\frac{\rho(\xi)-\rho(t^{-1}\xi)}{w\left(  \xi\right)  \left\vert
\xi\right\vert ^{2\theta}}\left(  i\xi\right)  ^{\alpha}d\xi & =\int\frac
{\rho(t\eta)-\rho(\eta)}{w\left(  t\eta\right)  \left\vert t\eta\right\vert
^{2\theta}}\left(  it\eta\right)  ^{\alpha}t^{d}d\eta\\
& =t^{d+\left\vert \alpha\right\vert -s-2\theta}\int\frac{\rho(t\eta
)-\rho(\eta)}{w\left(  \eta\right)  \left\vert \eta\right\vert ^{2\theta}%
}\left(  i\eta\right)  ^{\alpha}d\eta,
\end{align*}

so that%
\[
G_{\rho\left(  t^{-1}\xi\right)  }\left(  x\right)  =G_{\rho}\left(  x\right)
+\left(  2\pi\right)  ^{-d/2}\sum_{\left\vert \alpha\right\vert <2\theta
}\left(  t^{d+\left\vert \alpha\right\vert -s-2\theta}\int\frac{\rho
(t\eta)-\rho(\eta)}{w\left(  \eta\right)  \left\vert \eta\right\vert
^{2\theta}}\left(  i\eta\right)  ^{\alpha}d\eta\right)  \frac{x^{\alpha}%
}{\alpha!},
\]

and%
\begin{align*}
G_{\rho}\left(  tx\right)   & =t^{s+2\theta-d}G_{\rho\left(  t^{-1}\xi\right)
}\left(  x\right) \\
& =t^{s+2\theta-d}\left(  G_{\rho}\left(  x\right)  +\left(  2\pi\right)
^{-d/2}\sum_{\left\vert \alpha\right\vert <2\theta}\left(  t^{d+\left\vert
\alpha\right\vert -s-2\theta}\int\frac{\rho(t\eta)-\rho(\eta)}{w\left(
\eta\right)  \left\vert \eta\right\vert ^{2\theta}}\left(  i\eta\right)
^{\alpha}d\eta\right)  \frac{x^{\alpha}}{\alpha!}\right) \\
& =t^{s+2\theta-d}G_{\rho}\left(  x\right)  +\left(  2\pi\right)  ^{-d/2}%
\sum_{\left\vert \alpha\right\vert <2\theta}t^{\left\vert \alpha\right\vert
}\left(  \int\frac{\rho(t\eta)-\rho(\eta)}{w\left(  \eta\right)  \left\vert
\eta\right\vert ^{2\theta}}\left(  i\eta\right)  ^{\alpha}d\eta\right)
\frac{x^{\alpha}}{\alpha!}\\
& =t^{s+2\theta-d}G_{\rho}\left(  x\right)  +\left(  2\pi\right)  ^{-d/2}%
\sum_{\left\vert \alpha\right\vert <2\theta}\left(  \int\frac{\rho(t\eta
)-\rho(\eta)}{w\left(  \eta\right)  \left\vert \eta\right\vert ^{2\theta}%
}\left(  i\eta\right)  ^{\alpha}d\eta\right)  \frac{\left(  tx\right)
^{\alpha}}{\alpha!}.
\end{align*}

\end{proof}

\begin{corollary}
\label{Cor_Thm_Grho}Assume the weight function $w$ has properties W2.1 and
W3.2 for some order $\theta$ and parameter $\kappa$. Then:

\begin{enumerate}
\item There always exists a conjugate-even basis function.

\item If $w$ is even then there exists an even, real valued basis function.

\item If $w$ is radial then there exists a radial basis function.

\item If $w$ is homogeneous of order $s$ then there exists a basis function
$G$ which is homogeneous of order $s+2\theta-d$ modulo a polynomial of order
at most $2\theta$. More precisely
\[
G\left(  tx\right)  -t^{s+2\theta-d}G\left(  x\right)  =\sum_{\left\vert
\alpha\right\vert <2\theta}q_{\alpha}\left(  t\right)  x^{\alpha},\quad
x\in\mathbb{R}^{d},\text{ }t>0.
\]

\end{enumerate}
\end{corollary}

Regarding part 4 of the last theorem, we now give explicit formula for the
basis function $G_{\rho}$ when $w$ and $\rho$ are radial.

\begin{theorem}
If $w$ and $\rho$ are radial, say $w\left(  \xi\right)  =w_{\circ}\left(
\left\vert \xi\right\vert \right)  $ and $\rho\left(  \xi\right)  =\rho
_{\circ}\left(  \left\vert \xi\right\vert \right)  $, then the basis function
$G_{\rho}$ is radial, say $G_{\rho}\left(  x\right)  =\left(  G_{\rho}\right)
_{\circ}\left(  \left\vert x\right\vert \right)  $. In fact%
\begin{align*}
G_{\rho}\left(  x\right)   &  =\int\left(  e^{ix\xi}-\rho(\xi)\sum
_{k=0}^{2\theta}\frac{\left(  ix\xi\right)  ^{k}}{k!}\right)  \frac{d\xi
}{w\left(  \xi\right)  \left\vert \xi\right\vert ^{2\theta}}\\
&  =\Gamma\left(  \frac{d}{2}\right)  \int\left(  \frac{J_{\frac{d-2}{2}%
}\left(  \left\vert x\right\vert \left\vert \xi\right\vert \right)  }{\left(
\left\vert x\right\vert \left\vert \xi\right\vert /2\right)  ^{\frac{d-2}{2}}%
}-\rho\left(  \xi\right)  \sum_{j=0}^{\theta}\frac{\left(  -1\right)  ^{j}%
}{j!\Gamma\left(  j+\frac{d}{2}\right)  }\left(  \frac{\left\vert x\right\vert
\left\vert \xi\right\vert }{2}\right)  ^{2j}\right)  \frac{d\xi}{w\left(
\xi\right)  \left\vert \xi\right\vert ^{2\theta}},
\end{align*}

and%
\[
\left(  G_{\rho}\right)  _{\circ}\left(  s\right)  =\Gamma\left(  \frac{d}%
{2}\right)  \omega_{d}\int_{0}^{\infty}\left(  \frac{J_{\frac{d-2}{2}}\left(
st\right)  }{\left(  st/2\right)  ^{\frac{d-2}{2}}}-\rho_{\circ}\left(
t\right)  \sum_{j=0}^{\theta}\frac{\left(  -1\right)  ^{j}}{j!\Gamma\left(
j+\frac{d}{2}\right)  }\left(  \frac{st}{2}\right)  ^{2j}\right)
\frac{t^{d-1}dt}{w_{\circ}\left(  t\right)  t^{2\theta}}.
\]

\end{theorem}

\begin{proof}
From the definition in Theorem \ref{Thm_deriv_Grho},%
\begin{align*}
G_{\rho}\left(  x\right)   & =\left(  2\pi\right)  ^{-\frac{d}{2}}\int\left(
e^{ix\xi}-\rho(\xi)\sum_{k=0}^{2\theta}\frac{\left(  ix\xi\right)  ^{k}}%
{k!}\right)  \frac{d\xi}{w\left(  \xi\right)  \left\vert \xi\right\vert
^{2\theta}}\\
& =\left(  2\pi\right)  ^{-\frac{d}{2}}\lim_{\varepsilon\rightarrow0^{+}}%
\int\limits_{\left\vert \xi\right\vert \geq\varepsilon}\left(  e^{ix\xi}%
-\rho(\xi)\sum_{k=0}^{2\theta}\frac{\left(  ix\xi\right)  ^{k}}{k!}\right)
\frac{d\xi}{w\left(  \xi\right)  \left\vert \xi\right\vert ^{2\theta}}.
\end{align*}

From Theorem \ref{Cor_Thm_IntegBr_exp(ixy)f(|x|)dx},%
\begin{align*}
\int\limits_{\left\vert \xi\right\vert \geq\varepsilon} &  \left(  e^{ix\xi
}-\rho(\xi)\sum_{k=0}^{2\theta}\frac{\left(  ix\xi\right)  ^{k}}{k!}\right)
\frac{d\xi}{w\left(  \xi\right)  \left\vert \xi\right\vert ^{2\theta}}\\
&  =\int\limits_{\left\vert \xi\right\vert \geq\varepsilon}\frac{e^{ix\xi}%
}{w\left(  \xi\right)  \left\vert \xi\right\vert ^{2\theta}}-\int%
\limits_{\left\vert \xi\right\vert \geq\varepsilon}\rho\left(  \xi\right)
\sum_{k=0}^{2\theta}\frac{1}{k!}\frac{\left(  ix\xi\right)  ^{k}}{w\left(
\xi\right)  \left\vert \xi\right\vert ^{2\theta}}d\xi\\
&  =\int\limits_{\left\vert \xi\right\vert \geq\varepsilon}\frac{e^{ix\xi}%
}{w\left(  \xi\right)  \left\vert \xi\right\vert ^{2\theta}}-\sum
_{k=0}^{2\theta}\frac{1}{k!}\int\limits_{\left\vert \xi\right\vert
\geq\varepsilon}\frac{\rho\left(  \xi\right)  \left(  ix\xi\right)  ^{k}%
}{w\left(  \xi\right)  \left\vert \xi\right\vert ^{2\theta}}d\xi\\
&  =\int\limits_{\left\vert \xi\right\vert \geq\varepsilon}\frac{e^{ix\xi}%
}{w_{\circ}\left(  \xi\right)  \left\vert \xi\right\vert ^{2\theta}}%
-\sum_{k=0}^{2\theta}\int\limits_{\left\vert \xi\right\vert \geq\varepsilon
}\frac{\left(  ix\xi\right)  ^{k}}{k!}\frac{\rho_{\circ}\left(  \left\vert
\xi\right\vert \right)  }{w_{\circ}\left(  \left\vert \xi\right\vert \right)
\left\vert \xi\right\vert ^{2\theta}}d\xi\\
&  =\int\limits_{\left\vert \xi\right\vert \geq\varepsilon}\frac{e^{ix\xi}%
}{w_{\circ}\left(  \xi\right)  \left\vert \xi\right\vert ^{2\theta}}%
-\sum_{k=0}^{2\theta}\frac{1}{k!}\int\limits_{\left\vert \xi\right\vert
\geq\varepsilon}\left(  ix\widehat{\xi}\right)  ^{k}\frac{\rho_{\circ}\left(
\left\vert \xi\right\vert \right)  \left\vert \xi\right\vert ^{k}}{w_{\circ
}\left(  \left\vert \xi\right\vert \right)  \left\vert \xi\right\vert
^{2\theta}}d\xi.
\end{align*}

But from Corollary \ref{Cor_2_Thm_Integ_u(xy)f(|x|)dx} in the Appendix:%
\begin{align*}
\int_{\left\vert \xi\right\vert \geq\varepsilon}\left(  x\widehat{\xi}\right)
^{j}f\left(  \left\vert \xi\right\vert \right)  d\xi & =0,\quad j\text{
}is\text{ }odd,\\
\int_{\left\vert \xi\right\vert \geq\varepsilon}\left(  x\widehat{\xi}\right)
^{2j}f\left(  \left\vert \xi\right\vert \right)  d\xi & =\left\vert
x\right\vert ^{2j}B\left(  \frac{d-1}{2},n+\frac{1}{2}\right)  \omega
_{d-1}\int_{\varepsilon}^{\infty}f\left(  t\right)  t^{d-1}dt.
\end{align*}

Hence%
\begin{align*}
\int\limits_{\left\vert \xi\right\vert \geq\varepsilon} &  \left(  e^{ix\xi
}-\rho(\xi)\sum_{k=0}^{2\theta}\frac{\left(  ix\xi\right)  ^{k}}{k!}\right)
\frac{d\xi}{w\left(  \xi\right)  \left\vert \xi\right\vert ^{2\theta}}\\
&  =\int\limits_{\left\vert \xi\right\vert \geq\varepsilon}\frac{e^{ix\xi}%
}{w_{\circ}\left(  \xi\right)  \left\vert \xi\right\vert ^{2\theta}}%
-\sum_{j=0}^{\theta}\frac{1}{\left(  2j\right)  !}\int\limits_{\left\vert
\xi\right\vert \geq\varepsilon}\left(  ix\widehat{\xi}\right)  ^{2j}\frac
{\rho_{\circ}\left(  \left\vert \xi\right\vert \right)  \left\vert
\xi\right\vert ^{2j}}{w_{\circ}\left(  \left\vert \xi\right\vert \right)
\left\vert \xi\right\vert ^{2\theta}}d\xi\\
&  =\int\limits_{\left\vert \xi\right\vert \geq\varepsilon}\frac{e^{ix\xi}%
d\xi}{w_{\circ}\left(  \xi\right)  \left\vert \xi\right\vert ^{2\theta}}%
-\sum_{j=0}^{\theta}\frac{\left\vert x\right\vert ^{2j}}{\left(  2j\right)
!}B\left(  \frac{d-1}{2},j+\frac{1}{2}\right)  \omega_{d-1}\int_{\varepsilon
}^{\infty}\frac{\rho_{\circ}\left(  t\right)  t^{2j}}{w_{\circ}\left(
t\right)  t^{2\theta}}t^{d-1}dt.
\end{align*}

From Theorem \ref{Cor_Thm_IntegBr_exp(ixy)f(|x|)dx},%
\[
\int\limits_{\left\vert \xi\right\vert \leq r}e^{ix\xi}f\left(  \left\vert
\xi\right\vert \right)  d\xi=\Gamma\left(  \tfrac{d}{2}\right)  \omega_{d}%
\int_{0}^{r}\frac{J_{\frac{d-2}{2}}\left(  \left\vert x\right\vert t\right)
}{\left(  \left\vert x\right\vert t/2\right)  ^{\frac{d-2}{2}}}f\left(
t\right)  t^{d-1}dt,
\]

so%
\[
\int_{\left\vert \xi\right\vert \geq\varepsilon}\frac{e^{ix\xi}d\xi}{w\left(
\xi\right)  \left\vert \xi\right\vert ^{2\theta}}=\Gamma\left(  \tfrac{d}%
{2}\right)  \omega_{d}\int_{\varepsilon}^{\infty}\frac{J_{\frac{d-2}{2}%
}\left(  \left\vert x\right\vert t\right)  }{\left(  \left\vert x\right\vert
t/2\right)  ^{\frac{d-2}{2}}}\frac{t^{d-1}}{w_{\circ}\left(  t\right)
t^{2\theta}}dt,
\]

and we get%
\begin{align*}
& \int\limits_{\left\vert \xi\right\vert \geq\varepsilon}\left(  e^{ix\xi
}-\rho(\xi)\sum_{k=0}^{2\theta}\frac{\left(  ix\xi\right)  ^{k}}{k!}\right)
\frac{d\xi}{w\left(  \xi\right)  \left\vert \xi\right\vert ^{2\theta}}\\
& =\Gamma\left(  \tfrac{d}{2}\right)  \omega_{d}\int_{\varepsilon}^{\infty
}\frac{J_{\frac{d-2}{2}}\left(  \left\vert x\right\vert t\right)  }{\left(
\left\vert x\right\vert t/2\right)  ^{\frac{d-2}{2}}}\frac{t^{d-1}}{w_{\circ
}\left(  t\right)  t^{2\theta}}dt-\sum_{j=0}^{\theta}\left(  -1\right)
^{j}\frac{\left\vert x\right\vert ^{2j}}{\left(  2j\right)  !}B\left(
\frac{d-1}{2},j+\frac{1}{2}\right)  \omega_{d-1}\int_{\varepsilon}^{\infty
}\frac{\rho_{\circ}\left(  t\right)  t^{2j}}{w_{\circ}\left(  t\right)
t^{2\theta}}t^{d-1}dt\\
& =\int_{\varepsilon}^{\infty}\Gamma\left(  \tfrac{d}{2}\right)  \omega
_{d}\frac{J_{\frac{d-2}{2}}\left(  \left\vert x\right\vert t\right)  }{\left(
\left\vert x\right\vert t/2\right)  ^{\frac{d-2}{2}}}\frac{t^{d-1}}{w_{\circ
}\left(  t\right)  t^{2\theta}}dt-\int_{\varepsilon}^{\infty}\rho_{\circ
}\left(  t\right)  \sum_{j=0}^{\theta}\left(  -1\right)  ^{j}\frac{\left(
\left\vert x\right\vert t\right)  ^{2j}}{\left(  2j\right)  !}B\left(
\frac{d-1}{2},j+\frac{1}{2}\right)  \omega_{d-1}\frac{t^{d-1}dt}{w_{\circ
}\left(  t\right)  t^{2\theta}}\\
& =\int_{\varepsilon}^{\infty}\left(  \Gamma\left(  \tfrac{d}{2}\right)
\omega_{d}\frac{J_{\frac{d-2}{2}}\left(  \left\vert x\right\vert t\right)
}{\left(  \left\vert x\right\vert t/2\right)  ^{\frac{d-2}{2}}}-\rho_{\circ
}\left(  t\right)  \sum_{j=0}^{\theta}\left(  -1\right)  ^{j}\frac{\left(
\left\vert x\right\vert t\right)  ^{2j}}{\left(  2j\right)  !}B\left(
\frac{d-1}{2},j+\frac{1}{2}\right)  \omega_{d-1}\right)  \frac{t^{d-1}%
dt}{w_{\circ}\left(  t\right)  t^{2\theta}}\\
& =\Gamma\left(  \tfrac{d}{2}\right)  \omega_{d}\int_{\varepsilon}^{\infty
}\left(  \frac{J_{\frac{d-2}{2}}\left(  \left\vert x\right\vert t\right)
}{\left(  \left\vert x\right\vert t/2\right)  ^{\frac{d-2}{2}}}-\rho_{\circ
}\left(  t\right)  \sum_{j=0}^{\theta}\left(  -1\right)  ^{j}\frac{\left(
\left\vert x\right\vert t\right)  ^{2j}}{\left(  2j\right)  !}\frac{B\left(
\frac{d-1}{2},j+\frac{1}{2}\right)  }{\Gamma\left(  \tfrac{d}{2}\right)
}\frac{\omega_{d-1}}{\omega_{d}}\right)  \frac{t^{d-1}dt}{w_{\circ}\left(
t\right)  t^{2\theta}}\\
& =\Gamma\left(  \tfrac{d}{2}\right)  \omega_{d}\int_{\varepsilon}^{\infty
}\left(  \frac{J_{\frac{d-2}{2}}\left(  \left\vert x\right\vert t\right)
}{\left(  \left\vert x\right\vert t/2\right)  ^{\frac{d-2}{2}}}-\rho_{\circ
}\left(  t\right)  \sum_{j=0}^{\theta}\left(  -1\right)  ^{j}\frac{\left(
\left\vert x\right\vert t\right)  ^{2j}}{\left(  2j\right)  !}\frac{B\left(
\frac{d-1}{2},j+\frac{1}{2}\right)  }{\Gamma\left(  \tfrac{d}{2}\right)
}\frac{\frac{2\pi^{\left(  d-1\right)  /2}}{\Gamma\left(  \left(  d-1\right)
/2\right)  }}{\frac{2\pi^{d/2}}{\Gamma\left(  d/2\right)  }}\right)
\frac{t^{d-1}dt}{w_{\circ}\left(  t\right)  t^{2\theta}}\\
& =\Gamma\left(  \tfrac{d}{2}\right)  \omega_{d}\int_{\varepsilon}^{\infty
}\left(  \frac{J_{\frac{d-2}{2}}\left(  \left\vert x\right\vert t\right)
}{\left(  \left\vert x\right\vert t/2\right)  ^{\frac{d-2}{2}}}-\rho_{\circ
}\left(  t\right)  \sum_{j=0}^{\theta}\left(  -1\right)  ^{j}\frac{\left(
\left\vert x\right\vert t\right)  ^{2j}}{\left(  2j\right)  !}\frac{B\left(
\frac{d-1}{2},j+\frac{1}{2}\right)  }{\pi^{1/2}\Gamma\left(  \frac{d-1}%
{2}\right)  }\right)  \frac{t^{d-1}dt}{w_{\circ}\left(  t\right)  t^{2\theta}%
}.
\end{align*}

But%
\begin{align*}
\frac{B\left(  \frac{d-1}{2},j+\frac{1}{2}\right)  }{\pi^{1/2}\Gamma\left(
\frac{d-1}{2}\right)  }  & =\frac{\Gamma\left(  \frac{d-1}{2}\right)
\Gamma\left(  j+\frac{1}{2}\right)  }{\Gamma\left(  j+\frac{d}{2}\right)
\pi^{1/2}\Gamma\left(  \frac{d-1}{2}\right)  }=\frac{1}{\Gamma\left(
j+\frac{d}{2}\right)  }\frac{\Gamma\left(  j+\frac{1}{2}\right)  }%
{\Gamma\left(  \frac{1}{2}\right)  }=\frac{1}{\Gamma\left(  j+\frac{d}%
{2}\right)  }\frac{\Gamma\left(  2j\right)  }{2^{2j-1}\Gamma\left(  j\right)
}=\\
& =\frac{1}{\Gamma\left(  j+\frac{d}{2}\right)  }\frac{\left(  2j-1\right)
!}{2^{2j-1}\left(  j-1\right)  !}=\frac{1}{\Gamma\left(  j+\frac{d}{2}\right)
}\frac{\left(  2j\right)  !}{2^{2j}j!},
\end{align*}

so that%
\begin{align*}
\int\limits_{\left\vert \xi\right\vert \geq\varepsilon} &  \left(  e^{ix\xi
}-\rho(\xi)\sum_{k=0}^{2\theta}\frac{\left(  ix\xi\right)  ^{k}}{k!}\right)
\frac{d\xi}{w\left(  \xi\right)  \left\vert \xi\right\vert ^{2\theta}}\\
&  =\Gamma\left(  \tfrac{d}{2}\right)  \omega_{d}\int_{\varepsilon}^{\infty
}\left(  \frac{J_{\frac{d-2}{2}}\left(  \left\vert x\right\vert t\right)
}{\left(  \left\vert x\right\vert t/2\right)  ^{\frac{d-2}{2}}}-\rho_{\circ
}\left(  t\right)  \sum_{j=0}^{\theta}\left(  -1\right)  ^{j}\frac{\left(
\left\vert x\right\vert t\right)  ^{2j}}{\left(  2j\right)  !}\frac{1}%
{\Gamma\left(  j+\frac{d}{2}\right)  }\frac{\left(  2j\right)  !}{2^{2j}%
j!}\right)  \frac{t^{d-1}dt}{w_{\circ}\left(  t\right)  t^{2\theta}}\\
&  =\Gamma\left(  \tfrac{d}{2}\right)  \omega_{d}\int_{\varepsilon}^{\infty
}\left(  \frac{J_{\frac{d-2}{2}}\left(  \left\vert x\right\vert t\right)
}{\left(  \left\vert x\right\vert t/2\right)  ^{\frac{d-2}{2}}}-\rho_{\circ
}\left(  t\right)  \sum_{j=0}^{\theta}\frac{\left(  -1\right)  ^{j}}%
{j!\Gamma\left(  j+\frac{d}{2}\right)  }\left(  \frac{\left\vert x\right\vert
t}{2}\right)  ^{2j}\right)  \frac{t^{d-1}dt}{w_{\circ}\left(  t\right)
t^{2\theta}},
\end{align*}

and consequently%
\begin{align*}
G_{\rho}\left(  x\right)   &  =\int\left(  e^{ix\xi}-\rho(\xi)\sum
_{k=0}^{2\theta}\frac{\left(  ix\xi\right)  ^{k}}{k!}\right)  \frac{d\xi
}{w\left(  \xi\right)  \left\vert \xi\right\vert ^{2\theta}}\\
&  =\Gamma\left(  \frac{d}{2}\right)  \omega_{d}\int_{0}^{\infty}\left(
\frac{J_{\frac{d-2}{2}}\left(  \left\vert x\right\vert t\right)  }{\left(
\left\vert x\right\vert t/2\right)  ^{\frac{d-2}{2}}}-\rho_{\circ}\left(
t\right)  \sum_{j=0}^{\theta}\frac{\left(  -1\right)  ^{j}}{j!\Gamma\left(
j+\frac{d}{2}\right)  }\left(  \frac{\left\vert x\right\vert t}{2}\right)
^{2j}\right)  \frac{t^{d-1}dt}{w_{\circ}\left(  t\right)  t^{2\theta}}\\
&  =\Gamma\left(  \frac{d}{2}\right)  \int\left(  \frac{J_{\frac{d-2}{2}%
}\left(  \left\vert x\right\vert \left\vert \xi\right\vert \right)  }{\left(
\left\vert x\right\vert \left\vert \xi\right\vert /2\right)  ^{\frac{d-2}{2}}%
}-\rho\left(  \xi\right)  \sum_{j=0}^{\theta}\frac{\left(  -1\right)  ^{j}%
}{j!\Gamma\left(  j+\frac{d}{2}\right)  }\left(  \frac{\left\vert x\right\vert
\left\vert \xi\right\vert }{2}\right)  ^{2j}\right)  \frac{d\xi}{w\left(
\xi\right)  \left\vert \xi\right\vert ^{2\theta}}.
\end{align*}

Thus%
\[
G_{\circ}\left(  s\right)  =\Gamma\left(  \frac{d}{2}\right)  \omega_{d}%
\int_{0}^{\infty}\left(  \frac{J_{\frac{d-2}{2}}\left(  st\right)  }{\left(
st/2\right)  ^{\frac{d-2}{2}}}-\rho_{\circ}\left(  t\right)  \sum
_{j=0}^{\theta}\frac{\left(  -1\right)  ^{j}}{j!\Gamma\left(  j+\frac{d}%
{2}\right)  }\left(  \frac{st}{2}\right)  ^{2j}\right)  \frac{t^{d-1}%
dt}{w_{\circ}\left(  t\right)  t^{2\theta}}.
\]

From equation 11.6 of Arfken \cite{Arfken70}:%
\begin{align*}
J_{v}\left(  s\right)   & =\sum_{j=0}^{\infty}\frac{\left(  -1\right)  ^{j}%
}{j!\left(  j+v\right)  !}\left(  \frac{s}{2}\right)  ^{v+2j}=\sum
_{k=0}^{\infty}\frac{\left(  -1\right)  ^{j}}{j!\Gamma\left(  j+v+1\right)
}\left(  \frac{s}{2}\right)  ^{v+2j}.\\
\frac{J_{v}\left(  s\right)  }{\left(  s/2\right)  ^{v}}  & =\sum
_{j=0}^{\infty}\frac{\left(  -1\right)  ^{j}}{j!\Gamma\left(  j+v+1\right)
}\left(  \frac{s}{2}\right)  ^{2j}.\\
\frac{J_{\frac{d-2}{2}}\left(  s\right)  }{\left(  s/2\right)  ^{\frac{d-2}%
{2}}}  & =\sum_{j=0}^{\infty}\frac{\left(  -1\right)  ^{j}}{j!\Gamma\left(
j+\frac{d}{2}\right)  }\left(  \frac{s}{2}\right)  ^{2j}.
\end{align*}

\end{proof}

?? \textbf{Where to put this theorem}? ?? See the general estimate Lemma
\ref{Lem_bound_on_integ(|DQ(exp)|/w|x|)} above.

\begin{theorem}
\label{Thm_bnd_on_integ(|DQ(exp)|/w|x|)_radial}Suppose $0\leq\rho\leq1$ and
$\rho$ and $w$ are radial. Then, when $\left\vert \gamma\right\vert <2\theta$,%
\begin{align*}
&  \int\dfrac{\left\vert D_{x}^{\gamma}\mathcal{Q}_{\emptyset,2\theta,\xi
}\left(  e^{i\left(  x,\xi\right)  }\right)  \right\vert }{w\left(
\xi\right)  \left\vert \xi\right\vert ^{2\theta}}d\xi\\
&  \leq\tfrac{B\left(  \frac{1}{2}\left(  \gamma+\mathbf{1}\right)  \right)
}{B\left(  \frac{1}{2}\mathbf{1}\right)  }\left(  \frac{\left\Vert D^{2\theta
}\rho_{\circ}\right\Vert _{\infty;\leq r_{3}}}{\left(  2\theta\right)  !}%
\int_{\left\vert \cdot\right\vert \leq r_{3}}\frac{1}{w}+\int_{\left\vert
\cdot\right\vert \geq r_{3}}\dfrac{1}{w\left\vert \cdot\right\vert ^{2\theta}%
}\right)  +\\
&  +\frac{\left\vert x\right\vert ^{2\theta-\left\vert \gamma\right\vert }%
}{\left(  2\theta-\left\vert \gamma\right\vert +??1\right)  !}\min\left\{
\tfrac{B\left(  \frac{1}{2}\left(  \gamma+\mathbf{1}\right)  \right)
}{B\left(  \frac{1}{2}\mathbf{1}\right)  },\tfrac{B\left(  \frac{d-1}{2}%
,\frac{2\theta-\left\vert \gamma\right\vert +1}{2}\right)  }{B\left(
\frac{d-1}{2},\frac{1}{2}\right)  }\right\}  \left(  \int\limits_{\left\vert
\cdot\right\vert \leq r_{3}}\dfrac{1}{w}+\left\Vert \rho_{\circ}\left(
t\right)  t^{2\theta}\right\Vert _{\infty;\geq r_{3}}\int\limits_{\left\vert
\cdot\right\vert \geq r_{3}}\dfrac{1}{w\left\vert \cdot\right\vert ^{2\theta}%
}\right)  .
\end{align*}

\end{theorem}

\begin{proof}
From part 2 of the proof of Theorem \ref{Thm_bound_on_g(e,x)_2},%
\[
D_{x}^{\gamma}\mathcal{Q}_{\emptyset,2\theta,\xi}\left(  e^{ix\xi}\right)
=\left(  i\xi\right)  ^{\gamma}\left(  e^{ix\xi}\left(  1-\rho(\xi)\right)
+\rho(\xi)\left(  ix\xi\right)  ^{2\theta-\left\vert \gamma\right\vert }%
\mu_{2\theta-\left\vert \gamma\right\vert }(ix\xi)\right)  ,
\]

so that%
\begin{align*}
\int\dfrac{\left\vert D_{x}^{\gamma}\mathcal{Q}_{\emptyset,2\theta,\xi}\left(
e^{i\left(  x,\xi\right)  }\right)  \right\vert }{w\left(  \xi\right)
\left\vert \xi\right\vert ^{2\theta}}d\xi &  \leq\int\dfrac{\left\vert
\xi^{\gamma}\right\vert \left(  1-\rho(\xi)\right)  }{w\left(  \xi\right)
\left\vert \xi\right\vert ^{2\theta}}d\xi+\\
&  \qquad+\frac{\sqrt{2\pi}}{\left(  2\theta-\left\vert \gamma\right\vert
\right)  !}\int\dfrac{\left\vert \xi^{\gamma}\right\vert \left\vert \left(
x\xi\right)  ^{2\theta-\left\vert \gamma\right\vert }\widehat{g_{2\theta
-\left\vert \gamma\right\vert }}(ix\xi)\right\vert \rho(\xi)}{w\left(
\xi\right)  \left\vert \xi\right\vert ^{2\theta}}d\xi.
\end{align*}

From Theorem \ref{Thm_integ_X^a_f(absX)dX_R^d} and then Theorem
\ref{Thm_Tay_rem_zeros},%
\begin{align*}
\int\dfrac{\left\vert \xi^{\gamma}\right\vert \left(  1-\rho(\xi)\right)
}{w\left(  \xi\right)  \left\vert \xi\right\vert ^{2\theta}}d\xi &
=\tfrac{B\left(  \frac{1}{2}\left(  \gamma+\mathbf{1}\right)  \right)
}{B\left(  \frac{1}{2}\mathbf{1}\right)  }\int\dfrac{1-\rho}{w\left\vert
\cdot\right\vert ^{2\theta}}\\
& =\tfrac{B\left(  \frac{1}{2}\left(  \gamma+\mathbf{1}\right)  \right)
}{B\left(  \frac{1}{2}\mathbf{1}\right)  }\left(  \int_{\left\vert
\cdot\right\vert \leq r_{3}}\frac{1}{w}\dfrac{1-\rho}{\left\vert
\cdot\right\vert ^{2\theta}}+\int_{\left\vert \cdot\right\vert \geq r_{3}%
}\dfrac{1-\rho}{w\left\vert \cdot\right\vert ^{2\theta}}\right) \\
& \leq\tfrac{B\left(  \frac{1}{2}\left(  \gamma+\mathbf{1}\right)  \right)
}{B\left(  \frac{1}{2}\mathbf{1}\right)  }\left(  \int_{\left\vert
\cdot\right\vert \leq r_{3}}\frac{1}{w}\dfrac{1-\rho}{\left\vert
\cdot\right\vert ^{2\theta}}+\int_{\left\vert \cdot\right\vert \geq r_{3}%
}\dfrac{1}{w\left\vert \cdot\right\vert ^{2\theta}}\right) \\
& \leq\tfrac{B\left(  \frac{1}{2}\left(  \gamma+\mathbf{1}\right)  \right)
}{B\left(  \frac{1}{2}\mathbf{1}\right)  }\left(  \frac{\left\Vert D^{2\theta
}\rho_{\circ}\right\Vert _{\infty;\leq r_{3}}}{\left(  2\theta\right)  !}%
\int_{\left\vert \cdot\right\vert \leq r_{3}}\frac{1}{w}+\int_{\left\vert
\cdot\right\vert \geq r_{3}}\dfrac{1}{w\left\vert \cdot\right\vert ^{2\theta}%
}\right)
\end{align*}

Using \ref{Lem_gm_properties_2} and then Corollary
\ref{Cor_integ_inprod_to_1dim},%
\begin{align*}
\int &  \dfrac{\left\vert \xi^{\gamma}\right\vert \left\vert \left(
x\xi\right)  ^{2\theta-\left\vert \gamma\right\vert }\widehat{g_{2\theta
-\left\vert \gamma\right\vert }}(ix\xi)\right\vert \rho(\xi)}{w\left(
\xi\right)  \left\vert \xi\right\vert ^{2\theta}}d\xi\\
&  \leq\tfrac{1}{\sqrt{2\pi}}\tfrac{1}{2\theta-\left\vert \gamma\right\vert
+1}\int\dfrac{\left\vert \xi^{\gamma}\right\vert \left\vert x\xi\right\vert
^{2\theta-\left\vert \gamma\right\vert }\rho(\xi)}{w\left(  \xi\right)
\left\vert \xi\right\vert ^{2\theta}}d\xi=\\
&  =\tfrac{1}{\sqrt{2\pi}}\tfrac{\left\vert x\right\vert ^{2\theta-\left\vert
\gamma\right\vert }}{2\theta-\left\vert \gamma\right\vert +1}\int%
\dfrac{\left\vert \widehat{\xi}^{\gamma}\right\vert \left\vert \widehat
{x}\widehat{\xi}\right\vert ^{2\theta-\left\vert \gamma\right\vert }\rho(\xi
)}{w\left(  \xi\right)  }d\xi\\
&  \leq\tfrac{1}{\sqrt{2\pi}}\tfrac{\left\vert x\right\vert ^{2\theta
-\left\vert \gamma\right\vert }}{2\theta-\left\vert \gamma\right\vert +1}%
\min\left\{  \tfrac{B\left(  \frac{1}{2}\left(  \gamma+\mathbf{1}\right)
\right)  }{B\left(  \frac{1}{2}\mathbf{1}\right)  },\tfrac{B\left(  \frac
{d-1}{2},\frac{2\theta-\left\vert \gamma\right\vert +1}{2}\right)  }{B\left(
\frac{d-1}{2},\frac{1}{2}\right)  }\right\}  \int\dfrac{\rho}{w}\\
&  =\tfrac{1}{\sqrt{2\pi}}\tfrac{\left\vert x\right\vert ^{2\theta-\left\vert
\gamma\right\vert }}{2\theta-\left\vert \gamma\right\vert +1}\min\left\{
\ldots\right\}  \left(  \int_{\left\vert \cdot\right\vert \leq r_{3}}%
\dfrac{\rho}{w}+\int_{\left\vert \cdot\right\vert \geq r_{3}}\dfrac{\rho}%
{w}\right) \\
&  =\tfrac{1}{\sqrt{2\pi}}\tfrac{\left\vert x\right\vert ^{2\theta-\left\vert
\gamma\right\vert }}{2\theta-\left\vert \gamma\right\vert +1}\min\left\{
\ldots\right\}  \left(  \int_{\left\vert \cdot\right\vert \leq r_{3}}%
\dfrac{\rho}{w}+\int_{\left\vert \cdot\right\vert \geq r_{3}}\dfrac
{\rho\left\vert \cdot\right\vert ^{2\theta}}{w\left\vert \cdot\right\vert
^{2\theta}}\right) \\
&  \leq\tfrac{1}{\sqrt{2\pi}}\tfrac{\left\vert x\right\vert ^{2\theta
-\left\vert \gamma\right\vert }}{2\theta-\left\vert \gamma\right\vert +1}%
\min\left\{  \ldots\right\}  \left(  \int_{\left\vert \cdot\right\vert \leq
r_{3}}\dfrac{1}{w}+\left\Vert \rho\left\vert \cdot\right\vert ^{2\theta
}\right\Vert _{\infty;\geq r_{3}}\int_{\left\vert \cdot\right\vert \geq r_{3}%
}\dfrac{1}{w\left\vert \cdot\right\vert ^{2\theta}}\right) \\
&  =\tfrac{1}{\sqrt{2\pi}}\tfrac{\left\vert x\right\vert ^{2\theta-\left\vert
\gamma\right\vert }}{2\theta-\left\vert \gamma\right\vert +1}\min\left\{
\ldots\right\}  \left(  \int_{\left\vert \cdot\right\vert \leq r_{3}}\dfrac
{1}{w}+\left\Vert \rho_{\circ}\left(  t\right)  t^{2\theta}\right\Vert
_{\infty;\geq r_{3}}\int_{\left\vert \cdot\right\vert \geq r_{3}}\dfrac
{1}{w\left\vert \cdot\right\vert ^{2\theta}}\right)  ,
\end{align*}

so that%
\begin{align*}
&  \int\dfrac{\left\vert D_{x}^{\gamma}\mathcal{Q}_{\emptyset,2\theta,\xi
}\left(  e^{i\left(  x,\xi\right)  }\right)  \right\vert }{w\left(
\xi\right)  \left\vert \xi\right\vert ^{2\theta}}d\xi\\
&  \leq\tfrac{B\left(  \frac{1}{2}\left(  \gamma+\mathbf{1}\right)  \right)
}{B\left(  \frac{1}{2}\mathbf{1}\right)  }\left(  \frac{\left\Vert D^{2\theta
}\rho_{\circ}\right\Vert _{\infty;\leq r_{3}}}{\left(  2\theta\right)  !}%
\int_{\left\vert \cdot\right\vert \leq r_{3}}\frac{1}{w}+\int_{\left\vert
\cdot\right\vert \geq r_{3}}\dfrac{1}{w\left\vert \cdot\right\vert ^{2\theta}%
}\right)  +\\
&  +\frac{\left\vert x\right\vert ^{2\theta-\left\vert \gamma\right\vert }%
}{\left(  2\theta-\left\vert \gamma\right\vert +??1\right)  !}\min\left\{
\tfrac{B\left(  \frac{1}{2}\left(  \gamma+\mathbf{1}\right)  \right)
}{B\left(  \frac{1}{2}\mathbf{1}\right)  },\tfrac{B\left(  \frac{d-1}{2}%
,\frac{2\theta-\left\vert \gamma\right\vert +1}{2}\right)  }{B\left(
\frac{d-1}{2},\frac{1}{2}\right)  }\right\}  \left(  \int\limits_{\left\vert
\cdot\right\vert \leq r_{3}}\dfrac{1}{w}+\left\Vert \rho\left\vert
\cdot\right\vert ^{2\theta}\right\Vert _{\infty;\geq r_{3}}\int%
\limits_{\left\vert \cdot\right\vert \geq r_{3}}\dfrac{1}{w\left\vert
\cdot\right\vert ^{2\theta}}\right)  .
\end{align*}

\end{proof}

\section{Taylor expansion of basis functions when $w\in W3.2$%
\label{Sect_Tay_expan_basis_W3.2}}

?? \textbf{ADD BLURB}! ??

We begin by using the tempered distribution Taylor series expansion
\ref{a1.55}:%
\[
f\left(  \cdot+a\right)  -\sum_{\left\vert \beta\right\vert \leq n}%
\frac{a^{\beta}}{\beta!}D^{\beta}f=\left(  \mathcal{R}_{n+1}f\right)  \left(
\cdot,a\right)  ,\quad f\in S^{\prime},
\]

where%
\[
\left(  \mathcal{R}_{n+1}f\right)  \left(  \cdot,a\right)  =\frac{\sqrt{2\pi}%
}{n!}\left(  \left(  ia\xi\right)  ^{n+1}\overline{\widehat{g_{n}}}\left(
a\xi\right)  \widehat{f}\right)  ^{\vee}.
\]

Suppose $w\in W3.2$ is a weight function for order $\theta$ and smoothness
parameter $\kappa$. We want to\ estimate the remainder of this Taylor series
expansion. This makes sense since from part 4 of Theorem
\ref{Thm_basis_smth_W3.2_r3_pos} we have $G\subset C_{BP}^{\left(
\left\lfloor 2\kappa\right\rfloor \right)  }$ and accordingly $n$ will satisfy%
\[
n\leq\left\lfloor 2\kappa\right\rfloor .
\]

In this section I develop an analogue of the Taylor series expansions obtained
\textbf{below} for data functions in Section \ref{Sect_data_fn_Taylor_W3.2}.

From \ref{a1.55}, if $G\in S^{\prime}$ is a \textbf{basis distribution} then,%
\[
G\left(  \cdot+a\right)  -\sum_{k\leq n}\frac{\left(  aD\right)  ^{k}}%
{k!}G=\left(  \mathcal{R}_{n+1}G\right)  \left(  \cdot,a\right)  ,\text{\quad
}n\geq0,
\]

where%
\begin{equation}
\left(  \mathcal{R}_{n+1}G\right)  \left(  \cdot,a\right)  =\frac{\sqrt{2\pi}%
}{n!}\left(  \left(  ia\xi\right)  ^{n+1}\overline{\widehat{g_{n}}}\left(
a\xi\right)  \widehat{G}\right)  ^{\vee}.\label{a78}%
\end{equation}

Then from part 4 of Theorem \ref{Thm_basis_smth_W3.2_r3_pos}, $G\in
C_{BP}^{\left(  \left\lfloor 2\kappa\right\rfloor \right)  }$.

Use Theorem \ref{Thm_deriv_Grho}? which obtains formula \ref{a44}, namely%
\[
G=G_{\rho}+\left(  2\pi\right)  ^{-\frac{d}{2}}p_{2\theta-1,\widehat{G}},
\]

where $p_{2\theta-1,\widehat{G}}\in P_{2\theta-1}$ is defined using \ref{a45}
as%
\[
p_{2\theta-1,u}\left(  x\right)  =\sum_{\left\vert \alpha\right\vert <2\theta
}\frac{b_{2\theta-1,u,\alpha}}{\alpha!}x^{\alpha},\quad b_{2\theta-1,u,\alpha
}=\left[  u,\left(  -i\xi\right)  ^{\alpha}\rho\right]  ,\text{ }\rho\in
S_{1,2\theta},\text{ }u\in S^{\prime}.
\]

$G_{\rho}$ is given by \ref{a43} as
\[
G_{\rho}\left(  x\right)  =\left(  2\pi\right)  ^{-\frac{d}{2}}\int%
\frac{\mathcal{Q}_{\emptyset,2\theta,\xi}\left(  e^{i\left(  x,\xi\right)
}\right)  }{w\left(  \xi\right)  \left\vert \xi\right\vert ^{2\theta}}d\xi.
\]

Thus%
\begin{align}
&  \left(  \mathcal{R}_{n+1}G\right)  \left(  \cdot,a\right) \nonumber\\
&  =\frac{\sqrt{2\pi}}{n!}\left(  \left(  ia\xi\right)  ^{n+1}\overline
{\widehat{g_{n}}}\left(  a\xi\right)  \widehat{G_{\rho}}\right)  ^{\vee
}+\left(  \mathcal{R}_{n+1}\left(  \left(  2\pi\right)  ^{-\frac{d}{2}%
}p_{2\theta-1,\widehat{G}}\right)  \right)  \left(  \cdot,a\right) \nonumber\\
&  =\frac{\sqrt{2\pi}}{n!}\left(  \left(  ia\xi\right)  ^{n+1}\overline
{\widehat{g_{n}}}\left(  a\xi\right)  \widehat{G_{\rho}}\right)  ^{\vee
}+\left(  2\pi\right)  ^{-\frac{d}{2}}\left(  p_{2\theta-1,\widehat{G}}\left(
\cdot+a\right)  -\sum_{k\leq n}\frac{\left(  aD\right)  ^{k}}{k!}%
p_{2\theta-1,\widehat{G}}\right)  .\label{a62}%
\end{align}
\medskip

\underline{\textbf{Estimation of the }$G_{\rho}$\textbf{\ term of }%
\ref{a62}\textbf{\ when }$n\leq\left\lfloor 2\kappa\right\rfloor $}\medskip

Compare Section \ref{Sect_data_fn_Taylor_W3.2} which studies Taylor series
expansions of functions in $X_{0}^{\theta}$ where $w\in W3.2$.

We first introduce the partition of unity%
\[
1=\phi_{0}+\phi_{\infty},\text{\quad}\phi_{0}\in S_{1,2\theta},\text{\quad
}0\leq\phi_{0}\leq1,
\]

which implies that%
\[
\phi_{\infty}\in C_{\emptyset,2\theta}^{\infty}\cap C_{B}^{\infty}%
,\text{\quad}0\leq\phi_{\infty}\leq1.
\]

Observe that from Theorem \ref{Thm_Tay_rem_zeros}, for any $r>0$,%
\begin{equation}
\phi_{\infty}\left(  x\right)  \leq\frac{\left\vert x\right\vert ^{2\theta}%
}{\left(  2\theta\right)  !}\left\Vert \left(  \widehat{\cdot}D\right)
^{n}\phi_{\infty}\right\Vert _{\infty;\leq r},\text{\quad}x\in\overline{B}%
_{r}.\label{p85}%
\end{equation}

Set%
\begin{equation}
G_{F}:=\frac{1}{w\left\vert \cdot\right\vert ^{2\theta}}.\label{p86}%
\end{equation}

Also, since $\phi_{\infty}\in C_{\emptyset,2\theta}^{\infty}\cap C_{B}%
^{\infty}$ implies $\phi_{\infty}\phi\in S_{\emptyset,2\theta}$ when $\phi\in
S$, we have
\begin{align*}
\widehat{G_{\rho}}  & =\phi_{0}\widehat{G_{\rho}}+\phi_{\infty}\widehat
{G_{\rho}}=\phi_{0}\widehat{G_{\rho}}+\phi_{\infty}G_{F},\\
G_{\rho}  & =\overset{\vee}{\phi}_{0}\ast G_{\rho}+\overset{\vee}{\phi
_{\infty}}\ast G_{\rho},
\end{align*}

so that%
\begin{align*}
\left(  \left(  ia\xi\right)  ^{n+1}\overline{\widehat{g_{n}}}\left(
a\xi\right)  \widehat{G_{\rho}}\right)  ^{\vee}  & =\left(  \left(
ia\xi\right)  ^{n+1}\overline{\widehat{g_{n}}}\left(  a\xi\right)  \phi
_{0}\widehat{G_{\rho}}\right)  ^{\vee}+\left(  \left(  ia\xi\right)
^{n+1}\overline{\widehat{g_{n}}}\left(  a\xi\right)  \phi_{\infty}%
\widehat{G_{\rho}}\right)  ^{\vee}\\
& =\left(  \left(  ia\xi\right)  ^{n+1}\overline{\widehat{g_{n}}}\left(
a\xi\right)  \phi_{0}\widehat{G_{\rho}}\right)  ^{\vee}+\left(  \left(
ia\xi\right)  ^{n+1}\overline{\widehat{g_{n}}}\left(  a\xi\right)
\phi_{\infty}G_{F}\right)  ^{\vee}\\
& =\left(  \overline{\widehat{g_{n}}}\left(  a\xi\right)  \left(
ia\xi\right)  \phi_{0}\widehat{\left(  aD\right)  ^{n}G_{\rho}}\right)
^{\vee}+\left(  \left(  ia\xi\right)  ^{n+1}\overline{\widehat{g_{n}}}\left(
a\xi\right)  \phi_{\infty}G_{F}\right)  ^{\vee},
\end{align*}

and hence%
\begin{align}
& \left\Vert \left(  \left(  ia\xi\right)  ^{n+1}\overline{\widehat{g_{n}}%
}\left(  a\xi\right)  \widehat{G_{\rho}}\right)  ^{\vee}\right\Vert _{\infty
}\nonumber\\
& \leq\left\Vert \left(  \overline{\widehat{g_{n}}}\left(  a\xi\right)
\left(  ia\xi\right)  \phi_{0}\widehat{\left(  aD\right)  ^{n}G_{\rho}%
}\right)  ^{\vee}\right\Vert _{\infty}+\left\Vert \left(  \left(
ia\xi\right)  ^{n+1}\overline{\widehat{g_{n}}}\left(  a\xi\right)
\phi_{\infty}G_{F}\right)  ^{\vee}\right\Vert _{\infty}.\label{p28}%
\end{align}
\medskip

\underline{\textbf{Estimation of the }$G_{\rho}$\textbf{\ term on the RHS of
}\ref{p28}}\medskip

From Theorem \ref{Thm_deriv_Grho}, $G_{\rho}\in C_{BP}^{\left(  \left\lfloor
2\kappa\right\rfloor \right)  }$ and specifically from \ref{a43},
\[
D^{\gamma}G_{\rho}\left(  x\right)  =\left\{
\begin{array}
[c]{ll}%
\left(  2\pi\right)  ^{-\frac{d}{2}}\int\mathcal{Q}_{\emptyset,2\theta
-\left\vert \gamma\right\vert ,\xi}\left(  e^{ix\xi}\right)  \frac{\left(
i\xi\right)  ^{\gamma}}{w\left\vert \cdot\right\vert ^{2\theta}}d\xi, &
\left\vert \gamma\right\vert <2\theta,\\
\left(  2\pi\right)  ^{-\frac{d}{2}}\int e^{^{ix\xi}}\frac{\left(
i\xi\right)  ^{\gamma}}{w\left\vert \cdot\right\vert ^{2\theta}}d\xi, &
\left\vert \gamma\right\vert \geq2\theta,
\end{array}
\right.
\]

so that%
\begin{equation}
\left(  aD\right)  ^{n}G_{\rho}\left(  x\right)  =\left\{
\begin{array}
[c]{ll}%
\left(  2\pi\right)  ^{-\frac{d}{2}}\int\mathcal{Q}_{\emptyset,2\theta-n,\xi
}\left(  e^{ix\xi}\right)  \frac{\left(  ia\xi\right)  ^{n}}{w\left\vert
\cdot\right\vert ^{2\theta}}d\xi, & n<2\theta,\\
\left(  2\pi\right)  ^{-\frac{d}{2}}\int e^{^{ix\xi}}\frac{\left(
ia\xi\right)  ^{n}}{w\left\vert \cdot\right\vert ^{2\theta}}d\xi, &
n\geq2\theta.
\end{array}
\right. \label{p84}%
\end{equation}

From Lemma \ref{Lem_gm_properties_2} $\overline{\widehat{g_{n}}}\in
C_{B}^{\infty}\left(  \mathbb{R}^{1}\right)  $ which means that $\overline
{\widehat{g_{n}}}\left(  \left(  a,\cdot\right)  \right)  \in C_{B}^{\infty
}\left(  \mathbb{R}^{d}\right)  $. Define%
\begin{equation}
\left.
\begin{array}
[c]{l}%
\sigma_{a}\left(  \xi\right)  :=\frac{\sqrt{2\pi}}{n!}\overline{\widehat
{g_{n}}}\left(  a\xi\right)  \left(  ia\xi\right)  \phi_{0}\left(  \xi\right)
\in S,\\
\phi_{a}\left(  \xi\right)  :=\left(  ia\xi\right)  \phi_{0}\left(
\xi\right)  \in S,\\
\widehat{\phi_{a}}=aD\widehat{\phi_{0}},\text{ }\phi_{a}=\left\vert
a\right\vert \phi_{\widehat{a}},\quad\widehat{a}=a/\left\vert a\right\vert
\end{array}
\right\}  ,\label{a511}%
\end{equation}

and, since $\left(  aD\right)  ^{n}G_{\rho}\in C_{BP}^{\left(  0\right)
}\subset S^{\prime}$,\ by parts 1 and 2 of Definition \ref{Def_Convolution},%
\begin{align*}
\frac{\sqrt{2\pi}}{n!}\left(  \overline{\widehat{g_{n}}}\left(  a\xi\right)
\left(  ia\xi\right)  \phi_{0}\widehat{\left(  aD\right)  ^{n}G_{\rho}%
}\right)  ^{\vee}\left(  x\right)   & =\left(  \sigma_{a}\text{\thinspace
}\widehat{\left(  aD\right)  ^{n}G_{\rho}}\right)  ^{\vee}\left(  x\right) \\
& =\left(  \overset{\vee}{\sigma_{a}}\ast\left(  aD\right)  ^{n}G_{\rho
}\right)  \left(  x\right) \\
& =\left(  2\pi\right)  ^{-\frac{d}{2}}\int\overset{\vee}{\sigma_{a}}\left(
y\right)  \left(  \left(  aD\right)  ^{n}G_{\rho}\right)  \left(  x-y\right)
dy\\
& =\left(  2\pi\right)  ^{-\frac{d}{2}}\int\widehat{\sigma_{a}}\left(
-\xi\right)  \left(  \left(  aD\right)  ^{n}G_{\rho}\right)  \left(
x-\xi\right)  d\xi\\
& =\left(  2\pi\right)  ^{-\frac{d}{2}}\int\widehat{\sigma_{a}}\left(
\xi\right)  \left(  aD\right)  ^{n}G_{\rho}\left(  x+\xi\right)  d\xi\\
& =\left(  2\pi\right)  ^{-\frac{d}{2}}\left\vert a\right\vert ^{n}%
\int\widehat{\sigma_{a}}\left(  \xi\right)  \left(  \widehat{a}D\right)
^{n}G_{\rho}\left(  x+\xi\right)  d\xi.
\end{align*}

But because of the nice properties of $\phi_{a}$ and $g_{n}$ and
$\widehat{g_{n}}$ we can change the order of integration to obtain%
\begin{align*}
\widehat{\sigma_{a}}\left(  \xi\right)   & =\frac{\sqrt{2\pi}}{n!}\left(
2\pi\right)  ^{-\frac{d}{2}}\int e^{-ix\xi}\overline{\widehat{g_{n}}}\left(
ax\right)  \phi_{a}\left(  x\right)  dx\\
& =\frac{\sqrt{2\pi}}{n!}\left(  2\pi\right)  ^{-\frac{d}{2}}\int e^{-ix\xi
}\overline{\widehat{g_{n}}}\left(  ax\right)  \phi_{a}\left(  x\right)  dx\\
& =\frac{\sqrt{2\pi}}{n!}\left(  2\pi\right)  ^{-\frac{d+1}{2}}\int e^{-ix\xi
}\int_{0}^{1}e^{isax}g_{n}\left(  s\right)  ds\text{ }\phi_{a}\left(
x\right)  dx\\
& =\frac{\sqrt{2\pi}}{n!}\left(  2\pi\right)  ^{-\frac{d+1}{2}}\int_{0}%
^{1}\left(  \int e^{-i\left(  \xi-sa\right)  x}\phi_{a}\left(  x\right)
dx\right)  g_{n}\left(  s\right)  ds\\
& =\frac{\sqrt{2\pi}}{n!}\left(  2\pi\right)  ^{-\frac{1}{2}}\int_{0}%
^{1}\widehat{\phi_{a}}\left(  \xi-sa\right)  g_{n}\left(  s\right)  ds\\
& =\frac{1}{n!}\int_{0}^{1}\widehat{\phi_{a}}\left(  \xi-sa\right)
g_{n}\left(  s\right)  ds,
\end{align*}

and so%
\begin{align*}
&  \frac{\sqrt{2\pi}}{n!}\left(  \overline{\widehat{g_{n}}}\left(
a\xi\right)  \left(  ia\xi\right)  \phi_{0}\widehat{\left(  aD\right)
^{n}G_{\rho}}\right)  ^{\vee}\left(  x\right) \\
&  =\left(  2\pi\right)  ^{-\frac{d}{2}}\int\frac{\left\vert a\right\vert
^{n}}{n!}\int_{0}^{1}\widehat{\phi_{a}}\left(  \xi-sa\right)  g_{n}\left(
s\right)  ds\text{ }\left(  \widehat{a}D\right)  ^{n}G_{\rho}\left(
x+\xi\right)  d\xi\\
&  =\left(  2\pi\right)  ^{-\frac{d}{2}}\frac{\left\vert a\right\vert ^{n}%
}{n!}\int_{0}^{1}\int\widehat{\phi_{a}}\left(  \xi-sa\right)  \left(
\widehat{a}D\right)  ^{n}G_{\rho}\left(  x+\xi\right)  d\xi\text{ }%
g_{n}\left(  s\right)  ds\\
&  =\left(  2\pi\right)  ^{-\frac{d}{2}}\frac{\left\vert a\right\vert ^{n}%
}{n!}\int_{0}^{1}\int\widehat{\phi_{a}}\left(  \xi\right)  \left(  \left(
\widehat{a}D\right)  ^{n}G_{\rho}\right)  \left(  \xi+x+sa\right)  d\xi\text{
}g_{n}\left(  s\right)  ds\\
&  =\left(  2\pi\right)  ^{-\frac{d}{2}}\frac{\left\vert a\right\vert ^{n+1}%
}{n!}\int_{0}^{1}\int\widehat{\phi_{\widehat{a}}}\left(  \xi\right)  \left(
\left(  \widehat{a}D\right)  ^{n}G_{\rho}\right)  \left(  \xi+x+sa\right)
d\xi\text{ }g_{n}\left(  s\right)  ds,
\end{align*}

and hence%
\begin{align}
&  \left\vert \frac{\sqrt{2\pi}}{n!}\left(  \overline{\widehat{g_{n}}}\left(
a\xi\right)  \left(  ia\xi\right)  \phi_{0}\widehat{\left(  aD\right)
^{n}G_{\rho}}\right)  ^{\vee}\left(  x\right)  \right\vert \nonumber\\
&  \leq\left(  2\pi\right)  ^{-\frac{d}{2}}\frac{\left\vert a\right\vert
^{n+1}}{n!}\int_{0}^{1}\int\widehat{\phi_{\widehat{a}}}\left(  \xi\right)
\left(  \left(  \widehat{a}D\right)  ^{n}G_{\rho}\right)  \left(
\xi+x+sa\right)  d\xi\text{ }g_{n}\left(  s\right)  ds.\label{a801}%
\end{align}

Before proceeding we note that we have actually proved the more general
result: if $h\in C_{BP}^{\left(  0\right)  }$ then%
\begin{equation}
\left\vert \frac{\sqrt{2\pi}}{n!}\left(  \overline{\widehat{g_{n}}}\left(
a\xi\right)  \left(  ia\xi\right)  \phi_{0}\widehat{h}\right)  ^{\vee}\left(
x\right)  \right\vert \leq\left(  2\pi\right)  ^{-\frac{d}{2}}\frac{\left\vert
a\right\vert ^{n+1}}{n!}\int\limits_{0}^{1}\int\widehat{\phi_{\widehat{a}}%
}\left(  \xi\right)  h\left(  \xi+x+sa\right)  d\xi\text{\thinspace}%
g_{n}\left(  s\right)  ds.\label{a800}%
\end{equation}

The next step is to bound $\left\vert \left(  \widehat{a}D\right)  ^{n}%
G_{\rho}\left(  x\right)  \right\vert $ and there are two cases: $\left\lfloor
2\kappa\right\rfloor <2\theta$ and $\left\lfloor 2\kappa\right\rfloor
\geq2\theta$.\medskip

\fbox{\textbf{Case 1} $\left\lfloor 2\kappa\right\rfloor <2\theta$} From
\ref{p84},%
\begin{equation}
\left\vert \left(  \widehat{a}D\right)  ^{n}G_{\rho}\left(  x\right)
\right\vert \leq\left(  2\pi\right)  ^{-\frac{d}{2}}\int\left\vert
\mathcal{Q}_{\emptyset,2\theta-n,\xi}\left(  e^{ix\xi}\right)  \right\vert
\frac{\left\vert \widehat{a}\xi\right\vert ^{n}}{w\left\vert \cdot\right\vert
^{2\theta}}d\xi,\quad n\leq\left\lfloor 2\kappa\right\rfloor .\label{a161}%
\end{equation}

The next step is obtain an upper bound for the right side of this inequality
by first partitioning the domain of integration using the sphere $S_{r_{3}}$
and then using the upper bounds for $\mathcal{Q}_{\emptyset,2\theta-n,\xi
}\left(  e^{ix\xi}\right)  $ derived in Theorem \ref{Thm_bound_on_g(e,x)_2}.
From \ref{a702} with $\gamma=0$,%
\[
\left\vert \mathcal{Q}_{\emptyset,m,\xi}\left(  e^{ix\xi}\right)  \right\vert
\leq\underline{C}_{m,r}^{\left(  \rho\right)  }\left\vert \xi\right\vert
^{m}\left(  1+\frac{\left\vert x\widehat{\xi}\right\vert ^{m}}{m!}\right)
,\quad m\geq1,\text{ }\left\vert \xi\right\vert \leq r,
\]

so that%
\begin{align}
&  \int\limits_{\left\vert \cdot\right\vert \leq r_{3}}\left\vert
\mathcal{Q}_{\emptyset,2\theta-n,\xi}\left(  e^{ix\xi}\right)  \right\vert
\frac{\left\vert \widehat{a}\xi\right\vert ^{n}}{w\left\vert \cdot\right\vert
^{2\theta}}d\xi\nonumber\\
&  \leq\int\limits_{\left\vert \cdot\right\vert \leq r_{3}}\underline
{C}_{2\theta-n,r_{3}}^{\left(  \rho\right)  }\left\vert \cdot\right\vert
^{2\theta-n}\left(  1+\frac{\left\vert x\widehat{\xi}\right\vert ^{2\theta-n}%
}{\left(  2\theta-n\right)  !}\right)  \frac{\left\vert \widehat{a}%
\xi\right\vert ^{n}}{w\left\vert \cdot\right\vert ^{2\theta}}\nonumber\\
&  =\underline{C}_{2\theta-n,r_{3}}^{\left(  \rho\right)  }\int%
\limits_{\left\vert \cdot\right\vert \leq r_{3}}\left\vert \cdot\right\vert
^{2\theta-n}\frac{\left\vert \widehat{a}\xi\right\vert ^{n}}{w\left\vert
\cdot\right\vert ^{2\theta}}+\underline{C}_{\rho;2\theta-n,r_{3}}%
\int\limits_{\left\vert \cdot\right\vert \leq r_{3}}\left\vert \cdot
\right\vert ^{2\theta-n}\frac{\left\vert x\widehat{\xi}\right\vert
^{2\theta-n}}{\left(  2\theta-n\right)  !}\frac{\left\vert \widehat{a}%
\xi\right\vert ^{n}}{w\left\vert \cdot\right\vert ^{2\theta}}\nonumber\\
&  =\underline{C}_{2\theta-n,r_{3}}^{\left(  \rho\right)  }\left(
\int\limits_{\left\vert \cdot\right\vert \leq r_{3}}\frac{\left\vert
\widehat{a}\widehat{\xi}\right\vert ^{n}}{w}+\frac{\left\vert x\right\vert
^{2\theta-n}}{\left(  2\theta-n\right)  !}\int\limits_{\left\vert
\cdot\right\vert \leq r_{3}}\frac{\left\vert \widehat{x}\widehat{\xi
}\right\vert ^{2\theta-n}\left\vert \widehat{a}\widehat{\xi}\right\vert ^{n}%
}{w}\right) \nonumber\\
&  \leq\underline{C}_{2\theta-n,r_{3}}^{\left(  \rho\right)  }\left(
\int\limits_{\left\vert \cdot\right\vert \leq r_{3}}\frac{\left\vert
\widehat{a}\widehat{\xi}\right\vert ^{n}}{w}+\frac{\left\vert x\right\vert
^{2\theta-n}}{\left(  2\theta-n\right)  !}\min\left\{  \int\limits_{\left\vert
\cdot\right\vert \leq r_{3}}\frac{\left\vert \widehat{x}\widehat{\xi
}\right\vert ^{2\theta-n}}{w},\int\limits_{\left\vert \cdot\right\vert \leq
r_{3}}\frac{\left\vert \widehat{a}\widehat{\xi}\right\vert ^{n}}{w}\right\}
\right) \label{a631}\\
&  <\infty,\nonumber
\end{align}

since $w\in W2.1$. Note that the forms \ref{a63} and indeed \ref{a631} will be
useful when the weight function is radial.

From \ref{a27.2},%
\[
\left\vert \mathcal{Q}_{\emptyset,m,\xi}\left(  e^{ix\xi}\right)  \right\vert
\leq\overline{C}_{m,r}^{\left(  \rho\right)  }\left(  1+\sum\limits_{k<m}%
\frac{\left\vert x\widehat{\xi}\right\vert ^{k}}{k!}\right)  ,\quad
m\geq1,\text{ }\left\vert \xi\right\vert \geq r,
\]

so that%
\begin{align}
&  \int\limits_{\left\vert \cdot\right\vert \geq r_{3}}\left\vert
\mathcal{Q}_{\emptyset,2\theta-n,\xi}\left(  e^{ix\xi}\right)  \right\vert
\frac{\left\vert \widehat{a}\xi\right\vert ^{n}}{w\left\vert \cdot\right\vert
^{2\theta}}d\xi\nonumber\\
&  \leq\overline{C}_{2\theta-n,r_{3}}^{\left(  \rho\right)  }\int%
\limits_{\left\vert \cdot\right\vert \geq r_{3}}\left(  1+\sum
\limits_{k<2\theta-n}\frac{\left\vert x\widehat{\xi}\right\vert ^{k}}%
{k!}\right)  \frac{\left\vert \widehat{a}\xi\right\vert ^{n}}{w\left\vert
\cdot\right\vert ^{2\theta}}\nonumber\\
&  =\overline{C}_{2\theta-n,r_{3}}^{\left(  \rho\right)  }\left(
\int\limits_{\left\vert \cdot\right\vert \geq r_{3}}\frac{\left\vert
\widehat{a}\xi\right\vert ^{n}}{w\left\vert \cdot\right\vert ^{2\theta}}%
+\sum\limits_{k<2\theta-n}\frac{\left\vert x\right\vert ^{k}}{k!}%
\int\limits_{\left\vert \cdot\right\vert \geq r_{3}}\frac{\left\vert
\widehat{x}\widehat{\xi}\right\vert ^{k}\left\vert \widehat{a}\xi\right\vert
^{n}}{w\left\vert \cdot\right\vert ^{2\theta}}\right) \label{a67}\\
&  =\overline{C}_{2\theta-n,r_{3}}^{\left(  \rho\right)  }\left(
\int\limits_{\left\vert \cdot\right\vert \geq r_{3}}\frac{\left\vert
\widehat{a}\xi\right\vert ^{n}}{w\left\vert \cdot\right\vert ^{2\theta}}%
+\sum\limits_{k<2\theta-n}\frac{\left\vert x\right\vert ^{k}}{k!}%
\int\limits_{\left\vert \cdot\right\vert \geq r_{3}}\frac{\left\vert
\widehat{x}\widehat{\xi}\right\vert ^{k}\left\vert \widehat{a}\widehat{\xi
}\right\vert ^{n}\left\vert \cdot\right\vert ^{n}}{w\left\vert \cdot
\right\vert ^{2\theta}}\right) \nonumber\\
&  \leq\overline{C}_{2\theta-n,r_{3}}^{\left(  \rho\right)  }\left(
\int\limits_{\left\vert \cdot\right\vert \geq r_{3}}\frac{\left\vert
\widehat{a}\widehat{\xi}\right\vert ^{n}\left\vert \xi\right\vert ^{n}%
}{w\left\vert \cdot\right\vert ^{2\theta}}+\sum\limits_{k<2\theta-n}%
\frac{\left\vert x\right\vert ^{k}}{k!}\min\left\{  \int\limits_{\left\vert
\cdot\right\vert \geq r_{3}}\frac{\left\vert \widehat{x}\widehat{\xi
}\right\vert ^{k}\left\vert \cdot\right\vert ^{n}}{w\left\vert \cdot
\right\vert ^{2\theta}},\int\limits_{\left\vert \cdot\right\vert \geq r_{3}%
}\frac{\left\vert \widehat{a}\widehat{\xi}\right\vert ^{n}\left\vert
\cdot\right\vert ^{n}}{w\left\vert \cdot\right\vert ^{2\theta}}\right\}
\right) \nonumber\\
&  <\infty,\nonumber
\end{align}

since $w\in W3.2$. Note that the form \ref{a67} will be useful when the weight
function is radial.

Combining these two estimates with \ref{a161} yields%
\begin{align}
& \left\vert \left(  \widehat{a}D\right)  ^{n}G_{\rho}\left(  x\right)
\right\vert \nonumber\\
& \leq\left(  2\pi\right)  ^{-\frac{d}{2}}\int\left\vert \mathcal{Q}%
_{\emptyset,2\theta-n,\xi}\left(  e^{ix\xi}\right)  \right\vert \frac
{\left\vert \widehat{a}\xi\right\vert ^{n}}{w\left\vert \cdot\right\vert
^{2\theta}}d\xi\nonumber\\
& =\left(  2\pi\right)  ^{-\frac{d}{2}}\int\limits_{\left\vert \xi\right\vert
\leq r_{3}}\left\vert \mathcal{Q}_{\emptyset,2\theta-n,\xi}\left(  e^{ix\xi
}\right)  \right\vert \frac{\left\vert \widehat{a}\xi\right\vert ^{n}%
}{w\left\vert \cdot\right\vert ^{2\theta}}d\xi+\left(  2\pi\right)
^{-\frac{d}{2}}\int\limits_{\left\vert \xi\right\vert \geq r_{3}}\left\vert
\mathcal{Q}_{\emptyset,2\theta-n;\xi}\left(  e^{ix\xi}\right)  \right\vert
\frac{\left\vert \widehat{a}\xi\right\vert ^{n}}{w\left\vert \cdot\right\vert
^{2\theta}}d\xi\nonumber\\
& \leq\left(  2\pi\right)  ^{-\frac{d}{2}}\underline{C}_{2\theta-n,r_{3}%
}^{\left(  \rho\right)  }\left(  \int\limits_{\left\vert \cdot\right\vert \leq
r_{3}}\frac{\left\vert \widehat{a}\widehat{\xi}\right\vert ^{n}}{w}%
+\sum\limits_{\left\vert \beta\right\vert =2\theta-n}\frac{x_{+}^{\beta}%
}{\beta!}\int\limits_{\left\vert \cdot\right\vert \leq r_{3}}\frac{\left\vert
\widehat{\xi}^{\beta}\right\vert \left\vert \widehat{a}\widehat{\xi
}\right\vert ^{n}}{w}\right)  +\nonumber\\
& \qquad+\left(  2\pi\right)  ^{-\frac{d}{2}}\overline{C}_{2\theta-n,r_{3}%
}^{\left(  \rho\right)  }\left(  \int\limits_{\left\vert \cdot\right\vert \geq
r_{3}}\frac{\left\vert \widehat{a}\xi\right\vert ^{n}}{w\left\vert
\cdot\right\vert ^{2\theta}}+\sum\limits_{k<2\theta-n}\sum\limits_{\left\vert
\beta\right\vert =k}\frac{x_{+}^{\beta}}{\beta!}\int\limits_{\left\vert
\cdot\right\vert \geq r_{3}}\frac{\left\vert \widehat{\xi}^{\beta}\right\vert
\left\vert \widehat{a}\xi\right\vert ^{n}}{w\left\vert \cdot\right\vert
^{2\theta}}\right) \nonumber\\
& \leq\left(  2\pi\right)  ^{-\frac{d}{2}}\max\left\{  \underline{C}%
_{2\theta-n,r_{3}}^{\left(  \rho\right)  },\overline{C}_{2\theta-n,r_{3}%
}^{\left(  \rho\right)  }\right\}  \left(
\begin{array}
[c]{c}%
\int\limits_{\left\vert \cdot\right\vert \leq r_{3}}\frac{\left\vert
\widehat{a}\widehat{\xi}\right\vert ^{n}}{w}+\sum\limits_{\left\vert
\beta\right\vert =2\theta-n}\frac{x_{+}^{\beta}}{\beta!}\int%
\limits_{\left\vert \cdot\right\vert \leq r_{3}}\frac{\left\vert \widehat{\xi
}^{\beta}\right\vert \left\vert \widehat{a}\widehat{\xi}\right\vert ^{n}}%
{w}+\int\limits_{\left\vert \cdot\right\vert \geq r_{3}}\frac{\left\vert
\widehat{a}\xi\right\vert ^{n}}{w\left\vert \cdot\right\vert ^{2\theta}}+\\
+\sum\limits_{k<2\theta-n}\sum\limits_{\left\vert \beta\right\vert =k}%
\frac{x_{+}^{\beta}}{\beta!}\int\limits_{\left\vert \cdot\right\vert \geq
r_{3}}\frac{\left\vert \widehat{\xi}^{\beta}\right\vert \left\vert \widehat
{a}\xi\right\vert ^{n}}{w\left\vert \cdot\right\vert ^{2\theta}}%
\end{array}
\right) \nonumber\\
& =\left(  2\pi\right)  ^{-\frac{d}{2}}C_{2\theta-n,r_{3}}^{\left(
\rho\right)  }\left(
\begin{array}
[c]{c}%
\int\limits_{\left\vert \cdot\right\vert \leq r_{3}}\frac{\left\vert
\widehat{a}\widehat{\xi}\right\vert ^{n}}{w}+2\int\limits_{\left\vert
\cdot\right\vert \geq r_{3}}\frac{\left\vert \widehat{a}\xi\right\vert ^{n}%
}{w\left\vert \cdot\right\vert ^{2\theta}}+\sum\limits_{\substack{k<2\theta-n
\\k>0}}\sum\limits_{\left\vert \beta\right\vert =k}\frac{x_{+}^{\beta}}%
{\beta!}\int\limits_{\left\vert \cdot\right\vert \geq r_{3}}\frac{\left\vert
\widehat{\xi}^{\beta}\right\vert \left\vert \widehat{a}\xi\right\vert ^{n}%
}{w\left\vert \cdot\right\vert ^{2\theta}}+\\
+\sum\limits_{\left\vert \beta\right\vert =2\theta-n}\frac{x_{+}^{\beta}%
}{\beta!}\int\limits_{\left\vert \cdot\right\vert \leq r_{3}}\frac{\left\vert
\widehat{\xi}^{\beta}\right\vert \left\vert \widehat{a}\widehat{\xi
}\right\vert ^{n}}{w}%
\end{array}
\right) \nonumber\\
& =\left(  2\pi\right)  ^{-\frac{d}{2}}C_{2\theta-n,r_{3}}^{\left(
\rho\right)  }\sum\limits_{\left\vert \beta\right\vert \leq2\theta-n}%
\upsilon_{n,\beta}^{\left(  w\right)  }\frac{x_{+}^{\beta}}{\beta
!},\label{a1.09}%
\end{align}

where%
\begin{equation}
C_{2\theta-n,r_{3}}^{\left(  \rho\right)  }=\max\left\{  \underline
{C}_{2\theta-n,r_{3}}^{\left(  \rho\right)  },\overline{C}_{2\theta-n,r_{3}%
}^{\left(  \rho\right)  }\right\}  ,\label{a1.07}%
\end{equation}

and%
\begin{equation}
\upsilon_{n,\beta}^{\left(  w\right)  }=\left\{
\begin{array}
[c]{ll}%
\int\limits_{\left\vert \cdot\right\vert \leq r_{3}}\frac{\left\vert
\widehat{a}\widehat{\xi}\right\vert ^{n}}{w}+2\int\limits_{\left\vert
\cdot\right\vert \geq r_{3}}\frac{\left\vert \widehat{a}\xi\right\vert ^{n}%
}{w\left\vert \cdot\right\vert ^{2\theta}}, & \beta=0,\\
\int\limits_{\left\vert \cdot\right\vert \geq r_{3}}\frac{\left\vert
\widehat{\xi}^{\beta}\right\vert \left\vert \widehat{a}\xi\right\vert ^{n}%
}{w\left\vert \cdot\right\vert ^{2\theta}}, & 0<\left\vert \beta\right\vert
<2\theta-n,\\
\int\limits_{\left\vert \cdot\right\vert \leq r_{3}}\frac{\left\vert
\widehat{\xi}^{\beta}\right\vert \left\vert \widehat{a}\widehat{\xi
}\right\vert ^{n}}{w}, & \left\vert \beta\right\vert =2\theta-n.
\end{array}
\right. \label{a1.08}%
\end{equation}
\medskip

\fbox{\textbf{Case 2} $\left\lfloor 2\kappa\right\rfloor \geq2\theta$} From
\ref{p84},%
\begin{align}
\left\vert \left(  \widehat{a}D\right)  ^{n}G_{\rho}\left(  x\right)
\right\vert  & \leq\left(  2\pi\right)  ^{-\frac{d}{2}}\int\frac{\left\vert
\widehat{a}\xi\right\vert ^{n}}{w\left\vert \cdot\right\vert ^{2\theta}}%
d\xi\label{a561}\\
& =\left(  2\pi\right)  ^{-\frac{d}{2}}\left(  \int_{\left\vert \cdot
\right\vert \leq r_{3}}\frac{\left\vert \widehat{a}\xi\right\vert ^{n}%
}{w\left\vert \cdot\right\vert ^{2\theta}}+\int_{\left\vert \cdot\right\vert
\geq r_{3}}\frac{\left\vert \widehat{a}\xi\right\vert ^{n}}{w\left\vert
\cdot\right\vert ^{2\theta}}\right) \nonumber\\
& \leq\left(  2\pi\right)  ^{-\frac{d}{2}}\left(  \int_{\left\vert
\cdot\right\vert \leq r_{3}}\frac{\left\vert \xi\right\vert ^{n}}{w\left\vert
\cdot\right\vert ^{2\theta}}+\int_{\left\vert \cdot\right\vert \geq r_{3}%
}\frac{\left\vert \xi\right\vert ^{n}}{w\left\vert \cdot\right\vert ^{2\theta
}}\right) \nonumber\\
& =\left(  2\pi\right)  ^{-\frac{d}{2}}\left(  \int_{\left\vert \cdot
\right\vert \leq r_{3}}\frac{\left\vert \cdot\right\vert ^{n-2\theta}}{w}%
+\int_{\left\vert \cdot\right\vert \geq r_{3}}\frac{1}{\left\vert
\cdot\right\vert ^{2\kappa-n}}\frac{\left\vert \cdot\right\vert ^{2\kappa}%
}{w\left\vert \cdot\right\vert ^{2\theta}}\right) \nonumber\\
& \leq\left(  2\pi\right)  ^{-\frac{d}{2}}\left(  r_{3}^{n-2\theta}%
\int_{\left\vert \cdot\right\vert \leq r_{3}}\frac{1}{w}+\frac{1}%
{r_{3}^{2\kappa-n}}\int_{\left\vert \cdot\right\vert \geq r_{3}}%
\frac{\left\vert \cdot\right\vert ^{2\kappa}}{w\left\vert \cdot\right\vert
^{2\theta}}\right) \nonumber\\
& <\infty.\nonumber
\end{align}

This completes our upper bounds for $\left\vert \left(  \widehat{a}D\right)
^{n}G_{\rho}\left(  x\right)  \right\vert $ and now we have for $n\leq
\left\lfloor 2\kappa\right\rfloor $:%
\begin{equation}
\left\vert \left(  \widehat{a}D\right)  ^{n}G_{\rho}\left(  x\right)
\right\vert \leq\left\{
\begin{array}
[c]{ll}%
\left(  2\pi\right)  ^{-\frac{d}{2}}C_{2\theta-n,r_{3}}^{\left(  \rho\right)
}\sum\limits_{\left\vert \beta\right\vert \leq2\theta-n}\upsilon_{n,\beta
}^{\left(  w\right)  }\frac{\left\vert x^{\beta}\right\vert }{\beta!}, &
\left\lfloor 2\kappa\right\rfloor <2\theta,\\
\left(  2\pi\right)  ^{-\frac{d}{2}}\int\frac{\left\vert \widehat{a}%
\xi\right\vert ^{n}}{w\left\vert \cdot\right\vert ^{2\theta}}, & \left\lfloor
2\kappa\right\rfloor \geq2\theta,
\end{array}
\right. \label{a000}%
\end{equation}
\medskip

\underline{\textbf{Apply these estimates for }$\left\vert \left(  \widehat
{a}D\right)  ^{n}G_{\rho}\left(  x\right)  \right\vert $\textbf{\ to the right
side of \ref{a801}}}\medskip

We adapt Section \ref{Sect_data_fn_Taylor_W3.2} (data functions). There are
two cases: $\left\lfloor 2\kappa\right\rfloor <2\theta$ and $\left\lfloor
2\kappa\right\rfloor \geq2\theta$.\medskip

\fbox{\textbf{Case 1} $n\leq\left\lfloor 2\kappa\right\rfloor <2\theta$}
Estimating the right side of \ref{a801} using \ref{a000} yields%
\begin{align}
&  \left\vert \frac{\sqrt{2\pi}}{n!}\left(  \overline{\widehat{g_{n}}}\left(
a\xi\right)  \left(  ia\xi\right)  \phi_{0}\widehat{\left(  aD\right)
^{n}G_{\rho}}\right)  ^{\vee}\left(  x\right)  \right\vert \nonumber\\
&  \leq\left(  2\pi\right)  ^{-\frac{d}{2}}\frac{\left\vert a\right\vert
^{n+1}}{n!}\int_{0}^{1}\int\left\vert \widehat{\phi_{\widehat{a}}}\left(
\xi\right)  \right\vert \left\vert \left(  \left(  \widehat{a}D\right)
^{n}G_{\rho}\right)  \left(  \xi+x+sa\right)  \right\vert d\xi\text{ }%
g_{n}\left(  s\right)  ds\nonumber\\
&  \leq??\left(  2\pi\right)  ^{-\frac{d}{2}}\frac{\left\vert a\right\vert
^{n+1}}{n!}\int_{0}^{1}\int\left\vert \widehat{\phi_{\widehat{a}}}\left(
\xi\right)  \right\vert \left(  2\pi\right)  ^{-\frac{d}{2}}C_{2\theta
-n,r_{3}}^{\left(  \rho\right)  }\sum\limits_{\left\vert \alpha\right\vert
\leq2\theta-n}\upsilon_{n,\alpha}^{\left(  w\right)  }\frac{\left\vert \left(
\xi+x+sa\right)  ^{\alpha}\right\vert }{\alpha!}d\xi\text{ }g_{n}\left(
s\right)  ds\nonumber\\
&  =\frac{C_{2\theta-n,r_{3}}^{\left(  \rho\right)  }}{\left(  2\pi\right)
^{d}}\frac{\left\vert a\right\vert ^{n+1}}{n!}\int_{0}^{1}\int\left\vert
\widehat{\phi_{\widehat{a}}}\left(  \xi\right)  \right\vert \sum
\limits_{\left\vert \alpha\right\vert \leq2\theta-n}\upsilon_{n,\alpha
}^{\left(  w\right)  }\frac{\left\vert \left(  \xi+x+sa\right)  ^{\alpha
}\right\vert }{\alpha!}d\xi\text{ }g_{n}\left(  s\right)  ds\nonumber\\
&  =\frac{C_{2\theta-n,r_{3}}^{\left(  \rho\right)  }}{\left(  2\pi\right)
^{d}}\frac{\left\vert a\right\vert ^{n+1}}{n!}\sum\limits_{\left\vert
\alpha\right\vert \leq2\theta-n}\upsilon_{n,\alpha}^{\left(  w\right)  }%
\int_{0}^{1}\int\left\vert \frac{\left(  \xi+x+sa\right)  ^{\alpha}}{\alpha
!}\widehat{\phi_{\widehat{a}}}\left(  \xi\right)  \right\vert d\xi\text{
}g_{n}\left(  s\right)  ds\nonumber\\
&  =\frac{C_{2\theta-n,r_{3}}^{\left(  \rho\right)  }}{\left(  2\pi\right)
^{d}}\frac{\left\vert a\right\vert ^{n+1}}{n!}\sum\limits_{\left\vert
\alpha\right\vert \leq2\theta-n}\upsilon_{n,\alpha}^{\left(  w\right)  }%
\int_{0}^{1}\int\left\vert \frac{\left(  \xi+x+sa\right)  ^{\alpha}}{\alpha
!}\widehat{\phi_{\widehat{a}}}\left(  \xi\right)  \right\vert d\xi\text{
}g_{n}\left(  s\right)  ds.\label{a005}%
\end{align}

From part 9 of Summary \ref{Sum_Multi_index} i.e. by using the binomial
theorem,%
\[
\frac{\left(  x+sa+\xi\right)  ^{\alpha}}{\alpha!}=\sum_{\beta+\gamma
+\delta=\alpha}\frac{x^{\beta}\left(  sa\right)  ^{\gamma}\xi^{\delta}}%
{\beta!\gamma!\delta!}=\sum_{\beta+\gamma+\delta=\alpha}\frac{x^{\beta
}a^{\gamma}\xi^{\delta}}{\beta!\gamma!\delta!}s^{\left\vert \gamma\right\vert
},
\]

so by using the notation $\left(  \xi_{k}\right)  _{+}=\left(  \left\vert
\xi_{k}\right\vert \right)  $,%
\begin{align}
&  \int_{0}^{1}\int\left\vert \frac{\left(  \xi+x+sa\right)  ^{\alpha}}%
{\alpha!}\widehat{\phi_{\widehat{a}}}\left(  \xi\right)  \right\vert
d\xi\text{ }g_{n}\left(  s\right)  ds\nonumber\\
&  =\int_{0}^{1}\int\left\vert \sum_{\beta+\gamma+\delta=\alpha}\frac
{x^{\beta}a^{\gamma}\xi^{\delta}}{\beta!\gamma!\delta!}s^{\left\vert
\gamma\right\vert }\widehat{\phi_{\widehat{a}}}\left(  \xi\right)  \right\vert
d\xi ds\nonumber\\
&  \leq\int_{0}^{1}\int\sum_{\beta+\gamma+\delta=\alpha}\frac{x_{+}^{\beta
}a_{+}^{\gamma}\xi_{+}^{\delta}}{\beta!\gamma!\delta!}s^{\left\vert
\gamma\right\vert }\left\vert \widehat{\phi_{\widehat{a_{+}}}}\left(
\xi\right)  \right\vert d\xi\text{ }g_{n}\left(  s\right)  ds\nonumber\\
&  =\sum_{\beta+\gamma+\delta=\alpha}\frac{x_{+}^{\beta}a_{+}^{\gamma}}%
{\beta!\gamma!}\left(  \int\frac{\xi_{+}^{\delta}}{\delta!}\left\vert
\widehat{\phi_{\widehat{a}}}\left(  \xi\right)  \right\vert d\xi\right)
\int_{0}^{1}s^{\left\vert \gamma\right\vert }g_{n}\left(  s\right)
ds\nonumber\\
&  =\sum_{\beta+\gamma+\delta=\alpha}\frac{x_{+}^{\beta}a_{+}^{\gamma}}%
{\beta!\gamma!}\left(  \int\frac{\xi_{+}^{\delta}}{\delta!}\left\vert
\widehat{\phi_{\widehat{a}}}\left(  \xi\right)  \right\vert d\xi\right)
\int_{0}^{1}s^{\left\vert \gamma\right\vert }\left(  1-s\right)
^{n}ds\nonumber\\
&  \leq\frac{1}{n+1}\sum_{\beta+\gamma+\delta=\alpha}\frac{x_{+}^{\beta}%
a_{+}^{\gamma}}{\beta!\gamma!}\int\frac{\xi_{+}^{\delta}}{\delta!}\left\vert
\widehat{\phi_{\widehat{a}}}\left(  \xi\right)  \right\vert d\xi\nonumber\\
&  \leq\frac{1}{n+1}\sum\limits_{\delta\leq\alpha}\sum_{\beta+\gamma\leq
\alpha}\frac{x_{+}^{\beta}a_{+}^{\gamma}}{\beta!\gamma!}\int\frac{\xi
_{+}^{\delta}}{\delta!}\left\vert \widehat{\phi_{\widehat{a}}}\left(
\xi\right)  \right\vert d\xi\nonumber\\
&  =\frac{1}{n+1}\left(  \sum_{\beta+\gamma\leq\alpha}\frac{x_{+}^{\beta}%
a_{+}^{\gamma}}{\beta!\gamma!}\right)  \int\sum\limits_{\delta\leq\alpha}%
\frac{\xi_{+}^{\delta}}{\delta!}\left\vert \widehat{\phi_{\widehat{a}}}\left(
\xi\right)  \right\vert d\xi\nonumber\\
&  =\frac{1}{n+1}\left(  \sum\limits_{\sigma\leq\alpha}\sum_{\beta
+\gamma=\sigma}\frac{x_{+}^{\beta}a_{+}^{\gamma}}{\beta!\gamma!}\right)
\int\sum\limits_{\delta\leq\alpha}\frac{\xi_{+}^{\delta}}{\delta!}\left\vert
\widehat{\phi_{\widehat{a}}}\left(  \xi\right)  \right\vert d\xi\nonumber\\
&  =\frac{1}{n+1}\sum\limits_{\sigma\leq\alpha}\frac{\left(  x_{+}%
+a_{+}\right)  ^{\sigma}}{\sigma!}\int\sum\limits_{\delta\leq\alpha}\frac
{\xi_{+}^{\delta}}{\delta!}\left\vert \widehat{\phi_{\widehat{a}}}\left(
\xi\right)  \right\vert d\xi.\label{a1.19}%
\end{align}

and thus%
\begin{align}
\sum\limits_{\left\vert \alpha\right\vert \leq2\theta-n}\upsilon_{n,\alpha
}^{\left(  w\right)  }\int_{0}^{1}\int &  \left\vert \frac{\left(
\xi+x+sa\right)  ^{\alpha}}{\alpha!}\widehat{\phi_{\widehat{a}}}\left(
\xi\right)  \right\vert d\xi\text{ }g_{n}\left(  s\right)  ds\nonumber\\
&  \leq\sum\limits_{\left\vert \alpha\right\vert \leq2\theta-n}\frac
{\upsilon_{n,\alpha}^{\left(  w\right)  }}{n+1}\sum\limits_{\sigma\leq\alpha
}\frac{\left(  x_{+}+a_{+}\right)  ^{\sigma}}{\sigma!}\int\sum\limits_{\delta
\leq\alpha}\frac{\xi_{+}^{\delta}}{\delta!}\left\vert \widehat{\phi
_{\widehat{a}}}\left(  \xi\right)  \right\vert d\xi\nonumber\\
&  \leq\frac{1}{n+1}\sum\limits_{\left\vert \alpha\right\vert \leq2\theta
-n}\upsilon_{n,\alpha}^{\left(  w\right)  }\sum\limits_{\sigma\leq\alpha}%
\frac{\left(  x_{+}+a_{+}\right)  ^{\sigma}}{\sigma!}\max_{\left\vert
\delta\right\vert \leq2\theta-n}\int\sum\limits_{\delta\leq\alpha}\frac
{\xi_{+}^{\delta}}{\delta!}\left\vert \widehat{\phi_{\widehat{a}}}\left(
\xi\right)  \right\vert d\xi\nonumber\\
&  =\frac{1}{n+1}\sum\limits_{\left\vert \alpha\right\vert \leq2\theta
-n}\upsilon_{n,\alpha}^{\left(  w\right)  }\sum\limits_{\sigma\leq\alpha}%
\frac{\left(  x_{+}+a_{+}\right)  ^{\sigma}}{\sigma!}\int\sum
\limits_{\left\vert \delta\right\vert \leq2\theta-n}\frac{\xi_{+}^{\delta}%
}{\delta!}\left\vert \widehat{\phi_{\widehat{a}}}\right\vert \nonumber\\
&  =\frac{1}{n+1}\sum\limits_{\left\vert \alpha\right\vert \leq2\theta
-n}\upsilon_{n,\alpha}^{\left(  w\right)  }\sum\limits_{\sigma\leq\alpha}%
\frac{\left(  x_{+}+a_{+}\right)  ^{\sigma}}{\sigma!}\sum\limits_{j=0}%
^{2\theta-n}\int\frac{\left(  \xi_{+}\mathbf{1}\right)  ^{j}}{j!}\left\vert
\widehat{\phi_{\widehat{a}}}\right\vert .\label{a54}%
\end{align}

But $\phi_{\widehat{a}}\left(  \xi\right)  =\left(  i\widehat{a}\xi\right)
\phi_{0}\left(  \xi\right)  $ implies%
\begin{equation}
\widehat{\phi_{\widehat{a}}}\left(  \xi\right)  =\widehat{\widehat{a}x\phi
_{0}}\left(  \xi\right)  =\widehat{a}D\widehat{\phi_{0}}\left(  \xi\right)
,\label{a1.22}%
\end{equation}

so that
\begin{align}
\sum\limits_{\left\vert \alpha\right\vert \leq2\theta-n}\upsilon_{n,\alpha
}^{\left(  w\right)  }\int_{0}^{1}\int &  \left\vert \frac{\left(
\xi+x+sa\right)  ^{\alpha}}{\alpha!}\widehat{\phi_{\widehat{a}}}\left(
\xi\right)  \right\vert d\xi\text{ }g_{n}\left(  s\right)  ds\nonumber\\
&  \leq\frac{1}{n+1}\left(  \sum\limits_{\left\vert \alpha\right\vert
\leq2\theta-n}\upsilon_{n,\alpha}^{\left(  w\right)  }\sum\limits_{\sigma
\leq\alpha}\frac{\left(  x_{+}+a_{+}\right)  ^{\sigma}}{\sigma!}\right)
\sum\limits_{j=0}^{2\theta-n}\int\frac{\left(  \xi_{+}\mathbf{1}\right)  ^{j}%
}{j!}\left\vert \widehat{a}D\widehat{\phi_{0}}\left(  \xi\right)  \right\vert
.\label{a112}%
\end{align}

But the identity of part 5 of Theorem \ref{Thm_Ap_1} is,%
\[
\sum\limits_{\left\vert \alpha\right\vert \leq m}y^{\alpha}\sum\limits_{\beta
\leq\alpha}\frac{x^{\beta}}{\beta!}=\sum\limits_{j=0}^{m}\left(
\sum\limits_{k=0}^{m-j}h_{k}\left(  y\right)  \right)  \frac{\left(
xy\right)  ^{j}}{j!}=\sum\limits_{j=0}^{m}h_{j}\left(  y\right)
\sum\limits_{k=0}^{m-j}\frac{\left(  xy\right)  ^{k}}{k!},
\]

and by using \ref{a1.08} we get,%
\begin{align*}
& \sum\limits_{\left\vert \alpha\right\vert \leq2\theta-n}\upsilon_{n,\alpha
}^{\left(  w\right)  }\sum\limits_{\sigma\leq\alpha}\frac{\left(  x_{+}%
+a_{+}\right)  ^{\sigma}}{\sigma!}\\
& =\sum\limits_{\left\vert \alpha\right\vert \leq2\theta-n-1}\int%
\limits_{\left\vert \cdot\right\vert \geq r_{3}}\frac{\left\vert \widehat{\xi
}^{\alpha}\right\vert \left\vert \widehat{a}\xi\right\vert ^{n}}{w\left\vert
\cdot\right\vert ^{2\theta}}\sum\limits_{\sigma\leq\alpha}\frac{\left(
x_{+}+a_{+}\right)  ^{\sigma}}{\sigma!}+\ldots\\
& =\int\limits_{\left\vert \cdot\right\vert \geq r_{3}}\sum\limits_{\left\vert
\alpha\right\vert \leq2\theta-n-1}\widehat{\xi}_{+}^{\alpha}\sum
\limits_{\sigma\leq\alpha}\frac{\left(  x_{+}+a_{+}\right)  ^{\sigma}}%
{\sigma!}\frac{\left\vert \widehat{a}\xi\right\vert ^{n}}{w\left\vert
\cdot\right\vert ^{2\theta}}d\xi+\ldots\\
& =\int\limits_{\left\vert \cdot\right\vert \geq r_{3}}\sum\limits_{j=0}%
^{2\theta-n-1}\left(
{\textstyle\sum\limits_{k=0}^{2\theta-n-1-j}}
h_{k}\left(  \widehat{\xi}_{+}\right)  \right)  \tfrac{\left(  \left(
x_{+}+a_{+}\right)  \widehat{\xi}_{+}\right)  ^{j}}{j!}\tfrac{\left\vert
\widehat{a}\xi\right\vert ^{n}}{w\left\vert \cdot\right\vert ^{2\theta}%
}+\ldots\\
& =\sum\limits_{j=0}^{2\theta-n-1}\tfrac{1}{j!}\sum\limits_{k=0}%
^{2\theta-n-1-j}\sum\limits_{\left\vert \alpha\right\vert =k}\int%
\limits_{\left\vert \cdot\right\vert \geq r_{3}}\tfrac{\left\vert \widehat
{\xi}^{\alpha}\right\vert \left(  \left(  x_{+}+a_{+}\right)  \widehat{\xi
}_{+}\right)  ^{j}\left\vert \widehat{a}\xi\right\vert ^{n}}{w\left\vert
\cdot\right\vert ^{2\theta}}+\ldots\\
& =\sum\limits_{j=0}^{2\theta-n-1}\tfrac{\left(  x_{+}+a_{+}\right)  ^{j}}%
{j!}\sum\limits_{k=0}^{2\theta-n-1-j}\sum\limits_{\left\vert \alpha\right\vert
=k}\int\limits_{\left\vert \cdot\right\vert \geq r_{3}}\tfrac{\left\vert
\widehat{\xi}^{\alpha}\right\vert \left(  \left(  x_{+}+a_{+}\right)
^{\wedge}\widehat{\xi}_{+}\right)  ^{j}\left\vert \widehat{a}\widehat{\xi
}\right\vert ^{n}\left\vert \cdot\right\vert ^{n}}{w\left\vert \cdot
\right\vert ^{2\theta}}+\ldots\\
& \leq\sum\limits_{j=0}^{2\theta-n-1}\tfrac{\left(  x_{+}+a_{+}\right)  ^{j}%
}{j!}%
{\textstyle\sum\limits_{k=0}^{2\theta-n-1-j}}
{\textstyle\sum\limits_{\left\vert \alpha\right\vert =k}}
\min\left\{  \int\limits_{\left\vert \cdot\right\vert \geq r_{3}}%
\tfrac{\left\vert \widehat{\xi}^{\alpha}\right\vert \left\vert \cdot
\right\vert ^{n}}{w\left\vert \cdot\right\vert ^{2\theta}},\int%
\limits_{\left\vert \cdot\right\vert \geq r_{3}}\tfrac{\left(  \left(
x_{+}+a_{+}\right)  ^{\wedge}\widehat{\xi}_{+}\right)  ^{j}\left\vert
\cdot\right\vert ^{n}}{w\left\vert \cdot\right\vert ^{2\theta}},\int%
\limits_{\left\vert \cdot\right\vert \geq r_{3}}\tfrac{\left\vert \widehat
{a}\widehat{\xi}\right\vert ^{n}\left\vert \cdot\right\vert ^{n}}{w\left\vert
\cdot\right\vert ^{2\theta}}\right\}  +\ldots\\
& =\ldots+\left(  \int\limits_{\left\vert \cdot\right\vert \leq r_{3}}%
\tfrac{\left\vert \widehat{a}\widehat{\xi}\right\vert ^{n}}{w}+\int%
\limits_{\left\vert \cdot\right\vert \geq r_{3}}\tfrac{\left\vert \widehat
{a}\xi\right\vert ^{n}}{w\left\vert \cdot\right\vert ^{2\theta}}\right)
+\sum\limits_{\left\vert \alpha\right\vert =2\theta-n}\upsilon_{n,\alpha
}^{\left(  w\right)  }\sum\limits_{\sigma\leq\alpha}\frac{\left(  x_{+}%
+a_{+}\right)  ^{\sigma}}{\sigma!}.
\end{align*}

From \ref{a1.08} and then part 3 of Theorem \ref{Thm_Ap_1},%
\begin{align*}
& \sum\limits_{\left\vert \alpha\right\vert =2\theta-n}\upsilon_{n,\alpha
}^{\left(  w\right)  }\sum\limits_{\sigma\leq\alpha}\frac{\left(  x_{+}%
+a_{+}\right)  ^{\sigma}}{\sigma!}\\
& =\sum\limits_{\left\vert \alpha\right\vert =2\theta-n}\sum\limits_{\sigma
\leq\alpha}\frac{\left(  x_{+}+a_{+}\right)  ^{\sigma}}{\sigma!}%
\int\limits_{\left\vert \cdot\right\vert \leq r_{3}}\frac{\widehat{\xi}%
_{+}^{\alpha}\left\vert \widehat{a}\widehat{\xi}\right\vert ^{n}}{w}\\
& =\int\limits_{\left\vert \cdot\right\vert \leq r_{3}}\sum\limits_{\left\vert
\alpha\right\vert =2\theta-n}\left\vert \widehat{\xi}^{\alpha}\right\vert
\sum\limits_{\sigma\leq\alpha}\frac{\left(  x_{+}+a_{+}\right)  ^{\sigma}%
}{\sigma!}\frac{\left\vert \widehat{a}\widehat{\xi}\right\vert ^{n}}{w}\\
& =\int\limits_{\left\vert \cdot\right\vert \leq r_{3}}\sum\limits_{j=0}%
^{2\theta-n}h_{2\theta-n-j}\left(  \widehat{\xi}_{+}\right)  \frac{\left(
\left(  x_{+}+a_{+}\right)  \widehat{\xi}_{+}\right)  ^{j}}{j!}\frac
{\left\vert \widehat{a}\widehat{\xi}\right\vert ^{n}}{w}\\
& =\int\limits_{\left\vert \cdot\right\vert \leq r_{3}}\sum\limits_{j=0}%
^{2\theta-n}\sum\limits_{\left\vert \beta\right\vert =2\theta-n-j}\left\vert
\widehat{\xi}^{\beta}\right\vert \frac{\left(  \left(  x_{+}+a_{+}\right)
\widehat{\xi}_{+}\right)  ^{j}}{j!}\frac{\left\vert \widehat{a}\widehat{\xi
}\right\vert ^{n}}{w}\\
& =\sum\limits_{j=0}^{2\theta-n}\frac{\left\vert x_{+}+a_{+}\right\vert ^{j}%
}{j!}\sum\limits_{\left\vert \beta\right\vert =2\theta-n-j}\int%
\limits_{\left\vert \cdot\right\vert \leq r_{3}}\frac{\left\vert \widehat{\xi
}^{\beta}\right\vert \left(  \left(  x_{+}+a_{+}\right)  ^{\wedge}\widehat
{\xi}_{+}\right)  ^{j}\left\vert \widehat{a}\widehat{\xi}\right\vert ^{n}}%
{w}\\
& =\sum\limits_{j=0}^{2\theta-n}\frac{\left\vert x_{+}+a_{+}\right\vert ^{j}%
}{j!}\sum\limits_{\left\vert \beta\right\vert =2\theta-n-j}\min\left\{
\int\limits_{\left\vert \cdot\right\vert \leq r_{3}}\frac{\left\vert
\widehat{\xi}^{\beta}\right\vert }{w},\int\limits_{\left\vert \cdot\right\vert
\leq r_{3}}\frac{\left(  \left(  x_{+}+a_{+}\right)  ^{\wedge}\widehat{\xi
}_{+}\right)  ^{j}}{w},\int\limits_{\left\vert \cdot\right\vert \leq r_{3}%
}\frac{\left\vert \widehat{a}\widehat{\xi}\right\vert ^{n}}{w}\right\}  .
\end{align*}

Thus%
\begin{align}
& \sum\limits_{\left\vert \alpha\right\vert \leq2\theta-n}\upsilon_{n,\alpha
}^{\left(  w\right)  }\sum\limits_{\sigma\leq\alpha}\frac{\left(  x_{+}%
+a_{+}\right)  ^{\sigma}}{\sigma!}\nonumber\\
& \leq\int\limits_{\left\vert \cdot\right\vert \leq r_{3}}\frac{\left\vert
\widehat{a}\widehat{\xi}\right\vert ^{n}}{w}+\int\limits_{\left\vert
\cdot\right\vert \geq r_{3}}\frac{\left\vert \widehat{a}\xi\right\vert ^{n}%
}{w\left\vert \cdot\right\vert ^{2\theta}}+\nonumber\\
& +\sum\limits_{j=0}^{2\theta-n-1}\frac{\left\vert x_{+}+a_{+}\right\vert
^{j}}{j!}\sum\limits_{k=0}^{\substack{2\theta-n- \\-1-j}}\sum
\limits_{\left\vert \alpha\right\vert =k}\min\left\{  \int\limits_{\left\vert
\cdot\right\vert \geq r_{3}}\frac{\left\vert \widehat{\xi}^{\alpha}\right\vert
\left\vert \cdot\right\vert ^{n}}{w\left\vert \cdot\right\vert ^{2\theta}%
},\int\limits_{\left\vert \cdot\right\vert \geq r_{3}}\frac{\left(  \left(
x_{+}+a_{+}\right)  ^{\wedge}\widehat{\xi}_{+}\right)  ^{j}\left\vert
\cdot\right\vert ^{n}}{w\left\vert \cdot\right\vert ^{2\theta}},\int%
\limits_{\left\vert \cdot\right\vert \geq r_{3}}\frac{\left\vert \widehat
{a}\widehat{\xi}\right\vert ^{n}\left\vert \cdot\right\vert ^{n}}{w\left\vert
\cdot\right\vert ^{2\theta}}\right\}  +\nonumber\\
& +\sum\limits_{j=0}^{2\theta-n}\frac{\left\vert x_{+}+a_{+}\right\vert ^{j}%
}{j!}\sum\limits_{\left\vert \beta\right\vert =2\theta-n-j}\min\left\{
\int\limits_{\left\vert \cdot\right\vert \leq r_{3}}\frac{\left\vert
\widehat{\xi}^{\beta}\right\vert }{w},\int\limits_{\left\vert \cdot\right\vert
\leq r_{3}}\frac{\left(  \left(  x_{+}+a_{+}\right)  ^{\wedge}\widehat{\xi
}_{+}\right)  ^{j}}{w},\int\limits_{\left\vert \cdot\right\vert \leq r_{3}%
}\frac{\left\vert \widehat{a}\widehat{\xi}\right\vert ^{n}}{w}\right\}
\nonumber\\
& =\sum\limits_{j=0}^{2\theta-n}K_{j,n,r_{3}}^{\left(  w\right)  }%
\frac{\left\vert x_{+}+a_{+}\right\vert ^{j}}{j!},\label{a113}%
\end{align}

and so \ref{a112} becomes,%
\begin{multline*}
\sum\limits_{\left\vert \alpha\right\vert \leq2\theta-n}\upsilon_{n,\alpha
}^{\left(  w\right)  }\int_{0}^{1}\int\left\vert \frac{\left(  \xi
+x+sa\right)  ^{\alpha}}{\alpha!}\widehat{\phi_{\widehat{a}}}\left(
\xi\right)  \right\vert d\xi\text{ }g_{n}\left(  s\right)  ds\\
\leq\frac{1}{n+1}\left(  \sum\limits_{j=0}^{2\theta-n}K_{j,n,r_{3}}^{\left(
w\right)  }\frac{\left\vert x_{+}+a_{+}\right\vert ^{j}}{j!}\right)
\sum\limits_{j=0}^{2\theta-n}\int\frac{\left(  \xi_{+}\mathbf{1}\right)  ^{j}%
}{j!}\left\vert \widehat{a}D\widehat{\phi_{0}}\left(  \xi\right)  \right\vert
,
\end{multline*}

where%
\begin{align}
K_{0,n,r_{3}}^{\left(  w\right)  }  & =\int\limits_{\left\vert \cdot
\right\vert \leq r_{3}}\frac{\left\vert \widehat{a}\widehat{\xi}\right\vert
^{n}}{w}+\int\limits_{\left\vert \cdot\right\vert \geq r_{3}}\frac{\left\vert
\widehat{a}\xi\right\vert ^{n}}{w\left\vert \cdot\right\vert ^{2\theta}%
}+\nonumber\\
& +\sum\limits_{k=0}^{2\theta-n-1}\sum\limits_{\left\vert \alpha\right\vert
=k}\min\left\{  \int\limits_{\left\vert \cdot\right\vert \geq r_{3}}%
\frac{\left\vert \widehat{\xi}^{\alpha}\right\vert \left\vert \cdot\right\vert
^{n}}{w\left\vert \cdot\right\vert ^{2\theta}},\int\limits_{\left\vert
\cdot\right\vert \geq r_{3}}\frac{\left\vert \cdot\right\vert ^{n}%
}{w\left\vert \cdot\right\vert ^{2\theta}},\int\limits_{\left\vert
\cdot\right\vert \geq r_{3}}\frac{\left\vert \widehat{a}\widehat{\xi
}\right\vert ^{n}\left\vert \cdot\right\vert ^{n}}{w\left\vert \cdot
\right\vert ^{2\theta}}\right\}  +\nonumber\\
& +\sum\limits_{\left\vert \beta\right\vert =2\theta-n}\min\left\{
\int\limits_{\left\vert \cdot\right\vert \leq r_{3}}\frac{\left\vert
\widehat{\xi}^{\beta}\right\vert }{w},\int\limits_{\left\vert \cdot\right\vert
\leq r_{3}}\frac{1}{w},\int\limits_{\left\vert \cdot\right\vert \leq r_{3}%
}\frac{\left\vert \widehat{a}\widehat{\xi}\right\vert ^{n}}{w}\right\}
\nonumber\\
& =\int\limits_{\left\vert \cdot\right\vert \leq r_{3}}\frac{\left\vert
\widehat{a}\widehat{\xi}\right\vert ^{n}}{w}+\int\limits_{\left\vert
\cdot\right\vert \geq r_{3}}\frac{\left\vert \widehat{a}\xi\right\vert ^{n}%
}{w\left\vert \cdot\right\vert ^{2\theta}}+\sum\limits_{k=0}^{2\theta-n-1}%
\sum\limits_{\left\vert \alpha\right\vert =k}\min\left\{  \int%
\limits_{\left\vert \cdot\right\vert \geq r_{3}}\frac{\left\vert \widehat{\xi
}^{\alpha}\right\vert \left\vert \cdot\right\vert ^{n}}{w\left\vert
\cdot\right\vert ^{2\theta}},\int\limits_{\left\vert \cdot\right\vert \geq
r_{3}}\frac{\left\vert \widehat{a}\widehat{\xi}\right\vert ^{n}\left\vert
\cdot\right\vert ^{n}}{w\left\vert \cdot\right\vert ^{2\theta}}\right\}
+\nonumber\\
& +\sum\limits_{\left\vert \beta\right\vert =2\theta-n}\min\left\{
\int\limits_{\left\vert \cdot\right\vert \leq r_{3}}\frac{\left\vert
\widehat{\xi}^{\beta}\right\vert }{w},\int\limits_{\left\vert \cdot\right\vert
\leq r_{3}}\frac{\left\vert \widehat{a}\widehat{\xi}\right\vert ^{n}}%
{w}\right\}  ,\label{a1.13}%
\end{align}

and%
\begin{gather}
K_{j,n,r_{3}}^{\left(  w\right)  }=\sum\limits_{k=0}^{2\theta-n-1-j}%
\sum\limits_{\left\vert \alpha\right\vert =k}\min\left\{  \int%
\limits_{\left\vert \cdot\right\vert \geq r_{3}}\frac{\left\vert \widehat{\xi
}^{\alpha}\right\vert \left\vert \cdot\right\vert ^{n}}{w\left\vert
\cdot\right\vert ^{2\theta}},\int\limits_{\left\vert \cdot\right\vert \geq
r_{3}}\frac{\left(  \left(  x_{+}+a\right)  ^{\wedge}\widehat{\xi}_{+}\right)
^{j}\left\vert \cdot\right\vert ^{n}}{w\left\vert \cdot\right\vert ^{2\theta}%
},\int\limits_{\left\vert \cdot\right\vert \geq r_{3}}\frac{\left\vert
\widehat{a}\widehat{\xi}\right\vert ^{n}\left\vert \cdot\right\vert ^{n}%
}{w\left\vert \cdot\right\vert ^{2\theta}}\right\}  +\nonumber\\
\qquad+\sum\limits_{\left\vert \beta\right\vert =2\theta-n-j}\min\left\{
\int\limits_{\left\vert \cdot\right\vert \leq r_{3}}\frac{\left\vert
\widehat{\xi}^{\beta}\right\vert }{w},\int\limits_{\left\vert \cdot\right\vert
\leq r_{3}}\frac{\left(  \left(  x_{+}+a_{+}\right)  ^{\wedge}\widehat{\xi
}_{+}\right)  ^{j}}{w},\int\limits_{\left\vert \cdot\right\vert \leq r_{3}%
}\frac{\left\vert \widehat{a}\widehat{\xi}\right\vert ^{n}}{w}\right\}
,\nonumber\\
for\text{ }1\leq j\leq2\theta-n-1,\label{a1.14}%
\end{gather}

and%
\begin{align}
K_{2\theta-n,n,r_{3}}^{\left(  w\right)  }  & =\min\left\{  \int%
\limits_{\left\vert \cdot\right\vert \leq r_{3}}\frac{1}{w},\int%
\limits_{\left\vert \cdot\right\vert \leq r_{3}}\frac{\left(  \left(
x_{+}+a_{+}\right)  ^{\wedge}\widehat{\xi}_{+}\right)  ^{2\theta-n}}{w}%
,\int\limits_{\left\vert \cdot\right\vert \leq r_{3}}\frac{\left\vert
\widehat{a}\widehat{\xi}\right\vert ^{n}}{w}\right\} \nonumber\\
& =\min\left\{  \int\limits_{\left\vert \cdot\right\vert \leq r_{3}}%
\frac{\left(  \left(  x_{+}+a_{+}\right)  ^{\wedge}\widehat{\xi}_{+}\right)
^{2\theta-n}}{w},\int\limits_{\left\vert \cdot\right\vert \leq r_{3}}%
\frac{\left\vert \widehat{a}\widehat{\xi}\right\vert ^{n}}{w}\right\}
.\label{a1.17}%
\end{align}

Before proceeding we note for future use that we have actually proved but ??
not used yet? ??%
\begin{align}
&  \sum\limits_{\left\vert \alpha\right\vert \leq m}\int_{0}^{1}\int\left\vert
\frac{\left(  \xi+x+sa\right)  ^{\alpha}}{\alpha!}\widehat{\phi_{\widehat{a}}%
}\left(  \xi\right)  \right\vert d\xi\text{ }g_{n}\left(  s\right)
ds\nonumber\\
&  \leq\frac{1}{n+1}\left(  \sum\limits_{j=0}^{m}K_{w;j,n,r_{3}}%
\frac{\left\vert x_{+}+a_{+}\right\vert ^{j}}{j!}\right)  \sum_{k=1}^{d}%
\sum\limits_{j=0}^{m}\frac{1}{j!}\left\Vert \left\vert \cdot\right\vert
_{1}^{j}D_{k}\widehat{\phi_{0}}\right\Vert _{1},\quad m,n\geq0.\label{a512}%
\end{align}

Thus \ref{a005} becomes%
\begin{align}
&  \left\vert \frac{\sqrt{2\pi}}{n!}\left(  \overline{\widehat{g_{n}}}\left(
a\xi\right)  \left(  ia\xi\right)  \phi_{0}\widehat{\left(  aD\right)
^{n}G_{\rho}}\right)  ^{\vee}\left(  x\right)  \right\vert \nonumber\\
&  =\frac{C_{\rho;2\theta-n,r_{3}}}{\left(  2\pi\right)  ^{d}}\frac{\left\vert
a\right\vert ^{n+1}}{n!}\sum\limits_{\left\vert \alpha\right\vert \leq
2\theta-n}\upsilon_{n,\alpha}^{\left(  w\right)  }\int_{0}^{1}\int\left\vert
\frac{\left(  \xi+x+sa\right)  ^{\alpha}}{\alpha!}\widehat{\phi_{\widehat{a}}%
}\left(  \xi\right)  \right\vert d\xi\text{ }g_{n}\left(  s\right)
ds\nonumber\\
&  \leq\frac{C_{\rho;2\theta-n,r_{3}}}{\left(  2\pi\right)  ^{d}}%
\frac{\left\vert a\right\vert ^{n+1}}{\left(  n+1\right)  !}\left(
\sum\limits_{j=0}^{2\theta-n}K_{j,n,r_{3}}^{\left(  w\right)  }\frac
{\left\vert x_{+}+a_{+}\right\vert ^{j}}{j!}\right)  \sum\limits_{j=0}%
^{2\theta-n}\int\frac{\left(  \xi_{+}\mathbf{1}\right)  ^{j}}{j!}\left\vert
\widehat{a}D\widehat{\phi_{0}}\left(  \xi\right)  \right\vert ,\label{a1.18}\\
&
\begin{array}
[c]{c}%
when\text{ }n\leq\left\lfloor 2\kappa\right\rfloor <2\theta.
\end{array}
\nonumber
\end{align}
\medskip

\fbox{\textbf{Case 2} $\left\lfloor 2\kappa\right\rfloor \geq2\theta$,
$n\leq\left\lfloor 2\kappa\right\rfloor $} By \ref{a000},%
\[
\left\vert \left(  \widehat{a}D\right)  ^{n}G_{\rho}\left(  x\right)
\right\vert \leq\left(  2\pi\right)  ^{-\frac{d}{2}}\int\frac{\left\vert
\widehat{a}\eta\right\vert ^{n}}{w\left\vert \cdot\right\vert ^{2\theta}}%
d\eta,\quad x\in\mathbb{R}^{d},
\]

so using \ref{a1.22}, \ref{a801} becomes%
\begin{align}
&  \left\vert \frac{\sqrt{2\pi}}{n!}\left(  \overline{\widehat{g_{n}}}\left(
a\xi\right)  \left(  ia\xi\right)  \phi_{0}\widehat{\left(  aD\right)
^{n}G_{\rho}}\right)  ^{\vee}\left(  x\right)  \right\vert \nonumber\\
&  \leq\left(  2\pi\right)  ^{-\frac{d}{2}}\frac{\left\vert a\right\vert
^{n+1}}{n!}\int_{0}^{1}\int\left\vert \widehat{\phi_{\widehat{a}}}\left(
\xi\right)  \right\vert \left\vert \left(  \left(  \widehat{a}D\right)
^{n}G_{\rho}\right)  \left(  \xi+x+sa\right)  \right\vert d\xi\text{ }%
g_{n}\left(  s\right)  ds\nonumber\\
&  =\left(  2\pi\right)  ^{-d}\frac{\left\vert a\right\vert ^{n+1}}{n!}\left(
\int\left\vert \widehat{a}D\widehat{\phi_{0}}\right\vert \right)  \left(
\int\frac{\left\vert \widehat{a}\eta\right\vert ^{n}}{w\left\vert
\cdot\right\vert ^{2\theta}}d\eta\right)  \int_{0}^{1}g_{n}\left(  s\right)
ds\nonumber\\
&  =\left(  2\pi\right)  ^{-d}\frac{\left\vert a\right\vert ^{n+1}}{\left(
n+1\right)  !}\left(  \int\frac{\left\vert \widehat{a}\eta\right\vert ^{n}%
}{w\left\vert \cdot\right\vert ^{2\theta}}\right)  \int\left\vert \widehat
{a}D\widehat{\phi_{0}}\right\vert ,\label{a811}\\
&
\begin{array}
[c]{c}%
when\text{ }2\theta\leq n\leq\left\lfloor 2\kappa\right\rfloor
\end{array}
,\nonumber
\end{align}

Now combining cases \ref{a811} and \ref{a1.18} we get for $n\leq\left\lfloor
2\kappa\right\rfloor $,%
\begin{align}
& \left\vert \frac{\sqrt{2\pi}}{n!}\left(  \overline{\widehat{g_{n}}}\left(
a\xi\right)  \left(  ia\xi\right)  \phi_{0}\widehat{\left(  aD\right)
^{n}G_{\rho}}\right)  ^{\vee}\left(  x\right)  \right\vert \nonumber\\
& \leq\frac{1}{\left(  2\pi\right)  ^{d}}\frac{\left\vert a\right\vert ^{n+1}%
}{\left(  n+1\right)  !}\times\nonumber\\
& \times\left\{
\begin{array}
[c]{ll}%
\left(
\begin{array}
[c]{l}%
\underline{C}_{\rho;2\theta-n,0,r_{3}}\int\limits_{\left\vert \cdot\right\vert
\leq r_{3}}\frac{\left\vert \widehat{a}\widehat{\xi}\right\vert ^{n}}{w}+\\
\quad+\overline{C}_{2\theta-n,r_{3}}^{\left(  \rho\right)  }\int%
\limits_{\left\vert \cdot\right\vert \geq r_{3}}\frac{\left\vert \widehat
{a}\xi\right\vert ^{n}}{w\left\vert \cdot\right\vert ^{2\theta}}%
\end{array}
\right)  \left(  \sum\limits_{j=0}^{2\theta-n}K_{j,n,r_{3}}^{\left(  w\right)
}\frac{\left\vert x_{+}+a_{+}\right\vert ^{j}}{j!}\right)  \sum\limits_{j=0}%
^{2\theta-n}\int\frac{\left(  \xi_{+}\mathbf{1}\right)  ^{j}}{j!}\left\vert
\widehat{a}D\widehat{\phi_{0}}\right\vert , & n<2\theta,\\
\medskip & \\
\left(  \int\frac{\left\vert \widehat{a}\xi\right\vert ^{n}}{w\left\vert
\cdot\right\vert ^{2\theta}}\right)  \int\left\vert \widehat{a}D\widehat
{\phi_{0}}\right\vert , & n\geq2\theta,
\end{array}
\right. \label{a86}%
\end{align}

where $\overline{C}_{2\theta-n,r_{3}}^{\left(  \rho\right)  }$ is given by
\ref{a2.03}, $\underline{C}_{\rho;n,,r_{3}}$ is given by \ref{a75} and the
$\left\{  K_{j,n,r_{3}}^{\left(  w\right)  }\right\}  _{j=0}^{2\theta-n}$ are
given by \ref{a1.13}, \ref{a1.14} and \ref{a1.17}.\medskip

\underline{\textbf{Estimate the }$G_{F}$\textbf{\ term on RHS of }\ref{p28}%
}\label{SbSect_estim_Gf_term_Taylor_rem}\medskip

Compare Section \ref{Sect_data_fn_Taylor_W3.2} for data functions.

The goal here is to obtain the estimate the $G_{F}$ term of \ref{p28}. Now
write%
\[
\frac{\sqrt{2\pi}}{n!}\left(  \left(  ia\xi\right)  ^{n+1}\overline
{\widehat{g_{n}}}\left(  a\xi\right)  \phi_{\infty}G_{F}\right)  ^{\vee}%
=\frac{\sqrt{2\pi}}{n!}\left(  \overline{\widehat{g_{n}}}\left(  a\xi\right)
\phi_{\infty}\frac{\left(  ia\xi\right)  ^{n+1}}{w\left\vert \cdot\right\vert
^{2\theta}}\right)  ^{\vee},
\]

and observe that%
\begin{align}
& \left\Vert \overline{\widehat{g_{n}}}\left(  a\xi\right)  \phi_{\infty}%
\frac{\left(  ia\xi\right)  ^{n+1}}{w\left\vert \cdot\right\vert ^{2\theta}%
}\right\Vert _{1}\nonumber\\
& =\int\left\vert \overline{\widehat{g_{n}}}\left(  a\xi\right)  \phi_{\infty
}\frac{\left(  ia\xi\right)  ^{n+1}}{w\left\vert \cdot\right\vert ^{2\theta}%
}\right\vert d\xi\nonumber\\
& =\int\limits_{\left\vert \cdot\right\vert \leq r_{3}}\left\vert
\overline{\widehat{g_{n}}}\left(  a\xi\right)  \phi_{\infty}\frac{\left(
ia\xi\right)  ^{n+1}}{w\left\vert \cdot\right\vert ^{2\theta}}\right\vert
d\xi+\int\limits_{\left\vert \cdot\right\vert \geq r_{3}}\left\vert
\overline{\widehat{g_{n}}}\left(  a\xi\right)  \phi_{\infty}\frac{\left(
ia\xi\right)  ^{n+1}}{w\left\vert \cdot\right\vert ^{2\theta}}\right\vert
d\xi\nonumber\\
& =\int\limits_{\left\vert \cdot\right\vert \leq r_{3}}\left\vert
\widehat{g_{n}}\left(  a\xi\right)  \right\vert \phi_{\infty}\frac{\left\vert
a\xi\right\vert ^{n+1}}{w\left\vert \cdot\right\vert ^{2\theta}}d\xi
+\int\limits_{\left\vert \cdot\right\vert \geq r_{3}}\left\vert \widehat
{g_{n}}\left(  a\xi\right)  \right\vert \phi_{\infty}\frac{\left\vert
a\xi\right\vert ^{n+1}}{w\left\vert \cdot\right\vert ^{2\theta}}%
d\xi\nonumber\\
& =\left\vert a\right\vert ^{n+1}\left(  \int\limits_{\left\vert
\cdot\right\vert \leq r_{3}}\left\vert \widehat{g_{n}}\left(  a\xi\right)
\right\vert \phi_{\infty}\frac{\left\vert \widehat{a}\xi\right\vert ^{n+1}%
}{w\left\vert \cdot\right\vert ^{2\theta}}d\xi+\int\limits_{\left\vert
\cdot\right\vert \geq r_{3}}\left\vert \widehat{g_{n}}\left(  a\xi\right)
\right\vert \phi_{\infty}\frac{\left\vert \widehat{a}\xi\right\vert ^{n+1}%
}{w\left\vert \cdot\right\vert ^{2\theta}}d\xi\right)  .\label{a1.01}%
\end{align}

Now to estimate the two terms on the right side of \ref{a1.01}.\smallskip

\fbox{Show $\int\nolimits_{\left\vert \cdot\right\vert \leq r_{3}}\left\vert
\widehat{g_{n}}\left(  a\xi\right)  \right\vert \phi_{\infty}\frac{\left\vert
\widehat{a}\xi\right\vert ^{n+1}}{w\left\vert \cdot\right\vert ^{2\theta}}%
d\xi<\infty$} From part 2 of Lemma \ref{Lem_gm_properties_2},
\[
\left\vert \widehat{g_{n}}\left(  t\right)  \right\vert \leq\min\left\{
\frac{1}{\sqrt{2\pi}}\frac{1}{n+1},\frac{1}{\sqrt{2\pi}}\left(  1+\frac
{1}{n+1}\right)  \frac{1}{1+\left\vert t\right\vert }\right\}  ,\quad
t\in\mathbb{R}^{1},
\]

and from \ref{p85},%
\[
\phi_{\infty}\left(  x\right)  \leq\frac{\left\vert x\right\vert ^{2\theta}%
}{\left(  2\theta\right)  !}\left\Vert \left(  \widehat{\cdot}D\right)
^{2\theta}\phi_{\infty}\right\Vert _{\infty;\leq r},\text{\quad}x\in
\overline{B}_{r}.
\]

so
\[
\left\Vert \frac{\phi_{\infty}}{\left\vert \cdot\right\vert ^{2\theta}%
}\right\Vert _{\infty;B_{r_{3}}}\leq\frac{1}{\left(  2\theta\right)
!}\left\Vert \left(  \widehat{\cdot}D\right)  ^{2\theta}\phi_{\infty
}\right\Vert _{\infty;\leq r_{3}}%
\]

and%
\begin{align}
\int\limits_{\left\vert \cdot\right\vert \leq r_{3}}\left\vert \overline
{\widehat{g_{n}}}\left(  a\xi\right)  \phi_{\infty}\frac{\left(  i\widehat
{a}\xi\right)  ^{n+1}}{w\left\vert \cdot\right\vert ^{2\theta}}\right\vert
d\xi & \leq\int\limits_{\left\vert \cdot\right\vert \leq r_{3}}\left\vert
\widehat{g_{n}}\left(  a\xi\right)  \right\vert \frac{\phi_{\infty}%
}{\left\vert \cdot\right\vert ^{2\theta}}\frac{\left\vert \widehat{a}%
\xi\right\vert ^{n+1}}{w}\nonumber\\
& \leq\frac{1}{\sqrt{2\pi}}\frac{1}{n+1}\frac{\left\Vert \left(
\widehat{\cdot}D\right)  ^{2\theta}\phi_{\infty}\right\Vert _{\infty;\leq
r_{3}}}{\left(  2\theta\right)  !}\int\limits_{\left\vert \cdot\right\vert
\leq r_{3}}\frac{\left\vert \widehat{a}\xi\right\vert ^{n+1}}{w}%
\label{a2.12}\\
& <\infty.\nonumber
\end{align}
\medskip

\fbox{Show $\int\nolimits_{\left\vert \cdot\right\vert \geq r_{3}}\left\vert
\widehat{g_{n}}\left(  a\xi\right)  \right\vert \phi_{\infty}\frac{\left\vert
a\xi\right\vert ^{n+1}}{w\left\vert \cdot\right\vert ^{2\theta}}d\xi<\infty$%
.}\medskip

\qquad\fbox{If $n\leq\left\lfloor 2\kappa\right\rfloor -1$} then
$n+1\leq\left\lfloor 2\kappa\right\rfloor \leq2\kappa$ and so%
\begin{align*}
\int\limits_{\left\vert \cdot\right\vert \geq r_{3}}\left\vert \widehat{g_{n}%
}\left(  a\xi\right)  \right\vert \phi_{\infty}\frac{\left\vert a\xi
\right\vert ^{n+1}}{w\left\vert \cdot\right\vert ^{2\theta}}d\xi & \leq
\frac{1}{\sqrt{2\pi}}\frac{\left\Vert \phi_{\infty}\right\Vert _{\infty;\geq
r_{3}}}{n+1}\int\limits_{\left\vert \cdot\right\vert \geq r_{3}}%
\frac{\left\vert a\xi\right\vert ^{n+1}}{w\left\vert \cdot\right\vert
^{2\theta}}d\xi\\
& =\frac{\left\vert a\right\vert ^{n+1}}{\sqrt{2\pi}}\frac{\left\Vert
\phi_{\infty}\right\Vert _{\infty;\geq r_{3}}}{n+1}\int\limits_{\left\vert
\cdot\right\vert \geq r_{3}}\frac{\left\vert \widehat{a}\xi\right\vert ^{n+1}%
}{w\left\vert \cdot\right\vert ^{2\theta}}d\xi\\
& <\infty.
\end{align*}

\qquad\fbox{If $n=\left\lfloor 2\kappa\right\rfloor $} then $n+1=\left\lceil
2\kappa\right\rceil $ and%
\begin{align*}
\int\limits_{\left\vert \cdot\right\vert \geq r_{3}}\left\vert \widehat{g_{n}%
}\left(  a\xi\right)  \right\vert \phi_{\infty}\frac{\left\vert a\xi
\right\vert ^{n+1}}{w\left\vert \cdot\right\vert ^{2\theta}}d\xi &
=\int\limits_{\left\vert \cdot\right\vert \geq r_{3}}\left\vert \widehat
{g_{\left\lfloor 2\kappa\right\rfloor }}\left(  a\xi\right)  \right\vert
\phi_{\infty}\frac{\left\vert a\xi\right\vert ^{\left\lceil 2\kappa
\right\rceil }}{w\left\vert \cdot\right\vert ^{2\theta}}d\xi\\
& =\int\limits_{\left\vert \cdot\right\vert \geq r_{3}}\left\vert
\widehat{g_{\left\lfloor 2\kappa\right\rfloor }}\left(  a\xi\right)
\right\vert \phi_{\infty}\frac{\left\vert a\xi\right\vert ^{\left\lceil
2\kappa\right\rceil }}{w\left\vert \cdot\right\vert ^{2\theta}}d\xi\\
& =\int\limits_{\left\vert \cdot\right\vert \geq r_{3}}\left\vert
\widehat{g_{\left\lfloor 2\kappa\right\rfloor }}\left(  a\xi\right)
\right\vert \left\vert a\xi\right\vert ^{\left\lceil 2\kappa\right\rceil
-2\kappa}\phi_{\infty}\frac{\left\vert a\xi\right\vert ^{2\kappa}}{w\left\vert
\cdot\right\vert ^{2\theta}}d\xi\\
& \leq\int\limits_{\left\vert \cdot\right\vert \geq r_{3}}\frac{1}{\sqrt{2\pi
}}\left(  1+\frac{1}{\left\lfloor 2\kappa\right\rfloor +1}\right)
\frac{\left\vert a\xi\right\vert ^{\left\lceil 2\kappa\right\rceil -2\kappa}%
}{1+\left\vert a\xi\right\vert }\phi_{\infty}\frac{\left\vert a\xi\right\vert
^{2\kappa}}{w\left\vert \cdot\right\vert ^{2\theta}}d\xi\\
& =\frac{1}{\sqrt{2\pi}}\left(  1+\frac{1}{\left\lfloor 2\kappa\right\rfloor
+1}\right)  \int\limits_{\left\vert \cdot\right\vert \geq r_{3}}%
\frac{\left\vert a\xi\right\vert ^{\left\lceil 2\kappa\right\rceil -2\kappa}%
}{1+\left\vert a\xi\right\vert }\phi_{\infty}\frac{\left\vert a\xi\right\vert
^{2\kappa}}{w\left\vert \cdot\right\vert ^{2\theta}}d\xi\\
& \leq\frac{1}{\sqrt{2\pi}}\left(  1+\frac{1}{\left\lceil 2\kappa\right\rceil
}\right)  \left\Vert \phi_{\infty}\right\Vert _{\infty;\geq r_{3}}%
\int\limits_{\left\vert \cdot\right\vert \geq r_{3}}\frac{\left\vert
a\xi\right\vert ^{2\kappa}}{w\left\vert \cdot\right\vert ^{2\theta}}d\xi\\
& =\frac{\left\vert a\right\vert ^{2\kappa}}{\sqrt{2\pi}}\left(  1+\frac
{1}{\left\lceil 2\kappa\right\rceil }\right)  \left\Vert \phi_{\infty
}\right\Vert _{\infty;\geq r_{3}}\int\limits_{\left\vert \cdot\right\vert \geq
r_{3}}\frac{\left\vert \widehat{a}\xi\right\vert ^{2\kappa}}{w\left\vert
\cdot\right\vert ^{2\theta}}d\xi\\
& <\infty.
\end{align*}

Thus%
\begin{equation}
\int\limits_{\left\vert \cdot\right\vert \geq r_{3}}\left\vert \widehat{g_{n}%
}\left(  a\xi\right)  \right\vert \phi_{\infty}\frac{\left\vert a\xi
\right\vert ^{n+1}}{w\left\vert \cdot\right\vert ^{2\theta}}d\xi\leq\left\{
\begin{array}
[c]{ll}%
\frac{\left\vert a\right\vert ^{n+1}}{\sqrt{2\pi}}\frac{\left\Vert
\phi_{\infty}\right\Vert _{\infty;\geq r_{3}}}{n+1}\int\limits_{\left\vert
\cdot\right\vert \geq r_{3}}\frac{\left\vert \widehat{a}\xi\right\vert ^{n+1}%
}{w\left\vert \cdot\right\vert ^{2\theta}}d\xi, & n<\left\lfloor
2\kappa\right\rfloor ,\\
\frac{\left\vert a\right\vert ^{2\kappa}}{\sqrt{2\pi}}\left(  1+\frac
{1}{\left\lceil 2\kappa\right\rceil }\right)  \left\Vert \phi_{\infty
}\right\Vert _{\infty;\geq r_{3}}\int\limits_{\left\vert \cdot\right\vert \geq
r_{3}}\frac{\left\vert \widehat{a}\xi\right\vert ^{2\kappa}}{w\left\vert
\cdot\right\vert ^{2\theta}}d\xi, & n=\left\lfloor 2\kappa\right\rfloor ,
\end{array}
\right. \label{a2.00}%
\end{equation}

Substituting the estimates \ref{a2.12} and \ref{a2.00} into \ref{a1.01} yields%
\begin{align*}
&  \left\Vert \frac{\sqrt{2\pi}}{n!}\overline{\widehat{g_{n}}}\left(
a\xi\right)  \phi_{\infty}\frac{\left(  ia\xi\right)  ^{n+1}}{w\left\vert
\cdot\right\vert ^{2\theta}}\right\Vert _{1}\\
&  \leq\frac{\sqrt{2\pi}}{n!}\left\vert a\right\vert ^{n+1}\int%
\limits_{\left\vert \cdot\right\vert \leq r_{3}}\left\vert \widehat{g_{n}%
}\left(  a\xi\right)  \right\vert \phi_{\infty}\frac{\left\vert \widehat{a}%
\xi\right\vert ^{n+1}}{w\left\vert \cdot\right\vert ^{2\theta}}d\xi
+\frac{\sqrt{2\pi}}{n!}\int\limits_{\left\vert \cdot\right\vert \geq r_{3}%
}\left\vert \widehat{g_{n}}\left(  a\xi\right)  \right\vert \phi_{\infty}%
\frac{\left\vert a\xi\right\vert ^{n+1}}{w\left\vert \cdot\right\vert
^{2\theta}}d\xi\\
&  \leq\frac{\sqrt{2\pi}}{n!}\frac{1}{\sqrt{2\pi}}\frac{\left\vert
a\right\vert ^{n+1}}{n+1}\frac{\left\Vert \left(  \widehat{\cdot}D\right)
^{2\theta}\phi_{\infty}\right\Vert _{\infty;\leq r_{3}}}{\left(
2\theta\right)  !}\int\limits_{\left\vert \cdot\right\vert \leq r_{3}}%
\frac{\left\vert \widehat{a}\xi\right\vert ^{n+1}}{w}+\\
&  +\frac{\sqrt{2\pi}}{n!}\left\{
\begin{array}
[c]{ll}%
\frac{\left\vert a\right\vert ^{n+1}}{\sqrt{2\pi}}\frac{\left\Vert
\phi_{\infty}\right\Vert _{\infty;\geq r_{3}}}{n+1}\int\limits_{\left\vert
\cdot\right\vert \geq r_{3}}\frac{\left\vert \widehat{a}\xi\right\vert ^{n+1}%
}{w\left\vert \cdot\right\vert ^{2\theta}}d\xi, & n<\left\lfloor
2\kappa\right\rfloor ,\\
\frac{\left\vert a\right\vert ^{2\kappa}}{\sqrt{2\pi}}\left(  1+\frac
{1}{\left\lceil 2\kappa\right\rceil }\right)  \left\Vert \phi_{\infty
}\right\Vert _{\infty;\geq r_{3}}\int\limits_{\left\vert \cdot\right\vert \geq
r_{3}}\frac{\left\vert \widehat{a}\xi\right\vert ^{2\kappa}}{w\left\vert
\cdot\right\vert ^{2\theta}}d\xi, & n=\left\lfloor 2\kappa\right\rfloor ,
\end{array}
\right. \\
&  =\frac{\left\vert a\right\vert ^{n+1}}{\left(  n+1\right)  !}%
\frac{\left\Vert \left(  \widehat{\cdot}D\right)  ^{2\theta}\phi_{\infty
}\right\Vert _{\infty;\leq r_{3}}}{\left(  2\theta\right)  !}\int%
\limits_{\left\vert \cdot\right\vert \leq r_{3}}\frac{\left\vert \widehat
{a}\xi\right\vert ^{n+1}}{w}+\\
&  \qquad+\left\{
\begin{array}
[c]{ll}%
\frac{\left\vert a\right\vert ^{n+1}}{\left(  n+1\right)  !}\left\Vert
\phi_{\infty}\right\Vert _{\infty;\geq r_{3}}\int\limits_{\left\vert
\cdot\right\vert \geq r_{3}}\frac{\left\vert \widehat{a}\xi\right\vert ^{n+1}%
}{w\left\vert \cdot\right\vert ^{2\theta}}d\xi, & n<\left\lfloor
2\kappa\right\rfloor ,\\
\left\vert a\right\vert ^{2\kappa}\left(  \frac{1}{\left\lfloor 2\kappa
\right\rfloor !}+\frac{1}{\left\lceil 2\kappa\right\rceil !}\right)
\left\Vert \phi_{\infty}\right\Vert _{\infty;\geq r_{3}}\int%
\limits_{\left\vert \cdot\right\vert \geq r_{3}}\frac{\left\vert \widehat
{a}\xi\right\vert ^{2\kappa}}{w\left\vert \cdot\right\vert ^{2\theta}}d\xi, &
n=\left\lfloor 2\kappa\right\rfloor .
\end{array}
\right.
\end{align*}

Thus $\frac{\sqrt{2\pi}}{n!}\overline{\widehat{g_{n}}}\left(  a\xi\right)
\phi_{\infty}\frac{\left(  ia\xi\right)  ^{n+1}}{w\left\vert \cdot\right\vert
^{2\theta}}\in L^{1}$ means%
\[
\left\Vert \frac{\sqrt{2\pi}}{n!}\left(  \left(  ia\xi\right)  ^{n+1}%
\overline{\widehat{g_{n}}}\left(  a\xi\right)  \phi_{\infty}G_{F}\right)
^{\vee}\right\Vert _{\infty}\leq\frac{1}{\left(  2\pi\right)  ^{d/2}%
}\left\Vert \frac{\sqrt{2\pi}}{n!}\overline{\widehat{g_{n}}}\left(
a\xi\right)  \phi_{\infty}\frac{\left(  ia\xi\right)  ^{n+1}}{w\left\vert
\cdot\right\vert ^{2\theta}}\right\Vert _{1},
\]

and so%
\begin{align}
& \left\Vert \frac{\sqrt{2\pi}}{n!}\left(  \left(  ia\xi\right)
^{n+1}\overline{\widehat{g_{n}}}\left(  a\xi\right)  \phi_{\infty}%
G_{F}\right)  ^{\vee}\right\Vert _{\infty}\nonumber\\
& \leq\frac{1}{\left(  2\pi\right)  ^{\frac{d}{2}}}\frac{\left\vert
a\right\vert ^{n+1}}{\left(  n+1\right)  !}\frac{\left\Vert \left(
\widehat{\cdot}D\right)  ^{2\theta}\phi_{\infty}\right\Vert _{\infty;\leq
r_{3}}}{\left(  2\theta\right)  !}\int\limits_{\left\vert \cdot\right\vert
\leq r_{3}}\frac{\left\vert \widehat{a}\xi\right\vert ^{n+1}}{w}+\nonumber\\
& +\frac{1}{\left(  2\pi\right)  ^{\frac{d}{2}}}\left\{
\begin{array}
[c]{ll}%
\frac{\left\vert a\right\vert ^{n+1}}{\left(  n+1\right)  !}\left\Vert
\phi_{\infty}\right\Vert _{\infty;\geq r_{3}}\int\limits_{\left\vert
\cdot\right\vert \geq r_{3}}\frac{\left\vert \widehat{a}\xi\right\vert ^{n+1}%
}{w\left\vert \cdot\right\vert ^{2\theta}}d\xi, & n<\left\lfloor
2\kappa\right\rfloor ,\\
\left\vert a\right\vert ^{2\kappa}\left(  \frac{1}{\left\lfloor 2\kappa
\right\rfloor !}+\frac{1}{\left\lceil 2\kappa\right\rceil !}\right)
\left\Vert \phi_{\infty}\right\Vert _{\infty;\geq r_{3}}\int%
\limits_{\left\vert \cdot\right\vert \geq r_{3}}\frac{\left\vert \widehat
{a}\xi\right\vert ^{2\kappa}}{w\left\vert \cdot\right\vert ^{2\theta}}d\xi, &
n=\left\lfloor 2\kappa\right\rfloor .
\end{array}
\right. \label{a1.02}%
\end{align}

Substituting the estimates \ref{a1.02} for $G_{F}$ and \ref{a86} for $G_{\rho
}$ into \ref{p28} gives for $n\leq\left\lfloor 2\kappa\right\rfloor $,%
\begin{align}
&  \left\vert \frac{\sqrt{2\pi}}{n!}\left(  \left(  ia\xi\right)
^{n+1}\overline{\widehat{g_{n}}}\left(  a\xi\right)  \widehat{G_{\rho}%
}\right)  ^{\vee}\left(  x\right)  \right\vert \nonumber\\
&  \leq\left\vert \frac{\sqrt{2\pi}}{n!}\left(  \overline{\widehat{g_{n}}%
}\left(  a\xi\right)  \left(  ia\xi\right)  \phi_{0}\widehat{\left(
aD\right)  ^{n}G_{\rho}}\right)  ^{\vee}\left(  x\right)  \right\vert
+\left\vert \frac{\sqrt{2\pi}}{n!}\left(  \left(  ia\xi\right)  ^{n+1}%
\overline{\widehat{g_{n}}}\left(  a\xi\right)  \phi_{\infty}G_{F}\right)
^{\vee}\left(  x\right)  \right\vert \nonumber\\
&  \leq\left\{
\begin{array}
[c]{l}%
\frac{1}{\left(  2\pi\right)  ^{d}}\frac{\left\vert a\right\vert ^{n+1}%
}{\left(  n+1\right)  !}\times\\
\smallskip\\
\times\left\{
\begin{array}
[c]{l}%
\begin{array}
[c]{r}%
\left(
\begin{array}
[c]{l}%
\underline{C}_{\rho;2\theta-n,0,r_{3}}\int\limits_{\left\vert \cdot\right\vert
\leq r_{3}}\frac{\left\vert \widehat{a}\widehat{\xi}\right\vert ^{n}}{w}+\\
\smallskip\\
\quad+\overline{C}_{2\theta-n,r_{3}}^{\left(  \rho\right)  }\int%
\limits_{\left\vert \cdot\right\vert \geq r_{3}}\frac{\left\vert \widehat
{a}\xi\right\vert ^{n}}{w\left\vert \cdot\right\vert ^{2\theta}}%
\end{array}
\right)  \left(  \sum\limits_{j=0}^{2\theta-n}K_{j,n,r_{3}}^{\left(  w\right)
}\frac{\left\vert x_{+}+a_{+}\right\vert ^{j}}{j!}\right)  \sum\limits_{j=0}%
^{2\theta-n}\int\frac{\left(  \xi_{+}\mathbf{1}\right)  ^{j}}{j!}\left\vert
\widehat{a}D\widehat{\phi_{0}}\right\vert \\
if\text{ \ }\left\lfloor 2\kappa\right\rfloor <2\theta,
\end{array}
\\
\left(  \int\frac{\left\vert \widehat{a}\xi\right\vert ^{n}}{w\left\vert
\cdot\right\vert ^{2\theta}}\right)  \int\left\vert \widehat{a}D\widehat
{\phi_{0}}\right\vert ,\qquad\qquad\qquad\qquad\qquad\qquad\qquad\qquad\qquad
if\text{ \ }\left\lfloor 2\kappa\right\rfloor \geq2\theta,
\end{array}
\right\}  +\\
\medskip\\
+\frac{1}{\left(  2\pi\right)  ^{\frac{d}{2}}}\frac{\left\vert a\right\vert
^{n+1}}{\left(  n+1\right)  !}\frac{\left\Vert \left(  \widehat{\cdot
}D\right)  ^{2\theta}\phi_{\infty}\right\Vert _{\infty;\leq r_{3}}}{\left(
2\theta\right)  !}\int\limits_{\left\vert \cdot\right\vert \leq r_{3}}%
\frac{\left\vert \widehat{a}\xi\right\vert ^{n+1}}{w}+\\
\medskip\\
+\frac{1}{\left(  2\pi\right)  ^{\frac{d}{2}}}\left\{
\begin{array}
[c]{ll}%
\frac{\left\vert a\right\vert ^{n+1}}{\left(  n+1\right)  !}\left\Vert
\phi_{\infty}\right\Vert _{\infty;\geq r_{3}}\int\limits_{\left\vert
\cdot\right\vert \geq r_{3}}\frac{\left\vert \widehat{a}\xi\right\vert ^{n+1}%
}{w\left\vert \cdot\right\vert ^{2\theta}}d\xi, & n<\left\lfloor
2\kappa\right\rfloor ,\\
\smallskip & \\
\left\vert a\right\vert ^{2\kappa}\left(  \frac{1}{\left\lfloor 2\kappa
\right\rfloor !}+\frac{1}{\left\lceil 2\kappa\right\rceil !}\right)
\left\Vert \phi_{\infty}\right\Vert _{\infty;\geq r_{3}}\int%
\limits_{\left\vert \cdot\right\vert \geq r_{3}}\frac{\left\vert \widehat
{a}\xi\right\vert ^{2\kappa}}{w\left\vert \cdot\right\vert ^{2\theta}}d\xi, &
n=\left\lfloor 2\kappa\right\rfloor ,
\end{array}
\right\}
\end{array}
\right. \label{a1.04}%
\end{align}

where $\overline{C}_{2\theta-n,r_{3}}^{\left(  \rho\right)  }$ is given by
\ref{a2.03}, $\underline{C}_{n,,r_{3}}^{\left(  \rho\right)  }$ is given by
\ref{a75} and the $\left\{  K_{j,n,r_{3}}^{\left(  w\right)  }\right\}
_{j=0}^{2\theta-n}$ are given by \ref{a1.13}, \ref{a1.14} and \ref{a1.17}.

\begin{remark}
\label{Rem_Tay_expan_data_fns_W3.2}Concerning the estimates \ref{a1.04}:

\begin{enumerate}
\item When $n<\left\lfloor 2\kappa\right\rfloor $ the order of convergence is
$n+1\leq\left\lfloor 2\kappa\right\rfloor $ but that when $n=\left\lfloor
2\kappa\right\rfloor $ the order of convergence is (at least) $2\kappa$. The
$2\kappa$ convergence occurs in the estimate for $G_{F}$ when $n=\left\lfloor
2\kappa\right\rfloor $.

\item The only dependencies on $x$ are the terms $\frac{\left\vert x_{+}%
+a_{+}\right\vert ^{j}}{j!}$.

\item The dependency on $\rho\in S_{1,2\theta}$ is contained in $\underline
{C}_{2\theta-n,r_{3}}^{\left(  \rho\right)  }$ and $\overline{C}%
_{2\theta-n,r_{3}}^{\left(  \rho\right)  }$. The constant $\underline
{C}_{2\theta-n,r_{3}}^{\left(  \rho\right)  }$ can be made zero by choosing
$\operatorname*{supp}\rho\subset\overline{B}_{r_{3}\text{.}}$.

\item Expressions involving $\widehat{a}$ e.g. $\left\vert \widehat{a}%
\xi\right\vert ^{n+1}$ are retained because they allow more precise estimates
when the weight function is radial.

\item We can choose $\phi_{0}=\rho$.
\end{enumerate}
\end{remark}

\subsection{Radial weight functions and the estimate \ref{a1.04}}

Radial functions can be applied in three situations: to the function $\phi
_{0}$, to the the function $\rho$ and to the weight function $w$.\medskip

\fbox{\textbf{Part 1} $\phi_{0},\phi_{\infty}$ are \textbf{radial}.}\smallskip

Regarding the term $\sum\limits_{j=0}^{2\theta-n}\int\frac{\left(
\mathbf{1}\xi_{+}\right)  ^{j}}{j!}\left\vert \widehat{a}D\widehat{\phi_{0}%
}\left(  \xi\right)  \right\vert d\xi$: $\widehat{\phi_{0}}$ must be radial,
say $\widehat{\phi_{0}}\left(  \xi\right)  =\psi_{F}\left(  \left\vert
\xi\right\vert \right)  $. Part 2 of Lemma \ref{Lem_op_aD_estim} states that%
\[
\left(  \widehat{a}D\right)  \left(  f\left(  \left\vert x\right\vert \right)
\right)  =f^{\prime}\left(  \left\vert x\right\vert \right)  \left(
\widehat{a}D\right)  \left\vert x\right\vert =\widehat{a}\widehat{x}f^{\prime
}\left(  \left\vert x\right\vert \right)  ,
\]

so that%
\[
\widehat{a}D\widehat{\phi_{0}}\left(  \xi\right)  =\widehat{a}D\left(
\psi_{F}\left(  \left\vert \xi\right\vert \right)  \right)  =\widehat
{a}\widehat{\xi}\psi_{F}^{\prime}\left(  \left\vert \xi\right\vert \right)  ,
\]

and%
\begin{align*}
\sum\limits_{j=0}^{2\theta-n}\int\frac{\left(  \mathbf{1}\xi_{+}\right)  ^{j}%
}{j!}\left\vert \widehat{a}D\widehat{\phi_{0}}\left(  \xi\right)  \right\vert
d\xi & =\sum\limits_{j=0}^{2\theta-n}\int\frac{\left(  \mathbf{1}\xi
_{+}\right)  ^{j}}{j!}\left\vert \left(  \widehat{a}\widehat{\xi}\right)
\psi_{F}^{\prime}\left(  \left\vert \xi\right\vert \right)  \right\vert d\xi\\
& =\sum\limits_{j=0}^{2\theta-n}\frac{1}{j!}\int\left(  \mathbf{1}\widehat
{\xi}_{+}\right)  ^{j}\left\vert \widehat{a}\widehat{\xi}\right\vert
\left\vert \xi\right\vert ^{j}\left\vert \psi_{F}^{\prime}\left(  \left\vert
\xi\right\vert \right)  \right\vert d\xi\\
& \leq\sum\limits_{j=0}^{2\theta-n}\frac{1}{j!}\int\left(  \left\vert
\mathbf{1}\right\vert \left\vert \widehat{\xi}_{+}\right\vert \right)
^{j}\left\vert \widehat{a}\widehat{\xi}\right\vert \left\vert \xi\right\vert
^{j}\left\vert \psi_{F}^{\prime}\left(  \left\vert \xi\right\vert \right)
\right\vert d\xi\\
& =\sum\limits_{j=0}^{2\theta-n}\frac{d^{j/2}}{j!}\int\left\vert \widehat
{a}\widehat{\xi}\right\vert \left\vert \xi\right\vert ^{j}\left\vert \psi
_{F}^{\prime}\left(  \left\vert \xi\right\vert \right)  \right\vert d\xi.
\end{align*}

From \ref{Ap141},%
\[
\int\limits_{\left\vert \xi\right\vert \leq r}\left\vert \widehat{a}%
\xi\right\vert ^{p}f\left(  \left\vert \xi\right\vert \right)  d\xi
=\frac{B\left(  \frac{d}{2},\frac{p+1}{2}\right)  }{B\left(  \frac{d+p}%
{2},\frac{1}{2}\right)  }\int\limits_{\left\vert \xi\right\vert \leq
r}\left\vert \xi\right\vert ^{p}f\left(  \left\vert \xi\right\vert \right)
d\xi,
\]

which implies%
\begin{align*}
\int\left\vert \widehat{a}\widehat{\xi}\right\vert \left\vert \xi\right\vert
^{j}\left\vert \psi_{F}^{\prime}\left(  \left\vert \xi\right\vert \right)
\right\vert d\xi & =\frac{B\left(  \frac{d}{2},1\right)  }{B\left(  \frac
{d+1}{2},\frac{1}{2}\right)  }\int\left\vert \xi\right\vert ^{1}\left\vert
\xi\right\vert ^{j}\left\vert \psi_{F}^{\prime}\left(  \left\vert
\xi\right\vert \right)  \right\vert d\xi\\
& =\frac{B\left(  \frac{d}{2},1\right)  }{B\left(  \frac{d+1}{2},\frac{1}%
{2}\right)  }\int\left\vert \xi\right\vert ^{j+1}\left\vert \psi_{F}^{\prime
}\left(  \left\vert \xi\right\vert \right)  \right\vert d\xi,
\end{align*}

and so%
\begin{align*}
\sum\limits_{j=0}^{2\theta-n}\int\frac{\left(  \mathbf{1}\xi_{+}\right)  ^{j}%
}{j!}\left\vert \widehat{a}D\widehat{\phi_{0}}\left(  \xi\right)  \right\vert
d\xi & \leq\sum\limits_{j=0}^{2\theta-n}\frac{d^{j/2}}{j!}\int\left\vert
\widehat{a}\widehat{\xi}\right\vert \left\vert \xi\right\vert ^{j}\left\vert
\psi_{F}^{\prime}\left(  \left\vert \xi\right\vert \right)  \right\vert d\xi\\
& =\sum\limits_{j=0}^{2\theta-n}\frac{d^{j/2}}{j!}\frac{B\left(  \frac{d}%
{2},1\right)  }{B\left(  \frac{d+1}{2},\frac{1}{2}\right)  }\int\left\vert
\xi\right\vert ^{j+1}\left\vert \psi_{F}^{\prime}\left(  \left\vert
\xi\right\vert \right)  \right\vert d\xi\\
& =\frac{B\left(  \frac{d}{2},1\right)  }{B\left(  \frac{d+1}{2},\frac{1}%
{2}\right)  }\sum\limits_{j=0}^{2\theta-n}\frac{d^{j/2}}{j!}\int\left\vert
\xi\right\vert ^{j+1}\left\vert \psi_{F}^{\prime}\left(  \left\vert
\xi\right\vert \right)  \right\vert d\xi\\
& =\frac{B\left(  \frac{d}{2},1\right)  }{B\left(  \frac{d+1}{2},\frac{1}%
{2}\right)  }\omega_{d}\sum\limits_{j=0}^{2\theta-n}\frac{d^{j/2}}{j!}\int%
_{0}^{\infty}s^{j+d}\left\vert \psi_{F}^{\prime}\left(  s\right)  \right\vert
ds.
\end{align*}

But $\omega_{d}=\frac{2\pi^{d/2}}{\Gamma\left(  d/2\right)  }$ so%
\[
\frac{B\left(  \frac{d}{2},1\right)  }{B\left(  \frac{d+1}{2},\frac{1}%
{2}\right)  }\omega_{d}=\frac{\Gamma\left(  \frac{d}{2}\right)  }%
{\Gamma\left(  \frac{d+1}{2}\right)  \pi^{\frac{1}{2}}}\frac{2\pi^{\frac{d}%
{2}}}{\Gamma\left(  \frac{d}{2}\right)  }=\frac{1}{\pi}\frac{2\pi^{\frac
{d+1}{2}}}{\Gamma\left(  \frac{d+1}{2}\right)  }=\frac{\omega_{d+1}}{\pi},
\]

and as a consequence%
\begin{equation}
\sum\limits_{j=0}^{2\theta-n}\int\frac{\left(  \mathbf{1}\xi_{+}\right)  ^{j}%
}{j!}\left\vert \widehat{a}D\widehat{\phi_{0}}\left(  \xi\right)  \right\vert
d\xi\leq\frac{\omega_{d+1}}{\pi}\sum\limits_{j=0}^{2\theta-n}\frac{d^{j/2}%
}{j!}\int_{0}^{\infty}s^{j+d}\left\vert \psi_{F}^{\prime}\left(  s\right)
\right\vert ds.\label{a1.23}%
\end{equation}

Now to consider the expression $\left\Vert \left(  \widehat{\cdot}D\right)
^{2\theta}\phi_{\infty}\right\Vert _{\infty;\leq r_{3}}$. Since $\phi_{0}%
+\phi_{\infty}=1$,%
\[
\left\Vert \left(  \widehat{\cdot}D\right)  ^{2\theta}\phi_{\infty}\right\Vert
_{\infty;\leq r_{3}}=\left\Vert \left(  \widehat{\cdot}D\right)  ^{2\theta
}\phi_{0}\right\Vert _{\infty;\leq r_{3}},
\]

and suppose $\phi_{0}\left(  x\right)  =\left(  \phi_{0}\right)  _{\circ
}\left(  \left\vert x\right\vert \right)  $. Then%
\begin{equation}
\left\Vert \left(  \widehat{\cdot}D\right)  ^{2\theta}\phi_{\infty}\right\Vert
_{\infty;\leq r_{3}}=\left\Vert D^{2\theta}\left(  \phi_{0}\right)  _{\circ
}\right\Vert _{\infty;\left[  0,r_{3}\right]  }.\label{a2.01}%
\end{equation}
\medskip

\fbox{\textbf{Part 2} $\rho$ is radial}\smallskip

In \ref{a1.04} the dependency on $\rho\in S_{\emptyset,2\theta}$ is
encapsulated in the "constants" $\overline{C}_{2\theta-n,r_{3}}^{\left(
\rho\right)  }$ and $\underline{C}_{n,,r_{3}}^{\left(  \rho\right)  }$ where
$\overline{C}_{2\theta-n,r_{3}}^{\left(  \rho\right)  }$ is given by
\ref{a2.03}, $\underline{C}_{n,r_{3}}^{\left(  \rho\right)  }$ is given by
\ref{a75}. From \ref{a2.03},%
\[
\overline{C}_{m,r}^{\left(  \rho\right)  }=\max\left\{  1,\max_{j<m}\left\Vert
\left\vert \cdot\right\vert ^{j}\rho\left(  \cdot\right)  \right\Vert
_{\infty;\geq r}\right\}  ,\quad m\geq1,
\]

so simply%
\[
\overline{C}_{m,r}^{\left(  \rho\right)  }=\max\left\{  1,\max_{j<m}\left\Vert
t^{j}\rho_{\circ}\left(  t\right)  \right\Vert _{\infty;\left[  r,\infty
\right)  }\right\}  ,\quad m\geq1.
\]

From \ref{a75},%
\[
\underline{C}_{n,r}^{\left(  \rho\right)  }=\frac{1}{n!}\left\Vert \left(
\widehat{\cdot}D\right)  ^{n}\rho\right\Vert _{\infty;\leq r},\quad n\geq0.
\]

From part 1 of Lemma \ref{Lem_deriv_rad_funcs}, if $\phi\left(  x\right)
=\phi_{\circ}\left(  \left\vert x\right\vert \right)  $ then%
\[
\left(  \left(  \widehat{\cdot}D\right)  ^{n}\phi\right)  \left(  x\right)
=\left(  D^{n}\phi_{\circ}\right)  \left(  \left\vert x\right\vert \right)  ,
\]

so that in this case%
\[
\underline{C}_{n,r}^{\left(  \rho\right)  }=\frac{1}{n!}\left\Vert D^{n}%
\rho_{\circ}\right\Vert _{\infty;\left[  0,r\right]  }.
\]
\medskip

\fbox{\textbf{Part 3} $w$ is radial}\smallskip

In \ref{a1.04} the dependency on the weight function $w$ is encapsulated in
the "constants" $\left\{  K_{j,n,r_{3}}^{\left(  w\right)  }\right\}
_{j=0}^{2\theta-n}$ given by \ref{a1.13}, \ref{a1.14} and \ref{a1.17}. From
\ref{Ap017} and \ref{Ap123}, ?? CHECK! ??
\begin{align*}
\int_{\left\vert x\right\vert \leq r}\left\vert \widehat{\xi}\widehat
{x}\right\vert ^{n}f\left(  \left\vert x\right\vert \right)  dx  &
=\frac{B\left(  \frac{d}{2},\frac{n+1}{2}\right)  }{B\left(  \frac{1}{2}%
,\frac{n+d}{2}\right)  }\omega_{d}\int_{0}^{r}t^{d-1}f\left(  t\right)  dt,\\
& and\\
\int\limits_{\left\vert x\right\vert \geq r}\left\vert \widehat{x}^{\alpha
}\right\vert f\left(  \left\vert x\right\vert \right)  dx  & =2^{d+1}B\left(
\frac{\alpha+1}{2}\right)  \frac{1}{\omega_{d}}\int\limits_{\left\vert
x\right\vert \geq r,x\geq\mathbf{0}}f\left(  \left\vert x\right\vert \right)
dx\\
& =2B\left(  \frac{\alpha+1}{2}\right)  \frac{1}{\omega_{d}}\int%
\limits_{\left\vert x\right\vert \geq r}f\left(  \left\vert x\right\vert
\right)  dx\\
& =2B\left(  \frac{\alpha+1}{2}\right)  \int_{r}^{\infty}t^{d-1}f\left(
t\right)  dt.
\end{align*}

\begin{align*}
&  K_{0,n,r_{3}}^{\left(  w\right)  }\\
&  =\int\limits_{\left\vert \cdot\right\vert \leq r_{3}}\frac{\left\vert
\widehat{a}\widehat{\xi}\right\vert ^{n}}{w}+\int\limits_{\left\vert
\cdot\right\vert \geq r_{3}}\frac{\left\vert \widehat{a}\widehat{\xi
}\right\vert ^{n}\left\vert \cdot\right\vert ^{n}}{w\left\vert \cdot
\right\vert ^{2\theta}}+\sum\limits_{k=0}^{2\theta-n-1}\sum\limits_{\left\vert
\alpha\right\vert =k}\min\left\{  \int\limits_{\left\vert \cdot\right\vert
\geq r_{3}}\frac{\left\vert \widehat{\xi}^{\alpha}\right\vert \left\vert
\cdot\right\vert ^{n}}{w\left\vert \cdot\right\vert ^{2\theta}},\int%
\limits_{\left\vert \cdot\right\vert \geq r_{3}}\frac{\left\vert \widehat
{a}\widehat{\xi}\right\vert ^{n}\left\vert \cdot\right\vert ^{n}}{w\left\vert
\cdot\right\vert ^{2\theta}}\right\}  +\\
&  \qquad+\sum\limits_{\left\vert \beta\right\vert =2\theta-n}\min\left\{
\int\limits_{\left\vert \cdot\right\vert \leq r_{3}}\frac{\left\vert
\widehat{\xi}^{\beta}\right\vert }{w},\int\limits_{\left\vert \cdot\right\vert
\leq r_{3}}\frac{\left\vert \widehat{a}\widehat{\xi}\right\vert ^{n}}%
{w}\right\} \\
&  =\frac{B\left(  \frac{d}{2},\frac{n+1}{2}\right)  }{B\left(  \frac{1}%
{2},\frac{n+d}{2}\right)  }\omega_{d}\int_{0}^{r_{3}}\frac{t^{d-1}dt}%
{w_{\circ}\left(  t\right)  }+\frac{B\left(  \frac{d}{2},\frac{n+1}{2}\right)
}{B\left(  \frac{1}{2},\frac{n+d}{2}\right)  }\omega_{d}\int_{r_{3}}^{\infty
}\frac{t^{n-2\theta+d-1}}{w_{\circ}\left(  t\right)  }dt+\\
&  \qquad+\sum\limits_{k=0}^{2\theta-n-1}\sum\limits_{\left\vert
\alpha\right\vert =k}\min\left\{  \ldots\right\}  +\sum\limits_{\left\vert
\beta\right\vert =2\theta-n}\min\left\{  \ldots\right\} \\
&  =\frac{B\left(  \frac{d}{2},\frac{n+1}{2}\right)  }{B\left(  \frac{1}%
{2},\frac{n+d}{2}\right)  }\omega_{d}\left(  \int_{0}^{r_{3}}\frac{t^{d-1}%
dt}{w_{\circ}\left(  t\right)  }+\int_{r_{3}}^{\infty}\frac{t^{n-2\theta+d-1}%
}{w_{\circ}\left(  t\right)  }dt\right)  +\ldots\\
&  =\ldots+\sum\limits_{k=0}^{2\theta-n-1}\sum\limits_{\left\vert
\alpha\right\vert =k}\min\left\{  \int\limits_{\left\vert \cdot\right\vert
\geq r_{3}}\frac{\left\vert \widehat{\xi}^{\alpha}\right\vert \left\vert
\cdot\right\vert ^{n}}{w\left\vert \cdot\right\vert ^{2\theta}},\int%
\limits_{\left\vert \cdot\right\vert \geq r_{3}}\frac{\left\vert \widehat
{a}\widehat{\xi}\right\vert ^{n}\left\vert \cdot\right\vert ^{n}}{w\left\vert
\cdot\right\vert ^{2\theta}}\right\}  +\\
&  \qquad\qquad+\sum\limits_{\left\vert \beta\right\vert =2\theta-n}%
\min\left\{  \int\limits_{\left\vert \cdot\right\vert \leq r_{3}}%
\frac{\left\vert \widehat{\xi}^{\beta}\right\vert }{w},\int\limits_{\left\vert
\cdot\right\vert \leq r_{3}}\frac{\left\vert \widehat{a}\widehat{\xi
}\right\vert ^{n}}{w}\right\} \\
&  =\ldots+\sum\limits_{k=0}^{2\theta-n-1}\sum\limits_{\left\vert
\alpha\right\vert =k}\min\left\{  2B\left(  \frac{\alpha+1}{2}\right)
\int_{r_{3}}^{\infty}\frac{t^{n-2\theta+d-1}}{w_{\circ}\left(  t\right)
}dt,\frac{B\left(  \frac{d}{2},\frac{n+1}{2}\right)  }{B\left(  \frac{1}%
{2},\frac{n+d}{2}\right)  }\omega_{d}\int_{r_{3}}^{\infty}\frac{t^{n-2\theta
+d-1}}{w_{\circ}\left(  t\right)  }dt\right\}  +\\
&  \qquad+\sum\limits_{\left\vert \beta\right\vert =2\theta-n}\min\left\{
2B\left(  \frac{\beta+1}{2}\right)  \int_{0}^{r_{3}}\frac{t^{n+d-1}}{w_{\circ
}\left(  t\right)  }dt,\frac{B\left(  \frac{d}{2},\frac{n+1}{2}\right)
}{B\left(  \frac{1}{2},\frac{n+d}{2}\right)  }\omega_{d}\int_{0}^{r_{3}}%
\frac{t^{n+d-1}}{w_{\circ}\left(  t\right)  }dt\right\} \\
&  =\frac{B\left(  \frac{d}{2},\frac{n+1}{2}\right)  }{B\left(  \frac{1}%
{2},\frac{n+d}{2}\right)  }\omega_{d}\left(  \int_{0}^{r_{3}}\frac{t^{d-1}%
dt}{w_{\circ}\left(  t\right)  }+\int_{r_{3}}^{\infty}\frac{t^{n-2\theta+d-1}%
}{w_{\circ}\left(  t\right)  }dt\right)  +\\
&  \qquad+\sum\limits_{k=0}^{2\theta-n-1}\sum\limits_{\left\vert
\alpha\right\vert =k}\min\left\{  2B\left(  \frac{\alpha+1}{2}\right)
,\frac{B\left(  \frac{d}{2},\frac{n+1}{2}\right)  }{B\left(  \frac{1}{2}%
,\frac{n+d}{2}\right)  }\omega_{d}\right\}  \int_{r_{3}}^{\infty}%
\frac{t^{n-2\theta+d-1}}{w_{\circ}\left(  t\right)  }dt+\\
&  \qquad+\sum\limits_{\left\vert \beta\right\vert =2\theta-n}\min\left\{
2B\left(  \frac{\beta+1}{2}\right)  ,\frac{B\left(  \frac{d}{2},\frac{n+1}%
{2}\right)  }{B\left(  \frac{1}{2},\frac{n+d}{2}\right)  }\omega_{d}\right\}
\int_{0}^{r_{3}}\frac{t^{n+d-1}}{w_{\circ}\left(  t\right)  }dt.
\end{align*}

ETC. ??? Must consider $K_{j,n,r_{3}}^{\left(  w\right)  }$. Perhaps leave for now.

\begin{remark}
\textbf{Estimating the expression}%
\[
\int\limits_{\left\vert \xi\right\vert \leq r_{3}}\left\vert \overline
{\widehat{g_{\left\lfloor 2\kappa\right\rfloor }}}\left(  a\xi\right)
\right\vert \left(  \left\vert a\xi\right\vert \right)  ^{\left\lceil
2\kappa\right\rceil -2\theta}\frac{d\xi}{w\left(  \xi\right)  },
\]

when $w$ is radial.

From Theorem \ref{Thm_Integ_u(xy)f(|x|)dx},%
\[
\int_{\left\vert \xi\right\vert \leq r}u\left(  \widehat{a}\xi\right)
f\left(  \left\vert \xi\right\vert \right)  d\xi=\int_{\left\vert
\xi\right\vert \leq r}u\left(  \xi_{k}\right)  f\left(  \left\vert
\xi\right\vert \right)  d\xi,\quad1\leq k\leq d,
\]

so%
\[
\int_{\left\vert \xi\right\vert \leq r}u\left(  a\xi\right)  f\left(
\left\vert \xi\right\vert \right)  d\xi=\int_{\left\vert \xi\right\vert \leq
r}u\left(  \left\vert a\right\vert \xi_{k}\right)  f\left(  \left\vert
\xi\right\vert \right)  d\xi,\quad1\leq k\leq d.
\]

and%
\[
\int\limits_{\left\vert x\right\vert \leq r}u\left(  \left\vert a\right\vert
\xi_{k}\right)  f\left(  \left\vert \xi\right\vert \right)  d\xi=\omega
_{d-1}\int_{0}^{r}\left(  \int_{0}^{\pi}u\left(  \left\vert a\right\vert
\rho\cos t\right)  \sin^{d-2}tdt\right)  f\left(  \rho\right)  \rho^{d-1}%
d\rho,
\]

\[
\int_{0}^{\pi}u\left(  \left\vert a\right\vert \rho\cos t\right)  \sin
^{d-2}tdt=\int_{0}^{\pi}\left\vert \widehat{g_{\left\lfloor 2\kappa
\right\rfloor }}\left(  \left\vert a\right\vert \rho\cos t\right)  \right\vert
\left(  \left\vert a\right\vert \rho\left\vert \cos t\right\vert \right)
^{\left\lceil 2\kappa\right\rceil -2\theta}\sin^{d-2}tdt,
\]

Set $\tau=\left\vert a\right\vert \rho\cos t$ so that $\sin t=\left(
1-\frac{\tau^{2}}{\left(  \left\vert a\right\vert \rho\right)  ^{2}}\right)
^{1/2}$ and $dt=\frac{1}{\left\vert a\right\vert \rho}\left(  1-\frac{\tau
^{2}}{\left(  \left\vert a\right\vert \rho\right)  ^{2}}\right)  ^{-1/2}d\tau
$:%
\begin{align*}
& \int_{0}^{\pi}u\left(  \left\vert a\right\vert \rho\cos t\right)  \sin
^{d-2}tdt\\
& =\frac{1}{\left\vert a\right\vert \rho}\int_{-\left\vert a\right\vert \rho
}^{\left\vert a\right\vert \rho}\left\vert \widehat{g_{\left\lfloor
2\kappa\right\rfloor }}\left(  \tau\right)  \right\vert \left\vert
\tau\right\vert ^{\left\lceil 2\kappa\right\rceil -2\theta}\left(
\sqrt{1-\frac{\tau^{2}}{\left(  \left\vert a\right\vert \rho\right)  ^{2}}%
}\right)  ^{d-3}d\tau\\
& \leq\frac{1}{\left\vert a\right\vert \rho}\left(  \int_{-\left\vert
a\right\vert \rho}^{\left\vert a\right\vert \rho}\left\vert \widehat
{g_{\left\lfloor 2\kappa\right\rfloor }}\left(  \tau\right)  \right\vert
^{2}\right)  ^{1/2}\left(  \int_{-\left\vert a\right\vert \rho}^{\left\vert
a\right\vert \rho}\left\vert \tau\right\vert ^{2\left(  \left\lceil
2\kappa\right\rceil -2\theta\right)  }\left(  1-\frac{\tau^{2}}{\left(
\left\vert a\right\vert \rho\right)  ^{2}}\right)  ^{d-3}d\tau\right)
^{1/2}\\
& =??
\end{align*}

\end{remark}

\section{Taylor expansion of basis functions when $w\in W3.1$%
\label{SbSect_Tay_basis_W3.1}}

Assume the weight function $w$ has property $W3.1$ for order $\theta$ and
smoothness $\kappa$.

From \ref{a1.55}, if $G\in S^{\prime}$ is a \textbf{basis distribution} then,%
\[
G\left(  \cdot+a\right)  -\sum_{k\leq n}\frac{\left(  aD\right)  ^{k}}%
{k!}G=\left(  \mathcal{R}_{n+1}G\right)  \left(  \cdot,a\right)  ,\text{\quad
}n\geq0,
\]

where%
\begin{equation}
\left(  \mathcal{R}_{n+1}G\right)  \left(  \cdot,a\right)  =\frac{\sqrt{2\pi}%
}{n!}\left(  \left(  ia\xi\right)  ^{n+1}\overline{\widehat{g_{n}}}\left(
a\xi\right)  \widehat{G}\right)  ^{\vee}.\label{a781}%
\end{equation}

Using the expression \ref{a114} for $e^{ia\xi}$ i.e.%
\[
e^{ia\xi}=\sum_{k\leq n}\frac{\left(  ia\xi\right)  ^{k}}{k!}+\frac{\sqrt
{2\pi}}{n!}\left(  ia\xi\right)  ^{n+1}\overline{\widehat{g_{n}}}\left(
a\xi\right)  ,
\]

we can write%
\[
\left(  \mathcal{R}_{n+1}G\right)  \left(  \cdot,a\right)  =\left(  \left(
e^{ia\xi}-\sum\limits_{k=0}^{n}\frac{\left(  ia\xi\right)  ^{k}}{k!}\right)
\widehat{G}\right)  ^{\vee}.
\]

Thus%
\begin{align*}
\left\vert \left(  \mathcal{R}_{n+1}G\right)  \left(  \cdot,a\right)
\right\vert  & \leq\int\left\vert \left(  e^{ia\xi}-\sum\limits_{k=0}^{n}%
\frac{\left(  ia\xi\right)  ^{k}}{k!}\right)  \widehat{G}\right\vert \\
& =\int\left\vert \left(  e^{ia\xi}-\sum\limits_{k=0}^{n}\frac{\left(
ia\xi\right)  ^{k}}{k!}\right)  \sqrt{\widehat{G}}\sqrt{\widehat{G}%
}\right\vert \\
& \Rightarrow Cauchy-Schwartz\text{ }inequality\Rightarrow\\
& \leq\left(  \int\widehat{G}\right)  ^{1/2}\left(  \int\left\vert \left(
e^{ia\xi}-\sum\limits_{k=0}^{n}\frac{\left(  ia\xi\right)  ^{k}}{k!}\right)
\sqrt{\widehat{G}}\right\vert ^{2}\right)  ^{\frac{1}{2}}\\
& =\left(  \int\widehat{G}\right)  ^{1/2}\left(  \int\left\vert e^{ia\xi}%
-\sum\limits_{k=0}^{n}\frac{\left(  ia\xi\right)  ^{k}}{k!}\right\vert
^{2}\widehat{G}\right)  ^{1/2}\\
& =\left(  \left(  2\pi\right)  ^{\frac{d}{2}}\left(  2\pi\right)  ^{-\frac
{d}{2}}\int\widehat{G}\right)  ^{\frac{1}{2}}\left(  \int\left\vert e^{ia\xi
}-\sum\limits_{k=0}^{n}\frac{\left(  ia\xi\right)  ^{k}}{k!}\right\vert
^{2}\widehat{G}\right)  ^{\frac{1}{2}}\\
& =\left(  2\pi\right)  ^{\frac{d}{4}}\sqrt{G\left(  0\right)  }\left(
\int\left\vert e^{ia\xi}-\sum\limits_{k=0}^{n}\frac{\left(  ia\xi\right)
^{k}}{k!}\right\vert ^{2}\widehat{G}\right)  ^{\frac{1}{2}}.
\end{align*}

This estimate now enables us to directly use the approach of Subsection
\ref{SbSect_TaylorDataApproach2}.

??? FINISH! ??

\section{An inverse Fourier transform for data
functions\label{Sect_another_inv_Four_X}}

We now recall some properties of the semi-Hilbert data spaces $X_{w}^{\theta}$
(Definition \ref{Def_Xwth}) needed in this section. These results are taken
from Section \ref{Sect_Xwo_Xwth} and the relevant theorem is given in brackets.

\begin{summary}
\label{Sum_properties_Xwm_Ch3}Suppose $w$ is a weight function with property
W2. If $f\in X_{w}^{\theta}$ then $\widehat{f}\in L_{loc}^{1}\left(
\mathbb{R}^{d}\setminus0\right)  $ and we can define a.e. the function
$f_{F}:\mathbb{R}^{d}\rightarrow\mathbb{C}$ by: $f_{F}=\widehat{f}$ on
$\mathbb{R}^{d}\setminus0$. Further:

\begin{enumerate}
\item The seminorm and semi-inner product are given by%
\begin{equation}
\int w\left\vert \cdot\right\vert ^{2\theta}f_{F}\overline{g_{F}}=\left\langle
f,g\right\rangle _{w,\theta},\quad\int w\left\vert \cdot\right\vert ^{2\theta
}\left\vert f_{F}\right\vert ^{2}=\left\vert f\right\vert _{w,\theta}%
^{2}.\label{a50}%
\end{equation}

An alternative definition of $X_{w}^{\theta}$ is%
\begin{equation}
X_{w}^{\theta}=\left\{  f\in S^{\prime}:\xi^{\alpha}\widehat{f}\in L_{loc}%
^{1}\text{ }if\text{ }\left\vert \alpha\right\vert =\theta;\text{\thinspace
}\int w\left\vert \cdot\right\vert ^{2\theta}\left\vert f_{F}\right\vert
^{2}<\infty\right\}  .\label{a52}%
\end{equation}

Note that by part 2 Lemma \ref{Lem_properties_Xwm_distrib} the condition
$\widehat{f}\in L_{loc}^{1}\left(  \mathbb{R}^{d}\setminus0\right)  $ is
actually implied by $\xi^{\alpha}\widehat{f}\in L_{loc}^{1}$ for all
$\left\vert \alpha\right\vert =\theta$.\medskip

\item The functional $\left\vert \cdot\right\vert _{w,\theta}$ is a seminorm.
In fact, $\operatorname*{null}\left\vert \cdot\right\vert _{w,\theta
}=P_{\theta-1}$ and $P\cap X_{w}^{\theta}=P_{\theta-1}$ (part 3 Theorem
\ref{Thm_properties_Xwm}).\medskip

\item $f_{F}\in S_{\emptyset,\theta}^{\prime}$ with action $\int f_{F}\phi$ on
$S_{\emptyset,\theta}$. Also $f_{F}=\widehat{f}$ on $S_{\emptyset,\theta}$
(part 2 Theorem \ref{Thm_property_(g)_F}).

\item $X_{w}^{\theta}$ is complete in the seminorm sense (Theorem
\ref{Thm_Xw_complete}).

\item If $w$ has properties W2.1 and W3 for order $\theta$ and smoothness
parameter $\kappa$ then $X_{w}^{\theta}\subset C_{BP}^{\left(  \left\lfloor
\kappa\right\rfloor \right)  }$ (Theorem \ref{Thm_Xwth_W3_smooth}).
\end{enumerate}
\end{summary}

If a weight function $w$ has properties W2 and W3 for order $\theta$ and
smoothness parameter $\kappa$, we know from part 5 of Summary
\ref{Sum_properties_Xwm_Ch3} that $X_{w}^{\theta}\subset C_{BP}^{\left(
\left\lfloor \kappa\right\rfloor \right)  }$ i.e. functions in $X_{w}^{\theta
}$ have polynomial growth. The next lemma will allow us to use Corollary
\ref{Cor_Thm_integ(|DQ(exp)fF|)} to prove an inverse Fourier transform result
for the functions in $X_{w}^{\theta}$ and to derive upper bounds for the
growth rate of the derivatives near infinity.

\begin{lemma}
\label{Lem_integ(DQ(exp)|^2/w|.|_2}Suppose the weight function $w$ has
properties W2.1 and W3 for order $\theta$ and smoothness parameter $\kappa$.
Then if $\left\vert \gamma\right\vert \leq\underline{\kappa}$,%
\begin{equation}
\int\frac{\left\vert D_{x}^{\gamma}\mathcal{Q}_{\emptyset,\theta,\xi}\left(
e^{i\left(  x,\xi\right)  }\right)  \right\vert ^{2}}{w\left(  \xi\right)
\left\vert \xi\right\vert ^{2\theta}}d\xi\leq\left\{
\begin{array}
[c]{ll}%
\max\left\{
\begin{array}
[c]{c}%
\underline{C}_{\theta,r_{4}}^{\left(  \rho\right)  },\\
\overline{C}_{\theta-\left\vert \gamma\right\vert ,r_{4}}^{\left(
\rho\right)  }%
\end{array}
\right\}  ^{2}%
\begin{pmatrix}
K_{\gamma,r_{4}}^{\left(  w\right)  }\\
or\\
L_{\left\vert \gamma\right\vert ,r_{4}}^{\left(  w\right)  }%
\end{pmatrix}
\left(  1+\sum\limits_{k\leq\theta-\left\vert \gamma\right\vert }%
\frac{\left\vert x\right\vert ^{2k}}{k!}\right)  , & \left\vert \gamma
\right\vert <\theta,\\
\int\frac{\xi^{2\gamma}}{w\left\vert \cdot\right\vert ^{2\theta}}\text{
}or\text{ }\int\frac{\left\vert \cdot\right\vert ^{2\left\vert \gamma
\right\vert }}{w\left\vert \cdot\right\vert ^{2\theta}}, & \left\vert
\gamma\right\vert \geq\theta,
\end{array}
\right. \label{a062}%
\end{equation}

where $L_{\left\vert \gamma\right\vert ,r_{4}}^{\left(  w\right)  }$ is given
by \ref{a2.06} and $K_{\gamma,r_{4}}^{\left(  w\right)  }$ is given by
\ref{a2.06}, and both are independent of $x$.
\end{lemma}

\begin{proof}
There are two cases to be considered: $\left\vert \gamma\right\vert <\theta$
and $\left\vert \gamma\right\vert \geq\theta$. Define the constant $r_{4}$ by
$r_{4}=0$ if $w$ has property W3.1 and by $r_{4}=r_{3}$ if $w$ has property
W3.2. In both cases we will split the range of integration into the two
concentric regions defined by the sphere $S\left(  0;r_{4}\right)  $. We use
the estimate of part 2 of Theorem \ref{Thm_bound_on_g(e,x)_2} for $\left\vert
D_{x}^{\gamma}\mathcal{Q}_{\emptyset,\theta,\xi}\left(  e^{i\left(
x,\xi\right)  }\right)  \right\vert $ inside the sphere $S\left(
0;r_{4}\right)  $ and the estimate of part 1 of Theorem
\ref{Thm_bound_on_g(e,x)_2} outside this sphere.\medskip

\fbox{\textbf{Case 1} $\theta\geq\underline{\kappa}$ and $\left\vert
\gamma\right\vert \leq\underline{\kappa}$} By part 2 of Theorem
\ref{Thm_bound_on_g(e,x)_2},%
\begin{align}
&  \int\limits_{\left\vert \xi\right\vert \leq r_{4}}\frac{\left\vert
D_{x}^{\gamma}\mathcal{Q}_{\emptyset,\theta,\xi}\left(  e^{i\left(
x,\xi\right)  }\right)  \right\vert ^{2}d\xi}{w\left(  \xi\right)  \left\vert
\xi\right\vert ^{2\theta}}\nonumber\\
&  \leq\int\limits_{\left\vert \cdot\right\vert \leq r_{4}}\left(
\underline{C}_{\theta,r_{4}}^{\left(  \rho\right)  }\right)  ^{2}%
\frac{\left\vert \cdot\right\vert ^{2\theta}}{w\left\vert \cdot\right\vert
^{2\theta}}\left(  \left\vert \xi^{\gamma}\right\vert +\frac{\left\vert
x\widehat{\xi}\right\vert ^{\theta-\left\vert \gamma\right\vert }}{\left(
\theta-\left\vert \gamma\right\vert \right)  !}\right)  ^{2}\nonumber\\
&  =\left(  \underline{C}_{\theta,r_{4}}^{\left(  \rho\right)  }\right)
^{2}\int\limits_{\left\vert \cdot\right\vert \leq r_{4}}\frac{1}{w}\left(
\left\vert \widehat{\xi}^{\gamma}\right\vert \left\vert \cdot\right\vert
^{\left\vert \gamma\right\vert }+\left\vert x\right\vert ^{\theta-\left\vert
\gamma\right\vert }\frac{\left\vert \widehat{x}\widehat{\xi}\right\vert
^{\theta-\left\vert \gamma\right\vert }}{\left(  \theta-\left\vert
\gamma\right\vert \right)  !}\right)  ^{2}\nonumber\\
&  =\left(  \underline{C}_{\theta,r_{4}}^{\left(  \rho\right)  }\right)
^{2}\int\limits_{\left\vert \cdot\right\vert \leq r_{4}}\frac{1}{w}\left(
\left\vert \widehat{\xi}^{\gamma}\right\vert \left\vert \cdot\right\vert
^{\left\vert \gamma\right\vert }+\frac{\left\vert x\right\vert ^{\theta
-\left\vert \gamma\right\vert }}{\sqrt{\left(  \theta-\left\vert
\gamma\right\vert \right)  !}}\frac{\left\vert \widehat{x}\widehat{\xi
}\right\vert ^{\theta-\left\vert \gamma\right\vert }}{\sqrt{\left(
\theta-\left\vert \gamma\right\vert \right)  !}}\right)  ^{2}\nonumber\\
&  and\text{ }the\text{ }Cauchy-Schwartz\text{ }inequality\text{
}implies\nonumber\\
&  \leq\left(  \underline{C}_{\theta,r_{4}}^{\left(  \rho\right)  }\right)
^{2}\int\limits_{\left\vert \cdot\right\vert \leq r_{4}}\frac{1}{w}\left(
\left\vert \widehat{\xi}^{\gamma}\right\vert ^{2}\left\vert \cdot\right\vert
^{2\left\vert \gamma\right\vert }+\frac{\left\vert \widehat{x}\widehat{\xi
}\right\vert ^{2\left(  \theta-\left\vert \gamma\right\vert \right)  }%
}{\left(  \theta-\left\vert \gamma\right\vert \right)  !}\right)  \left(
1+\frac{\left\vert x\right\vert ^{2\left(  \theta-\left\vert \gamma\right\vert
\right)  }}{\left(  \theta-\left\vert \gamma\right\vert \right)  !}\right)
\nonumber\\
&  =\left(  \underline{C}_{\theta,r_{4}}^{\left(  \rho\right)  }\right)
^{2}\left(  \int\limits_{\left\vert \cdot\right\vert \leq r_{4}}\left\vert
\widehat{\xi}^{2\gamma}\right\vert \frac{\left\vert \cdot\right\vert
^{2\left\vert \gamma\right\vert }}{w}+\frac{1}{\left(  \theta-\left\vert
\gamma\right\vert \right)  !}\int\limits_{\left\vert \cdot\right\vert \leq
r_{4}}\frac{\left\vert \widehat{x}\widehat{\xi}\right\vert ^{2\left(
\theta-\left\vert \gamma\right\vert \right)  }}{w}\right)  \left(
1+\frac{\left\vert x\right\vert ^{2\left(  \theta-\left\vert \gamma\right\vert
\right)  }}{\left(  \theta-\left\vert \gamma\right\vert \right)  !}\right)
\label{a2.09}\\
&  \leq\left(  \underline{C}_{\theta,r_{4}}^{\left(  \rho\right)  }\right)
^{2}\left(  \int\limits_{\left\vert \cdot\right\vert \leq r_{4}}%
\frac{\left\vert \cdot\right\vert ^{2\left\vert \gamma\right\vert }}{w}%
+\frac{1}{\left(  \theta-\left\vert \gamma\right\vert \right)  !}%
\int\limits_{\left\vert \cdot\right\vert \leq r_{4}}\frac{1}{w}\right)
\left(  1+\frac{\left\vert x\right\vert ^{2\left(  \theta-\left\vert
\gamma\right\vert \right)  }}{\left(  \theta-\left\vert \gamma\right\vert
\right)  !}\right)  ,\label{a2.10}%
\end{align}

which exists since property W2.1 is $1/w\in L_{loc}^{1}$. Set%
\[
\underline{C}_{\gamma,r}^{\left(  w\right)  }:=\int\limits_{\left\vert
\cdot\right\vert \leq r_{4}}\left\vert \widehat{\xi}^{2\gamma}\right\vert
\frac{\left\vert \cdot\right\vert ^{2\left\vert \gamma\right\vert }}{w}%
+\frac{1}{\left(  \theta-\left\vert \gamma\right\vert \right)  !}%
\int\limits_{\left\vert \cdot\right\vert \leq r_{4}}\frac{\left\vert
\widehat{x}\widehat{\xi}\right\vert ^{2\left(  \theta-\left\vert
\gamma\right\vert \right)  }}{w},\quad\left\vert \gamma\right\vert <\theta.
\]

Further, by part 1 of Theorem \ref{Thm_bound_on_g(e,x)_2} when $n=\theta$,
\[
\left\vert D_{x}^{\gamma}\mathcal{Q}_{\emptyset,\theta,\xi}\left(  e^{ix\xi
}\right)  \right\vert \leq\overline{C}_{n-\left\vert \gamma\right\vert
,r}^{\left(  \rho\right)  }\left(  1+\sum\limits_{k<\theta-\left\vert
\gamma\right\vert }\frac{\left\vert x\widehat{\xi}\right\vert ^{k}}%
{k!}\right)  \left\vert \xi^{\gamma}\right\vert ,\quad\left\vert
\gamma\right\vert <\theta,\text{ }\left\vert \xi\right\vert \geq r,
\]

so%
\begin{align}
\int\limits_{\left\vert \xi\right\vert \geq r_{4}} &  \left\vert D_{x}%
^{\gamma}\mathcal{Q}_{\emptyset,\theta,\xi}\left(  e^{ix\xi}\right)
\right\vert ^{2}\frac{1}{w\left\vert \cdot\right\vert ^{2\theta}}\nonumber\\
&  \leq\int\limits_{\left\vert \cdot\right\vert \geq r_{4}}\left(
\overline{C}_{\theta-\left\vert \gamma\right\vert ,r}^{\left(  \rho\right)
}\right)  ^{2}\left(  1+\sum\limits_{k<\theta-\left\vert \gamma\right\vert
}\frac{\left\vert x\widehat{\xi}\right\vert ^{k}}{k!}\right)  ^{2}%
\frac{\left\vert \xi^{\gamma}\right\vert ^{2}}{w\left\vert \cdot\right\vert
^{2\theta}}\nonumber\\
&  =\left(  \overline{C}_{\theta-\left\vert \gamma\right\vert ,r_{4}}^{\left(
\rho\right)  }\right)  ^{2}\int\limits_{\left\vert \cdot\right\vert \geq
r_{4}}\left(  1+\sum\limits_{k<\theta-\left\vert \gamma\right\vert }%
\frac{\left\vert x\right\vert ^{k}}{\sqrt{k!}}\frac{\left\vert \widehat
{x}\widehat{\xi}\right\vert ^{k}}{\sqrt{k!}}\right)  ^{2}\frac{\left\vert
\xi^{\gamma}\right\vert ^{2}}{w\left\vert \cdot\right\vert ^{2\theta}%
},\nonumber\\
&  and\text{ }the\text{ }Cauchy-Schwartz\text{ }inequality\text{
}implies\nonumber\\
&  \leq\left(  \overline{C}_{\theta-\left\vert \gamma\right\vert ,r_{4}%
}^{\left(  \rho\right)  }\right)  ^{2}\int\limits_{\left\vert \cdot\right\vert
\geq r_{4}}\left(  1+\sum\limits_{k<\theta-\left\vert \gamma\right\vert }%
\frac{\left\vert x\right\vert ^{2k}}{k!}\right)  \left(  1+\sum
\limits_{k<\theta-\left\vert \gamma\right\vert }\frac{\left\vert \widehat
{x}\widehat{\xi}\right\vert ^{2k}}{k!}\right)  \frac{\left\vert \xi^{\gamma
}\right\vert ^{2}}{w\left\vert \cdot\right\vert ^{2\theta}}\nonumber\\
&  =\left(  \overline{C}_{\theta-\left\vert \gamma\right\vert ,r_{4}}^{\left(
\rho\right)  }\right)  ^{2}\left(  \int\limits_{\left\vert \cdot\right\vert
\geq r_{4}}\left(  1+\sum\limits_{k<\theta-\left\vert \gamma\right\vert }%
\frac{\left\vert \widehat{x}\widehat{\xi}\right\vert ^{2k}}{k!}\right)
\frac{\left\vert \xi^{\gamma}\right\vert ^{2}}{w\left\vert \cdot\right\vert
^{2\theta}}\right)  \left(  1+\sum\limits_{k<\theta-\left\vert \gamma
\right\vert }\frac{\left\vert x\right\vert ^{2k}}{k!}\right) \nonumber\\
&  =\left(  \overline{C}_{\theta-\left\vert \gamma\right\vert ,r_{4}}^{\left(
\rho\right)  }\right)  ^{2}\left(  \int\limits_{\left\vert \cdot\right\vert
\geq r_{4}}\left\vert \widehat{\xi}^{2\gamma}\right\vert \frac{\left\vert
\cdot\right\vert ^{2\left\vert \gamma\right\vert }}{w\left\vert \cdot
\right\vert ^{2\theta}}+\sum\limits_{k<\theta-\left\vert \gamma\right\vert
}\frac{1}{k!}\int\limits_{\left\vert \cdot\right\vert \geq r_{4}}\left\vert
\widehat{\xi}^{2\gamma}\right\vert \left\vert \widehat{x}\widehat{\xi
}\right\vert ^{2k}\frac{\left\vert \cdot\right\vert ^{2\left\vert
\gamma\right\vert }}{w\left\vert \cdot\right\vert ^{2\theta}}\right)
\times\nonumber\\
&  \qquad\times\left(  1+\sum\limits_{k<\theta-\left\vert \gamma\right\vert
}\frac{\left\vert x\right\vert ^{2k}}{k!}\right) \label{a2.05}\\
&  \leq\left(  \overline{C}_{\theta-\left\vert \gamma\right\vert ,r_{4}%
}^{\left(  \rho\right)  }\right)  ^{2}\left(  \int\limits_{\left\vert
\cdot\right\vert \geq r_{4}}\frac{\left\vert \cdot\right\vert ^{2\left\vert
\gamma\right\vert }}{w\left\vert \cdot\right\vert ^{2\theta}}+\sum
\limits_{k<\theta-\left\vert \gamma\right\vert }\frac{1}{k!}\int%
\limits_{\left\vert \cdot\right\vert \geq r_{4}}\frac{\left\vert
\cdot\right\vert ^{2\left\vert \gamma\right\vert }}{w\left\vert \cdot
\right\vert ^{2\theta}}\right)  \left(  1+\sum\limits_{k<\theta-\left\vert
\gamma\right\vert }\frac{\left\vert x\right\vert ^{2k}}{k!}\right) \nonumber\\
&  =\left(  \overline{C}_{\theta-\left\vert \gamma\right\vert ,r_{4}}^{\left(
\rho\right)  }\right)  ^{2}\left(  1+\sum\limits_{k<\theta-\left\vert
\gamma\right\vert }\frac{1}{k!}\right)  \int\limits_{\left\vert \cdot
\right\vert \geq r_{4}}\frac{\left\vert \cdot\right\vert ^{2\left\vert
\gamma\right\vert }}{w\left\vert \cdot\right\vert ^{2\theta}}\left(
1+\sum\limits_{k<\theta-\left\vert \gamma\right\vert }\frac{\left\vert
x\right\vert ^{2k}}{k!}\right) \nonumber\\
&  <\left(  \overline{C}_{\theta-\left\vert \gamma\right\vert ,r_{4}}^{\left(
\rho\right)  }\right)  ^{2}\left(  1+e_{\theta-\left\vert \gamma\right\vert
-1}\right)  \left(  \int\limits_{\left\vert \cdot\right\vert \geq r_{4}}%
\frac{\left\vert \cdot\right\vert ^{2\left\vert \gamma\right\vert }%
}{w\left\vert \cdot\right\vert ^{2\theta}}\right)  \left(  1+\sum
\limits_{k<\theta-\left\vert \gamma\right\vert }\frac{\left\vert x\right\vert
^{2k}}{k!}\right)  ,\label{a2.02}%
\end{align}

and since $w$ has property W3, the definition of $r_{4}$ and part 1 of Theorem
\ref{Thm_weight_property_relat}, imply the last integral exists.

Adding the estimates \ref{a2.10} and \ref{a2.02} we get%
\begin{align*}
& \int\frac{\left\vert D_{x}^{\gamma}\mathcal{Q}_{\emptyset,\theta,\xi}\left(
e^{i\left(  x,\xi\right)  }\right)  \right\vert ^{2}}{w\left(  \xi\right)
\left\vert \xi\right\vert ^{2\theta}}d\xi\\
& \leq\left(  \underline{C}_{\theta,r_{4}}^{\left(  \rho\right)  }\right)
^{2}\left(  \int\limits_{\left\vert \cdot\right\vert \leq r_{4}}%
\frac{\left\vert \cdot\right\vert ^{2\left\vert \gamma\right\vert }}{w}%
+\frac{1}{\left(  \theta-\left\vert \gamma\right\vert \right)  !}%
\int\limits_{\left\vert \cdot\right\vert \leq r_{4}}\frac{1}{w}\right)
\left(  1+\frac{\left\vert x\right\vert ^{2\left(  \theta-\left\vert
\gamma\right\vert \right)  }}{\left(  \theta-\left\vert \gamma\right\vert
\right)  !}\right)  +\\
& \qquad+\left(  \overline{C}_{\theta-\left\vert \gamma\right\vert ,r_{4}%
}^{\left(  \rho\right)  }\right)  ^{2}\left(  1+e_{??}\right)  \left(
\int\limits_{\left\vert \cdot\right\vert \geq r_{4}}\frac{\left\vert
\cdot\right\vert ^{2\left\vert \gamma\right\vert }}{w\left\vert \cdot
\right\vert ^{2\theta}}\right)  \left(  1+\sum\limits_{k<\theta-\left\vert
\gamma\right\vert }\frac{\left\vert x\right\vert ^{2k}}{k!}\right) \\
& \leq\left\{  \left(  \underline{C}_{\theta,r_{4}}^{\left(  \rho\right)
}\right)  ^{2}\left(  \int\limits_{\left\vert \cdot\right\vert \leq r_{4}%
}\frac{\left\vert \cdot\right\vert ^{2\left\vert \gamma\right\vert }}{w}%
+\frac{1}{\left(  \theta-\left\vert \gamma\right\vert \right)  !}%
\int\limits_{\left\vert \cdot\right\vert \leq r_{4}}\frac{1}{w}\right)
+\left(  \overline{C}_{\theta-\left\vert \gamma\right\vert ,r_{4}}^{\left(
\rho\right)  }\right)  ^{2}\left(  1+e_{??}\right)  \int\limits_{\left\vert
\cdot\right\vert \geq r_{4}}\frac{\left\vert \cdot\right\vert ^{2\left\vert
\gamma\right\vert }}{w\left\vert \cdot\right\vert ^{2\theta}}\right\}
\times\\
& \qquad\times\left(  1+\sum\limits_{k\leq\theta-\left\vert \gamma\right\vert
}\frac{\left\vert x\right\vert ^{2k}}{k!}\right) \\
& \leq\left(  \max\left\{
\begin{array}
[c]{l}%
\underline{C}_{\theta,r_{4}}^{\left(  \rho\right)  },\\
\overline{C}_{\theta-\left\vert \gamma\right\vert ,r_{4}}^{\left(
\rho\right)  }%
\end{array}
\right\}  \right)  ^{2}\left(  \int\limits_{\left\vert \cdot\right\vert \leq
r_{4}}\frac{\left\vert \cdot\right\vert ^{2\left\vert \gamma\right\vert }}%
{w}+\frac{1}{\left(  \theta-\left\vert \gamma\right\vert \right)  !}%
\int\limits_{\left\vert \cdot\right\vert \leq r_{4}}\frac{1}{w}+\left(
1+e\right)  \int\limits_{\left\vert \cdot\right\vert \geq r_{4}}%
\frac{\left\vert \cdot\right\vert ^{2\left\vert \gamma\right\vert }%
}{w\left\vert \cdot\right\vert ^{2\theta}}\right)  \times\\
& \qquad\times\left(  1+\sum\limits_{k\leq\theta-\left\vert \gamma\right\vert
}\frac{\left\vert x\right\vert ^{2k}}{k!}\right) \\
& =\left(  \max\left\{
\begin{array}
[c]{c}%
\underline{C}_{\theta,r_{4}}^{\left(  \rho\right)  },\\
\overline{C}_{\theta-\left\vert \gamma\right\vert ,r_{4}}^{\left(
\rho\right)  }%
\end{array}
\right\}  \right)  ^{2}C_{\left\vert \gamma\right\vert }^{\left(  w\right)
}\left(  1+\sum\limits_{k\leq\theta-\left\vert \gamma\right\vert }%
\frac{\left\vert x\right\vert ^{2k}}{k!}\right)  ,
\end{align*}

where%
\begin{equation}
C_{m}^{\left(  w\right)  }:=\int\limits_{\left\vert \cdot\right\vert \leq
r_{4}}\frac{\left\vert \cdot\right\vert ^{2m}}{w}+\frac{1}{\left(
\theta-m\right)  !}\int\limits_{\left\vert \cdot\right\vert \leq r_{4}}%
\frac{1}{w}+\left(  1+e\right)  \int\limits_{\left\vert \cdot\right\vert \geq
r_{4}}\frac{\left\vert \cdot\right\vert ^{2m}}{w\left\vert \cdot\right\vert
^{2\theta}},\quad m=0,1,2,\ldots,\min\left\{  \theta-1,\underline{\kappa
}\right\}  .\label{a2.14}%
\end{equation}

Adding the estimates \ref{a2.09} and \ref{a2.05} we get:%
\begin{align*}
\int &  \frac{\left\vert D_{x}^{\gamma}\mathcal{Q}_{\emptyset,\theta,\xi
}\left(  e^{i\left(  x,\xi\right)  }\right)  \right\vert ^{2}}{w\left(
\xi\right)  \left\vert \xi\right\vert ^{2\theta}}d\xi\\
&  \leq\left(  \underline{C}_{\theta,r_{4}}^{\left(  \rho\right)  }\right)
^{2}\left(  \int\limits_{\left\vert \cdot\right\vert \leq r_{4}}\left\vert
\widehat{\xi}^{2\gamma}\right\vert \frac{\left\vert \cdot\right\vert
^{2\left\vert \gamma\right\vert }}{w}+\frac{1}{\left(  \theta-\left\vert
\gamma\right\vert \right)  !}\int\limits_{\left\vert \cdot\right\vert \leq
r_{4}}\frac{\left\vert \widehat{x}\widehat{\xi}\right\vert ^{2\left(
\theta-\left\vert \gamma\right\vert \right)  }}{w}\right)  \left(
1+\frac{\left\vert x\right\vert ^{2\left(  \theta-\left\vert \gamma\right\vert
\right)  }}{\left(  \theta-\left\vert \gamma\right\vert \right)  !}\right)
+\\
&  +\left(  \overline{C}_{\theta-\left\vert \gamma\right\vert ,r_{4}}^{\left(
\rho\right)  }\right)  ^{2}\left(  \int\limits_{\left\vert \cdot\right\vert
\geq r_{4}}\left\vert \widehat{\xi}^{2\gamma}\right\vert \frac{\left\vert
\cdot\right\vert ^{2\left\vert \gamma\right\vert }}{w\left\vert \cdot
\right\vert ^{2\theta}}+\sum\limits_{k<\theta-\left\vert \gamma\right\vert
}\frac{1}{k!}\int\limits_{\left\vert \cdot\right\vert \geq r_{4}}\left\vert
\widehat{\xi}^{2\gamma}\right\vert \left\vert \widehat{x}\widehat{\xi
}\right\vert ^{2k}\frac{\left\vert \cdot\right\vert ^{2\left\vert
\gamma\right\vert }}{w\left\vert \cdot\right\vert ^{2\theta}}\right)  \left(
1+\sum\limits_{k<\theta-\left\vert \gamma\right\vert }\frac{\left\vert
x\right\vert ^{2k}}{k!}\right) \\
&  \leq\left(  \max\left\{
\begin{array}
[c]{c}%
\underline{C}_{\theta,r_{4}}^{\left(  \rho\right)  },\\
\overline{C}_{\theta-\left\vert \gamma\right\vert ,r_{4}}^{\left(
\rho\right)  }%
\end{array}
\right\}  \right)  ^{2}\left(
\begin{array}
[c]{c}%
\int\limits_{\left\vert \cdot\right\vert \leq r_{4}}\left\vert \widehat{\xi
}^{2\gamma}\right\vert \frac{\left\vert \cdot\right\vert ^{2\left\vert
\gamma\right\vert }}{w}+\frac{1}{\left(  \theta-\left\vert \gamma\right\vert
\right)  !}\int\limits_{\left\vert \cdot\right\vert \leq r_{4}}\frac
{\left\vert \widehat{x}\widehat{\xi}\right\vert ^{2\left(  \theta-\left\vert
\gamma\right\vert \right)  }}{w}+\\
+\int\limits_{\left\vert \cdot\right\vert \geq r_{4}}\left\vert \widehat{\xi
}^{2\gamma}\right\vert \frac{\left\vert \cdot\right\vert ^{2\left\vert
\gamma\right\vert }}{w\left\vert \cdot\right\vert ^{2\theta}}+\sum
\limits_{k<\theta-\left\vert \gamma\right\vert }\frac{1}{k!}\int%
\limits_{\left\vert \cdot\right\vert \geq r_{4}}\left\vert \widehat{\xi
}^{2\gamma}\right\vert \left\vert \widehat{x}\widehat{\xi}\right\vert
^{2k}\frac{\left\vert \cdot\right\vert ^{2\left\vert \gamma\right\vert }%
}{w\left\vert \cdot\right\vert ^{2\theta}}%
\end{array}
\right)  \times\\
&  \qquad\times\left(  1+\sum\limits_{k\leq\theta-\left\vert \gamma\right\vert
}\frac{\left\vert x\right\vert ^{2k}}{k!}\right) \\
&  =\left(  \max\left\{
\begin{array}
[c]{c}%
\underline{C}_{\theta,r_{4}}^{\left(  \rho\right)  },\\
\overline{C}_{\theta-\left\vert \gamma\right\vert ,r_{4}}^{\left(
\rho\right)  }%
\end{array}
\right\}  \right)  ^{2}L_{\gamma}^{\left(  w\right)  }\left(  1+\sum
\limits_{k\leq\theta-\left\vert \gamma\right\vert }\frac{\left\vert
x\right\vert ^{2k}}{k!}\right)  ,
\end{align*}

where%
\begin{equation}
L_{\gamma,r_{4}}^{\left(  w\right)  }:=\left(
\begin{array}
[c]{c}%
\int\limits_{\left\vert \cdot\right\vert \leq r_{4}}\left\vert \widehat{\xi
}^{2\gamma}\right\vert \frac{\left\vert \cdot\right\vert ^{2\left\vert
\gamma\right\vert }}{w}+\frac{1}{\left(  \theta-\left\vert \gamma\right\vert
\right)  !}\int\limits_{\left\vert \cdot\right\vert \leq r_{4}}\frac
{\left\vert \widehat{x}\widehat{\xi}\right\vert ^{2\left(  \theta-\left\vert
\gamma\right\vert \right)  }}{w}+\\
+\int\limits_{\left\vert \cdot\right\vert \geq r_{4}}\left\vert \widehat{\xi
}^{2\gamma}\right\vert \frac{\left\vert \cdot\right\vert ^{2\left\vert
\gamma\right\vert }}{w\left\vert \cdot\right\vert ^{2\theta}}+\sum
\limits_{k<\theta-\left\vert \gamma\right\vert }\frac{1}{k!}\int%
\limits_{\left\vert \cdot\right\vert \geq r_{4}}\left\vert \widehat{\xi
}^{2\gamma}\right\vert \left\vert \widehat{x}\widehat{\xi}\right\vert
^{2k}\frac{\left\vert \cdot\right\vert ^{2\left\vert \gamma\right\vert }%
}{w\left\vert \cdot\right\vert ^{2\theta}},
\end{array}
\right)  ,\quad\left\vert \gamma\right\vert \leq\min\left\{  \theta
-1,\underline{\kappa}\right\}  .\label{a2.06}%
\end{equation}
\medskip

\fbox{\textbf{Case 2} $\theta\leq\underline{\kappa}$ and $\left\vert
\gamma\right\vert \leq\underline{\kappa}$.} From \ref{a21},%
\begin{align*}
\int\frac{\left\vert D_{x}^{\gamma}\mathcal{Q}_{\emptyset,\theta,\xi}\left(
e^{ix\xi}\right)  \right\vert ^{2}}{w\left\vert \cdot\right\vert ^{2\theta}}
& =\int\frac{\left\vert \left(  i\xi\right)  ^{\gamma}\right\vert ^{2}%
}{w\left\vert \cdot\right\vert ^{2\theta}}=\int\frac{\xi^{2\gamma}%
}{w\left\vert \cdot\right\vert ^{2\theta}}\leq\int\frac{\left\vert
\cdot\right\vert ^{2\left\vert \gamma\right\vert }}{w\left\vert \cdot
\right\vert ^{2\theta}}=\\
& =\int\limits_{\left\vert \cdot\right\vert \leq r_{4}}\frac{\left\vert
\cdot\right\vert ^{2\left\vert \gamma\right\vert }}{w\left\vert \cdot
\right\vert ^{2\theta}}+\int\limits_{\left\vert \cdot\right\vert \geq r_{4}%
}\frac{\left\vert \cdot\right\vert ^{2\left\vert \gamma\right\vert }%
}{w\left\vert \cdot\right\vert ^{2\theta}}\\
& =\int\limits_{\left\vert \cdot\right\vert \leq r_{4}}\frac{\left\vert
\cdot\right\vert ^{2\left(  \left\vert \gamma\right\vert -\theta\right)  }}%
{w}+\int\limits_{\left\vert \cdot\right\vert \geq r_{4}}\frac{\left\vert
\cdot\right\vert ^{2\left\vert \gamma\right\vert }}{w\left\vert \cdot
\right\vert ^{2\theta}}\\
& \leq r_{4}^{2\left(  \left\vert \gamma\right\vert -\theta\right)  }%
\int\limits_{\left\vert \cdot\right\vert \leq r_{4}}\frac{1}{w}+\int%
\limits_{\left\vert \cdot\right\vert \geq r_{4}}\frac{\left\vert
\cdot\right\vert ^{2\left\vert \gamma\right\vert }}{w\left\vert \cdot
\right\vert ^{2\theta}}\\
& \leq\left(  1+r_{4}\right)  ^{2\left(  \left\vert \gamma\right\vert
-\theta\right)  }\left(  \int\limits_{\left\vert \cdot\right\vert \leq r_{4}%
}\frac{1}{w}+\int\limits_{\left\vert \cdot\right\vert \geq r_{4}}%
\frac{\left\vert \cdot\right\vert ^{2\left\vert \gamma\right\vert }%
}{w\left\vert \cdot\right\vert ^{2\theta}}\right) \\
& \leq\left(  1+r_{4}\right)  ^{2\underline{\kappa}}\max_{n\leq\underline
{\kappa}}\left(  \int\limits_{\left\vert \cdot\right\vert \leq r_{4}}\frac
{1}{w}+\int\limits_{\left\vert \cdot\right\vert \geq r_{4}}\frac{\left\vert
\cdot\right\vert ^{2n}}{w\left\vert \cdot\right\vert ^{2\theta}}\right) \\
& =\left(  1+r_{4}\right)  ^{2\underline{\kappa}}\left(  \int%
\limits_{\left\vert \cdot\right\vert \leq r_{4}}\frac{1}{w}+\max
_{n\leq\underline{\kappa}}\int\limits_{\left\vert \cdot\right\vert \geq r_{4}%
}\frac{\left\vert \cdot\right\vert ^{2n}}{w\left\vert \cdot\right\vert
^{2\theta}}\right)  .
\end{align*}

Since property W2.1 requires that $1/w\in L_{loc}^{1}$and property W3 implies
Theorem \ref{Thm_property_wt_fn_W3}, it follows that the last two integrals exist.
\end{proof}

\begin{example}
\label{Ex_Lem_integ(DQ(exp)|^2/w|.|_2}\textbf{Radial functions} ??
\end{example}

\begin{lemma}
\label{Lem_integ(DQ(exp)|^2/w|.|}Suppose the weight function $w$ has
properties W2.1 and W3 for order $\theta$ and smoothness parameter $\kappa$.
Then if $\left\vert \gamma\right\vert \leq\underline{\kappa}$ there exists a
constant $C_{w}$, independent of $x$, such that%
\begin{equation}
\int\frac{\left\vert D_{x}^{\gamma}\mathcal{Q}_{\emptyset,\theta,\xi}\left(
e^{ix\xi}\right)  \right\vert ^{2}}{w\left(  \xi\right)  \left\vert
\xi\right\vert ^{2\theta}}d\xi\leq\left\{
\begin{array}
[c]{ll}%
\left(  C_{w}\right)  ^{2}\left(  1+\left\vert x\right\vert \right)
^{2\left(  \theta-\left\vert \gamma\right\vert \right)  }, & \left\vert
\gamma\right\vert <\theta,\\
\left(  C_{w}\right)  ^{2}, & \left\vert \gamma\right\vert \geq\theta.
\end{array}
\right. \label{a06}%
\end{equation}

$C_{w}$ is given by \ref{a07} and only depends on the weight function $w$, the
function $\rho\in S_{1,2\theta}$ used to define $\mathcal{Q}_{\emptyset
,\theta}$ and on the parameters which define the weight function properties W3.
\end{lemma}

\begin{proof}
There are two cases to be considered: $\left\vert \gamma\right\vert <\theta$
and $\left\vert \gamma\right\vert \geq\theta$. Define the constant $r_{4}$ by
$r_{4}=0$ if $w$ has property W3.1 and by $r_{4}=r_{3}$ if $w$ has property
W3.2. In both cases we will split the range of integration into the two
concentric regions defined by the sphere $S\left(  0;r_{4}\right)  $. We use
the estimate of part 2 of Theorem \ref{Thm_bound_on_g(e,x)} for $\left\vert
D_{x}^{\gamma}\mathcal{Q}_{\emptyset,\theta,\xi}\left(  e^{i\left(
x,\xi\right)  }\right)  \right\vert $ inside the sphere $S\left(
0;r_{4}\right)  $ and the estimate of part 1 of Theorem ??
\ref{Thm_bound_on_g(e,x)} outside this sphere.\medskip

\fbox{\textbf{Case 1} ?? $\left\vert \gamma\right\vert <\theta$ and
$\left\vert \gamma\right\vert \leq\underline{\kappa}$} By part 2 of Theorem
\ref{Thm_bound_on_g(e,x)},
\begin{align*}
\int\limits_{\left\vert \xi\right\vert \leq r_{4}}\frac{\left\vert
D_{x}^{\gamma}\mathcal{Q}_{\emptyset,\theta,\xi}\left(  e^{i\left(
x,\xi\right)  }\right)  \right\vert ^{2}d\xi}{w\left(  \xi\right)  \left\vert
\xi\right\vert ^{2\theta}} &  \leq\int\limits_{\left\vert \cdot\right\vert
\leq r_{4}}\frac{\left(  C_{\theta,\left\vert \gamma\right\vert ,r_{4}%
}\right)  ^{2}\left\vert \cdot\right\vert ^{2\theta}\left(  1+\left\vert
x\right\vert \right)  ^{2\left(  \theta-\left\vert \gamma\right\vert \right)
}}{w\left\vert \cdot\right\vert ^{2\theta}}\\
&  =\left(  C_{\theta,\left\vert \gamma\right\vert ,r_{4}}\right)  ^{2}\left(
\int\limits_{\left\vert \cdot\right\vert \leq r_{4}}\frac{1}{w}\right)
\left(  1+\left\vert x\right\vert \right)  ^{2\left(  \theta-\left\vert
\gamma\right\vert \right)  }\\
&  \leq\left(  \max_{m<\theta}C_{\theta,m,r_{4}}\right)  ^{2}\left(
\int\limits_{\left\vert \cdot\right\vert \leq r_{4}}\frac{1}{w}\right)
\left(  1+\left\vert x\right\vert \right)  ^{2\left(  \theta-\left\vert
\gamma\right\vert \right)  },
\end{align*}

which exists since property W2.1 is $1/w\in L_{loc}^{1}$.

Further, by part 1 of Theorem \ref{Thm_bound_on_g(e,x)},%
\[
\left\vert D_{x}^{\gamma}\mathcal{Q}_{\emptyset,n,\xi}\left(  e^{i\left(
x,\xi\right)  }\right)  \right\vert \leq\left\{
\begin{array}
[c]{ll}%
C_{n-\left\vert \gamma\right\vert }\left\vert \xi\right\vert ^{\left\vert
\gamma\right\vert }\left(  1+\left\vert x\right\vert \right)  ^{n-\left\vert
\gamma\right\vert -1}, & \left\vert \gamma\right\vert <n,\\
\left\vert \xi\right\vert ^{\left\vert \gamma\right\vert }, & \left\vert
\gamma\right\vert \geq n,
\end{array}
\right.
\]

so that%
\begin{align*}
\int\limits_{\left\vert \xi\right\vert \geq r_{4}}\frac{\left\vert
D_{x}^{\gamma}\mathcal{Q}_{\emptyset,\theta,\xi}\left(  e^{i\left(
x,\xi\right)  }\right)  \right\vert ^{2}d\xi}{w\left(  \xi\right)  \left\vert
\xi\right\vert ^{2\theta}}  & \leq\int\limits_{\left\vert \cdot\right\vert
\geq r_{4}}\frac{\left(  C_{\theta-\left\vert \gamma\right\vert }\right)
^{2}\left\vert \cdot\right\vert ^{2\left\vert \gamma\right\vert }\left(
1+\left\vert x\right\vert \right)  ^{2\left(  \theta-\left\vert \gamma
\right\vert -1\right)  }}{w\left\vert \cdot\right\vert ^{2\theta}}\\
& =\left(  C_{\theta-\left\vert \gamma\right\vert }\right)  ^{2}\left(
\int\limits_{\left\vert \cdot\right\vert \geq r_{4}}\frac{\left\vert
\cdot\right\vert ^{2\left\vert \gamma\right\vert }}{w\left\vert \cdot
\right\vert ^{2\theta}}\right)  \left(  1+\left\vert x\right\vert \right)
^{2\left(  \theta-\left\vert \gamma\right\vert \right)  }\\
& \leq\left(  \max_{m<\theta}C_{m}\right)  ^{2}\left(  \max_{n<\underline
{\kappa}}\int\limits_{\left\vert \cdot\right\vert \geq r_{4}}\frac{\left\vert
\cdot\right\vert ^{2n}}{w\left\vert \cdot\right\vert ^{2\theta}}\right)
\left(  1+\left\vert x\right\vert \right)  ^{2\left(  \theta-\left\vert
\gamma\right\vert \right)  },
\end{align*}

and since $w$ has property W3, the definition of $r_{4}$ and part 1 of Theorem
\ref{Thm_weight_property_relat} implies the last integral exists. Adding the
last two estimates we get%
\begin{equation}
\int\frac{\left\vert D_{x}^{\gamma}\mathcal{Q}_{\emptyset,\theta,\xi}\left(
e^{i\left(  x,\xi\right)  }\right)  \right\vert ^{2}}{w\left(  \xi\right)
\left\vert \xi\right\vert ^{2\theta}}d\xi\leq\left(  C_{w}^{\prime}\right)
^{2}\left(  1+\left\vert x\right\vert \right)  ^{2\left(  \theta-\left\vert
\gamma\right\vert \right)  },\label{a39}%
\end{equation}

where
\begin{equation}
C_{w}^{\prime}=\max\left\{  \max_{m<\theta}C_{m},\max_{m<\theta}%
C_{\theta,m,r_{4}}\right\}  \left(  \int\limits_{\left\vert \cdot\right\vert
\leq r_{4}}\frac{1}{w}+\max_{n\leq\underline{\kappa}}\int\limits_{\left\vert
\cdot\right\vert \geq r_{4}}\frac{\left\vert \cdot\right\vert ^{2n}%
}{w\left\vert \cdot\right\vert ^{2\theta}}\right)  ^{\frac{1}{2}}.\label{a13}%
\end{equation}
\medskip

\fbox{\textbf{Case 2} ?? $\left\vert \gamma\right\vert \geq\theta$ and
$\left\vert \gamma\right\vert \leq\underline{\kappa}$}%
\begin{align}
\int\frac{\left\vert D_{x}^{\gamma}\mathcal{Q}_{\emptyset,\theta,\xi}\left(
e^{i\left(  x,\xi\right)  }\right)  \right\vert ^{2}}{w\left(  \xi\right)
\left\vert \xi\right\vert ^{2\theta}}d\xi & =\int\limits_{\left\vert
\cdot\right\vert \leq r_{4}}\frac{\left\vert \cdot\right\vert ^{2\left\vert
\gamma\right\vert }}{w\left\vert \cdot\right\vert ^{2\theta}}+\int%
\limits_{\left\vert \cdot\right\vert \geq r_{4}}\frac{\left\vert
\cdot\right\vert ^{2\left\vert \gamma\right\vert }}{w\left\vert \cdot
\right\vert ^{2\theta}}\nonumber\\
& =\int\limits_{\left\vert \cdot\right\vert \leq r_{4}}\left\vert
\cdot\right\vert ^{2\left(  \left\vert \gamma\right\vert -\theta\right)
}\frac{1}{w}+\int\limits_{\left\vert \cdot\right\vert \geq r_{4}}%
\frac{\left\vert \cdot\right\vert ^{2\left\vert \gamma\right\vert }%
}{w\left\vert \cdot\right\vert ^{2\theta}}\nonumber\\
& \leq r_{4}^{2\left(  \left\vert \gamma\right\vert -\theta\right)  }%
\int\limits_{\left\vert \cdot\right\vert \leq r_{4}}\frac{1}{w}+\int%
\limits_{\left\vert \cdot\right\vert \geq r_{4}}\frac{\left\vert
\cdot\right\vert ^{2\left\vert \gamma\right\vert }}{w\left\vert \cdot
\right\vert ^{2\theta}}\nonumber\\
& \leq\left(  1+r_{4}\right)  ^{2\left(  \left\vert \gamma\right\vert
-\theta\right)  }\left(  \int\limits_{\left\vert \cdot\right\vert \leq r_{4}%
}\frac{1}{w}+\int\limits_{\left\vert \cdot\right\vert \geq r_{4}}%
\frac{\left\vert \cdot\right\vert ^{2\left\vert \gamma\right\vert }%
}{w\left\vert \cdot\right\vert ^{2\theta}}\right) \nonumber\\
& \leq\left(  1+r_{4}\right)  ^{2\underline{\kappa}}\max_{n\leq\underline
{\kappa}}\left(  \int\limits_{\left\vert \cdot\right\vert \leq r_{4}}\frac
{1}{w}+\int\limits_{\left\vert \cdot\right\vert \geq r_{4}}\frac{\left\vert
\cdot\right\vert ^{2n}}{w\left\vert \cdot\right\vert ^{2\theta}}\right)
\nonumber\\
& =\left(  1+r_{4}\right)  ^{2\underline{\kappa}}\left(  \int%
\limits_{\left\vert \cdot\right\vert \leq r_{4}}\frac{1}{w}+\max
_{n\leq\underline{\kappa}}\int\limits_{\left\vert \cdot\right\vert \geq r_{4}%
}\frac{\left\vert \cdot\right\vert ^{2n}}{w\left\vert \cdot\right\vert
^{2\theta}}\right)  .\label{a05}%
\end{align}

Since property W2.1 requires that $1/w\in L_{loc}^{1}$and property W3 implies
Theorem \ref{Thm_property_wt_fn_W3}, it follows that the last two integrals
exist. We now combine both cases by setting%
\begin{equation}
C_{w}=\max\left\{  \max_{m<\theta}C_{m},\max_{m<\theta}C_{\theta,m,r_{4}%
},\left(  1+r_{4}\right)  ^{2\underline{\kappa}}\right\}  \left(
\int\limits_{\left\vert \cdot\right\vert \leq r_{3}}\frac{1}{w}+\max
_{n\leq\underline{\kappa}}\int\limits_{\left\vert \cdot\right\vert \geq r_{3}%
}\frac{\left\vert \cdot\right\vert ^{2n}}{w\left\vert \cdot\right\vert
^{2\theta}}\right) \label{a07}%
\end{equation}

so that inequalities \ref{a39} and \ref{a05} imply \ref{a06}.

The constants $C_{m}$ and $C_{\theta,m,r_{4}}$ are defined in the proof of
Theorem \ref{Thm_bound_on_g(e,x)}.
\end{proof}

The next theorem is our inverse Fourier transform result for distributions in
$X_{w}^{\theta}$ when the weight function has property W2.

\begin{theorem}
Suppose the weight function $w$ has property W2 and suppose $f\in
X_{w}^{\theta}$. Summary \ref{Sum_properties_Xwm_Ch3} allows us to define a.e.
the function $f_{F}:\mathbb{R}^{d}\rightarrow\mathbb{C}$ by $f_{F}=\widehat
{f}$ on $\mathbb{R}^{d}\setminus0$. Then for all multi-indexes $\gamma$,%
\begin{equation}
\left[  \widehat{D^{\gamma}f},\psi\right]  =\int\left(  i\xi\right)  ^{\gamma
}\left(  \mathcal{Q}_{\emptyset,\theta}\psi\right)  \left(  \xi\right)
f_{F}\left(  \xi\right)  d\xi+\left(  2\pi\right)  ^{-\frac{d}{2}}\left[
\widehat{p_{\widehat{D^{\gamma}f}}},\psi\right]  ,\quad\psi\in S,\label{a12}%
\end{equation}

where for $u\in S^{\prime}$, $p_{u}\in P_{\theta-1}$ and%
\begin{equation}
p_{u}\left(  x\right)  =\sum_{\left\vert \alpha\right\vert <\theta}%
\frac{b_{u,\alpha}}{\alpha!}x^{\alpha},\quad b_{u,\alpha}=\left[  u,\left(
-i\xi\right)  ^{\alpha}\rho\right]  .\label{a14}%
\end{equation}

Also for all multi-indexes $\gamma,$%
\begin{equation}
\left[  D^{\gamma}f,\psi\right]  =\int\left(  i\xi\right)  ^{\gamma}\left(
\mathcal{Q}_{\emptyset,\theta}\overset{\vee}{\psi}\right)  \left(  \xi\right)
f_{F}\left(  \xi\right)  d\xi+\left(  2\pi\right)  ^{-\frac{d}{2}}\int
p_{\widehat{D^{\gamma}f}}\psi,\quad\psi\in S,\label{a1.15}%
\end{equation}

and%
\begin{equation}
\left(  D^{\gamma}f\ast\psi\right)  \left(  x\right)  =\left(  2\pi\right)
^{-\frac{d}{2}}\left(  \int\left(  i\xi\right)  ^{\gamma}\left(
\mathcal{Q}_{\emptyset,\theta,\cdot}\left(  e^{i\left(  x,\cdot\right)
}\widehat{\psi}\right)  \right)  \left(  \xi\right)  f_{F}\left(  \xi\right)
d\xi+\left(  p_{\widehat{D^{\gamma}f}}\ast\psi\right)  \left(  x\right)
\right)  ,\quad\psi\in S,\label{a1.16}%
\end{equation}

where $f\ast\psi=\left[  f_{y},\psi\left(  \cdot-y\right)  \right]  $ for
$f\in S^{\prime}$ and $\psi\in S$.
\end{theorem}

\begin{proof}
If $f\in X_{w}^{\theta}$ then by Summary \ref{Sum_properties_Xwm_Ch3} $f\in
S^{\prime}$, $\widehat{f}\in L_{loc}^{1}\left(  \mathbb{R}^{d}\setminus
0\right)  $, $f_{F}\in S_{\emptyset,\theta}^{\prime}$ with action$\int
f_{F}\phi$ and $f_{F}=\widehat{f}$ on $S_{\emptyset,\theta}$. Hence $f$
satisfies the assumptions of part 1 of Theorem \ref{Thm_integ(|DQ(exp)fF|)}
and \ref{a19} and \ref{a20} implies \ref{a12} and \ref{a14}.

If a weight function also has property W3 for order $\theta$ then we have the
following modified inverse-Fourier transform theorem for the (continuous)
functions in $X_{w}^{\theta}$ as well as an upper bound for the growth rate
\ref{a28} of the derivatives near infinity. More precisely, a function in
$X_{w}^{\theta}$ is shown to be the sum of a function $f_{\rho}\in
X_{w}^{\theta}$ and a polynomial in $P_{\theta-1}$ and the function $f_{\rho}$
satisfies the modified inverse-Fourier transform equations \ref{a41}.

To prove \ref{a1.15} replace $\psi$ by $\overset{\vee}{\psi}$ in equation
\ref{a12}. We start with \ref{a1.15} and the convolution definition:%
\begin{align*}
\left(  D^{\gamma}f\ast\psi\right)  \left(  x\right)   & =\left(  2\pi\right)
^{-\frac{d}{2}}\left[  \left(  D^{\gamma}f\right)  _{y},\psi\left(
x-y\right)  \right] \\
& =\left(  2\pi\right)  ^{-\frac{d}{2}}\int\left(  i\xi\right)  ^{\gamma
}\left(  \mathcal{Q}_{\emptyset,\theta,\xi}\left(  F_{y}^{-1}\left[
\psi\left(  x-y\right)  \right]  \left(  \xi\right)  \right)  \right)
f_{F}\left(  \xi\right)  d\xi+\\
& \qquad+\left(  2\pi\right)  ^{-d}\int p_{\widehat{D^{\gamma}f}}\left(
y\right)  \psi\left(  x-y\right)  dy\\
& =\left(  2\pi\right)  ^{-\frac{d}{2}}\left(  \int\left(  i\xi\right)
^{\gamma}\left(  \mathcal{Q}_{\emptyset,\theta,\xi}F_{y}\left[  \psi\left(
x+y\right)  \right]  \left(  \xi\right)  \right)  f_{F}\left(  \xi\right)
d\xi+\left(  p_{\widehat{D^{\gamma}f}}\ast\psi\right)  \left(  x\right)
\right) \\
& =\left(  2\pi\right)  ^{-\frac{d}{2}}\left(  \int\left(  i\xi\right)
^{\gamma}\mathcal{Q}_{\emptyset,\theta,\cdot}\left(  e^{i\left(
x,\cdot\right)  }\widehat{\psi}\right)  \left(  \xi\right)  f_{F}\left(
\xi\right)  d\xi+\left(  p_{\widehat{D^{\gamma}f}}\ast\psi\right)  \left(
x\right)  \right)  .
\end{align*}

\end{proof}

\begin{theorem}
\label{Thm_Dfrho_W3}Now suppose the weight function $w$ also has property W2
or W3 for order $\theta$ and smoothing parameter $\kappa$. Set $\underline
{\kappa}=\min\kappa$. Then $X_{w}^{\theta}\subset C_{BP}^{\left(  \left\lfloor
\underline{\kappa}\right\rfloor \right)  }$. Now let $f_{\rho}=\left(
\mathcal{Q}_{\emptyset,n}^{\ast}\widehat{f}\right)  ^{\vee}$ for $f\in
X_{w}^{\theta}$ so that $f_{\rho}\in X_{w}^{\theta}$, $f_{F}=\left(  f_{\rho
}\right)  _{F}$ and%
\begin{equation}
f=f_{\rho}+\left(  2\pi\right)  ^{-\frac{d}{2}}p_{\widehat{f}},\quad f\in
X_{w}^{\theta},\label{a42}%
\end{equation}

where $p_{\widehat{f}}\in P_{\theta-1}$ is given by \ref{a14} and $\rho\in
S_{1,n}$ is used to define $\mathcal{P}_{\emptyset,\theta}$. Then $f_{\rho}\in
C_{BP}^{\left(  \left\lfloor \underline{\kappa}\right\rfloor \right)  }\cap
X_{w}^{\theta}$ and for $\left\vert \gamma\right\vert \leq\left\lfloor
\underline{\kappa}\right\rfloor $:%
\begin{align}
D^{\gamma}f_{\rho}\left(  x\right)   & =\left(  2\pi\right)  ^{-\frac{d}{2}%
}\int D_{x}^{\gamma}\mathcal{Q}_{\emptyset,\theta,\xi}\left(  e^{i\left(
x,\xi\right)  }\right)  f_{F}\left(  \xi\right)  d\xi\nonumber\\
& =\left\{
\begin{array}
[c]{ll}%
\left(  2\pi\right)  ^{-\frac{d}{2}}\int\mathcal{Q}_{\emptyset,\theta
-\left\vert \gamma\right\vert ,\xi}\left(  e^{i\left(  x,\xi\right)  }\right)
\left(  i\xi\right)  ^{\gamma}f_{F}\left(  \xi\right)  d\xi, & \left\vert
\gamma\right\vert <\theta,\\
\left(  2\pi\right)  ^{-\frac{d}{2}}\int e^{i\left(  x,\xi\right)  }\left(
i\xi\right)  ^{\gamma}f_{F}\left(  \xi\right)  d\xi, & \left\vert
\gamma\right\vert \geq\theta,
\end{array}
\right. \label{a48}\\
& =\left\{
\begin{array}
[c]{ll}%
\left(  2\pi\right)  ^{-\frac{d}{2}}\int\mathcal{Q}_{\emptyset,\theta
-\left\vert \gamma\right\vert ,\xi}\left(  e^{i\left(  x,\xi\right)  }\right)
\left(  D^{\gamma}f_{\rho}\right)  _{F}\left(  \xi\right)  d\xi, & \left\vert
\gamma\right\vert <\theta,\\
\left(  2\pi\right)  ^{-\frac{d}{2}}\int e^{i\left(  x,\xi\right)  }\left(
D^{\gamma}f_{\rho}\right)  _{F}\left(  \xi\right)  d\xi, & \left\vert
\gamma\right\vert \geq\theta,
\end{array}
\right. \label{a41}%
\end{align}

and $D^{\gamma}f_{\rho}$ satisfies the growth estimates%
\begin{equation}
\left\vert D^{\gamma}f_{\rho}\left(  x\right)  \right\vert \leq\left\{
\begin{array}
[c]{ll}%
C_{\rho,w}\left\vert f_{\rho}\right\vert _{w,\theta}\left(  1+\left\vert
x\right\vert \right)  ^{\theta-\left\vert \gamma\right\vert }, & \left\vert
\gamma\right\vert <\theta,\\
C_{\rho,w}\left\vert f_{\rho}\right\vert _{w,\theta}, & \left\vert
\gamma\right\vert \geq\theta,
\end{array}
\right. \label{a28}%
\end{equation}

where the constant $C_{\rho,w}$ is independent of $x$ and $f$.
\end{theorem}

\begin{proof}
This proof is an application of Corollary \ref{Cor_Thm_integ(|DQ(exp)fF|)}
with $n=\theta$ and $m=\left\lfloor \underline{\kappa}\right\rfloor $. If
$f\in X_{w}^{\theta}$ then by Summary \ref{Sum_properties_Xwm_Ch3}, $f\in
S^{\prime}$, $\widehat{f}\in L_{loc}^{1}\left(  \mathbb{R}^{d}\setminus
0\right)  $, $f_{F}\in S_{\emptyset,\theta}^{\prime}$ with action$\int
f_{F}\phi$ and $f_{F}=\widehat{f}$ on $S_{\emptyset,\theta}$ so $f$ satisfies
assumptions 1 and 2 of Corollary \ref{Cor_Thm_integ(|DQ(exp)fF|)}. Regarding
assumption 3: using the Cauchy-Schwartz inequality
\begin{align*}
\int\left\vert \left(  D_{x}^{\gamma}\mathcal{Q}_{\emptyset,\theta,\xi}\left(
e^{i\left(  x,\xi\right)  }\right)  \right)  f_{F}\left(  \xi\right)
\right\vert d\xi & =\int\left\vert \frac{D_{x}^{\gamma}\mathcal{Q}%
_{\emptyset,\theta,\xi}\left(  e^{i\left(  x,\xi\right)  }\right)  }%
{\sqrt{w\left(  \xi\right)  }\left\vert \xi\right\vert ^{\theta}}%
\sqrt{w\left(  \xi\right)  }\left\vert \xi\right\vert ^{\theta}f_{F}\left(
\xi\right)  \right\vert d\xi\\
& \leq\left(  \int\frac{\left\vert D_{x}^{\gamma}\mathcal{Q}_{\emptyset
,\theta,\xi}\left(  e^{i\left(  x,\xi\right)  }\right)  \right\vert ^{2}%
}{w\left(  \xi\right)  \left\vert \xi\right\vert ^{2\theta}}d\xi\right)
^{\frac{1}{2}}\left\vert f\right\vert _{w,\theta},
\end{align*}

since $\left\vert f\right\vert _{w,\theta}^{2}=\int w\left\vert \cdot
\right\vert ^{2\theta}\left\vert f_{F}\right\vert ^{2}d\xi$ by part 1 Summary
\ref{Sum_properties_Xwm_Ch3}. Applying the inequality derived in Lemma
\ref{Lem_integ(DQ(exp)|^2/w|.|_2} we get%
\[
\left(  \int\frac{\left\vert D_{x}^{\gamma}\mathcal{Q}_{\emptyset,\theta,\xi
}\left(  e^{i\left(  x,\xi\right)  }\right)  \right\vert ^{2}}{w\left(
\xi\right)  \left\vert \xi\right\vert ^{2\theta}}d\xi\right)  ^{\frac{1}{2}%
}\leq\left\{
\begin{array}
[c]{ll}%
\max\left\{
\begin{array}
[c]{c}%
\underline{C}_{\theta,r_{4}}^{\left(  \rho\right)  },\\
\overline{C}_{\theta-\left\vert \gamma\right\vert ,r_{4}}^{\left(
\rho\right)  }%
\end{array}
\right\}
\begin{pmatrix}
\sqrt{K_{\gamma,r_{4}}^{\left(  w\right)  }}\\
or\\
\sqrt{L_{\left\vert \gamma\right\vert ,r_{4}}^{\left(  w\right)  }}%
\end{pmatrix}
\left(  1+\sum\limits_{k\leq\theta-\left\vert \gamma\right\vert }%
\frac{\left\vert x\right\vert ^{2k}}{k!}\right)  ^{\frac{1}{2}}, & \left\vert
\gamma\right\vert <\theta,\\
\left(  \int\frac{\xi^{2\gamma}}{w\left\vert \cdot\right\vert ^{2\theta}%
}\right)  ^{1/2}\text{ }or\text{ }\left(  \int\frac{\left\vert \cdot
\right\vert ^{2\left\vert \gamma\right\vert }}{w\left\vert \cdot\right\vert
^{2\theta}}\right)  ^{1/2}, & \left\vert \gamma\right\vert \geq\theta.
\end{array}
\right.
\]

There now exists a constant $C_{\rho,w}$, independent of $x$, such that
inequality \ref{a32} of Corollary \ref{Cor_Thm_integ(|DQ(exp)fF|)} is
satisfied for $k_{\gamma}=C_{\rho,w}\left\vert f\right\vert _{w,\theta}$ and
$s_{\gamma}=\max\left\{  0,\theta-\left\vert \gamma\right\vert \right\}  $ and
so the rest of Corollary \ref{Cor_Thm_integ(|DQ(exp)fF|)} can be applied to
obtain $f\in C_{BP}^{\left(  m\right)  }$ and all equations and estimates for
$D^{\gamma}f_{\rho}$ except equation \ref{a41} where $f_{\rho}$ has replaced
$f$. But by \ref{a42}, $D^{\gamma}f=D^{\gamma}f_{\rho}+\left(  2\pi\right)
^{-\frac{d}{2}}D^{\gamma}p_{\widehat{f}}$ and since $D^{\gamma}p_{\widehat{f}%
}$ is a polynomial we have $\left(  D^{\gamma}f_{\rho}\right)  _{F}=\left(
D^{\gamma}f\right)  _{F}$.
\end{proof}

\begin{remark}
This result is closely related to Proposition 2.2 of Madych and Nelson
\cite{MadychNelson90}. However, in the comments following Theorem 2.1 Madych
and Nelson illustrate Theorem 2.1 by choosing $d\mu\left(  \xi\right)
=w\left(  \xi\right)  d\xi$ where $w$ corresponds to $1/w$ in this document
i.e. to the \textbf{reciprocal} of our weight function.
\end{remark}

\section{Taylor series expansions of data functions in $X_{w}^{\theta}$: $w\in
W3.2$\label{Sect_data_fn_Taylor_W3.2}}

We begin with the tempered distribution Taylor series expansion \ref{a1.55}:%
\[
f\left(  \cdot+a\right)  -\sum_{\left\vert \beta\right\vert \leq n}%
\frac{a^{\beta}}{\beta!}D^{\beta}f=\frac{\sqrt{2\pi}}{n!}\left(  \left(
ia\xi\right)  ^{n+1}\overline{\widehat{g_{n}}}\left(  a\xi\right)  \widehat
{f}\right)  ^{\vee},\quad f\in S^{\prime}.
\]

We now want to\ estimate the remainder of the Taylor series expansion when
$f\in X_{w}^{\theta}$, $w\in W3.2$ for order $\theta$ and smoothness
$\kappa\in\mathbb{R}^{1}$. This makes sense since from Theorem
\ref{Thm_Dfrho_W3} we have $X_{w}^{\theta}\subset C_{BP}^{\left(  \left\lfloor
\kappa\right\rfloor \right)  }$ and accordingly we choose%
\begin{equation}
n\leq\left\lfloor \kappa\right\rfloor .\label{a18}%
\end{equation}

As in Theorem \ref{Thm_Dfrho_W3}, set $f_{\rho}=\left(  \mathcal{Q}%
_{\emptyset,\theta}^{\ast}\widehat{f}\right)  ^{\vee}$ so that
\begin{align}
f_{\rho}  & =f-\left(  2\pi\right)  ^{-\frac{d}{2}}p_{\theta-1;\widehat{f}%
}=\left(  2\pi\right)  ^{-\frac{d}{2}}\int\mathcal{Q}_{\emptyset,\theta,\xi
}\left(  e^{i\left(  \cdot,\xi\right)  }\right)  f_{F}\left(  \xi\right)
d\xi,\label{a920}\\
& where\nonumber\\
p_{\theta-1;u}\left(  \xi\right)   & =\sum_{\left\vert \beta\right\vert
<\theta}\dfrac{b_{\beta;u}}{\beta!}\xi^{\beta}\in P_{\theta-1},\quad
b_{\beta;u}=\left[  u,\left(  -ix\right)  ^{\beta}\rho\right]  ,\label{a92}%
\end{align}

and hence%
\begin{equation}
f=f_{\rho}+\left(  2\pi\right)  ^{-\frac{d}{2}}p_{\theta-1;\widehat{f}}%
,\quad\deg p_{\theta-1;\widehat{f}}<\theta,\quad\left(  f_{\rho}\right)
_{F}=f_{F}\text{ }on\text{ }S_{\emptyset,\theta}.\label{a88}%
\end{equation}

We have added a $\theta$ to the subscript of $p_{\theta-1;\widehat{f}}$ for
later convenience.

Write for $f\in X_{w}^{\theta}$:%
\begin{align}
& f\left(  \cdot+a\right)  -\sum_{\left\vert \beta\right\vert \leq n}%
\frac{a^{\beta}}{\beta!}D^{\beta}f\nonumber\\
& =\left(  f_{\rho}\left(  \cdot+a\right)  -\sum_{\left\vert \beta\right\vert
\leq n}\frac{a^{\beta}}{\beta!}D^{\beta}f_{\rho}\right)  +\left(  2\pi\right)
^{-\frac{d}{2}}\left(  p_{\theta-1;\widehat{f}}\left(  \cdot+a\right)
-\sum_{\left\vert \beta\right\vert \leq n}\frac{a^{\beta}}{\beta!}D^{\beta
}p_{\theta-1;\widehat{f}}\right) \nonumber\\
& =\frac{\sqrt{2\pi}}{n!}\left(  \left(  ia\xi\right)  ^{n+1}\overline
{\widehat{g_{n}}}\left(  a\xi\right)  \widehat{f_{\rho}}\right)  ^{\vee
}+\left(  2\pi\right)  ^{-\frac{d}{2}}\left(  p_{\theta-1;\widehat{f}}\left(
\cdot+a\right)  -\sum_{\left\vert \beta\right\vert \leq n}\frac{a^{\beta}%
}{\beta!}D^{\beta}p_{\theta-1;\widehat{f}}\right)  .\label{a93}%
\end{align}
\medskip

\underline{\textbf{Estimation of the polynomial term in }\ref{a93}}\smallskip

\textbf{It will not be necessary to use this estimate in the applications to
interpolants and smoothers}.

Since
\[
p_{\theta-1;\widehat{f}}\left(  \cdot+a\right)  =\sum_{\left\vert
\beta\right\vert \leq\deg p_{\theta-1;\widehat{f}}}\frac{a^{\beta}}{\beta
!}D^{\beta}p_{\theta-1;\widehat{f}}=\sum_{\left\vert \beta\right\vert
\leq\theta-1}\frac{a^{\beta}}{\beta!}D^{\beta}p_{\theta-1;\widehat{f}},
\]

we have:\medskip

\fbox{If $n\geq\theta-1$} then $p_{\theta-1;\widehat{f}}\left(  \cdot
+a\right)  -\sum\limits_{\left\vert \beta\right\vert \leq n}\frac{a^{\beta}%
}{\beta!}D^{\beta}p_{\theta-1;\widehat{f}}=0$.\smallskip

?? \textbf{DEL}?

\fbox{If $n<\theta-1$,} from \ref{Ap147} and then part 3 of Theorem
\ref{Thm_Adel_poly_S1n},%
\begin{align*}
\left\vert p_{\widehat{f}}\left(  x+a\right)  -\sum_{\left\vert \beta
\right\vert \leq n}\frac{a^{\beta}}{\beta!}D^{\beta}p_{\theta-1;\widehat{f}%
}\left(  x\right)  \right\vert  & \leq\frac{\left\vert a\right\vert ^{n+1}%
}{\left(  n+1\right)  !}\left\Vert \left(  \widehat{a}D\right)  ^{n+1}%
p_{\theta-1;\widehat{f}}\right\Vert _{\infty;\left[  x,x+a\right]  }\\
& =\frac{\left\vert a\right\vert ^{n+1}}{\left(  n+1\right)  !}\left\Vert
p_{\theta-n-2;\widehat{v}}\right\Vert _{\infty;\left[  x,x+a\right]  },
\end{align*}

where
\begin{equation}
v=\left(  -i\widehat{a}D\right)  ^{n+1}f,\quad\widehat{v}=\left(  \widehat
{a}\xi\right)  ^{n+1}\widehat{f}.\label{a98}%
\end{equation}

Hence%
\begin{align}
&  \left\vert \left(  2\pi\right)  ^{-\frac{d}{2}}\left(  p_{\widehat{f}%
}\left(  x+a\right)  -\sum_{\left\vert \beta\right\vert \leq n}\frac{a^{\beta
}}{\beta!}D^{\beta}p_{\theta-1;\widehat{f}}\left(  x\right)  \right)
\right\vert \nonumber\\
&  \leq\left\{
\begin{array}
[c]{ll}%
\left(  2\pi\right)  ^{-\frac{d}{2}}\frac{\left\vert a\right\vert ^{n+1}%
}{\left(  n+1\right)  !}\left\Vert p_{\theta-n-2;\widehat{v}}\right\Vert
_{\infty;\left[  x,x+a\right]  }, & n<\theta-1,\\
0, & n\geq\theta-1.
\end{array}
\right. \label{a94}%
\end{align}

From \ref{a91},
\[
p_{\theta-n-1;u}\left(  \xi\right)  =\sum_{\left\vert \alpha\right\vert
<\theta-n-1}\dfrac{b_{\alpha,u}}{\alpha!}\xi^{\alpha}\in P_{\theta-n-2},\quad
b_{\alpha,u}=\left[  u,\left(  -ix\right)  ^{\alpha}\rho\right]  ,
\]

so%
\[
p_{\theta-n-1;\widehat{v}}\left(  \xi\right)  =\sum_{\left\vert \alpha
\right\vert <\theta-n-1}\dfrac{b_{\alpha,\widehat{v}}}{\alpha!}\xi^{\alpha
},\quad b_{\alpha,\widehat{v}}=\left[  \left(  -i\widehat{a}\eta\right)
^{n+1}\widehat{f},\left(  -i\eta\right)  ^{\alpha}\rho\right]  ,
\]

and if
\[
r=\max\left\{  \left\vert x\right\vert ,\left\vert x+a\right\vert \right\}
=\left(  \max\left\{  \left\vert x\right\vert ^{2},\left\vert x+a\right\vert
^{2}\right\}  \right)  ^{1/2},
\]

then%
\begin{align*}
\left\Vert p_{\theta-n-1;\widehat{v}}\right\Vert _{\infty;\left[
x,x+a\right]  }  & \leq\sum_{\left\vert \alpha\right\vert <\theta-n-1}%
\dfrac{\left\vert b_{\alpha,\widehat{v}}\right\vert }{\alpha!}\left\Vert
\xi^{\alpha}\right\Vert _{\infty;\left[  x,x+a\right]  }\leq\sum_{\left\vert
\alpha\right\vert <\theta-n-1}\dfrac{\left\vert b_{\alpha,\widehat{v}%
}\right\vert }{\alpha!}\left\Vert \left\vert \xi\right\vert ^{\left\vert
\alpha\right\vert }\right\Vert _{\infty;\left[  x,x+a\right]  }\leq\\
& \leq\sum_{\left\vert \alpha\right\vert <\theta-n-1}\dfrac{\left\vert
b_{\alpha,\widehat{v}}\right\vert }{\alpha!}r^{\left\vert \alpha\right\vert
}\leq\left(  \max_{\left\vert \alpha\right\vert <\theta-n-1}\left\vert
b_{\alpha,\widehat{v}}\right\vert \right)  \sum_{\left\vert \alpha\right\vert
<\theta-n-1}\dfrac{r^{\left\vert \alpha\right\vert }}{\alpha!}=\\
& =\left(  \max_{\left\vert \alpha\right\vert <\theta-n-1}\left\vert
b_{\alpha,\widehat{v}}\right\vert \right)  \sum\limits_{k=0}^{\theta-n-1}%
\sum_{\left\vert \alpha\right\vert =k}\dfrac{r^{\left\vert \alpha\right\vert
}}{\alpha!}\\
& =\left(  \max_{\left\vert \alpha\right\vert <\theta-n-1}\left\vert
b_{\alpha,\widehat{v}}\right\vert \right)  \sum\limits_{k=0}^{\theta-n-1}%
r^{k}\sum_{\left\vert \alpha\right\vert =k}\dfrac{\mathbf{1}^{\alpha
}\mathbf{1}^{\alpha}}{\alpha!}\\
& =\left(  \max_{\left\vert \alpha\right\vert <\theta-n-1}\left\vert
b_{\alpha,\widehat{v}}\right\vert \right)  \sum\limits_{k=0}^{\theta-n-1}%
\frac{r^{k}}{k!}\left(  \mathbf{1},\mathbf{1}\right)  ^{k}\\
& =\left(  \max_{\left\vert \alpha\right\vert <\theta-n-1}\left\vert
b_{\alpha,\widehat{v}}\right\vert \right)  \sum\limits_{k=0}^{\theta-n-1}%
\frac{\left(  rd\right)  ^{k}}{k!}.
\end{align*}

\underline{\textbf{Estimating the }$f_{\rho}$\textbf{\ term in }%
\ref{a93}\textbf{\ when ?? }$n<\theta$\textbf{, }$n\leq\left\lfloor
\kappa\right\rfloor $} We introduce the partition of unity%
\begin{equation}
1=\phi_{0}+\phi_{\infty},\text{\quad}\phi_{0}\in S_{1;\theta},\text{\quad
}0\leq\phi_{0}\leq1,\label{a15}%
\end{equation}

which implies that%
\begin{equation}
\phi_{\infty}\in C_{\emptyset,\theta}^{\infty}\cap C_{B}^{\infty},\text{\quad
}0\leq\phi_{\infty}\leq1.\label{a150}%
\end{equation}

Observe that from the remainder estimates of Theorem \ref{Thm_Tay_rem_zeros},
for any $r>0$,%
\begin{equation}
\left\vert \phi_{\infty}\left(  x\right)  \right\vert \leq\frac{1}{\theta
!}\left\vert x\right\vert ^{\theta}\left\Vert \left(  \widehat{\cdot}D\right)
^{\theta}\phi_{\infty}\right\Vert _{\infty;\leq r},\text{\quad}x\in
\overline{B}_{r}.\label{a11}%
\end{equation}

Also, since $\phi_{\infty}\in C_{\emptyset,\theta}^{\infty}\cap C_{B}^{\infty
}$ implies $\phi_{\infty}\phi\in S_{\emptyset,\theta}$ when $\phi\in S$, we
have
\begin{align*}
\widehat{f_{\rho}}  & =\phi_{0}\widehat{f_{\rho}}+\phi_{\infty}\widehat
{f_{\rho}}=\phi_{0}\widehat{f_{\rho}}+\phi_{\infty}f_{F},\\
f_{\rho}  & =\overset{\vee}{\phi}_{0}\ast f_{\rho}+\overset{\vee}{\phi
_{\infty}}\ast f_{\rho},
\end{align*}

and so%
\begin{align}
& f_{\rho}\left(  \cdot+a\right)  -\sum_{\left\vert \beta\right\vert \leq
n}\frac{a^{\beta}}{\beta!}D^{\beta}f_{\rho}\nonumber\\
& =\frac{\sqrt{2\pi}}{n!}\left(  \left(  ia\xi\right)  ^{n+1}\overline
{\widehat{g_{n}}}\left(  a\xi\right)  \widehat{f_{\rho}}\right)  ^{\vee
}\nonumber\\
& =\frac{\sqrt{2\pi}}{n!}\left(  \left(  ia\xi\right)  ^{n+1}\overline
{\widehat{g_{n}}}\left(  a\xi\right)  \phi_{0}\widehat{f_{\rho}}\right)
^{\vee}+\frac{\sqrt{2\pi}}{n!}\left(  \left(  ia\xi\right)  ^{n+1}%
\overline{\widehat{g_{n}}}\left(  a\xi\right)  \phi_{\infty}\widehat{f_{\rho}%
}\right)  ^{\vee}\nonumber\\
& =\frac{\sqrt{2\pi}}{n!}\left(  \left(  ia\xi\right)  ^{n+1}\overline
{\widehat{g_{n}}}\left(  a\xi\right)  \phi_{0}\widehat{f_{\rho}}\right)
^{\vee}+\frac{\sqrt{2\pi}}{n!}\left(  \left(  ia\xi\right)  ^{n+1}%
\overline{\widehat{g_{n}}}\left(  a\xi\right)  \phi_{\infty}f_{F}\right)
^{\vee}\nonumber\\
& =\frac{\sqrt{2\pi}}{n!}\left(  \overline{\widehat{g_{n}}}\left(
a\xi\right)  \left(  ia\xi\right)  \phi_{0}\widehat{\left(  aD\right)
^{n}f_{\rho}}\right)  ^{\vee}+\frac{\sqrt{2\pi}}{n!}\left(  \left(
ia\xi\right)  ^{n+1}\overline{\widehat{g_{n}}}\left(  a\xi\right)
\phi_{\infty}f_{F}\right)  ^{\vee}.\label{a35}%
\end{align}

\underline{\textbf{Estimating the }$\left(  aD\right)  ^{n}f_{\rho}%
$\textbf{\ term on the RHS of} \ref{a35}} From \ref{a48},%
\[
D^{\gamma}f_{\rho}\left(  x\right)  =\left\{
\begin{array}
[c]{ll}%
\left(  2\pi\right)  ^{-\frac{d}{2}}\int\mathcal{Q}_{\emptyset,\theta
-\left\vert \gamma\right\vert ,\xi}\left(  e^{ix\xi}\right)  \left(
i\xi\right)  ^{\gamma}f_{F}\left(  \xi\right)  d\xi, & \left\vert
\gamma\right\vert <\theta,\text{ }\left\vert \gamma\right\vert \leq
\left\lfloor \kappa\right\rfloor ,\\
\left(  2\pi\right)  ^{-\frac{d}{2}}\int e^{ix\xi}\left(  i\xi\right)
^{\gamma}f_{F}\left(  \xi\right)  d\xi, & \theta\leq\left\vert \gamma
\right\vert \leq\left\lfloor \kappa\right\rfloor ,
\end{array}
\right.
\]

so that%
\begin{equation}
\left(  aD\right)  ^{n}f_{\rho}\left(  x\right)  =\left\{
\begin{array}
[c]{ll}%
\left(  2\pi\right)  ^{-\frac{d}{2}}\int\mathcal{Q}_{\emptyset,\theta-n,\xi
}\left(  e^{ix\xi}\right)  \left(  ia\xi\right)  ^{n}f_{F}\left(  \xi\right)
d\xi, & n<\theta,\text{ }n\leq\left\lfloor \kappa\right\rfloor ,\\
\left(  2\pi\right)  ^{-\frac{d}{2}}\int e^{ix\xi}\left(  ia\xi\right)
^{n}f_{F}\left(  \xi\right)  d\xi, & \theta\leq n\leq\left\lfloor
\kappa\right\rfloor .
\end{array}
\right. \label{a49}%
\end{equation}

From Lemma \ref{Lem_gm_properties_2} $\overline{\widehat{g_{n}}}\in
C_{B}^{\infty}\left(  \mathbb{R}^{1}\right)  $ which means that $\overline
{\widehat{g_{n}}}\left(  \left(  a,\cdot\right)  \right)  \in C_{B}^{\infty
}\left(  \mathbb{R}^{d}\right)  $. Define%
\begin{equation}
\left.
\begin{array}
[c]{l}%
\sigma_{a}\left(  \xi\right)  :=\frac{\sqrt{2\pi}}{n!}\overline{\widehat
{g_{n}}}\left(  a\xi\right)  \left(  ia\xi\right)  \phi_{0}\left(  \xi\right)
\in S,\\
\phi_{a}\left(  \xi\right)  :=\left(  ia\xi\right)  \phi_{0}\left(
\xi\right)  \in S,\\
\widehat{\phi_{a}}=aD\widehat{\phi_{0}},\text{ }\phi_{a}=\left\vert
a\right\vert \phi_{\widehat{a}},\quad\widehat{a}=a/\left\vert a\right\vert .
\end{array}
\right\} \label{a51}%
\end{equation}

Since $\left(  aD\right)  ^{n}f_{\rho}\in C_{BP}^{\left(  0\right)  }$ when
$n\leq\left\lfloor \kappa\right\rfloor $ we can apply \ref{a800} to obtain the
inequality%
\begin{align}
& \left\vert \frac{\sqrt{2\pi}}{n!}\left(  \overline{\widehat{g_{n}}}\left(
a\xi\right)  \left(  ia\xi\right)  \phi_{0}\widehat{\left(  aD\right)
^{n}f_{\rho}}\right)  ^{\vee}\left(  x\right)  \right\vert \nonumber\\
& \leq\left(  2\pi\right)  ^{-\frac{d}{2}}\frac{\left\vert a\right\vert
^{n+1}}{n!}\int_{0}^{1}\int\left\vert \widehat{\phi_{\widehat{a}}}\left(
\xi\right)  \right\vert \left\vert \left(  \left(  \widehat{a}D\right)
^{n}f_{\rho}\right)  \left(  \xi+x+sa\right)  \right\vert d\xi\text{ }%
g_{n}\left(  s\right)  ds.\label{a80}%
\end{align}

The next step is to estimate $\left\vert \left(  \widehat{a}D\right)
^{n}f_{\rho}\left(  x\right)  \right\vert $ and there are two cases:\medskip

\fbox{\textbf{Case 1} $n\leq\left\lfloor \kappa\right\rfloor \leq\theta$.}
From \ref{a49},%
\begin{equation}
\left\vert \left(  \widehat{a}D\right)  ^{n}f_{\rho}\left(  x\right)
\right\vert \leq\left(  2\pi\right)  ^{-\frac{d}{2}}\int\left\vert
\mathcal{Q}_{\emptyset,\theta-n,\xi}\left(  e^{ix\xi}\right)  \right\vert
\left\vert \widehat{a}\xi\right\vert ^{n}\left\vert f_{F}\left(  \xi\right)
\right\vert d\xi.\label{a16}%
\end{equation}

Hence%
\begin{align}
\left\vert \left(  \widehat{a}D\right)  ^{n}f_{\rho}\left(  x\right)
\right\vert  &  \leq\left(  2\pi\right)  ^{-\frac{d}{2}}\int\left\vert
\mathcal{Q}_{\emptyset,\theta-n,\xi}\left(  e^{ix\xi}\right)  \right\vert
\frac{\left\vert \widehat{a}\xi\right\vert ^{n}}{\sqrt{w}\left\vert
\cdot\right\vert ^{\theta}}\sqrt{w}\left\vert \cdot\right\vert ^{\theta
}\left\vert f_{F}\right\vert \nonumber\\
&  \leq\left(  2\pi\right)  ^{-\frac{d}{2}}\left(  \int\left\vert
\mathcal{Q}_{\emptyset,\theta-n,\xi}\left(  e^{ix\xi}\right)  \right\vert
^{2}\frac{\left\vert \widehat{a}\xi\right\vert ^{2n}}{w\left\vert
\cdot\right\vert ^{2\theta}}\right)  ^{1/2}\left\Vert \sqrt{w}\left\vert
\cdot\right\vert ^{\theta}f_{F}\right\Vert _{2}\nonumber\\
&  =\left(  \int\left\vert \mathcal{Q}_{\emptyset,\theta-n,\xi}\left(
e^{ix\xi}\right)  \right\vert ^{2}\frac{\left\vert \widehat{a}\xi\right\vert
^{2n}}{w\left\vert \cdot\right\vert ^{2\theta}}\right)  ^{1/2}\left\vert
f\right\vert _{w,\theta}\nonumber\\
&  =\left(  \int\limits_{\left\vert \xi\right\vert \leq r_{3}}\left\vert
\mathcal{Q}_{\emptyset,\theta-n,\xi}\left(  e^{ix\xi}\right)  \right\vert
^{2}\frac{\left\vert \widehat{a}\xi\right\vert ^{2n}}{w\left\vert
\cdot\right\vert ^{2\theta}}+\int\limits_{\left\vert \xi\right\vert \geq
r_{3}}\left\vert \mathcal{Q}_{\emptyset,\theta-n,\xi}\left(  e^{ix\xi}\right)
\right\vert ^{2}\frac{\left\vert \widehat{a}\xi\right\vert ^{2n}}{w\left\vert
\cdot\right\vert ^{2\theta}}\right)  ^{1/2}\left\vert f\right\vert _{w,\theta
}.\label{a1.00}%
\end{align}

From \ref{a702} with $\gamma=0$ and $n=\theta-n$,%
\begin{align}
& \int\limits_{\left\vert \cdot\right\vert \leq r_{3}}\left\vert
\mathcal{Q}_{\emptyset,\theta-n,\xi}\left(  e^{ix\xi}\right)  \right\vert
^{2}\frac{\left\vert \widehat{a}\xi\right\vert ^{2n}}{w\left\vert
\cdot\right\vert ^{2\theta}}\nonumber\\
& \leq\int\limits_{\left\vert \cdot\right\vert \leq r_{3}}\left(
\underline{C}_{\theta-n,r_{3}}^{\left(  \rho\right)  }\left\vert
\xi\right\vert ^{\theta-n}\left(  1+\frac{\left\vert x\widehat{\xi}\right\vert
^{\theta-n}}{\left(  \theta-n\right)  !}\right)  \right)  ^{2}\frac{\left\vert
\widehat{a}\xi\right\vert ^{2n}}{w\left\vert \cdot\right\vert ^{2\theta}%
}\nonumber\\
& =\left(  \underline{C}_{\theta-n,r_{3}}^{\left(  \rho\right)  }\right)
^{2}\int\limits_{\left\vert \cdot\right\vert \leq r_{3}}\left\vert
\xi\right\vert ^{2\theta-2n}\left(  1+\frac{\left\vert x\widehat{\xi
}\right\vert ^{\theta-n}}{\left(  \theta-n\right)  !}\right)  ^{2}%
\frac{\left\vert \widehat{a}\xi\right\vert ^{2n}}{w\left\vert \cdot\right\vert
^{2\theta}}\nonumber\\
& =\left(  \underline{C}_{\theta-n,r_{3}}^{\left(  \rho\right)  }\right)
^{2}\int\limits_{\left\vert \cdot\right\vert \leq r_{3}}\left(  1+\frac
{\left\vert x\widehat{\xi}\right\vert ^{\theta-n}}{\left(  \theta-n\right)
!}\right)  ^{2}\frac{\left\vert \widehat{a}\widehat{\xi}\right\vert ^{2n}}%
{w}\nonumber\\
& =\left(  \underline{C}_{\theta-n,r_{3}}^{\left(  \rho\right)  }\right)
^{2}\int\limits_{\left\vert \cdot\right\vert \leq r_{3}}\left(  1+\frac
{\left\vert x\right\vert ^{\theta-n}}{\sqrt{\left(  \theta-n\right)  !}}%
\frac{\left\vert \widehat{x}\widehat{\xi}\right\vert ^{\theta-n}}%
{\sqrt{\left(  \theta-n\right)  !}}\right)  ^{2}\frac{\left\vert \widehat
{a}\widehat{\xi}\right\vert ^{2n}}{w}\nonumber\\
& \Rightarrow Cauchy-Schwartz\text{ }inequality\Rightarrow\nonumber\\
& =\left(  \underline{C}_{\theta-n,r_{3}}^{\left(  \rho\right)  }\right)
^{2}\int\limits_{\left\vert \cdot\right\vert \leq r_{3}}\left(  1+\frac
{\left\vert \widehat{x}\widehat{\xi}\right\vert ^{2\left(  \theta-n\right)  }%
}{\left(  \theta-n\right)  !}\right)  \frac{\left\vert \widehat{a}\widehat
{\xi}\right\vert ^{2n}}{w}\left(  1+\frac{\left\vert x\right\vert ^{2\left(
\theta-n\right)  }}{\left(  \theta-n\right)  !}\right) \nonumber\\
& =\left(  \underline{C}_{\theta-n,r_{3}}^{\left(  \rho\right)  }\right)
^{2}\left(  \int\limits_{\left\vert \cdot\right\vert \leq r_{3}}%
\frac{\left\vert \widehat{a}\widehat{\xi}\right\vert ^{2n}}{w}+\frac
{1}{\left(  \theta-n\right)  !}\int\limits_{\left\vert \cdot\right\vert \leq
r_{3}}\frac{\left\vert \widehat{x}\widehat{\xi}\right\vert ^{2\left(
\theta-n\right)  }\left\vert \widehat{a}\widehat{\xi}\right\vert ^{2n}}%
{w}\right)  \left(  1+\frac{\left\vert x\right\vert ^{2\left(  \theta
-n\right)  }}{\left(  \theta-n\right)  !}\right) \label{a611}\\
& \leq\left(  \underline{C}_{\theta-n,r_{3}}^{\left(  \rho\right)  }\right)
^{2}\left(  \int\limits_{\left\vert \cdot\right\vert \leq r_{3}}\frac{1}%
{w}+\frac{1}{\left(  \theta-n\right)  !}\int\limits_{\left\vert \cdot
\right\vert \leq r_{3}}\frac{1}{w}\right)  \left(  1+\frac{\left\vert
x\right\vert ^{2\left(  \theta-n\right)  }}{\left(  \theta-n\right)  !}\right)
\nonumber\\
& =\left(  \underline{C}_{\theta-n,r_{3}}^{\left(  \rho\right)  }\right)
^{2}\left(  1+\frac{1}{\left(  \theta-n\right)  !}\right)  \left(
\int\limits_{\left\vert \cdot\right\vert \leq r_{3}}\frac{1}{w}\right)
\left(  1+\frac{\left\vert x\right\vert ^{2\left(  \theta-n\right)  }}{\left(
\theta-n\right)  !}\right) \label{a60}\\
& <\infty,\nonumber
\end{align}

since $w\in W2.1$. Next from \ref{a27.2}, $\left\vert \mathcal{Q}%
_{\emptyset,\theta-n,\xi}\left(  e^{ix\xi}\right)  \right\vert \leq
\overline{C}_{\theta-n,r}^{\left(  \rho\right)  }\left(  1+\sum
\limits_{k<\theta-n}\frac{\left\vert x\widehat{\xi}\right\vert ^{k}}%
{k!}\right)  $, and so%
\begin{align}
\int\limits_{\left\vert \cdot\right\vert \geq r_{3}} &  \left\vert
\mathcal{Q}_{\emptyset,\theta-n,\xi}\left(  e^{ix\xi}\right)  \right\vert
^{2}\frac{\left\vert \widehat{a}\xi\right\vert ^{2n}}{w\left\vert
\cdot\right\vert ^{2\theta}}\nonumber\\
&  \leq\int\limits_{\left\vert \cdot\right\vert \geq r_{3}}\left(
\overline{C}_{\theta-n,r_{3}}^{\left(  \rho\right)  }\left(  1+\sum
\limits_{k<\theta-n}\frac{\left\vert x\widehat{\xi}\right\vert ^{k}}%
{k!}\right)  \right)  ^{2}\frac{\left\vert \widehat{a}\xi\right\vert ^{2n}%
}{w\left\vert \cdot\right\vert ^{2\theta}}\nonumber\\
&  =\left(  \overline{C}_{\theta-n,r_{3}}^{\left(  \rho\right)  }\right)
^{2}\int\limits_{\left\vert \cdot\right\vert \geq r_{3}}\left(  1+\sum
\limits_{k<\theta-n}\frac{\left\vert x\widehat{\xi}\right\vert ^{k}}%
{k!}\right)  ^{2}\frac{\left\vert \widehat{a}\xi\right\vert ^{2n}}{w\left\vert
\cdot\right\vert ^{2\theta}}\nonumber\\
&  =\left(  \overline{C}_{\theta-n,r}^{\left(  \rho\right)  }\right)  ^{2}%
\int\limits_{\left\vert \cdot\right\vert \geq r_{3}}\left(  1+\sum
\limits_{k<\theta-n}\frac{\left\vert x\right\vert ^{k}}{\sqrt{k!}}%
\frac{\left\vert \widehat{x}\widehat{\xi}\right\vert ^{k}}{\sqrt{k!}}\right)
^{2}\frac{\left\vert \widehat{a}\xi\right\vert ^{2n}}{w\left\vert
\cdot\right\vert ^{2\theta}}\nonumber\\
&  \Rightarrow Cauchy-Schwartz\text{ }inequality\Rightarrow\nonumber\\
&  \leq\left(  \overline{C}_{\theta-n,r_{3}}^{\left(  \rho\right)  }\right)
^{2}\int\limits_{\left\vert \cdot\right\vert \geq r_{3}}\left(  1+\sum
\limits_{k<\theta-n}\frac{\left\vert x\right\vert ^{2k}}{k!}\right)  \left(
1+\sum\limits_{k<\theta-n}\frac{\left\vert \widehat{x}\widehat{\xi}\right\vert
^{2k}}{k!}\right)  \frac{\left\vert \widehat{a}\xi\right\vert ^{2n}%
}{w\left\vert \cdot\right\vert ^{2\theta}}\nonumber\\
&  =\left(  \overline{C}_{\theta-n,r_{3}}^{\left(  \rho\right)  }\right)
^{2}\left(  \int\limits_{\left\vert \cdot\right\vert \geq r_{3}}\left(
1+\sum\limits_{k<\theta-n}\frac{\left\vert \widehat{x}\widehat{\xi}\right\vert
^{2k}}{k!}\right)  \frac{\left\vert \widehat{a}\xi\right\vert ^{2n}%
}{w\left\vert \cdot\right\vert ^{2\theta}}\right)  \left(  1+\sum
\limits_{k<\theta-n}\frac{\left\vert x\right\vert ^{2k}}{k!}\right)
\nonumber\\
&  =\left(  \overline{C}_{\theta-n,r_{3}}^{\left(  \rho\right)  }\right)
^{2}\left(  \int\limits_{\left\vert \cdot\right\vert \geq r_{3}}\left(
1+\sum\limits_{k<\theta-n}\frac{\left\vert \widehat{x}\widehat{\xi}\right\vert
^{2k}}{k!}\right)  \left\vert \widehat{a}\widehat{\xi}\right\vert ^{2n}%
\frac{\left\vert \cdot\right\vert ^{2n}}{w\left\vert \cdot\right\vert
^{2\theta}}\right)  \left(  1+\sum\limits_{k<\theta-n}\frac{\left\vert
x\right\vert ^{2k}}{k!}\right) \nonumber\\
&  =\left(  \overline{C}_{\theta-n,r_{3}}^{\left(  \rho\right)  }\right)
^{2}\left(
\begin{array}
[c]{c}%
\int\limits_{\left\vert \cdot\right\vert \geq r_{3}}\left\vert \widehat
{a}\widehat{\xi}\right\vert ^{2n}\frac{\left\vert \cdot\right\vert ^{2n}%
}{w\left\vert \cdot\right\vert ^{2\theta}}+\\
+\sum\limits_{k<\theta-n}\frac{1}{k!}\int\limits_{\left\vert \cdot\right\vert
\geq r_{3}}\left\vert \widehat{x}\widehat{\xi}\right\vert ^{2k}\left\vert
\widehat{a}\widehat{\xi}\right\vert ^{2n}\frac{\left\vert \cdot\right\vert
^{2n}}{w\left\vert \cdot\right\vert ^{2\theta}}%
\end{array}
\right)  \left(  1+\sum\limits_{k<\theta-n}\frac{\left\vert x\right\vert
^{2k}}{k!}\right) \label{a64}\\
&  <\left(  \overline{C}_{\theta-n,r_{3}}^{\left(  \rho\right)  }\right)
^{2}\left(  1+\sum\limits_{k<\theta-n}\frac{1}{k!}\right)  \left(
\int\limits_{\left\vert \cdot\right\vert \geq r_{3}}\frac{\left\vert
\cdot\right\vert ^{2n}}{w\left\vert \cdot\right\vert ^{2\theta}}\right)
\left(  1+\sum\limits_{k<\theta-n}\frac{\left\vert x\right\vert ^{2k}}%
{k!}\right)  .\label{a63}%
\end{align}

From \ref{a63},%
\begin{align}
\int\limits_{\left\vert \cdot\right\vert \geq r_{3}} &  \left\vert
\mathcal{Q}_{\emptyset,\theta-n,\xi}\left(  e^{ix\xi}\right)  \right\vert
^{2}\frac{\left\vert \widehat{a}\xi\right\vert ^{2n}}{w\left\vert
\cdot\right\vert ^{2\theta}}\nonumber\\
&  \leq\left(  \overline{C}_{\theta-n,r_{3}}^{\left(  \rho\right)  }\right)
^{2}\left(  1+e_{\theta-n-1}\right)  \left(  \int\limits_{\left\vert
\cdot\right\vert \geq r_{3}}\frac{\left\vert \cdot\right\vert ^{2n}%
}{w\left\vert \cdot\right\vert ^{2\theta}}\right)  \left(  1+\sum
\limits_{k<\theta-n}\frac{\left\vert x\right\vert ^{2k}}{k!}\right)
,\label{a1.24}%
\end{align}

where $e_{m}:=\sum\limits_{k\leq m}\frac{1}{k!}$, and from \ref{a64},%
\begin{align}
&  \int\limits_{\left\vert \cdot\right\vert \geq r_{3}}\left\vert
\mathcal{Q}_{\emptyset,\theta-n,\xi}\left(  e^{ix\xi}\right)  \right\vert
^{2}\frac{\left\vert \widehat{a}\xi\right\vert ^{2n}}{w\left\vert
\cdot\right\vert ^{2\theta}}\nonumber\\
&  \leq\left(  \overline{C}_{\theta-n,r_{3}}^{\left(  \rho\right)  }\right)
^{2}\left(
\begin{array}
[c]{c}%
\int\limits_{\left\vert \cdot\right\vert \geq r_{3}}\left\vert \widehat
{a}\widehat{\xi}\right\vert ^{2n}\frac{\left\vert \cdot\right\vert ^{2n}%
}{w\left\vert \cdot\right\vert ^{2\theta}}+\\
+e_{\theta-n-1}\int\limits_{\left\vert \cdot\right\vert \geq r_{3}}\left\vert
\widehat{x}\widehat{\xi}\right\vert ^{2k}\left\vert \widehat{a}\widehat{\xi
}\right\vert ^{2n}\frac{\left\vert \cdot\right\vert ^{2n}}{w\left\vert
\cdot\right\vert ^{2\theta}}%
\end{array}
\right)  \left(  1+\sum\limits_{k<\theta-n}\frac{\left\vert x\right\vert
^{2k}}{k!}\right) \label{a1.25}\\
&  <\infty.\nonumber
\end{align}

since $w\in W3.2$. Substituting the estimates \ref{a1.24} and \ref{a60} into
\ref{a1.00} yields
\begin{align}
&  \left\vert \left(  \widehat{a}D\right)  ^{n}f_{\rho}\left(  x\right)
\right\vert \nonumber\\
&  \leq\left(  \int\limits_{\left\vert \xi\right\vert \leq r_{3}}\left\vert
\mathcal{Q}_{\emptyset,\theta-n,\xi}\left(  e^{ix\xi}\right)  \right\vert
^{2}\frac{\left\vert \widehat{a}\xi\right\vert ^{2n}}{w\left\vert
\cdot\right\vert ^{2\theta}}+\int\limits_{\left\vert \xi\right\vert \geq
r_{3}}\left\vert \mathcal{Q}_{\emptyset,\theta-n,\xi}\left(  e^{ix\xi}\right)
\right\vert ^{2}\frac{\left\vert \widehat{a}\xi\right\vert ^{2n}}{w\left\vert
\cdot\right\vert ^{2\theta}}\right)  ^{1/2}\left\vert f\right\vert _{w,\theta
}\nonumber\\
&  \leq\left(
\begin{array}
[c]{l}%
\left(  \underline{C}_{\theta-n,r_{3}}^{\left(  \rho\right)  }\right)
^{2}\left(  1+\frac{1}{\left(  \theta-n\right)  !}\right)  \left(
\int\limits_{\left\vert \cdot\right\vert \leq r_{3}}\frac{1}{w}\right)
\left(  1+\frac{\left\vert x\right\vert ^{2\left(  \theta-n\right)  }}{\left(
\theta-n\right)  !}\right)  +\\
\qquad+\left(  \overline{C}_{\theta-n,r_{3}}^{\left(  \rho\right)  }\right)
^{2}\left(  1+e_{\theta-n-1}\right)  \left(  \int\limits_{\left\vert
\cdot\right\vert \geq r_{3}}\frac{\left\vert \cdot\right\vert ^{2n}%
}{w\left\vert \cdot\right\vert ^{2\theta}}\right)  \left(  1+\sum
\limits_{k<\theta-n}\frac{\left\vert x\right\vert ^{2k}}{k!}\right)
\end{array}
\right)  ^{1/2}\left\vert f\right\vert _{w,\theta}\nonumber\\
&  \leq\max\left\{
\begin{array}
[c]{c}%
\underline{C}_{\theta-n,r_{3}}^{\left(  \rho\right)  },\\
\overline{C}_{\theta-n,r_{3}}^{\left(  \rho\right)  }%
\end{array}
\right\}  \left(  1+\sum\limits_{k\leq\theta-n}\frac{1}{k!}\right)
^{1/2}\left(  \int\limits_{\left\vert \cdot\right\vert \leq r_{3}}\frac{1}%
{w}+\int\limits_{\left\vert \cdot\right\vert \geq r_{3}}\frac{\left\vert
\cdot\right\vert ^{2n}}{w\left\vert \cdot\right\vert ^{2\theta}}\right)
^{1/2}\left\vert f\right\vert _{w,\theta}\left(  1+\sum\limits_{k\leq\theta
-n}\frac{\left\vert x\right\vert ^{2k}}{k!}\right)  ^{1/2}\nonumber\\
&  <\max\left\{
\begin{array}
[c]{c}%
\underline{C}_{\theta-n,r_{3}}^{\left(  \rho\right)  },\\
\overline{C}_{\theta-n,r_{3}}^{\left(  \rho\right)  }%
\end{array}
\right\}  \sqrt{1+e}\left(  \int\limits_{\left\vert \cdot\right\vert \leq
r_{3}}\frac{1}{w}+\int\limits_{\left\vert \cdot\right\vert \geq r_{3}}%
\frac{\left\vert \cdot\right\vert ^{2n}}{w\left\vert \cdot\right\vert
^{2\theta}}\right)  ^{1/2}\left\vert f\right\vert _{w,\theta}\left(
1+\sum\limits_{k\leq\theta-n}\frac{\left\vert x\right\vert ^{2k}}{k!}\right)
^{1/2}.\label{a1.26}%
\end{align}

Substituting the estimates \ref{a1.25} and \ref{a611} into \ref{a1.00} yields%
\begin{align*}
&  \left\vert \left(  \widehat{a}D\right)  ^{n}f_{\rho}\left(  x\right)
\right\vert \\
&  \leq\left(  \int\limits_{\left\vert \xi\right\vert \leq r_{3}}\left\vert
\mathcal{Q}_{\emptyset,\theta-n,\xi}\left(  e^{ix\xi}\right)  \right\vert
^{2}\frac{\left\vert \widehat{a}\xi\right\vert ^{2n}}{w\left\vert
\cdot\right\vert ^{2\theta}}+\int\limits_{\left\vert \xi\right\vert \geq
r_{3}}\left\vert \mathcal{Q}_{\emptyset,\theta-n,\xi}\left(  e^{ix\xi}\right)
\right\vert ^{2}\frac{\left\vert \widehat{a}\xi\right\vert ^{2n}}{w\left\vert
\cdot\right\vert ^{2\theta}}\right)  ^{1/2}\left\vert f\right\vert _{w,\theta
}\\
&  \leq\left(
\begin{array}
[c]{l}%
\left(  \underline{C}_{\theta-n,r_{3}}^{\left(  \rho\right)  }\right)
^{2}\left(  \int\limits_{\left\vert \cdot\right\vert \leq r_{3}}%
\frac{\left\vert \widehat{a}\widehat{\xi}\right\vert ^{2n}}{w}+\frac
{1}{\left(  \theta-n\right)  !}\int\limits_{\left\vert \cdot\right\vert \leq
r_{3}}\frac{\left\vert \widehat{x}\widehat{\xi}\right\vert ^{2\left(
\theta-n\right)  }\left\vert \widehat{a}\widehat{\xi}\right\vert ^{2n}}%
{w}\right)  \left(  1+\frac{\left\vert x\right\vert ^{2\left(  \theta
-n\right)  }}{\left(  \theta-n\right)  !}\right)  +\\
+\left(  \overline{C}_{\theta-n,r_{3}}^{\left(  \rho\right)  }\right)
^{2}\left(
\begin{array}
[c]{c}%
\int\limits_{\left\vert \cdot\right\vert \geq r_{3}}\left\vert \widehat
{a}\widehat{\xi}\right\vert ^{2n}\frac{\left\vert \cdot\right\vert ^{2n}%
}{w\left\vert \cdot\right\vert ^{2\theta}}+\\
+\sum\limits_{k<\theta-n}\frac{1}{k!}\int\limits_{\left\vert \cdot\right\vert
\geq r_{3}}\left\vert \widehat{x}\widehat{\xi}\right\vert ^{2k}\left\vert
\widehat{a}\widehat{\xi}\right\vert ^{2n}\frac{\left\vert \cdot\right\vert
^{2n}}{w\left\vert \cdot\right\vert ^{2\theta}}%
\end{array}
\right)  \left(  1+\sum\limits_{k<\theta-n}\frac{\left\vert x\right\vert
^{2k}}{k!}\right)
\end{array}
\right)  ^{1/2}\left\vert f\right\vert _{w,\theta}\\
&  \leq\max\left\{
\begin{array}
[c]{c}%
\underline{C}_{\theta-n,r_{3}}^{\left(  \rho\right)  },\\
\overline{C}_{\theta-n,r_{3}}^{\left(  \rho\right)  }%
\end{array}
\right\}  \left(
\begin{array}
[c]{l}%
\int\limits_{\left\vert \cdot\right\vert \leq r_{3}}\frac{\left\vert
\widehat{a}\widehat{\xi}\right\vert ^{2n}}{w}+\frac{1}{\left(  \theta
-n\right)  !}\int\limits_{\left\vert \cdot\right\vert \leq r_{3}}%
\frac{\left\vert \widehat{x}\widehat{\xi}\right\vert ^{2\left(  \theta
-n\right)  }\left\vert \widehat{a}\widehat{\xi}\right\vert ^{2n}}{w}+\\
+\int\limits_{\left\vert \cdot\right\vert \geq r_{3}}\left\vert \widehat
{a}\widehat{\xi}\right\vert ^{2n}\frac{\left\vert \cdot\right\vert ^{2n}%
}{w\left\vert \cdot\right\vert ^{2\theta}}+\sum\limits_{k<\theta-n}\frac
{1}{k!}\int\limits_{\left\vert \cdot\right\vert \geq r_{3}}\left\vert
\widehat{x}\widehat{\xi}\right\vert ^{2k}\left\vert \widehat{a}\widehat{\xi
}\right\vert ^{2n}\frac{\left\vert \cdot\right\vert ^{2n}}{w\left\vert
\cdot\right\vert ^{2\theta}}%
\end{array}
\right)  ^{1/2}\times\\
&  \qquad\qquad\times\left(  1+\sum\limits_{k\leq\theta-n}\frac{\left\vert
x\right\vert ^{2k}}{k!}\right)  ^{1/2}\left\vert f\right\vert _{w,\theta},
\end{align*}

so we have the estimates%
\begin{align}
&  \left\vert \left(  \widehat{a}D\right)  ^{n}f_{\rho}\left(  x\right)
\right\vert \nonumber\\
&  \leq\left\{
\begin{array}
[c]{l}%
\max\left\{
\begin{array}
[c]{c}%
\underline{C}_{\theta-n,r_{3}}^{\left(  \rho\right)  },\\
\overline{C}_{\theta-n,r_{3}}^{\left(  \rho\right)  }%
\end{array}
\right\}  \sqrt{1+e}\left(  \int\limits_{\left\vert \cdot\right\vert \leq
r_{3}}\frac{1}{w}+\int\limits_{\left\vert \cdot\right\vert \geq r_{3}}%
\frac{\left\vert \cdot\right\vert ^{2n}}{w\left\vert \cdot\right\vert
^{2\theta}}\right)  ^{1/2}\left\vert f\right\vert _{w,\theta}\left(
1+\sum\limits_{k\leq\theta-n}\frac{\left\vert x\right\vert ^{2k}}{k!}\right)
^{1/2},\\
and\\
\max\left\{
\begin{array}
[c]{c}%
\underline{C}_{\theta-n,r_{3}}^{\left(  \rho\right)  },\\
\overline{C}_{\theta-n,r_{3}}^{\left(  \rho\right)  }%
\end{array}
\right\}  \left(
\begin{array}
[c]{l}%
\int\limits_{\left\vert \cdot\right\vert \leq r_{3}}\frac{\left\vert
\widehat{a}\widehat{\xi}\right\vert ^{2n}}{w}+\frac{1}{\left(  \theta
-n\right)  !}\int\limits_{\left\vert \cdot\right\vert \leq r_{3}}%
\frac{\left\vert \widehat{x}\widehat{\xi}\right\vert ^{2\left(  \theta
-n\right)  }\left\vert \widehat{a}\widehat{\xi}\right\vert ^{2n}}{w}+\\
+\int\limits_{\left\vert \cdot\right\vert \geq r_{3}}\left\vert \widehat
{a}\widehat{\xi}\right\vert ^{2n}\frac{\left\vert \cdot\right\vert ^{2n}%
}{w\left\vert \cdot\right\vert ^{2\theta}}+\sum\limits_{k<\theta-n}\frac
{1}{k!}\int\limits_{\left\vert \cdot\right\vert \geq r_{3}}\left\vert
\widehat{x}\widehat{\xi}\right\vert ^{2k}\left\vert \widehat{a}\widehat{\xi
}\right\vert ^{2n}\frac{\left\vert \cdot\right\vert ^{2n}}{w\left\vert
\cdot\right\vert ^{2\theta}}%
\end{array}
\right)  ^{\frac{1}{2}}\times\\
\qquad\qquad\times\left\vert f\right\vert _{w,\theta}\left(  1+\sum
\limits_{k\leq\theta-n}\frac{\left\vert x\right\vert ^{2k}}{k!}\right)
^{1/2},
\end{array}
\right. \label{a2.19}%
\end{align}

which can be written%
\begin{equation}
\left\vert \left(  \widehat{a}D\right)  ^{n}f_{\rho}\left(  x\right)
\right\vert \leq\max\left\{  \underline{C}_{\theta-n,r_{3}}^{\left(
\rho\right)  },\overline{C}_{\theta-n,r_{3}}^{\left(  \rho\right)  }\right\}
M_{n,r_{3}}^{\left(  w\right)  }\left\vert f\right\vert _{w,\theta}\left(
1+\sum\limits_{k\leq\theta-n}\frac{\left\vert x\right\vert ^{2k}}{k!}\right)
^{1/2},\text{ }n\leq\left\lfloor \kappa\right\rfloor \leq\theta,\label{a2.23}%
\end{equation}

where we have the series of "constants",%
\begin{equation}
\left.
\begin{array}
[c]{ll}%
M_{n,r_{3}}^{\left(  w\right)  }\left(  \widehat{a},\widehat{x}\right)  := &
\left(
\begin{array}
[c]{c}%
\int\limits_{\left\vert \cdot\right\vert \leq r_{3}}\frac{\left\vert
\widehat{a}\widehat{\xi}\right\vert ^{2n}}{w}+\frac{1}{\left(  \theta
-n\right)  !}\int\limits_{\left\vert \cdot\right\vert \leq r_{3}}%
\frac{\left\vert \widehat{x}\widehat{\xi}\right\vert ^{2\left(  \theta
-n\right)  }\left\vert \widehat{a}\widehat{\xi}\right\vert ^{2n}}{w}+\\
+\int\limits_{\left\vert \cdot\right\vert \geq r_{3}}\left\vert \widehat
{a}\widehat{\xi}\right\vert ^{2n}\frac{\left\vert \cdot\right\vert ^{2n}%
}{w\left\vert \cdot\right\vert ^{2\theta}}+\sum\limits_{k<\theta-n}\frac
{1}{k!}\int\limits_{\left\vert \cdot\right\vert \geq r_{3}}\left\vert
\widehat{x}\widehat{\xi}\right\vert ^{2k}\left\vert \widehat{a}\widehat{\xi
}\right\vert ^{2n}\frac{\left\vert \cdot\right\vert ^{2n}}{w\left\vert
\cdot\right\vert ^{2\theta}}%
\end{array}
\right)  ^{\frac{1}{2}},\\
& or\text{ }more\text{ }weakly\\
M_{n,r_{3}}^{\left(  w\right)  }\left(  \widehat{a},\widehat{x}\right)  := &
\left(
\begin{array}
[c]{c}%
\int\limits_{\left\vert \cdot\right\vert \leq r_{3}}\frac{\left\vert
\widehat{a}\widehat{\xi}\right\vert ^{2n}}{w}+\frac{1}{\left(  \theta
-n\right)  !}\min\left\{  \int\limits_{\left\vert \cdot\right\vert \leq r_{3}%
}\frac{\left\vert \widehat{x}\widehat{\xi}\right\vert ^{2\left(
\theta-n\right)  }}{w},\int\limits_{\left\vert \cdot\right\vert \leq r_{3}%
}\frac{\left\vert \widehat{a}\widehat{\xi}\right\vert ^{2n}}{w}\right\}  +\\
+\int\limits_{\left\vert \cdot\right\vert \geq r_{3}}\left\vert \widehat
{a}\widehat{\xi}\right\vert ^{2n}\frac{\left\vert \cdot\right\vert ^{2n}%
}{w\left\vert \cdot\right\vert ^{2\theta}}+\\
+\sum\limits_{k<\theta-n}\frac{1}{k!}\min\left\{  \int\limits_{\left\vert
\cdot\right\vert \geq r_{3}}\left\vert \widehat{a}\widehat{\xi}\right\vert
^{2n}\frac{\left\vert \cdot\right\vert ^{2n}}{w\left\vert \cdot\right\vert
^{2\theta}},\int\limits_{\left\vert \cdot\right\vert \geq r_{3}}\left\vert
\widehat{x}\widehat{\xi}\right\vert ^{2k}\frac{\left\vert \cdot\right\vert
^{2n}}{w\left\vert \cdot\right\vert ^{2\theta}}\right\}
\end{array}
\right)  ^{\frac{1}{2}},\\
& or\text{ }more\text{ }weakly\\
M_{n,r_{3}}^{\left(  w\right)  }\left(  \widehat{a}\right)  := & \sqrt
{1+e}\left(  \int\limits_{\left\vert \cdot\right\vert \leq r_{3}}\left\vert
\widehat{a}\widehat{\xi}\right\vert ^{2n}\frac{1}{w}+\int\limits_{\left\vert
\cdot\right\vert \geq r_{3}}\left\vert \widehat{a}\widehat{\xi}\right\vert
^{2n}\frac{\left\vert \cdot\right\vert ^{2n}}{w\left\vert \cdot\right\vert
^{2\theta}}\right)  ^{1/2},\\
& or\text{ }more\text{ }weakly\\
M_{n,r_{3}}^{\left(  w\right)  }:= & \sqrt{1+e}\left(  \int\limits_{\left\vert
\cdot\right\vert \leq r_{3}}\frac{1}{w}+\int\limits_{\left\vert \cdot
\right\vert \geq r_{3}}\frac{\left\vert \cdot\right\vert ^{2n}}{w\left\vert
\cdot\right\vert ^{2\theta}}\right)  ^{1/2}.
\end{array}
\right\} \label{a2.27}%
\end{equation}

Further, for compactness we will define%
\begin{equation}
C_{n}^{\left(  \rho,w\right)  }:=\max\left\{  \underline{C}_{\theta-n,r_{3}%
}^{\left(  \rho\right)  },\overline{C}_{\theta-n,r_{3}}^{\left(  \rho\right)
}\right\}  M_{n,r_{3}}^{\left(  w\right)  },\text{\quad}n\leq\left\lfloor
\kappa\right\rfloor \leq\theta.\label{a2.29}%
\end{equation}

\textbf{NOTE} that if $w$ is radial, $\widehat{a}$ and $\widehat{x}$ will be
eliminated from all but the strongest bound - see ??.\medskip

\fbox{\textbf{Case 2} $\left\lfloor \kappa\right\rfloor \geq\theta$ and
$n\leq\left\lfloor \kappa\right\rfloor $} From \ref{a49},%
\begin{align}
&  \left\vert \left(  \widehat{a}D\right)  ^{n}f_{\rho}\left(  x\right)
\right\vert \nonumber\\
&  \leq\left(  2\pi\right)  ^{-\frac{d}{2}}\int\left\vert \widehat{a}%
\xi\right\vert ^{n}\left\vert f_{F}\right\vert \nonumber\\
&  =\left(  2\pi\right)  ^{-\frac{d}{2}}\left(  \int_{\left\vert
\cdot\right\vert \leq r_{3}}\left\vert \widehat{a}\xi\right\vert
^{n}\left\vert f_{F}\right\vert +\int_{\left\vert \cdot\right\vert \geq r_{3}%
}\left\vert \widehat{a}\xi\right\vert ^{n}\left\vert f_{F}\right\vert \right)
\nonumber\\
&  =\left(  2\pi\right)  ^{-\frac{d}{2}}\left(  \int_{\left\vert
\cdot\right\vert \leq r_{3}}\frac{\left\vert \widehat{a}\xi\right\vert
^{n-\theta}\left\vert \widehat{a}\xi\right\vert ^{\theta}}{\sqrt{w}\left\vert
\cdot\right\vert ^{\theta}}\sqrt{w}\left\vert \cdot\right\vert ^{\theta
}\left\vert f_{F}\right\vert +\int_{\left\vert \cdot\right\vert \geq r_{3}%
}\frac{\left\vert \widehat{a}\xi\right\vert ^{n}}{\sqrt{w}\left\vert
\cdot\right\vert ^{\theta}}\sqrt{w}\left\vert \cdot\right\vert ^{\theta
}\left\vert f_{F}\right\vert \right) \nonumber\\
&  \leq\left(  2\pi\right)  ^{-\frac{d}{2}}\left(  \int_{\left\vert
\cdot\right\vert \leq r_{3}}\frac{\left\vert \widehat{a}\xi\right\vert
^{n-\theta}\left\vert \cdot\right\vert ^{\theta}}{\sqrt{w}\left\vert
\cdot\right\vert ^{\theta}}\sqrt{w}\left\vert \cdot\right\vert ^{\theta
}\left\vert f_{F}\right\vert +\int_{\left\vert \cdot\right\vert \geq r_{3}%
}\frac{\left\vert \widehat{a}\xi\right\vert ^{n}}{\sqrt{w}\left\vert
\cdot\right\vert ^{\theta}}\sqrt{w}\left\vert \cdot\right\vert ^{\theta
}\left\vert f_{F}\right\vert \right) \nonumber\\
&  =\left(  2\pi\right)  ^{-\frac{d}{2}}\left(  \int_{\left\vert
\cdot\right\vert \leq r_{3}}\frac{\left\vert \widehat{a}\xi\right\vert
^{n-\theta}}{\sqrt{w}}\sqrt{w}\left\vert \cdot\right\vert ^{\theta}\left\vert
f_{F}\right\vert +\int_{\left\vert \cdot\right\vert \geq r_{3}}\frac
{\left\vert \widehat{a}\xi\right\vert ^{n}}{\sqrt{w}\left\vert \cdot
\right\vert ^{\theta}}\sqrt{w}\left\vert \cdot\right\vert ^{\theta}\left\vert
f_{F}\right\vert \right) \nonumber\\
&  \leq\left(  2\pi\right)  ^{-\frac{d}{2}}\left\vert f_{\rho}\right\vert
_{w,\theta}\left(  \left(  \int_{\left\vert \cdot\right\vert \leq r_{3}}%
\frac{\left\vert \widehat{a}\xi\right\vert ^{2\left(  n-\theta\right)  }}%
{w}\right)  ^{\frac{1}{2}}+\left(  \int_{\left\vert \cdot\right\vert \geq
r_{3}}\frac{\left\vert \widehat{a}\xi\right\vert ^{2n}}{w\left\vert
\cdot\right\vert ^{2\theta}}\right)  ^{\frac{1}{2}}\right) \label{a56}\\
&  \leq\left(  2\pi\right)  ^{-\frac{d}{2}}\left\vert f_{\rho}\right\vert
_{w,\theta}\left(  \left(  \int_{\left\vert \cdot\right\vert \leq r_{3}}%
\frac{\left\vert \cdot\right\vert ^{2\left(  n-\theta\right)  }}{w}\right)
^{\frac{1}{2}}+\left(  \int_{\left\vert \cdot\right\vert \geq r_{3}}%
\frac{\left\vert \cdot\right\vert ^{2n}}{w\left\vert \cdot\right\vert
^{2\theta}}\right)  ^{\frac{1}{2}}\right)  .\nonumber\\
&  <\infty.\nonumber
\end{align}

Also, from the point of view of obtaining a bound%
\begin{align*}
\left\vert \left(  \widehat{a}D\right)  ^{n}f_{\rho}\left(  x\right)
\right\vert  & \leq\left(  2\pi\right)  ^{-\frac{d}{2}}\int\left\vert
\widehat{a}\xi\right\vert ^{n}\left\vert f_{F}\right\vert \\
& =\left(  2\pi\right)  ^{-\frac{d}{2}}\int\frac{\left\vert \widehat{a}%
\xi\right\vert ^{n}}{\sqrt{w}\left\vert \cdot\right\vert ^{\theta}}\sqrt
{w}\left\vert \cdot\right\vert ^{\theta}\left\vert f_{F}\right\vert \\
& \leq\left(  2\pi\right)  ^{-\frac{d}{2}}\left(  \int\frac{\left(
\widehat{a}\xi\right)  ^{2n}}{w\left\vert \cdot\right\vert ^{2\theta}}\right)
^{1/2}\left(  \int w\left\vert \cdot\right\vert ^{2\theta}\left\vert
f_{F}\right\vert ^{2}\right)  ^{1/2}\\
& =\left(  2\pi\right)  ^{-\frac{d}{2}}\left(  \int\frac{\left(  \widehat
{a}\xi\right)  ^{2n}}{w\left\vert \cdot\right\vert ^{2\theta}}\right)
^{1/2}\left\vert f_{\rho}\right\vert _{w,\theta}\\
& =\left(  2\pi\right)  ^{-\frac{d}{2}}\left(  \int\left(  \widehat{a}%
\widehat{\xi}\right)  ^{2n}\frac{\left\vert \cdot\right\vert ^{2n}%
}{w\left\vert \cdot\right\vert ^{2\theta}}\right)  ^{1/2}\left\vert f_{\rho
}\right\vert _{w,\theta}\\
& \leq\left(  2\pi\right)  ^{-\frac{d}{2}}\left(  \int\frac{\left\vert
\cdot\right\vert ^{2n}}{w\left\vert \cdot\right\vert ^{2\theta}}\right)
^{1/2}\left\vert f_{\rho}\right\vert _{w,\theta}\\
& <\infty.
\end{align*}

This completes our upper bounds for $\left\vert \left(  \widehat{a}D\right)
^{n}f_{\rho}\left(  x\right)  \right\vert $ and we have for $n\leq\left\lfloor
\kappa\right\rfloor $,
\begin{equation}
\left\vert \left(  \widehat{a}D\right)  ^{n}f_{\rho}\left(  x\right)
\right\vert \leq\left\{
\begin{array}
[c]{ll}%
\max\left\{  \underline{C}_{\theta-n,r_{3}}^{\left(  \rho\right)  }%
,\overline{C}_{\theta-n,r_{3}}^{\left(  \rho\right)  }\right\}  M_{n,r_{3}%
}^{\left(  w\right)  }\left\vert f_{\rho}\right\vert _{w,\theta}\left(
1+\sum\limits_{k\leq\theta-n}\frac{\left\vert x\right\vert ^{2k}}{k!}\right)
^{1/2}, & \theta>\left\lfloor \kappa\right\rfloor ,\\
\left(  2\pi\right)  ^{-\frac{d}{2}}\left(  \int\left(  \widehat{a}%
\widehat{\xi}\right)  ^{2n}\frac{\left\vert \cdot\right\vert ^{2n}%
}{w\left\vert \cdot\right\vert ^{2\theta}}\right)  ^{1/2}\left\vert f_{\rho
}\right\vert _{w,\theta}, & \theta\leq\left\lfloor \kappa\right\rfloor ,
\end{array}
\right. \label{a79}%
\end{equation}

where the $M_{n,r}^{\left(  w\right)  }$ are given \ by \ref{a2.27}.

\begin{remark}
We have $1+\sum\limits_{k\leq\theta-n}\frac{\left\vert x\right\vert ^{2k}}%
{k!}<e^{\left\vert x\right\vert ^{2}}$, $\left(  1+\sum\limits_{k\leq\theta
-n}\frac{\left\vert x\right\vert ^{2k}}{k!}\right)  ^{1/2}<e^{\frac{1}%
{2}\left\vert x\right\vert ^{2}}$ and $1+\sum\limits_{k<\theta-n}%
\frac{\left\vert x\widehat{\xi}\right\vert ^{k}}{k!}<e^{\left\vert
x\widehat{\xi}\right\vert ^{2}}$ etc. and I recall that the function
$e^{\left\vert x\right\vert ^{2}}$ \textbf{was important} in the chapter
concerning bounded linear functionals on $X_{w}^{\theta}$.
\end{remark}

\underline{\textbf{Apply the estimates for }$\left\vert \left(  \widehat
{a}D\right)  ^{n}f_{\rho}\left(  x\right)  \right\vert $\textbf{\ to the right
side of }\ref{a80}} There are two cases: $n\leq\left\lfloor \kappa
\right\rfloor <\theta$ and then $\theta\leq\left\lfloor \kappa\right\rfloor $
and $n\leq\left\lfloor \kappa\right\rfloor $.\medskip

\fbox{\textbf{Case 1} $n\leq\left\lfloor \kappa\right\rfloor <\theta$} Apply
\ref{a79} to \ref{a80}:%
\begin{align*}
&  \left\vert \frac{\sqrt{2\pi}}{n!}\left(  \overline{\widehat{g_{n}}}\left(
a\xi\right)  \left(  ia\xi\right)  \phi_{0}\widehat{\left(  aD\right)
^{n}f_{\rho}}\right)  ^{\vee}\left(  x\right)  \right\vert \\
&  \leq\frac{\left\vert a\right\vert ^{n+1}}{\left(  2\pi\right)  ^{\frac
{d}{2}}n!}\int_{0}^{1}\int\left\vert \widehat{\phi_{\widehat{a}}}\left(
\xi\right)  \right\vert \left\vert \left(  \left(  \widehat{a}D\right)
^{n}f_{\rho}\right)  \left(  \xi+x+sa\right)  \right\vert d\xi\text{ }%
g_{n}\left(  s\right)  ds\\
&  \leq\frac{\left\vert a\right\vert ^{n+1}}{\left(  2\pi\right)  ^{\frac
{d}{2}}n!}\int_{0}^{1}\int\left\vert \widehat{\phi_{\widehat{a}}}\left(
\xi\right)  \right\vert \left(  2\pi\right)  ^{-\frac{d}{2}}\left\vert
f_{\rho}\right\vert _{w,\theta}C_{n}^{\left(  \rho,w\right)  }\left(
1+\sum\limits_{k=0}^{\theta-n}\frac{\left(  \xi+x+sa\right)  ^{2k}}%
{k!}\right)  ^{\frac{1}{2}}d\xi\text{ }g_{n}\left(  s\right)  ds\\
&  =\left\vert f_{\rho}\right\vert _{w,\theta}\frac{C_{n}^{\left(
\rho,w\right)  }}{\left(  2\pi\right)  ^{\frac{d}{2}}}\frac{\left\vert
a\right\vert ^{n+1}}{n!}\int_{0}^{1}\int\left\vert \widehat{\phi_{\widehat{a}%
}}\left(  \xi\right)  \right\vert \left(  1+\sum\limits_{k=0}^{\theta-n}%
\frac{\left\vert \xi+x+sa\right\vert ^{2k}}{k!}\right)  ^{\frac{1}{2}}%
d\xi\text{ }g_{n}\left(  s\right)  ds\\
&  =\left\vert f_{\rho}\right\vert _{w,\theta}\frac{C_{n}^{\left(
\rho,w\right)  }}{\left(  2\pi\right)  ^{\frac{d}{2}}}\frac{\left\vert
a\right\vert ^{n+1}}{n!}\int_{0}^{1}\int\left\vert \widehat{\phi_{\widehat{a}%
}}\left(  \xi\right)  \right\vert ^{\frac{1}{2}}\left(  \left(  1+\sum
\limits_{k=0}^{\theta-n}\frac{\left(  \left\vert \xi+x+s\right\vert \right)
^{2k}}{k!}\right)  \left\vert \widehat{\phi_{\widehat{a}}}\left(  \xi\right)
\right\vert \right)  ^{\frac{1}{2}}d\xi\text{ }g_{n}\left(  s\right)  ds.
\end{align*}

The Cauchy-Schwartz inequality gives%
\begin{align}
&  \left\vert \frac{\sqrt{2\pi}}{n!}\left(  \overline{\widehat{g_{n}}}\left(
a\xi\right)  \left(  ia\xi\right)  \phi_{0}\widehat{\left(  aD\right)
^{n}f_{\rho}}\right)  ^{\vee}\left(  x\right)  \right\vert \nonumber\\
&  \leq\left\vert f_{\rho}\right\vert _{w,\theta}\frac{C_{n}^{\left(
\rho,w\right)  }}{\left(  2\pi\right)  ^{\frac{d}{2}}}\frac{\left\vert
a\right\vert ^{n+1}}{n!}\int_{0}^{1}\left\Vert \widehat{\phi_{\widehat{a}}%
}\right\Vert _{1}^{\frac{1}{2}}\left(  \int\left(  1+\sum\limits_{k=0}%
^{\theta-n}\frac{\left\vert \xi+x+sa\right\vert ^{2k}}{k!}\right)  \left\vert
\widehat{\phi_{\widehat{a}}}\left(  \xi\right)  \right\vert d\xi\right)
^{\frac{1}{2}}g_{n}\left(  s\right)  ds\nonumber\\
&  =\left\vert f_{\rho}\right\vert _{w,\theta}\frac{C_{n}^{\left(
\rho,w\right)  }}{\left(  2\pi\right)  ^{\frac{d}{2}}}\frac{\left\vert
a\right\vert ^{n+1}}{n!}\left\Vert \widehat{\phi_{\widehat{a}}}\right\Vert
_{1}^{\frac{1}{2}}\int_{0}^{1}\left(  \int\left(  1+\sum\limits_{k=0}%
^{\theta-n}\frac{\left\vert \xi+x+sa\right\vert ^{2k}}{k!}\right)  \left\vert
\widehat{\phi_{\widehat{a}}}\left(  \xi\right)  \right\vert d\xi\right)
^{\frac{1}{2}}g_{n}\left(  s\right)  ds\nonumber\\
&  =\left\vert f_{\rho}\right\vert _{w,\theta}\frac{C_{n}^{\left(
\rho,w\right)  }}{\left(  2\pi\right)  ^{\frac{d}{2}}}\frac{\left\vert
a\right\vert ^{n+1}}{n!}\left\Vert \widehat{\phi_{\widehat{a}}}\right\Vert
_{1}^{\frac{1}{2}}\int_{0}^{1}\left(  \left\Vert \widehat{\phi_{\widehat{a}}%
}\right\Vert _{1}+\sum\limits_{k=0}^{\theta-n}\int\frac{\left\vert
\xi+x+sa\right\vert ^{2k}}{\alpha!}\left\vert \widehat{\phi_{\widehat{a}}%
}\left(  \xi\right)  \right\vert d\xi\right)  ^{\frac{1}{2}}g_{n}\left(
s\right)  ds\nonumber\\
&  \leq\left\vert f_{\rho}\right\vert _{w,\theta}\frac{C_{n}^{\left(
\rho,w\right)  }}{\left(  2\pi\right)  ^{\frac{d}{2}}}\frac{\left\vert
a\right\vert ^{n+1}}{n!}\left\Vert \widehat{\phi_{\widehat{a}}}\right\Vert
_{1}^{\frac{1}{2}}\int_{0}^{1}\left(  \left(  \left\Vert \widehat
{\phi_{\widehat{a}}}\right\Vert _{1}\right)  ^{\frac{1}{2}}+\left(
\sum\limits_{k=0}^{\theta-n}\int\frac{\left\vert \xi+x+sa\right\vert ^{2k}%
}{k!}\left\vert \widehat{\phi_{\widehat{a}}}\left(  \xi\right)  \right\vert
d\xi\right)  ^{\frac{1}{2}}\right)  g_{n}\left(  s\right)  ds\nonumber\\
&  =\left\vert f_{\rho}\right\vert _{w,\theta}\frac{C_{n}^{\left(
\rho,w\right)  }}{\left(  2\pi\right)  ^{\frac{d}{2}}}\frac{\left\vert
a\right\vert ^{n+1}}{n!}\left\Vert \widehat{\phi_{\widehat{a}}}\right\Vert
_{1}\int_{0}^{1}g_{n}\left(  s\right)  ds+\nonumber\\
&  \qquad+\left\vert f_{\rho}\right\vert _{w,\theta}\frac{C_{n}^{\left(
\rho,w\right)  }}{\left(  2\pi\right)  ^{\frac{d}{2}}}\frac{\left\vert
a\right\vert ^{n+1}}{n!}\left\Vert \widehat{\phi_{\widehat{a}}}\right\Vert
_{1}\int_{0}^{1}\left(  \sum\limits_{k=0}^{\theta-n}\frac{1}{k!}\int\left\vert
\xi+x+sa\right\vert ^{2k}\left\vert \widehat{\phi_{\widehat{a}}}\left(
\xi\right)  \right\vert d\xi\right)  ^{\frac{1}{2}}g_{n}\left(  s\right)
ds\nonumber\\
&  =\left\vert f_{\rho}\right\vert _{w,\theta}\frac{C_{n}^{\left(
\rho,w\right)  }}{\left(  2\pi\right)  ^{\frac{d}{2}}}\frac{\left\vert
a\right\vert ^{n+1}}{\left(  n+1\right)  !}\left\Vert \widehat{\phi
_{\widehat{a}}}\right\Vert _{1}+\nonumber\\
&  \qquad+\left\vert f_{\rho}\right\vert _{w,\theta}\frac{C_{n}^{\left(
\rho,w\right)  }}{\left(  2\pi\right)  ^{\frac{d}{2}}}\frac{\left\vert
a\right\vert ^{n+1}}{n!}\left\Vert \widehat{\phi_{\widehat{a}}}\right\Vert
_{1}\int_{0}^{1}\left(  \sum\limits_{k=0}^{\theta-n}\frac{1}{k!}\int\left\vert
\xi+x+sa\right\vert ^{2k}\left\vert \widehat{\phi_{\widehat{a}}}\left(
\xi\right)  \right\vert d\xi\right)  ^{\frac{1}{2}}g_{n}\left(  s\right)
ds.\label{a81}%
\end{align}

Using the Cauchy-Schwartz estimate and then the identity \ref{Ap126} we get%
\begin{align*}
\frac{1}{k!}\left\vert x+y+z\right\vert ^{2k}  & \leq\frac{1}{k!}\left(
\left\vert x\right\vert +\left\vert y\right\vert +\left\vert z\right\vert
\right)  ^{2k}\\
& =\frac{1}{k!}\left(  \left(  \left\vert x\right\vert ,\left\vert
y\right\vert ,\left\vert z\right\vert \right)  \cdot\left(  1,1,1\right)
\right)  ^{2k}\\
& =\frac{3^{k}}{k!}\left(  \left\vert x\right\vert ^{2}+\left\vert
y\right\vert ^{2}+\left\vert z\right\vert ^{2}\right)  ^{k}\\
& =3^{k}\sum_{n_{1}+n_{2}+n_{3}=k}\frac{\left\vert x\right\vert ^{2n_{1}%
}\left\vert y\right\vert ^{2n_{2}}\left\vert z\right\vert ^{2n_{3}}}%
{n_{1}!n_{2}!n_{3}!},
\end{align*}

and so%
\begin{align*}
\int_{0}^{1} &  \left(  \sum\limits_{k=0}^{\theta-n}\frac{1}{k!}\int\left\vert
\xi+x+sa\right\vert ^{2k}\left\vert \widehat{\phi_{\widehat{a}}}\left(
\xi\right)  \right\vert d\xi\right)  ^{\frac{1}{2}}g_{n}\left(  s\right)  ds\\
&  =\int_{0}^{1}\left(  \sum\limits_{k=0}^{\theta-n}3^{k}\int\sum_{n_{1}%
+n_{2}+n_{3}=k}\frac{\left\vert \xi\right\vert ^{2n_{1}}\left\vert
x\right\vert ^{2n_{2}}\left\vert sa\right\vert ^{2n_{3}}}{n_{1}!n_{2}!n_{3}%
!}\left\vert \widehat{\phi_{\widehat{a}}}\left(  \xi\right)  \right\vert
d\xi\right)  ^{\frac{1}{2}}g_{n}\left(  s\right)  ds\\
&  \leq\int_{0}^{1}\left(  \sum\limits_{k=0}^{\theta-n}3^{k}\int\sum
_{n_{1}+n_{2}+n_{3}=k}\frac{\left\vert \xi\right\vert ^{2n_{1}}\left\vert
x\right\vert ^{2n_{2}}\left\vert a\right\vert ^{2n_{3}}}{n_{1}!n_{2}!n_{3}%
!}\left\vert \widehat{\phi_{\widehat{a}}}\left(  \xi\right)  \right\vert
d\xi\right)  ^{\frac{1}{2}}g_{n}\left(  s\right)  ds\\
&  =\frac{1}{n+1}\left(  \sum\limits_{k=0}^{\theta-n}3^{k}\int\sum
_{n_{1}+n_{2}+n_{3}=k}\frac{\left\vert \xi\right\vert ^{2n_{1}}\left\vert
x\right\vert ^{2n_{2}}\left\vert a\right\vert ^{2n_{3}}}{n_{1}!n_{2}!n_{3}%
!}\left\vert \widehat{\phi_{\widehat{a}}}\left(  \xi\right)  \right\vert
d\xi\right)  ^{\frac{1}{2}}\\
&  =\frac{1}{n+1}\left(  \sum\limits_{k=0}^{\theta-n}3^{k}\sum_{n_{1}%
+n_{2}+n_{3}=k}\frac{\left\vert x\right\vert ^{2n_{2}}\left\vert a\right\vert
^{2n_{3}}}{n_{2}!n_{3}!}\int\frac{\left\vert \xi\right\vert ^{2n_{1}}}{n_{1}%
!}\left\vert \widehat{\phi_{\widehat{a}}}\left(  \xi\right)  \right\vert
d\xi\right)  ^{\frac{1}{2}}\\
&  \leq\frac{1}{n+1}\left(  \sum\limits_{k=0}^{\theta-n}3^{k}\left(
\sum_{n_{2}+n_{3}\leq k}\frac{\left\vert x\right\vert ^{2n_{2}}\left\vert
a\right\vert ^{2n_{3}}}{n_{2}!n_{3}!}\right)  \sum_{n_{1}\leq k}\int%
\frac{\left\vert \xi\right\vert ^{2n_{1}}}{n_{1}!}\left\vert \widehat
{\phi_{\widehat{a}}}\left(  \xi\right)  \right\vert d\xi\right)  ^{\frac{1}%
{2}}\\
&  =\frac{1}{n+1}\left(  \sum\limits_{k=0}^{\theta-n}3^{k}\left(  \sum
_{j=0}^{k}\frac{\left(  \left\vert x\right\vert ^{2}+\left\vert a\right\vert
^{2}\right)  ^{j}}{j!}\right)  \sum_{l\leq k}\int\frac{\left\vert
\cdot\right\vert ^{2l}}{l!}\left\vert \widehat{\phi_{\widehat{a}}}\right\vert
\right)  ^{\frac{1}{2}}\\
&  \leq\frac{1}{n+1}\left(  \sum\limits_{k=0}^{\theta-n}3^{k}\right)
^{\frac{1}{2}}\left(  \sum_{j=0}^{\theta-n}\frac{\left(  \left\vert
x\right\vert ^{2}+\left\vert a\right\vert ^{2}\right)  ^{j}}{j!}\right)
^{\frac{1}{2}}\left(  \sum_{l=0}^{\theta-n}\int\frac{\left\vert \cdot
\right\vert ^{2l}}{l!}\left\vert \widehat{\phi_{\widehat{a}}}\right\vert
\right)  ^{\frac{1}{2}}\\
&  =\frac{1}{n+1}\left(  \frac{3^{\theta-n+1}-1}{2}\right)  ^{\frac{1}{2}%
}\left(  \sum_{j=0}^{\theta-n}\frac{\left(  \left\vert x\right\vert
^{2}+\left\vert a\right\vert ^{2}\right)  ^{j}}{j!}\right)  ^{\frac{1}{2}%
}\left(  \sum_{l=0}^{\theta-n}\int\frac{\left\vert \cdot\right\vert ^{2l}}%
{l!}\left\vert \widehat{\phi_{\widehat{a}}}\right\vert \right)  ^{\frac{1}{2}%
}.
\end{align*}

so that \ref{a81} becomes%
\begin{align*}
& \left\vert \frac{\sqrt{2\pi}}{n!}\left(  \overline{\widehat{g_{n}}}\left(
a\xi\right)  \left(  ia\xi\right)  \phi_{0}\widehat{\left(  aD\right)
^{n}f_{\rho}}\right)  ^{\vee}\left(  x\right)  \right\vert \\
& \leq\left\vert f_{\rho}\right\vert _{w,\theta}\frac{C_{n}^{\left(
\rho,w\right)  }}{\left(  2\pi\right)  ^{\frac{d}{2}}}\frac{\left\vert
a\right\vert ^{n+1}}{\left(  n+1\right)  !}\left\Vert \widehat{\phi
_{\widehat{a}}}\right\Vert _{1}+\\
& +\left\vert f_{\rho}\right\vert _{w,\theta}\frac{C_{n}^{\left(
\rho,w\right)  }}{\left(  2\pi\right)  ^{\frac{d}{2}}}\frac{\left\vert
a\right\vert ^{n+1}}{n!}\left\Vert \widehat{\phi_{\widehat{a}}}\right\Vert
_{1}^{1/2}\int_{0}^{1}\left(  \sum\limits_{k=0}^{\theta-n}\frac{1}{k!}%
\int\left\vert \xi+x+sa\right\vert ^{2k}\left\vert \widehat{\phi_{\widehat{a}%
}}\left(  \xi\right)  \right\vert d\xi\right)  ^{\frac{1}{2}}g_{n}\left(
s\right)  ds\\
& \leq\left\vert f_{\rho}\right\vert _{w,\theta}\frac{C_{n}^{\left(
\rho,w\right)  }}{\left(  2\pi\right)  ^{\frac{d}{2}}}\frac{\left\vert
a\right\vert ^{n+1}}{\left(  n+1\right)  !}\left\Vert \widehat{\phi
_{\widehat{a}}}\right\Vert _{1}+\\
& +\left\vert f_{\rho}\right\vert _{w,\theta}\frac{C_{n}^{\left(
\rho,w\right)  }}{\left(  2\pi\right)  ^{\frac{d}{2}}}\frac{\left\vert
a\right\vert ^{n+1}}{\left(  n+1\right)  !}\left\Vert \widehat{\phi
_{\widehat{a}}}\right\Vert _{1}^{1/2}\left(  \frac{3^{\theta-n+1}-1}%
{2}\right)  ^{\frac{1}{2}}\left(  \sum_{j=0}^{\theta-n}\frac{\left(
\left\vert x\right\vert ^{2}+\left\vert a\right\vert ^{2}\right)  ^{j}}%
{j!}\right)  ^{\frac{1}{2}}\left(  \sum_{l=0}^{\theta-n}\int\frac{\left\vert
\cdot\right\vert ^{2l}}{l!}\left\vert \widehat{\phi_{\widehat{a}}}\right\vert
\right)  ^{\frac{1}{2}}\\
& =\left\vert f_{\rho}\right\vert _{w,\theta}\frac{C_{n}^{\left(
\rho,w\right)  }}{\left(  2\pi\right)  ^{\frac{d}{2}}}\frac{\left\vert
a\right\vert ^{n+1}}{\left(  n+1\right)  !}\times\\
& \times\left(  \left\Vert \widehat{\phi_{\widehat{a}}}\right\Vert
_{1}+\left(  \frac{3^{\theta-n+1}-1}{2}\right)  ^{\frac{1}{2}}\left(
\sum\limits_{j=0}^{\theta-n}\frac{\left(  \left\vert x\right\vert
^{2}+\left\vert a\right\vert ^{2}\right)  ^{j}}{j!}\right)  ^{\frac{1}{2}%
}\left(  \left\Vert \widehat{\phi_{\widehat{a}}}\right\Vert _{1}%
\sum\limits_{l=0}^{\theta-n}\int\frac{\left\vert \cdot\right\vert ^{2l}}%
{l!}\left\vert \widehat{\phi_{\widehat{a}}}\right\vert \right)  ^{\frac{1}{2}%
}\right)  .
\end{align*}

Finally, from \ref{a1.22},%
\[
\widehat{\phi_{\widehat{a}}}\left(  \xi\right)  =\widehat{\widehat{a}x\phi
_{0}}\left(  \xi\right)  =\left(  \left(  \widehat{a}D\right)  \widehat
{\phi_{0}}\right)  \left(  \xi\right)  ,
\]

so that%
\begin{align}
&  \left\vert \frac{\sqrt{2\pi}}{n!}\left(  \overline{\widehat{g_{n}}}\left(
a\xi\right)  \left(  ia\xi\right)  \phi_{0}\widehat{\left(  aD\right)
^{n}f_{\rho}}\right)  ^{\vee}\left(  x\right)  \right\vert \nonumber\\
&  \leq\left\vert f_{\rho}\right\vert _{w,\theta}\frac{C_{n}^{\left(
\rho,w\right)  }}{\left(  2\pi\right)  ^{\frac{d}{2}}}\frac{\left\vert
a\right\vert ^{n+1}}{\left(  n+1\right)  !}\times\nonumber\\
&  \times\left(
\begin{array}
[c]{l}%
\left\Vert \left(  \widehat{a}D\right)  \widehat{\phi_{0}}\right\Vert _{1}+\\
+\left(  \frac{3^{\theta-n+1}-1}{2}\right)  ^{\frac{1}{2}}\left(
\sum\limits_{j=0}^{\theta-n}\frac{\left(  \left\vert x\right\vert
^{2}+\left\vert a\right\vert ^{2}\right)  ^{j}}{j!}\right)  ^{\frac{1}{2}%
}\left(  \left\Vert \left(  \widehat{a}D\right)  \widehat{\phi_{0}}\right\Vert
_{1}\sum\limits_{l=0}^{\theta-n}\frac{\left\Vert \left\vert \cdot\right\vert
^{2l}\widehat{a}D\widehat{\phi_{0}}\right\Vert _{1}}{l!}\right)  ^{\frac{1}%
{2}}%
\end{array}
\right) \label{a82}%
\end{align}

A special case is $\phi_{0}$ radial which will be treated in Subsection
\ref{SbSect_Tay_dat_rem_rad_fn}.\medskip

\fbox{\textbf{Case 2} $\left\lfloor \kappa\right\rfloor \geq\theta$ and
$n\leq\left\lfloor \kappa\right\rfloor $} By \ref{a79},%
\[
\left\vert \left(  \widehat{a}D\right)  ^{n}f_{\rho}\left(  x\right)
\right\vert \leq\left(  2\pi\right)  ^{-\frac{d}{2}}\left(  \int\left(
\widehat{a}\widehat{\xi}\right)  ^{2n}\frac{\left\vert \cdot\right\vert ^{2n}%
}{w\left\vert \cdot\right\vert ^{2\theta}}\right)  ^{\frac{1}{2}}\left\vert
f_{\rho}\right\vert _{w,\theta},
\]

so that by also using \ref{a1.05}, \ref{a80} becomes%
\begin{align}
&  \left\vert \frac{\sqrt{2\pi}}{n!}\left(  \overline{\widehat{g_{n}}}\left(
a\xi\right)  \left(  ia\xi\right)  \phi_{0}\widehat{\left(  aD\right)
^{n}f_{\rho}}\right)  ^{\vee}\left(  x\right)  \right\vert \nonumber\\
&  \leq\left(  2\pi\right)  ^{-\frac{d}{2}}\frac{\left\vert a\right\vert
^{n+1}}{n!}\int_{0}^{1}\int\left\vert \widehat{\phi_{\widehat{a}}}\left(
\xi\right)  \right\vert \left\vert \left(  \left(  \widehat{a}D\right)
^{n}f_{\rho}\right)  \left(  \xi+x+sa\right)  \right\vert d\xi\text{ }%
g_{n}\left(  s\right)  ds\nonumber\\
&  \leq\left(  2\pi\right)  ^{-\frac{d}{2}}\frac{\left\vert a\right\vert
^{n+1}}{n!}\int_{0}^{1}\left(  \int\left\vert \widehat{\phi_{\widehat{a}}%
}\left(  \xi\right)  \right\vert d\xi\right)  \left(  2\pi\right)  ^{-\frac
{d}{2}}\left(  \int\left(  \widehat{a}\widehat{\xi}\right)  ^{2n}%
\frac{\left\vert \cdot\right\vert ^{2n}}{w\left\vert \cdot\right\vert
^{2\theta}}\right)  ^{\frac{1}{2}}\left\vert f_{\rho}\right\vert _{w,\theta
}\text{ }g_{n}\left(  s\right)  ds\nonumber\\
&  =\frac{\left\vert f_{\rho}\right\vert _{w,\theta}}{\left(  2\pi\right)
^{d}}\frac{\left\vert a\right\vert ^{n+1}}{n!}\left(  \int\left\vert
\widehat{\phi_{\widehat{a}}}\right\vert \right)  \left(  \int\left(
\widehat{a}\widehat{\xi}\right)  ^{2n}\frac{\left\vert \cdot\right\vert ^{2n}%
}{w\left\vert \cdot\right\vert ^{2\theta}}\right)  ^{\frac{1}{2}}\int_{0}%
^{1}g_{n}\left(  s\right)  ds\nonumber\\
&  =\left\vert f_{\rho}\right\vert _{w,\theta}\frac{1}{\left(  2\pi\right)
^{d}}\left(  \int\left(  \widehat{a}\widehat{\xi}\right)  ^{2n}\frac
{\left\vert \cdot\right\vert ^{2n}}{w\left\vert \cdot\right\vert ^{2\theta}%
}\right)  ^{\frac{1}{2}}\frac{\left\vert a\right\vert ^{n+1}}{\left(
n+1\right)  !}\left\Vert \widehat{\phi_{\widehat{a}}}\right\Vert
_{1}\nonumber\\
&  =\left\vert f_{\rho}\right\vert _{w,\theta}\frac{1}{\left(  2\pi\right)
^{d}}\left(  \int\left(  \widehat{a}\widehat{\xi}\right)  ^{2n}\frac
{\left\vert \cdot\right\vert ^{2n}}{w\left\vert \cdot\right\vert ^{2\theta}%
}\right)  ^{\frac{1}{2}}\frac{\left\vert a\right\vert ^{n+1}}{\left(
n+1\right)  !}\left\Vert \widehat{a}D\widehat{\phi_{0}}\right\Vert
_{1}.\label{a85}%
\end{align}

Combining \ref{a82} and \ref{a85} we get for $n\leq\left\lfloor \kappa
\right\rfloor $,%
\begin{align}
&  \left\vert \frac{\sqrt{2\pi}}{n!}\left(  \overline{\widehat{g_{n}}}\left(
a\xi\right)  \left(  ia\xi\right)  \phi_{0}\widehat{\left(  aD\right)
^{n}f_{\rho}}\right)  ^{\vee}\left(  x\right)  \right\vert \nonumber\\
&  \leq\left\vert f_{\rho}\right\vert _{w,\theta}\frac{\left\vert a\right\vert
^{n+1}}{\left(  n+1\right)  !}\times\nonumber\\
&  \times\left\{
\begin{array}
[c]{ll}%
\frac{C_{n}^{\left(  \rho,w\right)  }}{\left(  2\pi\right)  ^{\frac{d}{2}}%
}\left(
\begin{array}
[c]{c}%
\left\Vert \widehat{a}D\widehat{\phi_{0}}\right\Vert _{1}+\\
+\left(  \frac{3^{\theta-n+1}-1}{2}\right)  ^{\frac{1}{2}}\left(
\sum\limits_{j=0}^{\theta-n}\frac{\left(  \left\vert x\right\vert
^{2}+\left\vert a\right\vert ^{2}\right)  ^{j}}{j!}\right)  ^{\frac{1}{2}%
}\left(  \left\Vert \widehat{a}D\widehat{\phi_{0}}\right\Vert _{1}%
\sum\limits_{l=0}^{\theta-n}\frac{\left\Vert \left\vert \cdot\right\vert
^{2l}\widehat{a}D\widehat{\phi_{0}}\right\Vert _{1}}{l!}\right)  ^{\frac{1}%
{2}}%
\end{array}
\right)  , & \theta\geq\left\lfloor \kappa\right\rfloor ,\\
\medskip & \\
\frac{1}{\left(  2\pi\right)  ^{d}}\left(  \int\left(  \widehat{a}\widehat
{\xi}\right)  ^{2n}\frac{\left\vert \cdot\right\vert ^{2n}}{w\left\vert
\cdot\right\vert ^{2\theta}}\right)  ^{\frac{1}{2}}\left\Vert \widehat
{a}D\widehat{\phi_{0}}\right\Vert _{1}, & \theta\leq\left\lfloor
\kappa\right\rfloor ,
\end{array}
\right. \label{a861}%
\end{align}

where $C_{n}^{\left(  \rho,w\right)  }$ is given by \ref{a2.29}.\medskip

\underline{\textbf{Estimate the }$f_{F}$\textbf{\ term in }\ref{a35}} The goal
here is to obtain the estimate \ref{a72}. To this end we use the step function
partition of unity $\left\{  \pi_{0},\pi_{\infty}\right\}  $ where
$\operatorname*{supp}\pi_{0}=\overline{B}_{r_{3}}$, $\operatorname*{supp}%
\pi_{\infty}=\mathbb{R}^{d}\setminus B_{r_{3}}$ where $r_{3}$ was introduced
in the definition of weight function property W3.2. Now write%
\begin{align*}
&  \left(  \left(  ia\xi\right)  ^{n+1}\overline{\widehat{g_{n}}}\left(
a\xi\right)  \phi_{\infty}f_{F}\right)  ^{\vee}\\
&  =\left(  \pi_{0}\left(  ia\xi\right)  ^{n+1}\overline{\widehat{g_{n}}%
}\left(  a\xi\right)  \phi_{\infty}f_{F}\right)  ^{\vee}+\left(  \pi_{\infty
}\left(  ia\xi\right)  ^{n+1}\overline{\widehat{g_{n}}}\left(  a\xi\right)
\phi_{\infty}f_{F}\right)  ^{\vee}\\
&  =\left(  \pi_{0}\frac{\left(  ia\xi\right)  ^{n+1}}{\sqrt{w}}%
\overline{\widehat{g_{n}}}\left(  a\xi\right)  \sqrt{w}\phi_{\infty}%
f_{F}\right)  ^{\vee}+\left(  \pi_{\infty}\phi_{\infty}\frac{\left(
ia\xi\right)  ^{n+1}}{\sqrt{w}\left\vert \cdot\right\vert ^{\theta}}%
\overline{\widehat{g_{n}}}\left(  a\xi\right)  \sqrt{w}\left\vert
\cdot\right\vert ^{\theta}f_{F}\right)  ^{\vee}\\
&  =\left(  \pi_{0}\frac{\left(  ia\xi\right)  ^{n+1}}{\sqrt{w}}%
\overline{\widehat{g_{n}}}\left(  a\xi\right)  \right)  ^{\vee}\ast\left(
\sqrt{w}\phi_{\infty}f_{F}\right)  ^{\vee}+\left(  \pi_{\infty}\phi_{\infty
}\frac{\left(  ia\xi\right)  ^{n+1}}{\sqrt{w}\left\vert \cdot\right\vert
^{\theta}}\overline{\widehat{g_{n}}}\left(  a\xi\right)  \right)  ^{\vee}%
\ast\left(  \sqrt{w}\left\vert \cdot\right\vert ^{\theta}f_{F}\right)  ^{\vee
},
\end{align*}

so that by Young's inequality \ref{1.056} and Parceval's theorem,%
\begin{align*}
&  \left\Vert \left(  \left(  ia\xi\right)  ^{n+1}\overline{\widehat{g_{n}}%
}\left(  a\xi\right)  \phi_{\infty}f_{F}\right)  ^{\vee}\right\Vert _{\infty
}\\
&  \leq\left\Vert \left(  \pi_{0}\frac{\left(  ia\xi\right)  ^{n+1}}{\sqrt{w}%
}\overline{\widehat{g_{n}}}\left(  a\xi\right)  \right)  ^{\vee}\right\Vert
_{2}\left\Vert \left(  \sqrt{w}\phi_{\infty}f_{F}\right)  ^{\vee}\right\Vert
_{2}+\\
&  \qquad\qquad+\left\Vert \left(  \pi_{\infty}\phi_{\infty}\frac{\left(
ia\xi\right)  ^{n+1}}{\sqrt{w}\left\vert \cdot\right\vert ^{\theta}}%
\overline{\widehat{g_{n}}}\left(  a\xi\right)  \right)  ^{\vee}\right\Vert
_{2}\left\Vert \left(  \sqrt{w}\left\vert \cdot\right\vert ^{\theta}%
f_{F}\right)  ^{\vee}\right\Vert _{2}\\
&  =\left\Vert \pi_{0}\frac{\left(  ia\xi\right)  ^{n+1}}{\sqrt{w}}%
\overline{\widehat{g_{n}}}\left(  a\xi\right)  \right\Vert _{2}\left\Vert
\sqrt{w}\phi_{\infty}f_{F}\right\Vert _{2}+\left\Vert \pi_{\infty}\phi
_{\infty}\frac{\left(  ia\xi\right)  ^{n+1}}{\sqrt{w}\left\vert \cdot
\right\vert ^{\theta}}\overline{\widehat{g_{n}}}\left(  a\xi\right)
\right\Vert _{2}\left\Vert \sqrt{w}\left\vert \cdot\right\vert ^{\theta}%
f_{F}\right\Vert _{2}\\
&  =\left\Vert \pi_{0}\frac{\left(  ia\xi\right)  ^{n+1}}{\sqrt{w}}%
\widehat{g_{n}}\left(  a\xi\right)  \right\Vert _{2}\left\Vert \sqrt{w}%
\phi_{\infty}f_{F}\right\Vert _{2}+\left\Vert \pi_{\infty}\phi_{\infty}%
\frac{\left(  ia\xi\right)  ^{n+1}}{\sqrt{w}\left\vert \cdot\right\vert
^{\theta}}\widehat{g_{n}}\left(  a\xi\right)  \right\Vert _{2}\left\vert
f\right\vert _{w,\theta}.
\end{align*}

By \ref{a11},
\[
\left\vert \phi_{\infty}\left(  x\right)  \right\vert \leq\left\vert
x\right\vert ^{\theta}\frac{\left\Vert \left(  \widehat{\cdot}D\right)
^{\theta}\phi_{\infty}\right\Vert _{\infty;\leq r_{3}}}{\theta!},\text{\quad
}x\in\overline{B}_{r_{3}},
\]

so%
\[
\left\Vert \sqrt{w}\phi_{\infty}f_{F}\right\Vert _{2}\leq\frac{\left\Vert
\left(  \widehat{\cdot}D\right)  ^{\theta}\phi_{\infty}\right\Vert
_{\infty;\leq r_{3}}}{\theta!}\left\Vert \sqrt{w}\left\vert \cdot\right\vert
^{\theta}f_{F}\right\Vert _{2}=\frac{\left\Vert \left(  \widehat{\cdot
}D\right)  ^{\theta}\phi_{\infty}\right\Vert _{\infty;\leq r_{3}}}{\theta
!}\left\vert f\right\vert _{w,\theta},
\]

which means that%
\begin{align}
&  \left\Vert \left(  \left(  ia\xi\right)  ^{n+1}\overline{\widehat{g_{n}}%
}\left(  a\xi\right)  \phi_{\infty}f_{F}\right)  ^{\vee}\right\Vert _{\infty
}\nonumber\\
&  \leq\left\Vert \pi_{0}\frac{\left(  ia\xi\right)  ^{n+1}}{\sqrt{w}}%
\widehat{g_{n}}\left(  a\xi\right)  \right\Vert _{2}\frac{\left\Vert \left(
\widehat{\cdot}D\right)  ^{\theta}\phi_{\infty}\right\Vert _{\infty;\leq
r_{3}}}{\theta!}\left\vert f\right\vert _{w,\theta}+\left\Vert \pi_{\infty
}\phi_{\infty}\frac{\left(  ia\xi\right)  ^{n+1}}{\sqrt{w}\left\vert
\cdot\right\vert ^{\theta}}\widehat{g_{n}}\left(  a\xi\right)  \right\Vert
_{2}\left\vert f\right\vert _{w,\theta}\nonumber\\
&  =\left(  \frac{\left\Vert \left(  \widehat{\cdot}D\right)  ^{\theta}%
\phi_{\infty}\right\Vert _{\infty;\leq r_{3}}}{\theta!}\left\Vert \pi_{0}%
\frac{\left(  ia\xi\right)  ^{n+1}}{\sqrt{w}}\widehat{g_{n}}\left(
a\xi\right)  \right\Vert _{2}+\left\Vert \pi_{\infty}\phi_{\infty}%
\frac{\left(  ia\xi\right)  ^{n+1}}{\sqrt{w}\left\vert \cdot\right\vert
^{\theta}}\widehat{g_{n}}\left(  a\xi\right)  \right\Vert _{2}\right)
\left\vert f\right\vert _{w,\theta}\nonumber\\
&  \leq\left(  \frac{\left\Vert \left(  \widehat{\cdot}D\right)  ^{\theta}%
\phi_{\infty}\right\Vert _{\infty;\leq r_{3}}}{\theta!}\left\Vert
\widehat{g_{n}}\right\Vert _{\infty}\left\Vert \pi_{0}\frac{\left(
ia\xi\right)  ^{n+1}}{\sqrt{w}}\right\Vert _{2}+\left\Vert \pi_{\infty}%
\frac{\left(  ia\xi\right)  ^{n+1}}{\sqrt{w}\left\vert \cdot\right\vert
^{\theta}}\widehat{g_{n}}\left(  a\xi\right)  \right\Vert _{2}\right)
\left\vert f\right\vert _{w,\theta}\nonumber\\
&  =\left(  \frac{\left\Vert \left(  \widehat{\cdot}D\right)  ^{\theta}%
\phi_{\infty}\right\Vert _{\infty;\leq r_{3}}}{\theta!}\left\Vert
\widehat{g_{n}}\right\Vert _{\infty}\left(  \int\limits_{\left\vert
\cdot\right\vert \leq r_{3}}\frac{\left(  a\xi\right)  ^{2\left(  n+1\right)
}}{w}\right)  ^{\frac{1}{2}}+\left(  \int\limits_{\left\vert \cdot\right\vert
\geq r_{3}}\frac{\left(  a\xi\right)  ^{2\left(  n+1\right)  }}{w\left\vert
\cdot\right\vert ^{2\theta}}\left\vert \widehat{g_{n}}\left(  a\xi\right)
\right\vert ^{2}\right)  ^{\frac{1}{2}}\right)  \left\vert f\right\vert
_{w,\theta}.\label{a59}%
\end{align}
\medskip

\fbox{If $n\leq\left\lfloor \kappa\right\rfloor -1$} then $2\left(
n+1\right)  \leq2\left\lfloor \kappa\right\rfloor \leq2\kappa$ and
so\smallskip%
\begin{align}
&  \left\Vert \left(  \left(  ia\xi\right)  ^{n+1}\overline{\widehat{g_{n}}%
}\left(  a\xi\right)  \phi_{\infty}f_{F}\right)  ^{\vee}\right\Vert _{\infty
}\nonumber\\
&  \leq\left\vert f\right\vert _{w,\theta}\left(  \frac{\left\Vert \left(
\widehat{\cdot}D\right)  ^{\theta}\phi_{\infty}\right\Vert _{\infty;\leq
r_{3}}}{\theta!}\left\Vert \widehat{g_{n}}\right\Vert _{\infty}\left(
\int\limits_{\left\vert \cdot\right\vert \leq r_{3}}\frac{\left(  a\xi\right)
^{2\left(  n+1\right)  }}{w}\right)  ^{\frac{1}{2}}+\left\Vert \widehat{g_{n}%
}\right\Vert _{\infty}\left(  \int\limits_{\left\vert \cdot\right\vert \geq
r_{3}}\frac{\left(  a\xi\right)  ^{2\left(  n+1\right)  }}{w\left\vert
\cdot\right\vert ^{2\theta}}\right)  ^{\frac{1}{2}}\right) \nonumber\\
&  \leq\left\vert f\right\vert _{w,\theta}\left\vert a\right\vert
^{n+1}\left\Vert \widehat{g_{n}}\right\Vert _{\infty}\left(  \frac{\left\Vert
\left(  \widehat{\cdot}D\right)  ^{\theta}\phi_{\infty}\right\Vert
_{\infty;\leq r_{3}}}{\theta!}\left(  \int\limits_{\left\vert \cdot\right\vert
\leq r_{3}}\frac{\left(  \widehat{a}\xi\right)  ^{2\left(  n+1\right)  }}%
{w}\right)  ^{\frac{1}{2}}+\left(  \int\limits_{\left\vert \cdot\right\vert
\geq r_{3}}\frac{\left(  \widehat{a}\xi\right)  ^{2\left(  n+1\right)  }%
}{w\left\vert \cdot\right\vert ^{2\theta}}\right)  ^{\frac{1}{2}}\right)
\nonumber\\
&  \leq\left\vert f\right\vert _{w,\theta}\frac{1}{\sqrt{2\pi}}\frac
{\left\vert a\right\vert ^{n+1}}{n+1}\left(  \frac{\left\Vert \left(
\widehat{\cdot}D\right)  ^{\theta}\phi_{\infty}\right\Vert _{\infty;\leq
r_{3}}}{\theta!}\left(  \int\limits_{\left\vert \cdot\right\vert \leq r_{3}%
}\frac{\left(  \widehat{a}\xi\right)  ^{2\left(  n+1\right)  }}{w}\right)
^{\frac{1}{2}}+\left(  \int\limits_{\left\vert \cdot\right\vert \geq r_{3}%
}\frac{\left(  \widehat{a}\xi\right)  ^{2\left(  n+1\right)  }}{w\left\vert
\cdot\right\vert ^{2\theta}}\right)  ^{\frac{1}{2}}\right) \label{a57}\\
&  <\infty,\nonumber
\end{align}

since $w\in W3.2$ for order $\theta$ and $\kappa$, and $w\in W2.1$.\medskip

\fbox{If $n=\left\lfloor \kappa\right\rfloor $} then $n+1=\left\lceil
\kappa\right\rceil $ and since from Lemma \ref{Lem_gm_properties_2},
$\left\vert \widehat{g_{n}}\left(  t\right)  \right\vert \leq\frac{1}%
{\sqrt{2\pi}}\min\left\{  \frac{1}{n+1},\frac{2+\frac{1}{n+1}}{1+\left\vert
t\right\vert }\right\}  $, we can proceed from \ref{a59} as follows:
\begin{align}
&  \left\Vert \left(  \left(  ia\xi\right)  ^{n+1}\overline{\widehat{g_{n}}%
}\left(  a\xi\right)  \phi_{\infty}f_{F}\right)  ^{\vee}\right\Vert _{\infty
}\frac{1}{\left\vert f\right\vert _{w,\theta}}\nonumber\\
&  \leq\frac{\left\Vert \left(  \widehat{\cdot}D\right)  ^{\theta}\phi
_{\infty}\right\Vert _{\infty;\leq r_{3}}}{\theta!}\left\Vert \widehat
{g_{\left\lfloor \kappa\right\rfloor }}\right\Vert _{\infty}\left(
\int\limits_{\left\vert \cdot\right\vert \leq r_{3}}\frac{\left(  a\xi\right)
^{2\left\lceil \kappa\right\rceil }}{w}\right)  ^{\frac{1}{2}}+\left(
\int\limits_{\left\vert \cdot\right\vert \geq r_{3}}\frac{\left(  a\xi\right)
^{2\left\lceil \kappa\right\rceil }}{w\left\vert \cdot\right\vert ^{2\theta}%
}\left\vert \widehat{g_{\left\lfloor \kappa\right\rfloor }}\left(
a\xi\right)  \right\vert ^{2}\right)  ^{\frac{1}{2}}\label{a66}\\
&  =\frac{\left\Vert \left(  \widehat{\cdot}D\right)  ^{\theta}\phi_{\infty
}\right\Vert _{\infty;\leq r_{3}}}{\theta!}\left\Vert \widehat{g_{\left\lfloor
\kappa\right\rfloor }}\right\Vert _{\infty}\left(  \int\limits_{\left\vert
\cdot\right\vert \leq r_{3}}\frac{\left(  a\xi\right)  ^{2\left\lceil
\kappa\right\rceil }}{w}\right)  ^{\frac{1}{2}}+\left(  \int%
\limits_{\left\vert \cdot\right\vert \geq r_{3}}\frac{\left(  a\xi\right)
^{2\kappa}}{w\left\vert \cdot\right\vert ^{2\theta}}\left\vert a\xi\right\vert
^{2\left(  \left\lceil \kappa\right\rceil -\kappa\right)  }\left\vert
\widehat{g_{\left\lfloor \kappa\right\rfloor }}\left(  a\xi\right)
\right\vert ^{2}\right)  ^{\frac{1}{2}}\nonumber\\
&  \leq\frac{\left\Vert \left(  \widehat{\cdot}D\right)  ^{\theta}\phi
_{\infty}\right\Vert _{\infty;\leq r_{3}}}{\theta!}\left\Vert \widehat
{g_{\left\lfloor \kappa\right\rfloor }}\right\Vert _{\infty}\left(
\int\limits_{\left\vert \cdot\right\vert \leq r_{3}}\frac{\left(  a\xi\right)
^{2\left\lceil \kappa\right\rceil }}{w}\right)  ^{\frac{1}{2}}+\left(
\int\limits_{\left\vert \cdot\right\vert \geq r_{3}}\frac{\left(  a\xi\right)
^{2\kappa}}{w\left\vert \cdot\right\vert ^{2\theta}}\left\vert \frac
{\left\vert a\xi\right\vert ^{\left(  \left\lceil \kappa\right\rceil
-\kappa\right)  }}{\sqrt{2\pi}}\frac{2+\frac{1}{\left\lceil \kappa\right\rceil
}}{1+\left\vert a\xi\right\vert }\right\vert ^{2}\right)  ^{\frac{1}{2}%
}\nonumber\\
&  \leq\frac{\left\Vert \left(  \widehat{\cdot}D\right)  ^{\theta}\phi
_{\infty}\right\Vert _{\infty;\leq r_{3}}}{\theta!}\left\Vert \widehat
{g_{\left\lfloor \kappa\right\rfloor }}\right\Vert _{\infty}\left(
\int\limits_{\left\vert \cdot\right\vert \leq r_{3}}\frac{\left(  a\xi\right)
^{2\left\lceil \kappa\right\rceil }}{w}\right)  ^{\frac{1}{2}}+\left(
\int\limits_{\left\vert \cdot\right\vert \geq r_{3}}\frac{\left\vert
a\xi\right\vert ^{2\kappa}}{w\left\vert \cdot\right\vert ^{2\theta}}\left\vert
\frac{1}{\sqrt{2\pi}}\left(  2+\frac{1}{\left\lceil \kappa\right\rceil
}\right)  \right\vert ^{2}\right)  ^{\frac{1}{2}}\nonumber\\
&  \leq\frac{\left\Vert \left(  \widehat{\cdot}D\right)  ^{\theta}\phi
_{\infty}\right\Vert _{\infty;\leq r_{3}}}{\theta!}\left\Vert \widehat
{g_{\left\lfloor \kappa\right\rfloor }}\right\Vert _{\infty}\left(
\int\limits_{\left\vert \cdot\right\vert \leq r_{3}}\frac{\left(  a\xi\right)
^{2\left\lceil \kappa\right\rceil }}{w}\right)  ^{\frac{1}{2}}+\frac{1}%
{\sqrt{2\pi}}\left(  2+\frac{1}{\left\lceil \kappa\right\rceil }\right)
\left(  \int\limits_{\left\vert \cdot\right\vert \geq r_{3}}\frac{\left\vert
a\xi\right\vert ^{2\kappa}}{w\left\vert \cdot\right\vert ^{2\theta}}\right)
^{\frac{1}{2}}\nonumber\\
&  \leq\left\vert a\right\vert ^{\left\lceil \kappa\right\rceil }%
\frac{\left\Vert \left(  \widehat{\cdot}D\right)  ^{\theta}\phi_{\infty
}\right\Vert _{\infty;\leq r_{3}}}{\theta!}\frac{1}{\sqrt{2\pi}\left\lceil
\kappa\right\rceil }\left(  \int\limits_{\left\vert \cdot\right\vert \leq
r_{3}}\frac{\left(  \widehat{a}\xi\right)  ^{2\left\lceil \kappa\right\rceil
}}{w}\right)  ^{\frac{1}{2}}+\frac{\left\vert a\right\vert ^{\kappa}}%
{\sqrt{2\pi}}\left(  2+\frac{1}{\left\lceil \kappa\right\rceil }\right)
\left(  \int\limits_{\left\vert \cdot\right\vert \geq r_{3}}\frac{\left\vert
\widehat{a}\xi\right\vert ^{2\kappa}}{w\left\vert \cdot\right\vert ^{2\theta}%
}\right)  ^{\frac{1}{2}}\label{a71}\\
&  \leq\left\vert a\right\vert ^{\left\lceil \kappa\right\rceil }%
\frac{\left\Vert \left(  \widehat{\cdot}D\right)  ^{\theta}\phi_{\infty
}\right\Vert _{\infty;\leq r_{3}}}{\theta!}\frac{r_{3}^{\left\lceil
\kappa\right\rceil }}{\sqrt{2\pi}\left\lceil \kappa\right\rceil }\left(
\int\limits_{\left\vert \cdot\right\vert \leq r_{3}}\frac{1}{w}\right)
^{\frac{1}{2}}+\frac{\left\vert a\right\vert ^{\kappa}}{\sqrt{2\pi}}\left(
2+\frac{1}{\left\lceil \kappa\right\rceil }\right)  \left(  \int%
\limits_{\left\vert \cdot\right\vert \geq r_{3}}\frac{\left\vert
\cdot\right\vert ^{2\kappa}}{w\left\vert \cdot\right\vert ^{2\theta}}\right)
^{\frac{1}{2}}\label{a661}\\
&  <\infty,\nonumber
\end{align}

which means that when $n=\left\lfloor \kappa\right\rfloor $,%
\begin{align}
&  \left\Vert \frac{\sqrt{2\pi}}{n!}\left(  \left(  ia\xi\right)
^{n+1}\overline{\widehat{g_{n}}}\left(  a\xi\right)  \phi_{\infty}%
f_{F}\right)  ^{\vee}\right\Vert _{\infty}\nonumber\\
&  \leq\left\vert f\right\vert _{w,\theta}\times\nonumber\\
&  \left(  \frac{\left\vert a\right\vert ^{\left\lceil \kappa\right\rceil }%
}{\left\lceil \kappa\right\rceil !}\frac{\left\Vert \left(  \widehat{\cdot
}D\right)  ^{\theta}\phi_{0}\right\Vert _{\infty}}{\theta!}\left(
\int\limits_{\left\vert \cdot\right\vert \leq r_{3}}\frac{\left(  \widehat
{a}\xi\right)  ^{2\left\lceil \kappa\right\rceil }}{w}\right)  ^{\frac{1}{2}%
}+\frac{\left\vert a\right\vert ^{\kappa}}{\sqrt{2\pi}}\left(  \frac
{2}{\left\lfloor \kappa\right\rfloor !}+\frac{1}{\left\lceil \kappa
\right\rceil !}\right)  \left(  \int\limits_{\left\vert \cdot\right\vert \geq
r_{3}}\frac{\left\vert \widehat{a}\xi\right\vert ^{2\kappa}}{w\left\vert
\cdot\right\vert ^{2\theta}}\right)  ^{\frac{1}{2}}\right) \label{a73}%
\end{align}

To summarize:%
\begin{align}
&  \left\Vert \frac{\sqrt{2\pi}}{n!}\left(  \left(  ia\xi\right)
^{n+1}\overline{\widehat{g_{n}}}\left(  a\xi\right)  \phi_{\infty}%
f_{F}\right)  ^{\vee}\right\Vert _{\infty}\nonumber\\
&  \leq\left\vert f\right\vert _{w,\theta}\left\{
\begin{array}
[c]{ll}%
\frac{\left\vert a\right\vert ^{n+1}}{\left(  n+1\right)  !}\left(
\frac{\left\Vert \left(  \widehat{\cdot}D\right)  ^{\theta}\phi_{\infty
}\right\Vert _{\infty;\leq r_{3}}}{\theta!}\left(  \int\limits_{\left\vert
\cdot\right\vert \leq r_{3}}\frac{\left(  \widehat{a}\xi\right)  ^{2\left(
n+1\right)  }}{w}\right)  ^{\frac{1}{2}}+\left(  \int\limits_{\left\vert
\cdot\right\vert \geq r_{3}}\frac{\left\vert \widehat{a}\xi\right\vert
^{2\left(  n+1\right)  }}{w\left\vert \cdot\right\vert ^{2\theta}}\right)
^{\frac{1}{2}}\right)  , & n\leq\left\lfloor \kappa\right\rfloor -1,\\
& \\
\frac{\left\vert a\right\vert ^{\left\lceil \kappa\right\rceil }}{\left\lceil
\kappa\right\rceil !}\frac{\left\Vert \left(  \widehat{\cdot}D\right)
^{\theta}\phi_{\infty}\right\Vert _{\infty;\leq r_{3}}}{\theta!}\left(
\int\limits_{\left\vert \cdot\right\vert \leq r_{3}}\frac{\left(  \widehat
{a}\xi\right)  ^{2\left\lceil \kappa\right\rceil }}{w}\right)  ^{\frac{1}{2}%
}+\left\vert a\right\vert ^{\kappa}\left(  \frac{2}{\left\lfloor
\kappa\right\rfloor !}+\frac{1}{\left\lceil \kappa\right\rceil !}\right)
\left(  \int\limits_{\left\vert \cdot\right\vert \geq r_{3}}\frac{\left\vert
\widehat{a}\xi\right\vert ^{2\kappa}}{w\left\vert \cdot\right\vert ^{2\theta}%
}\right)  ^{\frac{1}{2}}, & n=\left\lfloor \kappa\right\rfloor .
\end{array}
\right. \label{a72}%
\end{align}

\subsection{Remainder estimate for $f_{\rho}$\label{SbSect_rem_frho_W3.2}}

We have the following Taylor series remainder estimate for $f_{\rho}$: for
$n\leq\left\lfloor \kappa\right\rfloor $ the right side of \ref{a35} can be
estimated by \ref{a861} and \ref{a72} respectively as%
\begin{align}
&  \left\vert \mathcal{R}_{n+1}f_{\rho}\left(  x,a\right)  \right\vert
\nonumber\\
&  \leq\left\vert \frac{\sqrt{2\pi}}{n!}\left(  \overline{\widehat{g_{n}}%
}\left(  a\xi\right)  \left(  ia\xi\right)  \phi_{0}\widehat{\left(
aD\right)  ^{n}f_{\rho}}\right)  ^{\vee}\left(  x\right)  \right\vert
+\left\Vert \frac{\sqrt{2\pi}}{n!}\left(  \left(  ia\xi\right)  ^{n+1}%
\overline{\widehat{g_{n}}}\left(  a\xi\right)  \phi_{\infty}f_{F}\right)
^{\vee}\right\Vert _{\infty}\nonumber\\
&  \leq\left\vert f\right\vert _{w,\theta}\tfrac{\left\vert a\right\vert
^{n+1}}{\left(  n+1\right)  !}\left\{
\begin{array}
[c]{l}%
\frac{C_{n}^{\left(  \rho,w\right)  }}{\left(  2\pi\right)  ^{\frac{d}{2}}%
}\left(
\begin{array}
[c]{l}%
\left\Vert \widehat{a}D\widehat{\phi_{0}}\right\Vert _{1}+\\
\left(  \frac{3^{\theta-n+1}-1}{2}\right)  ^{\frac{1}{2}}\left(
\sum\limits_{j=0}^{\theta-n}\frac{\left(  \left\vert x\right\vert
^{2}+\left\vert a\right\vert ^{2}\right)  ^{j}}{j!}\right)  ^{\frac{1}{2}%
}\left(  \left\Vert \widehat{a}D\widehat{\phi_{0}}\right\Vert _{1}%
\sum\limits_{l=0}^{\theta-n}\frac{1}{l!}\left\Vert \left\vert \cdot\right\vert
^{2l}\widehat{a}D\widehat{\phi_{0}}\right\Vert _{1}\right)  ^{\frac{1}{2}}%
\end{array}
\right) \\
\qquad\qquad\qquad\qquad\qquad\qquad\qquad\qquad\qquad\qquad\qquad\qquad\qquad
if\text{ }n\leq\left\lfloor \kappa\right\rfloor <\theta,\\
\medskip\\
\frac{1}{\left(  2\pi\right)  ^{d}}\left(  \int\left(  \widehat{a}\widehat
{\xi}\right)  ^{2n}\frac{\left\vert \cdot\right\vert ^{2n}}{w\left\vert
\cdot\right\vert ^{2\theta}}\right)  ^{1/2}\left\Vert \widehat{a}%
D\widehat{\phi_{0}}\right\Vert _{1},\qquad\qquad\qquad\qquad\qquad if\text{
}\theta\leq n\leq\left\lfloor \kappa\right\rfloor ,
\end{array}
\right.  +\nonumber\\
&  \medskip\nonumber\\
&  +\left\vert f\right\vert _{w,\theta}\left\{
\begin{array}
[c]{r}%
\frac{\left\vert a\right\vert ^{n+1}}{\left(  n+1\right)  !}\left(
\frac{\left\Vert \left(  \widehat{\cdot}D\right)  ^{\theta}\phi_{\infty
}\right\Vert _{\infty;\leq r_{3}}}{\theta!}\left(  \int\limits_{\left\vert
\cdot\right\vert \leq r_{3}}\frac{\left(  \widehat{a}\xi\right)  ^{2\left(
n+1\right)  }}{w}\right)  ^{\frac{1}{2}}+\left(  \int\limits_{\left\vert
\cdot\right\vert \geq r_{3}}\frac{\left\vert \widehat{a}\xi\right\vert
^{2\left(  n+1\right)  }}{w\left\vert \cdot\right\vert ^{2\theta}}\right)
^{\frac{1}{2}}\right)  ,\qquad\qquad\\
if\text{ }n\leq\left\lfloor \kappa\right\rfloor -1,\\
\medskip\\
\frac{\left\vert a\right\vert ^{\left\lceil \kappa\right\rceil }}{\left\lceil
\kappa\right\rceil !}\frac{\left\Vert \left(  \widehat{\cdot}D\right)
^{\theta}\phi_{\infty}\right\Vert _{\infty;\leq r_{3}}}{\theta!}\left(
\int\limits_{\left\vert \cdot\right\vert \leq r_{3}}\frac{\left(  \widehat
{a}\xi\right)  ^{2\left\lceil \kappa\right\rceil }}{w}\right)  ^{\frac{1}{2}%
}+\left\vert a\right\vert ^{\kappa}\left(  \frac{2}{\left\lfloor
\kappa\right\rfloor !}+\frac{1}{\left\lceil \kappa\right\rceil !}\right)
\left(  \int\limits_{\left\vert \cdot\right\vert \geq r_{3}}\frac{\left\vert
\widehat{a}\xi\right\vert ^{2\kappa}}{w\left\vert \cdot\right\vert ^{2\theta}%
}\right)  ^{\frac{1}{2}},\\
if\text{ }n=\left\lfloor \kappa\right\rfloor ,
\end{array}
\right\}  ,\label{a65}%
\end{align}

where $C_{n}^{\left(  \rho,w\right)  }$ is given by \ref{a2.29}.

\subsection{Pointwise remainder estimate}

For $f\in X_{w}^{\theta}$, using \ref{a93} and then \ref{a35} we wrote%
\begin{align}
\left(  \mathcal{R}_{n+1}f\right)   &  \left(  x,a\right)  =\frac{\sqrt{2\pi}%
}{n!}\left(  \overline{\widehat{g_{n}}}\left(  a\xi\right)  \left(
ia\xi\right)  \phi_{0}\widehat{\left(  aD\right)  ^{n}f_{\rho}}\right)
^{\vee}\left(  x\right)  +\frac{\sqrt{2\pi}}{n!}\left(  \left(  ia\xi\right)
^{n+1}\overline{\widehat{g_{n}}}\left(  a\xi\right)  \phi_{\infty}%
f_{F}\right)  ^{\vee}\left(  x\right)  +\nonumber\\
&  \qquad\qquad+\left(  2\pi\right)  ^{-\frac{d}{2}}\left(  p_{\theta
-1;\widehat{f}}\left(  x+a\right)  -\sum_{\left\vert \beta\right\vert \leq
n}\frac{a^{\beta}}{\beta!}D^{\beta}p_{\theta-1;\widehat{f}}\left(  x\right)
\right)  ,\label{a53}%
\end{align}

and then estimated the \textbf{first term} on the right by \ref{a86}, the
\textbf{second term} by \ref{a72} and the \textbf{third (polynomial) term} by
\ref{a94}.

\begin{remark}
?? Base remarks on Remark \ref{Rem_Tay_expan_data_fns_W3.2} ??.

If $n<\left\lfloor \kappa\right\rfloor $ we have an estimate of order $n+1$ in
$\left\vert a\right\vert $ and when $n=\left\lfloor \kappa\right\rfloor $
\textbf{we always have a remainder estimate of order} $\kappa$ in $\left\vert
a\right\vert $.
\end{remark}

\subsection{The estimates \ref{a861} and \ref{a72} assuming radial
functions\label{SbSect_Tay_dat_rem_rad_fn}}

\smallskip

\fbox{\textbf{Estimate} \ref{a72}} In the light of Remark
\ref{Rem_Cor_integ_inprod_to_1dim} in the Appendix, \textbf{when the weight
function is radial}, say $w\left(  \xi\right)  =w_{\circ}\left(  \left\vert
\xi\right\vert \right)  $, we can use \ref{Ap141} i.e.
\[
\int\limits_{\left\vert x\right\vert \leq r}\left\vert \widehat{\xi
}x\right\vert ^{2p}f\left(  \left\vert x\right\vert \right)  dx=\frac{B\left(
\frac{d-1}{2},p+\frac{1}{2}\right)  }{B\left(  \frac{d}{2}+p,\frac{1}%
{2}\right)  }\int\limits_{\left\vert x\right\vert \leq r}\left\vert
x\right\vert ^{2p}f\left(  \left\vert x\right\vert \right)  dx,\quad
p>-\frac{1}{2},
\]

where $B$ is the beta function, to as%
\begin{align*}
&  \left\Vert \frac{\sqrt{2\pi}}{n!}\left(  \left(  ia\xi\right)
^{n+1}\overline{\widehat{g_{n}}}\left(  a\xi\right)  \phi_{\infty}%
f_{F}\right)  ^{\vee}\right\Vert _{\infty}\\
&  \leq\left\vert f\right\vert _{w,\theta}\left\{
\begin{array}
[c]{l}%
\frac{\left\vert a\right\vert ^{n+1}}{\left(  n+1\right)  !}\left(
\frac{\left\Vert \left(  \widehat{\cdot}D\right)  ^{\theta}\phi_{\infty
}\right\Vert _{\infty;\leq r_{3}}}{\theta!}\left(  \int\limits_{\left\vert
\cdot\right\vert \leq r_{3}}\frac{\left(  \widehat{a}\xi\right)  ^{2\left(
n+1\right)  }}{w}\right)  ^{\frac{1}{2}}+\left(  \int\limits_{\left\vert
\cdot\right\vert \geq r_{3}}\frac{\left\vert \widehat{a}\xi\right\vert
^{2\left(  n+1\right)  }}{w\left\vert \cdot\right\vert ^{2\theta}}\right)
^{\frac{1}{2}}\right)  ,\\
\qquad\qquad\qquad\qquad\qquad\qquad\qquad\qquad\qquad\qquad\qquad\qquad
if\text{ \ }n\leq\left\lfloor \kappa\right\rfloor -1,\\
\medskip\\
\frac{\left\vert a\right\vert ^{\left\lceil \kappa\right\rceil }}{\left\lceil
\kappa\right\rceil !}\frac{\left\Vert \left(  \widehat{\cdot}D\right)
^{\theta}\phi_{\infty}\right\Vert _{\infty;\leq r_{3}}}{\theta!}\left(
\int\limits_{\left\vert \cdot\right\vert \leq r_{3}}\frac{\left(  \widehat
{a}\xi\right)  ^{2\left\lceil \kappa\right\rceil }}{w}\right)  ^{\frac{1}{2}%
}+\left\vert a\right\vert ^{\kappa}\left(  \frac{2}{\left\lfloor
\kappa\right\rfloor !}+\frac{1}{\left\lceil \kappa\right\rceil !}\right)
\left(  \int\limits_{\left\vert \cdot\right\vert \geq r_{3}}\frac{\left\vert
\widehat{a}\xi\right\vert ^{2\kappa}}{w\left\vert \cdot\right\vert ^{2\theta}%
}\right)  ^{\frac{1}{2}},\\
\qquad\qquad\qquad\qquad\qquad\qquad\qquad\qquad\qquad\qquad\qquad\qquad
if\text{ \ }n=\left\lfloor \kappa\right\rfloor .
\end{array}
\right. \\
&  \medskip\\
&  =\left\vert f\right\vert _{w,\theta}\left\{
\begin{array}
[c]{l}%
\begin{array}
[c]{l}%
\frac{\left\vert a\right\vert ^{n+1}}{\left(  n+1\right)  !}\left(
\frac{B\left(  \frac{d}{2},\frac{1}{2}+n+1\right)  }{B\left(  \frac{1}%
{2},\frac{d}{2}+n+1\right)  }\right)  ^{\frac{1}{2}}\left(  \frac{\left\Vert
\left(  \widehat{\cdot}D\right)  ^{\theta}\phi_{\infty}\right\Vert
_{\infty;\leq r_{3}}}{\theta!}\left(  \int\limits_{\left\vert \cdot\right\vert
\leq r_{3}}\frac{\left\vert \xi\right\vert ^{2\left(  n+1\right)  }}%
{w}\right)  ^{\frac{1}{2}}+\left(  \int\limits_{\left\vert \cdot\right\vert
\geq r_{3}}\frac{\left\vert \xi\right\vert ^{2\left(  n+1\right)  }%
}{w\left\vert \cdot\right\vert ^{2\theta}}\right)  ^{\frac{1}{2}}\right)  ,\\
\qquad\qquad\qquad\qquad\qquad\qquad\qquad\qquad\qquad\qquad\qquad\qquad
if\text{ \ }n\leq\left\lfloor \kappa\right\rfloor -1,
\end{array}
\\
\medskip\\%
\begin{array}
[c]{l}%
\frac{\left\vert a\right\vert ^{\left\lceil \kappa\right\rceil }}{\left\lceil
\kappa\right\rceil !}\frac{\left\Vert \left(  \widehat{\cdot}D\right)
^{\theta}\phi_{\infty}\right\Vert _{\infty;\leq r_{3}}}{\theta!}\left(
\frac{B\left(  \frac{d}{2},\frac{1}{2}+\left\lceil \kappa\right\rceil \right)
}{B\left(  \frac{1}{2},\frac{d}{2}+\left\lceil \kappa\right\rceil \right)
}\right)  ^{\frac{1}{2}}\left(  \int\limits_{\left\vert \cdot\right\vert \leq
r_{3}}\frac{\left\vert \xi\right\vert ^{2\left\lceil \kappa\right\rceil }}%
{w}\right)  ^{\frac{1}{2}}+\\
\qquad+\left\vert a\right\vert ^{\kappa}\left(  \frac{2}{\left\lfloor
\kappa\right\rfloor !}+\frac{1}{\left\lceil \kappa\right\rceil !}\right)
\left(  \frac{B\left(  \frac{d}{2},\frac{1}{2}+\kappa\right)  }{B\left(
\frac{1}{2},\frac{d}{2}+\kappa\right)  }\right)  ^{\frac{1}{2}}\left(
\int\limits_{\left\vert \cdot\right\vert \geq r_{3}}\frac{\left\vert
\xi\right\vert ^{2\kappa}}{w\left\vert \cdot\right\vert ^{2\theta}}\right)
^{\frac{1}{2}},\quad n=\left\lfloor \kappa\right\rfloor .
\end{array}
\end{array}
\right.
\end{align*}

Further, by \ref{Ap148},%
\begin{align*}
\int_{\left\vert \cdot\right\vert \leq r_{3}}\frac{\left\vert \xi\right\vert
^{2\left(  n+1\right)  }}{w}  & =\omega_{d}\int_{0}^{r_{3}}\frac{t^{2n+d+1}%
}{w_{\circ}\left(  t\right)  }dt,\\
\int\limits_{\left\vert \cdot\right\vert \geq r_{3}}\frac{\left\vert
\xi\right\vert ^{2\left(  n+1\right)  }}{w\left\vert \cdot\right\vert
^{2\theta}}  & =\omega_{d}\int_{r_{3}}^{\infty}\frac{t^{2n+d+1}}{w_{\circ
}\left(  t\right)  \left\vert t\right\vert ^{2\theta}}dt,\\
\int\limits_{\left\vert \cdot\right\vert \leq r_{3}}\frac{\left\vert
\xi\right\vert ^{2\left\lceil \kappa\right\rceil }}{w}  & =\omega_{d}\int%
_{0}^{r_{3}}\frac{t^{2\left\lceil \kappa\right\rceil +d-1}}{w_{\circ}\left(
t\right)  }dt,\\
\int\limits_{\left\vert \cdot\right\vert \geq r_{3}}\frac{\left\vert
\xi\right\vert ^{2\kappa}}{w\left\vert \cdot\right\vert ^{2\theta}}  &
=\omega_{d}\int_{r_{3}}^{\infty}\frac{t^{2\kappa+d-1}}{w_{\circ}\left(
t\right)  \left\vert t\right\vert ^{2\theta}}dt.
\end{align*}

\textbf{Radial partition of unity} $\left\{  \phi_{0},\phi_{\infty}\right\}  $
given by \ref{a15}. Hence $\phi_{0}\in S_{1,\theta}$. Here $\phi_{\infty
}\left(  x\right)  =\left(  \phi_{\infty}\right)  _{\circ}\left(  \left\vert
x\right\vert \right)  $ and $\phi_{0}\left(  x\right)  =\left(  \phi
_{0}\right)  _{\circ}\left(  \left\vert x\right\vert \right)  $ and by part 1
of Lemma \ref{Lem_deriv_rad_funcs},
\begin{equation}
\left\Vert \left(  \widehat{\cdot}D\right)  ^{\theta}\phi_{\infty}\right\Vert
_{\infty;\leq r_{3}}=\left\Vert \left(  \widehat{\cdot}D\right)  ^{\theta}%
\phi_{0}\right\Vert _{\infty;\leq r_{3}}=\left\Vert D^{\theta}\left(  \phi
_{0}\right)  _{\circ}\left(  \left\vert x\right\vert \right)  \right\Vert
_{\infty;\leq r_{3}}=\max_{s\in\left[  0,r_{3}\right]  }\left\vert D^{\theta
}\left(  \phi_{0}\right)  _{\circ}\left(  s\right)  \right\vert .\label{a2.38}%
\end{equation}
\medskip

\fbox{\textbf{Estimate} \ref{a861}} Regarding $\left\Vert \widehat{a}%
D\widehat{\phi_{0}}\right\Vert _{1}$, since $\phi_{0}$ is radial,
$\widehat{\phi_{0}}$ is also radial, say
\begin{equation}
\widehat{\phi_{0}}\left(  \xi\right)  =\left(  \widehat{\phi_{0}}\right)
_{\circ}\left(  \left\vert \xi\right\vert \right)  :=\phi_{\mathcal{F}}\left(
\left\vert \xi\right\vert \right)  ,\label{a2.34}%
\end{equation}

and%
\[
\left(  \widehat{a}D\widehat{\phi_{0}}\right)  \left(  \xi\right)
=\widehat{a}D\left(  \phi_{\mathcal{F}}\left(  \left\vert \xi\right\vert
\right)  \right)  =\widehat{a}\widehat{\xi}\phi_{\mathcal{F}}^{\prime}\left(
\left\vert \xi\right\vert \right)  ,
\]

so that by \ref{Ap129}, Corollary \ref{Cor_integ_inprod_to_1dim},%
\begin{align*}
\left\Vert \widehat{a}D\widehat{\phi_{0}}\right\Vert _{1}  & =\int\left\vert
\widehat{a}\widehat{\xi}\right\vert \left\vert \phi_{\mathcal{F}}^{\prime
}\left(  \left\vert \xi\right\vert \right)  \right\vert d\xi=??\frac{B\left(
\frac{d}{2},1\right)  }{B\left(  \frac{1}{2},\frac{d+1}{2}\right)  }\omega
_{d}\int_{0}^{\infty}t^{d}\left\vert \phi_{\mathcal{F}}^{\prime}\left(
t\right)  \right\vert dt=\\
& =??\frac{\omega_{d}}{\pi}B\left(  \frac{d}{2},\frac{1}{2}\right)  \int%
_{0}^{\infty}t^{d}\left\vert \phi_{\mathcal{F}}^{\prime}\left(  t\right)
\right\vert dt.
\end{align*}

Regarding $\sum\limits_{l=0}^{\theta-n}\frac{1}{l!}\left\Vert \left\vert
\cdot\right\vert ^{2l}\widehat{a}D\widehat{\phi_{0}}\right\Vert _{1}$,%
\begin{align*}
\left\Vert \left\vert \cdot\right\vert ^{2l}\widehat{a}D\widehat{\phi_{0}%
}\right\Vert _{1}  & =\int\left\vert \cdot\right\vert ^{2l}\left\vert
\widehat{a}D\widehat{\phi_{0}}\right\vert =\int\widehat{a}\widehat{\xi
}\left\vert \cdot\right\vert ^{2l}\phi_{\mathcal{F}}^{\prime}\left(
\left\vert \xi\right\vert \right)  =\\
& =\frac{\omega_{d}}{\pi}B\left(  \frac{d}{2},\frac{1}{2}\right)  \int%
_{0}^{\infty}t^{2l+d}\left\vert \phi_{\mathcal{F}}^{\prime}\left(  t\right)
\right\vert dt,
\end{align*}

so that%
\[
\sum\limits_{l=0}^{\theta-n}\frac{1}{l!}\left\Vert \left\vert \cdot\right\vert
^{2l}\widehat{a}D\widehat{\phi_{0}}\right\Vert _{1}=\frac{\omega_{d}}{\pi
}B\left(  \frac{d}{2},\frac{1}{2}\right)  \sum\limits_{l=0}^{\theta-n}\frac
{1}{l!}\int_{0}^{\infty}t^{2l+d}\left\vert \phi_{\mathcal{F}}^{\prime}\left(
t\right)  \right\vert dt.
\]

Thus%
\begin{align}
&  \left(  \left\Vert \widehat{a}D\widehat{\phi_{0}}\right\Vert _{1}%
\sum\limits_{l=0}^{\theta-n}\frac{1}{l!}\left\Vert \left\vert \cdot\right\vert
^{2l}\widehat{a}D\widehat{\phi_{0}}\right\Vert _{1}\right)  ^{\frac{1}{2}%
}\nonumber\\
&  =\frac{\omega_{d}}{\pi}B\left(  \frac{d}{2},\frac{1}{2}\right)  \left(
\int_{0}^{\infty}t^{d}\left\vert \phi_{\mathcal{F}}^{\prime}\left(  t\right)
\right\vert dt\right)  ^{\frac{1}{2}}\left(  \sum\limits_{l=0}^{\theta-n}%
\frac{1}{l!}\int_{0}^{\infty}t^{2l+d}\left\vert \phi_{\mathcal{F}}^{\prime
}\left(  t\right)  \right\vert dt\right)  ^{\frac{1}{2}}.\label{a2.37}%
\end{align}
\medskip

\underline{\textbf{The constants} $C_{n}^{\left(  \rho,w\right)  }$} are given
by \ref{a2.29} as%
\[
C_{n}^{\left(  \rho,w\right)  }:=\max\left\{  \underline{C}_{\theta-n,r_{3}%
}^{\left(  \rho\right)  },\overline{C}_{\theta-n,r_{3}}^{\left(  \rho\right)
}\right\}  M_{n,r_{3}}^{\left(  w\right)  },
\]

where the $M_{n,r_{3}}^{\left(  w\right)  }$ are in turn given by \ref{a2.27}
as%
\[
\left.
\begin{array}
[c]{ll}%
M_{n,r_{3}}^{\left(  w\right)  }\left(  \widehat{a},\widehat{x}\right)  := &
\left(
\begin{array}
[c]{l}%
\int\limits_{\left\vert \cdot\right\vert \leq r_{3}}\frac{\left\vert
\widehat{a}\widehat{\xi}\right\vert ^{2n}}{w}+\frac{1}{\left(  \theta
-n\right)  !}\int\limits_{\left\vert \cdot\right\vert \leq r_{3}}%
\frac{\left\vert \widehat{x}\widehat{\xi}\right\vert ^{2\left(  \theta
-n\right)  }\left\vert \widehat{a}\widehat{\xi}\right\vert ^{2n}}{w}+\\
+\int\limits_{\left\vert \cdot\right\vert \geq r_{3}}\left\vert \widehat
{a}\widehat{\xi}\right\vert ^{2n}\frac{\left\vert \cdot\right\vert ^{2n}%
}{w\left\vert \cdot\right\vert ^{2\theta}}+\sum\limits_{k<\theta-n}\frac
{1}{k!}\int\limits_{\left\vert \cdot\right\vert \geq r_{3}}\left\vert
\widehat{x}\widehat{\xi}\right\vert ^{2k}\left\vert \widehat{a}\widehat{\xi
}\right\vert ^{2n}\frac{\left\vert \cdot\right\vert ^{2n}}{w\left\vert
\cdot\right\vert ^{2\theta}}%
\end{array}
\right)  ^{\frac{1}{2}},\\
& or\text{ }more\text{ }weakly\\
M_{n,r_{3}}^{\left(  w\right)  }\left(  \widehat{a},\widehat{x}\right)  := &
\left(
\begin{array}
[c]{c}%
\int\limits_{\left\vert \cdot\right\vert \leq r_{3}}\frac{\left\vert
\widehat{a}\widehat{\xi}\right\vert ^{2n}}{w}+\frac{1}{\left(  \theta
-n\right)  !}\min\left\{  \int\limits_{\left\vert \cdot\right\vert \leq r_{3}%
}\frac{\left\vert \widehat{x}\widehat{\xi}\right\vert ^{2\left(
\theta-n\right)  }}{w},\int\limits_{\left\vert \cdot\right\vert \leq r_{3}%
}\frac{\left\vert \widehat{a}\widehat{\xi}\right\vert ^{2n}}{w}\right\}  +\\
+\int\limits_{\left\vert \cdot\right\vert \geq r_{3}}\left\vert \widehat
{a}\widehat{\xi}\right\vert ^{2n}\frac{\left\vert \cdot\right\vert ^{2n}%
}{w\left\vert \cdot\right\vert ^{2\theta}}+\\
+\sum\limits_{k<\theta-n}\frac{1}{k!}\min\left\{  \int\limits_{\left\vert
\cdot\right\vert \geq r_{3}}\left\vert \widehat{a}\widehat{\xi}\right\vert
^{2n}\frac{\left\vert \cdot\right\vert ^{2n}}{w\left\vert \cdot\right\vert
^{2\theta}},\int\limits_{\left\vert \cdot\right\vert \geq r_{3}}\left\vert
\widehat{x}\widehat{\xi}\right\vert ^{2k}\frac{\left\vert \cdot\right\vert
^{2n}}{w\left\vert \cdot\right\vert ^{2\theta}}\right\}
\end{array}
\right)  ^{\frac{1}{2}},\\
& or\text{ }more\text{ }weakly\\
M_{n,r_{3}}^{\left(  w\right)  }\left(  \widehat{a}\right)  := & \sqrt
{1+e}\left(  \int\limits_{\left\vert \cdot\right\vert \leq r_{3}}\left\vert
\widehat{a}\widehat{\xi}\right\vert ^{2n}\frac{1}{w}+\int\limits_{\left\vert
\cdot\right\vert \geq r_{3}}\left\vert \widehat{a}\widehat{\xi}\right\vert
^{2n}\frac{\left\vert \cdot\right\vert ^{2n}}{w\left\vert \cdot\right\vert
^{2\theta}}\right)  ^{1/2},\\
& or\text{ }more\text{ }weakly\\
M_{n,r_{3}}^{\left(  w\right)  }:= & \sqrt{1+e}\left(  \int\limits_{\left\vert
\cdot\right\vert \leq r_{3}}\frac{1}{w}+\int\limits_{\left\vert \cdot
\right\vert \geq r_{3}}\frac{\left\vert \cdot\right\vert ^{2n}}{w\left\vert
\cdot\right\vert ^{2\theta}}\right)  ^{1/2}.
\end{array}
\right\}
\]

?? FINISH!

\underline{\textbf{Examples of radial }$\phi_{0}$} First we need:\medskip

\begin{theorem}
\label{Thm_C1n_from_C1_1}\textbf{Construction of }$C_{1;n}^{\infty}%
$\textbf{\ functions from }$C_{\mathbf{1};1}^{\infty}$\textbf{\ functions}

Assume $n\geq2$.

\begin{enumerate}
\item If $\psi\in C_{\mathbf{1};1}^{\infty}$ and $\sigma\in C_{1;n-1}^{\infty
}$ then $\phi\in C_{1;n}^{\infty}$, where $\phi$ is defined by%
\[
\phi\left(  x\right)  =\psi\left(  x\right)  -\sigma\left(  x\right)
\sum\limits_{\left\vert \alpha\right\vert =1}^{n-1}\frac{D^{\alpha}\psi\left(
0\right)  }{\alpha!}x^{\alpha}\in C_{1;n}^{\infty}.
\]

If $\psi\in S_{\mathbf{1};1}$ and $\sigma\in S_{1;n-1}$ then $\phi\in S_{1;n}$.

\item If $\psi$ and $\sigma$ are radial then $\phi$ is radial.
\end{enumerate}
\end{theorem}

\begin{proof}
?? \textbf{Part 1} \textbf{FINISH}! ??\medskip

\textbf{Part 2} From Theorem \ref{Thm_Taylor_rad_remain}, $\sum
\limits_{\left\vert \alpha\right\vert =k}\frac{D^{\alpha}\psi\left(  0\right)
}{\alpha!}x^{\alpha}$ is radial and so $\phi$ is radial.
\end{proof}

Using the last theorem \textbf{we now construct a member of }$C_{1;n}^{\infty
}\left(  \mathbb{R}^{d}\right)  $\textbf{\ from any function }$\psi
$\textbf{\ in} $C_{\mathbf{1};1}^{\infty}\left(  \mathbb{R}^{d}\right)  $.
This construction allows the Fourier transform of this member to be calculated
if the Fourier transform of $\psi$ can be calculated.

Now we know that if $\psi\in C_{\mathbf{1};1}^{\infty}$ then%
\[
\eta^{\left(  1\right)  }\left(  x\right)  =\psi\left(  x\right)  -\psi\left(
x\right)  \sum\limits_{\left\vert \alpha^{\left(  1\right)  }\right\vert
=1}\frac{D^{\alpha^{\left(  1\right)  }}\psi\left(  0\right)  }{\alpha
^{\left(  1\right)  }!}x^{\alpha^{\left(  1\right)  }}=\psi\left(  x\right)
\left(  1-\sum\limits_{\left\vert \alpha^{\left(  1\right)  }\right\vert
=1}\frac{D^{\alpha^{\left(  1\right)  }}\psi\left(  0\right)  }{\alpha
^{\left(  1\right)  }!}x^{\alpha^{\left(  1\right)  }}\right)  \in
C_{\mathbf{1};2}^{\infty}.
\]

Further, $\eta^{\left(  2\right)  }\in C_{\mathbf{1};3}^{\infty}$ where%
\begin{align*}
\eta^{\left(  2\right)  }\left(  x\right)   & =\psi\left(  x\right)
-\eta^{\left(  1\right)  }\left(  x\right)  \sum\limits_{\left\vert
\alpha^{\left(  2\right)  }\right\vert =1}^{2}\frac{D^{\alpha^{\left(
2\right)  }}\psi\left(  0\right)  }{\alpha^{\left(  2\right)  }!}%
x^{\alpha^{\left(  2\right)  }}\\
& =\psi\left(  x\right)  -\psi\left(  x\right)  \left(  1-\sum
\limits_{\left\vert \alpha^{\left(  1\right)  }\right\vert =1}\frac
{D^{\alpha^{\left(  1\right)  }}\psi\left(  0\right)  }{\alpha^{\left(
1\right)  }!}x^{\alpha^{\left(  1\right)  }}\right)  \sum\limits_{\left\vert
\alpha^{\left(  2\right)  }\right\vert =1}^{2}\frac{D^{\alpha^{\left(
2\right)  }}\psi\left(  0\right)  }{\alpha^{\left(  2\right)  }!}%
x^{\alpha^{\left(  2\right)  }}\\
& =\psi\left(  x\right)  \left(  1-\left(  1-\sum\limits_{\left\vert
\alpha^{\left(  1\right)  }\right\vert =1}\frac{D^{\alpha^{\left(  1\right)
}}\psi\left(  0\right)  }{\alpha^{\left(  1\right)  }!}x^{\alpha^{\left(
1\right)  }}\right)  \sum\limits_{\left\vert \alpha^{\left(  2\right)
}\right\vert =1}^{2}\frac{D^{\alpha^{\left(  2\right)  }}\psi\left(  0\right)
}{\alpha^{\left(  2\right)  }!}x^{\alpha^{\left(  2\right)  }}\right)  .
\end{align*}

Continuing, $\eta^{\left(  3\right)  }\in C_{\mathbf{1};4}^{\infty}$ where
\begin{align*}
&  \eta^{\left(  3\right)  }\left(  x\right) \\
&  =\psi\left(  x\right)  -\eta^{\left(  2\right)  }\left(  x\right)
\sum\limits_{\left\vert \alpha^{\left(  3\right)  }\right\vert =1}^{3}%
\frac{D^{\alpha^{\left(  3\right)  }}\psi\left(  0\right)  }{\alpha^{\left(
3\right)  }!}x^{\alpha^{\left(  3\right)  }}\\
&  =\psi\left(  x\right)  -\psi\left(  x\right)  \left(  1-\left(
1-\sum\limits_{\left\vert \alpha^{\left(  1\right)  }\right\vert =1}%
\frac{D^{\alpha^{\left(  1\right)  }}\psi\left(  0\right)  }{\alpha^{\left(
1\right)  }!}x^{\alpha^{\left(  1\right)  }}\right)  \sum\limits_{\left\vert
\alpha^{\left(  2\right)  }\right\vert =1}^{2}\frac{D^{\alpha^{\left(
2\right)  }}\psi\left(  0\right)  }{\alpha^{\left(  2\right)  }!}%
x^{\alpha^{\left(  2\right)  }}\right)  \times\\
&  \times\sum\limits_{\left\vert \alpha^{\left(  3\right)  }\right\vert
=1}^{3}\frac{D^{\alpha^{\left(  3\right)  }}\psi\left(  0\right)  }%
{\alpha^{\left(  3\right)  }!}x^{\alpha^{\left(  3\right)  }}\\
&  =\psi\left(  x\right)  \left(
\begin{array}
[c]{r}%
1-\left(  1-\left(  1-\sum\limits_{\left\vert \alpha^{\left(  1\right)
}\right\vert =1}\frac{D^{\alpha^{\left(  1\right)  }}\psi\left(  0\right)
}{\alpha^{\left(  1\right)  }!}x^{\alpha^{\left(  1\right)  }}\right)
\sum\limits_{\left\vert \alpha^{\left(  2\right)  }\right\vert =1}^{2}%
\frac{D^{\alpha^{\left(  2\right)  }}\psi\left(  0\right)  }{\alpha^{\left(
2\right)  }!}x^{\alpha^{\left(  2\right)  }}\right)  \times\\
\times\sum\limits_{\left\vert \alpha^{\left(  3\right)  }\right\vert =1}%
^{3}\frac{D^{\alpha^{\left(  3\right)  }}\psi\left(  0\right)  }%
{\alpha^{\left(  3\right)  }!}x^{\alpha^{\left(  3\right)  }}%
\end{array}
\right) \\
&  =\psi\left(  x\right)  \left(
\begin{array}
[c]{r}%
1-\sum\limits_{\left\vert \alpha^{\left(  3\right)  }\right\vert =1}^{3}%
\frac{D^{\alpha^{\left(  3\right)  }}\psi\left(  0\right)  }{\alpha^{\left(
3\right)  }!}x^{\alpha^{\left(  3\right)  }}+\sum\limits_{\left\vert
\alpha^{\left(  2\right)  }\right\vert =1}^{2}\sum\limits_{\left\vert
\alpha^{\left(  3\right)  }\right\vert =1}^{3}\frac{D^{\alpha^{\left(
2\right)  }}\psi\left(  0\right)  \text{ }D^{\alpha^{\left(  3\right)  }}%
\psi\left(  0\right)  }{\alpha^{\left(  2\right)  }!\alpha^{\left(  3\right)
}!}x^{\alpha^{\left(  2\right)  }+\alpha^{\left(  3\right)  }}-\\
-\sum\limits_{\left\vert \alpha^{\left(  1\right)  }\right\vert =1}%
\sum\limits_{\left\vert \alpha^{\left(  2\right)  }\right\vert =1}^{2}%
\sum\limits_{\left\vert \alpha^{\left(  3\right)  }\right\vert =1}^{3}%
\frac{D^{\alpha^{\left(  1\right)  }}\psi\left(  0\right)  \text{ }%
D^{\alpha^{\left(  2\right)  }}\psi\left(  0\right)  \text{ }D^{\alpha
^{\left(  3\right)  }}\psi\left(  0\right)  }{\alpha^{\left(  1\right)
}!\alpha^{\left(  2\right)  }!\alpha^{\left(  3\right)  }!}x^{\alpha^{\left(
1\right)  }+\alpha^{\left(  2\right)  }+\alpha^{\left(  3\right)  }}%
\end{array}
\right) \\
&  =\psi\left(  x\right)  q^{\left(  3\right)  }\left(  x\right)  ,
\end{align*}

where $q^{\left(  3\right)  }\left(  x\right)  $ is the polynomial%
\[
q^{\left(  3\right)  }\left(  x\right)  =\left\{
\begin{array}
[c]{r}%
1-\sum\limits_{\left\vert \alpha^{\left(  3\right)  }\right\vert =1}^{3}%
\frac{D^{\alpha^{\left(  3\right)  }}\psi\left(  0\right)  }{\alpha^{\left(
3\right)  }!}x^{\alpha^{\left(  3\right)  }}+\sum\limits_{\left\vert
\alpha^{\left(  2\right)  }\right\vert =1}^{2}\sum\limits_{\left\vert
\alpha^{\left(  3\right)  }\right\vert =1}^{3}\frac{D^{\alpha^{\left(
2\right)  }}\psi\left(  0\right)  \text{ }D^{\alpha^{\left(  3\right)  }}%
\psi\left(  0\right)  }{\alpha^{\left(  2\right)  }!\alpha^{\left(  3\right)
}!}x^{\alpha^{\left(  2\right)  }+\alpha^{\left(  3\right)  }}-\\
-\sum\limits_{\left\vert \alpha^{\left(  1\right)  }\right\vert =1}%
\sum\limits_{\left\vert \alpha^{\left(  2\right)  }\right\vert =1}^{2}%
\sum\limits_{\left\vert \alpha^{\left(  3\right)  }\right\vert =1}^{3}%
\frac{D^{\alpha^{\left(  1\right)  }}\psi\left(  0\right)  \text{ }%
D^{\alpha^{\left(  2\right)  }}\psi\left(  0\right)  \text{ }D^{\alpha
^{\left(  3\right)  }}\psi\left(  0\right)  }{\alpha^{\left(  1\right)
}!\alpha^{\left(  2\right)  }!\alpha^{\left(  3\right)  }!}x^{\alpha^{\left(
1\right)  }+\alpha^{\left(  2\right)  }+\alpha^{\left(  3\right)  }}%
\end{array}
\right.  ,
\]

with degree $\deg q^{\left(  3\right)  }\leq3$. Thus%
\[
\widehat{\eta^{\left(  3\right)  }}=q^{\left(  3\right)  }\left(  iD\right)
\widehat{\psi}.
\]

In general%
\[
\eta^{\left(  n\right)  }\left(  x\right)  =q^{\left(  n\right)  }\left(
x\right)  \psi\left(  x\right)  \in C_{\mathbf{1};n+1}^{\infty},
\]

where%
\begin{align}
&  q^{\left(  n\right)  }\left(  x\right) \nonumber\\
&  =\left\{
\begin{array}
[c]{l}%
1-\sum\limits_{\left\vert \alpha^{\left(  n\right)  }\right\vert =1}^{n}%
\frac{D^{\alpha^{\left(  n\right)  }}\psi\left(  0\right)  }{\alpha^{\left(
n\right)  }!}x^{\alpha^{\left(  n\right)  }}+\sum\limits_{\left\vert
\alpha^{\left(  n-1\right)  }\right\vert =1}^{n-1}\sum\limits_{\left\vert
\alpha^{\left(  n\right)  }\right\vert =1}^{n}\frac{D^{\alpha^{\left(
n-1\right)  }}\psi\left(  0\right)  \text{ }D^{\alpha^{\left(  n\right)  }%
}\psi\left(  0\right)  }{\alpha^{\left(  n-1\right)  }!\alpha^{\left(
n\right)  }!}x^{\alpha^{\left(  n-1\right)  }+\alpha^{\left(  n\right)  }%
}-\ldots\\
+\left(  -1\right)  ^{k+1}\sum\limits_{\left\vert \alpha^{\left(  n-k\right)
}\right\vert =1}\ldots\sum\limits_{\left\vert \alpha^{\left(  n-1\right)
}\right\vert =1}^{n-1}\sum\limits_{\left\vert \alpha^{\left(  n\right)
}\right\vert =1}^{n}\frac{D^{\alpha^{\left(  n-k\right)  }}\psi\left(
0\right)  \ldots D^{\alpha^{\left(  n-1\right)  }}\psi\left(  0\right)  \text{
}D^{\alpha^{\left(  n\right)  }}\psi\left(  0\right)  }{\alpha^{\left(
n-k\right)  }!\ldots\alpha^{\left(  n-1\right)  }!\alpha^{\left(  n\right)
}!}x^{\alpha^{\left(  n-k\right)  }+\ldots+\alpha^{\left(  n-1\right)
}+\alpha^{\left(  n\right)  }}+\ldots\\
+\left(  -1\right)  ^{n}\sum\limits_{\left\vert \alpha^{\left(  1\right)
}\right\vert =1}\ldots\sum\limits_{\left\vert \alpha^{\left(  n-1\right)
}\right\vert =1}^{n-1}\sum\limits_{\left\vert \alpha^{\left(  n\right)
}\right\vert =1}^{n}\frac{D^{\alpha^{\left(  1\right)  }}\psi\left(  0\right)
\ldots D^{\alpha^{\left(  n-1\right)  }}\psi\left(  0\right)  \text{
}D^{\alpha^{\left(  n\right)  }}\psi\left(  0\right)  }{\alpha^{\left(
1\right)  }!\ldots\alpha^{\left(  n-1\right)  }!\alpha^{\left(  n\right)  }%
!}x^{\alpha^{\left(  1\right)  }+\ldots+\alpha^{\left(  n-1\right)  }%
+\alpha^{\left(  n\right)  }},
\end{array}
\right. \label{a04}%
\end{align}

and for this polynomial $\deg q^{\left(  n\right)  }\leq\frac{1}{2}n\left(
n+1\right)  $.

We will now use part 3 of Theorem \ref{Thm_Taylor_rad_remain} of the Appendix
to rewrite $q^{\left(  n\right)  }\left(  x\right)  $ in terms of $\left\vert
x\right\vert $ \textbf{when }$\psi$\textbf{\ is radial}. From \ref{Ap133},
\[
\sum\limits_{\left\vert \beta\right\vert =m}\frac{D^{\beta}\psi\left(
0\right)  }{\beta!}x^{\beta}=c_{m}\left\vert x\right\vert ^{m},
\]

so%
\begin{equation}
\sum\limits_{\left\vert \beta\right\vert =1}^{j}\frac{D^{\beta}\psi\left(
0\right)  }{\beta!}x^{\beta}=\sum\limits_{m=1}^{j}c_{m}\left\vert x\right\vert
^{m}:=d_{j}\left(  \left\vert x\right\vert \right)  ,\label{a006}%
\end{equation}

so \ref{a04} becomes%
\begin{align}
q^{\left(  n\right)  }\left(  x\right)   & =1-d_{n}\left(  \left\vert
x\right\vert \right)  +d_{n-1}\left(  \left\vert x\right\vert \right)
d_{n}\left(  \left\vert x\right\vert \right)  -\ldots\nonumber\\
& \ldots+\left(  -1\right)  ^{k+1}\left(  d_{n-k}\left(  \left\vert
x\right\vert \right)  \ldots d_{n-1}\left(  \left\vert x\right\vert \right)
d_{n}\left(  \left\vert x\right\vert \right)  \right)  +\ldots\\
& \ldots+\left(  -1\right)  ^{n}\left(  d_{1}\left(  \left\vert x\right\vert
\right)  \ldots d_{n-1}\left(  \left\vert x\right\vert \right)  d_{n}\left(
\left\vert x\right\vert \right)  \right)  ,\label{a007}%
\end{align}

and hence \textbf{when }$\psi$\textbf{\ is radial}%
\begin{equation}
q_{\circ}^{\left(  n\right)  }=1-d_{n}+d_{n-1}d_{n}-\ldots+\left(  -1\right)
^{k+1}d_{n-k}\ldots d_{n-1}d_{n}+\ldots+\left(  -1\right)  ^{n}d_{1}\ldots
d_{n-1}d_{n}.\label{a016}%
\end{equation}

?? \textbf{SEE} \ref{a017} below ??\medskip

\underline{If $c_{m}=0$ when $m$ is odd} then%
\begin{align*}
\widehat{\eta^{\left(  n\right)  }}  & =q^{\left(  n\right)  }\left(
iD\right)  \widehat{\psi}\\
& =\left(
\begin{array}
[c]{c}%
1-d_{n}\left(  \left\vert iD\right\vert \right)  +d_{n-1}\left(  \left\vert
iD\right\vert \right)  d_{n}\left(  \left\vert iD\right\vert \right)  -\\
\vdots\\
+\left(  -1\right)  ^{k+1}d_{n-k}\left(  \left\vert iD\right\vert \right)
\ldots d_{n-1}\left(  \left\vert iD\right\vert \right)  d_{n}\left(
\left\vert iD\right\vert \right)  +\\
\vdots\\
+\left(  -1\right)  ^{n}d_{1}\left(  \left\vert iD\right\vert \right)  \ldots
d_{n-1}\left(  \left\vert iD\right\vert \right)  d_{n}\left(  \left\vert
iD\right\vert \right)
\end{array}
\right)  \widehat{\psi},
\end{align*}

where%
\[
d_{j}\left(  \left\vert iD\right\vert \right)  \widehat{\psi}=\sum
\limits_{\substack{m=1 \\m\text{ }even}}^{j}c_{m}\left\vert iD\right\vert
^{m}\widehat{\psi}=\sum\limits_{\substack{m=1 \\m\text{ }even}}^{j}\left(
-1\right)  ^{m/2}c_{m}\left(  \left\vert D\right\vert ^{2}\right)
^{m/2}\widehat{\psi}.
\]

\underline{Calculation of $\left\vert D\right\vert ^{2n}\left(  f\left(
\left\vert x\right\vert \right)  \right)  $} We can use Lemma
\ref{Lem_iterLaplacian_rad} i.e.%
\begin{equation}%
\begin{array}
[c]{c}%
\left\vert D\right\vert ^{2n}\left(  f\left(  \left\vert x\right\vert \right)
\right)  =\left(  \mathcal{L}^{n}f\right)  \left(  \left\vert x\right\vert
\right)  ,\\
where\\
\left(  \mathcal{L}f\right)  \left(  s\right)  =\left(  d-1\right)
s^{-1}f\left(  s\right)  +D^{2}f\left(  s\right)  .
\end{array}
\label{a020}%
\end{equation}

\begin{equation}
Hence\text{ }\left\vert D\right\vert ^{2n}\left(  f\left(  \left\vert
x\right\vert \right)  \right)  \text{ }is\text{ }a\text{ }radial\text{
}function.\label{a022}%
\end{equation}
\medskip

\fbox{\textbf{Example} $\psi=e^{-\left\vert \cdot\right\vert ^{2}/2}\in
S_{1;1}$} Here $\widehat{e^{-\left\vert x\right\vert ^{2}/2}}=e^{-\left\vert
x\right\vert ^{2}/2}\in S$; $c_{m}=0$ when $m$ is odd;%
\begin{equation}
\eta^{\left(  \theta-1\right)  }\left(  x\right)  =q^{\left(  \theta-1\right)
}\left(  x\right)  e^{-\left\vert x\right\vert ^{2}/2}\in S_{1;\theta
},\label{a011}%
\end{equation}

and%
\[
\widehat{\eta^{\left(  \theta-1\right)  }}=q^{\left(  \theta-1\right)
}\left(  iD\right)  \widehat{e^{-\left\vert \cdot\right\vert ^{2}/2}%
}=q^{\left(  \theta-1\right)  }\left(  iD\right)  e^{-\left\vert
\cdot\right\vert ^{2}/2}.
\]

Choose%
\begin{equation}
\phi_{0}\left(  x\right)  =\eta^{\left(  \theta-1\right)  }\left(  x\right)
=q_{\circ}^{\left(  \theta-1\right)  }\left(  \left\vert x\right\vert \right)
e^{-\left\vert x\right\vert ^{2}/2}\in S_{1;\theta}.\label{a010}%
\end{equation}

From the Taylor series expansion of $e^{-\left\vert x\right\vert ^{2}/2}$
about the origin%
\begin{equation}
\sum\limits_{\left\vert \beta\right\vert =1}^{j}\frac{D^{\beta}\psi\left(
0\right)  }{\beta!}x^{\beta}=\sum\limits_{k=1}^{\left\lfloor j/2\right\rfloor
}\frac{\left\vert x\right\vert ^{2k}}{2^{k}k!}=d_{j}\left(  \left\vert
x\right\vert \right)  ,\label{a012}%
\end{equation}

so that%
\[
c_{m}=\frac{1}{2^{m/2}\left(  m/2\right)  !},\text{ }m\text{ }even,
\]

and hence%
\[
d_{j}\left(  s\right)  =\sum\limits_{m=1}^{j}c_{m}s^{m}=\sum
\limits_{\substack{m=1 \\m\text{ }even}}^{j}c_{m}s^{m}=\sum\limits_{k=1}%
^{\left\lfloor j/2\right\rfloor }c_{2k}s^{2k}=\sum\limits_{k=1}^{\left\lfloor
j/2\right\rfloor }\frac{s^{2k}}{2^{k}k!}=\sum\limits_{k=1}^{\left\lfloor
j/2\right\rfloor }\frac{\left(  s/\sqrt{2}\right)  ^{2k}}{k!}<e^{s^{2}/2}.
\]

Thus%
\[
\phi_{0}\left(  x\right)  =??
\]

Thus%
\begin{equation}
\widehat{\phi_{0}}=\left(
\begin{array}
[c]{c}%
1-d_{\theta-1}\left(  \left\vert iD\right\vert \right)  +d_{\theta-2}\left(
\left\vert iD\right\vert \right)  d_{\theta-1}\left(  \left\vert iD\right\vert
\right)  -\\
\vdots\\
+\left(  -1\right)  ^{k+1}d_{\theta-1-k}\left(  \left\vert iD\right\vert
\right)  \ldots d_{\theta-2}\left(  \left\vert iD\right\vert \right)
d_{\theta-1}\left(  \left\vert iD\right\vert \right)  +\\
\vdots\\
+\left(  -1\right)  ^{\theta-1}d_{1}\left(  \left\vert iD\right\vert \right)
\ldots d_{\theta-2}\left(  \left\vert iD\right\vert \right)  d_{\theta
-1}\left(  \left\vert iD\right\vert \right)
\end{array}
\right)  e^{-\left\vert \cdot\right\vert ^{2}/2},\label{a009}%
\end{equation}

and%
\begin{align}
d_{j}\left(  \left\vert iD\right\vert \right)  e^{-\left\vert \cdot\right\vert
^{2}/2}=\sum\limits_{\substack{m=1 \\m\text{ }even}}^{j}\frac{\left\vert
iD\right\vert ^{m}e^{-\left\vert \cdot\right\vert ^{2}/2}}{\left(  m/2\right)
!2^{m/2}}  & =\sum\limits_{\substack{m=1 \\m\text{ }even }}^{j}\frac{\left(
-1\right)  ^{m/2}}{\left(  m/2\right)  !2^{m/2}}\left(  \left\vert
D\right\vert ^{2}\right)  ^{m/2}e^{-\left\vert \cdot\right\vert ^{2}%
/2}\nonumber\\
& =\sum\limits_{n=1}^{\left\lfloor j/2\right\rfloor }\frac{\left(  -1\right)
^{n}}{2^{n}n!}\left\vert D\right\vert ^{2n}e^{-\left\vert \cdot\right\vert
^{2}/2}.\label{a008}%
\end{align}

This enables us to write%
\begin{align*}
d_{k}\left(  \left\vert iD\right\vert \right)  d_{j}\left(  \left\vert
iD\right\vert \right)  e^{-\left\vert \cdot\right\vert ^{2}/2}  &
=d_{k}\left(  \left\vert iD\right\vert \right)  \sum\limits_{n=1}%
^{\left\lfloor j/2\right\rfloor }\frac{\left(  -1\right)  ^{n}}{2^{n}%
n!}\left\vert D\right\vert ^{2n}e^{-\left\vert \cdot\right\vert ^{2}/2}\\
& =\sum\limits_{n=1}^{\left\lfloor j/2\right\rfloor }\frac{\left(  -1\right)
^{n}}{2^{n}n!}\left\vert D\right\vert ^{2n}d_{k}\left(  \left\vert
iD\right\vert \right)  e^{-\left\vert \cdot\right\vert ^{2}/2}\\
& =\sum\limits_{n=1}^{\left\lfloor j/2\right\rfloor }\frac{\left(  -1\right)
^{n}}{2^{n}n!}\left\vert D\right\vert ^{2n}\left(  \sum\limits_{m=1}%
^{\left\lfloor k/2\right\rfloor }\frac{\left(  -1\right)  ^{m}}{2^{m}%
m!}\left\vert D\right\vert ^{2m}e^{-\left\vert \cdot\right\vert ^{2}/2}\right)
\\
& =\sum\limits_{n=1}^{\left\lfloor j/2\right\rfloor }\sum\limits_{m=1}%
^{\left\lfloor k/2\right\rfloor }\frac{\left(  -1\right)  ^{n+m}}{2^{n+m}%
n!m!}\left\vert D\right\vert ^{2\left(  n+m\right)  }e^{-\left\vert
\cdot\right\vert ^{2}/2}\\
& =\sum\limits_{\alpha=\mathbf{1}_{2}}^{\left\lfloor \left(  j,k\right)
/2\right\rfloor }\frac{\left(  -1\right)  ^{\left\vert \alpha\right\vert }%
}{2^{\left\vert \alpha\right\vert }\alpha!}\left\vert D\right\vert
^{2\left\vert \alpha\right\vert }e^{-\left\vert \cdot\right\vert ^{2}/2},
\end{align*}

and if $\theta^{\left(  k\right)  }:=\left(  \theta-1-k,\theta-k,\ldots
,\theta-1\right)  $ for $k=0,1,2,\ldots,\theta-2$, then%
\begin{equation}
\left(  -1\right)  ^{k+1}d_{\theta-1-k}\left(  \left\vert iD\right\vert
\right)  \ldots d_{\theta-2}\left(  \left\vert iD\right\vert \right)
d_{\theta-1}\left(  \left\vert iD\right\vert \right)  =\left(  -1\right)
^{k+1}\sum\limits_{\alpha=\mathbf{1}}^{\left\lfloor \theta^{\left(  k\right)
}/2\right\rfloor }\frac{\left(  -1\right)  ^{\left\vert \alpha\right\vert }%
}{2^{\left\vert \alpha\right\vert }\alpha!}\left\vert D\right\vert
^{2\left\vert \alpha\right\vert }e^{-\left\vert \cdot\right\vert ^{2}%
/2}.\label{a017}%
\end{equation}

\underline{Closed formula for $\left\vert D\right\vert ^{2n}e^{-\left\vert
\cdot\right\vert ^{2}/2}$ in \ref{a008}} The Hermite polynomials $H_{k}$ on
$\mathbb{R}^{1}$ are defined by%
\[
D_{t}^{k}e^{-t^{2}}=\left(  -1\right)  ^{k}H_{k}\left(  t\right)  e^{-t^{2}%
},\quad t\in\mathbb{R}^{1},
\]

and the Hermite polynomials on $\mathbb{R}^{d}$ are defined by%
\[
H_{\beta}\left(  x\right)  :=H_{\beta_{1}}\left(  x_{1}\right)  \ldots
H_{\beta_{d}}\left(  x_{d}\right)  ,\quad x\in\mathbb{R}^{d}.
\]

The generating function is ?? DEL?%
\[
g\left(  x,\tau\right)  =e^{-\left\vert \tau\right\vert ^{2}+2\tau x}%
=\sum\limits_{\beta\geq0}H_{\beta}\left(  x\right)  \frac{\tau^{\beta}}%
{\beta!},\quad x,\tau\in\mathbb{R}^{d}.
\]

We have%
\[
D_{t}^{k}e^{-t^{2}/2}=\left(  -1\right)  ^{k}2^{-k/2}H_{k}\left(  t/\sqrt
{2}\right)  e^{-t^{2}/2},
\]

and so on $\mathbb{R}^{d}$,%
\begin{align*}
\frac{1}{n!}\left\vert D\right\vert ^{2n}e^{-\left\vert x\right\vert ^{2}/2}
& =\sum\limits_{\left\vert \beta\right\vert =n}\frac{1}{\beta!}D^{2\beta
}e^{-\left\vert x\right\vert ^{2}/2}\\
& =\sum\limits_{\left\vert \beta\right\vert =n}\frac{1}{\beta!}\left(
D_{1}^{2\beta_{1}}e^{-x_{1}^{2}/2}\right)  \times\ldots\times\left(
D_{d}^{2\beta_{d}}e^{-x_{d}^{2}/2}\right) \\
& =\sum\limits_{\left\vert \beta\right\vert =n}\frac{1}{2^{\beta_{1}}\beta
!}H_{2\beta_{1}}\left(  \frac{x_{1}}{\sqrt{2}}\right)  e^{-x_{1}^{2}/2}%
\times\ldots\times2^{-\beta_{d}}H_{2\beta_{d}}\left(  \frac{x_{d}}{\sqrt{2}%
}\right)  e^{-x_{d}^{2}/2}\\
& =\sum\limits_{\left\vert \beta\right\vert =n}\frac{1}{2^{\beta_{1}}\beta
!}H_{2\beta}\left(  \frac{x}{\sqrt{2}}\right)  e^{-\left\vert x\right\vert
^{2}/2}\\
& =\frac{1}{2^{n}}\left(  \sum\limits_{\left\vert \beta\right\vert =n}\frac
{1}{\beta!}H_{2\beta}\right)  \left(  \frac{x}{\sqrt{2}}\right)
e^{-\left\vert x\right\vert ^{2}/2}.
\end{align*}

From Wikipedia [??] we have%
\[
H_{k}\left(  t\right)  =\left(  e^{-\frac{1}{2}D^{2}}s^{k}\right)  \left(
2t\right)  .
\]

Check:%
\[
H_{2k}\left(  t\right)  =\left(  e^{-\frac{1}{2}D^{2}}s^{2k}\right)  \left(
2t\right)  =\left(  s^{2k}-\frac{1}{2}D^{2}s^{2k}+\frac{1}{2!2^{2}}D^{4}%
s^{2k}-\ldots\right)  \left(  2t\right)  .
\]

Hence%
\begin{align*}
\sum\limits_{\left\vert \beta\right\vert =n}\frac{1}{\beta!}H_{2\beta}\left(
x\right)   & =\frac{1}{n!}\left(  e^{-\frac{1}{2}\left\vert D\right\vert ^{2}%
}\left\vert \cdot\right\vert ^{2n}\right)  \left(  2x\right) \\
& =\frac{1}{n!}\sum\limits_{n=0}^{\infty}\frac{1}{k!}\left(  -\frac{1}%
{2}\left\vert D\right\vert ^{2}\right)  ^{k}\left\vert \cdot\right\vert
^{2n}\\
& =\frac{1}{n!}\sum\limits_{k=0}^{\infty}\frac{\left(  -1\right)  ^{k}}%
{2^{k}k!}\left(  \left\vert D\right\vert ^{2k}\left\vert \cdot\right\vert
^{2n}\right)  \left(  2x\right)  ,
\end{align*}

but from part 8 of Corollary \ref{Lem_iterLaplacian_rad},%
\[
\left\vert D\right\vert ^{2k}\left(  \left\vert x\right\vert ^{2n}\right)
=\left\{
\begin{array}
[c]{ll}%
2^{2k}\left(  k!\right)  ^{2}\tbinom{n}{k}\tbinom{n+d/2-1}{k}\left\vert
x\right\vert ^{2n-2k} & k\leq n,\\
0, & k>n,
\end{array}
\right.
\]

so that when $k\leq n$,
\begin{align}
\sum\limits_{\left\vert \beta\right\vert =n}\frac{1}{\beta!}H_{2\beta}\left(
x\right)   & =\frac{1}{n!}\sum\limits_{k=0}^{n}\frac{\left(  -1\right)  ^{k}%
}{2^{k}k!}\left(  \left\vert D\right\vert ^{2k}\left\vert \cdot\right\vert
^{2n}\right)  \left(  2x\right) \nonumber\\
& =\frac{1}{n!}\sum\limits_{k=0}^{n}\frac{\left(  -1\right)  ^{k}}{k!2^{k}%
}2^{2k}\left(  k!\right)  ^{2}\tbinom{n}{k}\tbinom{n+d/2-1}{k}\left\vert
2x\right\vert ^{2n-2k}\nonumber\\
& =\frac{1}{n!}\sum\limits_{k=0}^{n}\left(  -1\right)  ^{k}2^{k}k!\tbinom
{n}{k}\tbinom{n+d/2-1}{k}\left\vert 2x\right\vert ^{2n-2k}\nonumber\\
& =\sum\limits_{k=0}^{n}\frac{\left(  -1\right)  ^{k}2^{k}}{\left(
n-k\right)  !}\tbinom{n+d/2-1}{k}\left\vert 2x\right\vert ^{2n-2k}\nonumber\\
& =\sum\limits_{k=0}^{n}\frac{\left(  -1\right)  ^{n-k}2^{n-k}}{k!}%
\tbinom{n-k+d/2-1}{n-k}\left\vert 2x\right\vert ^{2k}\nonumber\\
& =\left(  -1\right)  ^{n}2^{n}\sum\limits_{k=0}^{n}\frac{\left(  -1\right)
^{k}}{2^{k}k!}\tbinom{n-k+d/2-1}{n-k}\left\vert 2x\right\vert ^{2k}%
,\label{a014}%
\end{align}

and we can conclude that ?? \textbf{CHECK}! ??%
\begin{align}
\frac{1}{n!}\left\vert D\right\vert ^{2n}e^{-\left\vert \cdot\right\vert
^{2}/2}  & =2^{-n}\left(  \sum\limits_{\left\vert \beta\right\vert =n}\frac
{1}{\beta!}H_{2\beta}\right)  \left(  \frac{??x}{\sqrt{2}}\right)
e^{-\left\vert x\right\vert ^{2}/2}\nonumber\\
& =2^{-n}\left(  \sum\limits_{k=0}^{n}\frac{\left(  -1\right)  ^{k}2^{k}%
}{\left(  n-k\right)  !}\tbinom{n+d/2-1}{k}\left\vert \frac{x}{\sqrt{2}%
}\right\vert ^{2n-2k}\right)  e^{-\left\vert x\right\vert ^{2}/2}\nonumber\\
& =2^{-n}\left(  \sum\limits_{k=0}^{n}\frac{\left(  -1\right)  ^{k}2^{k}%
}{\left(  n-k\right)  !}\tbinom{n+d/2-1}{k}\frac{1}{2^{n-k}}\left\vert
x\right\vert ^{2n-2k}\right)  e^{-\left\vert x\right\vert ^{2}/2}\nonumber\\
& =\left(  \sum\limits_{k=0}^{n}\frac{\left(  -1\right)  ^{k}}{4^{n-k}}%
\tbinom{n+d/2-1}{k}\left\vert x\right\vert ^{2n-2k}\right)  e^{-\left\vert
x\right\vert ^{2}/2}\nonumber\\
& =\left(  \sum\limits_{k=0}^{n}\frac{\left(  -1\right)  ^{n-k}}{4^{k}}%
\tbinom{n+d/2-1}{n-k}\left\vert x\right\vert ^{2k}\right)  e^{-\left\vert
x\right\vert ^{2}/2}\nonumber\\
& =\left(  -1\right)  ^{n}\left(  \sum\limits_{k=0}^{n}\frac{\left(
-1\right)  ^{k}}{4^{k}}\tbinom{n+d/2-1}{n-k}\left\vert x\right\vert
^{2k}\right)  e^{-\left\vert x\right\vert ^{2}/2}.\label{a013}%
\end{align}

?? Check case $d=1$: ??

and from \ref{a008},%
\begin{align}
d_{j}\left(  \left\vert iD\right\vert \right)  e^{-\left\vert \cdot\right\vert
^{2}/2}  & =\sum\limits_{n=1}^{\left\lfloor j/2\right\rfloor }\frac{\left(
-1\right)  ^{n}}{2^{n}}\frac{1}{n!}\left\vert D\right\vert ^{2n}e^{-\left\vert
\cdot\right\vert ^{2}/2}\nonumber\\
& =\sum\limits_{n=1}^{\left\lfloor j/2\right\rfloor }\frac{\left(  -1\right)
^{n}}{2^{n}}\left(  -1\right)  ^{n}\left(  \sum\limits_{k=0}^{n}\frac{\left(
-1\right)  ^{k}}{4^{k}}\tbinom{n+d/2-1}{n-k}\left\vert x\right\vert
^{2k}\right)  e^{-\left\vert x\right\vert ^{2}/2}\nonumber\\
& =\left(  \sum\limits_{n=1}^{\left\lfloor j/2\right\rfloor }\sum
\limits_{k=0}^{n}\frac{1}{2^{n}}\frac{\left(  -1\right)  ^{k}}{4^{k}}%
\tbinom{n+d/2-1}{n-k}\left\vert x\right\vert ^{2k}\right)  e^{-\left\vert
x\right\vert ^{2}/2}\nonumber\\
& =\left(  \sum\limits_{k=0}^{\left\lfloor j/2\right\rfloor }\sum
\limits_{n=k}^{\left\lfloor j/2\right\rfloor }\frac{\left(  -1\right)  ^{k}%
}{4^{k}}\frac{1}{2^{n}}\tbinom{n+d/2-1}{n-k}\left\vert x\right\vert
^{2k}\right)  e^{-\left\vert x\right\vert ^{2}/2}\nonumber\\
& =\left(  \sum\limits_{k=0}^{\left\lfloor j/2\right\rfloor }\frac{\left(
-1\right)  ^{k}}{4^{k}}\left(  \sum\limits_{n=k}^{\left\lfloor
j/2\right\rfloor }\frac{1}{2^{n}}\tbinom{n+d/2-1}{n-k}\right)  \left\vert
x\right\vert ^{2k}\right)  e^{-\left\vert x\right\vert ^{2}/2}\nonumber\\
& =\left(  \sum\limits_{k=0}^{\left\lfloor j/2\right\rfloor }\frac{\left(
-1\right)  ^{k}}{4^{k}}\left(  \sum\limits_{m=0}^{\left\lfloor
j/2\right\rfloor -k}\frac{1}{2^{m+k}}\tbinom{m+k+d/2-1}{m}\right)  \left\vert
x\right\vert ^{2k}\right)  e^{-\left\vert x\right\vert ^{2}/2}\nonumber\\
& =\left(  \sum\limits_{k=0}^{\left\lfloor j/2\right\rfloor }\frac{\left(
-1\right)  ^{k}}{8^{k}}\left(  \sum\limits_{m=0}^{\left\lfloor
j/2\right\rfloor -k}\frac{1}{2^{m}}\tbinom{m+k+d/2-1}{m}\right)  \left\vert
x\right\vert ^{2k}\right)  e^{-\left\vert x\right\vert ^{2}/2}.\label{a1.27}%
\end{align}

and so \ref{a009} can be written%
\[
\widehat{\phi_{0}}\left(  x\right)  =q\left(  \left\vert x\right\vert
^{2}\right)  e^{-\left\vert x\right\vert ^{2}/2},
\]

where $q$ is a polynomial, and with reference to \ref{a2.34},%
\[
\phi_{\mathcal{F}}\left(  t\right)  =q\left(  t^{2}\right)  e^{-t^{2}/2}.
\]

\underline{Application to \ref{a2.37}}%
\[
\phi_{\mathcal{F}}^{\prime}\left(  t\right)  =2tq^{\prime}\left(
t^{2}\right)  e^{-t^{2}/2}-tq\left(  t^{2}\right)  e^{-t^{2}/2}=t\left(
2q^{\prime}\left(  t^{2}\right)  -q\left(  t^{2}\right)  \right)  e^{-t^{2}%
/2}.
\]

and \ref{a2.37} i.e.%
\[
\left(  \int_{0}^{\infty}t^{d}\left\vert \phi_{\mathcal{F}}^{\prime}\left(
t\right)  \right\vert dt\right)  ^{\frac{1}{2}}\left(  \sum\limits_{l=0}%
^{\theta-n}\frac{1}{l!}\int_{0}^{\infty}t^{2l+d}\left\vert \phi_{\mathcal{F}%
}^{\prime}\left(  t\right)  \right\vert dt\right)  ^{\frac{1}{2}},
\]

can be estimated numerically.\medskip

\underline{Estimate \ref{a2.38}} i.e. $\max\limits_{s\in\left[  0,r_{3}%
\right]  }\left\vert D_{s}^{\theta}\left(  \phi_{0}\right)  _{\circ}\left(
s\right)  \right\vert $ - see \ref{a011}, \ref{a010}, \ref{a012},
\ref{a007}.\medskip

\underline{We can set $\rho=\phi_{0}$}

\begin{remark}
From Hermite polynomial article in Wikipedia: ?? \textbf{CHECK}! ??%
\[
\left(  H_{k}\left(  t\right)  e^{-t^{2}/2}\right)  ^{\wedge}\left(  s\right)
=\left(  -i\right)  ^{k}H_{k}\left(  s\right)  e^{-s^{2}/2}.
\]

\end{remark}

??

\begin{remark}
\textbf{Part 1} From Exercise 13.1.5 of \cite{Arfken70},%
\[
H_{2n}\left(  t\right)  =\left(  -1\right)  ^{n}\sum\limits_{k=0}^{n}\left(
-1\right)  ^{k}\frac{\left(  2n\right)  !}{\left(  2k\right)  !\left(
n-k\right)  !}\left(  2t\right)  ^{2k},
\]

so that%
\begin{align*}
&  H_{2\beta}\left(  x\right) \\
&  =\left(  -1\right)  ^{\beta_{1}}\sum\limits_{\alpha_{1}\leq\beta_{1}}%
\frac{\left(  -1\right)  ^{\alpha_{1}}\left(  2\beta_{1}\right)  !}{\left(
2\alpha_{1}\right)  !\left(  \beta_{1}-\alpha_{1}\right)  !}\left(
2x_{1}\right)  ^{2\alpha_{1}}\times\ldots\times\sum\limits_{\alpha_{d}%
\leq\beta_{d}}\frac{\left(  -1\right)  ^{\alpha_{d}}\left(  2\beta_{d}\right)
!}{\left(  2\alpha_{d}\right)  !\left(  \beta_{d}-\alpha_{d}\right)  !}\left(
2x_{d}\right)  ^{2\alpha_{d}}\\
&  =\left(  -1\right)  ^{\left\vert \beta\right\vert }\sum\limits_{\alpha
\leq\beta}\frac{\left(  -1\right)  ^{\left\vert \alpha\right\vert }\left(
2\beta\right)  !}{\left(  2\alpha\right)  !\left(  \beta-\alpha\right)
!}\left(  2x\right)  ^{2\alpha},
\end{align*}

and hence%
\begin{align*}
\sum\limits_{\left\vert \beta\right\vert =n}\frac{1}{\beta!}H_{2\beta}\left(
x\right)   & =\sum\limits_{\left\vert \beta\right\vert =n}\frac{\left(
-1\right)  ^{\left\vert \beta\right\vert }}{\beta!}\sum\limits_{\alpha
\leq\beta}\frac{\left(  -1\right)  ^{\left\vert \alpha\right\vert }\left(
2\beta\right)  !}{\left(  2\alpha\right)  !\left(  \beta-\alpha\right)
!}\left(  2x\right)  ^{2\alpha}\\
& =\sum\limits_{\left\vert \beta\right\vert =n}\sum\limits_{\alpha\leq\beta
}\frac{\left(  -1\right)  ^{\left\vert \alpha\right\vert +\left\vert
\beta\right\vert }\alpha!\left(  2\beta\right)  !}{\beta!\left(
2\alpha\right)  !\left(  \beta-\alpha\right)  !}\frac{\left(  2x\right)
^{2\alpha}}{\alpha!}\\
& =\left(  -1\right)  ^{n}\sum\limits_{\left\vert \beta\right\vert =n}%
\sum\limits_{\alpha\leq\beta}\left(  -1\right)  ^{\left\vert \alpha\right\vert
}\frac{\alpha!\left(  2\beta\right)  !}{\left(  2\alpha\right)  !\beta!\left(
\beta-\alpha\right)  !}\frac{\left(  2x\right)  ^{2\alpha}}{\alpha!}.
\end{align*}

Set $b_{\alpha}^{\left(  \beta\right)  }=\left(  -1\right)  ^{\left\vert
\alpha\right\vert }\frac{\alpha!\left(  2\beta\right)  !}{\left(
2\alpha\right)  !\beta!\left(  \beta-\alpha\right)  !}$ so that from Theorem
\ref{Thm_double_sum},
\begin{align*}
\sum\limits_{\left\vert \beta\right\vert =n}\frac{1}{\beta!}H_{2\beta}\left(
x\right)   & =\left(  -1\right)  ^{n}\sum\limits_{\left\vert \beta\right\vert
=n}\sum\limits_{\alpha\leq\beta}b_{\alpha}^{\left(  \beta\right)  }%
\frac{\left(  2x\right)  ^{2\alpha}}{\alpha!}\\
& =\left(  -1\right)  ^{n}\sum\limits_{\left\vert \gamma\right\vert \leq
n}\sum\limits_{\left\vert \sigma\right\vert =n-\left\vert \gamma\right\vert
}b_{\gamma}^{\left(  \gamma+\sigma\right)  }\frac{\left(  2x\right)
^{2\gamma}}{\gamma!}\\
& =\left(  -1\right)  ^{n}\sum\limits_{\left\vert \gamma\right\vert \leq
n}\sum\limits_{\left\vert \sigma\right\vert =n-\left\vert \gamma\right\vert
}\left(  -1\right)  ^{\left\vert \gamma\right\vert }\frac{\left(  2\left(
\gamma+\sigma\right)  \right)  !}{\left(  2\gamma\right)  !\left(
\gamma+\sigma\right)  !\left(  \left(  \gamma+\sigma\right)  -\gamma\right)
!}\left(  2x\right)  ^{2\gamma}\\
& =\left(  -1\right)  ^{n}\sum\limits_{\left\vert \gamma\right\vert \leq
n}\sum\limits_{\left\vert \sigma\right\vert =n-\left\vert \gamma\right\vert
}\left(  -1\right)  ^{\left\vert \gamma\right\vert }\frac{\left(  2\left(
\gamma+\sigma\right)  \right)  !}{\left(  \gamma+\sigma\right)  !\sigma!}%
\frac{\left(  2x\right)  ^{2\gamma}}{\left(  2\gamma\right)  !}\\
& =\left(  -1\right)  ^{n}\sum\limits_{\left\vert \gamma\right\vert \leq
n}\left(  -1\right)  ^{\left\vert \gamma\right\vert }\left(  \sum
\limits_{\left\vert \sigma\right\vert =n-\left\vert \gamma\right\vert }%
\frac{\left(  2\gamma+2\sigma\right)  !}{\left(  \gamma+\sigma\right)
!\sigma!}\right)  \frac{\left(  2x\right)  ^{2\gamma}}{\left(  2\gamma\right)
!}\\
& =\left(  -1\right)  ^{n}\sum\limits_{k=0}^{n}\left(  -1\right)  ^{k}%
\sum\limits_{\left\vert \gamma\right\vert =k}\left(  \sum\limits_{\left\vert
\sigma\right\vert =n-k}\frac{\left(  2\gamma+2\sigma\right)  !}{\left(
\gamma+\sigma\right)  !\sigma!}\right)  \frac{\left(  2x\right)  ^{2\gamma}%
}{\left(  2\gamma\right)  !}\\
& =\left(  -1\right)  ^{n}\sum\limits_{k=0}^{n}\left(  -1\right)  ^{k}%
4^{k}\sum\limits_{\left\vert \gamma\right\vert =k}\left(  \sum
\limits_{\left\vert \sigma\right\vert =n-k}\frac{\left(  2\gamma
+2\sigma\right)  !}{\left(  \gamma+\sigma\right)  !\sigma!}\right)
\frac{x^{2\gamma}}{\left(  2\gamma\right)  !}\\
& =\left(  -1\right)  ^{n}\sum\limits_{k=0}^{n}\sum\limits_{\left\vert
\gamma\right\vert =k}\left(  -1\right)  ^{k}4^{k}\left(  \sum
\limits_{\left\vert \sigma\right\vert =n-k}\frac{\left(  2\gamma
+2\sigma\right)  !}{\left(  \gamma+\sigma\right)  !\sigma!}\right)
\frac{x^{2\gamma}}{\left(  2\gamma\right)  !}.
\end{align*}

Now compare with \ref{a014}:%
\begin{align}
\sum\limits_{\left\vert \beta\right\vert =n}\frac{1}{\beta!}H_{2\beta}\left(
x\right)   & =\left(  -1\right)  ^{n}2^{n}\sum\limits_{k=0}^{n}\frac{\left(
-1\right)  ^{k}}{2^{k}}\tbinom{n-k+d/2-1}{n-k}\frac{1}{k!}\left\vert
2x\right\vert ^{2k}\nonumber\\
& =\left(  -1\right)  ^{n}2^{n}\sum\limits_{k=0}^{n}\left(  -1\right)
^{k}2^{k}\tbinom{n-k+d/2-1}{n-k}\sum\limits_{\left\vert \gamma\right\vert
=k}\frac{x^{2\gamma}}{\gamma!}\nonumber\\
& =\left(  -1\right)  ^{n}\sum\limits_{k=0}^{n}\sum\limits_{\left\vert
\gamma\right\vert =k}\left(  -1\right)  ^{k}2^{n+k}\tbinom{n-k+d/2-1}%
{n-k}\frac{x^{2\gamma}}{\gamma!},\label{a015}%
\end{align}

so that%
\begin{align*}
\left(  -1\right)  ^{k}4^{k}\left(  \sum\limits_{\left\vert \sigma\right\vert
=n-k}\frac{\left(  2\gamma+2\sigma\right)  !}{\left(  \gamma+\sigma\right)
!\sigma!}\right)  \frac{1}{\left(  2\gamma\right)  !}  & =\left(  -1\right)
^{k}2^{n+k}\tbinom{n-k+d/2-1}{n-k}\frac{1}{\gamma!},\quad\left\vert
\gamma\right\vert =k,\text{ }n\geq k\geq0,\\
\sum\limits_{\left\vert \sigma\right\vert =n-k}\frac{\left(  2\gamma
+2\sigma\right)  !}{\left(  \gamma+\sigma\right)  !\sigma!}  & =2^{n-k}%
\tbinom{n-k+d/2-1}{n-k}\frac{\left(  2\gamma\right)  !}{\gamma!}%
,\quad\left\vert \gamma\right\vert =k,\text{ }n\geq k\geq0,\\
\sum\limits_{\left\vert \sigma\right\vert =m}\frac{\left(  2\gamma
+2\sigma\right)  !}{\left(  \gamma+\sigma\right)  !\sigma!}  & =2^{m}%
\tbinom{m+d/2-1}{m}\frac{\left(  2\gamma\right)  !}{\gamma!},\quad\gamma
\geq0,\text{ }m\geq0.
\end{align*}
\medskip

\textbf{Part 2} From Wikipedia the Hermite polynomials can be expressed in
terms of Laguerre polynomials as%
\[
\frac{1}{n!}H_{2n}\left(  t\right)  =\left(  -4\right)  ^{n}L_{n}^{\left(
-1/2\right)  }\left(  t^{2}\right)  =4^{n}\sum\limits_{k=0}^{n}\left(
-1\right)  ^{n-k}\tbinom{n-\frac{1}{2}}{n-k}\frac{t^{2k}}{k!}.
\]

so that%
\begin{align*}
\frac{1}{\beta!}H_{2\beta}\left(  x\right)   & =4^{\beta_{1}}\sum
\limits_{\alpha_{1}\leq\beta_{1}}\left(  -1\right)  ^{n-\alpha_{1}}%
\tbinom{\beta_{1}-\frac{1}{2}}{\beta_{1}-\alpha_{1}}\frac{x_{1}^{2\alpha_{1}}%
}{\alpha_{1}!}\times\ldots\times4^{\beta_{d}}\sum\limits_{\alpha_{d}\leq
\beta_{d}}\left(  -1\right)  ^{n-\alpha_{d}}\tbinom{\beta_{d}-\frac{1}{2}%
}{\beta_{d}-\alpha_{1}}\frac{x_{d}^{2\alpha_{d}}}{\alpha_{d}!}\\
& =4^{\left\vert \beta\right\vert }\sum\limits_{\alpha\leq\beta}\left(
-1\right)  ^{nd-\left\vert \alpha\right\vert }\tbinom{\beta-\frac{1}{2}}%
{\beta-\alpha}\frac{x^{2\alpha}}{\alpha!},
\end{align*}

and hence%
\begin{align*}
\sum\limits_{\left\vert \beta\right\vert =n}\frac{1}{\beta!}H_{2\beta}\left(
x\right)   & =\sum\limits_{\left\vert \beta\right\vert =n}\frac{1}{\beta
!}4^{\left\vert \beta\right\vert }\sum\limits_{\alpha\leq\beta}\left(
-1\right)  ^{nd-\left\vert \alpha\right\vert }\tbinom{\beta-\frac{1}{2}}%
{\beta-\alpha}\frac{x^{2\alpha}}{\alpha!}\\
& =\left(  -1\right)  ^{nd}4^{n}\sum\limits_{\left\vert \beta\right\vert
=n}\sum\limits_{\alpha\leq\beta}\frac{\left(  -1\right)  ^{\left\vert
\alpha\right\vert }}{\beta!}\tbinom{\beta-\frac{1}{2}}{\beta-\alpha}%
\frac{x^{2\alpha}}{\alpha!},
\end{align*}

If $b_{\alpha}^{\left(  \beta\right)  }=\frac{\left(  -1\right)  ^{\left\vert
\alpha\right\vert }}{\beta!}\tbinom{\beta-\frac{1}{2}}{\beta-\alpha}$ then
from Theorem \ref{Thm_double_sum},%
\begin{align*}
\sum\limits_{\left\vert \beta\right\vert =n}\frac{1}{\beta!}H_{2\beta}\left(
x\right)   & =\left(  -1\right)  ^{n}\sum\limits_{\left\vert \beta\right\vert
=n}\sum\limits_{\alpha\leq\beta}b_{\alpha}^{\left(  \beta\right)  }%
\frac{x^{2\alpha}}{\alpha!}\\
& =\left(  -1\right)  ^{n}\sum\limits_{\left\vert \gamma\right\vert \leq
n}\sum\limits_{\left\vert \sigma\right\vert =n-\left\vert \gamma\right\vert
}\frac{\left(  -1\right)  ^{\left\vert \gamma\right\vert }}{\left(
\gamma+\sigma\right)  !}\tbinom{\gamma+\sigma-\frac{1}{2}}{\gamma
+\sigma-\gamma}\frac{x^{2\gamma}}{\gamma!}\\
& =\left(  -1\right)  ^{n}\sum\limits_{\left\vert \gamma\right\vert \leq
n}\left(  -1\right)  ^{\left\vert \gamma\right\vert }\sum\limits_{\left\vert
\sigma\right\vert =n-\left\vert \gamma\right\vert }\frac{1}{\left(
\gamma+\sigma\right)  !}\tbinom{\gamma+\sigma-\frac{1}{2}}{\sigma}%
\frac{x^{2\gamma}}{\gamma!}\\
& =\left(  -1\right)  ^{n}\sum\limits_{\left\vert \gamma\right\vert \leq
n}\left(  -1\right)  ^{\left\vert \gamma\right\vert }\left(  \sum
\limits_{\left\vert \sigma\right\vert =n-\left\vert \gamma\right\vert }%
\frac{1}{\left(  \gamma+\sigma\right)  !}\tbinom{\gamma+\sigma-\frac{1}{2}%
}{\sigma}\right)  \frac{x^{2\gamma}}{\gamma!}\\
& =\left(  -1\right)  ^{n}\sum\limits_{k=0}^{n}\sum\limits_{\left\vert
\gamma\right\vert =k}\left(  -1\right)  ^{k}\left(  \sum\limits_{\left\vert
\sigma\right\vert =n-k}\frac{1}{\left(  \gamma+\sigma\right)  !}\tbinom
{\gamma+\sigma-\frac{1}{2}}{\sigma}\right)  \frac{x^{2\gamma}}{\gamma!}.
\end{align*}

Now compare with \ref{a014} using the form \ref{a015}: ?? \textbf{NEED }%
$2x$\textbf{\ instead of} $x$? ??%
\begin{align*}
\left(  -1\right)  ^{n}\sum\limits_{k=0}^{n}\sum\limits_{\left\vert
\gamma\right\vert =k}\left(  -1\right)  ^{k}\left(  \sum\limits_{\left\vert
\sigma\right\vert =n-k}\frac{\tbinom{\gamma+\sigma-\frac{1}{2}}{\sigma}%
}{\left(  \gamma+\sigma\right)  !}\right)  \frac{\left(  ??2x\right)
^{2\gamma}}{\gamma!}  & =\left(  -1\right)  ^{n}\sum\limits_{k=0}^{n}%
\sum\limits_{\left\vert \gamma\right\vert =k}\left(  -1\right)  ^{k}%
2^{n+k}\tbinom{n-k+\frac{d}{2}-1}{n-k}\frac{x^{2\gamma}}{\gamma!},\\
\sum\limits_{\left\vert \sigma\right\vert =n-k}\frac{1}{\left(  \gamma
+\sigma\right)  !}\tbinom{\gamma+\sigma-\frac{1}{2}}{\sigma}  & =2^{n-k}%
\tbinom{n-k+\frac{d}{2}-1}{n-k},\quad\left\vert \gamma\right\vert =k,\text{
}n\geq k\geq0,\\
\sum\limits_{\left\vert \sigma\right\vert =m}\frac{1}{\left(  \gamma
+\sigma\right)  !}\tbinom{\gamma+\sigma-\frac{1}{2}}{\sigma}  & =2^{m}%
\tbinom{m+\frac{d}{2}-1}{m},\quad\gamma\geq0,\text{ }m\geq0.
\end{align*}

\end{remark}

\section{Taylor series expansions for functions in $X_{w}^{\theta}$ when $w\in
W3.1$\label{Sect_data_fn_Taylor_W3.1}}

\subsection{Introduction}

Approach II below enables us to obtain significantly improved convergence
estimates for the extended B-splines compared to Approach I. In fact ??

\subsection{Approach 1\label{SbSect_TaylorDataApproach1}}

We now consider the case when the weight function has property W3.1
(\ref{p36}) for order $\theta$ and (uniform) parameter $\kappa$. Define
\[
\underline{\kappa}=\min\limits_{k}\kappa_{k}.
\]

The remainder term for the Taylor series for arbitrary $f\in X_{w}^{\theta}$
is given by \ref{a93} and the polynomial term was estimated by \ref{a94}. This
means we must estimate
\begin{equation}
\mathcal{R}_{n+1}f_{\rho}\left(  \cdot,a\right)  =\frac{\sqrt{2\pi}}%
{n!}\left(  \left(  ia\xi\right)  ^{n+1}\overline{\widehat{g_{n}}}\left(
a\xi\right)  \widehat{f_{\rho}}\right)  ^{\vee},\label{p76}%
\end{equation}

where $f_{\rho}\in X_{w}^{\theta}$ is given by \ref{a920}.

As with the property W3.2 we will consider two cases: $\left\lfloor
\underline{\kappa}\right\rfloor \geq\theta$ and $\left\lfloor \underline
{\kappa}\right\rfloor <\theta$. The former case is relatively simple and we
consider that first.\medskip

\fbox{\textbf{Case}: $\left\lfloor \underline{\kappa}\right\rfloor \geq\theta
$} Of course $n\leq\left\lfloor \underline{\kappa}\right\rfloor $. Since $f\in
X_{w}^{\theta}$ we have $\xi^{\beta}\widehat{f}=\xi^{\beta}f_{F}$ when
$\left\vert \beta\right\vert =\theta$ so
\[
\left(  a\xi\right)  ^{\theta}\widehat{f}=\sum_{\left\vert \beta\right\vert
=\theta}\frac{a^{\beta}\xi^{\beta}}{\beta!}\widehat{f}=\sum_{\left\vert
\beta\right\vert =\theta}\frac{a^{\beta}\xi^{\beta}}{\beta!}f_{F}=\left(
a\xi\right)  ^{\theta}f_{F},
\]

and hence%
\begin{align*}
\left(  \left(  ia\xi\right)  ^{n+1}\overline{\widehat{g_{n}}}\left(
a\xi\right)  \widehat{f_{\rho}}\right)  ^{\vee}  & =\left(  \left(
ia\xi\right)  ^{n-\theta+1}\overline{\widehat{g_{n}}}\left(  a\xi\right)
\left(  ia\xi\right)  ^{\theta}\widehat{f_{\rho}}\right)  ^{\vee}\\
& =\left(  \left(  ia\xi\right)  ^{n-\theta+1}\overline{\widehat{g_{n}}%
}\left(  a\xi\right)  \left(  ia\xi\right)  ^{\theta}f_{F}\right)  ^{\vee}\\
& =\left(  \left(  ia\xi\right)  ^{n+1}\overline{\widehat{g_{n}}}\left(
a\xi\right)  f_{F}\right)  ^{\vee}\\
& =\left(  \frac{\left(  ia\xi\right)  ^{n+1}\overline{\widehat{g_{n}}}\left(
a\xi\right)  }{\sqrt{w}\left\vert \cdot\right\vert ^{\theta}}\sqrt
{w}\left\vert \cdot\right\vert ^{\theta}f_{F}\right)  ^{\vee}\\
& =\left(  \frac{\left(  ia\xi\right)  ^{n+1}\overline{\widehat{g_{n}}}\left(
a\xi\right)  }{\sqrt{w}\left\vert \cdot\right\vert ^{\theta}}\right)  ^{\vee
}\ast\left(  \sqrt{w}\left\vert \cdot\right\vert ^{\theta}f_{F}\right)
^{\vee}.
\end{align*}

Applying Young's convolution estimate \ref{1.056} and then Parceval's theorem
we get%
\begin{align}
\left\vert \left(  \left(  ia\xi\right)  ^{n+1}\overline{\widehat{g_{n}}%
}\left(  a\xi\right)  \widehat{f_{\rho}}\right)  ^{\vee}\left(  x\right)
\right\vert  & \leq\left\Vert \left(  \frac{\left(  ia\xi\right)
^{n+1}\overline{\widehat{g_{n}}}\left(  a\xi\right)  }{\sqrt{w}\left\vert
\cdot\right\vert ^{\theta}}\right)  ^{\vee}\right\Vert _{2}\left\Vert \left(
\sqrt{w}\left\vert \cdot\right\vert ^{\theta}f_{F}\right)  ^{\vee}\right\Vert
_{2}\label{a2.35}\\
& =\left\Vert \frac{\left(  ia\xi\right)  ^{n+1}\overline{\widehat{g_{n}}%
}\left(  a\xi\right)  }{\sqrt{w}\left\vert \cdot\right\vert ^{\theta}%
}\right\Vert _{2}\left\Vert \sqrt{w}\left\vert \cdot\right\vert ^{\theta}%
f_{F}\right\Vert _{2}\nonumber\\
& \leq\left\Vert \frac{\left(  ia\xi\right)  ^{n+1}\overline{\widehat{g_{n}}%
}\left(  a\xi\right)  }{\sqrt{w}\left\vert \cdot\right\vert ^{\theta}%
}\right\Vert _{2}\left\vert f_{\rho}\right\vert _{w,\theta}\nonumber\\
& =\left(  \int\frac{\left(  a\xi\right)  ^{2\left(  n+1\right)  }\left\vert
\overline{\widehat{g_{n}}}\left(  a\xi\right)  \right\vert ^{2}}{w\left\vert
\cdot\right\vert ^{2\theta}}\right)  ^{\frac{1}{2}}\left\vert f_{\rho
}\right\vert _{w,\theta},\label{a87}%
\end{align}

so that%
\begin{equation}
\left\vert \mathcal{R}_{n+1}f_{\rho}\left(  x,a\right)  \right\vert \leq
\frac{\sqrt{2\pi}}{n!}\left(  \int\frac{\left(  a\xi\right)  ^{2\left(
n+1\right)  }\left\vert \overline{\widehat{g_{n}}}\left(  a\xi\right)
\right\vert ^{2}}{w\left\vert \cdot\right\vert ^{2\theta}}\right)  ^{\frac
{1}{2}}\left\vert f_{\rho}\right\vert _{w,\theta},\quad\left\{
\begin{array}
[c]{c}%
\theta\leq\left\lfloor \underline{\kappa}\right\rfloor ,\\
n\leq\left\lfloor \underline{\kappa}\right\rfloor .
\end{array}
\right. \label{a1011}%
\end{equation}
\medskip

\fbox{\textbf{Case}: $\left\lfloor \underline{\kappa}\right\rfloor <\theta$}
Since $n\leq\left\lfloor \underline{\kappa}\right\rfloor $ we have $n<\theta$,
and as above we employ the partition of unity \ref{a15} to obtain the sum
\ref{a35}.\medskip

\underline{\textbf{The first term} $\frac{\sqrt{2\pi}}{n!}\left(
\overline{\widehat{g_{n}}}\left(  a\xi\right)  \left(  ia\xi\right)  \phi
_{0}\widehat{\left(  aD\right)  ^{n}f_{\rho}}\right)  ^{\vee}$ \textbf{in
}\ref{a35}}\medskip

In an identical manner to what follows \ref{a35} (clunky reference) we obtain
the estimate \ref{a80}. The next step is to estimate $\left(  \widehat
{a}D\right)  ^{n}f_{\rho}\left(  x\right)  $ by applying the upper bounds for
$\mathcal{Q}_{\emptyset,\theta-n,\xi}\left(  e^{ix\xi}\right)  $ given in
Theorem \ref{Thm_bound_on_g(e,x)_2} to equation \ref{a49} for $\left(
\widehat{a}D\right)  ^{n}f_{\rho}$ i.e. to the formula%
\[
\left(  \widehat{a}D\right)  ^{n}f_{\rho}\left(  x\right)  =\left(
2\pi\right)  ^{-\frac{d}{2}}\int\left\vert \mathcal{Q}_{\emptyset,\theta
-n,\xi}\left(  e^{ix\xi}\right)  \right\vert \left\vert \widehat{a}%
\xi\right\vert ^{n}\left\vert f_{F}\left(  \xi\right)  \right\vert d\xi.
\]

But%
\begin{align*}
\int\left\vert \mathcal{Q}_{\emptyset,\theta-n,\xi}\left(  e^{ix\xi}\right)
\right\vert \left\vert \widehat{a}\xi\right\vert ^{n}\left\vert f_{F}\left(
\xi\right)  \right\vert d\xi & =\int\left\vert \mathcal{Q}_{\emptyset
,\theta-n,\xi}\left(  e^{ix\xi}\right)  \right\vert \frac{\left\vert
\widehat{a}\xi\right\vert ^{n}}{\sqrt{w}\left\vert \cdot\right\vert ^{\theta}%
}\sqrt{w}\left\vert \cdot\right\vert ^{\theta}\left\vert f_{F}\left(
\xi\right)  \right\vert d\xi\\
& \Rightarrow Cauchy-Schwartz\text{ }inequality\Rightarrow\\
& \leq\left(  \int\left\vert \mathcal{Q}_{\emptyset,\theta-n,\xi}\left(
e^{ix\xi}\right)  \right\vert ^{2}\frac{\left(  \widehat{a}\xi\right)  ^{2n}%
}{w\left\vert \cdot\right\vert ^{2\theta}}\right)  ^{1/2}\left\Vert \sqrt
{w}\left\vert \cdot\right\vert ^{\theta}f_{F}\right\Vert _{2}\\
& =\left(  \int\left\vert \mathcal{Q}_{\emptyset,\theta-n,\xi}\left(
e^{ix\xi}\right)  \right\vert ^{2}\frac{\left(  \widehat{a}\xi\right)  ^{2n}%
}{w\left\vert \cdot\right\vert ^{2\theta}}\right)  ^{1/2}\left\vert f_{\rho
}\right\vert _{w,\theta},
\end{align*}

and from \ref{a1.24} and \ref{a1.25} with $r_{3}=0$ we have%
\[
\int\left\vert \mathcal{Q}_{\emptyset,\theta-n,\xi}\left(  e^{ix\xi}\right)
\right\vert ^{2}\frac{\left\vert \widehat{a}\xi\right\vert ^{2n}}{w\left\vert
\cdot\right\vert ^{2\theta}}\leq\left(  \overline{C}_{\theta-n,0}^{\left(
\rho\right)  }\right)  ^{2}\left(  1+e_{\theta-n-1}\right)  \left(  \int%
\frac{\left\vert \cdot\right\vert ^{2n}}{w\left\vert \cdot\right\vert
^{2\theta}}\right)  \left(  1+\sum\limits_{k<\theta-n}\frac{\left\vert
x\right\vert ^{2k}}{k!}\right)  ,
\]

and%
\begin{align*}
&  \int\left\vert \mathcal{Q}_{\emptyset,\theta-n,\xi}\left(  e^{ix\xi
}\right)  \right\vert ^{2}\frac{\left\vert \widehat{a}\xi\right\vert ^{2n}%
}{w\left\vert \cdot\right\vert ^{2\theta}}\\
&  \leq\left(  \overline{C}_{\theta-n,0}^{\left(  \rho\right)  }\right)
^{2}\left(
\begin{array}
[c]{l}%
\int\left\vert \widehat{a}\widehat{\xi}\right\vert ^{2n}\frac{\left\vert
\cdot\right\vert ^{2n}}{w\left\vert \cdot\right\vert ^{2\theta}}+\\
+\sum\limits_{k<\theta-n}\frac{1}{k!}\int\left\vert \widehat{x}\widehat{\xi
}\right\vert ^{2k}\left\vert \widehat{a}\widehat{\xi}\right\vert ^{2n}%
\frac{\left\vert \cdot\right\vert ^{2n}}{w\left\vert \cdot\right\vert
^{2\theta}}%
\end{array}
\right)  \left(  1+\sum\limits_{k<\theta-n}\frac{\left\vert x\right\vert
^{2k}}{k!}\right)  ,
\end{align*}

which both exist by part 1 of Theorem \ref{Thm_property_wt_fn_W3}.

Hence, if $n\leq\left\lfloor \underline{\kappa}\right\rfloor <\theta$,%
\begin{align*}
&  \left\vert \left(  \widehat{a}D\right)  ^{n}f_{\rho}\left(  x\right)
\right\vert \\
&  \leq\frac{\overline{C}_{\theta-n,0}^{\left(  \rho\right)  }}{\left(
2\pi\right)  ^{\frac{d}{2}}}\left\vert f_{\rho}\right\vert _{w,\theta}\left(
1+\sum\limits_{k<\theta-n}\frac{\left\vert x\right\vert ^{2k}}{k!}\right)
^{\frac{1}{2}}\times\left\{
\begin{array}
[c]{l}%
\left(  1+e_{\theta-n-1}\right)  ^{\frac{1}{2}}\left(  \int\frac{\left\vert
\cdot\right\vert ^{2n}}{w\left\vert \cdot\right\vert ^{2\theta}}\right)
^{\frac{1}{2}},\\
\qquad and\\
\left(
\begin{array}
[c]{l}%
\int\left\vert \widehat{a}\widehat{\xi}\right\vert ^{2n}\frac{\left\vert
\cdot\right\vert ^{2n}}{w\left\vert \cdot\right\vert ^{2\theta}}+\\
+\sum\limits_{k<\theta-n}\frac{1}{k!}\int\left\vert \widehat{x}\widehat{\xi
}\right\vert ^{2k}\left\vert \widehat{a}\widehat{\xi}\right\vert ^{2n}%
\frac{\left\vert \cdot\right\vert ^{2n}}{w\left\vert \cdot\right\vert
^{2\theta}}%
\end{array}
\right)  ^{\frac{1}{2}},
\end{array}
\right.
\end{align*}

We now apply the estimate \ref{a2.23} for $\left\vert \left(  \widehat
{a}D\right)  ^{n}f_{\rho}\left(  x\right)  \right\vert $ to the right side of
\ref{a80}. In fact, using the approach of ??, it is clear that we obtain
\ref{a82} with $C_{n}^{\left(  \rho,w\right)  }$ replaced by $C_{n}^{\left(
\rho,w\right)  }\left(  r_{3}=0\right)  $ i.e. when $n\leq\left\lfloor
\underline{\kappa}\right\rfloor <\theta$,%
\begin{align}
&  \left\vert \frac{\sqrt{2\pi}}{n!}\left(  \overline{\widehat{g_{n}}}\left(
a\xi\right)  \left(  ia\xi\right)  \phi_{0}\widehat{\left(  aD\right)
^{n}f_{\rho}}\right)  ^{\vee}\left(  x\right)  \right\vert \nonumber\\
&  \leq\left\vert f_{\rho}\right\vert _{w,\theta}\frac{C_{n}^{\left(
\rho,w\right)  }\left(  r_{3}=0\right)  }{\left(  2\pi\right)  ^{\frac{d}{2}}%
}\frac{\left\vert a\right\vert ^{n+1}}{\left(  n+1\right)  !}\times\nonumber\\
&  \times\left(  \left\Vert \widehat{a}D\widehat{\phi_{0}}\right\Vert
_{1}+\left(  \frac{3^{\theta-n+1}-1}{2}\right)  ^{\frac{1}{2}}\left(
\sum\limits_{j=0}^{\theta-n}\frac{\left(  \left\vert x\right\vert
^{2}+\left\vert a\right\vert ^{2}\right)  ^{j}}{j!}\right)  ^{\frac{1}{2}%
}\left(  \left\Vert \widehat{a}D\widehat{\phi_{0}}\right\Vert _{1}%
\sum\limits_{l=0}^{\theta-n}\frac{\left\Vert \left\vert \cdot\right\vert
^{2l}\widehat{a}D\widehat{\phi_{0}}\right\Vert _{1}}{l!}\right)  ^{\frac{1}%
{2}}\right) \label{a841}%
\end{align}
\smallskip

\underline{\textbf{The second term} $\frac{\sqrt{2\pi}}{n!}\left(  \left(
ia\xi\right)  ^{n+1}\overline{\widehat{g_{n}}}\left(  a\xi\right)
\phi_{\infty}f_{F}\right)  ^{\vee}$ \textbf{in }\ref{a35}}\medskip

This corresponds to ?? \ref{a35} ?? for $w\in W3.2$ but the calculations are
much simpler:%
\begin{align*}
\left(  \left(  ia\xi\right)  ^{n+1}\overline{\widehat{g_{n}}}\left(
a\xi\right)  \phi_{\infty}f_{F}\right)  ^{\vee} &  =\left(  \phi_{\infty}%
\frac{\left(  ia\xi\right)  ^{n+1}}{\sqrt{w}\left\vert \cdot\right\vert
^{\theta}}\overline{\widehat{g_{n}}}\left(  a\xi\right)  \sqrt{w}\left\vert
\cdot\right\vert ^{\theta}f_{F}\right)  ^{\vee}\\
&  =\left(  \phi_{\infty}\frac{\left(  ia\xi\right)  ^{n+1}}{\sqrt
{w}\left\vert \cdot\right\vert ^{\theta}}\overline{\widehat{g_{n}}}\left(
a\xi\right)  \right)  ^{\vee}\ast\left(  \sqrt{w}\left\vert \cdot\right\vert
^{\theta}f_{F}\right)  ^{\vee},
\end{align*}

so that by Young's convolution estimate \ref{1.056} and Parceval's theorem we
get%
\begin{align}
&  \left\vert \left(  \left(  ia\xi\right)  ^{n+1}\overline{\widehat{g_{n}}%
}\left(  a\xi\right)  \phi_{\infty}f_{F}\right)  ^{\vee}\left(  x\right)
\right\vert \nonumber\\
&  =\left\vert \left(  \left(  \phi_{\infty}\frac{\left(  ia\xi\right)
^{n+1}}{\sqrt{w}\left\vert \cdot\right\vert ^{\theta}}\overline{\widehat
{g_{n}}}\left(  a\xi\right)  \right)  ^{\vee}\ast\left(  \sqrt{w}\left\vert
\cdot\right\vert ^{\theta}f_{F}\right)  ^{\vee}\right)  \left(  x\right)
\right\vert \nonumber\\
&  \leq\left\Vert \left(  \phi_{\infty}\frac{\left(  ia\xi\right)  ^{n+1}%
}{\sqrt{w}\left\vert \cdot\right\vert ^{\theta}}\overline{\widehat{g_{n}}%
}\left(  a\xi\right)  \right)  ^{\vee}\right\Vert _{2}\left\Vert \left(
\sqrt{w}\left\vert \cdot\right\vert ^{\theta}f_{F}\right)  ^{\vee}\right\Vert
_{2}\nonumber\\
&  =\left\Vert \phi_{\infty}\frac{\left(  ia\xi\right)  ^{n+1}}{\sqrt
{w}\left\vert \cdot\right\vert ^{\theta}}\overline{\widehat{g_{n}}}\left(
a\xi\right)  \right\Vert _{2}\left\Vert \sqrt{w}\left\vert \cdot\right\vert
^{\theta}f_{F}\right\Vert _{2}\nonumber\\
&  =\left\Vert \phi_{\infty}\frac{\left(  ia\xi\right)  ^{n+1}}{\sqrt
{w}\left\vert \cdot\right\vert ^{\theta}}\overline{\widehat{g_{n}}}\left(
a\xi\right)  \right\Vert _{2}\left\vert f\right\vert _{w,\theta}\nonumber\\
&  \leq\left\Vert \frac{\left(  ia\xi\right)  ^{n+1}}{\sqrt{w}\left\vert
\cdot\right\vert ^{\theta}}\overline{\widehat{g_{n}}}\left(  a\xi\right)
\right\Vert _{2}\left\vert f\right\vert _{w,\theta}\nonumber\\
&  =\left(  \int\frac{\left(  a\xi\right)  ^{2\left(  n+1\right)  }%
}{w\left\vert \cdot\right\vert ^{2\theta}}\left\vert \widehat{g_{n}}\left(
a\xi\right)  \right\vert ^{2}d\xi\right)  ^{1/2}\left\vert f\right\vert
_{w,\theta}.\label{a58}%
\end{align}
\smallskip

Applying \ref{a841} and \ref{a58} to \ref{a35} we get for the \fbox{case
$n\leq\left\lfloor \underline{\kappa}\right\rfloor <\theta$}:\smallskip%
\begin{align}
& \left\vert \mathcal{R}_{n+1}f_{\rho}\left(  x,a\right)  \right\vert
\nonumber\\
& \leq\frac{\sqrt{2\pi}}{n!}\left\vert \left(  \overline{\widehat{g_{n}}%
}\left(  a\xi\right)  \left(  ia\xi\right)  \phi_{0}\widehat{\left(
aD\right)  ^{n}f_{\rho}}\right)  ^{\vee}\left(  x\right)  \right\vert
+\frac{\sqrt{2\pi}}{n!}\left\vert \left(  \left(  ia\xi\right)  ^{n+1}%
\overline{\widehat{g_{n}}}\left(  a\xi\right)  \phi_{\infty}f_{F}\right)
^{\vee}\left(  x\right)  \right\vert \nonumber\\
& \leq\left(
\begin{array}
[c]{l}%
\frac{C_{n}^{\left(  \rho,w\right)  }\left(  r_{3}=0\right)  }{\left(
2\pi\right)  ^{d/2}}\frac{\left\vert a\right\vert ^{n+1}}{\left(  n+1\right)
!}\times\\
\times\left(  \left\Vert \widehat{a}D\widehat{\phi_{0}}\right\Vert
_{1}+\left(  \frac{3^{\theta-n+1}-1}{2}\right)  ^{\frac{1}{2}}\left(
\sum\limits_{j=0}^{\theta-n}\frac{\left(  \left\vert x\right\vert
^{2}+\left\vert a\right\vert ^{2}\right)  ^{j}}{j!}\right)  ^{\frac{1}{2}%
}\left(  \left\Vert \widehat{a}D\widehat{\phi_{0}}\right\Vert _{1}%
\sum\limits_{l=0}^{\theta-n}\frac{1}{l!}\left\Vert \left\vert \cdot\right\vert
^{2l}\widehat{a}D\widehat{\phi_{0}}\right\Vert _{1}\right)  ^{\frac{1}{2}%
}\right)
\end{array}
\right)  \left\vert f_{\rho}\right\vert _{w,\theta}+\nonumber\\
& \qquad\qquad+\frac{\sqrt{2\pi}}{n!}\left(  \int\frac{\left(  a\xi\right)
^{2\left(  n+1\right)  }}{w\left\vert \cdot\right\vert ^{2\theta}}\left\vert
\widehat{g_{n}}\left(  a\xi\right)  \right\vert ^{2}\right)  ^{1/2}\left\vert
f\right\vert _{w,\theta}.\label{a1001}%
\end{align}
\medskip

\underline{\textbf{A remainder estimate when} $w\in W3.1$}\medskip

Combining \ref{a101} and \ref{a100} we get for $n\leq\left\lfloor
\underline{\kappa}\right\rfloor $:%
\begin{multline}
\left\vert \mathcal{R}_{n+1}f_{\rho}\left(  x,a\right)  \right\vert
\nonumber\\
\leq\left\vert f_{\rho}\right\vert _{w,\theta}\left\{
\begin{array}
[c]{l}%
\frac{\sqrt{2\pi}}{n!}\left(  \int\frac{\left(  a\xi\right)  ^{2\left(
n+1\right)  }\left\vert \widehat{g_{n}}\left(  a\xi\right)  \right\vert ^{2}%
}{w\left\vert \cdot\right\vert ^{2\theta}}\right)  ^{\frac{1}{2}},\qquad
if\text{ \ }\theta\leq\left\lfloor \underline{\kappa}\right\rfloor ,\\
\medskip\\%
\begin{array}
[c]{l}%
\left.
\begin{array}
[c]{l}%
\frac{\sqrt{2\pi}}{n!}\left(  \int\frac{\left(  a\xi\right)  ^{2\left(
n+1\right)  }\left\vert \widehat{g_{n}}\left(  a\xi\right)  \right\vert ^{2}%
}{w\left\vert \cdot\right\vert ^{2\theta}}\right)  ^{\frac{1}{2}}+\\
+\frac{C_{n}^{\left(  \rho,w\right)  }\left(  r_{3}=0\right)  }{\left(
2\pi\right)  ^{d/2}}\frac{\left\vert a\right\vert ^{n+1}}{\left(  n+1\right)
!}\times\\
\times\left(  \left\Vert \widehat{a}D\widehat{\phi_{0}}\right\Vert
_{1}+\left(  \frac{3^{\theta-n+1}-1}{2}\right)  ^{\frac{1}{2}}\left(
\sum\limits_{j=0}^{\theta-n}\frac{\left(  \left\vert x\right\vert
^{2}+\left\vert a\right\vert ^{2}\right)  ^{j}}{j!}\right)  ^{\frac{1}{2}%
}\left(  \left\Vert \widehat{a}D\widehat{\phi_{0}}\right\Vert _{1}%
\sum\limits_{l=0}^{\theta-n}\frac{\left\Vert \left\vert \cdot\right\vert
^{2l}\widehat{a}D\widehat{\phi_{0}}\right\Vert _{1}}{l!}\right)  ^{\frac{1}%
{2}}\right)
\end{array}
\right\} \\
\qquad\qquad\qquad\qquad\qquad\qquad\qquad\qquad\qquad\qquad\qquad\qquad\qquad
if\text{ \ }\theta>\left\lfloor \underline{\kappa}\right\rfloor .
\end{array}
\end{array}
\right. \label{A102}%
\end{multline}
\medskip

\fbox{If $n\leq\left\lfloor \underline{\kappa}\right\rfloor -1$} then we have
an order of convergence of $\left\lfloor \underline{\kappa}\right\rfloor $
because $\left\vert \widehat{g_{n}}\left(  a\xi\right)  \right\vert \leq
\frac{1}{\sqrt{2\pi}}\frac{1}{n+1}$ implies,\smallskip%
\begin{align*}
&  \left\vert \mathcal{R}_{n+1}f_{\rho}\left(  x,a\right)  \right\vert \\
&  \leq\frac{\left\vert a\right\vert ^{n+1}}{\left(  n+1\right)  !}\left\vert
f_{\rho}\right\vert _{w,\theta}\times\\
&  \times\left\{
\begin{array}
[c]{l}%
\frac{\sqrt{2\pi}}{n!}\left(  \int\frac{\left(  \widehat{a}\xi\right)
^{2\left(  n+1\right)  }}{w\left\vert \cdot\right\vert ^{2\theta}}\right)
^{\frac{1}{2}},\qquad if\text{ \ }\theta\leq\left\lfloor \underline{\kappa
}\right\rfloor ,\\
\medskip\\%
\begin{array}
[c]{l}%
\begin{array}
[c]{l}%
\frac{\sqrt{2\pi}}{n!}\left(  \int\frac{\left(  \widehat{a}\xi\right)
^{2\left(  n+1\right)  }}{w\left\vert \cdot\right\vert ^{2\theta}}\right)
^{\frac{1}{2}}+\\
+\frac{C_{n}^{\left(  \rho,w\right)  }\left(  r_{3}=0\right)  }{\left(
2\pi\right)  ^{d/2}}\left(  \left\Vert \widehat{a}D\widehat{\phi_{0}%
}\right\Vert _{1}+\left(  \frac{3^{\theta-n+1}-1}{2}\right)  ^{\frac{1}{2}%
}\left(  \sum\limits_{j=0}^{\theta-n}\frac{\left(  \left\vert x\right\vert
^{2}+\left\vert a\right\vert ^{2}\right)  ^{j}}{j!}\right)  ^{\frac{1}{2}%
}\left(  \left\Vert \widehat{a}D\widehat{\phi_{0}}\right\Vert _{1}%
\sum\limits_{l=0}^{\theta-n}\frac{1}{l!}\left\Vert \left\vert \cdot\right\vert
^{2l}\widehat{a}D\widehat{\phi_{0}}\right\Vert _{1}\right)  ^{\frac{1}{2}%
}\right)
\end{array}
\\
\qquad\qquad\qquad\qquad\qquad\qquad\qquad\qquad\qquad\qquad\qquad\qquad
if\text{ \ }\theta>\left\lfloor \underline{\kappa}\right\rfloor .
\end{array}
\end{array}
\right.
\end{align*}
\medskip

\fbox{If $n=\left\lfloor \underline{\kappa}\right\rfloor $} we use the
technique used to prove inequality \ref{a661}, but with $r_{3}=0$ and $\kappa$
replaced by $\underline{\kappa}$, to obtain\smallskip%
\[
\sqrt{2\pi}\left(  \int\frac{\left(  a\xi\right)  ^{2\left\lceil
\kappa\right\rceil }}{w\left\vert \cdot\right\vert ^{2\theta}}\left\vert
\widehat{g_{\left\lfloor \kappa\right\rfloor }}\left(  a\xi\right)
\right\vert ^{2}\right)  ^{\frac{1}{2}}\leq\left\vert a\right\vert
^{\underline{\kappa}}\left(  2+\frac{1}{n+1}\right)  \left(  \int%
\frac{\left\vert \cdot\right\vert ^{2\underline{\kappa}}}{w\left\vert
\cdot\right\vert ^{2\theta}}\right)  ^{\frac{1}{2}},
\]

i.e.%
\[
\sqrt{2\pi}\left(  \int\frac{\left(  a\xi\right)  ^{2\left(  n+1\right)  }%
}{w\left\vert \cdot\right\vert ^{2\theta}}\left\vert \widehat{g_{n}}\left(
a\xi\right)  \right\vert ^{2}\right)  ^{\frac{1}{2}}\leq\left\vert
a\right\vert ^{\underline{\kappa}}\left(  2+\frac{1}{n+1}\right)  \left(
\int\frac{\left\vert \cdot\right\vert ^{2\underline{\kappa}}}{w\left\vert
\cdot\right\vert ^{2\theta}}\right)  ^{\frac{1}{2}},
\]

and this integral exists by Theorem \ref{Thm_property_wt_fn_W3}. Thus when
$n=\left\lfloor \underline{\kappa}\right\rfloor $, \ref{a102} becomes%
\begin{align}
&  \left\vert \mathcal{R}_{n+1}f_{\rho}\left(  x,a\right)  \right\vert
\nonumber\\
&  \leq\left\vert f_{\rho}\right\vert _{w,\theta}\left\{
\begin{array}
[c]{l}%
\left\vert a\right\vert ^{\underline{\kappa}}\left(  \frac{2}{n!}+\frac
{1}{\left(  n+1\right)  !}\right)  \left(  \int\frac{\left\vert \cdot
\right\vert ^{2\underline{\kappa}}}{w\left\vert \cdot\right\vert ^{2\theta}%
}\right)  ^{\frac{1}{2}},\text{\quad}if\text{ \ }\theta\leq\left\lfloor
\underline{\kappa}\right\rfloor ,\\
\medskip\\
\left\vert a\right\vert ^{\underline{\kappa}}\left(  \frac{2}{n!}+\frac
{1}{\left(  n+1\right)  !}\right)  \left(  \int\frac{\left\vert \cdot
\right\vert ^{2\underline{\kappa}}}{w\left\vert \cdot\right\vert ^{2\theta}%
}\right)  ^{\frac{1}{2}}+\frac{C_{n}^{\left(  \rho,w\right)  }\left(
r_{3}=0\right)  }{\left(  2\pi\right)  ^{d}}\frac{\left\vert a\right\vert
^{n+1}}{\left(  n+1\right)  !}\times\\
\left(  \left\Vert \widehat{a}D\widehat{\phi_{0}}\right\Vert _{1}+\left(
\frac{3^{\theta-n+1}-1}{2}\right)  ^{\frac{1}{2}}\left(  \sum\limits_{j=0}%
^{\theta-n}\frac{\left(  \left\vert x\right\vert ^{2}+\left\vert a\right\vert
^{2}\right)  ^{j}}{j!}\right)  ^{\frac{1}{2}}\left(  \left\Vert \widehat
{a}D\widehat{\phi_{0}}\right\Vert _{1}\sum\limits_{l=0}^{\theta-n}%
\frac{\left\Vert \left\vert \cdot\right\vert ^{2l}\widehat{a}D\widehat
{\phi_{0}}\right\Vert _{1}}{l!}\right)  ^{\frac{1}{2}}\right) \\
\qquad\qquad\qquad\qquad\qquad\qquad\qquad\qquad\qquad\qquad\qquad\qquad
\qquad\qquad if\text{ \ }\theta>\left\lfloor \underline{\kappa}\right\rfloor .
\end{array}
\right. \label{a001}%
\end{align}

Hence the order of convergence is at least $\underline{\kappa}=\min\kappa$.
However, there may be a value of $\kappa$, say $\kappa_{\infty}$, such that as
$\underline{\kappa}\rightarrow\underline{\kappa_{\infty}^{-}}$ we have
$\int\frac{\left\vert \cdot\right\vert ^{2\kappa}}{w\left\vert \cdot
\right\vert ^{2\theta}}\rightarrow\infty$. Thus the coefficient of the power
of $\left\vert a\right\vert $ varies.

\begin{remark}
?? Regarding \ref{a2.35}: calculation of $\left\Vert \left(  \frac{\left(
ia\xi\right)  ^{n+1}\overline{\widehat{g_{n}}}\left(  a\xi\right)  }{\sqrt
{w}\left\vert \cdot\right\vert ^{\theta}}\right)  ^{\vee}\right\Vert _{2}$
when $w$ is radial. From Theorem \ref{Thm_Integ_u(xy)f(|x|)dx},%
\[
\int_{\left\vert x\right\vert \leq r}u\left(  \xi x\right)  f\left(
\left\vert x\right\vert \right)  dx=\int_{\left\vert x\right\vert \leq
r}u\left(  \left\vert \xi\right\vert x_{k}\right)  f\left(  \left\vert
x\right\vert \right)  dx,\quad r\geq0.
\]

\end{remark}

\subsection{Approach 2\label{SbSect_TaylorDataApproach2}}

This approach enables us to obtain \textbf{significantly improved estimates}
for the extended splines compared to approach I above.

Suppose $G$ is the unique basis function given by $\widehat{G}=\frac
{1}{w\left\vert \cdot\right\vert ^{2\theta}}$ (see Theorem
\ref{Thm_basis_smth_W3.1}). Then using the expression \ref{a114} for
$e^{ia\xi}$ we can write the rewrite the estimates \ref{a101} and \ref{a58}
as,%
\begin{align}
\left\vert \mathcal{R}_{n+1}f_{\rho}\left(  x,a\right)  \right\vert  &
\leq\frac{\sqrt{2\pi}}{n!}\left(  \int\left\vert \left(  ia\xi\right)
^{n+1}\widehat{g_{n}}\left(  a\xi\right)  \right\vert ^{2}\widehat{G}\right)
^{\frac{1}{2}}\left\vert f_{\rho}\right\vert _{w,\theta}\nonumber\\
& =\frac{\sqrt{2\pi}}{n!}\left(  \int\left\vert \frac{n!}{\sqrt{2\pi}}\left(
e^{ia\xi}-\sum\limits_{k=0}^{n}\frac{\left(  ia\xi\right)  ^{k}}{k!}\right)
\right\vert ^{2}\widehat{G}\right)  ^{\frac{1}{2}}\left\vert f_{\rho
}\right\vert _{w,\theta}\nonumber\\
& =\left(  \int\left\vert e^{ia\xi}-\sum_{k=0}^{n}\frac{\left(  ia\xi\right)
^{k}}{k!}\right\vert ^{2}\widehat{G}\right)  ^{\frac{1}{2}}\left\vert f_{\rho
}\right\vert _{w,\theta},\label{a95}%
\end{align}

where $\mathcal{R}_{n+1}f$ is defined by \ref{a50.5} and $f_{\rho}$ is defined
before \ref{a31}.

We will next derive the following estimates for the integral in \ref{a95}:

\begin{theorem}
\label{Thm_Tayor_rem_estim_W3.1}Suppose the weight function $w$ has property
W2 and property W3.1 or W3.3 for order $\theta$ and parameter $\kappa$. Set
$m=\left\lfloor \underline{\kappa}\right\rfloor $. Then if $G$ is the basis
function given by $\widehat{G}=\frac{1}{w\left\vert \cdot\right\vert
^{2\theta}}$ we have $G\in C_{B}^{\left(  2m\right)  }$ and:

\begin{enumerate}
\item for all $a\in\mathbb{R}^{d}$ and $0\leq n\leq2m$,
\begin{align}
& \left(  2\pi\right)  ^{-\frac{d}{2}}\int\left\vert e^{ia\xi}-\sum_{k=0}%
^{n}\frac{\left(  ia\xi\right)  ^{k}}{k!}\right\vert ^{2}\widehat{G}\left(
\xi\right)  d\xi\label{a2.15}\\
& =2G\left(  0\right)  -\sum_{k=0}^{n}\frac{1}{k!}\left(  \left(  aD\right)
^{k}G\right)  \left(  -a\right)  -\sum_{k=0}^{n}\frac{1}{k!}\left(
-aD\right)  ^{k}G\left(  a\right)  +\nonumber\\
& \qquad\qquad+%
{\displaystyle\sum\limits_{\substack{j,k\leq n \\j+k>n}}}
\frac{\left(  -1\right)  ^{j}}{j!k!}\left(  \left(  aD\right)  ^{j+k}G\right)
\left(  0\right)  .\label{a2.31}%
\end{align}

\item Suppose further that $\left(  aD\right)  ^{2m+1}G\in L^{1}$ when
$\left\vert a\right\vert =1$. Then%
\begin{equation}
\left(  2\pi\right)  ^{-\frac{d}{2}}\int\left\vert e^{ia\xi}-\sum_{k=0}%
^{n}\frac{\left(  ia\xi\right)  ^{k}}{k!}\right\vert ^{2}\widehat{G}\left(
\xi\right)  d\xi\leq\mathcal{R}_{2n+1}G\left(  a\right)  ,\quad x,a\in
\mathbb{R}^{d},\label{a2.30}%
\end{equation}

where $\mathcal{R}_{2n+1}G$ is defined by%
\begin{align}
\left(  \mathcal{R}_{2n+1}G\right)  \left(  a\right)   & :=\sum_{k=n+1}%
^{2n}\frac{1}{k!}\left\{  \left(  \mathcal{R}_{2n+1-k}\left(  \left(
-aD\right)  ^{k}G\right)  \right)  \left(  0,a\right)  +\left(  \mathcal{R}%
_{2n+1-k}\left(  \left(  aD\right)  ^{k}G\right)  \right)  \left(
0,-a\right)  \right\}  +\nonumber\\
& \qquad\qquad+\left(  \mathcal{R}_{2n+1}G\right)  \left(  a,-a\right)
+\left(  \mathcal{R}_{2n+1}G\right)  \left(  -a,a\right) \label{a2.32}\\
& =2\int_{0}^{1}\left(  \sum\limits_{j=0}^{n-1}\frac{\left(  -1\right)
^{j}\left(  1-s\right)  ^{j}}{j!\left(  2n-j\right)  !}-\frac{s^{2n}}{\left(
2n\right)  !}\right)  \left(  \left(  aD\right)  ^{2n+1}\operatorname{Re}%
G\right)  \left(  sa\right)  ds.\label{a2.33}%
\end{align}

\item Finally,%
\begin{equation}
\int_{0}^{1}\left\vert \sum\limits_{j=0}^{n-1}\frac{\left(  -1\right)
^{j}\left(  1-s\right)  ^{j}}{j!\left(  2n-j\right)  !}-\frac{s^{2n}}{\left(
2n\right)  !}\right\vert ds=\frac{1}{\left(  2n+1\right)  \left(  n!\right)
^{2}},\quad n\geq1.\label{a96}%
\end{equation}

\end{enumerate}
\end{theorem}

\begin{proof}
\textbf{Part 1} Expanding gives%
\begin{align*}
\left\vert e^{ia\xi}-\sum_{k=0}^{n}\frac{\left(  ia\xi\right)  ^{k}}%
{k!}\right\vert ^{2} &  =\left(  \overline{e^{ia\xi}-\sum_{l=0}^{n}%
\frac{\left(  ia\xi\right)  ^{l}}{l!}}\right)  \left(  e^{ia\xi}-\sum
_{k=0}^{n}\frac{\left(  ia\xi\right)  ^{k}}{k!}\right) \\
&  =\left(  e^{-ia\xi}-\sum_{l=0}^{n}\frac{\left(  -ia\xi\right)  ^{l}}%
{l!}\right)  \left(  e^{ia\xi}-\sum_{k=0}^{n}\frac{\left(  ia\xi\right)  ^{k}%
}{k!}\right) \\
&  =1-e^{-ia\xi}\sum_{k=0}^{n}\frac{\left(  ia\xi\right)  ^{k}}{k!}-e^{ia\xi
}\sum_{l=0}^{n}\frac{\left(  -ia\xi\right)  ^{l}}{l!}+\sum_{k,l=0}^{n}%
\frac{\left(  -ia\xi\right)  ^{l}}{l!}\frac{\left(  ia\xi\right)  ^{k}}{k!}\\
&  =1-e^{-ia\xi}\sum_{k=0}^{n}\frac{\left(  ia\xi\right)  ^{k}}{k!}-e^{ia\xi
}\sum_{l=0}^{n}\frac{\left(  -1\right)  ^{l}\left(  ia\xi\right)  ^{l}}%
{l!}+\sum_{l,k=0}^{n}\frac{\left(  -1\right)  ^{l}\left(  ia\xi\right)
^{l+k}}{k!l!},
\end{align*}

and since we know that $e^{ia\xi}\widehat{G}\left(  \xi\right)  =\left(
G\left(  \cdot+a\right)  \right)  ^{\wedge}\left(  \xi\right)  $ and $\left(
ia\xi\right)  ^{k}\widehat{G}\left(  \xi\right)  =\left(  \left(  aD\right)
^{k}G\right)  ^{\wedge}\left(  \xi\right)  $ it follows that%
\begin{align*}
&  \left\vert e^{ia\xi}-\sum_{k=0}^{n}\frac{\left(  ia\xi\right)  ^{k}}%
{k!}\right\vert ^{2}\widehat{G}\\
&  =\left(  1-\sum_{k=0}^{n}\frac{\left(  ia\xi\right)  ^{k}}{k!}e^{-ia\xi
}-\sum_{l=0}^{n}\frac{\left(  -ia\xi\right)  ^{l}}{l!}e^{ia\xi}+\sum
_{l,k=0}^{n}\frac{\left(  -1\right)  ^{l}}{k!l!}\left(  ia\xi\right)
^{l+k}\right)  \widehat{G}\\
&  =\widehat{G}-\sum_{k=0}^{n}\frac{\left(  ia\xi\right)  ^{k}e^{-ia\xi
}\widehat{G}}{k!}-\sum_{l=0}^{n}\frac{\left(  -ia\xi\right)  ^{l}e^{ia\xi
}\widehat{G}}{l!}+\sum_{l,k=0}^{n}\frac{\left(  -1\right)  ^{l}}{k!l!}\left(
ia\xi\right)  ^{l+k}\widehat{G}\\
&  =\widehat{G}-\sum_{k=0}^{n}\frac{\left(  ia\xi\right)  ^{k}\widehat
{G\left(  \cdot-a\right)  }}{k!}-\sum_{l=0}^{n}\frac{\left(  -ia\xi\right)
^{l}\widehat{G\left(  \cdot+a\right)  }}{l!}+\sum_{l,k=0}^{n}\frac{\left(
-1\right)  ^{l}}{k!l!}\left(  ia\xi\right)  ^{l+k}\widehat{G}\\
&  =\widehat{G}-\sum_{k=0}^{n}\frac{\widehat{\left(  aD\right)  ^{k}G\left(
\cdot-a\right)  }}{k!}-\sum_{l=0}^{n}\frac{\widehat{\left(  -aD\right)
^{l}G\left(  \cdot+a\right)  }}{l!}+\sum_{l,k=0}^{n}\frac{\left(  -1\right)
^{l}}{k!l!}\widehat{\left(  aD\right)  ^{l+k}G}\\
&  =\left(  G-\sum_{k=0}^{n}\frac{\left(  \left(  aD\right)  ^{k}G\right)
\left(  \cdot-a\right)  }{k!}-\sum_{l=0}^{n}\frac{\left(  \left(  -aD\right)
^{l}G\right)  \left(  \cdot+a\right)  }{l!}+\sum_{l,k=0}^{n}\frac{\left(
-1\right)  ^{l}}{k!l!}\left(  aD\right)  ^{l+k}G\right)  ^{\wedge},
\end{align*}

and hence%
\begin{align}
\left(  2\pi\right)  ^{-\frac{d}{2}}\int &  \left\vert e^{ia\xi}-\sum
_{k=0}^{n}\frac{\left(  ia\xi\right)  ^{k}}{k!}\right\vert ^{2}\widehat
{G}\left(  \xi\right)  d\xi\nonumber\\
&  =G\left(  0\right)  -\sum_{k=0}^{n}\frac{\left(  \left(  aD\right)
^{k}G\right)  \left(  -a\right)  }{k!}-\sum_{k=0}^{n}\frac{\left(  \left(
-aD\right)  ^{k}G\right)  \left(  a\right)  }{k!}+\sum_{l,k=0}^{n}%
\frac{\left(  -1\right)  ^{l}}{k!l!}\left(  \left(  aD\right)  ^{l+k}G\right)
\left(  0\right)  ,\label{a2.20}%
\end{align}

which proves \ref{a2.31} and this part.\medskip

\fbox{\textbf{Part 2}} Regarding the last term in \ref{a2.20}:%
\begin{align*}
\sum_{l,k=0}^{n}\frac{\left(  -1\right)  ^{l}}{k!l!}\left(  \left(  aD\right)
^{l+k}G\right)  \left(  0\right)   & =%
{\displaystyle\sum\limits_{l+k\leq n}}
\frac{\left(  -1\right)  ^{l}}{k!l!}\left(  \left(  aD\right)  ^{l+k}G\right)
\left(  0\right)  +%
{\displaystyle\sum\limits_{\substack{l+k>n \\l,k\leq n}}}
\frac{\left(  -1\right)  ^{l}}{k!l!}\left(  \left(  aD\right)  ^{l+k}G\right)
\left(  0\right) \\
& =%
{\displaystyle\sum\limits_{q=0}^{n}}
{\displaystyle\sum\limits_{l+k=q}}
\frac{\left(  -1\right)  ^{l}}{k!l!}\left(  \left(  aD\right)  ^{l+k}G\right)
\left(  0\right)  +%
{\displaystyle\sum\limits_{\substack{l+k>n \\l,k\leq n}}}
\frac{\left(  -1\right)  ^{l}}{k!l!}\left(  \left(  aD\right)  ^{l+k}G\right)
\left(  0\right) \\
& =%
{\displaystyle\sum\limits_{q=0}^{n}}
\left(
{\displaystyle\sum\limits_{l+k=q}}
\frac{\left(  -1\right)  ^{l}}{k!l!}\right)  \left(  \left(  aD\right)
^{q}G\right)  \left(  0\right)  +%
{\displaystyle\sum\limits_{\substack{l+k>n \\l,k\leq n}}}
\frac{\left(  -1\right)  ^{l}}{k!l!}\left(  \left(  aD\right)  ^{l+k}G\right)
\left(  0\right)  ,
\end{align*}

and since: when $q>0$,%
\begin{align}%
{\displaystyle\sum\limits_{l+k=q}}
\frac{\left(  -1\right)  ^{l}}{k!l!}  & =%
{\displaystyle\sum\limits_{\left\vert \alpha\right\vert =q}}
\frac{\left(  \left(  -1\right)  ^{\alpha_{1}}1^{\alpha_{2}}\right)  \left(
1^{\alpha_{1}}1^{\alpha_{2}}\right)  }{\alpha!}=%
{\displaystyle\sum\limits_{\left\vert \alpha\right\vert =q}}
\frac{\left(  -1,1\right)  ^{\alpha}\left(  1,1\right)  ^{\alpha}}{\alpha
!}=\frac{\left(  \left(  -1,1\right)  \cdot\left(  1,1\right)  \right)  ^{q}%
}{q!}=\nonumber\\
& =0,\label{a2.13}%
\end{align}

we have%
\begin{align*}
\sum_{l,k=0}^{n}\frac{\left(  -1\right)  ^{l}}{k!l!}\left(  \left(  aD\right)
^{l+k}G\right)  \left(  0\right)   & =G\left(  0\right)  +%
{\displaystyle\sum\limits_{\substack{l+k>n \\l,k\leq n}}}
\frac{\left(  -1\right)  ^{l}}{k!l!}\left(  \left(  aD\right)  ^{l+k}G\right)
\left(  0\right) \\
& =G\left(  0\right)  +%
{\displaystyle\sum\limits_{q=n+1}^{2n}}
{\displaystyle\sum\limits_{\substack{l+k=q \\l,k\leq n}}}
\frac{\left(  -1\right)  ^{l}}{k!l!}\left(  \left(  aD\right)  ^{l+k}G\right)
\left(  0\right) \\
& =G\left(  0\right)  +%
{\displaystyle\sum\limits_{q=n+1}^{2n}}
\left(
{\displaystyle\sum\limits_{\substack{l+k=q \\l,k\leq n}}}
\frac{\left(  -1\right)  ^{l}}{k!l!}\right)  \left(  \left(  aD\right)
^{q}G\right)  \left(  0\right)  ,
\end{align*}

so that%
\begin{multline*}
\left(  2\pi\right)  ^{-\frac{d}{2}}\int\left\vert e^{ia\xi}-\sum_{k=0}%
^{n}\frac{\left(  ia\xi\right)  ^{k}}{k!}\right\vert ^{2}\widehat{G}\left(
\xi\right)  d\xi\\
=2G\left(  0\right)  -\sum_{k=0}^{n}\frac{\left(  \left(  aD\right)
^{k}G\right)  \left(  -a\right)  }{k!}-\sum_{k=0}^{n}\frac{\left(  \left(
-aD\right)  ^{k}G\right)  \left(  a\right)  }{k!}+%
{\displaystyle\sum\limits_{\substack{l+k>n \\l,k\leq n}}}
\frac{\left(  -1\right)  ^{l}}{k!l!}\left(  \left(  aD\right)  ^{l+k}G\right)
\left(  0\right)  .
\end{multline*}

But%
\begin{equation}%
{\displaystyle\sum\limits_{\substack{l+k>n \\l,k\leq n}}}
\frac{\left(  -1\right)  ^{l}}{k!l!}\left(  \left(  aD\right)  ^{l+k}G\right)
=%
{\displaystyle\sum\limits_{q=n+1}^{2n}}
{\displaystyle\sum\limits_{\substack{l+k=q \\l,k\leq n}}}
\frac{\left(  -1\right)  ^{l}}{k!l!}\left(  \left(  aD\right)  ^{l+k}G\right)
=%
{\displaystyle\sum\limits_{q=n+1}^{2n}}
\left(
{\displaystyle\sum\limits_{\substack{l+k=q \\l,k\leq n }}}
\frac{\left(  -1\right)  ^{l}}{k!l!}\right)  \left(  aD\right)  ^{q}%
G,\label{a331}%
\end{equation}

so that%
\begin{align}
\left(  2\pi\right)  ^{-\frac{d}{2}} &  \int\left\vert e^{ia\xi}-\sum
_{k=0}^{n}\frac{\left(  ia\xi\right)  ^{k}}{k!}\right\vert ^{2}\widehat
{G}\left(  \xi\right)  d\xi\nonumber\\
&  =2G\left(  0\right)  -\sum_{k=0}^{n}\frac{\left(  \left(  aD\right)
^{k}G\right)  \left(  -a\right)  }{k!}-\sum_{k=0}^{n}\frac{\left(  \left(
-aD\right)  ^{k}G\right)  \left(  a\right)  }{k!}+\nonumber\\
&  \qquad\qquad+%
{\displaystyle\sum\limits_{q=n+1}^{2n}}
\left(
{\displaystyle\sum\limits_{\substack{l+k=q \\l,k\leq n}}}
\frac{\left(  -1\right)  ^{l}}{k!l!}\right)  \left(  aD\right)  ^{q}G\left(
0\right)  .\label{a337}%
\end{align}

Since $D^{\alpha}G\in L^{1}$ when $\left\vert \alpha\right\vert =2m+1$ we can
apply \ref{a1.55} with remainder form \ref{a111} to get%
\begin{align*}
G\left(  0\right)   & =G\left(  a-a\right)  =\sum_{k=0}^{2n}\frac{\left(
\left(  -aD\right)  ^{k}G\right)  \left(  a\right)  }{k!}+\left(
\mathcal{R}_{2n+1}G\right)  \left(  a,-a\right)  ,\\
G\left(  0\right)   & =G\left(  -a+a\right)  =\sum_{k=0}^{2n}\frac{\left(
\left(  aD\right)  ^{k}G\right)  \left(  -a\right)  }{k!}+\left(
\mathcal{R}_{2n+1}G\right)  \left(  -a,a\right)  ,
\end{align*}

so that \ref{a337} becomes%
\begin{align}
&  \left(  2\pi\right)  ^{-\frac{d}{2}}\int\left\vert e^{ia\xi}-\sum_{k=0}%
^{n}\frac{\left(  ia\xi\right)  ^{k}}{k!}\right\vert ^{2}\widehat{G}\left(
\xi\right)  d\xi\nonumber\\
&  =\sum_{k=n+1}^{2n}\frac{\left(  \left(  -aD\right)  ^{k}G\right)  \left(
a\right)  }{k!}+\sum_{k=n+1}^{2n}\frac{\left(  \left(  aD\right)
^{k}G\right)  \left(  -a\right)  }{k!}+%
{\displaystyle\sum\limits_{q=n+1}^{2n}}
\left(
{\displaystyle\sum\limits_{\substack{l+k=q \\l,k\leq n}}}
\frac{\left(  -1\right)  ^{l}}{k!l!}\right)  \left(  aD\right)  ^{q}G\left(
0\right)  +\nonumber\\
&  \qquad\qquad+\left(  \mathcal{R}_{2n+1}G\right)  \left(  a,-a\right)
+\left(  \mathcal{R}_{2n+1}G\right)  \left(  -a,a\right)  .\label{a2.16}%
\end{align}

Since $D^{\alpha}G\in L^{1}$ when $\left\vert \alpha\right\vert =2n+1$ we can
apply \ref{a1.55} with remainder \ref{a111} to $\left(  \left(  -aD\right)
^{k}G\right)  \left(  b\right)  $ about $b=0$ to get%
\[
\left(  \left(  aD\right)  ^{k}G\right)  \left(  b\right)  =\sum_{l=0}%
^{2n-k}\frac{\left(  \left(  bD\right)  ^{l}\left(  aD\right)  ^{k}G\right)
\left(  0\right)  }{l!}+\mathcal{R}_{2n+1-k}\left(  \left(  aD\right)
^{k}G\right)  \left(  0,b\right)  ,\quad0\leq k\leq2n.
\]

When $b=-a$ \ we get%
\begin{equation}
\left(  \left(  aD\right)  ^{k}G\right)  \left(  -a\right)  =\sum_{l=0}%
^{2n-k}\frac{\left(  -1\right)  ^{l}}{l!}\left(  \left(  aD\right)
^{k+l}G\right)  \left(  0\right)  +\mathcal{R}_{2n+1-k}\left(  \left(
aD\right)  ^{k}G\right)  \left(  0,-a\right)  ,\quad0\leq k\leq
2n,\label{a2.17}%
\end{equation}

and when $a\rightarrow-a$ and $b=a$ we get%
\begin{align}
\left(  \left(  -aD\right)  ^{k}G\right)   &  \left(  a\right) \nonumber\\
&  =\sum_{l=0}^{2n-k}\frac{\left(  \left(  aD\right)  ^{l}\left(  -aD\right)
^{k}G\right)  \left(  0\right)  }{l!}+\mathcal{R}_{2n+1-k}\left(  \left(
-aD\right)  ^{k}G\right)  \left(  0,a\right) \nonumber\\
&  =\sum_{l=0}^{2n-k}\frac{\left(  -1\right)  ^{k}}{l!}\left(  \left(
aD\right)  ^{k+l}G\right)  \left(  0\right)  +\mathcal{R}_{2n+1-k}\left(
\left(  -aD\right)  ^{k}G\right)  \left(  0,a\right)  ,\quad0\leq
k\leq2n.\label{a2.18}%
\end{align}

Now substitute the Taylor series expansions \ref{a2.17} and \ref{a2.18}\ into
the terms of the first and second summations of \ref{a2.16} to get%
\begin{align}
\sum_{k=n+1}^{2n} &  \frac{\left(  \left(  -aD\right)  ^{k}G\right)  \left(
a\right)  }{k!}\nonumber\\
&  =\sum_{k=n+1}^{2n}\frac{1}{k!}\left(  \sum_{l=0}^{2n-k}\frac{\left(
-1\right)  ^{k}}{l!}\left(  aD\right)  ^{k+l}G\left(  0\right)  +\mathcal{R}%
_{2n+1-k}\left(  \left(  -aD\right)  ^{k}G\right)  \left(  0,a\right)  \right)
\nonumber\\
&  =\sum_{k=n+1}^{2n}\frac{\left(  -1\right)  ^{k}}{k!}\sum_{l=0}^{2n-k}%
\frac{\left(  aD\right)  ^{k+l}G\left(  0\right)  }{l!}+\sum_{k=n+1}^{2n}%
\frac{1}{k!}\mathcal{R}_{2n+1-k}\left(  \left(  -aD\right)  ^{k}G\right)
\left(  0,a\right) \nonumber\\
&  =\sum_{k=n+1}^{2n}\sum_{l=0}^{2n-k}\frac{\left(  -1\right)  ^{k}}%
{k!l!}\left(  aD\right)  ^{k+l}G\left(  0\right)  +\ldots\nonumber\\
&  =\sum_{q=n+1}^{2n}\sum_{\substack{k+l=q \\k\geq n+1}}\frac{\left(
-1\right)  ^{k}}{k!l!}\left(  aD\right)  ^{k+l}G\left(  0\right)
+\ldots\nonumber\\
&  =\sum_{q=n+1}^{2n}\left(  \sum_{\substack{k+l=q \\k\geq n+1}}\frac{\left(
-1\right)  ^{k}}{k!l!}\right)  \left(  aD\right)  ^{q}G\left(  0\right)
+\ldots\nonumber\\
&  =\sum_{q=n+1}^{2n}\left(  \sum_{\substack{k+l=q \\k\geq n+1}}\frac{\left(
-1\right)  ^{k}}{k!l!}\right)  \left(  aD\right)  ^{q}G\left(  0\right)
+\sum_{k=n+1}^{2n}\frac{1}{k!}\mathcal{R}_{2n+1-k}\left(  \left(  -aD\right)
^{k}G\right)  \left(  0,a\right)  ,\label{a2.21}%
\end{align}

and so%
\begin{align}
& \sum_{k=n+1}^{2n}\frac{\left(  \left(  aD\right)  ^{k}G\right)  \left(
-a\right)  }{k!}\nonumber\\
& =\sum_{q=n+1}^{2n}\left(  \sum_{\substack{k+l=q \\k\geq n+1}}\frac{\left(
-1\right)  ^{k}}{k!l!}\right)  \left(  \left(  -aD\right)  ^{q}G\right)
\left(  0\right)  +\sum_{k=n+1}^{2n}\frac{1}{k!}\mathcal{R}_{2n+1-k}\left(
\left(  aD\right)  ^{k}G\right)  \left(  0,-a\right) \nonumber\\
& =\sum_{q=n+1}^{2n}\left(  -1\right)  ^{q}\left(  \sum_{\substack{k+l=q
\\k\geq n+1}}\frac{\left(  -1\right)  ^{k}}{k!l!}\right)  \left(  aD\right)
^{q}G\left(  0\right)  +\ldots\nonumber\\
& =\sum_{q=n+1}^{2n}\left(  -1\right)  ^{q}\left(  \sum_{\substack{k+l=q
\\k\geq n+1}}\frac{\left(  -1\right)  ^{-k}}{k!l!}\right)  \left(  aD\right)
^{q}G\left(  0\right)  +\ldots\nonumber\\
& =\sum_{q=n+1}^{2n}\left(  \sum_{\substack{k+l=q \\k\geq n+1}}\frac{\left(
-1\right)  ^{q-k}}{k!l!}\right)  \left(  aD\right)  ^{q}G\left(  0\right)
+\ldots\nonumber\\
& =\sum_{q=n+1}^{2n}\left(  \sum_{\substack{k+l=q \\k\geq n+1}}\frac{\left(
-1\right)  ^{l}}{k!l!}\right)  \left(  aD\right)  ^{q}G\left(  0\right)
+\sum_{k=n+1}^{2n}\frac{1}{k!}\mathcal{R}_{2n+1-k}\left(  \left(  aD\right)
^{k}G\right)  \left(  0,-a\right)  .\label{a2.22}%
\end{align}

Substituting \ref{a2.21} and \ref{a2.22} into \ref{a2.16} gives%
\begin{align*}
& \left(  2\pi\right)  ^{-\frac{d}{2}}\int\left\vert e^{ia\xi}-\sum_{k\leq
n}\frac{\left(  ia\xi\right)  ^{k}}{k!}\right\vert ^{2}\widehat{G}\left(
\xi\right)  d\xi\\
& =\sum_{q=n+1}^{2n}\left(  \sum_{\substack{k+l=q \\k\geq n+1}}\frac{\left(
-1\right)  ^{k}}{k!l!}\right)  \left(  aD\right)  ^{q}G\left(  0\right)
+\sum_{k=n+1}^{2n}\frac{1}{k!}\left(  \mathcal{R}_{2n+1}\left(  \left(
-aD\right)  ^{k}G\right)  \right)  \left(  0,a\right)  +\\
& +\sum_{q=n+1}^{2n}\left(  \sum_{\substack{k+l=q \\k\geq n+1}}\frac{\left(
-1\right)  ^{l}}{k!l!}\right)  \left(  aD\right)  ^{q}G\left(  0\right)
+\sum_{k=n+1}^{2n}\frac{1}{k!}\left(  \mathcal{R}_{2n+1}\left(  \left(
aD\right)  ^{k}G\right)  \right)  \left(  0,-a\right)  +\\
& +%
{\displaystyle\sum\limits_{q=n+1}^{2n}}
\left(
{\displaystyle\sum\limits_{\substack{l+k=q \\l,k\leq n}}}
\frac{\left(  -1\right)  ^{l}}{k!l!}\right)  \left(  aD\right)  ^{q}G\left(
0\right)  +\left(  \mathcal{R}_{2n+1}G\right)  \left(  a,-a\right)  +\left(
\mathcal{R}_{2n+1}G\right)  \left(  -a,a\right) \\
& =\sum_{q=n+1}^{2n}\left(  \sum_{\substack{k+l=q \\k\geq n+1}}\frac{\left(
-1\right)  ^{k}}{k!l!}\right)  \left(  aD\right)  ^{q}G\left(  0\right)
+\sum_{q=n+1}^{2n}\left(  \sum_{\substack{k+l=q \\k\geq n+1}}\frac{\left(
-1\right)  ^{l}}{k!l!}\right)  \left(  aD\right)  ^{q}G\left(  0\right)  +\\
& +%
{\displaystyle\sum\limits_{q=n+1}^{2n}}
\left(
{\displaystyle\sum\limits_{\substack{l+k=q \\l,k\leq n}}}
\frac{\left(  -1\right)  ^{l}}{k!l!}\right)  \left(  aD\right)  ^{q}G\left(
0\right)  +\\
& +\sum_{k=n+1}^{2n}\frac{1}{k!}\left(  \mathcal{R}_{2n+1-k}\left(  \left(
-aD\right)  ^{k}G\right)  \right)  \left(  0,a\right)  +\sum_{k=n+1}^{2n}%
\frac{1}{k!}\left(  \mathcal{R}_{2n+1-k}\left(  \left(  aD\right)
^{k}G\right)  \right)  \left(  0,-a\right)  +\\
& +\left(  \mathcal{R}_{2n+1}G\right)  \left(  a,-a\right)  +\left(
\mathcal{R}_{2n+1}G\right)  \left(  -a,a\right) \\
& =\sum_{q=n+1}^{2n}\left(  \sum_{\substack{k+l=q \\k\geq n+1}}\frac{\left(
-1\right)  ^{k}}{k!l!}+\left(  -1\right)  ^{q}\sum_{\substack{k+l=q \\k\geq
n+1}}\frac{\left(  -1\right)  ^{k}}{k!l!}+%
{\displaystyle\sum\limits_{\substack{l+k=q \\l,k\leq n}}}
\frac{\left(  -1\right)  ^{l}}{k!l!}\right)  \left(  aD\right)  ^{q}G\left(
0\right)  +\left(  \mathcal{R}_{2n+1}G\right)  \left(  a\right)  ,
\end{align*}

i.e.%
\begin{align}
& \left(  2\pi\right)  ^{-\frac{d}{2}}\int\left\vert e^{ia\xi}-\sum_{k\leq
n}\frac{\left(  ia\xi\right)  ^{k}}{k!}\right\vert ^{2}\widehat{G}\left(
\xi\right)  d\xi\nonumber\\
& =\sum_{q=n+1}^{2n}\left(  \sum_{\substack{k+l=q \\k\geq n+1}}\frac{\left(
-1\right)  ^{k}}{k!l!}+\sum_{\substack{k+l=q \\k\geq n+1}}\frac{\left(
-1\right)  ^{l}}{k!l!}+%
{\displaystyle\sum\limits_{\substack{k+l=q \\k,l\leq n}}}
\frac{\left(  -1\right)  ^{l}}{k!l!}\right)  \left(  aD\right)  ^{q}G\left(
0\right)  +\mathcal{R}_{2n+1}G\left(  a\right)  ,\label{a2.24}%
\end{align}

where $\left(  \mathcal{R}_{2n+1}G\right)  \left(  a\right)  $ is defined by
\ref{a2.32} in the statement of part 4.

Using \ref{a2.13} we get for $n+1\leq q\leq2n$,%
\begin{align*}
\sum_{\substack{k+l=q \\k\geq n+1}}\frac{\left(  -1\right)  ^{k}}{k!l!}%
+\sum_{\substack{k+l=q \\k\geq n+1}}\frac{\left(  -1\right)  ^{l}}{k!l!}+%
{\displaystyle\sum\limits_{\substack{k+l=q \\k,l\leq n}}}
\frac{\left(  -1\right)  ^{l}}{k!l!}  & =\sum_{\substack{k+l=q \\l\geq
n+1}}\frac{\left(  -1\right)  ^{l}}{k!l!}+\sum_{\substack{k+l=q \\k\geq
n+1}}\frac{\left(  -1\right)  ^{l}}{k!l!}+%
{\displaystyle\sum\limits_{\substack{k+l=q \\k,l\leq n}}}
\frac{\left(  -1\right)  ^{l}}{k!l!}\\
& =\sum_{\substack{k+l=q \\l\geq n+1}}\frac{\left(  -1\right)  ^{l}}%
{k!l!}+\sum_{\substack{k+l=q \\k\geq n+1}}\frac{\left(  -1\right)  ^{l}}%
{k!l!}+%
{\displaystyle\sum\limits_{\substack{k+l=q \\k,l\leq n}}}
\frac{\left(  -1\right)  ^{l}}{k!l!}\\
& =%
{\displaystyle\sum\limits_{k+l=q}}
\frac{\left(  -1\right)  ^{l}}{k!l!}\\
& =0,
\end{align*}

so that \ref{a2.24} becomes%
\begin{equation}
\left(  2\pi\right)  ^{-\frac{d}{2}}\int\left\vert e^{ia\xi}-\sum_{k\leq
n}\frac{\left(  ia\xi\right)  ^{k}}{k!}\right\vert ^{2}\widehat{G}\left(
\xi\right)  d\xi\leq\mathcal{R}_{2n+1}G\left(  a\right)  ,\label{a1.59}%
\end{equation}

which proves \ref{a2.30}.

Since $G\in C_{B}^{\left(  2m\right)  }$ and $\left(  aD\right)  ^{2m+1}G\in
L^{1}$, Lemma \ref{Lem_Taylor_estim_L1_Linf} implies that%
\begin{align}
&  \sum_{k=n+1}^{2n}\frac{1}{k!}\left(  \mathcal{R}_{2n+1-k}\left(  \left(
-aD\right)  ^{k}G\right)  \right)  \left(  0,a\right) \nonumber\\
&  =\sum_{k=n+1}^{2n}\frac{1}{k!}\left(  \frac{1}{\left(  2n-k\right)  !}%
\int_{0}^{1}g_{2n-k}\left(  s\right)  \left(  \left(  aD\right)
^{2n+1-k}\left(  -aD\right)  ^{k}G\right)  \left(  sa\right)  ds\right)
\nonumber\\
&  =\sum_{k=n+1}^{2n}\frac{1}{k!\left(  2n-k\right)  !}\int_{0}^{1}\left(
-1\right)  ^{k}g_{2n-k}\left(  s\right)  \left(  \left(  aD\right)
^{2n+1}G\right)  \left(  sa\right)  ds\nonumber\\
&  =\int_{0}^{1}\sum_{k=n+1}^{2n}\frac{\left(  -1\right)  ^{k}\left(
1-s\right)  ^{2n-k}}{k!\left(  2n-k\right)  !}\left(  \left(  aD\right)
^{2n+1}G\right)  \left(  sa\right)  ds\nonumber\\
&  =\int_{0}^{1}\sum_{j=0}^{n-1}\frac{\left(  -1\right)  ^{j}\left(
1-s\right)  ^{j}}{j!\left(  2n-j\right)  !}\left(  \left(  aD\right)
^{2n+1}G\right)  \left(  sa\right)  ds\nonumber\\
&  =\int_{0}^{1}\phi_{2n}\left(  s\right)  \left(  \left(  aD\right)
^{2n+1}G\right)  \left(  sa\right)  ds,\label{1.075}%
\end{align}

where%
\begin{equation}
\phi_{2n}\left(  s\right)  :=\left\{
\begin{array}
[c]{ll}%
0, & n=0,\\
\sum\limits_{j=0}^{n-1}\frac{\left(  -1\right)  ^{j}\left(  1-s\right)  ^{j}%
}{j!\left(  2n-j\right)  !}, & n=1,2,3,\ldots,
\end{array}
\right. \label{1.074}%
\end{equation}

and so $a\rightarrow-a$ yields%
\begin{align*}
\sum_{k=n+1}^{2n}\frac{1}{k!}\left(  \mathcal{R}_{2n+1-k}\left(  \left(
aD\right)  ^{k}G\right)  \right)  \left(  0,-a\right)   & =\frac{1}{\left(
2n\right)  !}\int_{0}^{1}\phi_{2n}\left(  s\right)  \left(  \left(
-aD\right)  ^{2n+1}G\right)  \left(  -sa\right)  ds\\
& =\frac{-1}{\left(  2n\right)  !}\int_{0}^{1}\phi_{2n}\left(  s\right)
\left(  \left(  aD\right)  ^{2n+1}G\right)  \left(  -sa\right)  ds\\
& =\frac{1}{\left(  2n\right)  !}\int_{0}^{1}\phi_{2n}\left(  s\right)
\left(  \left(  aD\right)  ^{2n+1}\overline{G}\right)  \left(  sa\right)  ds,
\end{align*}

which means that%
\begin{align}
\sum_{k=n+1}^{2n} &  \frac{1}{k!}\left\{  \left(  \mathcal{R}_{2n+1}\left(
\left(  -aD\right)  ^{k}G\right)  \right)  \left(  0,a\right)  +\mathcal{R}%
_{2n+1}\left(  \left(  aD\right)  ^{k}G\right)  \left(  0,-a\right)  \right\}
\nonumber\\
&  =\frac{1}{\left(  2n\right)  !}\int_{0}^{1}\phi_{2n}\left(  s\right)
\left(  \left(  aD\right)  ^{2n+1}G\right)  \left(  sa\right)  ds+\frac
{1}{\left(  2n\right)  !}\int_{0}^{1}\phi_{2n}\left(  s\right)  \left(
\left(  aD\right)  ^{2n+1}\overline{G}\right)  \left(  sa\right)
ds\nonumber\\
&  =\frac{2}{\left(  2n\right)  !}\int_{0}^{1}\phi_{2n}\left(  s\right)
\left(  \left(  aD\right)  ^{2n+1}\operatorname{Re}G\right)  \left(
sa\right)  ds,\label{a2.26}%
\end{align}

Since $G\in C_{B}^{\left(  2m\right)  }$ and $\left(  aD\right)  ^{2m+1}G\in
L^{1}$, Lemma \ref{Lem_Taylor_estim_L1_Linf} again implies that%
\begin{align*}
\left(  \mathcal{R}_{2n+1}G\right)  \left(  a,-a\right)   & =\frac{1}{\left(
2n\right)  !}\int_{0}^{1}g_{2n}\left(  s\right)  \left(  \left(  -aD\right)
^{2n+1}G\right)  \left(  a-sa\right)  ds\\
& =\frac{-1}{\left(  2n\right)  !}\int_{0}^{1}g_{2n}\left(  s\right)  \left(
\left(  aD\right)  ^{2n+1}G\right)  \left(  \left(  1-s\right)  a\right)  ds\\
& =\frac{-1}{\left(  2n\right)  !}\int_{0}^{1}\left(  1-s\right)  ^{2n}\left(
\left(  aD\right)  ^{2n+1}G\right)  \left(  \left(  1-s\right)  a\right)  ds\\
& =\frac{1}{\left(  2n\right)  !}\int_{1}^{0}t^{2n}\left(  \left(  aD\right)
^{2n+1}G\right)  \left(  ta\right)  dt\\
& =\frac{-1}{\left(  2n\right)  !}\int_{0}^{1}t^{2n}\left(  \left(  aD\right)
^{2n+1}G\right)  \left(  ta\right)  dt,
\end{align*}

and so $a\rightarrow-a$ yields%
\begin{align}
\left(  \mathcal{R}_{2n+1}G\right)  \left(  -a,a\right)   & =\frac{-1}{\left(
2n\right)  !}\int_{0}^{1}t^{2n}\left(  \left(  -aD\right)  ^{2n+1}G\right)
\left(  -ta\right)  dt\nonumber\\
& =\frac{1}{\left(  2n\right)  !}\int_{0}^{1}t^{2n}\left(  \left(  aD\right)
^{2n+1}G\right)  \left(  -ta\right)  dt\label{a2.25}\\
& =\frac{-1}{\left(  2n\right)  !}\int_{0}^{1}t^{2n}\left(  \left(
aD_{x}\right)  ^{2n+1}\left(  G\left(  -x\right)  \right)  \right)  \left(
ta\right)  dt\nonumber\\
& =\frac{-1}{\left(  2n\right)  !}\int_{0}^{1}t^{2n}\left(  \left(  aD\right)
^{2n+1}\overline{G}\right)  \left(  ta\right)  dt,\nonumber
\end{align}

and hence
\begin{equation}
\left(  \mathcal{R}_{2n+1}G\right)  \left(  a,-a\right)  +\left(
\mathcal{R}_{2n+1}G\right)  \left(  -a,a\right)  =\frac{-2}{\left(  2n\right)
!}\int_{0}^{1}t^{2n}\left(  \left(  aD\right)  ^{2n+1}\operatorname{Re}%
G\right)  \left(  ta\right)  dt.\label{a2.28}%
\end{equation}

Now \ref{a2.26}, \ref{a2.28} substituted into \ref{a2.32}\ gives%
\begin{align*}
\left(  \mathcal{R}_{2n+1}G\right)  \left(  a\right)   & =\sum_{k=n+1}%
^{2n}\frac{1}{k!}\left\{  \left(  \mathcal{R}_{2n+1-k}\left(  \left(
-aD\right)  ^{k}G\right)  \right)  \left(  0,a\right)  +\left(  \mathcal{R}%
_{2n+1-k}\left(  \left(  aD\right)  ^{k}G\right)  \right)  \left(
0,-a\right)  \right\}  +\\
& \qquad\qquad+\left(  \mathcal{R}_{2n+1}G\right)  \left(  a,-a\right)
+\left(  \mathcal{R}_{2n+1}G\right)  \left(  -a,a\right) \\
& =2\int_{0}^{1}\phi_{2n}\left(  s\right)  \left(  \left(  aD\right)
^{2n+1}\operatorname{Re}G\right)  \left(  sa\right)  ds-\frac{2}{\left(
2n\right)  !}\int_{0}^{1}s^{2n}\left(  \left(  aD\right)  ^{2n+1}%
\operatorname{Re}G\right)  \left(  sa\right)  ds\\
& =2\int_{0}^{1}\left(  \sum\limits_{j=0}^{n-1}\frac{\left(  -1\right)
^{j}\left(  1-s\right)  ^{j}}{j!\left(  2n-j\right)  !}-\frac{s^{2n}}{\left(
2n\right)  !}\right)  \left(  \left(  aD\right)  ^{2n+1}\operatorname{Re}%
G\right)  \left(  sa\right)  ds,
\end{align*}

which is \ref{a2.33}.\medskip

\textbf{Part 3} We now calculate%
\[
\int_{0}^{1}\left\vert \sum\limits_{j=0}^{n-1}\frac{\left(  -1\right)
^{j}\left(  1-s\right)  ^{j}}{j!\left(  2n-j\right)  !}-\frac{s^{2n}}{\left(
2n\right)  !}\right\vert ds,\quad n\geq1.
\]

Clearly%
\[
\sum\limits_{j=0}^{n-1}\frac{\left(  -1\right)  ^{j}\left(  1-s\right)  ^{j}%
}{j!\left(  2n-j\right)  !}=\frac{1}{\left(  2n\right)  !}\sum\limits_{j=0}%
^{n-1}\tbinom{2n}{j}\left(  -1\right)  ^{j}\left(  1-s\right)  ^{j},
\]

but%
\begin{align*}
\sum\limits_{j=0}^{n-1}\tbinom{2n}{j}\left(  -1\right)  ^{j}\left(
1-s\right)  ^{j}  & =\sum\limits_{j=0}^{2n}\tbinom{2n}{j}\left(  -1\right)
^{2n-j}\left(  1-s\right)  ^{j}-\sum\limits_{j=n}^{2n}\tbinom{2n}{j}\left(
-1\right)  ^{2n-j}\left(  1-s\right)  ^{j}\\
& =\left(  -1+1-s\right)  ^{2n}-\sum\limits_{j=n}^{2n}\tbinom{2n}{j}\left(
-1\right)  ^{2n-j}\left(  1-s\right)  ^{j}\\
& =s^{2n}-\sum\limits_{j=n}^{2n}\tbinom{2n}{j}\left(  -1\right)
^{2n-j}\left(  1-s\right)  ^{j}\\
& =s^{2n}-\sum\limits_{k=n}^{0}\tbinom{2n}{2n-k}\left(  -1\right)  ^{k}\left(
1-s\right)  ^{2n-k}\\
& =s^{2n}-\sum\limits_{k=0}^{n}\tbinom{2n}{k}\left(  -1\right)  ^{k}\left(
1-s\right)  ^{2n-k}\\
& =s^{2n}-\left(  1-s\right)  ^{n}\sum\limits_{k=0}^{n}\tbinom{2n}{k}\left(
-1\right)  ^{k}\left(  1-s\right)  ^{n-k}\\
& =s^{2n}-\left(  1-s\right)  ^{n}\sum\limits_{k=0}^{n}\tbinom{2n}{n-k}\left(
-1\right)  ^{n-k}\left(  1-s\right)  ^{k}\\
& =s^{2n}+\left(  -1\right)  ^{n+1}\left(  1-s\right)  ^{n}\sum\limits_{k=0}%
^{n}\tbinom{2n}{n-k}\left(  -1\right)  ^{k}\left(  1-s\right)  ^{k},
\end{align*}

so that%
\begin{align*}
\sum\limits_{j=0}^{n-1}\frac{\left(  -1\right)  ^{j}\left(  1-s\right)  ^{j}%
}{j!\left(  2n-j\right)  !} &  -\frac{s^{2n}}{\left(  2n\right)  !}\\
&  =\frac{1}{\left(  2n\right)  !}\sum\limits_{j=0}^{n-1}\tbinom{2n}{j}\left(
-1\right)  ^{j}\left(  1-s\right)  ^{j}-\frac{s^{2n}}{\left(  2n\right)  !}\\
&  =\frac{1}{\left(  2n\right)  !}\left(  s^{2n}+\left(  -1\right)
^{n+1}\left(  1-s\right)  ^{n}\sum\limits_{k=0}^{n}\tbinom{2n}{n-k}\left(
-1\right)  ^{k}\left(  1-s\right)  ^{k}\right)  -\frac{s^{2n}}{\left(
2n\right)  !}\\
&  =\frac{\left(  -1\right)  ^{n+1}}{\left(  2n\right)  !}\left(  1-s\right)
^{n}\sum\limits_{k=0}^{n}\tbinom{2n}{n-k}\left(  -1\right)  ^{k}\left(
1-s\right)  ^{k},
\end{align*}

and hence%
\begin{align*}
\int_{0}^{1}\left\vert \sum\limits_{j=0}^{n-1}\frac{\left(  -1\right)
^{j}\left(  1-s\right)  ^{j}}{j!\left(  2n-j\right)  !}-\frac{s^{2n}}{\left(
2n\right)  !}\right\vert ds  & =\frac{1}{\left(  2n\right)  !}\int_{0}%
^{1}\left(  1-s\right)  ^{n}\sum\limits_{k=0}^{n}\tbinom{2n}{n-k}\left(
-1\right)  ^{k}\left(  1-s\right)  ^{k}ds\\
& =\frac{1}{\left(  2n\right)  !}\int_{0}^{1}t^{n}\sum\limits_{k=0}^{n}%
\tbinom{2n}{n-k}\left(  -1\right)  ^{k}t^{k}dt.
\end{align*}

We consider two cases: $n$ odd and $n$ even.\smallskip

\fbox{\textbf{Case} $n$ is odd} Regarding the integrand:%
\begin{align*}
&  t^{n}\sum\limits_{k=0}^{n}\tbinom{2n}{n-k}\left(  -1\right)  ^{k}t^{k}\\
&  =t^{n}\left(  \tbinom{2n}{n}-\tbinom{2n}{n-1}t+\tbinom{2n}{n-2}%
t^{2}-\tbinom{2n}{n-3}t^{3}+\ldots+\tbinom{2n}{1}t^{n-1}-\tbinom{2n}{0}%
t^{n}\right) \\
&  =t^{n}\left(  \left\{  \tbinom{2n}{n}-\tbinom{2n}{n-1}t\right\}  +\left\{
\tbinom{2n}{n-2}t^{2}-\tbinom{2n}{n-3}t^{3}\right\}  +\ldots+\left\{
\tbinom{2n}{1}t^{n-1}-\tbinom{2n}{0}t^{n}\right\}  \right) \\
&  =\left\{  \tbinom{2n}{n}t^{n}-\tbinom{2n}{n-1}t^{n+1}\right\}  +\left\{
\tbinom{2n}{n-2}t^{n+2}-\tbinom{2n}{n-3}t^{n+3}\right\}  +\ldots+\left\{
\tbinom{2n}{1}t^{2n-1}-\tbinom{2n}{0}t^{2n}\right\} \\
&  \geq0,
\end{align*}

since $0\leq t\leq1$ implies this expression is the sum of $n$ non-negative
terms. Further%
\begin{align*}
&  \frac{1}{\left(  2n\right)  !}\int_{0}^{1}t^{n}\sum\limits_{k=0}^{n}%
\tbinom{2n}{n-k}\left(  -1\right)  ^{k}t^{k}ds\\
&  =\frac{1}{\left(  2n\right)  !}\left(
\begin{array}
[c]{c}%
\left\{  \tbinom{2n}{n}\frac{1}{n+1}-\tbinom{2n}{n-1}\frac{1}{n+2}\right\}
+\left\{  \tbinom{2n}{n-2}\frac{1}{n+3}-\tbinom{2n}{n-3}\frac{1}{n+4}\right\}
+\ldots\\
+\left\{  \tbinom{2n}{1}\frac{1}{2n}-\tbinom{2n}{0}\frac{1}{2n+1}\right\}
\end{array}
\right) \\
&  =\left\{  \frac{1}{n!n!}\frac{1}{n+1}-\frac{1}{\left(  n-1\right)  !\left(
n+1\right)  !}\frac{1}{n+2}\right\}  +\\
&  \qquad\qquad+\left\{  \frac{1}{\left(  n-2\right)  !\left(  n+2\right)
!}\frac{1}{n+3}-\frac{1}{\left(  n+1\right)  !\left(  n+3\right)  !}\frac
{1}{n+4}\right\}  +\ldots\\
&  \qquad\qquad+\left\{  \frac{1}{1!\left(  2n-1\right)  !}\frac{1}{2n}%
-\frac{1}{0!\left(  2n\right)  !}\frac{1}{2n+1}\right\} \\
&  =\left\{  \frac{1}{n!\left(  n+1\right)  !}-\frac{1}{\left(  n-1\right)
!\left(  n+2\right)  !}\right\}  +\left\{  \frac{1}{\left(  n-2\right)
!\left(  n+3\right)  !}-\frac{1}{\left(  n-3\right)  !\left(  n+4\right)
!}\right\}  +\ldots\\
&  \qquad\qquad+\left\{  \frac{1}{1!\left(  2n\right)  !}-\frac{1}{0!\left(
2n+1\right)  !}\right\} \\
&  =\frac{1}{\left(  2n+1\right)  !}\left(  \left\{  \tbinom{2n+1}{n}%
-\tbinom{2n+1}{n-1}\right\}  +\left\{  \tbinom{2n+1}{n-2}-\tbinom{2n+1}%
{n-3}\right\}  +\ldots+\left\{  \tbinom{2n+1}{1}-\tbinom{2n+1}{0}\right\}
\right) \\
&  =\frac{1}{\left(  2n+1\right)  !}\left(
\begin{array}
[c]{c}%
\left\{  \tbinom{2n}{n}+\tbinom{2n}{n-1}\right\}  -\left\{  \tbinom{2n}%
{n-1}+\tbinom{2n}{n-2}\right\}  +\left\{  \tbinom{2n}{n-2}+\tbinom{2n}%
{n-3}\right\}  -\ldots\\
+\left\{  \tbinom{2n}{1}+\tbinom{2n}{0}\right\}  -\tbinom{2n}{0}%
\end{array}
\right) \\
&  =\frac{1}{\left(  2n+1\right)  !}\binom{2n}{n}=\frac{1}{\left(
2n+1\right)  \left(  n!\right)  ^{2}},
\end{align*}

which means that when $n$ is odd:%
\begin{align*}
\frac{1}{\left(  2n\right)  !}\int_{0}^{1}\left\vert \sum\limits_{j=0}%
^{n-1}\frac{\left(  -1\right)  ^{j}\left(  1-s\right)  ^{j}}{j!\left(
2n-j\right)  !}-\frac{s^{2n}}{\left(  2n\right)  !}\right\vert ds  & =\frac
{1}{\left(  2n+1\right)  \left(  n!\right)  ^{2}}.\\
& .
\end{align*}

\fbox{\textbf{Case} $n$ is even}%
\begin{align*}
&  t^{n}\sum\limits_{k=0}^{n}\tbinom{2n}{n-k}\left(  -1\right)  ^{k}t^{k}\\
&  =t^{n}\left(  \tbinom{2n}{n}-\tbinom{2n}{n-1}t+\tbinom{2n}{n-2}%
t^{2}-\tbinom{2n}{n-3}t^{3}+\ldots-\tbinom{2n}{1}t^{n-1}+\tbinom{2n}{0}%
t^{n}\right) \\
&  =t^{n}\left(
\begin{array}
[c]{c}%
\left\{  \tbinom{2n}{n}-\tbinom{2n}{n-1}t\right\}  +\left\{  \tbinom{2n}%
{n-2}t^{2}-\tbinom{2n}{n-3}t^{3}\right\}  +\ldots\\
+\left\{  \tbinom{2n}{2}t^{n-2}-\tbinom{2n}{1}t^{n-1}\right\}  +\tbinom{2n}%
{0}t^{n}%
\end{array}
\right) \\
&  =\left\{  \tbinom{2n}{n}t^{n}-\tbinom{2n}{n-1}t^{n+1}\right\}  +\left\{
\tbinom{2n}{n-2}t^{n+2}-\tbinom{2n}{n-3}t^{n+3}\right\}  +\ldots\\
&  \qquad\qquad+\left\{  \tbinom{2n}{2}t^{2n-2}-\tbinom{2n}{1}t^{2n-1}%
\right\}  +\tbinom{2n}{0}t^{2n}\\
&  \geq0,
\end{align*}

since $0\leq t\leq1$ implies this expression is the sum of $n+1$ non-negative
terms. Further%
\begin{align*}
&  \frac{1}{\left(  2n\right)  !}\int_{0}^{1}t^{n}\sum\limits_{k=0}^{n}%
\tbinom{2n}{n-k}\left(  -1\right)  ^{k}t^{k}ds\\
&  =\frac{1}{\left(  2n\right)  !}\left(
\begin{array}
[c]{c}%
\left\{  \tbinom{2n}{n}\frac{1}{n+1}-\tbinom{2n}{n-1}\frac{1}{n+2}\right\}
+\left\{  \tbinom{2n}{n-2}\frac{1}{n+3}-\tbinom{2n}{n-3}\frac{1}{n+4}\right\}
+\ldots\\
+\left\{  \tbinom{2n}{2}\frac{1}{2n-1}-\tbinom{2n}{1}\frac{1}{2n}\right\}
++\tbinom{2n}{0}\frac{1}{2n+1}%
\end{array}
\right) \\
&  =\left\{  \frac{1}{n!n!}\frac{1}{n+1}-\frac{1}{\left(  n-1\right)  !\left(
n+1\right)  !}\frac{1}{n+2}\right\}  +\\
&  \qquad\qquad+\left\{  \frac{1}{\left(  n-2\right)  !\left(  n+2\right)
!}\frac{1}{n+3}-\frac{1}{\left(  n+1\right)  !\left(  n+3\right)  !}\frac
{1}{n+4}\right\}  +\ldots\\
&  \qquad\qquad+\left\{  \frac{1}{2!\left(  2n-2\right)  !}\frac{1}%
{2n-1}-\frac{1}{1!\left(  2n-1\right)  !}\frac{1}{2n}\right\}  +\frac
{1}{0!\left(  2n\right)  !}\frac{1}{2n+1}\\
&  =\left\{  \frac{1}{n!\left(  n+1\right)  !}-\frac{1}{\left(  n-1\right)
!\left(  n+2\right)  !}\right\}  +\left\{  \frac{1}{\left(  n-2\right)
!\left(  n+3\right)  !}-\frac{1}{\left(  n-3\right)  !\left(  n+4\right)
!}\right\}  +\ldots\\
&  \qquad\qquad+\left\{  \frac{1}{2!\left(  2n-1\right)  !}-\frac{1}{1!\left(
2n\right)  !}\right\}  +\frac{1}{0!\left(  2n+1\right)  !}\\
&  =\frac{1}{\left(  2n+1\right)  !}\left(  \tbinom{2n+1}{n}-\tbinom
{2n+1}{n-1}+\tbinom{2n+1}{n-2}-\tbinom{2n+1}{n-3}+\ldots+\tbinom{2n+1}%
{2}-\tbinom{2n+1}{1}+\tbinom{2n+1}{0}\right) \\
&  =\frac{1}{\left(  2n+1\right)  !}\left(
\begin{array}
[c]{c}%
\left\{  \tbinom{2n}{n}+\tbinom{2n}{n-1}\right\}  -\left\{  \tbinom{2n}%
{n-1}+\tbinom{2n}{n-2}\right\}  +\left\{  \tbinom{2n}{n-2}+\tbinom{2n}%
{n-3}\right\}  -\ldots\\
-\left\{  \tbinom{2n}{1}+\tbinom{2n}{0}\right\}  -\tbinom{2n}{0}%
\end{array}
\right) \\
&  =\frac{1}{\left(  2n+1\right)  !}\binom{2n}{n}=\frac{1}{\left(
2n+1\right)  \left(  n!\right)  ^{2}},
\end{align*}

which means that when $n$ is even:%
\[
\frac{1}{\left(  2n\right)  !}\int_{0}^{1}\left\vert \sum\limits_{j=0}%
^{n-1}\frac{\left(  -1\right)  ^{j}\left(  1-s\right)  ^{j}}{j!\left(
2n-j\right)  !}-\frac{s^{2n}}{\left(  2n\right)  !}\right\vert ds=\frac
{1}{\left(  2n+1\right)  \left(  n!\right)  ^{2}}.
\]

\end{proof}

In general a basis function is complex-valued. However, the remainder term of
the last theorem can be bounded using only the real part of the basis function:

\begin{corollary}
\label{Cor_Thm_Tayor_rem_estim_W3.2}If $H=\operatorname{Re}G$ then for all
$a\in\mathbb{R}^{d}$,%
\begin{align}
\int &  \left\vert e^{ia\xi}-%
{\displaystyle\sum\limits_{k=0}^{m}}
\frac{\left(  ia\xi\right)  ^{k}}{k!}\right\vert ^{2}\widehat{G}\left(
\xi\right)  d\xi\nonumber\\
&  =\int\left\vert e^{ia\xi}-%
{\displaystyle\sum\limits_{k=0}^{m}}
\frac{\left(  ia\xi\right)  ^{k}}{k!}\right\vert ^{2}\widehat{H}\left(
\xi\right)  d\xi\nonumber\\
&  =\left(  2\pi\right)  ^{\frac{d}{2}}\left(  2H\left(  0\right)  -2%
{\displaystyle\sum\limits_{k=0}^{m}}
\frac{\left(  -1\right)  ^{k}}{k!}\left(  aD\right)  ^{k}H\left(  a\right)  +%
{\displaystyle\sum\limits_{\substack{j,k\leq m \\j+k>m}}}
\frac{\left(  -1\right)  ^{j}}{j!k!}\left(  aD\right)  ^{j+k}H\left(
0\right)  \right)  .\label{a1.66}%
\end{align}

\end{corollary}

\begin{proof}
Regarding part 2 of Theorem \ref{Thm_Tayor_rem_estim_W3.1} set%
\[
g\left(  a\right)  =\int\left\vert e^{ia\xi}-%
{\displaystyle\sum\limits_{k=0}^{m}}
\frac{\left(  ia\xi\right)  ^{k}}{k!}\right\vert ^{2}\widehat{G}\left(
\xi\right)  d\xi,
\]

so that%
\begin{align*}
g\left(  -a\right)  =\int\left\vert e^{-ia\xi}-%
{\displaystyle\sum\limits_{k=0}^{m}}
\frac{\left(  -ia\xi\right)  ^{k}}{k!}\right\vert ^{2}\widehat{G}\left(
\xi\right)  d\xi & =\int\left\vert e^{ia\eta}-%
{\displaystyle\sum\limits_{k=0}^{m}}
\frac{\left(  ia\eta\right)  ^{k}}{k!}\right\vert ^{2}\widehat{G}\left(
-\eta\right)  d\eta\\
& =\int\left\vert e^{ia\eta}-%
{\displaystyle\sum\limits_{k=0}^{m}}
\frac{\left(  ia\eta\right)  ^{k}}{k!}\right\vert ^{2}\widehat{G\left(
-x\right)  }\left(  \eta\right)  d\eta\\
& =\int\left\vert e^{ia\eta}-%
{\displaystyle\sum\limits_{k=0}^{m}}
\frac{\left(  ia\eta\right)  ^{k}}{k!}\right\vert ^{2}\widehat{\overline{G}%
}\left(  \eta\right)  d\eta.
\end{align*}

Hence%
\begin{equation}
\int\left\vert e^{ia\xi}-%
{\displaystyle\sum\limits_{k=0}^{m}}
\frac{\left(  ia\xi\right)  ^{k}}{k!}\right\vert ^{2}\widehat
{\operatorname{Im}G}\left(  \xi\right)  d\xi=0,\label{a1.64}%
\end{equation}

and%
\begin{equation}
\int\left\vert e^{ia\xi}-%
{\displaystyle\sum\limits_{k=0}^{m}}
\frac{\left(  ia\xi\right)  ^{k}}{k!}\right\vert ^{2}\widehat{G}\left(
\xi\right)  d\xi=\int\left\vert e^{ia\xi}-%
{\displaystyle\sum\limits_{k=0}^{m}}
\frac{\left(  ia\xi\right)  ^{k}}{k!}\right\vert ^{2}\widehat
{\operatorname{Re}G}\left(  \xi\right)  d\xi.\label{a1.65}%
\end{equation}

Because $G\left(  -\cdot\right)  =\overline{G}$, we know that as
distributions,
\[
\left(  \left(  aD\right)  ^{k}G\right)  \left(  -\cdot\right)  =\left(
-1\right)  ^{k}\left(  aD\right)  ^{k}\overline{G}=\left(  -aD\right)
^{k}\overline{G},\quad k\geq0,\text{ }a\in\mathbb{R}^{d}.
\]

Also $\operatorname{Re}D^{\alpha}G$ is an even distribution when $\left\vert
\alpha\right\vert $ is even and an odd distribution when $\left\vert
\alpha\right\vert $ is odd.

Consequently%
\begin{align*}%
{\displaystyle\sum\limits_{\substack{k\leq m \\\,}}}
\frac{1}{k!} &  \left(  \left(  aD\right)  ^{k}G\right)  \left(  -a\right)  +%
{\displaystyle\sum\limits_{k=0}^{m}}
\frac{1}{k!}\left(  -aD\right)  ^{k}G\left(  a\right) \\
&  =%
{\displaystyle\sum\limits_{k=0}^{m}}
\frac{\left(  -1\right)  ^{k}}{k!}\left(  \left(  aD\right)  ^{k}\overline
{G}\right)  \left(  a\right)  +%
{\displaystyle\sum\limits_{k=0}^{m}}
\frac{\left(  -1\right)  ^{k}}{k!}\left(  aD\right)  ^{k}G\left(  a\right) \\
&  =2%
{\displaystyle\sum\limits_{k=0}^{m}}
\frac{\left(  -1\right)  ^{k}}{k!}\left(  aD\right)  ^{k}\operatorname{Re}%
G\left(  a\right)
\end{align*}

and inequalities \ref{a2.15}\ and \ref{a2.31} become%
\begin{align*}
&  \int\left\vert e^{ia\xi}-%
{\displaystyle\sum\limits_{k=0}^{m}}
\frac{\left(  ia\xi\right)  ^{k}}{k!}\right\vert ^{2}\widehat
{\operatorname{Re}G}\left(  \xi\right)  d\xi\\
&  =\left(  2\pi\right)  ^{\frac{d}{2}}\left(  2G\left(  0\right)
-2\operatorname{Re}%
{\displaystyle\sum\limits_{k=0}^{m}}
\frac{\left(  -1\right)  ^{k}}{k!}\left(  aD\right)  ^{k}G\left(  a\right)
+\right. \\
&  \qquad\qquad\left.  +%
{\displaystyle\sum\limits_{j,k\leq m;\text{ }j+k>m}}
\frac{\left(  -1\right)  ^{j}}{j!k!}\left(  \left(  aD\right)  ^{j+k}G\right)
\left(  0\right)  \right) \\
&  =\operatorname{Re}\left(  2\pi\right)  ^{\frac{d}{2}}\left(  2G\left(
0\right)  -2\operatorname{Re}%
{\displaystyle\sum\limits_{k=0}^{m}}
\frac{\left(  -1\right)  ^{k}}{k!}\left(  aD\right)  ^{k}G\left(  a\right)
+\right. \\
&  \qquad\qquad\left.  +%
{\displaystyle\sum\limits_{j,k\leq m;\text{ }j+k>m\text{ }}}
\frac{\left(  -1\right)  ^{j}}{j!k!}\left(  \left(  aD\right)  ^{j+k}G\right)
\left(  0\right)  \right) \\
&  =\left(  2\pi\right)  ^{\frac{d}{2}}\left(  2\operatorname{Re}G\left(
0\right)  -2%
{\displaystyle\sum\limits_{k=0}^{m}}
\frac{\left(  -1\right)  ^{k}}{k!}\left(  aD\right)  ^{k}\operatorname{Re}%
G\left(  a\right)  +\right. \\
&  \qquad\qquad\left.  +%
{\displaystyle\sum\limits_{j,k\leq m;\text{ }j+k>m}}
\frac{\left(  -1\right)  ^{j}}{j!k!}\left(  \left(  aD\right)  ^{j+k}%
\operatorname{Re}G\right)  \left(  0\right)  \right)  ,
\end{align*}

i.e. if $H=\operatorname{Re}G$ then%
\begin{multline*}
\int\left\vert e^{ia\xi}-%
{\displaystyle\sum\limits_{k=0}^{m}}
\frac{\left(  ia\xi\right)  ^{k}}{k!}\right\vert ^{2}\widehat{H}\left(
\xi\right)  d\xi\\
=\left(  2\pi\right)  ^{\frac{d}{2}}\left(  2H\left(  0\right)  -2%
{\displaystyle\sum\limits_{k=0}^{m}}
\frac{\left(  -1\right)  ^{k}}{k!}\left(  aD\right)  ^{k}H\left(  a\right)  +%
{\displaystyle\sum\limits_{j,k\leq m;\text{ }j+k>m}}
\frac{\left(  -1\right)  ^{j}}{j!k!}\left(  aD\right)  ^{j+k}H\left(
0\right)  \right)  .
\end{multline*}

\end{proof}

We now want to obtain some upper bounds for $\mathcal{R}_{2m+1}\left(
a\right)  $ given by \ref{a2.32}.

\begin{corollary}
\label{Cor_Thm_Tayor_rem_estim_1}Suppose the weight function $w$ has property
W2 and property W3.1 or W3.3 for order $\theta$ and parameter $\kappa$. Set
$m=\left\lfloor \underline{\kappa}\right\rfloor $. Regarding the expression
$\mathcal{R}_{2m+1}G$ defined in part 2 of Theorem
\ref{Thm_Tayor_rem_estim_W3.1}, a selection of estimates is:

\begin{enumerate}
\item
\[
\mathcal{R}_{2m+1}G\left(  a\right)  \leq\frac{1}{\left(  2m+1\right)  \left(
m!\right)  ^{2}}\max_{0\leq s\leq1}\left\vert \left(  \left(  \widehat{a}%
D_{x}\right)  ^{2m+1}\operatorname{Re}G\left(  x\right)  \right)  \left(
sa\right)  \right\vert \left\vert a\right\vert ^{2m+1}.
\]

\item
\[
\mathcal{R}_{2m+1}G\left(  a\right)  \leq\frac{1}{\left(  2m+1\right)  \left(
m!\right)  ^{2}}\max_{0\leq s\leq1}\left\vert \left(  \left(  \widehat{\cdot
}D\right)  ^{2m+1}\operatorname{Re}G\right)  \left(  sa\right)  \right\vert
\left\vert a\right\vert ^{2m+1}.
\]

\item
\[
\mathcal{R}_{2m+1}G\left(  a\right)  \leq\tbinom{2m}{m}\left(  \max_{0\leq
s\leq1}\sum_{\left\vert \beta\right\vert =2m+1}\frac{1}{\beta!}\left\vert
\left(  D^{\beta}\operatorname{Re}G\right)  \left(  sa\right)  \right\vert
\right)  \left\vert a\right\vert ^{2m+1}.
\]

\item
\[
\mathcal{R}_{2m+1}G\left(  a\right)  \leq\tbinom{2m}{m}\left(  \sum
_{\left\vert \beta\right\vert =2m+1}\frac{1}{\beta!}\max\limits_{0\leq s\leq
1}\left\vert \left(  D^{\beta}\operatorname{Re}G\right)  \left(  sa\right)
\right\vert \right)  \left\vert a\right\vert ^{2m+1}.
\]

\item
\[
\mathcal{R}_{2m+1}G\left(  a\right)  \leq\frac{d^{m+1/2}}{\left(  2m+1\right)
\left(  m!\right)  ^{2}}\left(  \max_{\left\vert \beta\right\vert =2m+1}%
\max\limits_{0\leq s\leq1}\left\vert \left(  D^{\beta}\operatorname{Re}%
G\right)  \left(  sa\right)  \right\vert \right)  \left\vert a\right\vert
^{2m+1}.
\]
\medskip

If in addition $\max\limits_{\left\vert \beta\right\vert =2m+1}\left\Vert
D^{\beta}\operatorname{Re}G\right\Vert _{\infty,B_{r}}<\infty$ then:

\item
\[
\mathcal{R}_{2m+1}G\left(  a\right)  \leq\tbinom{2m}{m}\left(  \sum
_{\left\vert \beta\right\vert =2m+1}\frac{1}{\beta!}\left\Vert D^{\beta
}\operatorname{Re}G\right\Vert _{\infty,B_{r}}\right)  \left\vert a\right\vert
^{2m+1}.
\]

\item
\[
\mathcal{R}_{2m+1}G\left(  a\right)  \leq\frac{d^{m+1/2}}{\left(  2m+1\right)
\left(  m!\right)  ^{2}}\max_{\left\vert \beta\right\vert =2m+1}\left\Vert
D^{\beta}\operatorname{Re}G\right\Vert _{\infty,B_{r}}\left\vert a\right\vert
^{2m+1},\quad\left\vert a\right\vert \leq r.
\]

\end{enumerate}
\end{corollary}

\begin{proof}
\textbf{Part 1} From \ref{a2.33} and \ref{a96},%
\begin{align*}
\mathcal{R}_{2m+1}G\left(  a\right)   & =2\int_{0}^{1}\left(  \sum
\limits_{j=0}^{m-1}\frac{\left(  -1\right)  ^{j}\left(  1-s\right)  ^{j}%
}{j!\left(  2m-j\right)  !}-\frac{s^{2m}}{\left(  2m\right)  !}\right)
\left(  \left(  aD_{x}\right)  ^{2m+1}\operatorname{Re}G\left(  x\right)
\right)  \left(  sa\right)  ds\\
& \leq2\int_{0}^{1}\left\vert \sum\limits_{j=0}^{m-1}\frac{\left(  -1\right)
^{j}\left(  1-s\right)  ^{j}}{j!\left(  2m-j\right)  !}-\frac{s^{2m}}{\left(
2m\right)  !}\right\vert \left\vert \left(  \left(  aD_{x}\right)
^{2m+1}\operatorname{Re}G\left(  x\right)  \right)  \left(  sa\right)
\right\vert ds\\
& \leq2\int_{0}^{1}\left\vert \sum\limits_{j=0}^{m-1}\frac{\left(  -1\right)
^{j}\left(  1-s\right)  ^{j}}{j!\left(  2m-j\right)  !}-\frac{s^{2m}}{\left(
2m\right)  !}\right\vert ds\text{ }\max_{0\leq s\leq1}\left\vert \left(
\left(  aD_{x}\right)  ^{2m+1}\operatorname{Re}G\left(  x\right)  \right)
\left(  sa\right)  \right\vert \\
& =2\int_{0}^{1}\left\vert \sum\limits_{j=0}^{m-1}\frac{\left(  -1\right)
^{j}\left(  1-s\right)  ^{j}}{j!\left(  2m-j\right)  !}-\frac{s^{2m}}{\left(
2m\right)  !}\right\vert ds\text{ }\max_{0\leq s\leq1}\left\vert \left(
\left(  \widehat{a}D_{x}\right)  ^{2m+1}\operatorname{Re}G\left(  x\right)
\right)  \left(  sa\right)  \right\vert \left\vert a\right\vert ^{2m+1}\\
& =\frac{1}{\left(  2m+1\right)  \left(  m!\right)  ^{2}}\max_{0\leq s\leq
1}\left\vert \left(  \left(  \widehat{a}D_{x}\right)  ^{2m+1}\operatorname{Re}%
G\left(  x\right)  \right)  \left(  sa\right)  \right\vert \left\vert
a\right\vert ^{2m+1}.
\end{align*}
\smallskip

\textbf{Part 2} From part 1,%
\begin{align*}
\max_{0\leq s\leq1}\left\vert \left(  \left(  \widehat{a}D_{x}\right)
^{2m+1}\operatorname{Re}G\left(  x\right)  \right)  \left(  sa\right)
\right\vert  & =\max_{0\leq s\leq1}\left\vert \left(  \left(  \widehat
{sa}D\right)  ^{2m+1}\operatorname{Re}G\right)  \left(  sa\right)  \right\vert
\\
& =\max_{0\leq s\leq1}\left\vert \left(  \left(  \widehat{\cdot}D\right)
^{2m+1}\operatorname{Re}G\right)  \left(  sa\right)  \right\vert .
\end{align*}
\medskip

\textbf{Part 3} Using the expansion \ref{p08} we get%
\begin{align*}
\frac{1}{\left(  2m+1\right)  !}\left\vert \left(  aD\right)  ^{2m+1}G\left(
x\right)  \right\vert  & =\frac{\left\vert a\right\vert ^{2m+1}}{\left(
2m+1\right)  !}\left\vert \left(  \widehat{a}D\right)  ^{2m+1}G\left(
x\right)  \right\vert =\left\vert a\right\vert ^{2m+1}\left\vert
\sum_{\left\vert \beta\right\vert =2m+1}\frac{\widehat{a}^{\beta}}{\beta
!}D^{\beta}G\left(  x\right)  \right\vert \leq\\
& \leq\left\vert a\right\vert ^{2m+1}\sum_{\left\vert \beta\right\vert
=2m+1}\frac{\left\vert \widehat{a}^{\beta}\right\vert }{\beta!}\left\vert
D^{\beta}G\left(  x\right)  \right\vert \leq\left\vert a\right\vert
^{2m+1}\sum_{\left\vert \beta\right\vert =2m+1}\frac{\left\vert D^{\beta
}G\left(  x\right)  \right\vert }{\beta!},
\end{align*}

so that%
\begin{align*}
\mathcal{R}_{2m+1}G\left(  a\right)   & \leq\tbinom{2m}{m}\frac{1}{\left(
2m+1\right)  !}\max_{0\leq s\leq1}\left\vert \left(  \left(  aD\right)
^{2m+1}\operatorname{Re}G\right)  \left(  sa\right)  \right\vert \\
& \leq\tbinom{2m}{m}\left(  \max_{0\leq s\leq1}\sum_{\left\vert \beta
\right\vert =2m+1}\frac{\left\vert D^{\beta}G\left(  sa\right)  \right\vert
}{\beta!}\right)  \left\vert a\right\vert ^{2m+1}.
\end{align*}

\textbf{Part 4} Clearly%
\[
\mathcal{R}_{2m+1}G\left(  a\right)  \leq\tbinom{2m}{m}\left(  \sum
_{\left\vert \beta\right\vert =2m+1}\frac{1}{\beta!}\max\limits_{0\leq s\leq
1}\left\vert \left(  D^{\beta}\operatorname{Re}G\right)  \left(  sa\right)
\right\vert \right)  \left\vert a\right\vert ^{2m+1}.
\]

\textbf{Part 5} Further,%
\begin{align*}
\sum_{\left\vert \beta\right\vert =2m+1}\frac{1}{\beta!}\max\limits_{0\leq
s\leq1}\left\vert \left(  D^{\beta}\operatorname{Re}G\right)  \left(
sa\right)  \right\vert  & \leq\left(  \sum_{\left\vert \beta\right\vert
=2m+1}\frac{1}{\beta!}\right)  \max_{\left\vert \beta\right\vert =2m+1}%
\max\limits_{0\leq s\leq1}\left\vert \left(  D^{\beta}\operatorname{Re}%
G\right)  \left(  sa\right)  \right\vert \\
& =\left(  \sum_{\left\vert \beta\right\vert =2m+1}\frac{\mathbf{1}^{2\beta}%
}{\beta!}\right)  \max_{\left\vert \beta\right\vert =2m+1}\max\limits_{0\leq
s\leq1}\left\vert \left(  D^{\beta}\operatorname{Re}G\right)  \left(
sa\right)  \right\vert \\
& =\frac{\left\vert \mathbf{1}\right\vert ^{2m+1}}{\left(  2m+1\right)  !}%
\max_{\left\vert \beta\right\vert =2m+1}\max\limits_{0\leq s\leq1}\left\vert
\left(  D^{\beta}\operatorname{Re}G\right)  \left(  sa\right)  \right\vert \\
& =\frac{d^{m+1/2}}{\left(  2m+1\right)  !}\max_{\left\vert \beta\right\vert
=2m+1}\max\limits_{0\leq s\leq1}\left\vert \left(  D^{\beta}\operatorname{Re}%
G\right)  \left(  sa\right)  \right\vert ,
\end{align*}

and hence%
\begin{align*}
\mathcal{R}_{2m+1}G\left(  a\right)   & \leq\tbinom{2m}{m}\frac{d^{m+1/2}%
}{\left(  2m+1\right)  !}\left(  \max_{\left\vert \beta\right\vert =2m+1}%
\max\limits_{0\leq s\leq1}\left\vert \left(  D^{\beta}\operatorname{Re}%
G\right)  \left(  sa\right)  \right\vert \right)  \left\vert a\right\vert
^{2m+1}\\
& =\frac{d^{m+1/2}}{\left(  2m+1\right)  \left(  m!\right)  ^{2}}\left(
\max_{\left\vert \beta\right\vert =2m+1}\max\limits_{0\leq s\leq1}\left\vert
\left(  D^{\beta}\operatorname{Re}G\right)  \left(  sa\right)  \right\vert
\right)  \left\vert a\right\vert ^{2m+1}.
\end{align*}

\textbf{Part 6} If $\left\vert a\right\vert \leq r$ then from part 3,%
\begin{align*}
\mathcal{R}_{2m+1}G\left(  a\right)   & \leq\tbinom{2m}{m}\left(
\sum_{\left\vert \beta\right\vert =2m+1}\frac{1}{\beta!}\max\limits_{0\leq
s\leq1}\left\vert \left(  D^{\beta}\operatorname{Re}G\right)  \left(
sa\right)  \right\vert \right)  \left\vert a\right\vert ^{2m+1}\\
& \leq\tbinom{2m}{m}\left(  \sum_{\left\vert \beta\right\vert =2m+1}\frac
{1}{\beta!}\max\limits_{\left\vert a\right\vert \leq r}\left\vert \left(
D^{\beta}\operatorname{Re}G\right)  \left(  a\right)  \right\vert \right)
\left\vert a\right\vert ^{2m+1}\\
& =\tbinom{2m}{m}\left(  \sum_{\left\vert \beta\right\vert =2m+1}\frac
{1}{\beta!}\left\Vert D^{\beta}\operatorname{Re}G\right\Vert _{\infty,B_{r}%
}\right)  \left\vert a\right\vert ^{2m+1}.
\end{align*}

\textbf{Part 7} If $\left\vert a\right\vert \leq r$ then
\begin{align*}
\mathcal{R}_{2m+1}G\left(  a\right)   & \leq\frac{d^{m+1/2}}{\left(
2m+1\right)  \left(  m!\right)  ^{2}}\left(  \max_{\left\vert \beta\right\vert
=2m+1}\max\limits_{0\leq s\leq1}\left\vert \left(  D^{\beta}\operatorname{Re}%
G\right)  \left(  sa\right)  \right\vert \right)  \left\vert a\right\vert
^{2m+1}\\
& \leq\frac{d^{m+1/2}}{\left(  2m+1\right)  \left(  m!\right)  ^{2}}\left(
\max_{\left\vert \beta\right\vert =2m+1}\max\limits_{\left\vert a\right\vert
\leq r}\left\vert \left(  D^{\beta}\operatorname{Re}G\right)  \left(
a\right)  \right\vert \right)  \left\vert a\right\vert ^{2m+1}\\
& =\frac{d^{m+1/2}}{\left(  2m+1\right)  \left(  m!\right)  ^{2}}%
\max_{\left\vert \beta\right\vert =2m+1}\left\Vert D^{\beta}\operatorname{Re}%
G\right\Vert _{\infty,B_{r}}\left\vert a\right\vert ^{2m+1}%
\end{align*}

\end{proof}

Suppose that the weight function has property W3.1 ?? \textbf{or W3.3} ?? for
order $\theta$ and parameter $\kappa1_{d}=\kappa\mathbf{1}$ where $\kappa
\in\mathbb{R}^{1}$. The weight function assumptions used to obtain part 2 of
Theorem \ref{Thm_Tayor_rem_estim_W3.1} imply (Theorem
\ref{Thm_basis_smth_W3.1}) that the basis function satisfies $G\in
C_{B}^{\left(  \left\lfloor 2\kappa\mathbf{1}\right\rfloor \right)  }\subset
C_{B}^{\left(  2m\mathbf{1}\right)  }$ if $m\leq\kappa<m+1/2$, and $G\in
C_{B}^{\left(  \left\lfloor 2\kappa\mathbf{1}\right\rfloor \right)  }\subset
C_{B}^{\left(  \left(  2m+1\right)  \mathbf{1}\right)  }$ if $m+1/2\leq
\kappa<m+1$. If $m+1/2\leq\kappa<m+1$ then $G\subset C_{B}^{\left(  \left(
2m+1\right)  \mathbf{1}\right)  }$ means the additional assumptions of the
next corollary are satisfied. However, ?? Corollary
\ref{Cor_Thm_Tayor_rem_estim_1} implies that $\left\vert a\right\vert $ has
exponent $\kappa$ which may be greater that $m+\frac{1}{2}$ obtained in the
next corollary.

\begin{corollary}
\label{Cor_Thm_Tayor_rem_estim_2}Suppose that the weight function has property
W3.1 ?? \textbf{or W3.3} ?? for order $\theta$ and parameter $\kappa1_{d}$
where $\kappa\in\mathbb{R}^{1}$. Set $m=\left\lfloor \kappa\right\rfloor $ and
let $\left(  \mathcal{R}_{2n+1}G\right)  \left(  a\right)  $ be the expression
defined in part 2 of Theorem \ref{Thm_Tayor_rem_estim_W3.1}. Suppose further
that $\left(  aD\right)  ^{2m+1}G\in L^{1}$ when $\left\vert a\right\vert =1$.

Then the remainder $\mathcal{R}_{m+1}f_{\rho}\left(  x,a\right)  $ of
\ref{a95} can be estimated by
\begin{equation}
\left\vert \mathcal{R}_{m+1}f_{\rho}\left(  x,a\right)  \right\vert \leq
\sqrt{\mathcal{R}_{2m+1}G\left(  a\right)  }\left\vert f_{\rho}\right\vert
_{w,\theta},\quad x\in\mathbb{R}^{d}.\label{p72}%
\end{equation}

\begin{enumerate}
\item We can directly use the estimates of Corollary
\ref{Cor_Thm_Tayor_rem_estim_1} for $\mathcal{R}_{2m+1}G\left(  a\right)  $.

\item If $\left\Vert \left(  \widehat{a}D\right)  ^{2m+1}\operatorname{Re}%
G\right\Vert _{\infty,B_{r}}<\infty$ for some $r\leq\infty$ then%
\[
\mathcal{R}_{2m+1}G\left(  a\right)  \leq\frac{1}{\left(  2m+1\right)  \left(
m!\right)  ^{2}}\left\Vert \left(  \widehat{a}D\right)  ^{2m+1}%
\operatorname{Re}G\right\Vert _{\infty,B_{r}}\left\vert a\right\vert
^{2m+1},\quad\left\vert a\right\vert \leq r.
\]

\item If $\left\Vert \max\limits_{0\leq s\leq1}\left\vert \left(  \left(
\widehat{\cdot}D\right)  ^{2m+1}\operatorname{Re}G\right)  \left(  sa\right)
\right\vert \right\Vert _{\infty,B_{r}}<\infty$ for some $r\leq\infty$, then
for $\left\vert a\right\vert \leq r$,%
\begin{align*}
\mathcal{R}_{2m+1}G\left(  a\right)   & \leq\frac{1}{\left(  2m+1\right)
\left(  m!\right)  ^{2}}\left\Vert \max\limits_{0\leq s\leq1}\left\vert
\left(  \left(  \widehat{\cdot}D\right)  ^{2m+1}\operatorname{Re}G\right)
\left(  sa\right)  \right\vert \right\Vert _{\infty,B_{r}}\left\vert
a\right\vert ^{2m+1}\\
& =\frac{1}{\left(  2m+1\right)  \left(  m!\right)  ^{2}}\max\limits_{0\leq
s\leq1}\left\Vert \left(  \left(  \widehat{\cdot}D\right)  ^{2m+1}%
\operatorname{Re}G\right)  \left(  sa\right)  \right\Vert _{\infty,B_{r}%
}\left\vert a\right\vert ^{2m+1}.
\end{align*}

\item If $\max\limits_{\left\vert \beta\right\vert =2m+1}\left\Vert D^{\beta
}G\right\Vert _{\infty,B_{r}}<\infty$ for some $0<r\leq\infty$ then from parts
5 and 6 of Corollary \ref{Cor_Thm_Tayor_rem_estim_1},%
\[
\left\vert \mathcal{R}_{m+1}f_{\rho}\left(  x,a\right)  \right\vert
\leq\left(  2\pi\right)  ^{\frac{d}{4}}\binom{2m}{m}^{\frac{1}{2}}\left(
\sum_{\left\vert \beta\right\vert =2m+1}\frac{1}{\beta!}\left\Vert D^{\beta
}G\right\Vert _{\infty,B_{r}}\right)  ^{\frac{1}{2}}\left\vert a\right\vert
^{m+\frac{1}{2}},\quad\left\vert a\right\vert \leq r,\text{ }x\in
\mathbb{R}^{d},
\]

and%
\[
\left\vert \mathcal{R}_{m+1}f_{\rho}\left(  x,a\right)  \right\vert
\leq\left(  2\pi\right)  ^{\frac{d}{4}}\frac{\sqrt{d^{m+\frac{1}{2}}}}%
{\sqrt{2m+1}m!}\left(  \max_{\left\vert \beta\right\vert =2m+1}\left\Vert
D^{\beta}G\right\Vert _{\infty,B_{r}}\right)  ^{\frac{1}{2}}\left\vert
a\right\vert ^{m+\frac{1}{2}},\quad\left\vert a\right\vert \leq r,\text{ }%
x\in\mathbb{R}^{d}.
\]

\item If $m+1/2\leq\kappa<m+1$ then $G\in C_{B}^{\left(  2m+1\right)  }$ and%
\[
\left\vert \mathcal{R}_{m+1}f_{\rho}\left(  x,a\right)  \right\vert
\leq\left(  2\pi\right)  ^{\frac{d}{4}}\binom{2m}{m}^{\frac{1}{2}}\left(
\sum_{\left\vert \beta\right\vert =2m+1}\frac{1}{\beta!}\left\Vert D^{\beta
}G\right\Vert _{\infty,B_{r}}\right)  ^{\frac{1}{2}}\left\vert a\right\vert
^{m+\frac{1}{2}},\quad a,x\in\mathbb{R}^{d},
\]

and%
\[
\left\vert \mathcal{R}_{m+1}f_{\rho}\left(  x,a\right)  \right\vert
\leq\left(  2\pi\right)  ^{\frac{d}{4}}\frac{\sqrt{d^{m+\frac{1}{2}}}}%
{\sqrt{2m+1}m!}\left(  \max_{\left\vert \beta\right\vert =2m+1}\left\Vert
D^{\beta}G\right\Vert _{\infty,B_{r}}\right)  ^{\frac{1}{2}}\left\vert
a\right\vert ^{m+\frac{1}{2}},\quad a,x\in\mathbb{R}^{d}.
\]

\end{enumerate}
\end{corollary}

\begin{proof}
Since $\left(  aD\right)  ^{2m+1}G\in L^{1}$ when $\left\vert a\right\vert
=1$, part 2 of Theorem \ref{Thm_Tayor_rem_estim_W3.1} holds and can be applied
to estimate \ref{a95} for $\mathcal{R}_{n+1}f_{\rho}\left(  x,a\right)  $ to
obtain \ref{p72}.\medskip

\textbf{Part 2} Easy ??\medskip

\textbf{Part 3} Easy ??\medskip

\textbf{Part 4} Substitute the estimates for $\mathcal{R}_{2m+1}G$ derived in
parts 5 and 6 into Corollary \ref{Cor_Thm_Tayor_rem_estim_1} into the estimate
\ref{p72}.\medskip

\textbf{Part 5} From Theorem \ref{Thm_basis_smth_W3.1}, $G\in C_{B}^{\left(
\left\lfloor 2\kappa\right\rfloor \right)  }\subset C_{B}^{\left(
2m+1\right)  }$ so the additional assumptions of part 2 is satisfied for
$r=\infty$.
\end{proof}

\subsection{Example: the Gaussian}

?? The Gaussian is a tensor product radial function so part 3 of Corollary
\ref{Cor_Thm_Tayor_rem_estim_1} seems nice - see part 1 of Lemma
\ref{Lem_deriv_rad_funcs}.

\subsection{Example: tensor product B-spline weight functions in W3.1}

We will use part 2 of Corollary \ref{Cor_Thm_Tayor_rem_estim_1} which assumes
$\max\limits_{s\in\left[  0,1\right]  }\left\vert \left(  \left(
\widehat{\cdot}D\right)  ^{2m+1}G_{\theta}\right)  \left(  sa\right)
\right\vert <\infty$. We use this condition because the differential operators
$\left(  \widehat{\cdot}D\right)  ^{m}$ acts nicely on the radial functions
$\left\vert x\right\vert ^{n}$ and $\left\vert x\right\vert ^{n}\log\left\vert
x\right\vert $ as shown in part 1 of Lemma \ref{Lem_deriv_rad_funcs} and the
upper bounds are not based on results which are conjectures such as those of
Lemma \ref{Lem_op_aD_estim}.

Note however that the last corollary relies on part 2 of Theorem
\ref{Thm_Tayor_rem_estim_W3.1} which assumes $\left(  aD\right)
^{2n-1}G_{\theta}\in L^{1}$ when $a\in\mathbb{R}^{d}$, \textbf{so we first
prove}:

\begin{theorem}
\label{Thm_xD_W3.1_L1}Suppose $w$ is an \textbf{extended B-spline} weight
function \ref{p32} with parameters $n$ and $l$ and which has property W3.1 for
order $\theta$ and smoothness $\kappa$. Suppose $G_{\theta}=\left(  \frac
{1}{w\left\vert \cdot\right\vert ^{2\theta}}\right)  ^{\vee}$ is an order
$\theta$ basis function and $G_{0}=\left(  \frac{1}{w}\right)  ^{\vee}$.

Then if $n>\theta$ we have $\left(  aD\right)  ^{2n-1}G_{\theta}=\left\vert
a\right\vert ^{2n-1}\left(  \widehat{a}D\right)  ^{2n-1}G_{\theta}\in L^{1}$
for each $a\in\mathbb{R}^{d}$ where $\widehat{a}=a/\left\vert a\right\vert $.
\end{theorem}

\begin{proof}
We need a copy of equations \ref{p23}:%
\[
G_{\theta}=\left\{
\begin{array}
[c]{ll}%
\frac{1}{e\left(  \theta-d/2\right)  }G_{0}\ast T_{\theta}, & 2\theta>d,\\
\frac{\left(  -1\right)  ^{\frac{d+1}{2}-\theta}}{e\left(  1/2\right)
}\left(  \left\vert D\right\vert ^{2}\right)  ^{\frac{d+1}{2}-\theta}\left(
G_{0}\ast T_{\frac{d+1}{2}}\right)  , & 2\theta<d,\text{ }d\text{ }odd,\text{
}d\geq3,\\
\frac{\left(  -1\right)  ^{\frac{d}{2}+1-\theta}}{e\left(  1\right)  }\left(
\left\vert D\right\vert ^{2}\right)  ^{\frac{d}{2}+1-\theta}\left(  G_{0}\ast
T_{\frac{d}{2}+1}\right)  , & 2\theta\leq d,\text{ }d\text{ }even.
\end{array}
\right.
\]

We use a radial $C^{\infty}$ partition of unity $\left\{  \phi_{0}%
,\phi_{\infty}=1-\phi_{0}\right\}  $ where $0\leq\phi_{0},\phi_{\infty}\leq1
$, $\operatorname*{supp}\phi_{0}=\overline{B}_{1}$, $\phi_{0}=1$ on
$\overline{B}_{1/2}$. We consider four cases:\medskip

\fbox{\textbf{Case} $2\theta>d$} Here%
\begin{align*}
\left(  \widehat{a}D\right)  ^{2n-1}\left(  G_{0}\ast T_{\theta}\right)   &
=\left(  \widehat{a}D\right)  ^{2n-1}\left(  G_{0}\ast\phi_{0}T_{\theta
}\right)  +\left(  \widehat{a}D\right)  ^{2n-1}\left(  G_{0}\ast\phi_{\infty
}T_{\theta}\right) \\
& =\left(  \left(  \widehat{a}D\right)  ^{2n-1}G_{0}\right)  \ast\phi
_{0}T_{\theta}+\left(  \left(  \widehat{a}D\right)  ^{2n-1-\left(
2\theta+1\right)  }G_{0}\right)  \ast\left(  \widehat{a}D\right)  ^{2\theta
+1}\left(  \phi_{\infty}T_{\theta}\right) \\
& =\left(  \left(  \widehat{a}D\right)  ^{2n-1}G_{0}\right)  \ast\phi
_{0}T_{\theta}+\left(  \left(  \widehat{a}D\right)  ^{2\left(  n-\theta
-1\right)  }G_{0}\right)  \ast\left(  \widehat{a}D\right)  ^{2\theta+1}\left(
\phi_{\infty}T_{\theta}\right)  .
\end{align*}

Now if $\left\vert x\right\vert \geq1$ then by parts 3 and 4 of Lemma
\ref{Lem_thin_plate_splin},
\begin{align*}
\left\vert \left(  \widehat{a}D\right)  ^{2\theta+1}\left(  \phi_{\infty
}T_{\theta}\right)  \left(  x\right)  \right\vert =\left\vert \left(
\widehat{a}D\right)  ^{2\theta+1}T_{\theta}\left(  x\right)  \right\vert  &
\leq k_{2\theta+1,2\theta-d}\left\vert x\right\vert ^{2\theta-d-\left(
2\theta+1\right)  }\\
& =k_{2\theta+1,2\theta-d}\left\vert x\right\vert ^{-d-1},
\end{align*}

and this implies that $\left(  \widehat{a}D\right)  ^{2\theta+1}\left(
\phi_{\infty}T_{\theta}\right)  \in L^{1}$.

Also $\phi_{0}T_{\theta}=const\phi_{0}\left\vert \cdot\right\vert ^{2\theta
-d}\leq const\phi_{0}\in L^{1}$.

From Lemma \ref{Lem_basis_ord0_Bsplin_W3}, $\left(  \widehat{a}D\right)
^{k}G_{0}\in L^{1}$ when $k\leq2n-1$.

Thus $\phi_{0}T_{\theta},$ $\left(  \widehat{a}D\right)  ^{2\theta+1}\left(
\phi_{\infty}T_{\theta}\right)  ,$ $\left(  \widehat{a}D\right)  ^{2\left(
n-\theta-1\right)  }G_{0},$ $\left(  \widehat{a}D\right)  ^{2n-1}G_{0}$ all
lie in $L^{1}$ and Young's convolution estimate with $p=q=r=1$ now implies
\begin{align*}
&  \left\Vert \left(  \widehat{a}D\right)  ^{2n-1}\left(  G_{0}\ast T_{\theta
}\right)  \right\Vert _{1}\\
&  \leq\left\Vert \left(  \left(  \widehat{a}D\right)  ^{2n-1}G_{0}\right)
\right\Vert _{1}\left\Vert \phi_{0}T_{\theta}\right\Vert _{1}+\left\Vert
\left(  \widehat{a}D\right)  ^{2\left(  n-\theta-1\right)  }G_{0}\right\Vert
_{1}\left\Vert \left(  \widehat{a}D\right)  ^{2\theta+1}\left(  \phi_{\infty
}T_{\theta}\right)  \right\Vert _{1}\\
&  <\infty.
\end{align*}
\medskip

\fbox{\textbf{Case} $2\theta<d$ and $d$ odd} From Theorem
\ref{Thm_ExtNatSplin_convol}, $T_{\frac{d+1}{2}}=-\left\vert \cdot\right\vert
$ and if we let%
\[
k_{d,\theta}=\frac{\left(  -1\right)  ^{\frac{d+1}{2}-\theta}}{e\left(
1/2\right)  },
\]

then%
\begin{align}
&  \left(  \widehat{a}D\right)  ^{2n-1}G_{\theta}\nonumber\\
&  =k_{d,\theta}\left(  \widehat{a}D\right)  ^{2n-1}\left(  \left\vert
D\right\vert ^{2}\right)  ^{\frac{d+1}{2}-\theta}\left(  G_{0}\ast
T_{\frac{d+1}{2}}\right) \nonumber\\
&  =-k_{d,\theta}\left(  \widehat{a}D\right)  ^{2n-1}\left(  \left\vert
D\right\vert ^{2}\right)  ^{\frac{d+1}{2}-\theta}\left(  G_{0}\ast\left\vert
\cdot\right\vert \right) \nonumber\\
&  =-k_{d,\theta}\left(  \left\vert D\right\vert ^{2}\right)  ^{\frac{d+1}%
{2}-\theta}\left(  \widehat{a}D\right)  ^{2n-1}\left(  G_{0}\ast\left\vert
\cdot\right\vert \right) \nonumber\\
&  =-k_{d,\theta}\left(  \left\vert D\right\vert ^{2}\right)  ^{\frac{d+1}%
{2}-\theta}\left(  \widehat{a}D\right)  ^{2n-1}\left(  G_{0}\ast\phi
_{0}\left\vert \cdot\right\vert \right)  -k_{d,\theta}\left(  \left\vert
D\right\vert ^{2}\right)  ^{\frac{d+1}{2}-\theta}\left(  \widehat{a}D\right)
^{2n-1}\left(  G_{0}\ast\phi_{\infty}\left\vert \cdot\right\vert \right)
\nonumber\\
&  =-k_{d,\theta}\left(  \left(  \widehat{a}D\right)  ^{2n-1}G_{0}\right)
\ast\left(  \left\vert D\right\vert ^{2}\right)  ^{\frac{d+1}{2}-\theta
}\left(  \phi_{0}\left\vert \cdot\right\vert \right)  -k_{d,\theta}\left(
\widehat{a}D\right)  ^{2n-1}\left(  \left\vert D\right\vert ^{2}\right)
^{\frac{d+1}{2}-\theta}\left(  G_{0}\ast\phi_{\infty}\left\vert \cdot
\right\vert \right)  .\label{a1.10}%
\end{align}

Part 9 of Lemma \ref{Lem_iterLaplacian_rad} with $k=1$ implies%
\begin{align*}
\left\vert D\right\vert ^{2n}\left(  \left\vert x\right\vert \right)   &
=2^{2n}\left(  -k+\frac{1}{2}\right)  _{n}\left(  -k-\frac{d}{2}+\frac{3}%
{2}\right)  _{n}\left\vert x\right\vert ^{2k-1-2n}\\
& =2^{2n}\left(  -\frac{1}{2}\right)  _{n}\left(  -\frac{d-1}{2}\right)
_{n}\left\vert x\right\vert ^{1-2n}.
\end{align*}

Regarding the first term in \ref{a1.10}: near the origin
\[
\left(  \left\vert D\right\vert ^{2}\right)  ^{\frac{d+1}{2}-\theta}\left(
\phi_{0}\left\vert x\right\vert \right)  =\left(  \left\vert D\right\vert
^{2}\right)  ^{\frac{d+1}{2}-\theta}\left\vert x\right\vert =2^{2n}\left(
-\frac{1}{2}\right)  _{n}\left(  -\frac{d-1}{2}\right)  _{n}\frac
{1}{\left\vert x\right\vert ^{d-2\theta}}\in L_{loc}^{1},
\]

which implies $\left(  \left\vert D\right\vert ^{2}\right)  ^{\frac{d+1}%
{2}-\theta}\left(  \phi_{0}\left\vert x\right\vert \right)  \in L^{1}$. Since
$\left(  \widehat{a}D\right)  ^{2n-1}G_{0}$ is bounded with compact support
(Lemma \ref{Lem_basis_ord0_Bsplin_W3}) it follows that $\left(  \left(
\widehat{a}D\right)  ^{2n-1}G_{0}\right)  \ast\left(  \left\vert D\right\vert
^{2}\right)  ^{\frac{d+1}{2}-\theta}\left(  \phi_{0}\left\vert \cdot
\right\vert \right)  \in L^{1}$.

Regarding the second term of \ref{a1.10}: since%
\begin{align*}
\left(  \widehat{a}D\right)  ^{2n-1}\left(  \left\vert D\right\vert
^{2}\right)  ^{\frac{d+1}{2}-\theta}u  & =\left(  \frac{d+1}{2}-\theta\right)
!\left(  \widehat{a}D\right)  ^{2n-1}\sum\limits_{\left\vert \alpha\right\vert
=\frac{d+1}{2}-\theta}\frac{1}{\alpha!}D^{2\alpha}u\\
& =\left(  \frac{d+1}{2}-\theta\right)  !\left(  2n-1\right)  !\sum
\limits_{\left\vert \beta\right\vert =2n-1}\frac{1}{\beta!}\widehat{a}^{\beta
}D^{\beta}\sum\limits_{\left\vert \alpha\right\vert =\frac{d+1}{2}-\theta
}\frac{1}{\alpha!}D^{2\alpha}u\\
& =\left(  \frac{d+1}{2}-\theta\right)  !\left(  2n-1\right)  !\sum
\limits_{\left\vert \beta\right\vert =2n-1}\sum\limits_{\left\vert
\alpha\right\vert =\frac{d+1}{2}-\theta}\frac{\widehat{a}^{\beta}}%
{\alpha!\beta!}D^{\beta+2\alpha}u,
\end{align*}

we have%
\begin{align}
\left(  \widehat{a}D\right)  ^{2n-1}\left(  \left\vert D\right\vert
^{2}\right)   &  ^{\frac{d+1}{2}-\theta}\left(  G_{0}\ast\phi_{\infty
}\left\vert \cdot\right\vert \right) \nonumber\\
&  =\left(  \frac{d+1}{2}-\theta\right)  !\left(  2n-1\right)  !\sum
\limits_{\left\vert \beta\right\vert =2n-1}\sum\limits_{\left\vert
\alpha\right\vert =\frac{d+1}{2}-\theta}\frac{\widehat{a}^{\beta}}%
{\alpha!\beta!}D^{\beta+2\alpha}\left(  G_{0}\ast\phi_{\infty}\left\vert
\cdot\right\vert \right)  .\label{a1.12}%
\end{align}

But $n>\theta\Rightarrow n\geq\theta+1\Rightarrow2n\geq2\theta+2$ so that%
\[
\left\vert \beta+2\alpha\right\vert =\left(  2n-1\right)  +\left(
d+1-2\theta\right)  =2n+d-2\theta\geq d+2,
\]

and there exists $\gamma\leq2\alpha+\beta$ such that $\left\vert
\gamma\right\vert =d+2$. Now write%
\[
D^{\beta+2\alpha}\left(  G_{0}\ast\left(  \phi_{\infty}\left\vert
\cdot\right\vert \right)  \right)  =\left(  D^{2\alpha+\beta-\gamma}%
G_{0}\right)  \ast D^{\gamma}\left(  \phi_{\infty}\left\vert \cdot\right\vert
\right)  .
\]

Since $\left\vert \beta+2\alpha-\gamma\right\vert =2n+d-2\theta-d-2=2n-2\theta
-2$ and $0\leq\left\vert \beta+2\alpha-\gamma\right\vert \leq2n-4$, Lemma
\ref{Lem_basis_ord0_Bsplin_W3} implies $D^{\beta+2\alpha-\gamma}G_{0}\in
L^{1}$. Further, when $\left\vert x\right\vert \geq1$, by part 1 of Lemma
\ref{Lem_bnd_deriv_rad_func}, when $\left\vert x\right\vert \geq1$,
\[
\left\vert D^{\gamma}\left(  \phi_{\infty}\left\vert x\right\vert \right)
\right\vert =\left\vert D^{\gamma}\left\vert x\right\vert \right\vert \leq
c_{\gamma,1}\left\vert x\right\vert ^{1-\left\vert \gamma\right\vert
}=c_{\gamma,1}\left\vert x\right\vert ^{1-\left(  d+2\right)  }=c_{\gamma
,1}\left\vert x\right\vert ^{-\left(  d+1\right)  },
\]

which means that $\left\vert D^{\gamma}\left(  \phi_{\infty}\left\vert
\cdot\right\vert \right)  \right\vert \in L^{1}$.

We can now conclude using Young's convolution estimate that $D^{\beta+2\alpha
}\left(  G_{0}\ast\phi_{\infty}\left\vert \cdot\right\vert \right)  \in L^{1}$
and consequently from \ref{a1.12} that $\left(  \widehat{a}D\right)
^{2n-1}\left(  \left\vert D\right\vert ^{2}\right)  ^{\frac{d+1}{2}-\theta
}\left(  G_{0}\ast\phi_{\infty}\left\vert \cdot\right\vert \right)  \in L^{1}%
$.\medskip

\fbox{\textbf{Case} $2\theta<d$ and $d$ even} From the statement of Theorem
\ref{Thm_ExtNatSplin_convol}, $T_{\frac{d}{2}+1}=\left\vert \cdot\right\vert
^{2}\log\left\vert \cdot\right\vert $.%
\begin{align}
&  \left(  \widehat{a}D\right)  ^{2n-1}\left(  \left\vert D\right\vert
^{2}\right)  ^{\frac{d}{2}+1-\theta}G_{\theta}\nonumber\\
&  =\left(  \widehat{a}D\right)  ^{2n-1}\left(  \left\vert D\right\vert
^{2}\right)  ^{\frac{d}{2}+1-\theta}\left(  G_{0}\ast T_{\frac{d}{2}+1}\right)
\nonumber\\
&  =\left(  \widehat{a}D\right)  ^{2n-1}\left(  \left\vert D\right\vert
^{2}\right)  ^{\frac{d}{2}+1-\theta}\left(  G_{0}\ast\left(  \left\vert
\cdot\right\vert ^{2}\log\left\vert \cdot\right\vert \right)  \right)
\nonumber\\
&  =\left(  \widehat{a}D\right)  ^{2n-1}\left(  \left\vert D\right\vert
^{2}\right)  ^{\frac{d}{2}+1-\theta}\left(  G_{0}\ast\left(  \phi
_{0}\left\vert \cdot\right\vert ^{2}\log\left\vert \cdot\right\vert \right)
\right)  +\left(  \widehat{a}D\right)  ^{2n-1}\left(  \left\vert D\right\vert
^{2}\right)  ^{\frac{d}{2}+1-\theta}\left(  G_{0}\ast\left(  \phi_{\infty
}\left\vert \cdot\right\vert ^{2}\log\left\vert \cdot\right\vert \right)
\right) \nonumber\\
&  =\left(  \left(  \widehat{a}D\right)  ^{2n-1}G_{0}\right)  \ast\left(
\left\vert D\right\vert ^{2}\right)  ^{\frac{d}{2}+1-\theta}\left(  \phi
_{0}\left\vert \cdot\right\vert ^{2}\log\left\vert \cdot\right\vert \right)
+\nonumber\\
&  \qquad\qquad\qquad\qquad+\left(  \widehat{a}D\right)  ^{2n-1}\left(
\left\vert D\right\vert ^{2}\right)  ^{\frac{d}{2}+1-\theta}\left(  G_{0}%
\ast\left(  \phi_{\infty}\left\vert \cdot\right\vert ^{2}\log\left\vert
\cdot\right\vert \right)  \right)  .\label{a1.20}%
\end{align}

We first consider the \textbf{first convolution}. Since $\frac{d}{2}%
-\theta+1\geq2$ we will use part 7 of Lemma \ref{Lem_iterLaplacian_rad} i.e.%
\[
\left\vert D\right\vert ^{2n}\left(  \left\vert x\right\vert ^{2}%
\log\left\vert x\right\vert \right)  =2^{2n-1}\left(  n-2\right)  !\left(
-\frac{d}{2}\right)  _{n}\frac{1}{\left\vert x\right\vert ^{2n-2}},\quad
n\geq2.
\]

Thus when $n=\frac{d}{2}+1-\theta$ and $\left\vert x\right\vert \leq1/2$,%
\[
\left(  \left\vert D\right\vert ^{2}\right)  ^{\frac{d}{2}+1-\theta}\left(
\phi_{0}\left\vert \cdot\right\vert ^{2}\log\left\vert \cdot\right\vert
\right)  =\left\vert D\right\vert ^{2n}\left(  \left\vert \cdot\right\vert
^{2}\log\left\vert \cdot\right\vert \right)  =2^{2n-1}\left(  n-2\right)
!\left(  -\frac{d}{2}\right)  _{n}\frac{1}{\left\vert x\right\vert
^{d-2\theta}}\in L_{loc}^{1},
\]

which implies $\left(  \left\vert D\right\vert ^{2}\right)  ^{\frac{d+1}%
{2}-\theta}\left(  \phi_{0}\left\vert \cdot\right\vert ^{2}\log\left\vert
\cdot\right\vert \right)  \in L^{1}$. Since $\left(  \widehat{a}D\right)
^{2n-1}G_{0}$ is bounded with compact support (Lemma
\ref{Lem_basis_ord0_Bsplin_W3}) it follows that $\left(  \left(  \widehat
{a}D\right)  ^{2n-1}G_{0}\right)  \ast\left(  \left\vert D\right\vert
^{2}\right)  ^{\frac{d+1}{2}-\theta}\left(  \phi_{0}\left\vert \cdot
\right\vert \right)  \in L^{1}$.

Regarding the \textbf{second convolution} in \ref{a1.20} we again follow the
previous case to obtain%
\begin{align}
& \left(  \widehat{a}D\right)  ^{2n-1}\left(  \left\vert D\right\vert
^{2}\right)  ^{\frac{d}{2}+1-\theta}\left(  G_{0}\ast\phi_{\infty}\left\vert
\cdot\right\vert ^{2}\log\left\vert \cdot\right\vert \right) \nonumber\\
& =\left(  \frac{d}{2}+1-\theta\right)  !\left(  2n-1\right)  !\sum
\limits_{\left\vert \beta\right\vert =2n-1}\sum\limits_{\left\vert
\alpha\right\vert =\frac{d}{2}+1-\theta}\frac{\widehat{a}^{\beta}}%
{\alpha!\beta!}D^{\beta+2\alpha}\left(  G_{0}\ast\phi_{\infty}\left\vert
\cdot\right\vert ^{2}\log\left\vert \cdot\right\vert \right)  .\label{a1.21}%
\end{align}

But $n>\theta\Rightarrow n\geq\theta+1\Rightarrow2n\geq2\theta+2$ so that%
\[
\left\vert \beta+2\alpha\right\vert =2n-1+d+2-2\theta=2n+1+d-2\theta\geq d+3,
\]

and there exists $\gamma\leq\beta+2\alpha$ such that $\left\vert
\gamma\right\vert =d+3$. Now write%
\[
D^{\beta+2\alpha}\left(  G_{0}\ast\left(  \phi_{\infty}\left\vert
\cdot\right\vert ^{2}\log\left\vert \cdot\right\vert \right)  \right)
=\left(  D^{\beta+2\alpha-\gamma}G_{0}\right)  \ast D^{\gamma}\left(
\phi_{\infty}\left\vert \cdot\right\vert ^{2}\log\left\vert \cdot\right\vert
\right)  .
\]

Since $\left\vert \left(  \beta+2\alpha\right)  -\gamma\right\vert =\left(
2n+1+d-2\theta\right)  -\left(  d+3\right)  =2n-2\theta-2$ and $0\leq
\left\vert \beta+2\alpha-\gamma\right\vert \leq2n-4$, Lemma
\ref{Lem_basis_ord0_Bsplin_W3} implies $D^{\beta+2\alpha-\gamma}G_{0}\in L^{1}
$. Noting that $\left\vert \gamma\right\vert >2$, by part 2 of Lemma
\ref{Lem_bnd_deriv_rad_func} there exist constants $c_{\gamma,2}^{\prime}$
such that%
\[
\left\vert D^{\gamma}\left(  \left\vert x\right\vert ^{2}\log\left\vert
x\right\vert \right)  \right\vert \leq c_{\gamma,2}^{\prime}\left\vert
x\right\vert ^{2-\left\vert \gamma\right\vert }=c_{\gamma,2}^{\prime
}\left\vert x\right\vert ^{2-\left(  d+3\right)  }=c_{\gamma,2}^{\prime
}\left\vert x\right\vert ^{-\left(  d+1\right)  }.
\]

Thus, when $\left\vert x\right\vert \geq1$,
\[
\left\vert D^{\gamma}\left(  \phi_{\infty}\left\vert x\right\vert ^{2}%
\log\left\vert x\right\vert \right)  \right\vert =\left\vert D^{\gamma}\left(
\left\vert x\right\vert ^{2}\log\left\vert x\right\vert \right)  \right\vert
\leq c_{\gamma,2}^{\prime}\left\vert x\right\vert ^{-\left(  d+1\right)  },
\]

which means that $\left\vert D^{\gamma}\left(  \phi_{\infty}\left\vert
\cdot\right\vert ^{2}\log\left\vert \cdot\right\vert \right)  \right\vert \in
L^{1}$.

We can now conclude using Young's convolution estimate that $D^{\beta+2\alpha
}\left(  G_{0}\ast\phi_{\infty}\left\vert \cdot\right\vert \right)  \in L^{1}$
and consequently from \ref{a1.12} that $\left(  \widehat{a}D\right)
^{2n-1}\left(  \left\vert D\right\vert ^{2}\right)  ^{\frac{d}{2}+1-\theta
}\left(  G_{0}\ast\phi_{\infty}\left\vert \cdot\right\vert \right)  \in L^{1}%
$.\medskip

\fbox{\textbf{Case} $2\theta=d$}
\begin{align*}
\left(  \widehat{a}D\right)  ^{2n-1}\left(  \left\vert D\right\vert
^{2}\right)  ^{\frac{d}{2}+1-\theta}G_{\theta}  & =\left(  \widehat
{a}D\right)  ^{2n-1}\left\vert D\right\vert ^{2}G_{\theta}\\
& =\left(  \widehat{a}D\right)  ^{2n-1}\left(  G_{0}\ast\left\vert
D\right\vert ^{2}\left(  \left\vert \cdot\right\vert ^{2}\log\left\vert
\cdot\right\vert \right)  \right) \\
& =\left(  \widehat{a}D\right)  ^{2n-1}\left(  G_{0}\ast\left(  1+2d\log
\left\vert \cdot\right\vert \right)  \right) \\
& =2d\left(  \widehat{a}D\right)  ^{2n-1}\left(  G_{0}\ast\log\left\vert
\cdot\right\vert \right) \\
& =2d\left(  \widehat{a}D\right)  ^{2n-1}\left(  G_{0}\ast\phi_{0}%
\log\left\vert \cdot\right\vert \right)  +2d\left(  \widehat{a}D\right)
^{2n-1}\left(  G_{0}\ast\phi_{\infty}\log\left\vert \cdot\right\vert \right)
\\
& =2d\left(  \left(  \widehat{a}D\right)  ^{2n-1}G_{0}\right)  \ast\left(
\phi_{0}\log\left\vert \cdot\right\vert \right)  +2dG_{0}\ast\left(
\widehat{a}D\right)  ^{2n-1}\left(  \phi_{\infty}\log\left\vert \cdot
\right\vert \right)  ,
\end{align*}

so that%
\begin{align*}
&  \left\Vert \left(  \widehat{a}D\right)  ^{2n-1}\left(  \left\vert
D\right\vert ^{2}\right)  ^{\frac{d}{2}+1-\theta}G_{\theta}\right\Vert _{1}\\
&  \leq2d\left\Vert \left(  \widehat{a}D\right)  ^{2n-1}G_{0}\right\Vert
_{1}\left\Vert \phi_{0}\log\left\vert \cdot\right\vert \right\Vert
_{1}+2d\left\Vert G_{0}\right\Vert _{1}\left\Vert \left(  \widehat{a}D\right)
^{2n-1}\left(  \phi_{\infty}\log\left\vert \cdot\right\vert \right)
\right\Vert _{1}\\
&  =2d\left\Vert \left(  \widehat{a}D\right)  ^{2n-1}G_{0}\right\Vert _{1}%
\int_{\left\vert \cdot\right\vert \leq1}\left\vert \phi_{0}\log\left\vert
\cdot\right\vert \right\vert +2d\left\Vert G_{0}\right\Vert _{1}%
\int_{\left\vert \cdot\right\vert \geq\frac{1}{2}}\left\vert \left(
\widehat{a}D\right)  ^{2n-1}\left(  \phi_{\infty}\log\left\vert \cdot
\right\vert \right)  \right\vert \\
&  \leq2d\left\Vert \left(  \widehat{a}D\right)  ^{2n-1}G_{0}\right\Vert
_{1}\int_{\left\vert \cdot\right\vert \leq1}\left\vert \log\left\vert
\cdot\right\vert \right\vert +2d\left\Vert G_{0}\right\Vert _{1}%
\int_{\left\vert \cdot\right\vert \geq\frac{1}{2}}\left\vert \left(
\widehat{a}D\right)  ^{2n-1}\left(  \phi_{\infty}\log\left\vert \cdot
\right\vert \right)  \right\vert .
\end{align*}

But by part 1 of Corollary \ref{Cor_Thm_Integ_u(xy)f(|x|)dx} of the Appendix,
changing to spherical coordinates yields%
\[
\int_{\left\vert \cdot\right\vert \leq1}\left\vert \log\left\vert
\cdot\right\vert \right\vert =-\frac{2\pi^{d/2}}{\Gamma\left(  d/2\right)
}\int_{0}^{1}\rho^{d-1}\log\rho d\rho=-\frac{2\pi^{d/2}}{\Gamma\left(
d/2\right)  }\int_{0}^{1}\rho^{d-1}\log\rho d\rho=\frac{2\pi^{d/2}}%
{\Gamma\left(  d/2\right)  }\frac{1}{d^{2}},
\]

and by part 9 of \ Lemma \ref{Lem_op_aD_estim},
\[
\left\vert \left(  \widehat{a}D\right)  ^{2n-1}\log\left\vert x\right\vert
\right\vert \leq2\left(  2n-1\right)  !\left\vert x\right\vert ^{-\left(
2n-1\right)  },
\]

so the inequality $2n-1\geq2\theta+1=d+1$ implies that $\int_{\left\vert
\cdot\right\vert \geq1}\left\vert \left(  \widehat{a}D\right)  ^{2n-1}%
\log\left\vert \cdot\right\vert \right\vert <\infty$ and hence that

$\int_{\left\vert \cdot\right\vert \geq\frac{1}{2}}\left\vert \left(
\widehat{a}D\right)  ^{2n-1}\left(  \phi_{\infty}\log\left\vert \cdot
\right\vert \right)  \right\vert <\infty$.
\end{proof}

?? MORE\ BLURB? ?? In the case of the extended B-spline, part 1 of Corollary
\ref{Cor_Thm_Tayor_rem_estim_1} will yield a significantly better order for
the Taylor series remainder than the estimates of Approach 1 in Subsubsection
\ref{SbSect_TaylorDataApproach1} i.e. the estimate \ref{a001}.

We simply apply part 1 of Corollary \ref{Cor_Thm_Tayor_rem_estim_1}.

\begin{theorem}
\label{Thm_xhatD_G3.1*_bded}Suppose $w$ is an \textbf{extended B-spline}
weight function \ref{p32} with parameters $n$ and $l$ and which has property
W3.1 for order $\theta$ and smoothness $\kappa$. Suppose $G_{\theta}=\left(
\frac{1}{w\left\vert \cdot\right\vert ^{2\theta}}\right)  ^{\vee}$ is an order
$\theta$ basis function and $G_{0}=\left(  \frac{1}{w}\right)  ^{\vee}$. Also
suppose $\theta<n$.

Then $\left(  \widehat{a}D\right)  ^{2n-1}G_{\theta}$ is bounded on any ball
containing $\operatorname*{supp}G_{0}$ when $2\theta>d$, and $\left(
\widehat{a}D\right)  ^{2n-1}G_{\theta}$ is bounded everywhere when
$2\theta\leq d$.

\textbf{Specifically}: When $2\theta>d$ and $x\in B_{R}$ and $R\geq\frac{1}%
{2}\operatorname*{diam}\operatorname*{supp}G_{0}$:%
\begin{align*}
\left\vert \left(  \widehat{a}D\right)  ^{2n-1}G_{\theta}\left(  x\right)
\right\vert  & \leq\tfrac{\left(  2\pi\right)  ^{-\frac{d}{2}}}{e\left(
\theta-\frac{d}{2}\right)  }\min\limits_{m<2\theta-d}\left\Vert \left(
\widehat{a}D\right)  ^{2n-1-m}G_{0}\right\Vert _{1}\left(  2R\right)
^{2\theta-d-m}m!\times\\
& \qquad\times\left\{
\begin{array}
[c]{ll}%
\frac{k_{m,2\theta-d}}{m!}\log R+\sum\limits_{j=1}^{m}\left(  2-\frac{1}%
{j}\right)  \frac{k_{j,2\theta-d}}{j!}, & d\text{ }even,\\
\frac{k_{m,2\theta-d}}{m!}, & d\text{ }odd.
\end{array}
\right.
\end{align*}

When $2\theta<d$:%
\begin{align*}
\left\Vert \left(  \widehat{a}D\right)  ^{2n-1}G_{\theta}\right\Vert
_{\infty}  & \leq\left(  \frac{\omega_{d}}{2\theta}\left\Vert \left(
\widehat{a}D\right)  ^{2n-1}G_{0}\right\Vert _{\infty}+\left\Vert \left(
\widehat{a}D\right)  ^{2n-1}G_{0}\right\Vert _{1}\right)  \times\\
& \qquad\times2^{d+1-2\theta}\left\{
\begin{array}
[c]{ll}%
\frac{1}{e\left(  1\right)  }\left(  \frac{d}{2}-1-\theta\right)  !\left\vert
\left(  -\frac{d}{2}\right)  _{\frac{d}{2}+1-\theta}\right\vert , & d\text{
}even,\\
\frac{1}{e\left(  1/2\right)  }\left\vert \left(  -\frac{1}{2}\right)
_{\frac{d+1}{2}-\theta}\left(  -\frac{d-1}{2}\right)  _{\frac{d+1}{2}-\theta
}\right\vert , & d\text{ }odd.
\end{array}
\right.
\end{align*}

When $2\theta=d$:%
\[
\left\Vert \left(  \widehat{a}D\right)  ^{2n-1}G_{\theta}\right\Vert _{\infty
}\leq\frac{2\omega_{d}}{d}\left\Vert \left(  \widehat{a}D\right)  ^{2n-1}%
G_{0}\right\Vert _{\infty}+2d\left\Vert \left(  \widehat{a}D\right)
^{2n-2}G_{0}\right\Vert _{1}.
\]

Here $\omega_{d}=\frac{2\pi^{d/2}}{\Gamma\left(  d/2\right)  }$ and the values
of the function $e$ are given in Remark \ref{Rem_Lem_e(s)_pos}.
\end{theorem}

\begin{proof}
We will consider three cases: $2\theta>d$, $2\theta<d$ and $2\theta
=d$:\medskip

\fbox{\textbf{Case 1}: $2\theta>d$} Suppose $0\leq m\leq\min\left\{
2n-1,2\theta-d-1\right\}  =2\theta-d-1$. Then%
\begin{align*}
\left(  2\pi\right)  ^{\frac{d}{2}}\left\vert \left(  \widehat{a}D\right)
^{2n-1}G_{\theta}\left(  x\right)  \right\vert  & =\left(  2\pi\right)
^{\frac{d}{2}}\left\vert \left(  \widehat{a}D\right)  ^{2n-1}\left(  G_{0}\ast
T_{\theta}\right)  \left(  x\right)  \right\vert \\
& \leq\int\left\vert \left(  \left(  \widehat{a}D\right)  ^{2n-1-m}%
G_{0}\right)  \left(  x-y\right)  \right\vert \left\vert \left(  \widehat
{a}D\right)  ^{m}T_{\theta}\left(  y\right)  \right\vert dy\\
& =\int\limits_{\left\vert x-y\right\vert \leq R}\left\vert \left(  \left(
\widehat{a}D\right)  ^{2n-1-m}G_{0}\right)  \left(  x-y\right)  \right\vert
\left\vert \left(  \widehat{a}D\right)  ^{m}T_{\theta}\left(  y\right)
\right\vert dy+\\
& \quad+\int\limits_{\left\vert x-y\right\vert \geq R}\left\vert \left(
\left(  \widehat{a}D\right)  ^{2n-1-m}G_{0}\right)  \left(  x-y\right)
\right\vert \left\vert \left(  \widehat{a}D\right)  ^{m}T_{\theta}\left(
y\right)  \right\vert dy.
\end{align*}

Choose $R\geq\frac{1}{2}\operatorname*{diam}\operatorname*{supp}G_{0}=\frac
{1}{2}\operatorname*{diam}\left[  -l,l\right]  ^{d}=\frac{1}{2}\sqrt{d\left(
2l\right)  ^{2}}=l\sqrt{d}\geq1$.

Then \textbf{near infinity}%
\begin{equation}
\int\limits_{\left\vert x-y\right\vert \geq R}\left\vert \left(  \left(
\widehat{a}D\right)  ^{2n-1-m}G_{0}\right)  \left(  x-y\right)  \right\vert
\left\vert \left(  \widehat{a}D\right)  ^{m}T_{\theta}\left(  y\right)
\right\vert dy=0,\quad x\in B_{R}.\label{p041}%
\end{equation}

\textbf{Near the origin}%
\begin{align*}
\int\limits_{\left\vert x-y\right\vert \leq R} &  \left\vert \left(  \left(
\widehat{a}D\right)  ^{2n-1-m}G_{0}\right)  \left(  x-y\right)  \right\vert
\left\vert \left(  \widehat{a}D\right)  ^{m}T_{\theta}\left(  y\right)
\right\vert dy\\
&  \leq\left(  \max_{\left\vert y-x\right\vert \leq R}\left\vert \left(
\widehat{a}D\right)  ^{m}T_{\theta}\right\vert \right)  \int%
\limits_{\left\vert y-x\right\vert \leq R}\left\vert \left(  \left(
\widehat{a}D\right)  ^{2n-1-m}G_{0}\right)  \left(  x-y\right)  \right\vert
dy\\
&  \leq\left(  \max_{\left\vert \cdot\right\vert \leq2R}\left\vert \left(
\widehat{a}D\right)  ^{m}T_{\theta}\right\vert \right)  \left\Vert \left(
\widehat{a}D\right)  ^{2n-1-m}G_{0}\right\Vert _{1},
\end{align*}

so that%
\[
\left\vert \left(  \widehat{a}D\right)  ^{2n-1}G_{\theta}\left(  x\right)
\right\vert \leq\frac{\left(  2\pi\right)  ^{-\frac{d}{2}}}{e\left(
\theta-\frac{d}{2}\right)  }\min_{m<2\theta-d}\left(  \left\Vert \left(
\widehat{a}D\right)  ^{2n-1-m}G_{0}\right\Vert _{1}\max_{\left\vert
\cdot\right\vert \leq2R}\left\vert \left(  \widehat{a}D\right)  ^{m}T_{\theta
}\right\vert \right)  ,\quad x\in\operatorname*{supp}G_{0}.
\]

Regarding the expression $\max\limits_{\left\vert \cdot\right\vert \leq
2R}\left\vert \left(  \widehat{a}D\right)  ^{m}T_{\theta}\right\vert $, the
two formulas \ref{p24} for $T_{\theta}$ lead us to consider $d$ odd and $d$
even.\medskip

\underline{If $d$ is \textbf{odd}} then from part 3 of Lemma
\ref{Lem_thin_plate_splin},%
\[
\left\vert \left(  \widehat{a}D\right)  ^{m}T_{\theta}\left(  x\right)
\right\vert \leq k_{m,2\theta-d}\left\vert x\right\vert ^{2\theta-d-m},
\]

and thus when $x\in B_{R}$,
\begin{align*}
\left\vert \left(  \widehat{a}D\right)  ^{2n-1}G_{\theta}\left(  x\right)
\right\vert  & \leq\frac{\left(  2\pi\right)  ^{-\frac{d}{2}}}{e\left(
\theta-\frac{d}{2}\right)  }\min_{m<2\theta-d}\left(  \left\Vert \left(
\widehat{a}D\right)  ^{2n-1-m}G_{0}\right\Vert _{1}k_{m,2\theta-d}\left\vert
x\right\vert ^{2\theta-d-m}\right) \\
& \leq\frac{\left(  2\pi\right)  ^{-\frac{d}{2}}}{e\left(  \theta-\frac{d}%
{2}\right)  }\min_{m<2\theta-d}\left(  k_{m,2\theta-d}\left\Vert \left(
\widehat{a}D\right)  ^{2n-1-m}G_{0}\right\Vert _{1}\left(  2R\right)
^{2\theta-d-m}\right)  .
\end{align*}
\medskip

\underline{If $d$ is \textbf{even}} then from part 4 of Lemma
\ref{Lem_thin_plate_splin},%
\[
\left(  \widehat{a}D\right)  ^{m}T_{\theta}\left(  x\right)  \leq m!\left(
\frac{k_{m,2\theta-d}}{m!}\left\vert \log\left\vert x\right\vert \right\vert
+\sum\limits_{j=1}^{m}\left(  2-\frac{1}{j}\right)  \frac{k_{j,2\theta-d}}%
{j!}\right)  \left\vert x\right\vert ^{2\theta-d-m},\quad1\leq m\leq
2\theta-d-1,
\]

and since $R\geq1$,%
\begin{align*}
& \max_{\left\vert \cdot\right\vert \leq2R}\left\vert \left(  \widehat
{a}D\right)  ^{m}T_{\theta}\right\vert \\
& \leq m!\left(  \frac{k_{m,2\theta-d}}{m!}\max_{t\in\left[  0,2R\right]
}\left\{  t^{2\theta-d-m}\left\vert \log t\right\vert \right\}  +\sum
\limits_{j=1}^{m}\left(  2-\frac{1}{j}\right)  \frac{k_{j,2\theta-d}}%
{j!}\left(  2R\right)  ^{2\theta-d-m}\right) \\
& =m!\left(  \frac{k_{m,2\theta-d}}{m!}\left(  2R\right)  ^{2\theta-d-m}\log
R+\sum\limits_{j=1}^{m}\left(  2-\frac{1}{j}\right)  \frac{k_{j,2\theta-d}%
}{j!}\left(  2R\right)  ^{2\theta-d-m}\right) \\
& =m!\left(  2R\right)  ^{2\theta-d-m}\left(  \frac{k_{m,2\theta-d}}{m!}\log
R+\sum\limits_{j=1}^{m}\left(  2-\frac{1}{j}\right)  \frac{k_{j,2\theta-d}%
}{j!}\right)
\end{align*}

so that when $x\in B_{R}$,
\begin{align}
&  \left\vert \left(  \widehat{a}D\right)  ^{2n-1}G_{\theta}\left(  x\right)
\right\vert \nonumber\\
&  \leq\tfrac{\left(  2\pi\right)  ^{-\frac{d}{2}}}{e\left(  \theta-\frac
{d}{2}\right)  }\min_{m<2\theta-d}\left(  \left\Vert \left(  \widehat
{a}D\right)  ^{2n-1-m}G_{0}\right\Vert _{1}\max_{\left\vert x\right\vert
\leq2R}\left\vert \left(  \widehat{a}D\right)  ^{m}T_{\theta}\left(  x\right)
\right\vert \right) \nonumber\\
&  \leq\tfrac{\left(  2\pi\right)  ^{-\frac{d}{2}}}{e\left(  \theta-\frac
{d}{2}\right)  }\min_{m<2\theta-d}\left(  \left\Vert \left(  \widehat
{a}D\right)  ^{2n-1-m}G_{0}\right\Vert _{1}m!\left(  2R\right)  ^{2\theta
-d-m}\left(  \frac{k_{m,2\theta-d}}{m!}\log R+\sum\limits_{j=1}^{m}\left(
2-\frac{1}{j}\right)  \frac{k_{j,2\theta-d}}{j!}\right)  \right)  .\label{p81}%
\end{align}
\medskip

\fbox{\textbf{Case 2}: $2\theta<d$} \underline{If $d$\textbf{\ is odd}} then
from Theorem \ref{Thm_ExtNatSplin_convol}, $T_{\frac{d+1}{2}}=-\left\vert
\cdot\right\vert $ and so if $b_{\theta}:=\frac{\left(  -1\right)
^{\frac{d+1}{2}-\theta}}{e\left(  1/2\right)  }$,%
\[
\left(  \widehat{a}D\right)  ^{2n-1}G_{\theta}=-b_{\theta}\left(  \widehat
{a}D\right)  ^{2n-1}\left(  \left\vert D\right\vert ^{2}\right)  ^{\frac
{d+1}{2}-\theta}\left(  G_{0}\ast\left\vert \cdot\right\vert \right)
=-b_{\theta}\left(  \widehat{a}D\right)  ^{2n-1}G_{0}\ast\left(  \left\vert
D\right\vert ^{2}\right)  ^{\frac{d+1}{2}-\theta}\left\vert \cdot\right\vert .
\]

From part 9 of \ref{Lem_iterLaplacian_rad} with $k=1$,%
\[
\left\vert D\right\vert ^{2m}\left\vert x\right\vert =2^{2m}\left(  -\frac
{1}{2}\right)  _{m}\left(  -\frac{d-1}{2}\right)  _{m}\left\vert x\right\vert
^{1-2m},
\]

so that%
\begin{align*}
\left(  \widehat{a}D\right)  ^{2n-1}G_{\theta}  & =b_{\theta}2^{2\left(
\frac{d+1}{2}-\theta\right)  }\left(  -\frac{1}{2}\right)  _{\frac{d+1}%
{2}-\theta}\left(  -\frac{d-1}{2}\right)  _{\frac{d+1}{2}-\theta}\left(
\widehat{a}D\right)  ^{2n-1}G_{0}\ast\left\vert \cdot\right\vert ^{1-2\left(
\frac{d+1}{2}-\theta\right)  }\\
& =b_{\theta}2^{d-2\theta+1}\left(  -\frac{1}{2}\right)  _{\frac{d+1}%
{2}-\theta}\left(  -\frac{d-1}{2}\right)  _{\frac{d+1}{2}-\theta}\left(
\widehat{a}D\right)  ^{2n-1}G_{0}\ast\left\vert \cdot\right\vert ^{-d+2\theta
}.
\end{align*}

Let $\left\{  \phi_{0}\left(  \left\vert \cdot\right\vert \right)
,\phi_{\infty}\left(  \left\vert \cdot\right\vert \right)  \right\}  $ be the
step function radial partition of unity such that\textbf{\ }$\phi_{0}\left(
\left\vert x\right\vert \right)  =1$ when $\left\vert x\right\vert \leq1$ and
zero otherwise. Then%
\begin{align*}
\left(  \widehat{a}D\right)  ^{2n-1}G_{0}\ast\left\vert \cdot\right\vert
^{-d+2\theta}  & =\left(  \widehat{a}D\right)  ^{2n-1}G_{0}\ast\left(
\left\vert \cdot\right\vert ^{-d+2\theta}\phi_{0}+\left\vert \cdot\right\vert
^{-\left(  d-2\theta\right)  }\phi_{\infty}\right) \\
& =\left(  \widehat{a}D\right)  ^{2n-1}G_{0}\ast\left\vert \cdot\right\vert
^{-d+2\theta}\phi_{0}+\left(  \widehat{a}D\right)  ^{d-2\theta}G_{0}%
\ast\left\vert \cdot\right\vert ^{-\left(  d-2\theta\right)  }\phi_{\infty},
\end{align*}

so that by using Young's inequality \ref{1.056} and then the spherical
coordinate result of part 1 of Corollary \ref{Cor_Thm_Integ_u(xy)f(|x|)dx},%
\begin{align*}
& \left\vert \left(  \left(  \widehat{a}D\right)  ^{2n-1}G_{0}\ast\left\vert
\cdot\right\vert ^{-d+2\theta}\right)  \left(  x\right)  \right\vert \\
& \leq\left\vert \left(  \left(  \widehat{a}D\right)  ^{2n-1}G_{0}%
\ast\left\vert \cdot\right\vert ^{-d+2\theta}\phi_{0}\right)  \left(
x\right)  \right\vert +\left\vert \left(  \left(  \widehat{a}D\right)
^{2n-1}G_{0}\ast\left\vert \cdot\right\vert ^{-\left(  d-2\theta\right)  }%
\phi_{\infty}\right)  \left(  x\right)  \right\vert \\
& \leq\left\Vert \left(  \widehat{a}D\right)  ^{2n-1}G_{0}\right\Vert
_{\infty}\left\Vert \left\vert \cdot\right\vert ^{-d+2\theta}\phi
_{0}\right\Vert _{1}+\left\Vert \left(  \widehat{a}D\right)  ^{2n-1}%
G_{0}\right\Vert _{1}\left\Vert \left\vert \cdot\right\vert ^{-\left(
d-2\theta\right)  }\phi_{\infty}\right\Vert _{\infty}\\
& =\left\Vert \left(  \widehat{a}D\right)  ^{2n-1}G_{0}\right\Vert _{\infty
}\int_{\left\vert \cdot\right\vert \leq1}\left\vert \cdot\right\vert
^{-d+2\theta}\phi_{0}\left(  \left\vert \cdot\right\vert \right)  +\left\Vert
\left(  \widehat{a}D\right)  ^{2n-1}G_{0}\right\Vert _{1}\left\Vert \left\vert
\cdot\right\vert ^{-\left(  d-2\theta\right)  }\right\Vert _{\infty;\left\vert
\cdot\right\vert \geq1}\\
& =\left\Vert \left(  \widehat{a}D\right)  ^{2n-1}G_{0}\right\Vert _{\infty
}\omega_{d}\int_{0}^{1}\rho^{-d+2\theta}\rho^{d-1}d\rho+\left\Vert \left(
\widehat{a}D\right)  ^{2n-1}G_{0}\right\Vert _{1}\\
& =\frac{\omega_{d}}{2\theta}\left\Vert \left(  \widehat{a}D\right)
^{2n-1}G_{0}\right\Vert _{\infty}+\left\Vert \left(  \widehat{a}D\right)
^{2n-1}G_{0}\right\Vert _{1},
\end{align*}

where $\omega_{d}=\frac{2\pi^{d/2}}{\Gamma\left(  d/2\right)  }$ for $d\geq1$.
Thus%
\begin{align}
& \left\vert \left(  \widehat{a}D\right)  ^{2n-1}G_{\theta}\left(  x\right)
\right\vert \nonumber\\
& \leq\left\vert b_{\theta}\right\vert 2^{d-2\theta+1}\left\vert \left(
-\frac{1}{2}\right)  _{\frac{d+1}{2}-\theta}\left(  -\frac{d-1}{2}\right)
_{\frac{d+1}{2}-\theta}\right\vert \left(  \frac{\omega_{d}}{2\theta
-1}\left\Vert \left(  \widehat{a}D\right)  ^{2n-1}G_{0}\right\Vert _{\infty
}+\left\Vert \left(  \widehat{a}D\right)  ^{2n-1}G_{0}\right\Vert _{1}\right)
\nonumber\\
& =\frac{2^{d-2\theta+1}}{e\left(  1/2\right)  }\left\vert \left(  -\frac
{1}{2}\right)  _{\frac{d+1}{2}-\theta}\left(  -\frac{d-1}{2}\right)
_{\frac{d+1}{2}-\theta}\right\vert \left(  \frac{\omega_{d}}{2\theta
}\left\Vert \left(  \widehat{a}D\right)  ^{2n-1}G_{0}\right\Vert _{\infty
}+\left\Vert \left(  \widehat{a}D\right)  ^{2n-1}G_{0}\right\Vert _{1}\right)
.\label{p78}%
\end{align}
\medskip

\underline{\textbf{If }$d$\textbf{\ is even}} set $b_{\theta}^{\prime}%
=\frac{\left(  -1\right)  ^{\frac{d}{2}+1-\theta}}{e\left(  1\right)  }$ and
write $d=2\theta+2m$ where $m\geq1$. Then from Theorem
\ref{Thm_ExtNatSplin_convol}, $T_{\frac{d}{2}+1}=\left\vert \cdot\right\vert
^{2}\log\left\vert \cdot\right\vert $ and%
\begin{align*}
\left(  \widehat{a}D\right)  ^{2n-1}G_{\theta}  & =\frac{\left(  -1\right)
^{\frac{d}{2}+1-\theta}}{e\left(  1\right)  }\left(  \widehat{a}D\right)
^{2n-1}\left(  \left\vert D\right\vert ^{2}\right)  ^{\frac{d}{2}+1-\theta
}\left(  G_{0}\ast T_{\frac{d}{2}+1}\right) \\
& =b_{\theta}^{\prime}\left(  \widehat{a}D\right)  ^{2n-1}\left\vert
D\right\vert ^{2\left(  m+1\right)  }\left(  G_{0}\ast\left(  \left\vert
\cdot\right\vert ^{2}\log\left\vert \cdot\right\vert \right)  \right) \\
& =b_{\theta}^{\prime}\left(  \left(  \widehat{a}D\right)  ^{2n-1}%
G_{0}\right)  \ast\left\vert D\right\vert ^{2\left(  m+1\right)  }\left(
\left\vert \cdot\right\vert ^{2}\log\left\vert \cdot\right\vert \right)  .
\end{align*}

By part 7 of Lemma \ref{Lem_iterLaplacian_rad} i.e.%
\[
\left\vert D\right\vert ^{2n}\left(  \left\vert x\right\vert ^{2}%
\log\left\vert x\right\vert \right)  =2^{2n-1}\left(  n-2\right)  !\left(
-\frac{d}{2}\right)  _{n}\left\vert x\right\vert ^{2-2n},\quad n\geq2,
\]

so that%
\[
\left\vert D\right\vert ^{2\left(  m+1\right)  }\left(  \left\vert
x\right\vert ^{2}\log\left\vert x\right\vert \right)  =2^{2m+1}\left(
m-1\right)  !\left(  -\frac{d}{2}\right)  _{m+1}\left\vert x\right\vert
^{-2m},\quad m\geq1,
\]

and hence%
\begin{align*}
\left(  \widehat{a}D\right)  ^{2n-1}G_{\theta}  & =b_{\theta}^{\prime}\left(
\left(  \widehat{a}D\right)  ^{2n-1}G_{0}\right)  \ast2^{2m+1}\left(
m-1\right)  !\left(  -\frac{d}{2}\right)  _{m+1}\left\vert \cdot\right\vert
^{-2m}\\
& =b_{\theta}^{\prime}2^{2m+1}\left(  m-1\right)  !\left(  -\frac{d}%
{2}\right)  _{m+1}\left(  \left(  \widehat{a}D\right)  ^{2n-1}G_{0}\right)
\ast\left\vert \cdot\right\vert ^{-2m}\\
& =b_{\theta}^{\prime}2^{2m+1}\left(  m-1\right)  !\left(  -\frac{d}%
{2}\right)  _{m+1}\left(  \left(  \widehat{a}D\right)  ^{2n-1}G_{0}\right)
\ast\left\vert \cdot\right\vert ^{-d+2\theta}.
\end{align*}

Again applying the partition of unity $\left\{  \phi_{0},\phi_{\infty
}\right\}  $ and then using Young's inequality to estimate each term we get%
\begin{align*}
& \left\Vert \left(  \left(  \widehat{a}D\right)  ^{2n-1}G_{0}\right)
\ast\left\vert \cdot\right\vert ^{-2m}\right\Vert _{\infty}\\
& \leq\left\Vert \left(  \left(  \widehat{a}D\right)  ^{2n-1}G_{0}\right)
\ast\phi_{0}\left\vert \cdot\right\vert ^{-2m}\right\Vert _{\infty}+\left\Vert
\left(  \left(  \widehat{a}D\right)  ^{2n-1}G_{0}\right)  \ast\phi_{\infty
}\left\vert \cdot\right\vert ^{-2m}\right\Vert _{\infty}\\
& \leq\left\Vert \left(  \widehat{a}D\right)  ^{2n-1}G_{0}\right\Vert
_{\infty}\left\Vert \phi_{0}\left\vert \cdot\right\vert ^{-2m}\right\Vert
_{1}+\left\Vert \left(  \widehat{a}D\right)  ^{2n-1}G_{0}\right\Vert
_{1}\left\Vert \phi_{\infty}\left\vert \cdot\right\vert ^{-2m}\right\Vert
_{\infty}\\
& =\left\Vert \left(  \widehat{a}D\right)  ^{2n-1}G_{0}\right\Vert _{\infty
}\int_{\left\vert \cdot\right\vert \leq1}\left\vert \cdot\right\vert
^{-2m}+\left\Vert \left(  \widehat{a}D\right)  ^{2n-1}G_{0}\right\Vert
_{1}\sup_{\left\vert x\right\vert \geq1}\left\vert x\right\vert ^{-2m}\\
& =\left\Vert \left(  \widehat{a}D\right)  ^{2n-1}G_{0}\right\Vert _{\infty
}\omega_{d}\int_{0}^{1}\rho^{-2m}\rho^{d-1}d\rho+\left\Vert \left(
\widehat{a}D\right)  ^{2n-1}G_{0}\right\Vert _{1}\\
& =\left\Vert \left(  \widehat{a}D\right)  ^{2n-1}G_{0}\right\Vert _{\infty
}\omega_{d}\int_{0}^{1}\rho^{-\left(  d-2\theta\right)  }\rho^{d-1}%
d\rho+\left\Vert \left(  \widehat{a}D\right)  ^{2n-1}G_{0}\right\Vert _{1}\\
& =\left\Vert \left(  \widehat{a}D\right)  ^{2n-1}G_{0}\right\Vert _{\infty
}\omega_{d}\int_{0}^{1}\rho^{2\theta-1}d\rho+\left\Vert \left(  \widehat
{a}D\right)  ^{2n-1}G_{0}\right\Vert _{1}\\
& =\frac{\omega_{d}}{2\theta}\left\Vert \left(  \widehat{a}D\right)
^{2n-1}G_{0}\right\Vert _{\infty}+\left\Vert \left(  \widehat{a}D\right)
^{2n-1}G_{0}\right\Vert _{1},
\end{align*}

and since $2m=d-2\theta$ this means that%
\begin{align}
& \left\Vert \left(  \widehat{a}D\right)  ^{2n-1}G_{\theta}\right\Vert
_{\infty}\nonumber\\
& \leq\left\vert b_{\theta}^{\prime}\right\vert 2^{2m+1}\left(  m-1\right)
!\left\vert \left(  -\frac{d}{2}\right)  _{m+1}\right\vert \left(
\frac{\omega_{d}}{2\theta}\left\Vert \left(  \widehat{a}D\right)  ^{2n-1}%
G_{0}\right\Vert _{\infty}+\left\Vert \left(  \widehat{a}D\right)
^{2n-1}G_{0}\right\Vert _{1}\right) \nonumber\\
& =\frac{2^{d-2\theta+1}\left(  \frac{d}{2}-\theta-1\right)  !}{e\left(
1\right)  }\left\vert \left(  -\frac{d}{2}\right)  _{\frac{d}{2}-\theta
+1}\right\vert \left(  \frac{\omega_{d}}{2\theta}\left\Vert \left(
\widehat{a}D\right)  ^{2n-1}G_{0}\right\Vert _{\infty}+\left\Vert \left(
\widehat{a}D\right)  ^{2n-1}G_{0}\right\Vert _{1}\right)  .\label{p79}%
\end{align}
\medskip

\fbox{\textbf{Case} $2\theta=d$} Here we will use\textbf{\ }a radial partition
of unity $\left\{  \psi_{0}\left(  \left\vert \cdot\right\vert \right)
,\psi_{\infty}\left(  \left\vert \cdot\right\vert \right)  \right\}  $\ which
is continuous and piecewise linear such that
\[
\psi_{0}\left(  s\right)  =\left\{
\begin{array}
[c]{ll}%
1, & 0\leq s\leq1/2,\\
2-2s, & 1/2\leq s\leq1,\\
0 & 1\leq s.
\end{array}
\right.
\]

Using part 7 of Lemma \ref{Lem_iterLaplacian_rad},%
\begin{align*}
\left(  \widehat{a}D\right)  ^{2n-1}G_{\theta}  & =\left(  \widehat
{a}D\right)  ^{2n-1}\left\vert D\right\vert ^{2}\left(  G_{0}\ast T_{\frac
{d}{2}+1}\right) \\
& =\left(  \widehat{a}D\right)  ^{2n-1}\left\vert D\right\vert ^{2}\left(
G_{0}\ast\left(  \left\vert \cdot\right\vert ^{2}\log\left\vert \cdot
\right\vert \right)  \right) \\
& =\left(  \widehat{a}D\right)  ^{2n-1}\left(  G_{0}\ast\left\vert
D\right\vert ^{2}\left(  \left\vert \cdot\right\vert ^{2}\log\left\vert
\cdot\right\vert \right)  \right) \\
& =\left(  \widehat{a}D\right)  ^{2n-1}\left(  G_{0}\ast\left(  1+2d\log
\left\vert x\right\vert \right)  \right) \\
& =\left(  \widehat{a}D\right)  ^{2n-1}\left(  G_{0}\ast1+2dG_{0}\ast
\log\left\vert \cdot\right\vert \right) \\
& =G_{0}\ast\left(  \widehat{a}D\right)  ^{2n-1}1+2d\left(  aD\right)
^{2n-1}\left(  G_{0}\ast\log\left\vert \cdot\right\vert \right) \\
& =2d\left(  \widehat{a}D\right)  ^{2n-1}\left(  G_{0}\ast\log\left\vert
\cdot\right\vert \right) \\
& =2d\left(  \widehat{a}D\right)  ^{2n-1}\left(  G_{0}\ast\psi_{0}%
\log\left\vert \cdot\right\vert \right)  +2d\left(  \widehat{a}D\right)
^{2n-1}\left(  G_{0}\ast\psi_{\infty}\log\left\vert \cdot\right\vert \right)
\\
& =2d\left(  \left(  \widehat{a}D\right)  ^{2n-1}G_{0}\right)  \ast\psi
_{0}\log\left\vert \cdot\right\vert +2d\left(  \left(  \widehat{a}D\right)
^{2n-2}G_{0}\right)  \ast\left(  \widehat{a}D\right)  \left(  \psi_{\infty
}\log\left\vert \cdot\right\vert \right)  ,
\end{align*}

so that%
\begin{align}
& \left\Vert \left(  \widehat{a}D\right)  ^{2n-1}G_{\theta}\right\Vert
_{\infty}\nonumber\\
& \leq2d\left\Vert \left(  \widehat{a}D\right)  ^{2n-1}G_{0}\right\Vert
_{\infty}\left\Vert \psi_{0}\log\left\vert \cdot\right\vert \right\Vert
_{1}+2d\left\Vert \left(  \widehat{a}D\right)  ^{2n-2}G_{0}\right\Vert
_{1}\left\Vert \left(  \widehat{a}D\right)  \left(  \psi_{\infty}%
\log\left\vert \cdot\right\vert \right)  \right\Vert _{\infty},\label{p82}%
\end{align}

provided these norms exist.

From Lemma \ref{Lem_basis_ord0_Bsplin_W3} the norms relating to $G_{0}$ are
all finite. Also from \ref{Ap148} and formula 610.9 of Dwight \cite{Dwight61},%
\begin{align*}
\left\Vert \psi_{0}\log\left\vert x\right\vert \right\Vert _{1}=-\int%
_{\left\vert x\right\vert \leq1}\log\left\vert x\right\vert dx=-\omega_{d}%
\int_{0}^{1}\rho^{d-1}\log\rho d\rho & =-\omega_{d}\left[  \frac{\rho^{d}}%
{d}\log\rho-\frac{\rho^{d}}{d^{2}}\right]  _{0}^{1}\\
& =\frac{\omega_{d}}{d^{2}}<\infty,
\end{align*}

From part 2 of Lemma \ref{Lem_op_aD_estim},
\begin{align*}
\widehat{a}D_{x}\left(  \psi_{\infty}\left(  \left\vert x\right\vert \right)
\log\left\vert x\right\vert \right)   & =\widehat{a}\widehat{x}D_{s}\left(
\psi_{\infty}\left(  s\right)  \log s\right)  \left(  s=\left\vert
x\right\vert \right) \\
& =\widehat{a}\widehat{x}\left(  \psi_{\infty}^{\prime}\left(  s\right)  \log
s+\frac{\psi_{\infty}\left(  s\right)  }{s}\right)  \left(  s=\left\vert
x\right\vert \right) \\
& =\left(  \log\left\vert x\right\vert \right)  \widehat{a}D\psi_{\infty
}\left(  x\right)  +\psi_{\infty}\left(  x\right)  \widehat{a}D\log\left\vert
x\right\vert \\
& =\widehat{a}\widehat{x}\psi_{\infty}^{\prime}\log\left\vert x\right\vert
+\frac{\widehat{a}\widehat{x}}{\left\vert \cdot\right\vert }\psi_{\infty}\\
& =\widehat{a}\widehat{x}\left(  \psi_{\infty}^{\prime}\left(  \left\vert
x\right\vert \right)  \log\left\vert x\right\vert +\frac{1}{\left\vert
x\right\vert }\psi_{\infty}\left(  \left\vert x\right\vert \right)  \right)  ,
\end{align*}

so that%
\begin{align*}
\left\Vert \left(  \widehat{a}D\right)  \left(  \psi_{\infty}\log\left\vert
x\right\vert \right)  \right\Vert _{\infty}  & =\max_{t\geq\frac{1}{2}%
}\left\vert \psi_{\infty}^{\prime}\left(  t\right)  \log t+\frac{1}{t}%
\psi_{\infty}\left(  t\right)  \right\vert \\
& =\max\left\{  \max_{t\in\left[  \frac{1}{2},1\right]  }\left\vert
\psi_{\infty}^{\prime}\left(  t\right)  \log t+\frac{1}{t}\psi_{\infty}\left(
t\right)  \right\vert ,\max_{t\in\left[  1,\infty\right]  }\left\vert
\psi_{\infty}^{\prime}\left(  t\right)  \log t+\frac{1}{t}\psi_{\infty}\left(
t\right)  \right\vert \right\} \\
& =\max\left\{  \max_{t\in\left[  \frac{1}{2},1\right]  }\left\vert \log
t+\frac{1}{t}\left(  -1+2t\right)  \right\vert ,\max_{t\in\left[
1,\infty\right]  }\frac{1}{t}\right\} \\
& =\max\left\{  \max_{t\in\left[  \frac{1}{2},1\right]  }\left\vert \log
t-\frac{1}{t}+2\right\vert ,1\right\} \\
& =1.
\end{align*}

Thus \ref{p82} becomes%
\begin{align}
\left\Vert \left(  \widehat{a}D\right)  ^{2n-1}G_{\theta}\right\Vert
_{\infty}  & \leq2d\left\Vert \left(  \widehat{a}D\right)  ^{2n-1}%
G_{0}\right\Vert _{\infty}\frac{\omega_{d}}{d^{2}}+2d\left\Vert \left(
\widehat{a}D\right)  ^{2n-2}G_{0}\right\Vert _{1}\nonumber\\
& =\frac{2\omega_{d}}{d}\left\Vert \left(  \widehat{a}D\right)  ^{2n-1}%
G_{0}\right\Vert _{\infty}+2d\left\Vert \left(  \widehat{a}D\right)
^{2n-2}G_{0}\right\Vert _{1}.\label{p83}%
\end{align}

\end{proof}

\begin{remark}
\textbf{Part 1} Recalling that $b_{+}:=\left(  \left\vert b_{k}\right\vert
\right)  $ and $\left\vert b\right\vert _{1}:=%
{\textstyle\sum\limits_{k=1}^{d}}
\left\vert b_{k}\right\vert \leq\left\vert \mathbf{1}\right\vert \left\vert
b\right\vert =d^{1/2}\left\vert b\right\vert $ we have%
\begin{align*}
\frac{1}{k!}\left\Vert \left(  \widehat{a}D\right)  ^{k}G_{0}\right\Vert  &
=\left\Vert
{\textstyle\sum\limits_{\left\vert \beta\right\vert =k}}
\frac{\widehat{a}^{\beta}}{\beta!}D^{\beta}G_{0}\right\Vert \leq%
{\textstyle\sum\limits_{\left\vert \beta\right\vert =k}}
\frac{\widehat{a}_{+}^{\beta}}{\beta!}\left\Vert D^{\beta}G_{0}\right\Vert
\leq\\
& \leq\max_{\left\vert \beta\right\vert =k}\left\Vert D^{\beta}G_{0}%
\right\Vert
{\textstyle\sum\limits_{\left\vert \beta\right\vert =k}}
\frac{\left\vert \widehat{a}^{\beta}\right\vert }{\beta!}=\max_{\left\vert
\beta\right\vert =k}\left\Vert D^{\beta}G_{0}\right\Vert
{\textstyle\sum\limits_{\left\vert \beta\right\vert =k}}
\frac{\widehat{a}_{+}^{\beta}}{\beta!}=\\
& =\max_{\left\vert \beta\right\vert =k}\left\Vert D^{\beta}G_{0}\right\Vert
{\textstyle\sum\limits_{\left\vert \beta\right\vert =k}}
\frac{\widehat{a}_{+}^{\beta}\mathbf{1}^{\beta}}{\beta!}=\frac{\left(
\widehat{a}_{+}\mathbf{1}\right)  ^{k}}{k!}\max_{\left\vert \beta\right\vert
=k}\left\Vert D^{\beta}G_{0}\right\Vert =\\
& =\frac{\left\vert \widehat{a}\right\vert _{1}^{k}}{k!}\max_{\left\vert
\beta\right\vert =k}\left\Vert D^{\beta}G_{0}\right\Vert ,
\end{align*}

so that applying the Cauchy-Schwartz inequality yields%
\[
\left\Vert \left(  \widehat{a}D\right)  ^{k}G_{0}\right\Vert \leq\left\vert
\widehat{a}\right\vert _{1}^{k}\max_{\left\vert \beta\right\vert =k}\left\Vert
D^{\beta}G_{0}\right\Vert \leq d^{k/2}\max_{\left\vert \beta\right\vert
=k}\left\Vert D^{\beta}G_{0}\right\Vert .
\]
\medskip

\textbf{Part 2}%
\[
\left\Vert \left(  \widehat{a}D\right)  ^{k}G_{0}\right\Vert _{1}\leq\left(
\operatorname*{vol}\operatorname*{supp}G_{0}\right)  ^{\frac{1}{2}}\left\Vert
\left(  \widehat{a}D\right)  ^{k}G_{0}\right\Vert _{2},
\]

and%
\[
\left\Vert \left(  \widehat{a}D\right)  ^{k}G_{0}\right\Vert _{2}^{2}%
=\int\left\vert \left(  \widehat{a}D\right)  ^{k}G_{0}\right\vert ^{2}%
=\int\left\vert \left(  \left(  \widehat{a}D\right)  ^{k}G_{0}\right)
^{\wedge}\right\vert ^{2}=\int\left\vert \left(  \widehat{a}\xi\right)
^{k}\widehat{G}_{0}\right\vert ^{2}=\int\frac{\left(  \widehat{a}\xi\right)
^{2k}}{w_{0}^{2}}.
\]

Thus%
\[
\frac{1}{\left(  2k\right)  !}\left\Vert \left(  \widehat{a}D\right)
^{k}G_{0}\right\Vert _{2}^{2}=\int\frac{\left(  \widehat{a}\xi\right)  ^{2k}%
}{w_{0}^{2}}%
\]

\end{remark}

\subsection{?? Example: The tensor product function $\Lambda^{2}$ has order
$1$?}

??

\chapter{Calculating the weight function from the basis function without using
$S_{\emptyset,2n}$\label{Ch_basis_no_So2n}}

\section{Introduction}

In this chapter we will prove some more properties of the basis functions
studied in Chapter \ref{Ch_weight_fn_exten} where these objects were defined
in terms of a weight function $w$ and a positive integer order parameter
$\theta$. In Chapter \ref{Ch_weight_fn_exten} these weight function properties
were used to define reproducing kernel semi-Hilbert spaces of continuous
functions $X_{w}^{\theta}$ and continuous basis functions of order $\theta$.
In the Chapters \ref{Ch_Interpol}, \ref{Ch_ExactSmth}, \ref{Ch_Approx_smth}
these semi-Hilbert spaces will be used to formulate and study several
non-parametric, basis function interpolation and smoothing problems with the
basis functions being used to represent the solutions to these problems.

In this chapter we are only concerned with the basis functions and we prove a
result which will enable us to determine whether a given function is a basis
function without recourse to the basis function definition \ref{p01} used by
Light and Wayne in \cite{LightWayneX98Weight}. In practice their definition is
difficult to use and involves a weight function and the bounded linear
functionals on the subspace $S_{\emptyset,2\theta}=\left\{  \phi\in
S:D^{\alpha}\phi\left(  0\right)  =0,\text{\ }\left\vert \alpha\right\vert
<2\theta\right\}  $ of the test functions of the tempered distributions $S$
(Appendix \ref{Sect_tempered_distrib}). We give a simple test which can be
applied to the tempered distribution Fourier transform of a continuous function.

These results are then applied to several classes of well-known radial basis
functions, the choice here following Dyn \cite{Dyn89}: the thin-plate splines,
the shifted thin-plate splines, the multiquadric and inverse multiquadric
functions and the Gaussian. These classes of basis functions are well known in
the literature and details are given in Figure \ref{m24} below. In the last
section I will illustrate the method using a non-radial example: the
fundamental solutions of homogeneous elliptic differential operators of even order.

\section{Theory}

In this document we will prove that some of the important classes of functions
used to define basis function interpolants and smoothers are basis functions
in the Light sense i.e. generated by weight functions, without recourse to the
awkward Definition \ref{Def_basis_distrib} which uses the spaces
$S_{\emptyset,2\theta}$ and $S_{\emptyset,2\theta}^{\prime}$.

Our choice of basis functions is given in Table \ref{m24} below and follows
Dyn \cite{Dyn89}. These basis functions are well known in the literature.
Theorem \ref{Thm_rad_avoid_So2n} is our main result and it will be applied to
the various classes of radial and non-radial basis functions.

We will need the following basis distribution and weight function properties
which were proved in the previous chapter:

\begin{summary}
\label{Sum_basis_properties}\ Suppose the weight function $w$ has property W2
i.e. sub-properties W2.1 and W2.2. Then:

\begin{enumerate}
\item by part 2 of Appendix \ref{SbSect_property_S'}, $1/w$ is a regular
tempered distribution and thus $1/w\in L_{loc}^{1}\cap S^{\prime}$.

\item We can define a tempered basis distribution $G$ of order $\theta\geq1 $
generated by $w$.

\item $\left\vert \cdot\right\vert ^{2\theta}\widehat{G}=\frac{1}{w}$ as
tempered distributions (part 2 Theorem \ref{Thm_basis_smth_W3.2_r3_pos}).

\item $G+P_{2\theta-1}$ is the set of all basis distributions generated by $w$
(Theorem \ref{Thm_basis_fn_set}).\medskip

Now suppose the weight function $w$ only has properties W2.1 and \textbf{W3.2}
for order $\theta$ and $\kappa$. Then:\smallskip

\item $w$ also has property W2.2 (part 3 Theorem
\ref{Thm_weight_property_relat}).

\item $G\in C_{BP}^{\left(  \left\lfloor 2\kappa\right\rfloor \right)  }$,
where $\left\lfloor \cdot\right\rfloor $ is the floor function (part 4 Theorem
\ref{Thm_basis_smth_W3.2_r3_pos}). We say that $G$ is a basis function of
order $\theta$ generated by $w$.

\item If $w$ is radial then there exists a radial basis function (part 3
Corollary \ref{Cor_Thm_Grho}).

\item If $w$ is homogeneous of order $s$ then there exists a basis function
$G$ which is homogeneous of order $s+2\theta-d$ modulo a polynomial of order
at most $2\theta$ (part 4 Corollary \ref{Cor_Thm_Grho}). More precisely
\[
G\left(  tx\right)  -t^{s+2\theta-d}G\left(  x\right)  =\sum_{\left\vert
\alpha\right\vert <2\theta}q_{\alpha}\left(  t\right)  x^{\alpha},\quad
x\in\mathbb{R}^{d},\text{ }t>0.
\]

\end{enumerate}
\end{summary}

The next two theorems are the main results and will allow us to determine when
a $C_{BP}^{\left(  0\right)  }$ function is a basis function in the Light
sense by only studying its Fourier transform and so avoid the need to use the
awkward Definition \ref{Def_basis_distrib} of a basis function which uses
$S_{\emptyset,2n}$ subspaces. The first theorem deals with radial basis
functions and the second with homogeneous basis functions.

\begin{theorem}
\label{Thm_rad_avoid_So2n}\textbf{Radial case} Suppose $H\in S^{\prime}$ is
radial and that $\left\vert \cdot\right\vert ^{2\theta}\widehat{H}\in
L_{loc}^{1}$. Hence for some $\mathcal{B}\subseteq\left\{  \mathbf{0}\right\}
$ we can define the function $H_{F}\in L_{loc}^{1}\left(  \mathbb{R}%
^{d}\setminus\mathcal{B}\right)  $ by $H_{F}=\widehat{H}$ on $\mathbb{R}%
^{d}\setminus\mathcal{B}$.

Suppose further that $H_{F}\in C^{\left(  0\right)  }\left(  \mathbb{R}%
^{d}\setminus0\right)  $ and $H_{F}\left(  \xi\right)  >0$ on $\mathbb{R}%
^{d}\setminus0$. Now define the function $w$ by
\begin{equation}
w\left(  \xi\right)  =\frac{1}{\left\vert \xi\right\vert ^{2\theta}%
H_{F}\left(  \xi\right)  },\quad\xi\in\mathbb{R}^{d}\setminus0.\label{m22}%
\end{equation}

Then:

\begin{enumerate}
\item $w$ satisfies weight function property W1 w.r.t. the set $\mathcal{A}%
=\left\{  0\right\}  $.

\item Suppose $w$ also has weight function properties W2.1 and \textbf{W3.2}
for some order $\theta$ and $\kappa$. Then $H\in C_{BP}^{\left(  \left\lfloor
2\kappa\right\rfloor \right)  }$ and $H$ is a basis function of order $\theta$
generated by $w$.
\end{enumerate}
\end{theorem}

\begin{proof}
\textbf{Part 1} Clearly $\mathcal{A}$ is a closed set of measure zero and the
properties of $H$ imply that $w\in C^{\left(  0\right)  }\left(
\mathbb{R}^{d}\setminus\mathcal{A}\right)  $ and $w\left(  \xi\right)  >0$ on
$\mathbb{R}^{d}\setminus\mathcal{A}$. Hence $w$ has property W1 with w.r.t.
the set $\mathcal{A}=\left\{  0\right\}  $.\medskip

\textbf{Part 2}. Since $w$ has properties W2.1 and W3.2, part 5 of Summary
\ref{Sum_basis_properties} implies that $w$ has property W2. Part 1 of Summary
\ref{Sum_basis_properties} then implies that $1/w\in L_{loc}^{1}\cap
S^{\prime}$.

By definition of $H_{F}$, $\left\vert \cdot\right\vert ^{2\theta}\widehat
{H}=\left\vert \cdot\right\vert ^{2\theta}H_{F}$ as distributions on
$\mathbb{R}^{d}\setminus0$. By \ref{m22}, $\left\vert \cdot\right\vert
^{2\theta}H_{F}=1/w$ a.e. so that $\left\vert \cdot\right\vert ^{2\theta}%
H_{F}\in L_{loc}^{1}$ and since we have assumed that $\left\vert
\cdot\right\vert ^{2\theta}\widehat{H}\in L_{loc}^{1}$ it follows that
$\left\vert \cdot\right\vert ^{2\theta}H_{F}=\left\vert \cdot\right\vert
^{2\theta}\widehat{H}=1/w$ as distributions. But $\left\vert \cdot\right\vert
^{2\theta}\widehat{H}\in S^{\prime}$ and $1/w\in S^{\prime}$ so $\left\vert
\cdot\right\vert ^{2\theta}\widehat{H}=1/w,$ as tempered distributions. If $G$
is a basis distribution of order $\theta$ generated by $w$ then by part 3 of
Summary \ref{Sum_basis_properties}, $\left\vert \cdot\right\vert ^{2\theta
}\widehat{G}=1/w$ and the basis (tempered) distribution definition implies
$\widehat{G}=H_{F}$ on $\mathbb{R}^{d}\setminus0$ as distributions. From the
definition of $H_{F}$, $\widehat{H}=H_{F}$ on $\mathbb{R}^{d}\setminus0$, and
so $\operatorname*{supp}\left(  \widehat{G}-\widehat{H}\right)  \subseteq
\left\{  0\right\}  $, which implies $G-H=p$ where $p$ is some polynomial.

However, we have proved that $\left\vert \cdot\right\vert ^{2\theta}%
\widehat{G}=1/w$ and $\left\vert \cdot\right\vert ^{2\theta}\widehat{H}=1/w $
so, $0=\left\vert \cdot\right\vert ^{2\theta}\left(  \widehat{G}-\widehat
{H}\right)  =\left\vert \cdot\right\vert ^{2\theta}\widehat{p}$ i.e.
$\left\vert D\right\vert ^{2\theta}p=0$.

Since $H$ is radial it follows from \ref{m22} that $w$ is radial and hence
part 7 of Summary \ref{Sum_basis_properties} tells us that $w$ has a radial
basis function and so we can assume that $G$ is radial. Thus $p$ is a radial
polynomial and hence it must have the form $p=\sum_{k=0}^{m}a_{k}\left\vert
\cdot\right\vert ^{2k}$. But when $m\leq k$, $\left\vert D\right\vert
^{2m}\left(  \left\vert \cdot\right\vert ^{2k}\right)  =c_{m,k}\left\vert
\cdot\right\vert ^{2k-2m}$ for some constant $c_{m,k}$, so that
\[
0=\left\vert D\right\vert ^{2\theta}p=\sum_{k=0}^{m}a_{k}\left\vert
D\right\vert ^{2\theta}\left(  \left\vert \cdot\right\vert ^{2k}\right)
=\sum_{k=\theta}^{m}c_{\theta,k}a_{k}\left\vert \cdot\right\vert ^{2k-2\theta
},
\]

implies $a_{k}=0$ for $k\geq\theta$ and hence that $p\in P_{2\theta-1}$. Thus
$H\in G+P_{2\theta-1}$ and parts 4 and 6 of Summary \ref{Sum_basis_properties}
means that $H$ is a basis function and so $H\in C_{BP}^{\left(  \left\lfloor
2\kappa\right\rfloor \right)  }$.
\end{proof}

\begin{lemma}
\label{Lem_wt_homog_W2_W3.2}Suppose $w$ is a weight function w.r.t.
$\mathcal{A}=\left\{  0\right\}  $ and suppose $w$ is also homogeneous of
order $s$ on $\mathbb{R}^{d}\setminus\mathcal{A}$. Then $w$ has property W2
for $\sigma>0 $ and property W3.2 for order $\theta\geq1$ and smoothness
$\kappa$ iff the inequalities%
\begin{align}
d-2\theta+2\kappa & <s<d,\label{m00}\\
d-2\sigma & <s,\nonumber\\
\kappa & <\theta,\nonumber
\end{align}

hold.
\end{lemma}

\begin{proof}
The homogeneity of $w$ implies $w\left(  x\right)  =\left\vert x\right\vert
^{s}w\left(  \frac{x}{\left\vert x\right\vert }\right)  $ on $\mathbb{R}%
^{d}\setminus0$, and since $w\left(  \frac{x}{\left\vert x\right\vert
}\right)  $ is continuous and positive on $\mathbb{R}^{d}\setminus0$, there
exist constants $C_{i}>0$ such that $C_{1}\leq w\left(  \frac{x}{\left\vert
x\right\vert }\right)  \leq C_{2}$ when $x\neq0$.

Property W3.2 implies that for $0\leq t\leq\kappa$,%
\[
\infty>\int_{\left\vert x\right\vert \geq r_{3}}\frac{\left\vert x\right\vert
^{2t}}{w\left\vert x\right\vert ^{2\theta}}=\int_{\left\vert x\right\vert \geq
r_{3}}\frac{\left\vert x\right\vert ^{2t}}{w\left(  \frac{x}{\left\vert
x\right\vert }\right)  \left\vert x\right\vert ^{2\theta+s}}\geq\frac{1}%
{C_{2}}\int_{\left\vert x\right\vert \geq r_{3}}\frac{\left\vert x\right\vert
^{2t}}{\left\vert x\right\vert ^{2\theta+s}},
\]

which exists iff $2\theta+s-2\kappa>d$.

Property W2 comprises the conditions: $1/w\in L_{loc}^{1}$ and $\int%
_{\left\vert \cdot\right\vert \geq r_{2}}\frac{1}{w\left\vert \cdot\right\vert
^{2\sigma}}<\infty$ for some $\sigma>0$ and $r_{2}>0$.

Thus W2 holds iff $s<d$ and $s+2\sigma>d$. Further $d-2\theta+2\kappa<s<d$
implies $\kappa<\theta$.
\end{proof}

\begin{theorem}
\label{Thm_homog_avoid_So2n}\textbf{Homogeneous case} Suppose $H\in S^{\prime
}$ is homogeneous of order $t$ and that $\left\vert \cdot\right\vert
^{2\theta}\widehat{H}\in L_{loc}^{1}$. Hence for some $\mathcal{B}%
\subseteq\left\{  \mathbf{0}\right\}  $ we can define the function $H_{F}\in
L_{loc}^{1}\left(  \mathbb{R}^{d}\setminus\mathcal{B}\right)  $ by
$H_{F}=\widehat{H}$ on $\mathbb{R}^{d}\setminus\mathcal{B}$.

Suppose further that $H_{F}\in C^{\left(  0\right)  }\left(  \mathbb{R}%
^{d}\setminus0\right)  $ and $H_{F}\left(  \xi\right)  >0$ on $\mathbb{R}%
^{d}\setminus0$. Now define the function $w$ by
\begin{equation}
w\left(  \xi\right)  =\frac{1}{\left\vert \xi\right\vert ^{2\theta}%
H_{F}\left(  \xi\right)  },\quad\xi\in\mathbb{R}^{d}\setminus0.\label{m27}%
\end{equation}

Then:

\begin{enumerate}
\item $w$ satisfies weight function property W1 w.r.t. the set $\mathcal{A}%
=\left\{  0\right\}  $.

Suppose $w$ also has weight function properties W2.1 and W3.2 for order
$\theta$ and smoothness $\kappa$. Then:

\item $H\in C_{BP}^{\left(  \left\lfloor 2\kappa\right\rfloor \right)  }$ and
$H$ is a basis function of order $\theta$ generated by $w$.
\end{enumerate}
\end{theorem}

\begin{proof}
\textbf{Part 1} Clearly $\mathcal{A}$ is a closed set of measure zero and the
properties of $H$ imply that $w\in C^{\left(  0\right)  }\left(
\mathbb{R}^{d}\setminus\mathcal{A}\right)  $ and $w\left(  \xi\right)  >0$ on
$\mathbb{R}^{d}\setminus\mathcal{A}$. Hence $w$ has property W1 with w.r.t.
the set $\mathcal{A}=\left\{  0\right\}  $.\medskip

\textbf{Part 2}. Since $w$ has properties W2.1 and W3.2, part 5 of Summary
\ref{Sum_basis_properties} implies that $w$ has property W2. Part 1 of Summary
\ref{Sum_basis_properties} then implies that $1/w\in L_{loc}^{1}\cap
S^{\prime}$.

By definition of $H_{F}$, $\left\vert \cdot\right\vert ^{2\theta}\widehat
{H}=\left\vert \cdot\right\vert ^{2\theta}H_{F}$ as distributions on
$\mathbb{R}^{d}\setminus0$. By \ref{m22}, $\left\vert \cdot\right\vert
^{2\theta}H_{F}=1/w$ a.e. so that $\left\vert \cdot\right\vert ^{2\theta}%
H_{F}\in L_{loc}^{1}$ and since we have assumed that $\left\vert
\cdot\right\vert ^{2\theta}\widehat{H}\in L_{loc}^{1}$ it follows that
$\left\vert \cdot\right\vert ^{2\theta}H_{F}=\left\vert \cdot\right\vert
^{2\theta}\widehat{H}=1/w$ as distributions. But $\left\vert \cdot\right\vert
^{2\theta}\widehat{H}\in S^{\prime}$ and $1/w\in S^{\prime}$ so $\left\vert
\cdot\right\vert ^{2\theta}\widehat{H}=1/w,$ as tempered distributions. If $G$
is a basis distribution of order $\theta$ generated by $w$ then by part 3 of
Summary \ref{Sum_basis_properties}, $\left\vert \cdot\right\vert ^{2\theta
}\widehat{G}=1/w$ and the basis (tempered) distribution definition implies
$\widehat{G}=H_{F}$ on $\mathbb{R}^{d}\setminus0$ as distributions. From the
definition of $H_{F}$, $\widehat{H}=H_{F}$ on $\mathbb{R}^{d}\setminus0$, and
so $\operatorname*{supp}\left(  \widehat{G}-\widehat{H}\right)  \subseteq
\left\{  0\right\}  $, which implies $G-H=p$ where $p$ is some polynomial.

However, we have proved that $\left\vert \cdot\right\vert ^{2\theta}%
\widehat{G}=1/w$ and $\left\vert \cdot\right\vert ^{2\theta}\widehat{H}=1/w $
so, $0=\left\vert \cdot\right\vert ^{2\theta}\left(  \widehat{G}-\widehat
{H}\right)  =\left\vert \cdot\right\vert ^{2\theta}\widehat{p}$ i.e.
$\left\vert D\right\vert ^{2\theta}p=0$.

Since $H\in S^{\prime}$ has homogeneity $t$, $\widehat{H}$ has homogeneity
$-d-t$ and so $H_{F}$ has homogeneity $-d-t$ on $\mathbb{R}^{d}\setminus0$.
Equation \ref{m27} now implies that $w$ has homogeneity $t+d-2\theta$ on
$\mathbb{R}^{d}\setminus0$. But from part 8 of Summary
\ref{Sum_basis_properties} there exists a basis function $G$ which has
homogeneity $t^{s+2\theta-d}$ in the sense that%
\[
G\left(  \mu x\right)  =\mu^{t}G\left(  x\right)  +\sum_{\left\vert
\alpha\right\vert <2\theta}q_{\alpha}\left(  \mu\right)  x^{\alpha},\quad
x\in\mathbb{R}^{d},\text{ }\mu>0,
\]

and thus%
\begin{align*}
p\left(  \mu x\right)  =G\left(  \mu x\right)  -H\left(  \mu x\right)   &
=\mu^{t}G\left(  x\right)  +\sum_{\left\vert \alpha\right\vert <2\theta
}q_{\alpha}\left(  \mu\right)  x^{\alpha}-\mu^{t}H\left(  x\right) \\
& =\mu^{t}p\left(  x\right)  +\sum_{\left\vert \alpha\right\vert <2\theta
}q_{\alpha}\left(  \mu\right)  x^{\alpha}.
\end{align*}

Suppose $p\left(  x\right)  =\sum_{\alpha}p_{\alpha}x^{\alpha}$ where
$p_{\alpha}\neq0$ for some $\left\vert \alpha\right\vert \geq2\theta$. Then%
\[
\sum_{\left\vert \alpha\right\vert \geq2\theta}p_{\alpha}\left(  \mu x\right)
^{\alpha}=\mu^{t}\sum_{\left\vert \alpha\right\vert \geq2\theta}p_{\alpha
}x^{\alpha}.
\]

Comparing the monomial terms $x^{\alpha}$ with $\left\vert \alpha\right\vert
\geq2\theta$ implies that
\[
p_{\alpha}\mu^{\left\vert \alpha\right\vert -t}=p_{\alpha}\text{ when
}p_{\alpha}\neq0,\text{ }\left\vert \alpha\right\vert \geq2\theta,\text{ }%
\mu>0.
\]

Since $w$ has homogeneity $t+d-2\theta$ the inequalities \ref{m00} of Lemma
\ref{Lem_wt_homog_W2_W3.2} imply $t<2\theta$, and so $p_{\alpha}=0$ when
$\left\vert \alpha\right\vert \geq2\theta$, which in turn implies that $p\in
P_{2\theta-1}$. Thus $H\in G+P_{2\theta-1}$ and part 4 of Summary
\ref{Sum_basis_properties} means that $H$ is a basis function.
\end{proof}

\section{Examples: \textit{radial basis functions} generated by weight
functions in W3.2\label{Sect_rad_basis_fns}}

Table \ref{m24} below gives a list of the radial functions derived from Dyn
\cite{Dyn89} and which are used as examples. Both types of splines have been
generalized in the sense that the parameter $s$ has been allowed to take
non-integer values. Section 1 of Dyn \cite{Dyn89} describes splines for which
$s$ has integer values. However, latter in Section 2 after Theorem 4, more
general basis functions are introduced in equations 23 to 26. These are listed
in Table 4.1, with a significant change. The change is to multiply by $\left(
-1\right)  ^{\left\lceil s\right\rceil }$ or $\left(  -1\right)  ^{s+1}$ so
that the Fourier transform is positive.%

\begin{gather}%
\begin{tabular}
[c]{|l||ll|}\hline
Type & \multicolumn{2}{||l|}{Radial basis function}\\\hline\hline
Surface/thin- & $\left(  -1\right)  ^{\left\lceil s\right\rceil }r^{2s},$ &
$s>0,$ $s\neq1,2,3\ldots.$\\\cline{2-3}%
plate spline & $\left(  -1\right)  ^{s+1}r^{2s}\log r,$ & $s=1,2,3\ldots.$\\
(generalized $s$) &  & \\\hline
Shifted surface & $\left(  -1\right)  ^{\left\lceil s\right\rceil }\left(
a^{2}+r^{2}\right)  ^{s},$ & $a>0,$ $s>-d/2,$\\
/thin-plate &  & $s\neq1,2,3\ldots.$\\\cline{2-3}%
spline & $\frac{\left(  -1\right)  ^{s+1}}{2}\left(  a^{2}+r^{2}\right)
^{s}\log\left(  a^{2}+r^{2}\right)  $ & $a>0,$\\
(generalized $s$) &  & $s=1,2,3,\ldots.$\\\hline
Multiquadric & $-\left(  a^{2}+r^{2}\right)  ^{1/2},$ & $a>0,$ $d>1.$\\\hline
Inverse multi- & $\left(  a^{2}+r^{2}\right)  ^{-1/2},$ & $a>0.$\\
quadric &  & \\\hline
Gaussian & $\exp\left(  -\frac{r^{2}}{2}\right)  .$ & \\\hline
\end{tabular}
\label{m24}\\
\text{List of radial basis functions from Dyn \cite{Dyn89}.}\nonumber
\end{gather}

\subsection{Thin-plate spline or surface spline
functions\label{SbSect_thin_plate_basis}}

For arbitrary dimension $d$ the continuous thin-plate spline (or surface
spline) function $H$ can be defined by
\begin{equation}
H\left(  x\right)  =\left\{
\begin{array}
[c]{ll}%
\left(  -1\right)  ^{s+1}\left\vert x\right\vert ^{2s}\log\left\vert
x\right\vert , & s=1,2,3,\ldots,\\
\left(  -1\right)  ^{\left\lceil s\right\rceil }\left\vert x\right\vert
^{2s}, & s>0,\text{ }s\neq1,2,3\ldots.
\end{array}
\right. \label{m19}%
\end{equation}

Because $H$ is a regular tempered distribution, $\widehat{H}\in S^{\prime}$.
In fact, from equations 23 and 24 of Dyn \cite{Dyn89} we have that, as
distributions,%
\begin{equation}
\widehat{H}=e\left(  s\right)  \left\vert \cdot\right\vert ^{-2s-d}\text{
}on\text{ }\mathbb{R}^{d}\setminus0,\label{m23}%
\end{equation}

where $e\in C_{BP}^{\infty}$. Specifically
\begin{equation}
e\left(  s\right)  =\left\{
\begin{array}
[c]{ll}%
\left(  -1\right)  ^{s+1}c^{\prime}\left(  2s\right)  , & s=1,2,3,\ldots,\\
\left(  -1\right)  ^{\left\lceil s\right\rceil }c\left(  2s\right)  , &
s>0,\text{ }s\neq1,2,3\ldots,
\end{array}
\right. \label{m25}%
\end{equation}

where $c$ and $c^{\prime}$ are defined by
\begin{equation}
c\left(  t\right)  =\pi^{d/2}2^{t+d}\Gamma\left(  \frac{t+d}{2}\right)
/\Gamma\left(  -\frac{t}{2}\right)  ,\text{\quad}c^{\prime}=\frac{dc}%
{dt}.\label{m26}%
\end{equation}

The next theorem will require that $\widehat{H}\left(  \xi\right)  >0$ when
$\xi\neq0$ so the following lemma will be required.

\begin{lemma}
\label{Lem_e(s)_pos}$e\left(  s\right)  >0$ when $s>0$.
\end{lemma}

\begin{proof}
The proof is an application of the reflection formula
\[
\frac{1}{\Gamma\left(  -x\right)  }=-\frac{x\sin\pi x}{\pi}\Gamma\left(
x\right)  .
\]

First suppose that $s>0$ and $s\neq1,2,3\ldots$. Then from \ref{m25} and
\ref{m26}%

\begin{align*}
e\left(  s\right)  =\left(  -1\right)  ^{\left\lceil s\right\rceil }c\left(
2s\right)   &  =\left(  -1\right)  ^{\left\lceil s\right\rceil }\pi
^{d/2}2^{2s+d}\Gamma\left(  s+\frac{d}{2}\right)  /\Gamma\left(  -s\right) \\
&  =\left(  -1\right)  ^{\left\lceil s\right\rceil }\pi^{d/2}2^{2s+d}%
\Gamma\left(  s+\frac{d}{2}\right)  \Gamma\left(  s\right)  s\frac{-\sin\pi
s}{\pi}\\
&  =\pi^{d/2}2^{2s+d}\Gamma\left(  s+\frac{d}{2}\right)  \Gamma\left(
s\right)  s\frac{\left(  -1\right)  ^{1+\left\lceil s\right\rceil }\sin\pi
s}{\pi}\\
&  >0.
\end{align*}

Next assume $s=1,2,3,\ldots$. Then
\begin{align*}
e\left(  s\right)  =\left(  -1\right)  ^{s+1}c^{\prime}\left(  2s\right)  = &
\left(  -1\right)  ^{s+1}\frac{1}{2}\frac{d}{ds}c\left(  2s\right) \\
&  =\frac{\left(  -1\right)  ^{s+1}}{2}\frac{d}{ds}\left(  \pi^{d/2}%
2^{2s+d}\Gamma\left(  s+\frac{d}{2}\right)  \Gamma\left(  s\right)
s\frac{-\sin\pi s}{\pi}\right) \\
&  =\frac{\left(  -1\right)  ^{s+1}}{2}\pi^{d/2}2^{2s+d}\Gamma\left(
s+\frac{d}{2}\right)  \Gamma\left(  s\right)  s\,\left(  -\cos\pi s\right)
+\\
&  \qquad+\frac{\left(  -1\right)  ^{s+1}}{2}\frac{d}{ds}\left(  \pi
^{d/2}2^{2s+d}\Gamma\left(  s+\frac{d}{2}\right)  \Gamma\left(  s\right)
s\right)  \frac{-\sin\pi s}{\pi}\\
&  =\pi^{d/2}2^{2s+d}\Gamma\left(  s+\frac{d}{2}\right)  \Gamma\left(
s\right)  \frac{s}{2}\\
&  >0,
\end{align*}

as required.
\end{proof}

\begin{remark}
\label{Rem_Lem_e(s)_pos}$e\left(  s\right)  >0$ when $s>0$ and%
\begin{equation}
e\left(  s\right)  =\left\{
\begin{array}
[c]{ll}%
\pi^{\frac{d}{2}}2^{2s+d-1}\Gamma\left(  s+\frac{d}{2}\right)  \Gamma\left(
s\right)  s, & s=1,2,3,\ldots,\\
\left.
\begin{array}
[c]{l}%
\pi^{\frac{d}{2}}2^{2s+d}\Gamma\left(  s+\frac{d}{2}\right)  \Gamma\left(
s\right)  s\frac{\left(  -1\right)  ^{1+\left\lceil s\right\rceil }\sin\pi
s}{\pi},\\
\left(  -1\right)  ^{\left\lceil s\right\rceil }\pi^{\frac{d}{2}}%
2^{2s+d}\Gamma\left(  s+\frac{d}{2}\right)  /\Gamma\left(  -s\right)
\end{array}
\right\}  , & s>0,\text{ }s\neq1,2,3\ldots.
\end{array}
\right. \label{m29}%
\end{equation}

Specifically%
\begin{align*}
e\left(  \frac{1}{2}\right)   & =\pi^{\frac{d}{2}}2^{1+d}\Gamma\left(
\frac{1}{2}+\frac{d}{2}\right)  \Gamma\left(  \frac{1}{2}\right)  \frac{1}%
{2}\frac{\left(  -1\right)  ^{1+\left\lceil \frac{1}{2}\right\rceil }\sin
\frac{\pi}{2}}{\pi}\\
& =\pi^{\frac{d}{2}-1}2^{d}\Gamma\left(  \frac{d+1}{2}\right)  \Gamma\left(
\frac{1}{2}\right)  =\pi^{\frac{d-1}{2}}2^{d}\Gamma\left(  \frac{d+1}%
{2}\right)  =\\
& =\left\{
\begin{array}
[c]{ll}%
\pi^{\frac{d-1}{2}}2^{d}\left(  \frac{d-1}{2}\right)  !, & d\text{ }odd,\\
\pi^{\frac{d-1}{2}}2^{d}\Gamma\left(  \frac{d+1}{2}\right)  , & d\text{ }even,
\end{array}
\right. \\
& =\left\{
\begin{array}
[c]{ll}%
\pi^{\frac{d-1}{2}}2^{d}\left(  \frac{d-1}{2}\right)  !, & d\text{ }odd,\\
\pi^{\frac{d-1}{2}}2^{d}1\cdot3\cdot5\cdots\left(  d-1\right)  \frac
{\pi^{\frac{1}{2}}}{2^{\frac{d}{2}}}, & d\text{ }even,
\end{array}
\right. \\
& =\left\{
\begin{array}
[c]{ll}%
\pi^{\frac{d-1}{2}}2^{d}\left(  \frac{d-1}{2}\right)  !, & d\text{ }odd,\\
\left(  2\pi\right)  ^{\frac{d}{2}}\left(  d-1\right)  !!, & d\text{ }even,
\end{array}
\right.
\end{align*}

and%
\begin{align*}
e\left(  1\right)   & =\pi^{\frac{d}{2}}2^{d+2}\Gamma\left(  \frac{d}%
{2}+1\right)  \Gamma\left(  1\right)  \frac{1}{2}=\pi^{\frac{d}{2}}%
2^{d+1}\Gamma\left(  \frac{d}{2}+1\right)  =\\
& =\left\{
\begin{array}
[c]{ll}%
\pi^{\frac{d}{2}}2^{d+1}\frac{d}{2}\Gamma\left(  \frac{d}{2}\right)  , &
d\text{ }odd,\\
\pi^{\frac{d}{2}}2^{d+1}\left(  \frac{d}{2}\right)  !, & d\text{ }even.
\end{array}
\right.
\end{align*}

\end{remark}

\begin{theorem}
\label{Thm_surf_spline_ext_basis}Let $H$ be the generalized (radial)
thin-plate function defined by equation \ref{m19}. With reference to Theorem
\ref{Thm_rad_avoid_So2n}, equation \ref{m23} implies $H_{F}\left(  \xi\right)
=e\left(  s\right)  \left\vert \xi\right\vert ^{-2s-d}\in L_{loc}^{1}\left(
\mathbb{R}^{d}\setminus\mathcal{B}\right)  $ where $\mathcal{B}=\left\{
0\right\}  $. For integer order $\theta\geq1$ let $\mathcal{A}=\left\{
0\right\}  $ and define the function $w$ by \ref{m22} i.e.
\begin{equation}
w\left(  \xi\right)  =\frac{1}{\left\vert \xi\right\vert ^{2\theta}%
H_{F}\left(  \xi\right)  }=\frac{1}{e\left(  s\right)  }\left\vert
\xi\right\vert ^{-2\theta+2s+d},\quad\xi\neq0.\label{m16}%
\end{equation}

Then:

\begin{enumerate}
\item $w$ has weight function property W1 for the set $\mathcal{A}$.

\item $w$ also has weight function properties W2.1 and W3.2 for $\theta$ and
$\kappa\geq0$ iff $0\leq\kappa<s<\theta$.

\item If $0\leq\kappa<s<\theta$ then $H$ is a basis function of order $\theta$
generated by $w$.

\item Regarding the smoothness of functions in $X_{w}^{\theta}$,%
\begin{align*}
\max\left\lfloor \kappa\right\rfloor  & =\left\{
\begin{array}
[c]{ll}%
s-1, & s=1,2,3,\ldots,\\
\left\lfloor s\right\rfloor , & s>0,\text{ }s\neq1,2,3\ldots,
\end{array}
\right. \\
& =\left\{
\begin{array}
[c]{ll}%
n-1, & s=n,\\
n, & n<s<n+1,
\end{array}
\right.
\end{align*}

and the smoothness of the basis functions
\[
\max\left\lfloor 2\kappa\right\rfloor =\left\{
\begin{array}
[c]{ll}%
2s-1, & s=1,2,3,\ldots,\\
2\left\lfloor s\right\rfloor , & n<s\leq n+1/2,\text{ }n=0,1,2,\ldots,\\
2\left\lfloor s\right\rfloor +1, & n+1/2<s<n+1,\text{ }n=0,1,2,\ldots.
\end{array}
\right.
\]

The minimum order is%
\[
\min\theta=\left\{
\begin{array}
[c]{ll}%
s+1, & s=1,2,3,\ldots,\\
\left\lceil s\right\rceil , & s>0,\text{ }s\neq1,2,3\ldots.
\end{array}
\right.
\]

\item If $t>2s$ there exists a constant $c_{t}$ such that $\left\vert H\left(
x\right)  \right\vert \leq c_{t}\left(  1+\left\vert x\right\vert \right)
^{t}$ for all $x$.
\end{enumerate}
\end{theorem}

\begin{proof}
\textbf{Parts 1-3} This theorem is an application of Theorem
\ref{Thm_rad_avoid_So2n}. From \ref{m23} $\widehat{H}\left(  \xi\right)
=e\left(  s\right)  \left\vert \xi\right\vert ^{-2s-d}$ so the condition
$\left\vert \cdot\right\vert ^{2\theta}\widehat{H}\in L_{loc}^{1}$ is
satisfied iff $2\theta-2s-d>-d$ i.e.
\begin{equation}
s<\theta.\label{m15}%
\end{equation}

Lemma \ref{Lem_e(s)_pos} implies $H_{F}\in C^{\left(  0\right)  }\left(
\mathbb{R}^{d}\setminus\mathcal{A}\right)  $ and $H_{F}\left(  \xi\right)  >0$
on $\mathbb{R}^{d}\setminus\mathcal{A}$ so that if we define the function $w $
by \ref{m16}, part 1 of Theorem \ref{Thm_rad_avoid_So2n} implies $w$ has
Property \textbf{W1} w.r.t. the set $\mathcal{A}$.\smallskip

Property \textbf{W2.1} requires that $1/w\in L_{loc}^{1}$. But%
\begin{equation}
\frac{1}{w\left(  \xi\right)  }=e\left(  s\right)  \left\vert \xi\right\vert
^{2\theta-2s-d},\label{m20}%
\end{equation}

so that $1/w\in L_{loc}^{1}$ iff $2\theta-2s-d>-d$ i.e. iff $s<\theta$, which
is already implied by \ref{m15}.\smallskip

Property \textbf{W3.2} is true for order $\theta$ and $\kappa$ if there exists
$r_{3}>0$ such that $\int\limits_{\left\vert \cdot\right\vert \geq r_{3}%
}\dfrac{\left\vert \cdot\right\vert ^{2\kappa}}{w\left\vert \cdot\right\vert
^{2\theta}}<\infty$. But from \ref{m20}, $\int\limits_{\left\vert
\cdot\right\vert \geq r_{3}}\dfrac{\left\vert \cdot\right\vert ^{2\kappa}%
}{w\left\vert \cdot\right\vert ^{2\theta}}=\frac{1}{e\left(  s\right)  }%
\int\limits_{\left\vert \cdot\right\vert \geq r_{3}}\frac{1}{\left\vert
\cdot\right\vert ^{2s-2\kappa+d}}$ and this exists iff $2s-2\kappa+d>d$ i.e.
iff
\begin{equation}
\kappa<s.\label{m11}%
\end{equation}

Thus $w$ has properties W1, W2.1 and W3.2 for some $\theta$ and $\kappa$ if
and only if $\theta>s$ and $0\leq\kappa<s<\theta$. Now by part 2 of Theorem
\ref{Thm_rad_avoid_So2n}, $H$ is a basis function of order $\theta$ generated
by $w$.\smallskip

\textbf{Part 4} Follows from the conditions $0\leq\kappa<s<\theta$ and the
values of $s$ in \ref{m19}.\smallskip

\textbf{Part 5} From \ref{m19}
\[
\left(  1+\left\vert x\right\vert \right)  ^{-t}H\left(  x\right)  =\left\{
\begin{array}
[c]{ll}%
\left(  -1\right)  ^{s+1}\frac{\left\vert x\right\vert ^{2s}\log\left\vert
x\right\vert }{\left(  1+\left\vert x\right\vert \right)  ^{t}}, &
s=1,2,3,\ldots,\\
\left(  -1\right)  ^{\left\lceil s\right\rceil }\frac{\left\vert x\right\vert
^{2s}}{\left(  1+\left\vert x\right\vert \right)  ^{t}}, & s>0,\text{ }%
s\neq1,2,3\ldots,
\end{array}
\right.
\]

so that $t>2s$ implies $\left(  1+\left\vert x\right\vert \right)
^{-t}H\left(  x\right)  \rightarrow0$ as $\left\vert x\right\vert
\rightarrow\infty$. Thus $\left(  1+\left\vert x\right\vert \right)
^{-t}\left\vert H\left(  x\right)  \right\vert $\ is bounded.
\end{proof}

\subsection{Shifted thin-plate spline or shifted surface spline functions
\label{SbSect_shft_thin_plate_basis}}

For arbitrary dimension $d$ the shifted thin-plate spline functions $H$ are
defined for $a>0$ by
\begin{equation}
H\left(  x\right)  =\left\{
\begin{array}
[c]{ll}%
\frac{\left(  -1\right)  ^{s+1}}{2}\left(  a^{2}+\left\vert x\right\vert
^{2}\right)  ^{s}\log\left(  a^{2}+\left\vert x\right\vert ^{2}\right)  , &
s=1,2,3,\ldots,\\
\left(  -1\right)  ^{\left\lceil s\right\rceil }\left(  a^{2}+\left\vert
x\right\vert ^{2}\right)  ^{s}, & s>-d/2,\text{ }s\neq1,2,3\ldots.
\end{array}
\right. \label{m44}%
\end{equation}

Now $H\in S^{\prime}$ and from equations 25, 26 and 27 of Dyn \cite{Dyn89}, as
distributions,
\begin{equation}
\widehat{H}\left(  \xi\right)  =\widetilde{e}\left(  s\right)  \widetilde
{K}_{s+d/2}\left(  a\left\vert \xi\right\vert \right)  \left\vert
\xi\right\vert ^{-2s-d}\;on\text{ }\mathbb{R}^{d}\setminus0,\text{
}s>-d/2,\label{m18}%
\end{equation}

where%
\[
\widetilde{e}\left(  s\right)  =\left\{
\begin{array}
[c]{ll}%
\left(  -1\right)  ^{s+1}\widetilde{c}^{\prime}\left(  2s\right)  , &
s=1,2,3,\ldots,\\
\left(  -1\right)  ^{\left\lceil s\right\rceil }\widetilde{c}\left(
2s\right)  , & s>-d/2,\text{ }s\neq1,2,3\ldots,
\end{array}
\right.
\]

with
\[
\widetilde{c}\left(  t\right)  =\left(  2\pi\right)  ^{d/2}2^{\frac{t+2}{2}%
}/\Gamma\left(  -t/2\right)  ,\text{\quad}\widetilde{c}^{\prime}\left(
t\right)  =\dfrac{d\widetilde{c}\left(  t\right)  }{dt},\text{\quad}t>0,
\]

and
\[
\widetilde{K}_{\lambda}\left(  t\right)  =t^{\lambda}K_{\lambda}\left(
t\right)  ,\text{\quad}t,\text{ }\lambda\text{ }real,
\]

where $K_{\lambda}$ is the modified Bessel function of the second kind of
order $\lambda$. $\widetilde{K}_{\lambda}$ has the following properties
\begin{equation}
\widetilde{K}_{\lambda}\in C^{\left(  0\right)  }\left(  \mathbb{R}\right)
;\text{\quad}\widetilde{K}_{\lambda}\left(  t\right)  >0,\text{\quad}%
t\geq0;\text{\quad}\lim\limits_{t\rightarrow\infty}\widetilde{K}_{\lambda
}\left(  t\right)  =0\text{ }exponentially.\label{m14}%
\end{equation}

See for example Abramowitz and Stegun \cite{AbramowStegun70}. The next theorem
will require that $\widehat{H}\left(  \xi\right)  >0$ when $\xi\neq0 $, so we
will need the following lemma.

\begin{lemma}
\label{Lem_etild_Ktild_pos}$\widetilde{e}\left(  s\right)  >0$,$\quad\xi
\in\mathbb{R}^{d}$, $s>-d/2$.
\end{lemma}

\begin{proof}
The proof is very similar to the proof that $e\left(  s\right)  >0$ when $s>0$
(Lemma \ref{Lem_e(s)_pos}). It is again an application of the reflection
formula \ref{m36}.

First suppose that $s>-d/2$ and $s\neq1,2,3\ldots$ Then
\begin{align*}
\widetilde{e}\left(  s\right)  =\left(  -1\right)  ^{\left\lceil s\right\rceil
}\widetilde{c}\left(  2s\right)  =\left(  -1\right)  ^{\left\lceil
s\right\rceil }\left(  2\pi\right)  ^{d/2}2^{s+1}/\Gamma\left(  -s\right)   &
=\left(  -1\right)  ^{\left\lceil s\right\rceil }\left(  2\pi\right)
^{d/2}2^{s+1}\Gamma\left(  s\right)  s\frac{-\sin\pi s}{\pi}\\
&  =\left(  2\pi\right)  ^{d/2}2^{s+1}\Gamma\left(  s\right)  s\frac{\left(
-1\right)  ^{1+\left\lceil s\right\rceil }\sin\pi s}{\pi}\\
&  >0.
\end{align*}

Next assume $s=1,2,3,\ldots$ Then
\begin{align}
e\left(  s\right)  =\left(  -1\right)  ^{s+1}c^{\prime}\left(  2s\right)   &
=\left(  -1\right)  ^{s+1}\frac{1}{2}\frac{d}{ds}c\left(  2s\right)
\nonumber\\
&  =\frac{\left(  -1\right)  ^{s+1}}{2}\frac{d}{ds}\left(  \pi^{d/2}%
2^{2s+d}\Gamma\left(  s+\frac{d}{2}\right)  \Gamma\left(  s\right)
s\frac{-\sin\pi s}{\pi}\right) \nonumber\\
&  =\frac{\left(  -1\right)  ^{s+1}}{2}\pi^{d/2}2^{2s+d}\Gamma\left(
s+\frac{d}{2}\right)  \Gamma\left(  s\right)  s\,\left(  -\cos\pi s\right)
+\nonumber\\
&  \qquad+\frac{\left(  -1\right)  ^{s+1}}{2}\frac{d}{ds}\left(  \pi
^{d/2}2^{2s+d}\Gamma\left(  s+\frac{d}{2}\right)  \Gamma\left(  s\right)
s\right)  \frac{-\sin\pi s}{\pi}\nonumber\\
&  =\frac{1}{2}\pi^{d/2}2^{2s+d}\Gamma\left(  s+\frac{d}{2}\right)
\Gamma\left(  s\right)  s\label{p77}\\
&  >0,\nonumber
\end{align}

as required.
\end{proof}

\begin{theorem}
\label{Thm_shift_surf_spline_ext_basis}Let $H$ be the generalized (radial)
shifted thin-plate function defined by equation \ref{m44}. With reference to
Theorem \ref{Thm_rad_avoid_So2n}, equation \ref{m18} together with \ref{m14}
imply $H_{F}\left(  \xi\right)  =\widetilde{e}\left(  s\right)  \widetilde
{K}_{s+d/2}\left(  a\left\vert \xi\right\vert \right)  \left\vert
\xi\right\vert ^{-2\theta+2s+d}\in L_{loc}^{1}\left(  \mathbb{R}^{d}%
\setminus\mathcal{B}\right)  $ where $\mathcal{B}=\left\{  0\right\}  $. For
integer order $\theta\geq1$ let $\mathcal{A}=\left\{  0\right\}  $ and define
the function $w$ by \ref{m22} i.e.
\begin{equation}
w\left(  \xi\right)  =\frac{1}{\left\vert \xi\right\vert ^{2\theta}%
H_{F}\left(  \xi\right)  }=\frac{1}{\widetilde{e}\left(  s\right)
\widetilde{K}_{s+d/2}\left(  a\left\vert \xi\right\vert \right)  }\left\vert
\xi\right\vert ^{-2\theta+2s+d},\quad s>-d/2.\label{m33}%
\end{equation}

Then:

\begin{enumerate}
\item $w$ has weight function property W1 w.r.t. the set $\mathcal{A}$.

\item $w$ also has weight function properties W2.1 and W3.2 for $\theta$ and
all $\kappa\geq0$ iff $-d/2<s<\theta$. ?? \textbf{WHAT IF} $s\geq\theta$? ??

\item If $-d/2<s<\theta$ then $H$ is a basis function of order $\theta$
generated by $w$, and $H\in C_{BP}^{\infty}$.

\item If $t>2s$ there exists a constant $c_{t}$ such that $\left\vert H\left(
x\right)  \right\vert \leq c_{t}\left(  1+\left\vert x\right\vert \right)
^{t}$ for all $x$.
\end{enumerate}
\end{theorem}

\begin{proof}
This theorem is an application of Theorem \ref{Thm_rad_avoid_So2n}. From
\ref{m18},

$\widehat{H}\left(  \xi\right)  =\widetilde{e}\left(  s\right)  \widetilde
{K}_{s+d/2}\left(  a\left\vert \xi\right\vert \right)  \left\vert
\xi\right\vert ^{-2s-d},$\ $\xi\neq0,$ $s>-d/2$, so the condition $\left\vert
\cdot\right\vert ^{2\theta}\widehat{H}\in L_{loc}^{1}$ is satisfied iff
$2\theta-2s-d>-d$ and $s>-d/2$ i.e.
\begin{equation}
-d/2<s<\theta.\label{m17}%
\end{equation}

Lemma \ref{Lem_etild_Ktild_pos} and \ref{m14} imply $H_{F}\in C^{\left(
0\right)  }\left(  \mathbb{R}^{d}\setminus\mathcal{A}\right)  $ and
$H_{F}\left(  \xi\right)  >0$ on $\mathbb{R}^{d}\setminus\mathcal{A}$ so that
if we define the function $w$ by \ref{m33}, part 1 of Theorem
\ref{Thm_rad_avoid_So2n} implies $w$ has \textbf{Property W1} w.r.t. the set
$\mathcal{A}$.\smallskip

\textbf{Property W2.1} requires that $1/w\in L_{loc}^{1}$. But%
\begin{equation}
\frac{1}{w\left(  \xi\right)  }=\widetilde{e}\left(  s\right)  \widetilde
{K}_{s+d/2}\left(  a\left\vert \xi\right\vert \right)  \left\vert
\xi\right\vert ^{2\theta-2s-d},\;s>-d/2\label{m28}%
\end{equation}

and \ref{m14} implies $\widetilde{K}_{\lambda}\in C^{\left(  0\right)
}\left(  \mathbb{R}\right)  $ and $\widetilde{K}_{\lambda}\left(  t\right)
>0$ for $t\geq0$, so that $1/w\in L_{loc}^{1}$ iff $\left\vert \xi\right\vert
^{2\theta-2s-d}\in L_{loc}^{1}$ iff $2\theta-2s-d>-d$ i.e. $s<\theta$, which
is again constraint \ref{m17}.\smallskip

\textbf{Property W3.2} is true for order $\theta$ and $\kappa$ if there exists
$r_{3}>0$ such that $\int_{\left\vert \cdot\right\vert \geq r_{3}}%
\dfrac{\left\vert \cdot\right\vert ^{2\kappa}}{w\left\vert \cdot\right\vert
^{2\theta}}<\infty$. But by \ref{m28}, when $r_{3}>0$%
\[
\int_{\left\vert \cdot\right\vert \geq r_{3}}\dfrac{\left\vert \cdot
\right\vert ^{2\kappa}}{w\left\vert \cdot\right\vert ^{2\theta}}=\widetilde
{e}\left(  s\right)  \int_{\left\vert \xi\right\vert \geq r_{3}}\widetilde
{K}_{s+d/2}\left(  a\left\vert \xi\right\vert \right)  \left\vert
\xi\right\vert ^{2\kappa-2s-d}d\xi,
\]
and since by \ref{m14}, $\widetilde{K}_{s+d/2}\left(  a\left\vert
\xi\right\vert \right)  $ is continuous and $\lim\limits_{\left\vert
\xi\right\vert \rightarrow\infty}\widetilde{e}\left(  s\right)  \widetilde
{K}_{s+d/2}\left(  a\left\vert \xi\right\vert \right)  \rightarrow0$
exponentially, this integral always exists. Thus $w$ has properties W1, W2.1
and W3.2 for some $\theta$ and any $\kappa\geq0$ if and only \ref{m17} is
satisfied. Now by part 2 of Theorem \ref{Thm_rad_avoid_So2n}, $H$ is a basis
function of order $\theta$ generated by $w$ and by Theorem
\ref{Thm_basis_smth_W3.2_r3_pos}, $H\in C_{BP}^{\infty}$.

From \ref{m44}
\[
\left(  1+\left\vert x\right\vert \right)  ^{-t}H\left(  x\right)  =\left\{
\begin{array}
[c]{ll}%
\left(  -1\right)  ^{s+1}\frac{\left(  a^{2}+\left\vert \xi\right\vert
^{2}\right)  ^{s}\log\left(  a^{2}+\left\vert \xi\right\vert ^{2}\right)
}{\left(  1+\left\vert x\right\vert \right)  ^{t}}, & s=1,2,3,\ldots,\\
\left(  -1\right)  ^{\left\lceil s\right\rceil }\frac{\left(  a^{2}+\left\vert
\xi\right\vert ^{2}\right)  ^{s}}{\left(  1+\left\vert x\right\vert \right)
^{t}}, & s>0,\text{ }s\neq1,2,3\ldots,
\end{array}
\right.
\]

so that $t>2s$ implies $\left(  1+\left\vert x\right\vert \right)
^{-t}H\left(  x\right)  \rightarrow0$ as $\left\vert x\right\vert
\rightarrow\infty$. Thus $\left(  1+\left\vert x\right\vert \right)
^{-t}\left\vert H\left(  x\right)  \right\vert $\ is bounded.
\end{proof}

\subsection{Multiquadric and inverse multiquadric basis functions}

In arbitrary dimension the multiquadric function is defined by:
\[
H\left(  \xi\right)  =\left(  a^{2}+\left\vert \xi\right\vert ^{2}\right)
^{1/2},\text{\quad}a>0.
\]

In dimension $d\geq2$ the inverse multiquadric function is defined by:
\[
H\left(  \xi\right)  =\left(  a^{2}+\left\vert \xi\right\vert ^{2}\right)
^{-1/2},\text{\quad}a>0.
\]

Note that these are clearly specific cases of the shifted surface splines
discussed in the previous subsection. These functions were introduced by Hardy
\cite{Hardy71} for geophysical applications. See also the review paper by
Hardy \cite{Hardy90}.

\subsection{The Gaussian}

The Gaussian function is
\begin{equation}
H\left(  x\right)  =\exp\left(  -\left\vert x\right\vert ^{2}\right)
,\text{\quad}x\in\mathbb{R}^{d}.\label{m21}%
\end{equation}

and has Fourier transform
\begin{equation}
\widehat{H}\left(  \xi\right)  =\frac{\sqrt{\pi}}{2}\exp\left(  -\left\vert
\xi\right\vert ^{2}/4\right)  ,\text{\quad}\xi\in\mathbb{R}^{d}.\label{m12}%
\end{equation}

\begin{theorem}
Let $H$ be the Gaussian function defined by equation \ref{m21}. With reference
to Theorem \ref{Thm_rad_avoid_So2n}, equation \ref{m12} implies
\begin{equation}
H_{F}\left(  \xi\right)  =\widehat{H}\left(  \xi\right)  \in L_{loc}%
^{1}\left(  \mathbb{R}^{d}\right)  ,\label{m13}%
\end{equation}

so $\mathcal{B}=\left\{  {}\right\}  $. For integer order $\theta\geq1$ let
$\mathcal{A}=\left\{  0\right\}  $ and define the function $w$ by \ref{m22}
i.e.
\begin{equation}
w\left(  \xi\right)  =\frac{1}{\left\vert \xi\right\vert ^{2\theta}%
H_{F}\left(  \xi\right)  }=\frac{2}{\sqrt{\pi}}\left\vert \xi\right\vert
^{-2\theta}\exp\left(  \left\vert \xi\right\vert ^{2}/4\right)  ,\quad
\xi\notin\mathcal{A}.\label{m34}%
\end{equation}

Then:

\begin{enumerate}
\item $w$ has weight function property W1 w.r.t. the set $\mathcal{A}$.

\item $w$ also has weight function properties W2.1 and W3.2 for all $\theta$
and $\kappa\geq0$.

\item $H$ is a basis function of order $\theta$ generated by $w$.

\item For any real $t$ there exists a constant $c_{t}$ such that $\left\vert
H\left(  x\right)  \right\vert \leq c_{t}\left(  1+\left\vert x\right\vert
\right)  ^{t}$ for all $x$.
\end{enumerate}
\end{theorem}

\begin{proof}
This theorem is an application of Theorem \ref{Thm_rad_avoid_So2n}. From
\ref{m12} the condition $\left\vert \cdot\right\vert ^{2\theta}\widehat{H}\in
L_{loc}^{1}$ is always satisfied. Equations \ref{m12} and \ref{m13} imply
$H_{F}\in C^{\left(  0\right)  }\left(  \mathbb{R}^{d}\setminus\mathcal{A}%
\right)  $ and $H_{F}\left(  \xi\right)  >0$ on $\mathbb{R}^{d}\setminus
\mathcal{A}$ so that if we define the function $w$ by \ref{m16}, part 1 of
Theorem \ref{Thm_rad_avoid_So2n} means that $w$ has \textbf{Property W1}
w.r.t. the set $\mathcal{A}=\left\{  0\right\}  $.\smallskip

\textbf{Property W2.1} requires that $1/w\in L_{loc}^{1}$. But%
\begin{equation}
\frac{1}{w\left(  \xi\right)  }=\frac{2}{\sqrt{\pi}}\left\vert \xi\right\vert
^{2\theta}\exp\left(  -\left\vert \xi\right\vert ^{2}/4\right)  ,\label{m10}%
\end{equation}

so this is clearly true.\smallskip

\textbf{Property W3.2} is true for given order $\theta$ and $\kappa$ if there
exists $r_{3}>0$ such that $\int\limits_{\left\vert \cdot\right\vert \geq
r_{3}}\dfrac{\left\vert \cdot\right\vert ^{2\kappa}}{w\left\vert
\cdot\right\vert ^{2\theta}}<\infty$, but this is obviously true from
\ref{m10}.

Thus $w$ has properties W1, W2.1 and W3.2 for any order $\theta$ and $\kappa$
and by part 2 of Theorem \ref{Thm_rad_avoid_So2n}, $H$ is a basis function of
order $\theta$ generated by $w$.

Finally, it is clear that $\left(  1+\left\vert x\right\vert \right)
^{-t}\left\vert H\left(  x\right)  \right\vert $\ is bounded for any real $t$..
\end{proof}

\section{Examples: \textit{non-radial basis functions} generated by weight
functions in W3.2}

\subsection{Fundamental solutions of homogeneous elliptic differential
operators of even order}

In Section 2 of \cite{Dyn89} Dyn describes a large class of non-radial
functions on $\mathbb{R}^{d}$ which are positive definite, namely the
fundamental solutions of homogeneous elliptic differential operators of even
order $2\theta$, where $2\theta>d$. In the next result we will use the weight
function theory of Theorem \ref{Thm_homog_avoid_So2n} to prove that these
fundamental solutions are scalar multiples of basis functions of order
$\theta$ generated by a weight function.

\begin{theorem}
\label{Thm_fund_soln_of_p(D)}Suppose $p$ is a homogeneous polynomial of degree
$2\theta$ on $\mathbb{R}^{d}$, where $\theta\geq1$ and $2\theta>d$. Further,
suppose that $p\left(  x\right)  >0$ if $x\neq0$, and that $f\in S^{\prime}$
is a fundamental solution of the differential operator $p\left(  D\right)  $. Then:

\begin{enumerate}
\item The function $w$ defined by%
\begin{equation}
w\left(  x\right)  =p\left(  \dfrac{x}{\left\vert x\right\vert }\right)
,\quad x\neq0,\label{m36}%
\end{equation}

has property W1 of a weight function w.r.t. the set $\mathcal{A}=\left\{
0\right\}  $. It also has properties W2.1 and W3.2 for order $\theta$ and
$\kappa$ iff $0\leq2\kappa<2\theta-d$.

\item $H=\left(  -1\right)  ^{\theta}\left(  2\pi\right)  ^{d/2}f$ satisfies
$\left\vert \cdot\right\vert ^{2\theta}\widehat{H}=\frac{1}{w}\in L_{loc}^{1}$
as distributions.

\item $H\in C_{BP}^{\left(  2\theta-d-1\right)  }$ is a basis function of
order $\theta$ generated by the weight function $w$.
\end{enumerate}
\end{theorem}

\begin{proof}
\textbf{Part 1} Since the sphere $\left\vert x\right\vert =1$ is a compact
set
\[
c_{1}=\min\left\{  p\left(  x\right)  :\left\vert x\right\vert =1\right\}
,\text{\quad}c_{2}=\max\left\{  p\left(  x\right)  :\left\vert x\right\vert
=1\right\}  ,
\]

both exist. Also, $c_{1},c_{2}>0$ since $p\left(  x\right)  >0$ if $x\neq0$.
We conclude that
\[
0<c_{1}\leq w\left(  x\right)  \leq c_{2},\text{\quad}x\neq0,
\]

and clearly
\begin{equation}
w\in C_{B}^{\left(  0\right)  }\left(  \mathbb{R}^{d}\setminus0\right)
;\text{\quad}0<\frac{1}{c_{2}}\leq\frac{1}{w\left(  x\right)  }\leq\frac
{1}{c_{1}},\text{\quad}x\neq0,\label{m30}%
\end{equation}

and
\begin{equation}
\frac{1}{p\left(  x\right)  }=\frac{1}{w\left(  x\right)  \left\vert
x\right\vert ^{2\theta}},\text{\quad}x\neq0.\label{m31}%
\end{equation}

Thus $w$ has property W1 of a weight function for the set $\mathcal{A=}%
\left\{  0\right\}  $, as well as property W2.1. Property W3.2 holds for order
$\theta$ and $\kappa$ iff there exists $r_{3}>0$ such that $\int_{\left\vert
\cdot\right\vert \geq r_{3}}\dfrac{\left\vert \cdot\right\vert ^{2\kappa}%
}{w\left\vert \cdot\right\vert ^{2\theta}}<\infty$ which by \ref{m30} is true
iff $2\theta-2\kappa>d$. Thus $w$ satisfies properties W1, W2.1 and W3.2 for
order $\theta$ and $\kappa$ iff $0\leq2\kappa<2\theta-d$.\textbf{\smallskip}

\textbf{Part 2} A tempered distribution $f$ which satisfies the equation
$p\left(  D\right)  f=\delta$ is said to be a fundamental solution of the
differential operator $p\left(  D\right)  $. Taking the Fourier transform and
noting that $p$ is homogeneous of degree $2\theta$, $p\left(  D\right)
f=p\left(  -i\xi\right)  \widehat{f}=\left(  -1\right)  ^{\theta}p\widehat
{f}=\widehat{\delta}=\left(  2\pi\right)  ^{-d/2}$, so that $\left(
-1\right)  ^{\theta}p\widehat{f}=\left(  2\pi\right)  ^{-d/2}$ i.e.
$p\widehat{H}=1$ as tempered distributions.\textbf{\smallskip}

But \ref{m30} implies $\frac{1}{w}\in L_{loc}^{1}$ and since $\frac{1}%
{w}=\frac{\left\vert \cdot\right\vert ^{2\theta}}{w\left\vert \cdot\right\vert
^{2\theta}}=\frac{\left\vert \cdot\right\vert ^{2\theta}}{p}=\left\vert
\cdot\right\vert ^{2\theta}\widehat{H}$ a.e. it follows that $\frac{1}%
{w}=\left\vert \cdot\right\vert ^{2\theta}\widehat{H}$ as
distributions.\textbf{\smallskip}

\textbf{Part 3} Here we will make use of Theorem \ref{Thm_homog_avoid_So2n}.
Now $H\in S^{\prime}$ and by part 2, $\left\vert \cdot\right\vert ^{2\theta
}\widehat{H}=\frac{1}{w}\in L_{loc}^{1}$ as tempered distributions so for
$\mathcal{B}=\left\{  0\right\}  $ we can define $H_{F}\in L_{loc}^{1}\left(
\mathcal{B}\right)  $ by $H_{F}=\widehat{H}$ on $\mathbb{R}^{d}\setminus
\mathcal{B}$. From the proof of part 2, $p\widehat{H}=1$ as tempered
distributions. Thus if we set $\mathcal{A}=\mathcal{B}$ it follows that
$H_{F}=1/p$ on $\mathbb{R}^{d}\setminus\mathcal{A}$ and $H_{F}\in C^{\left(
0\right)  }\left(  \mathbb{R}^{d}\setminus\mathcal{A}\right)  $ and
$H_{F}\left(  \xi\right)  >0$ on $\mathbb{R}^{d}\setminus\mathcal{A}$.
Equation \ref{m31} implies the weight function $w$ of part 1 coincides with
that of \ref{m22} in Theorem \ref{Thm_homog_avoid_So2n}, and so the results of
part 1 allow us to use Theorem \ref{Thm_homog_avoid_So2n} to conclude that
$H\in C_{BP}^{\left(  \left\lfloor 2\kappa\right\rfloor \right)  }$ and $H$ is
a basis function of order $\theta$ generated by $w$. But $0\leq2\kappa
<2\theta-d$ implies $\left\lfloor 2\kappa\right\rfloor \leq2\theta-d-1$ so
$H\in C_{BP}^{\left(  2\theta-d-1\right)  }$.
\end{proof}

\chapter{The basis function interpolant and its convergence to the data
function\label{Ch_Interpol}}

\section{Introduction\label{Sect_intro}}

The goal of this chapter is to derive orders for the pointwise convergence of
an interpolant to its data function as the density of the independent data
increases.\medskip

\textbf{Section by section}:\medskip

\textbf{Section} \ref{Sect_unisolv} The concept of unisolvent data sets and
minimal unisolvent data sets is introduced together with the Lagrangian
interpolation operator $\mathcal{P}$ and the operator $\mathcal{Q}%
=I-\mathcal{P}$ which are defined using a minimal unisolvent subset of the
data.\smallskip

\textbf{Section} \ref{Sect_Q(exp)} Smoothness results and upper bounds are
obtained for the operators $D_{x}^{\gamma}\mathcal{Q}_{x}\left(  e^{i\left(
x,\xi\right)  }\right)  $. These are then used to estimate an upper bound for
$\int\frac{\left\vert D_{x}^{\gamma}\mathcal{Q}_{x}\left(  e^{i\left(
x,\xi\right)  }\right)  \right\vert ^{2}}{w\left(  \xi\right)  \left\vert
\xi\right\vert ^{2\theta}}d\xi$ where $w$ is the weight function.\smallskip

\textbf{Section} \ref{Sect_inverse_Four_Xwth} The properties of the function
$\mathcal{Q}_{x}\left(  e^{i\left(  x,\xi\right)  }\right)  $ are used to
derive an `inverse Fourier transform' theorem which expresses the value of a
function $f\in X_{w}^{\theta}$ in terms of its distribution Fourier transform,
which is a function in $L_{loc}^{1}\left(  \mathbb{R}^{d}\setminus0\right)
$.\smallskip

\textbf{Section} \ref{Sect_Riesz_rep} The data function space $X_{w}^{\theta}$
is endowed with the \textit{Light norm} which is based on a minimal unisolvent
set. We then derive explicit formulas for the Riesz representers of the
evaluation functionals $f\rightarrow D^{\gamma}f\left(  x\right)  $ where
$f\in X_{w}^{\theta}$. These formulas are expressed in terms of a basis
function. The existence of a Riesz representer of $f\rightarrow f\left(
x\right)  $ means that $X_{w}^{\theta}$ is a reproducing kernel Hilbert space
of continuous functions.\smallskip

\textbf{Section} \ref{Sect_mat_basis_reprod_ker} In this section we define the
\textit{unisolvency matrix} and the \textit{cardinal unisolvency matrix}, as
well as the \textit{basis function matrix} and the \textit{reproducing kernel
matrix}. These matrices will be used to construct the matrix equations for the
interpolant as well as the smoothers studied in later chapters.\smallskip

\textbf{Section} \ref{Sect_Wgx} The finite dimensional space $W_{G,X}$ is
introduced and some properties proved. This space will contain the solution to
the interpolation problem of this chapter and to the smoothing problems of
later chapters.\smallskip

\textbf{Section} \ref{Sect_Eval_&_T_ops} The vector-valued evaluation operator
$\widetilde{\mathcal{E}}_{X}f=\left(  f\left(  x^{\left(  k\right)  }\right)
\right)  $ is studied. This operator and its adjoint will be fundamental to
solving the interpolation and the smoothing problems.\smallskip

\textbf{Section} \ref{Sect_interpol_prob} The minimal seminorm and norm
interpolants are defined and shown to have the same basis function solution in
$W_{G,X}$. The concept of a data function is introduced and several matrix
equations for the smoother are derived.\smallskip

\textbf{Section} \ref{Sect_interpol_converg} We obtain estimates for the order
of pointwise convergence of the interpolant to its data function.\smallskip

\textbf{Section} \ref{Sect_better_results} The convergence estimates of the
previous section are improved using Taylor series expansions of distributions.
These results are applied to the thin-plate splines and the shifted thin-plate splines.

I have not included the results of any numerical experiments concerning the
interpolant error estimates.

\section{Unisolvency: sets, bases, operators and matrices\label{Sect_unisolv}}

The concept of unisolvent sets of data is fundamental to the theory of basis
function interpolation for orders $\theta\geq1$. The basic importance of
unisolvency is that it ensures that the interpolation problem has a unique
solution. Using minimal unisolvent sets we will then construct the Lagrange
polynomial interpolation operator $\mathcal{P}$ and the operator
$\mathcal{Q}=I-\mathcal{P}$ as well as the unisolvency matrices.

In Subsection \ref{SbSect_L&W_norm}, following Light and Wayne
\cite{LightWayneX98Weight}, a unisolvent set of points will be used to define
an inner product on the space $X_{w}^{\theta}$. With this inner product we
will show that $X_{w}^{\theta}$ is a reproducing kernel Hilbert space. But
first we need to define unisolvency and some related concepts.

\begin{definition}
\label{Def_unisolv}\textbf{Unisolvent sets and minimal unisolvent sets}

Recall that $P_{\theta-1}$ is the set of polynomials of degree $\theta-1$.

Then a finite (ordered) set of distinct points $X=\left\{  x^{\left(
i\right)  }\right\}  $ is said to be a \textbf{unisolvent set} with respect to
$P_{\theta-1}$ if: $p\in P_{\theta-1}$ and $p\left(  x^{\left(  i\right)
}\right)  =0$ for all $x^{\left(  i\right)  }\in X$ implies $p=0$.

Sometimes we say $X$ is unisolvent of order $\theta$ or that $X$ is $\theta$-unisolvent.

It is well known that any unisolvent set has at least $M=\dim P_{\theta
-1}=\binom{d+\theta-1}{d}$ points, and that any unisolvent set of more than
$M$ points has a unisolvent subset with $M$ points. Consequently, a unisolvent
set with $M$ points is called a \textbf{minimal unisolvent set}.
\end{definition}

It is convenient to introduce cardinal bases for polynomials and permutations together.

\begin{definition}
\label{Def_cardinal_basis}\textbf{Cardinal bases for polynomials,
permutations}:

\begin{enumerate}
\item A basis $\left\{  l_{i}\right\}  _{i=1}^{M}$ for $P_{\theta-1}$ is a
\textbf{cardinal basis} for the minimal unisolvent set $A=\left\{  a^{\left(
i\right)  }\right\}  _{i=1}^{M}$ if $l_{i}\left(  a^{\left(  j\right)
}\right)  =\delta_{i,j}$ and the $l_{i}$ are polynomials with \textbf{real
valued} coefficients.

\item We will call the column vector $\widetilde{l}=\left(  l_{i}\right)  $
the \textbf{cardinal basis (column) vector} and we will write $\widetilde
{l}_{A}$ to make the dependency on $A$ explicit. Evaluation at a point $x$
will be indicated by $\widetilde{l}_{A}\left(  x\right)  $.

\item It is well known that a set is minimally unisolvent iff there exists a
(unique) cardinal basis for the set. Also, re-ordering $A$ re-orders $\left\{
l_{i}\right\}  _{i=1}^{M}$ identically i.e. $\widetilde{l}_{\pi\left(
A\right)  }=\pi\left(  \widetilde{l}_{A}\right)  =\Pi\widetilde{l}_{A}$ where
$\pi$ is a permutation and $\Pi$ is the corresponding permutation matrix. A
\textbf{permutation} $\pi$ can also be thought of as a sequence of indexes
such that if $\left\{  y^{\left(  i\right)  }\right\}  =\pi\left(  X\right)  $
then $y^{\left(  i\right)  }=x^{\left(  \pi\left(  i\right)  \right)  } $.

\item When used as a vector an ordered set is regarded as a column vector by
default. Re-ordering by columns involves left-multiplication by $\Pi$.
Re-ordering a row vector involves right-multiplying by $\Pi^{T}$.

\item The inverse permutation is denoted $\pi^{-1}$ and has permutation matrix
$\Pi^{T}$. Thus $\Pi^{T}\Pi=\Pi\Pi^{T}=I$.
\end{enumerate}
\end{definition}

\begin{proposition}
\label{Prop_Lagr_interp_translat}\textbf{Translations and dilations of
unisolvent sets} Suppose $Q\subset\mathbb{R}^{d}$ and $f:\mathbb{R}%
^{d}\rightarrow\mathbb{C}$. Then we define translations of sets and functions
by: $\tau_{b}Q:=Q+b$ and $\tau_{b}f\left(  x\right)  :=f\left(  \tau
_{-b}\left(  x\right)  \right)  =f\left(  x-b\right)  $. Regarding Lemma
\ref{Lem_int_Lagrange_interpol}:

\begin{enumerate}
\item From Lemma \ref{Lem_unisolv_trans_dilat} below, translations and
dilations of unisolvent sets are unisolvent sets. The $\theta$-unisolvent
subsets of $\tau_{b}\Omega$ are the images of the $\theta$-unisolvent subsets
$A\subset\Omega$ under $\tau_{b}$.

\item $\left\vert \tau_{b}f\right\vert _{w,\theta}=\left\vert f\right\vert
_{w,\theta}$ and $\widetilde{l}_{\tau_{b}A}\left(  \tau_{b}x\right)
=\widetilde{l}_{A}\left(  x\right)  $, $x\in\Omega$. Also
$\operatorname*{diam}\left(  \tau_{b}A\right)  _{\tau_{b}x}%
=\operatorname*{diam}A_{x}$.

\item ?? \textbf{MOVE}? \textbf{regarding} \ref{p947} ?? $h_{\tau_{b}%
X,\tau_{b}\Omega}=h_{X,\Omega}$. ??

For any $b$ we can ?? choose ?? $c_{\tau_{b}\Omega,\theta}=c_{\Omega,\theta}$,
$h_{\tau_{b}\Omega,\theta}=h_{\Omega,\theta}$ and $K_{\tau_{b}\Omega,\theta
}^{\prime}=K_{\Omega,\theta}^{\prime}$.

\item ?? Finally, the unique \textbf{circumradius} has upper bound%
\[
\min\limits_{b}\max\limits_{x\in\overline{\Omega}}\left\vert x-b\right\vert
\leq\operatorname*{diam}\Omega\sqrt{\frac{d}{2\left(  d+1\right)  }}.
\]

\item Any open set contains a minimal unisolvent set.
\end{enumerate}
\end{proposition}

\begin{proof}
\textbf{Part 1} Proved in Lemma \ref{Lem_unisolv_trans_dilat}.\smallskip

\textbf{Part 2} With reference to ??, the first equation is true since
$\left(  \tau_{b}f\right)  _{F}=e^{-ib\xi}f_{F}$. The second ?? FINISH
??\smallskip

\textbf{Part 3} ?? \textbf{MOVE \& FINISH} ?? Suppose $h_{\tau_{b}X,\tau
_{b}\Omega}<h_{\Omega,\theta}$. Then $h_{\tau_{b}X,\tau_{b}\Omega}%
=h_{X,\Omega}<h_{\Omega,\theta}$ and so there exists a minimal unisolvent set
$A\subset X$ satisfying $\operatorname*{diam}A_{x}\leq c_{\Omega,\theta
}h_{X,\Omega}$. This implies that $\operatorname*{diam}\left(  \tau
_{b}A\right)  _{\tau_{b}x}=\operatorname*{diam}A_{x}\leq c_{\Omega,\theta
}h_{X,\Omega}=c_{\Omega,\theta}h_{\tau_{b}X,\tau_{b}\Omega}$.\smallskip

\textbf{Part 4} ?? \textbf{Use Jung's theorem}. Let $b$ be the centre of the
?? unique? ?? circumball of $\Omega$. Then by Jung's theorem [??],
$\max\limits_{x\in\overline{\Omega}}\left\vert \tau_{b}x\right\vert
\leq\operatorname*{diam}\Omega\sqrt{\frac{d}{2\left(  d+1\right)  }}%
$.\smallskip

\textbf{Part 5} An immediate consequence of part 1.
\end{proof}

\subsection{The Lagrange interpolation operator $\mathcal{P}$ and
$\mathcal{Q}=I-\mathcal{P}$\label{SbSect_P_Q}}

Now we have conditions for the continuity of functions in $X_{w}%
^{\mathbb{\theta}}$, we can introduce operators, norms etc. which involve
evaluating these functions at points. Following Light and Wayne, the first
step is to introduce the operators $\mathcal{P}$ and $\mathcal{Q}$, which will
play a pivotal role in the remainder of this document. Light and Wayne
\cite{LightWayneX98Weight} introduced the Lagrange polynomial interpolation
operator $\mathcal{P}$ as equation (7) following Lemma 3.4 as well as the
operator $\mathcal{Q}=I-\mathcal{P}$.

\begin{definition}
\label{Def_Aux_proj_operator}\textbf{The operators }$\mathcal{P}:C^{\left(
0\right)  }\rightarrow P_{\theta-1}$ \textbf{and} $\mathcal{Q}:C^{\left(
0\right)  }\rightarrow C^{\left(  0\right)  }$. $\mathcal{P}$\textbf{\ is the}
\textbf{Lagrange interpolation operator}.

Suppose the set $A=\left\{  a^{\left(  i\right)  }\right\}  _{i=1}^{M}$ is a
minimal unisolvent set with respect to the polynomials $P_{\theta-1}$. Suppose
$\left\{  l_{i}\right\}  _{i=1}^{M}$ is the cardinal basis of $P_{\theta-1}$
with respect to the unisolvent set of points $A$.

Then for any continuous function $f$ define the operators $\mathcal{P}$ and
$\mathcal{Q}$ by
\begin{equation}
\mathcal{P}f=\sum\limits_{i=1}^{M}f\left(  a^{\left(  i\right)  }\right)
l_{i},\qquad\mathcal{Q}f=f-\mathcal{P}f,\label{q82}%
\end{equation}

where the unisolvent set can be used as a subscript when necessary.

Sometimes we will indicate the dependence on the unisolvent set $A$ by
$\mathcal{Q}^{??A}\mathcal{Q}_{??A}$, $\mathcal{P}_{A}$, $l_{A;i}$.
\end{definition}

\begin{theorem}
\label{Thm_P_op_properties}The operators $\mathcal{P}$ and $\mathcal{Q}$ have
the following properties:

\begin{enumerate}
\item $p\in P_{\theta-1}$ implies $\mathcal{P}p=p$ and hence $\mathcal{Q}p=0 $.

\item For all $a^{\left(  i\right)  }\in A$, $\left(  \mathcal{P}f\right)
\left(  a^{\left(  i\right)  }\right)  =f\left(  a^{\left(  i\right)
}\right)  $ and $\left(  \mathcal{Q}f\right)  \left(  a^{\left(  i\right)
}\right)  =0$. The polynomial $\mathcal{P}f$ interpolates the data $\left\{
\text{ }\left(  a^{\left(  i\right)  },f\left(  a^{\left(  i\right)  }\right)
\right)  \,\right\}  _{i=1}^{M}$.

\item $\mathcal{P}^{2}=\mathcal{P}$, $\mathcal{PQ}=\mathcal{QP}=0$ and
$\mathcal{Q}^{2}=\mathcal{Q}$, so $\mathcal{P}$ and $\mathcal{Q}$ are projections.

\item $\operatorname*{null}\mathcal{Q}=\operatorname*{range}\mathcal{P}%
=P_{\theta-1}$.

\item $\mathcal{P}$ and $\mathcal{Q}$ are independent of the order of the
points in $A$.
\end{enumerate}
\end{theorem}

\begin{proof}
\textbf{Part 1} is true since each member of the cardinal basis satisfies
$\mathcal{P}l_{j}=\sum\limits_{i=1}^{M}l_{i}\left(  a^{\left(  i\right)
}\right)  l_{i}=l_{j}$ and \textbf{part 2} is true since $l_{i}\left(
a^{\left(  j\right)  }\right)  =\delta_{i,j}$. Regarding \textbf{part 3},
$\mathcal{P}$ is a projection since
\[
\mathcal{P}^{2}f=\mathcal{P}\left(  \sum\limits_{j=1}^{M}f\left(  a^{\left(
j\right)  }\right)  l_{j}\right)  =\sum\limits_{j=1}^{M}f\left(  a^{\left(
j\right)  }\right)  \mathcal{P}\left(  l_{j}\right)  =\sum\limits_{j=1}%
^{M}f\left(  a^{\left(  j\right)  }\right)  l_{j}=\mathcal{P}f.
\]

Thus $\mathcal{PQ}=\mathcal{P}\left(  I-\mathcal{P}\right)  =\mathcal{P}%
-\mathcal{P}^{2}=0=\left(  I-\mathcal{P}\right)  \mathcal{P}=\mathcal{QP}$,
and so $\mathcal{Q}^{2}=\mathcal{Q}\left(  I-\mathcal{P}\right)  =\mathcal{Q}%
$. Finally \textbf{part 4} is true since $\mathcal{P}+\mathcal{Q}=I$.
\textbf{Part 5} is true\textbf{\ }since from the definition of unisolvency
reordering $A$ reorders $\left\{  l_{i}\right\}  _{i=1}^{M}$ identically. In
terms of permutations (recall Definition \ref{Def_cardinal_basis}), if we
write $\mathcal{P}_{A}$ and define the \textit{vector-valued evaluation
function} by $\widetilde{\mathcal{E}}_{A}f=\left(  f\left(  a^{\left(
i\right)  }\right)  \right)  $ then%
\begin{align*}
\mathcal{P}_{\pi\left(  A\right)  }f\left(  x\right)  =\left(  \widetilde
{\mathcal{E}}_{\pi\left(  A\right)  }f\right)  ^{T}\widetilde{l}_{\pi\left(
A\right)  }\left(  x\right)  =\left(  \Pi\widetilde{\mathcal{E}}_{A}f\right)
^{T}\Pi\widetilde{l}_{A}\left(  x\right)  =\left(  \widetilde{\mathcal{E}}%
_{A}f\right)  ^{T}\Pi^{T}\Pi\widetilde{l}_{A}\left(  x\right)   &  =\left(
\widetilde{\mathcal{E}}_{A}f\right)  ^{T}\widetilde{l}_{A}\left(  x\right) \\
&  =\mathcal{P}_{A}f\left(  x\right)  .
\end{align*}

\end{proof}

\begin{theorem}
\label{Thm_translat_Qf}Noting Definition \ref{Def_Aux_proj_operator} we have:%
\[
\tau_{y}\mathcal{Q}_{??A}f=\mathcal{Q}_{\tau_{y}A}\tau_{y}f,\quad
y\in\mathbb{R}^{d}.
\]

\end{theorem}

\begin{proof}
?? Use Proposition \ref{Prop_Lagr_interp_translat}. ??
\end{proof}

\subsection{Unisolvency matrices\label{SbSect_unisolv_mat}}

This section builds on Section \ref{Sect_unisolv} which introduced the
concepts of unisolvent sets, minimal unisolvent sets and cardinal bases. In
this section we define the unisolvency matrix and the cardinal unisolvency
matrix which can be used to characterize unisolvent sets in a form suitable
for studying the interpolation and smoothing problems of this series of documents.

\begin{definition}
\label{Def_unisolv_matrix_Px}\textbf{Unisolvency matrix\ }and \textbf{cardinal
unisolvency matrix}.

Suppose $\left\{  p_{i}\right\}  _{i=1}^{M}$ is any basis for the polynomials
$P_{\theta-1}$ and we have an arbitrary ordered set of distinct points
$X=\left\{  x^{(i)}\right\}  _{i=1}^{N}\subset\mathbb{R}^{d}$. Clearly we must
have $N\geq M=\dim P_{\theta-1}$.

The $N\times M$ \textbf{unisolvency matrix} $P_{X}$ is now defined by
\begin{equation}
P_{X}=\left(  p_{j}\left(  x^{(i)}\right)  \right)  .\label{q01}%
\end{equation}

For example, for numerical work we may use the monomials $\left\{  x^{\alpha
}\right\}  _{\left\vert \alpha\right\vert <\theta}$ as a basis for
$P_{\theta-1}$.

Now suppose the subset $A\subset X$ is minimally unisolvent and has cardinal
basis $\left\{  l_{j}\right\}  _{j=1}^{M}$. Then we use the special notation
\begin{equation}
L_{X}=\left(  l_{j}\left(  x^{(i)}\right)  \right)  ,\label{q04}%
\end{equation}

for the unisolvency matrix which corresponds to the basis $\left\{
l_{j}\right\}  $. We call $L_{X}$ a \textbf{cardinal unisolvency matrix}.
\end{definition}

As well as proving some elementary properties of unisolvency matrices, the
next theorem introduces some matrix notation which will be used often later.
It will be noted that there is an emphasis on characterizing the set
$\operatorname*{null}P_{X}^{T}$. This is because the constraint $P_{X}^{T}v=0$
is part of the definition (Definition \ref{Def_Wgx}) of the finite dimensional
space $W_{G,X}$ which contains the interpolants studied in this document. We
now prove two theorems of properties of the unisolvency matrices. The first
theorem proves results which are true for any minimal unisolvent subset of $X
$. The second theorem assumes the first $M$ points of $X$ constitute a minimal
unisolvent subset.

\begin{theorem}
\label{Thm_Px_properties}\textbf{Properties of unisolvency matrices }Suppose
$X=\left\{  x^{(i)}\right\}  _{i=1}^{N}$ is any set of distinct points in
$\mathbb{R}^{d}$. Then if $M=\dim P_{\theta-1}$:

\begin{enumerate}
\item The null space of the unisolvency matrix $P_{X}$ can be used to
characterize unisolvent sets. In fact, $X$ is unisolvent if and only if
$\operatorname*{null}P_{X}=\left\{  0\right\}  $.\medskip

Now suppose the set $X$ is unisolvent and $A=\left\{  a^{\left(  i\right)
}\right\}  _{i=1}^{M}$ is any minimal unisolvent subset. Then:

\item $P_{X}^{T}\beta=0$ iff $\sum\limits_{j=1}^{N}\beta_{j}p\left(
x^{(j)}\right)  =0$ for all $p\in P_{\theta-1}$.

\item Suppose $P_{X}$ and $Q_{X}$ are the unisolvency matrices generated by
the bases $\widetilde{p}=\left\{  p_{i}\right\}  _{i=1}^{M}$ and
$\widetilde{q}=\left\{  q_{i}\right\}  _{i=1}^{M}$ of $P_{\theta-1}$. Let $R$
be the regular change of basis matrix i.e. $\widetilde{p}=R\widetilde{q}$.
Then%
\[
P_{X}=Q_{X}R^{T}=Q_{X}Q_{A}^{-1}P_{A},
\]

where $P_{A}=\left(  p_{j}\left(  a^{\left(  i\right)  }\right)  \right)
_{i,j=1}^{M}$ and $Q_{A}=\left(  q_{j}\left(  a^{\left(  i\right)  }\right)
\right)  _{i,j=1}^{M}$ are regular. Further%
\begin{equation}
P_{X}=L_{X}P_{A}.\label{q06}%
\end{equation}

\item $P_{X}^{T}\beta=0$ iff $L_{X}^{T}\beta=0$.

\item $\dim\operatorname*{null}P_{X}^{T}=\dim\operatorname*{null}L_{X}%
^{T}=N-M$.
\end{enumerate}
\end{theorem}

\begin{proof}
\textbf{Part 1} If $p\in P_{\theta-1}$ then $p=\sum\limits_{i=1}^{M}\beta
_{i}p_{i}$ for unique $\beta_{i}$. Thus the condition $p\left(  x^{(j)}%
\right)  =0$ for all $j=1,\ldots,N$ is equivalent to $\sum\limits_{i=1}%
^{M}\beta_{i}p_{i}\left(  x^{(j)}\right)  =0$ for all $j=1,\ldots,N$ which by
definition \ref{q01} is equivalent to $P_{X}\beta=0$.\medskip

\textbf{Part 2} Since $0=P_{X}^{T}\beta=\sum\limits_{j=1}^{N}\beta_{j}%
p_{i}\left(  x^{(j)}\right)  $ and $\left\{  p_{i}\right\}  _{i=1}^{M}$ is a
basis for $P_{\theta-1}$, it follows that $\sum\limits_{j=1}^{N}\beta
_{j}p\left(  x^{(j)}\right)  =0$ for any $p\in P_{\theta-1}$. Clearly the
argument is reversible.\medskip

\textbf{Part} \textbf{3} If $\widetilde{p}=\left(  p_{i}\right)  $ and
$\widetilde{q}=\left(  q_{i}\right)  $ are the basis vectors then there is a
regular change of basis matrix $R$ such that $\vec{p}=R\,\widetilde{q}$. Then
$\left(  \widetilde{p}\right)  ^{T}=\left(  \widetilde{q}\right)  ^{T}R^{T}$
so that $\left(  \widetilde{p}\right)  ^{T}\left(  x^{(j)}\right)  =\left(
\widetilde{q}\right)  ^{T}\left(  x^{(j)}\right)  R^{T}$ i.e. $P_{X}%
=Q_{X}R^{T}$ which implies $P_{A}=Q_{A}R^{T}$ and so $P_{X}=Q_{X}R^{T}%
=Q_{X}Q_{A}^{-1}P_{A}$. If we chose $\left(  \widetilde{q}\right)  $ to be the
cardinal basis i.e. $q_{i}=l_{i}$ it would follow from part 1 Definition
\ref{Def_cardinal_basis} that $Q_{A}=\left(  l_{j}\left(  a^{\left(  i\right)
}\right)  \right)  =I$ and so $P_{A}$ would be regular. Thus in general
$Q_{A}$ is regular and $P_{X}=Q_{X}Q_{A}^{-1}P_{A}$. By definition \ref{q04}
choosing $\left(  \widetilde{q}\right)  $ to be the cardinal basis would also
imply $Q_{X}=L_{X}$ and so $P_{X}=L_{X}P_{A}$.\medskip

\textbf{Part 4} By part 3, $P_{X}=L_{X}P_{A}$ and the regularity of $P_{A}$
implies $P_{X}^{T}\beta=0$ iff $L_{X}^{T}\beta=0$.\smallskip

\textbf{Part 5} That $\dim\operatorname*{null}P_{X}^{T}=\dim
\operatorname*{null}L_{X}^{T}$ is clear from part 4. Next observe that
$\dim\operatorname*{null}L_{X}^{T}=N-\operatorname*{rank}L_{X}^{T}%
=N-\operatorname*{rank}L_{X}$. Now recall the discussion of permutation
operators and matrices in Definition \ref{Def_cardinal_basis}. Let $\pi$
permute $X$ so the first $M$ elements of $X$ lie in $A$. Then
$\operatorname*{rank}L_{X}=\operatorname*{rank}L_{\pi\left(  X\right)  }$ and
since the first $M$ rows of $L_{\pi\left(  X\right)  }$ are $l_{i}\left(
a^{\left(  j\right)  }\right)  $ they compose the unit matrix and we can
conclude that $\operatorname*{rank}L_{\pi\left(  X\right)  }=M$ and
$\dim\operatorname*{null}L_{X}^{T}=N-M$, as claimed.
\end{proof}

\begin{theorem}
\label{Thm_Px_properties_2}\textbf{Properties of unisolvency matrices} Suppose
$X=\left\{  x^{(i)}\right\}  _{i=1}^{N}$ is a unisolvent set of points in
$\mathbb{R}^{d}$ and $X_{1}=\left\{  x^{(i)}\right\}  _{i=1}^{M}$ is a minimal
unisolvent subset. Set $X_{2}=\left\{  x^{(i)}\right\}  _{i=M+1}^{N}$. Then:

\begin{enumerate}
\item The cardinal unisolvency matrix w.r.t. $X_{1}$ is $L_{X}=\left(
\begin{array}
[c]{l}%
I_{M}\\
L_{X_{2}}%
\end{array}
\right)  $ where $L_{X_{2}}=\left(  l_{j}\left(  x^{(i)}\right)  \right)  $,
$i>M$. Also%
\[
P_{X}^{T}\beta=0\text{ iff }\beta_{k}=-\sum\limits_{j=M+1}^{N}\beta_{j}%
l_{k}\left(  x^{(j)}\right)  ,\quad k=1,\ldots,M.
\]

\item Suppose $\beta=%
\begin{pmatrix}
\beta^{\prime}\\
\beta^{\prime\prime}%
\end{pmatrix}
$ where $\beta^{\prime}\in\mathbb{R}^{M}$, $\beta^{\prime\prime}\in
\mathbb{R}^{N-M}$. Then%
\[
P_{X}^{T}\beta=0\text{ iff }\beta^{\prime}=-L_{X_{2}}^{T}\beta^{\prime\prime
}\text{ iff }\beta=%
\begin{pmatrix}
-L_{X_{2}}^{T}\\
I_{N-M}%
\end{pmatrix}
\beta^{\prime\prime}.
\]

\item Define the $N\times N$ matrix $L_{X;0}=%
\begin{pmatrix}
L_{X} & O_{N,N-M}%
\end{pmatrix}
$. Then $L_{X;0}L_{X}=L_{X}$ and $\left(  L_{X;0}\right)  ^{2}=L_{X;0}$.
\end{enumerate}
\end{theorem}

\begin{proof}
\textbf{Parts 1 and 2} Since by definition $L_{X}=\left(  l_{j}\left(
x^{(i)}\right)  \right)  $ and by definition $l_{j}\left(  x^{(i)}\right)
=\delta_{i,j}$ when $i\leq M$, we obtain the block form $L_{X}=\left(
\begin{array}
[c]{l}%
I_{M}\\
L_{X_{2}}%
\end{array}
\right)  $. The rest of the theorem follows easily by observing that
$L_{X}^{T}\beta=\beta^{\prime}+L_{X_{2}}^{T}\beta^{\prime\prime}$.

\textbf{Part 3}. $L_{X;0}=%
\begin{pmatrix}
L_{X} & O_{N,M-N}%
\end{pmatrix}
=%
\begin{pmatrix}
I_{M} & O\\
L_{X_{2}} & O_{N-M}%
\end{pmatrix}
$, and so
\[
L_{X;0}L_{X}=%
\begin{pmatrix}
I_{M} & O\\
L_{X_{2}} & O_{N-M}%
\end{pmatrix}
\left(
\begin{array}
[c]{l}%
I_{M}\\
L_{X_{2}}%
\end{array}
\right)  =L_{X}.
\]

Further%
\[
\left(  L_{X;0}\right)  ^{2}=%
\begin{pmatrix}
I_{M} & O\\
L_{X_{2}} & O_{N-M}%
\end{pmatrix}%
\begin{pmatrix}
I_{M} & O\\
L_{X_{2}} & O_{N-M}%
\end{pmatrix}
=%
\begin{pmatrix}
I_{M} & O\\
L_{X_{2}} & O_{N-M}%
\end{pmatrix}
=L_{X;0}.
\]

\end{proof}

\section{Properties of the function $\mathcal{Q}_{x}\left(  e^{ix\xi}\right)
$\label{Sect_Q(exp)}}

The operators $\mathcal{P}$ and $\mathcal{Q}$ were introduced in the previous
section and defined using a minimal unisolvent set of points in $\mathbb{R}%
^{d}$. In this section we study the functions $D_{x}^{\gamma}\mathcal{Q}%
_{x}\left(  e^{ix\xi}\right)  $ and obtain smoothness results and upper
bounds, both near the origin and near infinity. Finally, the function
$D_{x}^{\gamma}\mathcal{Q}_{x}\left(  e^{ix\xi}\right)  $ is related to the
weight function by an upper bound for $\int\frac{\left\vert D_{x}^{\gamma
}\mathcal{Q}_{x}\left(  e^{i\left(  x,\cdot\right)  }\right)  \right\vert
^{2}}{w\left\vert \cdot\right\vert ^{2\theta}}$. These results will be used in
the next section to derive an `inverse Fourier transform' theorem which
expresses the value of a function $f\in X_{w}^{\theta}$ in terms of its
Fourier transform on $\mathbb{R}^{d}\setminus0$.

The next theorem derives some smoothness properties of the function
$\mathcal{Q}_{x}\left(  e^{ix\xi}\right)  $.

\begin{theorem}
\label{Thm_Q(exp)_loc_property}Suppose $\theta\geq1$ is an integer and $M=\dim
P_{\theta-1}$. Let $\left\{  l_{i}\right\}  _{i=1}^{M}$ be the cardinal basis
polynomials of $P_{\theta-1}$ associated with the minimal unisolvent set
$A=\left\{  a^{\left(  i\right)  }\right\}  _{i=1}^{M}$. Then
\[
\mathcal{Q}_{x}\left(  e^{ix\xi}\right)  =e^{ix\xi}-\sum\limits_{i=1}^{M}%
l_{i}(x)e^{ia^{\left(  i\right)  }\xi},\text{\quad}\xi,x\in\mathbb{R}^{d},
\]

and $\mathcal{Q}_{x}\left(  e^{ix\xi}\right)  $ has the following properties.

\begin{enumerate}
\item $\mathcal{Q}_{x}\left(  e^{ix\xi}\right)  =0$ when $x\in A$.

\item If $\left\vert \beta\right\vert <\theta$, then for each $x$, $D_{\xi
}^{\beta}\mathcal{Q}_{x}\left(  e^{ix\xi}\right)  \in C_{\emptyset,\theta
}^{\infty}\cap C_{B}^{\infty}$.

\item If $\left\vert \beta\right\vert <\theta$, then for each $x$,
$D_{x}^{\beta}\mathcal{Q}_{x}\left(  e^{ix\xi}\right)  \in C_{\emptyset
,\theta}^{\infty}\cap C_{BP}^{\infty}$.

\item If $\left\vert \alpha\right\vert =\theta$ and $\left\vert \beta
\right\vert <\theta$, then for each $x$, $\left(  i\xi\right)  ^{\alpha}%
D_{x}^{\beta}\mathcal{Q}_{x}\left(  e^{ix\xi}\right)  \in C_{\emptyset
,2\theta}^{\infty}\cap C_{BP}^{\infty}$.
\end{enumerate}
\end{theorem}

\begin{proof}
\textbf{Part 1} Follows directly from the fact that $l_{i}\left(  a_{\beta
}\right)  =\delta_{\alpha,\beta}$.\smallskip

\textbf{Part 2} Clearly, for each $x$, $\mathcal{Q}_{x}\left(  e^{-ix\xi
}\right)  \in C^{\infty}$. Since
\begin{equation}
D_{\xi}^{\beta}\mathcal{Q}_{x}\left(  e^{-ix\xi}\right)  =\left(  -ix\right)
^{\beta}e^{-ix\xi}-\sum\limits_{i=1}^{M}\left(  -ia^{\left(  i\right)
}\right)  ^{\beta}e^{-ia^{\left(  i\right)  }\xi}l_{i}(x),\label{q74}%
\end{equation}

when $\xi=0$%
\[
D_{\xi}^{\beta}\mathcal{Q}_{x}\left(  e^{-ix\xi}\right)  =\left(  -i\right)
^{\left\vert \beta\right\vert }\left(  x^{\beta}-\sum\limits_{i=1}^{M}\left(
a^{\left(  i\right)  }\right)  ^{\beta}l_{i}(x)\right)  =\left(  -i\right)
^{\left\vert \beta\right\vert }\mathcal{Q}\left(  x^{\beta}\right)  =0,
\]

when $\left\vert \beta\right\vert <\theta$. Hence for each $x$, $D_{\xi
}^{\beta}\mathcal{Q}_{x}\left(  e^{-ix\xi}\right)  \in C_{\emptyset,\theta
}^{\infty}$.

\textbf{Part 3} Suppose $\left\vert \beta\right\vert <\theta$.%
\begin{align*}
D_{x}^{\beta}\mathcal{Q}_{x}\left(  e^{-ix\xi}\right)  =D_{x}^{\beta}\left(
e^{-ix\xi}-\sum\limits_{i=1}^{M}e^{-ia^{\left(  i\right)  }\xi}l_{i}%
(x)\right)   &  =\left(  -i\xi\right)  ^{\beta}e^{-ix\xi}-\sum\limits_{i=1}%
^{M}e^{-ia^{\left(  i\right)  }\xi}\left(  D^{\beta}l_{i}\right)  (x)\\
&  =e^{-ix\xi}\left(  \left(  -i\xi\right)  ^{\beta}-\sum\limits_{i=1}%
^{M}e^{i\left(  x-a^{\left(  i\right)  }\right)  \xi}\left(  D^{\beta}%
l_{i}\right)  (x)\right)  .
\end{align*}

Again noting Theorem \ref{Thm_product_of_Co,k_funcs}, since $e^{-ix\xi}$ is in
$C_{\emptyset,0}^{\infty}$ as a function of $\xi$, proving that the second
factor of the last term is in $C_{\emptyset,\theta}^{\infty}$ for any fixed
$x$ proves the result. In fact, if $\left\vert \gamma\right\vert <\theta$ and
$\gamma\neq\beta$, at $\xi=0$%
\begin{align*}
D_{\xi}^{\gamma}\left(  \left(  -i\xi\right)  ^{\beta}-\sum\limits_{i=1}%
^{M}e^{i\left(  x-a^{\left(  i\right)  }\right)  \xi}\left(  D^{\beta}%
l_{i}\right)  (x)\right)   & =-i^{\left\vert \gamma\right\vert }%
\sum\limits_{i=1}^{M}\left(  x-a^{\left(  i\right)  }\right)  ^{\gamma}\left(
D^{\beta}l_{i}\right)  (x)\\
& =-i^{\left\vert \gamma\right\vert }\left\{  D_{x}^{\beta}\mathcal{P}%
_{x}\left[  \left(  y-x\right)  ^{\gamma}\right]  \right\}  _{y=x}\\
& =-i^{\left\vert \gamma\right\vert }\left\{  D_{x}^{\beta}\left[  \left(
y-x\right)  ^{\gamma}\right]  \right\}  _{y=x}\\
& =0,
\end{align*}

since $\gamma\neq\beta$. There remains the case of $\gamma=\beta$, but this
follows easily using the same technique.\medskip

\textbf{Part 4} Suppose $\left\vert \alpha\right\vert =\theta$ and $\left\vert
\beta\right\vert <\theta$. First observe that by theorem
\ref{Thm_product_of_Co,k_funcs}, $\left(  i\xi\right)  ^{\alpha}\in
C_{\emptyset,\theta}^{\infty}$ and since it was shown in part 3 that for each
$x$, $D_{x}^{\beta}\mathcal{Q}_{x}\left(  e^{-ix\xi}\right)  \in
C_{\emptyset,\theta}^{\infty}$, this part is proved.
\end{proof}

To prove the next result we will require the following lemma concerning
differentiation under the integral sign. This result is Proposition 7.8.4 of
Malliavin \cite{Malliavin95}.

\begin{lemma}
\label{Lem_diff_under_integral}Suppose $f:\mathbb{R}^{m+n}\rightarrow
\mathbb{C}$ and we write $f\left(  \xi,x\right)  $ where $\xi\in\mathbb{R}%
^{m}$ and $x\in\mathbb{R}^{n}$. Suppose that:

\begin{enumerate}
\item For each $\xi$, $f\left(  \xi,\cdot\right)  \in C^{\left(  k\right)
}\left(  \mathbb{R}^{n}\right)  $.

\item For each $x$, $\int\left\vert D_{\xi}^{\alpha}f\left(  \xi,x\right)
\right\vert d\xi<\infty$ for $\left\vert \alpha\right\vert \leq k$.
\end{enumerate}

Then when $\left\vert \alpha\right\vert \leq k$ we have
\[
D_{x}^{\alpha}\int f\left(  \xi,x\right)  d\xi=\int D_{x}^{\alpha}f\left(
\xi,x\right)  d\xi,
\]

and $\int f\left(  \xi,\cdot\right)  d\xi\in C^{\left(  k\right)  }\left(
\mathbb{R}^{n}\right)  $.
\end{lemma}

In the next lemma we elucidate some of the structure of the function
$\mathcal{Q}_{x}\left(  e^{ix\xi}\right)  $ studied in the previous theorem.

??? WHAT ABOUT USING LEMMA \ref{Lem_Qf_estim_multipt} below? USE Lemma
\ref{Lem_Qf_estim_multipt} ??

\begin{lemma}
\label{Lem_estim_Q(exp)_multipt}???
\end{lemma}

\begin{proof}
?? From Lemma \ref{Lem_Qf_estim_multipt},
\[
\left\vert \mathcal{Q}f\left(  x\right)  \right\vert \leq\left(
\sum\limits_{k=1}^{M}\left\vert l_{k}\left(  x\right)  \right\vert \right)
\max_{k=1}^{M}\left\vert \left(  \mathcal{R}_{n+1}f\right)  \left(
x,a^{\left(  k\right)  }-x\right)  \right\vert ,\quad x\in\mathbb{R}^{d},
\]

where $\mathcal{R}_{n+1}f$ is the single point Taylor series remainder ??
\ref{a53}.
\end{proof}

\begin{lemma}
\label{Lem_estim_Q(exp)}Suppose the operator $\mathcal{Q}$ is generated by the
minimal $\theta$-unisolvent set $A=\left\{  a^{\left(  i\right)  }\right\}
_{i=1}^{M}$ and the cardinal basis $\left\{  l_{i}\right\}  _{i=1}^{M}$. Then:

\begin{enumerate}
\item
\begin{equation}
\mathcal{Q}_{x}\left(  e^{ix\xi}\right)  =\sum_{\left\vert \alpha\right\vert
=\theta}v_{\alpha}\left(  x,\xi\right)  \frac{\left(  i\xi\right)  ^{\alpha}%
}{\alpha!},\label{q91}%
\end{equation}

where
\[
\nu_{\alpha}\left(  x,\xi\right)  =\mathcal{Q}_{x}\left(  x^{\alpha}%
\mu_{\theta}\left(  \xi x\right)  \right)  ,
\]

and (see ?? )%
\[
\mu_{\theta}\left(  t\right)  =\frac{e^{it}}{\left(  \theta-1\right)  !}%
\int_{0}^{1}e^{-ist}s^{\theta-1}ds,\text{\quad}t\in\mathbb{R}.
\]

\item $\mu_{\theta}\in C_{B}^{\infty}\left(  \mathbb{R}\right)  $ and
$\left\Vert D^{k}\mu_{\theta}\right\Vert _{\infty}<\frac{1}{\theta!}$ for all
$k$.

\item For each $x$, $v_{\alpha}\left(  x,\cdot\right)  \in C_{BP}^{\infty
}\left(  \mathbb{R}^{d}\right)  $ and for each $\xi$, $v_{\alpha}\left(
\cdot,\xi\right)  \in C_{BP}^{\infty}\left(  \mathbb{R}^{d}\right)  $.

Also, for all $\left\vert \gamma\right\vert <\theta$ and $\left\vert
\alpha\right\vert =\theta$,%
\[
\left\vert D_{x}^{\gamma}v_{\alpha}\left(  x,\xi\right)  \right\vert \leq
\frac{2^{\left\vert \gamma\right\vert }}{\gamma!}\left(  1+\xi_{+}\right)
^{\gamma}\left\vert x\right\vert ^{\theta-\left\vert \gamma\right\vert
}\left(  1+\left\vert x\right\vert \right)  ^{\left\vert \gamma\right\vert
}+\frac{1}{\theta!}\sum_{i=1}^{M}\tfrac{\left(  a_{+}^{\left(  i\right)
}\right)  ^{\alpha}}{\alpha!}\left\vert D^{\gamma}l_{i}\left(  x\right)
\right\vert ,\text{\quad}x,\xi\in\mathbb{R}^{d}.
\]

\end{enumerate}
\end{lemma}

\begin{proof}
\textbf{Parts 1 and 2} We start by expanding $e^{it},$ $t\in\mathbb{R}$ about
zero using Taylor's theorem with remainder, as given in Appendix
\ref{Sect_taylor_expansion}, and then substitute $t=x\xi$ to obtain
\begin{equation}
e^{ix\xi}=p_{\theta-1}\left(  ix\xi\right)  +\frac{\left(  ix\xi\right)
^{\theta}e^{ix\xi}}{\left(  \theta-1\right)  !}\int_{0}^{1}e^{-isx\xi
}s^{\theta-1}ds,\label{p937}%
\end{equation}

where $p_{\theta-1}\in P_{\theta-1}$. For notational compactness set
\[
\mu_{\theta}\left(  t\right)  =\frac{e^{it}}{\left(  \theta-1\right)  !}%
\int_{0}^{1}e^{-ist}s^{\theta-1}ds,\text{\quad}t\in\mathbb{R}.
\]

Next observe that $\mu_{\theta}\in C_{BP}^{\infty}$, since $\mu_{\theta} $ is
$e^{it}$ times the Fourier transform of a bounded $L_{loc}^{1}$ function with
compact support. Further, since the derivatives of the integrand which defines
$\mu_{\theta}$ are $L^{1}$ functions on $\mathbb{R} $, Lemma
\ref{Lem_diff_under_integral} enables us to differentiate under the integral
sign and this implies that all the derivatives are bounded i.e. $\mu_{\theta
}\in C_{B}^{\infty}$. In fact, for any integer $k\geq0$
\begin{align*}
D^{k}\mu_{\theta}\left(  t\right)  =D^{k}\frac{e^{it}}{\left(  \theta
-1\right)  !}\int_{0}^{1}e^{-ist}s^{\theta-1}ds &  =D^{k}\frac{1}{\left(
\theta-1\right)  !}\int_{0}^{1}e^{i\left(  1-s\right)  t}s^{\theta-1}ds\\
&  =\frac{i^{k}}{\left(  \theta-1\right)  !}\int_{0}^{1}e^{i\left(
1-s\right)  t}\left(  1-s\right)  ^{k}s^{\theta-1}ds,
\end{align*}

so that by the elementary properties of the beta function%
\[
D^{k}\mu_{\theta}\left(  0\right)  =\frac{i^{k}}{\left(  \theta-1\right)
!}\frac{k!\left(  \theta-1\right)  !}{\left(  k+\theta\right)  !}=\frac
{i^{k}k!}{\left(  k+\theta\right)  !},
\]

and%
\[
\left\vert D^{k}\mu_{\theta}\left(  t\right)  \right\vert \leq\frac{1}{\left(
\theta-1\right)  !}\int_{0}^{1}\left(  1-s\right)  ^{k}s^{\theta-1}ds<\frac
{1}{\left(  \theta-1\right)  !}\int_{0}^{1}s^{\theta-1}ds=\frac{1}{\theta!},
\]

and more accurately%
\begin{equation}
\left\vert D^{k}\mu_{\theta}\left(  t\right)  \right\vert \leq\frac{1}{\left(
\theta-1\right)  !}\int_{0}^{1}\left(  1-s\right)  ^{k}s^{\theta-1}ds=\frac
{1}{\left(  \theta-1\right)  !}\frac{k!\left(  \theta-1\right)  !}{\left(
k+\theta\right)  !}=\frac{k!}{\left(  k+\theta\right)  !}.\label{a1.03}%
\end{equation}

We now rewrite \ref{p937} as
\begin{equation}
e^{ix\xi}=p_{\theta-1}\left(  ix\xi\right)  +\left(  ix\xi\right)  ^{\theta
}\mu_{\theta}\left(  x\xi\right)  .\label{q71}%
\end{equation}

But for given $\xi$, $p_{\theta-1}\left(  ix\xi\right)  $ is a polynomial in
$x$ of degree less than $\theta$, so by part 2 of Theorem
\ref{Thm_P_op_properties} $\mathcal{Q}_{x}\left(  p_{\theta-1}\left(
ix\xi\right)  \right)  =0$ for all $x$. Thus
\begin{align}
\mathcal{Q}_{x}\left(  e^{ix\xi}\right)   & =\mathcal{Q}_{x}\left(  \left(
i\xi x\right)  ^{\theta}\mu_{\theta}\left(  x\xi\right)  \right) \label{q70}\\
& =\mathcal{Q}_{x}\left(  \left(  \sum_{\left\vert \alpha\right\vert =\theta
}\frac{x^{\alpha}\left(  i\xi\right)  ^{\alpha}}{\alpha!}\right)  \mu_{\theta
}\left(  x\xi\right)  \right) \nonumber\\
& =\sum_{\left\vert \alpha\right\vert =\theta}\frac{1}{\alpha!}\mathcal{Q}%
_{x}\left(  x^{\alpha}\mu_{\theta}\left(  x\xi\right)  \right)  \left(
i\xi\right)  ^{\alpha}\nonumber\\
& =\sum_{\left\vert \alpha\right\vert =\theta}v_{\alpha}\left(  x,\xi\right)
\left(  i\xi\right)  ^{\alpha},\label{p933}%
\end{align}

where we have defined the functions $\left\{  v_{\alpha}\right\}  _{\left\vert
\alpha\right\vert =\theta}$ by
\begin{equation}
v_{\alpha}\left(  x,\xi\right)  =\frac{1}{\alpha!}\mathcal{Q}_{x}\left(
x^{\alpha}\mu_{\theta}\left(  x\xi\right)  \right)  ,\text{\quad}\left\vert
\alpha\right\vert =\theta.\label{q89}%
\end{equation}
\medskip

\textbf{Part 3} Since $\mu_{\theta}\in C_{B}^{\infty}\left(  \mathbb{R}%
^{1}\right)  $, for each $x$, $v_{\alpha}\left(  x,\cdot\right)  \in
C_{BP}^{\infty}\left(  \mathbb{R}^{d}\right)  $ and for each $\xi$,
$v_{\alpha}\left(  \cdot,\xi\right)  \in C_{BP}^{\infty}\left(  \mathbb{R}%
^{d}\right)  $. We now differentiate equation \ref{q89} of part 1 with respect
to $x$ so that for $\left\vert \gamma\right\vert <\theta$ and $\left\vert
\alpha\right\vert =\theta$,%
\begin{align*}
\left(  D_{x}^{\gamma}v_{\alpha}\right)  \left(  x,\xi\right)   & =\frac
{1}{\alpha!}D_{x}^{\gamma}\mathcal{Q}_{x}\left(  x^{\alpha}\mu_{\theta}\left(
x\xi\right)  \right) \\
& =\frac{1}{\alpha!}D_{x}^{\gamma}\left(  x^{\alpha}\mu_{\theta}\left(
x\xi\right)  -\sum_{i=1}^{M}\left(  a^{\left(  i\right)  }\right)  ^{\alpha
}\mu_{\theta}\left(  a^{\left(  i\right)  }\xi\right)  l_{i}\left(  x\right)
\right) \\
& =\frac{1}{\alpha!}D_{x}^{\gamma}\left(  x^{\alpha}\mu_{\theta}\left(
x\xi\right)  \right)  -\sum_{i=1}^{M}\left(  a^{\left(  i\right)  }\right)
^{\alpha}\mu_{\theta}\left(  a^{\left(  i\right)  }\xi\right)  D^{\gamma}%
l_{i}\left(  x\right) \\
& =\frac{1}{\alpha!}\sum_{\substack{\beta\leq\gamma\\\beta\leq\alpha}%
}\binom{\gamma}{\beta}D^{\beta}\left(  x^{\alpha}\right)  D_{x}^{\gamma-\beta
}\left(  \mu_{\theta}\left(  x\xi\right)  \right)  -\sum_{i=1}^{M}\left(
a^{\left(  i\right)  }\right)  ^{\alpha}\mu_{\theta}\left(  a^{\left(
i\right)  }\xi\right)  D^{\gamma}l_{i}\left(  x\right) \\
& =\frac{1}{\alpha!}\sum_{\substack{\beta\leq\gamma\\\beta\leq\alpha}%
}\binom{\gamma}{\beta}\frac{\alpha!}{\left(  \alpha-\beta\right)  !}%
x^{\alpha-\beta}\xi^{\gamma-\beta}\left(  D^{\left\vert \gamma-\beta
\right\vert }\mu_{\theta}\right)  \left(  x\xi\right)  -\\
& \qquad-\frac{1}{\alpha!}\sum_{i=1}^{M}\left(  a^{\left(  i\right)  }\right)
^{\alpha}\mu_{\theta}\left(  a^{\left(  i\right)  }\xi\right)  D^{\gamma}%
l_{i}\left(  x\right) \\
& =\frac{1}{\alpha!}\sum_{\substack{\beta\leq\gamma\\\beta\leq\alpha}%
}\binom{\gamma}{\beta}\binom{\alpha}{\beta}x^{\alpha-\beta}\xi^{\gamma-\beta
}\beta!\left(  D^{\left\vert \gamma-\beta\right\vert }\mu_{\theta}\right)
\left(  x\xi\right)  -\\
& \qquad-\frac{1}{\alpha!}\sum_{i=1}^{M}\left(  a^{\left(  i\right)  }\right)
^{\alpha}\mu_{\theta}\left(  a^{\left(  i\right)  }\xi\right)  D^{\gamma}%
l_{i}\left(  x\right)  .
\end{align*}

and hence using \ref{a1.03},%
\begin{align*}
&  \left\vert \left(  D_{x}^{\gamma}v_{\alpha}\right)  \left(  x,\xi\right)
\right\vert \\
&  \leq\frac{1}{\alpha!}\sum_{\beta\leq\gamma;\beta\leq\alpha}\tbinom{\gamma
}{\beta}\tbinom{\alpha}{\beta}x_{+}^{\alpha-\beta}\xi_{+}^{\gamma-\beta}%
\beta!\left\vert \left(  D_{x}^{\left\vert \gamma-\beta\right\vert }%
\mu_{\theta}\right)  \left(  x\xi\right)  \right\vert +\\
&  \qquad\qquad+\frac{1}{\alpha!}\sum_{i=1}^{M}\left\vert \left(  a^{\left(
i\right)  }\right)  ^{\alpha}\right\vert \left\vert \mu_{\theta}\left(
a^{\left(  i\right)  }\xi\right)  \right\vert \left\vert D^{\gamma}%
l_{i}\left(  x\right)  \right\vert \\
&  \leq\frac{1}{\alpha!}\sum_{\beta\leq\gamma;\beta\leq\alpha}\tbinom{\gamma
}{\beta}\tbinom{\alpha}{\beta}x_{+}^{\alpha-\beta}\xi_{+}^{\gamma-\beta}%
\tfrac{\beta!\left\vert \gamma-\beta\right\vert !}{\left(  \left\vert
\gamma-\beta\right\vert +\theta\right)  !}+\frac{1}{\alpha!}\frac{1}{\theta
!}\sum_{i=1}^{M}\left\vert \left(  a^{\left(  i\right)  }\right)  ^{\alpha
}\right\vert \left\vert D^{\gamma}l_{i}\left(  x\right)  \right\vert \\
&  \leq\frac{1}{\alpha!}\sum_{\beta\leq\gamma;\beta\leq\alpha}\tbinom{\gamma
}{\beta}\tbinom{\alpha}{\beta}x_{+}^{\alpha-\beta}\xi_{+}^{\gamma-\beta}%
+\frac{1}{\alpha!}\frac{1}{\theta!}\sum_{i=1}^{M}\left\vert \left(  a^{\left(
i\right)  }\right)  ^{\alpha}\right\vert \left\vert D^{\gamma}l_{i}\left(
x\right)  \right\vert .
\end{align*}

Now%
\begin{align*}
\frac{1}{\alpha!} &  \sum_{\beta\leq\gamma;\beta\leq\alpha}\tbinom{\gamma
}{\beta}\tbinom{\alpha}{\beta}x_{+}^{\alpha-\beta}\xi_{+}^{\gamma-\beta}\\
&  \leq\frac{1}{\alpha!}\sum_{\beta\leq\gamma;\beta\leq\alpha}\tbinom{\gamma
}{\beta}\tbinom{\alpha}{\beta}x_{+}^{\alpha-\beta}\xi_{+}^{\gamma-\beta}\\
&  \leq\frac{1}{\alpha!}\left(  \sum_{\beta\leq\gamma;\beta\leq\alpha}%
\tbinom{\gamma}{\beta}^{2}\xi^{2\left(  \gamma-\beta\right)  }\right)
^{\frac{1}{2}}\left(  \sum_{\beta\leq\gamma;\beta\leq\alpha}\tbinom{\alpha
}{\beta}^{2}x^{2\left(  \alpha-\beta\right)  }\right)  ^{\frac{1}{2}}\\
&  \leq\frac{1}{\alpha!}\left(  \sum_{\beta\leq\gamma}\tbinom{\gamma}{\beta
}\xi_{+}^{\gamma-\beta}\right)  \left(  \sum_{\beta\leq\gamma;\beta\leq\alpha
}\tbinom{\alpha}{\beta}x_{+}^{\alpha-\beta}\right) \\
&  =\frac{1}{\alpha!}\left(  \sum_{\beta\leq\gamma}\tbinom{\gamma}{\beta}%
\xi_{+}^{\beta}\right)  \left(  \sum_{\beta\leq\gamma;\beta\leq\alpha}%
\tbinom{\alpha}{\beta}x_{+}^{\alpha-\beta}\right) \\
&  =\left(  1+\xi_{+}\right)  ^{\gamma}\sum_{\substack{\beta\leq\gamma
\\\beta\leq\alpha}}\tfrac{1}{\alpha!}\tbinom{\alpha}{\beta}x_{+}^{\alpha
-\beta}\leq\left(  1+\xi_{+}\right)  ^{\gamma}\sum_{\substack{\beta\leq
\gamma\\\beta\leq\alpha}}\tfrac{1}{\alpha!}\tbinom{\alpha}{\beta}\left\vert
x\right\vert ^{\left\vert \alpha\right\vert -\left\vert \beta\right\vert }%
\leq\\
&  \leq\left(  1+\xi_{+}\right)  ^{\gamma}\sum_{\substack{\beta\leq
\gamma\\\beta\leq\alpha}}\tfrac{1}{\gamma!}\tbinom{\gamma}{\beta}\left\vert
x\right\vert ^{\theta-\left\vert \beta\right\vert }=\tfrac{\left(  1+\xi
_{+}\right)  ^{\gamma}}{\gamma!}\left\vert x\right\vert ^{\theta-\left\vert
\gamma\right\vert }\sum_{\substack{\beta\leq\gamma\\\beta\leq\alpha}%
}\tbinom{\gamma}{\beta}\left\vert x\right\vert ^{\left\vert \gamma\right\vert
-\left\vert \beta\right\vert }\leq\\
&  \leq\tfrac{\left(  1+\xi_{+}\right)  ^{\gamma}}{\gamma!}\left\vert
x\right\vert ^{\theta-\left\vert \gamma\right\vert }\sum_{\beta\leq\gamma
}\tbinom{\gamma}{\beta}\left\vert x\right\vert ^{\left\vert \gamma\right\vert
-\left\vert \beta\right\vert }=\tfrac{\left(  1+\xi_{+}\right)  ^{\gamma}%
}{\gamma!}\left\vert x\right\vert ^{\theta-\left\vert \gamma\right\vert }%
\sum_{\beta\leq\gamma}\tbinom{\gamma}{\beta}\left\vert x\right\vert
^{\left\vert \beta\right\vert }\leq\\
&  \leq\tfrac{\left(  1+\xi_{+}\right)  ^{\gamma}}{\gamma!}\left\vert
x\right\vert ^{\theta-\left\vert \gamma\right\vert }\sum_{\beta\leq\gamma
}\tbinom{\gamma}{\beta}\left(  1+\left\vert x\right\vert \right)  ^{\left\vert
\beta\right\vert }\leq\tfrac{\left(  1+\xi_{+}\right)  ^{\gamma}}{\gamma
!}\left\vert x\right\vert ^{\theta-\left\vert \gamma\right\vert }\left(
1+\left\vert x\right\vert \right)  ^{\left\vert \gamma\right\vert }\sum
_{\beta\leq\gamma}\tbinom{\gamma}{\beta}=\\
&  =\tfrac{2^{\left\vert \gamma\right\vert }}{\gamma!}\left(  1+\xi
_{+}\right)  ^{\gamma}\left\vert x\right\vert ^{\theta-\left\vert
\gamma\right\vert }\left(  1+\left\vert x\right\vert \right)  ^{\left\vert
\gamma\right\vert },
\end{align*}

so that%
\[
\left\vert \left(  D_{x}^{\gamma}v_{\alpha}\right)  \left(  x,\xi\right)
\right\vert \leq\frac{2^{\left\vert \gamma\right\vert }}{\gamma!}\left(
1+\xi_{+}\right)  ^{\gamma}\left\vert x\right\vert ^{\theta-\left\vert
\gamma\right\vert }\left(  1+\left\vert x\right\vert \right)  ^{\left\vert
\gamma\right\vert }+\frac{1}{\theta!}\sum_{i=1}^{M}\left\vert \tfrac{\left(
a^{\left(  i\right)  }\right)  ^{\alpha}}{\alpha!}\right\vert \left\vert
D^{\gamma}l_{i}\left(  x\right)  \right\vert .
\]

\end{proof}

\begin{theorem}
\label{Thm_bnd_DQ(exp)}Suppose the operator $\mathcal{Q}$ is generated by the
minimal $\theta$-unisolvent set $A=\left\{  a^{\left(  i\right)  }\right\}
_{i=1}^{M}$ and the cardinal basis $\left\{  l_{i}\right\}  _{i=1}^{M}$. Then
for $x,\xi\in\mathbb{R}^{d}$:

\begin{enumerate}
\item For all $\gamma$ we have the simple but useful estimate%
\[
\left\vert D_{x}^{\gamma}\mathcal{Q}_{x}\left(  e^{ix\xi}\right)  \right\vert
\leq\xi_{+}^{\gamma}+\sum\limits_{k=1}^{M}\left\vert \left(  D^{\gamma}%
l_{\tau_{y}A;k}\right)  (x+y)\right\vert ,\quad x,y,\xi\in\mathbb{R}^{d}.
\]

\item For all multi-indexes $\gamma$ and $x,y,z,\xi\in\mathbb{R}^{d}$,%
\[
\left\vert D_{x}^{\gamma}\mathcal{Q}_{x}\left(  e^{ix\xi}\right)  \right\vert
\leq\left\{
\begin{array}
[c]{ll}%
\left\vert \xi\right\vert ^{\theta}\left(
\begin{array}
[c]{r}%
\tbinom{d+\theta-1}{\theta}\frac{2^{\left\vert \gamma\right\vert }}{\gamma
!}\left(  1+\xi_{+}\right)  ^{\gamma}\left\vert x+y\right\vert ^{\theta
-\left\vert \gamma\right\vert }\left(  1+\left\vert x+y\right\vert \right)
^{\left\vert \gamma\right\vert }+\\
+\frac{\left\vert \tau_{y}A\right\vert _{\infty}^{\theta}}{\left(
\theta!\right)  ^{2}}c_{\tau_{z},\gamma}\left(  1+\left\vert x+z\right\vert
\right)  ^{\theta-\left\vert \gamma\right\vert -1}%
\end{array}
\right)  , & \left\vert \gamma\right\vert <\theta,\\
\xi_{+}^{\gamma}, & \left\vert \gamma\right\vert \geq\theta,
\end{array}
\right.
\]

where $c_{A,\gamma}$ satisfies \ref{p912} and $\left\vert \tau_{y}A\right\vert
_{\infty}=\max\limits_{k}\left\vert a^{\left(  k\right)  }+y\right\vert $.
This inequality is useful for $\xi$ near the origin.

Also%
\[
\left\vert D_{x}^{\gamma}\mathcal{Q}_{x}\left(  e^{ix\xi}\right)  \right\vert
\leq\left\{
\begin{array}
[c]{ll}%
\left\vert \xi\right\vert ^{\theta}\frac{\left\vert \tau_{-x}A\right\vert
_{\infty}^{\theta}}{\left(  \theta!\right)  ^{2}}c_{\tau_{z}A,\gamma}\left(
1+\left\vert x+z\right\vert \right)  ^{\theta-\left\vert \gamma\right\vert
-1}, & \left\vert \gamma\right\vert <\theta,\\
\xi_{+}^{\gamma}, & \left\vert \gamma\right\vert \geq\theta.
\end{array}
\right.
\]

\item For all $\gamma$ and all $z\in\mathbb{R}^{d}$,%
\[
\left\vert D_{x}^{\gamma}\mathcal{Q}_{x}\left(  e^{ix\xi}\right)  \right\vert
\leq\left\{
\begin{array}
[c]{ll}%
\xi_{+}^{\gamma}+c_{\tau_{z}A,\gamma}\left(  1+\left\vert x+z\right\vert
\right)  ^{\theta-\left\vert \gamma\right\vert -1}, & \left\vert
\gamma\right\vert <\theta,\\
\xi_{+}^{\gamma}, & \left\vert \gamma\right\vert \geq\theta.
\end{array}
\right.
\]

This inequality is useful for all $x$ or for when $\xi$ is large, as the
estimates have lower powers than those in part 2.
\end{enumerate}
\end{theorem}

\begin{proof}
\textbf{Part 1} By definition of $\mathcal{Q}_{x}$,%
\[
D_{x}^{\gamma}\mathcal{Q}_{x}\left(  e^{ix\xi}\right)  =\left(  i\xi\right)
^{\gamma}e^{ix\xi}-\sum\limits_{k=1}^{M}\left(  D^{\gamma}l_{k}\right)
(x)e^{ia^{\left(  k\right)  }\xi},
\]

so it follows that
\[
\left\vert D_{x}^{\gamma}\mathcal{Q}_{x}\left(  e^{ix\xi}\right)  \right\vert
\leq\xi_{+}^{\gamma}+\sum\limits_{k=1}^{M}\left\vert \left(  D^{\gamma}%
l_{k}\right)  (x)\right\vert .
\]

\textbf{Part 2} Since $\left\vert \gamma\right\vert <\theta$, we can apply the
estimate of part 3 Lemma \ref{Lem_estim_Q(exp)} to the formula \ref{q71} for
$\mathcal{Q}_{x}\left(  e^{ix\xi}\right)  $ of part 1 of Lemma
\ref{Lem_estim_Q(exp)}: ???%
\[
\left\vert D_{x}^{\gamma}v_{\alpha}\left(  x,\xi\right)  \right\vert \leq
\frac{2^{\left\vert \gamma\right\vert }}{\gamma!}\left(  1+\xi_{+}\right)
^{\gamma}\left\vert x\right\vert ^{\theta-\left\vert \gamma\right\vert
}\left(  1+\left\vert x\right\vert \right)  ^{\left\vert \gamma\right\vert
}+\frac{1}{\theta!}\sum_{i=1}^{M}\tfrac{\left(  a_{+}^{\left(  i\right)
}\right)  ^{\alpha}}{\alpha!}\left\vert D^{\gamma}l_{i}\left(  x\right)
\right\vert ,
\]

to get%
\begin{align}
&  \left\vert D_{x}^{\gamma}\mathcal{Q}_{x}\left(  e^{ix\xi}\right)
\right\vert \nonumber\\
&  \leq\sum_{\left\vert \alpha\right\vert =\theta}\left\vert D_{x}^{\gamma
}v_{\alpha}\left(  x,\xi\right)  \left(  i\xi\right)  ^{\alpha}\right\vert
=\sum_{\left\vert \alpha\right\vert =\theta}\left\vert D_{x}^{\gamma}%
v_{\alpha}\left(  x,\xi\right)  \xi_{+}^{\alpha}\right\vert \leq\nonumber\\
&  \leq\sum_{\left\vert \alpha\right\vert =\theta}\left(  \frac{2^{\left\vert
\gamma\right\vert }}{\gamma!}\left(  1+\xi_{+}\right)  ^{\gamma}\left\vert
x\right\vert ^{\theta-\left\vert \gamma\right\vert }\left(  1+\left\vert
x\right\vert \right)  ^{\left\vert \gamma\right\vert }\xi_{+}^{\alpha}\right)
+\frac{1}{\theta!}\sum_{i=1}^{M}\left(  \sum_{\left\vert \alpha\right\vert
=\theta}\tfrac{\left(  a_{+}^{\left(  i\right)  }\right)  ^{\alpha}\xi
_{+}^{\alpha}}{\alpha!}\right)  \left\vert D^{\gamma}l_{i}\left(  x\right)
\right\vert \nonumber\\
&  =\frac{2^{\left\vert \gamma\right\vert }}{\gamma!}\left(  1+\xi_{+}\right)
^{\gamma}\left\vert x\right\vert ^{\theta-\left\vert \gamma\right\vert
}\left(  1+\left\vert x\right\vert \right)  ^{\left\vert \gamma\right\vert
}\sum_{\left\vert \alpha\right\vert =\theta}\xi_{+}^{\alpha}+\frac{1}{\theta
!}\sum_{i=1}^{M}\tfrac{\left(  a_{+}^{\left(  i\right)  }\xi_{+}\right)
^{\theta}}{\theta!}\left\vert D^{\gamma}l_{i}\left(  x\right)  \right\vert
\nonumber\\
&  \leq\frac{2^{\left\vert \gamma\right\vert }}{\gamma!}\left(  1+\xi
_{+}\right)  ^{\gamma}\left\vert x\right\vert ^{\theta-\left\vert
\gamma\right\vert }\left(  1+\left\vert x\right\vert \right)  ^{\left\vert
\gamma\right\vert }\tbinom{d+\theta-1}{\theta}\left\vert \xi\right\vert
^{\theta}+\left(  \frac{1}{\theta!}\right)  ^{2}\sum_{i=1}^{M}\left\vert
a^{\left(  i\right)  }\right\vert ^{\theta}\left\vert \xi\right\vert ^{\theta
}\left\vert D^{\gamma}l_{i}\left(  x\right)  \right\vert \nonumber\\
&  =\left\vert \xi\right\vert ^{\theta}\left(  \tbinom{d+\theta-1}{\theta
}\frac{2^{\left\vert \gamma\right\vert }}{\gamma!}\left(  1+\xi_{+}\right)
^{\gamma}\left\vert x\right\vert ^{\theta-\left\vert \gamma\right\vert
}\left(  1+\left\vert x\right\vert \right)  ^{\left\vert \gamma\right\vert
}+\left(  \frac{1}{\theta!}\right)  ^{2}\sum_{i=1}^{M}\left\vert a^{\left(
i\right)  }\right\vert ^{\theta}\left\vert D^{\gamma}l_{i}\left(  x\right)
\right\vert \right) \nonumber\\
&  =\left\vert \xi\right\vert ^{\theta}\left(  \tbinom{d+\theta-1}{\theta
}\frac{2^{\left\vert \gamma\right\vert }}{\gamma!}\left(  1+\xi_{+}\right)
^{\gamma}\left\vert x\right\vert ^{\theta-\left\vert \gamma\right\vert
}\left(  1+\left\vert x\right\vert \right)  ^{\left\vert \gamma\right\vert
}+\frac{\left\vert A\right\vert _{\infty}^{\theta}}{\left(  \theta!\right)
^{2}}\sum_{i=1}^{M}\left\vert D^{\gamma}l_{i}\left(  x\right)  \right\vert
\right)  .\label{a029}%
\end{align}

From the translation result Lemma \ref{Lem_estim_Q(exp)},%
\[
\mathcal{Q}_{x}\left(  e^{ix\xi}\right)  =\mathcal{Q}_{x}\left(  e^{iy\xi
}e^{i\left(  x-y\right)  \xi}\right)  =e^{iy\xi}\mathcal{Q}_{x}\left(
e^{i\left(  x-y\right)  \xi}\right)  =e^{iy\xi}\mathcal{Q}_{x}\left(  \tau
_{y}e^{ix\xi}\right)  .
\]

But from Theorem \ref{Thm_translat_Qf} and then part 1 of Lemma
\ref{Lem_estim_Q(exp)},%
\begin{align*}
\mathcal{Q}_{x}\left(  \tau_{y}e^{ix\xi}\right)   & =\mathcal{Q}_{A;x}\left(
\tau_{y}e^{ix\xi}\right)  =\tau_{y}\mathcal{Q}_{\tau_{y}A;x}\left(  e^{ix\xi
}\right)  =\\
& =\tau_{y}\sum_{\left\vert \alpha\right\vert =\theta}v_{\tau_{y}A;\alpha
}\left(  x,\xi\right)  \frac{\left(  i\xi\right)  ^{\alpha}}{\alpha!}\\
& =\sum_{\left\vert \alpha\right\vert =\theta}v_{\tau_{y}A;\alpha}\left(
x-y,\xi\right)  \frac{\left(  i\xi\right)  ^{\alpha}}{\alpha!},
\end{align*}

and so%
\[
\left\vert D_{x}^{\gamma}\mathcal{Q}_{x}\left(  e^{ix\xi}\right)  \right\vert
=\left\vert \sum_{\left\vert \alpha\right\vert =\theta}D_{x}^{\gamma}%
v_{\tau_{y}A;\alpha}\left(  x+y,\xi\right)  \frac{\left(  i\xi\right)
^{\alpha}}{\alpha!}\right\vert .
\]

Thus for all $y\in\mathbb{R}^{d}$, \ref{a029} becomes%
\[
\left\vert D_{x}^{\gamma}\mathcal{Q}_{x}\left(  e^{ix\xi}\right)  \right\vert
\leq\left\vert \xi\right\vert ^{\theta}\left(
\begin{array}
[c]{r}%
\tbinom{d+\theta-1}{\theta}\frac{2^{\left\vert \gamma\right\vert }}{\gamma
!}\left(  1+\xi_{+}\right)  ^{\gamma}\left\vert x+y\right\vert ^{\theta
-\left\vert \gamma\right\vert }\left(  1+\left\vert x+y\right\vert \right)
^{\left\vert \gamma\right\vert }+\\
+\frac{\left\vert \tau_{y}A\right\vert _{\infty}^{\theta}}{\left(
\theta!\right)  ^{2}}\sum\limits_{i=1}^{M}\left\vert D_{x}^{\gamma}%
l_{i}\left(  x+y\right)  \right\vert
\end{array}
\right)  .
\]

Now $\left\{  l_{i}\right\}  $ is the cardinal basis corresponding to the
minimal unisolvent set $A=\left\{  a^{\left(  i\right)  }\right\}  _{i=1}^{M}%
$and there exist constants $c_{A,\gamma}$ such that ?? WHY? ??
\begin{equation}
\sum_{i=1}^{M}\left\vert D^{\gamma}l_{i}\left(  x\right)  \right\vert \leq
c_{A,\gamma}\left(  1+\left\vert x\right\vert \right)  ^{\theta-\left\vert
\gamma\right\vert -1},\quad\left\vert \gamma\right\vert <\theta,\text{ }%
x\in\mathbb{R}^{d}.\label{p912}%
\end{equation}

But for all $x,y,z\in\mathbb{R}^{d}$,%
\[
l_{\tau_{y}A;i}\left(  x+y\right)  =l_{A;i}\left(  x\right)  =l_{\tau_{z}%
A;i}\left(  x+z\right)  .
\]

Thus%
\[
\left\vert D_{x}^{\gamma}\mathcal{Q}_{x}\left(  e^{ix\xi}\right)  \right\vert
\leq\left\vert \xi\right\vert ^{\theta}\left(
\begin{array}
[c]{r}%
\tbinom{d+\theta-1}{\theta}\frac{2^{\left\vert \gamma\right\vert }}{\gamma
!}\left(  1+\xi_{+}\right)  ^{\gamma}\left\vert x+y\right\vert ^{\theta
-\left\vert \gamma\right\vert }\left(  1+\left\vert x+y\right\vert \right)
^{\left\vert \gamma\right\vert }+\\
+\frac{\left\vert \tau_{y}A\right\vert _{\infty}^{\theta}}{\left(
\theta!\right)  ^{2}}c_{\tau_{z}A,\gamma}\left(  1+\left\vert x+z\right\vert
\right)  ^{\theta-\left\vert \gamma\right\vert -1}%
\end{array}
\right)  .
\]

On the other hand, if $\left\vert \gamma\right\vert \geq\theta$, $\left\vert
D_{x}^{\gamma}\mathcal{Q}_{x}\left(  e^{ix\xi}\right)  \right\vert =\xi
_{+}^{\gamma}$.

The second set of inequalities is proved by setting $y=-x$.\medskip

\textbf{Part 3} If $\left\vert \gamma\right\vert <\theta$ and $z\in
\mathbb{R}^{d}$ then from part 1,
\[
\left\vert D_{x}^{\gamma}\mathcal{Q}_{x}\left(  e^{ix\xi}\right)  \right\vert
\leq\xi_{+}^{\gamma}+\sum\limits_{k=1}^{M}\left\vert \left(  D^{\gamma}%
l_{k}\right)  (x)\right\vert \leq\xi_{+}^{\gamma}+c_{\tau_{z}A,\gamma}\left(
1+\left\vert x+z\right\vert \right)  ^{\theta-1-\left\vert \gamma\right\vert
}.
\]

On the other hand, if $\left\vert \gamma\right\vert \geq\theta$, we use the
result of part 2.
\end{proof}

\begin{theorem}
\label{Thm_|DQ(exp)|^2/(w|x|^2theta)} Suppose that the weight function $w$ has
properties W2.1 and W3 for order $\theta$ and $\kappa$, and that the operator
$\mathcal{Q}$ is defined using a minimal unisolvent set $A$ of order $\theta$.
Define $r_{4}$ by $r_{4}=0$ if $w$ has property W3.1 or W3.3\ and $r_{4}%
=r_{3}$ if $w$ has property W3.2. Then for each $x\in\mathbb{R}^{d}$,
$\frac{D_{x}^{\gamma}\mathcal{Q}_{x}\left(  e^{i\left(  x,\cdot\right)
}\right)  }{\sqrt{w}\left\vert \cdot\right\vert ^{2\theta}}\in L^{2}$ whenever
$\left\vert \gamma\right\vert \leq\underline{\kappa}$. In fact, there exists a
positive constant $c_{w}$, independent of $x$, such that:%
\[
\int\frac{\left\vert D_{x}^{\gamma}\mathcal{Q}_{x}\left(  e^{i\left(
x,\cdot\right)  }\right)  \right\vert ^{2}}{w\left\vert \cdot\right\vert
^{2\theta}}\leq\left\{
\begin{array}
[c]{ll}%
\left(  c_{w,\gamma}\right)  ^{2}\left(  1+\left\vert x\right\vert \right)
^{2\theta}, & \left\vert \gamma\right\vert <\theta,\text{ }\left\vert
\gamma\right\vert \leq\underline{\kappa},\\
\left(  c_{w,\gamma}\right)  ^{2}, & \theta\leq\left\vert \gamma\right\vert
\leq\underline{\kappa}.
\end{array}
\right.
\]

\end{theorem}

\begin{proof}
\frame{ \textbf{Case 1} $\left\vert \gamma\right\vert <\theta$, $\left\vert
\gamma\right\vert \leq\underline{\kappa}$ } From part 2 of Theorem
\ref{Thm_bnd_DQ(exp)} with $z=y=0$:%
\begin{align}
\int\limits_{\left\vert \cdot\right\vert \leq r_{4}}\frac{\left\vert
D_{x}^{\gamma}\mathcal{Q}_{x}\left(  e^{i\left(  x,\cdot\right)  }\right)
\right\vert ^{2}}{w\left\vert \cdot\right\vert ^{2\theta}}  & \leq
Const\int\limits_{\left\vert \cdot\right\vert \leq r_{4}}\frac{\left(
\left\vert \cdot\right\vert ^{\theta}\left(  \tbinom{d+\theta-1}{\theta}%
\frac{2^{\left\vert \gamma\right\vert }}{\gamma!}\left(  1+\xi_{+}\right)
^{\gamma}+\frac{\left\vert A\right\vert _{\infty}^{\theta}}{\theta
!}c_{A,\gamma}\right)  \left(  1+\left\vert x\right\vert \right)  ^{\theta
}\right)  ^{2}}{w\left\vert \cdot\right\vert ^{2\theta}}\nonumber\\
& \leq\left(  \int\limits_{\left\vert \cdot\right\vert \leq r_{4}}%
\frac{\left(  \tbinom{d+\theta-1}{\theta}\frac{2^{\left\vert \gamma\right\vert
}}{\gamma!}\left(  1+\xi_{+}\right)  ^{\gamma}+\frac{\left\vert A\right\vert
_{\infty}^{\theta}}{\theta!}c_{A,\gamma}\right)  ^{2}}{w}\right)  \left(
1+\left\vert x\right\vert \right)  ^{2\theta}\label{a1.05}\\
& \leq2\left(  \int\limits_{\left\vert \cdot\right\vert \leq r_{4}}%
\frac{\left(  \tbinom{d+\theta-1}{\theta}\frac{2^{\left\vert \gamma\right\vert
}}{\gamma!}\right)  ^{2}\left(  1+\left\vert \xi\right\vert \right)
^{2\left\vert \gamma\right\vert }+\left(  \frac{\left\vert A\right\vert
_{\infty}^{\theta}}{\theta!}c_{A,\gamma}\right)  ^{2}}{w}\right)  \left(
1+\left\vert x\right\vert \right)  ^{2\theta},\nonumber
\end{align}

which exists since property W2.1 is $1/w\in L_{loc}^{1}$. If $w$ is radial
then the formulas of Theorem \ref{Thm_integ_X^a_f(absX)dX_R^d} can be used to
obtain tighter estimates.

Next we use part 3 of Theorem \ref{Thm_bnd_DQ(exp)} with $z=0$:%
\begin{align*}
\int\limits_{\left\vert \cdot\right\vert \geq r_{4}}\frac{\left\vert
D_{x}^{\gamma}\mathcal{Q}_{x}\left(  e^{i\left(  x,\cdot\right)  }\right)
\right\vert ^{2}}{w\left\vert \cdot\right\vert ^{2\theta}}  & \leq
\int\limits_{\left\vert \cdot\right\vert \geq r_{4}}\frac{\left(  \xi
_{+}^{\gamma}+c_{A,\gamma}\left(  1+\left\vert x\right\vert \right)
^{\theta-1-\left\vert \gamma\right\vert }\right)  ^{2}}{w\left\vert
\cdot\right\vert ^{2\theta}}\\
& \leq\int\limits_{\left\vert \cdot\right\vert \geq r_{4}}\frac{2\xi^{2\gamma
}+2\left(  c_{A,\gamma}\right)  ^{2}\left(  1+\left\vert x\right\vert \right)
^{2\left(  \theta-1-\left\vert \gamma\right\vert \right)  }}{w\left\vert
\cdot\right\vert ^{2\theta}}\\
& \leq2\int\limits_{\left\vert \cdot\right\vert \geq r_{4}}\frac{\xi^{2\gamma
}}{w\left\vert \cdot\right\vert ^{2\theta}}+2\left(  c_{A,\gamma}\right)
^{2}\left(  \int\limits_{\left\vert \cdot\right\vert \geq r_{4}}\frac
{1}{w\left\vert \cdot\right\vert ^{2\theta}}\right)  \left(  1+\left\vert
x\right\vert \right)  ^{2\left(  \theta-1-\left\vert \gamma\right\vert
\right)  }%
\end{align*}

Thus for some $c_{w,\gamma}>0$, independent of $x$,
\[
\int\frac{\left\vert D_{x}^{\gamma}\mathcal{Q}_{x}\left(  e^{i\left(
x,\cdot\right)  }\right)  \right\vert ^{2}}{w\left\vert \cdot\right\vert
^{2\theta}}<\left(  c_{w,\gamma}\right)  ^{2}\left(  1+\left\vert x\right\vert
\right)  ^{2\theta},
\]
\medskip

\frame{ \textbf{Case 2} $\theta\leq\left\vert \gamma\right\vert \leq
\underline{\kappa}$ } Using part 2 of Theorem \ref{Thm_bnd_DQ(exp)} with $z=0$
we obtain%
\[
\int\frac{\left\vert D_{x}^{\gamma}\mathcal{Q}_{x}\left(  e^{i\left(
x,\cdot\right)  }\right)  \right\vert ^{2}}{w\left\vert \cdot\right\vert
^{2\theta}}\leq\int\frac{\xi^{2\gamma}}{w\left\vert \cdot\right\vert
^{2\theta}}=\left(  c_{w,\gamma}\right)  ^{2}.
\]

\end{proof}

\section{ inverse Fourier transform result for functions in $X_{w}^{\theta}%
$\label{Sect_inverse_Four_Xwth}}

In this section we use the properties of the function $\mathcal{Q}_{x}\left(
e^{i\left(  x,\xi\right)  }\right)  $ derived in the previous section to
derive an inverse Fourier transform theorem which expresses the value of a
function $f\in X_{w}^{\theta}$ in terms of its distribution Fourier transform
which is a function in $L_{loc}^{1}\left(  \mathbb{R}^{d}\setminus0\right)  $.

Summary \ref{Sum_int_properties_Xwm} will now recall some results about the
semi-Hilbert data space%
\[
X_{w}^{\theta}=\left\{  f\in S^{\prime}:\widehat{D^{\alpha}f}\in L_{loc}%
^{1}\left(  \mathbb{R}^{d}\right)  ,\text{ }\int w\left\vert \widehat
{D^{\alpha}f}\right\vert ^{2}<\infty\text{ }for\text{ }all\text{ }\left\vert
\alpha\right\vert =\theta\right\}  ,
\]

with semi-inner product and seminorm%
\[
\left\langle f,g\right\rangle _{w,\theta}=\sum\limits_{\left\vert
\alpha\right\vert =\theta}\frac{\theta!}{\alpha!}\int w\widehat{D^{\alpha}%
f}\,\overline{\widehat{D^{\alpha}g}},\text{\quad}\left\vert f\right\vert
_{w,\theta}=\sqrt{\left\langle f,f\right\rangle _{w,\theta}},
\]
introduced as Definition \ref{Def_Xwth}.

\begin{summary}
\label{Sum_int_properties_Xwm}Suppose the weight function $w$ has property
\textbf{W2}. If $f\in X_{w}^{\theta}$ then $\widehat{f}\in L_{loc}^{1}\left(
\mathbb{R}^{d}\setminus0\right)  $ and we can define the function
$f_{F}:\mathbb{R}^{d}\rightarrow\mathbb{C}$ a.e. by $f_{F}=\widehat{f}$ on
$\mathbb{R}^{d}\setminus0$. Further:

\begin{enumerate}
\item The seminorm $\left\vert \cdot\right\vert _{w,\theta}$ satisfies (part 1
Theorem \ref{Thm_properties_Xwm})%
\begin{equation}
\int w\left\vert \cdot\right\vert ^{2\theta}\left\vert f_{F}\right\vert
^{2}=\left\vert f\right\vert _{w,\theta}^{2}.\label{q81}%
\end{equation}

\item The functional $\left\vert \cdot\right\vert _{w,\theta}$ is a seminorm.
In fact, $\operatorname*{null}\left\vert \cdot\right\vert _{w,\theta
}=P_{\theta-1}$ and $X_{w}^{\theta}\cap P=P_{\theta-1}$ (part 3 Theorem
\ref{Thm_properties_Xwm}).

\item $f_{F}=0$ iff $f\in P_{\theta-1}$ (part 1 Theorem
\ref{Thm_property_(g)_F}).

\item $f_{F}\in S_{\emptyset,\theta}^{\prime}$ and $\widehat{f}=f_{F}$ on
$S_{\emptyset,\theta}$ (part 2 Theorem \ref{Thm_property_(g)_F}).

\item Since the weight function has property W2 we can use the definition of
$X_{w}^{\theta}$ from Corollary \ref{Cor_2_Thm_property_(g)_F}, namely
\begin{align}
X_{w}^{\theta}  & =\left\{  f\in S^{\prime}:\widehat{f}\in L_{loc}^{1}\left(
\mathbb{R}^{d}\setminus0\right)  ,\text{ }\int w\left\vert \cdot\right\vert
^{2\theta}\left\vert f_{F}\right\vert ^{2}<\infty\text{, }\left\vert
\alpha\right\vert =\theta\text{ }implies\text{ }\xi^{\alpha}\widehat{f}%
=\xi^{\alpha}f_{F}\text{ }on\text{ }S\right\}  .\nonumber\\
& \label{p948}%
\end{align}

Note that this definition makes sense because the conditions

$\widehat{f}\in L_{loc}^{1}\left(  \mathbb{R}^{d}\setminus0\right)  $ and
$\int w\left\vert \cdot\right\vert ^{2\theta}\left\vert f_{F}\right\vert
^{2}<\infty$ imply (part 2 Corollary \ref{Lem_properties_Xwm_distrib}) that
when $\left\vert \alpha\right\vert =\theta$, $\xi^{\alpha}f_{F}$ is a regular
tempered distribution in the sense of part 2 Appendix \ref{SbSect_property_S'}%
.\medskip

\item $X_{w}^{\theta}$ is complete in the seminorm sense (Theorem
\ref{Thm_Xw_complete}).

\item If $w$ also has property W3 for order $\theta$ and smoothness parameter
$\kappa$ then $X_{w}^{\theta}\subset C_{BP}^{\left(  \left\lfloor
\kappa\right\rfloor \right)  }$ (Theorem \ref{Thm_Xwth_W3_smooth}).
\end{enumerate}
\end{summary}

\begin{lemma}
\label{Lem_DQ(e)fhat=DQ(e)fF}Suppose the weight function $w$ also has property W2.

Then for all $f\in X_{w}^{\mathbb{\theta}}$ and all multi-indexes $\gamma
\geq0$%
\[
\left(  D_{x}^{\gamma}\mathcal{Q}_{x}\left(  e^{i\left(  x,\cdot\right)
}\right)  \right)  \widehat{f}=\left(  D_{x}^{\gamma}\mathcal{Q}_{x}\left(
e^{i\left(  x,\cdot\right)  }\right)  \right)  f_{F}\in S^{\prime},
\]

where the function $f_{F}:\mathbb{R}^{d}\rightarrow\mathbb{C}$ is defined a.e.
by $f_{F}=\widehat{f}$ on $\mathbb{R}^{d}\setminus0$.
\end{lemma}

\begin{proof}
If $\left\vert \gamma\right\vert <\theta$ then from part 2 Theorem
\ref{Thm_Q(exp)_loc_property} we have $D_{x}^{\gamma}\mathcal{Q}_{x}\left(
e^{i\left(  x,\cdot\right)  }\right)  \in C_{\emptyset,\theta}^{\infty}\cap
C_{BP}^{\infty}$.

If $\left\vert \gamma\right\vert \geq\theta$ then $D_{x}^{\gamma}%
\mathcal{Q}_{x}\left(  e^{i\left(  x,\xi\right)  }\right)  =\left(
i\xi\right)  ^{\gamma}e^{i\left(  x,\cdot\right)  }$ and by part 3 Theorem
\ref{Thm_product_of_Co,k_funcs},

$D_{x}^{\gamma}\mathcal{Q}_{x}\left(  e^{i\left(  x,\xi\right)  }\right)  \in
C_{\emptyset,\theta}^{\infty}\cap C_{BP}^{\infty}$ for each $x$. Again by
Theorem \ref{Thm_product_of_Co,k_funcs}, $\phi D_{x}^{\gamma}\mathcal{Q}%
_{x}\left(  e^{i\left(  x,\cdot\right)  }\right)  \in S_{\emptyset,\theta}$ if
$\phi\in S$. In addition, from part 4 Summary \ref{Sum_int_properties_Xwm},
$f_{F}\in S_{\emptyset,\theta}^{\prime}$ and $\widehat{f}=f_{F}$ on
$S_{\emptyset,\theta}$. So for $\phi\in S$ we now have
\begin{align*}
\left[  \left(  D_{x}^{\gamma}\mathcal{Q}_{x}\left(  e^{i\left(
x,\cdot\right)  }\right)  \right)  \widehat{f},\phi\right]  =\left[
\widehat{f},\phi D_{x}^{\gamma}\mathcal{Q}_{x}\left(  e^{i\left(
x,\cdot\right)  }\right)  \right]   & =\left[  f_{F},\phi D_{x}^{\gamma
}\mathcal{Q}_{x}\left(  e^{i\left(  x,\cdot\right)  }\right)  \right] \\
& =\left[  \left(  D_{x}^{\gamma}\mathcal{Q}_{x}\left(  e^{i\left(
x,\cdot\right)  }\right)  \right)  f_{F},\phi\right]  ,
\end{align*}

which proves this lemma.
\end{proof}

The next theorem is an inverse Fourier transform result for functions in
$X_{w}^{\mathbb{\theta}}$.

\begin{theorem}
\label{Thm_inv_Fourier_tran_func}If the weight function $w$ has properties
W2.1 and W3 for order $\theta\geq1$ and $\kappa\geq0$ then $X_{w}^{\theta
}\subset C_{BP}^{\left(  \left\lfloor \kappa\right\rfloor \right)  }$.

Further, if $f\in X_{w}^{\theta}$, $f_{F}=\widehat{f}$ on $\mathbb{R}%
^{d}\setminus0$, $x\in\mathbb{R}^{d}$ and $\left\vert \gamma\right\vert
\leq\left\lfloor \min\kappa\right\rfloor $ then

$\left(  D_{x}^{\gamma}\mathcal{Q}_{x}\left(  e^{i\left(  x,\cdot\right)
}\right)  \right)  f_{F}\in L^{1}$ and%
\begin{equation}
D^{\gamma}f\left(  x\right)  =\left(  2\pi\right)  ^{-\frac{d}{2}}\int\left(
D_{x}^{\gamma}\mathcal{Q}_{x}\left(  e^{i\left(  x,\xi\right)  }\right)
\right)  f_{F}\left(  \xi\right)  d\xi+\left(  D^{\gamma}\mathcal{P}f\right)
\left(  x\right)  ,\label{p920}%
\end{equation}

and we have the growth estimates%
\begin{equation}
\left\vert D^{\gamma}f\left(  x\right)  \right\vert \leq\left\{
\begin{array}
[c]{ll}%
c_{w,\gamma}\left\vert f\right\vert _{w,\theta}\left(  1+\left\vert
x\right\vert \right)  ^{\theta}+\left\vert \left(  D^{\gamma}\mathcal{P}%
f\right)  \left(  x\right)  \right\vert , & \left\vert \gamma\right\vert
<\theta,\\
c_{w,\gamma}\left\vert f\right\vert _{w,\theta}, & \left\vert \gamma
\right\vert \geq\theta,
\end{array}
\right. \label{p930}%
\end{equation}

where $c_{w,\gamma}$ is defined in Theorem \ref{Thm_|DQ(exp)|^2/(w|x|^2theta)}.
\end{theorem}

\begin{proof}
First note that $f_{F}\in L_{loc}^{1}\left(  \mathbb{R}^{d}\setminus0\right)
$ implies

$\left(  D_{x}^{\gamma}\mathcal{Q}_{x}\left(  e^{i\left(  x,\cdot\right)
}\right)  \right)  f_{F}\in L_{loc}^{1}\left(  \mathbb{R}^{d}\setminus
0\right)  $. For each $x$, we apply the Cauchy-Schwartz inequality to obtain%
\begin{align}
\int\left\vert \left(  D_{x}^{\gamma}\mathcal{Q}_{x}\left(  e^{i\left(
x,\cdot\right)  }\right)  \right)  f_{F}\right\vert  & =\int\frac{\left\vert
D_{x}^{\gamma}\mathcal{Q}_{x}\left(  e^{i\left(  x,\cdot\right)  }\right)
\right\vert }{\sqrt{w}\left\vert \cdot\right\vert ^{\theta}}\sqrt{w}\left\vert
\cdot\right\vert ^{\theta}\left\vert f_{F}\right\vert \nonumber\\
& \leq\left(  \int\frac{\left\vert D_{x}^{\gamma}\mathcal{Q}_{x}\left(
e^{i\left(  x,\cdot\right)  }\right)  \right\vert ^{2}}{w\left\vert
\cdot\right\vert ^{2\theta}}\right)  ^{\frac{1}{2}}\left\vert f\right\vert
_{w,\theta},\label{p919}%
\end{align}

and the last integral exists by Theorem \ref{Thm_|DQ(exp)|^2/(w|x|^2theta)}.
From Lemma \ref{Lem_DQ(e)fhat=DQ(e)fF}, for each $x$,%
\begin{align*}
\left(  D_{x}^{\gamma}\mathcal{Q}_{x}\left(  e^{i\left(  x,\cdot\right)
}\right)  \right)  f_{F}=\left(  D_{x}^{\gamma}\mathcal{Q}_{x}\left(
e^{i\left(  x,\cdot\right)  }\right)  \right)  \widehat{f} &  =\left(  \left(
i\xi\right)  ^{\gamma}e^{i\left(  x,\cdot\right)  }-\sum\limits_{i=1}%
^{M}e^{i\left\langle a^{\left(  i\right)  },\cdot\right\rangle }D^{\gamma
}l_{i}\left(  x\right)  \right)  \widehat{f}\\
&  =\left(  D^{\gamma}f\left(  \cdot+x\right)  -\sum\limits_{i=1}^{M}f\left(
\cdot+a^{\left(  i\right)  }\right)  D^{\gamma}l_{i}\left(  x\right)  \right)
^{\wedge}.
\end{align*}

Now set $g_{x}=D^{\gamma}f\left(  \cdot+x\right)  -\sum\limits_{i=1}%
^{M}D^{\gamma}l_{i}\left(  x\right)  f\left(  \cdot+a^{\left(  i\right)
}\right)  $ so that for each $x$, $g_{x}\in S^{\prime}$ and $\widehat{g_{x}%
}=D_{x}^{\gamma}\mathcal{Q}_{x}\left(  e^{i\left(  x,\cdot\right)  }\right)
\widehat{f}\in L^{1}$. But by Lemma \ref{Lem_L1_Fourier_contin} Chapter
\ref{Ch_weight_fn_exten}, if $u\in S^{\prime}$ and $\widehat{u}\in L^{1}$,
then $u\in C_{B}^{\left(  0\right)  }$ and $u=\left(  2\pi\right)  ^{-d/2}\int
e^{i\left(  \cdot,\xi\right)  }\widehat{u}\left(  \xi\right)  d\xi$. Thus,
$D^{\gamma}f\left(  \cdot+x\right)  -\sum\limits_{i=1}^{M}f\left(
\cdot+a^{\left(  i\right)  }\right)  D^{\gamma}l_{i}\left(  x\right)  \in
C_{B}^{\left(  0\right)  }$ and
\[
D^{\gamma}f\left(  \cdot+x\right)  -\sum\limits_{i=1}^{M}f\left(
\cdot+a^{\left(  i\right)  }\right)  D^{\gamma}l_{i}\left(  x\right)  =\left(
2\pi\right)  ^{-\frac{d}{2}}\int e^{i\left(  \cdot,\xi\right)  }D_{x}^{\gamma
}\mathcal{Q}_{x}\left(  e^{i\left(  x,\xi\right)  }\right)  f_{F}\left(
\xi\right)  d\xi,
\]

which implies%
\[
D^{\gamma}f\left(  x\right)  =\left(  2\pi\right)  ^{-\frac{d}{2}}\int
D_{x}^{\gamma}\mathcal{Q}_{x}\left(  e^{i\left(  x,\xi\right)  }\right)
f_{F}\left(  \xi\right)  d\xi+D^{\gamma}\mathcal{P}f\left(  x\right)  ,
\]

which proves \ref{p920}. Now we use the Cauchy-Schwartz inequality to get%
\begin{align*}
\left\vert D^{\gamma}f\left(  x\right)  \right\vert  & \leq\left(
2\pi\right)  ^{-\frac{d}{2}}\int\left\vert D_{x}^{\gamma}\mathcal{Q}%
_{x}\left(  e^{i\left(  x,\xi\right)  }\right)  f_{F}\left(  \xi\right)
\right\vert d\xi+\left\vert D^{\gamma}\mathcal{P}f\left(  x\right)
\right\vert \\
& \leq\left(  2\pi\right)  ^{-\frac{d}{2}}\int\frac{\left\vert D_{x}^{\gamma
}\mathcal{Q}_{x}\left(  e^{i\left(  x,\xi\right)  }\right)  \right\vert
}{\sqrt{w}\left\vert \cdot\right\vert ^{\theta}}\sqrt{w}\left\vert
\cdot\right\vert ^{\theta}\left\vert f_{F}\left(  \xi\right)  \right\vert
d\xi+\left\vert D^{\gamma}\mathcal{P}f\left(  x\right)  \right\vert \\
& \leq\left(  2\pi\right)  ^{-\frac{d}{2}}\left(  \int\frac{\left\vert
D_{x}^{\gamma}\mathcal{Q}_{x}\left(  e^{i\left(  x,\cdot\right)  }\right)
\right\vert ^{2}}{w\left\vert \cdot\right\vert ^{2\theta}}\right)
^{1/2}\left\vert f\right\vert _{w,\theta}+\left\vert D^{\gamma}\mathcal{P}%
f\left(  x\right)  \right\vert .
\end{align*}

Substituting the inequality of Theorem \ref{Thm_|DQ(exp)|^2/(w|x|^2theta)}
yields
\[
\left\vert D^{\gamma}f\left(  x\right)  \right\vert \leq\left\{
\begin{array}
[c]{ll}%
c_{\gamma,w}\left\vert f\right\vert _{w,\theta}\left(  1+\left\vert
x\right\vert \right)  ^{\theta}+\left\vert D^{\gamma}\mathcal{P}f\left(
x\right)  \right\vert , & \left\vert \gamma\right\vert <\theta,\\
c_{\gamma,w}\left\vert f\right\vert _{w,\theta} & \left\vert \gamma\right\vert
\geq\theta,
\end{array}
\right.
\]

which proves inequality \ref{p930}.
\end{proof}

\section{Riesz representers for the functionals $u\rightarrow D^{\gamma
}u\left(  x\right)  $\label{Sect_Riesz_rep}}

Suppose that the weight function $w$ has properties W2 and W3 for order
$\theta$ and smoothness $\kappa$. Then\ the functions in $X_{w}^{\theta}$ are
continuous. In this section we will derive simple explicit formulas for the
Riesz representers of the evaluation functionals $f\rightarrow D^{\gamma
}f\left(  x\right)  $ where $f\in X_{w}^{\theta}$ and $\left\vert
\gamma\right\vert \leq\left\lfloor \min\kappa\right\rfloor $ and
$X_{w}^{\theta}$ is endowed with a special norm called the Light norm. When
$\gamma=0$ the representer is denoted by $R_{x}$ and $D_{x}^{\gamma}R_{x}$ is
shown to be the Riesz representer for the evaluation functional $D^{\gamma
}u\left(  x\right)  $ when $\left\vert \gamma\right\vert \leq\left\lfloor
\min\kappa\right\rfloor $. The formula for $R_{x}$ is exhibited in equation
\ref{p939} and is expressed in terms of a basis function and the cardinal
polynomial basis used to define the norm. The existence of $R_{x}\in
X_{w}^{\theta}$ means that $X_{w}^{\theta}$ is a reproducing kernel Hilbert
space of continuous functions.

The approach to deriving the Riesz representers is to first introduce a
special subspace of $X_{w}^{\theta}$ in the next subsection and then to deduce
the form $R_{x}$ must have from the definition of $R_{x}$ and its implications
for this subspace of functions. The operators $\mathcal{P}$ and $\mathcal{Q}$
introduced in Section \ref{SbSect_P_Q} are also useful.

\subsection{The functions $G\ast\widehat{S}_{\emptyset,\theta}$%
\label{SbSect_G*F[Son]}}

We will now prove some properties of the space of functions, denoted by
$G\ast\widehat{S}_{\emptyset,\theta}$, having the form $G\ast\varphi$ where
$\varphi\in\widehat{S}_{\emptyset,\theta}$. Here $\widehat{S}_{\emptyset
,\theta}$ denotes the Fourier transforms of the functions in $S_{\emptyset
,\theta}$. In the next section these functions will help us deduce the precise
form of the Riesz representer of the functionals $f\rightarrow D^{\alpha
}f\left(  x\right)  $.

\begin{corollary}
\label{Cor_f_G_convol_phi}Suppose the weight function $w$ has property W2, and
suppose $G$ is a basis distribution of order $\theta\geq1$ generated by $w$.
Then $\varphi\in\widehat{S}_{\emptyset,\theta}$ implies $G\ast\varphi\in
X_{w}^{\theta}$ and%
\begin{equation}
\left\vert G\ast\varphi\right\vert _{w,\theta}^{2}=\int\frac{\left\vert
\widehat{\varphi}\right\vert ^{2}}{w\left\vert \cdot\right\vert ^{2\theta}%
}=\left[  \widehat{G},\left\vert \widehat{\varphi}\right\vert ^{2}\right]
.\label{q69}%
\end{equation}

Also
\begin{equation}
\left\langle f,G\ast\varphi\right\rangle _{w,\theta}=\left[  f,\overline
{\varphi}\right]  ,\quad f\in X_{w}^{\theta}.\label{p926}%
\end{equation}

\end{corollary}

\begin{proof}
The first step is to show that $G\ast\varphi\subset X_{w}^{\theta}$. Let
$g=G\ast\varphi$ where $\varphi\in\widehat{S}_{\emptyset,\theta}$. To prove
that $g\in X_{w}^{\theta}$ we use definition \ref{p948} of $X_{w}^{\theta}$
i.e. we show that $g\in S^{\prime}$, $\widehat{g}\in L_{loc}^{1}\left(
\mathbb{R}^{d}\setminus0\right)  $, $\int w\left\vert \cdot\right\vert
^{2\theta}\left\vert g_{F}\right\vert ^{2}<\infty$ and $\xi^{\alpha}%
\widehat{g}=\xi^{\alpha}g_{F}$ on $S$.

Since $G\in S^{\prime}$ and $\varphi\in S$, $g=G\ast\varphi\in S^{\prime}$ and
$\widehat{g}=\widehat{\varphi}\,\widehat{G}$, as distributions. Further, from
Definition \ref{Def_basis_distrib} of the basis distribution $G $ we know that
$\left[  \widehat{G},\psi\right]  =\int\frac{\psi}{w\left\vert \cdot
\right\vert ^{2\theta}}$ for $\psi\in S_{\emptyset,2\theta}$, and so $\left[
\widehat{g},\psi\right]  =\int\frac{\widehat{\varphi}\psi}{w\left\vert
\cdot\right\vert ^{2\theta}}$ for $\psi\in S_{\emptyset,2\theta}$. But
property W2.1 implies $\frac{\widehat{\varphi}}{w\left\vert \cdot\right\vert
^{2\theta}}\in L_{loc}^{1}\left(  \mathbb{R}^{d}\setminus0\right)  $ which
means that $g_{F}=\frac{\widehat{\varphi}}{w\left\vert \cdot\right\vert
^{2\theta}}$ on $\mathbb{R}^{d}\setminus0$. But by definition $\widehat
{\varphi}\in S_{\emptyset,\theta}$ so that Theorem
\ref{Thm_product_of_Co,k_funcs} implies $\left\vert \widehat{\varphi
}\right\vert ^{2}=\widehat{\varphi}\overline{\widehat{\varphi}}\in
S_{\emptyset,2\theta}$ and
\begin{equation}
\int w\left\vert \cdot\right\vert ^{2\theta}\left\vert g_{F}\right\vert
^{2}=\int w\left\vert \cdot\right\vert ^{2\theta}\frac{\left\vert
\widehat{\varphi}\right\vert ^{2}}{w^{2}\left\vert \cdot\right\vert ^{4\theta
}}=\int\frac{\left\vert \widehat{\varphi}\right\vert ^{2}}{w\left\vert
\cdot\right\vert ^{2\theta}}=\left[  \widehat{G},\left\vert \widehat{\varphi
}\right\vert ^{2}\right]  <\infty.\label{p949}%
\end{equation}

To show $\xi^{\alpha}\widehat{g}=\xi^{\alpha}g_{F}$ on $S$, we need to show
that $\xi^{\alpha}\widehat{\varphi}\widehat{G}=\frac{\xi^{\alpha}%
\widehat{\varphi}}{w\left\vert \cdot\right\vert ^{2\theta}}$ on $S$. But by
Theorem \ref{Thm_product_of_Co,k_funcs}, $\xi^{\alpha}\in C_{\emptyset,\theta
}^{\infty}\cap C_{BP}^{\infty}$ and hence $\xi^{\alpha}\widehat{\varphi}\in
S_{\emptyset,2\theta}$ and hence for $\psi\in S$ the basis distribution
definition implies $\left[  \xi^{\alpha}\widehat{\varphi}\widehat{G}%
,\psi\right]  =\left[  \widehat{G},\xi^{\alpha}\widehat{\varphi}\psi\right]
=\int\frac{\xi^{\alpha}\widehat{\varphi}\psi}{w\left\vert \cdot\right\vert
^{2\theta}}$, as required. This means that \ref{p949} is true and so proves
\ref{q69}.

Now to prove formula \ref{p926}. From above we know that $f_{F}=\widehat{f}$
on $\mathbb{R}^{d}\setminus0$ then
\[
\left\langle f,G\ast\varphi\right\rangle _{w}=\int w\left\vert \cdot
\right\vert ^{2\theta}f_{F}\overline{g_{F}}=\int w\left\vert \cdot\right\vert
^{2\theta}f_{F}\frac{\overline{\widehat{\varphi}}}{w\left\vert \cdot
\right\vert ^{2\theta}}=\int f_{F}\overline{\widehat{\varphi}}=\int
f_{F}\widehat{\varphi^{\ast}},
\]

where $\varphi^{\ast}\left(  x\right)  =\overline{\varphi\left(  -x\right)  }$
and $\widehat{\phi^{\ast}}\in S_{\emptyset,\theta}$. From part 4 of Summary
\ref{Sum_int_properties_Xwm} $f_{F}\in S_{\emptyset,\theta}^{\prime}$ and
$\widehat{f}=f_{F}$ on $S_{\emptyset,\theta}$, so that
\[
\left\langle f,G\ast\varphi\right\rangle _{w}=\left[  \widehat{f}%
,\widehat{\varphi^{\ast}}\right]  =\left[  f,\overline{\varphi}\right]  ,
\]
as required.
\end{proof}

\subsection{The Light norm for $X_{w}^{\theta}$\label{SbSect_L&W_norm}}

Using the results of the previous section, a unisolvent set of points will be
used to create an inner product space from the semi-inner product space
$X_{w}^{\theta}$ introduced above in Definition \ref{Def_Xwth}. This
particular inner product was introduced by Light and Wayne in
\cite{LightWayneX98Weight} and involves evaluating functions in $X_{w}%
^{\theta}$ on the unisolvent set, which is OK since $X_{w}^{\theta}\subset
C_{BP}^{\left(  \left\lfloor \kappa\right\rfloor \right)  }$ by part 7 of
Summary \ref{Sum_int_properties_Xwm}.

\begin{definition}
\label{Def_inner_prod_Xpos}Suppose $\theta\geq1$ and suppose the set
$A=\left\{  a^{\left(  i\right)  }\right\}  $ is minimally unisolvent with
respect to the polynomial space $P_{\theta-1}$. Then following Light and Wayne
\cite{LightWayneX98Weight} we will introduce the following Light norm and
inner product for the space $X_{w}^{\theta}$:
\begin{equation}
\left(  u,v\right)  _{w,\theta}=\left\langle u,v\right\rangle _{w,\theta}%
+\sum\limits_{i}u\left(  a^{\left(  i\right)  }\right)  \overline{v\left(
a^{\left(  i\right)  }\right)  },\quad\left\Vert u\right\Vert _{w,\theta}%
^{2}=\left\vert u\right\vert _{w,\theta}^{2}+\sum\limits_{i}\left\vert
u\left(  a^{\left(  i\right)  }\right)  \right\vert ^{2},\label{p917}%
\end{equation}

where $\left\vert u\right\vert _{w,\theta}$ is defined by \ref{q60}.
$\left\Vert \cdot\right\Vert _{w,\theta}$ is actually a norm since $\left\Vert
\cdot\right\Vert _{w,\theta}=0$ implies that $\left\vert u\right\vert
_{w,\theta}=0$ and $u\left(  a^{\left(  i\right)  }\right)  =0$ for all
unisolvent $a^{\left(  i\right)  }$. By part 3 of Summary
\ref{Sum_int_properties_Xwm}, $\left\vert u\right\vert _{w,\theta}=0$ means
$u\in P_{\theta-1}$ and so from the definition of unisolvent sets $u=0$.
\end{definition}

\begin{theorem}
\label{Thm_Light_norm_property}Properties of the Light norm.

Suppose the Light norm and the Lagrangian interpolation operator $\mathcal{P}
$ are defined using the same minimal unisolvent set $\left\{  a^{\left(
i\right)  }\right\}  $. Then:

\begin{enumerate}
\item $\left(  \mathcal{Q}u,v\right)  _{w,\theta}=\left(  u,\mathcal{Q}%
v\right)  _{w,\theta}=\left\langle u,v\right\rangle _{w,\theta}$.

\item $\left(  \mathcal{P}u,v\right)  _{w,\theta}=\left(  u,\mathcal{P}%
v\right)  _{w,\theta}=\sum\limits_{i}u\left(  a^{\left(  i\right)  }\right)
\overline{v\left(  a^{\left(  i\right)  }\right)  }$.
\end{enumerate}
\end{theorem}

\begin{proof}
Define the Light norm using the minimal unisolvent set $A=\left\{  a^{\left(
i\right)  }\right\}  $. Then using results from Theorem
\ref{Thm_P_op_properties} concerning the operators $\mathcal{P}$ and
$\mathcal{Q} $:

\textbf{Part 1} $\left(  \mathcal{Q}u,v\right)  _{w,\theta}=\left\langle
\mathcal{Q}u,v\right\rangle _{w,\theta}+\sum\limits_{i}\mathcal{Q}u\left(
a^{\left(  i\right)  }\right)  \overline{v\left(  a^{\left(  i\right)
}\right)  }=\left\langle \mathcal{Q}u,v\right\rangle _{w,\theta}$, since
$\mathcal{Q}u\left(  a^{\left(  i\right)  }\right)  =0$. But $\left\langle
\mathcal{Q}u,v\right\rangle _{w,\theta}=\left\langle u-\mathcal{P}%
u,v\right\rangle _{w,\theta}=\left\langle u,v\right\rangle _{w,\theta}$ and so
$\left(  \mathcal{Q}u,v\right)  _{w,\theta}=\left\langle u,v\right\rangle
_{w,\theta} $. Similarly $\left(  u,\mathcal{Q}v\right)  _{w,\theta
}=\left\langle u,v\right\rangle _{w,\theta}$.

\textbf{Part 2} $\left(  \mathcal{P}u,v\right)  _{w,\theta}=\left\langle
\mathcal{P}u,v\right\rangle _{w,\theta}+\sum\limits_{i}\mathcal{P}u\left(
a^{\left(  i\right)  }\right)  \overline{v\left(  a^{\left(  i\right)
}\right)  }=\sum\limits_{i}\mathcal{P}u\left(  a^{\left(  i\right)  }\right)
\overline{v\left(  a^{\left(  i\right)  }\right)  }=$

$=\sum\limits_{i}u\left(  a^{\left(  i\right)  }\right)  \overline{v\left(
a^{\left(  i\right)  }\right)  }$, since $\mathcal{P}u$ interpolates $\left\{
a^{\left(  i\right)  },u\left(  a^{\left(  i\right)  }\right)  \right\}  $.
Similarly $\left(  u,\mathcal{P}v\right)  _{w,\theta}=\sum\limits_{i}u\left(
a^{\left(  i\right)  }\right)  \overline{v\left(  a^{\left(  i\right)
}\right)  }$.
\end{proof}

\subsection{Derivation of the Riesz representers $D_{x}^{\gamma}R_{x}%
$\label{SbSect_Riesz_rep_Xpos}}

In this subsection we derive simple explicit formulas for the Riesz
representers for the evaluation functionals where $D^{\gamma}u\in
X_{w}^{\theta}$ and $\left\vert \gamma\right\vert \leq\left\lfloor \min
\kappa\right\rfloor $. These formulas are exhibited in equation \ref{p939}
below and involve order $\gamma$ derivatives of the basis function and the
cardinal basis functions associated with the Light norm.

Suppose that the weight function $w$ has properties W2 and W3 for order
$\theta$ and parameter $\kappa$. Then\ by Summary \ref{Sum_int_properties_Xwm}
the functions in $X_{w}^{\theta}$ have order $\left\lfloor \kappa\right\rfloor
$ differentiability. We endow the space $X_{w}^{\theta}$ with the Light norm
generated by a minimal unisolvent set $A=\left\{  a^{\left(  i\right)
}\right\}  _{i=1}^{M}$ and its corresponding cardinal basis $\left\{
l_{i}\right\}  _{i=1}^{M}$ of $P_{\theta-1}$. We will also make good use of
the Lagrangian interpolation operator $\mathcal{P} $ and the operator
$\mathcal{Q}=I-\mathcal{P}$ introduced in Section \ref{SbSect_P_Q}.

The approach used here is to first assume the Riesz representers exist, and
then to derive the form of the representer by substituting for $u$ two classes
of functions. First we use the cardinal basis functions associated with the
unisolvent set used to define the Light norm and then the functions in
$G\ast\widehat{S}_{\emptyset,\theta}$ discussed in Subsection
\ref{SbSect_G*F[Son]}. By this means we obtain the form of the representers.
Once we have the candidates for Riesz representers it will be shown that they
have the properties initially assumed.\medskip

We start by assuming that for all $x\in\mathbb{R}^{d}$ there exists $R_{x}\in
X_{w}^{\theta}$ such that $\left(  u,R_{x}\right)  _{w,\theta}=u\left(
x\right)  $ for all $u\in X_{w}^{\theta}$. Now if $h\in\mathbb{R}^{d}$ then
$u\left(  x+h\right)  -u\left(  x\right)  =\left(  u,R_{x+h}\right)
_{w,\theta}-\left(  u,R_{x}\right)  _{w,\theta}=\left(  u,R_{x+h}%
-R_{x}\right)  _{w,\theta}$ and so we would expect the formula $D_{k}u\left(
x\right)  =\left(  u,\frac{\partial R_{x}}{\partial x_{k}}\right)  _{w,\theta
}$ to hold, and indeed we would expect the formulas%
\begin{equation}
D^{\gamma}u\left(  x\right)  =\left(  u,D_{x}^{\gamma}R_{x}\right)
_{w,\theta},\quad\left\vert \gamma\right\vert \leq\left\lfloor \min
\kappa\right\rfloor ,\label{p921}%
\end{equation}

to hold. Next observe that we have the anti-symmetric result%
\begin{equation}
R_{x}\left(  y\right)  =\left(  R_{x},R_{y}\right)  _{w,\theta}=\overline
{\left(  R_{y},R_{x}\right)  _{w,\theta}}=\overline{R_{y}\left(  x\right)
}.\label{p938}%
\end{equation}

By definition the cardinal basis $\left\{  l_{i}\right\}  _{i=1}^{M}$ of
$P_{\theta-1}$ has real-valued coefficients and satisfies $l_{i}\left(
a^{\left(  j\right)  }\right)  =\delta_{i,j}$. Thus%
\[
l_{j}\left(  x\right)  =\left(  l_{j},R_{x}\right)  _{w,\theta}=\left\langle
l_{j},R_{x}\right\rangle _{w,\theta}+\sum\limits_{i=1}^{M}l_{j}\left(
a^{\left(  i\right)  }\right)  \overline{R_{x}\left(  a^{\left(  i\right)
}\right)  }=\overline{R_{x}\left(  a^{\left(  j\right)  }\right)  }%
=R_{x}\left(  a^{\left(  j\right)  }\right)  ,
\]

so that
\begin{equation}
l_{j}\left(  x\right)  =R_{x}\left(  a^{\left(  j\right)  }\right)
.\label{p932}%
\end{equation}

Next we will substitute functions from the space $G\ast\widehat{S}%
_{\emptyset,\theta}=\left\{  G\ast\phi:\phi\in\widehat{S}_{\emptyset,\theta
}\right\}  $ introduced in Section \ref{SbSect_G*F[Son]}. Thus, if $\phi
\in\widehat{S}_{\emptyset,\theta}$ then by Corollary \ref{Cor_f_G_convol_phi}
$G\ast\phi\in X_{w}^{\theta}$ and so
\begin{align*}
\overline{\left(  G\ast\phi\right)  \left(  x\right)  }=\overline{\left(
G\ast\phi,R_{x}\right)  _{w,\theta}} &  =\left(  R_{x},G\ast\phi\right)
_{w,\theta}\\
&  =\left\langle R_{x},G\ast\phi\right\rangle _{w,\theta}+\sum\limits_{i=1}%
^{M}R_{x}\left(  a^{\left(  i\right)  }\right)  \overline{\left(  G\ast
\phi\right)  \left(  a^{\left(  i\right)  }\right)  }\\
&  =\left\langle R_{x},G\ast\phi\right\rangle _{w,\theta}+\sum\limits_{i=1}%
^{M}l_{i}\left(  x\right)  \overline{\left(  G\ast\phi\right)  \left(
a^{\left(  i\right)  }\right)  }\\
&  =\left[  R_{x},\overline{\phi}\right]  +\sum\limits_{i=1}^{M}l_{i}\left(
x\right)  \overline{\left(  G\ast\phi\right)  \left(  a^{\left(  i\right)
}\right)  },
\end{align*}

where the last line was derived using equation \ref{p926}. Rearrange this
equation to get%
\[
\left[  R_{x},\overline{\phi}\right]  =\overline{\left(  G\ast\phi\right)
\left(  x\right)  }-\sum\limits_{i=1}^{M}l_{i}\left(  x\right)  \overline
{\left(  G\ast\phi\right)  \left(  a^{\left(  i\right)  }\right)  }.
\]

But
\begin{align*}
\overline{\left(  G\ast\phi\right)  \left(  x\right)  }=\overline{\left(
G\ast\phi\right)  \left(  x\right)  }=\left(  2\pi\right)  ^{-\frac{d}{2}%
}\overline{\left[  G\left(  x-\cdot\right)  ,\phi\right]  } &  =\left(
2\pi\right)  ^{-\frac{d}{2}}\left[  \overline{G\left(  x-\cdot\right)
},\overline{\phi}\right] \\
&  =\left(  2\pi\right)  ^{-\frac{d}{2}}\left[  G\left(  \cdot-x\right)
,\overline{\phi}\right]  ,
\end{align*}

so%
\[
\left[  R_{x},\overline{\phi}\right]  =\left(  2\pi\right)  ^{-\frac{d}{2}%
}\left[  G\left(  \cdot-x\right)  ,\overline{\phi}\right]  -\left(
2\pi\right)  ^{-\frac{d}{2}}\sum\limits_{i=1}^{M}l_{i}\left(  x\right)
\left[  G\left(  \cdot-a^{\left(  i\right)  }\right)  ,\overline{\phi}\right]
,
\]

or%
\[
\left[  R_{x},\overline{\phi}\right]  =\left(  2\pi\right)  ^{-\frac{d}{2}%
}\left[  G\left(  \cdot-x\right)  -\sum\limits_{i=1}^{M}l_{i}\left(  x\right)
G\left(  \cdot-a^{\left(  i\right)  }\right)  ,\overline{\phi}\right]  ,
\]

for all $\phi\in\widehat{S}_{\emptyset,\theta}$. But by part 1 of Theorem
\ref{Thm_So,n_and_Pnhat}, $u\in\widehat{P}_{n}\ $iff$\ \left[  u,\phi\right]
=0$ for all $\phi\in S_{\emptyset,n}$. Hence%
\[
\left(  2\pi\right)  ^{\frac{d}{2}}R_{x}\in G\left(  \cdot-x\right)
-\sum\limits_{i=1}^{M}l_{i}\left(  x\right)  G\left(  \cdot-a^{\left(
i\right)  }\right)  +P_{\theta-1},
\]

so that for each $x$, $q_{x}\in P_{\theta-1}$ and%
\begin{equation}
\left(  2\pi\right)  ^{\frac{d}{2}}R_{x}=G\left(  \cdot-x\right)
-\sum\limits_{i=1}^{M}l_{i}\left(  x\right)  G\left(  \cdot-a^{\left(
i\right)  }\right)  +q_{x}.\label{p911}%
\end{equation}

We now calculate $q_{x}$ by applying the operator $\mathcal{P}$ to \ref{p911}.
This involves using the result $\mathcal{P}R_{x}=\sum\limits_{j=1}^{M}%
l_{j}\left(  x\right)  l_{j}$ and assuming that $q_{x}\in P_{\theta-1}$, so
that $\mathcal{P}q_{x}=q_{x}$. In this way we obtain%
\begin{equation}
q_{x}=\sum\limits_{i=1}^{M}l_{i}\left(  x\right)  \mathcal{P}G\left(
\cdot-a^{\left(  i\right)  }\right)  -\mathcal{P}G\left(  \cdot-x\right)
+\left(  2\pi\right)  ^{\frac{d}{2}}\sum\limits_{j=1}^{M}l_{j}\left(
x\right)  l_{j},\label{p936}%
\end{equation}

so that%
\begin{align}
\left(  2\pi\right)  ^{\frac{d}{2}}R_{x}\left(  y\right)   & =G\left(
y-x\right)  -\sum\limits_{i=1}^{M}l_{i}\left(  x\right)  G\left(  y-a^{\left(
i\right)  }\right)  -\sum\limits_{j=1}^{M}G\left(  a^{\left(  j\right)
}-x\right)  l_{j}\left(  y\right)  +\nonumber\\
& \qquad+\sum\limits_{i,j=1}^{M}l_{i}\left(  x\right)  G\left(  a^{\left(
j\right)  }-a^{\left(  i\right)  }\right)  l_{j}\left(  y\right)  +\left(
2\pi\right)  ^{\frac{d}{2}}\sum\limits_{j=1}^{M}l_{j}\left(  x\right)
l_{j}\left(  y\right)  .\label{p939}%
\end{align}

This equation will be used to define $R_{x}\left(  y\right)  $ and noting
\ref{p921}, $D_{x}^{\gamma}R_{x}$ will be our candidate for Riesz representer
for the evaluation functional $D^{\gamma}u\left(  x\right)  $ when $\left\vert
\gamma\right\vert \leq\left\lfloor \min\kappa\right\rfloor $.

\begin{theorem}
\label{Thm_Rx(y)=QyQxG(y-x)+Sum(lj(x)lj(y))}Suppose $G$ is a basis function of
order $\theta$ generated by a weight function with properties W2 and W3.
Suppose the operator $\mathcal{P}$ is defined with the same minimal unisolvent
set used to define $R_{x}$ in \ref{p939}. Then%
\begin{equation}
R_{x}\left(  y\right)  =\left(  2\pi\right)  ^{-\frac{d}{2}}\mathcal{Q}%
_{y}\mathcal{Q}_{x}G\left(  y-x\right)  +\sum\limits_{j=1}^{M}l_{j}\left(
x\right)  l_{j}\left(  y\right)  ,\label{p940}%
\end{equation}

where $\mathcal{Q}=I-\mathcal{P}$.
\end{theorem}

\begin{proof}
The basis function is continuous by Theorems \ref{Thm_basis_smth_W3.1} or
\ref{Thm_basis_smth_W3.2_r3_pos}. Expanding $\mathcal{Q}_{x}$ and then
$\mathcal{Q}_{y}$ yields equation \ref{p939} for $R_{x}$.
\end{proof}

\begin{theorem}
\label{Thm_QyQx(y-x)^a_=_0}$\mathcal{Q}_{y}\mathcal{Q}_{x}\left(  y+x\right)
^{\alpha}=0$\ when $\left\vert \alpha\right\vert <2\theta$.
\end{theorem}

\begin{proof}
$\mathcal{Q}_{y}\mathcal{Q}_{x}\left(  y+x\right)  ^{\alpha}=\mathcal{Q}%
_{y}\mathcal{Q}_{x}\sum\limits_{\beta\leq\alpha}\binom{\alpha}{\beta}y^{\beta
}\left(  -x\right)  ^{\alpha-\beta}=\sum\limits_{\beta\leq\alpha}\binom
{\alpha}{\beta}\left(  \mathcal{Q}_{y}y^{\beta}\right)  \left(  \mathcal{Q}%
_{x}x^{\alpha-\beta}\right)  $, and since $\left\vert \alpha\right\vert
=\left\vert \beta\right\vert +\left\vert \alpha-\beta\right\vert <2\theta$ it
follows that $\left\vert \beta\right\vert <\theta$ or $\left\vert \alpha
-\beta\right\vert <\theta$ and hence that $\left(  \mathcal{Q}_{y}y^{\beta
}\right)  \left(  \mathcal{Q}_{x}x^{\alpha-\beta}\right)  =0$ for all
$\beta\leq\alpha$.
\end{proof}

\begin{theorem}
\label{Thm_Rx_unique}For each $x$, as a function:

\begin{enumerate}
\item $R_{x}$ is independent of the basis function used to define it.

\item $R_{x}$ is independent of the order of $A$.
\end{enumerate}
\end{theorem}

\begin{proof}
\textbf{Part 1} From Definition \ref{Def_basis_distrib} of a basis function,
any two basis functions differ by a polynomial of order of at most $2\theta$.
Our result now follows easily from Theorem
\ref{Thm_Rx(y)=QyQxG(y-x)+Sum(lj(x)lj(y))} and Theorem
\ref{Thm_QyQx(y-x)^a_=_0}.\medskip

\textbf{Part 2} To prove this part we use equation \ref{p940} for
$R_{x}\left(  y\right)  $. From Theorem \ref{Thm_P_op_properties} the function
$\mathcal{Q}$ is independent of the order of the unisolvent set which
generates it. Thus $\left(  2\pi\right)  ^{-\frac{d}{2}}\mathcal{Q}%
_{y}\mathcal{Q}_{x}G\left(  y-x\right)  $ is also independent of the order of
the unisolvent set which generates it. Further, $\sum\limits_{j=1}^{M}%
l_{j}(x)l_{j}\left(  y\right)  $ is independent of the order of the cardinal
basis functions. Hence $R_{x}$ is independent of the order of $A$.
\end{proof}

\begin{theorem}
\label{Thm_DRx_in_Xw}Suppose the weight function $w$ has properties W2 and W3
for order $\theta$ and parameter $\kappa$. Suppose the distribution $R_{x}\in
S^{\prime}$ is defined by equation \ref{p939}. Then for $\left\vert
\gamma\right\vert \leq\left\lfloor \min\kappa\right\rfloor $ we have
$D_{x}^{\gamma}R_{x}\in X_{w}^{\theta}$ and%
\begin{equation}
\left\vert D_{x}^{\gamma}R_{x}\right\vert _{w,\theta}^{2}=\left(  2\pi\right)
^{-d}\int\frac{\left\vert D_{x}^{\gamma}\mathcal{Q}_{x}\left(  e^{-i\left(
x,\cdot\right)  }\right)  \right\vert ^{2}}{w\left\vert \cdot\right\vert
^{2\theta}},\quad x\in\mathbb{R}^{d},\label{p918}%
\end{equation}

where Theorem \ref{Thm_|DQ(exp)|^2/(w|x|^2theta)} provides an upper bound for
the right side of \ref{p918}.
\end{theorem}

\begin{proof}
We prove that $D_{x}^{\gamma}R_{x}\in X_{w}^{\theta}$ by using definition
\ref{p948} i.e. we show that $D_{x}^{\gamma}R_{x}\in S^{\prime}$,
$\widehat{D_{x}^{\gamma}R_{x}}\in L_{loc}^{1}\left(  \mathbb{R}^{d}%
\setminus0\right)  $, $\int w\left\vert \cdot\right\vert ^{2\theta}\left\vert
\left(  D_{x}^{\gamma}R_{x}\right)  _{F}\right\vert ^{2}<\infty$ and that
$\left\vert \alpha\right\vert =\theta$ implies $\xi^{\alpha}\widehat
{D_{x}^{\gamma}R_{x}}=\xi^{\alpha}\left(  D_{x}^{\gamma}R_{x}\right)  _{F}$ on
$S$. If these criteria are satisfied then $\left\vert D_{x}^{\gamma}%
R_{x}\right\vert _{w,0}^{2}=\int w\left\vert \cdot\right\vert ^{2\theta
}\left\vert \left(  D_{x}^{\gamma}R_{x}\right)  _{F}\right\vert ^{2}$.

From \ref{p911}, for each $x\in\mathbb{R}^{d}$,
\begin{align}
\left(  2\pi\right)  ^{d/2}\widehat{D_{x}^{\gamma}R_{x}}  & =e^{-i\left(
x,\cdot\right)  }\widehat{D^{\gamma}G}-\sum\limits_{i=1}^{M}D^{\gamma}%
l_{i}\left(  x\right)  e^{-i\left\langle a^{\left(  i\right)  },\cdot
\right\rangle }\widehat{G}+\widehat{q_{x}}\nonumber\\
& =e^{-i\left(  x,\cdot\right)  }\left(  -i\xi\right)  ^{\gamma}\widehat
{G}-\sum\limits_{i=1}^{M}D^{\gamma}l_{i}\left(  x\right)  e^{-i\left\langle
a^{\left(  i\right)  },\cdot\right\rangle }\widehat{G}+\widehat{q_{x}%
}\nonumber\\
& =\left(  e^{-i\left(  x,\cdot\right)  }\left(  -i\xi\right)  ^{\gamma}%
-\sum\limits_{i=1}^{M}D^{\gamma}l_{i}\left(  x\right)  e^{-i\left\langle
a^{\left(  i\right)  },\cdot\right\rangle }\right)  \widehat{G}+\widehat
{q_{x}}\nonumber\\
& =D_{x}^{\gamma}\left(  e^{-i\left(  x,\cdot\right)  }-\sum\limits_{i=1}%
^{M}l_{i}\left(  x\right)  e^{-i\left\langle a^{\left(  i\right)  }%
,\cdot\right\rangle }\right)  \widehat{G}+\widehat{q_{x}}\nonumber\\
& =\left(  D_{x}^{\gamma}\mathcal{Q}_{x}\left(  e^{-i\left(  x,\cdot\right)
}\right)  \right)  \widehat{G}+\widehat{q_{x}},\label{q67}%
\end{align}

where $q_{x}\in P_{\theta-1}$. Since W2 implies $1/w\in L_{loc}^{1}$,
Definition \ref{Def_basis_distrib} of a basis distribution implies
$\widehat{G}=\frac{1}{w\left\vert \cdot\right\vert ^{2\theta}}\in L_{loc}%
^{1}\left(  \mathbb{R}^{d}\setminus0\right)  $ and so \ref{q67} implies that
for each $x$
\begin{equation}
\left(  D_{x}^{\gamma}R_{x}\right)  _{F}=\left(  2\pi\right)  ^{-\frac{d}{2}%
}\frac{D_{x}^{\gamma}\mathcal{Q}_{x}\left(  e^{-i\left(  x,\cdot\right)
}\right)  }{w\left\vert \cdot\right\vert ^{2\theta}}\in L_{loc}^{1}\left(
\mathbb{R}^{d}\setminus0\right)  .\label{q55}%
\end{equation}

Again by \ref{q67} and Theorem \ref{Thm_|DQ(exp)|^2/(w|x|^2theta)},%
\begin{equation}
\int w\left\vert \cdot\right\vert ^{2\theta}\left\vert \left(  D_{x}^{\gamma
}R_{x}\right)  _{F}\right\vert ^{2}=\left(  2\pi\right)  ^{-d}\int%
\frac{\left\vert D_{x}^{\gamma}\mathcal{Q}_{x}\left(  e^{-i\left(
x,\cdot\right)  }\right)  \right\vert ^{2}}{w\left\vert \cdot\right\vert
^{2\theta}}<\infty,\label{q50}%
\end{equation}

Now use $\xi$ as an \textit{action }variable. It must be shown that
$\xi^{\alpha}\widehat{D_{x}^{\gamma}R_{x}}=\xi^{\alpha}\left(  D_{x}^{\gamma
}R_{x}\right)  _{F}$ on $S$ when $\left\vert \alpha\right\vert =\theta$, which
is true if $\xi^{\alpha}\widehat{R_{x}}=\xi^{\alpha}\left(  R_{x}\right)
_{F}$ on $S$ when $\left\vert \alpha\right\vert =\theta$. But from \ref{q67}
it follows that $\xi^{\alpha}\widehat{R_{x}}=\xi^{\alpha}\mathcal{Q}%
_{x}\left(  e^{-i\left(  x,\xi\right)  }\right)  \widehat{G}+\xi^{\alpha
}\widehat{q_{x}}=\xi^{\alpha}\mathcal{Q}_{x}\left(  e^{-i\left(  x,\xi\right)
}\right)  \widehat{G}$ since $q_{x}\in P_{\theta-1}$. From part 4 of Theorem
\ref{Thm_Q(exp)_loc_property}, if $\left\vert \alpha\right\vert =\theta$ then
for each $x$,

$\left(  i\xi\right)  ^{\alpha}\mathcal{Q}_{x}\left(  e^{i\left(
x,\xi\right)  }\right)  \in C_{\emptyset,2\theta}^{\infty}\cap C_{BP}^{\infty
}$ and so $\psi\in S$ implies $\left(  i\xi\right)  ^{\alpha}\mathcal{Q}%
_{x}\left(  e^{i\left(  x,\xi\right)  }\right)  \psi\in S_{\emptyset,2\theta}
$ by Theorem \ref{Thm_product_of_Co,k_funcs}. Thus by the definition of a
basis distribution
\begin{align*}
\left[  \xi^{\alpha}\widehat{R_{x}},\psi\right]  =\left[  \xi^{\alpha
}\mathcal{Q}_{x}\left(  e^{-i\left(  x,\xi\right)  }\right)  \widehat{G}%
,\psi\right]  =\left[  \widehat{G},\xi^{\alpha}\mathcal{Q}_{x}\left(
e^{-i\left(  x,\xi\right)  }\right)  \psi\right]   &  =\int\frac{\xi^{\alpha
}\mathcal{Q}_{x}\left(  e^{-i\left(  x,\xi\right)  }\right)  }{w\left(
\xi\right)  \left\vert \xi\right\vert ^{2\theta}}\psi\left(  \xi\right)
d\xi\\
&  =\int\xi^{\alpha}\left(  R_{x}\right)  _{F}\psi,
\end{align*}

and hence $\xi^{\alpha}\widehat{R_{x}}=\xi^{\alpha}\left(  R_{x}\right)  _{F}$
on $S$. We can now conclude that $D_{x}^{\gamma}R_{x}\in X_{w}^{\theta}$, that
\ref{q50} is true and so \ref{p918} is proven.
\end{proof}

\begin{corollary}
\label{Cor_Thm_DRx_in_Xw}Suppose the weight function $w$ has properties W2 and
W3 for order $\theta$ and smoothness parameter $\kappa$. Then for $\left\vert
\delta\right\vert ,\left\vert \gamma\right\vert \leq\left\lfloor \min
\kappa\right\rfloor $ we have%
\[
\left\langle D_{y}^{\delta}R_{y},D_{x}^{\gamma}R_{x}\right\rangle _{w,\theta
}=\left(  2\pi\right)  ^{-d}\int\frac{D_{y}^{\delta}\mathcal{Q}_{y}\left(
e^{-i\left(  y,\cdot\right)  }\right)  \overline{D_{x}^{\gamma}\mathcal{Q}%
_{x}\left(  e^{-i\left(  x,\cdot\right)  }\right)  }}{w\left\vert
\cdot\right\vert ^{2\theta}},\quad x,y\in\mathbb{R}^{d},
\]

\end{corollary}

\begin{proof}
Write $\left\langle D_{y}^{\delta}R_{y},D_{x}^{\gamma}R_{x}\right\rangle
_{w,\theta}=\int w\left\vert \cdot\right\vert ^{2\theta}\left(  D_{y}^{\delta
}R_{y}\right)  _{F}\overline{\left(  D_{x}^{\gamma}R_{x}\right)  _{F}}$ and
then use \ref{q55}.
\end{proof}

Now we prove the important result that when $\left\vert \gamma\right\vert
\leq\min\kappa$, $D_{x}^{\gamma}R_{x}$ is the Riesz representer for the
evaluation functional $u\rightarrow D^{\gamma}u\left(  x\right)  $.

\begin{theorem}
\label{Thm_Rx_is_Riesz_rep}If $u\in X_{w}^{\theta}$ and $\left\vert
\gamma\right\vert \leq\min\kappa$ then
\[
D^{\gamma}u\left(  x\right)  =\left(  u,D_{x}^{\gamma}R_{x}\right)
_{w,\theta},\quad x\in\mathbb{R}^{d},
\]

where $R_{x}$ and the Light norm \ref{p917} are defined using the same minimal
$\theta$-unisolvent set $A=\left\{  a^{\left(  k\right)  }\right\}  _{k=1}%
^{M}$.
\end{theorem}

\begin{proof}
Suppose $R_{x}$ and the Light norm are defined using the minimal unisolvent
set $A=\left\{  a^{\left(  k\right)  }\right\}  _{k=1}^{M}$ and that the
corresponding cardinal basis $\left\{  l_{k}\right\}  _{k=1}^{M}$ of
$P_{\theta-1}$. From the definition of the Light norm%
\[
\left(  u,D_{x}^{\gamma}R_{x}\right)  _{w,\theta}=\int w\left\vert
\cdot\right\vert ^{2\theta}u_{F}\overline{\left(  D_{x}^{\gamma}R_{x}\right)
_{F}}+\sum\limits_{k=1}^{M}u\left(  a^{\left(  k\right)  }\right)
\overline{\left(  D_{x}^{\gamma}R_{x}\right)  \left(  a^{\left(  k\right)
}\right)  },
\]

where the function $\left(  D_{x}^{\gamma}R_{x}\right)  _{F}$ is given by
\ref{q55}. Hence
\[
\int w\left\vert \cdot\right\vert ^{2\theta}u_{F}\text{ }\overline{\left(
D_{x}^{\gamma}R_{x}\right)  _{F}}=\left(  2\pi\right)  ^{-\frac{d}{2}}\int
w\left\vert \cdot\right\vert ^{2\theta}u_{F}\frac{D_{x}^{\gamma}%
\mathcal{Q}_{x}\left(  e^{i\left(  x,\cdot\right)  }\right)  }{w\left\vert
\cdot\right\vert ^{2\theta}}=\left(  2\pi\right)  ^{-\frac{d}{2}}\int
u_{F}D_{x}^{\gamma}\mathcal{Q}_{x}\left(  e^{i\left(  x,\cdot\right)
}\right)  ,
\]

and by \ref{p932}, $D^{\gamma}l_{k}\left(  x\right)  =\left(  D_{x}^{\gamma
}R_{x}\right)  \left(  a^{\left(  k\right)  }\right)  $ so that%
\[
\sum\limits_{k=1}^{M}u\left(  a^{\left(  k\right)  }\right)  \overline{\left(
D_{x}^{\gamma}R_{x}\right)  \left(  a^{\left(  k\right)  }\right)  }%
=\sum\limits_{k=1}^{M}u\left(  a^{\left(  k\right)  }\right)  D^{\gamma}%
l_{k}\left(  x\right)  =\left(  D^{\gamma}\mathcal{P}u\right)  \left(
x\right)  .
\]

Thus%
\[
\left(  u,D_{x}^{\gamma}R_{x}\right)  _{w,\theta}=\left(  2\pi\right)
^{-\frac{d}{2}}\int u_{F}D_{x}^{\gamma}\mathcal{Q}_{x}\left(  e^{i\left(
x,\cdot\right)  }\right)  +\left(  D^{\gamma}\mathcal{P}u\right)  \left(
x\right)  ,
\]

and comparison with equation \ref{p920} of Theorem
\ref{Thm_inv_Fourier_tran_func} for $D^{\gamma}u$ implies $D^{\gamma}u\left(
x\right)  =\left(  u,D_{x}^{\gamma}R_{x}\right)  _{w,\theta}$ as required.
\end{proof}

Showing that $X_{w}^{\theta}$ is a reproducing kernel Hilbert space now
becomes very simple.

\begin{corollary}
\label{Cor_Thm_Rx_is_Riesz_rep}Suppose the weight function $w$ has properties
W2 and W3 for order $\theta$ and parameter $\kappa$. Then $X_{w}^{\theta}$ is
a reproducing kernel Hilbert space when endowed with the Light norm.
\end{corollary}

\begin{proof}
By Theorem \ref{Thm_DRx_in_Xw} $R_{x}\in X_{w}^{\theta}$. Hence by Theorem
\ref{Thm_Rx_is_Riesz_rep} with $\gamma=0$ we have $\left\vert u\left(
x\right)  \right\vert \leq\left\vert \left(  u,R_{x}\right)  _{w,\theta
}\right\vert \leq\left\Vert R_{x}\right\Vert _{w,\theta}\left\Vert
u\right\Vert _{w,\theta}$ and $\left\Vert R_{x}\right\Vert _{w,\theta}$ is a
finite constant for each $x$. Hence $X_{w}^{\theta}$ is a reproducing kernel
Hilbert space e.g. Theorem III.9.1, Yosida \cite{Yosida58}.
\end{proof}

\begin{theorem}
\label{Thm_Riesz_property}\textbf{Properties of }$D_{x}^{\gamma}R_{x}$

Suppose the weight function $w$ satisfies properties W2 and W3 for order
$\theta$ and $\kappa$. Suppose $R_{x}$ is defined using the minimal unisolvent
set $A=\left\{  a^{\left(  j\right)  }\right\}  _{j=1}^{M}$. Then for all
$x,y\in\mathbb{R}^{d}$ and $\left\vert \gamma\right\vert \leq\left\lfloor
\min\kappa\right\rfloor $:

\begin{enumerate}
\item $R_{x}\left(  y\right)  =\overline{R_{y}\left(  x\right)  }$.

\item For $j=1,\ldots,M$, $\left(  D_{x}^{\gamma}R_{x}\right)  \left(
a^{\left(  j\right)  }\right)  =D^{\gamma}l_{j}\left(  x\right)  $, and each
cardinal basis polynomial $l_{j}$ has real coefficients.

\item $D^{\gamma}\mathcal{Q}u\left(  x\right)  =\left\langle u,D_{x}^{\gamma
}R_{x}\right\rangle _{w,\theta}$, $u\in X_{w}^{\theta}$.

\item $\mathcal{P}D_{x}^{\gamma}R_{x}=\sum\limits_{j=1}^{M}D^{\gamma}%
l_{j}\left(  x\right)  l_{j}=D_{x}^{\gamma}\mathcal{P}R_{x}$ and
$\mathcal{Q}D_{x}^{\gamma}R_{x}=D_{x}^{\gamma}\mathcal{Q}R_{x}$.

\item The Riesz representer is unique.
\end{enumerate}
\end{theorem}

\begin{proof}
\textbf{Parts 1} and \textbf{2} By Theorem \ref{Thm_Rx_is_Riesz_rep}
$D_{x}^{\gamma}R_{x}$ is the Riesz representer of the evaluation functional
$u\rightarrow D^{\gamma}u\left(  x\right)  $ and so equations \ref{p938} and
\ref{p932} are valid. By Definition \ref{Def_cardinal_basis} the cardinal
basis polynomials have real coefficients.\smallskip

\textbf{Part 3} If $u\in X_{w}^{\theta}$ then by part 1 of Theorem
\ref{Thm_Light_norm_property} $D^{\gamma}\mathcal{Q}u\left(  x\right)
=\left(  \mathcal{Q}u,D_{x}^{\gamma}R_{x}\right)  _{w,\theta}=\left\langle
u,D_{x}^{\gamma}R_{x}\right\rangle _{w,\theta}$.\smallskip

\textbf{Part 4} Using part 2%
\[
\mathcal{P}D_{x}^{\gamma}R_{x}=\sum\limits_{j=1}^{M}\left(  D_{x}^{\gamma
}R_{x}\right)  \left(  a^{\left(  j\right)  }\right)  l_{j}=\sum
\limits_{j=1}^{M}\left(  D^{\gamma}l_{j}\right)  \left(  x\right)  l_{j}%
=D_{x}^{\gamma}\sum\limits_{j=1}^{M}l_{j}\left(  x\right)  l_{j}=D_{x}%
^{\gamma}\mathcal{P}R_{x},
\]

so that%
\[
\mathcal{Q}D_{x}^{\gamma}R_{x}=D_{x}^{\gamma}R_{x}-\mathcal{P}D_{x}^{\gamma
}R_{x}=D_{x}^{\gamma}R_{x}-D_{x}^{\gamma}\mathcal{P}R_{x}=D_{x}^{\gamma
}\mathcal{Q}R_{x}.
\]
\smallskip

\textbf{Part 5} Suppose two representers exist, say $D_{x}^{\gamma}R_{x}$ and
$S_{x}^{\left(  \gamma\right)  }$. Then for all $u\in X_{w}^{\theta}$ and
$x\in\mathbb{R}^{d}$, $\left(  u,D_{x}^{\gamma}R_{x}\right)  _{w,\theta
}=\left(  u,S_{x}^{\left(  \gamma\right)  }\right)  _{w,\theta}=D^{\gamma
}u\left(  x\right)  $ and hence

$\left(  u,S_{x}^{\left(  \gamma\right)  }-D_{x}^{\gamma}R_{x}\right)
_{w,\theta}=0$ which implies $S_{x}^{\left(  \gamma\right)  }=D_{x}^{\gamma
}R_{x}$.
\end{proof}

\subsection{The semi-Riesz representer $r_{x}=\mathcal{Q}R_{x}$%
\label{SbSect_rx}}

Here we will introduce what I will call the \textit{semi-Riesz representer}
$r_{x}=\mathcal{Q}R_{x}$ for the evaluation functional $u\rightarrow u\left(
x\right)  $. This terminology is based on the equation $\mathcal{Q}u\left(
x\right)  =\left\langle u,r_{x}\right\rangle _{w,\theta}$, $u\in X_{w}%
^{\theta}$, derived in part 7 of the next theorem. A bound for the the
function $r_{x}\left(  x\right)  $ will be obtained in Subsection
\ref{SbSect_int_bound_rx(x)} which will then be used to estimate the pointwise
rate of convergence of the interpolant.

\begin{theorem}
\label{Thm_rx(y)_properties}Suppose the function $r_{x}=\mathcal{Q}R_{x}$ is
defined using the minimal unisolvent set $\left\{  a^{\left(  k\right)
}\right\}  _{k=1}^{M}$ and corresponding cardinal basis $\left\{
l_{k}\right\}  _{k=1}^{M}$ of $P_{\theta-1}$. Then $r_{x}$ has the properties:

\begin{enumerate}
\item $r_{x}\left(  y\right)  =R_{x}\left(  y\right)  -\sum\limits_{j=1}%
^{M}l_{j}(x)l_{j}\left(  y\right)  $.

\item $r_{x}\left(  y\right)  =\overline{r_{y}\left(  x\right)  }$.

\item $r_{x}\left(  a^{\left(  i\right)  }\right)  =r_{a^{\left(  i\right)  }%
}\left(  x\right)  =0$ for all $i$.

\item $\mathcal{P}r_{x}=0$ and $\mathcal{Q}r_{x}=r_{x}$.

\item $r_{x}\left(  y\right)  =\left\langle r_{x},r_{y}\right\rangle
_{w,\theta}$.

\item $r_{x}\left(  y\right)  =\left(  2\pi\right)  ^{-\frac{d}{2}}%
\mathcal{Q}_{y}\mathcal{Q}_{x}G\left(  y-x\right)  $.

\item $\mathcal{Q}u\left(  x\right)  =\left\langle u,r_{x}\right\rangle
_{w,\theta}$ when $u\in X_{w}^{\theta}$.
\end{enumerate}
\end{theorem}

\begin{proof}
\textbf{Part 1} Apply $\mathcal{Q}_{y}$ to equation \ref{p940}. \textbf{Part 2
}Follows from part 1 since $R_{x}\left(  y\right)  =\overline{R_{y}\left(
x\right)  }$ and the $l_{k}$ are real valued. \textbf{Part 3} Part 2 of
Theorem \ref{Thm_Riesz_property} with $\gamma=0$ implies $R_{x}\left(
a^{\left(  i\right)  }\right)  =l_{i}\left(  x\right)  $. Parts 1 and 2 now
give this result. \textbf{Part 4} Use part 3. \textbf{Part 5} Part 3 of
Theorem \ref{Thm_Riesz_property} with $u=r_{x}$ and $\gamma=0$ gives
$\mathcal{Q}r_{x}\left(  y\right)  =\left\langle r_{x},R_{y}\right\rangle
_{w,\theta}$. Since $\mathcal{Q}r_{x}=r_{x}$, Theorem
\ref{Thm_Light_norm_property} implies $r_{x}\left(  y\right)  =\mathcal{Q}%
r_{x}\left(  y\right)  =\left\langle r_{x},R_{y}\right\rangle _{w,\theta
}=\left\langle r_{x},\mathcal{Q}R_{y}\right\rangle _{w,\theta}=\left\langle
r_{x},r_{y}\right\rangle _{w,\theta}$.

\textbf{Part 6} Apply $\mathcal{Q}_{y}$ to formula \ref{p940}. \textbf{Part 7}
follows from part 3 of Theorem \ref{Thm_Riesz_property} with $\gamma=0$, and
then use part 4.
\end{proof}

\begin{remark}
The projection properties of the operators $\mathcal{P}$ and $\mathcal{Q}$ can
be used to show that $X_{w}^{\theta}=\mathcal{Q}\left(  X_{w}^{\theta}\right)
\oplus P_{\theta-1}$ and that $\mathcal{Q}\left(  X_{w}^{\theta}\right)  $ is
a Hilbert space when endowed with the semi-inner product $\left\langle
\cdot,\cdot\right\rangle _{w,\theta}$. In fact, on $\mathcal{Q}\left(
X_{w}^{\theta}\right)  $, $r_{x}$ is the Riesz representer of the functional
$f\rightarrow f\left(  x\right)  $ and $\mathcal{Q}\left(  X_{w}^{\theta
}\right)  $ is a reproducing kernel Hilbert space.
\end{remark}

\section{The basis function and reproducing kernel
matrices\label{Sect_mat_basis_reprod_ker}}

In this section we will introduce the \textit{basis function matrix} and the
\textit{reproducing kernel matrix}. Together with the unisolvency matrices
these\textit{\ }matrices will be used to construct matrix equations for the
basis function interpolants of this document.

\subsection{The basis function matrix $G_{X,X}$\label{SbSect_mat_basis}}

The \textit{basis function matrix }and the unisolvency matrix will be used to
construct the block matrix equation \ref{q64} for the variational
interpolation in Subsection \ref{SbSect_mat_eqn_interpol}.

\begin{definition}
\label{Def_matrices_from G}\textbf{The basis function matrix}

Let $X=\left\{  x^{\left(  n\right)  }\right\}  _{n=1}^{N}$ be $N$ distinct
points in $\mathbb{R}^{d}$ and suppose $G$ is a (continuous) basis function.

Then the $N\times N$ basis function matrix $G_{X,X}$ is defined by
\[
G_{X,X}=\left(  G\left(  x^{\left(  i\right)  }-x^{\left(  j\right)  }\right)
\right)  .
\]

\end{definition}

\subsection{The reproducing kernel matrix $R_{X,X}$\label{SbSect_mat_eval_rep}%
}

The \textit{reproducing kernel matrix }is derived from the Riesz representer
of the evaluation functional $u\rightarrow u\left(  x\right)  $ discussed in
Section \ref{SbSect_Riesz_rep_Xpos}. I call it the reproducing kernel matrix
because the existence of this special evaluation functional means that the
space $X_{w}^{\theta}$ is a reproducing kernel Hilbert space when endowed with
the Light norm. This matrix is closely linked to the basis function matrix of
the previous subsection will be used to construct a matrix equation for the
variational interpolation problems of Section \ref{Sect_interpol_prob}. In the
next chapter we will study the relationships between the reproducing kernel
matrix and the basis function matrix and derive a matrix equation for the
Exact smoother using the reproducing kernel matrix.

\begin{definition}
\label{Def_Matrices_from_Rx}\textbf{The reproducing kernel matrix }$R_{X,X}$

Suppose a basis function has order $\theta$ and suppose that $R_{x}(y)$ is the
Riesz representer of the evaluation functional $u\rightarrow u\left(
x\right)  $, $u\in X_{w}^{\theta}$ defined by \ref{p939}. Also suppose that
$X=\left\{  x^{\left(  n\right)  }\right\}  _{n=1}^{N}$ is a set of distinct
unisolvent data points in $\mathbb{R}^{d}$. Then define the $N\times N$ matrix
$R_{X,X}$ by
\[
R_{X,X}=\left(  R_{x^{\left(  j\right)  }}(x^{\left(  i\right)  })\right)  .
\]

\end{definition}

\begin{theorem}
\label{Thm_int_Rxx_property}Suppose the conditions imposed in the definition
of the reproducing kernel matrix $R_{X,X}$ hold. Then the following are true:

\begin{enumerate}
\item $R_{X,X}=\left(  (R_{x^{(j)}},R_{x^{(i)}})_{w,\theta}\right)  $ so it is
a Gram matrix i.e. has inner product elements $\left(  u_{i},u_{j}\right)
_{w,\theta}$.

\item The matrix $R_{X,X}$ is Hermitian, regular and positive definite.

\item The functions $\left\{  R_{x^{(i)}}\right\}  _{i=1}^{N}$ are linearly independent.
\end{enumerate}
\end{theorem}

\begin{proof}
\textbf{Part 1} By definition of $R_{x}$, $R_{X,X}=\left(  R_{x^{(j)}}\left(
x^{(i)}\right)  \right)  =\left(  R_{x^{(j)}},R_{x^{(i)}}\right)  _{w,\theta}%
$.\medskip

\textbf{Part 2} These are elementary properties of a Gram matrix.\medskip

\textbf{Part 3} This is a direct consequence of the regularity of $R_{X,X}$.
\end{proof}

\section{The basis function spaces $\protect\overset{\cdot}{W}_{G,X}$ and
$W_{G,X}$\label{Sect_Wgx}}

The importance of the finite dimensional spaces $\overset{\cdot}{W}_{G,X}$ and
$W_{G,X}$ is that they will contain the basis function interpolants of this
document and the basis function smoothers of the succeeding documents. Here
$X$ is the independent data and will be assumed to be unisolvent. The finite
dimensionality of these spaces means that matrix equations can be constructed
for the interpolant and the smoothers.

\begin{definition}
\label{Def_Wgx}\textbf{The basis function spaces }$W_{G,X}$\textbf{\ and
}$\overset{\cdot}{W}_{G,X}$

Suppose the weight function $w$ has properties W2 and W3 for order $\theta$
and smoothness parameter $\kappa$. Then the basis distributions of order
$\theta$ are continuous functions and let $G$ be a basis function. Let
$X=\left\{  x^{(i)}\right\}  _{i=1}^{N}$ be a $\theta$-unisolvent set of
distinct points in $\mathbb{R}^{d}$ and set $M=\dim P_{\theta-1}$. Next choose
a (real) basis $\left\{  p_{j}\right\}  _{j=1}^{M}$ of $P_{\theta-1}$ and
calculate the unisolvency matrix $P_{X}=\left(  p_{j}\left(  x^{(i)}\right)
\right)  $. We can now define
\begin{align}
\overset{\cdot}{W}_{G,X} &  =\left\{  \sum_{i=1}^{N}v_{i}G\left(
x-x^{(i)}\right)  :v_{i}\in\mathbb{C}\text{ }and\text{ }P_{X}^{T}v=0\right\}
,\label{q51}\\
W_{G,X} &  =\overset{\cdot}{W}_{G,X}+P_{\theta-1}.\nonumber
\end{align}

Clearly $\overset{\cdot}{W}_{G,X}$ and $W_{G,X}$ are vector spaces. When
convenient, functions in spaces of the form $\overset{\cdot}{W}_{G,X}$ will be
written, $f_{v}\left(  x\right)  =\sum\limits_{i=1}^{N}v_{i}G\left(
x-x^{(i)}\right)  $.

In words, $\overset{\cdot}{W}_{G,X}$ is a \textbf{subspace} of the complex
span of the $X-$\textit{translated} basis functions and $W_{G,X}$ is a
\textbf{subspace} of the complex span of the $X-$translated basis functions
plus the polynomials of at most basis function order. Sometimes we will say
data-translated basis functions.
\end{definition}

The next theorem shows that the set $\overset{\cdot}{W}_{G,X}$ is independent
of the order of the points in $X$ and of the basis of $P_{\theta-1}$ used to
calculate $P_{X}$.

\begin{theorem}
\label{Thm_Wgx_defin_indep_of_Px_chosen}Let $\left\{  p_{i}\right\}  $ and
$\left\{  q_{i}\right\}  $ be two bases for the polynomials $P_{\theta-1}$,
and set $P_{X}=\left(  p_{j}\left(  x^{(i)}\right)  \right)  $ and
$Q_{X}=\left(  q_{j}\left(  x^{(i)}\right)  \right)  $. Then
\begin{equation}
\left\{  \sum_{i=1}^{N}v_{i}G\left(  x-x^{(i)}\right)  :P_{X}^{T}v=0\right\}
=\left\{  \sum_{i=1}^{N}\mu_{i}G\left(  x-x^{(i)}\right)  :Q_{X}^{T}%
\mu=0\right\}  ,\label{q68}%
\end{equation}

and $\overset{\cdot}{W}_{G,X}$ is independent of the basis for $P_{\theta-1}$
used to define $P_{X}$.
\end{theorem}

\begin{proof}
From part 3 of Theorem \ref{Thm_Px_properties} there exists a regular matrix
$R$ such that $P_{X}=Q_{X}R^{T}$. Hence $P_{X}^{T}=RQ_{X}^{T}$ and then
equation \ref{q68} easily follows.
\end{proof}

Now recall Definition \ref{Def_cardinal_basis} regarding permutation theory.

\begin{theorem}
\label{Thm_Wgx_reorder}As a set $\overset{\cdot}{W}_{G,X}$ is independent of
the ordering of $X$ and independent of the basis function $G$ used to define it.
\end{theorem}

\begin{proof}
Using the notation $G_{x,X}=\left(  G\left(  x-x^{(j)}\right)  \right)  $ we
have%
\[
\overset{\cdot}{W}_{G,X}=\left\{  \sum_{i=1}^{N}v_{i}G\left(  x-x^{(i)}%
\right)  :P_{X}^{T}v=0\right\}  =\left\{  G_{x,X}v:P_{X}^{T}v=0\right\}  ,
\]

where $P_{X}^{T}=\left(  p_{i}\left(  x^{\left(  j\right)  }\right)  \right)
=\left(
\begin{array}
[c]{cccc}%
p_{1}\left(  x^{\left(  1\right)  }\right)  & p_{1}\left(  x^{\left(
2\right)  }\right)  & \cdots & p_{1}\left(  x^{\left(  N\right)  }\right) \\
p_{2}\left(  x^{\left(  1\right)  }\right)  & p_{2}\left(  x^{\left(
2\right)  }\right)  & \cdots & p_{2}\left(  x^{\left(  N\right)  }\right) \\
\vdots & \vdots & \ddots & \vdots\\
p_{M}\left(  x^{\left(  1\right)  }\right)  & p_{M}\left(  x^{\left(
2\right)  }\right)  & \cdots & p_{M}\left(  x^{\left(  N\right)  }\right)
\end{array}
\right)  $.

Suppose re-ordering $X$ uses the permutation $\pi$. Then re-ordering a row
vector involves right multiplication by the permutation matrix $\Pi^{T}$ so
that $G_{x,\pi\left(  X\right)  }=G_{x,X}\Pi^{T}$, $P_{\pi\left(  X\right)
}^{T}=P_{X}^{T}\Pi^{T}$ and hence%
\[
\overset{\cdot}{W}_{G,\pi\left(  X\right)  }=\left\{  G_{x,\pi\left(
X\right)  }:P_{\pi\left(  X\right)  }^{T}v=0\right\}  =\left\{  G_{x,\pi
\left(  X\right)  }\Pi^{T}v:P_{X}^{T}\Pi^{T}v=0\right\}  =\overset{\cdot
}{W}_{G,X}.
\]
\medskip

Set-wise independence of the basis function will hold if $p\in P_{\theta-1}$ implies

$\sum_{i=1}^{N}v_{i}p\left(  x-x^{(i)}\right)  =0$ and clearly we need only
show this for the basis functions $\left\{  x^{\alpha}\right\}  _{\left\vert
\alpha\right\vert <\theta}$. Since $\left\{  x^{\alpha}\right\}  _{\left\vert
\alpha\right\vert <\theta}$ is a basis, by Theorem
\ref{Thm_Wgx_defin_indep_of_Px_chosen} $P_{X}^{T}v=0$ iff $\sum_{i=1}^{N}%
v_{i}\left(  x^{(i)}\right)  ^{\alpha}=0$ for $\left\vert \alpha\right\vert
<\theta$. Thus
\[
\sum_{i=1}^{N}v_{i}\left(  x-x^{(i)}\right)  ^{\alpha}=\sum_{i=1}^{N}v_{i}%
\sum_{\beta\leq\alpha}\binom{\alpha}{\beta}x^{\beta}\left(  x^{(i)}\right)
^{\alpha-\beta}=\sum_{\beta\leq\alpha}\binom{\alpha}{\beta}x^{\beta}\sum
_{i=1}^{N}v_{i}\left(  x^{(i)}\right)  ^{\alpha-\beta}=0,
\]

confirming set-wise independence.
\end{proof}

\begin{lemma}
\label{Lem_sum(v.exp(-ixy))}Suppose $X=\left\{  x^{\left(  k\right)
}\right\}  _{k=1}^{N}$ is $\theta$-unisolvent and $P_{X}$ is a unisolvency
matrix of order $\theta$ (see Definition \ref{Def_unisolv_matrix_Px}).
Further, suppose $v=\left(  v_{k}\right)  \in\mathbb{R}^{N}$ satisfies
$P_{X}^{T}v=0$. Now define the function $a_{v}:\mathbb{R}^{d}\rightarrow
\mathbb{C}$ by
\begin{equation}
a_{v}\left(  z\right)  =\sum\limits_{k=1}^{N}v_{k}e^{-i\left(  x^{\left(
k\right)  },z\right)  },\quad z\in\mathbb{R}^{d}.\label{p943}%
\end{equation}

Then the function $a_{v}$ has the following properties:

\begin{enumerate}
\item $a_{v}\in C_{\emptyset,\theta}^{\infty}\cap C_{B}^{\infty}$.

\item If $\phi\in S$ then $\left\vert a_{v}\right\vert ^{2}\phi\in S$ and
\[
\widehat{\left\vert a_{v}\right\vert ^{2}\phi}=\sum\limits_{j,k=1}^{N}%
v_{j}v_{k}\widehat{\phi}\left(  \cdot-\left(  x^{\left(  k\right)
}-x^{\left(  j\right)  }\right)  \right)  .
\]

\item $a_{v}\left(  z\right)  =0$ a.e. implies $v_{k}=0$ for all $k$.
\end{enumerate}
\end{lemma}

\begin{proof}
\textbf{Part 1} Since
\[
D^{\beta}a_{v}\left(  z\right)  =D^{\beta}\sum_{k=1}^{N}v_{k}e^{-i\left(
x^{\left(  k\right)  },z\right)  }=\sum_{k=1}^{N}v_{k}\left(  -ix^{\left(
k\right)  }\right)  ^{\beta}e^{-i\left(  x^{\left(  k\right)  },z\right)  },
\]

it is clear that all derivatives are bounded and hence $a_{v}\in C_{B}%
^{\infty}$. If $\left\vert \beta\right\vert <\theta$ then
\[
D^{\beta}a_{v}\left(  0\right)  =\sum_{k=1}^{N}v_{k}\left(  -ix^{\left(
k\right)  }\right)  ^{\beta}=\left(  -i\right)  ^{\left\vert \beta\right\vert
}\sum_{k=1}^{N}v_{k}\left(  x^{\left(  k\right)  }\right)  ^{\beta}=0,
\]

since $P_{X}^{T}v=0.\medskip$

\textbf{Part 2}
\begin{align*}
\left\vert a_{v}\left(  x\right)  \right\vert ^{2}\phi\left(  x\right)
=\left\vert \sum_{j=1}^{N}v_{j}e^{-i\left(  x^{\left(  j\right)  },x\right)
}\right\vert ^{2}\phi\left(  x\right)   &  =\left(  \sum_{j=1}^{N}%
v_{j}e^{-i\left(  x^{\left(  j\right)  },x\right)  }\right)  \left(
\overline{\sum_{k=1}^{N}v_{k}e^{-i\left(  x^{\left(  k\right)  },x\right)  }%
}\right)  \phi\left(  x\right) \\
&  =\sum_{j,k=1}^{N}v_{j}v_{k}e^{-i\left(  x^{\left(  j\right)  }-x^{\left(
k\right)  },x\right)  }\phi\left(  x\right)  .
\end{align*}

Therefore
\begin{align*}
\widehat{\left\vert a_{v}\right\vert ^{2}\phi}\left(  \eta\right)
=\sum_{j,k=1}^{N}v_{j}v_{k}\left[  e^{-i\left(  x^{\left(  j\right)
}-x^{\left(  k\right)  },x\right)  }\phi\left(  x\right)  \right]  ^{\wedge
}\left(  \eta\right)   &  =\sum_{j,k=1}^{N}v_{j}v_{k}\widehat{\phi}\left(
x^{\left(  j\right)  }-x^{\left(  k\right)  }+\eta\right) \\
&  =\sum_{j,k=1}^{N}v_{j}v_{k}\widehat{\phi}\left(  \eta-\left(  x^{\left(
k\right)  }-x^{\left(  j\right)  }\right)  \right)  .
\end{align*}
\smallskip

\textbf{Part 3} Let $\Delta=\left\{  x^{\left(  k\right)  }-x^{\left(
j\right)  }:1\leq j,k\leq N,\text{ }k\neq j\right\}  $. Then $0\notin\Delta$
and so $\delta=\operatorname*{dist}\left(  0;\Delta\right)  >0$. Next let
$\omega_{\varepsilon}$ be the standard `cap-shaped' distribution test function%
\[
\omega_{\varepsilon}\left(  \xi\right)  =\left\{
\begin{array}
[c]{ll}%
C_{\varepsilon}\exp\left(  -\frac{\varepsilon^{2}}{\varepsilon^{2}-\left\vert
\xi\right\vert ^{2}}\right)  , & \left\vert \xi\right\vert \leq\varepsilon,\\
0, & \left\vert \xi\right\vert >\varepsilon,
\end{array}
\right.
\]

with $C_{\varepsilon}$ chosen so that $\omega_{\varepsilon}\left(  0\right)
=1$. Now $\omega_{\varepsilon}$ has support $\overline{B\left(  0;\varepsilon
\right)  }$ and if we set $\widehat{\phi}\left(  \xi\right)  =\omega_{\delta
}\left(  \xi\right)  $ it follows that $\widehat{\phi}\left(  0\right)  =1$
and $\widehat{\phi}\left(  \xi\right)  =0$ when $\xi\in\Delta$. Thus $a_{v}=0$
a.e. implies $\widehat{\left\vert a_{v}\right\vert ^{2}\phi}=0$ a.e. and so,
$0=\sum\limits_{j,k=1}^{N}v_{j}v_{k}\widehat{\phi}\left(  \xi-\left(
x^{\left(  k\right)  }-x^{\left(  j\right)  }\right)  \right)  =\sum
\limits_{k=1}^{N}v_{k}^{2}$ and therefore $v_{k}=0$ for all $k$.
\end{proof}

Next we obtain expressions for the semi-inner product and seminorm of
functions in $\overset{\cdot}{W}_{G,X}$.

\begin{theorem}
\label{Thm_Wgx_properties_1}Suppose $w$ is a weight function with properties
W2.1 and W3 for some order $\theta\geq1$. Suppose $v=\left(  v_{i}\right)
_{i=1}^{N}$ is a complex vector satisfying $P_{X}^{T}v=0$ where $P_{X}$ is the
unisolvency matrix. Then we have the following results:

\begin{enumerate}
\item If the functions $f_{v}$ and $a_{v}$ are defined by
\begin{equation}
f_{v}\left(  x\right)  =\sum\limits_{k=1}^{N}v_{k}G\left(  x-x^{\left(
k\right)  }\right)  \in\overset{\cdot}{W}_{G,X},\text{\quad}a_{v}\left(
\xi\right)  =\sum\limits_{j=1}^{N}v_{j}e^{-i\left(  x^{\left(  j\right)  }%
,\xi\right)  },\label{q58}%
\end{equation}

then $f_{v}\in X_{w}^{\theta}$ and
\begin{equation}
\left\vert f_{v}\right\vert _{w,\theta}^{2}=\int\frac{\left\vert
a_{v}\right\vert ^{2}}{w\left\vert \cdot\right\vert ^{2\theta}}=\left(
2\pi\right)  ^{\frac{d}{2}}\sum_{j,k=1}^{N}v_{j}\overline{v_{k}}G\left(
x^{\left(  j\right)  }-x^{\left(  k\right)  }\right)  =\left(  2\pi\right)
^{\frac{d}{2}}v^{T}G_{X,X}\overline{v}.\label{q54}%
\end{equation}

\item Each $f\in\overset{\cdot}{W}_{G,X}$ has a unique representation
\ref{q51} and $\dim\overset{\cdot}{W}_{G,X}=N-M$.

\item $G_{X,X}$ is conditionally positive definite on $\operatorname*{null}%
P_{X}^{T}$ i.e. $P_{X}^{T}v=0$ and $v\neq0$ implies $v^{T}G_{X,X}\overline
{v}>0$.

\item We have the direct sum $W_{G,X}=\overset{\cdot}{W}_{G,X}\oplus
P_{\theta-1}$ and $\dim W_{G,X}=N$.

\item $X_{w}^{\theta}=W_{G,X}\oplus W_{G,X}^{\perp}$ where
\[
W_{G,X}^{\perp}=\left\{  u\in X_{w}^{\theta}:u\left(  x^{\left(  k\right)
}\right)  =0\text{ }for\text{ }all\text{ }x^{\left(  k\right)  }\in X\right\}
.
\]

\end{enumerate}
\end{theorem}

\begin{proof}
\textbf{Part 1} The first step is to show that $f_{v}\in X_{w}^{\mathbb{\theta
}}$. Here we use the definition of $X_{w}^{\theta}$ given in equation
\ref{p948} i.e. $f_{v}\in S^{\prime}$, $\widehat{f_{v}}\in L_{loc}^{1}\left(
\mathbb{R}^{d}\setminus0\right)  $, $\int w\left\vert \cdot\right\vert
^{2\theta}\left\vert \left(  f_{v}\right)  _{F}\right\vert ^{2}<\infty$ and
$\left\vert \alpha\right\vert =\theta$ implies $\xi^{\alpha}\widehat{f_{v}%
}=\xi^{\alpha}\left(  f_{v}\right)  _{F}$ on $S$. Clearly since $G\in
S^{\prime}$ we have $f_{v}\in S^{\prime}$ and
\[
\widehat{f_{\nu}}=\sum_{j=1}^{N}v_{j}F_{x}\left[  G\left(  x-x^{\left(
j\right)  }\right)  \right]  =\sum_{j=1}^{N}v_{j}e^{-i\left(  x^{\left(
j\right)  },\xi\right)  }\widehat{G}=\left(  \sum_{j=1}^{N}v_{j}e^{-i\left(
x^{\left(  j\right)  },\xi\right)  }\right)  \widehat{G}=a_{v}\left(
\xi\right)  \widehat{G},
\]

where $a_{v}$ is given by \ref{q58}.

Thus $\widehat{f_{\nu}}=\frac{a_{v}}{w\left\vert \cdot\right\vert ^{2\theta}}$
on $\mathbb{R}^{d}\setminus0$ and since property W2.1 is $1/w\in L_{loc}^{1}$
it follows that $\widehat{f_{\nu}}\in L_{loc}^{1}\left(  \mathbb{R}%
^{d}\setminus0\right)  $ and
\begin{equation}
\left(  f_{\nu}\right)  _{F}=\frac{a_{v}}{w\left\vert \cdot\right\vert
^{2\theta}}.\label{q57}%
\end{equation}

By Theorem \ref{Lem_sum(v.exp(-ixy))}, $a_{v}\in C_{\emptyset,\theta}^{\infty
}\cap C_{B}^{\infty}$ and so by Theorem \ref{Thm_product_of_Co,k_funcs} we
have $\left\vert a_{v}\right\vert ^{2}\in C_{\emptyset,2\theta}^{\infty}\cap
C_{B}^{\infty}$. Further%
\begin{equation}
\left\vert f_{\nu}\right\vert _{w,\theta}^{2}=\int w\left\vert \cdot
\right\vert ^{2\theta}\left\vert \left(  f_{\nu}\right)  _{F}\right\vert
^{2}=\int w\left\vert \cdot\right\vert ^{2\theta}\left\vert \frac{a_{v}%
}{w\left\vert \cdot\right\vert ^{2\theta}}\right\vert ^{2}=\int\frac
{\left\vert a_{v}\right\vert ^{2}}{w\left\vert \cdot\right\vert ^{2\theta}%
}.\label{q93}%
\end{equation}

From Theorem \ref{Thm_weight_property_relat}, W3.1 implies property W3.2 and
W3.3 implies property W3.2, and by definition, property W3.2 implies
$\int_{\left\vert \cdot\right\vert \geq r_{3}}\frac{1}{w\left\vert
\cdot\right\vert ^{2\theta}}<\infty$ for some $r_{3}\geq0$. Choose $\phi\in
C_{0}^{\infty}$ such that $0\leq\phi\leq1$ and $\phi\left(  x\right)  =1$ in a
neighborhood of $0$, and define $\phi_{\varepsilon}$ by: $\phi_{\varepsilon
}\left(  \xi\right)  =\phi\left(  \varepsilon\xi\right)  $ for $\varepsilon
>0$. Then%
\begin{equation}
\int\frac{\left\vert a_{v}\right\vert ^{2}}{w\left\vert \cdot\right\vert
^{2\theta}}=\int\frac{\phi_{\varepsilon}\left\vert a_{v}\right\vert ^{2}%
}{w\left\vert \cdot\right\vert ^{2\theta}}+\int\frac{\left(  1-\phi
_{\varepsilon}\right)  \left\vert a_{v}\right\vert ^{2}}{w\left\vert
\cdot\right\vert ^{2\theta}},\label{q53}%
\end{equation}

and since $\phi_{\varepsilon}\left\vert a_{v}\right\vert ^{2}\in
S_{\emptyset,2\theta}$ we have by the basis function Definition
\ref{Def_basis_distrib} that $\int\frac{\phi_{\varepsilon}\left\vert
a_{v}\right\vert ^{2}}{w\left\vert \cdot\right\vert ^{2\theta}}=\left[
\widehat{G},\phi_{\varepsilon}\left\vert a_{v}\right\vert ^{2}\right]
<\infty$. Further, for sufficiently small $\varepsilon$ the support of
$\left(  1-\phi_{\varepsilon}\right)  \left\vert a_{v}\right\vert ^{2}$ lies
outside the sphere of radius $r_{3}$ so property W3.2 now ensures that the
integral \ref{q53} exists. To finish the proof that $f_{v}\in X_{w}^{\theta} $
we must show that $\xi^{\alpha}\left(  f_{\nu}\right)  _{F}=\xi^{\alpha
}\widehat{f_{\nu}}$ on $S$ i.e. $\xi^{\alpha}\frac{a_{v}}{w\left\vert
\cdot\right\vert ^{2\theta}}=\xi^{\alpha}a_{v}\widehat{G}$ on $S$.

If $\left\vert \alpha\right\vert =\theta$ then by Theorem
\ref{Thm_product_of_Co,k_funcs}, $\xi^{\alpha}\in C_{\emptyset,\theta}%
^{\infty}$, $\xi^{\alpha}a_{v}\in C_{\emptyset,2\theta}^{\infty}\cap
C_{BP}^{\infty}$ and $\xi^{\alpha}a_{v}\psi\in S_{\emptyset,2\theta}$ when
$\psi\in S$ so that
\[
\left[  \xi^{\alpha}\frac{a_{v}}{w\left\vert \cdot\right\vert ^{2\theta}}%
,\psi\right]  =\left[  \frac{1}{w\left\vert \cdot\right\vert ^{2\theta}}%
,\xi^{\alpha}a_{v}\psi\right]  =\left[  \widehat{G},\xi^{\alpha}a_{v}%
\psi\right]  =\left[  \xi^{\alpha}a_{v}\widehat{G},\psi\right]  ,
\]

as required. Thus $f_{v}\in X_{w}^{\theta}$ and \ref{q93} proves the first
equation of \ref{q54}.

We have already shown that $\left\vert a_{v}\right\vert ^{2}\in C_{\emptyset
,2\theta}^{\infty}\cap C_{B}^{\infty}$ so $\left\vert a_{v}\right\vert
^{2}\phi_{\varepsilon}\in S_{\emptyset,2\theta}$ and by definition of a basis
distribution%
\begin{align*}
\int\frac{\left\vert a_{v}\right\vert ^{2}}{w\left\vert \cdot\right\vert
^{2\theta}}\phi_{\varepsilon}=\left[  \frac{1}{w\left\vert \cdot\right\vert
^{2\theta}},\left\vert a_{v}\right\vert ^{2}\phi_{\varepsilon}\right]
=\left[  \widehat{G},\left\vert a_{v}\right\vert ^{2}\phi_{\varepsilon
}\right]   &  =\left[  G,\widehat{\left\vert a_{v}\right\vert ^{2}%
\phi_{\varepsilon}}\right] \\
&  =\sum_{j,k=1}^{N}v_{j}\overline{v_{k}}\left[  G,\widehat{\phi_{\varepsilon
}}\left(  \cdot-\left(  x^{\left(  k\right)  }-x^{\left(  j\right)  }\right)
\right)  \right] \\
&  =\sum_{j,k=1}^{N}v_{j}\overline{v_{k}}\left[  G,\overset{\vee
}{\phi_{\varepsilon}}\left(  \left(  x^{\left(  k\right)  }-x^{\left(
j\right)  }\right)  -\cdot\right)  \right]  ,
\end{align*}

by part 2 of Theorem \ref{Lem_sum(v.exp(-ixy))}. By the convolution definition
of part 2 Appendix \ref{Def_Convolution},
\[
\sum_{j,k=1}^{N}v_{j}\overline{v_{k}}\left[  G,\overset{\vee}{\phi
_{\varepsilon}}\left(  \left(  x^{\left(  k\right)  }-x^{\left(  j\right)
}\right)  -\cdot\right)  \right]  =\left(  2\pi\right)  ^{\frac{d}{2}}%
\sum_{j,k=1}^{N}v_{j}\overline{v_{k}}\left(  G\ast\overset{\vee}{\phi
_{\varepsilon}}\right)  \left(  x^{\left(  k\right)  }-x^{\left(  j\right)
}\right)  .
\]

Now in the tempered distribution sense as $\varepsilon\rightarrow0$,
$\phi_{\varepsilon}\rightarrow1$ and so $\overset{\vee}{\phi_{\varepsilon}%
}\rightarrow\overset{\vee}{1}=\left(  2\pi\right)  ^{\frac{d}{2}}\delta$ and
hence $G\ast\overset{\vee}{\phi_{\varepsilon}}\rightarrow G\ast\left(
2\pi\right)  ^{\frac{d}{2}}\delta=G$ pointwise because by Theorem
\ref{Thm_basis_smth_W3.1} and \ref{Thm_basis_smth_W3.2_r3_pos}, if $w$ has
properties W2.1 and W3 then $G\in C^{\left(  0\right)  }$. Here $\overset{\vee
}{\phi_{\varepsilon}}$ is called a mollifier e.g. Lemma 2.18 Adams
\cite{Adams75}. Thus%
\[
\lim\limits_{\varepsilon\rightarrow0}\int\frac{\left\vert a_{v}\right\vert
^{2}}{w\left\vert \cdot\right\vert ^{2\theta}}\phi_{\varepsilon}=\left(
2\pi\right)  ^{\frac{d}{2}}\sum_{j,k=1}^{N}v_{j}\overline{v_{k}}G\left(
x^{\left(  k\right)  }-x^{\left(  j\right)  }\right)  ,
\]

and finally we note that $\lim\limits_{\varepsilon\rightarrow0}\int%
\frac{\left\vert a_{v}\right\vert ^{2}}{w\left\vert \cdot\right\vert
^{2\theta}}\phi_{\varepsilon}=\int\frac{\left\vert a_{v}\right\vert ^{2}%
}{w\left\vert \cdot\right\vert ^{2\theta}}$.\medskip

\textbf{Part 2} Suppose $\sum\limits_{k=1}^{N}v_{k}G\left(  x-x^{\left(
k\right)  }\right)  =0$ for all $x$ and $P_{X}^{T}v=0$. Then
\[
0=\left\vert \sum_{k=1}^{N}v_{k}G\left(  x-x^{\left(  k\right)  }\right)
\right\vert _{w,\theta}^{2}=\int\frac{\left\vert a_{v}\right\vert ^{2}%
}{w\left\vert \cdot\right\vert ^{2\theta}}.
\]

Since $w>0$ a.e. we must conclude that $a_{v}=0$ a.e. and that by part 3 of
Lemma \ref{Lem_sum(v.exp(-ixy))}, $v_{k}=0$ for all $k$. Thus the operator
$\Phi:\operatorname*{null}P_{X}^{T}\rightarrow\overset{\cdot}{W}_{G,X}$
defined by $\Phi v=\sum_{k=1}^{N}v_{k}G\left(  \cdot-x^{\left(  k\right)
}\right)  $ is an isomorphism and so $\dim\overset{\cdot}{W}_{G,X}%
=\dim\operatorname*{null}P_{X}^{T}=N-M$ by part 5 Theorem
\ref{Thm_Px_properties}.\medskip

\textbf{Part 3} From part 1, $P_{X}^{T}v=0$ implies $v^{T}G_{X,X}\overline
{v}\geq0$. Suppose $v^{T}G_{X,X}\overline{v}=0$ and $P_{X}^{T}v=0$. By part 1
of this theorem
\[
0=v^{T}G_{X,X}\overline{v}=\left\vert \sum_{k=1}^{N}v_{k}G\left(  x-x^{\left(
k\right)  }\right)  \right\vert _{w,\theta},
\]

so that $\sum\limits_{k=1}^{N}v_{k}G\left(  x-x^{\left(  k\right)  }\right)
=0$. Part 2 of this theorem then implies $v=0$.\smallskip

\textbf{Part 4} We must show that $\overset{\cdot}{W}_{G,X}\cap P_{\theta
-1}=\left\{  0\right\}  $. Suppose $p\in\overset{\cdot}{W}_{G,X}\cap
P_{\theta-1}$. Then $p=\sum\limits_{k=1}^{N}v_{k}G\left(  x-x^{\left(
k\right)  }\right)  $ and $P_{X}^{T}v=0$ but by part 1, $0=\left\Vert
p\right\Vert _{w,\theta}=\left(  2\pi\right)  ^{\frac{d}{2}}v^{T}%
G_{X,X}\overline{v}$. Part 2 now implies $v=0$ and hence that $p=0$.\medskip

\textbf{Part 5}\ Since $\left\{  R_{x^{\left(  k\right)  }}:x^{\left(
k\right)  }\in X\right\}  $ is a basis for $W_{G,X}$%
\begin{align*}
W_{G,X}^{\perp}=\left\{  v\in X_{w}^{\theta}\mid\left(  v,u\right)
_{w,\theta}=0\text{ }if\text{ }u\in W_{G,X}\right\}   &  =\left\{  v\in
X_{w}^{\theta}\mid\left(  v,R_{x^{\left(  k\right)  }}\right)  _{w,\theta
}=0\text{ }if\text{ }x^{\left(  k\right)  }\in X\right\} \\
&  =\left\{  v\in X_{w}^{\theta}\mid v\left(  x^{\left(  k\right)  }\right)
=0\text{ }if\text{\ }x^{\left(  k\right)  }\in X\right\}  .
\end{align*}

\end{proof}

A consequence of the last theorem is the following representation result for
members of the set $W_{G,X}$:

\begin{corollary}
\label{Cor_uniq_rep_for_elt_of_Wgx}Suppose the definitions and assumptions of
Theorem \ref{Thm_Wgx_properties_1} hold. If $\left\{  p_{j}\right\}
_{j=1}^{M}$ is basis for $P_{\theta-1}$ then the representation
\[
W_{G,X}=\left\{  \sum\limits_{i=1}^{N}\alpha_{i}G\left(  \cdot-x^{\left(
i\right)  }\right)  +\sum\limits_{j=1}^{M}\beta_{j}p_{j}:P_{X}^{T}%
\alpha=0,\text{ }\alpha=\left(  \alpha_{i}\right)  ,\text{ }\alpha_{i}%
,\beta_{j}\in\mathbb{C}\right\}  ,
\]

is unique in terms of $\alpha_{i}$ and $\beta_{j}$.
\end{corollary}

\begin{proof}
This follows directly from parts 2 and 4 of Theorem \ref{Thm_Wgx_properties_1}.
\end{proof}

In Subsection \ref{SbSect_mat_eqn_interpol} a matrix equation will be derived
for the coefficients $\alpha_{i}$ and $\beta_{j}$. From Theorem
\ref{Thm_G(*-Xk)_lin_indep} we know that the data-translated basis functions
$\left\{  G\left(  \cdot-x^{\left(  i\right)  }\right)  \right\}  _{i=1}^{N}$
are linearly independent and from the last theorem we know that the dimension
of $\overset{\cdot}{W}_{G,X}$ is $N-M$. However no subset of the functions
$\left\{  G\left(  \cdot-x^{\left(  i\right)  }\right)  \right\}  _{i=1}^{N}$
containing $N-M$ of these functions spans $\overset{\cdot}{W}_{G,X}$. The
functions $\left\{  G\left(  \cdot-x^{\left(  i\right)  }\right)  \right\}
_{i=1}^{N}$ can be called basis functions w.r.t. $\overset{\cdot}{W}_{G,X}$ in
the sense that linear combinations of these functions can be used to construct
a basis for $\overset{\cdot}{W}_{G,X}$. In fact, in the next corollary we will
show that provided the minimal unisolvent subset $A$ used to construct the
Riesz representer $R_{x}$ lies in $X$ and constitutes the first $M$ points of
$X$, it will follow that $\left\{  R_{x^{\left(  i\right)  }}\right\}
_{i=M+1}^{N}$ is a basis for $\overset{\cdot}{W}_{G,X}$.

\begin{theorem}
\label{Thm_Wgx_properties_2}Suppose $w$ is a weight function with properties
W2 and W3 for order $\theta$. Next let $X=\left\{  x^{\left(  k\right)
}\right\}  _{k=1}^{N}$ be $\theta$-unisolvent and suppose $X_{1}=\left\{
x^{\left(  k\right)  }\right\}  _{k=1}^{M}$ is a minimal unisolvent set.
Define the Riesz representer function $R_{x}$ by \ref{p932} using $X_{1}$.

Then the spaces $W_{G,X}$ and $\overset{\cdot}{W}_{G,X}$have the following properties:

\begin{enumerate}
\item $W_{G,X}$ has basis $\left\{  R_{x^{\left(  i\right)  }}\right\}
_{i=1}^{N}$ over the complex numbers.

\item $\overset{\cdot}{W}_{G,X}$has basis $\left\{  R_{x^{\left(  i\right)  }%
}\right\}  _{i=M+1}^{N}$ over the complex numbers.
\end{enumerate}
\end{theorem}

\begin{proof}
It will first be shown that if $\beta=\left\{  \beta_{k}\right\}  _{k=M+1}%
^{N}$ are any complex numbers then $\sum\limits_{k=M+1}^{N}\beta
_{k}R_{x^{\left(  k\right)  }}\in W_{G,X}$. Let $g_{\beta}=\sum\limits_{k=M+1}%
^{N}\beta_{k}R_{x^{\left(  k\right)  }}$. Then each $R_{x^{\left(  k\right)
}}\in X_{w}^{\theta}$ and from equation \ref{q55}
\[
\left(  R_{x^{\left(  k\right)  }}\right)  _{F}=\left(  2\pi\right)
^{-d/2}\frac{\left(  \mathcal{Q}_{x}\left(  e^{-i\left(  x,\cdot\right)
}\right)  \right)  \left(  x=x^{\left(  k\right)  }\right)  }{w\left\vert
\cdot\right\vert ^{2\theta}}\text{ }on\text{ }\mathbb{R}^{d}\setminus0.
\]
Thus%
\begin{align*}
\left(  g_{\beta}\right)  _{F}  & =\sum\limits_{k=M+1}^{N}\beta_{k}\left(
R_{x^{\left(  k\right)  }}\right)  _{F}\\
& =\frac{\left(  2\pi\right)  ^{-d/2}}{w\left\vert \cdot\right\vert ^{2\theta
}}\sum\limits_{k=M+1}^{N}\beta_{k}\left(  \mathcal{Q}_{x}\left(  e^{-i\left(
x,\cdot\right)  }\right)  \right)  \left(  x=x^{\left(  k\right)  }\right) \\
& =\frac{\left(  2\pi\right)  ^{-d/2}}{w\left\vert \cdot\right\vert ^{2\theta
}}\sum\limits_{k=M+1}^{N}\beta_{k}\left(  e^{-i\left(  x,\cdot\right)  }%
-\sum\limits_{j=1}^{M}l_{j}\left(  x\right)  e^{-i\left(  x^{\left(  j\right)
},\cdot\right)  }\right)  \left(  x=x^{\left(  k\right)  }\right) \\
& =\frac{\left(  2\pi\right)  ^{-d/2}}{w\left\vert \cdot\right\vert ^{2\theta
}}\sum\limits_{k=M+1}^{N}\beta_{k}\left(  e^{-i\left(  x^{\left(  k\right)
},\cdot\right)  }-\sum\limits_{j=1}^{M}l_{j}\left(  x^{\left(  k\right)
}\right)  e^{-i\left(  x^{\left(  j\right)  },\cdot\right)  }\right) \\
& =\frac{\left(  2\pi\right)  ^{-d/2}}{w\left\vert \cdot\right\vert ^{2\theta
}}\left(  \sum\limits_{k=M+1}^{N}\beta_{k}e^{-i\left(  x^{\left(  k\right)
},\cdot\right)  }-\sum\limits_{j=1}^{M}\left(  \sum\limits_{k=M+1}^{N}%
\beta_{k}l_{j}\left(  x^{\left(  k\right)  }\right)  \right)  e^{-i\left(
x^{\left(  j\right)  },\cdot\right)  }\right)  .
\end{align*}

From \ref{q57} and \ref{q58} we see that if $f_{\gamma}=\sum\limits_{k=1}%
^{N}\gamma_{k}G\left(  \cdot-x^{\left(  k\right)  }\right)  $ then

$\left(  f_{\gamma}\right)  _{F}=\frac{\left(  2\pi\right)  ^{-d/2}%
}{w\left\vert \cdot\right\vert ^{2\theta}}\sum\limits_{k=1}^{N}\gamma
_{k}e^{-i\left(  x^{\left(  k\right)  },\cdot\right)  }$ and setting $\left(
g_{\beta}\right)  _{F}=\frac{\left(  2\pi\right)  ^{-d/2}}{w\left\vert
\cdot\right\vert ^{2\theta}}\sum\limits_{k=1}^{N}\gamma_{k}e^{-i\left(
x^{\left(  k\right)  },\cdot\right)  }$ yields%
\begin{equation}
\gamma_{j}=\left\{
\begin{array}
[c]{cc}%
-\left(  2\pi\right)  ^{d/2}\sum\limits_{k=M+1}^{N}l_{j}\left(  x^{\left(
k\right)  }\right)  \beta_{j}, & j\leq M,\\
\left(  2\pi\right)  ^{d/2}\beta_{j}, & j>M.
\end{array}
\right. \label{q73}%
\end{equation}

The criterion of part 1 of Theorem \ref{Thm_Px_properties_2} clearly implies
that $P_{X}^{T}\gamma=0$ and thus

$f_{\gamma}=\left(  2\pi\right)  ^{-d/2}\sum\limits_{j=1}^{N}\gamma
_{j}G\left(  x-x^{\left(  j\right)  }\right)  \in\overset{\cdot}{W}_{G,X}$ and
$\left(  g_{\beta}-f_{\gamma}\right)  _{F}=0$. But by part 4 Summary
\ref{Sum_int_properties_Xwm}, $g_{\beta}-f_{\gamma}\in P_{\theta-1}$ and thus
$g_{\beta}\in W_{G,X}$. Further, $g_{\beta}=\sum\limits_{k=M+1}^{N}\beta
_{k}R_{x^{\left(  k\right)  }}$ so $\sum\limits_{k=1}^{N}\beta_{k}R_{x^{(k)}%
}\in W_{G,X}$ and $\operatorname*{span}\left\{  R_{x^{(j)}}\right\}
_{j=1}^{N}\subset W_{G,X}$.

To prove the converse choose $u\in W_{G,X}$. Then $u=\sum\limits_{j=1}%
^{N}\beta_{j}G\left(  \cdot-x^{\left(  j\right)  }\right)  +q$ for some $q\in
P_{\theta-1}$ and $\beta\in\mathbb{C}^{N}$ such that $P_{X}^{T}\beta=0$. When
$\gamma=0$ and $x=x^{\left(  j\right)  }$ equation \ref{p911} for $R_{x}$ can
be rearranged to give%
\[
G\left(  \cdot-x^{\left(  j\right)  }\right)  =\left(  2\pi\right)
^{d/2}R_{x^{\left(  j\right)  }}+\sum\limits_{i=1}^{M}l_{i}\left(  x^{\left(
j\right)  }\right)  G\left(  \cdot-x^{\left(  i\right)  }\right)
-q_{x^{\left(  j\right)  }},
\]

where each $q_{x^{\left(  j\right)  }}\in P_{\theta-1}$. Hence, since
\ref{p936} implies $q_{x^{\left(  j\right)  }}=\left(  2\pi\right)  ^{\frac
{d}{2}}R_{x^{\left(  j\right)  }}$ for $j\leq M$,
\begin{align*}
\sum_{j=1}^{N}\beta_{j}G\left(  \cdot-x^{\left(  j\right)  }\right)   &
=\left(  2\pi\right)  ^{\frac{d}{2}}\sum_{j=1}^{N}\beta_{j}R_{x^{\left(
j\right)  }}+\sum_{j=1}^{N}\beta_{j}\sum\limits_{i=1}^{M}l_{i}\left(
x^{\left(  j\right)  }\right)  G\left(  \cdot-x^{\left(  i\right)  }\right)
-\sum_{j=1}^{N}\beta_{j}q_{x^{\left(  j\right)  }}\\
& =\left(  2\pi\right)  ^{\frac{d}{2}}\sum_{j=M+1}^{N}\beta_{j}R_{x^{\left(
j\right)  }}+\sum\limits_{i=1}^{M}\left(  \sum_{j=1}^{N}\beta_{j}l_{i}\left(
x^{\left(  j\right)  }\right)  \right)  G\left(  \cdot-x^{\left(  i\right)
}\right)  -\sum_{j=M+1}^{N}\beta_{j}q_{x^{\left(  j\right)  }},
\end{align*}

and so
\begin{align*}
u &  =\sum_{j=1}^{N}\beta_{j}G\left(  \cdot-x^{\left(  j\right)  }\right)
+q\\
&  =\left(  2\pi\right)  ^{\frac{d}{2}}\sum_{j=M+1}^{N}\beta_{j}R_{x^{\left(
j\right)  }}+\sum_{i=1}^{M}\left(  \sum_{j=1}^{N}\beta_{j}l_{i}(x^{\left(
j\right)  })\right)  G\left(  \cdot-x^{\left(  i\right)  }\right)
-\sum_{j=M+1}^{N}\beta_{j}q_{x^{\left(  j\right)  }}+q\\
&  =\left(  2\pi\right)  ^{\frac{d}{2}}\sum_{j=M+1}^{N}\beta_{j}R_{x^{\left(
j\right)  }}+\sum_{i=1}^{M}\left(  \beta_{i}+\sum_{j=M+1}^{N}\beta_{j}%
l_{i}(x^{\left(  j\right)  })\right)  G\left(  \cdot-x^{\left(  i\right)
}\right)  -\sum_{j=M+1}^{N}\beta_{j}q_{x^{\left(  j\right)  }}+q\\
&  =\left(  2\pi\right)  ^{\frac{d}{2}}\sum_{j=M+1}^{N}\beta_{j}R_{x^{\left(
j\right)  }}-\sum_{j=M+1}^{N}\beta_{j}q_{x^{\left(  j\right)  }}+q,
\end{align*}

since by part 1 Theorem \ref{Thm_Px_properties_2}, $P_{X}^{T}\beta=0$ iff
$\beta_{i}=-\sum\limits_{j=M+1}^{N}\beta_{j}l_{i}\left(  x^{\left(  j\right)
}\right)  $ for each $l_{i}$. Thus $W_{G,X}=\operatorname*{span}%
\limits_{\mathbb{C}}\left\{  R_{x^{\left(  j\right)  }}\right\}  _{j=1}^{N}$.
Further, by part 2 of Theorem \ref{Thm_Riesz_property}, we know that
$R_{x^{\left(  j\right)  }}=l_{j}$ for $1\leq j\leq M$ and so
$\operatorname*{span}\limits_{\mathbb{C}}\left\{  R_{x^{\left(  j\right)  }%
}\right\}  _{j=1}^{M}=P_{\theta-1}$. Therefore the result $W_{G,X}%
=\overset{\cdot}{W}_{G,X}\oplus P_{\theta-1}$ from Theorem
\ref{Thm_Wgx_properties_1} implies that $\overset{\cdot}{W}_{G,X}%
=\operatorname*{span}\limits_{\mathbb{C}}\left\{  R_{x^{\left(  j\right)  }%
}\right\}  _{j=M+1}^{N}$.
\end{proof}

The previous theorem assumed that the Riesz representers were constructed from
a minimal unisolvent subset $X_{1}$ of $X$ which consisted of the first $M$
points of $X$. The next corollary weakens this assumption by that we can
choose any minimal unisolvent subset of $X$. Now recall Definition
\ref{Def_cardinal_basis} regarding permutation theory.

\begin{corollary}
\label{Cor_Thm_Wgx_properties_2}Suppose $w$ is a weight function with
properties W2 and W3 for order $\theta$. Next let $X=\left\{  x^{\left(
k\right)  }\right\}  _{k=1}^{N}$ be $\theta$-unisolvent and suppose $A\subset
X$ is any minimal unisolvent subset. Denote its cardinal basis by $\left\{
l_{k}\right\}  _{k=1}^{M}$ and define the Riesz representer function $R_{x}$
by \ref{p939}. Then the spaces $W_{G,X}$ and $\overset{\cdot}{W}_{G,X}$ have
the following properties:

\begin{enumerate}
\item $W_{G,X}$ has basis $\left\{  R_{x^{\left(  i\right)  }}\right\}
_{i=1}^{N}$ over the complex numbers.

\item $\overset{\cdot}{W}_{G,X}$ has basis $\left\{  R_{x^{\left(  i\right)
}}:x^{\left(  i\right)  }\notin A\right\}  $ over the complex numbers.
\end{enumerate}
\end{corollary}

\begin{proof}
\textbf{Part 1} We want to use Theorem \ref{Thm_Wgx_properties_2} so start by
re-ordering $X$ using a permutation $\pi$ such that the first $M$ points in
$\pi X$ belong to $A$. By Theorem \ref{Thm_Wgx_properties_2}, $W_{G,\pi X} $
has basis $\left\{  R_{y^{\left(  i\right)  }}\left(  \pi A\right)  \right\}
_{i=1}^{N}$ where $y^{\left(  i\right)  }=x^{\left(  \pi\left(  i\right)
\right)  }$. By \ref{p939} the function $R_{x}$ only depends on $G$, $A$ and
$x $ and by Theorem \ref{Thm_Rx_unique}, $R_{x}$ is independent of the order
of $A$. Hence $W_{G,\pi X}$ has bases $\left\{  R_{y^{\left(  i\right)  }%
}\left(  \pi A\right)  \right\}  _{i=1}^{N}=\left\{  R_{y^{\left(  i\right)
}}\left(  A\right)  \right\}  _{i=1}^{N}=\left\{  R_{x^{\left(  i\right)  }%
}\right\}  _{i=1}^{N}$. By Theorem \ref{Thm_Wgx_reorder} the set $W_{G,X}$ is
independent of the order of the points in $X$. Hence $W_{G,\pi X}=W_{G,X}$ and
we have part 1.\medskip

\textbf{Part 2} We know from part 4 Theorem \ref{Thm_Wgx_properties_1} that
$W_{G,X}=\overset{\cdot}{W}_{G,X}\oplus P_{\theta-1}$, and from part 2 Theorem
\ref{Thm_Riesz_property} we know that $\left\{  R_{x^{\left(  i\right)  }%
}:x^{\left(  i\right)  }\in A\right\}  $ is the cardinal basis for
$P_{\theta-1}$ generated by $A$. Thus $\overset{\cdot}{W}_{G,X}$has basis
$\left\{  R_{x^{\left(  i\right)  }}:x^{\left(  i\right)  }\notin A\right\}  $.
\end{proof}

\section{The vector-valued evaluation operator $\widetilde{\mathcal{E}}_{X}%
$\label{Sect_Eval_&_T_ops}}

The vector-valued evaluation operator $\widetilde{\mathcal{E}}_{X}$ and its
adjoint will be fundamental to the study of the interpolant in this document
and to the smoothers in the later chapters. This function evaluates a function
on an ordered set of points to form a complex vector.

\begin{definition}
\label{Def_eval_operators}\textbf{The vector-valued evaluation operator
}$\widetilde{\mathcal{E}}_{X}$

Let $X=\left\{  x^{(i)}\right\}  _{i=1}^{N}$ be a set of $N$ distinct points
in $\mathbb{R}^{d}$. Let $u$ be a continuous function. Then the evaluation
operator $\widetilde{\mathcal{E}}_{X}$ is defined by
\[
\widetilde{\mathcal{E}}_{X}u=\left(  u\left(  x^{(i)}\right)  \right)
_{i=1}^{N}.
\]

Sometimes we will use the more compact notation $u_{X}$ for $\widetilde
{\mathcal{E}}_{X}u$ and when dealing with matrices $\widetilde{\mathcal{E}%
}_{X}u$ will be regarded as a column vector.
\end{definition}

We will now assume that the weight function $w$ satisfies properties W2 and W3
for order $\theta$ and smoothness $\kappa$. Then the functions in
$X_{w}^{\theta}$ are continuous. The Riesz representer $R_{x}$ of the
evaluation functional $f\rightarrow f(z)$ allows some important properties of
the evaluation operator $\widetilde{\mathcal{E}}_{X}$ to be proved. In the
next theorem we assume that the first $M$ points in $X$ are minimally
unisolvent and use these to generate $R_{x}$.

\begin{theorem}
\label{Thm_eval_op_properties}Suppose that $X=\left\{  x^{(i)}\right\}
_{i=1}^{N}$ is a $\theta$-unisolvent set of distinct points in $\mathbb{R}%
^{d}$. Assume that $A=\left\{  a^{\left(  i\right)  }\right\}  _{i=1}^{M}$ is
any minimal unisolvent subset with cardinal basis $\left\{  l_{k}\right\}
_{k=1}^{M}$ which we use to define the Riesz representer $R_{x}$, the Light
norm $\left\Vert \cdot\right\Vert _{w,\theta}$ and the Lagrangian
interpolation operator $\mathcal{P}$. Then:

\begin{enumerate}
\item The evaluation operator $\widetilde{\mathcal{E}}_{X}:\left(
X_{w}^{\mathbb{\theta}},\left\Vert \cdot\right\Vert _{w,\theta}\right)
\rightarrow\left(  \mathbb{C}^{N},\left\vert \cdot\right\vert \right)  $ is
continuous and onto. Also, $\operatorname*{null}\widetilde{\mathcal{E}}%
_{X}=W_{G,X}^{\bot}$.

\item The adjoint operator $\widetilde{\mathcal{E}}_{X}^{\ast}:\mathbb{C}%
^{N}\rightarrow X_{w}^{\theta}$, defined by $\left(  \widetilde{\mathcal{E}%
}_{X}f,\beta\right)  _{\mathbb{C}^{N}}=\left(  f,\widetilde{\mathcal{E}}%
_{X}^{\ast}\beta\right)  _{w,\theta}$, satisfies
\begin{equation}
\widetilde{\mathcal{E}}_{X}^{\ast}\beta=\sum\limits_{i=1}^{N}\beta
_{i}R_{x^{(i)}},\text{\quad}\beta\in\mathbb{C}^{N},\label{q66}%
\end{equation}

and $\widetilde{\mathcal{E}}_{X}^{\ast}$ is continuous, 1-1 and maps $\left(
\mathbb{C}^{N},\left\vert \cdot\right\vert \right)  $ onto $\left(
W_{G,X},\left\Vert \cdot\right\Vert _{w,\theta}\right)  $.

\item $\left\Vert \widetilde{\mathcal{E}}_{X}^{\ast}\right\Vert =\left\Vert
\widetilde{\mathcal{E}}_{X}\right\Vert =\left\Vert R_{X,X}\right\Vert $, where
$R_{X,X}=\left(  R_{x^{\left(  j\right)  }}(x^{\left(  i\right)  })\right)  $
is the\textbf{\ }(regular)\textbf{\ }reproducing kernel matrix discussed in
Subsection \ref{SbSect_mat_eval_rep} and $\left\Vert \cdot\right\Vert $ is the
matrix norm corresponding to the Euclidean vector norm.

\item Assuming $\beta$ is a column vector and $\widetilde{\mathcal{E}}_{X}$ is
a column vector we have%
\[
\widetilde{\mathcal{E}}_{X}\widetilde{\mathcal{E}}_{X}^{\ast}\beta
=R_{X,X}\beta,\quad\beta\in\mathbb{C}^{N}.
\]

\item Assuming that $\widetilde{\mathcal{E}}_{X}$ and $\gamma$ are column
vectors,%
\[
\widetilde{\mathcal{E}}_{X}\widetilde{\mathcal{E}}_{A}^{\ast}\gamma
=L_{X}\gamma,\quad\gamma\in\mathbb{C}^{M},
\]

where $L_{X}=\left(  l_{j}\left(  x^{\left(  i\right)  }\right)  \right)  $ is
the $N\times M$ cardinal unisolvent matrix (Definition
\ref{Def_unisolv_matrix_Px}) corresponding to the minimal unisolvent set $A$.

\item $\widetilde{\mathcal{E}}_{X}^{\ast}\widetilde{\mathcal{E}}_{X}%
:X_{w}^{\theta}\rightarrow W_{G,X}$ is onto and $\operatorname*{null}%
\widetilde{\mathcal{E}}_{X}^{\ast}\widetilde{\mathcal{E}}_{X}=W_{G,X}^{\perp}%
$. Also, for $u\in X_{w}^{\mathbb{\theta}}$%
\[
\left(  \widetilde{\mathcal{E}}_{X}^{\ast}\widetilde{\mathcal{E}}_{X}u\right)
\left(  x\right)  =\sum\limits_{i=1}^{N}u\left(  x^{\left(  i\right)
}\right)  R_{x^{\left(  i\right)  }}\left(  x\right)  =\left(  u,\sum
_{j=1}^{N}R_{x^{\left(  j\right)  }}\left(  x\right)  R_{x^{\left(  j\right)
}}\right)  _{w,\theta}.
\]

\item For $\alpha\in\mathbb{C}^{N}$, $\left(  \mathcal{P}\widetilde
{\mathcal{E}}_{X}^{\ast}\alpha\right)  \left(  x\right)  =\left(  \alpha
^{T}L_{X}\right)  \widetilde{l}\left(  x\right)  $ where $\widetilde{l}\left(
x\right)  =\left(  l_{i}\left(  x\right)  \right)  $.

\item $\widetilde{\mathcal{E}}_{X}\mathcal{P}f=L_{X}\widetilde{\mathcal{E}%
}_{A}f$.
\end{enumerate}
\end{theorem}

\begin{proof}
\textbf{Parts 1 and 2} We will first show that $\widetilde{\mathcal{E}}_{X}%
$\textbf{\ }is continuous\textbf{. }In fact
\begin{align*}
\left\vert \widetilde{\mathcal{E}}_{X}u\right\vert _{\mathbb{C}^{N}}^{2}%
=\sum_{i=1}^{N}\left\vert u(x^{(i)})\right\vert ^{2}=\sum_{i=1}^{N}\left\vert
\left(  u,R_{x^{(i)}}\right)  _{w,\theta}\right\vert ^{2} &  \leq\sum
_{i=1}^{N}\left\Vert u\right\Vert _{w,\theta}^{2}\left\Vert R_{x^{(i)}%
}\right\Vert _{w,\theta}^{2}\\
&  =\left(  \sum_{i=1}^{N}\left\Vert R_{x^{(i)}}\right\Vert _{w,\theta}%
^{2}\right)  \left\Vert u\right\Vert _{w,\theta}^{2},
\end{align*}

so that $\widetilde{\mathcal{E}}_{X}:X_{w}^{\mathbb{\theta}}\rightarrow
\mathbb{C}^{N}$ is continuous. Next we show that $\widetilde{\mathcal{E}}_{X}
$ is onto. The Hilbert space adjoint $\widetilde{\mathcal{E}}_{X}^{\ast}$ of
$\widetilde{\mathcal{E}}_{X}$ is defined by
\begin{equation}
\left(  \widetilde{\mathcal{E}}_{X}u,\beta\right)  _{\mathbb{C}^{N}}=\left(
u,\widetilde{\mathcal{E}}_{X}^{\ast}\beta\right)  _{w,\theta}.\label{p935}%
\end{equation}

We now calculate the adjoint by using the representer $R_{x}$ to good effect.
\[
\left(  \widetilde{\mathcal{E}}_{X}u,\beta\right)  _{\mathbb{C}^{N}}%
=\sum_{i=1}^{N}u\left(  x^{\left(  i\right)  }\right)  \overline{\beta_{i}%
}=\sum_{i=1}^{N}\left(  u,R_{x^{\left(  i\right)  }}\right)  _{w,\theta
}\overline{\beta_{i}}=\left(  u,\sum_{i=1}^{N}\beta_{i}R_{x^{\left(  i\right)
}}\right)  _{w,\theta},
\]

so that
\begin{equation}
\widetilde{\mathcal{E}}_{X}^{\ast}\beta=\sum_{i=1}^{N}\beta_{i}R_{x^{\left(
i\right)  }},\quad\beta\in\mathbb{C}^{N}.\label{q75}%
\end{equation}

In Corollary \ref{Cor_Thm_Wgx_properties_2} it was shown that the functions
$\left\{  R_{x^{\left(  i\right)  }}\right\}  _{i=1}^{N}$ form a basis for
$W_{G,X}$. Hence $\operatorname*{range}\widetilde{\mathcal{E}}_{X}^{\ast
}=W_{G,X}$ and $\operatorname*{null}\widetilde{\mathcal{E}}_{X}^{\ast
}=\left\{  0\right\}  $. We now recall the closed-range theorem e.g. Yosida
\cite{Yosida68}, which states that for a continuous linear operator
$\mathcal{S}$ the $\operatorname*{range}\mathcal{S}$ is closed iff
$\operatorname*{range}\mathcal{S}^{\ast}$ is closed. Since
$\operatorname*{range}\widetilde{\mathcal{E}}_{X}^{\ast}=W_{G,X}$ is finite
dimensional it is closed so we conclude that $\operatorname*{range}%
\widetilde{\mathcal{E}}_{X}$ is closed. Consequently, using the result that
$\overline{\operatorname*{range}\mathcal{S}}=\left(  \operatorname*{null}%
\mathcal{S}^{\ast}\right)  ^{\perp}$ for any continuous linear operator
$\mathcal{S}$, we have
\[
\operatorname*{range}\widetilde{\mathcal{E}}_{X}=\overline
{\operatorname*{range}\left(  \widetilde{\mathcal{E}}_{X}\right)  }=\left(
\operatorname*{null}\widetilde{\mathcal{E}}_{X}^{\ast}\right)  ^{\perp
}=\left\{  0\right\}  ^{\perp}=\mathbb{R}^{N}.
\]

Finally, from \ref{p935}, $\widetilde{\mathcal{E}}_{X}u=0$ iff $\left(
\widetilde{\mathcal{E}}_{X}u,\beta\right)  _{\mathbb{C}^{N}}=0$ for all
$\beta\in\mathbb{C}^{N}$ iff $\left(  u,\widetilde{\mathcal{E}}_{X}^{\ast
}\beta\right)  _{w,\theta}=0$ for all $\beta\in\mathbb{C}^{N}$. But
$\operatorname*{range}\widetilde{\mathcal{E}}_{X}^{\ast}=W_{G,X}$ so
$\widetilde{\mathcal{E}}_{X}u=0$ iff $\left(  u,f\right)  _{w,\theta}=0$ for
all $f\in W_{G,X}$ i.e. iff $u\in W_{G,X}^{\perp}$.\medskip

\textbf{Part 3} That $\left\Vert \widetilde{\mathcal{E}}_{X}^{\ast}\right\Vert
=\left\Vert \widetilde{\mathcal{E}}_{X}\right\Vert $ is an elementary property
of the adjoint. Now%
\begin{align*}
\left\Vert \widetilde{\mathcal{E}}_{X}^{\ast}\beta\right\Vert _{w,\theta}%
^{2}=\left(  \sum_{i=1}^{N}\beta_{i}R_{x^{\left(  i\right)  }},\sum_{j=1}%
^{N}\beta_{j}R_{x^{\left(  j\right)  }}\right)  _{w,\theta} &  =\sum
_{i,j=1}^{N}\beta_{i}\overline{\beta_{j}}\left(  R_{x^{\left(  i\right)  }%
},R_{x^{\left(  j\right)  }}\right)  _{w,\theta}\\
&  =\sum_{i,j=1}^{N}\beta_{i}\overline{\beta_{j}}R_{x^{\left(  i\right)  }%
}\left(  x^{\left(  j\right)  }\right) \\
&  =\sum_{i,j=1}^{N}\beta_{i}\overline{\beta_{j}}\overline{R_{x^{\left(
j\right)  }}\left(  x^{\left(  i\right)  }\right)  }=\beta^{T}\overline
{R_{X,X}}\overline{\beta}\\
&  =\overline{\beta}^{T}R_{X,X}\beta,
\end{align*}

so that, $\left\Vert \widetilde{\mathcal{E}}_{X}^{\ast}\right\Vert
=\max\limits_{\beta\in\mathbb{C}^{N}}\frac{\left\Vert \widetilde{\mathcal{E}%
}_{X}^{\ast}\beta\right\Vert _{w,\theta}}{\left\vert \beta\right\vert }%
=\max\limits_{\beta\in\mathbb{C}^{N}}\frac{\sqrt{\overline{\beta}^{T}%
R_{X,X}\beta}}{\left\vert \beta\right\vert }=\max\limits_{\beta\in
\mathbb{C}^{N}}\frac{\sqrt{\beta^{T}R_{X,X}\overline{\beta}}}{\left\vert
\beta\right\vert }$ and this expression is the largest (positive) eigenvalue
of the Hermitian matrix $R_{X,X}$. But this is also the value of $\left\Vert
R_{X,X}\right\Vert $ so $\left\Vert \widetilde{\mathcal{E}}_{X}^{\ast
}\right\Vert =\left\Vert R_{X,X}\right\Vert $.\medskip

\textbf{Part 4}
\[
\widetilde{\mathcal{E}}_{X}\widetilde{\mathcal{E}}_{X}^{\ast}\beta
=\widetilde{\mathcal{E}}_{X}\left(  \sum_{j=1}^{N}\beta_{j}R_{x^{\left(
j\right)  }}\right)  =\left(  \sum_{j=1}^{N}\beta_{j}R_{x^{\left(  j\right)
}}\left(  x^{\left(  i\right)  }\right)  \right)  _{j=1,N}=\left(
R_{x^{\left(  j\right)  }}\left(  x^{\left(  i\right)  }\right)  \right)
\beta=R_{X,X}\beta.
\]
\medskip

\textbf{Part 5} By part 2 Theorem \ref{Thm_Riesz_property}, $R_{a^{\left(
j\right)  }}=l_{_{j}}$ for $j\leq M$%
\[
\widetilde{\mathcal{E}}_{X}\widetilde{\mathcal{E}}_{A}^{\ast}\gamma
=\widetilde{\mathcal{E}}_{X}\left(  \sum_{j=1}^{M}\gamma_{j}R_{a^{\left(
j\right)  }}\right)  =\widetilde{\mathcal{E}}_{X}\left(  \sum_{j=1}^{M}%
\gamma_{j}l_{_{j}}\right)  =\left(  \sum_{j=1}^{M}\gamma_{j}l_{_{j}}\left(
x^{\left(  i\right)  }\right)  \right)  _{i=1,N}=L_{X}\gamma.
\]
\medskip

\textbf{Part 6} By part 1, $\widetilde{\mathcal{E}}_{X}:X_{w}^{\mathbb{\theta
}}\rightarrow\mathbb{C}^{N}$ is onto and $\operatorname*{null}\widetilde
{\mathcal{E}}_{X}=W_{G,X}^{\bot}$. By part 2, $\widetilde{\mathcal{E}}%
_{X}^{\ast}:\mathbb{C}^{N}\rightarrow W_{G,X}$ is onto and 1-1. These facts
directly imply $\widetilde{\mathcal{E}}_{X}^{\ast}\widetilde{\mathcal{E}}%
_{X}:X_{w}^{\theta}\mathbb{\rightarrow}W_{G,X}$ is onto and
$\operatorname*{null}\widetilde{\mathcal{E}}_{X}^{\ast}\widetilde{\mathcal{E}%
}_{X}=W_{G,X}^{\perp}$.%
\[
\left(  \widetilde{\mathcal{E}}_{X}^{\ast}\widetilde{\mathcal{E}}_{X}f\right)
\left(  x\right)  =\left(  \widetilde{\mathcal{E}}_{X}^{\ast}\widetilde
{\mathcal{E}}_{X}f,R_{x}\right)  _{w,\theta}=\left(  \widetilde{\mathcal{E}%
}_{X}f,\widetilde{\mathcal{E}}_{X}R_{x}\right)  _{\mathbb{C}^{N}}=\sum
_{j=1}^{N}f\left(  x^{\left(  j\right)  }\right)  \overline{R_{x}\left(
x^{\left(  j\right)  }\right)  }=\sum_{j=1}^{N}f\left(  x^{\left(  j\right)
}\right)  R_{x^{\left(  j\right)  }}\left(  x\right)  ,
\]

so that%
\begin{equation}
\widetilde{\mathcal{E}}_{X}^{\ast}\widetilde{\mathcal{E}}_{X}f=\sum_{j=1}%
^{N}f\left(  x^{\left(  j\right)  }\right)  R_{x^{\left(  j\right)  }%
}.\label{q65}%
\end{equation}

Also $\left(  \widetilde{\mathcal{E}}_{X}^{\ast}\widetilde{\mathcal{E}}%
_{X}f\right)  \left(  x\right)  =\left(  f,\widetilde{\mathcal{E}}_{X}^{\ast
}\widetilde{\mathcal{E}}_{X}R_{x}\right)  _{w,\theta}$, but by \ref{q65},
$\widetilde{\mathcal{E}}_{X}^{\ast}\widetilde{\mathcal{E}}_{X}R_{x}%
=\sum\limits_{j=1}^{N}R_{x^{\left(  j\right)  }}\left(  x\right)
R_{x^{\left(  j\right)  }}$, so that%
\[
\left(  \widetilde{\mathcal{E}}_{X}^{\ast}\widetilde{\mathcal{E}}_{X}f\right)
\left(  x\right)  =\left(  f,\sum_{j=1}^{N}R_{x^{\left(  j\right)  }}\left(
x\right)  R_{x^{\left(  j\right)  }}\right)  _{w,\theta}.
\]
\medskip

\textbf{Part 7} We must show that $\left(  \mathcal{P}\widetilde{\mathcal{E}%
}_{X}^{\ast}\alpha\right)  \left(  x\right)  =\alpha^{T}L_{X}\widetilde
{l}\left(  x\right)  $.%
\begin{align*}
\left(  \mathcal{P}\widetilde{\mathcal{E}}_{X}^{\ast}\alpha\right)  \left(
x\right)  =\mathcal{P}\sum\limits_{j=1}^{N}\alpha_{j}R_{x^{\left(  j\right)
}}\left(  x\right)  =\sum\limits_{j=1}^{N}\alpha_{j}\mathcal{P}R_{x^{\left(
j\right)  }}\left(  x\right)   &  =\sum\limits_{j=1}^{N}\alpha_{j}%
\sum\limits_{i=1}^{M}R_{x^{\left(  j\right)  }}\left(  a^{\left(  i\right)
}\right)  l_{i}\left(  x\right) \\
&  =\sum\limits_{j=1}^{N}\alpha_{j}\sum\limits_{i=1}^{M}R_{a^{\left(
i\right)  }}\left(  x^{\left(  j\right)  }\right)  l_{i}\left(  x\right) \\
&  =\sum\limits_{j=1}^{N}\alpha_{j}\sum\limits_{i=1}^{M}l_{i}\left(
x^{\left(  j\right)  }\right)  l_{i}\left(  x\right) \\
&  =\alpha^{T}L_{X}\widetilde{l}\left(  x\right)  .
\end{align*}
\medskip

\textbf{Part 8} Since $L_{X}=\left(  l_{j}\left(  x^{\left(  i\right)
}\right)  \right)  $%
\[
\widetilde{\mathcal{E}}_{X}\mathcal{P}f=\widetilde{\mathcal{E}}_{X}%
\sum\limits_{j=1}^{M}f\left(  a^{\left(  j\right)  }\right)  l_{j}%
=\sum\limits_{j=1}^{M}f\left(  a^{\left(  j\right)  }\right)  \widetilde
{\mathcal{E}}_{X}l_{j}=L_{X}\widetilde{\mathcal{E}}_{A}f.
\]

\end{proof}

\section{Variational interpolation with basis
functions\label{Sect_interpol_prob}}

The mathematical machinery we have developed above will now be applied to
studying two variational interpolation problems. These will be called the
minimal norm and minimal seminorm problems respectively. The minimal seminorm
problem will involve minimizing the seminorm of the space $X_{w}^{\theta}$
over the functions in $X_{w}^{\theta}$ which interpolate the data points
$\left\{  \left(  x^{\left(  i\right)  },y_{i}\right)  \right\}  _{i=1}^{N}$.
The minimal norm problem will involve minimizing the Light norm of the space
$X_{w}^{\theta}$ over the functions in $X_{w}^{\theta}$ which interpolate the
data. We show the minimal norm problem has a unique solution and then show the
minimal seminorm problem has the same, unique solution and that this solution
lies in the finite dimensional space $W_{G,X}$ where $G$ is an order $\theta$
basis function of the weight function $w$ and $X=\{x^{(i)}\}_{i=1}^{N}$ i.e.
the solution is a basis function interpolant. A matrix equation is also
derived for the coefficients of the basis functions.

To render the interpolation problem meaningful we suppose the weight function
$w$ has weight properties W2 and W3 for order $\theta$ and smoothness
parameter $\kappa$, so that $X_{w}^{\theta}$ is a space of continuous
functions. Further, the conditions on $w$ allow us to define a continuous
basis function $G$ of order $\theta$.

\subsection{The minimal seminorm interpolation problem}

The scattered data set to be interpolated consists of interpolating the $N$
data points
\[
\left\{  \left(  x^{(i)},y_{i}\right)  \right\}  _{i=1}^{N},\quad x^{(i)}%
\in\mathbb{R}^{d},\text{ }x^{(i)}\text{ }distinct,\text{ }y_{i}\in\mathbb{R}.
\]

We call $X=\{x^{(i)}\}_{i=1}^{N}$ the independent data and $y=\{y_{i}%
\}_{i=1}^{N}$ the dependent data, and the data will sometimes be referred to
as $\left[  X,y\right]  $. The points in $X$ are not completely unconstrained
because it will be assumed that $X$\textbf{\ }is unisolvent with respect to
the polynomials $P_{\theta-1}$ (Definition \ref{Def_unisolv}). We now require
that a \textit{minimal seminorm interpolant} $u_{I}$ satisfies
\[
\left\vert u_{I}\right\vert _{w,\theta}=\min\left\{  \left\vert u\right\vert
_{w,\theta}:u\in X_{w}^{\mathbb{\theta}};\text{ }u_{I}\left(  x^{(i)}\right)
=y_{i},\text{ }i=1,\ldots,N\right\}  .
\]

\subsection{The minimal norm interpolation problem}

Since the independent data $X$ is $\theta$-unisolvent, from Definition
\ref{Def_unisolv} there is at least one minimal unisolvent subset of $X$, say
$A$. Denote by $b$ the dependent data corresponding to $A$. The Light norm
$\left\Vert \cdot\right\Vert _{w,\theta}$ of \ref{p917} is now constructed
using $A$. The scattered data is the same as for the minimal seminorm
interpolation problem i.e. the \textit{minimal norm interpolant} $u_{I}$
satisfies
\begin{equation}
\left\Vert u_{I}\right\Vert _{w,\theta}=\min\left\{  \left\Vert u_{I}%
\right\Vert _{w,\theta}:u\in X_{w}^{\theta},\text{ }u_{I}\left(
x^{(i)}\right)  =y_{i}\text{ }for\text{ }i=1,\ldots,N\right\}  .\label{q59}%
\end{equation}

Using Hilbert space techniques, we now want to do the following:

\begin{enumerate}
\item Show their exists a unique minimal norm interpolant.

\item Show their exists a unique minimal seminorm interpolant.

\item Show these interpolants are identical.

\item Show the interpolant is a basis function interpolant lying in the space
$W_{G,X}$ i.e. it can be written as a linear combination of data-translated
basis functions plus a polynomial of degree less than $\theta$.
\end{enumerate}

\subsection{Solving the minimal interpolation problems}

In this subsection we show the minimal norm interpolation problem has a unique
solution which we call the minimal norm interpolant. It is then shown that
this solution is also the unique solution to the minimum seminorm
interpolation problem. We then derive a formula for this solution which
implies that the solution lies in the finite dimensional space $W_{G,X}$
introduced in Definition \ref{Def_Wgx}. The last result will allow us to
derive matrix results for the interpolant in the next subsection.

\begin{theorem}
\label{Thm_var_interpol_norm}Endow $X_{w}^{\theta}$ with the Light norm
$\left\Vert \cdot\right\Vert _{w,\theta}$ defined by \ref{p917} using $A$.
Then there exists a unique minimal norm interpolant. In fact, given
$y\in\mathbb{R}^{N}$ there is a unique interpolant $u_{I}\in X_{w}^{\theta}$
such that
\begin{equation}
\left\Vert u_{I}\right\Vert _{w,\theta}=\min\left\{  \left\Vert u\right\Vert
_{w,\theta}:u\in X_{w}^{\theta}\text{ }and\text{ }\widetilde{\mathcal{E}}%
_{X}u=y\right\}  .\label{p945}%
\end{equation}

If $v\in X_{w}^{\mathbb{\theta}}$ is any other interpolant satisfying
$\widetilde{\mathcal{E}}_{X}v=y$ then
\begin{equation}
\left\Vert u_{I}\right\Vert _{w,\theta}^{2}+\left\Vert v-u_{I}\right\Vert
_{w,\theta}^{2}=\left\Vert v\right\Vert _{w,\theta}^{2},\label{q60}%
\end{equation}

or equivalently
\begin{equation}
\left(  v-u_{I},u_{I}\right)  _{w,\theta}=0.\label{p922}%
\end{equation}

\end{theorem}

\begin{proof}
Now by part 2 of Theorem \ref{Thm_eval_op_properties} the evaluation operator
$\widetilde{\mathcal{E}}_{X}:X_{w}^{\mathbb{\theta}}\rightarrow\mathbb{R}^{N}$
is continuous and onto, and since the singleton set $\left\{  y\right\}  $ is
closed, it follows that the set
\[
\left\{  u:\widetilde{\mathcal{E}}_{X}u=y\right\}  =\widetilde{\mathcal{E}%
}_{X}^{-1}\left(  y\right)  ,
\]

is a non-empty, proper, closed subspace of the Hilbert space $X_{w}%
^{\mathbb{\theta}}$. Hence this subspace contains a unique element of smallest
norm, say $u_{I}$, which satisfies equation \ref{p945}. If $\widetilde
{\mathcal{E}}_{X}v=y$ then
\[
\left\{  u:u\in X_{w}^{\mathbb{\theta}}\text{ and }\widetilde{\mathcal{E}}%
_{X}u=y\right\}  =\left\{  v-s:s\in\operatorname*{null}\widetilde{\mathcal{E}%
}_{X}\right\}  .
\]

Now
\[
\min\left\{  \left\Vert v-s\right\Vert _{w,\theta}:s\in\operatorname*{null}%
\widetilde{\mathcal{E}}_{X}\right\}  =\operatorname*{dist}\left(
v,\operatorname*{null}\widetilde{\mathcal{E}}_{X}\right)  ,
\]

is the distance between $v$ and the closed subspace $\operatorname*{null}%
\widetilde{\mathcal{E}}_{X}$. Therefore there exists a unique $s_{I}%
\in\operatorname*{null}\widetilde{\mathcal{E}}_{X}$ such that
\begin{equation}
u_{I}=v-s_{I},\label{p923}%
\end{equation}

\[
\left\Vert v-s_{I}\right\Vert _{w,\theta}=\min\left\{  \left\Vert
v-s\right\Vert _{w,\theta}:s\in\operatorname*{null}\widetilde{\mathcal{E}}%
_{X}\right\}  ,
\]

and
\begin{equation}
\left\Vert s_{I}\right\Vert _{w,\theta}^{2}+\left\Vert v-s_{I}\right\Vert
_{w,\theta}^{2}=\left\Vert v\right\Vert _{w,\theta}^{2}.\label{p924}%
\end{equation}

Substituting for $s_{I}$ in \ref{p924} using \ref{p923} yields \ref{q60}.
Equation \ref{p922} follows since it is a necessary and sufficient condition
for \ref{q60} to be true.
\end{proof}

\begin{theorem}
\label{Thm_var_interpol_seminorm}The unique minimal norm interpolant $u_{I}$
of Theorem \ref{Thm_var_interpol_norm} is the unique minimal seminorm
interpolant. In fact, suppose $y\in\mathbb{R}^{N}$. Then
\begin{equation}
\left\vert u_{I}\right\vert _{w,\theta}=\min\left\{  \left\vert u\right\vert
_{w,\theta}:u\in X_{w}^{\mathbb{\theta}}\text{ }and\text{ }\widetilde
{\mathcal{E}}_{X}u=y\right\}  .\label{p934}%
\end{equation}

If $v$ is any other interpolant satisfying $\widetilde{\mathcal{E}}_{X}v=y$
then
\begin{equation}
\left\vert u_{I}\right\vert _{w,\theta}^{2}+\left\vert v-u_{I}\right\vert
_{w,\theta}^{2}=\left\vert v\right\vert _{w,\theta}^{2},\label{p941}%
\end{equation}

or equivalently
\begin{equation}
\left\langle v-u_{I},u_{I}\right\rangle _{w,\theta}=0.\label{p942}%
\end{equation}

\end{theorem}

\begin{proof}
Recall that the Light norm \ref{p917} for the minimal norm interpolant problem
\ref{q59} was constructed using the minimally unisolvent subset $A$ of
independent data which had corresponding dependent data $y^{\prime}$. Since
$\widetilde{\mathcal{E}}_{X}f=y$
\[
\left\Vert f\right\Vert _{w,\theta}^{2}=\left\vert f\right\vert _{w,\theta
}^{2}+\sum_{x^{\left(  i\right)  }\in A}\left\vert f(x^{(i)})\right\vert
^{2}=\left\vert f\right\vert _{w,\theta}^{2}+\left\vert y^{\prime}\right\vert
^{2},
\]

By Theorem \ref{Thm_var_interpol_norm} there exists a unique minimal norm
interpolant $u_{I}$ such that
\[
\left\Vert u_{I}\right\Vert _{w,\theta}=\min\left\{  \left\Vert u\right\Vert
_{w,\theta}:u\in X_{w}^{\mathbb{\theta}}\text{ }and\text{ }\widetilde
{\mathcal{E}}_{X}u=y\right\}  ,
\]

and hence
\begin{align*}
\left\vert u_{I}\right\vert _{w,\theta}^{2}+\left\vert b\right\vert ^{2} &
=\min\left\{  \left\vert u\right\vert _{w,\theta}^{2}+\left\vert y^{\prime
}\right\vert ^{2}:u\in X_{w}^{\theta}\text{ }and\text{ }\widetilde
{\mathcal{E}}_{X}u=y\right\} \\
&  =\min\left\{  \left\vert u\right\vert _{w,\theta}^{2}:u\in X_{w}%
^{\mathbb{\theta}}\text{ }and\text{ }\widetilde{\mathcal{E}}_{X}u=y\right\}
+\left\vert y^{\prime}\right\vert ^{2},
\end{align*}

so that
\[
\left\Vert u_{I}\right\Vert _{w,\theta}=\min\left\{  \left\vert u\right\vert
_{w,\theta}:u\in X_{w}^{\mathbb{\theta}}\text{ }and\text{ }\widetilde
{\mathcal{E}}_{X}u=y\right\}  .
\]

To prove equalities \ref{p941} and \ref{p942} first note that if
$\widetilde{\mathcal{E}}_{X}f=0$ then $\left\Vert f\right\Vert _{w,\theta}%
^{2}=\left\vert f\right\vert _{w,\theta}^{2}+\sum_{i=1}^{M}\left\vert
f(x^{(i)})\right\vert ^{2}=\left\vert f\right\vert _{w,\theta}^{2}$. Thus,
since $\widetilde{\mathcal{E}}_{X}\left(  v-u_{I}\right)  =0$,
\begin{align*}
\left\vert u_{I}\right\vert _{w,\theta}^{2}+\left\vert v-u_{I}\right\vert
_{w,\theta}^{2}=\left\vert u_{I}\right\vert _{w,\theta}^{2}+\left\Vert
v-u_{I}\right\Vert _{w,\theta}^{2} &  =\left\Vert u_{I}\right\Vert _{w,\theta
}^{2}-\left\vert y^{\prime}\right\vert ^{2}+\left\Vert v-u_{I}\right\Vert
_{w,\theta}^{2}\\
&  =\left\Vert v\right\Vert _{w,\theta}^{2}-\left\vert y^{\prime}\right\vert
^{2},\quad by\text{ }\ref{q60},\\
&  =\left\vert v\right\vert _{w,\theta}^{2},
\end{align*}

which proves \ref{p941}. To prove \ref{p942} we use equation \ref{p922} to
obtain:
\[
0=\left(  v-u_{I},u_{I}\right)  _{w,\theta}=\left\langle v-u_{I}%
,u_{I}\right\rangle _{w,\theta}+\sum_{i=1}^{M}\left(  v-u_{I}\right)  \left(
x^{(i)}\right)  \overline{u_{I}\left(  x^{(i)}\right)  }=\left\langle
v-u_{I},u_{I}\right\rangle _{w,\theta},
\]

since $v\left(  x^{(i)}\right)  =u_{I}\left(  x^{(i)}\right)  $ for
$i=1,\ldots,N$.
\end{proof}

There is a very close connection between the minimal norm interpolants, the
space $W_{G,X}$ and the adjoint $\widetilde{\mathcal{E}}_{X}^{\ast}$ \ref{q75}
of the evaluation operator $\widetilde{\mathcal{E}}_{X}$. The next theorem
describes this connection and characterizes the space of minimal norm
interpolants associated with a given unisolvent set of independent data points
$X$.

The minimal norm interpolation problem was formulated using a Light norm
$\left\Vert \cdot\right\Vert _{w,\theta}$ \ref{p917} constructed with a
minimal unisolvent set $A\subset X$ and the corresponding cardinal basis of
$A$. In the next theorem we will also construct the Riesz representer $R_{x}$
\ref{p939} using $A$ and its cardinal basis. Then by Definition
\ref{Def_Matrices_from_Rx} the\textbf{\ }reproducing kernel matrix\textbf{\ }%
$R_{X,X}$ is given by $R_{X,X}=\left(  R_{x^{\left(  j\right)  }}(x^{\left(
i\right)  })\right)  $.

\begin{theorem}
\label{Thm_min_norm_interpol_in_Wgx}We know from Theorem
\ref{Thm_var_interpol_norm} that for given dependent data $\left[  X;y\right]
$ we can associate a unique minimal norm interpolant $u_{I}\in X_{w}%
^{\mathbb{\theta}}$. In fact
\begin{equation}
u_{I}=\widetilde{\mathcal{E}}_{X}^{\ast}\left(  R_{X,X}\right)  ^{-1}y\in
W_{G,X},\label{p931}%
\end{equation}

and the space of all minimal norm interpolants is $W_{G,X}$. Here $R_{X,X}$ is
the\textbf{\ }reproducing kernel matrix.
\end{theorem}

\begin{proof}
From parts 3 and 4 of Theorem \ref{Thm_eval_op_properties} $\widetilde
{\mathcal{E}}_{X}\widetilde{\mathcal{E}}_{X}^{\ast}=R_{X,X}$ and so
\[
\widetilde{\mathcal{E}}_{X}\widetilde{\mathcal{E}}_{X}^{\ast}y=\left(  \left(
\widetilde{\mathcal{E}}_{X}^{\ast}y\right)  \left(  x^{\left(  i\right)
}\right)  \right)  =R_{X,X}y,
\]

where $R_{X,X}$ is regular by part 2 Theorem \ref{Thm_int_Rxx_property}. In
general $R_{X,X}\neq I$ so in general $\widetilde{\mathcal{E}}_{X}%
\widetilde{\mathcal{E}}_{X}^{\ast}y\neq y$ and $\widetilde{\mathcal{E}}%
_{X}^{\ast}y$ is not an interpolant of $y$.

However, $\widetilde{\mathcal{E}}_{X}\widetilde{\mathcal{E}}_{X}^{\ast}\left(
R_{X,X}\right)  ^{-1}=I$, so $\widetilde{\mathcal{E}}_{X}^{\ast}\left(
R_{X,X}\right)  ^{-1}y$ is an interpolant and

$\widetilde{\mathcal{E}}_{X}^{\ast}\left(  R_{X,X}\right)  ^{-1}y\in W_{G,X}$.
For convenience set $u=\widetilde{\mathcal{E}}_{X}^{\ast}\left(
R_{X,X}\right)  ^{-1}y$. We want to show that $u=u_{I}$. But since both $u_{I}
$ and $u$ interpolate $y$
\begin{align*}
\left(  u,u_{I}-u\right)  _{w,\theta}=\left(  \widetilde{\mathcal{E}}%
_{X}^{\ast}\left(  R_{X,X}\right)  ^{-1}y,u_{I}-u\right)  _{w,\theta} &
=\left(  \left(  R_{X,X}\right)  ^{-1}y,\widetilde{\mathcal{E}}_{X}\left(
u_{I}-u\right)  \right)  _{w,\theta}\\
&  =0,
\end{align*}

and so
\[
\left\Vert u_{I}\right\Vert _{w,\theta}^{2}=\left\Vert u_{I}-u+u\right\Vert
_{w,\theta}^{2}=\left\Vert u_{I}-u\right\Vert _{w,\theta}^{2}+\left\Vert
u\right\Vert _{w,\theta}^{2}.
\]

If $u\neq u_{I}$ then $\left\Vert u\right\Vert _{w,\theta}<\left\Vert
u_{I}\right\Vert _{w,\theta}$, contradicting the fact that $u_{I}$ is the
minimal norm interpolant. Because $R_{X,X}$ is regular, $\left(
R_{X,X}\right)  ^{-1}$ is an isomorphism from $\mathbb{R}^{N}$ to
$\mathbb{R}^{N}$, and because $\widetilde{\mathcal{E}}_{X}^{\ast}$ is an
isomorphism from $\mathbb{R}^{N}$ to $W_{G,X}$, the formula $u_{I}%
=\widetilde{\mathcal{E}}_{X}^{\ast}\left(  R_{X,X}\right)  ^{-1}y$ implies
that $\operatorname*{range}\left(  \widetilde{\mathcal{E}}_{X}^{\ast}\left(
R_{X,X}\right)  ^{-1}\right)  =W_{G,X}$.
\end{proof}

\subsection{Matrix equations for the minimal seminorm
interpolant\label{SbSect_mat_eqn_interpol}}

We now know from Theorem \ref{Thm_min_norm_interpol_in_Wgx}\ that the minimal
seminorm interpolant lies in the finite dimensional space $W_{G,X}$. Hence we
should now be able to derive a matrix equation for the coefficients of the
$G\left(  \cdot-x^{\left(  k\right)  }\right)  $. This equation will be
constructed using the basis function matrix and the unisolvent matrix. This
will require the following lemma which turns out to be very useful since it
provides key results for a block matrix which has a structure central to the
study of the matrix equations in both basis function interpolation and basis
function smoothing theory.

\begin{lemma}
\label{Lem_reg_cspd}Let $B$ be a complex-valued matrix and $C$ be a
real-valued matrix. Suppose the block matrix $\left(
\begin{array}
[c]{ll}%
B & C\\
C^{T} & O
\end{array}
\right)  $ is square.

Suppose further that for complex vectors $z$
\begin{equation}
z^{T}B\overline{z}=0\text{ }and\text{ }C^{T}z=0\text{ }implies\text{
}z=0.\label{q62}%
\end{equation}

\begin{enumerate}
\item Then the equation
\[
\left(
\begin{array}
[c]{ll}%
B & C\\
C^{T} & O
\end{array}
\right)  \left(
\begin{array}
[c]{l}%
u\\
v
\end{array}
\right)  =\left(
\begin{array}
[c]{l}%
0\\
0
\end{array}
\right)  ,
\]

implies $u=0$ and $v\in\operatorname*{null}C$.

\item If, further to part 1, $\operatorname*{null}C=\left\{  0\right\}  $ then
the block matrix is regular.
\end{enumerate}
\end{lemma}

\begin{proof}
\textbf{Part 1} We proceed by showing that%
\[
\left(
\begin{array}
[c]{ll}%
B & C\\
C^{T} & O
\end{array}
\right)  \left(
\begin{array}
[c]{l}%
u\\
v
\end{array}
\right)  =\left(
\begin{array}
[c]{l}%
0\\
0
\end{array}
\right)  ,
\]

implies $u=0$. Expanding this block matrix equation yields $Bu+Cv=0$ and
$C^{T}u=0$. Thus%
\[
0=u^{T}(\overline{B}\overline{u}+C\overline{v})=u^{T}\overline{B}\overline
{u}+u^{T}C\overline{v}=u^{T}\overline{B}\overline{u}+(C^{T}u)^{T}\overline
{v}=u^{T}\overline{B}\overline{u},
\]

and so $\overline{u}^{T}Bu=0$. Set $x=\overline{u}$. Then $x^{T}B\overline
{x}=0$ and $C^{T}x=0$, and so \ref{q62} implies $x=0$. Hence $u=0$ and
$Bv=0$.\medskip

\textbf{Part 2} Clearly if $\operatorname*{null}C=\left\{  0\right\}  $, by
part 1 the block matrix has null space zero and so is regular.
\end{proof}

We now derive the matrix equation for minimal seminorm interpolant:

\begin{theorem}
\label{Thm_interpol_matrix_eqn}Suppose $\left\{  p_{j}\right\}  _{j=1}^{M}$ is
any basis for $P_{\theta-1}$. The space $W_{G,X}$ contains only one minimal
seminorm interpolant $u_{I}$ to any given data vector $y\in\mathbb{R}^{N}$.
This minimal seminorm interpolant $u_{I}$ is given uniquely by
\begin{equation}
u_{I}\left(  x\right)  =\sum_{i=1}^{N}v_{i}G\left(  x-x^{(i)}\right)
+\sum_{j=1}^{M}\beta_{j}p_{j}(x),\label{q63}%
\end{equation}

where the coefficient vectors $v=\left(  v_{i}\right)  $ and $\beta=\left(
\beta_{j}\right)  $ satisfy the matrix equation
\begin{equation}
\left(
\begin{array}
[c]{ll}%
G_{X,X} & P_{X}\\
P_{X}^{T} & O
\end{array}
\right)  \left(
\begin{array}
[c]{l}%
v\\
\beta
\end{array}
\right)  =\left(
\begin{array}
[c]{l}%
y\\
0
\end{array}
\right)  .\label{q64}%
\end{equation}

Here $P_{X}$ is the unisolvency matrix $\left(  p_{j}\left(  x^{\left(
i\right)  }\right)  \right)  $ and $G_{X,X}$ is the basis function matrix
$\left(  G\left(  x^{\left(  i\right)  }-x^{\left(  j\right)  }\right)
\right)  $.

Also, the matrix $\left(
\begin{array}
[c]{ll}%
G_{X,X} & P_{X}\\
P_{X}^{T} & O
\end{array}
\right)  $ has size $\left(  N+M\right)  \times\left(  N+M\right)  $ and is regular.
\end{theorem}

\begin{proof}
By Theorem \ref{Thm_min_norm_interpol_in_Wgx}, $W_{G,X}$ contains only one
minimal seminorm interpolant. Using equation \ref{q63} and the interpolation
requirement $u_{I}\left(  x^{(i)}\right)  =y_{i}$ we obtain the pair of
equations
\[
G_{X,X}v+P_{X}\beta=y;\text{\quad}P_{X}^{T}v=0,
\]

or in block matrix form
\[
\left(
\begin{array}
[c]{ll}%
G_{X,X} & P_{X}\\
P_{X}^{T} & O
\end{array}
\right)  \left(
\begin{array}
[c]{l}%
v\\
\beta
\end{array}
\right)  =\left(
\begin{array}
[c]{l}%
y\\
0
\end{array}
\right)  .
\]

By part 3 Theorem \ref{Thm_Wgx_properties_1}, $v^{T}G_{X,X}\overline{v}=0$ and
$P_{X}^{T}v=0$ implies $v=0$, and since $X$ is unisolvent, part 1 Theorem
\ref{Thm_Px_properties} implies $\operatorname*{null}P_{X}=\left\{  0\right\}
$. Thus, by Lemma \ref{Lem_reg_cspd} above, the matrix of \ref{q64} is regular
and the interpolation problem has a unique solution $\left(  v^{T},\beta
^{T}\right)  $.
\end{proof}

From part 1 Corollary \ref{Cor_Thm_Wgx_properties_2} $\left\{  R_{x^{\left(
i\right)  }}\right\}  _{i=1}^{N}$ is a basis for $W_{G,X}$. In the next
theorem we give the corresponding matrix equation for the minimal norm interpolant.

\begin{theorem}
\label{Thm_interpol_mat_eqn_Rx}We know from Theorem
\ref{Thm_var_interpol_norm} that for given independent data $X=\left\{
x^{\left(  i\right)  }\right\}  $ and dependent data $y\in\mathbb{R}^{N}$ we
can associate a unique minimal norm interpolant $u_{I}\in X_{w}%
^{\mathbb{\theta}}$. In fact, if $R_{X,X}=\left\{  R_{x^{\left(  j\right)  }%
}\left(  x^{\left(  i\right)  }\right)  \right\}  $ is the reproducing kernel
matrix then
\begin{equation}
u_{I}=\sum_{i=1}^{N}v_{i}R_{x^{\left(  i\right)  }},\quad when\text{\ }%
R_{X,X}v=y,\label{p946}%
\end{equation}

where $R_{x}$ is the Riesz representer of the evaluation functional
$f\rightarrow f\left(  x\right)  $.
\end{theorem}

\begin{proof}
$u_{I}=\widetilde{\mathcal{E}}_{X}^{\ast}\left(  R_{X,X}\right)
^{-1}y=\widetilde{\mathcal{E}}_{X}^{\ast}v=\sum_{i=1}^{N}v_{i}R_{x^{\left(
i\right)  }}$.
\end{proof}

\subsection{The minimal interpolant mapping and data functions}

The last theorem allows us to define the mapping between a data function and
its corresponding interpolant. We call this the \textit{minimal interpolant
mapping}:

\begin{definition}
\label{Def_data_func_interpol_map}\textbf{Data functions and the minimal
interpolant mapping} $\mathcal{I}_{X}:X_{w}^{\theta}\rightarrow W_{G,X}$

Given an independent data set $X$, we shall assume that each member of
$X_{w}^{\theta}$ can act as a legitimate data function $f$ and generate the
data vector $\widetilde{\mathcal{E}}_{X}f$.

Equation \ref{p931} of Theorem \ref{Thm_min_norm_interpol_in_Wgx} enables us
to define the linear mapping $\mathcal{I}_{X}:X_{w}^{\theta}\rightarrow
W_{G,X}$ from the data functions to the corresponding unique minimal
interpolant $\mathcal{I}_{X}f=u_{I}$ given by%
\begin{equation}
\mathcal{I}_{X}f=\widetilde{\mathcal{E}}_{X}^{\ast}\left(  R_{X,X}\right)
^{-1}\widetilde{\mathcal{E}}_{X}f,\quad f\in X_{w}^{\theta}.\label{q61}%
\end{equation}

\end{definition}

The purpose of the minimal interpolant mapping will be to study the pointwise
convergence of the interpolant to its data function.

\begin{theorem}
\label{Thm_propert_interpol_map}The minimal interpolant mapping $\mathcal{I}%
_{X}:X_{w}^{\theta}\rightarrow W_{G,X}$ has the following properties:

\begin{enumerate}
\item $\left\Vert \mathcal{I}_{X}f\right\Vert _{w,\theta}\leq\left\Vert
f\right\Vert _{w,\theta}$ and $\left\Vert \left(  I-\mathcal{I}_{X}\right)
f\right\Vert _{w,\theta}\leq\left\Vert f\right\Vert _{w,\theta}$.

\item $\left\vert \mathcal{I}_{X}f\right\vert _{w,\theta}\leq\left\vert
f\right\vert _{w,\theta}$ and $\left\vert \left(  I-\mathcal{I}_{X}\right)
f\right\vert _{w,\theta}\leq\left\vert f\right\vert _{w,\theta}$.

\item $\mathcal{I}_{X}$ is a continuous projection onto $W_{G,X}$ with null
space $W_{G,X}^{\bot}$.

\item $\mathcal{I}_{X}$ is self-adjoint.
\end{enumerate}
\end{theorem}

\begin{proof}
\textbf{Part 1} Follows from \ref{q60} since the data function is an
interpolant.\smallskip

\textbf{Part 2} Follows from part 1 since $\mathcal{I}_{X}$ is an
interpolant.\smallskip

\textbf{Part 3} Part 1 implies continuity. Part 4 of Theorem
\ref{Thm_eval_op_properties} implies $\widetilde{\mathcal{E}}_{X}%
\widetilde{\mathcal{E}}_{X}^{\ast}=R_{X,X}$ and $R_{X,X}$ is Hermitian and
regular. Hence \ref{q61} implies $\mathcal{I}_{X}$ is a projection and thus
onto. Further, $\widetilde{\mathcal{E}}_{X}^{\ast}$ is 1-1 by part 2 of
Theorem \ref{Thm_eval_op_properties} so \ref{q61} implies $\mathcal{I}_{X}f=0$
iff $\widetilde{\mathcal{E}}_{X}f=0$ iff $f\in W_{G,X}^{\bot}$ by part 1 of
Theorem \ref{Thm_eval_op_properties}.\smallskip

\textbf{Part 4\ }By\textbf{\ }\ref{q61}, $\left(  \mathcal{I}_{X}f,g\right)
_{w,\theta}=\left(  \widetilde{\mathcal{E}}_{X}^{\ast}\left(  R_{X,X}\right)
^{-1}\widetilde{\mathcal{E}}_{X}f,g\right)  _{w,\theta}=$

$=\left(  \left(  R_{X,X}\right)  ^{-1}\widetilde{\mathcal{E}}_{X}%
f,\widetilde{\mathcal{E}}_{X}g\right)  _{w,\theta}=\left(  f,\widetilde
{\mathcal{E}}_{X}^{\ast}\left(  R_{X,X}\right)  ^{-1}\widetilde{\mathcal{E}%
}_{X}g\right)  _{w,\theta}$ since $R_{X,X}$ is Hermitian.
\end{proof}

\section{Convergence orders for the minimal
interpolant\label{Sect_interpol_converg}}

The estimates for the order of pointwise convergence obtained here are derived
using similar techniques to those of Light and Wayne in
\cite{LightWayne95ErrEst} and \cite{LightWayne98PowFunc}. We use the Lagrange
interpolation theory from \cite{LightWayne98PowFunc} and the function
$r_{x}\left(  x\right)  $.

\subsection{A unisolvency lemma}

To study the convergence of the minimal interpolant we will also need the
following lemma which supplies some elementary results from the theory of
Lagrange interpolation. These results are stated without proof. This lemma has
been created from Lemma 3.2, Lemma 3.5 and the first two paragraphs of the
proof of Theorem 3.6 of Light and Wayne \cite{LightWayne98PowFunc}. The
results of this lemma do not involve any reference to weight or basis
functions or functions in $X_{w}^{\theta}$, but consider the properties of the
set which contains the independent data points and the order of the
unisolvency used for the interpolation. Thus we have separated the part of the
proof that involves weight functions from the part that uses the detailed
theory of Lagrange interpolation operators.

\begin{lemma}
\label{Lem_int_Lagrange_interpol}Assume first that:

\begin{enumerate}
\item $\Omega$ is a bounded, open, connected subset of $\mathbb{R}^{d}$ having
the \textbf{cone property}.

\item $X$ is a unisolvent subset of $\Omega$ of order $\theta$.
\end{enumerate}

Suppose $\left(  l_{j}\right)  _{j=1}^{M}=\widetilde{l}_{A}$ is the
\textbf{cardinal basis} of $P_{\theta-1}$ with respect to a \textbf{minimal
unisolvent set} $A\subset\Omega$. Using \textbf{Lagrange interpolation}
techniques, it can be shown there exists a constant $K_{\Omega,\theta}%
^{\prime}>0$ such that
\[
\sum\limits_{j=1}^{M}\left\vert l_{j}\left(  x\right)  \right\vert =\left\vert
\widetilde{l}_{A}\left(  x\right)  \right\vert _{1}\leq K_{\Omega,\theta
}^{\prime},\quad x\in\overline{\Omega},
\]

for all minimal unisolvent subsets $A\subset\overline{\Omega}$. Now define the
\textbf{spherical cavity size} by
\begin{equation}
h_{X,\Omega}=\sup\limits_{\omega\in\Omega}\operatorname*{dist}\left(
\omega,X\right)  ,\label{p947}%
\end{equation}

and fix $x\in\Omega$. Using Lagrange interpolation techniques it can be shown
there are constants $c_{\Omega,\theta},h_{\Omega,\theta}>0$ such that when
$h_{X,\Omega}<h_{\Omega,\theta}$ there exists a minimal unisolvent set
$A\subset X$ satisfying
\[
\operatorname*{diam}A_{x}\leq c_{\Omega,\theta}h_{X,\Omega},
\]

where $A_{x}=A\cup\left\{  x\right\}  $.
\end{lemma}

\begin{remark}
If $\omega_{1}$ is one of the furthest points in $\overline{\Omega}$ from $X$
then $h_{X,\Omega}=\operatorname*{dist}\left(  \omega_{1},X\right)  $.
$h_{X,\Omega}$ can also be interpreted as the radius of the largest sphere
(cavity) centred in $\overline{\Omega}$ that can be placed between the points
in $X$.
\end{remark}

\subsection{An upper bound for $\sqrt{r_{x}\left(  x\right)  }$%
\label{SbSect_int_bound_rx(x)}}

\begin{theorem}
\label{Thm_int_rx(x)_bound}Suppose $w$ is a weight function with properties W2
and W3 for order $\theta$ and smoothness $\kappa$, and that $G$ is a basis
function of order $\theta$. Suppose $A=\left\{  a^{\left(  k\right)
}\right\}  _{k=1}^{M}$ is a minimal $\theta$-unisolvent set and that $\left\{
l_{k}\right\}  _{k=1}^{M}$ is the corresponding unique cardinal basis for
$P_{\theta-1}$. Now construct $\mathcal{P}$, $\mathcal{Q}$ and $R_{x}$ using
$A$ and $\left\{  l_{k}\right\}  _{k=1}^{M}$. Then if $r_{x}=\mathcal{Q}R_{x}$:

\begin{enumerate}
\item
\[
\sqrt{r_{x}\left(  x\right)  }\leq\frac{d^{\eta/2}}{\left(  2\pi\right)
^{\frac{d}{4}}}\frac{\left(  \operatorname*{diam}A_{x}\right)  ^{\eta}}%
{\sqrt{\left(  2\eta\right)  !}}\left(  1+\sum_{k=1}^{M}\left\vert
l_{k}\left(  x\right)  \right\vert \right)  \max\limits_{\substack{\left\vert
\beta\right\vert =2\eta\\\left\vert y\right\vert \leq\operatorname{diam}A_{x}%
}}\sqrt{\left\vert D^{\beta}G\left(  y\right)  \right\vert },\quad
x\in\mathbb{R}^{d},
\]

\item and suitable for radial basis functions:%
\[
\sqrt{r_{x}\left(  x\right)  }\leq\frac{1}{\left(  2\pi\right)  ^{\frac{d}{4}%
}}\frac{\left(  \operatorname*{diam}A_{x}\right)  ^{\eta}}{\sqrt{\left(
2\eta\right)  !}}\left(  1+\sum_{k=1}^{M}\left\vert l_{k}\left(  x\right)
\right\vert \right)  \max\limits_{\left\vert y\right\vert \leq
\operatorname{diam}A_{x}}\sqrt{\left\vert \left(  \left(  \widehat{\cdot
}D\right)  ^{2\eta}G\right)  \left(  y\right)  \right\vert },\quad
x\in\mathbb{R}^{d},
\]

\end{enumerate}

where $\eta=\min\left\{  \theta,\frac{1}{2}\left\lfloor \min2\kappa
\right\rfloor \right\}  $, $A_{x}=A\cup\left\{  x\right\}  $ and $\left(
\left(  \widehat{\cdot}D\right)  u\right)  \left(  x\right)  =\sum
\limits_{k=1}^{d}\widehat{x}_{k}D_{k}u\left(  x\right)  $, $\widehat
{x}=x/\left\vert x\right\vert $. For properties of the operator $\widehat
{\cdot}D$ and its action on radial functions see Lemma
\ref{Lem_deriv_rad_funcs} in the Appendix.
\end{theorem}

\begin{proof}
From part 4 of Theorem \ref{Thm_rx(y)_properties}, $r_{x}\left(  y\right)
=\left\langle r_{x},r_{y}\right\rangle _{w,\theta}$ so that $r_{x}\left(
x\right)  =\left\vert r_{x}\right\vert _{w,\theta}^{2}\geq0$.

Also from Theorem \ref{Thm_rx(y)_properties}, $r_{x}\left(  y\right)  =\left(
2\pi\right)  ^{-\frac{d}{2}}\mathcal{Q}_{y}\mathcal{Q}_{x}G\left(  y-x\right)
$ when $x\neq y$.

Now from Theorem \ref{Thm_basis_smth_W3.2_r3_pos} or \ref{Thm_basis_smth_W3.1}
or \ref{Thm_basis_smth_W3.3}, $G\in C_{BP}^{\left(  \left\lfloor
2\kappa\right\rfloor \right)  }\subset C_{BP}^{\left(  \left\lfloor
\min2\kappa\right\rfloor \right)  }$ and on expanding $G$ about the origin
using the Taylor's series with integral remainder in the Appendix
\ref{Sect_taylor_expansion} we get%
\begin{equation}
G\left(  y-x\right)  =\sum_{\left\vert \beta\right\vert <2\eta}\frac{\left(
y-x\right)  ^{\beta}}{\beta!}\left(  D^{\beta}G\right)  \left(  0\right)
+\mathcal{R}_{2\eta}\left(  0,y-x\right)  .\label{p925}%
\end{equation}
To calculate $r_{x}\left(  y\right)  $ apply the operator $\mathcal{Q}%
_{y}\mathcal{Q}_{x}$ to \ref{p925} and then, noting that $2\eta\leq2\theta$,
use Theorem \ref{Thm_QyQx(y-x)^a_=_0} to eliminate the power series and obtain

$r_{x}\left(  y\right)  =\left(  2\pi\right)  ^{-\frac{d}{2}}\mathcal{Q}%
_{y}\mathcal{Q}_{x}\mathcal{R}_{2\eta}\left(  0,y-x\right)  $. Expanding
$\mathcal{Q}_{y}$ and $\mathcal{Q}_{x}$ using $\mathcal{P}_{y}$ and
$\mathcal{P}_{x}$ now gives%
\begin{align*}
\left(  2\pi\right)  ^{\frac{d}{2}}r_{x}\left(  y\right)   & =\mathcal{R}%
_{2\eta}\left(  0,y-x\right)  -\mathcal{P}_{x}\left(  \mathcal{R}_{2\eta
}\left(  0,y-x\right)  \right)  -\mathcal{P}_{y}\left(  \mathcal{R}_{2\eta
}\left(  0,y-x\right)  \right)  +\mathcal{P}_{y}\mathcal{P}_{x}\left(
\mathcal{R}_{2\eta}\left(  0,y-x\right)  \right) \\
& =\mathcal{R}_{2\eta}\left(  0,y-x\right)  -\sum_{j=1}^{M}\mathcal{R}_{2\eta
}\left(  0,y-a^{\left(  j\right)  }\right)  l_{j}\left(  x\right)  -\sum
_{k=1}^{M}\mathcal{R}_{2\eta}\left(  0,a^{\left(  k\right)  }-x\right)
l_{k}\left(  y\right)  +\\
& +\sum_{j,k=1}^{M}\mathcal{R}_{2\eta}\left(  0,a^{\left(  k\right)
}-a^{\left(  j\right)  }\right)  l_{j}\left(  x\right)  l_{k}\left(  y\right)
,
\end{align*}

so that%
\begin{align*}
\left(  2\pi\right)  ^{\frac{d}{2}}r_{x}\left(  x\right)   & \leq\sum
_{j=1}^{M}\left\vert \mathcal{R}_{2\eta}\left(  0,x-a^{\left(  j\right)
}\right)  \right\vert \text{ }\left\vert l_{j}\left(  x\right)  \right\vert
+\sum_{k=1}^{M}\left\vert \mathcal{R}_{2\eta}\left(  0,a^{\left(  k\right)
}-x\right)  \right\vert \text{ }\left\vert l_{k}\left(  x\right)  \right\vert
+\\
& +\sum_{j,k=1}^{M}\left\vert \mathcal{R}_{2\eta}\left(  0,a^{\left(
k\right)  }-a^{\left(  j\right)  }\right)  \right\vert \text{ }\left\vert
l_{j}\left(  x\right)  \right\vert \text{ }\left\vert l_{k}\left(  x\right)
\right\vert \\
& \leq2\left(  \sum_{k=1}^{M}\left\vert l_{k}\left(  x\right)  \right\vert
\right)  \max_{k}\left\vert \mathcal{R}_{2\eta}\left(  0,x-a^{\left(
k\right)  }\right)  \right\vert +\left(  \sum_{k=1}^{M}\left\vert l_{k}\left(
x\right)  \right\vert \right)  ^{2}\max_{j,k}\left\vert \mathcal{R}_{2\eta
}\left(  0,a^{\left(  k\right)  }-a^{\left(  j\right)  }\right)  \right\vert .
\end{align*}

The two remainder estimates \ref{p34} and \ref{Ap144} for the Taylor series in
Appendix \ref{Sect_taylor_expansion}, namely:%
\[
\left\vert \mathcal{R}_{2\eta}\left(  0,b\right)  \right\vert \leq\left\{
\begin{array}
[c]{l}%
\frac{d^{\eta}\left\vert b\right\vert ^{2\eta}}{\left(  2\eta\right)  !}%
\max\limits_{\substack{\left\vert \beta\right\vert =2\eta\\s\in\left[
0,1\right]  }}\left\vert \left(  D^{\beta}G\right)  \left(  sb\right)
\right\vert ,\\
\medskip\\
\frac{\left\vert b\right\vert ^{2\eta}}{\left(  2\eta\right)  !}%
\max\limits_{s\in\left[  0,1\right]  }\left\vert \left(  \left(
\widehat{\cdot}D\right)  ^{2\eta}G\right)  \left(  sb\right)  \right\vert ,
\end{array}
\right.
\]

The first remainder estimate yields%
\begin{align*}
\left(  2\pi\right)  ^{\frac{d}{2}}r_{x}\left(  x\right)   & \leq\frac
{d^{\eta/2}}{\left(  2\eta\right)  !}2\left(  \sum_{k=1}^{M}\left\vert
l_{k}\left(  x\right)  \right\vert \right)  \left(  \max_{k}\left\vert
x-a^{\left(  k\right)  }\right\vert \right)  ^{2\eta}\max
\limits_{\substack{\left\vert \beta\right\vert =2\eta\\\left\vert y\right\vert
\leq\operatorname{diam}A_{x}}}\left\vert D^{\beta}G\left(  y\right)
\right\vert +\\
& \qquad+\frac{d^{\eta/2}}{\left(  2\eta\right)  !}\left(  \sum_{k=1}%
^{M}\left\vert l_{k}\left(  x\right)  \right\vert \right)  ^{2}\left(
\max_{j,k}\left\vert a^{\left(  j\right)  }-a^{\left(  k\right)  }\right\vert
\right)  ^{2\eta}\max\limits_{\substack{\left\vert \beta\right\vert
=2\eta\\\left\vert y\right\vert \leq\operatorname{diam}A_{x}}}\left\vert
D^{\beta}G\left(  y\right)  \right\vert \\
& <\frac{d^{\eta/2}}{\left(  2\eta\right)  !}\left(  1+\sum_{k=1}%
^{M}\left\vert l_{k}\left(  x\right)  \right\vert \right)  ^{2}\left(
\operatorname*{diam}A_{x}\right)  ^{2\eta}\max\limits_{\substack{\left\vert
\beta\right\vert =2\eta\\\left\vert y\right\vert \leq\operatorname{diam}A_{x}%
}}\left\vert D^{\beta}G\left(  y\right)  \right\vert ,
\end{align*}

so that%
\[
\sqrt{r_{x}\left(  x\right)  }\leq\frac{d^{\eta/2}}{\left(  2\pi\right)
^{\frac{d}{4}}}\frac{\left(  \operatorname*{diam}A_{x}\right)  ^{\eta}}%
{\sqrt{\left(  2\eta\right)  !}}\left(  1+\sum_{k=1}^{M}\left\vert
l_{k}\left(  x\right)  \right\vert \right)  \max\limits_{\substack{\left\vert
\beta\right\vert =2\eta\\\left\vert y\right\vert \leq\operatorname{diam}A_{x}%
}}\sqrt{\left\vert D^{\beta}G\left(  y\right)  \right\vert }.
\]

The second remainder estimate yields%
\begin{align*}
\left(  2\pi\right)  ^{\frac{d}{2}}r_{x}\left(  x\right)   & \leq2\left(
\sum_{k=1}^{M}\left\vert l_{k}\left(  x\right)  \right\vert \right)  \max
_{k}\left\vert \mathcal{R}_{2\eta}\left(  0,x-a^{\left(  k\right)  }\right)
\right\vert +\left(  \sum_{k=1}^{M}\left\vert l_{k}\left(  x\right)
\right\vert \right)  ^{2}\max_{j,k}\left\vert \mathcal{R}_{2\eta}\left(
0,a^{\left(  k\right)  }-a^{\left(  j\right)  }\right)  \right\vert \\
& \leq2\left(  \sum_{k=1}^{M}\left\vert l_{k}\left(  x\right)  \right\vert
\right)  \max_{k}\frac{\left\vert x-a^{\left(  k\right)  }\right\vert ^{2\eta
}}{\left(  2\eta\right)  !}\max\limits_{s\in\left[  0,1\right]  }\left\vert
\left(  \left(  \widehat{\cdot}D\right)  ^{2\eta}G\right)  \left(  s\left(
x-a^{\left(  k\right)  }\right)  \right)  \right\vert +\\
& \quad+\left(  \sum_{k=1}^{M}\left\vert l_{k}\left(  x\right)  \right\vert
\right)  ^{2}\max_{j,k}\frac{\left\vert a^{\left(  k\right)  }-a^{\left(
j\right)  }\right\vert ^{2\eta}}{\left(  2\eta\right)  !}\max\limits_{s\in
\left[  0,1\right]  }\left\vert \left(  \left(  \widehat{\cdot}D\right)
^{2\eta}G\right)  \left(  s\left(  a^{\left(  k\right)  }-a^{\left(  j\right)
}\right)  \right)  \right\vert \\
& \leq2\left(  \sum_{k=1}^{M}\left\vert l_{k}\left(  x\right)  \right\vert
\right)  \max_{k}\frac{\left\vert x-a^{\left(  k\right)  }\right\vert ^{2\eta
}}{\left(  2\eta\right)  !}\max\limits_{\left\vert y\right\vert \leq
\operatorname{diam}A_{x}}\left\vert \left(  \left(  \widehat{\cdot}D\right)
^{2\eta}G\right)  \left(  y\right)  \right\vert +\\
& +\left(  \sum_{k=1}^{M}\left\vert l_{k}\left(  x\right)  \right\vert
\right)  ^{2}\max_{j,k}\frac{\left\vert a^{\left(  k\right)  }-a^{\left(
j\right)  }\right\vert ^{2\eta}}{\left(  2\eta\right)  !}\max
\limits_{\left\vert y\right\vert \leq\operatorname{diam}A_{x}}\left\vert
\left(  \left(  \widehat{\cdot}D\right)  ^{2\eta}G\right)  \left(  y\right)
\right\vert \\
& \leq2\left(  \sum_{k=1}^{M}\left\vert l_{k}\left(  x\right)  \right\vert
\right)  \frac{\left(  \operatorname{diam}A_{x}\right)  ^{2\eta}}{\left(
2\eta\right)  !}\max\limits_{\left\vert y\right\vert \leq\operatorname{diam}%
A_{x}}\left\vert \left(  \left(  \widehat{\cdot}D\right)  ^{2\eta}G\right)
\left(  y\right)  \right\vert +\\
& +\left(  \sum_{k=1}^{M}\left\vert l_{k}\left(  x\right)  \right\vert
\right)  ^{2}\frac{\left(  \operatorname{diam}A_{x}\right)  ^{2\eta}}{\left(
2\eta\right)  !}\max\limits_{\left\vert y\right\vert \leq\operatorname{diam}%
A_{x}}\left\vert \left(  \left(  \widehat{\cdot}D\right)  ^{2\eta}G\right)
\left(  y\right)  \right\vert ,
\end{align*}

so that%
\[
\sqrt{r_{x}\left(  x\right)  }\leq\frac{\left(  \operatorname*{diam}%
A_{x}\right)  ^{\eta}}{\left(  2\pi\right)  ^{\frac{d}{4}}\sqrt{\left(
2\eta\right)  !}}\left(  1+\sum_{k=1}^{M}\left\vert l_{k}\left(  x\right)
\right\vert \right)  \max\limits_{\left\vert y\right\vert \leq
\operatorname{diam}A_{x}}\sqrt{\left\vert \left(  \left(  \widehat{\cdot
}D\right)  ^{2\eta}G\right)  \left(  y\right)  \right\vert }.
\]

\end{proof}

\subsection{The order of pointwise convergence of the interpolant to its data
function}

The result of the previous subsection now allow us to prove the pointwise
convergence of the minimal interpolant to its data function and to derive an
order of convergence. The concept of interpolation error and convergence go
hand in hand. The interpolation error is simply the difference between the
interpolant and the data function at a given point.

Before deriving the interpolation error estimate I will show that there exists
a nested sequence of independent data sets with sparsity tending to zero.

\begin{theorem}
\label{Thm_seq_data_regions_2}Suppose $\Omega$ is a bounded, open set
containing all the independent data sets. Then there exists a sequence of
independent data sets $X^{\left(  k\right)  }\subset\Omega$ such that
$X^{\left(  k\right)  }\subset X^{\left(  k+1\right)  }$ and $h_{X^{\left(
k\right)  }}\rightarrow0$ as $k\rightarrow\infty$.
\end{theorem}

\begin{proof}
For $k=1,2,3,\ldots$ there exists a finite covering of $\Omega$ by the balls

$\left\{  B\left(  a_{k}^{\left(  j\right)  };\frac{1}{k}\right)  \right\}
_{j=1}^{M_{k}}$. Construct $X^{\left(  1\right)  }$ by choosing points from
$\Omega$ so that one point lies in each ball $B\left(  a_{k}^{\left(
j\right)  };1\right)  $. Construct $X^{\left(  k+1\right)  }$ by first
choosing the points $X^{\left(  k\right)  }$ and then at least one extra point
so that $X^{\left(  k+1\right)  }$ contains points from each ball $B\left(
a_{k+1}^{\left(  j\right)  };\frac{1}{k+1}\right)  $.

Then $x\in\Omega\cap X^{\left(  k\right)  }$ implies $x\in B\left(
a_{k}^{\left(  j\right)  };\frac{1}{k}\right)  $ for some $j$ and hence
$\operatorname*{dist}\left(  x,X^{\left(  k\right)  }\right)  <\frac{1}{k}$.

Hence $h_{X^{\left(  k\right)  }}=\sup\limits_{x\in\Omega}\operatorname*{dist}%
\left(  x,X^{\left(  k\right)  }\right)  <\frac{1}{k}$ and $\lim
\limits_{k\rightarrow\infty}h_{X^{\left(  k\right)  }}=0$.
\end{proof}

The next theorem derives an upper bound for the uniform pointwise error of the
minimal interpolant in terms of the product of the seminorm of the data
function, a power of the sparsity of the independent data and various
constants derived from the theory of Lagrange interpolation.

The order of uniform pointwise convergence is defined to be the power of the
sparsity parameter, in this case $\min\left\{  \theta,\frac{1}{2}\left\lfloor
\min2\kappa\right\rfloor \right\}  $ where $\theta$ and $\kappa$ are the
weight function parameters.

\begin{theorem}
\label{Thm_interpol_converg}Let $w$ be a weight function with properties W2
and W3 for order $\theta$ and smoothness parameter $\kappa$. Set $\eta
=\min\left\{  \theta,\frac{1}{2}\left\lfloor 2\underline{\kappa}\right\rfloor
\right\}  $.

Suppose the notation and assumptions of Lemma \ref{Lem_int_Lagrange_interpol}
hold with the data point sparsity of $X$ in $\Omega$ defined by $h_{X}%
=\sup\limits_{\omega\in\Omega}\operatorname*{dist}\left(  \omega;X\right)  $.
Suppose also that $\mathcal{I}_{X}$ is the minimal interpolant on $X$ of the
data function $f_{d}\in X_{w}^{\theta}$.

Then there exist positive constants $c_{\Omega,\theta},h_{\Omega,\theta
},K_{\Omega,\theta}^{\prime}$ and

$k_{\Omega,\theta,\eta}=\frac{d^{\eta/2}}{\left(  2\pi\right)  ^{d/4}%
\sqrt{\left(  2\eta\right)  !}}\left(  1+K_{\Omega,\theta}^{\prime}\right)  $
such that
\begin{equation}
\left\vert f_{d}\left(  x\right)  -\mathcal{I}_{X}f_{d}\left(  x\right)
\right\vert \leq\left\vert f_{d}-\mathcal{I}_{X}f_{d}\right\vert _{w,\theta
}k_{\Omega,\theta,\eta}\left(  c_{\Omega,\theta}h_{X}\right)  ^{\eta}%
\max\limits_{\substack{\left\vert \beta\right\vert =2\eta\\\left\vert
y\right\vert \leq c_{\Omega,\theta}h_{X}}}\left\vert D^{\beta}G(y)\right\vert
,\quad x\in\overline{\Omega},\label{q52}%
\end{equation}

when $0<h_{X}\leq h_{\Omega,\theta}$. The Lagrange interpolation constants
$c_{\Omega,\theta},K_{\Omega,\theta}^{\prime}$ and $h_{\Omega,\theta}$ are
derived in Lemma \ref{Lem_int_Lagrange_interpol}.

Further, $\left\vert f_{d}-\mathcal{I}_{X}f_{d}\right\vert _{w,\theta}%
\leq\left\vert f_{d}\right\vert _{w,\theta}$ and the order of convergence is
$\eta$.

Finally we have the bound%
\[
\left\vert f_{d}\left(  x\right)  -\mathcal{I}_{X}f_{d}\left(  x\right)
\right\vert \leq\left\vert f_{d}-\mathcal{I}_{X}f_{d}\right\vert _{w,\theta
}k_{\Omega,\theta,\eta}\left(  \operatorname{diam}\Omega\right)  ^{\eta}%
\max\limits_{\substack{\left\vert \beta\right\vert =2\eta\\\left\vert
y\right\vert \leq\operatorname{diam}\Omega}}\left\vert D^{\beta}G\left(
y\right)  \right\vert ,\quad x\in\overline{\Omega}.
\]

\end{theorem}

\begin{proof}
Fix $x\in\Omega$. Since $\mathcal{I}_{X}$ interpolates the data from the
function $f_{d}$ it follows that $\mathcal{I}_{X}f\left(  x\right)  =f\left(
x\right)  $ when $x\in X$. Thus $\mathcal{P}\left(  f_{d}-\mathcal{I}_{X}%
f_{d}\right)  =0$ and hence $f_{d}-\mathcal{I}_{X}f_{d}=\mathcal{Q}\left(
f_{d}-\mathcal{I}_{X}f_{d}\right)  $.%
\[
f_{d}\left(  x\right)  -\mathcal{I}_{X}f_{d}\left(  x\right)  =\left(
f_{d}-\mathcal{I}_{X}f_{d},R_{x}\right)  _{w,\theta}=\left(  \mathcal{Q}%
\left(  f_{d}-\mathcal{I}_{X}f_{d}\right)  ,\mathcal{Q}R_{x}\right)
_{w,\theta}=\left\langle f_{d}-\mathcal{I}_{X}f_{d},r_{x}\right\rangle
_{w,\theta}.
\]

But by part 4 of Theorem \ref{Thm_rx(y)_properties}, $\left\vert
r_{x}\right\vert _{w,\theta}=\sqrt{r_{x}\left(  x\right)  }$ and thus
\begin{equation}
\left\vert f_{d}\left(  x\right)  -\mathcal{I}_{X}f_{d}\left(  x\right)
\right\vert \leq\left\vert f_{d}-\mathcal{I}_{X}f_{d}\right\vert _{w,\theta
}\left\vert r_{x}\right\vert _{w,\theta}=\left\vert f_{d}-\mathcal{I}_{X}%
f_{d}\right\vert _{w,\theta}\sqrt{r_{x}\left(  x\right)  }.\label{p9102}%
\end{equation}

The upper bound for $\sqrt{r_{x}\left(  x\right)  }$ given by\ Theorem
\ref{Thm_int_rx(x)_bound} is%
\[
\sqrt{r_{x}\left(  x\right)  }\leq\frac{1}{\left(  2\pi\right)  ^{\frac{d}{4}%
}}\frac{d^{\eta/2}}{\sqrt{\left(  2\eta\right)  !}}\left(  1+\sum_{k=1}%
^{M}\left\vert l_{k}\left(  x\right)  \right\vert \right)  \left(
\operatorname*{diam}A_{x}\right)  ^{\eta}\max\limits_{\substack{\left\vert
\beta\right\vert =2\eta\\\left\vert y\right\vert \leq\operatorname{diam}A_{x}%
}}\left\vert D^{\beta}G(y)\right\vert ,
\]

so that by Lemma \ref{Lem_int_Lagrange_interpol}%
\begin{equation}
\sqrt{r_{x}\left(  x\right)  }\leq\frac{d^{\eta/2}}{\left(  2\pi\right)
^{\frac{d}{4}}\sqrt{\left(  2\eta\right)  !}}\left(  1+K_{\Omega,\theta
}^{\prime}\right)  \left(  c_{\Omega,\theta}h_{X}\right)  ^{\eta}%
\max\limits_{\substack{\left\vert \beta\right\vert =2\eta\\\left\vert
y\right\vert \leq c_{\Omega,\theta}h_{X}}}\left\vert (D^{\beta}%
G)(y)\right\vert ,\label{q83}%
\end{equation}

when $h_{X}\leq h_{\Omega,\theta}$. Therefore%
\begin{align*}
\left\vert f_{d}\left(  x\right)  -\mathcal{I}_{X}f_{d}\left(  x\right)
\right\vert  & \leq\left\vert f_{d}-\mathcal{I}_{X}f_{d}\right\vert
_{w,\theta}\frac{d^{\frac{\eta}{2}}\left(  1+K_{\Omega,\theta}^{\prime
}\right)  }{\left(  2\pi\right)  ^{\frac{d}{4}}\sqrt{\left(  2\eta\right)  !}%
}\left(  c_{\Omega,\theta}h_{X}\right)  ^{\eta}\max
\limits_{\substack{\left\vert \beta\right\vert =2\eta\\\left\vert y\right\vert
\leq c_{\Omega,\theta}h_{X,}}}\left\vert (D^{\beta}G)(y)\right\vert \\
& \leq\left\vert f_{d}-\mathcal{I}_{X}f_{d}\right\vert _{w,\theta}%
k_{\Omega,\theta,\eta}\left(  c_{\Omega,\theta}h_{X}\right)  ^{\eta}%
\max\limits_{\substack{\left\vert \beta\right\vert =2\eta\\\left\vert
y\right\vert \leq c_{\Omega,\theta}h_{X}}}\left\vert (D^{\beta}%
G)(y)\right\vert ,\quad x\in\Omega,
\end{align*}

when $h_{X}\leq h_{\Omega,\theta}$. The last inequality can be extended
to\ $\overline{\Omega}$ because the inequality is uniform on $\Omega$,
$\Omega$ is bounded and $f_{d}-\mathcal{I}_{X}f_{d}$ is continuous on
$\mathbb{R}^{d}$. That $\left\vert f_{d}-\mathcal{I}_{X}f_{d}\right\vert
_{w,\theta}\leq\left\vert f_{d}\right\vert _{w,\theta}$ was shown in part 2 of
Theorem \ref{Thm_propert_interpol_map}.

The final bound is true since $\operatorname*{diam}A_{x}\leq
\operatorname*{diam}\Omega$.
\end{proof}

\begin{remark}
\label{Lem_Thm_interpol_converg}If $\kappa_{i}<1/2$ for some $i$ then
$\left\lfloor \min2\kappa\right\rfloor =0$ and so $\eta=0$, and the last
theorem only shows the interpolation error is bounded i.e. zero order of
convergence. \textbf{In the next section we show how to overcome this
limitation} and obtain positive convergence orders for these cases and
improved orders of convergence otherwise. In the process we will demonstrate
improved order of convergence estimates for almost all thin-plate splines and
shifted thin-plate splines.
\end{remark}

\section{Slightly improved convergence results\label{Sect_better_results}}

??

\subsection{Results using basis function Taylor series}

Noting Remark \ref{Lem_Thm_interpol_converg} to the last lemma in the previous
section, we will show how to obtain improved interpolant convergence estimates
using a Taylor series expansion result for distributions. Our results will be
illustrated using the thin-plate splines. The improvement in convergence order
for these examples is at most $1/2$.

To prove Lemma \ref{Lem_Taylor_extension} we will require the following result:

\begin{lemma}
\label{Lem_Lem_Taylor_extension}Suppose $u\in L_{loc}^{1}\left(
\mathbb{R}^{d}\right)  $ and $\phi\in C_{0}^{\infty}$. Then for all
$b\in\mathbb{R}^{d}$,%
\begin{equation}
\int_{0}^{1}\int\left\vert u\left(  x+tb\right)  \right\vert \text{
}\left\vert \phi\left(  x\right)  \right\vert dxdt<\infty.\label{p9108}%
\end{equation}

\end{lemma}

\begin{proof}
Suppose $\operatorname*{supp}\phi\subset B\left(  0;R\right)  $. Then
\begin{align*}
\int\limits_{0}^{1}\int\left\vert u\left(  x+tb\right)  \right\vert \text{
}\left\vert \phi\left(  x\right)  \right\vert dxdt=\int\limits_{0}^{1}%
\int\left\vert u\left(  y\right)  \right\vert \text{ }\left\vert \phi\left(
y-tb\right)  \right\vert dydt &  =\int\limits_{0}^{1}\int\limits_{\left\vert
y-tb\right\vert \leq R}\left\vert u\left(  y\right)  \right\vert \text{
}\left\vert \phi\left(  y-tb\right)  \right\vert dydt\\
&  \leq\int\limits_{0}^{1}\int\limits_{\left\vert y\right\vert \leq
R+\left\vert b\right\vert }\left\vert u\left(  y\right)  \right\vert \text{
}\left\vert \phi\left(  y-tb\right)  \right\vert dydt\\
&  \leq\left\Vert \phi\right\Vert _{\infty}\int\limits_{0}^{1}\int%
\limits_{\left\vert y\right\vert \leq R+\left\vert b\right\vert }\left\vert
u\left(  y\right)  \right\vert dydt\\
&  \leq\left\Vert \phi\right\Vert _{\infty}\int\limits_{\left\vert
y\right\vert \leq R+\left\vert b\right\vert }\left\vert u\left(  y\right)
\right\vert dy\\
&  <\infty.
\end{align*}

\end{proof}

Our Taylor series distribution result is:

\begin{lemma}
\label{Lem_Taylor_extension}This lemma uses the results and notation of Lemmas
\ref{Lem_deriv_rad_funcs} and \ref{Lem_op_aD_estim} which are in the Appendix.

Suppose $u\in C^{\left(  n-1\right)  }\left(  \mathbb{R}^{d}\right)  $ and the
distributional derivatives $\left\{  D^{\beta}u\right\}  _{\left\vert
\beta\right\vert =n}$ are $L_{loc}^{1}$ functions. Suppose also that for each
fixed $b\neq0$ the integrals%
\begin{equation}
\int_{0}^{1}\left(  1-t\right)  ^{n-1}\left\vert (D^{\beta}u)(z+tb)\right\vert
dt,\quad\left\vert \beta\right\vert =n;\text{ }z,b\in\mathbb{R}^{d}%
,\label{q07}%
\end{equation}

have polynomial growth in $z$. Then%
\begin{equation}
u(z+b)=\sum_{\left\vert \beta\right\vert <n}\frac{b^{\beta}}{\beta!}\left(
D^{\beta}u\right)  (z)+\left(  \mathcal{R}_{n}u\right)  \left(  z,b\right)
,\text{\quad}z,b\in\mathbb{R}^{d},\label{p91.87}%
\end{equation}

where $\mathcal{R}_{n}u$ is the integral remainder term
\begin{align}
\left(  \mathcal{R}_{n}u\right)  \left(  z,b\right)   & =n\sum_{\left\vert
\beta\right\vert =n}\frac{b^{\beta}}{\beta!}\int_{0}^{1}\left(  1-t\right)
^{n-1}(D^{\beta}u)(z+tb)dt\label{p915.47}\\
& =\frac{1}{\left(  n-1\right)  !}\int_{0}^{1}\left(  1-t\right)
^{n-1}(\left(  bD\right)  ^{n}u)(z+tb)dt.\label{p913}%
\end{align}

In particular, we have%
\begin{align}
\left(  \mathcal{R}_{n}u\right)  \left(  0,b\right)   & =\frac{1}{\left(
n-1\right)  !}\int_{0}^{1}\left(  1-t\right)  ^{n-1}\left(  \left(  \cdot
D\right)  ^{n}u\right)  \left(  tb\right)  dt\label{q11}\\
& =\frac{\left\vert b\right\vert ^{n}}{\left(  n-1\right)  !}\int_{0}%
^{1}\left(  1-t\right)  ^{n-1}\left(  \left(  \widehat{\cdot}D\right)
^{n}u\right)  \left(  tb\right)  dt,\label{q10}%
\end{align}

and the estimate%
\begin{equation}
\left\vert \left(  \mathcal{R}_{n}u\right)  \left(  0,b\right)  \right\vert
\leq\frac{\left\vert b\right\vert ^{n}}{n!}\max_{t\in\left[  0,1\right]
}\left\vert \left(  \widehat{\cdot}D\right)  ^{n}u\left(  tb\right)
\right\vert .\label{p915}%
\end{equation}

\end{lemma}

\begin{proof}
In order to overcome the fact that $D^{\beta}u$ may not be $C^{\left(
0\right)  }\left(  \mathbb{R}^{d}\right)  $ when $\left\vert \beta\right\vert
=n$, we will use a Taylor series expansion with remainder for distributions.
Suppose $\phi\in C_{0}^{\infty}$. Then the conditions on $u$ allow us to use
Lemma \ref{Lem_Lem_Taylor_extension} to show that the iterated integrals
\ref{q56} and \ref{q79} are absolutely convergent and thus apply Fubini's
theorem to swap the orders of integration in the following calculations: using
remainder form \ref{p88} we proceed as follows:%
\begin{align}
\left[  u\left(  z+b\right)  ,\phi\left(  z\right)  \right]   & =\left[
u\left(  z\right)  ,\phi\left(  z-b\right)  \right] \nonumber\\
& =\left[  u\left(  z\right)  ,\sum_{\left\vert \beta\right\vert \leq n}%
\frac{\left(  -b\right)  ^{\beta}}{\beta!}D^{\beta}\phi(z)\right]
+\nonumber\\
& \qquad+\left[  u\left(  z\right)  ,n\sum_{\left\vert \beta\right\vert
=n}\frac{\left(  -b\right)  ^{\beta}}{\beta!}\int_{0}^{1}\left(  1-t\right)
^{n-1}\left(  D^{\beta}\phi\right)  \left(  z+tb\right)  dt\right] \nonumber\\
& =\sum_{\left\vert \beta\right\vert <n}\frac{b^{\beta}}{\beta!}\left[
D^{\beta}u,\phi\right]  +\nonumber\\
& \qquad+n\sum_{\left\vert \beta\right\vert =n}\frac{\left(  -b\right)
^{\beta}}{\beta!}\left[  u,\int_{0}^{1}\left(  1-t\right)  ^{n-1}\left(
D^{\beta}\phi\right)  \left(  \cdot-tb\right)  dt\right]  .\label{p916}%
\end{align}

We now analyze the integral remainder term of \ref{p916}:
\begin{align}
&  \left[  u,\int_{0}^{1}\left(  1-t\right)  ^{n-1}\left(  D^{\beta}%
\phi\right)  \left(  \cdot-tb\right)  dt\right] \nonumber\\
&  =\int u\left(  z\right)  \int_{0}^{1}\left(  1-t\right)  ^{n-1}\left(
D^{\beta}\phi\right)  \left(  z-tb\right)  dt\text{ }dz\nonumber\\
&  =\int_{0}^{1}\left(  1-t\right)  ^{n-1}\int u\left(  z\right)  \left(
D^{\beta}\phi\right)  \left(  z-tb\right)  dz\text{ }dt\label{q56}\\
&  =\int_{0}^{1}\left(  1-t\right)  ^{n-1}\left[  u\left(  z\right)  ,\left(
D^{\beta}\phi\right)  \left(  z-tb\right)  \right]  \text{ }dt\nonumber\\
&  =\int_{0}^{1}\left(  1-t\right)  ^{n-1}\left[  u\left(  z+tb\right)
,\left(  D^{\beta}\phi\right)  \left(  z\right)  \right]  \text{
}dt\nonumber\\
&  =\left(  -1\right)  ^{\left\vert \beta\right\vert }\int_{0}^{1}\left(
1-t\right)  ^{n-1}\left[  \left(  D^{\beta}u\right)  \left(  z+tb\right)
,\phi\left(  z\right)  \right]  \text{ }dt\nonumber\\
&  =\left(  -1\right)  ^{\left\vert \beta\right\vert }\int_{0}^{1}\left(
1-t\right)  ^{n-1}\int\left(  D^{\beta}u\right)  \left(  z+tb\right)
\phi\left(  z\right)  dz\text{ }dt\nonumber\\
&  =\left(  -1\right)  ^{\left\vert \beta\right\vert }\int\int_{0}^{1}\left(
1-t\right)  ^{n-1}\left(  D^{\beta}u\right)  \left(  z+tb\right)  dt\text{
}\phi\left(  z\right)  dz\label{q79}\\
&  =\left(  -1\right)  ^{\left\vert \beta\right\vert }\left[  \int_{0}%
^{1}\left(  1-t\right)  ^{n-1}\left(  D^{\beta}u\right)  \left(  z+tb\right)
dt,\phi\left(  z\right)  \right]  ,\nonumber
\end{align}

so \ref{p916} now becomes%
\begin{align*}
\left[  u\left(  z+b\right)  ,\phi\left(  z\right)  \right]   & =\sum
_{\left\vert \beta\right\vert <n}\frac{b^{\beta}}{\beta!}\left[  D^{\beta
}u\left(  z\right)  ,\phi\left(  z\right)  \right]  +\\
& \qquad\qquad+n\sum_{\left\vert \beta\right\vert =n}\frac{b^{\beta}}{\beta
!}\left[  \int_{0}^{1}\left(  1-t\right)  ^{n-1}\left(  D^{\beta}u\right)
\left(  z+tb\right)  dt,\phi\left(  z\right)  \right]  ,
\end{align*}

for $\phi\in C_{0}^{\infty}$. Thus%
\begin{align*}
u\left(  z+b\right)   & =\sum_{\left\vert \beta\right\vert <n}\frac{b^{\beta}%
}{\beta!}\left(  D^{\beta}u\right)  \left(  z\right)  +n\sum_{\left\vert
\beta\right\vert =n}\frac{b^{\beta}}{\beta!}\int_{0}^{1}\left(  1-t\right)
^{n-1}\left(  D^{\beta}u\right)  \left(  z+tb\right)  dt\\
& =\sum_{\left\vert \beta\right\vert <n}\frac{b^{\beta}}{\beta!}\left(
D^{\beta}u\right)  \left(  z\right)  +n\int_{0}^{1}\left(  1-t\right)
^{n-1}\left(  \sum_{\left\vert \beta\right\vert =n}\frac{b^{\beta}D^{\beta}%
}{\beta!}u\right)  \left(  z+tb\right)  dt\\
& =\sum_{\left\vert \beta\right\vert <n}\frac{b^{\beta}}{\beta!}\left(
D^{\beta}u\right)  \left(  z\right)  +n\int_{0}^{1}\left(  1-t\right)
^{n-1}\left(  \frac{1}{n!}\left(  bD\right)  ^{n}u\right)  \left(
z+tb\right)  dt\\
& =\sum_{\left\vert \beta\right\vert <n}\frac{b^{\beta}}{\beta!}\left(
D^{\beta}u\right)  \left(  z\right)  +\frac{1}{\left(  n-1\right)  !}\int%
_{0}^{1}\left(  1-t\right)  ^{n-1}\left(  \left(  bD\right)  ^{n}u\right)
\left(  z+tb\right)  dt,
\end{align*}

as claimed. In particular, we have%
\begin{align*}
\left(  \mathcal{R}_{n}u\right)  \left(  0,b\right)   & =\frac{1}{\left(
n-1\right)  !}\int_{0}^{1}\left(  1-t\right)  ^{n-1}\left(  \left(  bD\right)
^{n}u\right)  \left(  tb\right)  dt\\
& =\frac{\left\vert b\right\vert ^{n}}{\left(  n-1\right)  !}\int_{0}%
^{1}\left(  1-t\right)  ^{n-1}\left(  \left(  \widehat{b}D\right)
^{n}u\right)  \left(  tb\right)  dt\\
& =\frac{\left\vert b\right\vert ^{n}}{\left(  n-1\right)  !}\int_{0}%
^{1}\left(  1-t\right)  ^{n-1}\left(  \left(  \widehat{tb}D\right)
^{n}u\right)  \left(  tb\right)  dt\\
& =\frac{\left\vert b\right\vert ^{n}}{\left(  n-1\right)  !}\int_{0}%
^{1}\left(  1-t\right)  ^{n-1}\left(  \left(  \widehat{\cdot}D\right)
^{n}u\right)  \left(  tb\right)  dt,
\end{align*}

and the estimate%
\[
\left\vert \left(  \mathcal{R}_{n}u\right)  \left(  0,b\right)  \right\vert
\leq\frac{\left\vert b\right\vert ^{n}}{n!}\max_{x\in\left[  0,b\right]
}\left\vert \left(  \widehat{\cdot}D\right)  ^{n}u\left(  x\right)
\right\vert .
\]

\end{proof}

We now consider the radial and homogeneous cases:

\begin{corollary}
\label{Cor_Lem_Taylor_extension}Two special cases of the remainder formula
\ref{q10} are:

\begin{enumerate}
\item If $u$ is \textbf{radial}, say $u\left(  x\right)  =u_{\circ}\left(
\left\vert x\right\vert \right)  $, then from Lemma \ref{Lem_deriv_rad_funcs},%
\[
\left(  \mathcal{R}_{n}u\right)  \left(  0,b\right)  =\frac{\left\vert
b\right\vert ^{n}}{\left(  n-1\right)  !}\int_{0}^{1}\left(  1-t\right)
^{n-1}\left(  D^{n}u_{\circ}\right)  \left(  t\left\vert b\right\vert \right)
dt,
\]

and%
\[
\left\vert \left(  \mathcal{R}_{n}u\right)  \left(  0,b\right)  \right\vert
\leq\frac{\left\vert b\right\vert ^{n}}{n!}\max_{t\in\left[  0,\left\vert
b\right\vert \right]  }\left\vert D^{n}u_{\circ}\left(  t\right)  \right\vert
.
\]

\item If $u\in C^{\left(  n\right)  }$ is homogeneous of order $s$ then
\[
\left(  \widehat{\cdot}D\right)  ^{n}u\left(  x\right)  =\left\vert
x\right\vert ^{s-n}\left(  \left(  \widehat{\cdot}D\right)  ^{n}u\right)
\left(  \widehat{x}\right)  ,
\]

and%
\[
\left(  \left(  \widehat{\cdot}D\right)  ^{n}u\right)  \left(  0\right)
=0,\text{\quad}n<s.
\]

\end{enumerate}
\end{corollary}

\begin{proof}
\textbf{Part 1} From Lemma \ref{Lem_deriv_rad_funcs}, $\left(  \left(
\widehat{\cdot}D\right)  ^{n}u\right)  \left(  x\right)  =\left(
D^{n}u_{\circ}\right)  \left(  \left\vert x\right\vert \right)  $.\smallskip

\textbf{Part 2} Since $u\left(  tx\right)  =t^{s}u\left(  x\right)  $ we have
$D^{\alpha}u\left(  tx\right)  =t^{\left\vert \alpha\right\vert }\left(
D^{\alpha}u\right)  \left(  tx\right)  =t^{s}D^{\alpha}u\left(  x\right)  $
and setting $t=\left\vert x\right\vert ^{-1}$ we get $D^{\alpha}u\left(
x\right)  =\left\vert x\right\vert ^{s-\left\vert \alpha\right\vert }\left(
D^{\alpha}u\right)  \left(  \widehat{x}\right)  $.

Thus
\begin{align*}
\frac{1}{k!}\left(  \widehat{\cdot}D\right)  ^{k}u\left(  x\right)
=\sum\limits_{\left\vert \alpha\right\vert =k}\frac{\widehat{x}^{\alpha}%
}{\alpha!}D^{\alpha}u\left(  x\right)   & =\sum\limits_{\left\vert
\alpha\right\vert =k}\frac{\widehat{x}^{\alpha}}{\alpha!}\left\vert
x\right\vert ^{s-\left\vert \alpha\right\vert }\left(  D^{\alpha}u\right)
\left(  \widehat{x}\right) \\
& =\left\vert x\right\vert ^{s-kn}\sum\limits_{\left\vert \alpha\right\vert
=n}\frac{\widehat{x}^{\alpha}}{\alpha!}\left(  D^{\alpha}u\right)  \left(
\widehat{x}\right) \\
& =\frac{1}{k!}\left\vert x\right\vert ^{s-k}\left(  \left(  \widehat{\cdot
}D\right)  ^{k}u\right)  \left(  \widehat{x}\right)  ,
\end{align*}

and thus $\left(  \widehat{\cdot}D\right)  ^{n}u\left(  0\right)  =0$ when
$n<s$.
\end{proof}

Using the previous lemma we now modify Theorem \ref{Thm_int_rx(x)_bound} to
obtain the following estimate for $\sqrt{r_{x}\left(  x\right)  }$ which
requires the assumptions \ref{q761} and \ref{q081} to hold.

\begin{theorem}
\label{Thm_int_rx(x)_bnd_better_order}Suppose $w$ is a weight function with
properties W2 and W3 for order $\theta$ and $\kappa$, and such that $\frac
{1}{2}\left\lfloor 2\underline{\kappa}\right\rfloor <\theta$. Set $\eta
=\frac{1}{2}\left\lfloor 2\underline{\kappa}\right\rfloor $. Also suppose $G$
is a basis function of order $\theta$ so we know that $G\in C_{BP}^{\left(
2\eta\right)  }$.

Now suppose that the distributions $\left\{  D^{\beta}G\right\}  _{\left\vert
\beta\right\vert =2\eta+1}$ are $L_{loc}^{1}$ functions and that for each
fixed $b\neq0$ the integrals%
\begin{equation}
\int_{0}^{1}\left(  1-t\right)  ^{2\eta}\left\vert \left(  D^{\beta}G\right)
\left(  z+tb\right)  \right\vert dt,\quad z,b\in\mathbb{R}^{d},\text{
}\left\vert \beta\right\vert =2\eta+1,\label{q761}%
\end{equation}

have polynomial growth in $z$.

Further, suppose there exist constants $r_{G},c_{G,\eta}>0$ and $\delta_{G}>0$
such that%
\begin{equation}
\left\vert \int_{0}^{1}\left(  1-t\right)  ^{2\eta}\left(  \left(  \cdot
D\right)  ^{2\eta+1}G\right)  \left(  tb\right)  dt\right\vert \leq\left(
2\pi\right)  ^{\frac{d}{2}}\left(  2\eta\right)  !c_{G,\eta}\left\vert
b\right\vert ^{2\left(  \eta+\delta_{G}\right)  },\quad\left\vert b\right\vert
\leq r_{G},\label{q081}%
\end{equation}

where $\left(  \left(  \cdot D\right)  G\right)  \left(  x\right)  =\left(
x_{1}D_{1}+\ldots+x_{d}D_{d}\right)  G\left(  x\right)  $.

Regarding unisolvency, assume $A=\left\{  a^{\left(  k\right)  }\right\}
_{k=1}^{M}$ is a minimal $\theta$-unisolvent set and that $\left\{
l_{k}\right\}  _{k=1}^{M}$ is the corresponding unique cardinal basis for
$P_{\theta-1}$. Now construct $\mathcal{P},\mathcal{Q},R_{x}$ using $A$ and
$\left\{  l_{k}\right\}  _{k=1}^{M}$.

Then if $r_{x}=\mathcal{Q}R_{x}$ we have the estimate%
\begin{equation}
\sqrt{r_{x}\left(  x\right)  }\leq\sqrt{c_{G,\eta}}\left(  1+\sum_{k=1}%
^{M}\left\vert l_{k}\left(  x\right)  \right\vert \right)  \left(
\operatorname*{diam}A_{x}\right)  ^{\eta+\delta_{G}},\quad\operatorname*{diam}%
A_{x}\leq r_{G},\text{ }x\in\Omega,\label{q87}%
\end{equation}

where $A_{x}=A\cup\left\{  x\right\}  $.
\end{theorem}

\begin{proof}
The proof will follow that of Theorem \ref{Thm_int_rx(x)_bound}. From part 4
of Theorem \ref{Thm_rx(y)_properties}, $r_{x}\left(  y\right)  =\left\langle
r_{x},r_{y}\right\rangle _{w,\theta}$ so that $r_{x}\left(  x\right)
=\left\vert r_{x}\right\vert _{w,\theta}^{2}\geq0$. Also from Theorem
\ref{Thm_rx(y)_properties}, $r_{x}\left(  y\right)  =\left(  2\pi\right)
^{-\frac{d}{2}}\mathcal{Q}_{y}\mathcal{Q}_{x}G\left(  y-x\right)  $ when
$x\neq y$.

Since $G\in C_{BP}^{\left(  2\eta\right)  }$ and the functions $\left\{
D^{\beta}G\right\}  _{\left\vert \beta\right\vert =2\eta+1}$ are $L_{loc}^{1}%
$, the conditions \ref{q761} mean we can use Lemma \ref{Lem_Taylor_extension}
to expand $G$ about the origin and get%
\begin{equation}
G\left(  y-x\right)  =\sum_{\left\vert \beta\right\vert \leq2\eta}%
\frac{\left(  y-x\right)  ^{\beta}}{\beta!}D^{\beta}G\left(  0\right)
+\mathcal{R}_{2\eta+1}\left(  0,y-x\right)  .\label{p9271}%
\end{equation}

To calculate $r_{x}\left(  y\right)  $ apply the operator $\mathcal{Q}%
_{y}\mathcal{Q}_{x}$ to \ref{p9271} and then noting that $\eta\leq\theta$, use
Theorem \ref{Thm_QyQx(y-x)^a_=_0} to eliminate the power series and so obtain%
\[
r_{x}\left(  y\right)  =\left(  2\pi\right)  ^{-\frac{d}{2}}\mathcal{Q}%
_{y}\mathcal{Q}_{x}\mathcal{R}_{2\eta+1}\left(  0,y-x\right)  .
\]

Expanding $\mathcal{Q}_{y}$ and $\mathcal{Q}_{x}$ using $\mathcal{P}_{y}$ and
$\mathcal{P}_{x}$ now gives%
\begin{align*}
\left(  2\pi\right)  ^{\frac{d}{2}}r_{x}\left(  y\right)   & =\mathcal{R}%
_{2\eta+1}\left(  0,y-x\right)  -\mathcal{P}_{x}\left(  \mathcal{R}_{2\eta
+1}\right)  -\mathcal{P}_{y}\left(  \mathcal{R}_{2\eta+1}\left(  0,y-x\right)
\right)  +\\
& \qquad\qquad+\mathcal{P}_{y}\mathcal{P}_{x}\left(  \mathcal{R}_{2\eta
+1}\left(  0,y-\left(  0,y-x\right)  x\right)  \right) \\
& =\mathcal{R}_{2\eta+1}\left(  0,y-x\right)  -\sum_{j=1}^{M}\mathcal{R}%
_{2\eta+1}\left(  0,y-a^{\left(  j\right)  }\right)  l_{j}\left(  x\right)
-\\
& \qquad-\sum_{k=1}^{M}\mathcal{R}_{2\eta+1}\left(  0,a^{\left(  k\right)
}-x\right)  l_{k}\left(  y\right)  +\\
& \qquad+\sum_{j,k=1}^{M}\mathcal{R}_{2\eta+1}\left(  0,a^{\left(  k\right)
}-a^{\left(  j\right)  }\right)  l_{j}\left(  x\right)  l_{k}\left(  y\right)
,
\end{align*}

so that%
\begin{align*}
\left(  2\pi\right)  ^{\frac{d}{2}}r_{x}\left(  x\right)   & \leq\sum
_{j=1}^{M}\left\vert \mathcal{R}_{2\eta+1}\left(  0,x-a^{\left(  j\right)
}\right)  \right\vert \left\vert l_{j}\left(  x\right)  \right\vert
+\sum_{k=1}^{M}\left\vert \mathcal{R}_{2\eta+1}\left(  0,a^{\left(  k\right)
}-x\right)  \right\vert \left\vert l_{k}\left(  x\right)  \right\vert +\\
& \qquad\qquad+\sum_{j,k=1}^{M}\left\vert \mathcal{R}_{2\eta+1}\left(
0,a^{\left(  k\right)  }-a^{\left(  j\right)  }\right)  \right\vert \left\vert
l_{j}\left(  x\right)  \right\vert \left\vert l_{k}\left(  x\right)
\right\vert \\
& \leq2\left(  \sum_{k=1}^{M}\left\vert l_{k}\left(  x\right)  \right\vert
\right)  \max_{k}\left\vert \mathcal{R}_{2\eta+1}\left(  0,x-a^{\left(
k\right)  }\right)  \right\vert +\left(  \sum_{k=1}^{M}\left\vert l_{k}\left(
x\right)  \right\vert \right)  ^{2}\max_{j,k}\left\vert \mathcal{R}_{2\eta
+1}\left(  0,a^{\left(  k\right)  }-a^{\left(  j\right)  }\right)  \right\vert
.
\end{align*}

But from \ref{q11} of Lemma \ref{Lem_Taylor_extension},%
\[
\mathcal{R}_{2\eta+1}\left(  0,b\right)  =\frac{1}{\left(  2\eta\right)
!}\int_{0}^{1}\left(  1-t\right)  ^{2\eta}\left(  \left(  \cdot D\right)
^{2\eta+1}G\right)  \left(  tb\right)  dt,
\]

and estimate \ref{q081} holds i.e.
\[
\left\vert \int_{0}^{1}\left(  1-t\right)  ^{2\eta}\left(  \left(  \cdot
D\right)  ^{2\eta+1}G\right)  (tb)dt\right\vert \leq\left(  2\pi\right)
^{\frac{d}{2}}c_{G,\eta}\left(  2\eta\right)  !\left\vert b\right\vert
^{2\left(  \eta+\delta_{G}\right)  },\quad\left\vert b\right\vert \leq r_{G}.
\]

Thus%
\[
\left\vert \mathcal{R}_{2\eta+1}\left(  0,b\right)  \right\vert \leq\left(
2\pi\right)  ^{\frac{d}{2}}c_{G,\eta}\left\vert b\right\vert ^{2\left(
\eta+\delta_{G}\right)  }.
\]

Hence if $\operatorname*{diam}A_{x}\leq r_{G}$
\begin{align*}
\left(  2\pi\right)  ^{\frac{d}{2}}r_{x}\left(  x\right)   & \leq2\left(
2\pi\right)  ^{\frac{d}{2}}c_{G,\eta}\left(  \sum_{k=1}^{M}\left\vert
l_{k}\left(  x\right)  \right\vert \right)  \left(  \max_{k}\left\vert
x-a^{\left(  k\right)  }\right\vert \right)  ^{2\left(  \eta+\delta
_{G}\right)  }+\\
& \qquad\qquad+\left(  2\pi\right)  ^{\frac{d}{2}}c_{G,\eta}\left(  \sum
_{k=1}^{M}\left\vert l_{k}\left(  x\right)  \right\vert \right)  ^{2}\left(
\max_{j,k}\left\vert a^{\left(  j\right)  }-a^{\left(  k\right)  }\right\vert
\right)  ^{2\left(  \eta+\delta_{G}\right)  }\\
& <\left(  2\pi\right)  ^{\frac{d}{2}}c_{G,\eta}\left(  1+\sum_{k=1}%
^{M}\left\vert l_{k}\left(  x\right)  \right\vert \right)  ^{2}\left(
\operatorname*{diam}A_{x}\right)  ^{2\left(  \eta+\delta_{G}\right)  },
\end{align*}

so that when $\operatorname*{diam}A_{x}\leq r_{G}$,%
\[
\sqrt{r_{x}\left(  x\right)  }\leq\sqrt{c_{G,\eta}}\left(  1+\sum_{k=1}%
^{M}\left\vert l_{k}\left(  x\right)  \right\vert \right)  \left(
\operatorname*{diam}A_{x}\right)  ^{\eta+\delta_{G}},
\]

as claimed.
\end{proof}

We can now prove our improved interpolant convergence estimate.

\begin{corollary}
\label{Cor_Thm_int_rx(x)_bnd_better_order}Suppose the notation and assumptions
of Lemma \ref{Lem_int_Lagrange_interpol} and Theorem
\ref{Thm_int_rx(x)_bnd_better_order} hold. Suppose also that $\mathcal{I}_{X}$
is the minimal interpolant on $X$ of the data function $f_{d}\in X_{w}%
^{\theta}$.

Then when $0<h_{X}\leq\min\left\{  h_{\Omega,\theta},r_{\Omega}\right\}  $ we
have%
\begin{equation}
\left\vert f_{d}\left(  x\right)  -\mathcal{I}_{X}f_{d}\left(  x\right)
\right\vert \leq\left\vert f_{d}-\mathcal{I}_{X}f_{d}\right\vert _{w,\theta
}k_{\Omega,\theta,\eta}\left(  c_{\Omega,\theta}h_{X}\right)  ^{\eta
+\delta_{G}},\quad x\in\overline{\Omega},\label{p9100}%
\end{equation}

where $k_{\Omega,\theta,\eta}=\sqrt{c_{G,\eta}}\left(  1+K_{\Omega,\theta
}^{\prime}\right)  $, the constants $\eta,c_{G,\eta},r_{\Omega},\delta_{G}$
come from Theorem \ref{Thm_int_rx(x)_bnd_better_order}, and the constants
$c_{\Omega,\theta},K_{\Omega,\theta}^{\prime}$, $h_{\Omega,\theta}$ come from
Lemma \ref{Lem_int_Lagrange_interpol}.

?? Next, since $\left\vert f_{d}-\mathcal{I}_{X}f_{d}\right\vert _{w,\theta
}\leq\left\vert f_{d}\right\vert _{w,\theta}$, the order of convergence is
$\eta+\delta_{G}$.

Finally, for all independent data $X$,%
\begin{equation}
\left\vert f_{d}\left(  x\right)  -\mathcal{I}_{X}f_{d}\left(  x\right)
\right\vert \leq\left\vert f_{d}-\mathcal{I}_{X}f_{d}\right\vert _{w,\theta
}k_{\Omega,\theta,\eta}\left(  \min\left\{  \operatorname{diam}\Omega
,r_{G}\right\}  \right)  ^{\eta+\delta_{G}},\quad x\in\overline{\Omega
}.\label{p9103}%
\end{equation}

\end{corollary}

\begin{proof}
If $x\in\Omega$ and $\operatorname*{diam}A_{x}\leq r_{G}$ then from
\ref{p9102}, \ref{q87} and Lemma \ref{Lem_int_Lagrange_interpol}%
\begin{align*}
\left\vert f_{d}\left(  x\right)  -\mathcal{I}_{X}f_{d}\left(  x\right)
\right\vert  & \leq\left\vert f_{d}-\mathcal{I}_{X}f_{d}\right\vert
_{w,\theta}\sqrt{r_{x}\left(  x\right)  }\\
& \leq\left\vert f_{d}-\mathcal{I}_{X}f_{d}\right\vert _{w,\theta}%
\sqrt{c_{G,\eta}}\left(  1+\sum_{k=1}^{M}\left\vert l_{k}\left(  x\right)
\right\vert \right)  \left(  \operatorname*{diam}A_{x}\right)  ^{\eta
+\delta_{G}}\\
& \leq\left\vert f_{d}-\mathcal{I}_{X}f_{d}\right\vert _{w,\theta}%
\sqrt{c_{G,\eta}}\left(  1+K_{\Omega,\theta}^{\prime}\right)  \left(
\operatorname*{diam}A_{x}\right)  ^{\eta+\delta_{G}}\\
& =\left\vert f_{d}-\mathcal{I}_{X}f_{d}\right\vert _{w,\theta}k_{\Omega
,\theta,\eta}\left(  \operatorname*{diam}A_{x}\right)  ^{\eta+\delta_{G}}.
\end{align*}

Now by Lemma \ref{Lem_int_Lagrange_interpol} if $h_{X}\leq h_{\Omega,\theta}$
then there exists a minimal unisolvent $A\subset\Omega$ such that
$\operatorname*{diam}A_{x}\leq c_{\Omega,\theta}h_{X}$ which proves
\ref{p9100}. Finally, the estimate \ref{p9103} is true since
$\operatorname*{diam}A_{x}\leq\operatorname{diam}\Omega$.
\end{proof}

We now assume that the weight function has \textbf{property W3.2}.

\begin{lemma}
\label{Lem_W3_range_parm_k}Suppose for a given order $\theta$ the weight
function $w$ has property \textbf{W3.2} for some $\kappa$. Then one of the
following must hold:

\begin{enumerate}
\item[Property 1] For some $s>\kappa$, $w$ has property W3.2 for all
$0\leq\kappa<s$ and does not have property W3.2 for $\kappa\geq s$;

\item[Property 2] For some $s\geq\kappa$, $w$ has property W3.2 for all
$0\leq\kappa\leq s$ and does not have property W3.2 for $\kappa>s$;

\item[Property 3] $w$ has property W3.2 for all $\kappa\geq0$ i.e. $s=\infty$.
\end{enumerate}
\end{lemma}

\begin{proof}
This result follows from the fact that if property W3.2 holds for
$\kappa=\kappa_{0}$ then W3.2 holds for all $0\leq\kappa\leq\kappa_{0}$.
\end{proof}

The examples of this section will have either property 1 or property 3 of the
last lemma and these examples will be concerned with finding the `maximum'
order of convergence given the constraint applied to $\kappa$. Now $\eta
=\max\limits_{0\leq\kappa<s}\min\left\{  \theta,\frac{1}{2}\left\lfloor
2\kappa\right\rfloor \right\}  $ is the maximum order of convergence defined
by Theorem \ref{Thm_int_rx(x)_bound}. To this is added an increment of
convergence order $\delta_{G}$ defined by Theorem
\ref{Thm_int_rx(x)_bnd_better_order}. The following result will be used to
calculate $\eta$ and $\sigma$ when the weight function has property 1 or 3 of
the last lemma.

\begin{theorem}
\label{Th_eta_sigma}Given $s>0$ suppose $w$ is a weight function with
properties W2 and \textbf{W3.2} for order $\theta$ and all $\kappa<s$. Set
$\eta=\max\limits_{0\leq\kappa<s}\min\left\{  \theta,\frac{1}{2}\left\lfloor
2\kappa\right\rfloor \right\}  $ and $\sigma=\max\limits_{0\leq\kappa<s}%
\min\left\{  \theta,\frac{1}{2}\left\lfloor 2\kappa+1\right\rfloor \right\}  $.

Then
\begin{equation}
\eta=\left\{
\begin{array}
[c]{ll}%
\min\left\{  \theta,\frac{1}{2}\left\lfloor 2s\right\rfloor \right\}  , &
if\text{ }2s\text{ }is\text{ }not\text{ }an\text{ }integer,\\
\min\left\{  \theta,s-\frac{1}{2}\right\}  , & if\text{ }2s\text{ }is\text{
}an\text{ }integer.
\end{array}
\right. \label{q88}%
\end{equation}

If $s\leq\theta$ then
\begin{equation}
\eta=\left\{
\begin{array}
[c]{ll}%
\frac{1}{2}\left\lfloor 2s\right\rfloor , & if\text{ }2s\text{ }is\text{
}not\text{ }an\text{ }integer,\\
s-\frac{1}{2}, & if\text{ }2s\text{ }is\text{ }an\text{ }integer.
\end{array}
\right. \label{p9107}%
\end{equation}

If $s>\theta$ then
\begin{equation}
\eta=\theta.\label{p801}%
\end{equation}

Regarding $\sigma$:
\begin{equation}
\sigma=\left\{
\begin{array}
[c]{ll}%
\min\left\{  \theta,\frac{1}{2}\left\lfloor 2s+1\right\rfloor \right\}  , &
if\text{ }2s\text{ }is\text{ }not\text{ }an\text{ }integer,\\
\min\left\{  \theta,s\right\}  , & if\text{ }2s\text{ }is\text{ }an\text{
}integer,
\end{array}
\right. \label{p9101}%
\end{equation}

and $\sigma$ is related to $\eta$ by%
\begin{equation}
\sigma=\left\{
\begin{array}
[c]{ll}%
\eta, & \left\lfloor 2s\right\rfloor \geq2\theta,\\
\eta+\frac{1}{2}, & \left\lfloor 2s\right\rfloor \leq2\theta-1.
\end{array}
\right. \label{p9104}%
\end{equation}

\end{theorem}

\begin{proof}
Let $\eta_{\kappa}=\min\left\{  \theta,\frac{1}{2}\left\lfloor 2\kappa
\right\rfloor \right\}  $. Suppose $2s$ is not an integer and $0<\varepsilon
<2s-\left\lfloor 2s\right\rfloor $. Then%
\begin{align*}
\eta=\max\left\{  \eta_{\kappa}:\left(  1-\frac{\varepsilon}{2s}\right)
s\leq\kappa<s\right\}  \geq\min\left\{  \theta,\frac{1}{2}\left\lfloor
2\left(  1-\frac{\varepsilon}{2s}\right)  s\right\rfloor \right\}   &
=\min\left\{  \theta,\frac{1}{2}\left\lfloor 2s-\varepsilon\right\rfloor
\right\} \\
&  \geq\min\left\{  \theta,\frac{1}{2}\left\lfloor 2s\right\rfloor \right\}  ,
\end{align*}

but $\eta=\max\left\{  \eta_{\kappa}:0\leq\kappa<s\right\}  \leq\max\left\{
\eta_{\kappa}:0\leq\kappa\leq s\right\}  \leq\min\left\{  \theta,\frac{1}%
{2}\left\lfloor 2s\right\rfloor \right\}  $ so when $2s$ is not an integer it
follows that $\eta=\min\left\{  \theta,\frac{1}{2}\left\lfloor 2s\right\rfloor
\right\}  $.

If $2s$ is an integer then
\begin{align*}
\eta=\max\left\{  \eta_{\kappa}:0\leq\kappa<s\right\}   &  =\max\left\{
\eta_{\kappa}:0\leq2\kappa<2s\right\} \\
&  =\max\left\{  \eta_{\kappa}:2s-1<2\kappa<2s\right\} \\
&  =\max\left\{  \min\left\{  \theta,\frac{1}{2}\left\lfloor 2\kappa
\right\rfloor \right\}  :2s-1<2\kappa<2s\right\} \\
&  =\max\left\{  \min\left\{  \theta,\frac{1}{2}\left(  2s-1\right)  \right\}
:2s-1<2\kappa<2s\right\} \\
&  =\min\left\{  \theta,s-\frac{1}{2}\right\}  .
\end{align*}
\medskip

To prove \ref{p9107} suppose $s\leq\theta$: If $2s$ is not an integer then by
\ref{q88}, $\eta=\min\left\{  \theta,\frac{1}{2}\left\lfloor 2s\right\rfloor
\right\}  $. But $\frac{1}{2}\left\lfloor 2s\right\rfloor \leq\frac{1}%
{2}\left\lfloor 2\theta\right\rfloor =\theta$ so $\eta=\frac{1}{2}\left\lfloor
2s\right\rfloor $. If $2s$ is an integer then again by \ref{q88}, $\eta
=\min\left\{  \theta,s-\frac{1}{2}\right\}  =s-\frac{1}{2}$.\medskip

To prove \ref{p801} suppose $s>\theta$: If $2s$ is not an integer then
$\eta=\min\left\{  \theta,\frac{1}{2}\left\lfloor 2s\right\rfloor \right\}
=\theta$. If $2s$ is an integer then $\eta=\min\left\{  \theta,s-\frac{1}%
{2}\right\}  $. But $2s>2\theta$ implies $2s-1\geq2\theta$ i.e. $s-\frac{1}%
{2}\geq\theta$ and $\eta=\theta$.\medskip

Let $\sigma_{\kappa}=\min\left\{  \theta,\frac{1}{2}\left\lfloor
2\kappa+1\right\rfloor \right\}  $. Suppose $2s$ is not an integer and
$0<\varepsilon<2s-\left\lfloor 2s\right\rfloor $. Then%
\begin{align*}
\sigma=\max\left\{  \sigma_{\kappa}:\left(  1-\frac{\varepsilon}{2s}\right)
s\leq\kappa<s\right\}   &  \geq\min\left\{  \theta,\frac{1}{2}\left\lfloor
2\left(  1-\frac{\varepsilon}{2s}\right)  s\right\rfloor +\frac{1}{2}\right\}
\\
&  =\min\left\{  \theta,\frac{1}{2}\left\lfloor 2s-\varepsilon\right\rfloor
+\frac{1}{2}\right\} \\
&  \geq\min\left\{  \theta,\frac{1}{2}\left\lfloor 2s\right\rfloor +\frac
{1}{2}\right\}  ,
\end{align*}

but $\sigma=\max\left\{  \sigma_{\kappa}:0\leq\kappa<s\right\}  \leq
\max\left\{  \sigma_{\kappa}:0\leq\kappa\leq s\right\}  \leq\min\left\{
\theta,\frac{1}{2}\left\lfloor 2s\right\rfloor +\frac{1}{2}\right\}  $ so when
$2s$ is not an integer it follows that $\eta=\min\left\{  \theta,\frac{1}%
{2}\left\lfloor 2s\right\rfloor +\frac{1}{2}\right\}  =\min\left\{
\theta,\frac{1}{2}\left\lfloor 2s+1\right\rfloor \right\}  $.

If $2s$ is an integer then
\begin{align*}
\sigma=\max\left\{  \sigma_{\kappa}:0\leq\kappa<s\right\}   &  =\max\left\{
\sigma_{\kappa}:0\leq2\kappa<2s\right\} \\
&  =\max\left\{  \sigma_{\kappa}:2s-1<2\kappa<2s\right\} \\
&  =\max\left\{  \min\left\{  \theta,\frac{1}{2}\left\lfloor 2\kappa
\right\rfloor +\frac{1}{2}\right\}  :2s-1<2\kappa<2s\right\} \\
&  =\max\left\{  \min\left\{  \theta,\frac{1}{2}\left(  2s-1\right)  +\frac
{1}{2}\right\}  :2s-1<2\kappa<2s\right\} \\
&  =\min\left\{  \theta,s\right\}  ,
\end{align*}

which proves \ref{p9101}.

Finally, \ref{p9104} is proved by using \ref{p9101} and \ref{q88} to calculate
$\sigma$ and $\eta$ for $\left\lfloor 2s\right\rfloor \geq2\theta$ and
$\left\lfloor 2s\right\rfloor \leq2\theta-1$.
\end{proof}

Next we will illustrate the theory of this subsection by using as the example
of the \textbf{radial thin-plate splines} which were among the positive order
basis functions studied in Chapter \ref{Ch_basis_no_So2n}.\medskip

\textbf{Example: The thin-plate splines} We will derive the combinations of
$\eta$ and $\delta_{G}$ for which the thin-plate spline\textbf{\ }basis
functions satisfy the requirements of Theorem
\ref{Thm_int_rx(x)_bnd_better_order}.

Corollary \ref{Cor_Thm_int_rx(x)_bnd_better_order} then tells us the orders of
convergence of the interpolant are $\eta+\delta_{G}$, which is an improvement
of $\delta_{G}$.

The thin-plate spline weight and basis functions were studied in Subsection
\ref{SbSect_thin_plate_basis} of Chapter \ref{Ch_basis_no_So2n}. By Theorem
\ref{Thm_surf_spline_ext_basis} the equation%
\begin{equation}
w\left(  \xi\right)  =\frac{1}{e\left(  s\right)  }\left\vert \xi\right\vert
^{-2\theta+2s+d},\label{q96}%
\end{equation}

defines a thin-plate spline weight function with properties W2.1 and W3.2 with
positive integer order $\theta$ and non-negative $\kappa\in\mathbb{R}^{1}$ iff
$\kappa<s<\theta$. The corresponding basis functions are defined by
\begin{equation}
G\left(  x\right)  =\left\{
\begin{array}
[c]{ll}%
\left(  -1\right)  ^{s+1}\left\vert x\right\vert ^{2s}\log\left\vert
x\right\vert , & s=1,2,3,\ldots,\\
\left(  -1\right)  ^{\left\lceil s\right\rceil }\left\vert x\right\vert
^{2s}, & s>0,\text{ }s\neq1,2,3,\ldots.
\end{array}
\right. \label{q95}%
\end{equation}

Now suppose $\kappa<s<\theta$.\medskip

\fbox{\textbf{Case 1:} $\kappa<s<\theta$, $s\notin\mathbb{Z}$} By \ref{q95},
$G\left(  x\right)  =\left(  -1\right)  ^{\left\lceil s\right\rceil
}\left\vert x\right\vert ^{2s}$. This is a homogeneous function of degree $2s
$ and thus%
\begin{equation}
\left\vert D^{\beta}G\left(  x\right)  \right\vert \leq k_{s,\left\vert
\beta\right\vert }\left\vert x\right\vert ^{2s-\left\vert \beta\right\vert
},\quad\beta\geq0,\label{q77}%
\end{equation}

where $k_{s,n}=\max\limits_{\left\vert \beta\right\vert =n}\max
\limits_{\left\vert x\right\vert =1}\left\vert D^{\beta}G\left(  x\right)
\right\vert <\infty$. Hence $D^{\beta}G\in L_{loc}^{1}$ iff
\begin{equation}
-d<2s-\left\vert \beta\right\vert ,\label{p92.1}%
\end{equation}

which implies that for all $\beta\geq0$,
\begin{equation}
\int_{0}^{1}\left(  1-t\right)  ^{2\eta}\left\vert \left(  D^{\beta}G\right)
(z+tb)\right\vert dt\leq k_{s,\left\vert \beta\right\vert }\int_{0}%
^{1}\left\vert z+tb\right\vert ^{2s-\left\vert \beta\right\vert }%
dt.\label{q78}%
\end{equation}

We want to apply Theorem \ref{Thm_int_rx(x)_bnd_better_order} which assumes
$\left\vert \beta\right\vert =2\eta+1=\left\lfloor 2\kappa\right\rfloor
+1=\left\lceil 2\kappa\right\rceil $. From \ref{p9107} of Theorem
\ref{Th_eta_sigma},
\begin{equation}
\eta=\frac{1}{2}\left\lfloor 2s\right\rfloor ,\label{q90}%
\end{equation}

so $2s-\left\vert \beta\right\vert =2s-2\eta-1=2s-\left\lfloor 2s\right\rfloor
-1$ i.e. $-1<2s-\left\vert \beta\right\vert \leq0$, which means \ref{p92.1} is
satisfied. Hence%
\begin{align*}
\int\limits_{0}^{1}\left\vert z+bt\right\vert ^{2s-\left\vert \beta\right\vert
}dt\leq\int\limits_{0}^{1}\left\vert \left\vert z\right\vert -\left\vert
b\right\vert t\right\vert ^{2s-\left\vert \beta\right\vert }dt &  =\left\vert
b\right\vert ^{2s-\left\vert \beta\right\vert }\int\limits_{0}^{1}\left\vert
\frac{\left\vert z\right\vert }{\left\vert b\right\vert }-t\right\vert
^{2s-\left\vert \beta\right\vert }dt\\
&  =\left\vert b\right\vert ^{2s-\left\vert \beta\right\vert }\int%
\limits_{\left\vert z\right\vert \left\vert b\right\vert ^{-1}-1}^{\left\vert
z\right\vert \left\vert b\right\vert ^{-1}}\left\vert u\right\vert
^{2s-\left\vert \beta\right\vert }du\\
&  <2\left\vert b\right\vert ^{2s-\left\vert \beta\right\vert }\int%
\limits_{0}^{1+\left\vert z\right\vert \left\vert b\right\vert ^{-1}%
}u^{2s-\left\vert \beta\right\vert }du\\
&  \leq\frac{2\left\vert b\right\vert ^{2s-\left\vert \beta\right\vert }%
}{2s-\left\vert \beta\right\vert +1}\left(  1+\frac{\left\vert z\right\vert
}{\left\vert b\right\vert }\right)  ^{2s-\left\vert \beta\right\vert +1}\\
&  <\infty,
\end{align*}

and \ref{q78} implies
\begin{align*}
\int_{0}^{1}\left(  1-t\right)  ^{2\eta}\left\vert \left(  D^{\beta}G\right)
(z+tb)\right\vert dt  & \leq\frac{2k_{s,\left\vert \beta\right\vert }%
}{2s-\left\vert \beta\right\vert +1}\left(  1+\frac{\left\vert z\right\vert
}{\left\vert b\right\vert }\right)  ^{2s-\left\vert \beta\right\vert +1}\\
& \leq\frac{2k_{s,\left\vert \beta\right\vert }}{2s-\left\vert \beta
\right\vert +1}\left(  1+\frac{\left\vert z\right\vert }{\left\vert
b\right\vert }\right)  ,
\end{align*}

which in turn implies polynomial growth in $\left\vert z\right\vert $.

Condition \ref{q081} is%
\[
\left\vert \int_{0}^{1}\left(  1-t\right)  ^{2\eta}\left(  \left(  \cdot
D\right)  ^{2\eta+1}G\right)  \left(  tb\right)  dt\right\vert \leq\left(
2\pi\right)  ^{\frac{d}{2}}c_{G,\eta}\left(  2\eta\right)  !\left\vert
b\right\vert ^{2\left(  \eta+\delta_{G}\right)  },\quad\left\vert b\right\vert
\leq r_{G}.
\]

But%
\[
\mathcal{R}_{2\eta+1}\left(  0,b\right)  =\frac{1}{\left(  2\eta\right)
!}\int_{0}^{1}\left(  1-t\right)  ^{2\eta}\left(  \left(  \cdot D\right)
^{2\eta+1}G\right)  \left(  tb\right)  dt,
\]

and $G\left(  x\right)  =\left(  -1\right)  ^{\left\lceil s\right\rceil
}\left\vert x\right\vert ^{2s}$. Thus by part 2 of Corollary
\ref{Cor_Lem_Taylor_extension}, $\left(  \left(  \widehat{\cdot}D\right)
^{k}G\right)  \left(  0\right)  =0$ when $0<k<2s$ i.e. when $0<k\leq
\left\lfloor 2s\right\rfloor =2\eta$ and this implies $\mathcal{R}_{2\eta
+1}\left(  0,b\right)  =G\left(  b\right)  =\left(  -1\right)  ^{\left\lceil
s\right\rceil }\left\vert b\right\vert ^{2s}$. Consequently%
\[
\left\vert \int_{0}^{1}\left(  1-t\right)  ^{2\eta}\left(  \left(
\widehat{\cdot}D\right)  ^{2\eta+1}G\right)  \left(  tb\right)  dt\right\vert
=\left(  2\eta\right)  !\left\vert b\right\vert ^{2s},\quad b\in\mathbb{R}%
^{d},
\]

For $\left\vert b\right\vert \leq r_{G}$ we want%
\begin{align*}
\left(  2\pi\right)  ^{\frac{d}{2}}c_{G,\eta}\left(  2\eta\right)  !\left\vert
b\right\vert ^{2\left(  \eta+\delta_{G}\right)  }  & \leq\left(  2\eta\right)
!\left\vert b\right\vert ^{2s},\\
\left(  2\pi\right)  ^{\frac{d}{2}}c_{G,\eta}\left\vert b\right\vert
^{2\left(  \eta+\delta_{G}\right)  }  & \leq\left\vert b\right\vert ^{2s},\\
c_{G,\eta}  & \leq\left(  2\pi\right)  ^{-\frac{d}{2}}\left\vert b\right\vert
^{2s-2\left(  \eta+\delta_{G}\right)  }.
\end{align*}

But from \ref{q90}, $\eta=\frac{1}{2}\left\lfloor 2s\right\rfloor $ so
$2s=2\eta+\left(  2s-\left\lfloor 2s\right\rfloor \right)  $ and we can write
$2s=2\eta+2\delta_{G}$ where $\delta_{G}=s-\frac{1}{2}\left\lfloor
2s\right\rfloor <\frac{1}{2}$. Hence we have
\[
c_{G,\eta}=\left(  2\pi\right)  ^{-\frac{d}{2}}\left(  r_{G}\right)
^{2s-\left\lfloor 2s\right\rfloor -2\delta_{G}}=\left(  2\pi\right)
^{-\frac{d}{2}},
\]

and%
\begin{equation}
r_{G}=\infty;\quad\delta_{G}=s-\frac{1}{2}\left\lfloor 2s\right\rfloor ;\quad
c_{G,\eta}=\left(  2\pi\right)  ^{-\frac{d}{2}}.\label{a018}%
\end{equation}

Note that $s-\frac{1}{2}\left\lfloor 2s\right\rfloor =0$ iff $s=n$ or
$n+\frac{1}{2}$ for some integer $n\geq0$. Hence if $s=n+\frac{1}{2}$ for some
integer $n$ then $\delta_{G}=0$ and there is no improvement.\medskip

\fbox{\textbf{Case 2:} $\kappa<s<\theta$, $s\in\mathbb{Z}.$} To apply Theorem
\ref{Thm_int_rx(x)_bnd_better_order} we need to estimate the derivatives
$\left\{  D^{\beta}G\right\}  _{\left\vert \beta\right\vert =2\eta+1}$. Now
$\eta=\frac{1}{2}\left\lfloor 2\kappa\right\rfloor $ so $2\eta+1=\left\lceil
2\kappa\right\rceil $. But from \ref{p9107} of Theorem \ref{Th_eta_sigma},
\[
\eta=s-\frac{1}{2}%
\]

so that $\left\vert \beta\right\vert =2\eta+1=2s=\left\lceil 2\kappa
\right\rceil $.

Next, \ref{q95} implies $G\left(  x\right)  =\left(  -1\right)  ^{s+1}%
\left\vert x\right\vert ^{2s}\log\left\vert x\right\vert $. But any first
order derivative of $\log\left\vert x\right\vert $ is a homogeneous function
of order $-1$ so there exist constants $k_{s,\left\vert \beta\right\vert }$
and $k_{s,\left\vert \beta\right\vert }^{\prime}$ such that%
\[
\left\vert D^{\beta}G\left(  x\right)  \right\vert \leq k_{s,\left\vert
\beta\right\vert }\left\vert x\right\vert ^{2s-\left\vert \beta\right\vert
}\left\vert \log\left\vert x\right\vert \right\vert +k_{s,\left\vert
\beta\right\vert }^{\prime}\left\vert x\right\vert ^{2s-\left\vert
\beta\right\vert },\quad0\leq\left\vert \beta\right\vert \leq2s,
\]

and when $\left\vert \beta\right\vert =2s$
\begin{equation}
\left\vert D^{\beta}G\left(  x\right)  \right\vert \leq k_{s,2s}\left\vert
\log\left\vert x\right\vert \right\vert +k_{s,2s}^{\prime},\quad\left\vert
\beta\right\vert =2s.\label{q84}%
\end{equation}

Now $\left\vert \log\left\vert x\right\vert \right\vert \in L_{loc}^{1}$ iff
$\int_{\left\vert \cdot\right\vert \leq1}\left\vert \log\left\vert
x\right\vert \right\vert dx<\infty$ and by using the spherical polar
coordinates $\left(  r,\phi\right)  $ of Section \ref{Sect_App_SpherCoords} it
is straight forward to show that this integral exists because $r^{d}\left\vert
\log r\right\vert \in L_{loc}^{1}\left(  \mathbb{R}^{1}\right)  $.

Further, the integrals \ref{q761} satisfy the inequalities%
\begin{align}
\int_{0}^{1}\left(  1-t\right)  ^{2\eta}\left\vert \left(  D^{\beta}G\right)
\left(  z+tb\right)  \right\vert dt  & \leq\int_{0}^{1}\left\vert \left(
D^{\beta}G\right)  (z+tb)\right\vert dt\nonumber\\
& \leq k_{s,2s}\int_{0}^{1}\left\vert \log\left\vert z+tb\right\vert
\right\vert dt+k_{s,2s}^{\prime}\left\vert b\right\vert ^{2s}\int_{0}%
^{1}dt\nonumber\\
& =k_{s,2s}\int_{0}^{1}\left\vert \log\left\vert z+tb\right\vert \right\vert
dt+k_{s,2s}^{\prime}\left\vert b\right\vert ^{2s}.\label{q86}%
\end{align}

Suppose $d=1$. Then if $\left\vert z+b\right\vert \leq1$ and $\left\vert
z\right\vert \leq1$ then $\int_{0}^{1}\left\vert \log\left\vert
z+tb\right\vert \right\vert dt=\frac{1}{b}\int_{z}^{z+b}\left\vert
\log\left\vert \tau\right\vert \right\vert d\tau\leq\frac{1}{b}\int_{-1}%
^{1}\left\vert \log\left\vert \tau\right\vert \right\vert d\tau=\frac{2}{b}$.
Otherwise, if $\left\vert z+b\right\vert \geq1$ or $\left\vert z\right\vert
\geq1$ then set $a_{1}=\max\left\{  \left\vert z+b\right\vert ,\left\vert
z\right\vert \right\}  $ so that
\begin{align*}
\int_{0}^{1}\left\vert \log\left\vert z+tb\right\vert \right\vert dt=\frac
{1}{b}\int_{z}^{z+b}\left\vert \log\left\vert \tau\right\vert \right\vert
d\tau\leq\frac{1}{b}\int_{-1}^{1}\left\vert \log\left\vert \tau\right\vert
\right\vert d\tau+\frac{2}{b}\int_{1}^{a_{1}}\left\vert \log\left\vert
\tau\right\vert \right\vert d\tau &  \leq\frac{2}{b}+\frac{2}{b}\int%
_{1}^{a_{1}}a_{1}d\tau\\
&  \leq\frac{2}{b}+\frac{2}{b}\left\vert a_{1}\right\vert ^{2}\\
&  \leq\frac{2}{b}+\frac{2}{b}\left\vert z+b\right\vert ^{2},
\end{align*}

proving that when $d=1$ the integrals \ref{q761} have polynomial increase in
$z$ for any given $b\neq0$.

Now assume $d>1$. To evaluate $\int_{0}^{1}\left\vert \log\left\vert
z+tb\right\vert \right\vert dt$ first observe that the interval $\left[
z,z+b\right]  $ can be rotated so that $b_{1}>0$ and $b_{i}=0$ for $i\geq2$.
To this end set $z=\left(  z_{1},z^{\prime}\right)  $ so that
\begin{align}
\int_{0}^{1}\left\vert \log\left\vert z+tb\right\vert \right\vert dt=\int%
_{0}^{1}\left\vert \log\left\vert \left(  z_{1}+b_{1}t,z^{\prime}\right)
\right\vert \right\vert dt &  =\frac{1}{2}\int_{0}^{1}\left\vert \log\left(
\left(  z_{1}+b_{1}t\right)  ^{2}+\left\vert z^{\prime}\right\vert
^{2}\right)  \right\vert dt.\nonumber\\
&  =\frac{1}{2}\int_{0}^{1}\left\vert \log\left(  \left(  z_{1}+b_{1}t\right)
^{2}+\left\vert z^{\prime}\right\vert ^{2}\right)  \right\vert dt.\label{q85}%
\end{align}

If $\left\vert z^{\prime}\right\vert \geq1$ then $z$ and $z+b$ both lie
outside the unit sphere and
\begin{equation}
\int_{0}^{1}\left\vert \log\left\vert z+tb\right\vert \right\vert dt\leq
\frac{1}{2}\int_{0}^{1}\left(  \left(  z_{1}+b_{1}\right)  ^{2}+z_{2}%
^{2}\right)  dt=\frac{1}{2}\left(  \left(  z_{1}+b_{1}\right)  ^{2}+z_{2}%
^{2}\right)  \leq\left\vert z\right\vert ^{2}+\left\vert b\right\vert
^{2}.\label{q72}%
\end{equation}

If $\left\vert z^{\prime}\right\vert <1$ and $z$ and $z+b$ both lie inside the
unit sphere we have $\left(  z_{1}+b_{1}t\right)  ^{2}+\left\vert z^{\prime
}\right\vert ^{2}<1$ and%
\begin{align}
\int_{0}^{1}\left\vert \log\left\vert z+tb\right\vert \right\vert dt\leq
-\frac{1}{2b_{1}}\int_{0}^{1}\log\left(  \left\vert z_{1}+b_{1}t\right\vert
^{2}\right)  dt  & =-\frac{1}{\left\vert b\right\vert }\int_{z_{1}}%
^{z_{1}+b_{1}}\log\left\vert s\right\vert ds\nonumber\\
& \leq-\frac{2}{\left\vert b\right\vert }\int_{0}^{1}\log sds\nonumber\\
& =\frac{2}{\left\vert b\right\vert }.\label{q80}%
\end{align}

If $\left\vert z^{\prime}\right\vert <1$ and $z$ and $z+b$ do not both lie
inside the unit sphere set $a_{1}=\max\left\{  \left\vert z_{1}\right\vert
,\left\vert z_{1}+b_{1}\right\vert \right\}  \geq1$ so that by \ref{q85}%
\begin{align*}
\int_{0}^{1}\left\vert \log\left\vert z+tb\right\vert \right\vert dt  &
\leq\frac{1}{2}\int_{0}^{1}\left\vert \log\left(  \left(  z_{1}+b_{1}t\right)
^{2}+\left\vert z^{\prime}\right\vert ^{2}\right)  \right\vert dt\\
& \leq\frac{1}{2b_{1}}\int_{z_{1}}^{z_{1}+b_{1}}\left\vert \log\left(
\tau^{2}+\left\vert z^{\prime}\right\vert ^{2}\right)  \right\vert d\tau\\
& =-\frac{1}{\left\vert b\right\vert }\int_{0}^{\sqrt{1-z_{2}^{2}}}\log\left(
\tau^{2}+\left\vert z^{\prime}\right\vert ^{2}\right)  d\tau+\frac
{1}{\left\vert b\right\vert }\int_{\sqrt{1-z_{2}^{2}}}^{a_{1}}\log\left(
\tau^{2}+\left\vert z^{\prime}\right\vert ^{2}\right)  d\tau\\
& \leq-\frac{2}{\left\vert b\right\vert }\int_{0}^{1}\log\tau d\tau+\frac
{1}{\left\vert b\right\vert }\int_{\sqrt{1-z_{2}^{2}}}^{a_{1}}\left(
a_{1}^{2}+\left\vert z^{\prime}\right\vert ^{2}\right)  d\tau\\
& \leq-\frac{2}{\left\vert b\right\vert }\int_{0}^{1}\log\tau d\tau+\frac
{1}{\left\vert b\right\vert }\int_{0}^{a_{1}}\left(  a_{1}^{2}+\left\vert
z^{\prime}\right\vert ^{2}\right)  d\tau\\
& =\frac{2}{\left\vert b\right\vert }+\frac{1}{\left\vert b\right\vert }%
a_{1}\left(  a_{1}^{2}+\left\vert z^{\prime}\right\vert ^{2}\right)  ,
\end{align*}

and since $\int_{0}^{1}\log\tau d\tau=1$, $a_{1}\leq\left\vert z\right\vert
+\left\vert b\right\vert $ and $\left\vert z^{\prime}\right\vert ^{2}%
\leq\left\vert z\right\vert ^{2}$ it follows that $\int_{0}^{1}\left\vert
\log\left\vert z+tb\right\vert \right\vert dt$ has \textbf{polynomial increase
in} $z$. Indeed, combining this result with the estimates \ref{q80} and
\ref{q72} we can conclude that the integrals \ref{q761} have polynomial
increase in $z$ when $d>1$.

Next, condition \ref{q081}:%
\[
\left\vert \int_{0}^{1}\left(  1-t\right)  ^{2\eta}\left(  \left(  \cdot
D\right)  ^{2\eta+1}G\right)  \left(  tb\right)  dt\right\vert \leq\left(
2\pi\right)  ^{\frac{d}{2}}c_{G,\eta}\left(  2\eta\right)  !\left\vert
b\right\vert ^{2\left(  \eta+\delta_{G}\right)  },\quad\left\vert b\right\vert
\leq r_{G},
\]

must be satisfied and here we use same technique as Case 1 i.e. show $\left(
\left(  \widehat{\cdot}D\right)  ^{k}G\right)  \left(  0\right)  =0$ when
$0<k\leq2\eta=2s-1$. ??

From part 1 of Lemma \ref{Lem_deriv_rad_funcs},
\begin{align*}
\left(  \left(  \cdot D\right)  ^{k}G\right)  \left(  x\right)   & =\left(
-1\right)  ^{s+1}\left(  \widehat{\cdot}D\right)  ^{k}\left(  \left\vert
x\right\vert ^{2s}\log\left\vert x\right\vert \right) \\
& =\left(  -1\right)  ^{s+1}\left(  D_{t}^{k}\left(  t^{2s}\log t\right)
\right)  \left(  \left\vert x\right\vert \right)  ,
\end{align*}

and from part 7 of Lemma \ref{Lem_deriv_rad_funcs},%
\[
D_{t}^{k}\left(  t^{2s}\log t\right)  =k!\left(  \tbinom{2s}{k}\log
t+\sum\limits_{j=1}^{k}\frac{\left(  -1\right)  ^{j+1}}{j}\tbinom{2s}%
{k-j}\right)  t^{2s-k},
\]

we means that%
\[
\left(  \left(  \cdot D\right)  ^{k}G\right)  \left(  x\right)  =\left(
-1\right)  ^{s+1}k!\left(  \tbinom{2s}{k}\log\left\vert x\right\vert
+\sum\limits_{j=1}^{k}\frac{\left(  -1\right)  ^{j+1}}{j}\tbinom{2s}%
{k-j}\right)  \left\vert x\right\vert ^{2s-k},
\]

and $\left(  \left(  \widehat{\cdot}D\right)  ^{k}G\right)  \left(  0\right)
=0$ when $0<k\leq2s-1$. Thus%
\[
\mathcal{R}_{2\eta+1}\left(  0,b\right)  =G\left(  b\right)  =\left(
-1\right)  ^{s+1}\left\vert b\right\vert ^{2s}\log\left\vert b\right\vert ,
\]

and hence%
\begin{align*}
\left\vert \frac{1}{\left(  2\eta\right)  !}\int_{0}^{1}\left(  1-t\right)
^{2\eta}\left(  \left(  \cdot D\right)  ^{2\eta+1}G\right)  \left(  tb\right)
dt\right\vert  & =\left\vert b\right\vert ^{2s}\left\vert \log\left\vert
b\right\vert \right\vert =\left\vert b\right\vert ^{2\eta+1}\left\vert
\log\left\vert b\right\vert \right\vert ,\\
\left\vert \int_{0}^{1}\left(  1-t\right)  ^{2\eta}\left(  \left(  \cdot
D\right)  ^{2\eta+1}G\right)  \left(  tb\right)  dt\right\vert  & =\left(
2\eta\right)  !\left\vert b\right\vert ^{2\eta+1}\left\vert \log\left\vert
b\right\vert \right\vert ,
\end{align*}

For $\left\vert b\right\vert \leq r_{G}$ we want%
\begin{align*}
\left(  2\pi\right)  ^{\frac{d}{2}}c_{G,\eta}\left(  2\eta\right)  !\left\vert
b\right\vert ^{2\left(  \eta+\delta_{G}\right)  }  & \leq\left(  2\eta\right)
!\left\vert b\right\vert ^{2\eta+1}\left\vert \log\left\vert b\right\vert
\right\vert ,\\
\left(  2\pi\right)  ^{\frac{d}{2}}c_{G,\eta}\left\vert b\right\vert
^{2\left(  \eta+\delta_{G}\right)  }  & \leq\left\vert b\right\vert ^{2\eta
+1}\left\vert \log\left\vert b\right\vert \right\vert ,\\
c_{G,\eta}  & \leq\left(  2\pi\right)  ^{-\frac{d}{2}}\left\vert b\right\vert
^{1-2\delta_{G}}\left\vert \log\left\vert b\right\vert \right\vert ,
\end{align*}

which is true when%
\begin{equation}
r_{G}>0;\quad0\leq\delta_{G}<\frac{1}{2};\quad c_{G,\eta}=\left(  2\pi\right)
^{-\frac{d}{2}}\max_{r\in\left[  0,r_{G}\right]  }\left\{  r^{1-2\delta_{G}%
}\left\vert \log r\right\vert \right\}  ,\quad\eta=s-\frac{1}{2}.\label{a019}%
\end{equation}

\textbf{Note} that as $\delta_{G}\rightarrow\frac{1}{2}$, $c_{G,\eta
}\rightarrow\infty$.

\begin{conclusion}
From \ref{a018} and \ref{a019}:

\begin{enumerate}
\item If $\kappa<s<\theta$, $s\notin\mathbb{Z}$ then we can choose
$r_{G}=\infty$, $\delta_{G}=s-\frac{1}{2}\left\lfloor 2s\right\rfloor $ and
$c_{G,\eta}=\left(  2\pi\right)  ^{-\frac{d}{2}}$.

\item If $\kappa<s<\theta$, $s\in\mathbb{Z}$ then we can choose $r_{G}>0$,
$0\leq\delta_{G}<\frac{1}{2}$ and $\left(  2\pi\right)  ^{-\frac{d}{2}}%
\max\limits_{r\in\left[  0,r_{G}\right]  }\left\{  r^{1-2\delta_{G}}\left\vert
\log r\right\vert \right\}  $.
\end{enumerate}
\end{conclusion}

\subsection{Using data function Taylor series}

We use the Taylor series expansions of functions in $X_{w}^{\theta}$ derived
in Sections \ref{Sect_data_fn_Taylor_W3.2} and \ref{Sect_data_fn_Taylor_W3.1}.

The pointwise convergence estimates of this section will be obtained by
modifying those of Subsection \ref{Sect_unisolv} which used multipoint Taylor
series expansions and Lagrange interpolation using minimal unisolvent sets of
independent data points. Recall that the operator $\mathcal{Q}$ was introduced
in Definition \ref{Def_Aux_proj_operator}.

\begin{lemma}
\label{Lem_Qf_estim_multipt}\textbf{Multipoint Taylor series expansion}
Suppose $w\in W2$ and also $w\in W3.2$ or $w\in W3.1$ i.e. $w\in W3.3$, for
order $\theta$ and parameter $\kappa$. Set $n=\min\left\{  \theta
-1,\left\lfloor \underline{\kappa}\right\rfloor \right\}  $. Suppose
$A=\left\{  a^{\left(  k\right)  }\right\}  _{k=1}^{M_{\theta}}$ is a minimal
unisolvent set of order $\theta$ with cardinal basis $\left(  l_{k}\right)  $.
Then $f\in X_{w}^{\theta}$ implies $f\in C_{BP}^{\left(  \left\lfloor
\underline{\kappa}\right\rfloor \right)  }$ and%
\[
\mathcal{Q}f\left(  x\right)  =-\sum_{k=1}^{M}\left(  \mathcal{R}%
_{n+1}f\right)  \left(  x,a^{\left(  k\right)  }-x\right)  l_{k}\left(
x\right)  ,\quad x\in\mathbb{R}^{d},
\]

and we have the upper bound%
\[
\left\vert \mathcal{Q}f\left(  x\right)  \right\vert \leq\left(
\sum\limits_{k=1}^{M}\left\vert l_{k}\left(  x\right)  \right\vert \right)
\max_{k=1}^{M}\left\vert \left(  \mathcal{R}_{n+1}f\right)  \left(
x,a^{\left(  k\right)  }-x\right)  \right\vert ,\quad x\in\mathbb{R}^{d},
\]

where $\mathcal{R}_{n+1}f$ is the single point Taylor series remainder
\ref{a53}.
\end{lemma}

\begin{proof}
From Definition \ref{Def_Aux_proj_operator} of the operators $\mathcal{P}$ and
$\mathcal{Q}$,%
\[
\mathcal{Q}f\left(  x\right)  =f\left(  x\right)  -\mathcal{P}f\left(
x\right)  =f\left(  x\right)  -\sum_{k=1}^{M}f\left(  a^{\left(  k\right)
}\right)  l_{k}\left(  x\right)  =f\left(  x\right)  -\sum_{k=1}^{M}f\left(
x+\left(  a^{\left(  k\right)  }-x\right)  \right)  l_{k}\left(  x\right)  .
\]

Using the single point Taylor series expansion formula with remainder
\ref{a53} we have for each $k$%
\[
f\left(  x+\left(  a^{\left(  k\right)  }-x\right)  \right)  =\sum_{\left\vert
\beta\right\vert \leq n}\frac{D^{\beta}f(x)}{\beta!}\left(  a^{\left(
k\right)  }-x\right)  ^{\beta}+\left(  \mathcal{R}_{n+1}f\right)  \left(
x,a^{\left(  k\right)  }-x\right)  ,
\]

so that%
\begin{align*}
\sum_{k=1}^{M}f\left(  x+\left(  a^{\left(  k\right)  }-x\right)  \right)
l_{k}\left(  x\right)   & =\sum_{k=1}^{M}\left(  \sum_{\left\vert
\beta\right\vert \leq n}\frac{D^{\beta}f(x)}{\beta!}\left(  a^{\left(
k\right)  }-x\right)  ^{\beta}+\left(  \mathcal{R}_{n+1}f\right)  \left(
x,a^{\left(  k\right)  }-x\right)  \right)  l_{k}\left(  x\right) \\
& =\sum_{k=1}^{M}\sum_{\left\vert \beta\right\vert \leq n}\frac{D^{\beta}%
f(x)}{\beta!}\left(  a^{\left(  k\right)  }-x\right)  ^{\beta}l_{k}\left(
x\right)  +\sum_{k=1}^{M}\left(  \mathcal{R}_{n+1}f\right)  \left(
x,a^{\left(  k\right)  }-x\right)  l_{k}\left(  x\right)  .
\end{align*}

But by part 1 of Theorem \ref{Thm_P_op_properties} the operator $\mathcal{P}$
preserves polynomials of degree $<\theta$. Hence%
\begin{align*}
\sum_{k=1}^{M}\sum_{\left\vert \beta\right\vert \leq n}\frac{D^{\beta}%
f(x)}{\beta!}\left(  a^{\left(  k\right)  }-x\right)  ^{\beta}l_{k}\left(
x\right)   & =\sum_{\left\vert \beta\right\vert \leq n}\frac{D^{\beta}%
f(x)}{\beta!}\mathcal{P}_{y}\left(  \left(  y-x\right)  ^{\beta}\right)
\left(  y=x\right) \\
& =\sum_{\left\vert \beta\right\vert \leq n}\frac{D^{\beta}f(x)}{\beta
!}\left(  \left(  y-x\right)  ^{\beta}\right)  \left(  y=x\right) \\
& =f\left(  x\right)  ,
\end{align*}

leaving us with%
\[
\sum_{k=1}^{M}f\left(  x+\left(  a^{\left(  k\right)  }-x\right)  \right)
l_{k}\left(  x\right)  =f\left(  x\right)  +\sum_{k=1}^{M}\left(
\mathcal{R}_{n+1}f\right)  \left(  x,a^{\left(  k\right)  }-x\right)
l_{k}\left(  x\right)  ,
\]

and%
\[
\mathcal{Q}f\left(  x\right)  =-\sum_{k=1}^{M}\left(  \mathcal{R}%
_{n+1}f\right)  \left(  x,a^{\left(  k\right)  }-x\right)  l_{k}\left(
x\right)  ,
\]

so that%
\[
\left\vert \mathcal{Q}f\left(  x\right)  \right\vert \leq\left(  \sum
_{k=1}^{M}\left\vert l_{k}\left(  x\right)  \right\vert \right)  \max
_{k=1}^{M}\left\vert \left(  \mathcal{R}_{n+1}f\right)  \left(  x,a^{\left(
k\right)  }-x\right)  \right\vert ,
\]

as claimed.
\end{proof}

To obtain the highest order of convergence from our results we
\textbf{assume}
\begin{equation}
n=\min\left\{  \theta-1,\left\lfloor \underline{\kappa}\right\rfloor \right\}
.\label{a55}%
\end{equation}

Since $\mathcal{Q}f=\mathcal{Q}f_{\rho}$ and $\mathcal{I}_{X}p=p$ when $p\in
P_{\theta-1}$,
\[
\mathcal{I}_{X}f-f=\mathcal{Q}\left(  \mathcal{I}_{X}f-f\right)
=\mathcal{Q}\left(  \mathcal{I}_{X}\left(  f_{\rho}+p\right)  -\left(
f_{\rho}+p\right)  \right)  =\mathcal{Q}\left(  \mathcal{I}_{X}f_{\rho
}-f_{\rho}\right)  .
\]

But by Lemma \ref{Lem_Qf_estim_multipt},%
\[
\left\vert \mathcal{Q}g\left(  x\right)  \right\vert \leq\left(  \sum
_{k=1}^{M}\left\vert l_{k}\left(  x\right)  \right\vert \right)  \max
_{k=1}^{M}\left\vert \left(  \mathcal{R}_{n+1}g\right)  \left(  x,a^{\left(
k\right)  }-x\right)  \right\vert ,\quad g\in X_{w}^{\theta},\text{ }%
x\in\mathbb{R}^{d},
\]

so that%
\begin{align}
\left\vert \mathcal{I}_{X}f\left(  x\right)  -f\left(  x\right)  \right\vert
& =\left\vert \mathcal{Q}\left(  \mathcal{I}_{X}f_{\rho}-f_{\rho}\right)
\left(  x\right)  \right\vert \nonumber\\
& \leq\left(  \sum_{k=1}^{M}\left\vert l_{k}\left(  x\right)  \right\vert
\right)  \max_{j=1}^{M}\left\vert \left(  \mathcal{R}_{n+1}f_{\rho}\right)
\left(  x,a^{\left(  j\right)  }-x\right)  \right\vert \nonumber\\
& =\left\vert \widetilde{l}_{A}\left(  x\right)  \right\vert _{1}\max
_{j=1}^{M}\left\vert \left(  \mathcal{R}_{n+1}f_{\rho}\right)  \left(
x,a^{\left(  j\right)  }-x\right)  \right\vert ,\label{a69}%
\end{align}

where $\widetilde{l}_{A}=\left(  l_{k}\right)  $. We estimate $\left(
\mathcal{R}_{n+1}f_{\rho}\right)  \left(  x,a^{\left(  j\right)  }-x\right)  $
using \ref{a65} i.e.%
\begin{align}
&  \left\vert \mathcal{R}_{n+1}f_{\rho}\left(  x,a\right)  \right\vert
\nonumber\\
&  \leq\left\vert f\right\vert _{w,\theta}\tfrac{\left\vert a\right\vert
^{n+1}}{\left(  n+1\right)  !}\left\{
\begin{array}
[c]{l}%
\frac{C_{n}^{\left(  \rho,w\right)  }}{\left(  2\pi\right)  ^{\frac{d}{2}}%
}\left(  \left\Vert \widehat{a}D\widehat{\phi_{0}}\right\Vert _{1}+\right. \\
\qquad\left.  +\left(  \frac{3^{\theta-n+1}-1}{2}\right)  ^{\frac{1}{2}%
}\left(  \sum\limits_{j=0}^{\theta-n}\frac{\left(  \left\vert x\right\vert
^{2}+\left\vert a\right\vert ^{2}\right)  ^{j}}{j!}\right)  ^{\frac{1}{2}%
}\left(  \left\Vert \widehat{a}D\widehat{\phi_{0}}\right\Vert _{1}%
\sum\limits_{l=0}^{\theta-n}\frac{1}{l!}\left\Vert \left\vert \cdot\right\vert
^{2l}\widehat{a}D\widehat{\phi_{0}}\right\Vert _{1}\right)  ^{\frac{1}{2}%
}\right)  ,\\
\qquad\qquad\qquad\qquad\qquad\qquad\qquad\qquad\qquad\qquad if\text{ \ }%
n\leq\left\lfloor \kappa\right\rfloor <\theta,\\
\medskip\\
\frac{1}{\left(  2\pi\right)  ^{d}}\left(  \int\left(  \widehat{a}\widehat
{\xi}\right)  ^{2n}\frac{\left\vert \cdot\right\vert ^{2n}}{w\left\vert
\cdot\right\vert ^{2\theta}}\right)  ^{\frac{1}{2}}\left\Vert \widehat
{a}D\widehat{\phi_{0}}\right\Vert _{1},\qquad\qquad\qquad if\text{ \ }%
\theta\leq n\leq\left\lfloor \kappa\right\rfloor ,
\end{array}
\right\}  +\nonumber\\
&  +\left\vert f\right\vert _{w,\theta}\left\{
\begin{array}
[c]{l}%
\frac{\left\vert a\right\vert ^{n+1}}{\left(  n+1\right)  !}\left(
\frac{\left\Vert \left(  \widehat{\cdot}D\right)  ^{\theta}\phi_{\infty
}\right\Vert _{\infty;\leq r_{3}}}{\theta!}\left(  \int\limits_{\left\vert
\cdot\right\vert \leq r_{3}}\frac{\left(  \widehat{a}\xi\right)  ^{2\left(
n+1\right)  }}{w}\right)  ^{\frac{1}{2}}+\left(  \int\limits_{\left\vert
\cdot\right\vert \geq r_{3}}\frac{\left\vert \widehat{a}\xi\right\vert
^{2\left(  n+1\right)  }}{w\left\vert \cdot\right\vert ^{2\theta}}\right)
^{\frac{1}{2}}\right)  ,\qquad n\leq\left\lfloor \kappa\right\rfloor -1.\\
\medskip\\
\frac{\left\vert a\right\vert ^{\left\lceil \kappa\right\rceil }}{\left\lceil
\kappa\right\rceil !}\frac{\left\Vert \left(  \widehat{\cdot}D\right)
^{\theta}\phi_{\infty}\right\Vert _{\infty;\leq r_{3}}}{\theta!}\left(
\int\limits_{\left\vert \cdot\right\vert \leq r_{3}}\frac{\left(  \widehat
{a}\xi\right)  ^{2\left\lceil \kappa\right\rceil }}{w}\right)  ^{\frac{1}{2}%
}+\left\vert a\right\vert ^{\kappa}\left(  \frac{2}{\left\lfloor
\kappa\right\rfloor !}+\frac{1}{\left\lceil \kappa\right\rceil !}\right)
\left(  \int\limits_{\left\vert \cdot\right\vert \geq r_{3}}\frac{\left\vert
\widehat{a}\xi\right\vert ^{2\kappa}}{w\left\vert \cdot\right\vert ^{2\theta}%
}\right)  ^{\frac{1}{2}},\text{ \ }n=\left\lfloor \kappa\right\rfloor .
\end{array}
\right. \label{a028}%
\end{align}

and we may be able to simplify our calculations by choosing $\rho$ and/or
$\phi_{0}$ to be radial, as discussed in Subsection
\ref{SbSect_Tay_dat_rem_rad_fn}.

Noting the form of \ref{a028}, we now assume an estimate of the form%
\[
\left\vert \mathcal{R}_{n+1}f_{\rho}\left(  x,a\right)  \right\vert
\leq\left\vert f\right\vert _{w,\theta}\Phi\left(  \left\vert x\right\vert
,\left\vert a\right\vert \right)  ,\quad a,x\in\mathbb{R}^{d},
\]

where%
\[
\Phi\left(  s,t\right)  \leq\Phi\left(  s+\varepsilon,t+\eta\right)  \text{
}when\text{ }s,t,\varepsilon,\eta\geq0.
\]

Set $A_{x}=A\cup\left\{  x\right\}  $ and note that $\left\vert a^{\left(
j\right)  }-x\right\vert \leq\operatorname*{diam}$ $A_{x}$.

Next we use Proposition \ref{Prop_Lagr_interp_translat} to write for any
$b\in\mathbb{R}^{d}$,%
\begin{align*}
\left\vert \mathcal{I}_{X}f\left(  x\right)  -f\left(  x\right)  \right\vert
& =\left\vert \mathcal{I}_{\tau_{b}X}\left(  \tau_{b}f\right)  \left(
\tau_{b}x\right)  -\left(  \tau_{b}f\right)  \left(  \tau_{b}x\right)
\right\vert \\
& \leq\left\vert \tau_{b}f\right\vert _{w,\theta}\left\vert \widetilde
{l}_{\tau_{b}A}\left(  \tau_{b}x\right)  \right\vert _{1}\Phi\left(
\left\vert \tau_{b}x\right\vert ,\left\vert \tau_{b}a^{\left(  j\right)
}-\tau_{b}x\right\vert \right) \\
& =\left\vert f\right\vert _{w,\theta}\left\vert \widetilde{l}_{A}\left(
x\right)  \right\vert _{1}\max_{j=1}^{M}\Phi\left(  \left\vert \tau
_{b}x\right\vert ,\left\vert a^{\left(  j\right)  }-x\right\vert \right) \\
& \leq\left\vert f\right\vert _{w,\theta}\left\vert \widetilde{l}_{A}\left(
x\right)  \right\vert _{1}\Phi\left(  \left\vert \tau_{b}x\right\vert
,\operatorname*{diam}A_{x}\right)  .
\end{align*}

By Lemma \ref{Lem_int_Lagrange_interpol} there exist constants $h_{\Omega
,\theta}$ and $c_{\Omega,\theta}$ such that for each $x\in\Omega$ there exists
$A_{x}$ such that $\operatorname*{diam}A_{x}\leq c_{\Omega,\theta}h_{X,\Omega
}$ when $h_{X,\Omega}<h_{\Omega,\theta}$. Also there exists $K_{\Omega,\theta
}^{\prime}$ such that $\left\vert \widetilde{l}_{A}\left(  x\right)
\right\vert _{1}\leq K_{\Omega,\theta}^{\prime}$ for all $A_{x}$.Thus%
\[
\left\vert \mathcal{I}_{X}f\left(  x\right)  -f\left(  x\right)  \right\vert
\leq\left\vert f\right\vert _{w,\theta}K_{\Omega,\theta}^{\prime}\Phi\left(
\left\vert \tau_{b}x\right\vert ,c_{\Omega,\theta}h_{X,\Omega}\right)  ,\quad
x\in\overline{\Omega},\text{ }b\in\mathbb{R}^{d}.
\]

Finally, by part 4 of Proposition \ref{Prop_Lagr_interp_translat}, if we
choose $b$ to be the centre $c_{\Omega}$ of the smallest ball containing
$\Omega$ then the circumradius has upper bound
\[
\left\vert x-c_{\Omega}\right\vert \leq\operatorname*{diam}\Omega\sqrt
{\frac{d}{2\left(  d+1\right)  }},\quad x\in\overline{\Omega},
\]

and so $h_{X,\Omega}<h_{\Omega,\theta}$ implies%
\[
\left\vert \mathcal{I}_{X}f\left(  x\right)  -f\left(  x\right)  \right\vert
\leq\left\vert f\right\vert _{w,\theta}K_{\Omega,\theta}^{\prime}\Phi\left(
\operatorname*{diam}\Omega\sqrt{\frac{d}{2\left(  d+1\right)  }}%
,c_{\Omega,\theta}h_{X,\Omega}\right)  ,\quad x\in\overline{\Omega}.
\]

To summarize:

\begin{theorem}
With reference to \ref{a55}: Set $n=\min\left\{  \theta-1,\left\lfloor
\underline{\kappa}\right\rfloor \right\}  $ and suppose%
\[
\left\vert \mathcal{R}_{n+1}f_{\rho}\left(  x,a\right)  \right\vert
\leq\left\vert f\right\vert _{w,\theta}\Phi\left(  \left\vert x\right\vert
,\left\vert a\right\vert \right)  ,\quad a,x\in\mathbb{R}^{d},
\]

where $\Phi$ has property:%
\begin{equation}
\Phi\left(  s,t\right)  \leq\Phi\left(  s+\varepsilon,t+\eta\right)  \text{
}when\text{ }s,t,\varepsilon,\eta\geq0.\label{a026}%
\end{equation}

Then when $h_{X,\Omega}<h_{\Omega,\theta}$,
\begin{equation}
\left\vert \mathcal{I}_{X}f\left(  x\right)  -f\left(  x\right)  \right\vert
\leq\left\vert f\right\vert _{w,\theta}K_{\Omega,\theta}^{\prime}\Phi\left(
\operatorname*{diam}\Omega\sqrt{\frac{d}{2\left(  d+1\right)  }}%
,c_{\Omega,\theta}h_{X,\Omega}\right)  ,\quad x\in\overline{\Omega
}.\label{a027}%
\end{equation}

\end{theorem}

??

\begin{remark}
\textbf{More estimates} Since $\kappa\leq n+1$ and $\kappa<\left\lceil
\kappa\right\rceil $ we write
\[
\left(  h_{X,\Omega}\right)  ^{\left\lceil \kappa\right\rceil }=\left(
h_{X,\Omega}\right)  ^{\left\lceil \kappa\right\rceil -\kappa}\left(
h_{X,\Omega}\right)  ^{\kappa}\leq\left(  h_{\Omega,\theta}\right)
^{\left\lceil \kappa\right\rceil -\kappa}\left(  h_{X,\Omega}\right)
^{\kappa},
\]

and so%
\begin{align*}
\left(  h_{X,\Omega}\right)  ^{n+1}  & =\left(  h_{X,\Omega}\right)
^{n+1-\kappa}\left(  h_{X,\Omega}\right)  ^{\kappa}\\
& \leq\left(  h_{\Omega,\theta}\right)  ^{n+1-\kappa}\left(  h_{X,\Omega
}\right)  ^{\kappa}.
\end{align*}

\end{remark}

\subsection{Improved convergence using basis function Taylor series}

?? Basis function Taylor series expansions were discussed in Subsection
\ref{SbSect_Tay_basis_W3.1} ??.

\chapter{The Exact smoother and its convergence to the data
function\label{Ch_ExactSmth}}

\section{Introduction}

\subsection{In brief}

\textbf{Section} \ref{Sect_ex_prepare} Some results from previous
chapters.\smallskip

\textbf{Section} \ref{Sect_soln_Exact_smth} Existence and uniqueness of the
solution to the Exact smoothing problem. Derivation of several well-known
identities. Some properties of the Exact smoother mapping.\smallskip

\textbf{Section} \ref{Sect_Rxx_in_terms_of_Gxx} We derive equations which
express the reproducing kernel matrix in terms of the basis function
matrix.\smallskip

\textbf{Section} \ref{Sect_mat_from_rx} Matrices, vectors and bases derived
from the semi-Riesz representer.\smallskip

\textbf{Section} \ref{Sect_Exact_smth_mat_eqn} Matrix equations are derived
for the Exact smoother.\smallskip

\textbf{Section} \ref{Sect_converg_Exact_smth} Estimates the order of
pointwise convergence of the Exact smoother to its data function on a bounded
region.\smallskip

\textbf{Section} \ref{Sect_improved_err_estim} We derive a slightly improved
order of convergence.

\subsection{In more detail}

Given independent data $\left\{  x^{(i)}\right\}  _{i=1}^{N}$ and dependent
data $\left\{  y_{i}\right\}  _{i=1}^{N}$ the \textit{Exact smoothing problem}
involves minimizing the functional $\rho\left\vert f\right\vert _{w,\theta
}^{2}+\frac{1}{N}\sum_{i=1}^{N}\left\vert f(x^{(i)})-y_{i}\right\vert ^{2}$
over the \textit{semi-Hilbert space} $X_{w}^{\theta}$ of continuous functions.
Here $\rho>0$ is termed the \textit{smoothing parameter}. The descriptor
\textit{Exact} is ours. \textbf{Note} that this smoothing functional is
different from that used in Narcowich, Ward and Wendland
\cite{NarcWardWend2004} which is%
\[
\rho\left\Vert f\right\Vert _{w,0}^{2}+\sum_{i=1}^{N}\left\vert f(x^{(i)}%
)-y_{i}\right\vert ^{2}.
\]

To obtain their error estimates you simply replace $\rho N$ by $\rho$ in our
error formulas.

We first prove the existence and uniqueness of a solution to this problem
using the theory of \textit{basis functions} and \textit{semi-Hilbert spaces}
generated by \textit{weight functions}. We then derive orders for the
pointwise convergence of the smoother to its data function as the density of
the data increases. Each function in $X_{w}^{\theta}$ is considered a
legitimate data function.

\textbf{Section} \ref{Sect_ex_prepare} summarizes some of the theory needed
from previous chapters.

The weight function can be considered to have two parameters, the
\textit{integer} \textit{order parameter} $\theta\geq1$ and the
\textit{smoothness parameter} $\kappa\geq0$. Such a weight function will
generate the spaces $X_{w}^{\theta}$ and the basis functions of order $\theta
$. We next state the results we require from Chapter \ref{Ch_Interpol} where
we studied the minimum seminorm interpolation problem. These results were
adapted from Light and Wayne. These results start with the fundamental concept
of a \textit{unisolvent set} of points. Using a \textit{minimal unisolvent
set} we define the Lagrange polynomial interpolation operator $\mathcal{P}$
and $\mathcal{Q}=I-\mathcal{P}$ as well as the unisolvency matrices and the
\textit{Light norm} which makes $X_{w}^{\theta}$ a reproducing kernel Hilbert
space. Thus the \textit{Riesz representer} $R_{x}$ of the evaluation
functional $f\rightarrow f\left(  x\right)  $ exists and can be expressed in
terms of any basis function. In fact, the Riesz representer of the evaluation
functional $f\rightarrow D^{\gamma}f\left(  x\right)  $ is $D_{x}^{\gamma
}R_{x}$ for $\left\vert \gamma\right\vert \leq\min\kappa$. In this document we
follow Light and Wayne and instead of using the reproducing kernel we use the
Riesz representer of the evaluation functional $f\rightarrow f\left(
x\right)  $ which can be readily generalized to the evaluation of derivatives.

We next present two types of matrix, the \textit{basis function matrix} and
the \textit{reproducing kernel matrix}. The basis function, reproducing kernel
and unisolvency\ matrices will be used to construct the block matrix equations
for the basis function smoother of this document. It will turn out that the
Exact smoother will lie in the same \textit{finite dimensional basis function
space} $W_{G,X}$ as the minimum seminorm interpolant. Here $X=\left\{
x^{\left(  i\right)  }\right\}  _{i=1}^{N}$ denotes a unisolvent set of
independent data points and $G$ is a basis function. A function in $W_{G,X}$
has the form
\begin{equation}
\sum\limits_{i=1}^{N}\phi_{i}G\left(  x-x^{\left(  i\right)  }\right)
+p:\phi_{i}\in\mathbb{C},\;P_{X}^{T}\phi=0,\text{\ }p\in P_{\theta
-1},\label{h44}%
\end{equation}

where $P_{X}$ is a \textit{unisolvency matrix} and $G\left(  x-x^{\left(
i\right)  }\right)  $ is termed a \textit{data-translated basis function}.
Since the interpolant will be shown to be in $W_{G,X}$ its
finite-dimensionality allows the matrix equations to be derived for the
$\phi_{i}$ and the coefficients of a basis for $P_{\theta-1}$. The
\textit{vector-valued evaluation operator} $\widetilde{\mathcal{E}}%
_{X}f=\left(  f\left(  x^{\left(  i\right)  }\right)  \right)  $ is
introduced. This operator and its adjoint under the Light norm will be
fundamental to solving the Exact smoothing problem.

In \textbf{Sections} \ref{Sect_soln_Exact_smth},
\ref{Sect_Rxx_in_terms_of_Gxx}, \ref{Sect_mat_from_rx} and
\ref{Sect_Exact_smth_mat_eqn} we will use this mathematical machinery to
define the Exact smoothing problem and then solve it using the Hilbert space
technique of \textit{orthogonal projection} based on the vector-valued
evaluation operator $\widetilde{\mathcal{E}}_{X}$ and its adjoint. The Exact
smoothing problem involves minimizing a functional over the space
$X_{w}^{\theta}$. More specifically, a special Hilbert space is constructed
based on the smoothing functional and then the Exact smoothing problem is
expressed as one of minimizing the distance between a point generated by the
dependent data and a hyperplane generated by the independent data. Standard
orthogonal projection results then yield the existence of a unique solution to
the Exact smoothing problem which lies in the finite dimensional subspace
$W_{G,X}$.

Using the unisolvency matrix and the reproducing kernel matrix a block matrix
equation is derived for the values of the Exact smoother on $X$, assuming the
Riesz representer basis $\left\{  R_{x^{\left(  i\right)  }}\right\}  $ for
$W_{G,X}$. Then using this equation and relationships between the basis
function matrix and the unisolvency matrix a block matrix equation is derived
for the Exact smoother, assuming the data-translated basis function basis for
$W_{G,X}$. The block matrix is constructed from a basis function matrix and a
unisolvency matrix.

The estimates for the order of pointwise convergence will use the Riesz
representer and the operators which were used to solve the smoothing problem.
We also employ a special lemma derived from the Lagrange interpolation theory
from Light and Wayne \cite{LightWayne98PowFunc}. The dependent data $y$ is
generated using the data functions in the data space $X_{w}^{\theta}$ so that
$y=\widetilde{\mathcal{E}}_{X}f$ for some $f\in X_{w}^{\theta}$. Given a
bounded open set $\Omega$, we are interested in the behavior of $\max
_{x\in\overline{\Omega}}\left\vert u_{I}\left(  x\right)  -f\left(  x\right)
\right\vert $ as some measure $h_{X}$ of the `density' of the data points $X$
increases. Following Light and Wayne \cite{LightWayne98PowFunc} we have used
the measure $h_{X}=\sup_{\omega\in\Omega}\operatorname*{dist}\left(
\omega,X\right)  $.

In Theorem \ref{Thm_Exact_smth_converg_at_x_1} we show that if the weight
function has order $\theta$ and smoothness parameter $\kappa$ then we have the
smoother error estimate%
\[
\left\vert s_{e}\left(  x\right)  -f_{d}\left(  x\right)  \right\vert
\leq\left\vert f_{d}\right\vert _{w,\theta}\left(  1+K_{\Omega,\theta}%
^{\prime}\right)  \left(  \sqrt{c_{G}}\left(  c_{\Omega,\theta}h_{X}\right)
^{\eta}+\sqrt{N_{X}\rho}\right)  ,\quad x\in\overline{\Omega},
\]

where $\eta=\min\left\{  \theta,\frac{1}{2}\left\lfloor \min2\kappa
\right\rfloor \right\}  $. This assumes $\sqrt{r_{x}\left(  x\right)  }$
satisfies the bound \ref{h103}. Here $c_{G},c_{\Omega,\theta},K_{\Omega
,\theta}^{\prime}$ are constants independent of $X$ and when $\rho=0$ we have
an estimate for the interpolant derived in Chapter \ref{Ch_Interpol}.

In \textbf{Section} \ref{Sect_improved_err_estim}, Theorem
\ref{Thm_Exact_smth_converg_at_x} demonstrates the slightly improved estimate
\begin{equation}
\left\vert s_{e}\left(  x\right)  -f_{d}\left(  x\right)  \right\vert
\leq\left\vert f_{d}\right\vert _{w,\theta}\left(  1+K_{\Omega,\theta}%
^{\prime}\right)  \left(  \sqrt{c_{G}}\left(  c_{\Omega,\theta}h_{X}\right)
^{\eta+\delta_{G}}+\sqrt{N\rho}\right)  ,\quad x\in\overline{\Omega
},\label{h90}%
\end{equation}

and some $\delta_{G}\geq0$. This result assumes that $\sqrt{r_{x}\left(
x\right)  }$ satisfies the weaker bound \ref{h102} which is \ref{h103} with an
additional exponent $\delta_{G}$. This result is illustrated using the radial
thin-plate and surface thin-plate splines.

Unfortunately, for given $h_{X}$, $N=\left\vert X\right\vert $ can be
arbitrarily large so to obtain an \textit{order of pointwise convergence} a
relationship between $N$ and $h_{X}$ is obtained by numerically constructing
several $X$ in $\left[  -1.5,1.5\right]  \subset\mathbb{R}^{1}$ using a
uniform statistical distribution and then using a least squares fit to obtain
$h_{X}\simeq3.09N^{-0.81}$. Substituting for $N$ in \ref{h90} and minimizing
the right side using $\rho$ we show that given $h_{X_{k}}\rightarrow0$ there
is a sequence $\rho_{k}\rightarrow0$ such that $\left\vert s_{e}^{\left(
k\right)  }\left(  x\right)  -f_{d}\left(  x\right)  \right\vert \leq
C_{1}\left(  h_{X_{k}}\right)  ^{\eta+\delta_{G}}\rightarrow0$. We say the
order of convergence is at least $\left(  h_{X_{k}}\right)  ^{\eta+\delta_{G}%
}$.

For special data functions that are linear combinations of Riesz representers
$R_{x_{i}^{\prime}}$ a \textit{doubled order of convergence} of $2\eta
+2\delta_{G}$ is demonstrated.

I have not included the results of any numerical experiments concerning the
smoother error estimates derived in this chapter.

\section{Preparation\label{Sect_ex_prepare}}

We will require the following results from previous chapters concerning the
data spaces $X_{w}^{\theta}$, the Riesz representer $R_{x}$ and the spaces
$W_{G,X}$.

\subsection{The data spaces $X_{w}^{\theta}$\label{SbSect_ex_Xwm}}

\begin{summary}
\label{Sum_ex_property_Xwm}Suppose the weight function $w$ has property W2. If
$f\in X_{w}^{\theta}$ then $\widehat{f}\in L_{loc}^{1}\left(  \mathbb{R}%
^{d}\setminus0\right)  $ and we can define the function $f_{F}:\mathbb{R}%
^{d}\rightarrow\mathbb{C}$ a.e. by: $f_{F}=\widehat{f}$ on $\mathbb{R}%
^{d}\setminus0$. Further:

\begin{enumerate}
\item The seminorm \ref{p61} of $X_{w}^{\theta}$ satisfies (Theorem
\ref{Thm_properties_Xwm})%
\begin{equation}
\int w\left\vert \cdot\right\vert ^{2\theta}\left\vert f_{F}\right\vert
^{2}=\left\vert f\right\vert _{w,\theta}^{2}.\label{h81}%
\end{equation}

\item An alternative definition of $X_{w}^{\theta}$ is (Theorem
\ref{Thm_properties_Xwm}):
\begin{equation}
X_{w}^{\theta}=\left\{  f\in S^{\prime}:\xi^{\alpha}\widehat{f}\in L_{loc}%
^{1}\text{ }if\text{ }\left\vert \alpha\right\vert =\theta;\text{\thinspace
}\int w\left\vert \cdot\right\vert ^{2\theta}\left\vert f_{F}\right\vert
^{2}<\infty\right\}  .\label{h65}%
\end{equation}

Note that the condition $\widehat{f}\in L_{loc}^{1}\left(  \mathbb{R}%
^{d}\setminus0\right)  $ is implied by $\xi^{\alpha}\widehat{f}\in L_{loc}%
^{1}$ for all $\left\vert \alpha\right\vert =\theta$.\medskip

\item The functional $\left\vert \cdot\right\vert _{w,\theta}$ is a seminorm
and $\operatorname*{null}\left\vert \cdot\right\vert _{w,\theta}=P_{\theta-1}%
$. Also $P\cap X_{w}^{\theta}=P_{\theta-1}$ (Theorem \ref{Thm_properties_Xwm}).

\item $X_{w}^{\theta}$ is complete in the seminorm sense for all orders
$\theta\geq1$ (Theorem \ref{Thm_Xw_complete}).

\item Suppose $w$ also has property W3 for order $\theta$ and $\kappa$. Then
$X_{w}^{\theta}\subset C_{BP}^{\left(  \left\lfloor \kappa\right\rfloor
\right)  }$ (Theorem \ref{Thm_Xwth_W3_smooth}).
\end{enumerate}
\end{summary}

\subsection{The Riesz representer $R_{x}$ of $u\rightarrow u\left(  x\right)
$\label{SbSect_ex_Riesz_rep}}

This summary combines results from Section \ref{Sect_Riesz_rep}.

\begin{summary}
\label{Sum_ex_Riesz_rep_property}\textbf{Properties of the Riesz representer}
$R_{x}$ \textbf{of the evaluation functional} $u\rightarrow u\left(  x\right)
$ \textbf{on} $X_{w}^{\theta}$\textbf{:}

\begin{enumerate}
\item Given any minimal $\theta$-unisolvent set of points $A=\left\{
a^{\left(  i\right)  }\right\}  _{i=1}^{M}$ we endow the semi-Hilbert space
$X_{w}^{\theta}$ with the Light inner product $\left(  \cdot,\cdot\right)
_{w,\theta}$ defined using the cardinal basis $\left\{  l_{i}\right\}
_{i=1}^{M}$ that corresponds to $A$. Then the representer is given by
\ref{p939} where $G$ is any basis function of order $\theta$ generated by the
weight function $w$ with properties W2 and W3 for order $\theta$.

\item For all $u\in X_{w}^{\theta}$, $u\left(  x\right)  =\left(
u,R_{x}\right)  _{w,\theta}$.

\item $R_{x}$ is unique and is independent of the basis function used to
define it.

\item $R_{x}\left(  y\right)  =\overline{R_{y}\left(  x\right)  }$.

\item $R_{x}\left(  a^{\left(  i\right)  }\right)  =l_{i}\left(  x\right)  $.

\item $\mathcal{Q}u\left(  x\right)  =\left\langle u,R_{x}\right\rangle
_{w,\theta}$ when $u\in X_{w}^{\theta}$.

\item $\mathcal{P}R_{x}=\sum\limits_{i=1}^{M}l_{i}\left(  x\right)  l_{i}$.
\end{enumerate}
\end{summary}

\subsection{The basis function spaces $\protect\overset{\cdot}{W}_{G,X}$ and
$W_{G,X}$\label{SbSect_ex_Wgx}}

The importance of the finite dimensional spaces $\overset{\cdot}{W}_{G,X}$ and
$W_{G,X}$ is that they contain the solutions to the variational interpolation
and smoothing problems studied in this document.

\begin{definition}
\label{Def_ex_Wgx}(Copy of Definition \ref{Def_Wgx})\textbf{The basis function
spaces }$W_{G,X}$\textbf{\ and }$\overset{\cdot}{W}_{G,X} $

Suppose the weight function $w$ has properties W2 and W3 for order $\theta
\geq1$ and $\kappa$. Then the basis distributions of order $\theta$ are
continuous functions and we let $G$ be a basis function. Let $X=\left\{
x^{(i)}\right\}  _{i=1}^{N}$ be a $\theta$-unisolvent set of distinct points
in $\mathbb{R}^{d}$ and set $M=\dim P_{\theta-1}$. Next choose a real-valued
basis $\left\{  p_{j}\right\}  _{j=1}^{M}$ of $P_{\theta-1}$ and calculate the
unisolvency matrix $P_{X}=\left(  p_{j}\left(  x^{(i)}\right)  \right)  $. We
can now define
\begin{align}
\overset{\cdot}{W}_{G,X} &  =\left\{  \sum_{i=1}^{N}v_{i}G\left(
x-x^{(i)}\right)  :\left(  v_{i}\right)  \in\mathbb{C}^{N}\text{ }and\text{
}P_{X}^{T}v=0\right\}  ,\nonumber\\
W_{G,X} &  =\overset{\cdot}{W}_{G,X}+P_{\theta-1}.\nonumber
\end{align}

\end{definition}

\begin{summary}
\label{Sum_ex_Wgx_properties_1}From Section \ref{Sect_Wgx} we know that set-wise:

\begin{enumerate}
\item $W_{G,X}$ is independent of the basis function of order $\theta$ used to
define it (Theorem \ref{Thm_Wgx_reorder}).

\item $\overset{\cdot}{W}_{G,X}$ and $W_{G,X}$ are independent of the basis
used to define $P_{\theta-1}$ (Theorem \ref{Thm_Wgx_defin_indep_of_Px_chosen}).

\item $W_{G,X}$ is independent of the ordering of the points in $X$ (Theorem
\ref{Thm_Wgx_reorder}).
\end{enumerate}
\end{summary}

The next theorem restates more results from Section \ref{Sect_Wgx}.

\begin{remark}
\begin{summary}
\label{Sum_ex_Wgx_properties_2}The spaces $\overset{\cdot}{W}_{G,X}$ and
$W_{G,X}$ have the following properties: from Theorem
\ref{Thm_Wgx_properties_1}:

\begin{enumerate}
\item If $f_{v}\left(  x\right)  =\sum\limits_{k=1}^{N}v_{k}G\left(
x-x^{\left(  k\right)  }\right)  $ then $\left\vert f_{v}\right\vert
_{w,\theta}^{2}=\left(  2\pi\right)  ^{\frac{d}{2}}v^{T}G_{X,X}\overline{v}$.

\item $G_{X,X}$ is conditionally positive definite on $\operatorname*{null}%
P_{X}^{T}$ i.e. when $P_{X}^{T}v=0$ and $v\neq0$ we have $v^{T}G_{X,X}%
\overline{v}>0$.

\item $W_{G,X}=\overset{\cdot}{W}_{G,X}\oplus P_{\theta-1}$, $\dim
\overset{\cdot}{W}_{G,X}=N-M$ and $\dim W_{G,X}=N$.

\item $X_{w}^{\theta}=W_{G,X}\oplus W_{G,X}^{\perp}$ where%
\[
W_{G,X}^{\perp}=\left\{  u\in X_{w}^{\theta}:u\left(  x^{\left(  k\right)
}\right)  =0\text{ }for\text{ }all\text{ }x^{\left(  k\right)  }\in X\right\}
.
\]
\medskip

From Corollary \ref{Cor_uniq_rep_for_elt_of_Wgx}:\smallskip

\item If $\left\{  p_{j}\right\}  _{j=1}^{M}$ is basis for $P_{\theta-1}$ then
the representation
\begin{equation}
W_{G,X}=\left\{  \sum\limits_{i=1}^{N}\alpha_{i}G\left(  \cdot-x^{\left(
i\right)  }\right)  +\sum\limits_{j=1}^{M}\beta_{j}p_{j}:P_{X}^{T}%
\alpha=0,\text{ }\alpha=\left(  \alpha_{i}\right)  ,\text{ }\alpha_{i}%
,\beta_{j}\in\mathbb{C}\right\}  ,\label{h14}%
\end{equation}

is unique in terms of $\alpha_{i}$ and $\beta_{j}$.\medskip

From Corollary \ref{Cor_Thm_Wgx_properties_2}:\smallskip

\item Suppose $A\subset X$ is a minimal unisolvent set and suppose $R_{x}$ is
the Riesz representer of the functional $f\rightarrow f\left(  x\right)  $
w.r.t. the Light norm, both being defined using $A$. Then $\overset{\cdot
}{W}_{G,X}$ has basis $\left\{  R_{x^{\left(  i\right)  }}:x^{\left(
i\right)  }\notin A\right\}  $ and $W_{G,X}$ has basis $\left\{  R_{x^{\left(
i\right)  }}\right\}  _{i=1}^{N}$.
\end{enumerate}
\end{summary}
\end{remark}

\section{Solution of the Exact smoothing problem and its
properties\label{Sect_soln_Exact_smth}}

The Exact smoothing problem involves minimizing a functional over the space of
continuous functions $X_{w}^{\theta}$. The functional is constructed using a
smoothing parameter $\rho>0$, the seminorm $\left\vert \cdot\right\vert
_{w,\theta}$ of $X_{w}^{\theta}$ and a set of distinct, scattered data points
$\left\{  \left(  x^{(i)},y_{i}\right)  \right\}  _{i=1}^{N}$, $x^{(i)}%
\in\mathbb{R}^{d}$ and $y_{i}\in\mathbb{R}$. As with the interpolation problem
we require the $x^{(i)}$ to be distinct. The data will sometimes be specified
using the notation $\left[  X,y\right]  $ where $X=\left\{  x^{(i)}\right\}
_{i=1}^{N}$ is termed the independent data and $y=\left\{  y_{i}\right\}
_{i=1}^{N}$ is termed the dependent data. The functional to be minimized is%
\begin{equation}
J_{e}\left[  f\right]  =\rho\left\vert f\right\vert _{w,\theta}^{2}+\frac
{1}{N}\sum_{i=1}^{N}\left\vert f\left(  x^{(i)}\right)  -y_{i}\right\vert
^{2},\quad f\in X_{w}^{\theta}.\label{h52}%
\end{equation}

The Exact smoothing problem is now stated as:

\begin{definition}
\label{Def_Exact_smth_problem}\textbf{The Exact smoothing problem and the
Exact smoother\medskip}

\fbox{Minimize$\text{ }$the$\text{ Exact }$smoothing$\text{ }$%
functional\ $J_{e}\text{ on the space }X_{w}^{\theta}$.}\textbf{\medskip}

A solution to this problem will be called an \textbf{Exact smoother} of the
data $\left[  X,y\right]  $\textbf{.}
\end{definition}

The first component $\rho\left\vert f\right\vert _{w,\theta}^{2}$ of the
functional can be regarded as a global smoothing component and the second
component as the localizing least squares component.

It is clear that when $\rho=0$ any interpolant of the data minimizes the Exact
smoother functional. Thus, from part 4 Theorem \ref{Sum_ex_Wgx_properties_2}
the set of solutions is $u_{I}+W_{G,X}^{\perp}$ where $u_{I}$ is the minimal
seminorm interpolant studied in Chapter \ref{Ch_Interpol}. In Remark
\ref{Rem_Thm_smooth_matrix_soln_1} we will note that as $\rho\rightarrow0$ the
matrix equation representing the solution to this problem becomes that of the
minimal norm interpolation problem so the limit "seeks out" one particular interpolant.

\subsection{Existence, uniqueness and formulas for the smoother}

Using the technique of orthogonal projection it will be shown below that a
solution to the smoothing problem exists and is unique. The proof will be
carried out within a Hilbert space framework by formulating the smoothing
functional in terms of a special inner product on the Hilbert product space
$V=X_{w}^{\theta}\otimes\mathbb{C}^{N}$. To this end I will introduce the
following definitions:

\begin{definition}
\label{Def_ex_Hilbert_smoothing}For order $\theta\geq1$:

\begin{enumerate}
\item Suppose $\rho>0$. Let $V=X_{w}^{\theta}\otimes\mathbb{C}^{N}$ be the
Hilbert product space with inner product
\[
\left(  \left(  u_{1},\widetilde{u}_{2}\right)  ,\left(  v_{1},\widetilde
{v}_{2}\right)  \right)  _{V}=\rho\left\langle u_{1},v_{1}\right\rangle
_{w,\theta}+\frac{1}{N}\left(  \widetilde{u}_{2},\widetilde{v}_{2}\right)
_{\mathbb{C}^{N}},
\]

and $\left\Vert f\right\Vert _{V}^{2}=\left(  f,f\right)  _{V}$.

\item Let $\mathcal{L}_{X}:X_{w}^{\theta}\rightarrow V$ be defined by
$\mathcal{L}_{X}f=\left(  f,\widetilde{\mathcal{E}}_{X}f\right)  $ where
$\widetilde{\mathcal{E}}_{X}$ is the evaluation operator of Definition
\ref{Def_eval_operators}.
\end{enumerate}
\end{definition}

\begin{remark}
\label{Rem_SmoothFunc1}The smoothing functional can be expressed in terms of
$\mathcal{L}_{X}$ and the data as follows: set $\varsigma=\left(  0,y\right)
\in V$ where $y=\left(  y_{i}\right)  $ is the dependent data given in the
Exact smoothing problem. Then for $f\in X_{w}^{\theta}$%
\begin{align*}
\left\Vert \mathcal{L}_{X}f-\varsigma\right\Vert _{V}^{2}=\left\Vert \left(
f,\widetilde{\mathcal{E}}_{X}f\right)  -\left(  0,y\right)  \right\Vert
_{V}^{2}=\left\Vert \left(  f,\widetilde{\mathcal{E}}_{X}f-y\right)
\right\Vert _{V}^{2} &  =\rho\left\vert f\right\vert _{w,\theta}^{2}+\frac
{1}{N}\left\vert \widetilde{\mathcal{E}}_{X}f-y\right\vert ^{2}\\
&  =J_{e}\left[  f\right]  .
\end{align*}

\end{remark}

The operator $\mathcal{L}_{X}:X_{w}^{\theta}\rightarrow V$ and its adjoint
$\mathcal{L}_{X}^{\ast}:V\rightarrow X_{w}^{\theta}$ have the following properties:

\begin{theorem}
\label{Thm_L_op_properties}Suppose $X$ is a $\theta$-unisolvent set and
$A\subset X$ is a minimal unisolvent set. Use $A$ to define the Light norm
$\left\Vert \cdot\right\Vert _{w,\theta}$ and the Lagrange interpolation
operator $\mathcal{P}$. Then when $\rho>0$:

\begin{enumerate}
\item $\left\Vert \mathcal{L}_{X}f\right\Vert _{V}$ and $\left\Vert
f\right\Vert _{w,\theta}$ are equivalent norms on $X_{w}^{\theta}$.

\item $\mathcal{L}_{X}:X_{w}^{\theta}\rightarrow V$ is continuous and 1-1.

\item $\mathcal{L}_{X}^{\ast}:V\rightarrow X_{w}^{\theta}$ is continuous
w.r.t. the Light norm and if $u=\left(  u_{1,}\widetilde{u}_{2}\right)  \in V$
then%
\[
\mathcal{L}_{X}^{\ast}u=\rho\mathcal{Q}u_{1}+\frac{1}{N}\widetilde
{\mathcal{E}}_{X}^{\ast}\widetilde{u}_{2},
\]

with $\operatorname*{range}\mathcal{L}_{X}^{\ast}=X_{w}^{\theta}$.

\item $\mathcal{L}_{X}^{\ast}\mathcal{L}_{X}:X_{w}^{\theta}\rightarrow
X_{w}^{\theta}$ and
\[
\mathcal{L}_{X}^{\ast}\mathcal{L}_{X}=\rho\mathcal{Q}+\frac{1}{N}%
\widetilde{\mathcal{E}}_{X}^{\ast}\widetilde{\mathcal{E}}_{X}.
\]

\item $\mathcal{L}_{X}^{\ast}\mathcal{L}_{X}$ is 1-1 on $X_{w}^{\theta}$.
Also, $\left(  \mathcal{L}_{X}^{\ast}\mathcal{L}_{X}\right)  ^{-1}$ is
continuous on $\operatorname*{range}\mathcal{L}_{X}^{\ast}\mathcal{L}_{X}$.
\end{enumerate}
\end{theorem}

\begin{proof}
\textbf{Part 1 }From the definition of $\mathcal{L}_{X}$, if $f\in
X_{w}^{\theta}$
\[
\left\Vert \mathcal{L}_{X}f\right\Vert _{V}^{2}=\left\Vert \left(
f,\widetilde{\mathcal{E}}_{X}f\right)  \right\Vert _{V}^{2}=\rho\left\vert
f\right\vert _{w,\theta}^{2}+\frac{1}{N}\left\vert \widetilde{\mathcal{E}}%
_{X}f\right\vert _{\mathbb{C}^{N}}^{2}\leq const\left\Vert f\right\Vert
_{w,\theta}^{2},
\]

since $\widetilde{\mathcal{E}}_{X}:X_{w}^{\theta}\rightarrow\mathbb{R}^{N}$ is
continuous by part 1 of Theorem \ref{Thm_eval_op_properties}. Also
\begin{align*}
\left\Vert f\right\Vert _{w,\theta}^{2}=\left\vert f\right\vert _{w,\theta
}^{2}+\left\vert \widetilde{\mathcal{E}}_{A}f\right\vert _{\mathbb{C}^{M}}^{2}
&  \leq\left(  \min\left\{  \rho,\frac{1}{N}\right\}  \right)  ^{-1}\left(
\rho\left\vert f\right\vert _{w,\theta}^{2}+\frac{1}{N}\left\vert
\widetilde{\mathcal{E}}_{X}f\right\vert _{\mathbb{C}^{N}}^{2}\right) \\
&  =\left(  \min\left\{  \rho,\frac{1}{N}\right\}  \right)  ^{-1}\left\Vert
\mathcal{L}_{X}f\right\Vert _{V}^{2},
\end{align*}

since $\rho>0$.\medskip

\textbf{Part 2} Part 1 implies $\mathcal{L}_{X}$ is continuous and that
$\mathcal{L}_{X}$ is 1-1 is clear from its definition.\medskip

\textbf{Part 3 }Since $\mathcal{L}_{X}$ is a continuous operator the adjoint
$\mathcal{L}_{X}^{\ast}:V\rightarrow X_{w}^{\theta}$ exists and is continuous.
Further, if $f\in X_{w}^{\theta}$ and $u=\left(  u_{1,}\widetilde{u}%
_{2}\right)  \in V$ then by the properties of the Light norm given in Theorem
\ref{Thm_Light_norm_property}
\begin{align*}
\left(  \mathcal{L}_{X}f,u\right)  _{V}=\left(  \left(  f,\widetilde
{\mathcal{E}}_{X}f\right)  ,\left(  u_{1},\widetilde{u}_{2}\right)  \right)
_{V} &  =\rho\,\left\langle f,u_{1}\right\rangle _{w,\theta}+\frac{1}%
{N}\left(  \widetilde{\mathcal{E}}_{X}f,\widetilde{u}_{2}\right)
_{\mathbb{R}^{N}}\\
&  =\rho\,\left(  f,\mathcal{Q}u_{1}\right)  _{w,\theta}+\frac{1}{N}\left(
f,\widetilde{\mathcal{E}}_{X}^{\ast}\widetilde{u}_{2}\right)  _{w,\theta}\\
&  =\left(  f,\rho\mathcal{Q}u_{1}+\frac{1}{N}\widetilde{\mathcal{E}}%
_{X}^{\ast}\widetilde{u}_{2}\right)  _{w,\theta}.
\end{align*}

Hence $\mathcal{L}_{X}^{\ast}u=\rho\mathcal{Q}u_{1}+\frac{1}{N}\widetilde
{\mathcal{E}}_{X}^{\ast}\widetilde{u}_{2}$.\medskip

\textbf{Part 4}
\[
\mathcal{L}_{X}^{\ast}\mathcal{L}_{X}f=\rho\mathcal{Q}\left(  \left(
\mathcal{L}_{X}f\right)  _{1}\right)  +\frac{1}{N}\widetilde{\mathcal{E}}%
_{X}^{\ast}\left(  \left(  \mathcal{L}_{X}f\right)  _{2}\right)
=\rho\mathcal{Q}f+\frac{1}{N}\widetilde{\mathcal{E}}_{X}^{\ast}\widetilde
{\mathcal{E}}_{X}f.
\]
\medskip

\textbf{Part 5} Suppose $\mathcal{L}_{X}^{\ast}\mathcal{L}_{X}f=0$. Then
$0=\left(  \mathcal{L}_{X}^{\ast}\mathcal{L}_{X}f,f\right)  _{w,\theta
}=\left(  \mathcal{L}_{X}f,\mathcal{L}_{X}f\right)  _{w,\theta}=\left\Vert
\mathcal{L}_{X}f\right\Vert _{w,\theta}^{2}$ so that $\mathcal{L}_{X}f=0$ and
$f=0$ since $\mathcal{L}_{X}$ is 1-1. Hence $\mathcal{L}_{X}^{\ast}%
\mathcal{L}_{X}$ is 1-1.

Since $\mathcal{L}_{X}^{\ast}\mathcal{L}_{X}$ is continuous it has closed
range. By the bounded inverse theorem $\left(  \mathcal{L}_{X}^{\ast
}\mathcal{L}_{X}\right)  ^{-1}$ is also continuous on $\operatorname*{range}%
\mathcal{L}_{X}^{\ast}\mathcal{L}_{X}$.
\end{proof}

Using the Hilbert space technique of orthogonal projection the next theorem
shows that when $\rho>0$ the Exact smoothing problem of Definition
\ref{Def_Exact_smth_problem} has a unique solution in $X_{w}^{\theta}$.

\begin{theorem}
\label{Thm_smooth_Exact}Fix $y\in\mathbb{R}^{N}$ and let $\varsigma=\left(
0,y\right)  \in V$. Then for $\rho>0$ there exists a unique function $s_{e}\in
X_{w}^{\theta}$ which solves the Exact smoothing problem with data $\left[
X,y\right]  $. This solution has the following properties:

\begin{enumerate}
\item $\left\Vert \mathcal{L}_{X}s_{e}-\varsigma\right\Vert _{V}<\left\Vert
\mathcal{L}_{X}f-\varsigma\right\Vert _{V}$ for all $f\in X_{w}^{\theta
}-\{s_{e}\}$.

\item $\left(  \mathcal{L}_{X}s_{e}-\varsigma,\mathcal{L}_{X}s_{e}%
-\mathcal{L}_{X}f\right)  _{V}=0$ for all $f\in X_{w}^{\theta}$.

\item $\left\Vert \mathcal{L}_{X}s_{e}-\varsigma\right\Vert _{V}%
^{2}+\left\Vert \mathcal{L}_{X}s_{e}-\mathcal{L}_{X}f\right\Vert _{V}%
^{2}=\left\Vert \mathcal{L}_{X}f-\varsigma\right\Vert _{V}^{2}$ for all $f\in
X_{w}^{\theta}$.

This equality is equivalent to that of part 2.

\item $s_{e}=\frac{1}{N}\left(  \mathcal{L}_{X}^{\ast}\mathcal{L}_{X}\right)
^{-1}\widetilde{\mathcal{E}}_{X}^{\ast}y$.

\item Since $\theta\geq1$, $p\in P_{\theta-1}$ implies $p=\frac{1}{N}\left(
\mathcal{L}_{X}^{\ast}\mathcal{L}_{X}\right)  ^{-1}\widetilde{\mathcal{E}}%
_{X}^{\ast}\widetilde{\mathcal{E}}_{X}p$. We say the Exact smoother preserves
polynomials up to order $\theta$.
\end{enumerate}
\end{theorem}

\begin{proof}
\textbf{Parts 1},\textbf{2},\textbf{3 }Since $\mathcal{L}_{X}$ is continuous,
we know that the hyperspace $\mathcal{L}_{X}\left(  X_{w}^{\theta}\right)  $
is closed and hence that the translated hyperspace $\mathcal{L}_{X}\left(
X_{w}^{\theta}\right)  -\varsigma$ is also closed. But from the remark
following Definition \ref{Def_ex_Hilbert_smoothing} $J_{e}\left[  f\right]
=\left\Vert \mathcal{L}_{X}f-\varsigma\right\Vert _{V}^{2}$ and the Exact
smoothing problem becomes minimize $\left\Vert \mathcal{L}_{X}f-\varsigma
\right\Vert _{V}$ over $X_{w}^{\theta}$. Thus by a well-known orthogonal
projection result concerning the distance between a point and a closed
subspace, there exists a unique element of $\mathcal{L}_{X}(X_{w}^{\theta
})-\varsigma$, call it $\mathcal{L}_{X}s_{e}-\varsigma$ such that $s_{e}$
satisfies parts 1, 2 and 3 of this theorem.\medskip

\textbf{Part 4 }Using the equation proved in part 2
\[
0=\left(  \mathcal{L}_{X}s_{e}-\varsigma,\mathcal{L}_{X}s_{e}-\mathcal{L}%
_{X}f\right)  _{V}=\left(  \mathcal{L}_{X}s_{e}-\varsigma,\mathcal{L}%
_{X}\left(  s_{e}-f\right)  \right)  _{V}=\left(  \mathcal{L}_{X}^{\ast
}\left(  \mathcal{\mathcal{L}}_{X}s_{e}-\varsigma\right)  ,s_{e}-f\right)
_{w,\theta},
\]

for all $f\in X_{w}^{\theta}$. Thus
\begin{equation}
\mathcal{L}_{X}^{\ast}\left(  \mathcal{L}_{X}s_{e}-\varsigma\right)
=0,\label{h36}%
\end{equation}

and therefore
\[
\mathcal{L}_{X}^{\ast}\mathcal{L}_{X}s_{e}=\mathcal{L}_{X}^{\ast}%
\varsigma=\mathcal{L}_{X}^{\ast}\left(  0,y\right)  =\frac{1}{N}%
\widetilde{\mathcal{E}}_{X}^{\ast}y.
\]

But by part 5 of Theorem \ref{Thm_L_op_properties}, $\mathcal{L}_{X}^{\ast
}\mathcal{L}_{X}$ is one-to-one and so $s_{e}=\frac{1}{N}\left(
\mathcal{L}_{X}^{\ast}\mathcal{L}_{X}\right)  ^{-1}\widetilde{\mathcal{E}}%
_{X}^{\ast}y$.\medskip

\textbf{Part 5} Substituting $y=\widetilde{\mathcal{E}}_{X}p$ into the result
proved in the previous part gives

$s_{e}=\frac{1}{N}\left(  \mathcal{L}_{X}^{\ast}\mathcal{L}_{X}\right)
^{-1}\widetilde{\mathcal{E}}_{X}^{\ast}\widetilde{\mathcal{E}}_{X}p$, or
equivalently $N\mathcal{L}_{X}^{\ast}\mathcal{L}_{X}s_{e}=\widetilde
{\mathcal{E}}_{X}^{\ast}\widetilde{\mathcal{E}}_{X}p$. But from part 4 of
Theorem \ref{Thm_L_op_properties}, $\mathcal{L}_{X}^{\ast}\mathcal{L}_{X}%
-\rho\mathcal{Q}=\frac{1}{N}\widetilde{\mathcal{E}}_{X}^{\ast}\widetilde
{\mathcal{E}}_{X}$.

Substituting this as an expression for $\widetilde{\mathcal{E}}_{X}^{\ast
}\widetilde{\mathcal{E}}_{X}$ into the last equation we get
\[
N\mathcal{L}_{X}^{\ast}\mathcal{L}_{X}s_{e}=N\left(  \mathcal{L}_{X}^{\ast
}\mathcal{L}_{X}p-\rho\mathcal{Q}p\right)  =N\mathcal{L}_{X}^{\ast}%
\mathcal{L}_{X}p,
\]

since $\mathcal{Q}p=\left(  I-\mathcal{P}\right)  p=0$. Finally, by part 5 of
Theorem \ref{Thm_L_op_properties}, $\mathcal{L}_{X}^{\ast}\mathcal{L}_{X}$ is
one-to-one and hence $s_{e}=p$.
\end{proof}

\subsection{Various identities}

The last theorem showed that the Exact smoothing problem has a unique
minimizer in $X_{w}^{\theta}$. In the next corollary we prove that the
smoother lies in the space $W_{G,X}$ introduced in Subsection
\ref{SbSect_ex_Wgx}. In the next two corollaries we will also prove several
well-known identities which involve all the data points. These identities
relate the Hilbert space properties and the pointwise properties of the data
and the Exact smoother.

\begin{corollary}
\label{Cor_propert_Exact_smth}Suppose $X=\left\{  x^{(k)}\right\}  _{k=1}^{N}$
is $\theta$-unisolvent so that $X$ contains a minimal unisolvent subset $A$.
Use this subset to define the operators $\mathcal{P}$, $\mathcal{Q}$ and the
Light norm. Then the unique solution $s_{e}\in X_{w}^{\theta}$ of the Exact
smoothing problem with data $\left[  X,y\right]  $ has the following properties:

\begin{enumerate}
\item $s_{e}\in W_{G,X}$ and $s_{e}=\mathcal{P}s_{e}-\frac{1}{N\rho}%
\widetilde{\mathcal{E}}_{X}^{\ast}\left(  \widetilde{\mathcal{E}}_{X}%
s_{e}-y\right)  $.

\item For all $f\in X_{w}^{\theta}$ we have
\begin{gather}
\rho\left\vert s_{e}\right\vert _{w,\theta}^{2}+\frac{1}{N}\sum_{k=1}%
^{N}\left\vert s_{e}\left(  x^{(k)}\right)  -y_{k}\right\vert ^{2}%
+\rho\left\vert s_{e}-f\right\vert _{w,\theta}^{2}+\frac{1}{N}\sum_{k=1}%
^{N}\left\vert s_{e}\left(  x^{(k)}\right)  -f\left(  x^{(k)}\right)
\right\vert ^{2}\label{h32}\\
=\rho\left\vert f\right\vert _{w,\theta}^{2}+\frac{1}{N}\sum_{k=1}%
^{N}\left\vert f\left(  x^{(k)}\right)  -y_{k}\right\vert ^{2},\nonumber
\end{gather}

or%
\begin{equation}
J_{e}\left[  s_{e}\right]  +\rho\left\vert s_{e}-f\right\vert _{w,\theta}%
^{2}+\frac{1}{N}\sum_{k=1}^{N}\left\vert s_{e}\left(  x^{(k)}\right)
-f\left(  x^{(k)}\right)  \right\vert ^{2}=J_{e}\left[  f\right]  .\label{h33}%
\end{equation}

\item $P_{X}^{T}\left(  \left(  s_{e}\right)  _{X}-y\right)  =0$.
\end{enumerate}
\end{corollary}

\begin{proof}
\textbf{Part 1}
\begin{align*}
0=\mathcal{L}_{X}^{\ast}\left(  \mathcal{L}_{X}s_{e}-\varsigma\right)
=\mathcal{L}_{X}^{\ast}\mathcal{L}_{X}s_{e}-\frac{1}{N}\widetilde{\mathcal{E}%
}_{X}^{\ast}y &  =\rho\mathcal{Q}s_{e}+\frac{1}{N}\widetilde{\mathcal{E}}%
_{X}^{\ast}\widetilde{\mathcal{E}}_{X}s_{e}-\frac{1}{N}\widetilde{\mathcal{E}%
}_{X}^{\ast}y\\
&  =\rho s_{e}-\rho\mathcal{P}s_{e}+\frac{1}{N}\widetilde{\mathcal{E}}%
_{X}^{\ast}\left(  \widetilde{\mathcal{E}}_{X}s_{e}-y\right)  ,
\end{align*}

so that%
\[
s_{e}=\mathcal{P}s_{e}-\frac{1}{N\rho}\widetilde{\mathcal{E}}_{X}^{\ast
}\left(  \widetilde{\mathcal{E}}_{X}s_{e}-y\right)  ,
\]

and $s_{e}\in W_{G,X}$ by part 2 of Theorem \ref{Thm_eval_op_properties}%
.\medskip

\textbf{Part 2} By part 3 of Theorem \ref{Thm_smooth_Exact}
\[
\left\Vert \mathcal{L}_{X}s_{e}-\varsigma\right\Vert _{V}^{2}+\left\Vert
\mathcal{L}_{X}s_{e}-\mathcal{L}_{X}f\right\Vert _{V}^{2}=\left\Vert
\mathcal{L}_{X}f-\varsigma\right\Vert _{V}^{2},\quad f\in X_{w}^{\theta}.
\]

Using Remark \ref{Rem_SmoothFunc1} and the definition of $\mathcal{L}_{X}$
this equation becomes
\[
\rho\left\vert s_{e}\right\vert _{w,\theta}^{2}+\frac{1}{N}\left\vert
\widetilde{\mathcal{E}}_{X}s_{e}-y\right\vert ^{2}+\rho\left\vert
s_{e}-f\right\vert _{w,\theta}^{2}+\frac{1}{N}\left\vert \widetilde
{\mathcal{E}}_{X}\left(  s_{e}-f\right)  \right\vert ^{2}=\rho\left\vert
f\right\vert _{w,\theta}^{2}+\frac{1}{N}\left\vert \widetilde{\mathcal{E}}%
_{X}f-y\right\vert ^{2},
\]

for all $f\in X_{w}^{\theta}$ and $y\in\mathbb{R}^{N}$, which proves this
part.\medskip

\textbf{Part 3} If $p\in P_{\theta-1}$, then from part 1
\begin{align*}
\frac{1}{N\rho}\left(  \widetilde{\mathcal{E}}_{X}^{\ast}\left(  \left(
s_{e}\right)  _{X}-y\right)  ,p\right)  _{w,\theta}=\frac{1}{N\rho}\left(
\widetilde{\mathcal{E}}_{X}^{\ast}\left(  \widetilde{\mathcal{E}}_{X}%
s_{e}-y\right)  ,p\right)  _{w,\theta} &  =\left(  \mathcal{P}s_{e}%
-s_{e},p\right)  _{w,\theta}\\
&  =-\left(  \mathcal{Q}s_{e},p\right)  _{w,\theta}\\
&  =-\left\langle s_{e},p\right\rangle _{w,\theta}\\
&  =0.
\end{align*}

But
\[
0=\left(  \widetilde{\mathcal{E}}_{X}^{\ast}\left(  \left(  s_{e}\right)
_{X}-y\right)  ,p\right)  _{w,\theta}=\left(  \left(  s_{e}\right)
_{X}-y,\widetilde{\mathcal{E}}_{X}p\right)  =\sum_{i=1}^{N}\left(  \left(
s_{e}\right)  _{X}-y\right)  _{i}\overline{p\left(  x^{\left(  i\right)
}\right)  },
\]

so by part 2 of Theorem \ref{Thm_Px_properties}, $P_{X}^{T}\left(  \left(
s_{e}\right)  _{X}-y\right)  =0$.
\end{proof}

We next prove some standard results which express the smoothing functional and
the seminorm of the smoother in terms of the dependent data and the values
taken by the smoother on the independent data.

\begin{corollary}
\label{Cor_propert_Exact_smth_2}Suppose $\left\{  \left(  x^{(k)}%
,y_{k}\right)  \right\}  _{k=1}^{N}$ is the data to be smoothed. Then the
Exact smoothing function $s_{e}$ has the following properties:

\begin{enumerate}
\item For all $p\in P_{\theta-1}$
\[
2\rho\left\vert s_{e}\right\vert _{w,\theta}^{2}+\frac{1}{N}\sum_{k=1}%
^{N}\left\vert s_{e}\left(  x^{(k)}\right)  -y_{k}\right\vert ^{2}+\frac{1}%
{N}\sum_{k=1}^{N}\left\vert s_{e}\left(  x^{(k)}\right)  -p\left(
x^{(k)}\right)  \right\vert ^{2}=\frac{1}{N}\sum_{k=1}^{N}\left\vert p\left(
x^{(k)}\right)  -y_{k}\right\vert ^{2}.
\]

\item $2\rho\left\vert s_{e}\right\vert _{w,\theta}^{2}+\frac{1}{N}%
\sum\limits_{k=1}^{N}\left\vert s_{e}\left(  x^{(k)}\right)  -y_{k}\right\vert
^{2}+\frac{1}{N}\sum\limits_{k=1}^{N}\left\vert s_{e}\left(  x^{(k)}\right)
\right\vert ^{2}=\frac{1}{N}\sum\limits_{k=1}^{N}\left\vert y_{k}\right\vert
^{2}.$

\item $\left\vert s_{e}\right\vert _{w,\theta}^{2}=\frac{1}{N\rho
}\operatorname{Re}\sum\limits_{k=1}^{N}\overline{s_{e}\left(  x^{(k)}\right)
}\left(  y_{k}-s_{e}\left(  x^{(k)}\right)  \right)  .$

\item $J_{e}\left[  s_{e}\right]  =\frac{1}{N}\operatorname{Re}\sum
\limits_{k=1}^{N}\left(  y_{k}-s_{e}\left(  x^{(k)}\right)  \right)
\overline{y_{k}}.$
\end{enumerate}
\end{corollary}

\begin{proof}
\textbf{Part 1} Substitute $f=p\in P_{\theta-1}$ in the equation of part 2 of
Corollary \ref{Cor_propert_Exact_smth} and use the fact that $\left\vert
g+p\right\vert _{w,\theta}=\left\vert g\right\vert _{w,\theta}$ when $g\in
X_{w}^{\theta}$.\medskip

\textbf{Part 2} Substitute $f=0$ in the equation of part 2 of Corollary
\ref{Cor_propert_Exact_smth}.\medskip

\textbf{Part 3} Substitute the expansion
\[
\left\vert s_{e}\left(  x^{(k)}\right)  -y_{k}\right\vert ^{2}=\left\vert
s_{e}\left(  x^{(k)}\right)  \right\vert ^{2}-2\operatorname{Re}s_{e}\left(
x^{(k)}\right)  \overline{y_{k}}+\left\vert y_{k}\right\vert ^{2},
\]

in part 2 of Corollary \ref{Cor_propert_Exact_smth} so that%
\begin{align*}
2\rho\left\vert s_{e}\right\vert _{w,\theta}^{2}  & =\frac{1}{N}\sum_{k=1}%
^{N}\left\vert y_{k}\right\vert ^{2}-\frac{1}{N}\sum_{k=1}^{N}\left\vert
s_{e}\left(  x^{(k)}\right)  -y_{k}\right\vert ^{2}-\frac{1}{N}\sum_{k=1}%
^{N}\left\vert s_{e}\left(  x^{(k)}\right)  \right\vert ^{2}\\
& =\frac{1}{N}\sum_{k=1}^{N}\left\vert y_{k}\right\vert ^{2}-\frac{1}{N}%
\sum_{k=1}^{N}\left(  \left\vert s_{e}\left(  x^{(k)}\right)  \right\vert
^{2}-2\operatorname{Re}s_{e}\left(  x^{(k)}\right)  \overline{y_{k}%
}+\left\vert y_{k}\right\vert ^{2}\right) \\
& \qquad\qquad-\frac{1}{N}\sum_{k=1}^{N}\left\vert s_{e}\left(  x^{(k)}%
\right)  \right\vert ^{2}\\
& =\frac{2}{N}\sum_{k=1}^{N}\operatorname{Re}s_{e}\left(  x^{(k)}\right)
\overline{y_{k}}-\frac{2}{N}\sum_{k=1}^{N}\left\vert s_{e}\left(
x^{(k)}\right)  \right\vert ^{2}\\
& =\frac{2}{N}\operatorname{Re}\sum_{k=1}^{N}\left(  s_{e}\left(
x^{(k)}\right)  \overline{y_{k}}-s_{e}\left(  x^{(k)}\right)  \overline
{s_{e}\left(  x^{(k)}\right)  }\right) \\
& =\frac{2}{N}\operatorname{Re}\sum_{k=1}^{N}s_{e}\left(  x^{(k)}\right)
\left(  \overline{y_{k}}-\overline{s_{e}\left(  x^{(k)}\right)  }\right) \\
& =\frac{2}{N}\operatorname{Re}\sum_{k=1}^{N}\overline{s_{e}\left(
x^{(k)}\right)  }\left(  y_{k}-s_{e}\left(  x^{(k)}\right)  \right)  .
\end{align*}
\medskip

\textbf{Part 4} Substitute the formula for $\left\vert s_{e}\right\vert
_{w,\theta}^{2}$ of part 3 into the definition of $J_{e}\left[  s_{e}\right]
$ to get%
\begin{align*}
J_{e}\left[  s_{e}\right]   &  =\rho\left\vert s_{e}\right\vert _{w,\theta
}^{2}+\frac{1}{N}\sum_{k=1}^{N}\left\vert s_{e}\left(  x^{(k)}\right)
-y_{k}\right\vert ^{2}\\
&  =\frac{1}{N}\operatorname{Re}\sum_{k=1}^{N}s_{e}\left(  x^{(k)}\right)
\left(  \overline{y_{k}}-\overline{s_{e}\left(  x^{(k)}\right)  }\right)  +\\
&  \qquad+\frac{1}{N}\sum_{k=1}^{N}\left(  \left\vert s_{e}\left(
x^{(k)}\right)  \right\vert ^{2}-2\operatorname{Re}s_{e}\left(  x^{(k)}%
\right)  \overline{y_{k}}+\left\vert y_{k}\right\vert ^{2}\right) \\
&  =\frac{1}{N}\sum_{k=1}^{N}\left\vert y_{k}\right\vert ^{2}-\frac{1}%
{N}\operatorname{Re}\sum_{k=1}^{N}s_{e}\left(  x^{(k)}\right)  \overline
{y_{k}}\\
&  =\frac{1}{N}\operatorname{Re}\sum_{k=1}^{N}\left(  y_{k}-s_{e}\left(
x^{(k)}\right)  \right)  \overline{y_{k}}.
\end{align*}

\end{proof}

\subsection{Data functions and the Exact smoother mapping}

Part 4 of Theorem \ref{Thm_smooth_Exact} allows us to define the mapping
between a data function and its corresponding Exact smoother. We call this the
Exact smoother mapping:

\begin{definition}
\label{Def_data_func_Exact smth_map}\textbf{Data functions and the Exact
smoother mapping}

Given an independent data set $X$, we shall assume that each member of
$X_{w}^{\theta}$ can act as a legitimate data function $f$ and generate the
data vector $\widetilde{\mathcal{E}}_{X}f$. By part 5 of Theorem
\ref{Thm_L_op_properties} the linear operator $\left(  \mathcal{L}_{X}^{\ast
}\mathcal{L}_{X}\right)  ^{-1}$ exists and is continuous. Thus part 4 of
Theorem \ref{Thm_smooth_Exact} enables us to define the continuous linear
mapping $\mathcal{S}_{X}^{e}:X_{w}^{\theta}\rightarrow W_{G,X}$ from the data
functions to the corresponding unique Exact smoother by
\begin{equation}
\mathcal{S}_{X}^{e}f=\frac{1}{N}\left(  \mathcal{L}_{X}^{\ast}\mathcal{L}%
_{X}\right)  ^{-1}\widetilde{\mathcal{E}}_{X}^{\ast}\widetilde{\mathcal{E}%
}_{X}f,\quad f\in X_{w}^{\theta}.\label{h10}%
\end{equation}

\end{definition}

We will now prove some properties of the operator $\mathcal{L}_{X}^{\ast
}\mathcal{L}_{X}$ which in turn will be used to prove some basic properties of
the Exact smoother mapping.

\begin{theorem}
\label{Thm_Lx*Lx_onto_Xw,th}\textbf{Properties of} $\mathcal{L}_{X}^{\ast
}\mathcal{L}_{X}$ \textbf{and} $\mathcal{S}_{X}^{e}$

When $\rho>0$ the composition $\mathcal{L}_{X}^{\ast}\mathcal{L}_{X}%
:X_{w}^{\theta}\rightarrow X_{w}^{\theta}$ is a homeomorphism and
$\mathcal{L}_{X}^{\ast}\mathcal{L}_{X}:W_{G,X}\rightarrow W_{G,X}$ is also a
homeomorphism. Regarding the properties of $\mathcal{S}_{X}^{e}$:%
\begin{equation}
\mathcal{S}_{X}^{e}f=f-\rho\,\left(  \mathcal{L}_{X}^{\ast}\mathcal{L}%
_{X}\right)  ^{-1}\mathcal{Q}f,\qquad f\in X_{w}^{\theta},\label{h51}%
\end{equation}

and the operator $\mathcal{S}_{X}^{e}$ is continuous. Further:

\begin{enumerate}
\item $\left\vert \mathcal{S}_{X}^{e}f\right\vert _{w,\theta}\leq\left\vert
f\right\vert _{w,\theta}$ and $\left\vert \left(  I-\mathcal{S}_{X}%
^{e}\right)  f\right\vert _{w,\theta}\leq\left\vert f\right\vert _{w,\theta}$
when $f\in X_{w}^{\theta}$ i.e. $\mathcal{S}_{X}^{e}$ and $I-\mathcal{S}%
_{X}^{e}$ are contractions in the seminorm sense.

\item $\mathcal{S}_{X}^{e}f=f$ iff $f\in P_{\theta-1}$.

\item $\mathcal{S}_{X}^{e}:X_{w}^{\theta}\rightarrow W_{G,X}$ is onto and
$\operatorname*{null}\mathcal{S}_{X}^{e}=W_{G,X}^{\bot}$.

\item The adjoint of the Exact smoother w.r.t. the Light norm is given by:
\[
\left(  \mathcal{S}_{X}^{e}\right)  ^{\ast}g=\frac{1}{N}\widetilde
{\mathcal{E}}_{X}^{\ast}\widetilde{\mathcal{E}}_{X}\left(  \mathcal{L}%
_{X}^{\ast}\mathcal{L}_{X}\right)  ^{-1}g=g-\rho\mathcal{Q}\left(
\mathcal{L}_{X}^{\ast}\mathcal{L}_{X}\right)  ^{-1}g,\quad g\in X_{w}^{\theta
}.
\]

\item $\mathcal{S}_{X}^{e}$ is self-adjoint iff $X$ is a minimal unisolvent
set iff $\mathcal{S}_{X}^{e}=\mathcal{P}$ i.e. the Lagrange polynomial
interpolation function \ref{q82}.
\end{enumerate}
\end{theorem}

\begin{proof}
Suppose $f\in X_{w}^{\theta}$ and let $s_{e}=\mathcal{S}_{X}^{e}f$. From part
5 of Theorem \ref{Thm_L_op_properties} is is known that $\mathcal{L}_{X}%
^{\ast}\mathcal{L}_{X}$ is 1-1. From part 4 of Theorem \ref{Thm_smooth_Exact},
$\mathcal{L}_{X}^{\ast}\mathcal{L}_{X}s_{e}=\frac{1}{N}\widetilde{\mathcal{E}%
}_{X}^{\ast}\widetilde{\mathcal{E}}_{X}f$ and from part 4 of Theorem
\ref{Thm_L_op_properties}, $\mathcal{L}_{X}^{\ast}\mathcal{L}_{X}%
f=\rho\mathcal{Q}f+\frac{1}{N}\widetilde{\mathcal{E}}_{X}^{\ast}%
\widetilde{\mathcal{E}}_{X}f$ so that $\mathcal{L}_{X}^{\ast}\mathcal{L}%
_{X}f=\rho\mathcal{Q}f+\mathcal{L}_{X}^{\ast}\mathcal{L}_{X}s_{e}$. Clearly
this equation implies equation \ref{h51} and also
\begin{equation}
\mathcal{Q}f=\mathcal{L}_{X}^{\ast}\mathcal{L}_{X}\left(  \frac{f-s_{e}}{\rho
}\right)  ,\qquad f\in X_{w}^{\theta}.\label{h50}%
\end{equation}

This last equation proves $\mathcal{Q}\left(  X_{w}^{\theta}\right)
\subset\mathcal{L}_{X}^{\ast}\mathcal{L}_{X}\left(  X_{w}^{\theta}\right)  $
and if we can show that $P_{\theta-1}\subset\mathcal{L}_{X}^{\ast}%
\mathcal{L}_{X}\left(  X_{w}^{\theta}\right)  $ then it follows that
$\mathcal{L}_{X}^{\ast}\mathcal{L}_{X}\left(  X_{w}^{\theta}\right)
=X_{w}^{\theta}$ and so by the \textit{open mapping} theorem $\mathcal{L}%
_{X}^{\ast}\mathcal{L}_{X}$ is a homeomorphism. Suppose $p\in P_{\theta-1}$
and for some $f\in X_{w}^{\theta}$%
\begin{equation}
p=\mathcal{L}_{X}^{\ast}\mathcal{L}_{X}f=\rho\mathcal{Q}f+\frac{1}%
{N}\widetilde{\mathcal{E}}_{X}^{\ast}\widetilde{\mathcal{E}}_{X}f.\label{h20}%
\end{equation}

By Theorem \ref{Thm_eval_op_properties}, $\widetilde{\mathcal{E}}_{X}^{\ast
}\widetilde{\mathcal{E}}_{X}:X_{w}^{\theta}\rightarrow W_{G,X}$ is onto and so
it follows that a solution $g$ to $\frac{1}{N}\widetilde{\mathcal{E}}%
_{X}^{\ast}\widetilde{\mathcal{E}}_{X}g=p$ exists. Now try a solution to
\ref{h20} of the form $f=g+u$ so that $p=\mathcal{L}_{X}^{\ast}\mathcal{L}%
_{X}\left(  g+u\right)  $ and $u$ must satisfy $\mathcal{L}_{X}^{\ast
}\mathcal{L}_{X}u=p-\mathcal{L}_{X}^{\ast}\mathcal{L}_{X}g=p-\left(
\rho\mathcal{Q}g+\frac{1}{N}\widetilde{\mathcal{E}}_{X}^{\ast}\widetilde
{\mathcal{E}}_{X}g\right)  =-\rho\mathcal{Q}g\in\mathcal{Q}\left(
X_{w}^{\theta}\right)  $. But we already know that $\mathcal{Q}\left(
X_{w}^{\theta}\right)  \subset\mathcal{L}_{X}^{\ast}\mathcal{L}_{X}\left(
X_{w}^{\theta}\right)  $ so a there exists $u$ such that $\mathcal{L}%
_{X}^{\ast}\mathcal{L}_{X}u=-\rho\mathcal{Q}g$ and hence a solution $f=g+u$ to
\ref{h20}. Consequently $P_{\theta-1}\subset\mathcal{L}_{X}^{\ast}%
\mathcal{L}_{X}\left(  X_{w}^{\theta}\right)  $ and $\mathcal{L}_{X}^{\ast
}\mathcal{L}_{X}:X_{w}^{\theta}\rightarrow X_{w}^{\theta}$ is a homeomorphism.

By part 3 of Theorem \ref{Sum_ex_Wgx_properties_2}, $W_{G,X}=\overset{\cdot
}{W}_{G,X}\oplus P_{\theta-1}$ and hence $\mathcal{P}\left(  W_{G,X}\right)
=P_{\theta-1}$, $W_{G,X}=\mathcal{Q}\left(  W_{G,X}\right)  \oplus
P_{\theta-1}$ and $\mathcal{Q}\left(  W_{G,X}\right)  =\overset{\cdot
}{W}_{G,X}$.

Thus \ref{p50} implies $\overset{\cdot}{W}_{G,X}=\mathcal{Q}\left(
W_{G,X}\right)  \subset\mathcal{L}_{X}^{\ast}\mathcal{L}_{X}\left(
W_{G,X}\right)  $ and so $\mathcal{L}_{X}^{\ast}\mathcal{L}_{X}\left(
W_{G,X}\right)  =W_{G,X}$. The \textit{open mapping} theorem then implies that
$\mathcal{L}_{X}^{\ast}\mathcal{L}_{X}:W_{G,X}\rightarrow W_{G,X}$ is a
homeomorphism.\medskip

\textbf{Part 1} In equation \ref{h32} let $f$ be the data function and note
that the data function interpolates the data.\smallskip

\textbf{Part 2} Since $\mathcal{L}_{X}^{\ast}\mathcal{L}_{X}$ is a
homeomorphism, from \ref{h51} $\mathcal{S}_{X}^{e}f=f$ iff $\left(
\mathcal{L}_{X}^{\ast}\mathcal{L}_{X}\right)  ^{-1}\mathcal{Q}f=0$ iff
$\mathcal{Q}f=0$ iff $f\in P_{\theta-1}$, where the last implication follows
from part 4 of Theorem \ref{Thm_P_op_properties}.\smallskip

\textbf{Part 3} By part 6 of Theorem \ref{Thm_eval_op_properties},
$\mathcal{S}_{X}^{e}f=0$ iff $\widetilde{\mathcal{E}}_{X}^{\ast}%
\widetilde{\mathcal{E}}_{X}f=0$ iff $f\in W_{G,X}^{\bot}$.

From \ref{h10}, $\mathcal{S}_{X}^{e}f=\frac{1}{N}\left(  \mathcal{L}_{X}%
^{\ast}\mathcal{L}_{X}\right)  ^{-1}\widetilde{\mathcal{E}}_{X}^{\ast
}\widetilde{\mathcal{E}}_{X}f$ and since we know from earlier in the proof
that $\widetilde{\mathcal{E}}_{X}^{\ast}\widetilde{\mathcal{E}}_{X}%
:X_{w}^{\theta}\rightarrow W_{G,X}$ and $\left(  \mathcal{L}_{X}^{\ast
}\mathcal{L}_{X}\right)  ^{-1}:W_{G,X}\rightarrow W_{G,X}$ are both onto we
have our result.\smallskip

\textbf{Part 4\ }A standard Hilbert space result is that if $\mathcal{K}$ is a
homeomorphism then $\left(  \mathcal{K}^{-1}\right)  ^{\ast}=\left(
\mathcal{K}^{\ast}\right)  ^{-1}$. Thus, since $\widetilde{\mathcal{E}}%
_{X}^{\ast}\widetilde{\mathcal{E}}_{X}$ and $\mathcal{L}_{X}^{\ast}%
\mathcal{L}_{X}$ are self-adjoint and by Theorem \ref{Thm_Light_norm_property}
$\mathcal{Q}$ is self-adjoint, taking the adjoint of equations \ref{h51} and
\ref{h50} easily yields the claimed formulas for $\left(  \mathcal{S}_{X}%
^{e}\right)  ^{\ast}$.\smallskip

\textbf{Part 5} Starting with the formulas of \ref{h51} and part 4 we have the
equivalences:%
\[%
\begin{array}
[c]{lll}%
\left(  \mathcal{S}_{X}^{e}\right)  ^{\ast}=\mathcal{S}_{X}^{e} & \text{iff} &
\mathcal{Q}\left(  \mathcal{L}_{X}^{\ast}\mathcal{L}_{X}\right)
^{-1}f=\left(  \mathcal{L}_{X}^{\ast}\mathcal{L}_{X}\right)  ^{-1}%
\mathcal{Q}f\text{,}\quad f\in X_{w}^{\theta}\text{,}\\
& \text{iff} & \mathcal{P}\left(  \mathcal{L}_{X}^{\ast}\mathcal{L}%
_{X}\right)  ^{-1}f=\left(  \mathcal{L}_{X}^{\ast}\mathcal{L}_{X}\right)
^{-1}\mathcal{P}f\text{,}\quad f\in X_{w}^{\theta}\text{.}%
\end{array}
\]

Since $\mathcal{L}_{X}^{\ast}\mathcal{L}_{X}:X_{w}^{\theta}\rightarrow
X_{w}^{\theta}$ is a homeomorphism, if we set $g=\left(  \mathcal{L}_{X}%
^{\ast}\mathcal{L}_{X}\right)  ^{-1}f$ then%
\[
\left(  \mathcal{S}_{X}^{e}\right)  ^{\ast}=\mathcal{S}_{X}^{e}\text{ iff
}\mathcal{L}_{X}^{\ast}\mathcal{L}_{X}\mathcal{P}g=\mathcal{PL}_{X}^{\ast
}\mathcal{L}_{X}g,\quad g\in X_{w}^{\theta}.
\]

Further, by Theorem \ref{Thm_L_op_properties}, $\mathcal{L}_{X}^{\ast
}\mathcal{L}_{X}=\rho\mathcal{Q}+\frac{1}{N}\widetilde{\mathcal{E}}_{X}^{\ast
}\widetilde{\mathcal{E}}_{X}$ and since $\mathcal{PQ}=0$
\[
\left(  \mathcal{S}_{X}^{e}\right)  ^{\ast}=\mathcal{S}_{X}^{e}\text{ iff
}\widetilde{\mathcal{E}}_{X}^{\ast}\widetilde{\mathcal{E}}_{X}\mathcal{P}%
g=\mathcal{P}\widetilde{\mathcal{E}}_{X}^{\ast}\widetilde{\mathcal{E}}%
_{X}g,\quad g\in X_{w}^{\theta}.
\]

Recall the properties of the Riesz representer $R_{x}$ given in Subsection
\ref{SbSect_Riesz_rep_Xpos}. Next suppose $X=\left\{  x^{\left(  k\right)
}\right\}  _{k=1}^{N}$ is not minimally unisolvent and that the minimal
unisolvent set $A\subset X$ was used to construct the smoother, $\mathcal{P}$,
$\mathcal{Q}$ etc. Then $\mathcal{P}\widetilde{\mathcal{E}}_{X}^{\ast
}\widetilde{\mathcal{E}}_{X}g\in P_{\theta-1}$ and so $\widetilde{\mathcal{E}%
}_{X}^{\ast}\widetilde{\mathcal{E}}_{X}\mathcal{P}g=\mathcal{P}g+\sum
\limits_{x^{\left(  i\right)  }\in X\setminus A}\left(  \mathcal{P}g\right)
\left(  x^{\left(  i\right)  }\right)  R_{x^{\left(  i\right)  }}\in
P_{\theta-1}$ with the independence of the functions $\left\{  R_{x^{\left(
i\right)  }}\right\}  _{i=1}^{N}$ implying that $\left(  \mathcal{P}g\right)
\left(  x^{\left(  i\right)  }\right)  =0$ when $x^{\left(  i\right)  }\in
X\setminus A$ and $g\in X_{w}^{\theta}$. Therefore, given $x^{\left(
i\right)  }\in X\setminus A$ we have $p\left(  x^{\left(  i\right)  }\right)
=0$ when $p\in P_{\theta-1}$ which is a contradiction. We conclude that
$\left(  \mathcal{S}_{X}^{e}\right)  ^{\ast}\neq\mathcal{S}_{X}^{e}$ or $X$ is
minimally unisolvent.

If $X$ is a minimal unisolvent set and $f$ is a data function then
$\mathcal{P}f\left(  x^{\left(  k\right)  }\right)  =f\left(  x^{\left(
k\right)  }\right)  $ and $\mathcal{P}f\in P_{\theta-1}$ so that the Exact
smoother functional satisfies $J_{e}\left[  \mathcal{P}f\right]  =0$ and thus
$\mathcal{S}_{X}^{e}f=\mathcal{P}f$ i.e. $\mathcal{S}_{X}^{e}=\mathcal{P}$.
Finally, if $\mathcal{S}_{X}^{e}=\mathcal{P}$ then $\mathcal{S}_{X}^{e}$ is
self-adjoint w.r.t. the Light norm by Theorem \ref{Thm_Light_norm_property}.
\end{proof}

\section{Expressing $R_{X,X}$ in terms of $G_{X,X}$%
\label{Sect_Rxx_in_terms_of_Gxx}}

In this section we will derive equations expressing the reproducing kernel
matrix in terms of the basis function matrix. I start with some notation:

\begin{definition}
\label{Def_Matrices_from_R_2}\textbf{Matrices and vectors derived from the
Riesz representer }$R_{x}$ \textbf{and the basis function} $G$.

If $X=\left\{  x^{\left(  k\right)  }\right\}  \subset\mathbb{R}^{d}$ and
$y\in\mathbb{R}^{d}$ then $R_{y,X}=\left(  R_{x^{\left(  j\right)  }}\left(
y\right)  \right)  $ is a row vector and $G_{y,X}=\left(  G\left(
y-x^{\left(  j\right)  }\right)  \right)  $ is a row vector.

Also $R_{X,y}=\left(  R_{y}\left(  x^{\left(  i\right)  }\right)  \right)  $
is a column vector and $G_{X,y}=\left(  G\left(  x^{\left(  i\right)
}-y\right)  \right)  $ is a column vector.
\end{definition}

Using the notation of Definition \ref{Def_Matrices_from_R_2} this theorem
derives some of the relationships between the matrices $R_{X,X}$ and $G_{X,X}
$.

\begin{theorem}
\label{Thm_relatnship_between_R_and_G_matrices}Suppose $A$ is any minimal
unisolvent set and let the corresponding cardinal basis be $\left\{
l_{i}\right\}  _{i=1}^{M}$. Define the Riesz representer $R_{x}$ using $A$ and
$\left\{  l_{i}\right\}  _{i=1}^{M}$. Then for any finite set $X\subset
\mathbb{R}^{d}$ and $y\in\mathbb{R}^{d}$
\begin{equation}
R_{y,X}=\left(  2\pi\right)  ^{-\frac{d}{2}}\left(  G_{y,X}-\widetilde
{l}\left(  y\right)  ^{T}G_{A,X}-G_{y,A}L_{X}^{T}+\widetilde{l}\left(
y\right)  ^{T}G_{A,A}L_{X}^{T}\right)  +\widetilde{l}\left(  y\right)
^{T}L_{X}^{T}.\nonumber
\end{equation}

The reproducing kernel\textbf{\ }matrix $R_{X,X}=\left(  R_{x^{(j)}}\left(
x^{(i)}\right)  \right)  $ and the basis function matrix $G_{X,X}=$

$\left(  G\left(  x^{(i)}-x^{(j)}\right)  \right)  $ are related by the
formulas%
\begin{equation}
R_{X,X}=\left(  2\pi\right)  ^{-\frac{d}{2}}\left(  G_{X,X}-L_{X}%
G_{A,X}-G_{X,A}L_{X}^{T}+L_{X}G_{A,A}L_{X}^{T}\right)  +L_{X}L_{X}%
^{T},\label{h21}%
\end{equation}

and%
\begin{equation}
R_{X,X}=\left(  2\pi\right)  ^{-\frac{d}{2}}\left(  I_{N}-L_{X;0}\right)
G_{X,X}\left(  I_{N}-L_{X;0}\right)  ^{T}+L_{X}L_{X}^{T},\label{h22}%
\end{equation}

where $L_{X}=\left(  l_{j}(x^{\left(  i\right)  })\right)  $ and $L_{X;0}=%
\begin{pmatrix}
L_{X} & O_{N,N-M}%
\end{pmatrix}
$.
\end{theorem}

\begin{proof}
From \ref{p939}
\begin{multline*}
\left(  2\pi\right)  ^{\frac{d}{2}}R_{x}\left(  y\right)  =G\left(
y-x\right)  -\sum_{j=1}^{M}l_{j}\left(  y\right)  G\left(  x^{\left(
j\right)  }-x\right)  -\sum_{i=1}^{M}G\left(  y-x^{\left(  i\right)  }\right)
l_{i}\left(  x\right)  +\\
+\sum_{i,j=1}^{M}l_{j}\left(  y\right)  G\left(  x^{\left(  j\right)
}-x^{\left(  i\right)  }\right)  l_{i}\left(  x\right)  +\left(  2\pi\right)
^{\frac{d}{2}}\sum\limits_{j=1}^{M}l_{j}(x)l_{j}\left(  y\right)  ,
\end{multline*}

or in the notation introduced in Definition \ref{Def_Matrices_from_R_2}%
\[
R_{x}\left(  y\right)  =\left(  2\pi\right)  ^{-\frac{d}{2}}\left(  G\left(
y-x\right)  -\widetilde{l}\left(  y\right)  ^{T}G_{A,x}-G_{y,A}\widetilde
{l}\left(  x\right)  +\widetilde{l}\left(  y\right)  ^{T}G_{A,A}\widetilde
{l}\left(  x\right)  \right)  +\widetilde{l}\left(  y\right)  ^{T}%
\widetilde{l}\left(  x\right)  ,
\]

Now $R_{y,X}$ is the row vector $\left(  R_{x^{\left(  j\right)  }}\left(
y\right)  \right)  $ and $L_{X}=\left(  l_{j}\left(  x^{(i)}\right)  \right)
$ so%
\[
R_{y,X}=\left(  2\pi\right)  ^{-\frac{d}{2}}\left(  G_{y,X}-\widetilde
{l}\left(  y\right)  ^{T}G_{A,X}-G_{y,A}L_{X}^{T}+\widetilde{l}\left(
y\right)  ^{T}G_{A,A}L_{X}^{T}\right)  +\widetilde{l}\left(  y\right)
^{T}L_{X}^{T},
\]

and hence, since $L_{X;0}=\left(  L_{X}\text{ }O_{N,N-M}\right)  $,
\begin{align*}
R_{X,X} &  =\left(  R_{x^{\left(  i\right)  },X}\right) \\
&  =\left(  2\pi\right)  ^{-\frac{d}{2}}\left(  G_{x^{\left(  i\right)  }%
,X}-\widetilde{l}\left(  x^{\left(  i\right)  }\right)  ^{T}G_{A,X}%
-G_{x^{\left(  i\right)  },A}L_{X}^{T}+\widetilde{l}\left(  x^{\left(
i\right)  }\right)  ^{T}G_{A,A}L_{X}^{T}\right)  +\widetilde{l}\left(
x^{\left(  i\right)  }\right)  ^{T}L_{X}^{T}\\
&  =\left(  2\pi\right)  ^{-\frac{d}{2}}\left(  G_{X,X}-L_{X}^{T}%
G_{A,X}-G_{X,A}L_{X}^{T}+L_{X}G_{A,A}L_{X}^{T}\right)  +L_{X}L_{X}^{T}\\
&  =\left(  2\pi\right)  ^{-\frac{d}{2}}\left(  G_{X,X}-L_{X;0}G_{X,X}%
-G_{X,X}L_{X;0}^{T}+L_{X;0}G_{X,X}L_{X;0}^{T}\right)  +L_{X;0}L_{X;0}^{T}\\
&  =\left(  2\pi\right)  ^{-\frac{d}{2}}\left(  I_{N}-L_{X;0}\right)
G_{X,X}\left(  I_{N}-L_{X;0}\right)  ^{T}+L_{X;0}L_{X;0}^{T}.
\end{align*}

\end{proof}

\section{Matrices, vectors and bases derived from the semi-Riesz representer
$r_{x}$\label{Sect_mat_from_rx}}

The \textit{semi-Riesz representer} $r_{x}=\mathcal{Q}R_{x}$ was introduced in
Chapter \ref{Ch_Interpol} where $r_{x}\left(  x\right)  $ was used to estimate
the convergence of the interpolant. In this document $r_{x}\left(  x\right)  $
will be used to estimate the convergence of the Exact smoother. Some basis
properties of $r_{x}$ were proved in Theorem \ref{Thm_rx(y)_properties} and we
will need these to prove the following result:

\begin{corollary}
\label{Cor_rxi_basis_Wgx}Suppose $X=\left\{  x^{\left(  i\right)  }\right\}
_{i=1}^{N}$ is unisolvent and $A=\left\{  a^{\left(  i\right)  }\right\}
_{i=1}^{M}$ is a minimal unisolvent subset. If $r_{x}$ is defined using this
set then the functions $\left\{  r_{x^{\left(  i\right)  }}:x^{\left(
i\right)  }\notin A\right\}  $ are independent and form a basis for
$\overset{\cdot}{W}_{G,X}$.
\end{corollary}

\begin{proof}
Firstly, by part 2 Theorem \ref{Thm_rx(y)_properties}, $r_{a^{\left(
i\right)  }}=0$ for $i=1,\ldots,M$. Next we prove independence. If
$\sum\limits_{x^{\left(  i\right)  }\notin A}\beta_{i}r_{x^{\left(  i\right)
}}=0$ then $\mathcal{Q}\sum\limits_{x^{\left(  i\right)  }\notin A}\beta
_{i}R_{x^{\left(  i\right)  }}=0$ so that $\sum\limits_{x^{\left(  i\right)
}\notin A}\beta_{i}R_{x^{\left(  i\right)  }}=\mathcal{P}\sum
\limits_{x^{\left(  i\right)  }\notin A}\beta_{i}R_{x^{\left(  i\right)  }}\in
P_{\theta-1} $. By part 6 of Summary \ref{Sum_ex_Wgx_properties_2},
$\sum\limits_{x^{\left(  i\right)  }\notin A}\beta_{i}R_{x^{\left(  i\right)
}}\in\overset{\cdot}{W}_{G,X}$. But by part 3 of Summary
\ref{Sum_ex_Wgx_properties_2}, $W_{G,X}=\overset{\cdot}{W}_{G,X}\oplus
P_{\theta-1}$ so all the $\beta_{i} $ are zero. Also by part 3 of Summary
\ref{Sum_ex_Wgx_properties_2}, $\overset{\cdot}{W}_{G,X}$ has dimension $N-M $
and so the $N-M$ functions $\left\{  r_{x^{\left(  i\right)  }}:x^{\left(
i\right)  }\notin A\right\}  $ form a basis.
\end{proof}

\begin{definition}
\label{Def_rX,x_rX,Y}\textbf{Matrices and vectors derived from }$r_{x}$

If $X=\left\{  x^{\left(  k\right)  }\right\}  \subset\mathbb{R}^{d}$ and
$y\in\mathbb{R}^{d}$ then $r_{X,X}=\left(  r_{x^{\left(  j\right)  }}\left(
x^{\left(  i\right)  }\right)  \right)  $, $r_{X,y}=\left(  r_{y}\left(
x^{\left(  i\right)  }\right)  \right)  $ is a column vector and
$r_{y,X}=\left(  r_{x^{\left(  j\right)  }}\left(  y\right)  \right)  $ is a
row vector.
\end{definition}

\begin{theorem}
\label{Thm_rx(y)_mat_vect}Suppose $X=\left\{  x^{\left(  i\right)  }\right\}
_{i=1}^{N}$ is unisolvent and $A$ is a minimal unisolvent subset with cardinal
basis $\left\{  l_{i}\right\}  _{i=1}^{M}$. Then if $r_{x}=\mathcal{Q}R_{x}$
is defined using $A$ we have:

\begin{enumerate}
\item $r_{X,y}=R_{X,y}-L_{X}\widetilde{l}\left(  y\right)  $.

\item $r_{X,X}=R_{X,X}-L_{X}L_{X}^{T}$.

\item The matrix $r_{A^{c},A^{c}}$ is positive definite, regular and Hermitian
where $A^{c}=X\setminus A$.
\end{enumerate}
\end{theorem}

\begin{proof}
\textbf{Parts 1 and 2} The definition of the cardinal unisolvency matrix is
$L_{X}=\left(  l_{j}\left(  x^{(i)}\right)  \right)  $ and by part 6 Theorem
\ref{Thm_rx(y)_properties}, $R_{y}=r_{y}+\sum\limits_{j=1}^{M}l_{j}\left(
y\right)  l_{j}$. Parts 1 and 2 then follow easily from the definitions of
$r_{X,y}$ and $r_{X,X}$.\medskip

\textbf{Part 3} From Theorem \ref{Thm_rx(y)_properties}, $r_{x^{\left(
j\right)  }}\left(  x^{\left(  i\right)  }\right)  =\left\langle r_{x^{\left(
j\right)  }},r_{x^{\left(  i\right)  }}\right\rangle _{w,\theta}$ and
$r_{x^{\left(  j\right)  }}\left(  a^{\left(  i\right)  }\right)  =0$ for
$a^{\left(  i\right)  }\in A$. Now let $\left(  \cdot,\cdot\right)
_{w,\theta} $ be the Light inner product \ref{p17}: $\left(  u,v\right)
_{w,\theta}=\left\langle u,v\right\rangle _{w,\theta}+\sum\limits_{i=1}%
^{M}u\left(  a^{\left(  i\right)  }\right)  \overline{v\left(  a^{\left(
i\right)  }\right)  }$ constructed using $A$. It now follows that
$r_{x^{\left(  j\right)  }}\left(  x^{\left(  i\right)  }\right)  =\left(
r_{x^{\left(  j\right)  }},r_{x^{\left(  i\right)  }}\right)  _{w,\theta}$ for
all $i$ and $j$, and so the matrix $r_{A^{c},A^{c}}=\left(  r_{x^{\left(
j\right)  }}\left(  x^{\left(  i\right)  }\right)  \right)  $ is a Gram matrix
and the properties stated in this theorem are well-known properties of Gram matrices.
\end{proof}

\section{Matrix equations for the Exact
smoother\label{Sect_Exact_smth_mat_eqn}}

In Corollary \ref{Cor_propert_Exact_smth} it was shown that the unique
solution of the Exact smoothing problem of the data $\left[  X,y\right]  $
lies in the finite dimensional space $W_{G,X}$ which has a unique
representation given by \ref{h14} i.e.%
\[
\sum\limits_{i=1}^{N}\alpha_{i}G\left(  \cdot-x^{\left(  i\right)  }\right)
+\sum\limits_{j=1}^{M}\beta_{j}p_{j},\quad\alpha_{i},\beta_{j}\in\mathbb{C},
\]

where $X=\left\{  x^{\left(  i\right)  }\right\}  _{i=1}^{N}$. The goal of
this section is to derive a matrix equation for the coefficients $\alpha
_{i},\beta_{j}$ of the basis functions of the space $W_{G,X}$. We start by
deriving a (singular) matrix equation for the values taken by the Exact
smoother at its independent data points $X$. Now recall the notation $f_{X}$
for $\widetilde{\mathcal{E}}_{X}f$ introduced in Definition
\ref{Def_eval_operators}.

\begin{theorem}
\label{Thm_mat_eqn_Sx}Let $X=\left\{  x^{(i)}\right\}  _{i=1}^{N}$ be a
unisolvent set of independent data with a minimally unisolvent subset $A$.
Construct $R_{x}$ using $A$ and its cardinal basis $\left\{  l_{i}\right\}
_{i=1}^{M}$, and denote the Exact smoother of the data $\left[  X,y\right]  $
by $s$. Then
\begin{equation}
\left(  N\rho\left(  I_{N}-L_{X;0}\right)  +R_{X,X}\right)  s_{X}%
=R_{X,X}y,\label{h66}%
\end{equation}

and%
\begin{equation}
L_{X}^{T}\left(  s_{X}-y\right)  =0,\label{h61}%
\end{equation}

where $L_{X}=\left(  l_{j}\left(  x^{(i)}\right)  \right)  $, $L_{X;0}=%
\begin{pmatrix}
L_{X} & O_{N,N-M}%
\end{pmatrix}
$ and $R_{X,X}=\left(  R_{x^{(j)}}(x^{(i)})\right)  $ is the reproducing
kernel matrix.
\end{theorem}

\begin{proof}
From parts 4 and 8 of Theorem \ref{Thm_eval_op_properties} we have
$\widetilde{\mathcal{E}}_{X}\widetilde{\mathcal{E}}_{X}^{\ast}\beta
=R_{X,X}\beta$ and $\widetilde{\mathcal{E}}_{X}\mathcal{P}s=L_{X}%
\widetilde{\mathcal{E}}_{A}s$. From part 1 of Corollary
\ref{Cor_propert_Exact_smth} we have $s=\mathcal{P}s-\frac{1}{N\rho}%
\widetilde{\mathcal{E}}_{X}^{\ast}\left(  \widetilde{\mathcal{E}}%
_{X}s-y\right)  $. Applying $\widetilde{\mathcal{E}}_{X}$ to the last equation
to get
\begin{align*}
s_{X}=\widetilde{\mathcal{E}}_{X}\mathcal{P}s-\frac{1}{N\rho}\widetilde
{\mathcal{E}}_{X}\widetilde{\mathcal{E}}_{X}^{\ast}\left(  s_{X}-y\right)   &
=L_{X}\widetilde{\mathcal{E}}_{A}s-\frac{1}{N\rho}R_{X,X}\left(
s_{X}-y\right) \\
&  =L_{X;0}s_{X}-\frac{1}{N\rho}R_{X,X}s_{X}+\frac{1}{N\rho}R_{X,X}y,
\end{align*}

or on rearranging%
\[
\left(  N\rho\left(  I_{N}-L_{X;0}\right)  +R_{X,X}\right)  s_{X}=R_{X,X}y,
\]

which is the desired matrix equation in $\widetilde{\mathcal{E}}_{X}s$.
Finally \ref{h61} follows from part 3 Corollary \ref{Cor_propert_Exact_smth}
and equation \ref{q06} of part 3 Theorem \ref{Thm_Px_properties}.
\end{proof}

Corollary \ref{Cor_propert_Exact_smth} established that the Exact smoother has
a unique solution in the finite dimensional space $W_{G,X}$. By Definition
\ref{Def_Wgx}, $W_{G,X}$ is independent of the basis function $G$, the order
of the points in $X$ and the basis of $P_{\theta-1}$ used to define $P_{X}$.
The next theorem derives the corresponding matrix equation for the
coefficients of the basis functions. This matrix equation is deduced using the
relationships between the reproducing kernel matrix $R_{X,X}$ and the basis
function matrix $G_{X,X}$ derived in Theorem
\ref{Thm_relatnship_between_R_and_G_matrices}. However, to prove the next
(well-known) result we need Lemma \ref{Lem_reg_cspd} from Chapter
\ref{Ch_Interpol}:

\begin{lemma}
\label{Lem_reg_cspd_2}(Lemma \ref{Lem_reg_cspd})Let $B$ be a complex-valued
matrix and $C$ be a real-valued matrix. Suppose the block matrix $\left(
\begin{array}
[c]{ll}%
B & C\\
C^{T} & O
\end{array}
\right)  $ is square and that for complex vectors $z$
\begin{equation}
z^{T}B\overline{z}=0\text{ }and\text{ }C^{T}z=0\text{ }implies\text{
}z=0.\label{h34}%
\end{equation}

\begin{enumerate}
\item Then the equation
\[
\left(
\begin{array}
[c]{ll}%
B & C\\
C^{T} & O
\end{array}
\right)  \left(
\begin{array}
[c]{l}%
u\\
v
\end{array}
\right)  =\left(
\begin{array}
[c]{l}%
0\\
0
\end{array}
\right)  ,
\]

implies $u=0$ and $v\in\operatorname*{null}C$.

\item If, in addition to \ref{h34} $\operatorname*{null}C=\left\{  0\right\}
$, then the block matrix is regular.
\end{enumerate}
\end{lemma}

\begin{theorem}
\label{Thm_smooth_matrix_soln_1}Suppose $\left[  X,y\right]  $ is the data for
the Exact smoothing problem of Definition \ref{Def_Exact_smth_problem} where
$X=\left\{  x^{(i)}\right\}  _{i=1}^{N}$ is $\theta$-unisolvent and
$y=\left\{  y_{i}\right\}  _{i=1}^{N}$. Then given a real-valued basis
$\left\{  p_{j}\right\}  _{j=1}^{M}$ for $P_{\theta-1}$ and a basis function
$G$ of order $\theta$, the Exact smoothing problem has a unique solution $s\in
W_{G,X}$ of the form%
\begin{equation}
s(x)=\sum_{i=1}^{N}v_{i}G\left(  x-x^{(i)}\right)  +\sum_{j=1}^{M}\beta
_{j}p_{j}\left(  x\right)  ,\label{h67}%
\end{equation}

where the coefficients $v=\left(  v_{i}\right)  \in\mathbb{C}^{N}$ and
$\beta=\left(  \beta_{j}\right)  \in\mathbb{C}^{N}$ satisfy the \textbf{Exact
smoothing matrix} \textbf{equation}
\begin{equation}
\left(
\begin{array}
[c]{ll}%
\left(  2\pi\right)  ^{\frac{d}{2}}N\rho I_{N}+G_{X,X} & P_{X}\\
P_{X}^{T} & O_{M}%
\end{array}
\right)  \left(
\begin{array}
[c]{l}%
v\\
\beta
\end{array}
\right)  =\left(
\begin{array}
[c]{l}%
y\\
0
\end{array}
\right)  .\label{h68}%
\end{equation}

Here $P_{X}=\left(  p_{j}\left(  x^{(i)}\right)  \right)  $ is a unisolvency
matrix and $G_{X,X}=\left(  G\left(  x^{(i)}-x^{\left(  j\right)  }\right)
\right)  $ is the basis function matrix.

The matrix equation \ref{h68} is independent of the ordering of the data
$\left[  X,y\right]  $.

Finally, the Exact smoothing\textbf{\ }matrix is regular and positive definite
but in general it is not Hermitian. However, there always exists a basis
function such that the Exact smoothing matrix is Hermitian.
\end{theorem}

\begin{proof}
\textbf{Step 1} Since $X=\left\{  x^{\left(  i\right)  }\right\}  _{i=1}^{N}$
is $\theta$-unisolvent there exists a minimal unisolvent subset, say $X_{1}$.
Assume that $X_{1}=\left\{  x^{\left(  i\right)  }\right\}  _{i=1}^{M}$. Now
construct $\mathcal{P}$, $\mathcal{Q}$, the Light norm and the Riesz
representer $R_{x}$ using $X_{1}$ and its cardinal basis $\left\{
l_{i}\right\}  _{i=1}^{M}$ for $P_{\theta-1}$.

Initially we will derive the Exact smoother matrix equation using the cardinal
basis associated with $X_{1}$ so that $P_{X}=L_{X}=\left(  l_{j}\left(
x^{\left(  i\right)  }\right)  \right)  $. Consequently, by the interpolation
result Theorem \ref{Thm_interpol_matrix_eqn} of Chapter \ref{Ch_Interpol}, for
given dependent data $s_{X}=\left(  s\left(  x^{(i)}\right)  \right)
_{i=1}^{N}$ there exist unique vectors $v$ and $\gamma$ such that
\[
\left(
\begin{array}
[c]{ll}%
G_{X,X} & L_{X}\\
L_{X}^{T} & O_{M}%
\end{array}
\right)  \left(
\begin{array}
[c]{l}%
v\\
\gamma
\end{array}
\right)  =\left(
\begin{array}
[c]{l}%
s_{X}\\
0
\end{array}
\right)  ,
\]

This matrix equation is equivalent to the two equations
\begin{align}
G_{X,X}v+L_{X}\gamma &  =s_{X}\label{h73}\\
L_{X}^{T}v &  =0.\label{h25}%
\end{align}

Multiplying equation \ref{h73} by $N\rho R_{X,X}^{-1}\left(  I_{N}%
-L_{X;0}\right)  +I_{N}$ and using equation \ref{h66} gives
\[
\left(  N\rho R_{X,X}^{-1}\left(  I_{N}-L_{X;0}\right)  +I_{N}\right)  \left(
G_{X,X}v+L_{X}\gamma\right)  =\left(  N\rho R_{X,X}^{-1}\left(  I_{N}%
-L_{X;0}\right)  +I_{N}\right)  s_{X}=y.
\]

From part 3 Theorem \ref{Thm_Px_properties_2}, $\left(  I_{N}-L_{X;0}\right)
L_{X}=O$ so the last equation simplifies to%
\begin{equation}
N\rho R_{X,X}^{-1}\left(  I_{N}-L_{X;0}\right)  G_{X,X}v+G_{X,X}v+L_{X}%
\gamma=y.\label{h30}%
\end{equation}

We now require an equation that expresses $R_{X,X}$ in terms of $G_{X,X}$. To
this end we use equation \ref{h22}, namely%
\[
R_{X,X}=\left(  2\pi\right)  ^{-\frac{d}{2}}\left(  I_{N}-L_{X;0}\right)
G_{X,X}\left(  I_{N}-L_{X;0}\right)  ^{T}+L_{X}L_{X}^{T}.
\]

Multiplying this equation on the left by $R_{X,X}^{-1}$ and by $v$ on the
right, and noting that $L_{X}^{T}v=0$, yields%
\begin{align*}
v  & =\left(  2\pi\right)  ^{-\frac{d}{2}}R_{X,X}^{-1}\left(  I_{N}%
-L_{X;0}\right)  G_{X,X}\left(  I_{N}-L_{X;0}\right)  ^{T}v+R_{X,X}^{-1}%
L_{X}L_{X}^{T}v\\
& =\left(  2\pi\right)  ^{-\frac{d}{2}}R_{X,X}^{-1}\left(  I_{N}%
-L_{X;0}\right)  G_{X,X}v,
\end{align*}

and thus \ref{h30} reduces to
\begin{equation}
\left(  \left(  2\pi\right)  ^{\frac{d}{2}}N\rho I_{N}+G_{X,X}\right)
v+L_{X}\gamma=y,\label{h72}%
\end{equation}

constrained by equation \ref{h25} i.e. $L_{X}^{T}v=0$.

Since $P_{X}=\left(  p_{j}\left(  x^{\left(  i\right)  }\right)  \right)  $
and $L_{X}=\left(  l_{j}\left(  x^{\left(  i\right)  }\right)  \right)  $,
where $\left\{  p_{j}\right\}  $ and $\left\{  l_{j}\right\}  $ are both
real-valued bases for $P_{\theta-1}$, it follows from part 3 Theorem
\ref{Thm_Px_properties}\ that there exists a regular matrix $C$ such that
$L_{X}=P_{X}C$. Substituting for $L_{X}$ in \ref{h72} and \ref{h25} and then
setting $\beta=C\gamma$ yields the Exact smoothing matrix equation
\ref{h68}.\medskip

\textbf{Step 2} We now show that re-ordering the data does not change the
equation for the smoother. Re-ordering a vector involves a permutation $\pi$.
Re-ordering a column vector involves left-multiplication by the permutation
matrix $\Pi$ and re-ordering a row vector involves right-multiplication by the
transpose of the permutation matrix. Thus the new matrix equation is%
\[
\left(
\begin{array}
[c]{ll}%
\left(  2\pi\right)  ^{\frac{d}{2}}N\rho I_{N}+G_{\pi\left(  X\right)
,\pi\left(  X\right)  } & P_{\pi\left(  X\right)  }\\
P_{\pi\left(  X\right)  }^{T} & O_{M}%
\end{array}
\right)  \left(
\begin{array}
[c]{l}%
v^{\prime}\\
\beta^{\prime}%
\end{array}
\right)  =\left(
\begin{array}
[c]{l}%
\Pi y\\
0
\end{array}
\right)  .
\]

Thus $G_{\pi\left(  X\right)  ,\pi\left(  X\right)  }=\Pi G_{X,X}\Pi^{T}$,
$P_{\pi\left(  X\right)  }=\Pi P_{X}$ and $\pi\left(  y\right)  =\Pi y$, or
since $\Pi\Pi^{T}=I_{N}$. Under the last set of transformations the Exact
smoothing matrix equation becomes
\[
\left(  \left(  2\pi\right)  ^{\frac{d}{2}}N\rho\Pi\Pi^{T}+\Pi G_{X,X}\Pi
^{T}\right)  v^{\prime}+\Pi P_{X}\beta^{\prime}=\Pi y,\quad P_{X}^{T}\Pi
^{T}v^{\prime}=0,
\]

or%
\begin{equation}
\left(  \left(  2\pi\right)  ^{\frac{d}{2}}N\rho+G_{X,X}\right)  \Pi
^{T}v^{\prime}+P_{X}\beta^{\prime}=y,\quad P_{X}^{T}\Pi^{T}v^{\prime
}=0,\label{h27}%
\end{equation}

and so the matrix equation is unchanged by re-ordering the data.\medskip

\textbf{Step 3} The next step is to show that the Exact smoothing matrix is
regular. To do this we use Lemma \ref{Lem_reg_cspd_2} with $B=\left(
2\pi\right)  ^{\frac{d}{2}}N\rho I_{N}+G_{X,X}$, $C=P_{X}$, and then show that
$z^{T}\left(  \left(  2\pi\right)  ^{\frac{d}{2}}N\rho I_{N}+G_{X,X}\right)
\overline{z}=0$ and $P_{X}^{T}z=0$ implies $z=0$. Now%
\begin{align*}
z^{T}\left(  \left(  2\pi\right)  ^{\frac{d}{2}}N\rho I_{N}+G_{X,X}\right)
\overline{z} &  =\left(  2\pi\right)  ^{\frac{d}{2}}z^{T}N\rho I_{N}%
\overline{z}+z^{T}G_{X,X}\overline{z}\\
&  =\left(  2\pi\right)  ^{\frac{d}{2}}N\rho\left\vert z\right\vert ^{2}%
+z^{T}G_{X,X}\overline{z}.
\end{align*}

But by part 2 Theorem \ref{Sum_ex_Wgx_properties_2}, $G_{X,X}$ is
conditionally positive definite on $\operatorname*{null}P_{X}^{T}$ i.e.
$z\in\operatorname*{null}P_{X}^{T}$ implies that $z^{T}G_{X,X}\overline{z}>0$
except when $z=0$. Thus $z\in\operatorname*{null}P_{X}^{T}$ and $z^{T}\left(
\left(  2\pi\right)  ^{\frac{d}{2}}N\rho I_{N}+G_{X,X}\right)  \overline{z}=0$
implies $z=0 $ and the Approximate smoother matrix is regular.\medskip

Finally, part 2 Theorem \ref{Thm_Grho} from Chapter \ref{Ch_MoreBasisTheory}
allows the basis function $G$ to be chosen so that $\overline{G\left(
x\right)  }=G\left(  -x\right)  $ and this implies $G_{X,X}$ is Hermitian.
\end{proof}

\begin{remark}
\label{Rem_Thm_smooth_matrix_soln_1}\ 

\begin{enumerate}
\item When $\rho=0$ the matrix equation for the Exact smoother becomes the
matrix equation for the minimal norm interpolant - see Theorem
\ref{Thm_interpol_matrix_eqn}.

\item The matrix \ref{h68} is $N\times N$ i.e. its size depends on the number
of data points, so this algorithm is \textbf{not} \textbf{scalable} i.e. the
time of execution is not linearly dependent on the number of data points. The
\textbf{Approximate smoother}, which overcomes this problem, will be derived
in Chapter \ref{Ch_Approx_smth}.
\end{enumerate}
\end{remark}

The basis function form of the Exact smoother matrix equation can also be
derived using Lagrange multipliers. The next result shows that the algebra can
be significantly simplified by assuming the basis function is real valued.

\begin{corollary}
\label{Cor_Exact_smth_real}If the dependent data $y$ is real-valued and the
basis function is real-valued, then the Exact smoother is real-valued. Also,
the smoother lies in the subspace of $W_{G,X}$ defined using the real scalars
instead of the complex scalars.
\end{corollary}

\begin{proof}
By Theorem \ref{Thm_relatnship_between_R_and_G_matrices}, if the basis
function $G$ is real-valued then the reproducing kernel matrix $R_{X,X}$ is
real-valued. Hence, since the cardinal unisolvency matrix $L_{X}$ is
real-valued, equations \ref{h67} and \ref{h68} imply that the smoother is real-valued.
\end{proof}

\section{Convergence to the data function - smoother error
\label{Sect_converg_Exact_smth}}

In this section we will prove that in the sense of Corollary
\ref{Cor_Exact_smth_converg} below\ the Exact smoother converges uniformly
pointwise to its data function on a bounded set. The rates of convergence are
shown to be\ the same as those obtained for the minimal seminorm interpolant.

\subsection{Convergence to the data function\label{SbSect_ex_Exact_smth_error}%
}

In Section \ref{Sect_soln_Exact_smth} the Exact smoother problem of Definition
\ref{Def_ex_Hilbert_smoothing} was studied using the Hilbert space
$V=X_{w}^{\theta}\otimes\mathbb{C}^{N}$ endowed with the inner product
$\left(  \left(  f,\widetilde{\alpha}\right)  ,\left(  g,\widetilde{\beta
}\right)  \right)  _{V}=\rho\left\langle f,g\right\rangle _{w,\theta}+\frac
{1}{N}\left(  \widetilde{\alpha},\widetilde{\beta}\right)  _{\mathbb{C}^{N}}$,
and the operator $\mathcal{L}_{X}:X_{w}^{\theta}\rightarrow V$ defined by
$\mathcal{L}_{X}f=\left(  f,\widetilde{\mathcal{E}}_{X}f\right)  $ for a set
of $N$ independent data points $X$.

In Theorem \ref{Thm_L_op_properties} it was proven that $\left\Vert
\mathcal{L}_{X}f\right\Vert _{V}$ is an equivalent norm to $\left\Vert
f\right\Vert _{w,\theta}$ which implies that $X_{w}^{\theta}$ is also a
reproducing kernel Hilbert space under the norm $\left\Vert \mathcal{L}%
_{X}f\right\Vert _{V}$. The space $V$ induces on $X_{w}^{\theta}$ the inner
product%
\[
\left(  f,g\right)  _{V,w,\theta}=\left(  \mathcal{L}_{X}f,\mathcal{L}%
_{X}g\right)  _{V},\quad f,\text{ }g\in X_{w}^{\theta}.
\]

Under this inner product $X_{w}^{\theta}$ is a reproducing kernel Hilbert
space with a unique reproducing kernel function, and consequently there is a
unique Riesz representer of the functional $f\rightarrow f\left(  x\right)  $
which we denote by $\mathfrak{R}_{V,x}$ i.e.%
\begin{equation}
f\left(  x\right)  =\left(  f,\mathfrak{R}_{V,x}\right)  _{V,w,\theta},\quad
f\in X_{w}^{\theta},\text{ }x\in\mathbb{R}^{d}.\label{h59}%
\end{equation}

The equations%
\begin{equation}
\left(  f,R_{x}\right)  _{w,\theta}=f\left(  x\right)  =\left(  f,\mathfrak{R}%
_{V,x}\right)  _{V,w,\theta}=\left(  \mathcal{L}_{X}f,\mathcal{L}%
_{X}\mathfrak{R}_{V,x}\right)  _{V}=\left(  f,\mathcal{L}_{X}^{\ast
}\mathcal{L}_{X}\mathfrak{R}_{V,x}\right)  _{w,\theta},\label{h24}%
\end{equation}

imply that $\mathcal{L}_{X}^{\ast}\mathcal{L}_{X}\mathfrak{R}_{V,x}=R_{x}$ and
by Theorem \ref{Thm_Lx*Lx_onto_Xw,th} the operator $\mathcal{L}_{X}^{\ast
}\mathcal{L}_{X}:X_{w}^{\theta}\rightarrow X_{w}^{\theta}$ is a homeomorphism
so we have%
\begin{equation}
\mathfrak{R}_{V,x}=\left(  \mathcal{L}_{X}^{\ast}\mathcal{L}_{X}\right)
^{-1}R_{x}.\label{h57}%
\end{equation}

The equations \ref{h24} also imply there exists a unique $R_{V,x}\in V$ such
that%
\begin{equation}
f\left(  x\right)  =\left(  \mathcal{L}_{X}f,R_{V,x}\right)  _{V},\quad
R_{V,x}=\mathcal{L}_{X}\mathfrak{R}_{V,x},\quad R_{x}=\mathcal{L}_{X}^{\ast
}R_{V,x}.\label{h58}%
\end{equation}

As mentioned above we are interested in estimating the pointwise error of the
Exact smoother $s_{e}$ with respect to its data function $f_{d}$, where the
independent data is $X$ and the dependent data is $y=\widetilde{\mathcal{E}%
}_{X}f_{d}$. By \ref{h58} the error is $s_{e}\left(  x\right)  -f_{d}\left(
x\right)  =\left(  \mathcal{L}_{X}\left(  s_{e}-f_{d}\right)  ,R_{V,x}\right)
_{V}$ and%
\[
\left\vert s_{e}\left(  x\right)  -f_{d}\left(  x\right)  \right\vert
=\left\vert \left(  \mathcal{L}_{X}\left(  s_{e}-f_{d}\right)  ,R_{V,x}%
\right)  _{V}\right\vert \leq\left\Vert \mathcal{L}_{X}\left(  s_{e}%
-f_{d}\right)  \right\Vert _{V}\left\Vert R_{V,x}\right\Vert _{V}.
\]

From part 3 of Theorem \ref{Thm_smooth_Exact}%
\[
\left\Vert \mathcal{L}_{X}s_{e}-\varsigma\right\Vert _{V}^{2}+\left\Vert
\mathcal{L}_{X}\left(  s_{e}-f\right)  \right\Vert _{V}^{2}=\left\Vert
\mathcal{L}_{X}f-\varsigma\right\Vert _{V}^{2},\quad f\in X_{w}^{\theta},
\]

where $\zeta=\left(  0,y\right)  =\left(  0,\widetilde{\mathcal{E}}_{X}%
f_{d}\right)  $. Thus when $f=f_{d}$%
\[
\left\Vert \mathcal{L}_{X}\left(  s_{e}-f_{d}\right)  \right\Vert _{V}%
\leq\left\Vert \mathcal{L}_{X}f_{d}-\left(  0,\widetilde{\mathcal{E}}_{X}%
f_{d}\right)  \right\Vert _{V}=\left\vert f_{d}\right\vert _{w,\theta}%
\sqrt{\rho},
\]

and%
\[
\left\vert s_{e}\left(  x\right)  -f_{d}\left(  x\right)  \right\vert
\leq\left\vert f_{d}\right\vert _{w,\theta}\sqrt{\rho}\left\Vert
R_{V,x}\right\Vert _{V}.
\]

Finally, using \ref{h24} we have
\[
\left\Vert R_{V,x}\right\Vert _{V}^{2}=\left(  R_{V,x},R_{V,x}\right)
_{V}=\left(  \mathcal{L}_{X}\mathfrak{R}_{V,x},\mathcal{L}_{X}\mathfrak{R}%
_{V,x}\right)  _{V}=\mathfrak{R}_{V,x}\left(  x\right)  .
\]

The above analysis is summarized as:

\begin{theorem}
\label{Thm_pt_estim_Exact_smth_Rvx}Suppose $s_{e}$ is the Exact smoother
generated by the independent data $X$ and the data function $f_{d}\in
X_{w}^{\theta}$. Then%
\begin{equation}
\left\vert s_{e}\left(  x\right)  -f_{d}\left(  x\right)  \right\vert
\leq\left\vert f_{d}\right\vert _{w,\theta}\sqrt{\rho}\left\Vert
R_{V,x}\right\Vert _{V},\quad x\in\mathbb{R}^{N},\label{h60}%
\end{equation}

where $\rho>0$ is the smoothing coefficient, $R_{V,x}\in V$ is given by
\ref{h57} and \ref{h58} and satisfies
\begin{equation}
\left\Vert R_{V,x}\right\Vert _{V}^{2}=\mathfrak{R}_{V,x}\left(  x\right)
.\label{h64}%
\end{equation}

\end{theorem}

It is clear from Definition \ref{Def_Exact_smth_problem} of the Exact
smoothing problem that the Exact smoother is independent of the order of the
points in the unisolvent independent data $X$. This allows the convenient
definition of a special minimal unisolvent set $X_{1}\subset X$: we assume
that $X_{1}=\left\{  x^{\left(  i\right)  }\right\}  _{i=1}^{M}$ is a minimal
unisolvent subset and let $X_{2}=\left\{  x^{\left(  i\right)  }\right\}
_{i=M+1}^{N}$. The set $X_{1}$ is then used to construct the Riesz representer
$R_{x}$, the semi-Riesz representer $r_{x}=\mathcal{Q}R_{x}$, the Lagrangian
operators $\mathcal{P}$, $\mathcal{Q}$, the basis function spaces $W_{G,X}%
$\ and $\overset{\cdot}{W}_{G,X}$, and the Light norm $\left(  \cdot
,\cdot\right)  _{w,\theta}$. This section will also use the
`prime-double-prime' notation to denote components related to $X_{1}$ and
$X_{2}$ e.g. $\alpha=\left(  \alpha^{\prime},\alpha^{\prime\prime}\right)  $
with $\alpha^{\prime}\in\mathbb{C}^{M}$ and $\alpha^{\prime\prime}%
\in\mathbb{C}^{N-M}$.

The next step is to calculate $\left\Vert R_{V,x}\right\Vert _{V}$ using
\ref{p64}. From \ref{h57}, $\mathcal{L}_{X}^{\ast}\mathcal{L}_{X}%
\mathfrak{R}_{V,x}=R_{x}$ and from part 4 Theorem \ref{Thm_L_op_properties},
$\mathcal{L}_{X}^{\ast}\mathcal{L}_{X}=\rho\mathcal{Q}+\frac{1}{N}%
\widetilde{\mathcal{E}}_{X}^{\ast}\widetilde{\mathcal{E}}_{X}$ so that\
\begin{equation}
\rho\mathcal{Q}\mathfrak{R}_{V,x}+\frac{1}{N}\widetilde{\mathcal{E}}_{X}%
^{\ast}\widetilde{\mathcal{E}}_{X}\mathfrak{R}_{V,x}=R_{x}.\label{h45}%
\end{equation}

We will solve equation \ref{h45} for $\mathfrak{R}_{V,x}$ by solving the
equivalent system%
\begin{align}
\mathcal{Q}v_{x}+\widetilde{\mathcal{E}}_{X}^{\ast}\alpha_{x}  &
=R_{x},\label{h78}\\
\frac{1}{N\rho}\widetilde{\mathcal{E}}_{X}v_{x}  & =\alpha_{x},\label{h26}%
\end{align}

for $\alpha_{x}\in\mathbb{C}^{N}$ and $v_{x}\in X_{w}^{\theta}$, so that%
\begin{equation}
\rho\mathfrak{R}_{V,x}=v_{x},\;R_{V,x}=\mathcal{L}_{X}\mathfrak{R}%
_{V,x}=\left(  \frac{1}{\rho}v_{x},N\alpha_{x}\right)  .\label{h46}%
\end{equation}

We first apply the operator $\mathcal{P}$ to \ref{h78}: by part 7 of Theorem
\ref{Thm_eval_op_properties}, $\mathcal{P}\widetilde{\mathcal{E}}_{X}^{\ast
}\alpha_{x}=\alpha_{x}^{T}L_{X}\widetilde{l}$ and applying $\mathcal{P}$ to
\ref{p940} gives $\mathcal{P}R_{x}=\widetilde{l}\left(  x\right)
^{T}\widetilde{l}$ so%
\[
\mathcal{P}\left(  \mathcal{Q}v_{x}+\widetilde{\mathcal{E}}_{X}^{\ast}%
\alpha_{x}\right)  =\alpha_{x}^{T}L_{X}\widetilde{l}=\widetilde{l}\left(
x\right)  ^{T}\widetilde{l},
\]

which implies the equivalent equations%
\begin{equation}
L_{X}^{T}\alpha_{x}=\widetilde{l}\left(  x\right)  ,\quad\alpha_{x}^{\prime
}=\widetilde{l}\left(  x\right)  -L_{X_{2}}^{T}\alpha_{x}^{\prime\prime
}.\label{h39}%
\end{equation}

Next apply the operator $\mathcal{Q}$ to \ref{h78}: since $\mathcal{Q}%
R_{x}=r_{x}$
\[
\mathcal{Q}\left(  v_{x}+\widetilde{\mathcal{E}}_{X}^{\ast}\alpha_{x}\right)
=\mathcal{Q}v_{x}+\mathcal{Q}\sum\limits_{k=1}^{N}\left(  \alpha_{x}\right)
_{k}R_{x^{\left(  k\right)  }}=\mathcal{Q}v_{x}+\sum\limits_{k=M+1}^{N}\left(
\alpha_{x}\right)  _{k}r_{x^{\left(  k\right)  }}=r_{x},
\]

and thus%
\begin{equation}
\mathcal{Q}v_{x}=r_{x}-\sum\limits_{k=M+1}^{N}\left(  \alpha_{x}\right)
_{k}r_{x^{\left(  k\right)  }}.\label{h63}%
\end{equation}

We now left-compose $\widetilde{\mathcal{E}}_{X_{2}}$ with \ref{h63}. From
part 8 Theorem \ref{Thm_eval_op_properties}, $\widetilde{\mathcal{E}}_{X_{2}%
}\mathcal{P}f=L_{X}\widetilde{\mathcal{E}}_{X_{1}}f$ so that applying
$\widetilde{\mathcal{E}}_{X_{2}}$ to the left of \ref{h63} and then using
\ref{p26} gives%
\[
\widetilde{\mathcal{E}}_{X_{2}}\mathcal{Q}v_{x}=\widetilde{\mathcal{E}}%
_{X_{2}}\left(  v_{x}-\mathcal{P}v_{x}\right)  =\widetilde{\mathcal{E}}%
_{X_{2}}v_{x}-L_{X_{2}}\widetilde{\mathcal{E}}_{X_{1}}v_{x}=N\rho\left(
\alpha_{x}^{\prime\prime}-L_{X_{2}}\alpha_{x}^{\prime}\right)  .
\]

Next left-compose $\widetilde{\mathcal{E}}_{X_{2}}$ with \ref{h63} and use the
notation described in Definition \ref{Def_rX,x_rX,Y} to obtain
\[
\widetilde{\mathcal{E}}_{X_{2}}\mathcal{Q}v_{x}=\widetilde{\mathcal{E}}%
_{X_{2}}\left(  r_{x}-\sum\limits_{k=M+1}^{N}\left(  \alpha_{x}\right)
_{k}r_{x^{\left(  k\right)  }}\right)  =r_{X_{2},x}-r_{X_{2},X_{2}}\alpha
_{x}^{\prime\prime}.
\]

The last two sequences of equations yield%
\[
\left(  N\rho I+r_{X_{2},X_{2}}\right)  \alpha_{x}^{\prime\prime}=N\rho
L_{X_{2}}\alpha_{x}^{\prime}+r_{X_{2},x},
\]

so that substituting for $\alpha_{x}^{\prime}$ in the last equation using
\ref{h39} and then rearranging gives%
\begin{equation}
\left(  N\rho\left(  I+L_{X_{2}}L_{X_{2}}^{T}\right)  +r_{X_{2},X_{2}}\right)
\alpha_{x}^{\prime\prime}=N\rho L_{X_{2}}\widetilde{l}\left(  x\right)
+r_{X_{2},x}.\label{h49}%
\end{equation}

We know from Theorem \ref{Thm_rx(y)_mat_vect} that $r_{X_{2},X_{2}}$ is
positive definite and since the cardinal unisolvency matrix $L_{X_{2}}$ is
real-valued it follows that $N\rho\left(  I+L_{X_{2}}L_{X_{2}}^{T}\right)
+r_{X_{2},X_{2}}$ is also positive definite and hence regular. Thus \ref{h49}
and \ref{h39} define $\alpha_{x}$ uniquely and which means \ref{h78} defines
$v$ up to a polynomial of order $\theta$. However \ref{h26} actually supplies
the information which defines $v_{x}$ uniquely. Indeed, from \ref{h26},
$\widetilde{\mathcal{E}}_{X_{1}}v_{x}=N\rho\alpha_{x}^{\prime}$ so that%
\begin{align*}
\mathcal{P}v_{x}\left(  y\right)  =N\rho\text{ }\left(  \alpha_{x}^{\prime
}\right)  ^{T}\widetilde{l}\left(  y\right)  =N\rho\left(  \widetilde
{l}\left(  x\right)  -L_{X_{2}}^{T}\alpha_{x}^{\prime\prime}\right)
^{T}\widetilde{l}\left(  y\right)   &  =N\rho\left(  \widetilde{l}\left(
x\right)  ^{T}-\left(  \alpha_{x}^{\prime\prime}\right)  ^{T}L_{X_{2}}\right)
\widetilde{l}\left(  y\right) \\
&  =N\rho\widetilde{l}\left(  x\right)  ^{T}\widetilde{l}\left(  y\right)
-N\rho\left(  \alpha_{x}^{\prime\prime}\right)  ^{T}L_{X_{2}}\widetilde
{l}\left(  y\right)  ,
\end{align*}

and as a consequence of \ref{h63}%
\begin{align}
v_{x}\left(  y\right)   &  =\mathcal{P}v_{x}\left(  y\right)  +\mathcal{Q}%
v_{x}\left(  y\right) \nonumber\\
&  =N\rho\widetilde{l}\left(  x\right)  ^{T}\widetilde{l}\left(  y\right)
-N\rho\left(  \alpha_{x}^{\prime\prime}\right)  ^{T}L_{X_{2}}\widetilde
{l}\left(  y\right)  +\mathcal{Q}v_{x}\left(  y\right) \nonumber\\
&  =N\rho\widetilde{l}\left(  x\right)  ^{T}\widetilde{l}\left(  y\right)
-N\rho\left(  \alpha_{x}^{\prime\prime}\right)  ^{T}L_{X_{2}}\widetilde
{l}\left(  y\right)  +r_{x}\left(  y\right)  -\sum\limits_{k=M+1}^{N}\left(
\alpha_{x}\right)  _{k}r_{x^{\left(  k\right)  }}\left(  y\right) \nonumber\\
&  =r_{x}\left(  y\right)  +N\rho\widetilde{l}\left(  x\right)  ^{T}%
\widetilde{l}\left(  y\right)  -N\rho\left(  \alpha_{x}^{\prime\prime}\right)
^{T}L_{X_{2}}\widetilde{l}\left(  y\right)  -\left(  \alpha_{x}^{\prime\prime
}\right)  ^{T}\overline{r_{X_{2},y}}\nonumber\\
&  =r_{x}\left(  y\right)  +N\rho\widetilde{l}\left(  x\right)  ^{T}%
\widetilde{l}\left(  y\right)  -\left(  \alpha_{x}^{\prime\prime}\right)
^{T}\left(  N\rho L_{X_{2}}\widetilde{l}\left(  y\right)  +\overline
{r_{X_{2},y}}\right) \nonumber\\
&  =r_{x}\left(  y\right)  +N\rho\widetilde{l}\left(  x\right)  ^{T}%
\widetilde{l}\left(  y\right)  -\left(  \alpha_{x}^{\prime\prime}\right)
^{T}\left(  N\rho\text{ }\left(  I+L_{X_{2}}L_{X_{2}}^{T}\right)
+\overline{r_{X_{2},X_{2}}}\right)  \overline{\alpha_{y}^{\prime\prime}%
}\label{h62}\\
&  =r_{x}\left(  y\right)  +N\rho\widetilde{l}\left(  x\right)  ^{T}%
\widetilde{l}\left(  y\right)  -\nonumber\\
&  -\left(  N\rho L_{X_{2}}\widetilde{l}\left(  x\right)  +r_{X_{2},x}\right)
^{T}\left(  N\rho\text{ }\left(  I+L_{X_{2}}L_{X_{2}}^{T}\right)
+\overline{r_{X_{2},X_{2}}}\right)  ^{-1}\left(  N\rho L_{X_{2}}\widetilde
{l}\left(  y\right)  +\overline{r_{X_{2},y}}\right)  ,\label{h56}%
\end{align}

where the last two equations were derived using \ref{h49}. Since $N\rho\left(
I+L_{X_{2}}L_{X_{2}}^{T}\right)  +r_{X_{2},X_{2}}$ is positive definite, its
complex conjugate inverse is positive definite and when $y=x$ equations
\ref{h64}, \ref{h62} and \ref{h46} imply
\[
\rho\left\Vert R_{V,x}\right\Vert _{V}^{2}=\rho\mathfrak{R}_{V,x}\left(
x\right)  =v_{x}\left(  x\right)  \leq r_{x}\left(  x\right)  +N\rho\left\vert
\widetilde{l}\left(  x\right)  \right\vert ^{2},
\]

The next theorem summarizes these results:

\begin{theorem}
\label{Thm_Rvx_solution}Suppose $s_{e}$ is the Exact smoother generated by the
independent data $X$ and the data function $f_{d}\in X_{w}^{\theta}$. Then%
\begin{equation}
\rho\left\Vert R_{V,x}\right\Vert _{V}^{2}\leq r_{x}\left(  x\right)
+N\rho\left\vert \widetilde{l}\left(  x\right)  \right\vert ^{2},\label{h69}%
\end{equation}

and%
\[
\left\vert s_{e}\left(  x\right)  -f_{d}\left(  x\right)  \right\vert
\leq\left\vert f_{d}\right\vert _{w,\theta}\sqrt{r_{x}\left(  x\right)
+N\rho\left\vert \widetilde{l}\left(  x\right)  \right\vert ^{2}},\quad
x\in\mathbb{R}^{N}.
\]

\end{theorem}

To study the convergence of the Exact smoother we will also need Lemma
\ref{Lem_int_Lagrange_interpol} which supplies some elementary results from
the theory of Lagrange interpolation and was used in Chapter 4 to derive
orders of convergence for the interpolant. These results are stated without
proof. This lemma has been created from Lemma 3.2, Lemma 3.5 and the first two
paragraphs of the proof of Theorem 3.6 of Light and Wayne
\cite{LightWayne98PowFunc}. The results of this lemma do not involve any
reference to weight or basis functions or functions in $X_{w}^{\theta}$, but
consider the properties of the set which contains the independent data points
and the order of the unisolvency used for the interpolation. Thus we have
separated the part of the proof that involves weight functions from the part
that uses the detailed theory of Lagrange interpolation operators.

\begin{lemma}
\label{Lem_Lagrange_interpol_2}(Copy of Lemma \ref{Lem_int_Lagrange_interpol})
Suppose that:

\begin{enumerate}
\item $\Omega$ is a bounded, open, connected subset of $\mathbb{R}^{d}$ having
the cone property.

\item $X$ is a unisolvent subset of $\Omega$ of order $\theta$.

\item $\left\{  l_{j}\right\}  _{j=1}^{M}$ is the cardinal basis of
$P_{\theta-1}$ with respect to a minimal unisolvent set of$\ \Omega$.

Then by using Lagrange interpolation techniques it can be shown there exists a
constant $K_{\Omega,\theta}^{\prime}>0$ such that
\begin{equation}
\left\vert \widetilde{l}\left(  x\right)  \right\vert _{1}=\sum\limits_{j=1}%
^{M}\left\vert l_{j}\left(  x\right)  \right\vert \leq K_{\Omega,\theta
}^{\prime},\quad x\in\overline{\Omega},\label{h09}%
\end{equation}

\end{enumerate}

and all minimal unisolvent subsets of $\Omega$. Now define
\[
h_{X}=\sup\limits_{\omega\in\Omega}\operatorname*{dist}\left(  \omega
,X\right)  ,
\]

and fix $x\in X$. Again using Lagrange interpolation techniques it can be
shown there are constants $c_{\Omega,\theta},h_{\Omega,\theta}>0$ such that
when $h_{X}<h_{\Omega,\theta}$ there exists a minimal unisolvent set $A\subset
X$ satisfying
\begin{equation}
\operatorname*{diam}A_{x}\leq c_{\Omega,\theta}h_{X},\label{h08}%
\end{equation}

where $A_{x}=A\cup\left\{  x\right\}  $.
\end{lemma}

We will now prove our first Exact smoother error estimate Theorem
\ref{Thm_Exact_smth_converg_at_x_1} which implies an order of convergence
$\eta$. In Theorem \ref{Thm_int_rx(x)_bound} conditions on the basis function
$G$ were supplied in order that the semi-Riesz representer $r_{x}\left(
y\right)  =\mathcal{Q}_{y}R_{x}\left(  y\right)  $ satisfy a condition of the
form \ref{h103} required below by Theorem \ref{Thm_Exact_smth_converg_at_x_1}.
Then in Theorem \ref{Thm_interpol_converg} an estimate for the interpolant
error was derived which is just the Exact smoother error estimate \ref{h104}
with $\rho=0$.

\begin{theorem}
\label{Thm_Exact_smth_converg_at_x_1}Let $w$ be a weight function with
properties W2 and W3 for order $\theta$ and parameter $\kappa$ and set
$\eta=\min\left\{  \theta,\frac{1}{2}\left\lfloor \min2\kappa\right\rfloor
\right\}  $. Suppose the notation and assumptions of Lemma
\ref{Lem_Lagrange_interpol_2} hold. Assume there exist constants $c_{G}%
,r_{G}>0$, independent of $A$ and $x\in\Omega$ such that%
\begin{equation}
\sqrt{r_{x}\left(  x\right)  }\leq\left(  1+\left\vert \widetilde{l}\left(
x\right)  \right\vert _{1}\right)  \sqrt{c_{G}}\left(  \operatorname*{diam}%
A_{x}\right)  ^{\eta},\quad\operatorname*{diam}A_{x}<r_{G},\text{ }x\in
\Omega,\label{h103}%
\end{equation}

where $A_{x}=A\cup\left\{  x\right\}  $.\smallskip

Then there exist constants $c_{\Omega,\theta},h_{\Omega,\theta},K_{\Omega
,\theta}^{\prime}>0$ such that%
\begin{equation}
\left\vert s_{e}\left(  x\right)  -f_{d}\left(  x\right)  \right\vert
\leq\left\vert f_{d}\right\vert _{w,\theta}\left(  1+K_{\Omega,\theta}%
^{\prime}\right)  \left(  \sqrt{c_{G}}\left(  c_{\Omega,\theta}h_{X}\right)
^{\eta}+\sqrt{N_{X}\rho}\right)  ,\quad x\in\overline{\Omega},\label{h104}%
\end{equation}

when $h_{X}\leq\min\left\{  h_{\Omega,\theta},r_{G}\right\}  $. Here
$N_{X}=\left\vert X\right\vert $ and the constants $c_{\Omega,\theta
},h_{\Omega,\theta},K_{\Omega,\theta}^{\prime}$ only depend on $\Omega
,\theta,\kappa$ and $d$.
\end{theorem}

\begin{proof}
From \ref{h69} of Theorem \ref{Thm_Rvx_solution}%
\begin{equation}
\sqrt{\rho}\left\Vert R_{V,x}\right\Vert _{V}\leq\sqrt{r_{x}\left(  x\right)
}+\left\vert \widetilde{l}\left(  x\right)  \right\vert _{1}\sqrt{N_{X}\rho
},\quad x\in\mathbb{R}^{d}.\label{h541}%
\end{equation}

Fix $x\in\mathbb{\Omega}$ and let $A$ be any minimal unisolvent subset of $X$.
Define $r_{x}$ using $A$ so that by \ref{h103} and Lemma
\ref{Lem_Lagrange_interpol_2} we have
\begin{equation}
\sqrt{r_{x}\left(  x\right)  }\leq\sqrt{c_{G}}\left(  1+\left\vert
\widetilde{l}\left(  x\right)  \right\vert _{1}\right)  \left(
\operatorname*{diam}A_{x}\right)  ^{\eta}\leq\left(  1+K_{\Omega,\theta
}^{\prime}\right)  \sqrt{c_{G}}\left(  \operatorname*{diam}A_{x}\right)
^{\eta+\delta_{G}},\label{h071}%
\end{equation}

when $\operatorname*{diam}A_{x}\leq r_{G}$. But $h_{X}\leq h_{\Omega,\theta}$
so that by Lemma \ref{Lem_Lagrange_interpol_2} there exists $A$ such that
$\operatorname*{diam}A_{x}\leq c_{\Omega,\theta}h_{X}$ and \ref{h071} implies%
\[
\sqrt{r_{x}\left(  x\right)  }\leq\left(  1+K_{\Omega,\theta}^{\prime}\right)
\sqrt{c_{G}}\left(  c_{\Omega,\theta}h_{X}\right)  ^{\eta}.
\]

Applying this estimate and the inequality \ref{h09} for $\left\vert
\widetilde{l}\left(  x\right)  \right\vert _{1}$ to the inequality \ref{h541}
for $R_{V,x}$ we get%
\[
\sqrt{\rho}\left\Vert R_{V,x}\right\Vert _{V}\leq\left(  1+K_{\Omega,\theta
}^{\prime}\right)  \left(  \sqrt{c_{G}}\left(  c_{\Omega,\theta}h_{X}\right)
^{\eta}+\sqrt{N_{X}\rho}\right)  ,\quad x\in\Omega,
\]

and the smoother error estimate \ref{h60} now becomes%
\begin{align*}
\left\vert s_{e}\left(  x\right)  -f_{d}\left(  x\right)  \right\vert  &
\leq\left\vert f_{d}\right\vert _{w,\theta}\sqrt{\rho}\left\Vert
R_{V,x}\right\Vert _{V}\\
& \leq\left\vert f_{d}\right\vert _{w,\theta}\left(  1+K_{\Omega,\theta
}^{\prime}\right)  \left(  \sqrt{c_{G}}\left(  c_{\Omega,\theta}h_{X}\right)
^{\eta}+\sqrt{N_{X}\rho}\right)  ,\quad x\in\Omega.
\end{align*}

Since $s_{e}$ and $f_{d}$ are continuous on $\mathbb{R}^{d}$ the last
inequality actually holds for $x\in\overline{\Omega}$ and \ref{h104} is true.
\end{proof}

The next result shows that in a certain sense the Exact smoother converges to
its data function.

\begin{corollary}
\label{Cor_Exact_smth_converg}The Exact smoother converges to its data
function $f_{d}$ in the sense that given $\varepsilon>0$ there exists a
positive integer $K\left(  f_{d};\varepsilon\right)  $, a nested sequence of
independent data sets $X^{\left(  k\right)  }\subset X^{\left(  k+1\right)
}\subset\Omega$ and a sequence of smoothing parameters $\rho_{k}>0$ such that
the corresponding sequence of Exact smoothers $s_{e}^{\left(  k\right)  }$
satisfies
\begin{equation}
\left\vert s_{e}^{\left(  k\right)  }\left(  x\right)  -f_{d}\left(  x\right)
\right\vert \leq\varepsilon,\quad x\in\overline{\Omega},\label{h13}%
\end{equation}

when $k\geq K\left(  f_{d};\varepsilon\right)  $.
\end{corollary}

\begin{proof}
From Theorem \ref{Thm_seq_data_regions_2} there exists a nested sequence of
independent data sets $X^{\left(  k\right)  }$ such that $X^{\left(  k\right)
}\subset X^{\left(  k+1\right)  }\subset\Omega$ and $h_{X^{\left(  k\right)
}}\rightarrow0$ as $k\rightarrow\infty$. To derive the bound \ref{h13} we
start with inequality \ref{h104} so that%
\[
\left\vert s_{e}^{\left(  k\right)  }\left(  x\right)  -f_{d}\left(  x\right)
\right\vert \leq\left\vert f_{d}\right\vert _{w,\theta}\left(  1+K_{\Omega
,\theta}^{\prime}\right)  \left(  \sqrt{c_{G}}\left(  c_{\Omega,\theta
}h_{X^{\left(  k\right)  }}\right)  ^{\eta_{G}}+\sqrt{N_{X^{\left(  k\right)
}}\rho_{k}}\right)  ,\quad x\in\overline{\Omega},
\]

when $h_{X^{\left(  k\right)  }}<\min\left\{  h_{\Omega,\theta},r_{G}\right\}
$. Now observe that for some positive integer $K_{0}\left(  f_{d}\right)  $,
$h_{X^{\left(  k\right)  }}<h_{\Omega,\theta}$ when $k\geq K_{0}\left(
f_{d}\right)  $. Next, choosing $\rho_{k}$ such that $\sqrt{c_{G}}\left(
c_{\Omega,\theta}h_{X^{\left(  k\right)  }}\right)  ^{\eta_{G}}=\sqrt
{N_{X^{\left(  k\right)  }}\rho_{k}}$ implies%
\[
\left\vert s_{e}^{\left(  k\right)  }\left(  x\right)  -f_{d}\left(  x\right)
\right\vert \leq2\left\vert f_{d}\right\vert _{w,\theta}\left(  1+K_{\Omega
,\theta}^{\prime}\right)  \sqrt{c_{G}}\left(  c_{\Omega,\theta}h_{X^{\left(
k\right)  }}\right)  ^{\eta_{G}},x\in\overline{\Omega},
\]
and as $h_{X^{\left(  k\right)  }}\rightarrow0$, $s_{e}\rightarrow f_{d}$
uniformly on $\overline{\Omega}$. The statement of this corollary now follows directly.
\end{proof}

\section{Improved error estimates\label{Sect_improved_err_estim}}

In this section we start by deriving an estimate which has a slightly improved
order of convergence: for the example of the shifted thin-plate spline this is
1/2. Then a double rate of convergence will be demonstrated for data functions
that are linear combinations of Riesz representers. As a function of the
smoothing parameter $\rho$ these convergence estimates are unbounded and we
end this chapter by deriving some error estimates bounded in $\rho$.

\subsection{A slightly increased order of convergence}

Our slightly improved Exact smoother error estimate will be Theorem
\ref{Thm_Exact_smth_converg_at_x}. In Theorem
\ref{Thm_int_rx(x)_bnd_better_order} conditions on the basis function $G$ were
supplied in order that the semi-Riesz representer $r_{x}\left(  y\right)
=\mathcal{Q}_{y}R_{x}\left(  y\right)  $ satisfy a condition of the form
\ref{h102} (below) which allows for an improved order of convergence of
$\eta+\delta_{G}$. Then in Corollary \ref{Cor_Thm_int_rx(x)_bnd_better_order}
an estimate for the interpolant error was derived which is just the improved
Exact smoother error estimate \ref{h40} (below) with $\rho=0$.

\begin{theorem}
\label{Thm_Exact_smth_converg_at_x}Let $w$ be a weight function with
properties W2 and W3 for order $\theta$ and parameter $\kappa$ and set
$\eta=\min\left\{  \theta,\frac{1}{2}\left\lfloor \min2\kappa\right\rfloor
\right\}  $. Suppose the notation and assumptions of Lemma
\ref{Lem_Lagrange_interpol_2} hold. Assume there exist constants $c_{G}%
,r_{G}>0$ and $\delta_{G}\geq0$, independent of $A$ and $x\in\Omega$, such
that%
\begin{equation}
\sqrt{r_{x}\left(  x\right)  }\leq\left(  1+\left\vert \widetilde{l}\left(
x\right)  \right\vert _{1}\right)  \sqrt{c_{G}}\left(  \operatorname*{diam}%
A_{x}\right)  ^{\eta+\delta_{G}},\quad\operatorname*{diam}A_{x}<r_{G},\text{
}x\in\Omega,\label{h102}%
\end{equation}

where $A_{x}=A\cup\left\{  x\right\}  $.\smallskip

Then there exist constants $c_{\Omega,\theta},h_{\Omega,\theta},K_{\Omega
,\theta}^{\prime}>0$ such that%
\begin{equation}
\left\vert s_{e}\left(  x\right)  -f_{d}\left(  x\right)  \right\vert
\leq\left\vert f_{d}\right\vert _{w,\theta}\left(  1+K_{\Omega,\theta}%
^{\prime}\right)  \left(  \sqrt{c_{G}}\left(  c_{\Omega,\theta}h_{X}\right)
^{\eta+\delta_{G}}+\sqrt{N_{X}\rho}\right)  ,\quad x\in\overline{\Omega
},\label{h40}%
\end{equation}

when $h_{X}\leq\min\left\{  h_{\Omega,\theta},r_{G}\right\}  $. Here
$N_{X}=\left\vert X\right\vert $ and the constants $c_{\Omega,\theta
},h_{\Omega,\theta},K_{\Omega,\theta}^{\prime}$ only depend on $\Omega
,\theta,\kappa$ and $d$.
\end{theorem}

\begin{proof}
The will follow closely that of Theorem \ref{Thm_Exact_smth_converg_at_x_1}.
From \ref{h69} of Theorem \ref{Thm_Rvx_solution}%
\begin{equation}
\sqrt{\rho}\left\Vert R_{V,x}\right\Vert _{V}\leq\sqrt{r_{x}\left(  x\right)
}+\left\vert \widetilde{l}\left(  x\right)  \right\vert _{1}\sqrt{N_{X}\rho
},\quad x\in\mathbb{R}^{d}.\label{h54}%
\end{equation}

Fix $x\in\mathbb{\Omega}$ and let $A$ be any minimal unisolvent subset of $X$.
Define $r_{x}$ using $A$ so that by \ref{h102} and Lemma
\ref{Lem_Lagrange_interpol_2} we have
\begin{equation}
\sqrt{r_{x}\left(  x\right)  }\leq\sqrt{c_{G}}\left(  1+\left\vert
\widetilde{l}\left(  x\right)  \right\vert _{1}\right)  \left(
\operatorname*{diam}A_{x}\right)  ^{\eta+\delta_{G}}\leq\left(  1+K_{\Omega
,\theta}^{\prime}\right)  \sqrt{c_{G}}\left(  \operatorname*{diam}%
A_{x}\right)  ^{\eta+\delta_{G}},\label{h07}%
\end{equation}

when $\operatorname*{diam}A_{x}\leq r_{G}$. But $h_{X}\leq h_{\Omega,\theta}$
so that by Lemma \ref{Lem_Lagrange_interpol_2} there exists $A$ such that
$\operatorname*{diam}A_{x}\leq c_{\Omega,\theta}h_{X}$ and \ref{h07} implies%
\begin{equation}
\sqrt{r_{x}\left(  x\right)  }\leq\left(  1+K_{\Omega,\theta}^{\prime}\right)
\sqrt{c_{G}}\left(  c_{\Omega,\theta}h_{X}\right)  ^{\eta+\delta_{G}%
}.\label{h03}%
\end{equation}

Applying this estimate and the inequality \ref{h09} for $\left\vert
\widetilde{l}\left(  x\right)  \right\vert _{1}$ to the inequality \ref{h54}
for $R_{V,x}$ we get%
\begin{equation}
\sqrt{\rho}\left\Vert R_{V,x}\right\Vert _{V}\leq\left(  1+K_{\Omega,\theta
}^{\prime}\right)  \left(  \sqrt{c_{G}}\left(  c_{\Omega,\theta}h_{X}\right)
^{\eta+\delta_{G}}+\sqrt{N_{X}\rho}\right)  ,\quad x\in\Omega,\label{h05}%
\end{equation}

and the smoother error estimate \ref{h60} now becomes%
\begin{align*}
\left\vert s_{e}\left(  x\right)  -f_{d}\left(  x\right)  \right\vert  &
\leq\left\vert f_{d}\right\vert _{w,\theta}\sqrt{\rho}\left\Vert
R_{V,x}\right\Vert _{V}\\
& \leq\left\vert f_{d}\right\vert _{w,\theta}\left(  1+K_{\Omega,\theta
}^{\prime}\right)  \left(  \sqrt{c_{G}}\left(  c_{\Omega,\theta}h_{X}\right)
^{\eta+\delta_{G}}+\sqrt{N_{X}\rho}\right)  ,\quad x\in\Omega.
\end{align*}

Since $s_{e}$ and $f_{d}$ are continuous on $\mathbb{R}^{d}$ the last
inequality actually holds for $x\in\overline{\Omega}$.
\end{proof}

In the examples below, to maximize the rate of convergence, we define
$\eta=\max\limits_{0\leq\kappa<s}\min\left\{  \theta,\frac{1}{2}\left\lfloor
\min2\kappa\right\rfloor \right\}  $.

\begin{example}
\textbf{Thin plate spline basis functions} By Theorem
\ref{Thm_surf_spline_ext_basis} the thin-plate spline weight functions are
given by%
\begin{equation}
w\left(  \xi\right)  =\frac{1}{e\left(  s\right)  }\left\vert \xi\right\vert
^{-2\theta+2s+d},\quad\xi\in\mathbb{R}^{d},\label{h76}%
\end{equation}

and have properties W2.1 and W3.2 with positive integer order $\theta$ and
non-negative $\kappa$ iff $\kappa<s<\theta$.

Now $\eta=\max\limits_{0\leq\kappa<s}\min\left\{  \theta,\frac{1}%
{2}\left\lfloor 2\kappa\right\rfloor \right\}  =\frac{1}{2}\left\lfloor
2\kappa\right\rfloor $ and from the conclusion of Example 1 of Section
\ref{Sect_better_results} Chapter \ref{Ch_Interpol}, equation \ref{h102} is
satisfied for $r_{G}=1$ and $\delta_{G}=\frac{1}{2}\left(  2s-\left\lfloor
2s\right\rfloor \right)  $ when $s>0$, $s\neq1,2,3,\ldots$ and for $r_{G}=1$
and any $0\leq\delta_{G}<\frac{1}{2}$ when $s=1,2,3,\ldots$.
\end{example}

\begin{example}
\textbf{Shifted thin plate spline basis functions} By Theorem
\ref{Thm_shift_surf_spline_ext_basis} the shifted thin-plate spline weight
functions are given by%
\[
w\left(  \xi\right)  =\frac{1}{\widetilde{e}\left(  s\right)  \widetilde
{K}_{s+d/2}\left(  a\left\vert \xi\right\vert \right)  }\left\vert
\xi\right\vert ^{-2\theta+2s+d},\quad s>-d/2,
\]

and have properties W2.1 and W3.2 for $\theta$ and all $\kappa\geq0$ iff
$-d/2<s<\theta$.

Now $\eta=\max\limits_{\kappa\geq0}\min\left\{  \theta,\frac{1}{2}\left\lfloor
2\kappa\right\rfloor \right\}  =\theta$ and from the conclusion of Example 2
of Section \ref{Sect_better_results} Chapter \ref{Ch_Interpol}, \ref{h102} is
satisfied for arbitrary $r_{G}>0$ and $\delta_{G}=\frac{1}{2}$.
\end{example}

\begin{remark}
\label{Rem_smooth_converg}\ 

\begin{enumerate}
\item \fbox{Relating $h_{X}$\ to $N_{X}$\ and the convergence of the smoother}

Theorem \ref{Thm_Exact_smth_converg_at_x} derives the smoother error estimate
\ref{h40}:
\[
\left\vert s_{e}\left(  x\right)  -f_{d}\left(  x\right)  \right\vert
\leq\left\vert f_{d}\right\vert _{w,\theta}\left(  1+K_{\Omega,\theta}%
^{\prime}\right)  \left(  \sqrt{c_{G}}\left(  c_{\Omega,\theta}h_{X}\right)
^{\eta+\delta_{G}}+\sqrt{N_{X}\rho}\right)  ,
\]

when $h_{X}\leq\min\left\{  h_{\Omega,\theta},r_{G}\right\}  $; with the
awkward term $N_{X}$ which is linked to $h_{X}$ in an indirect manner. In
order to have a meaningful concept of the convergence of a smoother a
relationship between $h_{X}$\ and $N_{X}$\ needs to be introduced. This
situation did not arise with interpolant.

\item In general the data $X$ is scattered and for a given value of $h_{X}$
the number of points $N_{X}$ can be arbitrarily large. However, for data on a
regular, rectangular grid we have the relation
\begin{equation}
d^{d/2}\operatorname*{vol}\left(  grid\right)  =N_{X}\left(  h_{X}\right)
^{d}.\label{h01}%
\end{equation}

\item In Williams \cite{WilliamsZeroOrdSmthV4} several \textbf{1-dimensional}
numerical experiments were run to compare the convergence of the \textbf{zero
order} Exact smoother with the predicted convergence. One-dimensional data
sets were constructed using a uniform distribution on the interval
$\Omega=\left[  -1.5,1.5\right]  $. Each of 20 data files were exponentially
sampled using a multiplier of approximately 1.2 and a maximum of 5000 points,
and then $\log_{10}h_{X}$ was plotted against $\log_{10}N$ where $N=\left\vert
X\right\vert $. It then seemed quite reasonable to use a least-squares linear
fit and in this case we obtained the relation%
\begin{equation}
h_{X}\simeq3.09N^{-0.81}.\label{h71}%
\end{equation}

For ease of calculation let
\begin{equation}
h_{X}=h_{1}\left(  N_{X}\right)  ^{-a},\text{\quad}h_{1}=3.09,\text{\quad
}a=0.81.\label{h74}%
\end{equation}

\item A barrier to the use of such a formula as \ref{h71} in higher dimensions
is the difficulty of actually calculating $h_{X}$ for a given data set. If a
sequence of\textbf{\ }independent test data sets was generated by a uniform
distribution in each dimension then the constants $a$ and $h_{1}$ might be
defined as the upper bound of the confidence interval of a statistical
distribution. Also, noting regular grid formula \ref{h01} we might hypothesize
a relationship of the form%
\[
h_{X}=h_{d}\left(  N_{X}\right)  ^{-a_{d}d},
\]

for \textbf{higher dimensions}.

\item \fbox{The order of convergence} Assuming \ref{h74} we will now show that
for a sequence of independent data points $X_{k}$ of increasing density there
exists a sequence $\rho_{k}$ of smoothing coefficients such that the smoother
error is of order $\eta_{G}=\eta+\delta_{G}$ in $h_{X_{k}}$.

By Theorem \ref{Thm_Exact_smth_converg_at_x}%
\[
\left\vert f_{d}\left(  x\right)  -\mathcal{S}_{X_{k}}^{e}f_{d}\left(
x\right)  \right\vert \leq\left\vert f_{d}\right\vert _{w,\theta}\left(
1+K_{\Omega,\theta}^{\prime}\right)  \left(  \sqrt{c_{G}}\left(
c_{\Omega,\theta}h_{X_{k}}\right)  ^{\eta_{G}}+\sqrt{N_{X_{k}}\rho_{k}%
}\right)  ,\quad x\in\overline{\Omega},
\]

when $h_{X_{k}}<\min\left\{  h_{\Omega,\theta},r_{G}\right\}  $. For clarity
define the constants%
\begin{equation}
A=\left\vert f_{d}\right\vert _{w,\theta}\left(  1+K_{\Omega,\theta}^{\prime
}\right)  \sqrt{c_{G}}\left(  c_{\Omega,\theta}\right)  ^{\eta_{G}%
},\text{\quad}B=\left\vert f_{d}\right\vert _{w,\theta}\left(  1+K_{\Omega
,\theta}^{\prime}\right)  ,\label{h42}%
\end{equation}

so that%
\[
\left\vert f_{d}\left(  x\right)  -\mathcal{S}_{X_{k}}^{e}f_{d}\left(
x\right)  \right\vert \leq A\left(  h_{X_{k}}\right)  ^{\eta_{G}}+B\sqrt
{\rho_{k}}\sqrt{N_{X_{k}}},\quad x\in\overline{\Omega}.
\]

Then condition \ref{h74} implies%
\begin{equation}
\left\vert f_{d}\left(  x\right)  -\mathcal{S}_{X_{k}}^{e}f_{d}\left(
x\right)  \right\vert \leq A\left(  h_{X_{k}}\right)  ^{\eta_{G}}+B\sqrt
{\rho_{k}}\left(  \frac{h_{X_{k}}}{h_{1}}\right)  ^{-\frac{1}{2a}},\quad
x\in\overline{\Omega},\label{h06}%
\end{equation}

and we want to minimize the right side as a function of $h_{X_{k}}$. This is
easily done by setting the derivative to zero so that%
\begin{align*}
\left.  D_{x}\left(  Ax^{\eta_{G}}+B\sqrt{\rho_{k}}\left(  \frac{x}{h_{1}%
}\right)  ^{-\frac{1}{2a}}\right)  \right\vert _{x=h_{X_{k}}}  & =A\eta
_{G}\left(  h_{X_{k}}\right)  ^{\eta_{G}-1}-\frac{1}{2ah_{1}}B\sqrt{\rho_{k}%
}\left(  \frac{h_{X_{k}}}{h_{1}}\right)  ^{-\frac{1}{2a}-1}\\
& =0,
\end{align*}

and a unique minimum is obtained when%
\begin{equation}
\sqrt{\rho_{k}}=\frac{A}{B}\frac{2a\eta_{G}}{\left(  h_{1}\right)  ^{\frac
{1}{2a}}}\left(  h_{X_{k}}\right)  ^{\eta_{G}+\frac{1}{2a}}.\label{h02}%
\end{equation}

By substituting for $\sqrt{\rho_{k}}$ in the right side of \ref{h06} we find
that the corresponding \textbf{minimum error} value is given by%
\begin{align*}
\left\vert f_{d}\left(  x\right)  -\mathcal{S}_{X_{k}}^{e}f_{d}\left(
x\right)  \right\vert  & \leq A\left(  h_{X_{k}}\right)  ^{\eta_{G}}%
+B\sqrt{\rho_{k}}\left(  \frac{h_{X_{k}}}{h_{1}}\right)  ^{-\frac{1}{2a}}\\
& =A\left(  h_{X_{k}}\right)  ^{\eta_{G}}+B\left(  \frac{A}{B}\frac{2a\eta
_{G}}{\left(  h_{1}\right)  ^{\frac{1}{2a}}}\left(  h_{X_{k}}\right)
^{\eta_{G}+\frac{1}{2a}}\right)  \left(  \frac{h_{X_{k}}}{h_{1}}\right)
^{-\frac{1}{2a}}\\
& =A\left(  h_{X_{k}}\right)  ^{\eta_{G}}+2a\eta_{G}A\left(  h_{X_{k}}\right)
^{\eta_{G}}\\
& =\left(  1+2a\eta_{G}\right)  A\left(  h_{X_{k}}\right)  ^{\eta_{G}}\\
& =\left(  1+2a\eta_{G}\right)  \left\vert f_{d}\right\vert _{w,\theta}\left(
1+K_{\Omega,\theta}^{\prime}\right)  \sqrt{c_{G}}\left(  c_{\Omega,\theta
}h_{X_{k}}\right)  ^{\eta_{G}}.
\end{align*}

Thus if $h_{X_{k}}\rightarrow0$ and $\rho_{k}$ is chosen to satisfy \ref{h02}
then the Exact smoother error always converges uniformly to zero and the order
of convergence is at least $\left(  h_{X_{k}}\right)  ^{\eta+\delta_{G}}$.
Note also that the minimum error is independent of $h_{1}$.

\item From \ref{h02} the convergence of $\rho_{k}$ to zero is of order
$2\eta_{G}+\frac{1}{a}$ in $h_{X}$. One might consider saying that the smaller
this order is the better in order to avoid stability problems when $\rho$ is small?

\item \fbox{Error behavior in terms of $\rho$} The $\rho$ error curve is well
known, always becoming constant with large $\rho$ and decreasing as $\rho$
decreases and then it may increase again to the interpolant error value. Its
observed slope behavior depends on the basis function: with decreasing $\rho$:

\begin{itemize}
\item \textbf{Thin-plate spline} slope increases like \ref{h40} then levels out.

\item \textbf{Shifted thin-plate spline} slope keeps increasing like \ref{p40}.

\item \textbf{Gaussian} slope oscillates with increasing positive amplitude.

\item \textbf{B-splines} - when $n=l=1$ the slope tends to zero. For other
splines it doe not.
\end{itemize}

Note that observing the behavior of the smoother for small $\rho$ may run into
stability problems and the limitations of accuracy.

If we assume our $X$ data lies in $\mathbb{R}^{1}$ and that $h_{X}$ and
$N_{X}$ satisfy \ref{h74} then we obtain the error estimate \ref{h06}:
\[
\left\vert f_{d}\left(  x\right)  -\mathcal{S}_{X_{k}}^{e}f_{d}\left(
x\right)  \right\vert \leq A\left(  h_{X_{k}}\right)  ^{\eta_{G}}+B\sqrt{\rho
}\left(  \frac{h_{X_{k}}}{h_{1}}\right)  ^{-\frac{1}{2a}},\quad x\in
\overline{\Omega}.
\]

Clearly the slope always tends to infinity for small $\rho$ and always tends
to infinity for large $\rho$, which leaves much to be desired. However, in
Subsection \ref{SbSect_ErrEstimBndedInSmthParm} we will establish error
estimates that, for given $X$, are constant for large $\rho$.
\end{enumerate}
\end{remark}

\subsection{Doubled order of convergence}

In this subsection a double rate of convergence in terms of $h_{X}$ will be
demonstrated for data functions that are linear combinations of Riesz
representers. We will again use the assumptions and notation of Lemma
\ref{Lem_Lagrange_interpol_2} regarding the open set which contains the
independent data points. The convergence results will be applied to the case
of semi-regular sequences of independent data points.

The operator $\mathcal{L}_{X}$ of Definition \ref{Def_ex_Hilbert_smoothing} is
not onto but by Theorem \ref{Thm_L_op_properties} it is 1-1, and since
$\mathcal{L}_{X}\left(  X_{w}^{\theta}\right)  $ is not dense in $V$ the usual
definition of the adjoint does not apply to $\mathcal{L}_{X}^{-1}$. However,
by Theorem \ref{Thm_Lx*Lx_onto_Xw,th}, $\mathcal{L}_{X}^{\ast}\mathcal{L}%
_{X}:X_{w}^{\theta}\rightarrow X_{w}^{\theta}$ is 1-1 and onto and this
property is used in the next lemma to extend $\mathcal{L}_{X}^{-1}$ to $V$ as
a continuous operator and allows a definition of $\left(  \mathcal{L}_{X}%
^{-1}\right)  ^{\ast}$:

\begin{lemma}
\label{Lem_adjoint_inverseLx}\textbf{The operator }$\left(  \mathcal{L}%
_{X}^{-1}\right)  ^{\ast}$: The operator $\left(  \mathcal{L}_{X}^{\ast
}\mathcal{L}_{X}\right)  ^{-1}\mathcal{L}_{X}^{\ast}:V\rightarrow
X_{w}^{\theta}$ is a continuous extension of $\mathcal{L}_{X}^{-1}$ to $V$.
Define the operator $\left(  \mathcal{L}_{X}^{-1}\right)  ^{\ast}$ by $\left(
\mathcal{L}_{X}^{-1}\right)  ^{\ast}=\left(  \left(  \mathcal{L}_{X}^{\ast
}\mathcal{L}_{X}\right)  ^{-1}\mathcal{L}_{X}^{\ast}\right)  ^{\ast}$. Then
$\left(  \mathcal{L}_{X}^{-1}\right)  ^{\ast}$ has the following properties:

\begin{enumerate}
\item $\left(  \mathcal{L}_{X}^{-1}\right)  ^{\ast}=\mathcal{L}_{X}\left(
\mathcal{L}_{X}^{\ast}\mathcal{L}_{X}\right)  ^{-1}$.

\item $\left(  \mathcal{L}_{X}^{\ast}\mathcal{L}_{X}\right)  ^{-1}%
=\mathcal{L}_{X}^{-1}\left(  \mathcal{L}_{X}^{-1}\right)  ^{\ast}$.

\item $\mathcal{L}_{X}^{\ast}\left(  \mathcal{L}_{X}^{-1}\right)  ^{\ast}=I$.

\item $\left(  \mathcal{L}_{X}^{-1}\right)  ^{\ast\ast}=\left(  \mathcal{L}%
_{X}^{\ast}\mathcal{L}_{X}\right)  ^{-1}\mathcal{L}_{X}^{\ast}$.
\end{enumerate}
\end{lemma}

\begin{proof}
The operator $\left(  \mathcal{L}_{X}^{\ast}\mathcal{L}_{X}\right)
^{-1}\mathcal{L}_{X}^{\ast}:X_{w}^{\theta}\rightarrow V$ is a continuous
extension of $\mathcal{L}_{X}^{-1}$ to $V$ because

$\left(  \mathcal{L}_{X}^{\ast}\mathcal{L}_{X}\right)  ^{-1}\mathcal{L}%
_{X}^{\ast}\mathcal{L}_{X}=I$.\medskip

\textbf{Part 1} Since $\mathcal{L}_{X}^{\ast}\mathcal{L}_{X}:X_{w}^{\theta
}\rightarrow X_{w}^{\theta}$ is 1-1 and onto, $\left(  \left(  \mathcal{L}%
_{X}^{\ast}\mathcal{L}_{X}\right)  ^{-1}\right)  ^{\ast}=\left(  \left(
\mathcal{L}_{X}^{\ast}\mathcal{L}_{X}\right)  ^{\ast}\right)  ^{-1}=\left(
\mathcal{L}_{X}^{\ast}\mathcal{L}_{X}\right)  ^{-1}$ and so, $\left(
\mathcal{L}_{X}^{-1}\right)  ^{\ast}=\left(  \left(  \mathcal{L}_{X}^{\ast
}\mathcal{L}_{X}\right)  ^{-1}\mathcal{L}_{X}^{\ast}\right)  ^{\ast
}=\mathcal{L}_{X}\left(  \mathcal{L}_{X}^{\ast}\mathcal{L}_{X}\right)  ^{-1}%
$.\medskip

\textbf{Part 2} Left-compose $\mathcal{L}_{X}^{-1}$ with the equation of part
1.\medskip

\textbf{Part 3} Left-compose $\mathcal{L}_{X}^{\ast}$ with the equation of
part 1.\medskip

\textbf{Part 4 }By definition $\left(  \mathcal{L}_{X}^{-1}\right)  ^{\ast
\ast}=\left(  \left(  \mathcal{L}_{X}^{\ast}\mathcal{L}_{X}\right)
^{-1}\mathcal{L}_{X}^{\ast}\right)  ^{\ast\ast}=\left(  \mathcal{L}_{X}^{\ast
}\mathcal{L}_{X}\right)  ^{-1}\mathcal{L}_{X}^{\ast}$.
\end{proof}

The next result gives more properties of the function $v_{x}$ studied in the
previous subsection, as well as defining and studying the object
$r_{V,x^{\prime}}=\left(  \mathcal{L}_{X}^{-1}\right)  ^{\ast}r_{x^{\prime}%
}\in V$ which is used to derive the double order convergence estimates.

\begin{theorem}
\label{Thm_rvx}Suppose $\mathcal{S}_{X}^{e}$ is the Exact smoother mapping and
$f_{d}\in X_{w}^{\theta}$ is a data function. Then:

\begin{enumerate}
\item $\mathcal{L}_{X}^{\ast}R_{V,x}=R_{x}$ and $R_{V,x}=\left(
\mathcal{L}_{X}^{-1}\right)  ^{\ast}R_{x}$.

\item $f_{d}\left(  x\right)  -\mathcal{S}_{X}^{e}f_{d}\left(  x\right)
=\rho$ $\left(  \left(  \mathcal{L}_{X}^{-1}\right)  ^{\ast}\mathcal{Q}%
f_{d},R_{V,x}\right)  _{V}=\left\langle f_{d},v_{x}\right\rangle _{w,\theta}$.

\item $R_{x^{\prime}}\left(  x\right)  -\mathcal{S}_{X}^{e}R_{x^{\prime}%
}\left(  x\right)  =\rho$ $\left(  r_{V,x^{\prime}},R_{V,x}\right)
_{V}=\mathcal{Q}_{x^{\prime}}\overline{v_{x^{\prime}}\left(  x\right)  }$,
where $r_{V,x^{\prime}}=\left(  \mathcal{L}_{X}^{-1}\right)  ^{\ast
}r_{x^{\prime}}$ and the subscript on the operator $\mathcal{Q}_{x^{\prime}}$
indicates the action variable.

\item $\rho\left\Vert r_{V,x}\right\Vert _{V}^{2}\leq r_{x}\left(  x\right)  $.
\end{enumerate}
\end{theorem}

\begin{proof}
\textbf{Part 1} From \ref{h57} and \ref{h58}, $R_{x}=\mathcal{L}_{X}^{\ast
}R_{V,x}$, $\mathfrak{R}_{V,x}=\left(  \mathcal{L}_{X}^{\ast}\mathcal{L}%
_{X}\right)  ^{-1}R_{x}$ and $R_{V,x}=\mathcal{L}_{X}\mathfrak{R}_{V,x}$. Thus
by part 1 Lemma \ref{Lem_adjoint_inverseLx}, $R_{V,x}=\mathcal{L}_{X}\left(
\mathcal{L}_{X}^{\ast}\mathcal{L}_{X}\right)  ^{-1}R_{x}=\left(
\mathcal{L}_{X}^{-1}\right)  ^{\ast}R_{x}$.\medskip

\textbf{Part 2} From \ref{h51}, $\mathcal{S}_{X}^{e}f_{d}=f_{d}-\rho\,\left(
\mathcal{L}_{X}^{\ast}\mathcal{L}_{X}\right)  ^{-1}\mathcal{Q}f_{d}$. Hence%
\begin{align*}
f_{d}\left(  x\right)  -\mathcal{S}_{X}^{e}f_{d}\left(  x\right)  =\left(
\mathcal{L}_{X}\left(  I-\mathcal{S}_{X}^{e}\right)  f_{d},R_{V,x}\right)
_{V} &  =\rho\text{ }\left(  \mathcal{L}_{X}\left(  \mathcal{L}_{X}^{\ast
}\mathcal{L}_{X}\right)  ^{-1}\mathcal{Q}f_{d},R_{V,x}\right)  _{V}\\
&  =\rho\text{ }\left(  \left(  \mathcal{L}_{X}^{-1}\right)  ^{\ast
}\mathcal{Q}f_{d},\mathcal{L}_{X}^{\ast}R_{V,x}\right)  _{w,\theta}\\
&  =\rho\left\langle f_{d},\left(  \mathcal{L}_{X}^{\ast}\mathcal{L}%
_{X}\right)  ^{-1}\mathcal{L}_{X}^{\ast}R_{V,x}\right\rangle _{w,\theta}.
\end{align*}

But the equations \ref{h58} imply $\mathfrak{R}_{V,x}=\left(  \mathcal{L}%
_{X}^{\ast}\mathcal{L}_{X}\right)  ^{-1}\mathcal{L}_{X}^{\ast}R_{V,x}$ and by
\ref{h46}, $\rho\mathfrak{R}_{V,x}=v_{x}$ so it follows that%
\[
f_{d}\left(  x\right)  -\mathcal{S}_{X}^{e}f_{d}\left(  x\right)
=\rho\left\langle f_{d},\mathfrak{R}_{V,x}\right\rangle _{w,\theta
}=\left\langle f_{d},v_{x}\right\rangle _{w,\theta}.
\]
\medskip

\textbf{Part 3} Substitute $f_{d}=R_{x^{\prime}}$ in part 2 so that
\[
R_{x^{\prime}}\left(  x\right)  -\mathcal{S}_{X}^{e}R_{x^{\prime}}\left(
x\right)  =\rho\text{ }\left(  r_{V,x^{\prime}},R_{V,x}\right)  _{V}%
=\mathcal{Q}v_{x}\left(  x^{\prime}\right)  =\mathcal{Q}_{x^{\prime}}%
\overline{v_{x^{\prime}}\left(  x\right)  }.
\]
\medskip

\textbf{Part 4} From part 3, $r_{V,x}=\left(  \mathcal{L}_{X}^{-1}\right)
^{\ast}r_{x}$ so that from part 1 Lemma \ref{Lem_adjoint_inverseLx}
\begin{align*}
\left\Vert r_{V,x}\right\Vert _{V}^{2}=\left(  r_{V,x},r_{V,x}\right)
_{V}=\left(  \left(  \mathcal{L}_{X}^{-1}\right)  ^{\ast}r_{x},\left(
\mathcal{L}_{X}^{-1}\right)  ^{\ast}r_{x}\right)  _{V} &  =\left(
\mathcal{L}_{X}\left(  \mathcal{L}_{X}^{\ast}\mathcal{L}_{X}\right)
^{-1}r_{x},\mathcal{L}_{X}\left(  \mathcal{L}_{X}^{\ast}\mathcal{L}%
_{X}\right)  ^{-1}r_{x}\right)  _{V}\\
&  =\left(  \left(  \mathcal{L}_{X}^{\ast}\mathcal{L}_{X}\right)  ^{-1}%
r_{x},r_{x}\right)  _{w,\theta}.
\end{align*}

Now let $z_{x}=\left(  \mathcal{L}_{X}^{\ast}\mathcal{L}_{X}\right)
^{-1}r_{x}$ so that Theorem \ref{Thm_Light_norm_property} and part 7 Theorem
\ref{Thm_rx(y)_properties} imply
\begin{equation}
\left\Vert r_{V,x}\right\Vert _{V}^{2}=\mathcal{Q}z_{x}\left(  x\right)
.\label{h41}%
\end{equation}

But by part 4 Theorem \ref{Thm_L_op_properties}%
\[
\mathcal{L}_{X}^{\ast}\mathcal{L}_{X}z_{x}=\rho\mathcal{Q}z_{x}+\frac{1}%
{N}\widetilde{\mathcal{E}}_{X}^{\ast}\widetilde{\mathcal{E}}_{X}z_{x}=r_{x},
\]

which is equivalent to the system%
\begin{align}
\mathcal{Q}w_{x}+\widetilde{\mathcal{E}}_{X}^{\ast}\beta_{x}  & =r_{x}%
,\label{h37}\\
\frac{1}{N\rho}\widetilde{\mathcal{E}}_{X}w_{x}  & =\beta_{x},\label{h35}%
\end{align}

where%
\begin{equation}
\rho z_{x}=w_{x}.\label{h38}%
\end{equation}

We now want to calculate $\left\Vert r_{V,x}\right\Vert _{V}$ by solving the
system \ref{h37}, \ref{h35}, \ref{h38} and then using \ref{h41}. To do this we
use the analogue of the method used to calculate $\left\Vert R_{V,x}%
\right\Vert _{V}=\mathfrak{R}_{V,x}\left(  x\right)  $ subsequent to Theorem
\ref{Thm_pt_estim_Exact_smth_Rvx}. This involved solving the system \ref{h78},
\ref{h26}, \ref{h46} for $\mathfrak{R}_{V,x}$. Using this technique yields%
\[
w_{x}\left(  y\right)  =r_{x}\left(  y\right)  -r_{X_{2},x}^{T}\left(
N\rho\text{ }\left(  I+L_{X_{2}}L_{X_{2}}^{T}\right)  +\overline
{r_{X_{2},X_{2}}}\right)  ^{-1}\overline{r_{X_{2},y}},
\]

so that%
\begin{align*}
\rho\left\Vert r_{V,x}\right\Vert _{V}^{2}  & =\rho\mathcal{Q}z_{x}\left(
x\right)  =\mathcal{Q}w_{x}\left(  x\right)  =w_{x}\left(  x\right) \\
& =r_{x}\left(  x\right)  -r_{X_{2},x}^{T}\left(  N\rho\text{ }\left(
I+L_{X_{2}}L_{X_{2}}^{T}\right)  +\overline{r_{X_{2},X_{2}}}\right)
^{-1}\overline{r_{X_{2},x}}\\
& \leq r_{x}\left(  x\right)  ,
\end{align*}

since it was shown after \ref{h49} that $\left(  N\rho\text{ }\left(
I+L_{X_{2}}L_{X_{2}}^{T}\right)  +\overline{r_{X_{2},X_{2}}}\right)  ^{-1}$ is
positive definite.
\end{proof}

We can now derive our double order convergence estimates.

\begin{theorem}
\label{Thm_double_ord_error}Suppose $\mathcal{S}_{X}^{e}$ is the Exact
smoother mapping of Definition \ref{Def_data_func_interpol_map},
$s_{e}=\mathcal{S}_{X}^{e}f_{d}$ and $R_{x}$ is the Riesz representer of the
functional $f\rightarrow f\left(  x\right)  $ defined using a minimal
unisolvent subset $A$. Set $\eta_{G}=\eta+\delta_{G}$. Then:

\begin{enumerate}
\item If $x,x^{\prime}\in\overline{\Omega}$ and $h_{X}<\min\left\{
h_{\Omega,\theta},r_{G}\right\}  $ we have the double order of convergence
estimate%
\[
\left\vert R_{x^{\prime}}\left(  x\right)  -\mathcal{S}_{X}^{e}R_{x^{\prime}%
}\left(  x\right)  \right\vert \leq\left(  1+K_{\Omega,\theta}^{\prime
}\right)  ^{2}\sqrt{c_{G}}\left(  c_{\Omega,\theta}h_{X}\right)  ^{\eta_{G}%
}\left(  \sqrt{c_{G}}\left(  c_{\Omega,\theta}h_{X}\right)  ^{\eta_{G}}%
+\sqrt{N\rho}\right)  .
\]

\item If $X^{\prime}=\left\{  x_{k}^{\prime}\right\}  _{k=1}^{N^{\prime}%
}\subset\Omega$, $h_{X}<\min\left\{  h_{\Omega,\theta},r_{G}\right\}  $ and
the data function has the form $f_{d}=\sum\limits_{k=1}^{N^{\prime}}\beta
_{k}R_{x_{k}^{\prime}}$ for some $\beta_{k}\in\mathbb{C}$ then%
\[
\left\vert f_{d}\left(  x\right)  -s_{e}\left(  x\right)  \right\vert
\leq\left(  \sum_{k=1}^{N^{\prime}}\left\vert \beta_{k}\right\vert \right)
\left(  1+K_{\Omega,\theta}^{\prime}\right)  ^{2}\sqrt{c_{G}}\left(
c_{\Omega,\theta}h_{X}\right)  ^{\eta_{G}}\left(  \sqrt{c_{G}}\left(
c_{\Omega,\theta}h_{X}\right)  ^{\eta_{G}}+\sqrt{N\rho}\right)  ,\quad
x\in\overline{\Omega}.
\]

\item If $A\subset X^{\prime}$ then the data function of part 2 lies in
$W_{G,X^{\prime}}$.
\end{enumerate}
\end{theorem}

\begin{proof}
\textbf{Part 1} If $x,x^{\prime}\in\mathbb{R}^{d}$ then by first using the
Cauchy-Schwartz inequality and then the estimates \ref{h69} and part 4 Theorem
\ref{Thm_rvx} for $\left\Vert R_{V,x}\right\Vert _{V}$ and $\left\Vert
r_{V,x^{\prime}}\right\Vert _{V}$ we obtain%
\begin{align}
\left\vert R_{x^{\prime}}\left(  x\right)  -\mathcal{S}_{X}^{e}R_{x^{\prime}%
}\left(  x\right)  \right\vert  & =\rho\left\vert \left(  r_{V,x^{\prime}%
},R_{V,x}\right)  _{V}\right\vert \leq\rho\left\Vert r_{V,x^{\prime}%
}\right\Vert _{V}\left\Vert R_{V,x}\right\Vert _{V}\nonumber\\
& \leq\rho\sqrt{\frac{r_{x^{\prime}}\left(  x^{\prime}\right)  }{\rho}}%
\sqrt{\frac{r_{x}\left(  x\right)  +N\rho\left\vert \widetilde{l}\left(
x\right)  \right\vert ^{2}}{\rho}}\nonumber\\
& \leq\sqrt{r_{x^{\prime}}\left(  x^{\prime}\right)  }\sqrt{r_{x}\left(
x\right)  +N\rho\left\vert \widetilde{l}\left(  x\right)  \right\vert ^{2}%
}\nonumber\\
& \leq\sqrt{r_{x^{\prime}}\left(  x^{\prime}\right)  }\left(  \sqrt
{r_{x}\left(  x\right)  }+\sqrt{N\rho}\left\vert \widetilde{l}\left(
x\right)  \right\vert \right) \nonumber\\
& \leq\sqrt{r_{x^{\prime}}\left(  x^{\prime}\right)  }\left(  \sqrt
{r_{x}\left(  x\right)  }+\sqrt{N\rho}\left\vert \widetilde{l}\left(
x\right)  \right\vert _{1}\right)  .\label{h53}%
\end{align}
\medskip

Regarding the last inequality, if $x,x^{\prime}\in\Omega$ and $h_{X}%
<h_{\Omega,\theta}$ then we can uniformly estimate $\sqrt{r_{x^{\prime}%
}\left(  x^{\prime}\right)  }$ and $\sqrt{r_{x}\left(  x\right)  } $using
\ref{h03} and uniformly estimate $\left\vert \widetilde{l}\left(  x^{\prime
}\right)  \right\vert _{1}$ using \ref{h09}. Hence
\begin{align*}
\left\vert R_{x^{\prime}}\left(  x\right)  -\mathcal{S}_{X}^{e}R_{x^{\prime}%
}\left(  x\right)  \right\vert  & \leq\sqrt{r_{x^{\prime}}\left(  x^{\prime
}\right)  }\left(  \sqrt{r_{x}\left(  x\right)  }+K_{\Omega,\theta}^{\prime
}\sqrt{N\rho}\right) \\
& \leq\left(  1+K_{\Omega,\theta}^{\prime}\right)  \sqrt{c_{G}}\left(
c_{\Omega,\theta}h_{X}\right)  ^{\eta_{G}}\left(  \left(  1+K_{\Omega,\theta
}^{\prime}\right)  \sqrt{c_{G}}\left(  c_{\Omega,\theta}h_{X}\right)
^{\eta_{G}}+K_{\Omega,\theta}^{\prime}\sqrt{N\rho}\right) \\
& =\left(  1+K_{\Omega,\theta}^{\prime}\right)  ^{2}\sqrt{c_{G}}\left(
c_{\Omega,\theta}h_{X}\right)  ^{\eta_{G}}\left(  \sqrt{c_{G}}\left(
c_{\Omega,\theta}h_{X}\right)  ^{\eta_{G}}+\sqrt{N\rho}\right)  .
\end{align*}
\medskip

\textbf{Part 2} The proof proceeds by applying the uniform estimate of part 1
to each term of $f_{d}=\sum\limits_{k=1}^{N^{\prime}}\beta_{k}R_{x_{k}%
^{\prime}}$. Thus%
\begin{align*}
\left\vert f_{d}\left(  x\right)  -\mathcal{S}_{X}^{e}f_{d}\left(  x\right)
\right\vert  & \leq\left\vert \sum\limits_{k=1}^{N^{\prime}}\beta_{k}\left(
R_{x_{k}^{\prime}}\left(  x\right)  -\mathcal{S}_{X}^{e}R_{x_{k}^{\prime}%
}\left(  x\right)  \right)  \right\vert \\
& \leq\left(  \sum\limits_{k=1}^{N^{\prime}}\left\vert \beta_{k}\right\vert
\right)  \max\limits_{\substack{x\in\overline{\Omega} \\1\leq k\leq N^{\prime
}}}\left\vert R_{x_{k}^{\prime}}\left(  x\right)  -\mathcal{S}_{X}^{e}%
R_{x_{k}^{\prime}}\left(  x\right)  \right\vert \\
& \leq\left(  \sum\limits_{k=1}^{N^{\prime}}\left\vert \beta_{k}\right\vert
\right)  \left(  1+K_{\Omega,\theta}^{\prime}\right)  ^{2}\sqrt{c_{G}}\left(
c_{\Omega,\theta}h_{X}\right)  ^{\eta_{G}}\left(  \sqrt{c_{G}}\left(
c_{\Omega,\theta}h_{X}\right)  ^{\eta_{G}}+\sqrt{N\rho}\right)  .
\end{align*}
\medskip

\textbf{Part 3} If $A\subset X^{\prime}$ then by part 6 Theorem
\ref{Sum_ex_Wgx_properties_2}, $R_{x_{k}^{\prime}}\in W_{G,X^{\prime}}$ for
each $k$, and so $f_{d}\in W_{G,X^{\prime}}$.
\end{proof}

\begin{remark}
\label{Rem_smooth_converg_2}Following the approach of Remark
\ref{Rem_smooth_converg} we will now show that if we assume $h_{X}$ and $N$
are related by \ref{h74} then for a sequence of independent data sets $X_{k}$
in $\mathbb{R}^{1}$ there exists a sequence $\rho_{k}$ of smoothing
coefficients such that the smoother error is of order $2\eta_{G}$ in
$h_{X_{k}}$.\smallskip

By part 1 Theorem \ref{Thm_double_ord_error}, when $h_{X_{k}}<\min\left\{
h_{\Omega,\theta},r_{G}\right\}  $%
\[
\left\vert f_{d}\left(  x\right)  -\mathcal{S}_{X_{k}}^{e}f_{d}\left(
x\right)  \right\vert \leq\left(  1+K_{\Omega,\theta}^{\prime}\right)
^{2}\sqrt{c_{G}}\left(  c_{\Omega,\theta}h_{X}\right)  ^{\eta_{G}}\left(
\sqrt{c_{G}}\left(  c_{\Omega,\theta}h_{X}\right)  ^{\eta_{G}}+\sqrt{N\rho
}\right)  ,\quad x\in\overline{\Omega}.
\]

For clarity define the constants%
\begin{equation}
A=\left(  1+K_{\Omega,\theta}^{\prime}\right)  \sqrt{c_{G}}\left(
c_{\Omega,\theta}\right)  ^{\eta_{G}},\quad B=1+K_{\Omega,\theta}^{\prime
},\label{h420}%
\end{equation}

so that%
\[
\left\vert f_{d}\left(  x\right)  -\mathcal{S}_{X_{k}}^{e}f_{d}\left(
x\right)  \right\vert \leq A^{2}\left(  h_{X_{k}}\right)  ^{2\eta_{G}}%
+AB\sqrt{\rho}\left(  h_{X_{k}}\right)  ^{\eta_{G}}\sqrt{N_{X_{k}}},\quad
x\in\overline{\Omega}.
\]

The condition \ref{h74} implies%
\begin{equation}
\left\vert f_{d}\left(  x\right)  -\mathcal{S}_{X_{k}}^{e}f_{d}\left(
x\right)  \right\vert \leq A^{2}\left(  h_{X_{k}}\right)  ^{2\eta_{G}}%
+AB\sqrt{\rho}\left(  h_{1}\right)  ^{\frac{1}{2a}}\left(  h_{X_{k}}\right)
^{\eta_{G}-\frac{1}{2a}},\quad x\in\overline{\Omega},\label{h060}%
\end{equation}

and we want to minimize the right side as a function of $h_{X_{k}}$. This is
easily done by setting the derivative to zero so that%
\begin{align*}
D_{h_{X_{k}}}\left(  A^{2}\left(  h_{X_{k}}\right)  ^{2\eta_{G}}+AB\sqrt{\rho
}\left(  h_{1}\right)  ^{\frac{1}{2a}}\left(  h_{X_{k}}\right)  ^{\eta
_{G}-\frac{1}{2a}}\right)   & =2\eta_{G}A^{2}\left(  h_{X_{k}}\right)
^{2\eta_{G}-1}+\\
& +\left(  \eta_{G}-\frac{1}{2a}\right)  AB\sqrt{\rho}\left(  h_{X_{k}%
}\right)  ^{\eta_{G}-\frac{1}{2a}-1}\\
& =0,
\end{align*}

and a unique minimum is obtained if $\eta_{G}<\frac{1}{2a}$. In this case
$h_{X_{k}}$ satisfies
\begin{equation}
2\eta_{G}A\left(  h_{X_{k}}\right)  ^{\eta_{G}+\frac{1}{2a}}=\left(  \frac
{1}{2a}-\eta_{G}\right)  B\sqrt{\rho},\label{h040}%
\end{equation}

which can be written%
\begin{equation}
\sqrt{\rho_{k}}=\frac{2\eta_{G}}{\frac{1}{2a}-\eta_{G}}\frac{A}{B}\left(
h_{X_{k}}\right)  ^{\eta_{G}+\frac{1}{2a}}.\label{h020}%
\end{equation}

By substituting for $\rho$ in the right side of \ref{h060} we find that the
corresponding minimum error value is given by%
\begin{align*}
A^{2}\left(  h_{X_{k}}\right)  ^{2\eta_{G}}+AB\sqrt{\rho}\left(  h_{1}\right)
^{\frac{1}{2a}}\left(  h_{X_{k}}\right)  ^{\eta_{G}-\frac{1}{2a}}  &
=A^{2}\left(  h_{X_{k}}\right)  ^{2\eta_{G}}+\\
& +AB\left(  \frac{2\eta_{G}}{\frac{1}{2a}-\eta_{G}}\frac{A}{B}\left(
h_{X_{k}}\right)  ^{\eta_{G}+\frac{1}{2a}}\right)  \left(  h_{X_{k}}\right)
^{\eta_{G}-\frac{1}{2a}}\\
& =A^{2}\left(  h_{X_{k}}\right)  ^{2\eta_{G}}+A^{2}\frac{2\eta_{G}}{\frac
{1}{2a}-\eta_{G}}\left(  h_{X_{k}}\right)  ^{2\eta_{G}}\\
& =\frac{\frac{1}{2a}+\eta_{G}}{\frac{1}{2a}-\eta_{G}}A^{2}\left(  h_{X_{k}%
}\right)  ^{2\eta_{G}},
\end{align*}

so that%
\[
\left\vert f_{d}\left(  x\right)  -\mathcal{S}_{X_{k}}^{e}f_{d}\left(
x\right)  \right\vert \leq\frac{\frac{1}{2a}+\eta_{G}}{\frac{1}{2a}-\eta_{G}%
}A^{2}\left(  h_{X_{k}}\right)  ^{2\eta_{G}},\quad x\in\overline{\Omega}.
\]

Thus if $\eta_{G}<\frac{1}{2a}$ and $h_{X_{k}}\rightarrow0$ and $\rho_{k}$ is
chosen to satisfy \ref{h020}, then the Exact smoother error will converge
uniformly to zero and the order of convergence will be at least $\left(
h_{X_{k}}\right)  ^{2\eta_{G}}$.

On the other hand if $\eta_{G}\geq\frac{1}{2a}$ then noting \ref{h060} choose
$\rho_{k}$ so that $\sqrt{\rho_{k}}\left(  h_{X_{k}}\right)  ^{\eta_{G}%
-\frac{1}{2a}}=\left(  h_{X_{k}}\right)  ^{2\eta_{G}}$ i.e. so that%
\[
\sqrt{\rho_{k}}=\left(  h_{X_{k}}\right)  ^{\eta_{G}+\frac{1}{2a}},
\]

and hence \ref{h060} becomes
\[
\left\vert f_{d}\left(  x\right)  -\mathcal{S}_{X_{k}}^{e}f_{d}\left(
x\right)  \right\vert \leq\left(  A^{2}+AB\left(  h_{1}\right)  ^{\frac{1}%
{2a}}\right)  \left(  h_{X_{k}}\right)  ^{2\eta_{G}},\quad x\in\overline
{\Omega},
\]

and again the Exact smoother error converges uniformly to zero and the order
of convergence is at least $\left(  h_{X_{k}}\right)  ^{2\eta_{G}}$.
\end{remark}

\subsection{Error estimates bounded in the smoothing
parameter\label{SbSect_ErrEstimBndedInSmthParm}}

As a function of the smoothing parameter $\rho$ the convergence estimates of
the previous two subsections tend to infinity as $\rho\rightarrow\infty$.
However, numerical experiments indicate that the error is a bounded function
of $\rho$ and we will now derive some convergence estimates with this
property. The first step is to derive an estimate for the Exact smoother error
which involves the Exact smoother functional $J_{e}$. First recall the
notation of Lemma \ref{Lem_Lagrange_interpol_2}.

\begin{theorem}
\label{Thm_err_Exact_smth}Let $w$ be a weight function with properties W2 and
W3 for order $\theta$ and parameter $\kappa$ and set $\eta=\max\limits_{\kappa
}\min\left\{  \theta,\frac{1}{2}\left\lfloor \min2\kappa\right\rfloor
\right\}  $. Assume $G$ is a basis function of order $\theta$ such that there
exist constants $c_{G}>0$ and $\delta_{G}\geq0$ such that \ref{h102} holds.
Set $\eta_{G}=\eta+\delta_{G}$.

Suppose $s_{e}$ is the Exact smoother of the data function $f_{d}$ on the
unisolvent independent data $X$ with $N_{X}$ points. Suppose $A\subset X$ is a
minimal $\theta$-unisolvent set and that $\left\{  l_{k}\right\}  _{k=1}^{M} $
is the corresponding unique cardinal basis for $P_{\theta-1}$. Set
$A_{x}=A\cup\left\{  x\right\}  $.\smallskip

Then if $x\in\overline{\Omega}$ and $\operatorname*{diam}A_{x}\leq r_{G}$ then%
\begin{equation}
\left\vert f_{d}\left(  x\right)  -s_{e}\left(  x\right)  \right\vert
\leq\left\vert f_{d}\right\vert _{w,\theta}\left(  1+K_{\Omega,\theta}%
^{\prime}\right)  \sqrt{c_{G}}\left(  \operatorname*{diam}A_{x}\right)
^{\eta_{G}}+K_{\Omega,\theta}^{\prime}\sqrt{N_{X}J_{e}\left[  s_{e}\right]
},\label{h19}%
\end{equation}

and further if $h_{X}\leq\left\{  h_{\Omega,\theta},r_{G}\right\}  $ then%
\begin{equation}
\left\vert f_{d}\left(  x\right)  -s_{e}\left(  x\right)  \right\vert
\leq\left\vert f_{d}\right\vert _{w,\theta}\left(  1+K_{\Omega,\theta}%
^{\prime}\right)  \sqrt{c_{G}}\left(  c_{\Omega,\theta}h_{X}\right)
^{\eta_{G}}+K_{\Omega,\theta}^{\prime}\sqrt{N_{X}J_{e}\left[  s_{e}\right]
}.\label{h31}%
\end{equation}

Here $J_{e}$ is the Exact smoother functional.
\end{theorem}

\begin{proof}
Fix $x\in\Omega$ and suppose $A=\left\{  a^{\left(  k\right)  }\right\}
_{k=1}^{M}$. Then%
\[
f_{d}\left(  x\right)  -s_{e}\left(  x\right)  =\left(  f_{d}-s_{e}%
,R_{x}\right)  _{w,\theta}=\left\langle f_{d}-s_{e},r_{x}\right\rangle
_{w,\theta}+\sum\limits_{k=1}^{M}\left(  f_{d}\left(  a^{\left(  k\right)
}\right)  -s_{e}\left(  a^{\left(  k\right)  }\right)  \right)  l_{k}\left(
x\right)  ,
\]

so that
\begin{align*}
\left\vert f_{d}\left(  x\right)  -s_{e}\left(  x\right)  \right\vert  &
\leq\left\vert \left\langle f_{d}-s_{e},r_{x}\right\rangle _{w,\theta
}\right\vert +\sum\limits_{k=1}^{M}\left\vert f_{d}\left(  a^{\left(
k\right)  }\right)  -s_{e}\left(  a^{\left(  k\right)  }\right)  \right\vert
\left\vert l_{k}\left(  x\right)  \right\vert \\
& \leq\left\vert f_{d}-s_{e}\right\vert _{w,\theta}\left\vert r_{x}\right\vert
_{w,\theta}+\left(  \max_{k}\left\vert f_{d}\left(  a^{\left(  k\right)
}\right)  -s_{e}\left(  a^{\left(  k\right)  }\right)  \right\vert \right)
\sum\limits_{k=1}^{M}\left\vert l_{k}\left(  x\right)  \right\vert \\
& =\left\vert f_{d}-s_{e}\right\vert _{w,\theta}\sqrt{r_{x}\left(  x\right)
}+\left(  \max_{k}\left\vert f_{d}\left(  a^{\left(  k\right)  }\right)
-s_{e}\left(  a^{\left(  k\right)  }\right)  \right\vert \right)
\sum\limits_{k=1}^{M}\left\vert l_{k}\left(  x\right)  \right\vert .
\end{align*}

Since $x\in\Omega$ the estimates \ref{h09} and \ref{h07} imply%
\[
\left\vert f_{d}\left(  x\right)  -s_{e}\left(  x\right)  \right\vert
\leq\left\vert f_{d}-s_{a}\right\vert _{w,\theta}\left(  1+K_{\Omega,\theta
}^{\prime}\right)  \sqrt{c_{G}}\left(  \operatorname*{diam}A_{x}\right)
^{\eta_{G}}+K_{\Omega,\theta}^{\prime}\left(  \max_{k}\left\vert f_{d}\left(
a^{\left(  k\right)  }\right)  -s_{e}\left(  a^{\left(  k\right)  }\right)
\right\vert \right)  .
\]

From \ref{h58}, $f_{d}\left(  x\right)  -s_{e}\left(  x\right)  =\left(
\mathcal{L}_{X}\left(  f_{d}-s_{e}\right)  ,R_{V,x}\right)  _{V}$ so that%
\[
\left\vert f_{d}\left(  a^{\left(  k\right)  }\right)  -s_{e}\left(
a^{\left(  k\right)  }\right)  \right\vert \leq\left\Vert \mathcal{L}%
_{X}\left(  f_{d}-s_{e}\right)  \right\Vert _{V}\left\Vert R_{V,a^{\left(
k\right)  }}\right\Vert _{V},
\]

but Remark \ref{Rem_SmoothFunc1} with $f=s_{e}$ and part 3 of Theorem
\ref{Thm_smooth_Exact} with $f=f_{d}$ yields $\left\Vert \mathcal{L}%
_{X}\left(  f_{d}-s_{e}\right)  \right\Vert _{V}\leq\sqrt{J_{e}\left[
s_{e}\right]  }$ so that%
\[
\left\vert f_{d}\left(  x\right)  -s_{e}\left(  x\right)  \right\vert
\leq\left\vert f_{d}-s_{e}\right\vert _{w,\theta}\left(  1+K_{\Omega,\theta
}^{\prime}\right)  \sqrt{c_{G}}\left(  \operatorname*{diam}A_{x}\right)
^{\eta_{G}}+K_{\Omega,\theta}^{\prime}\sqrt{J_{e}\left[  s_{e}\right]  }%
\max_{k}\left\Vert R_{V,a^{\left(  k\right)  }}\right\Vert _{V}.
\]

From \ref{h69}, $\left\Vert R_{V,a^{\left(  k\right)  }}\right\Vert _{V}%
\leq\sqrt{N_{X}}$ and by \ref{h32} with $f=f_{d}$, $\left\vert f_{d}%
-s_{e}\right\vert _{w,\theta}\leq\left\vert f_{d}\right\vert _{w,\theta}$ so
that%
\begin{equation}
\left\vert f_{d}\left(  x\right)  -s_{e}\left(  x\right)  \right\vert
\leq\left\vert f_{d}\right\vert _{w,\theta}\left(  1+K_{\Omega,\theta}%
^{\prime}\right)  \sqrt{c_{G}}\left(  \operatorname*{diam}A_{x}\right)
^{\eta_{G}}+K_{\Omega,\theta}^{\prime}\sqrt{N_{X}J_{e}\left[  s_{e}\right]
},\label{h79}%
\end{equation}

which proves \ref{h19}. Finally, since $h_{X}\leq\left\{  h_{\Omega,\theta
},r_{G}\right\}  $, Lemma \ref{Lem_Lagrange_interpol_2} implies
$\operatorname*{diam}A_{x}\leq c_{\Omega,\theta}h_{X}$ and so \ref{h79}
implies \ref{h31}.
\end{proof}

Next we estimate the term $J_{e}\left[  s_{e}\right]  $ which occurs in both
inequalities of the last theorem.

\begin{theorem}
\label{Thm_estim_Je}Let $w$ be a weight function with properties W2 and W3 for
order $\theta$ and parameter $\kappa$ and set $\eta=\max\limits_{\kappa}%
\min\left\{  \theta,\frac{1}{2}\left\lfloor 2\kappa\right\rfloor \right\}  $.
Assume $G$ is a basis function of order $\theta$ such that there exist
constants $c_{G}>0$ and $\delta_{G}\geq0$ such that \ref{h102} holds. Set
$\eta_{G}=\eta+\delta_{G}$.

Suppose $s_{e}$ is the Exact smoother of the data function $f_{d}$ on the
independent data $X$.\smallskip

Then if $\Omega$ is a bounded, open, connected subset of $\mathbb{R}^{d}$
having the cone property:
\begin{equation}
\sqrt{J_{e}\left[  s_{e}\right]  }\leq\left\vert f_{d}\right\vert _{w,\theta
}\min\left\{  \sqrt{\rho},\left(  1+K_{\Omega,\theta}^{\prime}\right)
\sqrt{c_{G}}\left(  \operatorname*{diam}\Omega\right)  ^{\eta_{G}}\right\}
,\label{h47}%
\end{equation}

where $\rho>0$ is the smoothing parameter.
\end{theorem}

\begin{proof}
Choose a unisolvent set $A\subset X$ to define the operators $\mathcal{P}$ and
$\mathcal{Q}$.

From the definition of the smoother problem $J_{e}\left[  s_{e}\right]  \leq
J_{e}\left[  f_{d}\right]  =\rho\left\vert f_{d}\right\vert _{w,\theta}^{2}.$
Also%
\begin{align*}
J_{e}\left[  s_{e}\right]  \leq J_{e}\left[  \mathcal{P}f_{d}\right]
=\rho\left\vert \mathcal{P}f_{d}\right\vert _{w,\theta}^{2}+\frac{1}{N}%
\sum_{k=1}^{N}\left\vert \mathcal{P}f_{d}\left(  x^{\left(  k\right)
}\right)  -f_{d}\left(  x^{\left(  k\right)  }\right)  \right\vert ^{2} &
=\frac{1}{N}\sum_{k=1}^{N}\left\vert \mathcal{Q}f_{d}\left(  x^{\left(
k\right)  }\right)  \right\vert ^{2}\\
&  \leq\max_{x\in\overline{\Omega}}\left\vert \mathcal{Q}f_{d}\left(
x\right)  \right\vert ^{2}.
\end{align*}

Fix $x\in\overline{\Omega}$ and using the properties of $\mathcal{Q}$ given in
Theorem \ref{Thm_Light_norm_property} and the properties of the semi-Riesz
representer $r_{x}$ given in Theorem \ref{Thm_rx(y)_properties} we have%
\[
\left\vert \mathcal{Q}f_{d}\left(  x\right)  \right\vert =\left\vert
\left\langle f_{d},r_{x}\right\rangle _{w,\theta}\right\vert \leq\left\vert
f_{d}\right\vert _{w,\theta}\left\vert r_{x}\right\vert _{w,\theta}=\left\vert
f_{d}\right\vert _{w,\theta}\sqrt{r_{x}\left(  x\right)  }.
\]

But from \ref{h102}
\[
\sqrt{r_{x}\left(  x\right)  }\leq\left(  1+\left\vert \widetilde{l}\left(
x\right)  \right\vert _{1}\right)  \sqrt{c_{G}}\left(  \operatorname*{diam}%
A_{x}\right)  ^{\eta_{G}}\leq\left(  1+K_{\Omega,\theta}^{\prime}\right)
\sqrt{c_{G}}\left(  \operatorname*{diam}\Omega\right)  ^{\eta_{G}},
\]

so that%
\[
\sqrt{J_{e}\left[  s_{e}\right]  }\leq\max_{x\in\overline{\Omega}}\left\vert
\mathcal{Q}f_{d}\left(  x\right)  \right\vert \leq\left\vert f_{d}\right\vert
_{w,\theta}\max_{x\in\overline{\Omega}}\sqrt{r_{x}\left(  x\right)  }%
\leq\left\vert f_{d}\right\vert _{w,\theta}\left(  1+K_{\Omega,\theta}%
^{\prime}\right)  \sqrt{c_{G}}\left(  \operatorname*{diam}\Omega\right)
^{\eta_{G}}.
\]

\end{proof}

We now combine the last two results to obtain our improved error estimates
which are bounded in $\rho$.

\begin{theorem}
\label{Thm_Exact_smth_err_bded_smth_parm}Let $w$ be a weight function with
properties W2 and W3 for order $\theta$ and parameter $\kappa$ and set
$\eta=\min\left\{  \theta,\frac{1}{2}\left\lfloor \min2\kappa\right\rfloor
\right\}  $. Assume $G$ is a basis function of order $\theta$ such that there
exist constants $c_{G}>0$ and $\delta_{G}\geq0$ such that \ref{h102} holds.
Set $\eta_{G}=\eta+\delta_{G}$.

Suppose $s_{e}$ is the Exact smoother of the data function $f_{d}$ on the
unisolvent independent data $X$ with $N_{X}$ points. Suppose the notation and
assumptions of Lemma \ref{Lem_Lagrange_interpol_2} hold so that the data point
density is $h_{X}=\sup\limits_{\omega\in\Omega}\operatorname*{dist}\left(
\omega;X\right)  $.\smallskip

Then there exist constants $c_{\Omega,\theta},h_{\Omega,\theta},K_{\Omega
,\theta}^{\prime}>0$ such that
\begin{align}
\left\vert f_{d}\left(  x\right)  -s_{e}\left(  x\right)  \right\vert  &
\leq\left\vert f_{d}\right\vert _{w,\theta}\left(  1+K_{\Omega,\theta}%
^{\prime}\right)  \sqrt{c_{G}}\left(  \operatorname*{diam}A_{x}\right)
^{\eta_{G}}+\nonumber\\
&  +\left\vert f_{d}\right\vert _{w,\theta}\sqrt{N_{X}}\min\left\{  \sqrt
{\rho},\left(  1+K_{\Omega,\theta}^{\prime}\right)  \sqrt{c_{G}}\left(
\operatorname*{diam}\Omega\right)  ^{\eta_{G}}\right\}  ,\label{h80}%
\end{align}

and%
\begin{align*}
\left\vert f_{d}\left(  x\right)  -s_{e}\left(  x\right)  \right\vert  &
\leq\left\vert f_{d}\right\vert _{w,\theta}\left(  1+K_{\Omega,\theta}%
^{\prime}\right)  \sqrt{c_{G}}\left(  c_{\Omega,\theta}h_{X}\right)
^{\eta_{G}}+\\
& +\left\vert f_{d}\right\vert _{w,\theta}K_{\Omega,\theta}^{\prime}%
\sqrt{N_{X}}\min\left\{  \sqrt{\rho},\left(  1+K_{\Omega,\theta}^{\prime
}\right)  \sqrt{c_{G}}\left(  \operatorname*{diam}\Omega\right)  ^{\eta_{G}%
}\right\}  ,
\end{align*}

for $x\in\Omega$.

Here the constants $c_{\Omega,\theta},h_{\Omega,\theta},K_{\Omega,\theta
}^{\prime}$ only depend on $\Omega,\theta,\kappa$ and $d$.
\end{theorem}

\begin{proof}
The inequalities of this theorem are obtained by applying the estimates for
$\sqrt{J_{e}\left[  s_{e}\right]  }$ proved in Theorem \ref{Thm_estim_Je} to
the estimates of Theorem \ref{Thm_err_Exact_smth}.
\end{proof}

\begin{remark}
\label{Rem_smooth_converg_3}For a given independent data set $X$ the bound
\ref{h80} is a bounded function of $\rho$. However, numerical experiments
indicate that there should be a bound that is also independent of $N_{X}$.
\end{remark}

\chapter{The Approximate smoother and its convergence to the Exact smoother
and it's data function\label{Ch_Approx_smth}}

\section{Introduction\label{Sect_ap_introd_Approx_smth}}

The Approximate smoother problem is derived from the Exact smoother problem by
restricting the range of the minimizing functions from $X_{w}^{\theta}$ to a
space $W_{G,X^{\prime}}$, where $X^{\prime}$ is an arbitrary set of unisolvent
points in $\mathbb{R}^{d}$.

The Approximate smoother problem is solved twice, first using Hilbert space
methods and then using matrix methods. A matrix equation is derived and the
construction and solution of this equation is shown to be scalable w.r.t. the
number of data points i.e. it depends linearly on the number of data points.

Estimates are first derived for the pointwise convergence of the Approximate
smoother to the Exact smoother and these are then added to the Exact smoother
error estimates to obtain the error of the Approximate smoother. These
estimates involve uniform pointwise convergence and derive orders of
convergence. However they are basically unsatisfactory and are very rough
estimates when compared with numerical results. I have not included the
results of any numerical experiments in this chapter.\medskip

\textbf{Section by section:}\medskip

\textbf{Section \ref{Sect_Jg} The convolution spaces }$\overset{\cdot}{J}_{G}%
$\textbf{\ and }$J_{G}$: A study of the convolution spaces $J_{G}%
=G\ast\widehat{S}_{\emptyset,\theta}\oplus P_{\theta-1}$ and $\overset{\cdot
}{J}_{G}=G\ast\widehat{S}_{\emptyset,\theta}$. The hat indicates the Fourier
transform.\smallskip

\textbf{Section \ref{Sect_Apprx_smth_problem}} \textbf{Formulation of the
Approximate smoothing problem} The Approximate smoother problem is derived
from the Exact smoother problem studied in Chapter \ref{Ch_ExactSmth} by
restricting the range of the minimizing functions from $X_{w}^{\theta}$ to a
space $W_{G,X^{\prime}}$, where $X^{\prime}$ is an arbitrary set of $\theta
$-unisolvent points in $\mathbb{R}^{d}$. The first step is to assume that
$X^{\prime}$ is a regular grid containing the Exact smoother data and to
approximate or discretize the functions in $X_{w}^{\theta}$ using this grid.
This leads to the finite dimensional subspace $W_{G,X^{\prime}}$ and the
Approximate smoother problem, namely $\min\limits_{f\in W_{G,X^{\prime}}}%
J_{e}\left[  f\right]  $, where $J_{e}$ is the Exact smoothing functional. We
then generalize this problem to an arbitrary unisolvent $X^{\prime}%
$.\smallskip

\textbf{Section \ref{Sect_Apprx_smth_Hilbert_sp}} \textbf{Studying the
Approximate smoothing problem using Hilbert space techniques} We now know the
Approximate smoother exists, is unique and is a member of $W_{G,X^{\prime}}$.
The next step is to derive a matrix equation for the coefficients of the
$X^{\prime}$-translated basis functions and the basis polynomials. This proof
is quite similar to that of Theorem \ref{Thm_smooth_matrix_soln_1} which
derives the Exact smoother matrix equation.

Hilbert space techniques are used to show that the Approximate smoothing
problem has a unique solution. We obtain a matrix equation for the Approximate
smoother.\smallskip

\textbf{Section \ref{Sect_Appr_smth_mat} Solving the Approximate smoothing
problem using matrix techniques} Matrix techniques are used to derive the
matrix equation. We write $f=\sum\limits_{i=1}^{N^{\prime}}\alpha_{i}G\left(
\cdot-x_{i}^{\prime}\right)  +\sum\limits_{j=1}^{M}\beta_{j}p_{j}\in
W_{G,X^{\prime}}$ and calculate $J\left[  f\right]  $ as a quadratic form in
terms of $\left(  \alpha^{T}\text{ }\beta^{T}\right)  $ restrained by
$P_{X^{\prime}}^{T}\alpha=0$. We next show the quadratic component is is
positive definite and then obtain the matrix equation for the Exact smoother
using Lagrange multipliers.\smallskip

\textbf{Section \ref{Sect_Appr_smth_conv} Convergence of the Approximate
smoother to the Exact smoother} For a bounded data region $\Omega$ we first
derive some uniform, pointwise convergence results which do not involve orders
of convergence e.g. in Corollary \ref{Cor_Thm_Jsd[sig(Zk)]toJsd(sig(Z))_2} it
is shown that as the points in $X^{\prime}$ get closer to those in $X$ the
Approximate smoother converges uniformly to the Exact smoother on
$\overline{\Omega}$. We next derive some orders of pointwise convergence for
the Approximate smoother to the Exact smoother and then add these to the Exact
smoother estimates of Chapter \ref{Ch_ExactSmth} to obtain the Approximate
smoother error.\smallskip

\textbf{Section} \ref{Sect_numer_implem_Approx_smth} \textbf{A numerical
implementation of the Approximate smoother} ???

\section{The convolution spaces $\protect\overset{\cdot}{J}_{G}$ and $J_{G}%
$\label{Sect_Jg}}

In this section we define the convolution spaces $J_{G}$ which will be used to
discretize the space $X_{w}^{\theta}$ in the derivation of the Approximate
smoothing problem in Section \ref{Sect_Apprx_smth_problem}. The notation
$J_{G}$ comes from Section 3 of Nira Dyn's review paper \cite{Dyn89} and
similar results are proved here.

\begin{theorem}
\label{Thm_f_G_convol_phi}Suppose the weight function $w$ has property W2, and
suppose $G$ is a basis distribution of order $\theta\geq1$ generated by $w$.
Then $\varphi\in\widehat{S}_{\emptyset,\theta}$ implies $G\ast\varphi\in
X_{w}^{\theta}$ and if $\psi\in\widehat{S}_{\emptyset,\theta}$%
\begin{equation}
\left\langle G\ast\varphi,G\ast\psi\right\rangle _{w,\theta}=\int%
\frac{\widehat{\varphi}\overline{\widehat{\psi}}}{w\left\vert \cdot\right\vert
^{2\theta}}=\left[  G,\varphi_{\_}\ast\overline{\psi}\right]  ,\label{s8}%
\end{equation}

where $\varphi_{\_}\left(  x\right)  =\varphi\left(  -x\right)  $. Also, if
$f\in X_{w}^{\theta}$ then%
\begin{equation}
\left\langle f,G\ast\varphi\right\rangle _{w,\theta}=\left[  f,\overline
{\varphi}\right]  .\label{s6}%
\end{equation}

\end{theorem}

\begin{proof}
Corollary \ref{Cor_f_G_convol_phi} implies that $\varphi\in\widehat
{S}_{\emptyset,\theta}$ implies $G\ast\varphi\in X_{w}^{\theta}$ for
$\theta\geq1$ and that \ref{s6} is true for $\theta\geq1$. It remains to prove
\ref{s8}.

From Corollary \ref{Cor_f_G_convol_phi}, $\left\vert G\ast\varphi\right\vert
_{w,\theta}^{2}=\int\frac{\left\vert \widehat{\varphi}\right\vert ^{2}%
}{w\left\vert \cdot\right\vert ^{2\theta}}$ when $\varphi\in\widehat
{S}_{\emptyset,\theta}$. Since $X_{w}^{\theta}$ is a semi-inner product vector
space with complex scalars the semi-inner product can be recovered from the
seminorm by%
\begin{align*}
\left\langle G\ast\varphi,G\ast\psi\right\rangle _{w,\theta} &  =\frac{1}%
{4}\left(  \left\vert G\ast\left(  \phi+\psi\right)  \right\vert _{w,\theta
}^{2}-\left\vert G\ast\left(  \phi-\psi\right)  \right\vert _{w,\theta}%
^{2}\right)  +\\
&  \qquad+\frac{i}{4}\left(  \left\vert G\ast\left(  \psi-i\phi\right)
\right\vert _{w,\theta}^{2}-\left\vert G\ast\left(  \psi+i\phi\right)
\right\vert _{w,\theta}^{2}\right) \\
&  =\int\frac{f\left(  \phi,\psi\right)  }{w\left\vert \cdot\right\vert
^{2\theta}},
\end{align*}

where%
\[
f\left(  \phi,\psi\right)  =\frac{1}{4}\left(  \left\vert \widehat{\varphi
}+\widehat{\psi}\right\vert ^{2}-\left\vert \widehat{\varphi}-\widehat{\psi
}\right\vert ^{2}\right)  +\frac{i}{4}\left(  \left\vert \widehat{\psi
}-i\widehat{\phi}\right\vert -\left\vert \widehat{\psi}+i\widehat{\phi
}\right\vert \right)  =\widehat{\varphi}\overline{\widehat{\psi}},
\]

and so the first equation of \ref{s8} holds.

Now if $\varphi,\psi\in\widehat{S}_{\emptyset,\theta}$ then by part 2 of
Theorem \ref{Thm_product_of_Co,k_funcs}, $\widehat{\varphi}\overline
{\widehat{\psi}}\in S_{\emptyset,2\theta}$. Thus by Definition
\ref{Def_basis_distrib} of a basis distribution of order $\theta$ we have%
\[
\left\langle G\ast\varphi,G\ast\psi\right\rangle _{w,\theta}=\int%
\frac{\widehat{\varphi}\overline{\widehat{\psi}}}{w\left\vert \cdot\right\vert
^{2\theta}}=\left[  \widehat{G},\widehat{\varphi}\overline{\widehat{\psi}%
}\right]  =\left[  G,\left(  \widehat{\varphi}\overline{\widehat{\psi}%
}\right)  ^{\wedge}\right]  =\left[  G,\varphi_{\_}\ast\overline{\psi}\right]
,
\]

where we have used the Fourier transform and convolution identities listed in
the Appendix.
\end{proof}

\begin{definition}
\label{sDef_Jg_dotJg}\textbf{The spaces }$\overset{\cdot}{J}_{G}$\textbf{and
}$J_{G}$. Suppose the weight function $w$ has property W2. Suppose $G$ is a
basis distribution of order $\theta\geq0$ generated by a weight function $w$.
Then the spaces $\overset{\cdot}{J}_{G}$and $J_{G}$ are defined by:
\[
\overset{\cdot}{J}_{G}=G\ast\widehat{S}_{\emptyset,\theta},\qquad
J_{G}=\overset{\cdot}{J}_{G}\oplus P_{\theta-1}.
\]

Here $\widehat{S}_{\emptyset,\theta}$ denotes the Fourier transform of the
functions in $S_{\emptyset,\theta}$ and $G\ast\widehat{S}_{\emptyset,\theta}$
denotes the convolution of the tempered distribution $G$ with the functions in
$\widehat{S}_{\emptyset,\theta}$.
\end{definition}

The use of the direct product $\oplus$ in the definition of $J_{G}$ must be
justified and we must show that the definition is independent of the
particular basis function chosen.

\begin{theorem}
\label{Thm_Jg_property}Suppose the weight function $w$ has property W2.
Suppose $G$ is a basis distribution of order $\theta\geq1$ generated by a
weight function $w$. Then:

\begin{enumerate}
\item $\overset{\cdot}{J}_{G}\cap P=\left\{  0\right\}  .$

\item $P_{2\theta-1}\ast\widehat{S}_{\emptyset,\theta}\subset P_{\theta-1}$.

\item $\overset{\cdot}{J}_{G_{1}}=\overset{\cdot}{J}_{G_{2}}$iff $G_{1}%
-G_{2}\in P_{\theta-1}$.

\item Set-wise, $J_{G}$ is independent of the basis function $G$ of order
$\theta$ used to define it.

\item $\overset{\cdot}{J}_{G}$ is dense in $X_{w}^{\theta}$ in the seminorm
sense and $J_{G}$ is dense in $X_{w}^{\theta}$ when endowed with the Light norm.
\end{enumerate}
\end{theorem}

\begin{proof}
\textbf{Part 1} Suppose $\overset{\cdot}{J}_{G}\cap P\neq\{0\}$. Then we can
choose $p\in P$, $p\neq0$ and $\varphi\in\widehat{S}_{\emptyset,\theta}$ such
that $p=G\ast\varphi$. But $p\in X_{w}^{\theta}$ implies $p\in P_{\theta-1}$
so from equation \ref{s8}%
\[
0=\left\vert p\right\vert _{w,\theta}^{2}=\left\vert G\ast\varphi\right\vert
_{w,\theta}^{2}=\int\frac{\left\vert \widehat{\varphi}\right\vert ^{2}%
}{w\left\vert \cdot\right\vert ^{2\theta}},
\]

which implies $\widehat{\varphi}=0$, and so $\varphi=0$ and $p=0$.\medskip

\textbf{Part 2} Suppose $q\in P_{2\theta-1}$. Then by part 1 of Definition
\ref{Def_basis_distrib} $q=G_{1}-G_{2}$ for some basis distributions $G_{1}$
and $G_{2}$. Thus by Theorem \ref{Thm_f_G_convol_phi}, $q\ast\widehat
{S}_{\emptyset,\theta}=G_{1}\ast\widehat{S}_{\emptyset,\theta}-G_{2}%
\ast\widehat{S}_{\emptyset,\theta}\subset X_{w}^{\theta}$. But $q\ast
\widehat{S}_{\emptyset,\theta}$ is a space of polynomials so $q\ast\widehat
{S}_{\emptyset,\theta}\subset P_{\theta-1}$.\medskip

\textbf{Part 3} Suppose $\overset{\cdot}{J}_{G_{1}}=\overset{\cdot}{J}_{G_{2}%
}$. Then given $\phi\in S_{\emptyset,\theta}$ there exists $\psi\in
S_{\emptyset,\theta}$ such that $G_{1}\ast\widehat{\phi}=G_{2}\ast
\widehat{\psi}$.

Thus $G_{2}\ast\left(  \widehat{\psi}-\widehat{\phi}\right)  =\left(
G_{1}-G_{2}\right)  \ast\widehat{\phi}=q\ast\widehat{\phi}\in P$, for some
$q\in P_{2\theta-1}$.

Part 1 of this theorem now implies that $\psi=\phi$ and $q\ast\widehat{\phi
}=0$. Indeed, we can conclude that $q\ast\widehat{\phi}=0$ for all $\phi\in
S_{\emptyset,\theta}$, or equivalently, $\phi\widehat{q}=0$ for all $\phi\in
S_{\emptyset,\theta}$. But part 2 of Theorem \ref{Thm_So,n_and_Pnhat} implies
that $q\in P_{\theta-1}$ and so $G_{1}-G_{2}\in P_{\theta-1}$.

Conversely, if $G_{1}-G_{2}\in P_{\theta-1}$ and $\phi\in S_{\emptyset,\theta
}$ then $\phi$ $\left(  G_{1}-G_{2}\right)  ^{\vee}=0$ for all $\phi\in
S_{\emptyset,\theta}$, or $\left(  G_{1}-G_{2}\right)  \ast\phi=0$ for all
$\phi\in S_{\emptyset,\theta}$. Hence $\overset{\cdot}{J}_{G_{1}%
}=\overset{\cdot}{J}_{G_{2}}$.\medskip

\textbf{Part 4} By Definition \ref{Def_basis_distrib}, $G_{1}-G_{2}\in
P_{2\theta-1}$. Set $q=G_{1}-G_{2}$. Now suppose $u\in J_{G_{1}}$, so that
$u=G_{1}\ast\phi+p$ for some $\phi\in\widehat{S}_{\emptyset,\theta}$ and $p\in
P_{\theta-1}$. Consequently
\[
u=G_{1}\ast\phi+p=G_{2}\ast\phi+\left(  G_{1}-G_{2}\right)  \ast\phi
+p=G_{2}\ast\phi+q\ast\phi+p,
\]

and by part 2 the last two terms are members of $J_{G_{2}}$. Repeating the
argument for $u\in J_{G_{1}}$ proves this part.\medskip

\textbf{Part 5} ?? \textbf{ADD proof of first part of this proof}.??

A standard result is that a subspace of a Hilbert space is dense iff its
orthogonal complement is $\left\{  0\right\}  $.

Suppose $\left(  f,G\ast\phi+p\right)  _{w,\theta}=0$ for all $\phi\in
\widehat{S}_{\emptyset,\theta}$ and $p\in P_{\theta-1}$. Now when
$p=-\mathcal{P}\left(  G\ast\phi\right)  $
\begin{align*}
0=\left(  f,G\ast\phi+p\right)  _{w,\theta} &  =\left\langle f,G\ast
\phi+p\right\rangle _{w,\theta}+\sum\limits_{k}f\left(  a^{\left(  k\right)
}\right)  \overline{\left(  G\ast\phi+p\right)  \left(  a^{\left(  k\right)
}\right)  }\\
&  =\left\langle f,G\ast\phi\right\rangle _{w,\theta}\\
&  =\left[  f,\overline{\phi}\right]  ,
\end{align*}

where the last step used \ref{s6}. Hence $\left[  f,\overline{\phi}\right]  =0
$ for all $\phi\in\widehat{S}_{\emptyset,\theta}$ and by Theorem
\ref{Thm_So,n_and_Pnhat} $f\in P_{\theta-1}$. Now we have $\left(  f,p\right)
_{w,\theta}=0$ for all $p\in P_{\theta-1}$, which implies $\left\Vert
f\right\Vert _{w,\theta}^{2}=0$ and so $f=0$.
\end{proof}

\section{Formulation of the Approximate smoothing
problem\label{Sect_Apprx_smth_problem}}

The Approximate smoother problem is derived from the Exact smoother problem
studied in Chapter \ref{Ch_ExactSmth} by restricting the range of the
minimizing functions from $X_{w}^{\theta}$ to a space $W_{G,X^{\prime}}$,
where $X^{\prime}$ is an arbitrary set of $\theta$-unisolvent points in
$\mathbb{R}^{d}$. The first step is to assume that $X^{\prime}$ is a regular
grid containing the Exact smoother data and to approximate or discretize the
functions in $X_{w}^{\theta}$ using this grid. This leads to the finite
dimensional subspace $W_{G,X^{\prime}}$ and the Approximate smoother problem,
namely $\min\limits_{f\in W_{G,X^{\prime}}}J_{e}\left[  f\right]  $, where
$J_{e}$ is the Exact smoothing functional. We then generalize this problem to
an arbitrary unisolvent $X^{\prime}$.

In this section we will provide some justification for approximating the
infinite dimensional Hilbert space $X_{w}^{\theta}$ by a finite dimensional
subspace $W_{G,X^{\prime}}$, where $X^{\prime}$ is a regular, rectangular grid
of points in $\mathbb{R}^{d}$. The space $W_{G,X^{\prime}}$ will be used to
define the Approximate smoothing problem. The set $X^{\prime}$ will then be
generalized to include any finite set of distinct points.

\begin{definition}
\label{Def_grid}\textbf{A regular, rectangular grid in }$\mathbb{R}^{d}$

Let the grid occupy a rectangle $R\left(  a;b\right)  $, which has left-most
point $a\in\mathbb{R}^{d}$ and right-most point $b$. Suppose the grid has
$\mathcal{N}^{\prime}=\left(  N_{1}^{\prime},N_{2}^{\prime},\ldots
,N_{d}^{\prime}\right)  $ points in each dimension and let $h\in\mathbb{R}%
^{d}$ denote the grid sizes.

Then $X^{\prime}=\left\{  x_{\alpha}^{\prime}=a+h\alpha\mid\alpha\in
\mathbb{Z}^{d}\text{ and }0\leq\alpha<\mathcal{N}^{\prime}\right\}  $ is the
set of grid points.

Let $N^{^{\prime}}$ be the number of grid points so that $N^{^{\prime}%
}=\left(  \mathcal{N}^{\prime}\right)  ^{\mathbf{1}}=\mathcal{N}_{1}^{\prime
}\mathcal{N}_{2}^{\prime}\ldots\mathcal{N}_{d}^{\prime}$, and of course we
have the constraint $\mathcal{N}^{^{\prime}}h=b-a$.
\end{definition}

The definition of the space $W_{G,X^{\prime}}$ requires that $X^{\prime}$ is
unisolvent. If $X^{\prime}$ is a regular, rectangular grid the next theorem
shows that if the grid is made finer in all dimensions the grid eventually
becomes unisolvent, no matter what order of unisolvency is required. We will
need the following lemma:

\begin{lemma}
\label{Lem_unisolv_trans_dilat}We have the following unisolvency results:

\begin{enumerate}
\item The set $\left\{  \gamma\in\mathbb{Z}^{d}:0\leq\gamma<n\right\}  $ is
unisolvent w.r.t. $P_{n}$.

\item Translations of minimal unisolvent sets are minimal unisolvent sets.

\item Dilations of minimal unisolvent sets are minimal unisolvent sets.
\end{enumerate}
\end{lemma}

\begin{proof}
\textbf{Part 1}. From the definition of unisolvency, Definition
\ref{Def_unisolv}, we must show that for each $p\in P_{n}$, $p\left(
\gamma\right)  =0$ for $0\leq\gamma<n$ implies $p=0$. The proof will be by
induction on the order of the polynomial.

Clearly the lemma is true for $n=1$ since $P_{1}$ is the constant polynomials.

Now assume that $n\geq2$ and that if $p\in P_{n}$ and $p\left(  \gamma\right)
=0$ for $0\leq\gamma<n$ then $p=0$. Set $p\left(  x\right)  =\sum
\limits_{\left\vert \beta\right\vert <n}c_{\beta}x^{\beta}$. Then if
$\gamma_{k}=0$ and $\gamma_{i}=1$ when $i\neq k$ then $0=\sum
\limits_{\left\vert \beta\right\vert <n}c_{\beta}\gamma^{\beta}$ implies
$c_{\beta}=0$ when $\beta_{k}=0$. Thus $p\left(  \gamma\right)  =0$ for
$0\leq\gamma<n$ implies $c_{\beta}=0$ when $\beta_{i}=0$ for some $i$.
Consequently, $p=0$ if $n\leq d$ else $p\left(  x\right)  =\sum
\limits_{\substack{\left\vert \beta\right\vert <n \\\beta>0}}c_{\beta}%
x^{\beta}$. If $n>d$ we can write $p\left(  x\right)  =x^{\mathbf{1}}q\left(
x\right)  $ where $q\in P_{n-1}$, so that $q\left(  \gamma\right)  =0$ for
$1\leq\gamma<n$ must imply $q=0$. Finally, if we define $r\in P_{n-1}$ by
$r\left(  x\right)  =q\left(  x+1\right)  $ then $r\left(  \gamma\right)  =0$
for $0\leq\gamma<n-1$ must imply $r=0$, and so the truth of our lemma for
$n-1$ implies the truth of the lemma for $n$, and the lemma is proved.\medskip

\textbf{Parts 2 and 3}. From the definition of a cardinal basis, Definition
\ref{Def_cardinal_basis}, there is a unique cardinal basis $\left\{
l_{i}\right\}  $ of $P_{n}$ associated with a set $A=\left\{  a^{\left(
i\right)  }\right\}  $ iff $A$ is minimally unisolvent of order $n$. Now by
definition $l_{i}\left(  a^{\left(  j\right)  }\right)  =\delta_{i,j}$. Hence,
if $\tau,\delta\in\mathbb{R}^{d}$ and $\delta$ has positive components, then
the cardinal basis associated with the translation $A+\tau$ is $\left\{
l_{i}\left(  \cdot-\tau\right)  \right\}  $ and the cardinal basis associated
with the dilation $\delta A$ is $\left\{  l_{i}\left(  \cdot/\delta\right)
\right\}  $.
\end{proof}

\begin{theorem}
\label{Lem_reg_grid_unisolv}Suppose $X^{\prime}=\left\{  x_{\alpha}^{\prime
}=a+h\alpha\mid\alpha\in\mathbb{Z}^{d}\text{ and }0\leq\alpha<\mathcal{N}%
^{\prime}\right\}  $ is the regular, rectangular grid introduced in Definition
\ref{Def_grid}. Then $X^{\prime}$ is $\theta$-unisolvent if $\mathcal{N}%
^{\prime}\geq\theta$.
\end{theorem}

\begin{proof}
Since $\left(  X^{\prime}-a\right)  /h=\left\{  \alpha\mid\alpha\in
\mathbb{Z}^{d}\text{ and }0\leq\alpha<\mathcal{N}^{\prime}\right\}  $ and
$\mathcal{N}^{\prime}\geq\theta$, part 1 of the lemma implies $\left(
X^{\prime}-a\right)  /h$ is $\theta$-unisolvent and thus from the definition
of unisolvency, Definition \ref{Def_unisolv}, $\left(  X^{\prime}-a\right)
/h$ must contain a minimal unisolvent subset. Parts 2 and 3 of the lemma imply
that $X^{\prime}$ contains a minimal unisolvent subset and so $X^{\prime}$ is unisolvent.
\end{proof}

By Theorem \ref{Thm_Jg_property} the space $J_{G}=G\ast\widehat{S}%
_{\emptyset,\theta}+P_{\theta-1}$ is dense in $X_{w}^{\theta}$ under the Light
norm sense. So we will approximate functions in $S_{\emptyset,\theta} $ using
the grid $X^{\prime}$ defined above. Our analysis will be matrix-based so
order the grid points and write $X^{^{\prime}}=\left\{  x_{n}^{\prime
}\right\}  _{n=1}^{N^{^{\prime}}}$. We will approximate integrals on the grid
region using the trapezoidal rule i.e.
\begin{equation}
\int\limits_{grid}f\left(  x\right)  dx\simeq h^{1}\sum\limits_{n=1}%
^{N^{^{\prime}}}f\left(  x_{n}^{\prime}\right)  ,\label{s80}%
\end{equation}

where $h^{\mathbf{1}}=h_{1}h_{2}\times\ldots\times h_{d}$.\medskip

\textbf{Step1} Approximation of the functions in $\overset{\cdot}{J}_{G}%
=G\ast\widehat{S}_{\emptyset,\theta}$ by functions in $\overset{\cdot
}{W}_{G,X^{\prime}}$.\medskip

We now need a basis for $P_{\theta-1}$, say $\left\{  q_{l}\right\}
_{l=1}^{M}$ where $M=\dim P_{\theta-1}$. Then $\widehat{S}_{\emptyset,\theta}$
has the following characterization which we will not prove here:%
\begin{equation}
\widehat{S}_{\emptyset,\theta}=\left\{  \phi\in S:\int q_{l}\left(  x\right)
\phi\left(  x\right)  dx=0,\text{ }for\text{ }all\text{ }l\right\}
.\label{s18}%
\end{equation}

Suppose $\phi\in\widehat{S}_{\emptyset,\theta}$. Then the trapezoidal
approximation \ref{s80} on the grid $X^{\prime}$ we have%
\begin{equation}
G\ast\phi=\int\limits_{\mathbb{R}^{d}}G\left(  x-y\right)  \phi\left(
y\right)  dy\simeq\sum_{n=1}^{N^{\prime}}G\left(  x-x_{n}^{\prime}\right)
\left(  h^{1}\phi\left(  x_{n}^{\prime}\right)  \right)  ,\label{s20}%
\end{equation}

and%
\begin{equation}
\int\limits_{\mathbb{R}^{d}}q_{l}\left(  x\right)  \phi\left(  x\right)
dx\simeq\sum_{n=1}^{N^{^{\prime}}}q_{l}\left(  x_{n}^{\prime}\right)  \left(
h^{1}\phi\left(  x_{n}^{\prime}\right)  \right)  .\label{s62}%
\end{equation}

Considering the equations and approximations \ref{s18}, \ref{s20} and
\ref{s62} we will approximate the space $G\ast\widehat{S}_{\emptyset,\theta}$
by functions of the form $\sum\limits_{n=1}^{N^{\prime}}G\left(
x-x_{n}^{\prime}\right)  \alpha_{n}$, $\alpha_{n}\in\mathbb{C}$ constrained by
$\sum\limits_{n=1}^{N^{^{\prime}}}q_{l}\left(  x_{n}^{\prime}\right)
\alpha_{n}=0$ for $l=1,\ldots,M$. In matrix terms these constraints become
$P_{X^{\prime}}^{T}\alpha=0$, where $\alpha=\left(  \alpha_{n}\right)  $ and
$P_{X^{\prime}}=\left(  q_{j}\left(  x_{i}^{\prime}\right)  \right)  $ is the
unisolvent matrix introduced in Definition \ref{Def_unisolv_matrix_Px}. But
noting Definition \ref{Def_ex_Wgx} we see that our approximating space is just
$\overset{\cdot}{W}_{G,X^{\prime}}$, provided $X^{\prime}$ is unisolvent. But
by Lemma \ref{Lem_reg_grid_unisolv} this is true if $N^{\prime}>\dim
P_{\theta-1}$.\textbf{\medskip}

\textbf{Step 2} Approximate $J_{G}=G\ast\widehat{S}_{\emptyset,\theta
}+P_{\theta-1}$ by $W_{G,X^{\prime}}=\overset{\cdot}{W}_{G,X^{\prime}%
}+P_{\theta-1}$.\medskip

\textbf{Step 3} With this motivation we could now specify a smoothing problem,
which we will call an Approximate smoothing problem, which involves minimizing
the Exact smoothing functional $J_{e}$ given by \ref{h52} over $W_{G,X^{\prime
}}$ where $X^{\prime}$ is a rectangular grid. However, since the space
$W_{G,X^{\prime}}$ is defined when $X^{\prime}$ is any $\theta$-unisolvent set
of distinct points we will define the following more general problem:

\begin{definition}
\label{Def_Approx_smth_prob}\textbf{The Approximate smoothing problem}

Minimize the Exact smoothing functional $J_{e}\left[  f\right]  $\ for $f\in
W_{G,X^{\prime}}$, where $X^{\prime}$ is a $\theta$-unisolvent set of distinct
points in $\mathbb{R}^{d}$. More concisely we can write $\min\limits_{f\in
W_{G,X^{\prime}}}J_{e}\left[  f\right]  $.

The Exact smoothing functional $J_{e}\left[  f\right]  $\ is given by
\ref{h52} and is defined using the scattered data $\left[  X,y\right]  $ where
$X=\left\{  x^{(i)}\right\}  _{i=1}^{N}$ is the $\theta$-unisolvent
independent data and $y=\left\{  y_{i}\right\}  _{i=1}^{N}$ is termed the
dependent data.
\end{definition}

\section{Solving the Approximate smoothing problem using Hilbert space
techniques\label{Sect_Apprx_smth_Hilbert_sp}}

In this section Hilbert space techniques are used to show the Approximate
smoothing problem has a unique solution in the finite dimensional space
$W_{G,X^{\prime}}$. We then prove several identities satisfied by the
smoothing function and obtain a matrix equation for the coefficients of the
basis functions of the space $W_{G,X^{\prime}}$.

\subsection{Summary of Exact smoother properties}

To start with we will require the following properties of the mapping
$\mathcal{L}_{X}$ and the Exact smoother which were proved in Section
\ref{Sect_soln_Exact_smth}:

\begin{summary}
\label{Sum_ap_property_L_Exact_smth}Suppose $\left[  X,y\right]  $ is the data
for the \textbf{Exact smoothing problem} and $X$ is a $\theta$-unisolvent set.
Assume the operators $\mathcal{P}$, $\mathcal{Q}$ and the Light norm
$\left\Vert \cdot\right\Vert _{w,\theta}$ are all constructed using the same
minimal unisolvent subset of $X$. Then:

\begin{enumerate}
\item If $\zeta=\left(  0,y\right)  $ then $\left\Vert \mathcal{L}%
_{X}f-\varsigma\right\Vert _{V}^{2}=J_{e}[f]$ for $f\in X_{w}^{\theta}$
(Remark \ref{Rem_SmoothFunc1}).\medskip

From Theorem \ref{Thm_L_op_properties}:\smallskip

\item $\mathcal{L}_{X}:X_{w}^{\theta}\rightarrow V$ is continuous and 1-1.

\item If $u=\left(  u_{1,}\widetilde{u}_{2}\right)  \in V$ then $\mathcal{L}%
_{X}^{\ast}u=\rho\mathcal{Q}u_{1}+\frac{1}{N}\widetilde{\mathcal{E}}_{X}%
^{\ast}\widetilde{u}_{2}$ w.r.t. the Light norm constructed from a minimal
unisolvent subset of $X$.

\item $\mathcal{L}_{X}^{\ast}\mathcal{L}_{X}f=\rho\mathcal{Q}f+\frac{1}%
{N}\widetilde{\mathcal{E}}_{X}^{\ast}\widetilde{\mathcal{E}}_{X}f$ when $f\in
X_{w}^{\theta}$. Also, $\mathcal{L}_{X}^{\ast}\mathcal{L}_{X}:X_{w}^{\theta
}\rightarrow X_{w}^{\theta}$ is a homeomorphism and $\mathcal{L}_{X}^{\ast
}\mathcal{L}_{X}:W_{G,X}\rightarrow W_{G,X}$ is a homeomorphism.\medskip

From Theorem \ref{Thm_smooth_Exact}:\smallskip

\item $\left\Vert \mathcal{L}_{X}s_{e}-\varsigma\right\Vert _{V}%
^{2}+\left\Vert \mathcal{L}_{X}s_{e}-\mathcal{L}_{X}f\right\Vert _{V}%
^{2}=\left\Vert \mathcal{L}_{X}f-\varsigma\right\Vert _{V}^{2}$ for all $f\in
X_{w}^{\theta}$.

\item The Exact smoother problem has a unique solution $s_{e}$ in $W_{G,X}$
and this preserves polynomials in $P_{\theta-1}$.

\item The Exact smoother of the data $\left[  X,y\right]  $ is given by
$s_{e}=\frac{1}{N}\left(  \mathcal{L}_{X}^{\ast}\mathcal{L}_{X}\right)
^{-1}\widetilde{\mathcal{E}}_{X}^{\ast}y$.
\end{enumerate}
\end{summary}

\subsection{The existence and uniqueness of the Approximate smoother}

The geometric Hilbert space technique of orthogonal projection is used to show
the Approximate smoothing problem has a unique solution in $W_{G,X^{\prime}}$.

\begin{theorem}
\label{Thm_Appr_smth_prop1}\textbf{The Approximate smoothing problem} of
Definition \ref{Def_Approx_smth_prob} has a unique solution in $W_{G,X^{\prime
}}$, say $s_{a}$, which satisfies:

\begin{enumerate}
\item $J_{e}\left[  s_{a}\right]  <J_{e}\left[  f\right]  $ for all $f\in
W_{G,X^{\prime}}$ and $f\neq s_{a}$.

\item $\left(  \mathcal{L}_{X}s_{a}-\varsigma,\mathcal{L}_{X}s_{a}%
-\mathcal{L}_{X}f\right)  _{V}=0$ for all $f\in W_{G,X^{\prime}}$, where
$\varsigma=\left(  0,y\right)  \in V$.

\item $\left\Vert \mathcal{L}_{X}s_{a}-\varsigma\right\Vert _{V}%
^{2}+\left\Vert \mathcal{L}_{X}s_{a}-\mathcal{L}_{X}f\right\Vert _{V}%
^{2}=\left\Vert \mathcal{L}_{X}f-\varsigma\right\Vert _{V}^{2}$ for all $f\in
W_{G,X^{\prime}}$. This is equivalent to part 2.

\item $\left(  \mathcal{L}_{X}^{\ast}\mathcal{L}_{X}s_{a}-\frac{1}%
{N}\widetilde{\mathcal{E}}_{X}^{\ast}y,\,f\right)  _{w,\theta}=0$ for all
$f\in W_{G,X^{\prime}}$. This is equivalent to part 2.

\item The Approximate smoother is independent of the basis function $G$ chosen
to construct $W_{G,X^{\prime}}$.

\item If $s_{e}$ is the Exact smoother of the data $y$ then $s_{a}%
-s_{e}=\left(  \mathcal{L}_{X}^{\ast}\mathcal{L}_{X}\right)  ^{-1}g$ for some
unique $g\in W_{G,X^{\prime}}^{\bot}$.

\item $\rho\left\vert s_{e}-s_{a}\right\vert _{w,\theta}^{2}+\frac{1}{N}%
\sum\limits_{k=1}^{N}\left\vert s_{e}\left(  x^{(k)}\right)  -s_{a}\left(
x^{(k)}\right)  \right\vert ^{2}=J_{e}\left[  s_{a}\right]  -J_{e}\left[
s_{e}\right]  .$
\end{enumerate}
\end{theorem}

\begin{proof}
\textbf{Part 1} From part 1 of Summary \ref{Sum_ap_property_L_Exact_smth}
$J_{e}[f]=\left\Vert \mathcal{L}_{X}f-\varsigma\right\Vert _{V}^{2}$. So now
we want to show that there is a unique function $f\in W_{G,X^{\prime}}$ which
minimizes the functional $\left\Vert \mathcal{L}_{X}f-\varsigma\right\Vert
_{V}^{2}$ over $W_{G,X^{\prime}}$. By part 3 Theorem
\ref{Sum_ex_Wgx_properties_2}, $W_{G,X^{\prime}}$ is a finite dimensional
subspace of $X_{w}^{\theta}$ so $\mathcal{L}_{X}\left(  W_{G,X^{\prime}%
}\right)  $ must be a finite dimensional subspace of $V$ and hence a closed
subspace of $V$. Consequently there exists a unique element of $\mathcal{L}%
_{X}\left(  W_{G,X^{\prime}}\right)  $, say $v$, which is the orthogonal
projection of $\zeta$ onto $\mathcal{L}_{X}\left(  W_{G,X^{\prime}}\right)  $
such that $\left\Vert v-\varsigma\right\Vert _{V}<\left\Vert \mathcal{L}%
_{X}f-\varsigma\right\Vert _{V}$ for all $f\in W_{G,X^{\prime}}$ and
$\mathcal{L}_{X}\left(  f\right)  \neq v$.

Since $\mathcal{L}_{X}$ is 1-1 on $X_{w}^{\theta}$ there exists a unique
element of $W_{G,X^{\prime}}$, call it $s_{a}$, such that $v=\mathcal{L}%
_{X}\left(  s_{a}\right)  $.

In terms of $J_{e}$ we have $J_{e}\left[  s_{a}\right]  <J_{e}\left[
f\right]  $ for all $f\in W_{G,X^{\prime}}$ and $f\neq s_{a}$.\medskip

\textbf{Parts 2} and \textbf{3} Since $v$ is the projection of $\zeta$ onto
$\mathcal{L}_{X}\left(  W_{G,X^{\prime}}\right)  $ we have the equivalent
equations of parts 2 and 3.\medskip

\textbf{Part 4} We study the equation of part 2 using the properties of the
operators $\mathcal{L}_{X}f=\left(  f,\widetilde{\mathcal{E}}_{X}f\right)  \in
V$ and $\mathcal{L}_{X}^{\ast}$ given in Theorem
\ref{Sum_ap_property_L_Exact_smth}. In fact, for all $g\in W_{G,X^{\prime}}$
\begin{align*}
0=\left(  \mathcal{L}_{X}s_{a}-\varsigma,\mathcal{L}_{X}s_{a}-\mathcal{L}%
_{X}g\right)  _{V}=\left(  \mathcal{L}_{X}s_{a}-\varsigma,\mathcal{L}%
_{X}\left(  s_{a}-g\right)  \right)  _{V} &  =\left(  \mathcal{L}_{X}^{\ast
}\mathcal{L}_{X}s_{a}-\mathcal{L}_{X}^{\ast}\varsigma,s_{a}-g\right)
_{w,\theta}\\
&  =\left(  \mathcal{L}_{X}^{\ast}\mathcal{L}_{X}s_{a}-\frac{1}{N}%
\widetilde{\mathcal{E}}_{X}^{\ast}y,s_{a}-g\right)  _{w,\theta}.
\end{align*}

But $s_{a}\in W_{G,X^{\prime}}$ so%
\[
\left(  \mathcal{L}_{X}^{\ast}\mathcal{L}_{X}s_{a}-\frac{1}{N}\widetilde
{\mathcal{E}}_{X}^{\ast}y,f\right)  _{w,\theta}=0,\quad f\in W_{G,X^{\prime}},
\]

and this is clearly equivalent to part 2.\medskip

\textbf{Part 5} From part 1 of Definition \ref{Def_ex_Wgx} we know that the
set $W_{G,X^{\prime}}$ is independent of the basis function used in its
construction, and since the Exact smoother functional \ref{h52} is independent
of the basis function the result follows.\medskip

\textbf{Part 6} Part 4 implies directly that $\mathcal{L}_{X}^{\ast
}\mathcal{L}_{X}s_{a}-\frac{1}{N}\widetilde{\mathcal{E}}_{X}^{\ast}y\in
W_{G,X^{\prime}}^{\bot}$, say $\mathcal{L}_{X}^{\ast}\mathcal{L}_{X}%
s_{a}-\frac{1}{N}\widetilde{\mathcal{E}}_{X}^{\ast}y=g$, and this part follows
since by part 7 of Summary \ref{Sum_ap_property_L_Exact_smth}, $s_{e}=\frac
{1}{N}\left(  \mathcal{L}_{X}^{\ast}\mathcal{L}_{X}\right)  ^{-1}%
\widetilde{\mathcal{E}}_{X}^{\ast}y$.\medskip

\textbf{Part 7} Substitute $f=s_{a}$ in the identity \ref{h33}:%
\[
J_{e}\left[  s_{e}\right]  +\rho\left\vert s_{e}-f\right\vert _{w,\theta}%
^{2}+\frac{1}{N}\sum_{k=1}^{N}\left\vert s_{e}\left(  x^{(k)}\right)
-f\left(  x^{(k)}\right)  \right\vert ^{2}=J_{e}\left[  f\right]  .
\]

\end{proof}

\begin{remark}
\label{Rem_Thm_Appr_smth_prop1}In part 4 above we proved that $\mathcal{L}%
_{X}^{\ast}\mathcal{L}_{X}s_{a}-\frac{1}{N}\widetilde{\mathcal{E}}_{X}^{\ast
}y\in W_{G,X^{\prime}}^{\bot}$. It was necessary to assume that $X$ was
$\theta$-unisolvent and that $X^{\prime}$ was $\theta$-unisolvent. However,
the minimal unisolvent subsets of $X$ and $X^{\prime}$ used to construct
$\widetilde{\mathcal{E}}_{X}^{\ast}$, $\mathcal{P}$, $W_{G,X^{\prime}}$ etc.
are not required to have any points in common. We must be careful to avoid any
calculations which are dependent on both $X$ and $X^{\prime}$ but which may
only be valid when it is assumed that $X$ and $X^{\prime}$ share a minimal
unisolvent subset and this set is used to calculate $\mathcal{P}$,
$\mathcal{Q}$, the Light norm etc. For example, consider the equation
$\widetilde{\mathcal{E}}_{X^{\prime}}\widetilde{\mathcal{E}}_{X}^{\ast}%
\alpha=R_{X^{\prime},X}\alpha$. What minimal unisolvent subset is used to
calculate $R_{x}$? In this example a minimal unisolvent subset of $X$ is used
to construct $R_{x}$ which is then evaluated using $X^{\prime}$ and $X$.
\end{remark}

\subsection{Various identities\label{SbSect_ap_identities}}

The identities of this subsection are similar to the Exact smoother identities
of Corollary \ref{Cor_propert_Exact_smth} and Corollary
\ref{Cor_propert_Exact_smth_2}. They relate the Hilbert space properties and
the pointwise properties of the data and the Approximate smoother.

\begin{corollary}
Suppose $s_{a}$ is the Approximate smoother of the data $X=\left\{  x^{\left(
i\right)  }\right\}  $ and $y=\left\{  y_{i}\right\}  $ induced by the points
$X^{\prime}$, and $f_{d}\in X_{w}^{\theta}$ is a data function for the Exact
smoother. Then for all $f\in W_{G,X^{\prime}}$:

\begin{enumerate}
\item
\begin{align*}
\rho\left\vert s_{a}\right\vert _{w,\theta}^{2}+\frac{1}{N}\sum\limits_{i=1}%
^{N}\left\vert s_{a}\left(  x^{\left(  i\right)  }\right)  -y_{i}\right\vert
^{2}+\rho\left\vert s_{a}-f\right\vert _{w,\theta}^{2}+ &  \frac{1}{N}%
\sum\limits_{i=1}^{N}\left\vert s_{a}\left(  x^{\left(  i\right)  }\right)
-f\left(  x^{\left(  i\right)  }\right)  \right\vert ^{2}\\
&  =\rho\left\vert f\right\vert _{w,\theta}^{2}+\frac{1}{N}\sum\limits_{i=1}%
^{N}\left\vert f\left(  x^{\left(  i\right)  }\right)  -y_{i}\right\vert ^{2}.
\end{align*}

\item If $f_{d}\in W_{G,X^{\prime}}$ then%
\[
\left\vert s_{a}\right\vert _{w,\theta}^{2}+\frac{2}{N\rho}\sum\limits_{i=1}%
^{N}\left\vert s_{a}\left(  x^{\left(  i\right)  }\right)  -f_{d}\left(
x^{\left(  i\right)  }\right)  \right\vert ^{2}+\left\vert s_{a}%
-f_{d}\right\vert _{w,\theta}^{2}=\left\vert f_{d}\right\vert _{w,\theta}^{2}.
\]

\end{enumerate}
\end{corollary}

\begin{proof}
\textbf{Part 1} Expand the equation of part 3 of Theorem
\ref{Thm_Appr_smth_prop1} using the formula $\left\Vert \mathcal{L}%
_{X}u\right\Vert _{V}^{2}=\rho\left\vert u\right\vert _{w,\theta}^{2}+\frac
{1}{N}\sum\limits_{i=1}^{N}\left\vert u\left(  x^{\left(  i\right)  }\right)
\right\vert ^{2}$ implied by the definition of $\mathcal{L}_{X}$ and $V$
(Definition \ref{Def_ex_Hilbert_smoothing}).\medskip

\textbf{Part 2} Set $f=f_{d}$ in part 1.
\end{proof}

The result of part 4 of Theorem \ref{Thm_Appr_smth_prop1} is used to prove the
next corollary.

\begin{corollary}
\label{Cor_Appr_smth_prop2}Suppose $s_{a}\in W_{G,X^{\prime}}$ is the (unique)
Approximate smoother of the data $\left[  X,y\right]  $ induced by the points
$X^{\prime}$. Then:

\begin{enumerate}
\item $\left\langle s_{a},f\right\rangle _{w,\theta}=\frac{1}{N\rho}\left(
y-\left(  s_{a}\right)  _{X}\right)  ^{T}\overline{f_{X}}$,\quad$f\in
W_{G,X^{\prime}}$.

\item $P_{X}^{T}\left(  y-\left(  s_{a}\right)  _{X}\right)  =0$ where $P_{X}$
is any unisolvency matrix defined using $X$.

\item $\left\vert s\right\vert _{w,\theta}^{2}=\frac{1}{N\rho}\left(
y-\left(  s_{a}\right)  _{X}\right)  ^{T}\left(  \overline{s_{a}}\right)
_{X}^{T} $.

\item $J_{e}\left(  s_{a}\right)  =\frac{1}{N}\left(  y-\left(  s_{a}\right)
_{X}\right)  ^{T}\overline{y}$.

\item $\left\vert \left(  s_{a}\right)  _{X}\right\vert ^{2}\leq\left(
s_{a}\right)  _{X}^{T}\overline{y}=\left(  \overline{s_{a}}\right)  _{X}%
^{T}y\leq\left\vert y\right\vert ^{2}$.
\end{enumerate}
\end{corollary}

\begin{proof}
\textbf{Part 1} By part 4 Summary \ref{Sum_ap_property_L_Exact_smth},
$\mathcal{L}_{X}^{\ast}\mathcal{L}_{X}s_{a}=\rho\mathcal{Q}s_{a}+\frac{1}%
{N}\widetilde{\mathcal{E}}_{X}^{\ast}\widetilde{\mathcal{E}}_{X}s_{a}$ and by
part 1 Theorem \ref{Thm_Light_norm_property}, $\left(  \mathcal{Q}%
s_{a},f\right)  _{w,\theta}=\left\langle s_{a},f\right\rangle _{w,\theta}$ so
that if $f\in W_{G,X^{\prime}}$,%
\begin{align}
0=\left(  \mathcal{L}_{X}^{\ast}\mathcal{L}_{X}s_{a}-\frac{1}{N}%
\widetilde{\mathcal{E}}_{X}^{\ast}y,f\right)  _{w,\theta} &  =\left(
\rho\mathcal{Q}s_{a}+\frac{1}{N}\widetilde{\mathcal{E}}_{X}^{\ast}%
\widetilde{\mathcal{E}}_{X}s_{a}-\frac{1}{N}\widetilde{\mathcal{E}}_{X}^{\ast
}y,f\right)  _{w,\theta}\nonumber\\
&  =\left(  \rho\mathcal{Q}s_{a}+\frac{1}{N}\widetilde{\mathcal{E}}_{X}^{\ast
}\left(  \widetilde{\mathcal{E}}_{X}s_{a}-y\right)  ,f\right)  _{w,\theta
}\nonumber\\
&  =\rho\left\langle s_{a},f\right\rangle _{w,\theta}+\frac{1}{N}\left(
\widetilde{\mathcal{E}}_{X}s_{a}-y,\widetilde{\mathcal{E}}_{X}f\right)
_{\mathbb{C}^{N}}\label{s09}\\
&  =\rho\left\langle s_{a},f\right\rangle _{w,\theta}+\frac{1}{N}\left(
s-y,f_{X}\right)  _{\mathbb{C}^{N}}\nonumber\\
&  =\rho\left\langle s_{a},f\right\rangle _{w,\theta}+\frac{1}{N}\left(
\left(  s_{a}\right)  _{X}-y\right)  ^{T}\overline{f_{X}},\nonumber
\end{align}

and hence%
\[
\left\langle s_{a},f\right\rangle _{w,\theta}=\frac{1}{N\rho}\left(  y-\left(
s_{a}\right)  _{X}\right)  ^{T}\overline{f_{X}}.
\]
\medskip

\textbf{Part 2} Suppose $\left\{  p_{i}\right\}  $ is a basis for
$P_{\theta-1}$. Then if $f=p_{i}\in P_{\theta-1}$ the equation proved in part
1 becomes $\left(  p_{i}\right)  _{X}^{T}\left(  y-\left(  s_{a}\right)
_{X}\right)  =0$ and so $\left(  P_{X}\right)  ^{T}\left(  \left(
s_{a}\right)  _{X}-y\right)  =0$ since $P_{X}=\left(  p_{j}\left(  x^{\left(
i\right)  }\right)  \right)  $.\medskip

\textbf{Part 3} Let $f=s_{a}$ in the equation proved in part 1.\medskip

\textbf{Part 4} By part 3, $\left\vert s_{a}\right\vert _{w,\theta}^{2}%
=\frac{1}{N\rho}\left(  y-\left(  s_{a}\right)  _{X}\right)  ^{T}\left(
\overline{s_{a}}\right)  _{X}^{T}$ so
\begin{align*}
J_{e}\left(  s_{a}\right)   & =\rho\left\vert s_{a}\right\vert _{w,\theta}%
^{2}+\frac{1}{N}\left\vert \left(  s_{a}\right)  _{X}-y\right\vert ^{2}\\
& =\rho\left\vert s_{a}\right\vert _{w,\theta}^{2}+\frac{1}{N}\left(  \left(
s_{a}\right)  _{X}-y^{T}\right)  \left(  \left(  \overline{s_{a}}\right)
_{X}^{T}-\overline{y}\right) \\
& =\rho\left\vert s_{a}\right\vert _{w,\theta}^{2}+\frac{1}{N}\left(  \left(
s_{a}\right)  _{X}^{T}\left(  \overline{s_{a}}\right)  _{X}^{T}-\left(
s_{a}\right)  _{X}^{T}\overline{y}-y^{T}\left(  \overline{s_{a}}\right)
_{X}^{T}+y^{T}\overline{y}\right) \\
& =\frac{1}{N}\left(  y-\left(  s_{a}\right)  _{X}\right)  ^{T}\left(
\overline{s_{a}}\right)  _{X}^{T}+\frac{1}{N}\left(  \left(  s_{a}\right)
_{X}^{T}\left(  \overline{s_{a}}\right)  _{X}^{T}-\left(  s_{a}\right)
_{X}^{T}\overline{y}-y^{T}\left(  \overline{s_{a}}\right)  _{X}^{T}%
+y^{T}\overline{y}\right) \\
& =\frac{1}{N}\left(  y^{T}\left(  \overline{s_{a}}\right)  _{X}^{T}-\left(
s_{a}\right)  _{X}^{T}\left(  \overline{s_{a}}\right)  _{X}^{T}+\left(
s_{a}\right)  _{X}^{T}\left(  \overline{s_{a}}\right)  _{X}^{T}-\left(
s_{a}\right)  _{X}^{T}\overline{y}-y^{T}\left(  \overline{s_{a}}\right)
_{X}^{T}+y^{T}\overline{y}\right) \\
& =\frac{1}{N}\left(  -\left(  s_{a}\right)  _{X}^{T}\overline{y}%
+y^{T}\overline{y}\right) \\
& =\frac{1}{N}\left(  y-\left(  s_{a}\right)  _{X}\right)  ^{T}\overline{y}.
\end{align*}
\medskip

\textbf{Part 5} Part 3 implies $\left\vert \left(  s_{a}\right)
_{X}\right\vert ^{2}\leq\left(  s_{a}\right)  _{X}^{T}\overline{y}=\left(
\overline{s_{a}}\right)  _{X}^{T}y$ and part 4 implies $\left(  s_{a}\right)
_{X}^{T}\overline{y}\leq\left\vert y\right\vert ^{2}$.
\end{proof}

\subsection{Matrices and vectors derived from the Riesz representer and the
basis function}

The interpolation and Exact smoother problems only involve a single
independent data set $X$ and this leads to matrix equations which only use
matrices of the form $R_{X,X}=\left(  R_{x^{\left(  j\right)  }}\left(
x^{\left(  i\right)  }\right)  \right)  $ and $G_{X,.X}=\left(  G\left(
x^{\left(  i\right)  }-x^{\left(  j\right)  }\right)  \right)  $. However, the
Approximate smoother problem (Definition \ref{Def_Approx_smth_prob}) involves
two independent data sets $X$ and $X^{\prime}$ which will require the
following definitions:

\begin{definition}
\label{Def_Matrices_from_R}\textbf{Matrices and vectors derived from the Riesz
representer }$R_{x}$ \textbf{and the basis function} $G$

Suppose $Y=\left\{  y^{\left(  k\right)  }\right\}  $ and $Z=\left\{
z^{\left(  k\right)  }\right\}  $ are arbitrary sets of points in
$\mathbb{R}^{d}$ and $y,z\in\mathbb{R}^{d}$. Then:

\begin{enumerate}
\item $R_{Y,Z}=\left(  R_{z^{\left(  j\right)  }}\left(  y^{\left(  i\right)
}\right)  \right)  $ and $G_{Y,Z}=\left(  G\left(  y^{\left(  i\right)
}-z^{\left(  j\right)  }\right)  \right)  $.

\item $R_{Y,z}=\left(  R_{z}\left(  y^{\left(  i\right)  }\right)  \right)  $
and $G_{Y,z}=\left(  G\left(  y^{\left(  i\right)  }-z\right)  \right)  $.

\item $R_{y,Z}=\left(  R_{z^{\left(  j\right)  }}\left(  y\right)  \right)  $
and $G_{y,Z}=\left(  G\left(  y-z^{\left(  j\right)  }\right)  \right)  $.

$R_{Y,Z}$ is called a (asymmetric) \textbf{reproducing kernel matrix} and
$G_{Y,Z}$ is called a (asymmetric) \textbf{basis function matrix}.
\end{enumerate}
\end{definition}

The next result derives an important relationship between the reproducing
kernel matrix and the basis function matrix.

\begin{theorem}
\label{Thm_ap_Ryz_property}Suppose now that $A$ is a minimal unisolvent set
with cardinal basis by $\left\{  l_{i}\right\}  _{i=1}^{M}$. Define the Riesz
representer $R_{x}$ using $A$ and $\left\{  l_{i}\right\}  _{i=1}^{M}$. Then
if $Y=\left\{  y^{\left(  k\right)  }\right\}  $ and $Z=\left\{  z^{\left(
k\right)  }\right\}  $ are arbitrary sets of points in $\mathbb{R}^{d}$:

\begin{enumerate}
\item
\begin{equation}
R_{Y,Z}=\left(  2\pi\right)  ^{-\frac{d}{2}}\left(  G_{Y,Z}-L_{Y}%
G_{A,Z}-G_{Y,A}L_{Z}^{T}+L_{Y}G_{A,A}L_{Z}^{T}\right)  +L_{Y}L_{Z}%
^{T}.\label{s33}%
\end{equation}

\item If $Y$ and $Z$ are sets of distinct points in $\mathbb{R}^{d}$ then
$\widetilde{\mathcal{E}}_{Y}\widetilde{\mathcal{E}}_{Z}^{\ast}=R_{Y,Z}$.
\end{enumerate}
\end{theorem}

\begin{proof}
\textbf{Part 1} From \ref{p939}
\begin{align*}
\left(  2\pi\right)  ^{\frac{d}{2}}R_{z}\left(  y\right)  =G\left(
y-z\right)   &  -\sum_{j=1}^{M}l_{j}\left(  y\right)  G\left(  a^{\left(
j\right)  }-z\right)  -\sum_{i=1}^{M}G\left(  y-a^{\left(  i\right)  }\right)
l_{i}\left(  z\right)  +\\
&  +\sum_{i,j=1}^{M}l_{j}\left(  y\right)  G\left(  a^{\left(  j\right)
}-a^{\left(  i\right)  }\right)  l_{i}\left(  z\right)  +\left(  2\pi\right)
^{\frac{d}{2}}\sum\limits_{j=1}^{M}l_{j}(z)l_{j}\left(  y\right)  ,
\end{align*}

or in the notation introduced in Definition \ref{Def_Matrices_from_R}%
\[
R_{z}\left(  y\right)  =\left(  2\pi\right)  ^{-\frac{d}{2}}\left(  G\left(
y-z\right)  -\widetilde{l}\left(  y\right)  ^{T}G_{A,z}-G_{y,A}\widetilde
{l}\left(  z\right)  +\widetilde{l}\left(  y\right)  ^{T}G_{A,A}\widetilde
{l}\left(  z\right)  \right)  +\widetilde{l}\left(  y\right)  ^{T}%
\widetilde{l}\left(  z\right)  ,
\]

Now $R_{y,X}$ is the row vector $\left(  R_{x^{\left(  j\right)  }}\left(
y\right)  \right)  $ and $L_{X}=\left(  l_{j}\left(  x^{(i)}\right)  \right)
$ so%
\[
R_{y,Z}=\left(  2\pi\right)  ^{-\frac{d}{2}}\left(  G_{y,Z}-\widetilde
{l}\left(  y\right)  ^{T}G_{A,Z}-G_{y,A}L_{Z}^{T}+\widetilde{l}\left(
y\right)  ^{T}G_{A,A}L_{Z}^{T}\right)  +\widetilde{l}\left(  y\right)
^{T}L_{Z}^{T},
\]

and hence%
\[
R_{Y,Z}=\left(  2\pi\right)  ^{-\frac{d}{2}}\left(  G_{Y,Z}-L_{Y}%
G_{A,X}-G_{Y,A}L_{Z}^{T}+L_{Y}G_{X,X}L_{Z}^{T}\right)  +L_{Y}L_{Z}^{T}%
\]
\medskip

\textbf{Part 2} Suppose $Y=\left\{  y^{\left(  i\right)  }\right\}  $,
$Z=\left\{  z^{\left(  j\right)  }\right\}  $ and $\beta=\left(  \beta
_{j}\right)  $. Then from Definition \ref{Def_eval_operators} and \ref{q66},%
\[
\widetilde{\mathcal{E}}_{Y}\widetilde{\mathcal{E}}_{Z}^{\ast}\beta
=\sum\limits_{j=1}^{N}\beta_{j}\widetilde{\mathcal{E}}_{Y}R_{z^{(j)}}%
=\sum\limits_{j=1}^{N}\beta_{j}\left(  R_{z^{(j)}}\left(  y^{\left(  i\right)
}\right)  \right)  =R_{Y,Z}\beta,
\]
where the last step used the definition of $R_{Y,Z}$ from Definition
\ref{Def_Matrices_from_R}.
\end{proof}

\subsection{The Approximate smoother matrix
equation\label{SbSect_ap_Hilb_Appr_mat_eqn}}

We now know the Approximate smoother exists, is unique and is a member of
$W_{G,X^{\prime}}$. The next step is to derive a matrix equation for the
coefficients of the data-translated basis functions and the basis polynomials.
This proof makes good use of the identities of Subsection
\ref{SbSect_ap_identities} and the properties of unisolvency matrices, and is
quite similar to the proof of Theorem \ref{Thm_smooth_matrix_soln_1} which
derives the Exact smoother matrix equation. This proof will require notation
to deal with inner products, Riesz representers etc. which are generated by
two minimal unisolvent sets $A$ and $A^{\prime}$:

\begin{notation}
\label{Not_A'}Suppose $A$ and $A^{\prime}$ are two minimal unisolvent sets and
$A$ is the "default" notation. Let $\left(  \cdot,\cdot\right)  _{w,\theta
}^{\prime}$ be the Light inner product generated by $A^{\prime}$,
$R_{x}^{\prime}$ be the Riesz representer generated by $A^{\prime}$ and
$\widetilde{\mathcal{E}}_{A^{\prime};X}^{\ast}$ be the adjoint of
$\widetilde{\mathcal{E}}_{X}$ w.r.t. $\left(  \cdot,\cdot\right)  _{w,\theta
}^{\prime}$.

Regarding unisolvency matrix notation, if the $P_{\theta-1}$ basis used for
$P_{X^{\prime}}$ is different to that used for $P_{X}$ then use the notation
$P_{X^{\prime}}^{\prime}$ instead of $P_{X^{\prime}}$. In contrast, no `prime'
notation is required for the cardinal unisolvency matrix $L_{X^{\prime}}$
because the (cardinal) basis used to calculate $L_{X^{\prime}}$ is uniquely
determined by $A^{\prime}$.
\end{notation}

\begin{theorem}
\label{Thm_Appr_smth_mat_eqn_2}\textbf{The Approximate smoother matrix
equation}

Fix a basis function $G$ of order $\theta$. Choose a minimal unisolvent set
$A^{\prime}\subset X^{\prime}$ and a $P_{\theta-1}$ basis $\left\{
p_{i}\right\}  $ to define the unisolvency matrix $P_{X^{\prime}}$, the Riesz
representer $R_{x}^{\prime}$ and the space $W_{G,X^{\prime}}$. Choose a
minimal unisolvent set $A\subset X$ and the $P_{\theta-1}$ basis $\left\{
p_{i}\right\}  $ to define the unisolvency matrix $P_{X}$, the Riesz
representer $R_{x}$ and $W_{G,X}$. Now suppose $s\in W_{G,X^{\prime}}$ is the
(unique) Approximate smoother of the data $\left[  X=\left\{  x^{\left(
i\right)  }\right\}  _{i=1}^{N},y=\left(  y_{i}\right)  \right]  $ induced by
the set of points $X^{\prime}=\left\{  x_{i}^{\prime}\right\}  _{i=1}%
^{N^{\prime}}$. Since $s\in W_{G,X^{\prime}}$%
\begin{equation}
s\left(  x\right)  =\sum\limits_{i=1}^{N^{\prime}}\alpha_{i}^{\prime}G\left(
x-x_{i}^{\prime}\right)  +\sum\limits_{i=1}^{M}\beta_{i}^{\prime}p_{i}\left(
x\right)  ,\label{s29}%
\end{equation}

for some $\alpha=\left(  \alpha_{i}^{\prime}\right)  \in\mathbb{C}^{N}$ and
$\beta^{\prime}=\left(  \beta_{i}^{\prime}\right)  \in\mathbb{C}^{M}$. In
fact, $\alpha^{\prime}$ and $\beta^{\prime}$ satisfy the \textbf{Approximate
smoother matrix equation}
\begin{equation}
\left(
\begin{array}
[c]{lll}%
\left(  2\pi\right)  ^{\frac{d}{2}}N\rho G_{X^{\prime},X^{\prime}%
}+G_{X^{\prime},X}G_{X,X^{\prime}} & G_{X^{\prime},X}P_{X} & P_{X^{\prime}}\\
P_{X}^{T}G_{X,X^{\prime}} & P_{X}^{T}P_{X} & O_{M}\\
P_{X^{\prime}}^{T} & O_{M} & O_{M}%
\end{array}
\right)  \left(
\begin{array}
[c]{c}%
\alpha^{\prime}\\
\beta^{\prime}\\
\gamma^{\prime}%
\end{array}
\right)  =\left(
\begin{array}
[c]{l}%
G_{X^{\prime},X}\\
P_{X}^{T}\\
O_{M,N}%
\end{array}
\right)  y,\label{s16}%
\end{equation}

where $G_{X^{\prime},X^{\prime}}=G\left(  x_{i}^{\prime}-x_{j}^{\prime
}\right)  $, $G_{X^{\prime},X}=G\left(  x_{i}^{\prime}-x^{\left(  j\right)
}\right)  $, $G_{X^{\prime},X}=G\left(  x^{\left(  i\right)  }-x_{j}^{\prime
}\right)  $ are basis function matrices and $P_{X}=\left(  p_{j}\left(
x^{\left(  i\right)  }\right)  \right)  $ and $P_{X^{\prime}}=\left(
p_{j}\left(  x_{i}^{\prime}\right)  \right)  $ are unisolvency matrices
(Definition \ref{Def_unisolv_matrix_Px}).

The matrix on the left of this equation will be called the \textbf{Approximate
smoother matrix }and will usually be denoted by the symbol $\Psi$.
\end{theorem}

\begin{proof}
From part 3 of Definition \ref{Def_ex_Wgx} the space $W_{G,X^{\prime}}$ is
independent of the ordering of the points in $X^{\prime}$. Thus the
Approximate smoothing problem is also independent of the ordering of the
points in $X^{\prime}$ and we now take advantage of this to order $X^{\prime}$
so that the first $M$ points of $X^{\prime}$ lie in $A^{\prime}$. Equation
\ref{s09} of the proof of Corollary \ref{Cor_Appr_smth_prop2} is%
\[
\left\langle s,f\right\rangle _{w,\theta}=\frac{1}{N\rho}\left(
y-\widetilde{\mathcal{E}}_{X}s,\widetilde{\mathcal{E}}_{X}f\right)  ,\quad
f\in W_{G,X^{\prime}},
\]

so the Light norm \ref{p917} corresponding to $A^{\prime}$ is%
\begin{align*}
\left(  s,f\right)  _{w,\theta}^{\prime}  & =\frac{1}{N\rho}\left(
y-\widetilde{\mathcal{E}}_{X}s,\widetilde{\mathcal{E}}_{X}f\right)  +\left(
\widetilde{\mathcal{E}}_{A^{\prime}}s,\widetilde{\mathcal{E}}_{A^{\prime}%
}f\right) \\
& =\frac{1}{N\rho}\left(  \widetilde{\mathcal{E}}_{A^{\prime};X}^{\ast}\left(
y-\widetilde{\mathcal{E}}_{X}s\right)  ,f\right)  _{w,\theta}^{\prime}+\left(
\widetilde{\mathcal{E}}_{A^{\prime};A^{\prime}}^{\ast}\widetilde{\mathcal{E}%
}_{A^{\prime}}s,f\right)  _{w,\theta}^{\prime},
\end{align*}

where we have used Notation \ref{Not_A'}. By part 6 Theorem
\ref{Sum_ex_Wgx_properties_2} the Riesz representer $R_{x}^{\prime}\in
W_{G,X^{\prime}}$ when $x\in X^{\prime}$ so the last equation implies%
\[
\widetilde{\mathcal{E}}_{X^{\prime}}\left(  s-\widetilde{\mathcal{E}%
}_{A^{\prime};A^{\prime}}^{\ast}\widetilde{\mathcal{E}}_{A^{\prime}}s+\frac
{1}{N\rho}\widetilde{\mathcal{E}}_{A^{\prime};X}^{\ast}\left(  s-\widetilde
{\mathcal{E}}_{X}y\right)  \right)  =0,
\]

or%
\[
N\rho\widetilde{\mathcal{E}}_{X^{\prime}}s-N\rho\widetilde{\mathcal{E}%
}_{X^{\prime}}\widetilde{\mathcal{E}}_{A^{\prime};A^{\prime}}^{\ast}%
\widetilde{\mathcal{E}}_{A^{\prime}}s+\widetilde{\mathcal{E}}_{X^{\prime}%
}\widetilde{\mathcal{E}}_{A^{\prime};X}^{\ast}\widetilde{\mathcal{E}}%
_{X}s-\widetilde{\mathcal{E}}_{X^{\prime}}\widetilde{\mathcal{E}}_{A^{\prime
};X}^{\ast}y=0.
\]

By part 5 Theorem \ref{Thm_eval_op_properties} and Definition
\ref{Def_Matrices_from_R}, $\widetilde{\mathcal{E}}_{X^{\prime}}%
\widetilde{\mathcal{E}}_{A^{\prime};A^{\prime}}^{\ast}=L_{X^{\prime}}$ and
$\widetilde{\mathcal{E}}_{X^{\prime}}\widetilde{\mathcal{E}}_{A^{\prime}%
;X}^{\ast}=R_{X^{\prime};X}^{\prime}$ so%
\[
N\rho\widetilde{\mathcal{E}}_{X^{\prime}}s-N\rho L_{X^{\prime}}\mathcal{E}%
_{A^{\prime}}s+R_{X^{\prime};X}^{\prime}\widetilde{\mathcal{E}}_{X}%
s-R_{X^{\prime};X}^{\prime}y=0,
\]

and for clarity we revert to the subscript notation for evaluations:%
\[
N\rho s_{X^{\prime}}-N\rho L_{X^{\prime}}s_{A^{\prime}}+R_{X^{\prime}%
;X}^{\prime}s_{X}-R_{X^{\prime};X}^{\prime}y=0,
\]

or%
\begin{equation}
N\rho\left(  I_{N^{\prime}}-L_{X^{\prime}:0}\right)  s_{X^{\prime}%
}+R_{X^{\prime},X}^{\prime}\left(  s_{X}-y\right)  =0\label{s35}%
\end{equation}

where the augmented square matrix $L_{X^{\prime}:0}=%
\begin{pmatrix}
L_{X^{\prime}} & O_{N,N-M}%
\end{pmatrix}
$ was introduced in part 3 Theorem \ref{Thm_Px_properties_2}. Setting
$Y=X^{\prime}$ and $Z=X$ in formula \ref{s33} and using the special ordering
of $X^{\prime}$ mentioned at the start of the proof we obtain%
\begin{equation}
\left(  2\pi\right)  ^{\frac{d}{2}}R_{X^{\prime},X}^{\prime}=G_{X^{\prime}%
,X}-L_{X^{\prime}}G_{A^{\prime},X}-G_{X^{\prime},A^{\prime}}L_{X}%
^{T}+L_{X^{\prime}}G_{A^{\prime},A^{\prime}}L_{X}^{T}+\left(  2\pi\right)
^{\frac{d}{2}}L_{X^{\prime}}L_{X}^{T}.\label{s40}%
\end{equation}

Part 2 of Corollary \ref{Cor_Appr_smth_prop2} implies $L_{X}^{T}\left(
s_{X}-y\right)  =0$. Using this equation and \ref{s40}, \ref{s35} becomes%
\begin{align}
0  & =\left(  2\pi\right)  ^{\frac{d}{2}}N\rho\left(  I_{N^{\prime}%
}-L_{X^{\prime}:0}\right)  s_{X^{\prime}}+\left(  2\pi\right)  ^{\frac{d}{2}%
}R_{X^{\prime},X}^{\prime}\left(  s_{X}-y\right) \nonumber\\
& =\left(  2\pi\right)  ^{\frac{d}{2}}N\rho\left(  I_{N^{\prime}}%
-L_{X^{\prime}:0}\right)  s_{X^{\prime}}+\left(  G_{X^{\prime},X}%
-L_{X^{\prime}}G_{A^{\prime},X}-G_{X^{\prime},A^{\prime}}L_{X}^{T}%
+L_{X^{\prime}}G_{A^{\prime},A^{\prime}}L_{X}^{T}\right)  \left(
s_{X}-y\right)  +\nonumber\\
& \qquad\qquad\qquad+\left(  2\pi\right)  ^{\frac{d}{2}}L_{X^{\prime}}%
L_{X}^{T}\left(  s_{X}-y\right) \nonumber\\
& =\left(  2\pi\right)  ^{\frac{d}{2}}N\rho\left(  I_{N^{\prime}}%
-L_{X^{\prime}:0}\right)  s_{X^{\prime}}+\left(  G_{X^{\prime},X}%
-L_{X^{\prime}}G_{A^{\prime},X}\right)  \left(  s_{X}-y\right) \nonumber\\
& =\left(  2\pi\right)  ^{\frac{d}{2}}N\rho\left(  I_{N^{\prime}}%
-L_{X^{\prime}:0}\right)  s_{X^{\prime}}+\left(  G_{X^{\prime},X}%
-L_{X^{\prime}:0}G_{X^{\prime},X}\right)  s_{X}-\left(  G_{X^{\prime}%
,X}-L_{X^{\prime}:0}G_{X^{\prime},X}\right)  y\nonumber\\
& =\left(  2\pi\right)  ^{\frac{d}{2}}N\rho\left(  I_{N^{\prime}}%
-L_{X^{\prime}:0}\right)  s_{X^{\prime}}+\left(  I_{N^{\prime}}-L_{X^{\prime
}:0}\right)  G_{X^{\prime},X}s_{X}-\left(  I_{N^{\prime}}-L_{X^{\prime}%
:0}\right)  G_{X^{\prime},X}y.\label{s38}%
\end{align}

The next step is to express the factors $s_{X^{\prime}}$ and $s_{X}$ in basis
function terms. But from the statement of this theorem, $s\left(  x\right)
=\sum\limits_{i=1}^{N^{\prime}}\alpha_{i}^{\prime}G\left(  x-x_{i}^{\prime
}\right)  +\sum\limits_{j=1}^{M}\beta_{j}^{\prime}p_{j}\left(  x\right)  $
constrained by
\begin{equation}
P_{X^{\prime}}^{T}\alpha^{\prime}=0,\label{s37}%
\end{equation}

so that%
\begin{equation}
s_{X}=G_{X,X^{\prime}}\alpha^{\prime}+P_{X}\beta^{\prime},\quad s_{X^{\prime}%
}=G_{X^{\prime},X^{\prime}}\alpha^{\prime}+P_{X^{\prime}}\beta^{\prime
}.\label{s10}%
\end{equation}

From part 3 of Theorem \ref{Thm_Px_properties} we know that $P_{X^{\prime}%
}=L_{X^{\prime}}P_{A^{\prime}}$, where $P_{A^{\prime}}$ is regular. Also from
part 3 Theorem \ref{Thm_Px_properties_2} $L_{X^{\prime}:0}L_{X^{\prime}%
}=L_{X^{\prime}}$. Hence%
\begin{align*}
\left(  2\pi\right)  ^{\frac{d}{2}} &  N\rho\left(  I_{N^{\prime}%
}-L_{X^{\prime}:0}\right)  s_{X^{\prime}}+\left(  I_{N^{\prime}}-L_{X^{\prime
}:0}\right)  G_{X^{\prime},X}s_{X}\\
&  =\left(  2\pi\right)  ^{\frac{d}{2}}N\rho\left(  I_{N^{\prime}%
}-L_{X^{\prime}:0}\right)  \left(  G_{X^{\prime},X^{\prime}}\alpha^{\prime
}+P_{X^{\prime}}\beta^{\prime}\right)  +\left(  I_{N^{\prime}}-L_{X^{\prime
}:0}\right)  G_{X^{\prime},X}s_{X}\\
&  =\left(  2\pi\right)  ^{\frac{d}{2}}N\rho\left(  I_{N^{\prime}%
}-L_{X^{\prime}:0}\right)  \left(  G_{X^{\prime},X^{\prime}}\alpha^{\prime
}+L_{X^{\prime}}P_{A^{\prime}}\beta^{\prime}\right)  +\left(  I_{N^{\prime}%
}-L_{X^{\prime}:0}\right)  G_{X^{\prime},X}s_{X}\\
&  =\left(  2\pi\right)  ^{\frac{d}{2}}N\rho\left(  I_{N^{\prime}%
}-L_{X^{\prime}:0}\right)  G_{X^{\prime},X^{\prime}}\alpha^{\prime}+\left(
I_{N^{\prime}}-L_{X^{\prime}:0}\right)  G_{X^{\prime},X}s_{X}\\
&  =\left(  2\pi\right)  ^{\frac{d}{2}}N\rho\left(  I_{N^{\prime}%
}-L_{X^{\prime}:0}\right)  G_{X^{\prime},X^{\prime}}\alpha^{\prime}+\left(
I_{N^{\prime}}-L_{X^{\prime}:0}\right)  G_{X^{\prime},X}\left(  G_{X,X^{\prime
}}\alpha^{\prime}+P_{X}\beta^{\prime}\right) \\
&  =\left(  I_{N^{\prime}}-L_{X^{\prime}:0}\right)  \left(  \left(  \left(
2\pi\right)  ^{\frac{d}{2}}N\rho G_{X^{\prime},X^{\prime}}+G_{X^{\prime}%
,X}G_{X,X^{\prime}}\right)  \alpha^{\prime}+G_{X^{\prime},X}P_{X}\beta
^{\prime}\right) \\
&  =\left(  \left(  2\pi\right)  ^{\frac{d}{2}}N\rho G_{X^{\prime},X^{\prime}%
}+G_{X^{\prime},X}G_{X,X^{\prime}}\right)  \alpha^{\prime}+G_{X^{\prime}%
,X}P_{X}\beta^{\prime}-\\
&  -L_{X^{\prime}:0}\left(  \left(  2\pi\right)  ^{\frac{d}{2}}N\rho
G_{X^{\prime},X^{\prime}}\alpha^{\prime}+G_{X^{\prime},X}G_{X,X^{\prime}%
}\alpha^{\prime}+G_{X^{\prime},X}P_{X}\beta^{\prime}\right)  ,
\end{align*}

and \ref{s38} becomes%
\begin{align*}
0  & =\left(  \left(  2\pi\right)  ^{\frac{d}{2}}N\rho G_{X^{\prime}%
,X^{\prime}}+G_{X^{\prime},X}G_{X,X^{\prime}}\right)  \alpha^{\prime
}+G_{X^{\prime},X}P_{X}\beta^{\prime}-\\
& \qquad-L_{X^{\prime}:0}\left(  \left(  2\pi\right)  ^{\frac{d}{2}}N\rho
G_{X^{\prime},X^{\prime}}\alpha^{\prime}+G_{X^{\prime},X}G_{X,X^{\prime}%
}\alpha^{\prime}+G_{X^{\prime},X}P_{X}\beta^{\prime}\right) \\
& \qquad-\left(  I_{N^{\prime}}-L_{X^{\prime}:0}\right)  G_{X^{\prime},X}y\\
& =\left(  \left(  2\pi\right)  ^{\frac{d}{2}}N\rho G_{X^{\prime},X^{\prime}%
}+G_{X^{\prime},X}G_{X,X^{\prime}}\right)  \alpha^{\prime}+G_{X^{\prime}%
,X}P_{X}\beta^{\prime}-G_{X^{\prime},X}y-\\
& \qquad-L_{X^{\prime}:0}\left(  \left(  2\pi\right)  ^{\frac{d}{2}}N\rho
G_{X^{\prime},X^{\prime}}\alpha^{\prime}+G_{X^{\prime},X}G_{X,X^{\prime}%
}\alpha^{\prime}+G_{X^{\prime},X}P_{X}\beta^{\prime}+G_{X^{\prime},X}y\right)
\\
& =\left(  \left(  2\pi\right)  ^{\frac{d}{2}}N\rho G_{X^{\prime},X^{\prime}%
}+G_{X^{\prime},X}G_{X,X^{\prime}}\right)  \alpha^{\prime}+G_{X^{\prime}%
,X}P_{X}\beta^{\prime}-G_{X^{\prime},X}y-\\
& \qquad-L_{X^{\prime}}\left(  \left(  2\pi\right)  ^{\frac{d}{2}}N\rho
G_{A^{\prime},X^{\prime}}\alpha^{\prime}+G_{A^{\prime},X}G_{X,X^{\prime}%
}\alpha^{\prime}+G_{A^{\prime},X}P_{X}\beta^{\prime}+G_{A^{\prime},X}y\right)
,
\end{align*}

From part 3 of Theorem \ref{Thm_Px_properties} we know that $L_{X^{\prime}%
}=P_{X^{\prime}}P_{A^{\prime}}^{-1}$, where $P_{A^{\prime}}$ is regular, so%
\begin{align*}
0= &  \left(  \left(  2\pi\right)  ^{\frac{d}{2}}N\rho G_{X^{\prime}%
,X^{\prime}}+G_{X^{\prime},X}G_{X,X^{\prime}}\right)  \alpha^{\prime
}+G_{X^{\prime},X}P_{X}\beta^{\prime}-G_{X^{\prime},X}y-\\
&  -P_{X^{\prime}}P_{A^{\prime}}^{-1}\left(  \left(  2\pi\right)  ^{\frac
{d}{2}}N\rho G_{A^{\prime},X^{\prime}}\alpha^{\prime}+G_{A^{\prime}%
,X}G_{X,X^{\prime}}\alpha^{\prime}+G_{A^{\prime},X}P_{X}\beta^{\prime
}+G_{A^{\prime},X}y\right)  ,
\end{align*}

which we write as%
\begin{equation}
\left(  \left(  2\pi\right)  ^{\frac{d}{2}}N\rho G_{X^{\prime},X^{\prime}%
}+G_{X^{\prime},X}G_{X,X^{\prime}}\right)  \alpha^{\prime}+G_{X^{\prime}%
,X}P_{X}\beta^{\prime}+P_{X^{\prime}}\gamma^{\prime}=G_{X^{\prime}%
,X}y,\label{s39}%
\end{equation}

where%
\[
\gamma^{\prime}=P_{A^{\prime}}^{-1}\left(  \left(  2\pi\right)  ^{\frac{d}{2}%
}N\rho G_{A^{\prime},X^{\prime}}\alpha^{\prime}+G_{A^{\prime},X}%
G_{X,X^{\prime}}\alpha^{\prime}+G_{A^{\prime},X}P_{X}\beta^{\prime
}+G_{A^{\prime},X}y\right)  .
\]

Finally, from part 2 of Corollary \ref{Cor_Appr_smth_prop2} it follows that
$P_{X}^{T}\left(  s_{X}-y\right)  =0$ and thus substituting for $s_{X}$ from
\ref{s10} we obtain%
\begin{equation}
P_{X}^{T}G_{X,X^{\prime}}\alpha^{\prime}+P_{X}^{T}P_{X}\beta^{\prime}%
=P_{X}^{T}y.\label{s23}%
\end{equation}

Equations \ref{s37}, \ref{s39} and \ref{s23} now combine to give the
Approximate smoother (block) matrix equation \ref{s16}.

At the start of this proof we assumed a special order for the points in the
ordered set $X^{\prime}$. To prove that the Approximate smoother is unchanged
by reordering the points $X^{\prime}$ we use the permutation matrix approach.
In fact we will show that for the permutation $\pi$ of $X^{\prime}$%
\begin{equation}
\left(
\begin{array}
[c]{lll}%
\begin{array}
[c]{l}%
\left(  2\pi\right)  ^{\frac{d}{2}}N\rho G_{\pi\left(  X^{\prime}\right)
,\pi\left(  X^{\prime}\right)  }+\\
\qquad+G_{\pi\left(  X^{\prime}\right)  ,X}G_{X,\pi\left(  X^{\prime}\right)
}%
\end{array}
& G_{\pi\left(  X^{\prime}\right)  ,X}P_{X} & P_{\pi\left(  X^{\prime}\right)
}\\
P_{X}^{T}G_{X,\pi\left(  X^{\prime}\right)  } & P_{X}^{T}P_{X} & O_{M}\\
P_{\pi\left(  X^{\prime}\right)  }^{T} & O_{M} & O_{M}%
\end{array}
\right)  \left(
\begin{array}
[c]{c}%
\pi\left(  \alpha^{\prime}\right) \\
\beta^{\prime}\\
\gamma^{\prime}%
\end{array}
\right)  =\left(
\begin{array}
[c]{l}%
G_{\pi\left(  X^{\prime}\right)  ,X}\\
P_{X}^{T}\\
O_{M,N}%
\end{array}
\right)  y,\label{s36}%
\end{equation}

and%
\begin{equation}
s\left(  x\right)  =G_{x,\pi\left(  X^{\prime}\right)  }\pi\left(
\alpha^{\prime}\right)  +\sum\limits_{i=1}^{M}\beta_{i}p_{i}\left(  x\right)
.\label{s41}%
\end{equation}

Denote the permutation matrix of $\pi$ by $\Pi$. To reorder the rows of a
matrix we left-multiply by $\Pi$ and to reorder the columns we right-multiply
by $\Pi^{T}$. Also $\Pi^{T}\Pi=\Pi\Pi^{T}=I$. Hence, $G_{\pi\left(  X^{\prime
}\right)  ,\pi\left(  X^{\prime}\right)  }=\Pi G_{X^{\prime},X^{\prime}}%
\Pi^{T}$, $G_{X,\pi\left(  X^{\prime}\right)  }=G_{X,X^{\prime}}\Pi^{T}$,
$G_{\pi\left(  X^{\prime}\right)  ,X}=\Pi G_{X^{\prime},X}$ and $P_{\pi\left(
X^{\prime}\right)  }=\Pi P_{X^{\prime}}$. The left side of the equation of the
first row of \ref{s36} now becomes%
\begin{align*}
&  \left(  \left(  2\pi\right)  ^{\frac{d}{2}}N\rho G_{\pi\left(  X^{\prime
}\right)  ,\pi\left(  X^{\prime}\right)  }+G_{\pi\left(  X^{\prime}\right)
,X}G_{X,\pi\left(  X^{\prime}\right)  }\right)  \pi\left(  \alpha^{\prime
}\right)  +G_{\pi\left(  X^{\prime}\right)  ,X}P_{X}\beta^{\prime}%
+P_{\pi\left(  X^{\prime}\right)  }\gamma^{\prime}\\
&  =\left(  \left(  2\pi\right)  ^{\frac{d}{2}}N\rho\Pi G_{X^{\prime
},X^{\prime}}\Pi^{T}+\Pi G_{X^{\prime},X}G_{X,X^{\prime}}\Pi^{T}\right)
\Pi\alpha^{\prime}+\Pi G_{X^{\prime},X}P_{X}\beta^{\prime}+\Pi P_{X^{\prime}%
}\gamma^{\prime}\\
&  =\Pi\left(  \left(  \left(  2\pi\right)  ^{\frac{d}{2}}N\rho G_{X^{\prime
},X^{\prime}}+G_{X^{\prime},X}G_{X,X^{\prime}}\right)  \Pi^{T}\Pi
\alpha^{\prime}+G_{X^{\prime},X}P_{X}\beta^{\prime}+P_{X^{\prime}}%
\gamma^{\prime}\right) \\
&  =\Pi\left(  \left(  \left(  2\pi\right)  ^{\frac{d}{2}}N\rho G_{X^{\prime
},X^{\prime}}+G_{X^{\prime},X}G_{X,X^{\prime}}\right)  \alpha^{\prime
}+G_{X^{\prime},X}P_{X}\beta^{\prime}+P_{X^{\prime}}\gamma^{\prime}\right) \\
&  =\Pi G_{X^{\prime},X}y\\
&  =G_{\pi\left(  X^{\prime}\right)  ,X},
\end{align*}

where the last line is implied by \ref{s16}. Thus the equation of the first
row of \ref{s36} is true, and the other equations are proved in a similar manner.

Using the notation of Definition \ref{Def_Matrices_from_R} we can write
\ref{s29} as
\[
s\left(  x\right)  =G_{x,X^{\prime}}\alpha^{\prime}+\sum\limits_{i=1}^{M}%
\beta_{i}^{\prime}p_{i}\left(  x\right)  =G_{x,X^{\prime}}\Pi^{T}\Pi
\alpha^{\prime}+\sum\limits_{i=1}^{M}\beta_{i}^{\prime}p_{i}\left(  x\right)
=G_{X,\pi\left(  X^{\prime}\right)  }\pi\left(  \alpha^{\prime}\right)
+\sum\limits_{i=1}^{M}\beta_{i}^{\prime}p_{i}\left(  x\right)  ,
\]

which proves \ref{s41} and finishes the proof.
\end{proof}

\begin{remark}
A different basis for $P_{\theta-1}$, say $\left\{  p_{j}^{\prime}\right\}  $,
could be used to define the unisolvency matrix $P_{X^{\prime}}$, the Riesz
representer $R_{x}^{\prime}$ and the space $W_{G,X^{\prime}}$. Using the
conventions of Notation \ref{Not_A'} we would write $P_{X^{\prime}}^{\prime}$
instead of $P_{X^{\prime}}$. All the results proved below would still hold.
\end{remark}

Recall from Definition \ref{Def_basis_distrib} that the set of basis functions
of order $\theta$ is $G+P_{2\theta-1}$, where $G$ is a particular basis
function. However, from part 5 of Theorem \ref{Thm_Appr_smth_prop1} we know
that the Approximate smoother function is independent of the basis function
chosen. The next corollary tells us that we can always choose a basis function
such that the Approximate smoother matrix is Hermitian and indeed this
particular matrix allows a simple proof of regularity which is given in
Theorem \ref{Thm_mat_psi}.

\begin{corollary}
\label{Cor_Thm_Appr_smth_mat_eqn_1}Suppose the weight function $w$ has
properties W2 and W3 for parameter $\theta$.

Then a (possibly complex-valued) basis function $G$ of order $\theta$ can
always be chosen so that $G_{Y,Z}^{T}=\overline{G_{Z,Y}}$ for all finite
subsets $Y$ and $Z$ of $\mathbb{R}^{d}$ i.e. so that $G\left(  -x\right)
=\overline{G\left(  x\right)  }$ for $x\in\mathbb{R}^{d}$. For such a basis
function the Approximate smoother matrix $\Psi$ of \ref{s16} is Hermitian and
we can write%
\begin{equation}
\Psi=\left(
\begin{array}
[c]{lll}%
\left(  2\pi\right)  ^{\frac{d}{2}}N\rho G_{X^{\prime},X^{\prime}}%
+\overline{G_{X,X^{\prime}}^{T}}G_{X,X^{\prime}} & \overline{G_{X,X^{\prime}%
}^{T}}P_{X} & P_{X^{\prime}}\\
P_{X}^{T}G_{X,X^{\prime}} & P_{X}^{T}P_{X} & O_{M}\\
P_{X^{\prime}}^{T} & O_{M} & O_{M}%
\end{array}
\right)  ,\label{s26}%
\end{equation}

where $M=\dim P_{\theta-1}$.
\end{corollary}

\begin{proof}
From part 2 Theorem \ref{Thm_Grho} there exists a basis function $G$ which
satisfies $G\left(  -x\right)  =\overline{G\left(  x\right)  }$ for
$x\in\mathbb{R}^{d}$. The condition $G_{Y,Z}^{T}=\overline{G_{Z,Y}}$ implies
that $G_{X^{\prime},X^{\prime}}^{T}=G_{X^{\prime},X^{\prime}}$ and
$G_{X,X^{\prime}}^{T}=G_{X^{\prime},X}$ and since $P_{X}$ and $P_{X^{\prime}}$
are real valued it follows that $\Psi$ has the required Hermitian form
\ref{s26}.
\end{proof}

In the next theorem we will prove the Approximate smoother matrix \ref{s26} is regular.

\begin{theorem}
\label{Thm_mat_psi}The Approximate smoother matrix $\Psi$ specified in
Corollary \ref{Cor_Thm_Appr_smth_mat_eqn_1} has the following properties:

\begin{enumerate}
\item $\Psi$ is a regular Hermitian matrix.

\item $\Psi$ is square with $N^{\prime}+2M$ rows. Hence the size of $\Psi$ is
independent of the number $N$ of (scattered) data points.
\end{enumerate}
\end{theorem}

\begin{proof}
\textbf{Part 1} We first write $\Psi$ in the form.
\[
\Psi=\left(
\begin{array}
[c]{cc}%
\left(  2\pi\right)  ^{\frac{d}{2}}N\rho G_{X^{\prime};0}+\overline{G_{P}^{T}%
}G_{P} &
\begin{array}
[c]{l}%
P_{X^{\prime}}\\
O_{M,N^{\prime}}%
\end{array}
\\%
\begin{array}
[c]{cc}%
P_{X^{\prime}}^{T} & O_{N^{\prime},M}%
\end{array}
& O_{M}%
\end{array}
\right)  ,
\]

where
\[
G_{X^{\prime};0}=\left(
\begin{array}
[c]{cc}%
G_{X^{\prime},X^{\prime}} & O\\
O & O_{M}%
\end{array}
\right)  ,\quad G_{P}=\left(
\begin{array}
[c]{cc}%
G_{X,X^{\prime}} & P_{X}%
\end{array}
\right)  .
\]

Since $G_{X,X}$ is Hermitian $\Psi$ is Hermitian. To prove that $\Psi$ is
regular we use Lemma \ref{Lem_reg_cspd}. To this end set
\begin{align*}
B  & =\left(  2\pi\right)  ^{\frac{d}{2}}N\rho G_{X^{\prime};0}+\overline
{G_{P}^{T}}G_{P},\\
C  & =\left(
\begin{array}
[c]{l}%
P_{X^{\prime}}\\
O_{M,N^{\prime}}%
\end{array}
\right)  .
\end{align*}

It must be shown that $\operatorname*{null}C=\left\{  0\right\}  $ and that
$z^{T}B\overline{z}=0$ and $C^{T}z=0$ implies $z=0$.

Firstly, $C\lambda=0$ implies $P_{X^{\prime}}\lambda=0$ which implies
$\lambda=0$ since $\operatorname*{null}P_{X^{\prime}}=\left\{  0\right\}  $.
Therefore $\operatorname*{null}C=\left\{  0\right\}  $. Next, set
$z^{T}=\left(  v^{T},\beta^{T}\right)  $. Hence, since
\[
C^{T}z=\left(
\begin{array}
[c]{cc}%
P_{X^{\prime}}^{T} & O
\end{array}
\right)  \left(
\begin{array}
[c]{c}%
v\\
\beta
\end{array}
\right)  =P_{X^{\prime}}^{T}v,
\]

we conclude that $C^{T}z=0$ implies $P_{X^{\prime}}^{T}v=0$. Now assume that
$C^{T}z=0$. Then
\begin{align*}
z^{T}B\overline{z}=z^{T}\left(  \left(  2\pi\right)  ^{\frac{d}{2}}N\rho
G_{X^{\prime};0}+\overline{G_{P}^{T}}G_{P}\right)  \overline{z} &  =\left(
2\pi\right)  ^{\frac{d}{2}}N\rho\,v^{T}G_{X^{\prime},X^{\prime}}\overline
{v}^{T}+z^{T}\overline{G_{P}^{T}}G_{P}\overline{z}\\
&  =\left(  2\pi\right)  ^{\frac{d}{2}}N\rho\,v^{T}G_{X^{\prime},X^{\prime}%
}\overline{v}+\left\vert G_{P}\overline{z}\right\vert _{\mathbb{C}^{N}}^{2}.
\end{align*}

By part 4 of Theorem \ref{Sum_ex_Wgx_properties_2}, $G_{X^{\prime},X^{\prime}%
}$ is conditionally positive definite on $\operatorname*{null}P_{X^{\prime}%
}^{T}$ i.e. $P_{X^{\prime}}^{T}v=0$ and $v\neq0$ implies $v^{T}G_{X^{\prime
},X^{\prime}}\overline{v}>0$. Therefore, $z^{T}B\overline{z}=0$ implies $v=0 $.

In addition, we have $G_{P}\overline{z}=0$ i.e. $0=G_{X,X^{\prime}}%
\overline{v}+P_{X}\overline{\beta}=P_{X}\overline{\beta}$, so that $\beta=0$
since $\operatorname*{null}P_{X}=\left\{  0\right\}  $ by part 1 Theorem
\ref{Thm_Px_properties}. We conclude that $z=0$.\medskip

\textbf{Part 2} We observe that $G_{X^{\prime},X^{\prime}}$ has size
$N^{\prime}\times N^{\prime}$, $G_{X,X^{\prime}}$ has size $N\times N^{\prime
}$,

$P_{X}$ has size $N\times M$, and $P_{X^{\prime}}$ has size $N^{\prime}\times
M$. From the block sizes it is clear that $\Psi$ is square with $N^{\prime
}+2M$ rows. Hence the size of $\Psi$ is independent of the number of
(scattered) data points $N$.
\end{proof}

Our next result shows the Approximate smoother algorithm is \textit{scalable}
i.e. the time of execution of construction and solution is linearly dependent
on the number of data points. This is in contrast with the Exact smoother
which is not scalable but which has quadratic dependency on the number of data points.

\begin{corollary}
\label{Cor_scalable}The Approximate smoother algorithm is scalable.
\end{corollary}

\begin{proof}
The Hermitian matrix equation \ref{s26} is
\[
\left(
\begin{array}
[c]{lll}%
\left(  2\pi\right)  ^{\frac{d}{2}}N\rho G_{X^{\prime},X^{\prime}}%
+\overline{G_{X,X^{\prime}}^{T}}G_{X,X^{\prime}} & \overline{G_{X,X^{\prime}%
}^{T}}P_{X} & P_{X^{\prime}}\\
P_{X}^{T}G_{X^{\prime},X} & P_{X}^{T}P_{X} & O_{M}\\
P_{X^{\prime}}^{T} & O_{M} & O_{M}%
\end{array}
\right)  \left(
\begin{array}
[c]{c}%
\alpha\\
\beta\\
\lambda
\end{array}
\right)  =\left(
\begin{array}
[c]{l}%
\overline{G_{X,X^{\prime}}^{T}}\\
P_{X}^{T}\\
O_{M,N}%
\end{array}
\right)  y,
\]

where $G_{X^{\prime},X^{\prime}}$ is $N^{\prime}\times N^{\prime}$,
$G_{X,X^{\prime}}$ is $N\times N^{\prime}$, $P_{X^{\prime}}$ is $N^{\prime
}\times M$ and $P_{X}$ is $N\times M$.

Suppose the \textit{evaluation cost} for $G\left(  x\right)  $ is $m_{G}$
multiplications and that $N\gg N^{\prime}$, $N\gg M$ and $N\gg m_{G}$. The
\textit{construction costs} for component matrices and matrix multiplications are:%

\begin{tabular}
[c]{|l|l|l|l|l|l|}\hline
\multicolumn{6}{|c|}{Construction costs for Approximate smoother
matrix}\\\hline
$G_{X^{\prime},X^{\prime}}$ & $G_{X,X^{\prime}}$ & $\overline{G_{X,X^{\prime}%
}^{T}}G_{X,X^{\prime}}$ & $P_{X}^{T}G_{X^{\prime},X}$ & $\Psi\left(
\begin{array}
[c]{c}%
\alpha\\
\beta\\
\lambda
\end{array}
\right)  $ & $\left(
\begin{array}
[c]{l}%
\overline{G_{X,X^{\prime}}^{T}}\\
P_{X}^{T}\\
O_{M,N}%
\end{array}
\right)  y$\\\hline
$\left(  N^{\prime}\right)  ^{2}m_{G}$ & $N^{\prime}Nm_{G}$ & $\left(
N^{\prime}\right)  ^{2}N$ & $N^{\prime}NM$ & $\left(  N^{\prime}+2M\right)
^{2}$ & $N\left(  N^{\prime}+2M\right)  $\\\hline
\end{tabular}

From the table the dominant \textit{construction costs} are $N^{\prime
}N\left(  N^{\prime}+m_{G}+M\right)  $. The \textit{solution cost} of an
$N^{\prime}\times N^{\prime}$ matrix equation is $\frac{1}{3}\left(
N^{\prime}\right)  ^{3}$ multiplications for a \textit{dense} matrix. Thus the
total cost is
\begin{align*}
N^{\prime}N\left(  N^{\prime}+m_{G}+M\right)  +\frac{1}{3}\left(  N^{\prime
}\right)  ^{3}  & =\left(  \left(  m_{G}+M\right)  N^{\prime}+\left(
N^{\prime}\right)  ^{2}\right)  N+\frac{1}{3}\left(  N^{\prime}\right)  ^{3}\\
& \gg\left(  \left(  m_{G}+M\right)  N^{\prime}+\left(  N^{\prime}\right)
^{2}\right)  N,
\end{align*}

which is linearly dependent on the number of data points. However, by the use
of a basis function with support containing only several points in $X^{\prime
}$ e.g. the extended natural spline basis functions of Lemma
\ref{Lem_basis_ord0_Bsplin_W3}, the construction and solution costs can be
reduced significantly. However we still have linear dependency on $N$.
\end{proof}

\section{Solving the Approximate smoothing problem using matrix techniques
\label{Sect_Appr_smth_mat}}

In the last section, by using Hilbert space orthogonal projection techniques,
the Approximate smoother matrix equation was derived for complex data and all
conjugate-even basis functions. In this section we present an alternative
derivation which uses matrix techniques. A key step in the Hilbert space
approach was to choose a congugate-even basis function. This was sufficiently
general and enabled us to prove that the Approximate smoother matrix was
regular. This derivation is less general and imposes the \textbf{extra
conditions} that the \textbf{basis function and data are real valued}. This is
to simplify the use of the method of Lagrange multipliers and we note that the
Approximate smoother matrix will now be symmetric and the basis function will
be even.

The Approximate smoothing problem is $\min\limits_{f\in W_{G,X^{\prime}}}%
J_{e}\left[  f\right]  $ where $J_{e}\left[  f\right]  $ is the Exact smoother
functional. The first step will be to calculate $J_{e}\left[  f\right]  $ for
$f\in W_{G,X^{\prime}}$. The result is that $J_{e}\left[  f\right]  $ is a
constrained quadratic form \ref{s76} in terms of the coefficients of the
data-translated basis functions and the polynomial basis functions. Then it is
shown that if the data is real and the basis function is real valued then the
solution is unique and real valued. Finally, Lagrange multipliers are used to
minimize the constrained quadratic form and derive the matrix equation as well
as the identities derived using the Hilbert space method.

\subsection{The Exact smoother functional for functions in $W_{G,X^{\prime}}
$}

In this subsection we will calculate the Exact smoothing functional
$J_{e}\left[  f\right]  $ for $f\in W_{G,X^{\prime}}$.

\begin{theorem}
\label{Thm_Je_on_Wgz}Suppose the basis function $G$ is a \textbf{real valued}
\textbf{even function} and $f\in W_{G,X^{\prime}}$, where $X^{\prime}=\left\{
x_{i}^{\prime}\right\}  _{i=1}^{N^{\prime}}$ is unisolvent. Then if
\[
f=\sum_{i=1}^{N^{\prime}}\alpha_{i}G\left(  \cdot-x_{i}^{\prime}\right)
+\sum_{j=1}^{M}\beta_{j}p_{j}\in W_{G,X^{\prime}},\quad\alpha_{i},\beta_{j}%
\in\mathbb{C},
\]

we have%
\begin{align}
J_{e}\left[  f\right]  =\left(  \alpha^{T}\,\beta^{T}\right)   &  \left(
\left(  2\pi\right)  ^{\frac{d}{2}}\rho G_{X^{\prime};0}+\frac{1}{N}G_{P}%
^{T}G_{P}\right)  \left(
\begin{array}
[c]{c}%
\overline{\alpha}\\
\overline{\beta}%
\end{array}
\right)  -\frac{1}{N}y^{T}G_{P}\left(
\begin{array}
[c]{c}%
\alpha\\
\beta
\end{array}
\right)  -\nonumber\\
&  -\frac{1}{N}y^{T}G_{P}\left(
\begin{array}
[c]{c}%
\overline{\alpha}\\
\overline{\beta}%
\end{array}
\right)  +\frac{1}{N}\left\vert y\right\vert ^{2}.\label{s76}%
\end{align}

where $J_{e}[\cdot]$ is the Exact smoother functional \ref{h52} and
\[
G_{X^{\prime};0}=\left(
\begin{array}
[c]{cc}%
G_{X^{\prime},X^{\prime}} & O\\
O & O_{M}%
\end{array}
\right)  ,\quad G_{P}=\left(
\begin{array}
[c]{ll}%
G_{X,X^{\prime}} & P_{X}%
\end{array}
\right)  .
\]

\end{theorem}

\begin{proof}
The Exact smoothing functional is
\[
J_{e}\left[  f\right]  =\rho\left\vert f\right\vert _{w,\theta}^{2}+\frac
{1}{N}\sum\limits_{i=1}^{N}\left\vert f\left(  x^{\left(  i\right)  }\right)
-y^{\left(  i\right)  }\right\vert ^{2},
\]

and from part 1 of Theorem \ref{Sum_ex_Wgx_properties_2}%
\begin{equation}
\rho\left\vert f\right\vert _{w,\theta}^{2}=\rho_{\pi}\alpha^{T}G_{X^{\prime
},X^{\prime}}\alpha=\rho_{\pi}\left(  \alpha^{T}\,\beta^{T}\right)
G_{X^{\prime};0}\left(
\begin{array}
[c]{c}%
\overline{\alpha}\\
\overline{\beta}%
\end{array}
\right)  .\label{s28}%
\end{equation}

The next step is to approximate the second `least squares' term $\sum
\limits_{i=1}^{N}\left\vert f\left(  x^{\left(  i\right)  }\right)
-y^{\left(  i\right)  }\right\vert ^{2}$. Let $f_{X}=\left(  f\left(
x^{\left(  i\right)  }\right)  \right)  _{i=1}^{N}$. Then
\begin{equation}
\sum_{i=1}^{N}\left\vert f(x^{\left(  i\right)  })-y^{\left(  i\right)
}\right\vert ^{2}=\left(  f_{X}-y\right)  ^{T}\left(  \overline{f_{X}%
}-y\right)  =f_{X}^{T}\overline{f_{X}}-y^{T}f_{X}-y^{T}\overline{f_{X}%
}+\left\vert y\right\vert ^{2}.\label{s27}%
\end{equation}

But in matrix form
\[
f_{X}=G_{X,X^{\prime}}\alpha+P_{X}\beta,
\]

where $G_{X,X^{\prime}}=\left(  G\left(  x^{\left(  i\right)  }-x_{j}^{\prime
}\right)  \right)  $ and $P_{X}=\left(  p_{j}\left(  x^{\left(  i\right)
}\right)  \right)  $, or in block form
\[
f_{X}=\left(
\begin{array}
[c]{ll}%
G_{X,X^{\prime}} & P_{X}%
\end{array}
\right)  \left(
\begin{array}
[c]{c}%
\alpha\\
\beta
\end{array}
\right)  =G_{P}\left(
\begin{array}
[c]{c}%
\alpha\\
\beta
\end{array}
\right)  ,
\]

where $G_{P}=\left(
\begin{array}
[c]{ll}%
G_{X,X^{\prime}} & P_{X}%
\end{array}
\right)  $. Thus
\[
f_{X}^{T}\overline{f_{X}}=\left(  \alpha^{T}\,\beta^{T}\right)  G_{P}^{T}%
G_{P}\left(
\begin{array}
[c]{c}%
\overline{\alpha}\\
\overline{\beta}%
\end{array}
\right)  ,
\]

and by \ref{s27}
\begin{equation}
\sum_{i=1}^{N}\left\vert f\left(  x^{\left(  i\right)  }\right)  -y^{\left(
i\right)  }\right\vert ^{2}=\left(  \alpha^{T}\,\beta^{T}\right)  G_{P}%
^{T}G_{P}\left(
\begin{array}
[c]{c}%
\overline{\alpha}\\
\overline{\beta}%
\end{array}
\right)  -y^{T}G_{P}\left(
\begin{array}
[c]{c}%
\alpha\\
\beta
\end{array}
\right)  -y^{T}G_{P}\left(
\begin{array}
[c]{c}%
\overline{\alpha}\\
\overline{\beta}%
\end{array}
\right)  +\left\vert y\right\vert ^{2},\label{s30}%
\end{equation}

so that combining \ref{s28} and \ref{s30}
\begin{align*}
J_{e}\left[  f\right]   &  =\rho\left\vert f\right\vert _{w,\theta}^{2}%
+\frac{1}{N}\sum_{i=1}^{N}\left\vert f\left(  x^{\left(  i\right)  }\right)
-y^{\left(  i\right)  }\right\vert ^{2}\\
&  =\rho_{\pi}\left(  \alpha^{T}\,\beta^{T}\right)  G_{X^{\prime};0}\left(
\begin{array}
[c]{c}%
\overline{\alpha}\\
\overline{\beta}%
\end{array}
\right)  +\frac{1}{N}\left(  \alpha^{T}\,\beta^{T}\right)  G_{P}^{T}%
G_{P}\left(
\begin{array}
[c]{c}%
\overline{\alpha}\\
\overline{\beta}%
\end{array}
\right)  -\\
&  \qquad\qquad-\frac{1}{N}y^{T}G_{P}\left(
\begin{array}
[c]{c}%
\alpha\\
\beta
\end{array}
\right)  -\frac{1}{N}y^{T}G_{P}\left(
\begin{array}
[c]{c}%
\overline{\alpha}\\
\overline{\beta}%
\end{array}
\right)  +\frac{1}{N}\left\vert y\right\vert ^{2}\\
&  =\left(  \alpha^{T}\,\beta^{T}\right)  \left(  \rho_{\pi}G_{X^{\prime}%
;0}+\frac{1}{N}G_{P}^{T}G_{P}\right)  \left(
\begin{array}
[c]{c}%
\overline{\alpha}\\
\overline{\beta}%
\end{array}
\right)  -\\
&  \qquad\qquad-\frac{1}{N}y^{T}G_{P}\left(
\begin{array}
[c]{c}%
\alpha\\
\beta
\end{array}
\right)  -\frac{1}{N}y^{T}G_{P}\left(
\begin{array}
[c]{c}%
\overline{\alpha}\\
\overline{\beta}%
\end{array}
\right)  +\frac{1}{N}\left\vert y\right\vert ^{2}.
\end{align*}

\end{proof}

\subsection{Proof that the smoother is unique}

We will now supply another proof, this time based on matrices, that there
\textit{exists} a unique solution in $W_{G,X^{\prime}}$ for the Approximate
smoother problem. After that we will obtain a matrix equation for this solution.

\begin{theorem}
\label{Thm_Appr_smth_unique}The Approximate smoothing problem has a unique
solution in $W_{G,X^{\prime}}$.
\end{theorem}

\begin{proof}
The Approximate smoother problem involves minimizing the quadratic form
\ref{s76} constrained by $P_{X^{\prime}}^{T}\alpha=0$. But from part 2 of
Theorem \ref{Thm_Px_properties_2} the null space of $P_{X^{\prime}}^{T}$ has
the form $\alpha=%
\begin{pmatrix}
-L_{X_{2}^{\prime}}^{T}\\
I_{N^{\prime}-M}%
\end{pmatrix}
\alpha^{\prime\prime}$ where $\alpha^{\prime\prime}=\left(  \alpha_{i}\right)
_{i=M+1}^{N^{\prime}}\in\mathbb{R}^{\left(  N^{\prime}-M\right)  }$ and
$X_{2}^{\prime}=\left\{  x_{i}^{\prime}\right\}  _{i=M+1}^{N^{\prime}}$. Set
$A=%
\begin{pmatrix}
-L_{X_{2}^{\prime}}^{T}\\
I_{N^{\prime}-M}%
\end{pmatrix}
$ so that $\left(
\begin{array}
[c]{c}%
\alpha\\
\beta
\end{array}
\right)  =%
\begin{pmatrix}
A & O\\
O & I_{N^{\prime}}%
\end{pmatrix}
\left(
\begin{array}
[c]{c}%
\alpha^{\prime\prime}\\
\beta
\end{array}
\right)  $ and the first term on the right of \ref{s76} becomes%
\begin{align}
&  \left(  \alpha^{T}\,\beta^{T}\right)  \left(  \rho_{\pi}G_{X^{\prime}%
;0}+\frac{1}{N}G_{P}^{T}G_{P}\right)  \left(
\begin{array}
[c]{c}%
\overline{\alpha}\\
\overline{\beta}%
\end{array}
\right) \label{s75}\\
&  =\left(  \alpha^{\prime\prime T}\,\beta^{T}\right)
\begin{pmatrix}
A^{T} & O\\
O & I_{N^{\prime}}%
\end{pmatrix}
\left(  \rho_{\pi}G_{X^{\prime};0}+\frac{1}{N}G_{P}^{T}G_{P}\right)
\begin{pmatrix}
A & O\\
O & I_{N^{\prime}}%
\end{pmatrix}
\left(
\begin{array}
[c]{c}%
\overline{\alpha^{\prime\prime}}\\
\overline{\beta}%
\end{array}
\right)  .\nonumber
\end{align}

Thus the Approximate smoother problem is equivalent to minimizing \ref{s75}
for all $\alpha^{\prime\prime}\in\mathbb{R}^{\left(  N^{\prime}-M\right)  }$
and $\beta\in\mathbb{R}^{N^{\prime}}$. But if the matrix
\begin{equation}%
\begin{pmatrix}
A^{T} & O\\
O & I_{N^{\prime}}%
\end{pmatrix}
\left(  \rho_{\pi}G_{X^{\prime};0}+\frac{1}{N}G_{P}^{T}G_{P}\right)
\begin{pmatrix}
A & O\\
O & I_{N^{\prime}}%
\end{pmatrix}
,\label{s19}%
\end{equation}

is positive definite then the quadratic form \ref{s75} has a unique stationary
point and this is a minimum point. Indeed%
\begin{align}
\left(  \alpha^{\prime\prime T}\,\beta^{T}\right)
\begin{pmatrix}
A^{T} & O\\
O & I_{N^{\prime}}%
\end{pmatrix}
&  \left(  \rho_{\pi}G_{X^{\prime};0}+\frac{1}{N}G_{P}^{T}G_{P}\right)
\begin{pmatrix}
A & O\\
O & I_{N^{\prime}}%
\end{pmatrix}
\left(
\begin{array}
[c]{c}%
\overline{\alpha^{\prime\prime}}\\
\overline{\beta}%
\end{array}
\right) \nonumber\\
&  =\left(  \alpha^{T}\,\beta^{T}\right)  \left(  \rho_{\pi}G_{X^{\prime}%
;0}+\frac{1}{N}G_{P}^{T}G_{P}\right)  \left(
\begin{array}
[c]{c}%
\overline{\alpha}\\
\overline{\beta}%
\end{array}
\right) \nonumber\\
&  =\rho_{\pi}\alpha^{T}G_{X^{\prime},X^{\prime}}\overline{\alpha}+\frac{1}%
{N}\left\vert \left(
\begin{array}
[c]{ll}%
G_{X,X^{\prime}} & P_{X}%
\end{array}
\right)  \left(
\begin{array}
[c]{c}%
\alpha\\
\beta
\end{array}
\right)  \right\vert _{\mathbb{C}^{N}}^{2}\nonumber\\
&  =\rho_{\pi}\alpha^{T}G_{X^{\prime},X^{\prime}}\overline{\alpha}+\frac{1}%
{N}\left\vert G_{X,X^{\prime}}\alpha+P_{X}\beta\right\vert _{\mathbb{C}^{N}%
}^{2},\label{s13}%
\end{align}

where $P_{X^{\prime}}^{T}\alpha=0$. So the matrix \ref{s19} is positive
definite if the matrix \ref{s13} is conditionally positive definite on
$\operatorname*{null}P_{X^{\prime}}^{T}$. So suppose%
\begin{equation}
\rho_{\pi}\alpha^{T}G_{X^{\prime},X^{\prime}}\overline{\alpha}+\frac{1}%
{N}\left\vert G_{X,X^{\prime}}\alpha+P_{X}\beta\right\vert _{\mathbb{C}^{N}%
}^{2}=0,\quad P_{X^{\prime}}^{T}\alpha=0.\label{s50}%
\end{equation}

Now by part 2 of Theorem \ref{Sum_ex_Wgx_properties_2}, $G_{X^{\prime
},X^{\prime}}$ is conditionally positive definite on $P_{X^{\prime}}^{T}%
\alpha=0$ i.e. when $P_{X^{\prime}}^{T}\alpha=0$ we have $\alpha
^{T}G_{X^{\prime},X^{\prime}}\overline{\alpha}\geq0$, and $P_{X^{\prime}}%
^{T}\alpha=0$ and $\alpha^{T}G_{X^{\prime},X^{\prime}}\overline{\alpha}=0$
implies $\alpha=0$. Thus $\alpha^{T}G_{X^{\prime},X^{\prime}}\overline{\alpha
}=0$ and $G_{X,X^{\prime}}\alpha+P_{X}\beta=0$ i.e. $\alpha=0$ and $P_{X}%
\beta=0$. But from part 1 of Theorem \ref{Thm_Px_properties}
$\operatorname*{null}P_{X}=\left\{  0\right\}  $, and this implies $\beta=0$.
\end{proof}

The next result will simplify the use of Lagrange multipliers to solve the
Approximate smoothing problem.

\begin{theorem}
\label{Thm_Approx_smth_real_val}Suppose the basis function $G$ is real-valued
with order $\theta\geq1$. Then the solution to the Approximate smoothing
problem with real valued data $\left[  X,y\right]  $ induced by the points
$X^{\prime}=\left\{  x_{i}^{\prime}\right\}  _{i=1}^{N^{\prime}}$ lies in the
subspace of $W_{G,X^{\prime}}$ defined by the real scalars.
\end{theorem}

\begin{proof}
The Approximate smoothing problem is $\min\limits_{f\in W_{G,X^{\prime}}}%
J_{e}\left[  f\right]  $ and by the previous theorem this problem has a unique
solution. Now let $f=\sum\limits_{i=1}^{N^{\prime}}\alpha_{i}G\left(
\cdot-x_{i}^{\prime}\right)  +\sum\limits_{j=1}^{M}\beta_{j}p_{j}\in
W_{G,X^{\prime}}$, where $\alpha_{i},\beta_{j}\in\mathbb{C}$. Then by Theorem
\ref{Thm_Je_on_Wgz}
\[
J_{e}\left[  f\right]  =\left(  \alpha^{T}\,\beta^{T}\right)  \left(
\rho_{\pi}G_{X^{\prime};0}+\frac{1}{N}G_{P}^{T}G_{P}\right)  \left(
\begin{array}
[c]{c}%
\overline{\alpha}\\
\overline{\beta}%
\end{array}
\right)  -\frac{1}{N}y^{T}G_{P}\left(
\begin{array}
[c]{c}%
\alpha\\
\beta
\end{array}
\right)  -\frac{1}{N}y^{T}G_{P}\left(
\begin{array}
[c]{c}%
\overline{\alpha}\\
\overline{\beta}%
\end{array}
\right)  +\frac{1}{N}\left\vert y\right\vert ^{2},
\]

and so%
\begin{align*}
J_{e}\left[  \overline{f}\right]   &  =\left(  \overline{\alpha}%
^{T}\,\overline{\beta}^{T}\right)  \left(  \rho_{\pi}G_{X^{\prime};0}+\frac
{1}{N}G_{P}^{T}G_{P}\right)  \left(
\begin{array}
[c]{c}%
\alpha\\
\beta
\end{array}
\right)  -\frac{1}{N}y^{T}G_{P}\left(
\begin{array}
[c]{c}%
\overline{\alpha}\\
\overline{\beta}%
\end{array}
\right)  -\frac{1}{N}y^{T}G_{P}\left(
\begin{array}
[c]{c}%
\alpha\\
\beta
\end{array}
\right)  +\\
&  \qquad\qquad\qquad+\frac{1}{N}\left\vert y\right\vert ^{2}\\
&  =\left(  \alpha^{T}\,\beta^{T}\right)  \left(  \rho_{\pi}G_{X^{\prime}%
;0}+\frac{1}{N}G_{P}^{T}G_{P}\right)  ^{T}\left(
\begin{array}
[c]{c}%
\overline{\alpha}\\
\overline{\beta}%
\end{array}
\right)  -\frac{1}{N}y^{T}G_{P}\left(
\begin{array}
[c]{c}%
\alpha\\
\beta
\end{array}
\right)  -\frac{1}{N}y^{T}G_{P}\left(
\begin{array}
[c]{c}%
\overline{\alpha}\\
\overline{\beta}%
\end{array}
\right)  +\\
&  \qquad\qquad\qquad+\frac{1}{N}\left\vert y\right\vert ^{2}\\
&  =\left(  \alpha^{T}\,\beta^{T}\right)  \left(  \rho_{\pi}G_{X^{\prime}%
;0}+\frac{1}{N}G_{P}^{T}G_{P}^{T}\right)  \left(
\begin{array}
[c]{c}%
\overline{\alpha}\\
\overline{\beta}%
\end{array}
\right)  -\frac{1}{N}y^{T}G_{P}\left(
\begin{array}
[c]{c}%
\alpha\\
\beta
\end{array}
\right)  -\frac{1}{N}y^{T}G_{P}\left(
\begin{array}
[c]{c}%
\overline{\alpha}\\
\overline{\beta}%
\end{array}
\right)  +\\
&  \qquad\qquad\qquad+\frac{1}{N}\left\vert y\right\vert ^{2}\\
&  =J_{e}\left[  f\right]  .
\end{align*}

\end{proof}

Thus, if $s$ is the Approximate smoother then $J_{e}\left[  s\right]
=J_{e}\left[  \overline{s}\right]  $ and $\overline{s}$ is also a solution to
the Approximate smoothing problem. But the solution is unique so
$s=\overline{s}$ and $s$ is real-valued. Hence if $s=\sum\limits_{i=1}%
^{N^{\prime}}\alpha_{i}G\left(  \cdot-x_{i}^{\prime}\right)  +\sum
\limits_{j=1}^{M}\beta_{j}p_{j}$ then $s=\overline{s}=\sum\limits_{i=1}%
^{N^{\prime}}\overline{\alpha_{i}}G\left(  \cdot-x_{i}^{\prime}\right)
+\sum\limits_{j=1}^{M}\overline{\beta_{j}}p_{j}$ and so $\sum\limits_{i=1}%
^{N^{\prime}}\left(  \alpha_{i}-\overline{\alpha_{i}}\right)  G\left(
\cdot-x_{i}^{\prime}\right)  +\sum\limits_{j=1}^{M}\left(  \beta_{j}%
-\overline{\beta_{j}}\right)  p_{j}=0$. The unique representation result of
part 5 Theorem \ref{Sum_ex_Wgx_properties_2} now implies that for all $i$ and
$j$, $\alpha_{i}=\overline{\alpha_{i}}$ and $\beta_{j}=\overline{\beta_{j}}$
i.e. $\alpha_{i}$ and $\beta_{j}$ are real. This means that $\sum
\limits_{j=1}^{M}\beta_{j}p_{j}\in P_{\theta-1}$ and we are done.

\subsection{The Approximate smoother matrix equation - a Lagrange multiplier
derivation\label{SbSect_approx_smth_Lagrange}}

Recall that the Approximate smoothing problem is $\min\limits_{f\in
W_{G,X^{\prime}}}J_{e}\left[  f\right]  $ and so this involves minimizing the
quadratic form \ref{s76} constrained by $P_{X^{\prime}}^{T}\alpha=0$. In
Theorem \ref{Thm_Appr_smth_unique} it was proved that a solution to the
Approximate smoothing problem exists and is unique. The proof also
demonstrated that the constrained problem only has one stationary point which
means the technique of Lagrange multipliers will always yield the Approximate
smoother. For convenience denote the quadratic form \ref{s76} by $Q_{a}\left(
\alpha,\beta\right)  $ and denote the Lagrange multiplier by the vector
$\lambda$ with length $N^{\prime}$.\ Since we have assumed the basis function
is real-valued Theorem \ref{Thm_Approx_smth_real_val} tells us that the
Approximate smoother can be expressed in terms of real $\alpha$ and $\beta$.

In preparation we expand the matrix equation \ref{s76} for $Q_{a}\left(
\alpha,\beta\right)  $ to include $\lambda$ and for algebraic simplicity we
will minimize $NQ_{a}$ instead of $Q_{a}$. Thus for real vectors $\alpha$ and
$\beta$, and $\rho_{\pi}=\left(  2\pi\right)  ^{\frac{d}{2}}\rho$
\begin{multline*}
NQ_{a}\left(  \alpha,\beta\right)  =\left(  \alpha^{T}\,\beta^{T}%
\,\vdots\,\lambda^{T}\right)  \left(
\begin{array}
[c]{cc}%
N\rho_{\pi}G_{X^{\prime};0}+G_{P}^{T}G_{P} & O\\
O & O_{N^{\prime}}%
\end{array}
\right)  \left(
\begin{array}
[c]{c}%
\alpha\\
\beta\\
\cdots\\
\lambda
\end{array}
\right)  -\\
-2y^{T}\left(
\begin{array}
[c]{cc}%
G_{P} & O_{N,N^{\prime}}%
\end{array}
\right)  \left(
\begin{array}
[c]{c}%
\alpha\\
\beta\\
\cdots\\
\lambda
\end{array}
\right)  +\left\vert y\right\vert ^{2},\quad\alpha\in\mathbb{R}^{N^{\prime}%
},\beta\in\mathbb{R}^{M}.
\end{multline*}

Let the Lagrangian be denoted by $L(\alpha,\beta,\lambda)$ so that
\begin{align}
L(\alpha,\beta,\lambda) &  =NQ_{a}\left(  \alpha,\beta\right)  +2\lambda
^{T}P_{X^{\prime}}^{T}\beta\nonumber\\
&  =NQ_{a}\left(  \alpha,\beta\right)  +\left(  \alpha^{T}\,\beta^{T}%
\,\lambda^{T}\right)  \left(
\begin{array}
[c]{ccc}%
O & O & P_{X^{\prime}}\\
O & O & O\\
P_{X^{\prime}}^{T} & O & O_{N^{\prime}}%
\end{array}
\right)  \left(
\begin{array}
[c]{c}%
\alpha\\
\beta\\
\lambda
\end{array}
\right) \nonumber\\
&  =\left(  \alpha^{T}\,\beta^{T}\,\vdots\,\lambda^{T}\right)  \left(
\begin{array}
[c]{cc}%
N\rho_{\pi}G_{X^{\prime};0}+G_{P}^{T}G_{P} & O\\
O & O_{N^{\prime}}%
\end{array}
\right)  \left(
\begin{array}
[c]{c}%
\alpha\\
\beta\\
\cdots\\
\lambda
\end{array}
\right)  -\nonumber\\
&  \qquad\qquad-2y^{T}\left(
\begin{array}
[c]{cc}%
G_{P} & O_{N,N^{\prime}}%
\end{array}
\right)  \left(
\begin{array}
[c]{c}%
\alpha\\
\beta\\
\cdots\\
\lambda
\end{array}
\right)  +\left\vert y\right\vert ^{2}+\nonumber\\
&  \qquad\qquad+\left(  \alpha^{T}\,\beta^{T}\,\lambda^{T}\right)  \left(
\begin{array}
[c]{ccc}%
O_{N^{\prime}} & O & P_{X^{\prime}}\\
O & O_{M} & O\\
P_{X^{\prime}}^{T} & O & O_{N^{\prime}}%
\end{array}
\right)  \left(
\begin{array}
[c]{c}%
\alpha\\
\beta\\
\lambda
\end{array}
\right) \nonumber\\
&  =\left(  \alpha^{T}\,\beta^{T}\,\vdots\,\lambda^{T}\right)  \left(
\begin{array}
[c]{cc}%
N\rho_{\pi}G_{X^{\prime};0}+G_{P}^{T}G_{P} &
\begin{array}
[c]{l}%
P_{X^{\prime}}\\
O
\end{array}
\\%
\begin{array}
[c]{cc}%
P_{X^{\prime}}^{T} & O
\end{array}
& O_{N^{\prime}}%
\end{array}
\right)  \left(
\begin{array}
[c]{c}%
\alpha\\
\beta\\
\cdots\\
\lambda
\end{array}
\right)  -\nonumber\\
&  \qquad\qquad-2y^{T}\left(
\begin{array}
[c]{cc}%
G_{P}^{T} & O_{N,N^{\prime}}%
\end{array}
\right)  \left(
\begin{array}
[c]{c}%
\alpha\\
\beta\\
\lambda
\end{array}
\right)  +|y|^{2}\label{s92}%
\end{align}

A necessary condition for a minimum is that the derivative is zero i.e.
$DL(\alpha,\beta,\lambda)=0$. Differentiating equation \ref{s92} gives%
\begin{align*}
DL &  =\left(  \alpha^{T}\,\beta^{T}\,\vdots\,\lambda^{T}\right)  \left(
\begin{array}
[c]{cc}%
N\rho_{\pi}G_{X^{\prime};0}+G_{P}^{T}G_{P} &
\begin{array}
[c]{l}%
P_{X^{\prime}}\\
O
\end{array}
\\%
\begin{array}
[c]{cc}%
P_{X^{\prime}}^{T} & O
\end{array}
& O_{N^{\prime}}%
\end{array}
\right)  +\\
&  \qquad\qquad+\left(  \alpha^{T}\,\beta^{T}\,\vdots\,\lambda^{T}\right)
\left(
\begin{array}
[c]{cc}%
N\rho_{\pi}G_{X^{\prime};0}^{T}+G_{P}^{T}G_{P} &
\begin{array}
[c]{l}%
P_{X^{\prime}}\\
O
\end{array}
\\%
\begin{array}
[c]{cc}%
P_{X^{\prime}}^{T} & O
\end{array}
& O_{N^{\prime}}%
\end{array}
\right)  ^{T}-2y^{T}\left(
\begin{array}
[c]{cc}%
G_{P} & O
\end{array}
\right) \\
&  =2\left(  \alpha^{T}\,\beta^{T}\,\vdots\,\lambda^{T}\right)  \left(
\begin{array}
[c]{cc}%
N\rho_{\pi}G_{X^{\prime};0}+G_{P}^{T}G_{P} &
\begin{array}
[c]{l}%
P_{X^{\prime}}\\
O
\end{array}
\\%
\begin{array}
[c]{cc}%
P_{X^{\prime}}^{T} & O
\end{array}
& O_{N^{\prime}}%
\end{array}
\right)  -2y^{T}\left(
\begin{array}
[c]{cc}%
G_{P} & O
\end{array}
\right)
\end{align*}

At the minimum, $(DL)^{T}=0$
\begin{equation}
\left(
\begin{array}
[c]{cc}%
N\rho_{\pi}G_{X^{\prime};0}+G_{P}^{T}G_{P} &
\begin{array}
[c]{l}%
P_{X^{\prime}}\\
O
\end{array}
\\%
\begin{array}
[c]{cc}%
P_{X^{\prime}}^{T} & O
\end{array}
& O_{N^{\prime}}%
\end{array}
\right)  \left(
\begin{array}
[c]{c}%
\alpha\\
\beta\\
\cdots\\
\lambda
\end{array}
\right)  =\left(
\begin{array}
[c]{c}%
G_{P}^{T}\\
O
\end{array}
\right)  y,\label{s93}%
\end{equation}

and since $G_{P}=\left(
\begin{array}
[c]{cc}%
G_{X,X^{\prime}} & P_{X}%
\end{array}
\right)  $,
\begin{align*}
N\rho_{\pi}G_{X^{\prime};0}+G_{P}^{T}G_{P} &  =N\rho_{\pi}G_{X^{\prime}%
;0}+\left(
\begin{array}
[c]{cc}%
G_{X,X^{\prime}} & P_{X}%
\end{array}
\right)  ^{T}\left(
\begin{array}
[c]{cc}%
G_{X,X^{\prime}} & P_{X}%
\end{array}
\right) \\
&  =N\rho_{\pi}G_{X^{\prime};0}+\left(
\begin{array}
[c]{l}%
G_{X,X^{\prime}}^{T}\\
P_{X}^{T}%
\end{array}
\right)  \left(
\begin{array}
[c]{cc}%
G_{X,X^{\prime}} & P_{X}%
\end{array}
\right) \\
&  =N\rho_{\pi}G_{X^{\prime};0}+\left(
\begin{array}
[c]{ll}%
G_{X,X^{\prime}}^{T}G_{X,X^{\prime}} & G_{X,X^{\prime}}^{T}P_{X}\\
P_{X}^{T}G_{X,X^{\prime}} & P_{X}^{T}P_{X}%
\end{array}
\right) \\
&  =\left(
\begin{array}
[c]{cc}%
N\rho_{\pi}G_{X^{\prime},X^{\prime}} & O\\
O & O
\end{array}
\right)  +\left(
\begin{array}
[c]{ll}%
G_{X,X^{\prime}}^{T}G_{X,X^{\prime}} & G_{X,X^{\prime}}^{T}P_{X}\\
P_{X}^{T}G_{X,X^{\prime}} & P_{X}^{T}P_{X}%
\end{array}
\right) \\
&  =\left(
\begin{array}
[c]{ll}%
N\rho_{\pi}G_{X^{\prime},X^{\prime}}+G_{X,X^{\prime}}^{T}G_{X,X^{\prime}} &
G_{X,X^{\prime}}^{T}P_{X}\\
P_{X}^{T}G_{X,X^{\prime}} & P_{X}^{T}P_{X}%
\end{array}
\right)
\end{align*}

Thus
\begin{equation}
\left(
\begin{array}
[c]{cc}%
N\rho_{\pi}G_{X^{\prime};0}+G_{P}^{T}G_{P} &
\begin{array}
[c]{l}%
P_{X^{\prime}}\\
O
\end{array}
\\%
\begin{array}
[c]{cc}%
P_{X^{\prime}}^{T} & O
\end{array}
& O_{N^{\prime}}%
\end{array}
\right)  =\left(
\begin{array}
[c]{lll}%
N\rho_{\pi}G_{X^{\prime},X^{\prime}}+G_{X,X^{\prime}}^{T}G_{X,X^{\prime}} &
G_{X,X^{\prime}}^{T}P_{X} & P_{X^{\prime}}\\
P_{X}^{T}G_{X,X^{\prime}} & P_{X}^{T}P_{X} & O_{M}\\
P_{X^{\prime}}^{T} & O_{M} & O_{M}%
\end{array}
\right)  ,\label{s94}%
\end{equation}

and noting that $\rho_{\pi}=\left(  2\pi\right)  ^{\frac{d}{2}}\rho$, equation
\ref{s93} becomes
\[
\left(
\begin{array}
[c]{lll}%
\left(  2\pi\right)  ^{\frac{d}{2}}N\rho G_{X^{\prime},X^{\prime}%
}+G_{X,X^{\prime}}^{T}G_{X,X^{\prime}} & G_{X,X^{\prime}}^{T}P_{X} &
P_{X^{\prime}}\\
P_{X}^{T}G_{X,X^{\prime}} & P_{X}^{T}P_{X} & O_{M}\\
P_{X^{\prime}}^{T} & O_{M} & O_{M}%
\end{array}
\right)  \left(
\begin{array}
[c]{c}%
\alpha\\
\beta\\
\lambda
\end{array}
\right)  =\left(
\begin{array}
[c]{l}%
G_{X,X^{\prime}}^{T}\\
P_{X}^{T}\\
O
\end{array}
\right)  y,
\]

which matches the Approximate matrix equation \ref{s26} which was derived
using Hilbert space techniques and assumed a complex valued basis function.

\section{Convergence of the Approximate smoother to the Exact
smoother\label{Sect_Appr_smth_conv}}

In this section we will study the uniform pointwise convergence on a bounded
region of the Approximate smoother $s_{a}$ to the Exact smoother $s_{e}$ and
to the data function $f_{d}$.

In Subsection \ref{SbSect_App_to_Exact_no_order} we study the uniform
convergence of the Approximate smoother to the Exact smoother. Here we make no
use of Lagrange interpolation theory and obtain no order of convergence results.

In Subsection \ref{SbSect_App_minus_Exact} we again consider the convergence
of the Approximate smoother to the Exact smoother but here we use Lagrange
interpolation theory and order of convergence results are derived.

Subsection \ref{SbSect_AppSmth_to_DatFunc_Lagran} deals with the convergence
of the Approximate smoother to the data function using the simple inequality
\[
\left\vert f_{d}\left(  x\right)  -s_{a}\left(  x\right)  \right\vert
\leq\left\vert f_{d}\left(  x\right)  -s_{e}\left(  x\right)  \right\vert
+\left\vert s_{e}\left(  x\right)  -s_{a}\left(  x\right)  \right\vert ,
\]

which combines the results of the previous subsection with those derived for
the convergence of the Exact smoother to the data function in Section
\ref{Sect_converg_Exact_smth}.

\subsection{Convergence results not involving
order\label{SbSect_App_to_Exact_no_order}}

The convergence results of this subsection do not involve data functions and
the dependent data is just given as an arbitrary vector. Also, no data
densities are estimated and no orders of convergence are calculated, as will
be done in Subsection \ref{SbSect_App_minus_Exact}. We simply establish
uniform pointwise convergence on a bounded data region. To start with we will
introduce a concept of convergence for finite subsets of $\mathbb{R}^{d} $.

\begin{definition}
\label{Def_ap_XDatk_to_XDat}\textbf{Convergence of finite sets }A sequence of
finite sets $X_{n}=\left\{  x_{n}^{\left(  i\right)  }\right\}  _{i=1}^{N_{n}}
$ is said to converge to the finite set $X=\left\{  x^{\left(  i\right)
}\right\}  _{i=1}^{N}$, denoted $X_{n}\rightarrow X$, if there exists an
integer $K$ such that when $n\geq K$, $N_{n}=N$ and for each $i$, $\left\vert
x_{n}^{\left(  i\right)  }-x^{\left(  i\right)  }\right\vert \rightarrow0$ as
$n\rightarrow\infty$.

When $n\geq K$, $X_{n}$ and $X$ can be regarded as members of $\mathbb{R}%
^{Nd}$ and convergence as convergence in $\mathbb{R}^{Nd}$ under the Euclidean norm.
\end{definition}

\begin{theorem}
\label{Thm_Jsd[sig(Zk)]_to_Jsd(sig(Z))}Suppose $G$ is a basis function and
$s_{e}$ is the Exact smoother generated by the data $\left[  X,y\right]  $.

Suppose $X_{n}^{\prime}$ is a sequence of independent data sets which
converges to $X_{0}$ in the sense of Definition \ref{Def_ap_XDatk_to_XDat} and
that $s_{a}^{\left(  n\right)  }$ is the Approximate smoother generated by
$X_{n}^{\prime}$ and $\left[  X,y\right]  $.

Then the Approximate smoothers satisfy $J_{e}\left[  s_{a}^{\left(  n\right)
}\right]  \rightarrow J_{e}\left[  s_{e}\right]  $ as $n\rightarrow\infty$,
where $J_{e}$ is the functional used to define the Exact smoother.
\end{theorem}

\begin{proof}
We first note that the definition of the convergence of independent data sets
allows us to assume that the $X_{n}^{\prime}$ have the same number of points
as $X$.

Now suppose $X^{\prime}$ is an arbitrary independent data set with the same
number of points as $X_{0}$ and let $s_{a}=s_{a}\left(  X^{\prime}\right)  $
be the corresponding Approximate smoother. If it can be shown that as a
function of $X^{\prime}$, $J_{e}\left[  s_{a}\left(  X^{\prime}\right)
\right]  $ is continuous everywhere this theorem holds since $J_{e}\left[
s_{e}\right]  =J_{e}\left[  s_{a}\left(  X\right)  \right]  $. In fact, by
Theorem \ref{Thm_Appr_smth_mat_eqn_2}
\[
s_{a}\left(  X^{\prime}\right)  \left(  x\right)  =\sum\limits_{i=1}^{N}%
\alpha_{i}G\left(  x-x_{i}^{\prime}\right)  +\sum\limits_{j=1}^{M}\beta
_{j}p_{j}\left(  x\right)  ,
\]

where $\alpha=\left(  \alpha_{i}\right)  $ and $\beta=\left(  \beta
_{j}\right)  $ satisfy the matrix equations
\[
\left(
\begin{array}
[c]{l}%
\alpha\\
\beta\\
\lambda
\end{array}
\right)  =\Psi^{-1}\left(
\begin{array}
[c]{l}%
G_{X,X^{\prime}}^{T}\\
P_{X}^{T}\\
O
\end{array}
\right)  y.
\]

and the Approximate smoothing matrix is given by%
\[
\Psi=\left(
\begin{array}
[c]{lll}%
\left(  2\pi\right)  ^{\frac{d}{2}}N\rho\,G_{X^{\prime},X^{\prime}%
}+G_{X,X^{\prime}}^{T}G_{X,X^{\prime}} & G_{X,X^{\prime}}^{T}P_{X} &
P_{X^{\prime}}\\
P_{X}^{T}G_{X,X^{\prime}} & P_{X}^{T}P_{X} & O_{M}\\
P_{X^{\prime}}^{T} & O_{M} & O_{M}%
\end{array}
\right)  .
\]

Also, from part 4 Corollary \ref{Cor_Appr_smth_prop2} and then \ref{s10} we
have%
\begin{align*}
J_{e}\left[  s_{a}\left(  X^{\prime}\right)  \right]   & =\frac{1}{N}%
y^{T}\left(  y-\left(  s_{a}\left(  X^{\prime}\right)  \right)  _{X}\right) \\
& =\frac{1}{N}y^{T}\left(  y-%
\begin{pmatrix}
G_{X,X^{\prime}} & P_{X}%
\end{pmatrix}
\left(
\begin{array}
[c]{c}%
\alpha\\
\beta
\end{array}
\right)  \right) \\
& =\frac{1}{N}y^{T}\left(  y-%
\begin{pmatrix}
G_{X,X^{\prime}} & P_{X} & O
\end{pmatrix}
\left(
\begin{array}
[c]{c}%
\alpha\\
\beta\\
\lambda
\end{array}
\right)  \right) \\
& =\frac{1}{N}y^{T}\left(  I_{N}-%
\begin{pmatrix}
G_{X,X^{\prime}} & P_{X} & O
\end{pmatrix}
\Psi^{-1}\left(
\begin{array}
[c]{l}%
G_{X,X^{\prime}}^{T}\\
P_{X}^{T}\\
O
\end{array}
\right)  \right)  y.
\end{align*}

If we can show that $G_{X,X^{\prime}}$ and $\Psi_{X^{\prime}}^{-1}$ are
continuous functions of $X^{\prime}$ in a neighborhood of $X_{0}$ then
$J_{e}\left[  s_{a}\left(  X^{\prime}\right)  \right]  $ is a continuous
function of $X^{\prime}$. But $G_{X^{\prime},X^{\prime}}$ and $G_{X,X^{\prime
}}$ are clearly continuous for all $X^{\prime}$ so $\Psi_{X^{\prime}}$ and
$\det\Psi_{X^{\prime}}$ are continuous for all $X^{\prime}$. Further, since by
Theorem \ref{Thm_mat_psi}, $\Psi_{X^{\prime}}$ is positive definite and
regular for all $X^{\prime}$, $\det\Psi_{X^{\prime}}>0$ for all $X^{\prime}$
and it is clear from Cramer's rule that $\Psi_{X^{\prime}}^{-1}$ is continuous
everywhere. Thus $J_{e}\left[  s_{a}\left(  X^{\prime}\right)  \right]  $ is
continuous everywhere and the proof is complete.
\end{proof}

The next corollary shows that for given data the Approximate smoother
converges uniformly to the Exact smoother on $\mathbb{R}^{d}$ for the case
that $X^{\prime}$ is a grid containing a data region and the grid size goes to zero.

\begin{corollary}
\label{Cor_Thm_Jsd[sig(Zk)]toJsd(sig(Z))_1}\ 

\begin{enumerate}
\item Suppose $s_{e}$ is the Exact smoother generated by the data $\left[
X,y\right]  $ and that we have a sequence of regular, rectangular grids with
grid nodes $X_{n}^{\prime}$, grid sizes $h_{n}^{\prime}$ and a fixed common
grid region $\mathcal{R}$ - a closed rectangle - containing $X$.

\item Let $s_{a}\left(  X_{n}^{\prime}\right)  $ denote the Approximate
smoother generated by $X_{n}^{\prime}$ and data $\left[  X,y\right]  $.

\item Assume that $X\subset\Omega\subseteq\mathcal{R}$ for all $n$, where the
data region $\Omega$ is a bounded, open, connected subset of $\mathbb{R}^{d}$
having the cone property.

\item Assume that the conditions of Theorem \ref{Thm_ap_rx(x)_bound_better} hold.
\end{enumerate}

Then $\left\vert h_{n}^{\prime}\right\vert \rightarrow0$ implies $\left\Vert
s_{a}\left(  X_{n}^{\prime}\right)  -s_{e}\right\Vert _{w,\theta}\rightarrow0$
and also that $s_{a}\left(  X_{n}^{\prime}\right)  $ converges uniformly
pointwise to $s_{e}$ on $\overline{\Omega}$.
\end{corollary}

\begin{proof}
For each data point $x^{\left(  k\right)  }\in X$ there exists a sequence of
distinct points $\left(  z_{n}^{\left(  k\right)  }\right)  _{n=1}^{\infty}$
such that $z_{n}^{\left(  k\right)  }\in X_{n}^{\prime}$ and $z_{n}^{\left(
k\right)  }\rightarrow x^{\left(  k\right)  }$ in $\mathbb{R}^{d}$ as
$n\rightarrow\infty$.

Set $Z_{n}=\left\{  z_{n}^{\left(  k\right)  }\right\}  _{k=1}^{N}$ and let
$s\left(  Z_{n}\right)  $ be the Approximate smoother generated by $Z_{n}$.
Then $Z_{n}\rightarrow X$ as independent data and by Theorem
\ref{Thm_Jsd[sig(Zk)]_to_Jsd(sig(Z))}, $J_{e}\left[  s_{a}\left(
Z_{n}\right)  \right]  \rightarrow J_{e}\left[  s_{e}\right]  $.

But since $Z_{n}\subset X_{n}^{\prime}$ we have $J_{e}\left[  s_{a}\left(
X_{n}^{\prime}\right)  \right]  \leq J_{e}\left[  s_{a}\left(  Z_{n}\right)
\right]  $.

Also from the definition of $s_{e}$ we have $J_{e}\left[  s_{e}\right]  \leq
J_{e}\left[  s_{a}\left(  X_{n}^{\prime}\right)  \right]  $.

Thus
\[
J_{e}\left[  s_{e}\right]  \leq J_{e}\left[  s_{a}\left(  X_{n}^{\prime
}\right)  \right]  \leq J_{e}\left[  s_{a}\left(  Z_{n}\right)  \right]  ,
\]

and so $J_{e}\left[  s_{a}\left(  X_{n}^{\prime}\right)  \right]  \rightarrow
J_{e}\left[  s_{e}\right]  $.

To establish the convergence of $s_{a}\left(  X_{n}^{\prime}\right)  $ we use
part 2 of Corollary \ref{Cor_propert_Exact_smth}, namely that%
\[
J_{e}\left[  s_{e}\right]  +\rho\left\vert s_{e}-f\right\vert _{w,\theta}%
^{2}+\frac{1}{N}\sum_{k=1}^{N}\left\vert s_{e}\left(  x^{(k)}\right)
-f\left(  x^{(k)}\right)  \right\vert ^{2}=J_{e}\left[  f\right]  ,
\]

for all $f\in X_{w}^{\theta}$. Clearly $f=s_{a}\left(  X_{n}^{\prime}\right)
$ implies that as $n\rightarrow\infty$, $\left\vert s_{e}-s_{a}\left(
X_{n}^{\prime}\right)  \right\vert _{w,\theta}\rightarrow0$ and $s_{a}\left(
X_{n}^{\prime}\right)  \left(  x^{(k)}\right)  \rightarrow s_{e}\left(
x^{(k)}\right)  $ for each $k$, and so
\[
\left\Vert s_{e}-s_{a}\left(  X_{n}^{\prime}\right)  \right\Vert _{w,\theta
}^{2}=\left\vert s_{e}-s_{a}\left(  X_{n}^{\prime}\right)  \right\vert
_{w,\theta}^{2}+\sum_{m=1}^{M}\left\vert s_{e}\left(  a_{m}\right)
-s_{a}\left(  X_{n}^{\prime}\right)  \left(  a_{m}\right)  \right\vert
^{2}\rightarrow0,
\]

where the minimally unisolvent set of points $A=\left\{  a_{m}\right\}
_{m=1}^{M}\subset X$ define the Light norm $\left\Vert \cdot\right\Vert
_{w,\theta}$. Finally, if $R_{x}$ is the Riesz representer of the functional
$f\rightarrow f\left(  x\right)  $ defined by $A$ then for $x\in\Omega$
\begin{align*}
\left\vert s_{e}\left(  x\right)  -s_{a}\left(  X_{n}^{\prime}\right)  \left(
x\right)  \right\vert  & =\left\vert \left(  s_{e}-s_{a}\left(  X_{n}^{\prime
}\right)  ,R_{x}\right)  _{w,\theta}\right\vert \\
& \leq\left\Vert s_{e}-s_{a}\left(  X_{n}^{\prime}\right)  \right\Vert
_{w,\theta}\left\Vert R_{x}\right\Vert _{w,\theta}\\
& =\left\Vert s_{e}-s_{a}\left(  X_{n}^{\prime}\right)  \right\Vert
_{w,\theta}\sqrt{\left\vert r_{x}\right\vert _{w,\theta}^{2}+\left\vert
\widetilde{l}\left(  x\right)  \right\vert ^{2}}\\
& \leq\left\Vert s_{e}-s_{a}\left(  X_{n}^{\prime}\right)  \right\Vert
_{w,\theta}\left(  \left\vert r_{x}\right\vert _{w,\theta}+\left\vert
\widetilde{l}\left(  x\right)  \right\vert _{1}\right) \\
& =\left\Vert s_{e}-s_{a}\left(  X_{n}^{\prime}\right)  \right\Vert
_{w,\theta}\left(  \sqrt{r_{x}\left(  x\right)  }+\sum\limits_{j=1}%
^{M}\left\vert l_{j}\left(  x\right)  \right\vert \right)  .
\end{align*}

Our assumptions now allow us to choose $A$ such that we can estimate
$\sqrt{r_{x}\left(  x\right)  }$ and $\sum\limits_{j=1}^{M}\left\vert
l_{j}\left(  x\right)  \right\vert $ using \ref{s03} and \ref{s07}, and so
obtain%
\[
\left\vert s_{e}\left(  x\right)  -s_{a}\left(  X_{n}^{\prime}\right)  \left(
x\right)  \right\vert \leq\sqrt{c_{G,\eta,\sigma}}\left(  1+K_{\Omega,\theta
}^{\prime}\right)  \left(  \operatorname*{diam}\Omega\right)  ^{\eta
+\delta_{G}}\left\Vert s_{e}-s_{a}\left(  X_{n}^{\prime}\right)  \right\Vert
_{w,\theta},\quad x\in\Omega,
\]

where $c_{G,\eta,\sigma}$ and $K_{\Omega,\theta}^{\prime}$ are independent of
$x\in\Omega$. Finally, since $s_{e}$ and $s_{a}$ are continuous on
$\mathbb{R}^{d}$ the estimate is actually valid on $\overline{\Omega}$. Hence
$s_{a}\left(  X_{n}^{\prime}\right)  \rightarrow s_{e}$ uniformly, pointwise
on $\overline{\Omega}$ as $n\rightarrow\infty$.
\end{proof}

The next corollary replaces the rectangular grids by scattered sets of points
and thus it can be applied to sparse grids, for example. This corollary shows
that for given data the Approximate smoother converges to the Exact smoother
on $\mathbb{R}^{d}$ if the scattered sets converge to the independent data in
the sense of Definition \ref{Def_ap_XDatk_to_XDat}.

\begin{corollary}
\label{Cor_Thm_Jsd[sig(Zk)]toJsd(sig(Z))_2}\ 

\begin{enumerate}
\item Suppose $s_{e}$ is the Exact smoother generated by the data $\left[
X,y\right]  $ and that $X\subset\Omega$ where the data region $\Omega$ is a
bounded, open, connected subset of $\mathbb{R}^{d}$ having the cone property.

\item Suppose that $X_{n}^{\prime}$ is a sequence of finite sets and that
there exist a sequence of data sets $X_{n}^{\prime\prime}$ such that
$X_{n}^{\prime\prime}\subset X_{n}^{\prime}$ and $X_{n}^{\prime\prime
}\rightarrow X$ in the sense of Definition \ref{Def_ap_XDatk_to_XDat}.

\item Assume that the conditions of Theorem \ref{Thm_ap_rx(x)_bound_better} hold.
\end{enumerate}

Then, if $s_{a}\left(  X_{n}^{\prime\prime}\right)  $ is the Approximate
smoother generated by $\left[  X,y\right]  $ and $X_{n}^{\prime\prime}$, it
follows that

$\left\Vert s_{a}\left(  X_{n}^{\prime\prime}\right)  -s_{e}\right\Vert
_{w,\theta}\rightarrow0$ and $s_{a}\left(  X_{n}^{\prime\prime}\right)  $
converges uniformly pointwise to $s_{e}$ on $\overline{\Omega}$.
\end{corollary}

\begin{proof}
A simple modification of the proof of the previous Corollary
\ref{Cor_Thm_Jsd[sig(Zk)]toJsd(sig(Z))_1}.
\end{proof}

\subsection{Convergence to the Exact smoother\label{SbSect_App_minus_Exact}}

In Subsection \ref{SbSect_ex_Exact_smth_error} we proved results concerning
the pointwise convergence of the Exact smoother to its data function using
results from Lagrange interpolation theory and the `cavity' measure of
independent data set density $h_{X}=\sup\limits_{\omega\in\Omega
}\operatorname*{dist}\left(  \omega;X\right)  $. In this subsection we will
employ the same techniques to derive conditions under which the Approximate
smoother converges pointwise to the Exact smoother - conditions expressed in
terms of the smoothing coefficient $\rho$ and the densities of the data sets
$X$ and $X^{\prime}$.

We will require the following lemma which is derived from equation \ref{h58}
and Theorem \ref{Thm_Rvx_solution} of Chapter \ref{Ch_ExactSmth}.

\begin{lemma}
\label{Lem_app_Rvx}Suppose $X=\left\{  x^{\left(  i\right)  }\right\}
_{i=1}^{N} $ is a $\theta$-unisolvent set and let $A=\left\{  a^{\left(
i\right)  }\right\}  _{i=1}^{M}$ be a minimally $\theta$-unisolvent subset.
Construct $R_{x}$, $\mathcal{P}$ and $\mathcal{Q}=I-\mathcal{P}$ using $A$.

Then for each $x\in\mathbb{R}^{d}$ there exists a unique element $R_{V,x}\in
V$ such that%
\[
f\left(  x\right)  =\left(  \mathcal{L}_{X}f,R_{V,x}\right)  _{V},\quad\text{
}f\in X_{w}^{\theta},
\]

and%
\begin{equation}
\rho\left\Vert R_{V,x}\right\Vert _{V}^{2}\leq r_{x}\left(  x\right)
+N_{X}\rho\left\vert \widetilde{l}\left(  x\right)  \right\vert ^{2}%
.\label{s14}%
\end{equation}

\end{lemma}

We will also need the following lemma which supplies some results from the
theory of Lagrange interpolation. This lemma has been created from Lemma 3.2,
Lemma 3.5 and the first two paragraphs of the proof of Theorem 3.6 of Light
and Wayne \cite{LightWayne98PowFunc}. The results of this lemma do not involve
any reference to weight or basis functions or to functions in $X_{w}^{\theta}%
$, but consider the properties of the region $\Omega$ which contains the
independent data points $X$ and the order of the unisolvency used for the
interpolation. Thus we have separated the part of the proof that involves
basis functions from the part that uses the detailed theory of Lagrange
interpolation operators.

\begin{lemma}
\label{Lem_app_Lagrange_interpol}(A copy of Lemma
\ref{Lem_int_Lagrange_interpol}) Assume first that:

\begin{enumerate}
\item $\Omega$ is a bounded, open, connected subset of $\mathbb{R}^{d}$ having
the cone property.

\item $X$ is a unisolvent subset of $\Omega$ of order $\theta$.
\end{enumerate}

Suppose $\left\{  l_{j}\right\}  _{j=1}^{M}$ is the cardinal basis of
$P_{\theta-1}$ with respect to a minimal unisolvent set of $\Omega$. Using
Lagrange interpolation techniques, it can be shown there exists a constant
$K_{\Omega,\theta}^{\prime}>0$ such that
\begin{equation}
\sum\limits_{j=1}^{M}\left\vert l_{j}\left(  x\right)  \right\vert \leq
K_{\Omega,\theta}^{\prime},\quad x\in\overline{\Omega},\label{s03}%
\end{equation}

and all minimal unisolvent subsets of $\Omega$. Now define
\[
h_{X}=\sup\limits_{x\in\Omega}\operatorname*{dist}\left(  x,X\right)  ,
\]

and fix $x\in X$. Using Lagrange interpolation techniques it can be shown
there are constants $c_{\Omega,\theta},h_{\Omega,\theta}>0$ such that when
$h_{X}<h_{\Omega,\theta}$ there exists a minimal unisolvent set $A\subset X$
satisfying
\begin{equation}
\operatorname*{diam}A_{x}\leq c_{\Omega,\theta}h_{X},\label{s05}%
\end{equation}

where $A_{x}=A\cup\left\{  x\right\}  $.
\end{lemma}

The next result supplies conditions on the basis function which may yield a
higher order estimate for $\sqrt{r_{x}\left(  x\right)  }$ than Theorem
\ref{Thm_int_rx(x)_bound}.

\begin{theorem}
\label{Thm_ap_rx(x)_bound_better}(Copy of Theorem
\ref{Thm_int_rx(x)_bnd_better_order}) Suppose $w$ is a weight function with
properties W2 and W3 for order $\theta$ and $\kappa$. Set $\eta=\min\left\{
\theta,\frac{1}{2}\left\lfloor 2\kappa\right\rfloor \right\}  $. Also suppose
$G$ is a basis function of order $\theta$ such that the distributions
$\left\{  D^{\beta}G\right\}  _{\left\vert \beta\right\vert =2\eta+1}$ are
$L_{loc}^{1}$ functions such that for each fixed $b\neq0$ the integrals%
\begin{equation}
\int_{0}^{1}\left(  1-t\right)  ^{2\eta}\left\vert \left(  D^{\beta}G\right)
\left(  z+tb\right)  \right\vert dt,\quad z,b\in\mathbb{R}^{d},\text{
}\left\vert \beta\right\vert =2\eta+1,\label{s01}%
\end{equation}

have polynomial growth in $z$.

Further, there exist constants $r_{G},c_{G,\eta}>0$ and $\delta_{G}\geq0$ such
that%
\begin{equation}
\left\vert b^{\beta}\right\vert \int_{0}^{1}\left(  1-t\right)  ^{2\eta
}\left\vert \left(  D^{\beta}G\right)  \left(  tb\right)  \right\vert
dt\leq\frac{c_{G,\eta}}{2\sigma}\left\vert b\right\vert ^{2\left(  \eta
+\delta_{G}\right)  },\quad\left\vert b\right\vert \leq r_{G},\text{
}\left\vert \beta\right\vert =2\eta+1.\label{s02}%
\end{equation}

Regarding unisolvency, assume $A=\left\{  a^{\left(  k\right)  }\right\}
_{k=1}^{M}$ is a minimal $\theta$-unisolvent set and that $\left\{
l_{k}\right\}  _{k=1}^{M}$ is the corresponding unique cardinal basis for
$P_{\theta-1}$. Now construct $\mathcal{P},\mathcal{Q},R_{x}$ using $A$ and
$\left\{  l_{k}\right\}  _{k=1}^{M}$.

Now if $r_{x}=\mathcal{Q}R_{x}$ we have the estimate%
\begin{equation}
\sqrt{r_{x}\left(  x\right)  }\leq\sqrt{c_{G,\eta,\sigma}}\left(  1+\sum
_{k=1}^{M}\left\vert l_{k}\left(  x\right)  \right\vert \right)  \left(
\operatorname*{diam}A_{x}\right)  ^{\eta+\delta_{G}},\quad\operatorname*{diam}%
A_{x}\leq r_{G},\text{ }x\in\Omega,\label{s07}%
\end{equation}

where $\sigma=\min\left\{  \theta,\frac{1}{2}\left\lfloor \min2\kappa
+1\right\rfloor \right\}  $, $A_{x}=A\cup\left\{  x\right\}  $ and
$c_{G,\eta,\sigma}=\frac{d^{\left\lceil \sigma\right\rceil }}{\left(
2\pi\right)  ^{d/2}\left\lceil \sigma\right\rceil !}c_{G,\eta}$.
\end{theorem}

\begin{remark}
In Section \ref{Sect_better_results}, $\eta$ and $\delta_{G}$ were calculated
for the thin-plate and shifted thin-plate splines by way of examples.
\end{remark}

We now derive our estimates for the pointwise convergence of the Approximate
smoother to the Exact smoother.

\begin{theorem}
\label{Thm_ap_Se(x)minusSa(x)_bound_RVx}We will need the following assumptions
and notation:

\begin{enumerate}
\item[(a)] Suppose $w$ is a weight function with properties W2 and W3 for
order $\theta$ and $\kappa$, and set $\eta=\min\left\{  \theta,\frac{1}%
{2}\left\lfloor \min2\kappa\right\rfloor \right\}  $. Assume $G$ is a basis
function of order $\theta$ such that there exist constants $c_{G},r_{G}>0$ and
$\delta_{G}\geq0$ satisfying \ref{s07}. Set $\eta_{G}=\eta+\delta_{G}$.

\item[(b)] Denote by $s_{e}$ the Exact smoother generated by the the smoothing
parameter $\rho$, the independent data $X$ and data function $f_{d}\in
X_{w}^{\theta}$.

\item[(c)] Assume $X\subset\Omega$ where the data region $\Omega$ has the
properties given in the Lagrange interpolation Lemma
\ref{Lem_app_Lagrange_interpol} and the constants $h_{\Omega,\theta}%
,c_{\Omega,\theta},K_{\Omega,\theta}^{\prime}$ are as given in Lemma
\ref{Lem_app_Lagrange_interpol}. Regarding Theorem
\ref{Thm_ap_rx(x)_bound_better} set%
\begin{equation}
c_{G,\eta,\sigma}=\frac{d^{\left\lceil \sigma\right\rceil }}{\left(
2\pi\right)  ^{d/2}\left\lceil \sigma\right\rceil !}c_{G,\eta}.\label{s64}%
\end{equation}

\item[(d)] Denote by $s_{a}$ the Approximate smoother generated by the
unisolvent set of points $X^{\prime}\subset\Omega$. Let $h_{X}=\sup
\limits_{\omega\in\Omega}\operatorname*{dist}\left(  \omega,X\right)  $,
$h_{X^{\prime}}=\sup\limits_{\omega\in\Omega}\operatorname*{dist}\left(
\omega,X^{\prime}\right)  $ measure the density of the point sets $X$ and
$X^{\prime}$.
\end{enumerate}

Then when $h_{X},h_{X^{\prime}}\leq\min\left\{  h_{\Omega,\theta}%
,r_{G}\right\}  $ we have for $x\in\overline{\Omega}$:

\begin{enumerate}
\item If $f_{d}\in W_{G,X^{\prime}}$%
\[
\left\vert s_{e}\left(  x\right)  -s_{a}\left(  x\right)  \right\vert
\leq\left\vert f_{d}\right\vert _{w,\theta}\left(  1+K_{\Omega,\theta}%
^{\prime}\right)  \left(  \sqrt{c_{G,\eta,\sigma}}\left(  c_{\Omega,\theta
}h_{X}\right)  ^{\eta_{G}}+\sqrt{N_{X}\rho}\right)  .
\]

\item If $f_{d}\in X_{w}^{\theta}$%
\begin{align}
\left\vert s_{e}\left(  x\right)  -s_{a}\left(  x\right)  \right\vert  &
\leq\left\vert f_{d}\right\vert _{w,\theta}\left(  1+\left(  1+K_{\Omega
,\theta}^{\prime}\right)  \sqrt{c_{G,\eta,\sigma}}\frac{\left(  c_{\Omega
,\theta}h_{X^{\prime}}\right)  ^{\eta_{G}}}{\sqrt{\rho}}\right)
\times\nonumber\\
& \qquad\times\left(  1+K_{\Omega,\theta}^{\prime}\right)  \left(
\sqrt{c_{G,\eta,\sigma}}\left(  c_{\Omega,\theta}h_{X}\right)  ^{\eta_{G}%
}+\sqrt{N_{X}\rho}\right)  .\label{s46}%
\end{align}

??? \textbf{Note}: The term $\sqrt{\rho}$ in the denominator does not appear
in the zero order case - suspicious.\medskip

Two similar bounds on the error for large $\rho$ are:\medskip

\item If $f_{d}\in X_{w}^{\theta}$%
\[
\left\vert s_{e}\left(  x\right)  -s_{a}\left(  x\right)  \right\vert
\leq\left\vert f_{d}\right\vert _{w,\theta}\sqrt{c_{G,\eta,\sigma}}\left(
1+K_{\Omega,\theta}^{\prime}\right)  ^{2}\left(  \operatorname*{diam}%
\Omega\right)  ^{\eta_{G}}\left(  \sqrt{c_{G,\eta,\sigma}}\frac{\left(
c_{\Omega,\theta}h_{X}\right)  ^{\eta_{G}}}{\sqrt{\rho}}+\sqrt{N_{X}}\right)
.
\]

\item If $f_{d}\in X_{w}^{\theta}$
\[
\left\vert s_{e}\left(  x\right)  -s_{a}\left(  x\right)  \right\vert
\leq\left(  \max\limits_{\omega\in\overline{\Omega}}\left\vert f_{d}\left(
\omega\right)  \right\vert \right)  \left(  1+K_{\Omega,\theta}^{\prime
}\right)  \left(  \sqrt{c_{G,\eta,\sigma}}\frac{\left(  c_{\Omega,\theta}%
h_{X}\right)  ^{\eta_{G}}}{\sqrt{\rho}}+\sqrt{N_{X}}\right)  .
\]

\item If $f_{d}\in P_{\theta-1}$ then $f_{d}=s_{e}=s_{a}$.
\end{enumerate}
\end{theorem}

\begin{proof}
Fix $x\in\Omega$ and construct $r_{x}=\mathcal{Q}R_{x}$, $s_{e}$ and $s_{a} $
from a minimal $\theta$-unisolvent set $A\subset X$. From Lemma
\ref{Lem_app_Rvx}
\[
\left\vert s_{e}\left(  x\right)  -s_{a}\left(  x\right)  \right\vert
=\left\vert \left(  \mathcal{L}_{X}\left(  s_{e}-s_{a}\right)  ,R_{V,x}%
\right)  _{V}\right\vert \leq\left\Vert \mathcal{L}_{X}\left(  s_{e}%
-s_{a}\right)  \right\Vert _{V}\left\Vert R_{V,x}\right\Vert _{V}.
\]

But from part 5 Summary \ref{Sum_ap_property_L_Exact_smth}, $\left\Vert
\mathcal{L}_{X}\left(  s_{e}-f\right)  \right\Vert _{V}\leq\left\Vert
\mathcal{L}_{X}f-\varsigma\right\Vert _{V}$ for all $f\in W_{G,X^{\prime}}$.
Choose $f=s_{a}$ so that by part 1 Summary \ref{Sum_ap_property_L_Exact_smth}%
\[
\left\Vert \mathcal{L}_{X}\left(  s_{e}-s_{a}\right)  \right\Vert _{V}^{2}%
\leq\left\Vert \mathcal{L}_{X}s_{a}-\varsigma\right\Vert _{V}^{2}=J_{e}\left[
s_{a}\right]  ,
\]

and hence%
\[
\left\vert s_{e}\left(  x\right)  -s_{a}\left(  x\right)  \right\vert
\leq\sqrt{J_{e}\left[  s_{a}\right]  }\left\Vert R_{V,x}\right\Vert _{V}.
\]

But by \ref{s14}%
\[
\rho\left\Vert R_{V,x}\right\Vert _{V}^{2}\leq r_{x}\left(  x\right)
+N_{X}\rho\left\vert \widetilde{l}\left(  x\right)  \right\vert ^{2},
\]

so that%
\[
\left\vert s_{e}\left(  x\right)  -s_{a}\left(  x\right)  \right\vert
\leq\sqrt{J_{e}\left[  s_{a}\right]  }\frac{1}{\sqrt{\rho}}\left(  \sqrt
{r_{x}\left(  x\right)  }+\sqrt{N_{X}\rho}\left\vert \widetilde{l}\left(
x\right)  \right\vert \right)  ,\quad x\in\mathbb{R}^{d}.
\]

By Lemma \ref{Lem_app_Lagrange_interpol} there exists a constant
$K_{\Omega,\theta}^{\prime}$, independent of $x\in\Omega$, such that
$\left\vert \widetilde{l}\left(  x\right)  \right\vert \leq\sum\limits_{k=1}%
^{M}\left\vert l_{k}\left(  x\right)  \right\vert \leq K_{\Omega,\theta
}^{\prime}$ so that%
\[
\left\vert s_{e}\left(  x\right)  -s_{a}\left(  x\right)  \right\vert
\leq\sqrt{J_{e}\left[  s_{a}\right]  }\frac{1}{\sqrt{\rho}}\left(  \sqrt
{r_{x}\left(  x\right)  }+K_{\Omega,\theta}^{\prime}\sqrt{N_{X}\rho}\right)
,\quad x\in\Omega,
\]

and the estimate \ref{s07} for $\sqrt{r_{x}\left(  x\right)  }$ on $\Omega$
implies%
\[
\sqrt{r_{x}\left(  x\right)  }\leq\sqrt{c_{G,\eta,\sigma}}\left(  1+\sum
_{k=1}^{M}\left\vert l_{k}\left(  x\right)  \right\vert \right)  \left(
\operatorname*{diam}A_{x}\right)  ^{\eta_{G}}\leq\sqrt{c_{G,\eta,\sigma}%
}\left(  1+K_{\Omega,\theta}^{\prime}\right)  \left(  \operatorname*{diam}%
A_{x}\right)  ^{\eta_{G}},
\]

provided $\operatorname*{diam}A_{x}\leq r_{G}$. But an assumption of this
theorem is $h_{X}\leq\min\left\{  h_{\Omega,\theta},r_{G}\right\}  $ so by
Lemma \ref{Lem_app_Lagrange_interpol}, $A$ can be chosen so that
$\operatorname*{diam}A_{x}\leq c_{\Omega,\theta}h_{X}$ so that
\begin{equation}
\sqrt{r_{x}\left(  x\right)  }\leq\sqrt{c_{G,\eta,\sigma}}\left(
1+K_{\Omega,\theta}^{\prime}\right)  \left(  c_{\Omega,\theta}h_{X}\right)
^{\eta_{G}},\label{s08}%
\end{equation}

and%
\begin{align*}
\left\vert s_{e}\left(  x\right)  -s_{a}\left(  x\right)  \right\vert  &
\leq\frac{\sqrt{J_{e}\left[  s_{a}\right]  }}{\sqrt{\rho}}\left(
1+K_{\Omega,\theta}^{\prime}\right)  \left(  \sqrt{c_{G,\eta,\sigma}}\left(
c_{\Omega,\theta}h_{X}\right)  ^{\eta_{G}}+\sqrt{N_{X}\rho}\right) \\
& =\left(  1+K_{\Omega,\theta}^{\prime}\right)  \frac{\sqrt{J_{e}\left[
s_{a}\right]  }}{\sqrt{\rho}}\left(  \sqrt{c_{G,\eta,\sigma}}\left(
c_{\Omega,\theta}h_{X}\right)  ^{\eta_{G}}+\sqrt{N_{X}\rho}\right)  .
\end{align*}

Thus when $f\in W_{G,X^{\prime}}$ and $f\neq s_{a}$, the definition of $s_{a}$
implies $J_{e}\left[  s_{a}\right]  <J_{e}\left[  f\right]  $ and%
\begin{equation}
\left\vert s_{e}\left(  x\right)  -s_{a}\left(  x\right)  \right\vert
\leq\left(  1+K_{\Omega,\theta}^{\prime}\right)  \sqrt{J_{e}\left[  f\right]
}\left(  \sqrt{c_{G,\eta,\sigma}}\frac{\left(  c_{\Omega,\theta}h_{X}\right)
^{\eta_{G}}}{\sqrt{\rho}}+\sqrt{N_{X}}\right)  .\label{s04}%
\end{equation}

The estimates for $\left\vert s_{e}\left(  x\right)  -s_{a}\left(  x\right)
\right\vert $ stated in the theorem will be proved by substituting various
$f\in W_{G,X^{\prime}}$ into \ref{s04}.\medskip

\textbf{Part 1} We choose $f=f_{d}\in W_{G,X^{\prime}}$. Then $J_{e}\left[
f\right]  =\rho\left\vert f_{d}\right\vert _{w,\theta}^{2}$ and \ref{s04}
implies
\[
\left\vert s_{e}\left(  x\right)  -s_{a}\left(  x\right)  \right\vert
\leq\left\vert f_{d}\right\vert _{w,\theta}\left(  1+K_{\Omega,\theta}%
^{\prime}\right)  \left(  \sqrt{c_{G,\eta,\sigma}}\left(  c_{\Omega,\theta
}h_{X}\right)  ^{\eta_{G}}+\sqrt{N_{X}\rho}\right)  ,\quad x\in\Omega,
\]

the extension to $\overline{\Omega}$ being valid since both $s_{e}$ and
$s_{a}$ are continuous on $\mathbb{R}^{d}$.\medskip

\textbf{Part 2} Choose $f=\mathcal{I}_{X^{\prime}}^{\prime}f_{d}$ where
$\mathcal{I}_{X^{\prime}}^{\prime}f_{d}$ is the minimal seminorm interpolant
of $f_{d}$ on the unisolvent set $X^{\prime}$. Then by Corollary
\ref{Cor_Thm_int_rx(x)_bnd_better_order}, when $h_{X^{\prime}}<\min\left\{
h_{\Omega,\theta},r_{G}\right\}  $ it follows that
\[
\left\vert f_{d}\left(  x\right)  -\mathcal{I}_{X^{\prime}}^{\prime}%
f_{d}\left(  x\right)  \right\vert \leq\left\vert f_{d}-\mathcal{I}%
_{X^{\prime}}^{\prime}f_{d}\right\vert _{w,\theta}\left(  1+K_{\Omega,\theta
}^{\prime}\right)  \sqrt{c_{G,\eta,\sigma}}\left(  c_{\Omega,\theta
}h_{X^{\prime}}\right)  ^{\eta_{G}},\quad x\in\overline{\Omega},
\]

and so, by using the properties: $\left\vert \mathcal{I}_{X^{\prime}}^{\prime
}f\right\vert _{w,\theta}\leq\left\vert f\right\vert _{w,\theta}$ and
$\left\vert \left(  I-\mathcal{I}_{X^{\prime}}^{\prime}\right)  f\right\vert
_{w,\theta}\leq\left\vert f\right\vert _{w,\theta}$, of part 2 Theorem
\ref{Thm_propert_interpol_map}%
\begin{align*}
J_{e}\left[  f\right]   & =J_{e}\left[  \mathcal{I}_{X^{\prime}}^{\prime}%
f_{d}\right] \\
& =\rho\left\vert \mathcal{I}_{X^{\prime}}^{\prime}f_{d}\right\vert
_{w,\theta}^{2}+\frac{1}{N}\sum\limits_{k=1}^{N}\left\vert \left(
\mathcal{I}_{X^{\prime}}^{\prime}f_{d}\right)  \left(  x^{\left(  k\right)
}\right)  -f_{d}\left(  x^{\left(  k\right)  }\right)  \right\vert ^{2}\\
& \leq\rho\left\vert \mathcal{I}_{X^{\prime}}^{\prime}f_{d}\right\vert
_{w,\theta}^{2}+\left\vert f_{d}-\mathcal{I}_{X^{\prime}}^{\prime}%
f_{d}\right\vert _{w,\theta}^{2}\left(  1+K_{\Omega,\theta}^{\prime}\right)
^{2}c_{G,\eta,\sigma}\left(  c_{\Omega,\theta}h_{X^{\prime}}\right)
^{2\eta_{G}}\\
& \leq\rho\left\vert f_{d}\right\vert _{w,\theta}^{2}+\left\vert
f_{d}\right\vert _{w,\theta}^{2}\left(  1+K_{\Omega,\theta}^{\prime}\right)
^{2}c_{G,\eta,\sigma}\left(  c_{\Omega,\theta}h_{X^{\prime}}\right)
^{2\eta_{G}}\\
& =\left\vert f_{d}\right\vert _{w,\theta}^{2}\left(  \rho+\left(
1+K_{\Omega,\theta}^{\prime}\right)  ^{2}c_{G,\eta,\sigma}\left(
c_{\Omega,\theta}h_{X^{\prime}}\right)  ^{2\eta_{G}}\right)  ,
\end{align*}

i.e.%
\[
\sqrt{J_{e}\left[  f\right]  }\leq\left\vert f_{d}\right\vert _{w,\theta
}\left(  \sqrt{\rho}+\left(  1+K_{\Omega,\theta}^{\prime}\right)
\sqrt{c_{G,\eta,\sigma}}\left(  c_{\Omega,\theta}h_{X^{\prime}}\right)
^{\eta_{G}}\right)  .
\]

Hence by \ref{s04}%
\begin{align*}
&  \left\vert s_{e}\left(  x\right)  -s_{a}\left(  x\right)  \right\vert \\
&  \leq\left(  1+K_{\Omega,\theta}^{\prime}\right)  \sqrt{J_{e}\left[
f\right]  }\left(  \sqrt{c_{G,\eta,\sigma}}\frac{\left(  c_{\Omega,\theta
}h_{X}\right)  ^{\eta_{G}}}{\sqrt{\rho}}+\sqrt{N_{X}}\right) \\
&  \leq\left\vert f_{d}\right\vert _{w,\theta}\left(  \sqrt{\rho}+\left(
1+K_{\Omega,\theta}^{\prime}\right)  \sqrt{c_{G,\eta,\sigma}}\left(
c_{\Omega,\theta}h_{X^{\prime}}\right)  ^{\eta_{G}}\right)  \left(
1+K_{\Omega,\theta}^{\prime}\right)  \left(  \sqrt{c_{G,\eta,\sigma}}%
\frac{\left(  c_{\Omega,\theta}h_{X}\right)  ^{\eta_{G}}}{\sqrt{\rho}}%
+\sqrt{N_{X}}\right) \\
&  =\left\vert f_{d}\right\vert _{w,\theta}\left(  1+\left(  1+K_{\Omega
,\theta}^{\prime}\right)  \sqrt{c_{G,\eta,\sigma}}\frac{\left(  c_{\Omega
,\theta}h_{X^{\prime}}\right)  ^{\eta_{G}}}{\sqrt{\rho}}\right)  \left(
1+K_{\Omega,\theta}^{\prime}\right)  \left(  \sqrt{c_{G,\eta,\sigma}}\left(
c_{\Omega,\theta}h_{X}\right)  ^{\eta_{G}}+\sqrt{N_{X}\rho}\right)  .
\end{align*}
\medskip

\textbf{Part 3} Suppose $\mathcal{P}^{\prime}$, $\mathcal{Q}^{\prime}$ are
$r_{x}^{\prime}$ are defined using a minimal unisolvent subset $A^{\prime
}\subset X^{\prime}$. Choose $f=\mathcal{P}^{\prime}f_{d}$. Then%
\begin{align*}
J_{e}\left[  f\right]  =J_{e}\left[  \mathcal{P}^{\prime}f_{d}\right]   &
=\rho\left\vert \mathcal{P}^{\prime}f_{d}\right\vert _{w,\theta}^{2}+\frac
{1}{N}\sum\limits_{k=1}^{N}\left\vert \left(  \mathcal{P}^{\prime}%
f_{d}\right)  \left(  x^{\left(  k\right)  }\right)  -f_{d}\left(  x^{\left(
k\right)  }\right)  \right\vert ^{2}\\
&  =\frac{1}{N}\sum\limits_{k=1}^{N}\left\vert \left(  \mathcal{Q}^{\prime
}f_{d}\right)  \left(  x^{\left(  k\right)  }\right)  \right\vert ^{2}\\
&  =\frac{1}{N}\sum\limits_{k=1}^{N}\left\vert \left\langle f_{d}%
,r_{x^{\left(  k\right)  }}^{\prime}\right\rangle _{w,\theta}\right\vert
^{2}\\
&  \leq\frac{1}{N}\sum\limits_{k=1}^{N}\left\vert f_{d}\right\vert _{w,\theta
}^{2}\left\vert r_{x^{\left(  k\right)  }}^{\prime}\right\vert _{w,\theta}%
^{2}\\
&  =\left\vert f_{d}\right\vert _{w,\theta}^{2}\frac{1}{N}\sum\limits_{k=1}%
^{N}\left\vert r_{x^{\left(  k\right)  }}^{\prime}\right\vert _{w,\theta}%
^{2}\\
&  \leq\left\vert f_{d}\right\vert _{w,\theta}^{2}\max_{\omega\in
\overline{\Omega}}\left\vert r_{\omega}^{\prime}\right\vert _{w,\theta}^{2}\\
&  =\left\vert f_{d}\right\vert _{w,\theta}^{2}\max_{\omega\in\overline
{\Omega}}r_{\omega}^{\prime}\left(  \omega\right)  .
\end{align*}

From \ref{s07}%
\[
\sqrt{r_{x}\left(  x\right)  }\leq\sqrt{c_{G,\eta,\sigma}}\left(
1+K_{\Omega,\theta}^{\prime}\right)  \left(  \operatorname*{diam}%
\Omega\right)  ^{\eta_{G}},\quad x\in\Omega,
\]

and since this is still true with $A$ replaced by $A^{\prime}$ we have
\[
\sqrt{J_{e}\left[  f\right]  }\leq\left\vert f_{d}\right\vert _{w,\theta}%
\max_{\omega\in\overline{\Omega}}\sqrt{r_{\omega}^{\prime}\left(
\omega\right)  }\leq\left\vert f_{d}\right\vert _{w,\theta}\sqrt
{c_{G,\eta,\sigma}}\left(  1+K_{\Omega,\theta}^{\prime}\right)  \left(
\operatorname*{diam}\Omega\right)  ^{\eta_{G}},
\]

and substitution into \ref{s04}%
\begin{align*}
\left\vert s_{e}\left(  x\right)  -s_{a}\left(  x\right)  \right\vert  &
\leq\left(  1+K_{\Omega,\theta}^{\prime}\right)  \sqrt{J_{e}\left[  f\right]
}\left(  \sqrt{c_{G,\eta,\sigma}}\frac{\left(  c_{\Omega,\theta}h_{X}\right)
^{\eta_{G}}}{\sqrt{\rho}}+\sqrt{N_{X}}\right) \\
& \leq\left\vert f_{d}\right\vert _{w,\theta}\sqrt{c_{G,\eta,\sigma}}\left(
1+K_{\Omega,\theta}^{\prime}\right)  ^{2}\left(  \operatorname*{diam}%
\Omega\right)  ^{\eta_{G}}\left(  \sqrt{c_{G,\eta,\sigma}}\frac{\left(
c_{\Omega,\theta}h_{X}\right)  ^{\eta_{G}}}{\sqrt{\rho}}+\sqrt{N_{X}}\right)
,
\end{align*}

which proves this part for $x\in\Omega$. Continuity on $\overline{\Omega}$ is
valid since $s_{e}$ and $s_{a}$ are continuous on $\mathbb{R}^{d}$.\medskip

\textbf{Part 4} Choose $f=0$ so that $J_{e}\left[  f\right]  =J_{e}\left[
0\right]  \leq\left(  \max\limits_{\omega\in\overline{\Omega}}\left\vert
f_{d}\left(  \omega\right)  \right\vert \right)  ^{2}$ and substitution into
\ref{s04} yields this part for $x\in\Omega$. Continuity on $\overline{\Omega}$
is valid since $s_{e}$ and $s_{a}$ are continuous on $\mathbb{R}^{d}$.
\medskip

\textbf{Part 5} If $f_{d}\in P_{\theta-1}$ then by part 4, $J_{e}\left[
s_{a}\right]  =0$ and so $\left\vert s_{a}\right\vert _{w,\theta}=0$ and
$s_{a}\left(  x^{\left(  k\right)  }\right)  =f_{d}\left(  x^{\left(
k\right)  }\right)  $ for all $x^{\left(  k\right)  }\in X$. Thus $s_{a}\in
P_{\theta-1}$ and so $s_{a}=f_{d}$. But by property 2 of Theorem
\ref{Thm_Lx*Lx_onto_Xw,th}, $\mathcal{S}_{X}^{e}f=f$ iff $f\in P_{\theta-1}$
and consequently $s_{e}=f_{d}$, proving this part.
\end{proof}

\begin{remark}
\label{Rem_Thm_Se(x)minusSa(x)_bound_RVx}\ 

\begin{enumerate}
\item The error formula of part 1 above is the same as the general Exact
smoother error estimate of Theorem \ref{Thm_ap_Exact_smth_converg_at_x}
Chapter \ref{Ch_ExactSmth}. However, here it only applies to the special
finite-dimensional subspace of data functions $W_{G,X^{\prime}}$. Thus a
sequence of independent data points $X$ and smoothing coefficients $\rho$ can
be chosen so that the corresponding sequence of Approximate smoothers
converges uniformly pointwise to the sequence of Exact smoothers,
independently of the chosen $X^{\prime}$.

\item The right-most factor of the estimate derived in part 2 of Theorem
\ref{Thm_ap_Se(x)minusSa(x)_bound_RVx} is the estimate for the Exact smoother
given below in Theorem \ref{Thm_ap_Exact_smth_converg_at_x}. If we choose
$\rho$ such that

$\sqrt{N_{X}\rho}=\sqrt{c_{G,\eta,\sigma}}\left(  c_{\Omega,\theta}%
h_{X}\right)  ^{\eta_{G}}$ i.e.
\begin{equation}
\sqrt{\rho}=\frac{\sqrt{c_{G,\eta,\sigma}}\left(  c_{\Omega,\theta}%
h_{X}\right)  ^{\eta_{G}}}{\sqrt{N_{X}}},\label{s42}%
\end{equation}

and then require that $\left(  1+K_{\Omega,\theta}^{\prime}\right)
^{2}c_{G,\eta,\sigma}\frac{\left(  c_{\Omega,\theta}h_{X^{\prime}}\right)
^{\eta_{G}}}{\sqrt{\rho}}=1$ we get%
\[
\left(  c_{\Omega,\theta}h_{X^{\prime}}\right)  ^{\eta_{G}}=\frac{\sqrt{\rho}%
}{\left(  1+K_{\Omega,\theta}^{\prime}\right)  ^{2}c_{G,\eta,\sigma}}%
=\frac{\frac{\sqrt{c_{G,\eta,\sigma}}\left(  c_{\Omega,\theta}h_{X}\right)
^{\eta_{G}}}{\sqrt{N_{X}}}}{\left(  1+K_{\Omega,\theta}^{\prime}\right)
^{2}c_{G,\eta,\sigma}}=\frac{\left(  c_{\Omega,\theta}h_{X}\right)  ^{\eta
_{G}}}{\left(  1+K_{\Omega,\theta}^{\prime}\right)  ^{2}\sqrt{c_{G,\eta
,\sigma}N_{X}}}.
\]

So if in addition%
\begin{equation}
h_{X^{\prime}}\leq\frac{h_{X}}{\left(  \left(  1+K_{\Omega,\theta}^{\prime
}\right)  ^{2}\sqrt{c_{G,\eta,\sigma}N_{X}}\right)  ^{1/\eta_{G}}},\label{s06}%
\end{equation}

then
\begin{equation}
\left\vert s_{a}\left(  x\right)  -s_{e}\left(  x\right)  \right\vert
\leq4\left\vert f_{d}\right\vert _{w,\theta}\left(  1+K_{\Omega,\theta
}^{\prime}\right)  \sqrt{c_{G,\eta,\sigma}}\left(  c_{\Omega,\theta}%
h_{X}\right)  ^{\eta_{G}},\text{\quad}x\in\overline{\Omega}.\label{s24}%
\end{equation}

Hence if the sequences $X_{k}$ and $X_{k}^{\prime}$ are such that $h_{X_{k}%
}\rightarrow0$ and $h_{X_{k}^{\prime}}\rightarrow0$ and $\rho_{k} $ is
constrained by \ref{s42} then there is a subsequence $X_{I_{k}}^{\prime}$ of
$X_{k}^{\prime}$ such that $h_{X_{I_{k}}^{\prime}}\leq\frac{h_{X_{k}}}{\left(
\left(  1+K_{\Omega,\theta}^{\prime}\right)  ^{2}\sqrt{c_{G,\eta,\sigma
}N_{X_{k}}}\right)  ^{1/\eta_{G}}}$ and the sequence of Approximate smoothers
converges to the sequence of Exact smoothers in the sense implied by \ref{s24}.

If $X^{\prime}$ is a regular grid then $h_{X^{\prime}}$ can be calculated
exactly:%
\begin{equation}
d^{\frac{d}{2}}\operatorname*{vol}\left(  grid\right)  =N_{X^{\prime}}\left(
h_{X^{\prime}}\right)  ^{d}.\label{s81}%
\end{equation}

\item The approach of part 2 can be augmented by assuming a relationship
between $N_{X}$ and $h_{X}$. This was done in Section ?? 4.8 of Williams
\cite{WilliamsZeroOrdSmthV4} for the zero order Exact smoother and in Remark
\ref{Rem_smooth_converg} for the positive order Exact smoother. Several
\textbf{1-dimensional} numerical experiments were run to compare the
convergence of the zero order Exact smoother with the predicted convergence.
1-dimensional \textbf{test data sets} were constructed using a uniform
distribution on the interval $\Omega=\left[  -1.5,1.5\right]  $. Each of 20
data files were exponentially sampled using a multiplier of approximately 1.2
and a maximum of 5000 points, and then $\log_{10}h_{X}$ was plotted against
$\log_{10}N_{X}$ where $N_{X}=\left\vert X\right\vert $. It then seemed quite
reasonable to use a least-squares linear fit and in this case we obtained the
relation%
\begin{equation}
h_{X}\simeq3.09N^{-0.81}.\label{s84}%
\end{equation}

For ease of calculation let
\begin{equation}
h_{X}=h_{1}\left(  N_{X}\right)  ^{-a},\text{\quad}h_{1}=3.09,\text{\quad
}a=0.81.\label{s85}%
\end{equation}

A barrier to the use of such a formula as \ref{s84} in higher dimensions is
the difficulty of calculating $h_{X}$ for such data sets. If a sequence
of\textbf{\ }independent test data sets was generated by a uniform
distribution in each dimension then the constants $a$ and $h_{1}$ might be
defined as the upper bound of the confidence interval of a statistical
distribution. Also, noting the regular grid formula \ref{s81}, we might
hypothesize a relationship of the form%
\[
h_{X}=h_{d}\left(  N_{X}\right)  ^{-a_{d}d},
\]

for \textbf{higher dimensions}.

\item A similar approach, following Chapter \ref{Ch_Interpol} and Chapter
\ref{Ch_ExactSmth} for the Exact smoother, is to substitute for $N_{X}$ and
$N_{X^{\prime}}$ in \ref{s46} and minimize the estimator for $\rho$.
\end{enumerate}
\end{remark}

\subsection{Convergence to the data
function\label{SbSect_AppSmth_to_DatFunc_Lagran}}

In the previous subsection we studied the convergence of the Approximate
smoother to the Exact smoother. In this subsection we will combine these
results with the Exact smoother error and thus estimate the error of the
Approximate smoother. The relevant result concerning the Exact smoother error
is Theorem \ref{Thm_Exact_smth_converg_at_x} and our next result combines this
result with Theorem \ref{Thm_ap_rx(x)_bound_better} above:

\begin{theorem}
\label{Thm_ap_Exact_smth_converg_at_x}\ 

\begin{enumerate}
\item The notation and assumptions of the Lagrange interpolation lemma (Lemma
\ref{Lem_app_Lagrange_interpol}) hold.

\item Let $w$ be a weight function with properties W2 and W3 for order
$\theta$ and parm. $\kappa$ and set $\eta=\min\left\{  \theta,\frac{1}%
{2}\left\lfloor \min2\kappa\right\rfloor \right\}  $. Assume $G$ is a basis
function of order $\theta$ such that there exist constants $c_{G},r_{G}>0$ and
$\delta_{G}\geq0$ such that the estimate of the form \ref{s07} holds i.e.
\[
\sqrt{r_{x}\left(  x\right)  }\leq\sqrt{c_{G}}\left(  1+\sum_{k=1}%
^{M}\left\vert l_{k}\left(  x\right)  \right\vert \right)  \left(
\operatorname*{diam}A_{x}\right)  ^{\eta_{G}},\quad\operatorname*{diam}%
A_{x}\leq r_{G},\text{ }x\in\Omega,
\]

where $\eta_{G}=\eta+\delta_{G}$.

\item Denote by $s_{e}$ the Exact smoother generated by the the smoothing
parameter $\rho$, the unisolvent independent data $X\subset\Omega$ and data
function $f_{d}\in X_{w}^{\theta}$.
\end{enumerate}

Then when $h_{X}\leq\min\left\{  h_{\Omega,\theta},r_{G}\right\}  $,
\begin{equation}
\left\vert s_{e}\left(  x\right)  -f_{d}\left(  x\right)  \right\vert
\leq\left\vert f_{d}\right\vert _{w,\theta}\left(  1+K_{\Omega,\theta}%
^{\prime}\right)  \left(  \sqrt{c_{G}}\left(  c_{\Omega,\theta}h_{X}\right)
^{\eta_{G}}+\sqrt{N_{X}\rho}\right)  ,\quad x\in\overline{\Omega}.\label{s45}%
\end{equation}

\end{theorem}

?? We now present our main result for the convergence of the Approximate
smoother to the data function. However, the estimate \ref{s25} has a factor of
$\sqrt{\rho}$ in the denominator which is not in the corresponding zero order
result. \textbf{Perhaps this factor can be eliminated}?

\begin{theorem}
\label{Thm_Se(x)minusFd(x)_bound_RVx}We will use the assumptions and notation
of Theorem \ref{Thm_ap_Exact_smth_converg_at_x} and Lemma
\ref{Lem_app_Lagrange_interpol}. Now denote by $s_{a}$ the Approximate
smoother of the data function $f_{d}$ generated by the data set $X$ and the
points $X^{\prime}$. Let $h_{X^{\prime}}=\sup\limits_{\omega\in\Omega
}\operatorname*{dist}\left(  \omega,X^{\prime}\right)  $ measure the density
of the point set $X^{\prime}$.

Then when $h_{X},h_{X^{\prime}}\leq\min\left\{  h_{\Omega,\theta}%
,r_{G}\right\}  $, $x\in\overline{\Omega}$ and $f_{d}\in X_{w}^{\theta}$ the
estimate%
\begin{align}
& \left\vert f_{d}\left(  x\right)  -s_{a}\left(  x\right)  \right\vert
\nonumber\\
& \leq\left\vert f_{d}\right\vert _{w,\theta}\left(  1+K_{\Omega,\theta
}^{\prime}\right)  \left(  \sqrt{c_{G}}\left(  c_{\Omega,\theta}h_{X}\right)
^{\eta_{G}}+\sqrt{N_{X}\rho}\right)  \left(  2+\left(  1+K_{\Omega,\theta
}^{\prime}\right)  ^{2}c_{G}\frac{\left(  c_{\Omega,\theta}h_{X^{\prime}%
}\right)  ^{\eta_{G}}}{\sqrt{\rho}}\right)  ,\label{s25}%
\end{align}

holds.
\end{theorem}

\begin{proof}
From \ref{s46},%
\begin{align*}
\left\vert s_{e}\left(  x\right)  -s_{a}\left(  x\right)  \right\vert  &
\leq\left\vert f_{d}\right\vert _{w,\theta}\left(  1+\left(  1+K_{\Omega
,\theta}^{\prime}\right)  \sqrt{c_{G}}\frac{\left(  c_{\Omega,\theta
}h_{X^{\prime}}\right)  ^{\eta_{G}}}{\sqrt{\rho}}\right)  \times\\
& \qquad\times\left(  1+K_{\Omega,\theta}^{\prime}\right)  \left(  \sqrt
{c_{G}}\left(  c_{\Omega,\theta}h_{X}\right)  ^{\eta_{G}}+\sqrt{N_{X}\rho
}\right)  ,\quad x\in\overline{\Omega},
\end{align*}

and from \ref{s45}%
\[
\left\vert f_{d}\left(  x\right)  -s_{e}\left(  x\right)  \right\vert
\leq\left\vert f_{d}\right\vert _{w,\theta}\left(  1+K_{\Omega,\theta}%
^{\prime}\right)  \left(  \sqrt{c_{G}}\left(  c_{\Omega,\theta}h_{X}\right)
^{\eta_{G}}+\sqrt{N_{X}\rho}\right)  ,\quad x\in\overline{\Omega},
\]

so that%
\begin{align*}
& \left\vert f_{d}\left(  x\right)  -s_{a}\left(  x\right)  \right\vert \\
& \leq\left\vert f_{d}\left(  x\right)  -s_{e}\left(  x\right)  \right\vert
+\left\vert s_{e}\left(  x\right)  -s_{a}\left(  x\right)  \right\vert \\
& \leq\left\vert f_{d}\right\vert _{w,\theta}\left(  1+K_{\Omega,\theta
}^{\prime}\right)  \left(  \sqrt{c_{G}}\left(  c_{\Omega,\theta}h_{X}\right)
^{\eta_{G}}+\sqrt{N_{X}\rho}\right)  +\\
& \qquad+\left\vert f_{d}\right\vert _{w,\theta}\left(  1+\left(
1+K_{\Omega,\theta}^{\prime}\right)  \sqrt{c_{G}}\frac{\left(  c_{\Omega
,\theta}h_{X^{\prime}}\right)  ^{\eta_{G}}}{\sqrt{\rho}}\right)  \left(
1+K_{\Omega,\theta}^{\prime}\right)  \left(  \sqrt{c_{G}}\left(
c_{\Omega,\theta}h_{X}\right)  ^{\eta_{G}}+\sqrt{N_{X}\rho}\right) \\
& =\left\vert f_{d}\right\vert _{w,\theta}\left(  1+K_{\Omega,\theta}^{\prime
}\right)  \left(  \sqrt{c_{G}}\left(  c_{\Omega,\theta}h_{X}\right)
^{\eta_{G}}+\sqrt{N_{X}\rho}\right)  \left(  2+\left(  1+K_{\Omega,\theta
}^{\prime}\right)  \sqrt{c_{G}}\frac{\left(  c_{\Omega,\theta}h_{X^{\prime}%
}\right)  ^{\eta_{G}}}{\sqrt{\rho}}\right)  .
\end{align*}

\end{proof}

\begin{remark}
\label{Rem_Thm_Se(x)minusFd(x)_bound_RVx}\ 

\begin{enumerate}
\item Comparison of the Approximate smoother error estimate proved in the last
theorem with that of the estimate of part 2 of Theorem
\ref{Thm_ap_Se(x)minusSa(x)_bound_RVx} shows that the convergence analysis
given in Remark \ref{Rem_Thm_Se(x)minusSa(x)_bound_RVx} can be applied to the
estimate \ref{s25}.

\item The Approximate smoother error estimate \ref{s25} is not bounded in
$\rho$ near zero or near infinity. However, estimates that are bounded in
$\rho$ near infinity can be derived as follows: First estimate $\left\vert
s_{e}\left(  x\right)  -s_{a}\left(  x\right)  \right\vert $ by combining the
estimate \ref{s46} with those of either part 3 or part 4 of the same theorem.
Then add an Exact smoother error estimate from Theorem
\ref{Thm_Exact_smth_err_bded_smth_parm} to obtain the Approximate smoother
error bounded in the smoothing parameter.

\item We now show that there exist sequences $X_{k}^{\prime}$, $X_{k}$ and
$\rho_{k}$ such that $\left\vert f_{d}\left(  x\right)  -s_{a}^{\left(
k\right)  }\left(  x\right)  \right\vert \rightarrow0$ uniformly on
$\overline{\Omega}$:

The estimator \ref{s25} has the form $\left(  A+Bx\right)  \left(
C+D/x\right)  $ where $x=\sqrt{\rho}$, $A=\sqrt{c_{G}}\left(  c_{\Omega
,\theta}h_{X}\right)  ^{\eta_{G}}$, $B=\sqrt{N_{X}}$, $C=2$ and $D=\left(
1+K_{\Omega,\theta}^{\prime}\right)  \sqrt{c_{G}}\left(  c_{\Omega,\theta
}h_{X^{\prime}}\right)  ^{\eta_{G}}$.

The estimator is minimized to $\left(  \sqrt{AC}+\sqrt{BD}\right)  ^{2}$ when
$\rho=AD/BC$ i.e.
\begin{align*}
\left\vert f_{d}\left(  x\right)  -s_{a}\left(  x\right)  \right\vert  &
\leq\left(  \sqrt{AC}+\sqrt{BD}\right)  ^{2}\\
& \leq AC+BD\\
& =2\sqrt{c_{G}}\left(  c_{\Omega,\theta}h_{X}\right)  ^{\eta_{G}}+\left(
1+K_{\Omega,\theta}^{\prime}\right)  \sqrt{c_{G}}\left(  c_{\Omega,\theta
}h_{X^{\prime}}\right)  ^{\eta_{G}}\sqrt{N_{X}},
\end{align*}

when%
\begin{align*}
\rho & =\frac{\sqrt{c_{G}}\left(  c_{\Omega,\theta}h_{X^{\prime}}\right)
^{\eta_{G}}\left(  1+K_{\Omega,\theta}^{\prime}\right)  \sqrt{c_{G}}\left(
c_{\Omega,\theta}h_{X}\right)  ^{\eta_{G}}}{2\sqrt{N_{X}}}\\
& =\frac{1}{2}\left(  1+K_{\Omega,\theta}^{\prime}\right)  \sqrt{c_{G}}%
\frac{\left(  c_{\Omega,\theta}h_{X^{\prime}}\right)  ^{\eta_{G}}\left(
c_{\Omega,\theta}h_{X}\right)  ^{\eta_{G}}}{\sqrt{N_{X}}}.
\end{align*}

Hence if $\left(  h_{X_{k}^{\prime}}\right)  ^{\eta_{G}}\sqrt{N_{X_{k}}%
}\rightarrow0$ then $\left\vert f_{d}\left(  x\right)  -s_{a}^{\left(
k\right)  }\left(  x\right)  \right\vert \rightarrow0$.
\end{enumerate}
\end{remark}

\section{A numerical implementation of the Approximate
smoother\label{Sect_numer_implem_Approx_smth}}

\subsection{The \textit{SmoothOperator} software (freeware)}

In this section we discuss a numerical implementation for the construction and
solution of the positive order Approximate smoother matrix equation \ref{s16}
i.e.%
\[
\left(
\begin{array}
[c]{lll}%
\left(  2\pi\right)  ^{\frac{d}{2}}N\rho G_{X^{\prime},X^{\prime}%
}+G_{X^{\prime},X}G_{X,X^{\prime}} & G_{X^{\prime},X}P_{X} & P_{X^{\prime}}\\
P_{X}^{T}G_{X,X^{\prime}} & P_{X}^{T}P_{X} & O_{M}\\
P_{X^{\prime}}^{T} & O_{M} & O_{M}%
\end{array}
\right)  \left(
\begin{array}
[c]{c}%
\alpha^{\prime}\\
\beta^{\prime}\\
\gamma^{\prime}%
\end{array}
\right)  =\left(
\begin{array}
[c]{l}%
G_{X^{\prime},X}\\
P_{X}^{T}\\
O_{M,N}%
\end{array}
\right)  y.
\]

I called this software \textit{SmoothOperator}. In Corollary
\ref{Cor_scalable} it was shown that the construction and solution of this
matrix equation is a scalable process and thus worthy of numerical
implementation. This software also implements the positive order version of
equation derived above and the tutorials concentrate on the positive order
basis functions.

The algorithm has been implemented in \textbf{Matlab 6.0} with a GUI interface
but has only been tested on Windows. \textit{SmoothOperator} can be obtained
by emailing the author. However there is a \textbf{short user document }(4
pages) and the potential user can read this document first to decide whether
they want the software. The top-level directory of the software contains the
file \textit{read\_me.txt} which can also be downloaded separately. A
\textbf{full user manual} comes with the software. The main features of the
full user manual are:\medskip

\textbf{1} \textbf{Tutorials and data experiments }To learn about the system
and the behavior of the algorithm, I have prepared three tutorials and five
data experiments.\medskip

\textbf{2} \textbf{Context-sensitive help }Each dialog box incorporates
context-sensitive help which is invoked using the right mouse button. An F1
key facility could be implemented. The actual help text is contained in the
text file \textbf{%
$\backslash$%
Help%
$\backslash$%
ContextHelpText.m} that can be easily edited.\medskip

\textbf{3} \textbf{Matlab diary facility} When the system is started the
Matlab diary facility is invoked. This means that most user information
generated by \textit{SmoothOperator} and written to the Matlab command line is
also written to the text diary file. Each dialog box has a drop-down diary
menu which allows the user to view the diary file using Notepad. There is also
a facility for you to choose another editor/browser. The file can also be
emptied, if, for example, it gets too big. It can also be disabled.\medskip

\textbf{4} \textbf{Tools for viewing data files} Before generating a smoother
you can view the contents of the data file.

For ASCII delimited text data files, use the (slow) \textbf{View records}
facility to display records and the file header, and then with this
information you can use the high speed \textbf{Study records} facility to
check the records and then obtain detailed information about single fields
e.g. a histogram, and multiple fields e.g. correlation coefficients and
scatter plots.

For binary \textbf{test data files} only the \textbf{View records} facility is
needed. The parameters which generated the file can also be viewed using the
\textbf{Make or View data} option.\medskip

\textbf{5} \textbf{The output data} The output from the experiments and
tutorials mentioned in point 1 above consists of well-documented Matlab one
and two-dimensional plots, command line output and diary output. There is
currently no file output, but this can be implemented on request.\medskip

\textbf{6} \textbf{Reading delimited text files} A MEX C file allows ASCII
delimited text files to be read very quickly. You can specify:

\textbf{6.1} that the file be read in chunks of records, and not in one go.

\textbf{6.2} The ids of the fields to be read. This means that the file can
contain non-numeric fields.

\textbf{6.3} Fields can be checked to ensure they are numeric.\medskip

The tutorials which create smoothers allow the data, smoothed data, and
related functions to be viewed using scatter plots and plots along lines and
planes.\medskip

\textbf{7} Please note that there is no explicit suite of functions -
application programmer interface or API - supplied by \textit{SmoothOperator}
for immediate use by the user. After reflecting on how I would create such an
API, I decided that I lacked the experience to produce it, and that, anyhow,
there were too many possibilities to anticipate. I would expect that possible
users of this software, designed to be applied to perhaps millions of records,
would need to familiarize themselves thoroughly with how this system works. I
urge them to contact me and discuss their application.

\subsection{Algorithms}

The following three algorithms are used in the \textit{SmoothOperator}
software package to calculate the Approximate smoother.

\textbf{Algorithm 1} uses data generated internally according to
specifications supplied by the user using the interface. The smoothing
parameter can be either specified or calculated using an error grid. This
algorithm is \textbf{scalable}.

\textbf{Algorithm 2} uses a `small' subset of actual data to get an idea of a
suitable smoothing parameter to use for the full data set. This algorithm is
\textbf{not scalable}.

\textbf{Algorithm 3} uses all the actual data. The smoothing parameter can be
either specified or calculated using an \textit{error grid}. This algorithm is
\textbf{scalable}.

We will now explain these algorithms in more detail.

\subsubsection{\underline{Algorithm 1: using experimental data}}

\begin{enumerate}
\item Generate the experimental data $\left[  X,y\right]  $. The independent
data\textbf{\ }$X$ is generated by uniformly distributed random numbers on $X
$. The dependent data $y$ is generated by uniformly perturbing an analytic
data function $g_{dat}$.

\item Choose a smoothing grid $X^{\prime}$ whose boundary contains $X$. Choose
an error grid $X_{err}^{\prime}$ which will be used to estimate the optimal
smoothing parameter $\rho$.

\item Read the data $\left[  X,y\right]  $ and construct the matrix equation.
If the matrix $G_{X,X^{\prime}}$ is dense, the matrix $G_{X,X^{\prime}}%
^{T}G_{X,X^{\prime}}$ is constructed using a Matlab MEX file (a compiled C
file). If $G_{X,X^{\prime}}$ is sparse the usual matrix multiplication is used.

\item Given a value for the smoothing parameter $\rho$ we can solve the matrix
equation and evaluate the smoother.

We want to estimate the value of $\rho$ which minimizes the `sum of squares'
error between the smoother and the data function
\[
\delta_{1}\left(  \rho\right)  =\sum\limits_{x^{\prime\prime}\in
X_{err}^{\prime}}\left(  \sigma_{\rho}\left(  x^{\prime\prime}\right)
-g_{dat}\left(  x^{\prime\prime}\right)  \right)  ^{2},\text{\quad}\rho>0.
\]

Empirical work indicates a standard shape for $\delta_{1}\left(  \rho\right)
$, namely decreasing from right to left, reaching a minimum and then
increasing at a decreasing rate. To find the minimum we basically use the
standard iterative algorithm of dividing and multiplying by a factor e.g. 10
and choosing the smallest value. The process is stopped when the percentage
change of one or both of $\delta_{1}\left(  \rho\right)  $ and $\rho$ are less
than prescribed values.
\end{enumerate}

\subsubsection{\underline{Algorithm 2: using a `test' subset of actual data}}

\begin{enumerate}
\item Perhaps based on results using Algorithm 1, choose an initial value for
smoothing\textbf{\ }parameter $\rho$ and choose a smoothing grid $X^{\prime}$.

\item Step 2 of Algorithm 1. Denote the data by $\left[  X,y\right]  $ where
$X=\left(  x^{\left(  i\right)  }\right)  $ and $y=\left(  y_{i}\right)  $.

\item This is the same as step 4 of Algorithm 2 except we now minimize
\begin{equation}
\delta_{2}\left(  \rho\right)  =\sum\limits_{i=1}^{N}\left(  \sigma_{\rho
}\left(  x^{\left(  i\right)  }\right)  -y_{i}\right)  ^{2},\text{\quad}%
\rho>0,\label{1.018}%
\end{equation}

because we do not have a data function.
\end{enumerate}

\subsubsection{\underline{Algorithm 3: using all the actual data}}

\begin{enumerate}
\item Choose a value for smoothing parameter $\rho$, based on experiments
using Algorithm 2. Choose a smoothing grid $X^{\prime}$.

\item Read all the data $\left[  X,y\right]  $ and construct the matrix
equation. If the matrix $G_{X,X^{\prime}}$ is dense, the matrix
$G_{X,X^{\prime}}^{T}G_{X,X^{\prime}}$ is constructed using a Matlab MEX file
(a compiled C file). If $G_{X,X^{\prime}}$ is sparse the usual multiplication
is used.

\item Solve the matrix equation to obtain the basis function coefficients and
evaluate the smoother at the desired points.
\end{enumerate}

\subsection{Features of the smoothing algorithm and its
implementation\label{SbSect_AlgorImplem}}

\textbf{1} The Short user manual and the User manual contains a lot of detail
regarding the \textit{SmoothOperator} system and algorithms. So we will
content ourselves here with just some key points.\medskip

\textbf{2} Although the algorithm is scalable there can still be a problem
with rapidly increasing memory usage as the grid size decreases and the
dimension increases. The classical radial basis functions, such as the thin
plate spline functions, have support everywhere. Hence the smoothing matrix is
completely full. To significantly reduce this problem we do the following :

\textbf{a)} We use basis functions with bounded support.

\textbf{b)} We shrink the basis function support to the magnitude of the grid
cells. This makes $G_{X,X^{\prime}}^{T}G_{X,X^{\prime}}$ a \textbf{very
sparse} banded diagonal matrix, and $G_{X^{\prime},X^{\prime}}+G_{X^{\prime
},X^{\prime}}^{T}$ has a small number of non-zero diagonals. We say this basis
function has \textit{small support}.\medskip

\textbf{3} Instead of using the memory devouring Matlab \texttt{repmat}
function to calculate $G_{X,X^{\prime}}$ and $G_{X^{\prime},X^{\prime}}$, we
directly calculate the arguments of the Matlab function \texttt{sparse}. This
function takes three arrays, namely the row ids, the column ids and the
corresponding matrix elements, as well as the matrix dimensions, and converts
them to the Matlab sparse internal representation. I must mention that
although this is quick and space efficient, the algorithm is much more
complicated than using \texttt{repmat}.

Matlab's sparse multiplication facility is then used to quickly and
efficiently calculate $G_{X,X^{\prime}}^{T}G_{X,X^{\prime}}$.\medskip

\textbf{4} This software has implemented the above techniques for the
smoothing matrix, and a selection of basis functions is supplied.\medskip

\textbf{5} ?? The tutorials and exercises are based around the zero order
\textbf{tensor product hat (triangle) function, }denoted by $\Lambda$. The hat
function, also known as the triangle function, is used because of its
simplicity and its analytic properties. The higher dimensional hat functions
are defined as the tensor product of the one dimensional hat function :
$\Lambda\left(  x_{1}\right)  =1-\left\vert x_{1}\right\vert $, when
$\left\vert x_{1}\right\vert \leq1$, and zero otherwise. This function has
zero order so that the smoothing matrix simplifies to%
\begin{equation}
\Psi_{G}=\rho NG_{X^{\prime},X^{\prime}}+G_{X,X^{\prime}}^{T}G_{X,X^{\prime}%
},\text{\quad}G_{X^{\prime},X^{\prime}}=I,\label{2.39}%
\end{equation}

with matrix equation%
\[
\Psi_{G}\alpha=G_{X^{\prime},X}y.
\]

Smoothers constructed from the hat function are continuous but not smooth. I
have also included the \textbf{B-spline of order 3}, namely $\Lambda
\ast\Lambda$, which is a basis function of order zero or one. This can be used
when a smoother is required to be smooth i.e. continuously
differentiable.\medskip

\textbf{6} \textbf{Note that this software was written before I embarked on
the error analysis contained in this document} so \textit{SmoothOperator}
concentrates on the Approximate smoother and the user cannot study the error
of the Exact smoother or the Approximate smoother or compare the Approximate
smoother with the Exact smoother.

\subsection{An application - predictive modelling of forest cover type}

In this section we demonstrate how the method developed in the previous
sections can be used in data mining for predictive modelling. This application
smooths \textbf{binary-valued} data - I was unable to obtain a large
`continuous-valued' data set for distribution on the web. The source of our
data is the web site file:%
\[%
\begin{tabular}
[c]{l}%
\textit{http://kdd.ics.uci.edu/databases/covertype/covertype.html,}%
\end{tabular}
\]

in the UCI KDD Archive, Information and Computer Science, University of
California, Urvine.

The data gives the forest cover type in $30\times30$ meter cells as a function
of the following cartographic parameters:%

\begin{table}[htbp] \centering
$%
\begin{tabular}
[c]{|l|l|l|}\hline
Id & Independent variable & Description\\\hline\hline
1 & ELEVATION & Altitude above sea level\\
2 & ASPECT & Azimuth\\
3 & SLOPE & Inclination\\
4 & HORIZ\_HYDRO & Horizontal distance to water\\
5 & VERT\_HYDRO & Vertical distance to water\\
6 & HORIZ\_ROAD & Horizontal distance to roadways\\
7 & HILL\_SHADE\_9 & Hill shade at 9am\\
8 & HILL\_SHADE\_12 & Hill shade at noon\\
9 & HILL\_SHADE\_15 & Hill shade at 3pm\\
10 & HORIZ\_FIRE & Horizontal distance to fire points\\\hline
\end{tabular}
$\caption{}\label{Tbl_cart_parms}%
\end{table}%

Forest cover type is the dependent variable\ $y$ and it takes on one of seven values:%

\begin{table}[htbp] \centering
$%
\begin{tabular}
[c]{|l|l|}\hline
Forest cover type & Id\\\hline\hline
Spruce fir & 1\\
Lodge-pole pine & 2\\
Ponderosa pine & 3\\
Cottonwood/Willow & 4\\
Aspen & 5\\
Douglas fir & 6\\
Krummholtz & 7\\\hline
\end{tabular}
$\caption{}\label{Tbl_forest_cover}%
\end{table}%

\subsubsection{\underline{The data file}}

In this study we will use the data to train a model predicting on the
\textbf{presence or absence of the Ponderosa pine forest cover (}id = 3), but
the results will be similar for the other forest types. To this end we have
created from the full web site file a file called \texttt{%
$\backslash$%
UserData%
$\backslash$%
forest\_1\_to\_10\_pondpin.dat} which contains the ten independent variables
of Table \ref{Tbl_cart_parms} and then a binary dependent variable derived
from the variable \textbf{Ponderosa pine} of Table \ref{Tbl_forest_cover}.
This variable is 1 if the cover is Ponderosa pine and zero otherwise.

\subsubsection{\underline{Methodology}}

We recommend you first use the user interface of the \textit{SmoothOperator
}software\textit{\ }to construct artificial data sets to understand the
behavior of the smoother and run the experiments to get a feel for the
influence of the parameters.\medskip

\textbf{1} We chose a small\ subset of the forest-cover data and selected
various smoothing parameters and smoothing grid sizes to study the performance
of the smoother using plots and the value of the error.

Note that \textit{two subsets of data} could be used here, often called the
\textit{training set} and the \textit{test set}. The training set would be
used to calculate the smoother for the initial value of the smoothing
coefficient, and then the test set could be used to determine the smoothing
coefficient which minimizes the least squares\ error.\medskip

\textbf{2} Having chosen our parameters we run the smoother program on the
full data set.

\chapter{The bounded linear functionals on the data space $X_{w}^{\theta}%
$\label{Ch_bnd_lin_fnal_Xwth}}

In this chapter we study the spaces $X_{1/w}^{-\theta}$ where $\theta\geq1 $.
Note that we use the notation $X_{1/w}^{-\theta}$ instead of $X_{w}^{-\theta}$
so we can distinguish between the space $X_{1/w}^{0}$ and the space $X_{w}%
^{0}$. The spaces $X_{1/w}^{-\theta}$ are thus not the negative order versions
of $X_{w}^{\theta}$ but are actually isometrically isomorphic to the spaces of
bounded linear functionals on $X_{w}^{\theta}$, denoted $\left(  X_{w}%
^{\theta}\right)  ^{\prime}$. The properties of the spaces $X_{1/w}^{-\theta}$
are established by constructing inverse isometric mappings between
$X_{1/w}^{-\theta}$ and $X_{w}^{\theta}$. These mapping is intimately related
to the basis distribution of order $\theta$.

We start by constructing $X_{1/w}^{-\theta}$ as a space of tempered
distributions $S^{\prime}$ and construct isometric mappings $\mathcal{L}$ and
$\mathcal{M}$ which are 1-1 and onto in the seminorm sense - see the figure
below.%
\begin{gather*}%
\begin{array}
[c]{ccccc}%
X_{1/w}^{-\theta} &
\begin{array}
[c]{l}%
\mathcal{L}\\
\mathbf{\longleftarrow}\\
\longrightarrow\\
\mathcal{M}%
\end{array}
& L^{2} &
\begin{array}
[c]{l}%
\mathcal{I}\\
\mathbf{\longleftarrow}\\
\longrightarrow\\
\mathcal{J}%
\end{array}
& X_{w}^{\theta}%
\end{array}
\\
\medskip\\
\fbox{Mappings between $X_{1/w}^{-\theta}$, $L^{2}$ and $X_{w}^{\theta}$.}%
\end{gather*}

To define $\mathcal{M}$ conditions are placed on the weight function which are
only satisfied by the thin-plate spline basis functions discussed in Section
\ref{Sect_rad_basis_fns} - ??? SHOW THIS!. To overcome this limitation the
space $S^{\prime}$ is replaced by the bounded linear functionals on a subspace
of $S$ defined by the weight function.

??? I need to ADD some examples as was done for the zero order case. Examples
should include the central difference and extended B-spline weight functions.

\section{The space $X_{1/w}^{-\theta},$ $\theta\geq1$}

This result will be needed to show that the space $X_{1/w}^{-\theta}$ defined
below is not empty.

\begin{lemma}
\label{Lem_Co,inf(Rd/C)}Let $\mathcal{B}$ be a closed set of measure zero.
Define
\[
C_{0}^{\infty}\left(  \mathbb{R}^{d}\setminus\mathcal{B}\right)  =\left\{
\psi\in C_{0}^{\infty}:\operatorname*{supp}\psi\subset\mathbb{R}^{d}%
\setminus\mathcal{B}\right\}  .
\]

Then $C_{0}^{\infty}\left(  \mathbb{R}^{d}\setminus\mathcal{B}\right)  \subset
L_{loc}^{1}\left(  \mathbb{R}^{d}\setminus\mathcal{B}\right)  \cap
S_{w,\theta}$ and $\int\frac{\left\vert \widehat{\phi}\right\vert ^{2}%
}{w\left\vert \cdot\right\vert ^{2\theta}}<\infty$ when $\phi\in C_{0}%
^{\infty}\left(  \mathbb{R}^{d}\setminus\mathcal{B}\right)  $.

Note that the space $S_{w,\theta}$ is introduced in Definition
\ref{Def_Sw,th_and_invF(Sw,th)} \textbf{below}.
\end{lemma}

\begin{proof}
Since $\mathcal{B}$ is a closed set $C_{0}^{\infty}\left(  \mathbb{R}%
^{d}\setminus\mathcal{B}\right)  \subset L_{loc}^{1}\left(  \mathbb{R}%
^{d}\setminus\mathcal{B}\right)  $.

If $\phi\in C_{0}^{\infty}\left(  \mathbb{R}^{d}\setminus\mathcal{B}\right)  $
then $\left(  \int w\left\vert \cdot\right\vert ^{2\theta}\left\vert
\phi\right\vert ^{2}\right)  ^{1/2}<\infty$ and so $\phi\in S_{w,\theta}$ from
the Definition \ref{Def_Sw,th_and_invF(Sw,th)} of $S_{w,\theta}$.

Since $w\left\vert \cdot\right\vert ^{2\theta}>0$ on $\mathbb{R}^{d}%
\setminus\mathcal{B}$ we have $w\left\vert \cdot\right\vert ^{2\theta}>0$ on
the compact set $\operatorname*{supp}\phi$ and hence%
\[
\int\frac{\left\vert \widehat{\phi}\right\vert ^{2}}{w\left\vert
\cdot\right\vert ^{2\theta}}\leq\frac{1}{\min\limits_{\operatorname*{supp}%
\phi}w\left\vert \cdot\right\vert ^{2\theta}}\int\left\vert \phi\right\vert
^{2}<\infty.
\]

\end{proof}

\begin{definition}
\label{Def_Xneg_order}\textbf{The spaces }$X_{1/w}^{-\theta}$, $\theta\geq1$.

Suppose $w$ is a weight function i.e. it has property W1. Since $w\in
C^{\left(  0\right)  }\left(  \mathbb{R}^{d}\setminus\mathcal{A}\right)  $ we
have $w\left\vert \cdot\right\vert ^{2\theta}\in C^{\left(  0\right)  }\left(
\mathbb{R}^{d}\setminus\mathcal{B}\right)  $\ where $\mathcal{B}=\mathcal{A}$
or $\mathcal{B}=\mathcal{A}\setminus0$, whichever is the smallest valid set.
We want $\mathcal{B}$ to be a closed set so if $\mathcal{B}=\mathcal{A}%
\setminus0$ we assume $0$ is an isolated point. This is because $\mathcal{A}%
\setminus0$ is a closed set iff $0$ is an isolated point. Clearly $w\left\vert
\cdot\right\vert ^{2\theta}>0$ on $\mathcal{B}$.

Now define the semi-inner product space%
\begin{equation}
X_{1/w}^{-\theta}=\left\{  u\in S^{\prime}:u_{F}\in L_{loc}^{1}\left(
\mathbb{R}^{d}\setminus\mathcal{B}\right)  \text{ }and\text{ }\int%
\frac{\left\vert u_{F}\right\vert ^{2}}{w\left\vert \cdot\right\vert
^{2\theta}}<\infty\right\}  ,\label{p1.016}%
\end{equation}

where $u_{F}=\widehat{u}$ on $\mathbb{R}^{d}\setminus\mathcal{B}$ and
$u_{F}=0$ on $\mathcal{B}$. We endow $X_{1/w}^{-\theta}$ with the semi-inner
product%
\[
\left\langle u,v\right\rangle _{1/w,-\theta}=\left(  \int\frac{u_{F}%
\overline{v_{F}}}{w\left\vert \cdot\right\vert ^{2\theta}}\right)  ^{1/2},
\]

and the norm will be denoted by $\left\vert u\right\vert _{1/w,-\theta}$.
Observe that the seminorm has null space $\left(  S_{\mathcal{B}}^{\prime
}\right)  ^{\vee}$ where
\[
S_{\mathcal{B}}^{\prime}=\left\{  u\in S^{\prime}:\operatorname*{supp}%
u\subset\mathcal{B}\right\}  .
\]

Also note that Lemma \ref{Lem_Co,inf(Rd/C)} implies $\left(  C_{0}^{\infty
}\left(  \mathbb{R}^{d}\setminus\mathcal{B}\right)  \right)  ^{\vee}\subset
X_{1/w}^{-\theta}$, and so $X_{1/w}^{-\theta}$ is not empty.
\end{definition}

\section{The operator $\mathcal{M}:X_{1/w}^{-\theta}\rightarrow L^{2}$}

In this section we study the analogue of the mapping $\mathcal{J}%
:L^{2}\rightarrow X_{w}^{\theta}$ from $X_{1/w}^{-\theta}$ to $L^{2}$, which
we denote by $\mathcal{M}$.

\begin{definition}
\label{Def_op_V}\textbf{The operator} $\mathcal{M}:X_{1/w}^{-\theta
}\rightarrow L^{2}$

Suppose the weight function $w$ has property W1. From the definition of
$X_{1/w}^{-\theta}$, $u\in X_{1/w}^{-\theta}$ implies $\frac{u_{F}}{\sqrt
{w}\left\vert \cdot\right\vert ^{\theta}}\in L^{2}$. We can now define the
linear mapping $\mathcal{M}:X_{1/w}^{-\theta}\rightarrow L^{2}$ by%
\begin{equation}
\mathcal{M}u=\left(  \frac{u_{F}}{\sqrt{w}\left\vert \cdot\right\vert
^{\theta}}\right)  ^{\vee},\text{\quad}u\in X_{1/w}^{-\theta}.\label{p1.067}%
\end{equation}

\end{definition}

The operator $\mathcal{M}:X_{1/w}^{-\theta}\rightarrow L^{2}$ has the
following properties:

\begin{theorem}
\label{Thm_prop_Mth} The linear operator $\mathcal{M}:X_{1/w}^{-\theta
}\rightarrow L^{2}$ is isometric with null space $\left(  S_{\mathcal{B}%
}^{\prime}\right)  ^{\vee}$.
\end{theorem}

\begin{proof}
That $\mathcal{M}$ is an isometry is clear from the definition of
$X_{1/w}^{-\theta}$. Since $\mathcal{M}$ is an isometry the null space of
$\mathcal{M}$ is the null space of the seminorm $\left\vert u\right\vert
_{1/w,\theta}$, namely $\left(  S_{\mathcal{B}}^{\prime}\right)  ^{\vee}$.
\end{proof}

\section{The operator $\mathcal{L}:L^{2}\rightarrow X_{1/w}^{-\theta}$}

Making certain assumptions about the weight function we introduce the operator
$\mathcal{L}:L^{2}\rightarrow S^{\prime}$. This is the analogue of the mapping
$\mathcal{I}:X_{w}^{\theta}\rightarrow L^{2}$. We show that $\mathcal{L}%
:L^{2}\rightarrow X_{1/w}^{\theta}$ and that it has nice properties in the
seminorm sense.

\begin{lemma}
\label{Lem_wt_fn_extra_property_1}Suppose the weight function $w$ has property
W1. Suppose also that $w$ satisfies the conditions
\begin{equation}
w\left\vert \cdot\right\vert ^{2\theta}\in L_{loc}^{1}\text{ }and\text{ }%
\int\limits_{\left\vert \cdot\right\vert \geq r}\frac{w\left\vert
\cdot\right\vert ^{2\theta}}{\left\vert \cdot\right\vert ^{2\tau}}%
<\infty\text{ }for\text{ }some\text{ }\tau\geq0\text{ }and\text{
}r>0.\label{p1.026}%
\end{equation}

Then $g\in L^{2}$ implies $\sqrt{w}\left\vert \cdot\right\vert ^{\theta}g\in
L_{loc}^{1}\cap S^{\prime}$.
\end{lemma}

\begin{proof}
Suppose $K$ is compact. Then by the Cauchy-Schwartz inequality%
\[
\int\limits_{K}\sqrt{w}\left\vert \cdot\right\vert ^{\theta}\left\vert
g\right\vert \leq\left(  \int\limits_{K}w\left\vert \cdot\right\vert
^{2\theta}\right)  ^{1/2}\left\Vert g\right\Vert _{2}<\infty.
\]

Also%
\[
\int_{\left\vert \cdot\right\vert \geq r}\frac{\sqrt{w}\left\vert
\cdot\right\vert ^{\theta}\left\vert g\right\vert }{\left\vert \cdot
\right\vert ^{\tau}}\leq\left(  \int_{\left\vert \cdot\right\vert \geq r}%
\frac{w\left\vert \cdot\right\vert ^{2\theta}}{\left\vert \cdot\right\vert
^{2\tau}}\right)  ^{1/2}\left\Vert g\right\Vert _{2}<\infty.
\]

Together, these two inequalities imply that $\sqrt{w}\left\vert \cdot
\right\vert ^{\theta}\widehat{g}\in L_{loc}^{1}\cap S^{\prime}$.
\end{proof}

\begin{theorem}
Weight function condition \ref{p1.026} implies that $X_{1/w}^{-\theta}$ is an
inner product space.
\end{theorem}

\begin{proof}
Suppose $u\in X_{1/w}^{-\theta}$. Then, if $K$ is compact
\[
\int\limits_{K}\left\vert u_{F}\right\vert =\int\limits_{K}\sqrt{w}\left\vert
\cdot\right\vert ^{\theta}\frac{\left\vert u_{F}\right\vert }{\sqrt
{w}\left\vert \cdot\right\vert ^{\theta}}\leq\left(  \int\limits_{K}%
w\left\vert \cdot\right\vert ^{2\theta}\right)  ^{1/2}\left\vert u\right\vert
_{1/w,-\theta}<\infty,
\]

and $\widehat{u}=u_{F}\in L_{loc}^{1}$. Thus
\[
X_{1/w}^{-\theta}=\left\{  u\in S^{\prime}:\widehat{u}\in L_{loc}^{1}\text{
}and\text{ }\int\frac{\left\vert \widehat{u}\right\vert ^{2}}{w\left\vert
\cdot\right\vert ^{2\theta}}<\infty\right\}  ,
\]

and $\left\vert u\right\vert _{1/w,-\theta}=0$ implies $u=0$.
\end{proof}

We now define the operator $\mathcal{L}:L^{2}\rightarrow S^{\prime}$. It will
then be shown that $\mathcal{L}:L^{2}\rightarrow X_{1/w}^{-\theta}$ and that
$\mathcal{L}$ is an inverse of $\mathcal{M}$. This will allow us to prove that
$\mathcal{M}$ is onto and hence that $X_{1/w}^{-\theta}$ is complete.

\begin{definition}
\label{Def_Lth_I}\textbf{The operator} $\mathcal{L}:L^{2}\rightarrow
S^{\prime}$

Suppose the weight function $w$ has the properties assumed in Lemma
\ref{Lem_wt_fn_extra_property_1}. Then we can define the operator
$\mathcal{L}:L^{2}\rightarrow S^{\prime}$ by
\[
\mathcal{L}g=\left(  \sqrt{w}\left\vert \cdot\right\vert ^{\theta}\widehat
{g}\right)  ^{\vee},\quad g\in L^{2}.
\]

\end{definition}

The next theorem shows that $\mathcal{L}$ maps $L^{2}$ to $X_{1/w}^{-\theta} $.

\begin{theorem}
\label{Thm_prop_Lth}Suppose the weight function $w$ has the properties assumed
in Lemma \ref{Lem_wt_fn_extra_property_1}.

Then $\mathcal{L}:L^{2}\rightarrow X_{1/w}^{-\theta}$, $\mathcal{L}$ is an
isometry and 1-1.
\end{theorem}

\begin{proof}
From Lemma \ref{Lem_wt_fn_extra_property_1} we know that if $g\in L^{2}$ then
$\mathcal{L}g\in S^{\prime}$, $\widehat{\mathcal{L}g}\in L_{loc}^{1}$. From
the definition of $\mathcal{L}$%
\[
\left\vert \mathcal{L}g\right\vert _{1/w,-\theta}^{2}=\int\frac{\left\vert
\widehat{\mathcal{L}g}\right\vert ^{2}}{w\left\vert \cdot\right\vert
^{2\theta}}=\int\frac{\left\vert \sqrt{w}\left\vert \cdot\right\vert ^{\theta
}\widehat{g}\right\vert ^{2}}{w\left\vert \cdot\right\vert ^{2\theta}%
}=\left\Vert g\right\Vert _{2}^{2}<\infty.
\]

and so $\mathcal{L}g\in X_{1/w}^{-\theta}$ and $\mathcal{L}$ is an isometry.
Clearly $\mathcal{L}g=0$ implies $g=0$.
\end{proof}

\begin{theorem}
Suppose the weight function $w$ has the properties \ref{p1.026} assumed in
Lemma \ref{Lem_wt_fn_extra_property_1}.

Then $\mathcal{ML}=I$ and $\mathcal{LM}=I$.
\end{theorem}

\begin{proof}
From the proof of Theorem \ref{Thm_prop_Lth}, $u\in X_{1/w}^{-\theta}$ implies
$\widehat{u}=u_{F}$. Hence $\mathcal{M}u=\left(  \frac{\widehat{u}}{\sqrt
{w}\left\vert \cdot\right\vert ^{\theta}}\right)  ^{\vee}$ for $u\in
X_{1/w}^{-\theta}$ and, since $\mathcal{L}g=\left(  \sqrt{w}\left\vert
\cdot\right\vert ^{\theta}\widehat{g}\right)  ^{\vee}$ for $g\in L^{2}$, it
follows immediately that $\mathcal{ML}=I$ and $\mathcal{LM}=I$. Also, we know
from Theorems \ref{Thm_prop_Lth} and \ref{Thm_prop_Mth} that $\mathcal{L}$ and
$\mathcal{M}$ are isometries.
\end{proof}

Since $L^{2}$ is complete, the mappings of the previous theorem directly yield
the following important result.

\begin{corollary}
Suppose the weight function $w$ has the properties assumed in Lemma
\ref{Lem_wt_fn_extra_property_1}. Then $X_{1/w}^{-\theta}$ is a semi-Hilbert space.
\end{corollary}

In Section \ref{Sect_rad_basis_fns} weight functions were generated by the
power function $\left\vert \cdot\right\vert ^{2\theta}$ and three standard
types of radial basis function: the Gaussian, the thin-plate splines and the
shifted thin-plate splines (the multiquadrics are specific cases of the
shifted thin-plate splines). Conditions were placed on the order $\theta$ and
the parameters which defined each basis function so that they were basis
functions in the sense of Light and Wayne. In fact, from the definition of a
basis function, we have $\frac{1}{\widehat{G}}=w\left\vert \cdot\right\vert
^{2\theta}$ outside the set $\mathcal{B}$. Thus the conditions of Lemma
\ref{Lem_wt_fn_extra_property_1} are satisfied for a weight function generated
by a thin-plate spline function. Here $\frac{1}{\widehat{G}}=w\left\vert
\cdot\right\vert ^{2\theta}$ increases at a polynomial rate so $\int%
\limits_{\left\vert \cdot\right\vert \geq r}\frac{w\left\vert \cdot\right\vert
^{2\theta}}{\left\vert \cdot\right\vert ^{2\tau}}<\infty$ is satisfied for
some $\tau$. But for the Gaussian basis function $w\left\vert \cdot\right\vert
^{2\theta}=\frac{2}{\sqrt{\pi}}e^{\frac{1}{4}\left\vert \cdot\right\vert ^{2}%
}$, and for a shifted thin-plate spline basis function $w\left\vert
\cdot\right\vert ^{2\theta}=e^{\beta\left\vert \cdot\right\vert }$ where
$\beta>0$, so the conditions of Lemma \ref{Lem_wt_fn_extra_property_1} are not satisfied.

The extended natural spline weight functions have powers of sines in the
denominator and clearly do not satisfy property \ref{p1.026}.

We summarize these observations in:

\begin{theorem}
Regarding the radial basis functions discussed in Section
\ref{Sect_rad_basis_fns} and the extended natural spline weight functions of
Theorem \ref{Thm_ExtNatSplin_wt_fn_W3.1*}:

\begin{enumerate}
\item The thin-plate spline weight functions satisfy property \ref{p1.026}.

\item The Gaussian weight functions satisfy $w\left\vert \cdot\right\vert
^{2\theta}=e^{\frac{1}{4}\left\vert \cdot\right\vert ^{2}}$ and do not satisfy
property \ref{p1.026}.

\item The shifted thin-plate splines weight functions satisfy $w\left\vert
\cdot\right\vert ^{2\theta}=e^{\beta\left\vert \cdot\right\vert }$ for some
$\beta>0$ and do not satisfy property \ref{p1.026}.

\item The extended natural spline weight functions do not satisfy property
\ref{p1.026}.
\end{enumerate}
\end{theorem}

In the next section we tackle these difficulties with the Gaussian, the
shifted thin-plate splines and the extended natural splines weight functions
by replacing $S^{\prime}$ by a larger space and then enlarging the space
$X_{1/w}^{-\theta}$. This will cope with all weight functions which have
property W1.

\section{The space $S_{w,\theta}$ and the operator $\mathcal{L}_{2}%
:L^{2}\rightarrow S^{\prime}$}

We start by defining a special subspace of $S$ and denote it by $S_{w,\theta}%
$, where $w$ is a weight function with $\theta\geq1$. Then we construct
analogues of $\mathcal{L}$ and $\mathcal{M}$ which we denote by $\mathcal{L}%
_{2}$ and $\mathcal{M}_{2}$ respectively.

\begin{definition}
\label{Def_Sw,th_and_invF(Sw,th)}\textbf{The spaces} $S_{w,\theta}$,
$\theta\geq1$.

Suppose the weight function $w$ has property W1. Then
\[
S_{w,\theta}=\left\{  \phi\in S:\int w\left\vert \cdot\right\vert ^{2\theta
}\left\vert \phi\right\vert ^{2}<\infty\right\}  .
\]

We note that $S_{w,\theta}$ is not empty because $w\left\vert \cdot\right\vert
^{2\theta}\in C^{\left(  0\right)  }\left(  \mathbb{R}^{d}\setminus
\mathcal{B}\right)  $ implies $C_{0}^{\infty}\left(  \mathbb{R}^{d}%
\setminus\mathcal{B}\right)  \subset S_{w,\theta}$.
\end{definition}

\begin{remark}
?? Denote the inverse Fourier transforms of functions in $S_{w,\theta}$ by
$\overset{\vee}{S}_{w,\theta}$. Then $\overset{\vee}{S}_{w,\theta}%
=X_{w}^{\theta}\cap S$. If $w$ also has property W2 then by Theorem
\ref{Thm_Co_inf_dense_Xwth}, $\overset{\vee}{S}_{w,\theta}$is dense in
$X_{w}^{\theta}$.
\end{remark}

\begin{definition}
\label{Def_Lth_II}\textbf{The space }$W_{S;\theta}$\textbf{\ and the operator}
$\mathcal{L}_{2}$, $\theta\geq1$.

Suppose the weight function $w$ satisfies property W1.

Suppose the non-linear functional $\left(  \int w\left\vert \cdot\right\vert
^{2\theta}\left\vert \phi\right\vert ^{2}\right)  ^{1/2}$ defined on
$S_{w,\theta}$ is a member of $S_{w,\theta}^{\prime}$ with respect to $S$
i.e.
\begin{equation}
\left(  \int w\left\vert \cdot\right\vert ^{2\theta}\left\vert \phi\right\vert
^{2}\right)  ^{1/2}\leq\left\vert \phi\right\vert _{\sigma},\quad\phi\in
S_{w,\theta},\label{p1.043}%
\end{equation}

where $\left\vert \cdot\right\vert _{\sigma}$ is a finite, positive linear
combination of the seminorms used to define the topology on $S$.

We say $w\in W_{S;\theta}$ if there exists some positive, linear combination
$\left\vert \cdot\right\vert _{\sigma}$ of seminorms used to define the
topology on $S$ such that the \textbf{non-linear} functional $\left(  \int
w\left\vert \cdot\right\vert ^{2\theta}\left\vert \phi\right\vert ^{2}\right)
^{1/2}$ satisfies%
\[
\left(  \int w\left\vert \cdot\right\vert ^{2\theta}\left\vert \phi\right\vert
^{2}\right)  ^{1/2}\leq\left\vert \phi\right\vert _{\sigma},\quad\phi\in
S_{w,\theta}.
\]

Then for $\phi\in S_{w,\theta}$ and $g\in L^{2}$%
\[
\left\vert \int\sqrt{w}\left\vert \cdot\right\vert ^{\theta}\widehat{g}%
\phi\right\vert =\left\vert \int\left(  \sqrt{w}\left\vert \cdot\right\vert
^{\theta}\phi\right)  \widehat{g}\right\vert \leq\left(  \int w\left\vert
\cdot\right\vert ^{2\theta}\left\vert \phi\right\vert ^{2}\right)
^{1/2}\left\Vert g\right\Vert _{2},
\]

so that $\sqrt{w}\left\vert \cdot\right\vert ^{\theta}\widehat{g}\in
S_{w,\theta}^{\prime}$. We can then (non-uniquely) extend $\sqrt{w}\left\vert
\cdot\right\vert ^{\theta}\widehat{g}$ to $S$ as a member of $S^{\prime}$.
Denote such an extension by $\left(  \sqrt{w}\left\vert \cdot\right\vert
^{\theta}\widehat{g}\right)  ^{e}\in S^{\prime}$. We can now define the class
of operators $\mathcal{L}_{2}:L^{2}\rightarrow S^{\prime}$ by%
\[
\mathcal{L}_{2}g=\left(  \left(  \sqrt{w}\left\vert \cdot\right\vert ^{\theta
}\widehat{g}\right)  ^{e}\right)  ^{\vee},\quad g\in L^{2}.
\]

In general $\mathcal{L}_{2}$ is not unique and not linear.
\end{definition}

\begin{remark}
???

When the basis function is the Gaussian or a shifted thin-plate spline the
functional $\left(  \int w\left\vert \cdot\right\vert ^{2\theta}\left\vert
\phi\right\vert ^{2}\right)  ^{1/2}$ defined on $S_{w,\theta}$ is \textbf{not}
a member of $S_{w,\theta}^{\prime}$ with respect to $S$ (Prove this remark -
see Subsection 9.3.6 of zero order document Williams
\cite{WilliamsZeroOrdSmthV4}).

Also show that the extended natural spline weight functions introduced in
Theorem \ref{Thm_ExtNatSplin_wt_fn_W3.1*} satisfy property \ref{p1.043}.
\end{remark}

The operator $\mathcal{L}_{2}$ has the following properties analogous to those
of the operator $\mathcal{J}$ - see Theorem \ref{Thm_properties_op_Jtheta}.

\begin{theorem}
\label{Thm_prop_Lth_II}The operator $\mathcal{L}_{2}$ has the following properties:

\begin{enumerate}
\item $\mathcal{L}_{2}:L^{2}\rightarrow X_{1/w}^{-\theta}$.

\item $\mathcal{L}_{2}$ is an isometry in the seminorm sense.

\item $\mathcal{L}_{2}$ is 1-1.

\item $\mathcal{L}_{2}$ is linear in the sense that%
\[
\mathcal{L}_{2}\left(  \lambda_{1}g_{1}+\lambda_{2}g_{2}\right)
-\mathcal{L}_{2}\left(  \lambda_{1}g_{1}\right)  -\mathcal{L}_{2}\left(
\lambda_{2}g_{2}\right)  \in\left(  S_{\mathcal{B}}^{\prime}\right)  ^{\vee}.
\]

\end{enumerate}
\end{theorem}

\begin{proof}
\textbf{Part 1}. From the definition of $\mathcal{L}_{2}$, $\widehat
{\mathcal{L}_{2}g}=\sqrt{w}\left\vert \cdot\right\vert ^{\theta}\widehat{g}$
on $S_{w,\theta}$. Next observe that because $w\left\vert \cdot\right\vert
^{2\theta}\in C^{\left(  0\right)  }\left(  \mathbb{R}^{d}\setminus
\mathcal{B}\right)  $ we have $C_{0}^{\infty}\left(  \mathbb{R}^{d}%
\setminus\mathcal{B}\right)  \subset S_{w,\theta}$ and hence
\[
\left[  \widehat{\mathcal{L}_{2}g},\phi\right]  =\left[  \sqrt{w}\left\vert
\cdot\right\vert ^{\theta}\widehat{g},\phi\right]  ,\quad\phi\in C_{0}%
^{\infty}\left(  \mathbb{R}^{d}\setminus\mathcal{B}\right)  .
\]

If we can show $\sqrt{w}\left\vert \cdot\right\vert ^{\theta}\widehat{g}\in
L_{loc}^{1}\left(  \mathbb{R}^{d}\setminus\mathcal{B}\right)  $ it follows
that $\widehat{\mathcal{L}_{2}g}\in L_{loc}^{1}\left(  \mathbb{R}^{d}%
\setminus\mathcal{B}\right)  $.

But if $K\subset\mathbb{R}^{d}\setminus\mathcal{B}$ is compact, $w\left\vert
\cdot\right\vert ^{2\theta}\in C^{\left(  0\right)  }\left(  \mathbb{R}%
^{d}\setminus\mathcal{B}\right)  $ implies
\[
\int_{K}\sqrt{w}\left\vert \cdot\right\vert ^{\theta}\left\vert \widehat
{g}\right\vert \leq\left(  \int_{K}w\left\vert \cdot\right\vert ^{2\theta
}\right)  ^{1/2}\left\Vert g\right\Vert _{2}\leq\max_{K}\left(  w\left\vert
\cdot\right\vert ^{2\theta}\right)  \left(  \int_{K}1\right)  \left\Vert
g\right\Vert _{2}.
\]

Thus $\widehat{\mathcal{L}_{2}g}\in L_{loc}^{1}\left(  \mathbb{R}^{d}%
\setminus\mathcal{B}\right)  $ and $\left(  \mathcal{L}_{2}g\right)
_{F}=\sqrt{w}\left\vert \cdot\right\vert ^{\theta}\widehat{g}$.\medskip

\textbf{Part 2}. $\left\vert \mathcal{L}_{2}g\right\vert _{1/w,-\theta
}=\left(  \int\frac{\left\vert \left(  \mathcal{L}_{2}g\right)  _{F}%
\right\vert ^{2}}{w\left\vert \cdot\right\vert ^{2\theta}}\right)
^{1/2}=\left\Vert g\right\Vert _{2}$.\medskip

\textbf{Part 3}. Follows directly from isometry.\medskip

\textbf{Part 4}. Since
\[
\left(  \mathcal{L}_{2}\left(  \lambda_{1}g_{1}+\lambda_{2}g_{2}\right)
-\mathcal{L}_{2}\left(  \lambda_{1}g_{1}\right)  -\mathcal{L}_{2}\left(
\lambda_{2}g_{2}\right)  \right)  _{F}=0,
\]

we have
\[
\left\vert \mathcal{L}_{2}\left(  \lambda_{1}g_{1}+\lambda_{2}g_{2}\right)
-\mathcal{L}_{2}\left(  \lambda_{1}g_{1}\right)  -\mathcal{L}_{2}\left(
\lambda_{2}g_{2}\right)  \right\vert _{1/w,-\theta}^{2}=0.
\]

\end{proof}

The following theorem shows how the operators $\mathcal{L}_{2}:L^{2}%
\rightarrow X_{1/w}^{-\theta}$ and $\mathcal{M}:X_{1/w}^{-\theta}\rightarrow
L^{2}$ interact.

\begin{theorem}
\label{Thm_L2M_ML2}Suppose the weight function $w$ has the properties assumed
in Definition \ref{Def_Lth_II}. Then $\mathcal{M}$ and $\mathcal{L}_{2}$ are
inverses in the following sense: $\mathcal{ML}_{2}=I$ on $X_{1/w}^{-\theta}$
and $\mathcal{L}_{2}\mathcal{M}u-u\in\left(  S_{\mathcal{B}}^{\prime}\right)
^{\vee}$ when $u\in L^{2}$.
\end{theorem}

\begin{proof}
By definition,$\mathcal{M}u=\left(  \frac{u_{F}}{\sqrt{w}\left\vert
\cdot\right\vert ^{\theta}}\right)  ^{\vee}$ when $u\in X_{1/w}^{-\theta}$.

From the proof of Theorem \ref{Thm_prop_Lth_II} $\left(  \mathcal{L}%
_{2}g\right)  _{F}=\sqrt{w}\left\vert \cdot\right\vert ^{\theta}\widehat{g}$.
Thus from the definition of $\mathcal{M}$, for $g\in L^{2}$
\[
\left(  \mathcal{ML}_{2}g\right)  ^{\wedge}=\frac{\left(  \mathcal{L}%
_{2}g\right)  _{F}}{\sqrt{w}\left\vert \cdot\right\vert ^{\theta}}=\frac
{\sqrt{w}\left\vert \cdot\right\vert ^{\theta}\widehat{g}}{\sqrt{w}\left\vert
\cdot\right\vert ^{\theta}}=\widehat{g},
\]

and thus $\mathcal{ML}_{2}=I$.

Next we show $\mathcal{L}_{2}\mathcal{M}u-u\in\left(  S_{\mathcal{B}}^{\prime
}\right)  ^{\vee}$ when $u\in X_{1/w}^{-\theta}$. In fact, on $S_{w,\theta}$
$\widehat{\mathcal{L}_{2}\mathcal{M}u}=\sqrt{w}\left\vert \cdot\right\vert
^{\theta}\widehat{\mathcal{M}u}=u_{F}$ so that%
\[
\left\vert \mathcal{L}_{2}\mathcal{M}u-u\right\vert _{1/w,-\theta}^{2}%
=\int\frac{\left\vert \left(  \mathcal{L}_{2}\mathcal{M}u-u\right)
_{F}\right\vert ^{2}}{w\left\vert \cdot\right\vert ^{2\theta}}=0,
\]

and $\mathcal{L}_{2}\mathcal{M}u-u\in\left(  S_{\mathcal{B}}^{\prime}\right)
^{\vee}$.
\end{proof}

Since $L^{2}$ is complete, the mappings of the previous theorem will yield the
following important result.

\begin{corollary}
\label{Cor_X1/w,-th_semiHilb_III}If the weight function has property
\ref{p1.043} then $X_{1/w}^{-\theta}$ is a semi-Hilbert space.
\end{corollary}

\begin{proof}
If $\left\{  u_{k}\right\}  $ is Cauchy in $X_{1/w}^{-\theta}$ then $\left\{
\mathcal{M}u_{k}\right\}  $ is Cauchy in the Hilbert space $L^{2}$. Thus
$\mathcal{M}u_{k}\rightarrow f$ for some $f\in L^{2}$ and so there is a
sequence $\left\{  w_{k}\right\}  $ in $\left(  S_{\mathcal{B}}^{\prime
}\right)  ^{\vee}$ such that $\mathcal{L}_{2}\mathcal{M}u_{k}=u_{k}%
+w_{k}\rightarrow\mathcal{L}_{2}f\in X_{1/w}^{-\theta}$ and so $u_{k}%
\rightarrow\mathcal{L}_{2}f$.
\end{proof}

\section{The general case: the space $\widetilde{X}_{1/w}^{-\theta}$}

\section{Overview}

For the cases of the shifted thin-plate functions and the Gaussian, $\left(
\int w\left\vert \phi\right\vert ^{2}\right)  ^{1/2}\notin W_{S;\theta}$ (???
Prove this: see Subsection 9.3.6 of the zero order document Williams
\cite{WilliamsZeroOrdSmthV4}).

To handle this case we try defining the space $X_{1/w}^{\theta}$ using the
space $\left(  \widehat{S}_{w,\theta}\right)  ^{\prime}$ instead of
$S^{\prime}$, where $\left(  \widehat{S}_{w,\theta}\right)  ^{\prime}$ will
denote the continuous linear functionals on the space $\widehat{S}_{w,\theta}$
of Fourier transforms of functions in $S_{w,\theta}$. The space $S_{w,\theta}$
was introduced in Definition \ref{Def_Sw,th_and_invF(Sw,th)}. Here
$S_{w,\theta}$ will be considered not as a subspace of $S$ but as a space in
itself and will be endowed with the topology using the $S$ seminorms
\textbf{and} the norm $\int w\left\vert \cdot\right\vert ^{2\theta}\left\vert
\phi\right\vert ^{2}$. We will call the modified space $\widetilde{X}%
_{1/w}^{\theta}$ so that%
\[
\widetilde{X}_{1/w}^{\theta}=\left\{  u\in\left(  \widehat{S}_{w,\theta
}\right)  ^{\prime}:\widehat{u}\in L_{loc}^{1}\left(  \mathbb{R}^{d}%
\setminus\mathcal{B}\right)  \text{ }and\text{ }\int\frac{\left\vert
u_{F}\right\vert ^{2}}{w\left\vert \cdot\right\vert ^{2\theta}}<\infty
\right\}  ,
\]

where $u_{F}=\widehat{u}$ on $\mathbb{R}^{d}\setminus\mathcal{B}$ and
$u_{F}\left(  \mathcal{B}\right)  =\left\{  0\right\}  $.

The space $\overset{\wedge}{S}_{w,0}$will be endowed with the topology that
makes the inverse-Fourier transform a homeomorphism.

Now if $g\in L^{2}$ we have $\mathcal{L}g=\left(  \sqrt{w}\left\vert
\cdot\right\vert ^{\theta}g\right)  ^{\vee}\in\left(  \widehat{S}_{w,\theta
}\right)  ^{\prime}$ and so it must be shown that $\mathcal{L}g\in L_{loc}%
^{1}\left(  \mathbb{R}^{d}\setminus\mathcal{A}\right)  $.\ But if
$K\subset\mathbb{R}^{d}\setminus\mathcal{B}$ is compact then $\sqrt
{w}\left\vert \cdot\right\vert ^{\theta}\in C^{\left(  0\right)  }\left(
\mathbb{R}^{d}\setminus\mathcal{B}\right)  $ and
\[
\int\limits_{K}\sqrt{w}\left\vert \cdot\right\vert ^{\theta}\left\vert
\widehat{g}\right\vert \leq\left\Vert \sqrt{w}\left\vert \cdot\right\vert
^{\theta}\right\Vert _{\infty}\int\limits_{K}\left\vert \widehat{g}\right\vert
\leq\max_{K}\left(  \sqrt{w}\left\vert \cdot\right\vert ^{\theta}\right)
\left\Vert g\right\Vert _{2}\int\limits_{K}1,
\]

proving $\sqrt{w}\left\vert \cdot\right\vert ^{\theta}\widehat{g}\in
L_{loc}^{1}\left(  \mathbb{R}^{d}\setminus\mathcal{B}\right)  $.

But what implications does this have for the `old' operator $\mathcal{M}%
:X_{1/w}^{-\theta}\rightarrow L^{2}$? It is the `same' operator i.e. if
$u\in\widetilde{X}_{1/w}^{-\theta}$ then $\mathcal{M}u=\frac{u_{F}}{\sqrt
{w}\left\vert \cdot\right\vert ^{\theta}}$. Clearly it isometric since
$\left\Vert \mathcal{M}u\right\Vert _{2}=\left\vert u\right\vert _{1/w,\theta
}$.

The space $S_{w,\theta}$ was introduced in Definition
\ref{Def_Sw,th_and_invF(Sw,th)}.

\begin{definition}
\label{Def_Swtheta}\textbf{The spaces }$\widehat{S}_{w,\theta}$,
$\overset{\vee}{S}_{w,\theta}$\textbf{and the spaces of functionals}
$S_{w,\theta}^{\prime}$, $\left(  \widehat{S}_{w,\theta}\right)  ^{\prime}$
\textbf{and} $\left(  \overset{\vee}{S}_{w,\theta}\right)  ^{\prime}$.

The space $S_{w,\theta}\subset S$ was introduced in Definition
\ref{Def_Sw,th_and_invF(Sw,th)} where it had the subspace topology induced by
$S$. Now we\textbf{\ }will topologize\textbf{\ }$S_{w,\theta}$\textbf{\ }using
the usual seminorms used to topologize\textbf{\ }$S$\textbf{\ }as well as the
seminorm used to define\textbf{\ }$S_{w,\theta}$, namely $\left(  \int
w\left\vert \cdot\right\vert ^{2\theta}\left\vert \phi\right\vert ^{2}\right)
^{1/2}$\textbf{.} Thus $f\in S_{w,\theta}^{\prime}$ iff $f$ is a linear
functional on $S_{w,\theta}$ and $\left\vert \left[  f,\phi\right]
\right\vert $ is bounded by a finite, positive linear combination of the norm
$\left(  \int w\left\vert \cdot\right\vert ^{2\theta}\left\vert \phi
\right\vert ^{2}\right)  ^{1/2}$ and the seminorms which define the topology
of $S$. Of course the linear combination is independent of $\phi$.

Next, the spaces $\widehat{S}_{w,\theta}$ and $\overset{\vee}{S}_{w,\theta}$
are defined by
\[
\widehat{S}_{w,\theta}=\left(  S_{w,\theta}\right)  ^{\wedge};\quad
\overset{\vee}{S}_{w,\theta}=\left(  S_{w,\theta}\right)  ^{\vee}.
\]

$\widehat{S}_{w,\theta}$ is topologized by composing the seminorms which
define the topology of $S_{w,\theta}$ with the Fourier transform.
$\overset{\vee}{S}_{w,\theta}$ is topologized by composing the seminorms which
define the topology of $S_{w,\theta}$ with the inverse-Fourier transform.

With these topologies the Fourier transform is now a homeomorphism from
$S_{w,\theta}$ to $\widehat{S}_{w,\theta}$ and the inverse-Fourier transform
is now a homeomorphism from $S_{w,\theta}$ to $\overset{\vee}{S}_{w,\theta}$.

Now define $S_{w,\theta}^{\prime}$, $\left(  \widehat{S}_{w,\theta}\right)
^{\prime}$ and $\left(  \overset{\vee}{S}_{w,\theta}\right)  ^{\prime}$ to be
the spaces of continuous, linear functionals on $S_{w,\theta}$, $\widehat
{S}_{w,\theta}$ and $\overset{\vee}{S}_{w,\theta}$ respectively.

We define the Fourier transform on $S_{w,\theta}^{\prime}$ and inverse-Fourier
transform on $S_{w,\theta}^{\prime}$ by%
\begin{align*}
\left[  \widehat{f},\phi\right]   & =\left[  f,\widehat{\phi}\right]  ,\quad
f\in S_{w,\theta}^{\prime},\phi\in\overset{\vee}{S}_{w,\theta},\\
\left[  \overset{\vee}{f},\phi\right]   & =\left[  f,\overset{\vee}{\phi
}\right]  ,\quad f\in S_{w,\theta}^{\prime},\phi\in\widehat{S}_{w,\theta},
\end{align*}

so that $\left(  S_{w,\theta}^{\prime}\right)  ^{\vee}=\left(  \widehat
{S}_{w,\theta}\right)  ^{\prime}$ and $\left(  S_{w,\theta}^{\prime}\right)
^{\wedge}=\left(  \overset{\vee}{S}_{w,\theta}\right)  ^{\prime}$. These
definitions imply that the Fourier transform is a homeomorphism from
$S_{w,\theta}^{\prime}$ to $\left(  \overset{\vee}{S}_{w,\theta}\right)
^{\prime}$, and that the inverse-Fourier transform is a homeomorphism from
$S_{w,\theta}^{\prime}$ to $\left(  \overset{\wedge}{S}_{w,\theta}\right)
^{\prime}$.
\end{definition}

We now `officially' define our more general version of $X_{1/w}^{-\theta}$
where $S^{\prime}$ is replaced by $\left(  \widehat{S}_{w,\theta}\right)
^{\prime}$.

\begin{definition}
\label{Def_X1/w,-th_II}\textbf{The space} $\widetilde{X}_{1/w}^{-\theta}$,
$\theta=1,2,3,\ldots$.

Suppose the weight function $w$ satisfies property W1. Then for each integer
$\theta\geq1$%
\[
\widetilde{X}_{1/w}^{-\theta}=\left\{  u\in\left(  \widehat{S}_{w,\theta
}\right)  ^{\prime}:\widehat{u}\in L_{loc}^{1}\left(  \mathbb{R}^{d}%
\setminus\mathcal{B}\right)  \text{ }and\text{ }\int\frac{\left\vert
u_{F}\right\vert ^{2}}{w\left\vert \cdot\right\vert ^{2\theta}}<\infty
\right\}  ,
\]

where $u_{F}=\widehat{u}$ on $\mathbb{R}^{d}\setminus\mathcal{B}$ and
$u_{F}\left(  \mathcal{B}\right)  =\left\{  0\right\}  $.

Note that $u\in\left(  \widehat{S}_{w,\theta}\right)  ^{\prime}$ implies
$\widehat{u}\in S_{w,\theta}^{\prime}$ and so $\widehat{u}\in L_{loc}%
^{1}\left(  \mathbb{R}^{d}\setminus\mathcal{B}\right)  $ means there exists
$f\in L_{loc}^{1}\left(  \mathbb{R}^{d}\setminus\mathcal{B}\right)  $ such
that $\left[  \widehat{u},\phi\right]  =\int f\phi$ for all $\phi\in
S_{w,\theta}$.

Also observe that Lemma \ref{Lem_Co,inf(Rd/C)} implies $\left(  C_{0}^{\infty
}\left(  \mathbb{R}^{d}\setminus\mathcal{B}\right)  \right)  ^{\vee}%
\subset\widetilde{X}_{1/w}^{-\theta}$ and so $\widetilde{X}_{1/w}^{-\theta}$
is not empty.

Endow $\widetilde{X}_{1/w}^{-\theta}$ with semi-inner product $\left\langle
u,v\right\rangle _{1/w,-\theta}=\int\frac{u_{F}\overline{v_{F}}}{w\left\vert
\cdot\right\vert ^{2\theta}}$ and seminorm $\left\vert u\right\vert
_{1/w,-\theta}=\left(  \int\frac{\left\vert u_{F}\right\vert ^{2}}{w\left\vert
\cdot\right\vert ^{2\theta}}\right)  ^{\frac{1}{2}}$.

Clearly%
\[
\operatorname*{null}\left\vert \cdot\right\vert _{1/w,-\theta}=\left(
\widetilde{S}_{\mathcal{B}}^{\prime}\right)  ^{\vee},
\]

where we have defined $\widetilde{S}_{\mathcal{B}}^{\prime}$ by%
\begin{equation}
\widetilde{S}_{\mathcal{B}}^{\prime}=\left\{  v\in S_{w,\theta}^{\prime
}:\operatorname*{supp}v\subset\mathcal{B}\right\}  .\label{p1.042}%
\end{equation}

\end{definition}

\begin{theorem}
$\widetilde{X}_{1/w}^{-\theta}=X_{1/w}^{-\theta}$ if inequality \ref{p1.043}
holds i.e. $w\left\vert \cdot\right\vert ^{2\theta}\in S_{w,\theta}^{\prime}$
when $S_{w,\theta}$ has the subspace topology induced by $S$.
\end{theorem}

\begin{proof}
Suppose $w\left\vert \cdot\right\vert ^{2\theta}\in S_{w,\theta}^{\prime} $
when $S_{w,\theta}$ has the subspace topology induced by $S$. Then the
topology on $S_{w,\theta}$ used to define $\widetilde{X}_{1/w}^{-\theta}$ is
identical to the topology on $S_{w,\theta}$ used to define $X_{1/w}^{-\theta}%
$. Thus $\widetilde{X}_{1/w}^{-\theta}$ and $X_{1/w}^{-\theta}$ are identical.
\end{proof}

\section{The operators $\widetilde{\mathcal{L}}$ and $\widetilde{\mathcal{M}}
$}

We now construct the equivalent of the mappings $\mathcal{L}$ and
$\mathcal{L}_{2}$.

\begin{definition}
\textbf{The linear operator} $\widetilde{\mathcal{L}}:L^{2}\rightarrow
S_{w,\theta}^{\prime}$, $\theta\geq1$.

To be consistent with previous definitions for $g\in L^{2}$ we try
$\widetilde{\mathcal{L}}g=\left(  \sqrt{w}\left\vert \cdot\right\vert
^{\theta}\widehat{g}\right)  ^{\vee}\in\left(  \widehat{S}_{w,\theta}\right)
^{\prime}$ i.e. $\sqrt{w}\left\vert \cdot\right\vert ^{\theta}\widehat{g}\in
S_{w,\theta}^{\prime}$. So suppose $\phi\in S_{w,\theta}$. Then%
\[
\left\vert \int\sqrt{w}\left\vert \cdot\right\vert ^{\theta}\widehat{g}%
\phi\right\vert \leq\int\left(  \sqrt{w}\left\vert \cdot\right\vert ^{\theta
}\left\vert \phi\right\vert \right)  \left\vert \widehat{g}\right\vert
\leq\left(  \int w\left\vert \cdot\right\vert ^{2\theta}\left\vert
\phi\right\vert ^{2}\right)  ^{1/2}\left\Vert g\right\Vert _{2},
\]

and the topology of $S_{w,\theta}^{\prime}$ means $\sqrt{w}\left\vert
\cdot\right\vert ^{\theta}\widehat{g}\in S_{w,\theta}^{\prime}$.

So we define $\widetilde{\mathcal{L}}$ by%
\[
\widetilde{\mathcal{L}}g=\left(  \sqrt{w}\left\vert \cdot\right\vert ^{\theta
}\widehat{g}\right)  ^{\vee},\quad g\in L^{2}.
\]

\end{definition}

We now prove some properties of $\widetilde{\mathcal{L}}$ which relate to the
space $\widetilde{X}_{1/w}^{-\theta}$. This theorem accounts for the condition
$\widehat{u}\in L_{loc}^{1}\left(  \mathbb{R}^{d}\setminus\mathcal{B}\right)
$ in the definition of $X_{1/w}^{-\theta}$. Here we need to prove that
$\widetilde{\mathcal{L}}g\in\widetilde{X}_{1/w}^{-\theta}$ and this requires
that $\left(  \widetilde{\mathcal{L}}g\right)  ^{\wedge}\in L_{loc}^{1}\left(
\mathbb{R}^{d}\setminus\mathcal{B}\right)  $ instead of the possibly stronger
$\left(  \widetilde{\mathcal{L}}g\right)  ^{\wedge}\in L_{loc}^{1}\left(
\mathbb{R}^{d}\setminus0\right)  $.

\begin{theorem}
\label{Thm_Lth_III} We have $\widetilde{\mathcal{L}}:L^{2}\rightarrow
\widetilde{X}_{1/w}^{-\theta}$, and $\widetilde{\mathcal{L}}$ is a linear isometry.

Also, $\left(  \widetilde{\mathcal{L}}g\right)  _{F}=\sqrt{w}\left\vert
\cdot\right\vert ^{\theta}\widehat{g}$.
\end{theorem}

\begin{proof}
From the definition of $\widetilde{\mathcal{L}}$, $\left(  \widetilde
{\mathcal{L}}g\right)  ^{\wedge}=\sqrt{w}\left\vert \cdot\right\vert ^{\theta
}\widehat{g}\in S_{w,\theta}^{\prime}$. Next observe that because $w\left\vert
\cdot\right\vert ^{2\theta}\in C^{\left(  0\right)  }\left(  \mathbb{R}%
^{d}\setminus\mathcal{B}\right)  $ we have $C_{0}^{\infty}\left(
\mathbb{R}^{d}\setminus\mathcal{B}\right)  \subset S_{w,\theta}$ and hence
that
\[
\left[  \left(  \widetilde{\mathcal{L}}g\right)  ^{\wedge},\phi\right]
=\left[  \sqrt{w}\left\vert \cdot\right\vert ^{\theta}\widehat{g},\phi\right]
,\quad\phi\in C_{0}^{\infty}\left(  \mathbb{R}^{d}\setminus\mathcal{B}\right)
.
\]

If we can show $\sqrt{w}\left\vert \cdot\right\vert ^{\theta}g\in L_{loc}%
^{1}\left(  \mathbb{R}^{d}\setminus\mathcal{B}\right)  $ it follows that
$\left(  \widetilde{\mathcal{L}}g\right)  ^{\wedge}\in L_{loc}^{1}\left(
\mathbb{R}^{d}\setminus\mathcal{B}\right)  $. Indeed, if $K\subset
\mathbb{R}^{d}\setminus\mathcal{B}$ is compact, $w\left\vert \cdot\right\vert
^{2\theta}\in C^{\left(  0\right)  }\left(  \mathbb{R}^{d}\setminus
\mathcal{B}\right)  $ implies
\[
\int_{K}\sqrt{w}\left\vert \cdot\right\vert ^{\theta}\left\vert g\right\vert
\leq\left(  \int_{K}w\left\vert \cdot\right\vert ^{2\theta}\right)
^{1/2}\left\Vert g\right\Vert _{2}\leq\max_{K}\left(  w\left\vert
\cdot\right\vert ^{2\theta}\right)  \left(  \int_{K}1\right)  \left\Vert
g\right\Vert _{2}.
\]

Thus $\left(  \widetilde{\mathcal{L}}g\right)  ^{\wedge}\in L_{loc}^{1}\left(
\mathbb{R}^{d}\setminus\mathcal{B}\right)  $ and $\left(  \widetilde
{\mathcal{L}}g\right)  _{F}=\sqrt{w}\left\vert \cdot\right\vert ^{\theta}g$ a.e.

Finally, $\left\vert \widetilde{\mathcal{L}}g\right\vert _{1/w,-\theta
}=\left(  \int\frac{\left\vert \left(  \widetilde{\mathcal{L}}g\right)
_{F}\right\vert ^{2}}{w\left\vert \cdot\right\vert ^{2\theta}}\right)
^{1/2}=\left\Vert g\right\Vert _{2}$.
\end{proof}

Our definition of the operator $\widetilde{\mathcal{M}}:\widetilde{X}%
_{1/w}^{-\theta}\rightarrow L^{2}$ will formally the same as $\mathcal{M}$ of
Definition \ref{Def_op_V}.

\begin{definition}
\label{Def_Mtilda}\textbf{The operator} $\widetilde{\mathcal{M}}%
:X_{1/w}^{-\theta}\rightarrow L^{2}$

Suppose the weight function $w$ has property W1. From the definition of
$\widetilde{X}_{1/w}^{-\theta}$, $u\in\widetilde{X}_{1/w}^{-\theta}$ implies
$\frac{u_{F}}{\sqrt{w}\left\vert \cdot\right\vert ^{\theta}}\in L^{2}$. We can
now define the linear mapping $\widetilde{\mathcal{M}}:\widetilde{X}%
_{1/w}^{-\theta}\rightarrow L^{2}$ by%
\[
\widetilde{\mathcal{M}}u=\left(  \frac{u_{F}}{\sqrt{w}\left\vert
\cdot\right\vert ^{\theta}}\right)  ^{\vee},\text{\quad}u\in\widetilde
{X}_{1/w}^{-\theta}.
\]

\end{definition}

$\widetilde{\mathcal{M}}$ has the following properties:

\begin{theorem}
The operator $\widetilde{\mathcal{M}}:\widetilde{X}_{1/w}^{-\theta}\rightarrow
L^{2}$ is linear and isometric with null space $\left(  \widetilde
{S}_{\mathcal{B}}^{\prime}\right)  ^{\vee}$.
\end{theorem}

\begin{proof}
That $\widetilde{\mathcal{M}}$ is an isometry is clear from the definition of
$\widetilde{X}_{1/w}^{-\theta}$. Since $\widetilde{\mathcal{M}}$ is an
isometry the null space of $\widetilde{\mathcal{M}}$ is the null space of the
seminorm $\left\vert \cdot\right\vert _{1/w,\theta}$, namely $\left(
\widetilde{S}_{\mathcal{B}}^{\prime}\right)  ^{\vee}$.
\end{proof}

The following theorem indicates how the operators $\widetilde{\mathcal{L}}$
and $\widetilde{\mathcal{M}}$ interact.

\begin{theorem}
\label{Thm_Lth_Mth_III_links}Suppose the weight function $w$ has the
properties assumed in Definition \ref{Def_Lth_II}. Then the operators
$\widetilde{\mathcal{L}}:L^{2}\rightarrow\widetilde{X}_{1/w}^{-\theta}$ and
$\widetilde{\mathcal{M}}:\widetilde{X}_{1/w}^{-\theta}\rightarrow L^{2}$ satisfy:

\begin{enumerate}
\item $\widetilde{\mathcal{M}}\widetilde{\mathcal{L}}=I$ on $L^{2}$;
$\widetilde{\mathcal{L}}\widetilde{\mathcal{M}}u-u\in\left(  \widetilde
{S}_{\mathcal{B}}^{\prime}\right)  ^{\vee}$ when $u\in\widetilde{X}%
_{1/w}^{-\theta}$.

Thus $\widetilde{\mathcal{M}}$ and $\widetilde{\mathcal{L}}$ are inverses in
the seminorm sense.

\item $\widetilde{\mathcal{L}}$ is 1-1. $\widetilde{\mathcal{L}}$ is also onto
in the seminorm sense.

\item $\widetilde{\mathcal{M}}$ is onto. It is also 1-1 in the seminorm sense.

\item $\widetilde{\mathcal{M}}$ and $\widetilde{\mathcal{L}}$ are adjoints in
the sense that $\left\langle \widetilde{\mathcal{L}}g,f\right\rangle
_{1/w,-\theta}=\left(  g,\widetilde{\mathcal{M}}f\right)  _{2}$.
\end{enumerate}
\end{theorem}

\begin{proof}
\textbf{Property 1}. From the definition of $\widetilde{\mathcal{L}}$,
$\left(  \widetilde{\mathcal{L}}g\right)  _{F}=\sqrt{w}\left\vert
\cdot\right\vert ^{\theta}\widehat{g}$. Thus, from the definition of
$\widetilde{\mathcal{M}}$ and Theorem \ref{Thm_Lth_III}, for $g\in L^{2}$
\[
\left(  \widetilde{\mathcal{M}}\widetilde{\mathcal{L}}g\right)  ^{\wedge
}=\frac{\left(  \widetilde{\mathcal{L}}g\right)  _{F}}{\sqrt{w}\left\vert
\cdot\right\vert ^{\theta}}=\frac{\sqrt{w}\left\vert \cdot\right\vert
^{\theta}\widehat{g}}{\sqrt{w}\left\vert \cdot\right\vert ^{\theta}}%
=\widehat{g},
\]

and thus $\widetilde{\mathcal{M}}\widetilde{\mathcal{L}}=I$.

Next we show $\widetilde{\mathcal{L}}\widetilde{\mathcal{M}}u-u\in\left(
\widetilde{S}_{\mathcal{B}}^{\prime}\right)  ^{\vee}$ on $\widetilde{X}%
_{1/w}^{-\theta}$.

In fact by Theorem \ref{Thm_Lth_III}, $u\in\widetilde{X}_{1/w}^{-\theta}$
implies $\left(  \widetilde{\mathcal{L}}\widetilde{\mathcal{M}}u\right)
_{F}=\sqrt{w}\left\vert \cdot\right\vert ^{\theta}\left(  \widetilde
{\mathcal{M}}u\right)  ^{\wedge}=u_{F}$ so that
\[
\left\vert \widetilde{\mathcal{L}}\widetilde{\mathcal{M}}u-u\right\vert
_{1/w,-\theta}^{2}=\int w\left\vert \cdot\right\vert ^{2\theta}\left\vert
\left(  \widetilde{\mathcal{L}}\widetilde{\mathcal{M}}u-u\right)
_{F}\right\vert ^{2}=0,
\]

so that $\widetilde{\mathcal{L}}\widetilde{\mathcal{M}}u-u\in\left(
\widetilde{S}_{\mathcal{B}}^{\prime}\right)  ^{\vee}$ i.e. $\widetilde
{\mathcal{L}}\widetilde{\mathcal{M}}u-u\in\left(  \widetilde{S}_{\mathcal{B}%
}^{\prime}\right)  ^{\vee}$, where $S_{\mathcal{B}}^{\prime}$ was defined in
\ref{p1.042}.

Since $\left(  \widetilde{S}_{\mathcal{B}}^{\prime}\right)  ^{\vee}$ is the
null space of the seminorm on $\widetilde{X}_{1/w}^{-\theta}$, $\widetilde
{\mathcal{M}}$ and $\widetilde{\mathcal{L}}$ are inverses in the seminorm
sense.\medskip

\textbf{Property 2}. That $\widetilde{\mathcal{L}}$ is 1-1 follows from
$\widetilde{\mathcal{M}}\widetilde{\mathcal{L}}=I$. That $\widetilde
{\mathcal{L}}$ is onto follows from $\widetilde{\mathcal{L}}\widetilde
{\mathcal{M}}u-u\in\left(  \widetilde{S}_{\mathcal{B}}^{\prime}\right)
^{\vee}$ when $u\in\widetilde{X}_{1/w}^{-\theta}$.\medskip

\textbf{Property 3}. That $\widetilde{\mathcal{M}}$ is onto follows from
$\widetilde{\mathcal{M}}\widetilde{\mathcal{L}}=I$. That $\widetilde
{\mathcal{L}}$ is onto follows from $\widetilde{\mathcal{L}}\widetilde
{\mathcal{M}}u-u\in\left(  \widetilde{S}_{\mathcal{B}}^{\prime}\right)
^{\vee}$ when $u\in\widetilde{X}_{1/w}^{-\theta}$.\medskip

\textbf{Property 4}. If $g\in L^{2}$ and $f\in\widetilde{X}_{1/w}^{-\theta}
$,
\begin{align*}
\left\langle \widetilde{\mathcal{L}}g,f\right\rangle _{1/w,-\theta}=\int%
\frac{\left(  \widetilde{\mathcal{L}}g\right)  _{F}\overline{f_{F}}%
}{w\left\vert \cdot\right\vert ^{2\theta}}=\int\frac{\sqrt{w}\left\vert
\cdot\right\vert ^{\theta}\widehat{g}\overline{f_{F}}}{w\left\vert
\cdot\right\vert ^{2\theta}}=\int\widehat{g}\frac{\overline{f_{F}}}{\sqrt
{w}\left\vert \cdot\right\vert ^{\theta}}  & =\left(  \widehat{g}%
,\overline{\left(  \widetilde{\mathcal{M}}f\right)  ^{\wedge}}\right)  _{2}\\
& =\left(  g,\overline{\widetilde{\mathcal{M}}f}\right)  _{2}.
\end{align*}

\end{proof}

Since $L^{2}$ is complete, the mappings of the previous theorem will yield the
following important result.

\begin{corollary}
\label{Cor_X1/w,-th_semiHilb_II}Suppose the weight function $w$ has the
properties assumed in Definition \ref{Def_Lth_II} for some order $\theta\geq
1$. Then $\widetilde{X}_{1/w}^{-\theta}$ is a semi-Hilbert space.
\end{corollary}

\begin{proof}
If $\left\{  u_{k}\right\}  $ is Cauchy in $\widetilde{X}_{1/w}^{-\theta}$
then $\left\{  \widetilde{\mathcal{M}}u_{k}\right\}  $ is Cauchy in $L^{2}$
since $\widetilde{\mathcal{M}}$ is isometric. Thus $\widetilde{\mathcal{M}%
}u_{k}\rightarrow u$ for some $u\in L^{2}$ and since $\widetilde{\mathcal{M}%
}\widetilde{\mathcal{L}}=I$ by Theorem \ref{Thm_Lth_Mth_III_links}%
\[
\left\vert u_{k}-\widetilde{\mathcal{L}}u\right\vert _{1/w,-\theta}=\left\Vert
\widetilde{\mathcal{M}}\left(  u_{k}-\widetilde{\mathcal{L}}u\right)
\right\Vert _{2}=\left\Vert \widetilde{\mathcal{M}}u_{k}-u\right\Vert
_{2}\rightarrow0.
\]

\end{proof}

\section{The operator $\widetilde{\mathcal{V}}:\widetilde{X}_{1/w}^{-\theta
}\rightarrow X_{w}^{\theta}$}

In this section the operator $\widetilde{\mathcal{V}}=\mathcal{J}%
\widetilde{\mathcal{M}}$ will be studied and used to prove various properties
of $\widetilde{X}_{1/w}^{-\theta},$ including that $\widetilde{X}%
_{1/w}^{-\theta}$ is a semi-Hilbert space and that $\widehat{S}_{\emptyset
,\theta}$ is dense in $\widetilde{X}_{1/w}^{-\theta}$.

\begin{definition}
\textbf{The operator} $\widetilde{\mathcal{V}}:\widetilde{X}_{1/w}^{-\theta
}\rightarrow X_{w}^{\theta}$.

Suppose the weight function $w$ has properties W1 and W2. These conditions
ensure that both $\widetilde{\mathcal{M}}$ and $\mathcal{J}$ are defined.

The operator $\widetilde{\mathcal{V}}$ is now defined by $\widetilde
{\mathcal{V}}=\mathcal{J}\widetilde{\mathcal{M}}$.
\end{definition}

Amongst other things the next result shows $\widetilde{\mathcal{V}}$ and
$\widetilde{\mathcal{L}}\mathcal{I}$ are inverses in the seminorm sense.

\begin{theorem}
\label{Thm_prop_op_V_inverse}Suppose the weight function $w$ has properties W1
and W2 for order $\theta\geq0$. In the \textbf{seminorm sense}:

\begin{enumerate}
\item $\widetilde{\mathcal{V}}:\widetilde{X}_{1/w}^{-\theta}\rightarrow
X_{w}^{\theta}$ is an isometry but is in general non-linear.

\item $\widetilde{\mathcal{L}}\mathcal{I}:X_{w}^{\theta}\rightarrow
\widetilde{X}_{1/w}^{-\theta}$ is a linear isometry.

\item $\widetilde{\mathcal{V}}$ and $\widetilde{\mathcal{L}}\mathcal{I}$ are inverses.

\item $\widetilde{\mathcal{V}}$ and $\widetilde{\mathcal{L}}\mathcal{I}$ are onto.

\item $\widetilde{\mathcal{V}}$ and $\widetilde{\mathcal{L}}\mathcal{I}$ are 1-1.

\item $\widetilde{\mathcal{V}}$ and $\widetilde{\mathcal{L}}\mathcal{I}$ are adjoints.

\item If $f\in X_{w}^{\theta}$ then $\widetilde{\mathcal{L}}\mathcal{I}%
f=\left(  w\left\vert \cdot\right\vert ^{2\theta}f_{F}\right)  ^{\vee}$, where
$f_{F}=\widehat{f}$ on $\mathbb{R}^{d}\setminus0$.
\end{enumerate}
\end{theorem}

\begin{proof}
\textbf{Parts 1 and 2}. True since $\mathcal{I}$, $\mathcal{J}$,
$\widetilde{\mathcal{L}}$ and $\widetilde{\mathcal{M}}$ are isometric.\medskip

\textbf{Part 3}. From Theorem \ref{Thm_Lth_Mth_III_links} $\widetilde
{\mathcal{M}}\widetilde{\mathcal{L}}=I$\ and Theorem \ref{Thm_Itheta_Jtheta}
$\mathcal{JI}f-f\in P_{\theta}$. Hence $\widetilde{\mathcal{V}}\widetilde
{\mathcal{L}}\mathcal{I}=\mathcal{J}\widetilde{\mathcal{M}}\widetilde
{\mathcal{L}}\mathcal{I}=\mathcal{JI}$ and so $\widetilde{\mathcal{V}%
}\widetilde{\mathcal{L}}\mathcal{I}f-f\in P_{\theta}$.

From Theorem \ref{Thm_Itheta_Jtheta} $\mathcal{IJ}=I$ and from Theorem
\ref{Thm_Lth_Mth_III_links} $\widetilde{\mathcal{L}}\widetilde{\mathcal{M}%
}g-g\in\left(  \widetilde{S}_{\mathcal{B}}^{\prime}\right)  ^{\vee}$. Hence
$\widetilde{\mathcal{L}}\mathcal{I}\widetilde{\mathcal{V}}=\widetilde
{\mathcal{L}}\mathcal{IJ}\widetilde{\mathcal{M}}=\widetilde{\mathcal{L}%
}\widetilde{\mathcal{M}}$ and $\widetilde{\mathcal{L}}\mathcal{I}%
\widetilde{\mathcal{V}}g-g\in\left(  \widetilde{S}_{\mathcal{B}}^{\prime
}\right)  ^{\vee}$.\medskip

\textbf{Parts 4 and 5}. follow directly from part 3.\medskip

\textbf{Part 6}. Theorem \ref{Thm_Lth_Mth_III_links}\ and Theorem
\ref{Thm_Itheta_Jtheta} imply $\widetilde{\mathcal{L}}$, $\widetilde
{\mathcal{M}} $ are adjoints and $\mathcal{I}$, $\mathcal{J}$ are adjoints.
Thus
\[
\left\langle g,\widetilde{\mathcal{L}}\mathcal{I}f\right\rangle _{1/w,-\theta
}=\left(  \widetilde{\mathcal{M}}g,\mathcal{I}f\right)  _{2}=\left\langle
\mathcal{J}\widetilde{\mathcal{M}}g,f\right\rangle _{w,\theta}=\left\langle
\widetilde{\mathcal{V}}g,f\right\rangle _{w,\theta}.
\]
\medskip

\textbf{Part 7}. If $f\in X_{w}^{\theta}$ then
\[
\widehat{\widetilde{\mathcal{L}}\mathcal{I}f}=\sqrt{w}\left\vert
\cdot\right\vert ^{\theta}\widehat{\mathcal{I}f}=\sqrt{w}\left\vert
\cdot\right\vert ^{\theta}\left(  \sqrt{w}\left\vert \cdot\right\vert
^{\theta}f_{F}\right)  =w\left\vert \cdot\right\vert ^{2\theta}f_{F}.
\]

\end{proof}

The next result evaluates $\widetilde{\mathcal{V}}^{\ast}\left(  D^{\gamma
}R_{x}\right)  $ where $R_{x}$ is the Riesz representer of the evaluation
functional $f\rightarrow f\left(  x\right)  $ discussed in Section
\ref{Sect_Riesz_rep}.

\begin{theorem}
\label{Thm_V*(DRx)=delta}Suppose the weight function $w$ has properties W1, W2
and W3 for parameters $\theta$ and $\kappa$. Then
\[
\left(  2\pi\right)  ^{d/2}\widetilde{\mathcal{V}}^{\ast}\left(  D^{\gamma
}R_{x}\right)  =\left(  -1\right)  ^{\left\vert \gamma\right\vert }D^{\gamma
}\delta\left(  \cdot-x\right)  -\sum\limits_{i=1}^{M}D^{\gamma}l_{i}%
(x)\delta\left(  \cdot-a_{i}\right)  ,\quad\left\vert \gamma\right\vert
\leq\kappa.
\]

\end{theorem}

\begin{proof}
From the proof of Theorem \ref{Thm_DRx_in_Xw}%
\[
\left(  2\pi\right)  ^{d/2}\widehat{D^{\gamma}R_{x}}=\frac{D_{x}^{\gamma
}\mathcal{Q}_{x}\left(  e^{-i\left\langle x,\cdot\right\rangle }\right)
}{w\left\vert \cdot\right\vert ^{2\theta}}\in L_{loc}^{1}\left(
\mathbb{R}^{d}\setminus0\right)  .
\]

Thus by parts 6 and 7 of Theorem \ref{Thm_prop_op_V_inverse}, for fixed $x$%
\begin{align*}
\left(  2\pi\right)  ^{d/2}\widetilde{\mathcal{V}}^{\ast}\left(  D^{\gamma
}R_{x}\right)   & =\left(  D_{x}^{\gamma}\mathcal{Q}_{x}\left(
e^{-i\left\langle x,\xi\right\rangle }\right)  \right)  ^{\vee}\\
& =\left(  D_{x}^{\gamma}\left(  e^{-i\left\langle x,\cdot\right\rangle }%
-\sum\limits_{i=1}^{M}l_{i}(x)e^{-i\left\langle a_{i},\xi\right\rangle
}\right)  \right)  ^{\vee}\\
& =\left(  \left(  -i\xi\right)  ^{\gamma}e^{-i\left\langle x,\xi\right\rangle
}-\sum\limits_{i=1}^{M}\left(  D^{\gamma}l_{i}\right)  (x)e^{-i\left\langle
a_{i},\xi\right\rangle }\right)  ^{\vee}\\
& =\left(  -1\right)  ^{\left\vert \gamma\right\vert }D^{\gamma}\left(
e^{-i\left\langle x,\xi\right\rangle }\right)  ^{\vee}-\sum\limits_{i=1}%
^{M}\left(  D^{\gamma}l_{i}\right)  (x)\left(  e^{-i\left\langle a_{i}%
,\xi\right\rangle }\right)  ^{\vee}\\
& =\left(  -1\right)  ^{\left\vert \gamma\right\vert }D^{\gamma}\delta\left(
\cdot-x\right)  -\sum\limits_{i=1}^{M}D^{\gamma}l_{i}(x)\delta\left(
\cdot-a_{i}\right)  .
\end{align*}

\end{proof}

The next result shows that the operator $\widetilde{\mathcal{V}}$ behaves very
nicely on the subspace $\widehat{S}_{\emptyset,\theta}$.

\begin{theorem}
\label{Thm_Vth_and_F(So,th)}Suppose the weight function $w$ satisfies
properties W1, W2.

Then $\widetilde{\mathcal{V}}\phi=G\ast\phi$ when $\phi\in\widehat
{S_{\emptyset,\theta}}$ and so $\widetilde{\mathcal{V}}:\widehat{S}%
_{\emptyset,\theta}\rightarrow G\ast\widehat{S}_{\emptyset,\theta}$. Indeed,
the restriction of $\widetilde{\mathcal{V}}$ to $\widehat{S}_{\emptyset
,\theta}$ is 1-1 and onto $G\ast\widehat{S}_{\emptyset,\theta}$.
\end{theorem}

\begin{proof}
From the definitions of $\mathcal{I}$ and $\widetilde{\mathcal{L}}$
\begin{equation}
\widehat{\mathcal{I}\left(  G\ast\phi\right)  }=\sqrt{w}\left\vert
\cdot\right\vert ^{\theta}\left(  G\ast\phi\right)  _{F}=\frac{\widehat{\phi}%
}{\sqrt{w}\left\vert \cdot\right\vert ^{\theta}},\label{p1.081}%
\end{equation}

and%
\begin{align*}
\widetilde{\mathcal{L}}\mathcal{I}\left(  G\ast\phi\right)  =\left(  \sqrt
{w}\left\vert \cdot\right\vert ^{\theta}\left(  \mathcal{I}\left(  G\ast
\phi\right)  \right)  _{F}\right)  ^{\vee}  & =\left(  \sqrt{w}\left\vert
\cdot\right\vert ^{\theta}\left(  \frac{\widehat{\phi}}{\sqrt{w}\left\vert
\cdot\right\vert ^{\theta}}\right)  \right)  ^{\vee}\\
& =\phi.
\end{align*}

Since $\mathcal{J}\widetilde{\mathcal{M}}\widetilde{\mathcal{L}}\mathcal{I}=I
$ we have $G\ast\phi=\mathcal{J}\widetilde{\mathcal{M}}\phi=\widetilde
{\mathcal{V}}\phi$. Clearly $\widetilde{\mathcal{V}}$ is onto. Since
$\widetilde{\mathcal{V}}$ is isometric $\widetilde{\mathcal{V}}\phi=0$ implies
$\left\vert \widetilde{\mathcal{V}}\phi\right\vert _{w,\theta}=\left\vert
\phi\right\vert _{1/w,-\theta}=\left(  \int\frac{\left\vert \widehat{\phi
}\right\vert ^{2}}{w\left\vert \cdot\right\vert ^{2\theta}}\right)  ^{\vee}=0$
which implies $\phi=0$.
\end{proof}

We now prove a density result for $X_{1/w}^{-\theta}$.

\begin{corollary}
\label{Cor_F[So,th]_dense_in_X1/w,-th} Suppose the weight function $w$ has
properties W1 and W2 for order $\theta$.

Then $\widehat{S}_{\emptyset,\theta}$ is dense in $\widetilde{X}%
_{1/w}^{-\theta}$ and $\widetilde{X}_{1/w}^{-\theta}$ is the completion of
$\widehat{S}_{\emptyset,\theta}$.
\end{corollary}

\begin{proof}
Choose $u\in\widetilde{X}_{1/w}^{-\theta}$ so that $\widetilde{\mathcal{V}%
}u\in X_{w}^{\theta}$. Since by Theorem \ref{Thm_Jg_property} $G\ast
\widehat{S}_{\emptyset,\theta}$ is dense in $X_{w}^{\theta}$ in the seminorm
sense, given $\varepsilon>0$ there exists $\phi_{\varepsilon}\in G\ast
\widehat{S}_{\emptyset,\theta}$ such that $\left\vert \widetilde{\mathcal{V}%
}u-\phi_{\varepsilon}\right\vert _{w,\theta}<\varepsilon$. But $\widetilde
{\mathcal{V}}:\widehat{S}_{\emptyset,\theta}\rightarrow G\ast\widehat
{S}_{\emptyset,\theta}$ is an isomorphism so $\widetilde{\mathcal{V}}^{-1}%
\phi_{\varepsilon}\in\widehat{S}_{\emptyset,\theta}$. Since $\widetilde
{\mathcal{V}}:X_{1/w}^{-\theta}\rightarrow X_{w}^{\theta}$ is an isometry%
\[
\left\vert u-\widetilde{\mathcal{V}}^{-1}\phi_{\varepsilon}\right\vert
_{1/w,-\theta}=\left\vert \widetilde{\mathcal{V}}\left(  u-\widetilde
{\mathcal{V}}^{-1}\phi_{\varepsilon}\right)  \right\vert _{w,\theta
}=\left\vert \widetilde{\mathcal{V}}u-\phi_{\varepsilon}\right\vert
_{w,\theta}<\varepsilon,
\]

which proves this density result.

Since $\widetilde{X}_{1/w}^{-\theta}$ is complete and $\widehat{S}%
_{\emptyset,\theta}$ is dense in $\widetilde{X}_{1/w}^{-\theta}$ it follows
that $\widetilde{X}_{1/w}^{-\theta}$ is the completion of $\widehat
{S}_{\emptyset,\theta}$.
\end{proof}

\begin{theorem}
Suppose the weight function $w$ has properties W1 and W2 for order $\theta
\geq1$.

Denote by $\left(  X_{w}^{\theta}\right)  ^{\prime}$ the space of bounded
linear functionals on $X_{w}^{\theta}$.

Then there exists a isomorphic linear operator $\Phi:\widetilde{X}%
_{1/w}^{-\theta}\rightarrow\left(  X_{w}^{\theta}\right)  ^{\prime}$ defined
by%
\[
\left(  \Phi u\right)  \left(  g\right)  =\left\langle g,\widetilde
{\mathcal{V}}u\right\rangle _{w,\theta},\quad g\in X_{w}^{\theta}%
,u\in\widetilde{X}_{1/w}^{-\theta},
\]

which is 1-1 and onto in the seminorm sense.
\end{theorem}

\begin{proof}
Since $\widetilde{\mathcal{V}}:\widetilde{X}_{1/w}^{-\theta}\rightarrow
X_{w}^{\theta}$ is an isometry, given $u\in\widetilde{X}_{1/w}^{-\theta}$ and
$g\in X_{w}^{\theta}$ we have
\[
\left\vert \left\langle g,\widetilde{\mathcal{V}}u\right\rangle _{w,\theta
}\right\vert \leq\left\vert g\right\vert _{w,\theta}\left\vert \widetilde
{\mathcal{V}}u\right\vert _{w,\theta}=\left\vert g\right\vert _{w,\theta
}\left\Vert u\right\Vert _{1/w,-\theta},
\]

and for each $u\in\widetilde{X}_{1/w}^{-\theta}$, the expression $\left\langle
g,\widetilde{\mathcal{V}}u\right\rangle _{w,\theta}$, $g\in X_{w}^{\theta}$
defines a bounded linear functional on $X_{w}^{\theta}$. Denote this
functional by $\Phi u$ so that
\begin{equation}
\left(  \Phi u\right)  \left(  g\right)  =\left\langle g,\widetilde
{\mathcal{V}}u\right\rangle _{w,\theta},\quad g\in X_{w}^{\theta
},\label{p1.076}%
\end{equation}

and the operator norm is
\[
\left\Vert \Phi u\right\Vert _{op}=\sup\limits_{g\in X_{w}^{\theta}}%
\frac{\left\vert \left\langle g,\widetilde{\mathcal{V}}u\right\rangle
_{w,\theta}\right\vert }{\left\vert g\right\vert _{w,\theta}}\leq\left\vert
u\right\vert _{1/w,-\theta}.
\]

Thus $\Phi:\widetilde{X}_{1/w}^{-\theta}\rightarrow\left(  X_{w}^{\theta
}\right)  ^{\prime}$. In fact, we can easily prove that $\left\Vert \Phi
u\right\Vert _{op}=\left\vert u\right\vert _{1/w,-\theta}$ by noting that when
$g=\widetilde{\mathcal{V}}u$, $\frac{\left\vert \left\langle g,\widetilde
{\mathcal{V}}u\right\rangle _{w,\theta}\right\vert }{\left\vert g\right\vert
_{w,\theta}}=\left\vert \widetilde{\mathcal{V}}u\right\vert _{w,\theta
}=\left\vert u\right\vert _{1/w,-\theta}$. Clearly $\Phi u=0$ implies
$\left\vert u\right\vert _{1/w,-\theta}=0$ so $\Phi$ is 1-1.

To prove $\Phi$ is onto choose $\mathcal{Y}\in\left(  X_{w}^{\theta}\right)
^{\prime}$. Since $X_{w}^{\theta}$ is an inner product space there exists
$v\in X_{w}^{\theta}$ such that $\mathcal{Y}g=\left\langle g,v\right\rangle
_{w,\theta}$ when $g\in X_{w}^{\theta}$. But by Theorem
\ref{Thm_prop_op_V_inverse} $\mathcal{Y}g=\left\langle g,v\right\rangle
_{w,\theta}=\left\langle g,\widetilde{\mathcal{V}}\widetilde{\mathcal{V}%
}^{\ast}v\right\rangle _{w,\theta}$ and comparison with \ref{p1.076} yields
$\left(  \Phi\left(  \widetilde{\mathcal{V}}^{\ast}v\right)  \right)  \left(
g\right)  =\mathcal{Y}g$ and $\Phi\left(  \widetilde{\mathcal{V}}^{\ast
}v\right)  =\mathcal{Y}$.
\end{proof}

\section{The operators $\mathcal{I}$, $\mathcal{J}$, $\mathcal{J}_{\pm}$and
the spaces $\protect\overset{\wedge}{S}_{\emptyset,\theta}$ and $G\ast
\protect\overset{\wedge}{S}_{\emptyset,\theta}$}

The spaces $\widehat{S}_{\emptyset,\theta}$ and $G\ast\widehat{S}%
_{\emptyset,\theta}=\overset{\cdot}{J}_{G}$ were introduced and discussed in
Section \ref{Sect_Jg}. The operators $\mathcal{I}$ and the class of operators
$\mathcal{J}$ were introduced in Subsection \ref{SbSect_I_J} and it was shown
that $\mathcal{I}:X_{w}^{\theta}\rightarrow L^{2}$ and $\mathcal{J}%
:L^{2}\rightarrow X_{w}^{\theta}$ are isometries and adjoints in the seminorm
sense. It was also shown that $\mathcal{I}$ is linear and $\mathcal{J}$ is
linear modulo a polynomial in $P_{\theta}$. From Corollary
\ref{Cor_F[So,th]_dense_in_X1/w,-th} $\widehat{S}_{\emptyset,\theta}$ is dense
in $\widetilde{X}_{1/w}^{-\theta}$ and from Theorem \ref{Thm_Jg_property},
$G\ast\widehat{S}_{\emptyset,\theta}$ is dense in $X_{w}^{\theta}$ in the
seminorm sense.

We start by showing that $\left(  \frac{S_{\emptyset,\theta}}{\sqrt
{w}\left\vert \cdot\right\vert ^{\theta}}\right)  ^{\vee}\subset X_{w}%
^{\theta}\cap L^{2}$ and then restrict $\mathcal{I}$ to $\left(
\frac{S_{\emptyset,\theta}}{\sqrt{w}\left\vert \cdot\right\vert ^{\theta}%
}\right)  ^{\vee}$ and also to $G\ast\widehat{S}_{\emptyset,\theta}\subset
X_{w}^{\theta}$. Next, inverses $\mathcal{J}_{-}$ and $\mathcal{J}_{+}$ of
these restrictions are constructed which are restrictions of members of the
class $\mathcal{J}$.

We will then show that $\mathcal{I}=\widetilde{\mathcal{L}}$ on $\left(
\frac{S_{\emptyset,\theta}}{\sqrt{w}\left\vert \cdot\right\vert ^{\theta}%
}\right)  ^{\vee}$ and $\mathcal{J}_{-}=\widetilde{\mathcal{M}}$ on
$\widehat{S}_{\emptyset,\theta}$. Finally it is demonstrated that
$\mathcal{I}^{2}:G\ast\widehat{S}_{\emptyset,\theta}\rightarrow\widehat
{S}_{\emptyset,\theta}$ has inverse $\mathcal{J}_{+}\mathcal{J}_{-}$ and
$\mathcal{I}^{2}=\widetilde{\mathcal{V}}$. Figure \ref{Fig_I_Jplus_Jminus}
below summarizes these results.%
\begin{equation}%
\begin{array}
[c]{ccccc}%
\widetilde{X}_{1/w}^{-\theta} &
\begin{array}
[c]{c}%
\widetilde{\mathcal{L}}\\
\mathbf{\longleftarrow-}\\
-\longrightarrow\\
\widetilde{\mathcal{M}}%
\end{array}
& L^{2} &
\begin{array}
[c]{c}%
\mathcal{I}\\
\mathbf{\longleftarrow-}\\
-\longrightarrow\\
\mathcal{J}%
\end{array}
& X_{w}^{\theta}\\
\uparrow d &  & \uparrow d &  & \uparrow d\\
\widehat{S}_{\emptyset,\theta} &
\begin{array}
[c]{c}%
\\
\mathcal{I}\\
\mathbf{\longleftarrow-}\\
-\longrightarrow\\
\mathcal{J}_{-}\\
isomet,isom
\end{array}
& \left(  \frac{S_{\emptyset,\theta}}{\sqrt{w}\left\vert \cdot\right\vert
^{\theta}}\right)  ^{\vee} &
\begin{array}
[c]{c}%
\\
\mathcal{I}\\
\mathbf{\longleftarrow-}\\
-\longrightarrow\\
\mathcal{J}_{+}\\
isomet,isom
\end{array}
& G\ast\widehat{S}_{\emptyset,\theta}%
\end{array}
\label{Fig_I_Jplus_Jminus}%
\end{equation}

\begin{lemma}
\label{Lem_F(S/sqrt(w|.|^th))_in_Xthw}Suppose the weight function $w$ has
properties W1 and W2 for order $\theta\geq1$. Then
\[
\left(  \frac{S_{\emptyset,\theta}}{\sqrt{w}\left\vert \cdot\right\vert
^{\theta}}\right)  ^{\vee}\subset X_{w}^{\theta}\cap L^{2},
\]

and%
\begin{equation}
\left(  \mathcal{I}\left(  G\ast\phi\right)  \right)  ^{\wedge}=\frac
{\widehat{\phi}}{\sqrt{w}\left\vert \cdot\right\vert ^{\theta}},\text{\quad
}\phi\in\widehat{S}_{\emptyset,\theta}.\label{p1.7}%
\end{equation}

\end{lemma}

\begin{proof}
If $\psi\in S_{\emptyset,\theta}$ then $\left\vert \psi\right\vert ^{2}\in
S_{\emptyset,2\theta}$ and so by Theorem \ref{Thm_1/w|.|^2u_in_So,2u'},%
\[
\int\left\vert \frac{\psi}{\sqrt{w}\left\vert \cdot\right\vert ^{\theta}%
}\right\vert ^{2}=\int\frac{\left\vert \psi\right\vert ^{2}}{w\left\vert
\cdot\right\vert ^{2\theta}}<\infty.
\]

Thus $\frac{\psi}{\sqrt{w}\left\vert \cdot\right\vert ^{\theta}}\in
L^{2}\subset S^{\prime}$ and hence $\left(  \frac{S_{\emptyset,\theta}}%
{\sqrt{w}\left\vert \cdot\right\vert ^{\theta}}\right)  ^{\vee}\in
L^{2}\subset S^{\prime}$.

To show that $\left(  \frac{S_{\emptyset,\theta}}{\sqrt{w}\left\vert
\cdot\right\vert ^{\theta}}\right)  ^{\vee}\subset X_{w}^{\theta}$ we will use
the definition of $X_{w}^{\theta}$ given by equation \ref{p40} i.e.%
\[
X_{w}^{\theta}=\left\{  g\in S^{\prime}:\xi^{\alpha}\widehat{g}\in L_{loc}%
^{1}\text{ }if\text{ }\left\vert \alpha\right\vert =\theta;\text{\thinspace
}\int w\left\vert \cdot\right\vert ^{2\theta}\left\vert g_{F}\right\vert
^{2}<\infty\right\}  .
\]

Set $g=\left(  \frac{\psi}{\sqrt{w}\left\vert \cdot\right\vert ^{\theta}%
}\right)  ^{\vee}$. Then $g_{F}=\widehat{g}=\frac{\psi}{\sqrt{w}\left\vert
\cdot\right\vert ^{\theta}}$ and so $\int w\left\vert \cdot\right\vert
^{2\theta}\left\vert g_{F}\right\vert ^{2}=\int\left\vert \psi\right\vert
^{2}<\infty$. Finally, if $K\subset\mathbb{R}^{d}$ is a compact then%
\[
\int_{K}\left\vert \xi^{\alpha}\widehat{g}\right\vert =\int_{K}\frac
{\left\vert \xi^{\alpha}\psi\right\vert }{\sqrt{w}\left\vert \cdot\right\vert
^{\theta}}\leq\int_{K}\frac{\left\vert \psi\right\vert }{\sqrt{w}}\leq\left(
\int_{K}\left\vert \psi\right\vert ^{2}\right)  ^{1/2}\left(  \int_{K}\frac
{1}{w}\right)  ^{1/2}<\infty,
\]

since property W1 requires $1/w\in L_{loc}^{1}$.\bigskip

By the Definition \ref{Def_Itheta} of $\mathcal{I}:X_{w}^{\theta}\rightarrow
L^{2}$, $f\in X_{w}^{\theta}$ implies $\mathcal{I}f=\left(  \sqrt{w}\left\vert
\cdot\right\vert ^{\theta}f_{F}\right)  ^{\vee}\ $where $f_{F}\in L_{loc}%
^{1}\left(  \mathbb{R}^{d}\setminus0\right)  $, $f_{F}=\widehat{f}$ on
$\mathbb{R}^{d}\setminus0$ and $f_{F}\left(  0\right)  =0$.

Suppose $\phi\in\widehat{S}_{\emptyset,\theta}$. Then $G\ast\widehat
{S}_{\emptyset,\theta}\subset X_{w}^{\theta}$ implies $\mathcal{I}\left(
G\ast\phi\right)  \in L^{2}$ and
\[
\left(  \mathcal{I}\left(  G\ast\phi\right)  \right)  ^{\wedge}=\sqrt
{w}\left\vert \cdot\right\vert ^{\theta}\left(  G\ast\phi\right)  _{F}%
=\frac{\widehat{\phi}}{\sqrt{w}\left\vert \cdot\right\vert ^{\theta}}\in
\frac{S_{\emptyset,\theta}}{\sqrt{w}\left\vert \cdot\right\vert ^{\theta}}.
\]

\end{proof}

\begin{definition}
\label{Def_operator_Ttheta_minus}\textbf{The linear operator} $\mathcal{J}%
_{-}:\widehat{S}_{\emptyset,\theta}\rightarrow\left(  \frac{S_{\emptyset
,\theta}}{\sqrt{w}\left\vert \cdot\right\vert ^{\theta}}\right)  ^{\vee}$.

If $\phi\in\widehat{S}_{\emptyset,\theta}$ then $\widehat{\phi}\in
S_{\emptyset,\theta}\subset L^{2}$ and by the previous lemma $\left(
\frac{\widehat{\phi}}{\sqrt{w}\left\vert \cdot\right\vert ^{\theta}}\right)
^{\vee}\in X_{w}^{\theta}$. We define $\mathcal{J}_{-}$ by:%
\begin{equation}
\mathcal{J}_{-}\phi=\left(  \frac{\widehat{\phi}}{\sqrt{w}\left\vert
\cdot\right\vert ^{\theta}}\right)  ^{\vee},\quad\phi\in\widehat{S}%
_{\emptyset,\theta}.\label{p1.013}%
\end{equation}

\end{definition}

We now introduce the linear operator $\mathcal{J}_{+}:L^{2}\rightarrow
X_{w}^{\theta}$.

\begin{theorem}
\label{Thm_operator_Ttheta_plus}\textbf{The operator} $\mathcal{J}_{+}$ If
$\phi\in\widehat{S}_{\emptyset,\theta}$ then $\left(  \frac{\widehat{\phi}%
}{\sqrt{w}\left\vert \cdot\right\vert ^{\theta}}\right)  ^{\vee}\in L^{2}$ and
there exists a \textbf{linear} member $\mathcal{J}_{+}$ of the class of
operators $\mathcal{J}:L^{2}\rightarrow X_{w}^{\theta}$ of Definition
\ref{Def_operator_Ttheta} which satisfies
\begin{equation}
\mathcal{J}_{+}\left(  \left(  \frac{\widehat{\phi}}{\sqrt{w}\left\vert
\cdot\right\vert ^{\theta}}\right)  ^{\vee}\right)  =G\ast\phi,\quad\phi
\in\overset{\wedge}{S}_{\emptyset,\theta}.\label{p1.019}%
\end{equation}

\end{theorem}

\begin{proof}
Suppose $g=\left(  \frac{\widehat{\phi}}{\sqrt{w}\left\vert \cdot\right\vert
^{\theta}}\right)  ^{\vee}$. It is sufficient to show that $\frac{\widehat{g}%
}{\sqrt{w}\left\vert \cdot\right\vert ^{\theta}}=G\ast\phi$ on $S_{\emptyset
,\theta}$, since this would imply that $\frac{\widehat{g}}{\sqrt{w}\left\vert
\cdot\right\vert ^{\theta}}\in S_{\emptyset,\theta}^{\prime}$.

But $\frac{\widehat{g}}{\sqrt{w}\left\vert \cdot\right\vert ^{\theta}}%
=\frac{\widehat{\phi}}{w\left\vert \cdot\right\vert ^{2\theta}}$ so for
$\psi\in S_{\emptyset,\theta}$%
\begin{align*}
\left[  \frac{\widehat{g}}{\sqrt{w}\left\vert \cdot\right\vert ^{\theta}}%
,\psi\right]  =\left[  \frac{\widehat{\phi}}{w\left\vert \cdot\right\vert
^{2\theta}},\psi\right]  =\int\frac{\widehat{\phi}}{w\left\vert \cdot
\right\vert ^{2\theta}}\overset{\vee}{\psi}  & =\int\frac{1}{w\left\vert
\cdot\right\vert ^{2\theta}}\widehat{\phi}\overset{\vee}{\psi}\\
& =\left[  \widehat{G},\widehat{\phi}\overset{\vee}{\psi}\right]
,\quad\widehat{\phi}\overset{\vee}{\psi}\in S_{\emptyset,2\theta},\\
& =\left[  \widehat{G\ast\phi},\overset{\vee}{\psi}\right] \\
& =\left[  G\ast\phi,\psi\right]  .
\end{align*}

By the Hahn-Banach theorem we can extend the linear operator defined by
\ref{p1.019} to a linear operator $\mathcal{J}_{+}$ in the class $\mathcal{J}$.
\end{proof}

\begin{definition}
\label{Def_operator_Ttheta_plus}\textbf{The linear operator} $\mathcal{J}%
_{+}:L^{2}\rightarrow X_{w}^{\theta}$

Let $\mathcal{J}_{+}$ be the operator introduced in Theorem
\ref{Thm_operator_Ttheta_plus}.
\end{definition}

We next establish our results of Figure \ref{Fig_I_Jplus_Jminus}\ concerning
$\mathcal{I}$, $\mathcal{J}_{+}$ and $\mathcal{J}_{-}$.

\begin{theorem}
\label{Thm_I_J_S_G*S}Suppose the weight function $w$ has properties W1 and W2,
and suppose $G$ is a basis distribution of order $\theta$ generated by $w$.

Then we have:

\begin{enumerate}
\item As tempered distributions%
\[
\mathcal{I}\left(  G\ast\phi\right)  =\left(  \frac{\widehat{\phi}}{\sqrt
{w}\left\vert \cdot\right\vert ^{\theta}}\right)  ^{\vee}=\mathcal{J}_{-}%
\phi,\quad\phi\in\widehat{S}_{\emptyset,\theta}.
\]

\item As tempered distributions%
\[
\mathcal{I}\left(  \left(  \frac{\widehat{\phi}}{\sqrt{w}\left\vert
\cdot\right\vert ^{\theta}}\right)  ^{\vee}\right)  =\phi,\quad\phi\in
\widehat{S}_{\emptyset,\theta}.
\]

\item $\mathcal{I}:G\ast\widehat{S}_{\emptyset,\theta}\rightarrow\left(
\frac{S_{\emptyset,\theta}}{\sqrt{w}\left\vert \cdot\right\vert ^{\theta}%
}\right)  ^{\vee}$ and $\mathcal{J}_{+}:\left(  \frac{S_{\emptyset,\theta}%
}{\sqrt{w}\left\vert \cdot\right\vert ^{\theta}}\right)  ^{\vee}\rightarrow
G\ast\widehat{S}_{\emptyset,\theta}$.

$\mathcal{I}:G\ast\widehat{S}_{\emptyset,\theta}\rightarrow\left(
\frac{S_{\emptyset,\theta}}{\sqrt{w}\left\vert \cdot\right\vert ^{\theta}%
}\right)  ^{\vee}$ is an isometric isomorphism with inverse $\mathcal{J}_{+}%
$when $G\ast\widehat{S}_{\emptyset,\theta}$ is regarded as a subspace of
$X_{w}^{\theta}$ and $\left(  \frac{S_{\emptyset,\theta}}{\sqrt{w}\left\vert
\cdot\right\vert ^{\theta}}\right)  ^{\vee}$ is regarded as a subspace of
$L^{2}$. $\mathcal{I}$ has inverse $\mathcal{J}_{+}$.

\item $\mathcal{I}:\left(  \frac{S_{\emptyset,\theta}}{\sqrt{w}\left\vert
\cdot\right\vert ^{\theta}}\right)  ^{\vee}\rightarrow\widehat{S}%
_{\emptyset,\theta}$ and $\mathcal{J}_{-}:\overset{\wedge}{S}_{\emptyset
,\theta}\rightarrow\left(  \frac{S_{\emptyset,\theta}}{\sqrt{w}\left\vert
\cdot\right\vert ^{\theta}}\right)  ^{\vee}$.

$\mathcal{I}:\left(  \frac{S_{\emptyset,\theta}}{\sqrt{w}\left\vert
\cdot\right\vert ^{\theta}}\right)  ^{\vee}\rightarrow\widehat{S}%
_{\emptyset,\theta}$ is an isometric isomorphism with inverse $\mathcal{J}%
_{-}$ when $\overset{\wedge}{S}_{\emptyset,\theta}$ is regarded as a subspace
of $X_{1/w}^{-\theta}$ and $\left(  \frac{S_{\emptyset,\theta}}{\sqrt
{w}\left\vert \cdot\right\vert ^{\theta}}\right)  ^{\vee}$ is regarded as a
subspace of $L^{2}$.
\end{enumerate}
\end{theorem}

\begin{proof}
\textbf{Part 1}. The first equation of part 2 follows immediately from
expression \ref{p1.7} of Lemma \ref{Lem_F(S/sqrt(w|.|^th))_in_Xthw}. The
second equation of this part follows from the definition of $\mathcal{J}_{-}%
$.\medskip

\textbf{Part 2}. Apply $\mathcal{I}$ to the second equation of part 1.\medskip

\textbf{Part 3}. The first line follows from part 1 and the definition of
$\mathcal{J}_{+}$.

Applying $\mathcal{J}_{+}$ to the first equation of part 1 and then using
\ref{p1.019} gives $\mathcal{J}_{+}\mathcal{I}=I$. We already know that
$\mathcal{IJ}_{+}=I$ so it follows that $\mathcal{I}$ is an isomorphism with
inverse $\mathcal{J}_{+}$. Since $\mathcal{I}:X_{w}^{\theta}\rightarrow L^{2}$
is isometric this part is true.\medskip

\textbf{Part 4}. The first line follows from part 2 and the second equation of
part 1. It is also clear that $\mathcal{J}_{-}\mathcal{I}=I=\mathcal{IJ}_{-}$
and so $\mathcal{I}$ is an isomorphism with inverse $\mathcal{J}_{+}$. From
part 2, $\mathcal{I}\left(  \left(  \frac{\widehat{\phi}}{\sqrt{w}\left\vert
\cdot\right\vert ^{\theta}}\right)  ^{\vee}\right)  =\phi$ when $\phi
\in\widehat{S}_{\emptyset,\theta}$. Using definition \ref{Def_X1/w,-th_II} of
the semi-inner product vector space $X_{1/w}^{-\theta}$ and setting $g=\left(
\frac{\widehat{\phi}}{\sqrt{w}\left\vert \cdot\right\vert ^{\theta}}\right)
^{\vee}\in L^{2}$, we get the isometry
\[
\left\vert \mathcal{I}g\right\vert _{1/w,-\theta}^{2}=\int\frac{\left\vert
\widehat{\phi}\right\vert ^{2}}{w\left\vert \cdot\right\vert ^{2\theta}%
}=\left\Vert g\right\Vert _{2}^{2}.
\]

\end{proof}

\begin{corollary}
$\left(  \frac{S_{\emptyset,\theta}}{\sqrt{w}\left\vert \cdot\right\vert
^{\theta}}\right)  ^{\vee}$ is dense in $L^{2}$.
\end{corollary}

\begin{proof}
Use Theorem \ref{Thm_I_J_S_G*S} and Corollary
\ref{Cor_F[So,th]_dense_in_X1/w,-th}.
\end{proof}

\begin{corollary}
$\mathcal{I}=\widetilde{\mathcal{L}}$ on $\left(  \frac{S_{\emptyset,\theta}%
}{\sqrt{w}\left\vert \cdot\right\vert ^{\theta}}\right)  ^{\vee}$;
$\mathcal{J}_{-}=\widetilde{\mathcal{M}}$ on $\widehat{S}_{\emptyset,\theta}
$; $\mathcal{J}_{+}=\mathcal{J}$ on $\left(  \frac{S_{\emptyset,\theta}}%
{\sqrt{w}\left\vert \cdot\right\vert ^{\theta}}\right)  ^{\vee}$.
\end{corollary}

\begin{proof}
If $f\in\left(  \frac{S_{\emptyset,\theta}}{\sqrt{w}\left\vert \cdot
\right\vert ^{\theta}}\right)  ^{\vee}\subset L^{2}$ then $\widehat{f}%
=\frac{\widehat{\phi}}{\sqrt{w}\left\vert \cdot\right\vert ^{\theta}}$ for
some $\phi\in\widehat{S}_{\emptyset,\theta}$.

Then $\widehat{\mathcal{J}\phi}=\frac{\widehat{\phi}}{\sqrt{w}\left\vert
\cdot\right\vert ^{\theta}}=\widehat{\mathcal{J}_{-}\phi}$.
\end{proof}

\begin{corollary}
The class $\mathcal{J}$ contains a linear operator which equals $\mathcal{J}%
_{+}$ on $\left(  \frac{S_{\emptyset,\theta}}{\sqrt{w}\left\vert
\cdot\right\vert ^{\theta}}\right)  ^{\vee}$.
\end{corollary}

\begin{proof}
Apply the Hahn-Banach theorem.
\end{proof}

\begin{corollary}
\label{Cor_Thm_I_J_S_G*S}$\mathcal{I}^{2}$ and $\mathcal{J}_{+}\mathcal{J}%
_{-}$ have the following properties:

\begin{enumerate}
\item $\mathcal{I}^{2}\left(  G\ast\phi\right)  =\phi$,$\quad\phi\in
\widehat{S}_{\emptyset,\theta}$.

\item $\mathcal{I}^{2}:G\ast\widehat{S}_{\emptyset,\theta}\rightarrow
\widehat{S}_{\emptyset,\theta}$ is an isometric isomorphism with inverse
$\mathcal{J}_{+}\mathcal{J}_{-}$ when $G\ast\widehat{S}_{\emptyset,\theta}$ is
regarded as a subspace of the Hilbert space $X_{w}^{\theta}$ and $\widehat
{S}_{\emptyset,\theta}$ is regarded as a subspace of the Hilbert space
$X_{1/w}^{-\theta}$.

\item $\mathcal{V}=\mathcal{J}_{+}\mathcal{J}_{-}$.
\end{enumerate}
\end{corollary}

\begin{proof}
\textbf{Part 1}. Apply operator $\mathcal{I}$ to the first equation of part 1
of Theorem \ref{Thm_I_J_S_G*S} and then use part 2 of the same
theorem.\medskip

\textbf{Part 2}. Direct consequence of parts 3 and 4 of Theorem
\ref{Thm_I_J_S_G*S}.\medskip

\textbf{Part 3}. From Theorem \ref{Thm_Vth_and_F(So,th)} $\mathcal{V}%
\phi=G\ast\phi$ when $\phi\in\widehat{S}_{\emptyset,\theta}$. From equations
\ref{p1.013} and \ref{p1.019} $\mathcal{J}_{+}\mathcal{J}_{-}\phi=G\ast\phi$
when $\phi\in\widehat{S}_{\emptyset,\theta}$.
\end{proof}%

\appendix

\chapter{Basic notation, definitions and symbols\label{App_basic_notation}}

\section{Basic function and distribution spaces}

\begin{definition}
\label{Def_Some_basic_spaces}\textbf{Basic function and distribution spaces}

All spaces below consist of\textbf{\ complex-valued} functions:

\begin{enumerate}
\item $P_{0}=\left\{  0\right\}  $. For $n\geq1$, $P_{n}$ denotes the
polynomials of degree at most $n$. These polynomials have the form
$\sum\limits_{\left\vert \alpha\right\vert <n}a_{\alpha}\xi^{\alpha}$, where
$a_{\alpha}\in\mathbb{C}$ and $\xi\in\mathbb{R}^{d}$. The space of all
polynomials will be denoted by $P$.

\item $C^{\left(  0\right)  }$ is the space of continuous functions.

\noindent$C_{B}^{\left(  0\right)  }$ is the space of bounded continuous functions.

\noindent$C_{BP}^{\left(  0\right)  }$ is the space of continuous functions
bounded by a polynomial.

\noindent$C^{\left(  m\right)  }=\left\{  f\in C^{\left(  0\right)
}:D^{\alpha}f\in C^{\left(  0\right)  },\text{ }when\text{ }\left\vert
\alpha\right\vert =m\right\}  $.

\noindent$C_{B}^{\left(  m\right)  }=\left\{  f\in C_{B}^{\left(  0\right)
}:D^{\alpha}f\in C_{B}^{\left(  0\right)  },\text{ }when\text{ }\left\vert
\alpha\right\vert \leq m\right\}  $.

\noindent$C_{BP}^{\left(  m\right)  }=\left\{  f\in C_{BP}^{\left(  0\right)
}:D^{\alpha}f\in C_{BP}^{\left(  0\right)  },\text{ }when\text{ }\left\vert
\alpha\right\vert \leq m\right\}  $.

\noindent$C^{\infty}=\bigcap\limits_{m\geq0}C^{\left(  m\right)  }$%
\textbf{;}$\quad C_{B}^{\infty}=\bigcap\limits_{m\geq0}C_{B}^{\left(
m\right)  }$\textbf{;}\quad$C_{BP}^{\infty}=\bigcap\limits_{m\geq0}%
C_{BP}^{\left(  m\right)  }$.

\item $L_{loc}^{1}$ is the space of measurable functions which are absolutely
integrable on any compact set i.e. any closed, bounded set.

$L^{1}$ is the Hilbert space of measurable functions $f$ such that
$\int\left\vert f\right\vert <\infty$. Norm is $\left\Vert f\right\Vert
_{1}=\int\left\vert f\right\vert $ and inner product is $\left(  f,g\right)
_{1}$.

$L^{2}$ is the Hilbert space of measurable functions $f$ such that
$\int\left\vert f\right\vert ^{2}<\infty$. Norm is $\left\Vert f\right\Vert
_{2}=\left(  \int\left\vert f\right\vert ^{2}\right)  ^{1/2}$ and inner
product is $\left(  f,g\right)  _{2}$.

$L^{\infty}$ is the space of a.e. bounded functions. Norm is $\left\Vert
f\right\Vert _{\infty}$.and inner product is $\left(  f,g\right)  _{\infty}$.

\item $C_{0}^{\infty}$ is the space of $C^{\infty}$ functions that have
compact support. These are the test functions for the space of distributions
defined on $\mathbb{R}^{d}$, sometimes denoted by $\mathcal{D}$.

$S$ is the $C^{\infty}$ space of rapidly decreasing functions. A function
$f\in C^{\infty}$ is in $S$ if given any multi-index $\alpha$ and integer
$n\geq0$, there exists a constant $k_{\alpha,n}$ such that, $\left\vert
D^{\alpha}f\left(  x\right)  \right\vert \leq k_{\alpha,n}\left(  1+\left\vert
x\right\vert \right)  ^{-n}$, $x\in\mathbb{R}^{d}$. These are the test
functions for the tempered distributions of Definition \ref{Def_Distributions} below.

$\mathcal{D}^{\prime}$ is the space of distributions (generalized functions).

$\mathcal{E}^{\prime}$ is the space of distributions with compact (bounded) support.
\end{enumerate}
\end{definition}

\section{Vector notation}

\begin{definition}
\label{Def_vector_notation}\textbf{Vector notation}

Suppose $v,w\in\mathbb{R}^{d}$. Then:

\begin{itemize}
\item $\mathbf{0}=\left(  0,0,\ldots,0\right)  =0_{d}$ and $\mathbf{1}=\left(
1,1,\ldots,1\right)  =1_{d}$.

\item Suppose $\symbol{126}$ is one of the binary operations $<,\leq,=,>,\geq$.

Then we write $v\symbol{126}w$ if $v_{i}\symbol{126}w_{i}$ for all $i$.

For $x\in\mathbb{R}$ we write $v\symbol{126}x$ if $v_{i}\symbol{126}x$ for all
$i$.

\item $v^{w}=\left(  v_{1}\right)  ^{w_{1}}\left(  v_{2}\right)  ^{w_{2}%
}\ldots\left(  v_{d}\right)  ^{w_{d}}$.

\item If $s\in\mathbb{R}$ then $v^{s}=v^{s\mathbf{1}}$.
\end{itemize}
\end{definition}

\section{Topology}

\begin{definition}
\label{Def_topol_notation}\textbf{Topology on} $R^{d}$

\begin{itemize}
\item The Euclidean norm an inner product are denoted by $\left\vert
x\right\vert $ and $\left(  x,y\right)  $.

\item $\left\vert x\right\vert _{\infty}=\max\limits_{i}\left\vert
x_{i}\right\vert $, $\left\vert x\right\vert _{1}=\sum\limits_{i}\left\vert
x_{i}\right\vert $ and in general for $p\geq1$: $\left\vert x\right\vert
_{p}=\left(  \sum\limits_{i}\left\vert x_{i}\right\vert ^{p}\right)  ^{1/p}$.

Note that $\left\vert x\right\vert _{\infty}=\lim\limits_{p\rightarrow\infty
}\left\vert x\right\vert _{p}$ and $\left\vert x\right\vert =\left\vert
x\right\vert _{2}$.

\item Ball $B\left(  x;r\right)  =\left\{  y:\left\vert x-y\right\vert
<r\right\}  $.

\item $\varepsilon\mathbf{-}$neighborhood of a set - for $\varepsilon>0$, the
$\varepsilon-$neighborhood of a set $S$ is denoted $S_{\varepsilon}$ and
$S_{\varepsilon}=\bigcup\limits_{x\in S}B\left(  x;\varepsilon\right)  $.
\end{itemize}
\end{definition}

\section{Multi-indexes}

\begin{summary}
\label{Sum_Multi_index}\textbf{Multi-indexes, definitions and identities}

\begin{enumerate}
\item Multi-indexes are vectors with non-negative integer components.

Let $\alpha=\left(  \alpha_{1},\alpha_{2},\ldots,\alpha_{d}\right)  $ and
$\beta=\left(  \beta_{1},\beta_{2},\ldots,\beta_{d}\right)  $ be multi-indexes.

Suppose $\symbol{126}$ is one of the binary operations $<$, $\leq$, $=$, $>$,
$\geq$.

Write $\beta\symbol{126}a$ if $\beta_{i}\symbol{126}\alpha_{i}$ for all $i $.

For $x\in\mathbb{R}$ write $\beta\symbol{126}x$ if $\beta_{i}\symbol{126}x $
for all $i$.

\item Denote $\left\vert \alpha\right\vert =\sum\limits_{i=1}^{d}\alpha_{i} $.
Then $D^{\alpha}f\left(  x\right)  $ is the derivative of the function $f$ of
degree $\alpha$
\[
D^{0}f\left(  x\right)  =f\left(  x\right)  ,\text{\quad}D^{\alpha}f\left(
x\right)  =\frac{D^{\left\vert \alpha\right\vert }f\left(  x_{1},x_{2}%
,\ldots,x_{d}\right)  }{D_{1}^{\alpha_{1}}x_{1}D_{2}^{\alpha_{2}}x_{2}\ldots
D_{d}^{\alpha_{d}}x_{d}}.
\]

\[
\frac{1}{\beta!}D^{\beta}x^{\alpha}=\left\{
\begin{array}
[c]{ll}%
\binom{\alpha}{\beta}x^{\alpha-\beta}, & \beta\leq\alpha,\\
0, & otherwise.
\end{array}
\right.
\]

\item We shall also use the notation
\[%
\begin{array}
[c]{ll}%
\text{monomial} & x^{\alpha}=x_{1}^{\alpha_{1}}x_{2}^{\alpha_{2}}\ldots
x_{d}^{\alpha_{d}},\\
\text{factorial} & \alpha!=\alpha_{1}!\alpha_{2}!\ldots\alpha_{d}!\text{
}and\text{ }0!=1,\\
\text{binomial} & \dbinom{\alpha}{\beta}=\dfrac{\alpha!}{\left(  \alpha
-\beta\right)  !\beta!},\quad if\text{ }\beta\leq\alpha,
\end{array}
\]

and we have%
\[
\alpha!\leq\left\vert \alpha\right\vert !
\]

\item The inequalities $\left\vert x^{\alpha}\right\vert \leq\left\vert
x\right\vert ^{\left\vert \alpha\right\vert }$ and $\left\vert \left(
1+x\right)  ^{\alpha}\right\vert \leq\left(  1+\left\vert x\right\vert
\right)  ^{\left\vert \alpha\right\vert }$ are used often. These are proved
using the fact that
\[
\text{weighted\_geometric\_mean}\leq\text{weighted\_arithmetic\_mean.}%
\]

Similarly one can prove, using the generalized means of \S 4.6 Archbold
\cite{Archbold64}, that for $0<p\leq\infty$,
\begin{align*}
\left\vert x^{\alpha}\right\vert  & \leq\left(  \frac{\alpha}{\left\vert
\alpha\right\vert }x_{+}^{p\mathbf{1}}\right)  ^{\left\vert \alpha\right\vert
/p}\leq\left(  \frac{\left\vert \alpha\right\vert _{\max}}{\left\vert
\alpha\right\vert }\right)  ^{\left\vert \alpha\right\vert /p}\left\vert
x\right\vert _{p}^{\left\vert \alpha\right\vert }\leq\left\vert x\right\vert
_{p}^{\left\vert \alpha\right\vert },\\
\left\vert \left(  1+x\right)  ^{\alpha}\right\vert  & \leq\left(
1+\left\vert x\right\vert _{p}\right)  ^{\left\vert \alpha\right\vert },\\
\left\vert \left(  1+x\right)  ^{\alpha}\right\vert  & \leq\left(
1+\frac{\alpha}{\left\vert \alpha\right\vert }x_{+}^{p\mathbf{1}}\right)
^{\left\vert \alpha\right\vert }\leq\left(  1+\frac{\left\vert \alpha
\right\vert _{\max}}{\left\vert \alpha\right\vert }\left\vert x\right\vert
_{p}\right)  ^{\left\vert \alpha\right\vert }\leq\left(  1+\left\vert
x\right\vert _{p}\right)  ^{\left\vert \alpha\right\vert }.
\end{align*}

\item Important identities are%
\begin{equation}
\frac{1}{k!}\left(  x,y\right)  ^{k}=\sum_{\left\vert \alpha\right\vert
=k}\frac{1}{\alpha!}x^{\alpha}y^{\alpha},\text{\quad}\frac{1}{k!}\left\vert
x\right\vert ^{2k}=\sum_{\left\vert \alpha\right\vert =k}\frac{1}{\alpha
!}x^{2\alpha}.\label{p08}%
\end{equation}

and a direct consequence is%
\[
\sum_{\beta\leq\alpha}\binom{\alpha}{\beta}=2^{\left\vert \alpha\right\vert }.
\]

Also%
\begin{align}
\frac{1}{k!}\left(  \left\vert x\right\vert ^{2}+\left\vert y\right\vert
^{2}\right)  ^{k}  & =\sum_{\left\vert \alpha+\beta\right\vert =k}%
\frac{x^{2\alpha}y^{2\beta}}{\alpha!\beta!},\label{Ap125}\\
\frac{1}{k!}\left\vert \left\vert x\right\vert ^{2}+\left\vert y\right\vert
^{2}+\left\vert z\right\vert ^{2}\right\vert ^{k}  & =\sum_{\left\vert
\alpha+\beta+\gamma\right\vert =k}\frac{x^{2\alpha}y^{2\beta}z^{2\gamma}%
}{\alpha!\beta!\gamma!}.\label{Ap126}%
\end{align}

\item If $\binom{k}{\alpha}=\frac{k!}{\alpha!\left(  k-\left\vert
\alpha\right\vert \right)  !}$ then
\begin{equation}
\left(  1+xy\right)  ^{k}=\sum_{\left\vert \alpha\right\vert \leq k}\binom
{k}{\alpha}x^{\alpha}y^{\alpha}.\label{p021}%
\end{equation}

\item Useful identities are%
\begin{equation}
\sum\limits_{\left\vert \alpha\right\vert =k}1=\binom{d+k-1}{k},\quad
k\geq1,\label{Ap023}%
\end{equation}

and the dimension of a polynomial of degree $k$ is given by%
\begin{equation}
\sum\limits_{\left\vert \alpha\right\vert \leq k}1=\binom{d+k}{d},\quad
k\geq0.\label{a1.11}%
\end{equation}

\item \textbf{Leibniz's rule} is formally%
\begin{equation}
D^{\alpha}\left(  uv\right)  =\sum_{\beta\leq\alpha}\dbinom{\alpha}{\beta
}D^{\beta}uD^{\alpha-\beta}v.\label{p68}%
\end{equation}

Examples are: $u\in\mathcal{D}$, $v\in C_{0}^{\infty}$ and $u\in S^{\prime}$,
$v\in C_{BP}^{\infty}$ and $u,v\in C^{\left(  \alpha\right)  }=\left\{  f\in
C^{\left(  0\right)  }:D^{\beta}f\in C^{\left(  0\right)  }\text{ }\beta
\leq\alpha\right\}  $.

\item \textbf{Binomial expansions}%
\[
\left(  x+y\right)  ^{\gamma}=\sum_{\alpha\leq\gamma}\tbinom{\gamma}{\alpha
}x^{\alpha}y^{\gamma-\alpha},
\]

or equivalently%
\[
\frac{1}{\gamma!}\left(  x+y\right)  ^{\gamma}=\sum_{\alpha+\beta=\gamma}%
\frac{x^{\alpha}y^{\beta}}{\alpha!\beta!}.
\]

In the same vein%
\[
\frac{1}{\delta!}\left(  x+y+z\right)  ^{\delta}=\sum_{\alpha+\beta
+\gamma=\delta}\frac{x^{\alpha}y^{\beta}z^{\gamma}}{\alpha!\beta!\gamma!}.
\]

\end{enumerate}
\end{summary}

We will need:

\begin{theorem}
\label{Thm_Ap_1}??

\begin{enumerate}
\item
\[
\sum\limits_{\left\vert \alpha\right\vert \leq m}\sum\limits_{\beta\leq\alpha
}\frac{x^{\beta}}{\beta!}=\sum\limits_{\left\vert \gamma\right\vert \leq
m}\tbinom{d+m-\left\vert \gamma\right\vert }{d}\frac{x^{\gamma}}{\gamma!}%
=\sum\limits_{k=0}^{m}\tbinom{d+m-k}{d}\frac{\left(  x\mathbf{1}\right)  ^{k}%
}{k!},
\]

and%
\[
\sum\limits_{\left\vert \alpha\right\vert \leq m}\sum\limits_{\beta\leq\alpha
}\frac{\left(  x.\ast y\right)  ^{\beta}}{\beta!}=\sum\limits_{k=0}^{m}%
\tbinom{d+m-k}{d}\frac{\left(  xy\right)  ^{k}}{k!}.
\]

Also, if $xy\geq0$,%
\[
\left\vert \sum\limits_{\left\vert \alpha\right\vert \leq m}\sum
\limits_{\beta\leq\alpha}\frac{\left(  x.\ast y\right)  ^{\beta}}{\beta
!}\right\vert \leq\tbinom{d+m}{d}\left(  1+xy\right)  ^{m}.
\]

\item
\[
\sum\limits_{\left\vert \alpha\right\vert =k}\sum\limits_{\beta\leq\alpha
}\frac{x^{\beta}}{\beta!}=\sum\limits_{\left\vert \gamma\right\vert \leq
k}\tbinom{d+k-\left\vert \gamma\right\vert -1}{k-\left\vert \gamma\right\vert
}\frac{x^{\gamma}}{\gamma!}=\sum\limits_{j=0}^{k}\tbinom{d+k-j-1}{k-j}%
\frac{\left(  x\mathbf{1}\right)  ^{j}}{j!},
\]

and%
\[
\sum\limits_{\left\vert \alpha\right\vert =k}\sum\limits_{\beta\leq\alpha
}\frac{\left(  x.\ast y\right)  ^{\beta}}{\beta!}=\sum\limits_{j=0}^{k}%
\tbinom{d+k-j-1}{k-j}\frac{\left(  xy\right)  ^{j}}{j!}.
\]

Also, if $xy\geq0$,%
\[
\sum\limits_{\left\vert \alpha\right\vert =k}\sum\limits_{\beta\leq\alpha
}\frac{\left(  x.\ast y\right)  ^{\beta}}{\beta!}\leq\tbinom{d+k-1}%
{d-1}\left(  1+xy\right)  ^{k}.
\]

\item
\[
\sum\limits_{\left\vert \alpha\right\vert =k}y^{\alpha}\sum\limits_{\beta
\leq\alpha}\frac{x^{\beta}}{\beta!}=\sum\limits_{j=0}^{k}h_{k-j}\left(
y\right)  \frac{\left(  xy\right)  ^{j}}{j!}=\sum\limits_{j=0}^{k}%
h_{k-j}\left(  \widehat{y}\right)  \left\vert y\right\vert ^{k-j}\frac{\left(
xy\right)  ^{j}}{j!}=\left\vert y\right\vert ^{k}\sum\limits_{j=0}^{k}%
h_{k-j}\left(  \widehat{y}\right)  \frac{\left(  x\widehat{y}\right)  ^{j}%
}{j!},
\]

where $h_{k-j}\left(  y\right)  =\sum\limits_{\left\vert \gamma\right\vert
=k-j}y^{\gamma}$ is the \textbf{complete homogeneous symmetric polynomial} of
order $k-j$.

\item If $x\widehat{y}\geq0$ then%
\[
\sum\limits_{\left\vert \alpha\right\vert =k}y^{\alpha}\sum\limits_{\beta
\leq\alpha}\frac{x^{\beta}}{\beta!}\leq\binom{d+k-1}{k}\left\vert y\right\vert
^{k}\left(  1+x\widehat{y}\right)  ^{k}.
\]

\item ?? \textbf{CHECK}! ??%
\[
\sum\limits_{\left\vert \alpha\right\vert \leq m}y^{\alpha}\sum\limits_{\beta
\leq\alpha}\frac{x^{\beta}}{\beta!}=\sum\limits_{j=0}^{m}\left(
\sum\limits_{k=0}^{m-j}h_{k}\left(  y\right)  \right)  \frac{\left(
xy\right)  ^{j}}{j!}=\sum\limits_{j=0}^{m}h_{j}\left(  y\right)
\sum\limits_{k=0}^{m-j}\frac{\left(  xy\right)  ^{k}}{k!}.
\]

\item If $y\geq0_{d}$ then ??%
\[
h_{k}\left(  \widehat{y}\right)  \leq h_{k}\left(  \frac{1}{\sqrt{d}%
}\mathbf{1}\right)  =d^{-k/2}\tbinom{d+k-1}{k}.
\]

\end{enumerate}
\end{theorem}

\begin{proof}
\textbf{Part 1} Since $\left\vert \beta\right\vert \leq\left\vert
\alpha\right\vert $, there exist constants $\mu_{\gamma}$ such that%
\[
\sum\limits_{\left\vert \alpha\right\vert \leq m}\sum\limits_{\beta\leq\alpha
}\frac{x^{\beta}}{\beta!}=\sum\limits_{\left\vert \gamma\right\vert \leq m}%
\mu_{\gamma}\frac{x^{\gamma}}{\gamma!}.
\]

In fact ?? OBSCURE - FIX! ??%
\begin{align*}
\mu_{\gamma}  & =\#\left\{  \delta:\left\vert \delta\right\vert \leq m\text{
}and\text{ }\gamma\leq\delta\right\} \\
& =\#\left\{  \gamma+\sigma:\left\vert \gamma+\sigma\right\vert \leq
m\right\}  =\\
& =\#\left\{  \sigma:\left\vert \gamma+\sigma\right\vert \leq m\right\} \\
& =\#\left\{  \sigma:\left\vert \sigma\right\vert \leq m-\left\vert
\gamma\right\vert \right\} \\
& =\sum\limits_{\left\vert \sigma\right\vert \leq m-\left\vert \gamma
\right\vert }\mathbf{1}.
\end{align*}

But by \ref{a1.11}, $\mu_{\gamma}=\tbinom{d+m-\left\vert \gamma\right\vert
}{d}$ and so%
\begin{align*}
\sum\limits_{\left\vert \alpha\right\vert \leq m}\sum\limits_{\beta\leq\alpha
}\frac{x^{\beta}}{\beta!}  & =\sum\limits_{\left\vert \gamma\right\vert \leq
m}\tbinom{d+m-\left\vert \gamma\right\vert }{d}\frac{x^{\gamma}}{\gamma!}%
=\sum\limits_{k=0}^{m}\sum\limits_{\left\vert \gamma\right\vert =k}%
\tbinom{d+m-\left\vert \gamma\right\vert }{d}\frac{x^{\gamma}}{\gamma!}=\\
& =\sum\limits_{k=0}^{m}\sum\limits_{\left\vert \gamma\right\vert =k}%
\tbinom{d+m-k}{d}\frac{x^{\gamma}}{\gamma!}=\sum\limits_{k=0}^{m}%
\tbinom{d+m-k}{d}\sum\limits_{\left\vert \gamma\right\vert =k}\frac{x^{\gamma
}}{\gamma!}=\\
& =\sum\limits_{k=0}^{m}\tbinom{d+m-k}{d}\sum\limits_{\left\vert
\gamma\right\vert =k}\frac{x^{\gamma}\mathbf{1}^{\gamma}}{\gamma!}%
=\sum\limits_{k=0}^{m}\tbinom{d+m-k}{d}\frac{\left(  x\mathbf{1}\right)  ^{k}%
}{k!}.
\end{align*}

\begin{align*}
\sum\limits_{k=0}^{m}\tbinom{d+m-k}{d}\frac{\left(  xy\right)  ^{k}}{k!}  &
=\sum\limits_{k=0}^{m}\tbinom{d+m-k}{d}\frac{\left(  xy\right)  ^{k}}{k!}%
=\sum\limits_{k=0}^{m}\tbinom{d+m-k}{d}\sum\limits_{\left\vert \gamma
\right\vert =k}\frac{x^{\gamma}y^{\gamma}}{\gamma!}=\\
& =\sum\limits_{k=0}^{m}\sum\limits_{\left\vert \gamma\right\vert =k}%
\tbinom{d+m-k}{d}\frac{x^{\gamma}y^{\gamma}}{\gamma!}=\sum\limits_{k=0}%
^{m}\sum\limits_{\left\vert \gamma\right\vert =k}\tbinom{d+m-\left\vert
\gamma\right\vert }{d}\frac{x^{\gamma}y^{\gamma}}{\gamma!}=\\
& =\sum\limits_{\left\vert \gamma\right\vert \leq m}\tbinom{d+m-\left\vert
\gamma\right\vert }{d}\frac{\left(  x.\ast y\right)  ^{\gamma}}{\gamma!}%
=\sum\limits_{\left\vert \alpha\right\vert \leq m}\sum\limits_{\beta\leq
\alpha}\frac{\left(  x.\ast y\right)  ^{\beta}}{\beta!}.
\end{align*}

\begin{align*}
\sum\limits_{k=0}^{m}\tbinom{d+m-k}{d}\frac{\left(  xy\right)  ^{k}}{k!}  &
=\sum\limits_{k=0}^{m}\frac{\left(  d+m-k\right)  !}{d!\left(  m-k\right)
!}\frac{1}{k!}\left(  xy\right)  ^{k}=\sum\limits_{k=0}^{m}\frac{\left(
d+m-k\right)  !}{m!d!}\frac{m!}{\left(  m-k\right)  !k!}\left(  xy\right)
^{k}\leq\\
& \leq\frac{\left(  d+m\right)  !}{m!d!}\sum\limits_{k=0}^{m}\frac{m!}{\left(
m-k\right)  !k!}\left(  xy\right)  ^{k}=\tbinom{d+m}{d}\left(  1+xy\right)
^{m}.
\end{align*}

\textbf{Part 2} Since $\left\vert \beta\right\vert \leq\left\vert
\alpha\right\vert $, there exist constants $\eta_{\gamma}$ such that%
\[
\sum\limits_{\left\vert \alpha\right\vert =k}\sum\limits_{\beta\leq\alpha
}\frac{x^{\beta}}{\beta!}=\sum\limits_{\left\vert \gamma\right\vert \leq
k}\eta_{\gamma}\frac{x^{\gamma}}{\gamma!}.
\]

In fact%
\begin{align*}
\eta_{\gamma}  & =\#\left\{  \delta:\left\vert \delta\right\vert =k\text{
}and\text{ }\gamma\leq\delta\right\}  =\#\left\{  \gamma+\sigma:\left\vert
\gamma+\sigma\right\vert =k\right\}  =\#\left\{  \sigma:\left\vert
\gamma+\sigma\right\vert =k\right\}  =\\
& =\#\left\{  \sigma:\left\vert \sigma\right\vert =k-\left\vert \gamma
\right\vert \right\}  =\sum\limits_{\left\vert \sigma\right\vert =k-\left\vert
\gamma\right\vert }\mathbf{1}.
\end{align*}

But from \ref{Ap023}, $\eta_{\gamma}=\binom{d+k-\left\vert \gamma\right\vert
-1}{k-\left\vert \gamma\right\vert }$ so%
\begin{align*}
\sum\limits_{\left\vert \alpha\right\vert =k}\sum\limits_{\beta\leq\alpha
}\frac{x^{\beta}}{\beta!}  & =\sum\limits_{\left\vert \gamma\right\vert \leq
k}\tbinom{d+k-\left\vert \gamma\right\vert -1}{k-\left\vert \gamma\right\vert
}\frac{x^{\gamma}}{\gamma!}=\sum\limits_{j=0}^{k}\sum\limits_{\left\vert
\gamma\right\vert =j}\tbinom{d+k-\left\vert \gamma\right\vert -1}{k-\left\vert
\gamma\right\vert }\frac{x^{\gamma}}{\gamma!}=\\
& =\sum\limits_{j=0}^{k}\sum\limits_{\left\vert \gamma\right\vert =j}%
\tbinom{d+k-j-1}{k-j}\frac{x^{\gamma}}{\gamma!}=\sum\limits_{j=0}^{k}%
\tbinom{d+k-j-1}{k-j}\sum\limits_{\left\vert \gamma\right\vert =j}%
\frac{x^{\gamma}}{\gamma!}=\\
& =\sum\limits_{j=0}^{k}\tbinom{d+k-j-1}{k-j}\sum\limits_{\left\vert
\gamma\right\vert =j}\frac{x^{\gamma}\mathbf{1}^{\gamma}}{\gamma!}%
=\sum\limits_{j=0}^{k}\tbinom{d+k-j-1}{k-j}\frac{\left(  x\mathbf{1}\right)
^{j}}{j!}.
\end{align*}
\medskip

\textbf{Part 3} We first write%
\[
\sum\limits_{\left\vert \alpha\right\vert =k}y^{\alpha}\sum\limits_{\beta
\leq\alpha}\frac{x^{\beta}}{\beta!}=\sum\limits_{\left\vert \alpha\right\vert
=k}\sum\limits_{\left\vert \beta\right\vert \leq k}q_{\alpha,\beta}%
\frac{y^{\alpha}x^{\beta}}{\beta!},
\]

where%
\[
q_{\alpha,\beta}=\left\{
\begin{array}
[c]{ll}%
1, & \beta\leq\alpha,\\
0, & otherwise,
\end{array}
\right.
\]

and then change the order of summation so that%
\[
\sum\limits_{\left\vert \alpha\right\vert =k}y^{\alpha}\sum\limits_{\beta
\leq\alpha}\frac{x^{\beta}}{\beta!}=\sum\limits_{\left\vert \beta\right\vert
\leq k}\frac{x^{\beta}}{\beta!}\sum\limits_{\left\vert \alpha\right\vert
=k}q_{\alpha,\beta}y^{\alpha}.
\]

Now given a value of $\beta$ for what values of $\alpha$ is $q_{\alpha,\beta
}=1$? This is answered by the calculation%
\begin{align*}
\left\{  \alpha:q_{\alpha,\beta}=1\right\}   & =\left\{  \alpha:\left\vert
\alpha\right\vert =k,\text{ }\alpha\geq\beta\right\}  =\left\{  \beta
+\gamma:\left\vert \beta+\gamma\right\vert =k\right\}  =\\
& =\left\{  \beta+\gamma:\left\vert \beta\right\vert +\left\vert
\gamma\right\vert =k\right\}  =\left\{  \beta+\gamma:\left\vert \gamma
\right\vert =k-\left\vert \beta\right\vert \right\}  =\\
& =\beta+\left\{  \gamma:\left\vert \gamma\right\vert =k-\left\vert
\beta\right\vert \right\}  ,
\end{align*}

so that%
\[
\sum\limits_{\left\vert \alpha\right\vert =k}y^{\alpha}\sum\limits_{\beta
\leq\alpha}\frac{x^{\beta}}{\beta!}=\sum\limits_{\left\vert \beta\right\vert
\leq k}\frac{x^{\beta}}{\beta!}\sum\limits_{\left\vert \gamma\right\vert
=k-\left\vert \beta\right\vert }y^{\beta+\gamma}=\sum\limits_{\left\vert
\beta\right\vert \leq k}\frac{x^{\beta}y^{\beta}}{\beta!}\sum
\limits_{\left\vert \gamma\right\vert =k-\left\vert \beta\right\vert
}y^{\gamma}=\sum\limits_{\left\vert \beta\right\vert \leq k}\frac{x^{\beta
}y^{\beta}}{\beta!}h_{k-\left\vert \beta\right\vert }\left(  y\right)  .
\]

Further%
\[
\sum\limits_{\left\vert \beta\right\vert \leq k}\frac{x^{\beta}y^{\beta}%
}{\beta!}h_{k-\left\vert \beta\right\vert }\left(  y\right)  =\sum
\limits_{j=0}^{k}\sum\limits_{\left\vert \beta\right\vert =j}\frac{x^{\beta
}y^{\beta}}{\beta!}h_{k-\left\vert \beta\right\vert }\left(  y\right)
=\sum\limits_{j=0}^{k}h_{k-j}\left(  y\right)  \sum\limits_{\left\vert
\beta\right\vert =j}\frac{x^{\beta}y^{\beta}}{\beta!}=\sum\limits_{j=0}%
^{k}h_{k-j}\left(  y\right)  \frac{\left(  xy\right)  ^{j}}{j!},
\]

and so%
\[
\sum\limits_{\left\vert \alpha\right\vert =k}y^{\alpha}\sum\limits_{\beta
\leq\alpha}\frac{x^{\beta}}{\beta!}=\sum\limits_{j=0}^{k}h_{k-j}\left(
y\right)  \frac{\left(  xy\right)  ^{j}}{j!}.
\]

\textbf{Part 4} If $x\widehat{y}\geq0$ then from part 3%
\begin{align*}
\sum\limits_{\left\vert \alpha\right\vert =k}y^{\alpha}\sum\limits_{\beta
\leq\alpha}\frac{x^{\beta}}{\beta!}  & =\left\vert y\right\vert ^{k}%
\sum\limits_{j=0}^{k}h_{k-j}\left(  \widehat{y}\right)  \frac{\left(
x\widehat{y}\right)  ^{j}}{j!}\\
& \leq\left\vert y\right\vert ^{k}\sum\limits_{j=0}^{k}\left(  \sum
\limits_{\left\vert \alpha\right\vert =k-j}1\right)  \frac{\left(
x\widehat{y}\right)  ^{j}}{j!}\\
& =\left\vert y\right\vert ^{k}\sum\limits_{j=0}^{k}\binom{d+k-j-1}{k-j}%
\frac{\left(  x\widehat{y}\right)  ^{j}}{j!}\\
& =\left\vert y\right\vert ^{k}\sum\limits_{j=0}^{k}\frac{\left(
d+k-j-1\right)  !}{\left(  d-1\right)  !\left(  k-j\right)  !}\frac{\left(
x\widehat{y}\right)  ^{j}}{j!}\\
& =\left\vert y\right\vert ^{k}\sum\limits_{j=0}^{k}\frac{\left(
d+k-j-1\right)  !}{\left(  d-1\right)  !k!}\frac{k!}{\left(  k-j\right)
!j!}\left(  x\widehat{y}\right)  ^{j}\\
& \leq\left\vert y\right\vert ^{k}\sum\limits_{j=0}^{k}\frac{\left(
d+k-1\right)  !}{\left(  d-1\right)  !k!}\frac{k!}{\left(  k-j\right)
!j!}\left(  x\widehat{y}\right)  ^{j}\\
& =\binom{d+k-1}{k}\left\vert y\right\vert ^{k}\sum\limits_{j=0}^{k}\binom
{k}{j}\left(  x\widehat{y}\right)  ^{j}\\
& =\binom{d+k-1}{k}\left\vert y\right\vert ^{k}\left(  1+x\widehat{y}\right)
^{k}.
\end{align*}
\medskip

\textbf{Part 5} From part 3,%
\begin{align*}
\sum\limits_{\left\vert \alpha\right\vert \leq m}y^{\alpha}\sum\limits_{\beta
\leq\alpha}\frac{x^{\beta}}{\beta!}=\sum\limits_{k=0}^{m}\sum
\limits_{\left\vert \alpha\right\vert =k}y^{\alpha}\sum\limits_{\beta
\leq\alpha}\frac{x^{\beta}}{\beta!} &  =\sum\limits_{k=0}^{m}\sum
\limits_{j=0}^{k}h_{k-j}\left(  y\right)  \frac{\left(  xy\right)  ^{j}}{j!}\\
&  =\sum\limits_{j=0}^{m}\sum\limits_{k=j}^{m}h_{k-j}\left(  y\right)
\frac{\left(  xy\right)  ^{j}}{j!}\\
&  =\sum\limits_{j=0}^{m}\sum\limits_{k=j}^{m}h_{j}\left(  y\right)
\frac{\left(  xy\right)  ^{k-j}}{\left(  k-j\right)  !}.
\end{align*}

Thus%
\[
\sum\limits_{\left\vert \alpha\right\vert \leq m}y^{\alpha}\sum\limits_{\beta
\leq\alpha}\frac{x^{\beta}}{\beta!}=\sum\limits_{j=0}^{m}\left(
\sum\limits_{k=0}^{m-j}h_{k}\left(  y\right)  \right)  \frac{\left(
xy\right)  ^{j}}{j!}=\sum\limits_{j=0}^{m}h_{j}\left(  y\right)
\sum\limits_{k=0}^{m-j}\frac{\left(  xy\right)  ^{k}}{k!}.
\]
\medskip

\textbf{Part 6} ?? No component of $\widehat{y}$ dominates so
intuitively\ldots
\end{proof}

\begin{theorem}
\label{Thm_double_sum} ?? \textbf{CHECK}!%
\[
\sum\limits_{\left\vert \alpha\right\vert =k}\sum\limits_{\beta\leq\alpha
}b_{\beta}^{\left(  \alpha\right)  }\frac{x^{2\beta}}{\beta!}=\sum
\limits_{\left\vert \gamma\right\vert \leq k}\left(  \sum\limits_{\left\vert
\sigma\right\vert =k-\left\vert \gamma\right\vert }b_{\gamma}^{\left(
\gamma+\sigma\right)  }\right)  \frac{x^{2\gamma}}{\gamma!}.
\]

\end{theorem}

\begin{proof}
Since $\left\vert \beta\right\vert \leq\left\vert \alpha\right\vert $, there
exist constants $\eta_{\gamma}$ such that%
\[
\sum\limits_{\left\vert \alpha\right\vert =k}\sum\limits_{\beta\leq\alpha
}b_{\beta}^{\left(  \alpha\right)  }\frac{x^{2\beta}}{\beta!}=\sum
\limits_{\left\vert \gamma\right\vert \leq k}\eta_{\gamma}\frac{x^{2\gamma}%
}{\gamma!}.
\]

Now%
\begin{align*}
\left\{  \alpha:\left\vert \alpha\right\vert =k\text{ }and\text{ }\gamma
\leq\alpha\right\}   & =\left\{  \alpha:\left\vert \alpha\right\vert =k\text{
}and\text{ }\alpha\geq\gamma\right\} \\
& =\left\{  \gamma+\sigma:\left\vert \gamma+\sigma\right\vert =k\right\} \\
& =\gamma+\left\{  \sigma:\left\vert \gamma+\sigma\right\vert =k\right\} \\
& =\gamma+\left\{  \sigma:\left\vert \gamma\right\vert +\left\vert
\sigma\right\vert =k\right\} \\
& =\gamma+\left\{  \sigma:\left\vert \sigma\right\vert =k-\left\vert
\gamma\right\vert \right\}  ,
\end{align*}

and so%
\[
\sum\limits_{\left\vert \alpha\right\vert =k}\sum\limits_{\beta\leq\alpha
}b_{\beta}^{\left(  \alpha\right)  }\frac{x^{2\beta}}{\beta!}=\sum
\limits_{\left\vert \gamma\right\vert \leq k}\left(  \sum\limits_{\left\vert
\sigma\right\vert =k-\left\vert \gamma\right\vert }b_{\gamma}^{\left(
\gamma+\sigma\right)  }\right)  \frac{x^{2\gamma}}{\gamma!},
\]

which implies that%
\[
\eta_{\gamma}=\sum\limits_{\left\vert \sigma\right\vert =k-\left\vert
\gamma\right\vert }b_{\gamma}^{\left(  \gamma+\sigma\right)  }.
\]

\end{proof}

\section{Tempered distributions\label{Sect_tempered_distrib}}

\begin{definition}
\label{Def_Distributions}\textbf{Tempered distributions} $S^{\prime}$\ (also
temperate or generalized functions of slow growth)

\begin{enumerate}
\item $S$ is the space of rapidly decreasing functions. We endow $S$ with the
topology defined using the countable set of seminorms $p_{n,\alpha}$ given by%
\begin{equation}
p_{n,\alpha}\left(  \psi\right)  =\left\Vert \left(  1+\left\vert
\cdot\right\vert \right)  ^{n}D^{\alpha}\psi\right\Vert _{\infty},\text{\quad
}n=0,1,2,\ldots;\text{ }\alpha\geq0.\label{a68}%
\end{equation}

\item $S^{\prime}$ is the space of tempered distributions or generalized
functions of slow growth. It is the set of all continuous linear functionals
on $S$ under the seminorm topology of part 1 of this definition.

\item A linear functional on $S$ is continuous if its absolute value is
bounded by a finite positive linear combination of seminorms.

\item If $f\in S^{\prime}$ and $\phi\in S$ then $\left[  f,\phi\right]
\in\mathbb{C}$ will represent the action of $f$ on the test function $\phi$.
\end{enumerate}
\end{definition}

\subsection{Properties\label{SbSect_property_S'}\ }

\begin{enumerate}
\item The continuous embeddings $S\subset C_{BP}^{\infty}\subset S^{\prime
}\subset\mathcal{D}^{\prime}$ are dense.

\item If $f\in L_{loc}^{1}$ and, $\int\left\vert f\left(  x\right)
\right\vert \left(  1+\left\vert x\right\vert \right)  ^{-\lambda}dx$ or
$\int_{\left\vert \cdot\right\vert \geq r}\left\vert f\left(  x\right)
\right\vert $ $\left\vert x\right\vert ^{-\lambda}dx$ exist for some
$\lambda,r\geq0$, then $f\in S^{\prime}$ with action $\left[  f,\phi\right]
=\int f\left(  x\right)  \phi\left(  x\right)  dx$, $\phi\in S$. We call $f$ a
\textit{regular tempered distribution (function)}\textbf{\ }e.g. Example
2.8.3(a) Vladimirov \cite{Vladimirov}.

\item If $f\in S^{\prime}$ and $a\in C_{BP}^{\infty}$ then $af\in S^{\prime}$
e.g. Example 2.8.3(e) Vladimirov \cite{Vladimirov} but note that Vladimirov
denotes $C_{BP}^{\infty}$ by $\theta_{M}$.
\end{enumerate}

\section{Fourier Transforms}

\begin{definition}
\label{Def_Fourier_transform}\textbf{Fourier and Inverse transforms on }%
$S$\textbf{\ and }$S^{\prime}$\textbf{.}

\begin{enumerate}
\item If $f\in S$ then let the Fourier transform be
\[
\widehat{f}(\xi)=\left(  2\pi\right)  ^{-d/2}\int\limits_{\mathbb{R}^{d}%
}e^{-ix\xi}f(x)dx.
\]

Where convenient the alternative notation $F\left[  f\right]  =\widehat{f}$
will be used.

The mapping $f\rightarrow\widehat{f}$ is continuous from $S\rightarrow S$.

The inverse Fourier transform is
\[
\overset{\vee}{f}(\xi)=\left(  2\pi\right)  ^{-d/2}\int\limits_{\mathbb{R}%
^{d}}e^{ix\xi}f(x)dx.
\]

Where convenient the alternative notation $F^{-1}\left[  f\right]
=\overset{\vee}{f}$ will be used.

We have the important property that $\left(  \widehat{f}\right)  ^{\vee}=f$.

\item If $f\in S^{\prime}$ and $\phi\in S$ then $\left[  \widehat{f}%
,\phi\right]  =\left[  f,\overset{\vee}{\phi}\right]  $ and $\left[
\overset{\vee}{f},\phi\right]  =\left[  f,\overset{\vee}{\phi}\right]  $.

The mappings $f\rightarrow\widehat{f}$ and $f\rightarrow\overset{\vee}{f}$ are
continuous from $S^{\prime}\rightarrow S^{\prime}$.

We have the important property that $\left(  \widehat{f}\right)  ^{\vee}=f$.
\end{enumerate}
\end{definition}

\subsection{Fourier transform properties on $S$ and $S^{\prime}$}

\begin{enumerate}
\item If $b$ is a complex constant then $\left(  f\left(  x+b\right)  \right)
^{\wedge}=e^{ib\xi}\widehat{f}$.

\item $\overset{\vee}{f}(\xi)=\widehat{f}(-\xi)$, $\overline{\widehat{f}}%
(\xi)=\overset{\vee}{\bar{f}}(\xi)$ and $\widehat{\widehat{f}}(\xi)=f\left(
-\xi\right)  $.

\item If $b$ is a complex constant then $\left(  e^{ibx}\right)  ^{\wedge
}=\left(  2\pi\right)  ^{d/2}\delta\left(  \xi+b\right)  $.

\item $D^{\alpha}\widehat{f}=\left(  -i\right)  ^{\left\vert \alpha\right\vert
}\widehat{x^{\alpha}f}$ and $D^{\alpha}\overset{\vee}{f}=i^{\left\vert
\alpha\right\vert }\left(  x^{\alpha}f\right)  ^{\vee}$.

\item $\widehat{x^{\alpha}f}=i^{\left\vert \alpha\right\vert }D^{\alpha
}\widehat{f}$.

\item $\widehat{D^{\alpha}f}=i^{\left\vert \alpha\right\vert }\xi^{\alpha
}\widehat{f}=\left(  i\xi\right)  ^{\alpha}\widehat{f}$.

\item $\xi^{\alpha}\widehat{f}=\left(  -i\right)  ^{\left\vert \alpha
\right\vert }\widehat{D^{\alpha}f}$.

\item $\widehat{\delta}=\left(  2\pi\right)  ^{-d/2}$ and $\widehat{1}=\left(
2\pi\right)  ^{d/2}\delta$.

\item $\widehat{x^{\alpha}}=\left(  2\pi\right)  ^{d/2}\left(  iD\right)
^{\alpha}\delta$ and $\widehat{D^{\alpha}\delta}=\left(  2\pi\right)
^{-d/2}\left(  -i\right)  ^{\left\vert \alpha\right\vert }\xi^{\alpha}$.

\item If $p$ is a polynomial then $\widehat{p}=\left(  2\pi\right)
^{d/2}p(iD)\delta$, and $\widehat{p\left(  D\right)  f}=p\left(  -i\xi\right)
\widehat{f}$.
\end{enumerate}

\subsection{Fourier transform properties on $\mathcal{E}^{\prime}%
$\label{SbSect_property_E'}}

\begin{enumerate}
\item Suppose $f\in\mathcal{E}^{\prime}$, $\eta\in C_{0}^{\infty}$ and
$\eta=1$ in a neighborhood of the support of $f$ - denoted
$\operatorname*{supp}f$. Then $f$ can be extended uniquely to $S^{\prime}$ by
$\left[  f,\phi\right]  =\left[  f,\eta\phi\right]  $, where $\phi\in S$. e.g.
Example 2.8.3(b) Vladimirov \cite{Vladimirov} but note that Vladimirov denotes
$C_{BP}^{\infty}$ by $\theta_{M}$.

\item If $f\in\mathcal{E}^{\prime}$ then $f\in S^{\prime}$ in the sense of
part 1and $\widehat{f}\in C_{BP}^{\infty}$ e.g. Theorem 2.9.4 Vladimirov
\cite{Vladimirov}.
\end{enumerate}

\section{Convolutions}

\begin{definition}
\label{Def_Convolution}\textbf{Convolution}

\begin{enumerate}
\item If $f\in C_{BP}^{\left(  0\right)  }$ and $\phi\in S$ then $f\ast\phi\in
C_{BP}^{\infty}$ where
\[
\left(  f\ast\phi\right)  \left(  x\right)  =\left(  \phi\ast f\right)
\left(  x\right)  =\left(  2\pi\right)  ^{-\frac{d}{2}}\int f\left(
x-y\right)  \phi\left(  y\right)  dy=\left(  2\pi\right)  ^{-\frac{d}{2}}\int
f\left(  y\right)  \phi\left(  x-y\right)  dy.
\]

\item The previous definition can be extended to $S^{\prime}$ by the formulas
\[
f\ast\phi=\phi\ast f=\left(  \widehat{\phi}\widehat{f}\right)  ^{\vee}=\left(
2\pi\right)  ^{-\frac{d}{2}}\left[  f_{y},\phi\left(  \cdot-y\right)  \right]
,\quad f\in S^{\prime},\text{ }\phi\in S.
\]

Here $f\ast\phi\in C_{BP}^{\infty}$. Further, if $f\in S^{\prime}$ is defined
by the seminorms $\left\{  p_{n_{k},\alpha^{\left(  k\right)  }}\right\}  $
then (Theorem 6.2 Petersen \cite{Petersen83})%
\[
\left\vert \left(  \phi\ast f\right)  \left(  x\right)  \right\vert \leq
C\left(  1+\left\vert x\right\vert \right)  ^{\max n_{k}},\quad\phi\in
S,\text{ }x\in\mathbb{R}^{d}.
\]

\item Noting parts 1 and 3 of Appendix \ref{SbSect_property_S'}, the following
definition is consistent with parts 1 and 2:%
\[
f\ast g=g\ast f=\left(  \widehat{g}\widehat{f}\right)  ^{\vee},\quad f\in
S^{\prime},\text{ }g\in\left(  C_{BP}^{\infty}\right)  ^{\wedge}.
\]

Here $f\ast\phi\in S^{\prime}$.

\item Also, noting the results of Appendix \ref{SbSect_property_E'}%
\[
f\ast g=g\ast f=\left(  \widehat{g}\widehat{f}\right)  ^{\vee},\quad f\in
S^{\prime},\text{ }g\in\mathcal{E}^{\prime}.
\]

Here $f\ast g\in S^{\prime}$.

\item If $f\in C_{BP}^{\left(  0\right)  }$ and $g\in C_{0}^{\left(  0\right)
}$ then $f\ast g$ is a regular tempered distribution and%
\[
\left(  f\ast g\right)  \left(  x\right)  =\left(  g\ast f\right)  \left(
x\right)  =\left(  2\pi\right)  ^{-\frac{d}{2}}\int f\left(  x-y\right)
g\left(  y\right)  dy=\left(  2\pi\right)  ^{-\frac{d}{2}}\int f\left(
y\right)  g\left(  x-y\right)  dy.
\]

This can be proved using Theorem 2.7.5 of Vladimirov \cite{Vladimirov}.

\item If $f\in L^{p}$, $g\in L^{q}$ satisfy $\frac{1}{p}+\frac{1}{q}%
=1+\frac{1}{r}$ and $1\leq p,q,r\leq\infty$, then $f\ast g\in L^{r}$ where the
convolution is given by the formulas of part 5.

Further we have \textbf{Young's inequality}%
\begin{equation}
\left\Vert f\ast g\right\Vert _{r}\leq\left\Vert f\right\Vert _{p}\left\Vert
g\right\Vert _{q}.\label{1.056}%
\end{equation}

\end{enumerate}
\end{definition}

\subsection{Convolution properties\label{SbSect_property_Convol}}

For the convolutions defined in parts 2 and 4 above:

\begin{enumerate}
\item $f\ast g=g\ast f$.

\item $D^{\alpha}\left(  f\ast g\right)  =f\ast D^{\alpha}g=D^{\alpha}f\ast g$
for all $\alpha$.

\item $\left(  f\ast g\right)  ^{\wedge}=\widehat{g}\widehat{f}$,$\quad$where
$\widehat{f}\in S^{\prime}$.

\item The mapping $f\rightarrow f\ast g$ is continuous from $S^{\prime}$ to
$S^{\prime}$.
\end{enumerate}

\section{Taylor series expansion\label{Sect_taylor_expansion}}

Suppose $u:\mathbb{R}^{d}\rightarrow\mathbb{C}$ and $u\in C^{\left(  n\right)
}\left(  \mathbb{R}^{d}\right)  $ for some $n\geq1$. Then the Taylor series
expansion about $z\in\mathbb{R}^{d}$ is given by $\frac{\left\vert
b\right\vert ^{n}}{n!}$%
\begin{align}
u\left(  z+b\right)   & =\sum_{\left\vert \beta\right\vert <n}\frac{b^{\beta}%
}{\beta!}\left(  D^{\beta}u\right)  \left(  z\right)  +\left(  \mathcal{R}%
_{n}u\right)  \left(  z,b\right) \label{p87}\\
& =\sum_{k<n}\frac{1}{k!}\left(  \left(  bD\right)  ^{k}u\right)  \left(
z\right)  +\left(  \mathcal{R}_{n}u\right)  \left(  z,b\right) \nonumber\\
& =\sum_{k<n}\frac{\left\vert b\right\vert ^{k}}{k!}\left(  \left(
\widehat{b}D\right)  ^{k}u\right)  \left(  z\right)  +\left(  \mathcal{R}%
_{n}u\right)  \left(  z,b\right)  ,\nonumber
\end{align}

where $\widehat{b}=\left\vert b\right\vert ^{-1}b$ and $\mathcal{R}_{n}u$ is
the integral remainder term
\begin{align}
\left(  \mathcal{R}_{n}u\right)  \left(  z,b\right)   & =n\sum_{\left\vert
\beta\right\vert =n}\frac{b^{\beta}}{\beta!}\int_{0}^{1}s^{n-1}\left(
D^{\beta}u\right)  \left(  z+\left(  1-s\right)  b\right)  ds\label{p31}\\
& =n\sum_{\left\vert \beta\right\vert =n}\frac{b^{\beta}}{\beta!}\int_{0}%
^{1}\left(  1-s\right)  ^{n-1}\left(  D^{\beta}u\right)  \left(  z+sb\right)
ds\label{p88}\\
& =\frac{1}{\left(  n-1\right)  !}\int_{0}^{1}\left(  1-s\right)
^{n-1}\left(  \left(  bD\right)  ^{n}u\right)  \left(  z+sb\right)
ds\nonumber\\
& =\frac{\left\vert b\right\vert ^{n}}{\left(  n-1\right)  !}\int_{0}%
^{1}\left(  1-s\right)  ^{n-1}\left(  \left(  \widehat{b}D\right)
^{n}u\right)  \left(  z+sb\right)  ds.\label{Ap130}%
\end{align}

We have the estimates%
\begin{equation}
\left\vert \left(  \mathcal{R}_{n}u\right)  \left(  z,b\right)  \right\vert
\leq\frac{\left\vert b\right\vert ^{n}}{n!}\max_{\substack{\left\vert
\beta\right\vert =n \\y\in\left[  z,z+b\right]  }}\left\vert \left(
\widehat{b}D\right)  ^{n}u\left(  y\right)  \right\vert ,\label{Ap147}%
\end{equation}

\begin{equation}
\left\vert \left(  \mathcal{R}_{n}u\right)  \left(  z,b\right)  \right\vert
\leq\frac{\left\vert b\right\vert _{1}^{n}}{n!}\max_{\substack{\left\vert
\beta\right\vert =n \\y\in\left[  z,z+b\right]  }}\left\vert \left(  D^{\beta
}u\right)  \left(  y\right)  \right\vert \leq d^{\frac{n}{2}}\frac{\left\vert
b\right\vert ^{n}}{n!}\max_{\substack{\left\vert \beta\right\vert =n
\\y\in\left[  z,z+b\right]  }}\left\vert \left(  D^{\beta}u\right)  \left(
y\right)  \right\vert ,\label{p34}%
\end{equation}

\begin{equation}
\left\vert \left(  \mathcal{R}_{n}u\right)  \left(  z,b\right)  \right\vert
\leq\sum_{\left\vert \beta\right\vert =n}\frac{\left\vert b^{\beta}\right\vert
}{\beta!}\max_{y\in\left[  z,z+b\right]  }\left\vert \left(  D^{\beta
}u\right)  \left(  y\right)  \right\vert \leq\left\vert b\right\vert ^{n}%
\sum_{\left\vert \beta\right\vert =n}\frac{1}{\beta!}\max_{y\in\left[
z,z+b\right]  }\left\vert \left(  D^{\beta}u\right)  \left(  y\right)
\right\vert ,\label{p20}%
\end{equation}

where we have used the notation $\left[  z,z+b\right]  =\left\{
z+sb:s\in\left[  0,1\right]  \right\}  $.

\begin{theorem}
\label{Thm_Taylor_rad_remain}Regarding \ref{Ap130}:

\begin{enumerate}
\item When $z=0$ we can write%
\begin{equation}
\left(  \mathcal{R}_{n}u\right)  \left(  0,b\right)  =\frac{\left\vert
b\right\vert ^{n}}{\left(  n-1\right)  !}\int_{0}^{1}\left(  1-s\right)
^{n-1}\left(  \left(  \widehat{\cdot}D\right)  ^{n}u\right)  \left(
sb\right)  ds,\label{Ap134}%
\end{equation}

with the bound%
\begin{equation}
\left\vert \left(  \mathcal{R}_{n}u\right)  \left(  0,b\right)  \right\vert
\leq\frac{\left\vert b\right\vert ^{n}}{n!}\max_{s\in\left[  0,1\right]
}\left\vert \left(  \left(  \widehat{\cdot}D\right)  ^{n}u\right)  \left(
sb\right)  \right\vert .\label{Ap144}%
\end{equation}

When $u$ is radial $\left(  \mathcal{R}_{n}u\right)  \left(  0,\cdot\right)  $
is also radial. In fact, if $u_{\odot}\left(  \left\vert \cdot\right\vert
\right)  =u\left(  \cdot\right)  $ then%
\begin{equation}
\left(  \mathcal{R}_{n}u\right)  \left(  0,b\right)  =\frac{\left\vert
b\right\vert ^{n}}{\left(  n-1\right)  !}\int_{0}^{1}\left(  1-s\right)
^{n-1}\left(  D^{n}u_{\odot}\right)  \left(  s\left\vert b\right\vert \right)
ds,\label{Ap135}%
\end{equation}

with upper bound%
\begin{equation}
\left\vert \left(  \mathcal{R}_{n}u\right)  \left(  0,b\right)  \right\vert
\leq\frac{\left\vert b\right\vert ^{n}}{n!}\max_{s\in\left[  0,1\right]
}\left\vert \left(  D^{n}u_{\odot}\right)  \left(  sb\right)  \right\vert
.\label{Ap146}%
\end{equation}

\item If $u$ is radial then $\frac{1}{n!}\left(  \left(  bD\right)
^{n}u\right)  \left(  0\right)  =\sum\limits_{\left\vert \alpha\right\vert
=n}\frac{b^{\alpha}}{\alpha!}D^{\alpha}u\left(  0\right)  $ is a radial
function of $b$ for $n\geq0$.

\item In fact, if $u$ is radial then
\begin{equation}
\sum\limits_{\left\vert \alpha\right\vert =n}\frac{b^{\alpha}}{\alpha
!}D^{\alpha}u\left(  0\right)  =\frac{\left\vert b\right\vert ^{n}}{n!}%
D^{n}u_{\odot}\left(  0\right)  =\frac{\left\vert b\right\vert ^{n}}%
{n!}\left(  \left(  \widehat{\cdot}D\right)  ^{n}u\right)  \left(  0\right)
.\label{Ap133}%
\end{equation}

\end{enumerate}
\end{theorem}

\begin{proof}
\textbf{Part 1} When $z=0$ \ref{Ap130} becomes%
\begin{align*}
\left(  \mathcal{R}_{n}u\right)  \left(  0,b\right)   & =\frac{\left\vert
b\right\vert ^{n}}{\left(  n-1\right)  !}\int_{0}^{1}\left(  1-t\right)
^{n-1}\left(  \left(  \widehat{b}D\right)  ^{n}u\right)  \left(  tb\right)
ds\\
& =\frac{\left\vert b\right\vert ^{n}}{\left(  n-1\right)  !}\int_{0}%
^{1}\left(  1-s\right)  ^{n-1}\left(  \left(  \widehat{sb}D\right)
^{n}u\right)  \left(  sb\right)  ds\\
& =\frac{\left\vert b\right\vert ^{n}}{\left(  n-1\right)  !}\int_{0}%
^{1}\left(  1-s\right)  ^{n-1}\left(  \left(  \widehat{\cdot}D\right)
^{n}u\right)  \left(  sb\right)  ds.
\end{align*}

Now suppose $u\left(  \cdot\right)  =u_{\odot}\left(  \left\vert
\cdot\right\vert \right)  $. Then from \ref{Ap131}, $\left(  \left(
\widehat{\cdot}D\right)  ^{n}u\right)  \left(  x\right)  =\left(
D^{n}u_{\odot}\right)  \left(  \left\vert x\right\vert \right)  $ and so%
\[
\left(  \mathcal{R}_{n}u\right)  \left(  0,b\right)  =\frac{\left\vert
b\right\vert ^{n}}{\left(  n-1\right)  !}\int_{0}^{1}\left(  1-s\right)
^{n-1}\left(  D^{n}u_{\odot}\right)  \left(  s\left\vert b\right\vert \right)
ds
\]
$\medskip$

\textbf{Part 2} From \ref{p87} we have\textbf{\ }%
\begin{align*}
u\left(  b\right)   & =\sum_{\left\vert \alpha\right\vert \leq n-1}%
\frac{b^{\alpha}}{\alpha!}\left(  D^{\alpha}u\right)  \left(  0\right)
+\left(  \mathcal{R}_{n}u\right)  \left(  0,b\right)  ,\\
u\left(  b\right)   & =\sum_{\left\vert \alpha\right\vert \leq n-1}%
\frac{b^{\alpha}}{\alpha!}\left(  D^{\alpha}u\right)  \left(  0\right)
+\sum_{\left\vert \alpha\right\vert =n}\frac{b^{\alpha}}{\alpha!}\left(
D^{\alpha}u\right)  \left(  0\right)  +\left(  \mathcal{R}_{n+1}u\right)
\left(  0,b\right)  ,
\end{align*}

so that by inspection%
\begin{equation}
\sum_{\left\vert \alpha\right\vert =n}\frac{b^{\alpha}}{\alpha!}\left(
D^{\alpha}u\right)  \left(  0\right)  =\left(  \mathcal{R}_{n}u\right)
\left(  0,b\right)  -\left(  \mathcal{R}_{n+1}u\right)  \left(  0,b\right)
.\label{Ap132}%
\end{equation}

This result follows since from part 1 we know that $\left(  \mathcal{R}%
_{n+1}u\right)  \left(  0,b\right)  $ and $\left(  \mathcal{R}_{n}u\right)
\left(  0,b\right)  $ are both radial.\medskip

\textbf{Part 3} Equation \ref{Ap132} can be written as%
\[
\frac{1}{n!}\left(  \left(  bD\right)  ^{n}u\right)  \left(  0\right)
=\left(  \mathcal{R}_{n}u\right)  \left(  0,b\right)  -\left(  \mathcal{R}%
_{n+1}u\right)  \left(  0,b\right)  ,
\]

and using \ref{Ap135} and applying integration by parts gives
\begin{align*}
& \left(  \mathcal{R}_{n}u\right)  \left(  0,b\right)  -\left(  \mathcal{R}%
_{n+1}u\right)  \left(  0,b\right) \\
& =\frac{1}{\left(  n-1\right)  !}\int_{0}^{1}\left(  1-s\right)
^{n-1}\left(  D^{n}u_{\odot}\right)  \left(  s\right)  ds-\frac{1}{n!}\int%
_{0}^{1}\left(  1-s\right)  ^{n}\left(  D^{n+1}u_{\odot}\right)  \left(
s\right)  ds\\
& =\frac{1}{\left(  n-1\right)  !}\int_{0}^{1}\left(  1-s\right)
^{n-1}\left(  D^{n}u_{\odot}\right)  \left(  s\right)  ds-\\
& \qquad\qquad-\frac{1}{n!}\left(  \left[  \left(  1-s\right)  ^{n}%
D^{n}u_{\odot}\left(  s\right)  \right]  _{0}^{1}-\int_{0}^{1}\left(
D_{s}\left(  1-s\right)  ^{n}\right)  D^{n}u_{\odot}\left(  s\right)
ds\right) \\
& =\frac{1}{\left(  n-1\right)  !}\int_{0}^{1}\left(  1-s\right)
^{n-1}\left(  D^{n}u_{\odot}\right)  \left(  s\right)  ds-\\
& \qquad\qquad-\frac{1}{n!}\left(  -D^{n}u_{\odot}\left(  0\right)  +\int%
_{0}^{1}\left(  1-s\right)  ^{n-1}D^{n}u_{\odot}\left(  s\right)  ds\right) \\
& =\frac{1}{n!}D^{n}u_{\odot}\left(  0\right)  .
\end{align*}

The second equality of \ref{Ap133} follows from part 1 of Lemma
\ref{Lem_deriv_rad_funcs}.
\end{proof}

\begin{theorem}
\label{Thm_Tay_rem_zeros}Suppose $u\in C_{B}^{\left(  n\right)  }\left(
\mathbb{R}^{d}\right)  $ and $D^{\alpha}u\left(  0\right)  =0$ when
$\left\vert \alpha\right\vert <n$. Then for all $r>0$,%
\begin{equation}
\left\vert u\left(  x\right)  \right\vert \leq\frac{\left\vert x\right\vert
^{n}}{n!}\left\Vert \left(  \widehat{\cdot}D\right)  ^{n}u\right\Vert
_{\infty;\overline{B}_{r}},\quad\left\vert x\right\vert \leq r,\label{Ap019}%
\end{equation}

where $\left(  \left(  \widehat{\cdot}D\right)  u\right)  \left(  x\right)
:=\sum\limits_{k=1}^{d}\widehat{x}_{k}D_{k}u\left(  x\right)  =\frac
{1}{\left\vert x\right\vert }\sum\limits_{k=1}^{d}x_{k}D_{k}u\left(  x\right)
$. Further%
\begin{equation}
\left\vert u\left(  x\right)  \right\vert \leq\frac{\left\vert x\right\vert
^{n}}{n!}\left\Vert \left(  \widehat{\cdot}D\right)  ^{n}u\right\Vert
_{\infty},\quad x\in\mathbb{R}^{d}.\label{Ap01}%
\end{equation}

Things are especially nice when $u$ is radial. In fact, if $u\left(  x\right)
=u_{\odot}\left(  \left\vert x\right\vert \right)  $ then $\left(
\widehat{\cdot}D\right)  ^{n}u$ is radial and%
\begin{equation}
\left(  \left(  \widehat{\cdot}D\right)  ^{n}u\right)  \left(  x\right)
=\left(  D^{n}u_{\odot}\right)  \left(  \left\vert x\right\vert \right)
.\label{Ap027}%
\end{equation}

\end{theorem}

\begin{proof}
Inequality \ref{Ap144} with $b=x$ and $\left\vert x\right\vert \leq r$ gives%
\begin{align*}
\left\vert u\left(  x\right)  \right\vert  & \leq\frac{\left\vert x\right\vert
^{n}}{n!}\max_{t\in\left[  0,1\right]  }\left\vert \left(  \left(  \widehat
{x}D\right)  ^{n}u\right)  \left(  tx\right)  \right\vert =\frac{\left\vert
x\right\vert ^{n}}{n!}\max_{t\in\left[  0,1\right]  }\left\vert \left(
\left(  \widehat{tx}D\right)  ^{n}u\right)  \left(  tx\right)  \right\vert =\\
& =\frac{\left\vert x\right\vert ^{n}}{n!}\max_{t\in\left[  0,1\right]
}\left\vert \left(  \left(  \widehat{\cdot}D\right)  ^{n}u\right)  \left(
tx\right)  \right\vert \leq\frac{\left\vert x\right\vert ^{n}}{n!}\max
_{z\in\overline{B}_{r}}\left\vert \left(  \left(  \widehat{\cdot}D\right)
^{n}u\right)  \left(  z\right)  \right\vert \leq\\
& =\frac{\left\vert x\right\vert ^{n}}{n!}\left\Vert \left(  \widehat{\cdot
}D\right)  ^{n}u\right\Vert _{\infty;\leq r}\leq\frac{\left\vert x\right\vert
^{n}}{n!}\left\Vert \left(  \widehat{\cdot}D\right)  ^{n}u\right\Vert
_{\infty}.
\end{align*}

Equation \ref{Ap027} follows directly from \ref{Ap131}.
\end{proof}

The next result may or may not be useful.

\begin{theorem}
\label{Thm_Taylor_rem_1divBk}Set%
\[
f_{n}\left(  \beta\right)  =\max_{y\in\left[  z,z+b\right]  }\left\vert
\left(  D^{\beta}u\right)  \left(  y\right)  \right\vert .
\]

Then%
\[
\left\vert \left(  \mathcal{R}_{n}u\right)  \left(  z,b\right)  \right\vert
\leq\left\vert b\right\vert ^{n}\frac{\left(  \chi_{n}d\right)  ^{n}}{n!},
\]

where%
\[
\chi_{n}=\max_{\left\vert \beta\right\vert =n}\max_{k}\left\vert f_{n}\left(
\beta\right)  \right\vert ^{\frac{1}{d\beta_{k}}}=\left(  \max_{\left\vert
\beta\right\vert =n}\max_{k}\left\vert f_{n}\left(  \beta\right)  \right\vert
^{\frac{1}{\beta_{k}}}\right)  ^{\frac{1}{d}}.
\]

\end{theorem}

\begin{proof}
From \ref{p20},%
\begin{align*}
\left\vert \left(  \mathcal{R}_{n}u\right)  \left(  z,b\right)  \right\vert  &
\leq\left\vert b\right\vert ^{n}\sum\limits_{\left\vert \beta\right\vert
=n}\frac{1}{\beta!}\max_{y\in\left[  z,z+b\right]  }\left\vert \left(
D^{\beta}u\right)  \left(  y\right)  \right\vert =\left\vert b\right\vert
^{n}\sum\limits_{\left\vert \beta\right\vert =n}\frac{1}{\beta!}f_{n}\left(
\beta\right)  =\\
& =\left\vert b\right\vert ^{n}\sum\limits_{\left\vert \beta\right\vert
=n}\frac{1}{\beta!}\mathbf{1}^{\beta}\left(  \left\vert f_{n}\left(
\beta\right)  \right\vert ^{\frac{1}{d\beta_{k}}}\right)  ^{\beta}%
\leq\left\vert b\right\vert ^{n}\sum\limits_{\left\vert \beta\right\vert
=n}\frac{1}{\beta!}\mathbf{1}^{\beta}\left(  \chi_{n}\mathbf{1}\right)
^{\beta}=\\
& =\left\vert b\right\vert ^{n}\left(  \chi_{n}\right)  ^{n}\sum
\limits_{\left\vert \beta\right\vert =n}\frac{\mathbf{1}^{2\beta}}{\beta
!}=\left\vert b\right\vert ^{n}\left(  \chi_{n}\right)  ^{n}\frac{1}%
{n!}\left\vert \mathbf{1}\right\vert ^{2n}=\left\vert b\right\vert ^{n}%
d^{n}\frac{\left(  \chi_{n}\right)  ^{n}}{n!}=\\
& =\left\vert b\right\vert ^{n}\frac{\left(  \chi_{n}d\right)  ^{n}}{n!}.
\end{align*}

\end{proof}

\section{Approximating the gamma function}

A well known gamma function approximation is
\begin{equation}
\sqrt{2\pi}e^{-x+\frac{1}{1+12x}}x^{x+\frac{1}{2}}<\Gamma\left(  x+1\right)
<\sqrt{2\pi}e^{-x+\frac{1}{12x}}x^{x+\frac{1}{2}},\quad x\geq0,\label{Ap029}%
\end{equation}

\begin{align*}
1  & \leq\frac{\sin x}{x}\leq\frac{2}{\pi},\quad\quad\quad0\leq x\leq\pi/2.\\
1-\frac{1}{6}x^{2}  & \leq\frac{\sin x}{x}\leq1-\frac{1}{7}x^{2},\quad0\leq
x\leq\pi/2.
\end{align*}

and also that%
\begin{equation}
\Gamma\left(  x\right)  \Gamma\left(  1-x\right)  =\frac{\pi}{\sin\pi x}%
,\quad\Gamma\left(  1/2\right)  =\sqrt{\pi}.\label{Ap039}%
\end{equation}

Now $\Gamma\left(  x+1\right)  =x\Gamma\left(  x\right)  $ so \ref{Ap029}
becomes%
\begin{equation}
\sqrt{2\pi}e^{-x+\frac{1}{1+12x}}x^{x-\frac{1}{2}}<\Gamma\left(  x\right)
<\sqrt{2\pi}e^{-x+\frac{1}{12x}}x^{x-\frac{1}{2}},\quad x>0.\label{Ap049}%
\end{equation}

The estimates \ref{Ap049} are OK for $x\geq1/2$ ??
SUPPLY\ NUM\ EVIDENCE?.\medskip

\fbox{For $0<x<1/2$} we can use \ref{Ap039}. In fact%
\begin{align*}
x\Gamma\left(  x\right)   & =\frac{\pi x}{\sin\pi x}\frac{1}{\Gamma\left(
1-x\right)  }\leq\frac{\pi}{2}\frac{1}{\Gamma\left(  1-x\right)  }\leq
\frac{\pi}{2}\frac{1}{\Gamma\left(  1\right)  }=\frac{\pi}{2}.\\
x\Gamma\left(  x\right)   & =\frac{\pi x}{\sin\pi x}\frac{1}{\Gamma\left(
1-x\right)  }\geq\frac{1}{\Gamma\left(  1-x\right)  }\geq\frac{1}{\sqrt{\pi}}.
\end{align*}

From \ref{Ap049},%
\[
\frac{1}{\Gamma\left(  1-x\right)  }\leq\frac{1}{\sqrt{2\pi}e^{-\left(
1-x\right)  +\frac{1}{1+12\left(  1-x\right)  }}\left(  1-x\right)  ^{\left(
1-x\right)  -\frac{1}{2}}}\leq\frac{e^{1-x}}{\sqrt{2\pi}\left(  1-x\right)
^{\frac{1}{2}-x}}%
\]

Numerically%
\begin{align*}
x\Gamma\left(  x\right)  \leq\left(  x\Gamma\left(  x\right)  \right)  \left(
0\right)  -\frac{\left(  x\Gamma\left(  x\right)  \right)  \left(  0\right)
-\left(  x\Gamma\left(  x\right)  \right)  \left(  1/2\right)  }{1/2}x  &
=1-\frac{1-\sqrt{\pi}/2}{1/2}x\\
& =1-\left(  2-\sqrt{\pi}\right)  x.
\end{align*}

Also, by adjusting the formula%
\[
1-\gamma x+0.5\left(  \frac{1}{2}\gamma^{2}+\frac{\pi^{2}}{12}\right)
x^{2}\leq x\Gamma\left(  x\right)  \leq1-\gamma x+\left(  \frac{1}{2}%
\gamma^{2}+\frac{\pi^{2}}{12}\right)  x^{2},
\]

by numerical observation we obtain the tighter%
\[
1-\gamma x+0.7\left(  \frac{1}{2}\gamma^{2}+\frac{\pi^{2}}{12}\right)
x^{2}\leq x\Gamma\left(  x\right)  \leq1-\gamma x+\left(  \frac{1}{2}%
\gamma^{2}+\frac{\pi^{2}}{12}\right)  x^{2}-0.525x^{3}.
\]

\section{Application of general multivariate spherical
coordinates\label{Sect_App_SpherCoords}}

In this section multivariate spherical coordinates are applied to obtain
one-dimensional integrals for $\int_{\left\vert x\right\vert \leq r}f\left(
\left\vert x\right\vert \right)  dx$ (see \ref{Ap148}) and $\int_{\left\vert
x\right\vert \leq r}\left\vert x_{k}\right\vert ^{p}f\left(  \left\vert
x\right\vert \right)  dx$ (see \ref{Ap016}) and then derive the formula
\ref{Ap017} i.e.%
\[
\int_{\left\vert x\right\vert \leq r}\left\vert \widehat{\xi}x\right\vert
^{p}f\left(  \left\vert x\right\vert \right)  dx=\frac{B\left(  \frac{d}%
{2},\frac{1+p}{2}\right)  }{B\left(  \frac{1}{2},\frac{d+p}{2}\right)  }%
\int\limits_{\left\vert x\right\vert \leq r}\left\vert x\right\vert
^{p}f\left(  \left\vert x\right\vert \right)  dx,\quad p>-1,
\]

and then obtain the upper bounds \ref{Ap220} and \ref{Ap225} for
$\frac{B\left(  \frac{d}{2},\frac{1+p}{2}\right)  }{B\left(  \frac{1}{2}%
,\frac{d+p}{2}\right)  }$. These results are useful for obtaining upper bounds
of certain expressions involving radial weight functions e.g. Subsection
\ref{SbSect_Tay_dat_rem_rad_fn}.

Following Section 5.44 of Adams \cite{Adams75} we will first describe
spherical polar coordinates and then give two applications that will be used
in this document.

Suppose $x\in\mathbb{R}^{d}$ and $d\geq2$. The \textbf{spherical polar
coordinate representation} is%
\[
x=x\left(  \rho,\phi\right)  =x\left(  \rho,\phi_{1},\phi_{2},\ldots
,\phi_{d-1}\right)  ,
\]

where%
\begin{equation}%
\begin{array}
[c]{rr}%
x_{1}= & \rho\sin\phi_{1}\sin\phi_{2}\ldots\sin\phi_{d-1},\\
x_{2}= & \rho\cos\phi_{1}\sin\phi_{2}\ldots\sin\phi_{d-1},\\
x_{3}= & \rho\cos\phi_{2}\ldots\sin\phi_{d-1},\\
\vdots\quad & \vdots\qquad\\
x_{d}= & \rho\cos\phi_{d-1},
\end{array}
\label{Ap051}%
\end{equation}

and on $\mathbb{R}^{d}$,%
\[
\rho\geq0,\text{ }-\pi\leq\phi_{1}\leq\pi,\text{ }0\leq\phi_{2},\ldots
,\phi_{d-1}\leq\pi.
\]

The volume element is%
\[
dx=dx_{1}dx_{2}\ldots dx_{d}=\rho^{d-1}%
{\textstyle\prod\limits_{j=1}^{d-1}}
\sin^{j-1}\phi_{j}d\rho d\phi,
\]

where $d\phi=d\phi_{1}\cdots d\phi_{d-1}$.\medskip

\textbf{Integration in the ball }$B_{r}$:%
\begin{align}
&  \int\limits_{\left\vert x\right\vert \leq r}g\left(  x\right)  f\left(
\left\vert x\right\vert \right)  dx\nonumber\\
& \nonumber\\
&  =\left\{
\begin{array}
[c]{ll}%
\int\limits_{-r}^{r}g\left(  \rho\right)  f\left(  \left\vert \rho\right\vert
\right)  d\rho, & d=1,\\
& \\
\int_{-\pi}^{\pi}\int_{0}^{r}g\left(  \rho\sin\phi_{1},\rho\cos\phi
_{1}\right)  f\left(  \rho\right)  \rho d\rho d\phi_{1}, & d=2,\\
& \\
\int\limits_{0}^{\pi}\int\limits_{-\pi}^{\pi}\int\limits_{0}^{r}g\left(
\rho\sin\phi_{1}\sin\phi_{2},\rho\cos\phi_{1}\sin\phi_{2},\rho\cos\phi
_{2}\right)  f\left(  \rho\right)  \rho^{2}d\rho\sin\phi_{2}d\phi_{1}d\phi
_{2}, & d=3,\\
& \\
\underset{d-2\times}{\underbrace{\int\limits_{0}^{\pi}}}\int\limits_{-\pi
}^{\pi}\int\limits_{0}^{r}g\left(  x\right)  f\left(  \rho\right)  \rho
^{d-1}d\rho%
{\textstyle\prod\limits_{j=1}^{d-1}}
\sin^{j-1}\phi_{j}d\phi_{1}\ldots d\phi_{d-1}, & d\geq3,
\end{array}
\right. \label{Ap086}%
\end{align}

where%
\[
0\leq\rho\leq r,\text{ }-\pi\leq\phi_{1}\leq\pi,\text{ }0\leq\phi_{2}%
,\ldots,\phi_{d-1}\leq\pi.
\]

\textbf{Integration in the portion of the orthant} $\overline{B}_{r}%
\cap\left\{  x:x\geq\mathbf{0}\right\}  $:%
\begin{equation}%
\begin{array}
[c]{l}%
\int\limits_{\left\vert x\right\vert \leq r;x\geq\mathbf{0}}g\left(  x\right)
f\left(  \left\vert x\right\vert \right)  dx\\
\\
=\left\{
\begin{array}
[c]{ll}%
\int_{0}^{r}g\left(  \rho\right)  f\left(  \left\vert \rho\right\vert \right)
d\rho, & d=1,\\
& \\
\int_{0}^{\pi/2}\int_{0}^{r}g\left(  \rho\sin\phi_{1},\rho\cos\phi_{1}\right)
f\left(  \rho\right)  \rho d\rho d\phi_{1}, & d=2,\\
& \\
\int\limits_{0}^{\pi/2}\int\limits_{0}^{\pi/2}\int\limits_{0}^{r}g\left(
\rho\sin\phi_{1}\sin\phi_{2},\rho\cos\phi_{1}\sin\phi_{2},\rho\cos\phi
_{2}\right)  f\left(  \rho\right)  \rho^{2}d\rho\sin\phi_{2}d\phi_{1}d\phi
_{2}, & d=3,\\
& \\
\underset{d-1\times}{\underbrace{\int_{0}^{\pi/2}}}\int_{0}^{r}g\left(
x\right)  f\left(  \rho\right)  \rho^{d-1}d\rho%
{\textstyle\prod\limits_{j=1}^{d-1}}
\sin^{j-1}\phi_{j}d\phi_{1}\ldots d\phi_{d-1}, & d\geq2,
\end{array}
\right. \\
\smallskip\\
0\leq\rho\leq r,\text{ }0\leq\phi_{1},\phi_{2},\ldots,\phi_{d-1}\leq\frac{\pi
}{2}.
\end{array}
\label{Ap138}%
\end{equation}

?? MOST OF\ THE RESULTS\ OF\ THIS\ SECTION can be proved using this result
stated in Prudnikov \cite{PrudBryMar86}:

\begin{theorem}
\label{Thm_Integ_u(xy,|x|)dx}When $d\geq2$, from Prudnikov \cite{PrudBryMar86}%
,%
\begin{align*}
\int_{\left\vert x\right\vert \leq r}\Phi\left(  \xi x,\left\vert x\right\vert
\right)  dx  & =\omega_{d-1}\int_{0}^{r}\rho^{d-1}\int_{0}^{\pi}\Phi\left(
\left\vert \xi\right\vert \rho\cos\theta,\rho\right)  \sin^{d-2}\theta d\theta
d\rho,\\
\omega_{t}  & =\frac{2\pi^{t/2}}{\Gamma\left(  t/2\right)  },\text{ }t>0.
\end{align*}

Formula 858.515 of Dwight \cite{Dwight61} is%
\[%
\begin{array}
[c]{c}%
\int\limits_{0}^{\pi/2}\sin^{p}\theta\cos^{q}\theta d\theta=\frac{1}{2}%
\frac{\Gamma\left(  \frac{p+1}{2}\right)  \Gamma\left(  \frac{q+1}{2}\right)
}{\Gamma\left(  \frac{p+q}{2}+1\right)  }=\frac{1}{2}B\left(  \frac{p+1}%
{2},\frac{q+1}{2}\right)  ,\\
any\text{ }real\text{ }p,q>-\frac{1}{2}.
\end{array}
\]

so that%
\[
\int_{0}^{\pi}\sin^{d-2}\theta d\theta=\frac{\omega_{d}}{\omega_{d-1}}.
\]

Also from Prudnikov \cite{PrudBryMar86},%
\[
\int_{\left\vert x\right\vert \leq r}\Phi\left(  \left\vert x\right\vert
^{2},\xi x,\eta x\right)  dx=\omega_{d-2}\int\limits_{-r}^{r}\int%
\limits_{-\sqrt{r^{2}-t^{2}}}^{\sqrt{r^{2}-t^{2}}}\int\limits_{0}^{\sqrt
{r^{2}-t^{2}-u^{2}}}v^{d-3}\Phi\left(  t^{2}+u^{2}+v^{2},\left\vert
\xi\right\vert t,\beta t+\gamma u\right)  \text{ }dvdudt,
\]

where%
\begin{align*}
\beta & =\widehat{\xi}\eta=\left\vert \eta\right\vert \widehat{\xi}%
\widehat{\eta},\\
\gamma & =\frac{\left\vert \xi\right\vert ^{2}\left\vert \eta\right\vert
^{2}-\left\vert \xi\eta\right\vert ^{2}}{\left\vert \xi\right\vert ^{2}%
}=\left\vert \eta\right\vert ^{2}\left(  1-\left\vert \widehat{\xi}%
\widehat{\eta}\right\vert ^{2}\right)  =\left\vert \eta\right\vert ^{2}%
-\beta^{2}.
\end{align*}

\end{theorem}

Basic result:

\begin{theorem}
\label{Thm_Integ_u(xy)f(|x|)dx}Suppose $\xi,x\in\mathbb{R}^{d}$ and $\xi
x=\left(  \xi,x\right)  $ denotes the Euclidean inner product. Then:

\begin{enumerate}
\item $\int_{\left\vert x\right\vert \leq r}u\left(  \xi x\right)  f\left(
\left\vert x\right\vert \right)  dx$ is a radial function of $\xi$ and%
\begin{equation}
\int_{\left\vert x\right\vert \leq r}u\left(  \xi x\right)  f\left(
\left\vert x\right\vert \right)  dx=\int_{\left\vert x\right\vert \leq
r}u\left(  \left\vert \xi\right\vert x_{k}\right)  f\left(  \left\vert
x\right\vert \right)  dx,\quad1\leq k\leq d.\label{Ap016}%
\end{equation}

\item When $d\geq2$, $s\in\mathbb{R}^{1}$ and $1\leq k\leq d$,%
\begin{equation}
\int\limits_{\left\vert x\right\vert \leq r}u\left(  sx_{k}\right)  f\left(
\left\vert x\right\vert \right)  dx=\omega_{d-1}\int_{0}^{r}\left(  \int%
_{0}^{\pi}u\left(  s\rho\cos t\right)  \sin^{d-2}tdt\right)  f\left(
\rho\right)  \rho^{d-1}d\rho,\label{Ap001}%
\end{equation}

where $\omega_{k}=\frac{2\pi^{k/2}}{\Gamma\left(  k/2\right)  }$ for $k\geq1$.
\end{enumerate}
\end{theorem}

\begin{proof}
\textbf{Part 1} This lemma will employ the same technique as the previous
lemma: Section 4.1 Stein and Weiss \cite{SteinWeiss71} which defines radial
functions in terms of orthogonal transformations. Stein and Weiss observed
that a function $g$ is radial if and only if $g\left(  \mathcal{O}x\right)
=g\left(  x\right)  $ for any linear, orthogonal transformation $\mathcal{O}%
:B_{r}\rightarrow B_{r}$ and any $x\in\mathbb{R}^{d}$. Now an orthogonal
transformation $\mathcal{O}$ satisfies $\mathcal{O}^{T}=\mathcal{O}^{-1}$
where $\mathcal{O}xy=x\mathcal{O}^{T}y$ for the Euclidean inner product, and
an orthogonal transformation has a Jacobian of one.

Now
\begin{align*}
\int\limits_{\left\vert x\right\vert \leq r}u\left(  \mathcal{O}\xi x\right)
f\left(  \left\vert x\right\vert \right)  dx=\int\limits_{\left\vert
x\right\vert \leq r}u\left(  \xi\mathcal{O}^{-1}x\right)  f\left(  \left\vert
x\right\vert \right)  dx  & =\int\limits_{\left\vert y\right\vert \leq
r}u\left(  \xi y\right)  f\left(  \left\vert \mathcal{O}y\right\vert \right)
\left\vert J%
\genfrac{(}{)}{}{}{x}{y}%
\right\vert dy\\
& =\int\limits_{\left\vert y\right\vert \leq r}u\left(  \xi y\right)  f\left(
\left\vert y\right\vert \right)  dy,
\end{align*}

and so $\int_{\left\vert x\right\vert \leq r}u\left(  \xi x\right)  f\left(
\left\vert x\right\vert \right)  dx$ is a radial function of $\xi$ and we can
define the radial function
\[
\mu\left(  \left\vert \eta\right\vert \right)  =\int_{\left\vert x\right\vert
\leq r}u\left(  \eta x\right)  f\left(  \left\vert x\right\vert \right)  dx.
\]

Set $\eta=\left(  \mathbf{0}_{1,k-1},\left\vert \xi\right\vert ,\mathbf{0}%
_{k+1,d}\right)  $ so that%
\[
\mu\left(  \left\vert \xi\right\vert \right)  =\int_{\left\vert x\right\vert
\leq r}u\left(  \left\vert \xi\right\vert x_{k}\right)  f\left(  \left\vert
x\right\vert \right)  dx=\int_{\left\vert x\right\vert \leq r}u\left(
\left\vert \xi\right\vert x\right)  f\left(  \left\vert x\right\vert \right)
dx.
\]
\medskip

\textbf{Part 2} From formula 858.515 of the tables \cite{Dwight61}, so noting
that $\omega_{k}=\frac{2\pi^{k/2}}{\Gamma\left(  k/2\right)  }$ we get%
\begin{equation}
\int_{0}^{\pi}\sin^{j}\theta d\theta=\frac{\Gamma\left(  \frac{j+1}{2}\right)
\Gamma\left(  \frac{1}{2}\right)  }{\Gamma\left(  \frac{j+2}{2}\right)
}=\sqrt{\pi}\frac{\Gamma\left(  \frac{j+1}{2}\right)  }{\Gamma\left(
\frac{j+2}{2}\right)  }=\frac{\omega_{j+2}}{\omega_{j+1}},\label{Ap009}%
\end{equation}

so that when $d\geq3$:%
\begin{equation}%
{\textstyle\prod\limits_{j=1}^{d-2}}
\int_{0}^{\pi}\sin^{j}\theta d\theta=\frac{\omega_{d}}{\omega_{2}%
},\label{Ap137}%
\end{equation}

and%
\begin{equation}
2%
{\textstyle\prod\limits_{j=0}^{d-2}}
\int_{0}^{\pi}\sin^{j}\theta d\theta=2\pi\frac{\omega_{d}}{\omega_{2}}%
=\omega_{d}.\label{Ap1371}%
\end{equation}
\medskip

\fbox{\textbf{Case:} $d=2$} From \ref{Ap086},%
\begin{align*}
\int\limits_{\left\vert x\right\vert \leq r}u\left(  sx_{2}\right)  f\left(
\left\vert x\right\vert \right)  dx  & =\int_{-\pi}^{\pi}\int_{0}^{r}u\left(
s\rho\cos\phi_{1}\right)  f\left(  \rho\right)  \rho d\rho d\phi_{1}\\
& =\int_{0}^{r}\left(  \int_{-\pi}^{\pi}u\left(  s\rho\cos t\right)
dt\right)  f\left(  \rho\right)  \rho d\rho\\
& =2\int_{0}^{r}\left(  \int_{0}^{\pi}u\left(  s\rho\cos t\right)  dt\right)
f\left(  \rho\right)  \rho d\rho\\
& =\omega_{1}\int_{0}^{r}\left(  \int_{0}^{\pi}u\left(  s\rho\cos t\right)
dt\right)  f\left(  \rho\right)  \rho d\rho.
\end{align*}
\medskip

\fbox{\textbf{Case:} $d\geq3$} From \ref{Ap086},%
\begin{align*}
\int\limits_{\left\vert x\right\vert \leq r}u\left(  sx_{d}\right)  f\left(
\left\vert x\right\vert \right)  dx  & =\underset{d-2\times}{\underbrace{\int%
\limits_{0}^{\pi}}}\int\limits_{-\pi}^{\pi}\int\limits_{0}^{r}u\left(
sx_{d}\right)  f\left(  \rho\right)  \rho^{d-1}d\rho%
{\textstyle\prod\limits_{j=1}^{d-1}}
\sin^{j-1}\phi_{j}d\phi_{1}\ldots d\phi_{d-1}\\
& =2\underset{d-1\times}{\underbrace{\int\limits_{0}^{\pi}}}\int%
\limits_{0}^{r}u\left(  s\rho\cos\phi_{d-1}\right)  f\left(  \rho\right)
\rho^{d-1}d\rho%
{\textstyle\prod\limits_{j=1}^{d-1}}
\sin^{j-1}\phi_{j}d\phi_{1}\ldots d\phi_{d-1}\\
& =\int_{0}^{r}\left(  \int_{0}^{\pi}u\left(  s\rho\cos\phi_{d-1}\right)
d\phi_{d-1}\right)  f\left(  \rho\right)  \rho^{d-1}d\rho\times\\
& \times2\underset{d-2\times}{\underbrace{\int_{0}^{\pi}}}%
{\textstyle\prod\limits_{j=1}^{d-2}}
\sin^{j-1}\phi_{j}d\phi_{1}d\phi_{2}\ldots d\phi_{d-2}.
\end{align*}

But using \ref{Ap1371}: $2%
{\textstyle\prod\limits_{j=0}^{d-3}}
\int_{0}^{\pi}\sin^{j}\theta d\theta=\omega_{d-1}$,%
\[
\underset{d-2\times}{2\underbrace{\int_{0}^{\pi}}}%
{\textstyle\prod\limits_{j=1}^{d-2}}
\sin^{j-1}\phi_{j}d\phi_{1}d\phi_{2}\ldots d\phi_{d-2}=2%
{\textstyle\prod\limits_{j=0}^{d-3}}
\int_{0}^{\pi}\sin^{j}\theta d\theta=\omega_{d-1},
\]

so%
\[
\int_{\left\vert x\right\vert \leq r}u\left(  sx_{d}\right)  f\left(
\left\vert x\right\vert \right)  dx=\omega_{d-1}\int_{0}^{r}\left(  \int%
_{0}^{\pi}u\left(  s\rho\cos\phi_{d-1}\right)  d\phi_{d-1}\right)  f\left(
\rho\right)  \rho^{d-1}d\rho.
\]

Now apply a permutation to $x$.
\end{proof}

??

\begin{theorem}
\label{Thm_Integ_u(xyhat)f(|x|)dx}?? Suppose $\xi,x\in\mathbb{R}^{d}$ and $\xi
x=\left(  \xi,x\right)  $ denotes the Euclidean inner product. Then:

\begin{enumerate}
\item $\int_{\left\vert x\right\vert \leq r}u\left(  \xi x\right)  f\left(
\left\vert x\right\vert \right)  dx$ is a radial function of $\xi$ and%
\[
\int_{\left\vert x\right\vert \leq r}u\left(  \xi\widehat{x}\right)  f\left(
\left\vert x\right\vert \right)  dx=\int_{\left\vert x\right\vert \leq
r}u\left(  \left\vert \xi\right\vert \widehat{x}_{k}\right)  f\left(
\left\vert x\right\vert \right)  dx,\quad1\leq k\leq d.
\]

\item When $d\geq2$, $s\in\mathbb{R}^{1}$ and $1\leq k\leq d$,%
\[
\int\limits_{\left\vert x\right\vert \leq r}u\left(  s\widehat{x}_{k}\right)
f\left(  \left\vert x\right\vert \right)  dx=??\omega_{d-1}\left(  \int%
_{0}^{\pi}u\left(  s\cos t\right)  \sin^{d-2}tdt\right)  \int_{0}^{r}f\left(
\rho\right)  \rho^{d-1}d\rho,
\]

where $\omega_{k}=\frac{2\pi^{k/2}}{\Gamma\left(  k/2\right)  }$ for $k\geq1$.
\end{enumerate}
\end{theorem}

\begin{proof}
?? ADD proof.
\end{proof}

\begin{corollary}
\label{Cor_2_Thm_Integ_u(xy)f(|x|)dx}%
\[
??\int_{\left\vert x\right\vert \leq r}\left(  \xi x\right)  ^{k}f\left(
\left\vert x\right\vert \right)  dx=0,\quad k\text{ }is\text{ }odd,
\]

and%
\begin{align*}
\int_{\left\vert x\right\vert \leq r}\left(  \xi x\right)  ^{2n}f\left(
\left\vert x\right\vert \right)  dx  & =\left\vert \xi\right\vert
^{2n}B\left(  \frac{d-1}{2},n+\frac{1}{2}\right)  \omega_{d-1}\int_{0}%
^{r}f\left(  t\right)  t^{2n+d-1}dt\\
& =\left\vert \xi\right\vert ^{2n}B\left(  \frac{d-1}{2},n+\frac{1}{2}\right)
\int_{0}^{r}\left\vert x\right\vert ^{2n}f\left(  \left\vert x\right\vert
\right)  dx.
\end{align*}

\end{corollary}

\begin{proof}
?? \textbf{FINISH}! Prove first equation. ??

If $u\left(  s\right)  =s^{2n}$ then from \ref{Ap016},%
\begin{align*}
\int_{\left\vert x\right\vert \leq r}u\left(  \xi x\right)  f\left(
\left\vert x\right\vert \right)  dx  & =\int_{\left\vert x\right\vert \leq
r}u\left(  \left\vert \xi\right\vert x_{k}\right)  f\left(  \left\vert
x\right\vert \right)  dx.\\
\int_{\left\vert x\right\vert \leq r}\left(  \xi x\right)  ^{2n}f\left(
\left\vert x\right\vert \right)  dx  & =\int_{\left\vert x\right\vert \leq
r}\left(  \left\vert \xi\right\vert x_{d}\right)  ^{2n}f\left(  \left\vert
x\right\vert \right)  dx\\
& =\left\vert \xi\right\vert ^{2n}\int_{\left\vert x\right\vert \leq r}%
x_{d}^{2n}f\left(  \left\vert x\right\vert \right)  dx.
\end{align*}

From \ref{Ap001},%
\begin{align*}
\int_{\left\vert x\right\vert \leq r}x_{d}^{2n}f\left(  \left\vert
x\right\vert \right)  dx  & =\omega_{d-1}\int_{0}^{r}\left(  \int_{0}^{\pi
}\left(  t\cos\theta\right)  ^{2n}\sin^{d-2}\theta d\theta\right)  f\left(
t\right)  t^{d-1}dt\\
& =\omega_{d-1}\int_{0}^{r}\left(  \int_{0}^{\pi}\cos^{2n}\theta\sin
^{d-2}\theta d\theta\right)  f\left(  t\right)  t^{2n+d-1}dt\\
& =\left(  \int_{0}^{\pi}\cos^{2n}\theta\sin^{d-2}\theta d\theta\right)
\omega_{d-1}\int_{0}^{r}f\left(  t\right)  t^{2n+d-1}dt,
\end{align*}

and from \ref{Ap010}, $\int_{0}^{\pi}\cos^{2n}\theta\sin^{d-2}\theta
d\theta=B\left(  \frac{d-1}{2},n+\frac{1}{2}\right)  $, so%
\[
\int_{\left\vert x\right\vert \leq r}x_{d}^{2n}f\left(  \left\vert
x\right\vert \right)  dx=B\left(  \frac{d-1}{2},n+\frac{1}{2}\right)
\omega_{d-1}\int_{0}^{r}f\left(  t\right)  t^{2n+d-1}dt,
\]

which implies%
\begin{align*}
\int_{\left\vert x\right\vert \leq r}\left(  \xi x\right)  ^{2n}f\left(
\left\vert x\right\vert \right)  dx  & =\left\vert \xi\right\vert
^{2n}B\left(  \frac{d-1}{2},n+\frac{1}{2}\right)  \omega_{d-1}\int_{0}%
^{r}f\left(  t\right)  t^{2n+d-1}dt\\
& =\left\vert \xi\right\vert ^{2n}B\left(  \frac{d-1}{2},n+\frac{1}{2}\right)
\int_{0}^{r}\left\vert x\right\vert ^{2n}f\left(  \left\vert x\right\vert
\right)  dx.
\end{align*}

\end{proof}

\begin{corollary}
\label{Cor_Thm_Integ_u(xy)f(|x|)dx}If the integral exists then for $r\geq0$
and $d\geq2$:

\begin{enumerate}
\item
\begin{equation}
\int\limits_{\left\vert x\right\vert \leq r}f\left(  \left\vert x\right\vert
\right)  dx=\omega_{d}\int_{0}^{r}t^{d-1}f\left(  t\right)  dt,\label{Ap148}%
\end{equation}

where $\omega_{d}=\frac{2\pi^{d/2}}{\Gamma\left(  d/2\right)  }$.

\item Assuming $p>-1/2$,%
\begin{align}
\int\limits_{\left\vert x\right\vert \leq r}\left\vert \widehat{\xi
}x\right\vert ^{p}f\left(  \left\vert x\right\vert \right)  dx  &
=\int\limits_{\left\vert x\right\vert \leq r}\left\vert x_{k}\right\vert
^{p}f\left(  \left\vert x\right\vert \right)  dx\nonumber\\
& =\frac{B\left(  \frac{d-1}{2},\frac{p+1}{2}\right)  }{B\left(  \frac{d-1}%
{2},\frac{1}{2}\right)  }\omega_{d}\int_{0}^{r}t^{p+d-1}f\left(  t\right)
dt\label{Ap017}\\
& =\frac{B\left(  \frac{d}{2},\frac{p+1}{2}\right)  }{B\left(  \frac{d+p}%
{2},\frac{1}{2}\right)  }\omega_{d}\int_{0}^{r}t^{p+d-1}f\left(  t\right)
dt,\label{Ap142}%
\end{align}

when $1\leq k\leq d$.

\item
\begin{equation}
\int\limits_{\left\vert x\right\vert \leq r}\left\vert \widehat{\xi
}x\right\vert ^{p}f\left(  \left\vert x\right\vert \right)  dx=\frac{B\left(
\frac{d}{2},\frac{p+1}{2}\right)  }{B\left(  \frac{d+p}{2},\frac{1}{2}\right)
}\int\limits_{\left\vert x\right\vert \leq r}\left\vert x\right\vert
^{p}f\left(  \left\vert x\right\vert \right)  dx.\label{Ap141}%
\end{equation}

\item If $x_{+}:=\left(  \left\vert x_{k}\right\vert \right)  $ for
$x\in\mathbb{R}^{d}$ then%
\[
\int\limits_{\left\vert x\right\vert \leq r}\left(  \widehat{\xi}_{+}%
x_{+}\right)  ^{m}f\left(  \left\vert x\right\vert \right)  dx=2^{d}%
\int\limits_{\left\vert x\right\vert \leq r;x\geq\mathbf{0}}\left(
\widehat{\xi}_{+}x\right)  ^{m}f\left(  \left\vert x\right\vert \right)  dx.
\]

\end{enumerate}
\end{corollary}

\begin{proof}
\textbf{Part 1} When $u=1$, noting \ref{Ap009}, part 2 of Theorem
\ref{Thm_Integ_u(xy)f(|x|)dx} becomes%
\begin{align*}
\int_{\left\vert x\right\vert \leq r}f\left(  \left\vert x\right\vert \right)
dx  & =\int_{\left\vert x\right\vert \leq r}u\left(  sx_{k}\right)  f\left(
\left\vert x\right\vert \right)  dx\\
& =\omega_{d-1}\int_{0}^{r}\left(  \int_{0}^{\pi}\sin^{d-2}tdt\right)
f\left(  \rho\right)  \rho^{d-1}d\rho\\
& =\omega_{d-1}\int_{0}^{r}\left(  \int_{0}^{\pi}\sin^{d-2}tdt\right)
f\left(  \rho\right)  \rho^{d-1}d\rho\\
& =2\left(  \int_{0}^{\pi/2}\sin^{d-2}tdt\right)  \omega_{d-1}\int_{0}%
^{r}f\left(  \rho\right)  \rho^{d-1}d\rho\\
& =2\frac{\sqrt{\pi}}{2}\frac{\Gamma\left(  \frac{d-1}{2}\right)  }%
{\Gamma\left(  \frac{d}{2}\right)  }\omega_{d-1}\int_{0}^{r}f\left(
\rho\right)  \rho^{d-1}d\rho\\
& =\sqrt{\pi}\frac{\Gamma\left(  \frac{d-1}{2}\right)  }{\Gamma\left(
\frac{d}{2}\right)  }\frac{2\pi^{\frac{d-1}{2}}}{\Gamma\left(  \frac{d-1}%
{2}\right)  }\int_{0}^{r}f\left(  \rho\right)  \rho^{d-1}d\rho\\
& =\omega_{d}\int_{0}^{r}f\left(  \rho\right)  \rho^{d-1}d\rho.
\end{align*}
\smallskip

\textbf{Part 2} Set $u\left(  t\right)  =\left\vert t\right\vert ^{p}$. Then
from Theorem \ref{Thm_Integ_u(xy)f(|x|)dx}, for any $k$,%

\[
\int_{\left\vert x\right\vert \leq r}\left\vert \xi x\right\vert ^{p}f\left(
\left\vert x\right\vert \right)  dx=\int_{\left\vert x\right\vert \leq
r}\left(  \left\vert \xi\right\vert \left\vert x_{k}\right\vert \right)
^{p}f\left(  \left\vert x\right\vert \right)  dx,
\]

and%
\begin{align*}
\int_{\left\vert x\right\vert \leq r}\left\vert \xi x\right\vert ^{p}f\left(
\left\vert x\right\vert \right)  dx  & =\omega_{d-1}\int_{0}^{r}\left(
\int_{0}^{\pi}\left\vert \left\vert \xi\right\vert \rho\cos t\right\vert
^{p}\sin^{d-2}tdt\right)  f\left(  \rho\right)  \rho^{d-1}d\rho\\
& =\left\vert \xi\right\vert ^{p}\left(  \int_{0}^{\pi}\left\vert \cos
t\right\vert ^{p}\sin^{d-2}tdt\right)  \omega_{d-1}\int_{0}^{r}f\left(
\rho\right)  \rho^{p+d-1}d\rho,
\end{align*}

so that%
\begin{align*}
\int_{\left\vert x\right\vert \leq r}\left\vert \widehat{\xi}x\right\vert
^{p}f\left(  \left\vert x\right\vert \right)  dx  & =\int_{\left\vert
x\right\vert \leq r}\left\vert x_{k}\right\vert ^{p}f\left(  \left\vert
x\right\vert \right)  dx\\
& =\left(  \int_{0}^{\pi}\left\vert \cos\theta\right\vert ^{p}\sin^{d-2}\theta
d\theta\right)  \omega_{d-1}\int_{0}^{r}f\left(  \rho\right)  \rho
^{p+d-1}d\rho.
\end{align*}

Formula 858.515 of Dwight \cite{Dwight61} is%
\begin{equation}%
\begin{array}
[c]{c}%
\int\limits_{0}^{\pi/2}\sin^{p}\theta\cos^{q}\theta d\theta=\frac{1}{2}%
\frac{\Gamma\left(  \frac{p+1}{2}\right)  \Gamma\left(  \frac{q+1}{2}\right)
}{\Gamma\left(  \frac{p+q}{2}+1\right)  }=\frac{1}{2}B\left(  \frac{p+1}%
{2},\frac{q+1}{2}\right)  =\int\limits_{0}^{\pi/2}\cos^{p}\theta\sin^{q}\theta
d\theta,\\
where\text{ }p,q>-\frac{1}{2},
\end{array}
\label{Ap005}%
\end{equation}

so that%
\begin{equation}
\int_{0}^{\pi/2}\cos^{p}\theta\sin^{d-2}\theta d\theta=\frac{1}{2}\frac
{\Gamma\left(  \frac{d-1}{2}\right)  \Gamma\left(  \frac{p+1}{2}\right)
}{\Gamma\left(  \frac{p+d}{2}\right)  }=\frac{1}{2}B\left(  \frac{d-1}%
{2},\frac{p+1}{2}\right)  .\label{Ap006}%
\end{equation}

Thus%
\begin{align}
\int_{0}^{\pi}\left\vert \cos\theta\right\vert ^{p}\sin^{d-2}\theta d\theta &
=\int_{0}^{\pi/2}\left\vert \cos\theta\right\vert ^{p}\sin^{d-2}\theta
d\theta+\int_{\pi/2}^{\pi}\left\vert \cos\theta\right\vert ^{p}\sin
^{d-2}\theta d\theta\nonumber\\
& =\int_{0}^{\pi/2}\cos^{p}\theta\sin^{d-2}\theta d\theta+\int_{0}^{\pi/2}%
\sin^{p}\theta\cos^{d-2}\theta d\theta\nonumber\\
& =2\int_{0}^{\pi/2}\cos^{p}\theta\sin^{d-2}\theta d\theta\nonumber\\
& =B\left(  \frac{d-1}{2},\frac{p+1}{2}\right)  ,\label{Ap010}%
\end{align}

and%
\begin{align*}
\int_{\left\vert x\right\vert \leq r}\left\vert \widehat{\xi}x\right\vert
^{p}f\left(  \left\vert x\right\vert \right)  dx  & =\int_{\left\vert
x\right\vert \leq r}\left\vert x_{k}\right\vert ^{p}f\left(  \left\vert
x\right\vert \right)  dx\\
& =B\left(  \frac{d-1}{2},\frac{p+1}{2}\right)  \omega_{d-1}\int_{0}%
^{r}f\left(  \rho\right)  \rho^{p+d-1}d\rho.
\end{align*}

But%
\begin{align}
\omega_{d}  & =\frac{2\pi^{\frac{d}{2}}}{\Gamma\left(  \frac{d}{2}\right)
}=\frac{2\pi^{\frac{d-1}{2}}}{\Gamma\left(  \frac{d-1}{2}\right)  }\pi
^{\frac{1}{2}}\frac{\Gamma\left(  \frac{d-1}{2}\right)  }{\Gamma\left(
\frac{d}{2}\right)  }=\frac{2\pi^{\frac{d-1}{2}}}{\Gamma\left(  \frac{d-1}%
{2}\right)  }\frac{\Gamma\left(  \frac{1}{2}\right)  \Gamma\left(  \frac
{d-1}{2}\right)  }{\Gamma\left(  \frac{d}{2}\right)  }=\nonumber\\
& =\omega_{d-1}B\left(  \frac{d-1}{2},\frac{1}{2}\right)  ,\label{Ap139}%
\end{align}

so%
\[
\int_{\left\vert x\right\vert \leq r}\left\vert \widehat{\xi}x\right\vert
^{p}f\left(  \left\vert x\right\vert \right)  dx=\frac{B\left(  \frac{d-1}%
{2},\frac{p+1}{2}\right)  }{B\left(  \frac{d-1}{2},\frac{1}{2}\right)  }%
\omega_{d}\int_{0}^{r}f\left(  \rho\right)  \rho^{p+d-1}d\rho.
\]

Finally%
\begin{align}
\frac{B\left(  \frac{d-1}{2},\frac{p+1}{2}\right)  }{B\left(  \frac{d-1}%
{2},\frac{1}{2}\right)  }  & =\frac{\Gamma\left(  \frac{d-1}{2}\right)
\Gamma\left(  \frac{p+1}{2}\right)  }{\Gamma\left(  \frac{d+p}{2}\right)
}\frac{\Gamma\left(  \frac{d}{2}\right)  }{\Gamma\left(  \frac{1}{2}\right)
\Gamma\left(  \frac{d-1}{2}\right)  }=\frac{\Gamma\left(  \frac{p+1}%
{2}\right)  }{\Gamma\left(  \frac{d+p}{2}\right)  }\frac{\Gamma\left(
\frac{d}{2}\right)  }{\Gamma\left(  \frac{1}{2}\right)  }=\nonumber\\
& =\frac{B\left(  \frac{d}{2},\frac{p+1}{2}\right)  }{B\left(  \frac{d+p}%
{2},\frac{1}{2}\right)  }.\label{Ap140}%
\end{align}
\medskip

\textbf{Part 3} Apply part 1 to part 2.\medskip

\textbf{Part 4} The first equation is obtained by splitting the domain of
integration $\left\vert x\right\vert \leq r$ into the $2^{d}$ orthants and
then reflecting each orthants into the orthant $x\geq\mathbf{0}$.
\end{proof}

??

\begin{theorem}
\label{Thm_integ_X^a_f(absX)dX_R^d}For $d\geq2$:

\begin{enumerate}
\item
\[
\int\limits_{\left\vert x\right\vert \leq r}\left\vert \widehat{x}^{\alpha
}\right\vert f\left(  \left\vert x\right\vert \right)  dx=2^{d}\int%
\limits_{\left\vert x\right\vert \leq r;x\geq\mathbf{0}}\widehat{x}^{\alpha
}f\left(  \left\vert x\right\vert \right)  dx,\quad\alpha\geq0.
\]

\item
\begin{equation}
\int\limits_{\left\vert x\right\vert \leq r;x\geq\mathbf{0}}\widehat
{x}^{\alpha}f\left(  \left\vert x\right\vert \right)  dx=\frac{B\left(
\frac{1}{2}\left(  \alpha+\mathbf{1}\right)  \right)  }{B\left(  \frac{1}%
{2}\mathbf{1}\right)  }\int\limits_{\left\vert x\right\vert \leq
r;x\geq\mathbf{0}}f\left(  \left\vert x\right\vert \right)  dx,\label{Ap123}%
\end{equation}

where the multivariate beta function $B\left(  \cdot\right)  $ is defined by:%
\begin{equation}
B\left(  x\right)  :=\frac{\Gamma\left(  x_{1}\right)  \Gamma\left(
x_{2}\right)  \Gamma\left(  x_{3}\right)  \ldots\Gamma\left(  x_{d}\right)
}{\Gamma\left(  x_{1}+x_{2}+\ldots+x_{d}\right)  },\quad x\in\mathbb{R}%
^{d}.\label{Ap128}%
\end{equation}

\item
\[
\frac{B\left(  \frac{1}{2}\left(  \alpha+\mathbf{1}\right)  \right)
}{B\left(  \frac{1}{2}\mathbf{1}\right)  }\leq1,\quad\alpha\geq0;\quad
\sum\limits_{\left\vert \alpha\right\vert =n}\frac{B\left(  \alpha+\frac{1}%
{2}\mathbf{1}\right)  }{\alpha!}=\frac{B\left(  \frac{1}{2}\mathbf{1}\right)
}{n!}.
\]

\item
\begin{equation}
\int\limits_{\left\vert x\right\vert \leq r;x\geq\mathbf{0}}\left(
a\widehat{x}\right)  ^{m}f\left(  \left\vert x\right\vert \right)
dx=\frac{m!\Gamma\left(  \frac{d}{2}\right)  }{2^{m}\Gamma\left(  \frac
{m+d}{2}\right)  }\left(  \sum\limits_{\left\vert \beta\right\vert =m}%
\frac{a^{\beta}}{\left(  \beta/2\right)  !}\right)  \int\limits_{\left\vert
x\right\vert \leq r;x\geq\mathbf{0}}f\left(  \left\vert x\right\vert \right)
dx.\label{Ap143}%
\end{equation}

\item
\[
\int\limits_{\left\vert x\right\vert \leq r}\left(  \widehat{a}_{+}\widehat
{x}_{+}\right)  ^{m}f\left(  \left\vert x\right\vert \right)  dx=\frac
{m!\Gamma\left(  \frac{d}{2}\right)  }{2^{m}\Gamma\left(  \frac{m+d}%
{2}\right)  }\left(  \sum\limits_{\left\vert \beta\right\vert =m}\frac
{a_{+}^{\beta}}{\left(  \beta/2\right)  !}\right)  \int\limits_{\left\vert
x\right\vert \leq r}f\left(  \left\vert x\right\vert \right)  dx.
\]

\textbf{See remark to this theorem}: specifically inequality \ref{Ap127}.

\item ??
\[
\left\vert \int\limits_{\left\vert x\right\vert \leq r}\left(  \widehat
{a}\widehat{x}\right)  ^{m}f\left(  \left\vert x\right\vert \right)
dx\right\vert \leq\frac{m!\Gamma\left(  \frac{d}{2}\right)  }{2^{m}%
\Gamma\left(  \frac{m+d}{2}\right)  }\left(  \sum\limits_{\left\vert
\beta\right\vert =m}\frac{a_{+}^{\beta}}{\left(  \beta/2\right)  !}\right)
\int\limits_{\left\vert x\right\vert \leq r}\left\vert f\left(  \left\vert
x\right\vert \right)  \right\vert dx.
\]

??%
\[
\int\limits_{\left\vert x\right\vert \leq r}\left(  \widehat{a}\widehat
{x}\right)  ^{m}f\left(  \left\vert x\right\vert \right)  dx=??\omega
_{d-1}B\left(  \frac{m+1}{2},\frac{d-2}{2}\right)  \int_{0}^{r}f\left(
\rho\right)  \rho^{d-1}d\rho.
\]

\end{enumerate}
\end{theorem}

\begin{proof}
\textbf{Part 1} By splitting the domain of integration $\left\vert
x\right\vert \leq r$ into the $2^{d}$ orthants and reflecting each orthants
into the orthant $x\geq\mathbf{0}$ we easily obtain
\[
\int\limits_{\left\vert x\right\vert \leq r}\left\vert \widehat{x}^{\alpha
}\right\vert f\left(  \left\vert x\right\vert \right)  dx=2^{d}\int%
\limits_{\left\vert x\right\vert \leq r;x\geq\mathbf{0}}\widehat{x}^{\alpha
}f\left(  \left\vert x\right\vert \right)  dx.
\]
\smallskip

\textbf{Part 2} From \ref{Ap138} the relevant formula is%
\[
\int\limits_{\left\vert x\right\vert \leq r;x\geq\mathbf{0}}g\left(  x\right)
f\left(  \left\vert x\right\vert \right)  dx=\underset{d-1\times
}{\underbrace{\int_{0}^{\pi/2}}}\int_{0}^{r}g\left(  x\right)  f\left(
\rho\right)  \rho^{d-1}d\rho%
{\textstyle\prod\limits_{j=1}^{d-1}}
\sin^{j-1}\phi_{j}d\phi_{1}\ldots d\phi_{d-1},
\]

and here $g\left(  x\right)  =\widehat{x}^{\alpha}$. From \ref{Ap051},\
\begin{align*}
\widehat{x}^{\alpha}  & =\left(  \frac{x}{\left\vert x\right\vert }\right)
^{\alpha}=\frac{x^{\alpha}}{\left\vert x\right\vert ^{\left\vert
\alpha\right\vert }}=\frac{x^{\alpha}}{\rho^{\left\vert \alpha\right\vert }%
}=\\
& =\left(  \rho\sin\phi_{1}\sin\phi_{2}\ldots\sin\phi_{d-1}\right)
^{\alpha_{1}}\left(  \rho\cos\phi_{1}\sin\phi_{2}\ldots\sin\phi_{d-1}\right)
^{\alpha_{2}}\times\\
& \qquad\times\left(  \rho\cos\phi_{2}\ldots\sin\phi_{d-1}\right)
^{\alpha_{3}}\times\ldots\times\left(  \rho\cos\phi_{d-1}\right)  ^{\alpha
_{d}}/\rho^{\left\vert \alpha\right\vert }\\
& =\sin^{\alpha_{1}}\phi_{1}\sin^{\alpha_{1}+\alpha_{2}}\phi_{2}\times
\ldots\times\sin^{\alpha_{1}+\alpha_{2}+\ldots+\alpha_{d-1}}\phi_{d-1}\times\\
& \qquad\times\left(  \cos^{\alpha_{2}}\phi_{1}\cos^{\alpha_{3}}\phi_{2}%
\cos^{\alpha_{4}}\phi_{3}\times\ldots\times\cos^{\alpha_{d}}\phi_{d-1}\right)
,
\end{align*}

and hence%
\[
\int\limits_{\substack{\left\vert x\right\vert \leq r \\x\geq\mathbf{0}%
}}\widehat{x}^{\alpha}f\left(  \left\vert x\right\vert \right)  dx=\left(
\int_{0}^{r}\rho^{d-1}f\left(  \rho\right)  d\rho\right)  \underset{\times
d-1}{\underbrace{\int_{0}^{\pi/2}}}\widehat{x}^{\alpha}%
{\textstyle\prod\limits_{j=1}^{d-1}}
\sin^{j-1}\phi_{j}d\phi_{1}\cdots d\phi_{d-1}.
\]

Regarding the angular term%
\begin{align*}
&  \underset{\times d-1}{\underbrace{\int_{0}^{\frac{\pi}{2}}\ldots\int%
_{0}^{\frac{\pi}{2}}}}\widehat{x}^{\alpha}%
{\textstyle\prod\limits_{j=1}^{d-1}}
\sin^{j-1}\phi_{j}d\phi_{1}\cdots d\phi_{d-1}\\
&  =\underset{\times d-1}{\underbrace{\int\limits_{0}^{\frac{\pi}{2}}%
\ldots\int\limits_{0}^{\frac{\pi}{2}}}}\left(
\begin{array}
[c]{c}%
\sin^{\alpha_{1}}\phi_{1}\sin^{\alpha_{1}+\alpha_{2}}\phi_{2}\times
\ldots\times\sin^{\alpha_{1}+\alpha_{2}+\ldots+\alpha_{d-1}}\phi_{d-1}\times\\
\times\cos^{\alpha_{2}}\phi_{1}\cos^{\alpha_{3}}\phi_{2}\cos^{\alpha_{4}}%
\phi_{3}\times\ldots\times\cos^{\alpha_{d}}\phi_{d-1}%
\end{array}
\right)
{\textstyle\prod\limits_{j=1}^{d-1}}
\sin^{j-1}\phi_{j}d\phi_{1}\ldots d\phi_{d-1}\\
&  =\left(  \int\limits_{0}^{\frac{\pi}{2}}\sin^{\alpha_{1}}\phi_{1}%
\cos^{\alpha_{2}}\phi_{1}d\phi_{1}\right)  \underset{\times
d-2}{\underbrace{\int\limits_{0}^{\frac{\pi}{2}}\ldots\int\limits_{0}%
^{\frac{\pi}{2}}}}\left(
\begin{array}
[c]{c}%
\sin^{\alpha_{1}+\alpha_{2}}\phi_{2}\times\ldots\times\sin^{\alpha_{1}%
+\alpha_{2}+\ldots+\alpha_{d-1}}\phi_{d-1}\times\\
\times\cos^{\alpha_{3}}\phi_{2}\cos^{\alpha_{4}}\phi_{3}\times\ldots\times
\cos^{\alpha_{d}}\phi_{d-1}%
\end{array}
\right)  \times\\
&  \qquad\qquad\times%
{\textstyle\prod\limits_{j=2}^{d-1}}
\sin^{j-1}\phi_{j}d\phi_{2}\ldots d\phi_{d-1}\\
&  =\left(  \int\limits_{0}^{\frac{\pi}{2}}\sin^{\alpha_{1}}\phi_{1}%
\cos^{\alpha_{2}}\phi_{1}d\phi_{1}\right)  \left(  \int\limits_{0}^{\frac{\pi
}{2}}\sin^{\alpha_{1}+\alpha_{2}+1}\phi_{2}\cos^{\alpha_{3}}\phi_{2}d\phi
_{2}\right)  \times\\
&  \qquad\qquad\qquad\times\left(  \int\limits_{0}^{\frac{\pi}{2}}\sin
^{\alpha_{1}+\alpha_{2}+\alpha_{3}+2}\phi_{3}\cos^{\alpha_{4}}\phi_{3}%
d\phi_{3}\right)  \times\ldots\\
&  \ldots\times\left(  \int\limits_{0}^{\frac{\pi}{2}}\sin^{\alpha_{1}%
+\ldots+\alpha_{k}+k-1}\phi_{k}\cos^{\alpha_{k+1}}\phi_{k}d\phi_{k}\right)
\times\ldots\times\left(  \int\limits_{0}^{\frac{\pi}{2}}\sin^{\alpha
_{1}+\ldots+\alpha_{d-1}+d-2}\phi_{d-1}\cos^{\alpha_{d}}\phi_{d-1}d\phi
_{d-1}\right)  .
\end{align*}

From \ref{Ap005},%
\[
\int_{0}^{\frac{\pi}{2}}\sin^{p}\theta\cos^{q}\theta d\theta=\frac{1}{2}%
\frac{\Gamma\left(  \frac{p+1}{2}\right)  \Gamma\left(  \frac{q+1}{2}\right)
}{\Gamma\left(  \frac{p+q}{2}+1\right)  },\quad p,q>-\frac{1}{2},
\]

so that%
\begin{align*}
&  \underset{\times d-1}{\underbrace{\int_{0}^{\pi/2}\ldots\int_{0}^{\pi/2}%
}\widehat{x}^{\alpha}}%
{\textstyle\prod\limits_{j=1}^{d-1}}
\sin^{j-1}\phi_{j}d\phi_{1}\cdots d\phi_{d-1}\\
&  =\left(  \int\limits_{0}^{\frac{\pi}{2}}\sin^{\alpha_{1}}\phi_{1}%
\cos^{\alpha_{2}}\phi_{1}d\phi_{1}\right)  \left(  \int\limits_{0}^{\frac{\pi
}{2}}\sin^{\alpha_{1}+\alpha_{2}+1}\phi_{2}\cos^{\alpha_{3}}\phi_{2}d\phi
_{2}\right)  \left(  \int\limits_{0}^{\frac{\pi}{2}}\sin^{\alpha_{1}%
+\alpha_{2}+\alpha_{3}+2}\phi_{3}\cos^{\alpha_{4}}\phi_{3}d\phi_{3}\right)
\times\ldots\\
&  \ldots\times\left(  \int\limits_{0}^{\frac{\pi}{2}}\sin^{\alpha_{1}%
+\ldots+\alpha_{k}+k-1}\phi_{k}\cos^{\alpha_{k+1}}\phi_{k}d\phi_{k}\right)
\times\ldots\times\left(  \int\limits_{0}^{\frac{\pi}{2}}\sin^{\alpha
_{1}+\ldots+\alpha_{d-1}+d-2}\phi_{d-1}\cos^{\alpha_{d}}\phi_{d-1}d\phi
_{d-1}\right) \\
&  =\frac{1}{2^{d-1}}\frac{\Gamma\left(  \frac{\alpha_{1}+1}{2}\right)
\Gamma\left(  \frac{\alpha_{2}+1}{2}\right)  }{\Gamma\left(  \frac{\alpha
_{1}+\alpha_{2}}{2}+1\right)  }\frac{\Gamma\left(  \frac{\alpha_{1}+\alpha
_{2}}{2}+1\right)  \Gamma\left(  \frac{\alpha_{3}+1}{2}\right)  }%
{\Gamma\left(  \frac{\alpha_{1}+\alpha_{2}+\alpha_{3}+1}{2}+1\right)  }%
\times\ldots\\
&  \ldots\times\frac{\Gamma\left(  \frac{\alpha_{1}+\ldots+\alpha_{k}+k}%
{2}\right)  \Gamma\left(  \frac{\alpha_{k+1}+1}{2}\right)  }{\Gamma\left(
\frac{\alpha_{1}+\ldots+\alpha_{k+1}+k-1}{2}+1\right)  }\times\ldots
\times\frac{\Gamma\left(  \frac{\alpha_{1}+\ldots+\alpha_{d-1}+d-1}{2}\right)
\Gamma\left(  \frac{\alpha_{d}+1}{2}\right)  }{\Gamma\left(  \frac{\alpha
_{1}+\ldots+\alpha_{d-1}+\alpha_{d}+d-2}{2}+1\right)  }\\
&  =\frac{1}{2^{d-1}}\frac{\Gamma\left(  \frac{\alpha_{1}+1}{2}\right)
\Gamma\left(  \frac{\alpha_{2}+1}{2}\right)  \ldots\Gamma\left(  \frac
{\alpha_{d}+1}{2}\right)  }{\Gamma\left(  \frac{\alpha_{1}+\ldots+\alpha
_{d-1}+\alpha_{d}+d-2}{2}+1\right)  }\\
&  =\frac{1}{2^{d-1}}\frac{\Gamma\left(  \frac{\alpha_{1}+1}{2}\right)
\Gamma\left(  \frac{\alpha_{2}+1}{2}\right)  \ldots\Gamma\left(  \frac
{\alpha_{d}+1}{2}\right)  }{\Gamma\left(  \frac{\left\vert \alpha\right\vert
+d}{2}\right)  }\\
&  =\frac{1}{2^{d-1}}B\left(  \frac{\alpha+1}{2}\right)  ,
\end{align*}

Since $\omega_{d}=\frac{2\pi^{d/2}}{\Gamma\left(  d/2\right)  }$ $\left(
=2B\left(  \frac{1}{2}\mathbf{1}\right)  \right)  $ we can now write
\begin{align*}
\int\limits_{\left\vert x\right\vert \leq r;x\geq\mathbf{0}}\widehat
{x}^{\alpha}f\left(  \left\vert x\right\vert \right)  dx  & =\frac{1}{2^{d-1}%
}B\left(  \frac{\alpha+1}{2}\right)  \int_{0}^{r}\rho^{d-1}f\left(
\rho\right)  d\rho\\
& =\frac{1}{2^{d-1}}B\left(  \frac{\alpha+1}{2}\right)  \frac{1}{\omega_{d}%
}\int\limits_{\left\vert x\right\vert \leq r}f\left(  \left\vert x\right\vert
\right)  dx\\
& =2B\left(  \frac{\alpha+1}{2}\right)  \frac{1}{\omega_{d}}\int%
\limits_{\left\vert x\right\vert \leq r;x\geq\mathbf{0}}f\left(  \left\vert
x\right\vert \right)  dx\\
& =\frac{B\left(  \frac{1}{2}\left(  \alpha+\mathbf{1}\right)  \right)
}{B\left(  \frac{1}{2}\mathbf{1}\right)  }\int\limits_{\left\vert x\right\vert
\leq r;x\geq\mathbf{0}}f\left(  \left\vert x\right\vert \right)  dx.
\end{align*}
\medskip

\textbf{Part 3} Choose $f$ such that $f\left(  x\right)  >0$. Then from part 1
and part 2,%
\begin{align*}
\frac{B\left(  \frac{1}{2}\left(  \alpha+\mathbf{1}\right)  \right)
}{B\left(  \frac{1}{2}\mathbf{1}\right)  }\int_{\substack{\left\vert
x\right\vert \leq r \\x\geq\mathbf{0}}}f\left(  \left\vert x\right\vert
\right)  dx=\int_{\substack{\left\vert x\right\vert \leq r \\x\geq\mathbf{0}%
}}\widehat{x}^{\alpha}f\left(  \left\vert x\right\vert \right)  dx  & \leq
\int_{\substack{\left\vert x\right\vert \leq r \\x\geq\mathbf{0}}}\left\vert
\widehat{x}\right\vert ^{\left\vert \alpha\right\vert }f\left(  \left\vert
x\right\vert \right)  dx\\
& =\int_{\substack{\left\vert x\right\vert \leq r \\x\geq\mathbf{0}}}f\left(
\left\vert x\right\vert \right)  dx,
\end{align*}

so%
\[
\frac{B\left(  \frac{1}{2}\left(  \alpha+\mathbf{1}\right)  \right)
}{B\left(  \frac{1}{2}\mathbf{1}\right)  }\leq1.
\]

From part 2%
\[
\int\limits_{\left\vert x\right\vert \leq r}\frac{\widehat{x}^{2\alpha}%
}{\alpha!}f\left(  \left\vert x\right\vert \right)  dx=\frac{B\left(  \frac
{1}{2}\left(  2\alpha+\mathbf{1}\right)  \right)  }{\alpha!B\left(  \frac
{1}{2}\mathbf{1}\right)  }\int\limits_{\left\vert x\right\vert \leq r}f\left(
\left\vert x\right\vert \right)  dx.
\]

The result follows since%
\[
\sum\limits_{\left\vert \alpha\right\vert =n}\int\limits_{\left\vert
x\right\vert \leq r}\frac{\widehat{x}^{2\alpha}}{\alpha!}f\left(  \left\vert
x\right\vert \right)  dx=\int\limits_{\left\vert x\right\vert \leq r}\left(
\sum\limits_{\left\vert \alpha\right\vert =n}\frac{\widehat{x}^{2\alpha}%
}{\alpha!}\right)  f\left(  \left\vert x\right\vert \right)  dx=\frac{1}%
{n!}\int\limits_{\left\vert x\right\vert \leq r}f\left(  \left\vert
x\right\vert \right)  dx,
\]

and%
\begin{align*}
\sum\limits_{\left\vert \alpha\right\vert =n}\frac{B\left(  \frac{1}{2}\left(
2\alpha+\mathbf{1}\right)  \right)  }{\alpha!B\left(  \frac{1}{2}%
\mathbf{1}\right)  }\int\limits_{\left\vert x\right\vert \leq r}f\left(
\left\vert x\right\vert \right)  dx  & =\sum\limits_{\left\vert \alpha
\right\vert =n}\frac{B\left(  \alpha+\frac{1}{2}\mathbf{1}\right)  }%
{\alpha!B\left(  \frac{1}{2}\mathbf{1}\right)  }\int\limits_{\left\vert
x\right\vert \leq r}f\left(  \left\vert x\right\vert \right)  dx\\
& =\left(  \sum\limits_{\left\vert \alpha\right\vert =n}\frac{B\left(
\alpha+\frac{1}{2}\mathbf{1}\right)  }{\alpha!}\right)  \frac{1}{B\left(
\frac{1}{2}\mathbf{1}\right)  }\int\limits_{\left\vert x\right\vert \leq
r}f\left(  \left\vert x\right\vert \right)  dx.
\end{align*}
\medskip

\textbf{Part 4} Using part 2,%
\begin{align}
\int\limits_{\substack{\left\vert x\right\vert \leq r, \\x\geq\mathbf{0}
}}\frac{\left(  a\widehat{x}\right)  ^{m}}{m!}f\left(  \left\vert x\right\vert
\right)  dx  & =\int\limits_{\substack{\left\vert x\right\vert \leq r,
\\x\geq\mathbf{0}}}\sum\limits_{\left\vert \alpha\right\vert =m}%
\frac{a^{\alpha}\widehat{x}^{\alpha}}{\alpha!}f\left(  \left\vert x\right\vert
\right)  dx=\sum\limits_{\left\vert \alpha\right\vert =m}\frac{a^{\alpha}%
}{\alpha!}\int\limits_{\substack{\left\vert x\right\vert \leq r \\x\geq
\mathbf{0}}}\widehat{x}^{\alpha}f\left(  \left\vert x\right\vert \right)
dx=\nonumber\\
& =\left(  2\sum\limits_{\left\vert \alpha\right\vert =m}\frac{a^{\alpha}%
}{\alpha!}B\left(  \frac{\alpha+1}{2}\right)  \right)  \frac{1}{\omega_{d}%
}\int\limits_{\left\vert x\right\vert \leq r;x\geq\mathbf{0}}f\left(
\left\vert x\right\vert \right)  dx.\label{Ap230}%
\end{align}

Next we use the duplication formulas%
\begin{align}
\Gamma\left(  x\right)  \Gamma\left(  x+\frac{1}{2}\right)   & =2^{1-2x}%
\sqrt{\pi}\Gamma\left(  2x\right)  ,\label{Ap120}\\
\Gamma\left(  \frac{x}{2}\right)  \Gamma\left(  \frac{x+1}{2}\right)   &
=2^{1-x}\sqrt{\pi}\Gamma\left(  x\right)  ,\label{Ap121}%
\end{align}

to write%
\begin{align*}
\frac{1}{\alpha!}B\left(  \frac{\alpha+1}{2}\right)   & =\frac{1}{\alpha
!}\frac{\Gamma\left(  \frac{\alpha_{1}+1}{2}\right)  \Gamma\left(
\frac{\alpha_{2}+1}{2}\right)  \ldots\Gamma\left(  \frac{\alpha_{d}+1}%
{2}\right)  }{\Gamma\left(  \frac{\left\vert \alpha\right\vert +d}{2}\right)
}\\
& =\frac{\Gamma\left(  \frac{\alpha_{1}+1}{2}\right)  \Gamma\left(
\frac{\alpha_{2}+1}{2}\right)  \ldots\Gamma\left(  \frac{\alpha_{d}+1}%
{2}\right)  }{\Gamma\left(  \alpha_{1}+1\right)  \Gamma\left(  \alpha
_{2}+1\right)  \ldots\Gamma\left(  \alpha_{d}+1\right)  }\frac{1}%
{\Gamma\left(  \frac{\left\vert \alpha\right\vert +d}{2}\right)  }\\
& =\frac{2^{1-2\left(  \frac{\alpha_{1}+1}{2}\right)  }\sqrt{\pi}}%
{\Gamma\left(  \frac{\alpha_{1}}{2}+1\right)  }\frac{2^{1-2\left(
\frac{\alpha_{2}+1}{2}\right)  }\sqrt{\pi}}{\Gamma\left(  \frac{\alpha_{2}}%
{2}+1\right)  }\ldots\frac{2^{1-2\left(  \frac{\alpha_{d}+1}{2}\right)  }%
\sqrt{\pi}}{\Gamma\left(  \frac{\alpha_{d}}{2}+1\right)  }\frac{1}%
{\Gamma\left(  \frac{\left\vert \alpha\right\vert +d}{2}\right)  }\\
& =\pi^{d/2}\frac{2^{-\alpha_{1}}}{\Gamma\left(  \frac{\alpha_{1}}%
{2}+1\right)  }\frac{2^{-\alpha_{2}}}{\Gamma\left(  \frac{\alpha_{2}}%
{2}+1\right)  }\ldots\frac{2^{-\alpha_{d}}}{\Gamma\left(  \frac{\alpha_{d}}%
{2}+1\right)  }\frac{1}{\Gamma\left(  \frac{\left\vert \alpha\right\vert
+d}{2}\right)  }\\
& =\frac{\pi^{d/2}}{2^{\left\vert \alpha\right\vert }}\frac{1}{\Gamma\left(
\frac{\alpha_{1}}{2}+1\right)  }\frac{1}{\Gamma\left(  \frac{\alpha_{2}}%
{2}+1\right)  }\ldots\frac{1}{\Gamma\left(  \frac{\alpha_{d}}{2}+1\right)
}\frac{1}{\Gamma\left(  \frac{\left\vert \alpha\right\vert +d}{2}\right)  },
\end{align*}

so that if we use the factorial (or pi) function%
\begin{align*}
2\sum\limits_{\left\vert \alpha\right\vert =m} &  \frac{a^{\alpha}}{\alpha
!}B\left(  \frac{\alpha+1}{2}\right) \\
&  =\sum\limits_{\left\vert \alpha\right\vert =m}\frac{\pi^{d/2}%
}{2^{\left\vert \alpha\right\vert -1}}a^{\alpha}\frac{1}{\Gamma\left(
\frac{\alpha_{1}}{2}+1\right)  }\frac{1}{\Gamma\left(  \frac{\alpha_{2}}%
{2}+1\right)  }\ldots\frac{1}{\Gamma\left(  \frac{\alpha_{d}}{2}+1\right)
}\frac{1}{\Gamma\left(  \frac{\left\vert \alpha\right\vert +d}{2}\right)  }\\
&  =\frac{\pi^{d/2}}{2^{m-1}\Gamma\left(  \frac{m+d}{2}\right)  }%
\sum\limits_{\left\vert \alpha\right\vert =m}a^{\alpha}\frac{1}{\Gamma\left(
\frac{\alpha_{1}}{2}+1\right)  }\frac{1}{\Gamma\left(  \frac{\alpha_{2}}%
{2}+1\right)  }\ldots\frac{1}{\Gamma\left(  \frac{\alpha_{d}}{2}+1\right)  }\\
&  =\frac{\pi^{d/2}}{2^{m-1}\Gamma\left(  \frac{m+d}{2}\right)  }%
\sum\limits_{\left\vert \alpha\right\vert =m}\frac{a^{\alpha}}{\Gamma\left(
\frac{\alpha}{2}+1\right)  }\\
&  =\frac{\pi^{d/2}}{2^{m-1}\Gamma\left(  \frac{m+d}{2}\right)  }%
\sum\limits_{\left\vert \alpha\right\vert =m}\frac{a^{\alpha}}{\left(
\alpha/2\right)  !},
\end{align*}

and because $\omega_{d}=\frac{2\pi^{d/2}}{\Gamma\left(  d/2\right)  }$,
\ref{Ap230} becomes%
\begin{align*}
\int\limits_{\left\vert x\right\vert \leq r;x\geq\mathbf{0}}\frac{\left(
a\widehat{x}\right)  ^{m}}{m!}f\left(  \left\vert x\right\vert \right)  dx  &
=2\sum\limits_{\left\vert \alpha\right\vert =m}\frac{a^{\alpha}}{\alpha
!}B\left(  \frac{\alpha+1}{2}\right)  \frac{1}{\omega_{d}}\int%
\limits_{\left\vert x\right\vert \leq r;x\geq\mathbf{0}}f\left(  \left\vert
x\right\vert \right)  dx\\
& =\frac{\pi^{d/2}}{2^{m-1}\Gamma\left(  \frac{m+d}{2}\right)  }\left(
\sum\limits_{\left\vert \alpha\right\vert =m}\frac{a^{\alpha}}{\left(
\alpha/2\right)  !}\right)  \frac{1}{\omega_{d}}\int\limits_{\left\vert
x\right\vert \leq r;x\geq\mathbf{0}}f\left(  \left\vert x\right\vert \right)
dx\\
& =\frac{\pi^{d/2}}{2^{m-1}\Gamma\left(  \frac{m+d}{2}\right)  }\left(
\sum\limits_{\left\vert \alpha\right\vert =m}\frac{a^{\alpha}}{\left(
\alpha/2\right)  !}\right)  \frac{1}{\frac{2\pi^{d/2}}{\Gamma\left(
d/2\right)  }}\int\limits_{\left\vert x\right\vert \leq r;x\geq\mathbf{0}%
}f\left(  \left\vert x\right\vert \right)  dx\\
& =\frac{\Gamma\left(  \frac{d}{2}\right)  }{2^{m}\Gamma\left(  \frac{m+d}%
{2}\right)  }\left(  \sum\limits_{\left\vert \alpha\right\vert =m}%
\frac{a^{\alpha}}{\left(  \alpha/2\right)  !}\right)  \int\limits_{\left\vert
x\right\vert \leq r;x\geq\mathbf{0}}f\left(  \left\vert x\right\vert \right)
dx,
\end{align*}

which proves part 4.\medskip

\textbf{Part 5} From part 4 of Corollary \ref{Cor_Thm_Integ_u(xy)f(|x|)dx} and
then part 4,%
\begin{align*}
\int\limits_{\left\vert x\right\vert \leq r}\left(  \widehat{a}_{+}\widehat
{x}_{+}\right)  ^{m}f\left(  \left\vert x\right\vert \right)  dx  & =2^{d}%
\int\limits_{\left\vert x\right\vert \leq r;x\geq\mathbf{0}}\left(
\widehat{a}_{+}\widehat{x}\right)  ^{m}f\left(  \left\vert x\right\vert
\right)  dx\\
& =\frac{m!\Gamma\left(  \frac{d}{2}\right)  }{2^{m}\Gamma\left(  \frac
{m+d}{2}\right)  }\left(  \sum\limits_{\left\vert \beta\right\vert =m}%
\frac{a_{+}^{\beta}}{\left(  \beta/2\right)  !}\right)  \int%
\limits_{\left\vert x\right\vert \leq r}f\left(  \left\vert x\right\vert
\right)  dx.
\end{align*}
\medskip

\textbf{Part 6} ?? From part 5,%
\begin{align*}
\left\vert \int\limits_{\left\vert x\right\vert \leq r}\left(  \widehat
{a}\widehat{x}\right)  ^{m}f\left(  \left\vert x\right\vert \right)
dx\right\vert  & \leq\int\limits_{\left\vert x\right\vert \leq r}\left\vert
\widehat{a}\widehat{x}\right\vert ^{m}\left\vert f\left(  \left\vert
x\right\vert \right)  \right\vert dx\leq\int\limits_{\left\vert x\right\vert
\leq r}\left(  \widehat{a}_{+}\widehat{x_{+}}\right)  ^{m}\left\vert f\left(
\left\vert x\right\vert \right)  \right\vert dx=\\
& =\frac{m!\Gamma\left(  \frac{d}{2}\right)  }{2^{m}\Gamma\left(  \frac
{m+d}{2}\right)  }\left(  \sum\limits_{\left\vert \beta\right\vert =m}%
\frac{a_{+}^{\beta}}{\left(  \beta/2\right)  !}\right)  \int%
\limits_{\left\vert x\right\vert \leq r}\left\vert f\left(  \left\vert
x\right\vert \right)  \right\vert dx.
\end{align*}

*****************************

Since%
\begin{equation}
\int_{\left\vert x\right\vert \leq r}u\left(  \xi x\right)  f\left(
\left\vert x\right\vert \right)  dx=\omega_{d-1}\int_{0}^{r}\left(  \int%
_{0}^{\pi}u\left(  \left\vert \xi\right\vert \rho\cos t\right)  \sin
^{d-2}tdt\right)  f\left(  \rho\right)  \rho^{d-1}d\rho.
\end{equation}

\begin{align*}
\int\limits_{\left\vert x\right\vert \leq r}\left(  \widehat{a}\widehat
{x}\right)  ^{m}g\left(  \left\vert x\right\vert \right)  dx  & =\int%
\limits_{\left\vert x\right\vert \leq r}\left(  \widehat{a}x\right)  ^{m}%
\frac{g\left(  \left\vert x\right\vert \right)  }{\left\vert x\right\vert
^{m}}dx\\
& =\omega_{d-1}\int_{0}^{r}\left(  \int_{0}^{\pi}\left(  \rho\cos
\theta\right)  ^{m}\sin^{d-2}\theta d\theta\right)  \frac{g\left(
\rho\right)  }{\rho^{m}}\rho^{d-1}d\rho\\
& =\omega_{d-1}\int_{0}^{r}\left(  \int_{0}^{\pi}\cos^{m}\theta\sin
^{d-2}\theta d\theta\right)  g\left(  \rho\right)  \rho^{d-1}d\rho\\
& =\omega_{d-1}\left(  \int_{0}^{\pi}\cos^{m}\theta\sin^{d-2}\theta
d\theta\right)  \int_{0}^{r}g\left(  \rho\right)  \rho^{d-1}d\rho.
\end{align*}

If $m$ even then%
\begin{align*}
\int_{0}^{\pi}\cos^{m}\theta\sin^{d-2}\theta d\theta & =\int_{0}^{\pi/2}%
\cos^{m}\theta\sin^{d-2}\theta d\theta+\int_{\pi/2}^{\pi}\cos^{m}\theta
\sin^{d-2}\theta d\theta\\
& =\int_{0}^{\pi/2}\cos^{m}\theta\sin^{d-2}\theta d\theta+\int_{0}^{\pi/2}%
\cos^{m}\left(  \phi+\frac{\pi}{2}\right)  \sin^{d-2}\left(  \phi+\frac{\pi
}{2}\right)  d\phi\\
& =\int_{0}^{\pi/2}\cos^{m}\theta\sin^{d-2}\theta d\theta+\left(  -1\right)
^{m}\int_{0}^{\pi/2}\sin^{m}\phi\cos^{d-2}\phi d\phi\\
& =2\int_{0}^{\pi/2}\cos^{m}\theta\sin^{d-2}\theta d\theta\\
& =B\left(  \frac{m+1}{2},\frac{d-2}{2}\right)  ,
\end{align*}

so that%
\[
\int\limits_{\left\vert x\right\vert \leq r}\left(  \widehat{a}\widehat
{x}\right)  ^{m}g\left(  \left\vert x\right\vert \right)  dx=\omega
_{d-1}B\left(  \frac{m+1}{2},\frac{d-2}{2}\right)  \int_{0}^{r}g\left(
\rho\right)  \rho^{d-1}d\rho.
\]

\end{proof}

\begin{remark}
?? \textbf{CHECK!} ?? We use the estimates \ref{Ap029} i.e.%
\[
\sqrt{2\pi}e^{-x}x^{x+\frac{1}{2}}<\sqrt{2\pi}e^{-x+\frac{1}{1+12x}}%
x^{x+\frac{1}{2}}<\Gamma\left(  x+1\right)  <\sqrt{2\pi}e^{-x+\frac{1}{12x}%
}x^{x+\frac{1}{2}},\quad x\geq0.
\]

Hence%
\begin{align*}
\frac{\beta_{k}!}{\left(  \beta_{k}/2\right)  !}  & =\frac{\Gamma\left(
\beta_{k}+1\right)  }{\Gamma\left(  \frac{\beta_{k}}{2}+1\right)  }%
<\frac{\sqrt{2\pi}e^{-\beta_{k}+\frac{1}{12\beta_{k}}}\beta_{k}^{\beta
_{k}+\frac{1}{2}}}{\sqrt{2\pi}e^{-\beta_{k}/2+\frac{1}{1+6\beta_{k}}}\left(
\frac{\beta_{k}}{2}\right)  ^{\frac{\beta_{k}}{2}+\frac{1}{2}}}=\\
& =\frac{e^{-\beta_{k}+\frac{1}{12\beta_{k}}}\beta_{k}^{\beta_{k}+\frac{1}{2}%
}}{e^{-\beta_{k}/2+\frac{1}{1+6\beta_{k}}}\left(  \frac{\beta_{k}}{2}\right)
^{\frac{\beta_{k}}{2}+\frac{1}{2}}}=\frac{e^{-\beta_{k}+\frac{1}{12\beta_{k}}%
}\beta_{k}^{\beta_{k}+\frac{1}{2}}2^{\frac{\beta_{k}}{2}+\frac{1}{2}}%
}{e^{-\beta_{k}/2+\frac{1}{1+6\beta_{k}}}\beta_{k}^{\frac{\beta_{k}}{2}%
+\frac{1}{2}}}=\\
& =\left(  e^{-\beta_{k}+\frac{1}{12\beta_{k}}}\beta_{k}^{\beta_{k}+\frac
{1}{2}}2^{\frac{\beta_{k}}{2}+\frac{1}{2}}\right)  \left(  e^{\beta
_{k}/2-\frac{1}{1+6\beta_{k}}}\beta_{k}^{-\frac{\beta_{k}}{2}-\frac{1}{2}%
}\right) \\
& =2^{\frac{\beta_{k}}{2}+\frac{1}{2}}\times\left(  e^{-\beta_{k}+\frac
{1}{12\beta_{k}}}e^{\beta_{k}/2-\frac{1}{1+6\beta_{k}}}\right)  \left(
\beta_{k}^{\beta_{k}+\frac{1}{2}}\beta_{k}^{-\frac{\beta_{k}}{2}-\frac{1}{2}%
}\right) \\
& =2^{\frac{\beta_{k}}{2}+\frac{1}{2}}e^{-\frac{\beta_{k}}{2}+\frac{1}%
{12\beta_{k}}-\frac{1}{1+6\beta_{k}}}\beta_{k}^{\frac{\beta_{k}}{2}}\\
& <\sqrt{2}2^{\frac{\beta_{k}}{2}}e^{-\frac{\beta_{k}}{2}}\beta_{k}%
^{\frac{\beta_{k}}{2}}\text{ }when\text{ }\beta_{k}\geq1.
\end{align*}

Inequality also holds when $\beta_{k}=0$ so%
\[
\frac{\beta_{k}!}{\left(  \beta_{k}/2\right)  !}\leq\sqrt{2}2^{\frac{\beta
_{k}}{2}}e^{-\frac{\beta_{k}}{2}}\beta_{k}^{\frac{\beta_{k}}{2}}\text{
}when\text{ }\beta_{k}\geq0,
\]

and so again from \ref{Ap029},%
\begin{align*}
\frac{\beta!}{\left(  \beta/2\right)  !}<2^{\frac{d}{2}}2^{\frac{\left\vert
\beta\right\vert }{2}}e^{-\frac{\left\vert \beta\right\vert }{2}}\left\vert
\beta\right\vert ^{\frac{\left\vert \beta\right\vert }{2}}  & =\frac
{2^{\frac{d}{2}}2^{\left\vert \beta\right\vert }}{\sqrt{2\pi}\sqrt{\left\vert
\beta\right\vert /2}}\sqrt{2\pi}e^{-\frac{\left\vert \beta\right\vert }{2}%
}\left(  \left\vert \beta\right\vert /2\right)  ^{\frac{\left\vert
\beta\right\vert }{2}}\left(  \left\vert \beta\right\vert /2\right)
^{\frac{1}{2}}\\
& <\frac{2^{\frac{d}{2}+\left\vert \beta\right\vert }}{\sqrt{\left\vert
\beta\right\vert \pi}}\left(  \left\vert \beta\right\vert /2\right)
!=\frac{2^{\frac{d}{2}+\left\vert \beta\right\vert }}{\sqrt{\left\vert
\beta\right\vert \pi}}\left(  \left\vert \beta\right\vert /2\right)  !.
\end{align*}

Thus%
\begin{align*}
\sum_{\left\vert \beta\right\vert =m}\frac{a_{+}^{\beta}}{\left(
\beta/2\right)  !}  & =\sum_{\left\vert \beta\right\vert =m}\frac{\beta
!}{\left(  \beta/2\right)  !}\frac{a_{+}^{\beta}}{\beta!}<\frac{2^{\frac{d}%
{2}+m}}{\sqrt{m\pi}}\left(  m/2\right)  !\sum_{\left\vert \beta\right\vert
=m}\frac{a_{+}^{\beta}}{\beta!}=\\
& =\frac{2^{\frac{d}{2}+m}}{\sqrt{m\pi}}\left(  m/2\right)  !\sum_{\left\vert
\beta\right\vert =m}\frac{a_{+}^{\beta}\mathbf{1}^{\beta}}{\beta!}%
=\frac{2^{\frac{d}{2}+m}}{\sqrt{m\pi}}\frac{\left(  m/2\right)  !}%
{m!}\left\vert a\right\vert _{1}^{m}%
\end{align*}

\begin{align*}
\frac{m!\Gamma\left(  \frac{d}{2}\right)  }{2^{m}\Gamma\left(  \frac{m+d}%
{2}\right)  }\sum_{\left\vert \beta\right\vert =m}\frac{a_{+}^{\beta}}{\left(
\beta/2\right)  !}  & <\frac{m!\Gamma\left(  \frac{d}{2}\right)  }{2^{m}%
\Gamma\left(  \frac{m+d}{2}\right)  }\frac{2^{\frac{d}{2}+m}}{\sqrt{m\pi}%
}\frac{\left(  m/2\right)  !}{m!}\left\vert a\right\vert _{1}^{m}\\
& =\frac{\Gamma\left(  \frac{d}{2}\right)  }{\Gamma\left(  \frac{m+d}%
{2}\right)  }\frac{2^{\frac{d}{2}}}{\sqrt{m\pi}}\left(  \frac{m}{2}\right)
!\left\vert a\right\vert _{1}^{m}\\
& =\frac{\Gamma\left(  \frac{d}{2}\right)  }{\Gamma\left(  \frac{m+d}%
{2}\right)  }\frac{2^{\frac{d}{2}}}{\sqrt{m\pi}}\Gamma\left(  \frac{m}%
{2}+1\right)  \left\vert a\right\vert _{1}^{m}\\
& =\frac{\Gamma\left(  \frac{d}{2}\right)  }{\Gamma\left(  \frac{m+d}%
{2}\right)  }\frac{2^{\frac{d}{2}}}{\sqrt{m\pi}}\frac{m}{2}\Gamma\left(
\frac{m}{2}\right)  \left\vert a\right\vert _{1}^{m}\\
& =B\left(  \frac{m}{2},\frac{d}{2}\right)  \frac{2^{\frac{d}{2}-1}}%
{\sqrt{m\pi}}m\left\vert a\right\vert _{1}^{m}\\
& =\sqrt{\frac{2^{d-2}m}{\pi}}B\left(  \frac{m}{2},\frac{d}{2}\right)
\left\vert a\right\vert _{1}^{m}.
\end{align*}

Finally note that $\left\vert a\right\vert _{1}=\mathbf{1}a_{+}\mathbf{\leq
}\left\vert \mathbf{1}\right\vert \left\vert a\right\vert =d^{1/2}\left\vert
a\right\vert $ so that%
\begin{equation}
\frac{m!\Gamma\left(  \frac{d}{2}\right)  }{2^{m}\Gamma\left(  \frac{m+d}%
{2}\right)  }\sum_{\left\vert \beta\right\vert =m}\frac{a_{+}^{\beta}}{\left(
\beta/2\right)  !}<\sqrt{\frac{2^{d-2}md^{m}}{\pi}}B\left(  \frac{m}{2}%
,\frac{d}{2}\right)  \left\vert a\right\vert ^{m}.\label{Ap127}%
\end{equation}

?? \textbf{Examine this inequality numerically}? ??
\end{remark}

??

\begin{corollary}
\label{Cor_integ_inprod_to_1dim}Assuming $p>-1$,%
\begin{equation}
\int_{\left\vert x\right\vert \leq r}\left\vert \widehat{\xi}x\right\vert
^{p}\left\vert \widehat{\eta}x\right\vert ^{q}f\left(  \left\vert x\right\vert
\right)  dx\leq\frac{B\left(  \frac{d-1}{2},\frac{\max\left\{  p,q\right\}
+1}{2}\right)  }{B\left(  \frac{d-1}{2},\frac{1}{2}\right)  }\int%
\limits_{\left\vert x\right\vert \leq r}\left\vert x\right\vert ^{p+q}f\left(
\left\vert x\right\vert \right)  dx,\label{Ap122}%
\end{equation}

and%
\begin{align*}
\int\limits_{\left\vert x\right\vert \leq r}\left\vert x^{\beta}\right\vert
&  \left\vert \widehat{\xi}x\right\vert ^{p}\left\vert \widehat{\eta
}x\right\vert ^{q}f\left(  \left\vert x\right\vert \right)  dx\\
&  \leq\min\left\{  \frac{2}{\omega_{d}}B\left(  \frac{\beta+1}{2}\right)
,\frac{B\left(  \frac{d-1}{2},\frac{\max\left\{  p,q\right\}  +1}{2}\right)
}{B\left(  \frac{d-1}{2},\frac{1}{2}\right)  }\right\}  \int%
\limits_{\left\vert x\right\vert \leq r}\left\vert x\right\vert ^{\left\vert
\beta\right\vert +p+q}f\left(  \left\vert x\right\vert \right)  dx,
\end{align*}

where $B\left(  \frac{\beta+1}{2}\right)  $ is defined by \ref{Ap128}.
\end{corollary}

\begin{proof}
We prove \ref{Ap122} as follows: from \ref{Ap141},%
\[
\int\limits_{\left\vert x\right\vert \leq r}\left(  \widehat{\xi}\widehat
{x}\right)  ^{s}f\left(  \left\vert x\right\vert \right)  dx=\frac{B\left(
\frac{d-1}{2},\frac{s+1}{2}\right)  }{B\left(  \frac{d-1}{2},\frac{1}%
{2}\right)  }\int\limits_{\left\vert x\right\vert \leq r}f\left(  \left\vert
x\right\vert \right)  dx,
\]

so%
\begin{align*}
\int\limits_{\left\vert x\right\vert \leq r}\left\vert \widehat{\xi
}x\right\vert ^{p} &  \left\vert \widehat{\eta}x\right\vert ^{q}f\left(
\left\vert x\right\vert \right)  dx\\
&  =\int\limits_{\left\vert x\right\vert \leq r}\left\vert \widehat{\xi
}\widehat{x}\right\vert ^{p}\left\vert \widehat{\eta}\widehat{x}\right\vert
^{q}\left\vert x\right\vert ^{p+q}f\left(  \left\vert x\right\vert \right)
dx\\
&  \leq\min\left\{  \int\limits_{\left\vert x\right\vert \leq r}\left\vert
\widehat{\xi}\widehat{x}\right\vert ^{p}\left\vert x\right\vert ^{p+q}f\left(
\left\vert x\right\vert \right)  dx,\int_{\left\vert x\right\vert \leq
r}\left\vert \widehat{\eta}\widehat{x}\right\vert ^{q}\left\vert x\right\vert
^{p+q}f\left(  \left\vert x\right\vert \right)  dx\right\} \\
&  =\min\left\{  \frac{B\left(  \frac{d-1}{2},\frac{p+1}{2}\right)  }{B\left(
\frac{d-1}{2},\frac{1}{2}\right)  }\int\limits_{\left\vert x\right\vert \leq
r}\left\vert x\right\vert ^{p+q}f\left(  \left\vert x\right\vert \right)
dx,\frac{B\left(  \frac{d-1}{2},\frac{q+1}{2}\right)  }{B\left(  \frac{d-1}%
{2},\frac{1}{2}\right)  }\int\limits_{\left\vert x\right\vert \leq
r}\left\vert x\right\vert ^{p+q}f\left(  \left\vert x\right\vert \right)
dx\right\} \\
&  =\frac{B\left(  \frac{d-1}{2},\frac{\max\left\{  p,q\right\}  +1}%
{2}\right)  }{B\left(  \frac{d-1}{2},\frac{1}{2}\right)  }\int%
\limits_{\left\vert x\right\vert \leq r}\left\vert x\right\vert ^{p+q}f\left(
\left\vert x\right\vert \right)  dx.
\end{align*}

Next%
\begin{align*}
\int\limits_{\left\vert x\right\vert \leq r} &  \left\vert x^{\beta
}\right\vert \left\vert \widehat{\xi}x\right\vert ^{p}\left\vert \widehat
{\eta}x\right\vert ^{q}f\left(  \left\vert x\right\vert \right)  dx\\
&  =\int\limits_{\left\vert x\right\vert \leq r}\left\vert \widehat{x}^{\beta
}\right\vert \left\vert \widehat{\xi}\widehat{x}\right\vert ^{p}\left\vert
\widehat{\eta}\widehat{x}\right\vert ^{q}\left\vert x\right\vert ^{\left\vert
\beta\right\vert +p+q}f\left(  \left\vert x\right\vert \right)  dx\\
&  \leq\min\left\{  \int\limits_{\left\vert x\right\vert \leq r}\left\vert
\widehat{x}^{\beta}\right\vert \left\vert x\right\vert ^{\left\vert
\beta\right\vert +p+q}f\left(  \left\vert x\right\vert \right)  dx,\int%
\limits_{\left\vert x\right\vert \leq r}\left\vert \widehat{\xi}\widehat
{x}\right\vert ^{p}\left\vert \widehat{\eta}\widehat{x}\right\vert
^{q}\left\vert x\right\vert ^{\left\vert \beta\right\vert +p+q}f\left(
\left\vert x\right\vert \right)  dx\right\} \\
&  \leq\min\left\{  \frac{2}{\omega_{d}}B\left(  \frac{\beta+1}{2}\right)
\int\limits_{\left\vert x\right\vert \leq r}\left\vert x\right\vert
^{\left\vert \beta\right\vert +p+q}f\left(  \left\vert x\right\vert \right)
dx,\frac{B\left(  \frac{d-1}{2},\frac{\max\left\{  p,q\right\}  +1}{2}\right)
}{B\left(  \frac{d-1}{2},\frac{1}{2}\right)  }\int\limits_{\left\vert
x\right\vert \leq r}\left\vert x\right\vert ^{\left\vert \beta\right\vert
+p+q}f\left(  \left\vert x\right\vert \right)  dx\right\} \\
&  \leq\min\left\{  \frac{2}{\omega_{d}}B\left(  \frac{\beta+1}{2}\right)
,\frac{B\left(  \frac{d-1}{2},\frac{\max\left\{  p,q\right\}  +1}{2}\right)
}{B\left(  \frac{d-1}{2},\frac{1}{2}\right)  }\right\}  \int%
\limits_{\left\vert x\right\vert \leq r}\left\vert x\right\vert ^{\left\vert
\beta\right\vert +p+q}f\left(  \left\vert x\right\vert \right)  dx
\end{align*}

\end{proof}

\begin{corollary}
\label{Cor_Thm_IntegBr_exp(ixy)f(|x|)dx}Suppose $d\geq2$. Then for $\xi
\in\mathbb{R}^{d}$ and $1\leq k\leq d$,
\begin{align*}
\int\limits_{\left\vert x\right\vert \leq r}e^{i\xi x}f\left(  \left\vert
x\right\vert \right)  dx  & =\int\limits_{\left\vert x\right\vert \leq
r}e^{i\left\vert \xi\right\vert x_{k}}f\left(  \left\vert x\right\vert
\right)  dx\\
& =\Gamma\left(  \tfrac{d}{2}\right)  \omega_{d}\int_{0}^{r}\frac
{J_{\frac{d-2}{2}}\left(  s\rho\right)  }{\left(  s\rho/2\right)  ^{\frac
{d-2}{2}}}f\left(  \rho\right)  \rho^{d-1}d\rho\\
& =\Gamma\left(  \frac{d}{2}\right)  \int\limits_{\left\vert x\right\vert \leq
r}\frac{J_{\frac{d-2}{2}}\left(  \left\vert \xi\right\vert \left\vert
x\right\vert \right)  }{\left(  \frac{\left\vert \xi\right\vert \left\vert
x\right\vert }{2}\right)  ^{\frac{d-2}{2}}}f\left(  \left\vert x\right\vert
\right)  dx,
\end{align*}

where $\omega_{d}=\frac{2\pi^{d/2}}{\Gamma\left(  d/2\right)  }$.
\end{corollary}

\begin{proof}
Here we apply both parts of Theorem \ref{Thm_Integ_u(xy)f(|x|)dx} with
$u\left(  \mu\right)  =e^{i\mu}$ and $s=\left\vert \xi\right\vert $. Indeed%
\begin{equation}
\int\limits_{\left\vert x\right\vert \leq r}e^{isx_{d}}f\left(  \left\vert
x\right\vert \right)  dx=\frac{2\pi^{\frac{d-1}{2}}}{\Gamma\left(  \frac
{d-1}{2}\right)  }\int_{0}^{r}\left(  \int_{0}^{\pi}e^{is\rho\cos t}\sin
^{d-2}tdt\right)  f\left(  \rho\right)  \rho^{d-1}d\rho.\label{Ap118}%
\end{equation}

But the change of variables: $q=\cos t$, $dt=-\left(  1-q^{2}\right)
^{-1/2}dq $ we obtain,%
\begin{align*}
\int_{0}^{\pi}e^{is\rho\cos t}\sin^{d-2}tdt  & =-\int_{1}^{-1}e^{is\rho
p}\left(  1-p^{2}\right)  ^{\frac{d-2}{2}}\left(  1-p^{2}\right)  ^{-1/2}dp\\
& =\int_{-1}^{1}e^{is\rho p}\left(  1-p^{2}\right)  ^{\frac{d-2}{2}-\frac
{1}{2}}dp.
\end{align*}

From Exercise 11.7.3 of \cite{Arfken70}:%
\[
J_{v}\left(  x\right)  =\frac{\left(  x/2\right)  ^{v}}{\sqrt{\pi}%
\Gamma\left(  v+\frac{1}{2}\right)  }\int_{-1}^{1}e^{ixp}\left(
1-p^{2}\right)  ^{v-\frac{1}{2}}dp,\quad v>\frac{1}{2},\text{ }x\in
\mathbb{R}^{1},
\]

so when $v=\frac{d-2}{2}$ and $x=s\rho$,%
\[
\int_{-1}^{1}e^{is\rho p}\left(  1-p^{2}\right)  ^{\frac{d-2}{2}-\frac{1}{2}%
}dp=\frac{\sqrt{\pi}\Gamma\left(  \frac{d-1}{2}\right)  }{\left(
s\rho/2\right)  ^{v}}J_{\frac{d-2}{2}}\left(  s\rho\right)  =\frac{\sqrt{\pi
}\Gamma\left(  \frac{d-1}{2}\right)  }{\left(  s\rho/2\right)  ^{v}}%
J_{\frac{d-2}{2}}\left(  s\rho\right)  ,
\]

which implies%
\begin{align*}
\int_{0}^{\pi}e^{is\rho\cos t}\sin^{d-2}tdt  & =\frac{\sqrt{\pi}\Gamma\left(
\frac{d-1}{2}\right)  }{\left(  s\rho/2\right)  ^{v}}J_{\frac{d-2}{2}}\left(
s\rho\right)  ,\\
\int_{\left\vert x\right\vert \leq r}e^{isx_{d}}f\left(  \left\vert
x\right\vert \right)  dx  & =\frac{2\pi^{\frac{d-1}{2}}}{\Gamma\left(
\frac{d-1}{2}\right)  }\int_{0}^{r}\frac{\sqrt{\pi}\Gamma\left(  \frac{d-1}%
{2}\right)  }{\left(  s\rho/2\right)  ^{\frac{d-2}{2}}}J_{\frac{d-2}{2}%
}\left(  s\rho\right)  f\left(  \rho\right)  \rho^{d-1}d\rho\\
& =2\pi^{\frac{d}{2}}\int_{0}^{r}\frac{J_{\frac{d-2}{2}}\left(  s\rho\right)
}{\left(  s\rho/2\right)  ^{\frac{d-2}{2}}}f\left(  \rho\right)  \rho
^{d-1}d\rho\\
& =\Gamma\left(  \tfrac{d}{2}\right)  \frac{2\pi^{\frac{d}{2}}}{\Gamma\left(
\frac{d}{2}\right)  }\int_{0}^{r}\frac{J_{\frac{d-2}{2}}\left(  s\rho\right)
}{\left(  s\rho/2\right)  ^{\frac{d-2}{2}}}f\left(  \rho\right)  \rho
^{d-1}d\rho\\
& =\Gamma\left(  \tfrac{d}{2}\right)  \omega_{d}\int_{0}^{r}\frac
{J_{\frac{d-2}{2}}\left(  s\rho\right)  }{\left(  s\rho/2\right)  ^{\frac
{d-2}{2}}}f\left(  \rho\right)  \rho^{d-1}d\rho.
\end{align*}

Finally, using a permutation change of variables that only swaps $x_{k}$ and
$x_{d}$, we have%
\[
\int\limits_{\left\vert x\right\vert \leq r}e^{isx_{k}}f\left(  \left\vert
x\right\vert \right)  dx=\int\limits_{\left\vert x\right\vert \leq
r}e^{isx_{d}}f\left(  \left\vert x\right\vert \right)  dx,\quad1\leq k\leq d.
\]

\end{proof}

\begin{remark}
\label{Rem_Cor_integ_inprod_to_1dim}If $B$ is the beta function, noting
Corollary \ref{Cor_Thm_Integ_u(xy)f(|x|)dx} we define the function $C$ by
\[
C\left(  d,p\right)  :=\frac{B\left(  \frac{d-1}{2},\frac{p+1}{2}\right)
}{B\left(  \frac{d-1}{2},\frac{1}{2}\right)  }=\frac{B\left(  \frac{d}%
{2},\frac{p+1}{2}\right)  }{B\left(  \frac{1}{2},\frac{p+d}{2}\right)  },\quad
d\geq1,\text{ }p>-1.
\]

We have:

\begin{enumerate}
\item $0<C\left(  d,p\right)  \leq1$ and $C\left(  1,p\right)  =1$ for all
$p\geq0$.

\item When $d>1$, $C\left(  d,p\right)  \leq1$ with equality iff $p=0$.

\item $\frac{\partial}{\partial d}C\left(  d,p\right)  <0$ and $\frac
{\partial}{\partial p}C\left(  d,p\right)  <0$ when $d>1$ and $p\geq0$.

\item $\lim\limits_{p\rightarrow\infty}C\left(  d,p\right)  =0$ and
$\lim\limits_{d\rightarrow\infty}C\left(  d,p\right)  =0$.
\end{enumerate}

Regarding the equation ??\ref{Ap017}, when $p>-\frac{1}{2}$ and $d\geq2$ it is
clear that
\[
\int_{\left\vert x\right\vert \leq r}\left\vert \widehat{\xi}x\right\vert
^{p}f\left(  \left\vert x\right\vert \right)  dx<\int_{\left\vert x\right\vert
\leq r}\left(  \left\vert \widehat{\xi}\right\vert \left\vert x\right\vert
\right)  ^{p}f\left(  \left\vert x\right\vert \right)  dx=\int_{\left\vert
x\right\vert \leq r}\left\vert x\right\vert ^{p}f\left(  \left\vert
x\right\vert \right)  dx,
\]

so \ref{Ap017} is always better.
\end{remark}

\begin{remark}
?? What about the approximation:%
\begin{align*}
\left\vert \int\limits_{\left\vert x\right\vert \leq r}u\left(  x_{k}\right)
f\left(  \left\vert x\right\vert \right)  dx\right\vert  & \leq\int%
_{\left\vert x\right\vert \leq r}\left\vert u\left(  x_{k}\right)  \right\vert
\left\vert f\left(  \left\vert x\right\vert \right)  \right\vert dx\leq
\int\limits_{\left\vert x^{\prime}\right\vert \leq r;\text{ }\left\vert
x_{k}\right\vert \leq r}\left\vert u\left(  x_{k}\right)  \right\vert
\left\vert f\left(  \left\vert x\right\vert \right)  \right\vert dx=\\
& =\int_{-r}^{r}\left\vert u\left(  x_{k}\right)  \right\vert \int_{\left\vert
x^{\prime}\right\vert \leq r}\left\vert f\left(  \left\vert x\right\vert
\right)  \right\vert dx^{\prime}dx_{k}\\
& =\int_{-r}^{r}\left\vert u\left(  x_{k}\right)  \right\vert \int_{\left\vert
x^{\prime}\right\vert \leq r}\left\vert f\left(  \sqrt{x_{k}^{2}+\left\vert
x^{\prime}\right\vert ^{2}}\right)  \right\vert dx^{\prime}dx_{k}\\
& =\omega_{d-1}\int_{-r}^{r}\left\vert u\left(  x_{k}\right)  \right\vert
\int_{0}^{r}\rho^{d-2}\left\vert f\left(  \sqrt{x_{k}^{2}+\rho^{2}}\right)
\right\vert d\rho dx_{k}\\
& =\omega_{d-1}\int_{0}^{r}\left(  \left\vert u\left(  x_{k}\right)
\right\vert +\left\vert u\left(  -x_{k}\right)  \right\vert \right)  \int%
_{0}^{r}\rho^{d-2}\left\vert f\left(  \sqrt{x_{k}^{2}+\rho^{2}}\right)
\right\vert d\rho dx_{k}\\
& =\omega_{d-1}\int_{0}^{r}\int_{0}^{r}\left(  \left\vert u\left(
y_{1}\right)  \right\vert +\left\vert u\left(  -y_{1}\right)  \right\vert
\right)  \left(  y_{2}\right)  ^{d-2}\left\vert f\left(  \left\vert
y\right\vert \right)  \right\vert dy_{1}dy_{2}\\
& \leq\omega_{d-1}\int_{\left\vert y\right\vert \leq r}\left(  \left\vert
u\left(  y_{1}\right)  \right\vert +\left\vert u\left(  -y_{1}\right)
\right\vert \right)  \left(  y_{2}\right)  ^{d-2}\left\vert f\left(
\left\vert y\right\vert \right)  \right\vert dy_{1}dy_{2}\\
& =\omega_{d-1}\int_{0}^{r}\int_{0}^{2\pi}\left(  \left\vert u\left(  \rho
\cos\phi\right)  \right\vert +\left\vert u\left(  -\rho\cos\phi\right)
\right\vert \right)  \left(  \rho\sin\phi\right)  ^{d-2}\left\vert f\left(
\rho\right)  \right\vert \rho\text{ }d\rho d\phi\\
& =\omega_{d-1}\int_{0}^{r}\int_{0}^{2\pi}\left(  \left\vert u\left(  \rho
\cos\phi\right)  \right\vert +\left\vert u\left(  -\rho\cos\phi\right)
\right\vert \right)  \left(  \rho\sin\phi\right)  ^{d-2}\left\vert f\left(
\rho\right)  \right\vert \rho d\phi d\rho\\
& =\omega_{d-1}\int_{0}^{r}\int_{0}^{2\pi}\rho^{d-1}\left\vert f\left(
\rho\right)  \right\vert \left\vert u\left(  \rho\cos\phi\right)  \right\vert
\sin^{d-2}\phi d\phi d\rho+\\
& +\omega_{d-1}\int_{0}^{r}\int_{0}^{2\pi}\rho^{d-1}\left\vert f\left(
\rho\right)  \right\vert \left\vert u\left(  -\rho\cos\phi\right)  \right\vert
\sin^{d-2}\phi d\phi d\rho?
\end{align*}

?? Application:%
\begin{align*}
\int_{\left\vert \cdot\right\vert \leq r}\frac{\left\vert a\xi\right\vert
^{2\left(  n+1\right)  }}{w_{\circ}\left(  \left\vert \cdot\right\vert
\right)  }\left\vert \widehat{g_{n}}\left(  a\xi\right)  \right\vert ^{2}  &
=\int_{\left\vert \xi\right\vert \leq r}\frac{1}{w_{\circ}\left(  \left\vert
\cdot\right\vert \right)  }\left(  \frac{\xi_{k}}{\left\vert a\right\vert
}\right)  ^{2\left(  n+1\right)  }\left\vert \widehat{g_{n}}\left(  \frac
{\xi_{k}}{\left\vert a\right\vert }\right)  \right\vert ^{2}d\xi\\
& =\int_{\left\vert \xi\right\vert \leq r}\frac{d\xi^{\prime}}{w_{\circ
}\left(  \left\vert \cdot\right\vert \right)  }\left(  \frac{\xi_{k}%
}{\left\vert a\right\vert }\right)  ^{2\left(  n+1\right)  }\left\vert
\widehat{g_{n}}\left(  \frac{\xi_{k}}{\left\vert a\right\vert }\right)
\right\vert ^{2}d\xi_{k}\\
& \leq\int_{\left\vert \xi_{k}\right\vert \leq r}\int_{\left\vert \xi^{\prime
}\right\vert \leq r}\frac{d\xi^{\prime}}{w_{\circ}\left(  \left\vert
\cdot\right\vert \right)  }\left(  \frac{\xi_{k}}{\left\vert a\right\vert
}\right)  ^{2\left(  n+1\right)  }\left\vert \widehat{g_{n}}\left(  \frac
{\xi_{k}}{\left\vert a\right\vert }\right)  \right\vert ^{2}d\xi_{k}\\
& ETC.??
\end{align*}

\end{remark}

\subsection{Expressions for $C\left(  d,1\right)  $}%

\begin{align}
C\left(  d,1\right)   & =\frac{B\left(  \frac{d}{2},\frac{2}{2}\right)
}{B\left(  \frac{1}{2},\frac{1+d}{2}\right)  }=\frac{B\left(  \frac{d}%
{2},1\right)  }{B\left(  \frac{1}{2},\frac{1+d}{2}\right)  }=\frac
{\Gamma\left(  \frac{d}{2}\right)  \Gamma\left(  1\right)  }{\Gamma\left(
\frac{1}{2}\right)  \Gamma\left(  \frac{1+d}{2}\right)  }=\label{Ap129}\\
& =\frac{\Gamma\left(  \frac{d}{2}\right)  \Gamma\left(  \frac{1}{2}\right)
}{\Gamma\left(  \frac{1}{2}\right)  ^{2}\Gamma\left(  \frac{1+d}{2}\right)
}=\frac{1}{\pi}B\left(  \frac{d}{2},\frac{1}{2}\right)  .\nonumber
\end{align}

Also%
\[
C\left(  d,1\right)  =\frac{\Gamma\left(  \frac{d}{2}\right)  \Gamma\left(
1\right)  }{\Gamma\left(  \frac{1}{2}\right)  \Gamma\left(  \frac{1+d}%
{2}\right)  }=\frac{1}{\sqrt{\pi}}\frac{\Gamma\left(  \frac{d}{2}\right)
}{\Gamma\left(  \frac{1+d}{2}\right)  }.
\]

From \ref{Ap127},%
\[
\Gamma\left(  \frac{1+m}{2}\right)  \Gamma\left(  \frac{m}{2}\right)
=\sqrt{\pi}\frac{\left(  m-1\right)  !}{2^{m-1}},
\]

so%
\[
C\left(  d,1\right)  =\frac{1}{\sqrt{\pi}}\frac{\Gamma\left(  \frac{d}%
{2}\right)  }{\Gamma\left(  \frac{1+d}{2}\right)  }=\frac{\left(  d-1\right)
!}{2^{d-1}}\frac{1}{\Gamma\left(  \frac{1+d}{2}\right)  ^{2}},
\]

and%
\[
C\left(  d,1\right)  =\frac{1}{\sqrt{\pi}}\frac{\Gamma\left(  \frac{d}%
{2}\right)  }{\Gamma\left(  \frac{1+d}{2}\right)  }=\frac{1}{\pi}\frac
{2^{d-1}}{\left(  d-1\right)  !}\Gamma\left(  \frac{d}{2}\right)  ^{2}.
\]

Thus \fbox{\textbf{when }$d$\textbf{\ is even}}%
\begin{align*}
C\left(  d,1\right)   & =\frac{1}{\pi}\frac{2^{d-1}}{\left(  d-1\right)
!}\Gamma\left(  \frac{d}{2}\right)  ^{2}=\frac{1}{\pi}\frac{2^{d-1}}{\left(
d-1\right)  !}\left(  \left(  \frac{d}{2}-1\right)  !\right)  ^{2}=\\
& =\frac{1}{\pi}\frac{2^{d-1}}{\left(  d-1\right)  !}\left(  \left(
\frac{d-2}{2}\right)  !\right)  ^{2}\\
& =\frac{1}{\pi}\frac{2^{d-1}}{d-1}\frac{1}{\left(  d-2\right)  !}\left(
\left(  \frac{d-2}{2}\right)  !\right)  ^{2}=\\
& =\frac{1}{\pi}\frac{2^{d-1}}{d-1}\frac{1}{\tbinom{d-2}{\frac{d-2}{2}}}.
\end{align*}

and \fbox{\textbf{when }$d$\textbf{\ is odd}}%
\begin{align*}
C\left(  d,1\right)   & =\frac{\left(  d-1\right)  !}{2^{d-1}}\frac{1}%
{\Gamma\left(  \frac{1+d}{2}\right)  ^{2}}=\frac{\left(  d-1\right)
!}{2^{d-1}}\frac{1}{\left(  \left(  \frac{1+d}{2}-1\right)  !\right)  ^{2}}=\\
& =\frac{1}{2^{d-1}}\frac{\left(  d-1\right)  !}{\left(  \left(  \frac{d-1}%
{2}\right)  !\right)  ^{2}}=\frac{1}{2^{d-1}}\tbinom{d-1}{\frac{d-1}{2}}.
\end{align*}

\subsection{Upper bounds for $C\left(  d,p\right)  $ when $p\geq0$}%

\begin{equation}
C\left(  d,p\right)  =\frac{B\left(  \frac{d}{2},\frac{p+1}{2}\right)
}{B\left(  \frac{1}{2},\frac{p+d}{2}\right)  }=\frac{\Gamma\left(  \frac{d}%
{2}\right)  \Gamma\left(  \frac{p+1}{2}\right)  }{\Gamma\left(  \frac{1}%
{2}\right)  \Gamma\left(  \frac{p+d}{2}\right)  },\text{ }p\geq0.\label{Ap222}%
\end{equation}

Two cases: $d$ even and $d$ odd.\medskip

\fbox{\textbf{Case} $d$ even} Set $d=2n\geq2$. Then%
\[
C\left(  d,p\right)  =\frac{\Gamma\left(  n\right)  \Gamma\left(  \frac
{p+1}{2}\right)  }{\Gamma\left(  \frac{1}{2}\right)  \Gamma\left(  \frac{p}%
{2}+n\right)  }.
\]

But if $n\geq2$,%
\[
\Gamma\left(  \frac{d}{2}\right)  =\Gamma\left(  n\right)  =\left(
n-1\right)  !
\]

Also%
\begin{align*}
\Gamma\left(  2+\frac{p}{2}\right)   & =\left(  1+\frac{p}{2}\right)
\Gamma\left(  1+\frac{p}{2}\right)  ,\\
\Gamma\left(  1+\frac{p}{2}\right)   & =\frac{p}{2}\Gamma\left(  \frac{p}%
{2}\right)  ,
\end{align*}

and in general%
\begin{equation}
\Gamma\left(  n+\frac{p}{2}\right)  =\left(  n-1+\frac{p}{2}\right)
\ldots\left(  1+\frac{p}{2}\right)  \Gamma\left(  1+\frac{p}{2}\right)  ,\quad
n\geq2,\text{ }p\geq0,\label{Ap221}%
\end{equation}

and%
\begin{equation}
\Gamma\left(  n+\frac{p}{2}\right)  =\left(  n-1+\frac{p}{2}\right)
\ldots\left(  1+\frac{p}{2}\right)  \frac{p}{2}\Gamma\left(  \frac{p}%
{2}\right)  ,\quad n\geq1,\text{ }p>0,\label{Ap224}%
\end{equation}

so that%
\begin{align}
C\left(  d,p\right)   & =\frac{\Gamma\left(  n\right)  \Gamma\left(
\frac{1+p}{2}\right)  }{\Gamma\left(  \frac{1}{2}\right)  \Gamma\left(
n+\frac{p}{2}\right)  }\nonumber\\
& =\frac{\left(  n-1\right)  !}{\left(  1+\frac{p}{2}\right)  \ldots\left(
n-1+\frac{p}{2}\right)  }\frac{\Gamma\left(  \frac{1+p}{2}\right)  }%
{\Gamma\left(  \frac{1}{2}\right)  \Gamma\left(  1+\frac{p}{2}\right)  }\\
& =\frac{\left(  n-1\right)  !}{\left(  1+\frac{p}{2}\right)  \ldots\left(
n-1+\frac{p}{2}\right)  }\frac{1}{\Gamma\left(  \frac{1}{2}\right)  ^{2}}%
\frac{\Gamma\left(  \frac{1}{2}\right)  \Gamma\left(  \frac{1+p}{2}\right)
}{\Gamma\left(  1+\frac{p}{2}\right)  }\nonumber\\
& =\left\{
\begin{array}
[c]{ll}%
1, & n=1,\\
\frac{1}{1+\frac{p}{2}}\frac{2}{2+\frac{p}{2}}\ldots\frac{n-1}{n-1+\frac{p}%
{2}}, & n\geq2,
\end{array}
\right\}  \frac{1}{\pi}B\left(  \frac{1}{2},\frac{1+p}{2}\right)
.\label{Ap223}%
\end{align}

The function
\[
f\left(  x\right)  =\frac{x}{x+\frac{p}{2}}%
\]

has positive, decreasing slope so%
\[
\frac{f\left(  1\right)  +f\left(  2\right)  +\ldots+f\left(  m\right)  }%
{m}<f\left(  \frac{m+1}{2}\right)  .
\]

But $GeomMean\leq ArithMean$ so%
\[
f\left(  1\right)  f\left(  2\right)  \ldots f\left(  m\right)  \leq\left(
\frac{f\left(  1\right)  +f\left(  2\right)  +\ldots+f\left(  m\right)  }%
{m}\right)  ^{m}<\left(  f\left(  \frac{m+1}{2}\right)  \right)  ^{m},
\]

and thus%
\[
\frac{1}{1+\frac{p}{2}}\frac{2}{2+\frac{p}{2}}\ldots\frac{m}{m+\frac{p}{2}%
}<\left(  \frac{\frac{m+1}{2}}{\frac{m+1}{2}+\frac{p}{2}}\right)  ^{m}%
=\frac{1}{\left(  1+\frac{p}{m+1}\right)  ^{m}}=\frac{1+\frac{p}{m+1}}{\left(
1+\frac{p}{m+1}\right)  ^{m+1}}.
\]

\textbf{Numerical experiments} show that for $s\geq1$,
\[
\frac{1}{\left(  1+\frac{q}{s}\right)  ^{s}}\leq e^{-q}\left(  1+\frac{q}%
{s}\right)  ^{1/2},\quad0\leq q\leq1,
\]

so that%
\[
\frac{1}{1+\frac{p}{2}}\frac{2}{2+\frac{p}{2}}\ldots\frac{m}{m+\frac{p}{2}%
}<\left\{
\begin{array}
[c]{ll}%
e^{-p}\left(  1+\frac{p}{m+1}\right)  ^{3/2}, & 0\leq p\leq1,\\
\left(  1+\frac{p}{m+1}\right)  ^{-m}, & p\geq1,
\end{array}
\right.
\]

and%
\begin{equation}
C\left(  d,p\right)  <\left\{
\begin{array}
[c]{ll}%
\left\{
\begin{array}
[c]{ll}%
e^{-p}\left(  1+\frac{p}{n}\right)  ^{3/2}, & 0\leq p\leq1,\\
\frac{1}{\left(  1+\frac{p}{n}\right)  ^{n-1}}, & p\geq1,
\end{array}
\right\}  , & n\geq2,\\
1, & n=1.
\end{array}
\right\}  \frac{1}{\pi}B\left(  \frac{1}{2},\frac{1+p}{2}\right)
.\label{Ap218}%
\end{equation}

From \textbf{numerical experiments}:%
\[
\frac{1.5708}{\sqrt{\frac{1}{4}+p}}<B\left(  \frac{1}{2},\frac{1+p}{2}\right)
<\frac{1.7725}{\sqrt{\frac{1}{4}+\frac{p}{2}}},
\]

and \textbf{I suspect that}%
\[
\frac{\pi/2}{\sqrt{\frac{1}{4}+p}}<B\left(  \frac{1}{2},\frac{1+p}{2}\right)
<\frac{\sqrt{\pi}}{\sqrt{\frac{1}{4}+\frac{p}{2}}}.
\]

Thus%
\begin{equation}
\frac{1/2}{\sqrt{\frac{1}{4}+\frac{p}{2}}}<\frac{1}{\pi}B\left(  \frac{1}%
{2},\frac{1+p}{2}\right)  <\frac{1/\sqrt{\pi}}{\sqrt{\frac{1}{4}+\frac{p}{2}}%
},\label{Ap219}%
\end{equation}

and \ref{Ap218} becomes%
\begin{align}
C\left(  d,p\right)   & <\left\{
\begin{array}
[c]{ll}%
\left\{
\begin{array}
[c]{ll}%
e^{-p}\left(  1+\frac{p}{n}\right)  ^{3/2}, & 0\leq p\leq1,\\
\frac{1}{\left(  1+\frac{p}{n}\right)  ^{n-1}}, & p\geq1,
\end{array}
\right\}  , & n\geq2,\\
1, & n=1.
\end{array}
\right\}  \frac{1/\sqrt{\pi}}{\sqrt{\frac{1}{4}+\frac{p}{2}}}\nonumber\\
& =\left\{
\begin{array}
[c]{ll}%
\left\{
\begin{array}
[c]{ll}%
e^{-p}\left(  1+\frac{2p}{d}\right)  ^{3/2}, & 0\leq p\leq1,\\
\frac{1}{\left(  1+\frac{2p}{d}\right)  ^{\frac{d}{2}-1}}, & p\geq1,
\end{array}
\right\}  , & d=4,6,8,\ldots,\\
1, & d=2.
\end{array}
\right\}  \frac{\sqrt{4/\pi}}{\sqrt{1+2p}}.\label{Ap220}%
\end{align}
\medskip

\fbox{\textbf{Case} $d$ odd} Set $d=2n+1$. Then from \ref{Ap224},
$\Gamma\left(  n+\frac{p}{2}\right)  =\left(  n-1+\frac{p}{2}\right)
\ldots\left(  1+\frac{p}{2}\right)  \frac{p}{2}\Gamma\left(  \frac{p}%
{2}\right)  $ when $n\geq1,$ $p>0$,%
\begin{align*}
C\left(  d,p\right)   & =\frac{\Gamma\left(  \frac{d}{2}\right)  \Gamma\left(
\frac{1}{2}+\frac{p}{2}\right)  }{\Gamma\left(  \frac{1}{2}\right)
\Gamma\left(  \frac{d}{2}+\frac{p}{2}\right)  }=\frac{\Gamma\left(  n+\frac
{1}{2}\right)  \Gamma\left(  \frac{1}{2}+\frac{p}{2}\right)  }{\Gamma\left(
\frac{1}{2}\right)  \Gamma\left(  n+\frac{1}{2}+\frac{p}{2}\right)  }=\\
& =\frac{\left(  n-1+\frac{1}{2}\right)  \ldots\left(  1+\frac{1}{2}\right)
\frac{1}{2}\Gamma\left(  \frac{1}{2}\right)  \Gamma\left(  \frac{1}{2}%
+\frac{p}{2}\right)  }{\Gamma\left(  \frac{1}{2}\right)  \left(  n-1+\frac
{1}{2}+\frac{p}{2}\right)  \ldots\left(  1+\frac{1}{2}+\frac{p}{2}\right)
\left(  \frac{1}{2}+\frac{p}{2}\right)  \Gamma\left(  \frac{1}{2}+\frac{p}%
{2}\right)  }\\
& =\frac{\left(  n-1+\frac{1}{2}\right)  \ldots\left(  1+\frac{1}{2}\right)
\frac{1}{2}}{\left(  n-1+\frac{1}{2}+\frac{p}{2}\right)  \ldots\left(
1+\frac{1}{2}+\frac{p}{2}\right)  \left(  \frac{1}{2}+\frac{p}{2}\right)  }\\
& =\frac{\frac{1}{2}}{\frac{1}{2}+\frac{p}{2}}\frac{\frac{3}{2}}{\frac{3}%
{2}+\frac{p}{2}}\ldots\frac{n-\frac{3}{2}}{n-\frac{3}{2}+\frac{p}{2}}%
\frac{n-\frac{1}{2}}{n-\frac{1}{2}+\frac{p}{2}}.
\end{align*}

Following the case of $d$ even:%
\begin{align*}
C\left(  d,p\right)   & =f\left(  \frac{1}{2}\right)  f\left(  \frac{3}%
{2}\right)  \ldots f\left(  n-\frac{1}{2}\right) \\
& \leq\left(  \frac{f\left(  \frac{1}{2}\right)  +f\left(  \frac{3}{2}\right)
+f\left(  \frac{5}{2}\right)  +\ldots+f\left(  n-\frac{1}{2}\right)  }%
{n}\right)  ^{n}<\\
& <\left(  f\left(  \frac{n}{2}\right)  \right)  ^{n}=\left(  \frac{\frac
{n}{2}}{\frac{n}{2}+\frac{p}{2}}\right)  ^{n}=\left(  \frac{1}{1+\frac{p}{n}%
}\right)  ^{n}=\left(  1+\frac{p}{n}\right)  ^{-n}\leq\\
& \leq\left\{
\begin{array}
[c]{ll}%
e^{-p}\left(  1+\frac{p}{n}\right)  ^{1/2}, & 0\leq p\leq1,\\
\left(  1+\frac{p}{n}\right)  ^{-n}, & p\geq1,
\end{array}
\right. \\
& =\left\{
\begin{array}
[c]{ll}%
e^{-p}\left(  1+\frac{2p}{d-1}\right)  ^{1/2}, & 0\leq p\leq1,\\
\left(  1+\frac{2p}{d-1}\right)  ^{-\frac{d-1}{2}}, & p\geq1,
\end{array}
\right.
\end{align*}

so that%
\begin{equation}
C\left(  d,p\right)  <\left\{
\begin{array}
[c]{ll}%
e^{-p}\left(  1+\frac{2p}{d-1}\right)  ^{1/2}, & 0\leq p\leq1,\\
\left(  1+\frac{2p}{d-1}\right)  ^{-\frac{d-1}{2}}, & p\geq1,
\end{array}
\right\}  ,\quad d=3,5,7,\ldots.\label{Ap225}%
\end{equation}

\subsection{The convolution of radial functions}

?? \textbf{THE\ PROOF\ OF} \textbf{THIS theorem is wrong! FIX}! USE Theorem
\ref{Thm_Integ_u(xy,|x|)dx} by writing:

Since%
\[
\left\vert y-x\right\vert =\sqrt{\left\vert y\right\vert ^{2}-2yx+\left\vert
x\right\vert ^{2}},
\]

we can write%
\[
\int\limits_{\left\vert x\right\vert \leq R}f\left(  \left\vert x\right\vert
\right)  g\left(  \left\vert y-x\right\vert \right)  dx=\int%
\limits_{\left\vert x\right\vert \leq R}f\left(  \left\vert x\right\vert
\right)  g\left(  \sqrt{\left\vert y\right\vert ^{2}-2yx+\left\vert
x\right\vert ^{2}}\right)  dx,
\]

and apply Theorem \ref{Thm_Integ_u(xy,|x|)dx}:%
\[
\int_{\left\vert x\right\vert \leq R}\Phi\left(  yx,\left\vert x\right\vert
\right)  dx=\omega_{d-1}\int_{0}^{R}\rho^{d-1}\int_{0}^{\pi}\Phi\left(
\left\vert y\right\vert \rho\cos\theta,\rho\right)  \sin^{d-2}\theta d\theta
d\rho,
\]

with%
\[
\Phi\left(  u,v\right)  =f\left(  v\right)  g\left(  \sqrt{\left\vert
y\right\vert ^{2}-2u+v^{2}}\right)  ,
\]

to obtain%
\begin{align*}
\int\limits_{\left\vert x\right\vert \leq R}f\left(  \left\vert x\right\vert
\right)  g\left(  \left\vert y-x\right\vert \right)  dx  & =\omega_{d-1}%
\int_{0}^{R}\rho^{d-1}\int_{0}^{\pi}f\left(  \rho\right)  g\left(
\sqrt{\left\vert y\right\vert ^{2}-2\left\vert y\right\vert \rho\cos
\theta+\rho^{2}}\right)  \sin^{d-2}\theta d\theta d\rho\\
& =\omega_{d-1}\int_{0}^{R}\rho^{d-1}f\left(  \rho\right)  \int_{0}^{\pi
}g\left(  \sqrt{\left\vert y\right\vert ^{2}-2\left\vert y\right\vert \rho
\cos\theta+\rho^{2}}\right)  \sin^{d-2}\theta d\theta d\rho.
\end{align*}

\begin{theorem}
\label{Thm_convol_rad_fn}\textbf{Convolution of radial functions} The
expression $\int\limits_{\left\vert x\right\vert \leq R}f\left(  \left\vert
x\right\vert \right)  g\left(  \left\vert y-x\right\vert \right)  dx$, is an
even function of $y$ when $d=1$ and a radial function of $y$ when $d>1$. In
fact,
\begin{align*}
&  \int\limits_{\left\vert x\right\vert \leq R}f\left(  \left\vert
x\right\vert \right)  g\left(  \left\vert y-x\right\vert \right)  dx\\
&  =\left\{
\begin{array}
[c]{ll}%
\int\limits_{0}^{R}f\left(  x\right)  \left(  g\left(  \left\vert \left\vert
y\right\vert +x\right\vert \right)  +g\left(  \left\vert \left\vert
y\right\vert -x\right\vert \right)  \right)  dx, & d=1,\\
\left\{
\begin{array}
[c]{l}%
\omega_{d-1}\int\limits_{0}^{\pi}\int\limits_{0}^{R}\rho^{d-1}f\left(
\rho\right)  g\left(  \sqrt{\left\vert y\right\vert ^{2}-2\left\vert
y\right\vert \rho\cos\theta+\rho^{2}}\right)  \sin^{d-2}\theta\text{ }d\rho
d\theta,\\
\omega_{d-1}\int\limits_{-1}^{1}\int\limits_{0}^{R}\rho^{d-1}f\left(
\rho\right)  g\left(  \sqrt{\left\vert y\right\vert ^{2}-2\left\vert
y\right\vert \rho t+\rho^{2}}\right)  d\rho\text{ }\left(  1-t^{2}\right)
^{\frac{d-3}{2}}dt,\\
\omega_{d-1}\left\vert y\right\vert ^{d}\int_{-1}^{1}\int_{0}^{R/\left\vert
y\right\vert }\rho^{d-1}f\left(  \left\vert y\right\vert \rho\right)  g\left(
\left\vert y\right\vert \sqrt{1-2\rho t+\rho^{2}}\right)  d\rho\text{ }\left(
1-t^{2}\right)  ^{\frac{d-3}{2}}dt,
\end{array}
\right.  & d\geq2.
\end{array}
\right.
\end{align*}

where%
\[
\omega_{k}=\frac{2\pi^{k/2}}{\Gamma\left(  k/2\right)  },\quad k\geq1.
\]

\end{theorem}

\begin{proof}
?? \textbf{PROOF IS WRONG}! It obtains $\frac{\omega_{d}}{\pi}$ instead of
$\omega_{d-1}$.??

??

If $d=1$ use the change of variables $y\rightarrow-y$. If $d\geq2$ an
orthogonal transformation argument (see the proof of part 1 Theorem
\ref{Thm_Integ_u(xy)f(|x|)dx}) shows that $\int\limits_{\left\vert
x\right\vert \leq r}f\left(  \left\vert x\right\vert \right)  g\left(
\left\vert y-x\right\vert \right)  dx$ is a radial function. Thus we can set%
\begin{align*}
r  & :=\left\vert y\right\vert ,\\
\kappa\left(  r\right)   & :=\int\limits_{\left\vert x\right\vert \leq
R}f\left(  \left\vert x\right\vert \right)  g\left(  \left\vert y-x\right\vert
\right)  dx,
\end{align*}

and choosing $y=\left(  0^{\prime},r\right)  $, $r\geq0$ we have%
\[
\kappa\left(  r\right)  :=\int\limits_{\left\vert x\right\vert \leq R}f\left(
\left\vert x\right\vert \right)  g\left(  \sqrt{\left\vert \left(  0^{\prime
},r\right)  -x\right\vert ^{2}}\right)  dx=\int\limits_{\left\vert
x\right\vert \leq R}f\left(  \left\vert x\right\vert \right)  g\left(
\sqrt{r^{2}-2rx_{d}+\left\vert x\right\vert ^{2}}\right)  dx.
\]

We consider three cases:\medskip

\fbox{Case $d=1$}\medskip\
\begin{align*}
\int\limits_{\left\vert x\right\vert \leq R}f\left(  \left\vert x\right\vert
\right)  g\left(  \left\vert y-x\right\vert \right)  dx  & =\int%
\limits_{-R}^{R}f\left(  \left\vert x\right\vert \right)  g\left(  \left\vert
\left\vert y\right\vert -x\right\vert \right)  dx\\
& =\int\limits_{-R}^{0}f\left(  \left\vert x\right\vert \right)  g\left(
\left\vert \left\vert y\right\vert -x\right\vert \right)  dx+\int%
\limits_{0}^{R}f\left(  x\right)  g\left(  \left\vert \left\vert y\right\vert
-x\right\vert \right)  dx\\
& =\int\limits_{0}^{R}f\left(  x\right)  g\left(  \left\vert \left\vert
y\right\vert +x\right\vert \right)  dx+\int\limits_{0}^{R}f\left(  x\right)
g\left(  \left\vert \left\vert y\right\vert -x\right\vert \right)  dx\\
& =\int\limits_{0}^{R}f\left(  x\right)  \left\{  g\left(  \left\vert
\left\vert y\right\vert +x\right\vert \right)  +g\left(  \left\vert \left\vert
y\right\vert -x\right\vert \right)  \right\}  dx.
\end{align*}
\medskip

\fbox{Case $d=2$}\medskip\
\begin{align*}
\kappa\left(  r\right)   & =\int_{-\pi}^{\pi}\int_{0}^{R}f\left(  \rho\right)
g\left(  \sqrt{r^{2}-2r\rho\cos\phi_{1}+\rho^{2}}\right)  \rho d\rho d\phi
_{1}\\
& =2\int_{0}^{\pi}\int_{0}^{R}\rho f\left(  \rho\right)  g\left(  \sqrt
{r^{2}-2r\rho\cos\theta+\rho^{2}}\right)  d\rho d\theta\\
& =2\int_{-1}^{1}\int_{0}^{R}\rho f\left(  \rho\right)  \frac{g\left(
\sqrt{r^{2}-2r\rho t+\rho^{2}}\right)  }{\sqrt{1-t^{2}}}d\rho dt.
\end{align*}

Now apply the change variables $\rho^{\prime}=\frac{\rho}{r},$ $d\rho
=rd\rho^{\prime}$ to get%
\[
\kappa\left(  r\right)  =2r\int_{-1}^{1}\int_{0}^{R/r}\rho^{\prime}f\left(
r\rho^{\prime}\right)  \frac{g\left(  r\sqrt{1-2\rho^{\prime}t+\left(
\rho^{\prime}\right)  ^{2}}\right)  }{\sqrt{1-t^{2}}}d\rho^{\prime}d\theta.
\]
\medskip

\fbox{Case $d\geq3$} Changing to spherical coordinates using \ref{Ap083}%
:\medskip%
\begin{align*}
\kappa\left(  r\right)   & =\underset{d-3\times}{\underbrace{\int%
\limits_{0}^{\pi}}}\int\limits_{-\pi}^{\pi}\int\limits_{0}^{R}f\left(
\rho\right)  g\left(  \sqrt{r^{2}-2r\rho\cos\phi_{d-1}+\rho^{2}}\right)
\rho^{d-1}%
{\textstyle\prod\limits_{j=1}^{d-1}}
\sin^{j-1}\phi_{j}d\rho d\phi\\
& =\int\limits_{0}^{\pi}\int\limits_{0}^{R}\rho^{d-1}f\left(  \rho\right)
g\left(  \sqrt{r^{2}-2r\rho\cos\phi_{d-1}+\rho^{2}}\right)  d\rho\sin
^{d-2}\phi_{d-1}d\phi_{d-1}\times\\
& \times\underset{d-3\times}{\underbrace{\int_{0}^{\pi}}}\int_{-\pi}^{\pi}%
{\textstyle\prod\limits_{j=1}^{d-2}}
\sin^{j-1}\phi_{j}d\phi^{\prime}\\
& =\int_{0}^{\pi}\int_{0}^{R}\rho^{d-1}f\left(  \rho\right)  g\left(
\sqrt{r^{2}-2r\rho\cos\theta+\rho^{2}}\right)  \sin^{d-2}\theta d\rho
d\theta\times\\
& \times\underset{d-2}{2\underbrace{\int_{0}^{\pi}}}%
{\textstyle\prod\limits_{j=1}^{d-2}}
\sin^{j-1}\phi_{j}d\phi^{\prime}.
\end{align*}

From \ref{Ap137},
\[
\underset{d-2}{\underbrace{\int_{0}^{\pi}}}%
{\textstyle\prod\limits_{j=1}^{d-2}}
\sin^{j-1}\phi_{j}d\phi^{\prime}=\frac{\omega_{d}}{\omega_{2}}=\frac
{\omega_{d}}{\frac{2\pi}{\Gamma\left(  1\right)  }}=\frac{\omega_{d}}{2\pi},
\]

so%
\begin{align*}
\kappa\left(  r\right)   & =\frac{\omega_{d}}{\pi}\int_{0}^{\pi}\int_{0}%
^{R}\rho^{d-1}f\left(  \rho\right)  g\left(  \sqrt{r^{2}-2r\rho\cos\theta
+\rho^{2}}\right)  \sin^{d-2}\theta d\rho d\theta\\
& =\frac{\omega_{d}}{\pi}\int_{-1}^{1}\int_{0}^{R}\rho^{d-1}f\left(
\rho\right)  g\left(  \sqrt{r^{2}-2r\rho t+\rho^{2}}\right)  \left(
1-t^{2}\right)  ^{\frac{d-3}{2}}d\rho dt\\
& =\frac{\omega_{d}}{\pi}\int_{-1}^{1}\int_{0}^{R}\rho^{d-1}f\left(
\rho\right)  g\left(  r\sqrt{1-2\frac{\rho}{r}t+\left(  \frac{\rho}{r}\right)
^{2}}\right)  \left(  1-t^{2}\right)  ^{\frac{d-3}{2}}d\rho dt.
\end{align*}

Now apply the change variables $\rho^{\prime}=\frac{\rho}{r},$ $d\rho
=rd\rho^{\prime}$ to get%
\[
\kappa\left(  r\right)  =\frac{\omega_{d}}{\pi}r^{d}\int_{-1}^{1}\int%
_{0}^{R/r}\left(  \rho^{\prime}\right)  ^{d-1}f\left(  r\rho^{\prime}\right)
g\left(  r\sqrt{1-2\rho^{\prime}t+\left(  \rho^{\prime}\right)  ^{2}}\right)
\left(  1-t^{2}\right)  ^{\frac{d-3}{2}}d\rho^{\prime}dt.
\]

\end{proof}

\begin{corollary}
\label{Cor_Thm_convol_rad_fn}Under the assumptions of Theorem
\ref{Thm_convol_rad_fn}, if $f,g\geq0$ then%
\begin{align*}
&  \int\limits_{\left\vert x\right\vert \leq R}f\left(  \left\vert
x\right\vert \right)  g\left(  \left\vert y-x\right\vert \right)  dx\\
&  =\frac{\omega_{d}}{2^{d-3}\pi}\frac{1}{\left\vert y\right\vert ^{d-2}%
}\times\left\{
\begin{array}
[c]{c}%
\int\limits_{0}^{R}\rho f\left(  \rho\right)  \int\limits_{\left(  \left\vert
y\right\vert -\rho\right)  ^{2}}^{\left(  \left\vert y\right\vert
+\rho\right)  ^{2}}g\left(  \sqrt{s}\right)  \left(  s-\left(  \left\vert
y\right\vert -\rho\right)  ^{2}\right)  ^{\frac{d-3}{2}}\left(  \left(
\left\vert y\right\vert +\rho\right)  ^{2}-s\right)  ^{\frac{d-3}{2}}ds\text{
}d\rho,\\
\int\limits_{0}^{R}\rho f\left(  \rho\right)  \int\limits_{\left\vert
\left\vert y\right\vert -\rho\right\vert }^{\left\vert y\right\vert +\rho
}g\left(  u\right)  \left(  u^{2}-\left(  \left\vert y\right\vert
-\rho\right)  ^{2}\right)  ^{\frac{d-3}{2}}\left(  \left(  \left\vert
y\right\vert +\rho\right)  ^{2}-u^{2}\right)  ^{\frac{d-3}{2}}udu\text{ }%
d\rho.
\end{array}
\right.
\end{align*}

\end{corollary}

\begin{proof}
Because $f,g\geq0$ we can change the order of integration if the new integral
exists:%
\begin{align}
\int\limits_{\left\vert x\right\vert \leq R}f\left(  \left\vert x\right\vert
\right)  g\left(  \left\vert y-x\right\vert \right)  dx  & =\frac{\omega_{d}%
}{\pi}\int\limits_{-1}^{1}\int\limits_{0}^{R}\rho^{d-1}f\left(  \rho\right)
g\left(  \sqrt{r^{2}-2r\rho t+\rho^{2}}\right)  \left(  1-t^{2}\right)
^{\frac{d-3}{2}}d\rho dt\nonumber\\
& =\frac{\omega_{d}}{\pi}\int\limits_{0}^{R}\rho^{d-1}f\left(  \rho\right)
\int\limits_{-1}^{1}g\left(  \sqrt{r^{2}-2r\rho t+\rho^{2}}\right)  \left(
1-t^{2}\right)  ^{\frac{d-3}{2}}dt\text{ }d\rho\label{a032}%
\end{align}

Set $a_{\rho}=r^{2}+\rho^{2}$ and $b_{\rho}=2r\rho$. Then $a_{\rho}+b_{\rho
}=\left(  r+\rho\right)  ^{2}$ and $a_{\rho}-b_{\rho}=\left(  r-\rho\right)
^{2}$ so\ that the change of variables: $s=a_{\rho}-b_{\rho}t$, $t=\frac
{a_{\rho}-s}{b_{\rho}}$, $dt=-\frac{1}{b_{\rho}}ds$, yields%
\begin{align}
\int_{-1}^{1} &  g\left(  \sqrt{r^{2}-2r\rho t+\rho^{2}}\right)  \left(
1-t^{2}\right)  ^{\frac{d-3}{2}}dt\nonumber\\
&  =\int_{-1}^{1}g\left(  \sqrt{a_{\rho}-b_{\rho}t}\right)  \left(
1-t^{2}\right)  ^{\frac{d-3}{2}}dt\nonumber\\
&  =\int_{\left(  r+\rho\right)  ^{2}}^{\left(  r-\rho\right)  ^{2}}g\left(
\sqrt{s}\right)  \left(  1-\left(  \frac{a_{\rho}-s}{b_{\rho}}\right)
^{2}\right)  ^{\frac{d-3}{2}}\left(  -\frac{ds}{b_{\rho}}\right) \nonumber\\
&  =\frac{1}{b_{\rho}^{d-2}}\int_{\left(  r-\rho\right)  ^{2}}^{\left(
r+\rho\right)  ^{2}}g\left(  \sqrt{s}\right)  \left(  b_{\rho}^{2}-\left(
a_{\rho}-s\right)  ^{2}\right)  ^{\frac{d-3}{2}}ds\nonumber\\
&  =\frac{1}{\left(  2r\rho\right)  ^{d-2}}\int_{\left(  r-\rho\right)  ^{2}%
}^{\left(  r+\rho\right)  ^{2}}g\left(  \sqrt{s}\right)  \left(  b_{\rho}%
^{2}-\left(  a_{\rho}-s\right)  ^{2}\right)  ^{\frac{d-3}{2}}ds\nonumber\\
&  =\frac{1}{\left(  2r\rho\right)  ^{d-2}}\int_{\left(  r-\rho\right)  ^{2}%
}^{\left(  r+\rho\right)  ^{2}}g\left(  \sqrt{s}\right)  \left(  b_{\rho
}-a_{\rho}+s\right)  ^{\frac{d-3}{2}}\left(  b_{\rho}+a_{\rho}-s\right)
^{\frac{d-3}{2}}ds\nonumber\\
&  =\frac{1}{\left(  2r\rho\right)  ^{d-2}}\int_{\left(  r-\rho\right)  ^{2}%
}^{\left(  r+\rho\right)  ^{2}}g\left(  \sqrt{s}\right)  \left(  s-\left(
r-\rho\right)  ^{2}\right)  ^{\frac{d-3}{2}}\left(  \left(  r+\rho\right)
^{2}-s\right)  ^{\frac{d-3}{2}}ds.\label{a030}%
\end{align}

Applying the change of variables $u=\sqrt{s}$, $s=u^{2}$, $ds=2udu$ gives%
\begin{align}
\int_{-1}^{1}g &  \left(  \sqrt{r^{2}-2r\rho t+\rho^{2}}\right)  \left(
1-t^{2}\right)  ^{\frac{d-3}{2}}dt\nonumber\\
&  =\frac{2}{\left(  2r\rho\right)  ^{d-2}}\int_{\left\vert r-\rho\right\vert
}^{r+\rho}g\left(  u\right)  \left(  u^{2}-\left(  r-\rho\right)  ^{2}\right)
^{\frac{d-3}{2}}\left(  \left(  r+\rho\right)  ^{2}-u^{2}\right)  ^{\frac
{d-3}{2}}udu.\label{a031}%
\end{align}

$\times$%
\begin{align*}
\int\limits_{\left\vert x\right\vert \leq R} &  f\left(  \left\vert
x\right\vert \right)  g\left(  \left\vert y-x\right\vert \right)  dx\\
&  =\frac{\omega_{d}}{\pi}\int\limits_{0}^{R}\frac{\rho^{d-1}f\left(
\rho\right)  }{\left(  2r\rho\right)  ^{d-2}}\int\limits_{\left(
r-\rho\right)  ^{2}}^{\left(  r+\rho\right)  ^{2}}g\left(  \sqrt{s}\right)
\left(  s-\left(  r-\rho\right)  ^{2}\right)  ^{\frac{d-3}{2}}\left(  \left(
r+\rho\right)  ^{2}-s\right)  ^{\frac{d-3}{2}}ds\text{ }d\rho\\
&  =\frac{\omega_{d}}{2^{d-3}\pi}\frac{1}{\left\vert y\right\vert ^{d-2}}%
\int\limits_{0}^{R}\rho f\left(  \rho\right)  \int\limits_{\left(  \left\vert
y\right\vert -\rho\right)  ^{2}}^{\left(  \left\vert y\right\vert
+\rho\right)  ^{2}}g\left(  \sqrt{s}\right)  \left(  s-\left(  \left\vert
y\right\vert -\rho\right)  ^{2}\right)  ^{\frac{d-3}{2}}\left(  \left(
\left\vert y\right\vert +\rho\right)  ^{2}-s\right)  ^{\frac{d-3}{2}}ds\text{
}d\rho.
\end{align*}

Substituting \ref{a031} into \ref{a032} we get%
\begin{align*}
\int\limits_{\left\vert x\right\vert \leq R} &  f\left(  \left\vert
x\right\vert \right)  g\left(  \left\vert y-x\right\vert \right)  dx\\
&  =\frac{\omega_{d}}{\pi}\int\limits_{0}^{R}\frac{2\rho^{d-1}f\left(
\rho\right)  }{\left(  2r\rho\right)  ^{d-2}}\int\limits_{\left\vert
r-\rho\right\vert }^{r+\rho}g\left(  u\right)  \left(  u^{2}-\left(
r-\rho\right)  ^{2}\right)  ^{\frac{d-3}{2}}\left(  \left(  r+\rho\right)
^{2}-u^{2}\right)  ^{\frac{d-3}{2}}udu\text{ }d\rho\\
&  =\frac{\omega_{d}}{\pi}\int\limits_{0}^{R}\frac{2\rho^{d-1}f\left(
\rho\right)  }{\left(  2\left\vert y\right\vert \rho\right)  ^{d-2}}%
\int\limits_{\left\vert \left\vert y\right\vert -\rho\right\vert }^{\left\vert
y\right\vert +\rho}g\left(  u\right)  \left(  u^{2}-\left(  \left\vert
y\right\vert -\rho\right)  ^{2}\right)  ^{\frac{d-3}{2}}\left(  \left(
\left\vert y\right\vert +\rho\right)  ^{2}-u^{2}\right)  ^{\frac{d-3}{2}%
}udu\text{ }d\rho\\
&  =\frac{\omega_{d}}{2^{d-3}\pi}\frac{1}{\left\vert y\right\vert ^{d-2}}%
\int\limits_{0}^{R}\rho f\left(  \rho\right)  \int\limits_{\left\vert
\left\vert y\right\vert -\rho\right\vert }^{\left\vert y\right\vert +\rho
}g\left(  u\right)  \left(  u^{2}-\left(  \left\vert y\right\vert
-\rho\right)  ^{2}\right)  ^{\frac{d-3}{2}}\left(  \left(  \left\vert
y\right\vert +\rho\right)  ^{2}-u^{2}\right)  ^{\frac{d-3}{2}}udu\text{ }%
d\rho.
\end{align*}

\end{proof}

\chapter{Estimates for $D^{\alpha}f$, $\left(  \widehat{a}D\right)  ^{n}f$,
$\left\vert D\right\vert ^{2n}f$, $\left(  \widehat{\cdot}D\right)  ^{n}f$
where $f=\left\vert \cdot\right\vert ^{t}$, $\left\vert \cdot\right\vert
^{t}log\left\vert \cdot\right\vert $}

\section{Introduction}

\begin{enumerate}
\item Motivated by some examples of formulas for $\left(  \widehat{a}D\right)
^{k}\left\vert x\right\vert ^{-s}$ we will first derive the equations
\ref{Ap044}, \ref{Ap046} and \ref{Ap053}, \ref{Ap041} which have the form
\begin{align*}
\left(  \widehat{a}D\right)  ^{2n}\left\vert x\right\vert ^{-s}  &
=2^{n}\left(  \frac{s}{2}\right)  _{n}p_{s}^{2n}\left(  \widehat{a}\widehat
{x}\right)  \left\vert x\right\vert ^{-\left(  s+2n\right)  },\\
\left(  \widehat{a}D\right)  ^{2n+1}\left\vert x\right\vert ^{-s}  &
=2^{n+1}\left(  \frac{s}{2}\right)  _{n+1}\left(  \widehat{a}\widehat
{x}\right)  p_{s}^{2n+1}\left(  \widehat{a}\widehat{x}\right)  \left\vert
x\right\vert ^{-s-\left(  2n+1\right)  },
\end{align*}

where $p_{s}^{2n}$ and $p_{s}^{2n+1}$ are polynomials and $\widehat
{a}=a/\left\vert a\right\vert $. We make extensive use of the pochhammer or
rising factorial notation $\left(  s\right)  _{k}$ introduced in Definition
\ref{Def_pochhammer_sym}. Iterative formulas are derived for the coefficients
of the polynomials, some involve matrices and are suitable for numerical use.

\item Matlab experiments suggested several \textbf{conjectures} concerning the
properties of the polynomials $p_{s}^{k}$ e.g. Conjecture
\ref{Conj_argmax_p^(2m+1)_neg2k}, and this enables us to conjecture the bounds
for $\left(  \widehat{a}D\right)  ^{m}\left\vert x\right\vert ^{k}$, where
$m,k=0,1,2,\ldots$, described in Theorem \ref{Thm_(aD)^m_|x|^n}.

\item Estimates of $D_{k}^{m}\left\vert x\right\vert ^{-s}$ follow from those
obtained for $\left(  \widehat{a}D\right)  ^{m}\left\vert x\right\vert ^{k}$.

\item Motivated by more Matlab experiments, further \textbf{conjectures} in
Subsubsection \ref{SbSbSect_bound_fn} finally enable us to obtain (conjectures
for) the upper bounds \ref{Ap113} and then \ref{Ap115} for $\left\vert
D^{\alpha}\left\vert x\right\vert ^{-s}\right\vert $ when $\left\vert
\alpha\right\vert +s>0$.
\end{enumerate}

\section{Upper bounds for $\left(  \widehat{a}D\right)  ^{k}\left\vert
x\right\vert ^{-s}$}

??

\subsection{Some expansions of $\left(  \widehat{a}D\right)  ^{k}\left\vert
x\right\vert ^{-s}$}

We know that%
\begin{equation}
\left(  \widehat{a}D\right)  \left\vert x\right\vert ^{-t}=-t\left(
\widehat{a}x\right)  \left\vert x\right\vert ^{-\left(  t+2\right)  }%
,\quad\left(  \widehat{a}D\right)  \widehat{a}x=1.\label{Ap048}%
\end{equation}

Hence%
\[
\left(  \widehat{a}D\right)  \left\vert x\right\vert ^{-s}=-s\left(
\widehat{a}x\right)  \left\vert x\right\vert ^{-\left(  s+2\right)  }.
\]%
\begin{align*}
\left(  \widehat{a}D\right)  ^{2}\left\vert x\right\vert ^{-s}  &
=-s\left\vert x\right\vert ^{-\left(  s+2\right)  }+s\left(  s+2\right)
\left(  \widehat{a}x\right)  ^{2}\left\vert x\right\vert ^{-\left(
s+4\right)  }\\
& =s\left(  -1+\left(  s+2\right)  \left(  \widehat{a}\widehat{x}\right)
^{2}\right)  \left\vert x\right\vert ^{-\left(  s+2\right)  }.
\end{align*}

\begin{align*}
\left(  \widehat{a}D\right)  ^{3}\left\vert x\right\vert ^{-s}  & =3s\left(
s+2\right)  \left(  \widehat{a}x\right)  \left\vert x\right\vert ^{-\left(
s+4\right)  }-s\left(  s+2\right)  \left(  s+4\right)  \left(  \widehat
{a}x\right)  ^{3}\left\vert x\right\vert ^{-\left(  s+6\right)  }\\
& =s\left(  s+2\right)  \left(  \widehat{a}\widehat{x}\right)  \left\{
3-\left(  s+4\right)  \left(  \widehat{a}\widehat{x}\right)  ^{2}\right\}
\left\vert x\right\vert ^{-\left(  s+3\right)  }.
\end{align*}

\begin{align*}
\left(  \widehat{a}D\right)  ^{4}\left\vert x\right\vert ^{-s}  & =3s\left(
s+2\right)  \left\vert x\right\vert ^{-\left(  s+4\right)  }-6s\left(
s+2\right)  \left(  s+4\right)  \left(  \widehat{a}x\right)  ^{2}\left\vert
x\right\vert ^{-\left(  s+6\right)  }+\\
& +s\left(  s+2\right)  \left(  s+4\right)  \left(  s+6\right)  \left(
\widehat{a}x\right)  ^{4}\left\vert x\right\vert ^{-\left(  s+8\right)  }\\
& =s\left(  s+2\right)  \left\{  3-6\left(  s+4\right)  \left(  \widehat
{a}\widehat{x}\right)  ^{2}+\left(  s+4\right)  \left(  s+6\right)  \left(
\widehat{a}\widehat{x}\right)  ^{4}\right\}  \left\vert x\right\vert
^{-\left(  s+4\right)  }.
\end{align*}

\begin{align*}
&  \left(  \widehat{a}D\right)  ^{5}\left\vert x\right\vert ^{-s}\\
&  =3s\left(  s+2\right)  \left(  \widehat{a}D\right)  \left\vert x\right\vert
^{-\left(  s+4\right)  }-6s\left(  s+2\right)  \left(  s+4\right)  \left(
\widehat{a}D\right)  \left(  \left(  \widehat{a}x\right)  ^{2}\left\vert
x\right\vert ^{-\left(  s+6\right)  }\right)  +\\
&  \quad+s\left(  s+2\right)  \left(  s+4\right)  \left(  s+6\right)  \left(
\widehat{a}D\right)  \left(  \left(  \widehat{a}x\right)  ^{4}\left\vert
x\right\vert ^{-\left(  s+8\right)  }\right) \\
&  =-3s\left(  s+2\right)  \left(  s+4\right)  \left(  \widehat{a}x\right)
\left\vert x\right\vert ^{-\left(  s+6\right)  }-12s\left(  s+2\right)
\left(  s+4\right)  \left(  \widehat{a}x\right)  \left\vert x\right\vert
^{-\left(  s+6\right)  }+\\
&  \quad+6s\left(  s+2\right)  \left(  s+4\right)  \left(  s+6\right)  \left(
\widehat{a}x\right)  ^{3}\left\vert x\right\vert ^{-\left(  s+8\right)  }+\\
&  \quad+4s\left(  s+2\right)  \left(  s+4\right)  \left(  s+6\right)  \left(
\widehat{a}x\right)  ^{3}\left\vert x\right\vert ^{-\left(  s+8\right)  }-\\
&  \quad-s\left(  s+2\right)  \left(  s+4\right)  \left(  s+6\right)  \left(
s+8\right)  \left(  \widehat{a}x\right)  ^{5}\left\vert x\right\vert
^{-\left(  s+10\right)  }\\
&  =-15s\left(  s+2\right)  \left(  s+4\right)  \left(  \widehat{a}x\right)
\left\vert x\right\vert ^{-\left(  s+6\right)  }+10s\left(  s+2\right)
\left(  s+4\right)  \left(  s+6\right)  \left(  \widehat{a}x\right)
^{3}\left\vert x\right\vert ^{-\left(  s+8\right)  }-\\
&  \quad-s\left(  s+2\right)  \left(  s+4\right)  \left(  s+6\right)  \left(
s+8\right)  \left(  \widehat{a}x\right)  ^{5}\left\vert x\right\vert
^{-\left(  s+10\right)  }\\
&  =s\left(  s+2\right)  \left(  s+4\right)  \left(  \widehat{a}\widehat
{x}\right)  \left\{  -15+10\left(  s+6\right)  \left(  \widehat{a}\widehat
{x}\right)  ^{2}-\left(  s+6\right)  \left(  s+8\right)  \left(  \widehat
{a}\widehat{x}\right)  ^{4}\right\}  \left\vert x\right\vert ^{-\left(
s+6\right)  }.
\end{align*}

\begin{align*}
&  \frac{\left(  \widehat{a}D\right)  ^{6}\left\vert x\right\vert ^{-s}%
}{s\left(  s+2\right)  \left(  s+4\right)  }\\
&  =-15\left(  \widehat{a}D\right)  \left\{  \left(  \widehat{a}x\right)
\left\vert x\right\vert ^{-\left(  s+6\right)  }\right\}  +10\left(
s+6\right)  \left(  \widehat{a}D\right)  \left\{  \left(  \widehat{a}x\right)
^{3}\left\vert x\right\vert ^{-\left(  s+8\right)  }\right\}  -\\
&  -\left(  s+6\right)  \left(  s+8\right)  \left(  \widehat{a}D\right)
\left\{  \left(  \widehat{a}x\right)  ^{5}\left\vert x\right\vert ^{-\left(
s+10\right)  }\right\} \\
&  =-15\left(  \widehat{a}D\right)  \left\{  \left(  \widehat{a}x\right)
\left\vert x\right\vert ^{-\left(  s+6\right)  }\right\}  +10\left(
s+6\right)  \left(  \widehat{a}D\right)  \left\{  \left(  \widehat{a}x\right)
^{3}\left\vert x\right\vert ^{-\left(  s+8\right)  }\right\}  -\\
&  \quad-\left(  s+6\right)  \left(  s+8\right)  \left(  \widehat{a}D\right)
\left\{  \left(  \widehat{a}x\right)  ^{5}\left\vert x\right\vert ^{-\left(
s+10\right)  }\right\} \\
&  =-15\left\{  \left\vert x\right\vert ^{-\left(  s+6\right)  }-\left(
s+6\right)  \left(  \widehat{a}x\right)  ^{2}\left\vert x\right\vert
^{-\left(  s+8\right)  }\right\}  +\\
&  \quad+10\left(  s+6\right)  \left\{  3\left(  \widehat{a}x\right)
^{2}\left\vert x\right\vert ^{-\left(  s+8\right)  }-\left(  s+8\right)
\left(  \widehat{a}x\right)  ^{4}\left\vert x\right\vert ^{-\left(
s+10\right)  }\right\}  -\\
&  \quad-\left(  s+6\right)  \left(  s+8\right)  \left\{  5\left(  \widehat
{a}x\right)  ^{4}\left\vert x\right\vert ^{-\left(  s+10\right)  }-\left(
s+10\right)  \left(  \widehat{a}x\right)  ^{6}\left\vert x\right\vert
^{-\left(  s+12\right)  }\right\} \\
&  =-15\left\vert x\right\vert ^{-\left(  s+6\right)  }+15s\left(  s+2\right)
\left(  s+4\right)  \left(  s+6\right)  \left(  \widehat{a}x\right)
^{2}\left\vert x\right\vert ^{-\left(  s+8\right)  }+30\left(  s+6\right)
\left(  \widehat{a}x\right)  ^{2}\left\vert x\right\vert ^{-\left(
s+8\right)  }-\\
&  -10\left(  s+6\right)  \left(  s+8\right)  \left(  \widehat{a}x\right)
^{4}\left\vert x\right\vert ^{-\left(  s+10\right)  }-5\left(  s+6\right)
\left(  s+8\right)  \left(  \widehat{a}x\right)  ^{4}\left\vert x\right\vert
^{-\left(  s+10\right)  }+\\
&  \quad+\left(  s+6\right)  \left(  s+8\right)  \left(  s+10\right)  \left(
\widehat{a}x\right)  ^{6}\left\vert x\right\vert ^{-\left(  s+12\right)  }\\
&  =-15\left\vert x\right\vert ^{-\left(  s+6\right)  }+45\left(  s+6\right)
\left(  \widehat{a}x\right)  ^{2}\left\vert x\right\vert ^{-\left(
s+8\right)  }-15\left(  s+6\right)  \left(  s+8\right)  \left(  \widehat
{a}x\right)  ^{4}\left\vert x\right\vert ^{-\left(  s+10\right)  }+\\
&  \quad+\left(  s+6\right)  \left(  s+8\right)  \left(  s+10\right)  \left(
\widehat{a}x\right)  ^{6}\left\vert x\right\vert ^{-\left(  s+12\right)  }\\
&  =-15\left\vert x\right\vert ^{-\left(  s+6\right)  }+45\left(  s+6\right)
\left(  \widehat{a}x\right)  ^{2}\left\vert x\right\vert ^{-\left(
s+8\right)  }-15\left(  s+6\right)  \left(  s+8\right)  \left(  \widehat
{a}x\right)  ^{4}\left\vert x\right\vert ^{-\left(  s+10\right)  }+\\
&  \quad+\left(  s+6\right)  \left(  s+8\right)  \left(  s+10\right)  \left(
\widehat{a}x\right)  ^{6}\left\vert x\right\vert ^{-\left(  s+12\right)  }\\
&  =\left\{
\begin{array}
[c]{c}%
-15+45\left(  s+6\right)  \left(  \widehat{a}\widehat{x}\right)
^{2}-15\left(  s+6\right)  \left(  s+8\right)  \left(  \widehat{a}\widehat
{x}\right)  ^{4}+\\
+\left(  s+6\right)  \left(  s+8\right)  \left(  s+10\right)  \left(
\widehat{a}\widehat{x}\right)  ^{6}%
\end{array}
\right\}  \left\vert x\right\vert ^{-\left(  s+6\right)  }.
\end{align*}

\begin{align*}
&  \frac{\left(  \widehat{a}D\right)  ^{7}\left\vert x\right\vert ^{-s}%
}{s\left(  s+2\right)  \left(  s+4\right)  \left(  s+6\right)  }\\
&  =-\frac{15}{s+6}\left(  \widehat{a}D\right)  \left\{  \left\vert
x\right\vert ^{-\left(  s+6\right)  }\right\}  +45\left(  \widehat{a}D\right)
\left\{  \left(  \widehat{a}x\right)  ^{2}\left\vert x\right\vert ^{-\left(
s+8\right)  }\right\}  -\\
&  \quad-15\left(  s+8\right)  \left(  \widehat{a}D\right)  \left\{  \left(
\widehat{a}x\right)  ^{4}\left\vert x\right\vert ^{-\left(  s+10\right)
}\right\}  +\\
&  +\left(  s+8\right)  \left(  s+10\right)  \left(  \widehat{a}D\right)
\left\{  \left(  \widehat{a}x\right)  ^{6}\left\vert x\right\vert ^{-\left(
s+12\right)  }\right\} \\
&  =15\left(  \widehat{a}x\right)  \left\vert x\right\vert ^{-\left(
s+8\right)  }+45\left\{  2\left(  \widehat{a}x\right)  \left\vert x\right\vert
^{-\left(  s+8\right)  }-\left(  s+8\right)  \left(  \widehat{a}x\right)
^{3}\left\vert x\right\vert ^{-\left(  s+10\right)  }\right\}  -\\
&  \quad-15\left(  s+8\right)  \left\{  4\left(  \widehat{a}x\right)
^{3}\left\vert x\right\vert ^{-\left(  s+10\right)  }-\left(  s+10\right)
\left(  \widehat{a}x\right)  ^{5}\left\vert x\right\vert ^{-\left(
s+12\right)  }\right\}  +\\
&  \quad+\left(  s+8\right)  \left(  s+10\right)  \left\{  6\left(
\widehat{a}x\right)  ^{5}\left\vert x\right\vert ^{-\left(  s+12\right)
}-\left(  s+12\right)  \left(  \widehat{a}x\right)  ^{7}\left\vert
x\right\vert ^{-\left(  s+14\right)  }\right\} \\
&  =15\left(  \widehat{a}x\right)  \left\vert x\right\vert ^{-\left(
s+8\right)  }+90\left(  \widehat{a}x\right)  \left\vert x\right\vert
^{-\left(  s+8\right)  }-45\left(  s+8\right)  \left(  \widehat{a}x\right)
^{3}\left\vert x\right\vert ^{-\left(  s+10\right)  }-\\
&  \quad-60\left(  s+8\right)  \left(  \widehat{a}x\right)  ^{3}\left\vert
x\right\vert ^{-\left(  s+10\right)  }+15\left(  s+8\right)  \left(
s+10\right)  \left(  \widehat{a}x\right)  ^{5}\left\vert x\right\vert
^{-\left(  s+12\right)  }+\\
&  \quad+6\left(  s+8\right)  \left(  s+10\right)  \left(  \widehat
{a}x\right)  ^{5}\left\vert x\right\vert ^{-\left(  s+12\right)  }-\left(
s+8\right)  \left(  s+10\right)  \left(  s+12\right)  \left(  \widehat
{a}x\right)  ^{7}\left\vert x\right\vert ^{-\left(  s+14\right)  }\\
&  =105\left(  \widehat{a}x\right)  \left\vert x\right\vert ^{-\left(
s+8\right)  }-105\left(  s+8\right)  \left(  \widehat{a}x\right)
^{3}\left\vert x\right\vert ^{-\left(  s+10\right)  }+21\left(  s+8\right)
\left(  s+10\right)  \left(  \widehat{a}x\right)  ^{5}\left\vert x\right\vert
^{-\left(  s+12\right)  }-\\
&  \quad-\left(  s+8\right)  \left(  s+10\right)  \left(  s+12\right)  \left(
\widehat{a}x\right)  ^{7}\left\vert x\right\vert ^{-\left(  s+14\right)  }\\
&  =\left(  \widehat{a}\widehat{x}\right)  \left\{
\begin{array}
[c]{r}%
105-105\left(  s+8\right)  \left(  \widehat{a}\widehat{x}\right)
^{2}+21\left(  s+8\right)  \left(  s+10\right)  \left(  \widehat{a}\widehat
{x}\right)  ^{4}-\\
-\left(  s+8\right)  \left(  s+10\right)  \left(  s+12\right)  \left(
\widehat{a}\widehat{x}\right)  ^{6}%
\end{array}
\right\}  \left\vert x\right\vert ^{-\left(  s+8\right)  }.
\end{align*}

\subsection{General form of $\left(  \widehat{a}D\right)  ^{m}\left\vert
x\right\vert ^{-s}$ and related polynomials $p_{s}^{\left(  m\right)  }$}

We first introduce the Pochhammer or rising factorial symbol:

\begin{definition}
\label{Def_pochhammer_sym}\textbf{Pochhammer or rising factorial symbol
}$\left(  s\right)  _{k}$ It is convenient to define the symbol%
\begin{equation}
\left(  s\right)  _{k}:=s\left(  s+1\right)  \left(  s+2\right)  \ldots\left(
s+k-1\right)  ,\quad k=1,2,3,\ldots;\text{ }s\in\mathbb{R}^{1},\label{Ap047}%
\end{equation}

and observe that since the gamma function satisfies $\Gamma\left(  z+1\right)
=z\Gamma\left(  z\right)  $ we have%
\begin{equation}
\left(  s\right)  _{k}:=\frac{\Gamma\left(  s+k\right)  }{\Gamma\left(
s\right)  }.\label{Ap065}%
\end{equation}

This suggests we define%
\[
\left(  s\right)  _{0}=1.
\]

Indeed, we might use \ref{Ap065} to extend the definition of $\left(
s\right)  _{k}$ to negative integer $k$. In this case we have%
\begin{equation}
\left(  s\right)  _{-k}:=\frac{1}{\left(  s-1\right)  \left(  s-2\right)
\ldots\left(  s-k\right)  },\quad k=1,2,3,\ldots;\text{ }s\neq1,2,3,\ldots
,k,\label{Ap066}%
\end{equation}

which leads to the identity%
\begin{equation}
\left(  s\right)  _{-k}\left(  1-s\right)  _{k}=\left(  -1\right)
^{k}.\label{Ap067}%
\end{equation}

We note the important identity%
\[
\Gamma\left(  s\right)  \Gamma\left(  1-s\right)  =\frac{\pi}{\sin\pi s}.
\]

Finally we have the \textbf{further generalization}%
\begin{equation}
\left(  s\right)  _{t}:=\frac{\Gamma\left(  s+t\right)  }{\Gamma\left(
s\right)  },\quad s,t\in\mathbb{R}^{1},\text{ }s\neq0,-1,-2,-3,\ldots
.\label{Ap085}%
\end{equation}

\end{definition}

From these calculations we \textbf{hypothesize} that for some constants
$\left\{  b_{2k}^{2n}\right\}  _{k=0}^{n}$,%
\begin{equation}
\left(  \widehat{a}D\right)  ^{2n}\left\vert x\right\vert ^{-s}=\left(
\sum_{k=0}^{n}b_{2k}^{2n}s\left(  s+2\right)  \ldots\left(  s+2n+2k-2\right)
\left(  \widehat{a}\widehat{x}\right)  ^{2k}\right)  \left\vert x\right\vert
^{-s-2n},\quad n=1,2,3,\ldots,\label{Ap052}%
\end{equation}

which can be written%
\begin{align}
\left(  \widehat{a}D\right)  ^{2n}\left\vert x\right\vert ^{-s}  & =s\left(
s+2\right)  \left(  s+4\right)  \ldots\left(  s+2n-2\right)  p_{s}^{2n}\left(
\widehat{a}\widehat{x}\right)  \left\vert x\right\vert ^{-\left(  s+2n\right)
}\nonumber\\
& =2^{n}\left(  \frac{s}{2}\right)  _{n}p_{s}^{2n}\left(  \widehat{a}%
\widehat{x}\right)  \left\vert x\right\vert ^{-\left(  s+2n\right)
},\label{Ap044}%
\end{align}

where%
\begin{align}
p_{s}^{\left(  2n\right)  }\left(  t\right)   & =b_{0}^{2n}+\sum_{k=1}%
^{n-1}b_{2k}^{2n}\left\{  \left(  s+2n\right)  \left(  s+2n+2\right)
\ldots\left(  s+2n+2k-2\right)  \right\}  t^{2k}+\nonumber\\
& \qquad+b_{2n}^{2n}\left\{  \left(  s+2n\right)  \left(  s+2n+2\right)
\ldots\left(  s+4n-2\right)  \right\}  t^{2n}\nonumber\\
& =b_{0}^{2n}+\sum_{k=1}^{n-1}b_{2k}^{2n}2^{k}\left\{  \left(  \frac{s}%
{2}+n\right)  \left(  \frac{s}{2}+n+1\right)  \ldots\left(  \frac{s}%
{2}+n+k-1\right)  \right\}  \left(  \widehat{a}\widehat{x}\right)
^{2k}+\nonumber\\
& \qquad+b_{2n}^{2n}\left\{  \left(  \frac{s}{2}+n\right)  \left(  \frac{s}%
{2}+n+1\right)  \ldots\left(  \frac{s}{2}+n+n-1\right)  \right\}
t^{2n}\nonumber\\
& =b_{0}^{2n}+\sum_{k=1}^{n}b_{2k}^{2n}2^{k}\left(  \frac{s}{2}+n\right)
_{k}t^{2k}\nonumber\\
& =\sum_{k=0}^{n}b_{2k}^{2n}2^{k}\left(  \frac{s}{2}+n\right)  _{k}%
t^{2k},\label{Ap046}%
\end{align}

We also \textbf{hypothesize} that for some constants $\left\{  b_{2k+1}%
^{2n+1}\right\}  _{k=0}^{n}$,%
\begin{equation}
\left(  \widehat{a}D\right)  ^{2n+1}\left\vert x\right\vert ^{-s}=\left(
\sum_{k=0}^{n}b_{2k+1}^{2n+1}\left\{  s\left(  s+2\right)  \ldots\left(
s+2n+2k\right)  \right\}  \left(  \widehat{a}\widehat{x}\right)
^{2k+1}\right)  \left\vert x\right\vert ^{-s-\left(  2n+1\right)  },\quad
n=0,1,2,\ldots,\label{Ap045}%
\end{equation}

so that when $n=0$,%
\begin{equation}
\left(  \widehat{a}D\right)  \left\vert x\right\vert ^{-s}=b_{1}^{1}s\left(
\widehat{a}\widehat{x}\right)  \left\vert x\right\vert ^{-s-1},\label{Ap043}%
\end{equation}

and when $n\geq0$,%
\begin{align}
&  \left(  \widehat{a}D\right)  ^{2n+1}\left\vert x\right\vert ^{-s}%
\nonumber\\
&  =\left(  \sum_{k=0}^{n}b_{2k+1}^{2n+1}\left\{  s\left(  s+2\right)
\ldots\left(  s+2n+2k\right)  \right\}  \left(  \widehat{a}\widehat{x}\right)
^{2k+1}\right)  \left\vert x\right\vert ^{-s-\left(  2n+1\right)
}\label{Ap050}\\
&  =\left\{  s\left(  s+2\right)  \left(  s+4\right)  \ldots\left(
s+2n\right)  \right\}  \left(  \widehat{a}\widehat{x}\right)  \times
\nonumber\\
&  \quad\times\left(  b_{1}^{2n+1}+\sum_{k=1}^{n}b_{2k+1}^{2n+1}\left\{
\left(  s+2n+2\right)  \left(  s+2n+4\right)  \ldots\left(  s+2n+2k\right)
\right\}  \left(  \widehat{a}\widehat{x}\right)  ^{2k}\right)  \left\vert
x\right\vert ^{-s-\left(  2n+1\right)  }\nonumber\\
&  =2^{n}\frac{s}{2}\left(  \frac{s}{2}+1\right)  \ldots\left(  \frac{s}%
{2}+n\right)  \left(  \widehat{a}\widehat{x}\right)  \times\nonumber\\
&  \quad\times\left(  b_{1}^{2n+1}+\sum_{k=1}^{n}b_{2k+1}^{2n+1}2^{k}\left\{
\left(  \frac{s}{2}+n+1\right)  \left(  \frac{s}{2}+n+2\right)  \ldots\left(
\frac{s}{2}+n+k\right)  \right\}  \left(  \widehat{a}\widehat{x}\right)
^{2k}\right)  \left\vert x\right\vert ^{-s-\left(  2n+1\right)  }\nonumber\\
&  =2^{n+1}\left(  \frac{s}{2}\right)  _{n+1}\left(  \widehat{a}\widehat
{x}\right)  \left(  b_{1}^{2n+1}+\sum_{k=1}^{n}b_{2k+1}^{2n+1}2^{k}\left(
\frac{s}{2}+n+1\right)  _{k}\left(  \widehat{a}\widehat{x}\right)
^{2k}\right)  \left\vert x\right\vert ^{-s-\left(  2n+1\right)  }\nonumber\\
&  =2^{n+1}\left(  \frac{s}{2}\right)  _{n+1}\left(  \widehat{a}\widehat
{x}\right)  \left(  \sum_{k=0}^{n}b_{2k+1}^{2n+1}2^{k}\left(  \frac{s}%
{2}+n+1\right)  _{k}\left(  \widehat{a}\widehat{x}\right)  ^{2k}\right)
\left\vert x\right\vert ^{-s-\left(  2n+1\right)  },\label{Ap042}%
\end{align}

and we define the polynomial%
\begin{equation}
p_{s}^{2n+1}\left(  t\right)  =\sum_{k=0}^{n}b_{2k+1}^{2n+1}2^{k}\left(
\frac{s}{2}+n+1\right)  _{k}t^{2k},\quad n=0,1,2,\ldots,\label{Ap041}%
\end{equation}

so that%
\begin{equation}
\left(  \widehat{a}D\right)  ^{2n+1}\left\vert x\right\vert ^{-s}%
=2^{n+1}\left(  \frac{s}{2}\right)  _{n+1}\left(  \widehat{a}\widehat
{x}\right)  p_{s}^{2n+1}\left(  \widehat{a}\widehat{x}\right)  \left\vert
x\right\vert ^{-s-\left(  2n+1\right)  }.\label{Ap053}%
\end{equation}

From the calculations at the start of this subsection:%
\begin{equation}
\left.
\begin{array}
[c]{ll}%
p_{s}^{1}\left(  t\right)  =-1, & b^{1}=-1,\\
p_{s}^{2}\left(  t\right)  =-1+\left(  s+2\right)  t^{2}, & b_{even}%
^{2}=-1,1,\\
& \\
p_{s}^{3}\left(  t\right)  =3-\left(  s+4\right)  t^{2}, & b_{odd}^{3}=3,-1,\\
p_{s}^{4}\left(  t\right)  =3-6\left(  s+4\right)  t^{2}+\left(  s+4\right)
\left(  s+6\right)  t^{4}, & b_{even}^{4}=3,-6,1,\\
& \\
p_{s}^{5}\left(  t\right)  =-15+10\left(  s+6\right)  t^{2}-\left(
s+6\right)  \left(  s+8\right)  t^{4}, & b_{odd}^{5}=-15,10,-1,\\
\left.
\begin{array}
[c]{r}%
p_{s}^{6}\left(  t\right)  =-15+45\left(  s+6\right)  t^{2}-15\left(
s+6\right)  \left(  s+8\right)  t^{4}+\\
+\left(  s+6\right)  \left(  s+8\right)  \left(  s+10\right)  t^{6},
\end{array}
\right\}  & b_{even}^{6}=-15,45,-15,1,\\
& \\
\left.
\begin{array}
[c]{r}%
p_{s}^{7}\left(  t\right)  =105-105\left(  s+8\right)  t^{2}+21\left(
s+8\right)  \left(  s+10\right)  t^{4}-\\
-\left(  s+8\right)  \left(  s+10\right)  \left(  s+12\right)  t^{6},
\end{array}
\right\}  & b_{odd}^{7}=105,-105,21,-1.
\end{array}
\right\} \label{Ap040}%
\end{equation}

\subsection{Calculating $p_{s}^{m}\left(  0\right)  $}

From the equations \ref{Ap040} we guess that:

\begin{conjecture}
\label{Conj_p^m_s_at_zero}For $n=1,2,3,\ldots$,%
\begin{align*}
p_{s}^{2n}\left(  0\right)   & =p_{s}^{2n-1}\left(  0\right)  .\\
p_{s}^{2n+1}\left(  0\right)   & =-\left(  2n+1\right)  p_{s}^{2n-1}\left(
0\right)  .
\end{align*}

which is equivalent to%
\begin{align*}
b_{0}^{2n}  & =b_{1}^{2n-1}.\\
b_{1}^{2n+1}  & =-\left(  2n+1\right)  b_{1}^{2n-1}.
\end{align*}

Also, we include here%
\begin{equation}
b_{2}^{2n}=-\left(  n-1\right)  b_{0}^{2n},\quad n\geq2.\label{Ap084}%
\end{equation}

\end{conjecture}

This conjecture implies%
\begin{align*}
p_{s}^{2n+1}\left(  0\right)   & =\left(  2n+1\right)  \left(  2n-1\right)
p_{s}^{2n-3}\left(  0\right)  =\ldots=\\
& =\left(  -1\right)  ^{n}\left(  2n+1\right)  \left(  2n-1\right)
\ldots\left(  3\right)  \left(  1\right)  p_{s}^{1}\left(  0\right) \\
& =\left(  -1\right)  ^{n+1}\left(  2n+1\right)  \left(  2n-1\right)
\ldots\left(  3\right)  \left(  1\right) \\
& =\left(  -1\right)  ^{n+1}\frac{\left(  2n+1\right)  !}{\left(  2n\right)
\left(  2n-2\right)  \ldots2}\\
& =\left(  -1\right)  ^{n+1}\frac{\left(  2n+1\right)  !}{2^{n}n!}\\
& =\left(  -1\right)  ^{n+1}\left(  2n+1\right)  !!,
\end{align*}

or better numerically,%
\begin{align}
p_{s}^{2n+1}\left(  0\right)   & =\left(  -1\right)  ^{n+1}\left(  1\right)
\left(  1+2\right)  \ldots\left(  1+2n-2\right)  \left(  1+2n\right)
\nonumber\\
& =\left(  -1\right)  ^{n+1}2^{n+1}\left(  \frac{1}{2}\right)  \left(
\frac{1}{2}+1\right)  \ldots\left(  \frac{1}{2}+n-1\right)  \left(  \frac
{1}{2}+n\right) \nonumber\\
& =\left(  -1\right)  ^{n+1}2^{n+1}\left(  \frac{1}{2}\right)  _{n+1}%
,\label{Ap038}%
\end{align}

and so we \textbf{guess} that%
\begin{equation}
p_{s}^{2n}\left(  0\right)  =p_{s}^{2n-1}\left(  0\right)  =\left(  -1\right)
^{n}2^{n}\left(  \frac{1}{2}\right)  _{n},\quad n\geq1;\text{ }s\in
\mathbb{R}^{1}.\label{Ap031}%
\end{equation}

I have \textbf{no proof} of \ref{Ap031} but this implies that%
\begin{equation}
b_{0}^{2n}=b_{1}^{2n-1}=\left(  -1\right)  ^{n}2^{n}\left(  \frac{1}%
{2}\right)  _{n},\quad n\geq1.\label{Ap083}%
\end{equation}

and so guess \ref{Ap084} implies%
\[
b_{2}^{2n}=\left(  -1\right)  ^{n+1}2^{n}\left(  n-1\right)  \left(  \frac
{1}{2}\right)  _{n},\quad n\geq2.
\]

\begin{align*}
b_{2k}^{2n}  & =b_{2k-2}^{2n-2}-\left(  4k+1\right)  b_{2k}^{2n-2}+\left(
2k+1\right)  \left(  2k+2\right)  b_{2k+2}^{2n-2},\Rightarrow\\
b_{2}^{2n+2}  & =b_{0}^{2n}-5b_{2}^{2n}+12b_{4}^{2n},\\
-nb_{0}^{2n+2}  & =b_{0}^{2n}+5\left(  n-1\right)  b_{0}^{2n}+12b_{4}^{2n},\\
-nb_{0}^{2n+2}  & =\left(  5n-4\right)  b_{0}^{2n}+12b_{4}^{2n}\\
\left(  -1\right)  ^{n}2^{n+1}n\left(  \frac{1}{2}\right)  _{n+1}  & =\left(
-1\right)  ^{n}2^{n}\left(  5n-4\right)  \left(  \frac{1}{2}\right)
_{n}+12b_{4}^{2n}\\
\left(  -1\right)  ^{n}2^{n}n\left(  2n+1\right)  \left(  \frac{1}{2}\right)
_{n}  & =\left(  -1\right)  ^{n}2^{n}\left(  5n-4\right)  \left(  \frac{1}%
{2}\right)  _{n}+12b_{4}^{2n}\\
\left(  -1\right)  ^{n}2^{n}\left(  2n^{2}-4n+4\right)  \left(  \frac{1}%
{2}\right)  _{n}  & =12b_{4}^{2n}%
\end{align*}

\subsection{A bound for $\left(  \widehat{a}D\right)  ^{2n+1}\left\vert
x\right\vert ^{-s}$, $s\geq-4$.}

For $s\geq-4$ we will obtain the bounds \ref{Ap057} and \ref{Ap062} which have
the form
\[
\left\vert \left(  \widehat{a}D\right)  ^{2n+1}\left\vert x\right\vert
^{-s}\right\vert \leq C_{s,n}\left\vert \widehat{a}\widehat{x}\right\vert
\left\vert x\right\vert ^{-s-\left(  2n+1\right)  }.
\]

We will need the identities:%
\begin{align}
\left(  \widehat{a}D\right)  \left(  \widehat{a}\widehat{x}\right)   &
=\left(  1-\left(  \widehat{a}\widehat{x}\right)  ^{2}\right)  \left\vert
x\right\vert ^{-1},\nonumber\\
\left(  \widehat{a}D\right)  \left(  \widehat{a}\widehat{x}\right)  ^{k}  &
=k\left(  \left(  \widehat{a}\widehat{x}\right)  ^{k-1}-\left(  \widehat
{a}\widehat{x}\right)  ^{k+1}\right)  \left\vert x\right\vert ^{-1}%
,\nonumber\\
\left(  \widehat{a}D\right)  \left\vert x\right\vert ^{-t}  & =-t\left(
\widehat{a}\widehat{x}\right)  \left\vert x\right\vert ^{-\left(  t+1\right)
},\nonumber\\
\left(  \widehat{a}D\right)  \left(  \left(  \widehat{a}\widehat{x}\right)
^{k}\left\vert x\right\vert ^{-t}\right)   & =\left(  \left(  \widehat
{a}D\right)  \left(  \widehat{a}\widehat{x}\right)  ^{k}\right)  \left\vert
x\right\vert ^{-t}+\left(  \widehat{a}\widehat{x}\right)  ^{k}\left(
\widehat{a}D\right)  \left\vert x\right\vert ^{-t}\nonumber\\
& =k\left(  \left(  \widehat{a}\widehat{x}\right)  ^{k-1}-\left(  \widehat
{a}\widehat{x}\right)  ^{k+1}\right)  \left\vert x\right\vert ^{-\left(
t+1\right)  }-\left(  \widehat{a}\widehat{x}\right)  ^{k}t\left(  \widehat
{a}\widehat{x}\right)  \left\vert x\right\vert ^{-\left(  t+1\right)
}\nonumber\\
& =\left\{  k\left(  \left(  \widehat{a}\widehat{x}\right)  ^{k-1}-\left(
\widehat{a}\widehat{x}\right)  ^{k+1}\right)  -t\left(  \widehat{a}\widehat
{x}\right)  ^{k+1}\right\}  \left\vert x\right\vert ^{-\left(  t+1\right)
}\nonumber\\
& =\left\{  k\left(  \widehat{a}\widehat{x}\right)  ^{k-1}-\left(  t+k\right)
\left(  \widehat{a}\widehat{x}\right)  ^{k+1}\right\}  \left\vert x\right\vert
^{-\left(  t+1\right)  }.\nonumber\\
i.e.\text{ }\left(  \widehat{a}D\right)  \left(  \left(  \widehat{a}%
\widehat{x}\right)  ^{k}\left\vert x\right\vert ^{-t}\right)   & =\left\{
k\left(  \widehat{a}\widehat{x}\right)  ^{k-1}-\left(  t+k\right)  \left(
\widehat{a}\widehat{x}\right)  ^{k+1}\right\}  \left\vert x\right\vert
^{-\left(  t+1\right)  }.\label{Ap058}%
\end{align}

Then from \ref{Ap052}\ and then the identity \ref{Ap058},%
\begin{align*}
&  \left(  \widehat{a}D\right)  ^{2n+1}\left\vert x\right\vert ^{-s}\\
&  =\left(  \widehat{a}D\right)  \left(  \widehat{a}D\right)  ^{2n}\left\vert
x\right\vert ^{-s}\\
&  =\sum_{k=0}^{n}b_{2k}^{2n}s\left(  s+2\right)  \ldots\left(
s+2n+2k-2\right)  \left(  \widehat{a}D\right)  \left\{  \left(  \widehat
{a}\widehat{x}\right)  ^{2k}\left\vert x\right\vert ^{-s-2n}\right\} \\
&  =b_{0}^{2n}s\left(  s+2\right)  \ldots\left(  s+2n-2\right)  \left(
\widehat{a}D\right)  \left\vert x\right\vert ^{-\left(  s+2n\right)  }+\\
&  \qquad+\sum_{k=1}^{n}b_{2k}^{2n}s\left(  s+2\right)  \ldots\left(
s+2n+2k-2\right)  \left(  \widehat{a}D\right)  \left(  \left(  \widehat
{a}\widehat{x}\right)  ^{2k}\left\vert x\right\vert ^{-\left(  s+2n\right)
}\right) \\
&  =\left(
\begin{array}
[c]{l}%
-b_{0}^{2n}s\left(  s+2\right)  \ldots\left(  s+2n-2\right)  \left(
s+2n\right)  \left(  \widehat{a}\widehat{x}\right)  +\\
+\sum\limits_{k=1}^{n}b_{2k}^{2n}s\left(  s+2\right)  \ldots\left(
s+2n+2k-2\right)  \left\{  2k\left(  \widehat{a}\widehat{x}\right)
^{2k-1}-\left(  s+2n+2k\right)  \left(  \widehat{a}\widehat{x}\right)
^{2k+1}\right\}
\end{array}
\right)  \left\vert x\right\vert ^{-\left(  s+2n+1\right)  }\\
&  =\left(
\begin{array}
[c]{l}%
-b_{0}^{2n}s\left(  s+2\right)  \ldots\left(  s+2n-2\right)  \left(
s+2n\right)  \left(  \widehat{a}\widehat{x}\right)  +\\
+\sum\limits_{k=1}^{n}b_{2k}^{2n}s\left(  s+2\right)  \ldots\left(
s+2n+2k-2\right)  2k\left(  \widehat{a}\widehat{x}\right)  ^{2k-1}-\\
-\sum\limits_{k=1}^{n}b_{2k}^{2n}s\left(  s+2\right)  \ldots\left(
s+2n+2k\right)  \left(  \widehat{a}\widehat{x}\right)  ^{2k+1}%
\end{array}
\right)  \left\vert x\right\vert ^{-\left(  s+2n+1\right)  }\\
&  =\left(
\begin{array}
[c]{l}%
\sum\limits_{k=0}^{n-1}\left(  2k+2\right)  b_{2k+2}^{2n}s\left(  s+2\right)
\ldots\left(  s+2n+2k\right)  \left(  \widehat{a}\widehat{x}\right)
^{2k+1}-\\
-\sum\limits_{k=0}^{n}b_{2k}^{2n}s\left(  s+2\right)  \ldots\left(
s+2n+2k\right)  \left(  \widehat{a}\widehat{x}\right)  ^{2k+1}%
\end{array}
\right)  \left\vert x\right\vert ^{-\left(  s+2n+1\right)  }\\
&  =\left(
\begin{array}
[c]{l}%
\sum\limits_{k=0}^{n-1}\left(  2k+2\right)  b_{2k+2}^{2n}s\left(  s+2\right)
\ldots\left(  s+2n+2k\right)  \left(  \widehat{a}\widehat{x}\right)
^{2k+1}-\\
-\sum\limits_{k=0}^{n-1}b_{2k}^{2n}s\left(  s+2\right)  \ldots\left(
s+2n+2k\right)  \left(  \widehat{a}\widehat{x}\right)  ^{2k+1}-\\
-b_{2n}^{2n}s\left(  s+2\right)  \ldots\left(  s+4n-2\right)  \left(
s+4n\right)  \left(  \widehat{a}\widehat{x}\right)  ^{2n+1}%
\end{array}
\right)  \left\vert x\right\vert ^{-\left(  s+2n+1\right)  }\\
&  =\left(
\begin{array}
[c]{l}%
\sum\limits_{k=0}^{n-1}s\left(  s+2\right)  \ldots\left(  s+2n+2k\right)
\left(  -b_{2k}^{2n}+\left(  2k+2\right)  b_{2k+2}^{2n}\right)  \left(
\widehat{a}\widehat{x}\right)  ^{2k+1}-\\
-b_{2n}^{2n}s\left(  s+2\right)  \ldots\left(  s+4n\right)  \left(
\widehat{a}\widehat{x}\right)  ^{2n+1}%
\end{array}
\right)  \left\vert x\right\vert ^{-\left(  s+2n+1\right)  },
\end{align*}

but from \ref{Ap045},%
\[
\left(  \widehat{a}D\right)  ^{2n+1}\left\vert x\right\vert ^{-s}=\left(
\sum_{k=0}^{n}b_{2k+1}^{2n+1}\left\{  s\left(  s+2\right)  \ldots\left(
s+2n+2k\right)  \right\}  \left(  \widehat{a}\widehat{x}\right)
^{2k+1}\right)  \left\vert x\right\vert ^{-s-\left(  2n+1\right)  },
\]

so comparing powers of $\widehat{a}\widehat{x}$ we get
\begin{equation}%
\begin{array}
[c]{ll}%
b_{2k+1}^{2n+1}=-b_{2k}^{2n}+\left(  2k+2\right)  b_{2k+2}^{2n}, &
k=0,1,2,\ldots,n-1.\\
b_{2n+1}^{2n+1}=-b_{2n}^{2n}. & \\
b_{1}^{2n+1}=-b_{0}^{2n}+2b_{2}^{2n}. &
\end{array}
\label{Ap028}%
\end{equation}

From \ref{Ap041} and \ref{Ap040},%
\begin{equation}
p_{s}^{\left(  2n+1\right)  }\left(  t\right)  =\left\{
\begin{array}
[c]{ll}%
b_{1}^{1}=-1, & n=0,\\
b_{1}^{2n+1}+\sum\limits_{k=1}^{n}b_{2k+1}^{2n+1}2^{k}\left(  \frac{s}%
{2}+n+1\right)  _{k}t^{2k}, & n\geq1,
\end{array}
\right. \label{Ap030}%
\end{equation}

and from \ref{Ap053},%
\begin{equation}
\left(  \widehat{a}D\right)  ^{2n+1}\left\vert x\right\vert ^{-s}%
=2^{n+1}\left(  \frac{s}{2}\right)  _{n+1}\left(  \widehat{a}\widehat
{x}\right)  p_{s}^{\left(  2n+1\right)  }\left(  \widehat{a}\widehat
{x}\right)  \left\vert x\right\vert ^{-s-\left(  2n+1\right)  },\quad
n\geq0,\label{Ap022}%
\end{equation}

so that when $x=a$,%
\[
\left(  \widehat{x}D\right)  ^{2n+1}\left\vert x\right\vert ^{-s}%
=2^{n+1}\left(  \frac{s}{2}\right)  _{n+1}p_{s}^{\left(  2n+1\right)  }\left(
1\right)  \left\vert x\right\vert ^{-s-\left(  2n+1\right)  },\quad n\geq0.
\]

But part 5 of Lemma \ref{Lem_deriv_rad_funcs} means that%
\[
\left(  \widehat{x}D\right)  ^{2n+1}\left\vert x\right\vert ^{-s}=-\left(
s\right)  _{2n+1}\left\vert x\right\vert ^{-s-\left(  2n+1\right)  },\quad
n=0,1,2,3,\ldots,
\]

so that%
\begin{align*}
2^{n+1}\left(  \frac{s}{2}\right)  _{n+1}p_{s}^{\left(  2n+1\right)  }\left(
1\right)  \left\vert x\right\vert ^{-s-\left(  2n+1\right)  }  & =-\left(
s\right)  _{2n+1}\left\vert x\right\vert ^{-s-\left(  2n+1\right)
}\Rightarrow\\
2^{n+1}\left(  \frac{s}{2}\right)  _{n+1}p_{s}^{\left(  2n+1\right)  }\left(
1\right)   & =-\left(  s\right)  _{2n+1}\Rightarrow\\
p_{s}^{\left(  2n+1\right)  }\left(  1\right)   & =-\frac{\left(  s\right)
_{2n+1}}{2^{n+1}\left(  \frac{s}{2}\right)  _{n+1}},
\end{align*}

and we simplify the right side as follows:%
\begin{align*}
\frac{\left(  s\right)  _{2n+1}}{2^{n+1}\left(  \frac{s}{2}\right)  _{n+1}}  &
=\frac{1}{2^{n+1}}\frac{s\left(  s+1\right)  \left(  s+2\right)  \ldots\left(
s+2n\right)  }{\left(  \frac{s}{2}\right)  \left(  \frac{s}{2}+1\right)
\ldots\left(  \frac{s}{2}+n\right)  }\\
& =\frac{s\left(  s+1\right)  \left(  s+2\right)  \ldots\left(  s+2n\right)
}{s\left(  s+2\right)  \ldots\left(  s+2n\right)  }\\
& =\left(  s+1\right)  \left(  s+3\right)  \ldots\left(  s+2n-1\right) \\
& =\left(  s+1\right)  \left(  s+1+2\right)  \ldots\left(  s+1+2k\right)
\ldots\left(  s+1+2n-2\right) \\
& =2^{n}\left(  \frac{s+1}{2}\right)  \left(  \frac{s+1}{2}+1\right)
\ldots\left(  \frac{s+1}{2}+k\right)  \ldots\left(  \frac{s+1}{2}+n-1\right)
\\
& =2^{n}\left(  \frac{s+1}{2}\right)  _{n},
\end{align*}

so that%
\begin{equation}
p_{s}^{\left(  2n+1\right)  }\left(  1\right)  =-\frac{\left(  s\right)
_{2n+1}}{2^{n+1}\left(  \frac{s}{2}\right)  _{n+1}}=2^{n}\left(  \frac{s+1}%
{2}\right)  _{n}.\label{Ap026}%
\end{equation}

We now need:

\begin{conjecture}
\label{Conj_argmax_p_odd}From the results of Matlab numerical experiments up
to $n=20$ we conjecture that for $n\geq0$,%
\[
0\in\operatorname*{argmax}_{t\in\left[  -1,1\right]  }\left\vert
p_{s}^{\left(  2n+1\right)  }\left(  t\right)  \right\vert ,\quad-4\leq
s\leq2,
\]

and%
\[
1\in\operatorname*{argmax}_{t\in\left[  -1,1\right]  }\left\vert
p_{s}^{\left(  2n+1\right)  }\left(  t\right)  \right\vert ,\quad s\geq2.
\]

\end{conjecture}

\textbf{Provided this conjecture is true} and in the light of \ref{Ap026} and
\ref{Ap038},%
\begin{align}
\max_{t\in\left[  -1,1\right]  }\left\vert p_{s}^{\left(  2n+1\right)
}\left(  t\right)  \right\vert  & =\left\{
\begin{array}
[c]{ll}%
\left\vert p_{s}^{\left(  2n+1\right)  }\left(  0\right)  \right\vert , &
-4\leq s\leq2,\\
\left\vert p_{s}^{\left(  2n+1\right)  }\left(  1\right)  \right\vert , &
s\geq2,
\end{array}
\right. \nonumber\\
& =\left\{
\begin{array}
[c]{ll}%
2^{n+1}\left(  \frac{1}{2}\right)  _{n+1}, & -4\leq s\leq2,\\
2^{n}\left(  \frac{s+1}{2}\right)  _{n}=\frac{\left(  s\right)  _{2n+1}%
}{2^{n+1}\left(  \frac{s}{2}\right)  _{n+1}}, & s\geq2.
\end{array}
\right. \label{Ap054}%
\end{align}

Regarding \ref{Ap054}, we note that when $s=2$:%
\begin{align*}
\left(  \frac{s+1}{2}\right)  _{n}=\left(  \frac{2+1}{2}\right)  _{n}  &
=2\left(  \frac{2+1}{2}\right)  \left(  \frac{2+1}{2}+1\right)  \ldots\left(
\frac{2+1}{2}+n-1\right) \\
& =2\left(  \frac{1}{2}+1\right)  \left(  \frac{1}{2}+2\right)  \ldots\left(
\frac{1}{2}+n\right) \\
& =2\left(  \frac{1}{2}\right)  _{n+1},
\end{align*}

which is consistent.

We now apply our conjecture to \ref{Ap022} using \ref{Ap054},%
\begin{align}
\left\vert \left(  \widehat{a}D\right)  ^{2n+1}\left\vert x\right\vert
^{-s}\right\vert  & =2^{n+1}\left\vert \left(  \frac{s}{2}\right)
_{n+1}\right\vert \left\vert \widehat{a}\widehat{x}\right\vert \left\vert
p_{s}^{\left(  2n+1\right)  }\left(  \widehat{a}\widehat{x}\right)
\right\vert \left\vert x\right\vert ^{-s-\left(  2n+1\right)  }\nonumber\\
& \leq2^{n+1}\left\vert \left(  \frac{s}{2}\right)  _{n+1}\right\vert
\left\vert p_{s}^{\left(  2n+1\right)  }\left(  \widehat{a}\widehat{x}\right)
\right\vert \left\vert \widehat{a}\widehat{x}\right\vert \left\vert
x\right\vert ^{-s-\left(  2n+1\right)  }\nonumber\\
& \leq2^{n+1}\left\vert \left(  \frac{s}{2}\right)  _{n+1}\right\vert
\max_{t\in\left[  -1,1\right]  }\left\vert p_{s}^{\left(  2n+1\right)
}\left(  t\right)  \right\vert \left\vert \widehat{a}\widehat{x}\right\vert
\left\vert x\right\vert ^{-s-\left(  2n+1\right)  }\nonumber\\
& =\left\{
\begin{array}
[c]{ll}%
2^{2n+2}\left\vert \left(  \frac{s}{2}\right)  _{n+1}\right\vert \left(
\frac{1}{2}\right)  _{n+1}, & -4\leq s\leq2,\\
2^{2n+1}\left(  s\right)  _{2n+1}, & s\geq2,
\end{array}
\right\}  \left\vert \widehat{a}\widehat{x}\right\vert \left\vert x\right\vert
^{-s-\left(  2n+1\right)  }\label{Ap062}\\
& =\left\{
\begin{array}
[c]{ll}%
2^{2n+2}\left\vert \left(  \frac{s}{2}\right)  _{n+1}\right\vert \left(
\frac{1}{2}\right)  _{n+1}, & -4\leq s\leq2,\\
2^{2n+1}\left(  \frac{s}{2}\right)  _{n+1}\left(  \frac{s+1}{2}\right)
_{n}, & s\geq2,
\end{array}
\right\}  \left\vert \widehat{a}\widehat{x}\right\vert \left\vert x\right\vert
^{-s-\left(  2n+1\right)  }\label{Ap057}\\
for\text{ }n  & \geq0.\nonumber
\end{align}

?? Equality when $a=\lambda x$? ??

\subsection{Estimating $\left(  \widehat{a}D\right)  ^{2n}\left\vert
x\right\vert ^{-s}$, $s\geq-2$.}

For $s\geq-2$ we will obtain the bounds \ref{Ap055} and \ref{Ap061} which have
the form $\left\vert \left(  \widehat{a}D\right)  ^{2n}\left\vert x\right\vert
^{-s}\right\vert \leq C_{s,2n}^{\prime}\left\vert x\right\vert ^{-s-2n}$.

From \ref{Ap045}\ and then employing the identity \ref{Ap058},%
\begin{align*}
&  \left(  \widehat{a}D\right)  ^{2n}\left\vert x\right\vert ^{-s}\\
&  =\widehat{a}D\left(  \widehat{a}D\right)  ^{2n-1}\left(  \left\vert
x\right\vert ^{-s}\right) \\
&  =\widehat{a}D\left(  \left(  \sum_{k=0}^{n-1}b_{2k+1}^{2n-1}\left\{
s\left(  s+2\right)  \ldots\left(  s+2n+2k-2\right)  \right\}  \left(
\widehat{a}\widehat{x}\right)  ^{2k+1}\right)  \left\vert x\right\vert
^{-s-\left(  2n-1\right)  }\right) \\
&  =\sum_{k=0}^{n-1}b_{2k+1}^{2n-1}\left\{  s\left(  s+2\right)  \ldots\left(
s+2n+2k-2\right)  \right\}  \widehat{a}D\left(  \left(  \widehat{a}\widehat
{x}\right)  ^{2k+1}\left\vert x\right\vert ^{-\left(  s+2n-1\right)  }\right)
\\
&  =\sum_{k=0}^{n-1}b_{2k+1}^{2n-1}\times\left(
\begin{array}
[c]{l}%
\left(  s\left(  s+2\right)  \ldots\left(  s+2n+2k-2\right)  \right)  \times\\
\times\left(  \left(  2k+1\right)  \left(  \widehat{a}\widehat{x}\right)
^{2k}-\left(  s+2n+2k\right)  \left(  \widehat{a}\widehat{x}\right)
^{2k+2}\right)
\end{array}
\right)  \left\vert x\right\vert ^{-\left(  s+2n\right)  }\\
&  =\left(
\begin{array}
[c]{r}%
\sum\limits_{k=0}^{n-1}\left(  2k+1\right)  b_{2k+1}^{2n-1}\left\{  s\left(
s+2\right)  \ldots\left(  s+2n+2k-2\right)  \right\}  \left(  \widehat
{a}\widehat{x}\right)  ^{2k}-\\
-\sum\limits_{k=0}^{n-1}b_{2k+1}^{2n-1}\left\{  s\left(  s+2\right)
\ldots\left(  s+2n+2k\right)  \right\}  \left(  \widehat{a}\widehat{x}\right)
^{2k+2}%
\end{array}
\right)  \left\vert x\right\vert ^{-\left(  s+2n\right)  }\\
&  =\left(
\begin{array}
[c]{r}%
\sum\limits_{k=0}^{n-1}\left(  2k+1\right)  b_{2k+1}^{2n-1}\left\{  s\left(
s+2\right)  \ldots\left(  s+2n+2k-2\right)  \right\}  \left(  \widehat
{a}\widehat{x}\right)  ^{2k}-\\
-\sum\limits_{k=1}^{n}b_{2k-1}^{2n-1}\left\{  s\left(  s+2\right)
\ldots\left(  s+2n+2k-2\right)  \right\}  \left(  \widehat{a}\widehat
{x}\right)  ^{2k}%
\end{array}
\right)  \left\vert x\right\vert ^{-\left(  s+2n\right)  }\\
&  =\left(
\begin{array}
[c]{l}%
b_{1}^{2n-1}\left\{  s\left(  s+2\right)  \ldots\left(  s+2n-2\right)
\right\}  +\\
+\sum\limits_{k=1}^{n-1}\left(  2k+1\right)  b_{2k+1}^{2n-1}\left\{  s\left(
s+2\right)  \ldots\left(  s+2n+2k-2\right)  \right\}  \left(  \widehat
{a}\widehat{x}\right)  ^{2k}-\\
-\sum\limits_{k=1}^{n-1}b_{2k-1}^{2n-1}\left\{  s\left(  s+2\right)
\ldots\left(  s+2n+2k-2\right)  \right\}  \left(  \widehat{a}\widehat
{x}\right)  ^{2k}-\\
-b_{2n-1}^{2n-1}\left\{  s\left(  s+2\right)  \ldots\left(  s+4n-2\right)
\right\}  \left(  \widehat{a}\widehat{x}\right)  ^{2n}%
\end{array}
\right)  \left\vert x\right\vert ^{-\left(  s+2n\right)  }\\
&  =\left(
\begin{array}
[c]{l}%
b_{1}^{2n-1}\left\{  s\left(  s+2\right)  \ldots\left(  s+2n-2\right)
\right\}  +\\
+\sum\limits_{k=1}^{n-1}\left\{  -b_{2k-1}^{2n-1}+\left(  2k+1\right)
b_{2k+1}^{2n-1}\right\}  \left\{  s\left(  s+2\right)  \ldots\left(
s+2n+2k\right)  \right\}  \left(  \widehat{a}\widehat{x}\right)  ^{2k}-\\
-b_{2n-1}^{2n-1}\left\{  s\left(  s+2\right)  \ldots\left(  s+4n-2\right)
\right\}  \left(  \widehat{a}\widehat{x}\right)  ^{2n}%
\end{array}
\right)  \left\vert x\right\vert ^{-\left(  s+2n\right)  },
\end{align*}

but from \ref{Ap052},%
\[
\left(  \widehat{a}D\right)  ^{2n}\left\vert x\right\vert ^{-s}=\left(
\begin{array}
[c]{l}%
b_{0}^{2n}\left\{  s\left(  s+2\right)  \ldots\left(  s+2n-2\right)  \right\}
+\\
+\sum\limits_{k=1}^{n-1}b_{2k}^{2n}\left\{  s\left(  s+2\right)  \ldots\left(
s+2n+2k-2\right)  \right\}  \left(  \widehat{a}\widehat{x}\right)  ^{2k}+\\
+b_{2n}^{2n}\left\{  s\left(  s+2\right)  \ldots\left(  s+4n-2\right)
\right\}  \left(  \widehat{a}\widehat{x}\right)  ^{2n}%
\end{array}
\right)  \times\left\vert x\right\vert ^{-\left(  s+2n\right)  },
\]

so comparing powers of $\widehat{a}\widehat{x}$ we conclude that%
\begin{equation}
\left.
\begin{array}
[c]{l}%
b_{0}^{2n}=b_{1}^{2n-1},\\
b_{2k}^{2n}=-b_{2k-1}^{2n-1}+\left(  2k+1\right)  b_{2k+1}^{2n-1},\\
b_{2n}^{2n}=-b_{2n-1}^{2n-1}.
\end{array}
\right\} \label{Ap033}%
\end{equation}

Now define%
\begin{align}
p_{s}^{\left(  2n\right)  }\left(  t\right)   &  =b_{0}^{2n}+\sum_{k=1}%
^{n-1}b_{2k}^{2n}\left\{  \left(  s+2n\right)  \left(  s+2n+2\right)
\ldots\left(  s+2n+2k-2\right)  \right\}  t^{2k}+\nonumber\\
&  \quad+b_{2n}^{2n}\left\{  \left(  s+2n\right)  \left(  s+2n+2\right)
\ldots\left(  s+4n-2\right)  \right\}  t^{2n}\nonumber\\
&  =b_{0}^{2n}+\sum_{k=1}^{n-1}b_{2k}^{2n}2^{k}\left\{  \left(  \frac{s}%
{2}+n\right)  \left(  \frac{s}{2}+n+1\right)  \ldots\left(  \frac{s}%
{2}+n+k-1\right)  \right\}  \left(  \widehat{a}\widehat{x}\right)
^{2k}+\nonumber\\
&  \quad+b_{2n}^{2n}\left\{  \left(  \frac{s}{2}+n\right)  \left(  \frac{s}%
{2}+n+1\right)  \ldots\left(  \frac{s}{2}+n+n-1\right)  \right\}
t^{2n}\nonumber\\
&  =b_{0}^{2n}+\sum_{k=1}^{n}b_{2k}^{2n}2^{k}\left(  \frac{s}{2}+n\right)
_{k}t^{2k}\nonumber\\
&  =\sum_{k=0}^{n}b_{2k}^{2n}2^{k}\left(  \frac{s}{2}+n\right)  _{k}%
t^{2k},\label{Ap072}%
\end{align}

so that%
\begin{align}
&  \left(  \widehat{a}D\right)  ^{2n}\left\vert x\right\vert ^{-s}\nonumber\\
&  =\left\{  s\left(  s+2\right)  \ldots\left(  s+2n-2\right)  \right\}
\times\nonumber\\
&  \times\left(  b_{0}^{2n}+\sum_{k=1}^{n-1}b_{2k}^{2n}\left(
\begin{array}
[c]{r}%
\left\{  \left(  s+2n\right)  \ldots\left(  s+2n+2k-2\right)  \right\}
\left(  \widehat{a}\widehat{x}\right)  ^{2k}+\\
+b_{2n}^{2n}\left\{  \left(  s+2n\right)  \ldots\left(  s+4n-2\right)
\right\}
\end{array}
\right)  \left(  \widehat{a}\widehat{x}\right)  ^{2n}\right)  \left\vert
x\right\vert ^{-\left(  s+2n\right)  }\nonumber\\
&  =\left\{  s\left(  s+2\right)  \left(  s+4\right)  \ldots\left(
s+2n-2\right)  \right\}  p_{s}^{\left(  2n\right)  }\left(  \widehat
{a}\widehat{x}\right)  \left\vert x\right\vert ^{-\left(  s+2n\right)
}\nonumber\\
&  =2^{n}\left\{  \frac{s}{2}\left(  \frac{s}{2}+1\right)  \left(  \frac{s}%
{2}+2\right)  \ldots\left(  \frac{s}{2}+n-1\right)  \right\}  p_{s}^{\left(
2n\right)  }\left(  \widehat{a}\widehat{x}\right)  \left\vert x\right\vert
^{-\left(  s+2n\right)  }\nonumber\\
&  =2^{n}\left(  \frac{s}{2}\right)  _{n}p_{s}^{\left(  2n\right)  }\left(
\widehat{a}\widehat{x}\right)  \left\vert x\right\vert ^{-\left(  s+2n\right)
},\label{Ap059}%
\end{align}

and hence when $a=x$,%
\begin{equation}
\left(  \widehat{x}D\right)  ^{2n}\left\vert x\right\vert ^{-s}=2^{n}\left(
\frac{s}{2}\right)  _{n}p_{s}^{\left(  2n\right)  }\left(  1\right)
\left\vert x\right\vert ^{-s-2n}.\label{Ap056}%
\end{equation}

Now from part 5 of Lemma \ref{Lem_deriv_rad_funcs} (below),%
\begin{equation}
\left(  \widehat{x}D\right)  ^{2n}\left\vert x\right\vert ^{-s}=\left(
s\right)  _{2n}\left\vert x\right\vert ^{-s-2n},\quad n=0,1,2,\ldots
,\label{Ap060}%
\end{equation}

and comparing \ref{Ap056} with \ref{Ap060} we get%
\[
p_{s}^{\left(  2n\right)  }\left(  1\right)  =\frac{\left(  s\right)  _{2n}%
}{2^{n}\left(  \frac{s}{2}\right)  _{n}}.
\]

Also, when $n\geq0$,%
\begin{align*}
\frac{\left(  s\right)  _{2n}}{2^{n}\left(  \frac{s}{2}\right)  _{n}}  &
=\frac{s\left(  s+1\right)  \ldots\left(  s+2n-1\right)  }{2^{n}\frac{s}%
{2}\left(  \frac{s}{2}+1\right)  \ldots\left(  \frac{s}{2}+n-1\right)  }%
=\frac{s\left(  s+1\right)  \ldots\left(  s+2n-1\right)  }{2s\left(
s+2\right)  \ldots\left(  s+2n-2\right)  }=\\
& =\left(  s+1\right)  \left(  s+3\right)  \ldots\left(  s+2n-3\right)
\left(  s+2n-1\right) \\
& =\left(  s+1\right)  \left(  s+1+2\right)  \left(  s+1+4\right)
\ldots\left(  s+1+2n-4\right)  \left(  s+1+2n-2\right) \\
& =2^{n}\left(  \frac{s+1}{2}\right)  \left(  \frac{s+1}{2}+1\right)  \left(
\frac{s+1}{2}+2\right)  \ldots\left(  \frac{s+1}{2}+n-2\right)  \left(
\frac{s+1}{2}+n-1\right) \\
& =2^{n}\left(  \frac{s+1}{2}\right)  _{n},
\end{align*}

so that%
\begin{equation}
p_{s}^{\left(  2n\right)  }\left(  1\right)  =\frac{\left(  s\right)  _{2n}%
}{2^{n}\left(  \frac{s}{2}\right)  _{n}}=2^{n}\left(  \frac{s+1}{2}\right)
_{n},\quad n=0,1,2,\ldots\label{Ap032}%
\end{equation}

Now we need:

\begin{conjecture}
\label{Conj_argmax_p_even}From the results of Matlab numerical experiments up
to $n=20$ we conjecture that for $n\geq1$,%
\[
\operatorname*{argmax}_{t\in\left[  -1,1\right]  }\left\vert p_{s}^{\left(
2n\right)  }\left(  t\right)  \right\vert =\left\{
\begin{array}
[c]{ll}%
0, & -2\leq s\leq0,\\
1, & s\geq0.
\end{array}
\right.
\]

\end{conjecture}

This conjecture can be used as follows: clearly if $-2\leq s\leq0$ then
$\max\limits_{t\in\left[  -1,1\right]  }\left\vert p_{s}^{\left(  2n\right)
}\left(  t\right)  \right\vert =\left\vert p_{s}^{\left(  2n\right)  }\left(
0\right)  \right\vert $, and if $s\geq0$ then $\max\limits_{t\in\left[
-1,1\right]  }\left\vert p_{s}^{\left(  2n\right)  }\left(  t\right)
\right\vert =\left\vert p_{s}^{\left(  2n\right)  }\left(  1\right)
\right\vert $. Hence by \ref{Ap031} and \ref{Ap032},%
\[
\max_{t\in\left[  -1,1\right]  }\left\vert p_{s}^{\left(  2n\right)  }\left(
t\right)  \right\vert =\left\{
\begin{array}
[c]{ll}%
2^{n}\left(  \frac{1}{2}\right)  _{n}, & -2\leq s\leq0,\\
\frac{\left(  s\right)  _{2n}}{2^{n}\left(  \frac{s}{2}\right)  _{n}}%
=2^{n}\left(  \frac{s+1}{2}\right)  _{n}, & s\geq0,
\end{array}
\right.
\]

and from \ref{Ap059},%
\begin{align}
\left\vert \left(  \widehat{a}D\right)  ^{2n}\left\vert x\right\vert
^{-s}\right\vert  & \leq2^{n}\left\vert \left(  \frac{s}{2}\right)
_{n}\right\vert \max_{t\in\left[  -1,1\right]  }\left\vert p_{s}^{\left(
2n\right)  }\left(  t\right)  \right\vert \left\vert x\right\vert
^{-s-2n}\nonumber\\
& =\left\{
\begin{array}
[c]{ll}%
2^{n}\left\vert \left(  \frac{s}{2}\right)  _{n}\right\vert 2^{n}\left(
\frac{1}{2}\right)  _{n}, & -2\leq s\leq0,\\
2^{n}\left(  \frac{s}{2}\right)  _{n}\frac{\left(  s\right)  _{2n}}%
{2^{n}\left(  \frac{s}{2}\right)  _{n}}, & s\geq0,
\end{array}
\right\}  \left\vert x\right\vert ^{-s-2n}\nonumber\\
& =\left\{
\begin{array}
[c]{ll}%
2^{2n}\left\vert \left(  \frac{s}{2}\right)  _{n}\right\vert \left(  \frac
{1}{2}\right)  _{n}, & -2\leq s\leq0,\\
\left(  s\right)  _{2n}, & s\geq0,
\end{array}
\right\}  \left\vert x\right\vert ^{-s-2n}\label{Ap061}\\
& =\left\{
\begin{array}
[c]{ll}%
2^{2n}\left\vert \left(  \frac{s}{2}\right)  _{n}\right\vert \left(  \frac
{1}{2}\right)  _{n}, & -2\leq s\leq0,\\
2^{2n}\left(  \frac{s}{2}\right)  _{n}\left(  \frac{s+1}{2}\right)  _{n}, &
s\geq0,
\end{array}
\right\}  \left\vert x\right\vert ^{-s-2n}\label{Ap055}\\
for\text{ }n  & \geq0.\nonumber
\end{align}

\subsection{Summary - upper bounds for $\left(  \widehat{a}D\right)
^{k}\left\vert x\right\vert ^{-s}$\label{SbSect_bnd_(aD)^k(absX^(-s))}}

From \ref{Ap055},%
\[
\left\vert \left(  \widehat{a}D\right)  ^{2n}\left\vert x\right\vert
^{-s}\right\vert \leq\left\{
\begin{array}
[c]{ll}%
2^{2n}\left\vert \left(  \frac{s}{2}\right)  _{n}\right\vert \left(  \frac
{1}{2}\right)  _{n}, & -2\leq s\leq0,\\
2^{2n}\left(  \frac{s}{2}\right)  _{n}\left(  \frac{s+1}{2}\right)  _{n}, &
s\geq0,
\end{array}
\right\}  \left\vert x\right\vert ^{-s-2n},
\]

and from \ref{Ap057},%
\[
\left\vert \left(  \widehat{a}D\right)  ^{2n+1}\left\vert x\right\vert
^{-s}\right\vert \leq\left\{
\begin{array}
[c]{ll}%
2^{2n+2}\left\vert \left(  \frac{s}{2}\right)  _{n+1}\right\vert \left(
\frac{1}{2}\right)  _{n+1}, & -4\leq s\leq2,\\
2^{2n+1}\left(  \frac{s}{2}\right)  _{n+1}\left(  \frac{s+1}{2}\right)
_{n}, & s\geq2,
\end{array}
\right\}  \left\vert \widehat{a}\widehat{x}\right\vert \left\vert x\right\vert
^{-s-\left(  2n+1\right)  }.
\]

\subsection{Calculating $p_{s}^{\left(  2n\right)  }$ - for numerical use}

This subsection is needed to calculate the iterative formulas \ref{Ap036} for
the constants $b_{2k}^{2n}$ of the polynomials $p_{s}^{\left(  2n\right)  } $
and so obtain Conjecture \ref{Conj_argmax_p_even} above. These results can be
also used to prove \ref{Ap045} and \ref{Ap052} by induction.

From \ref{Ap033},%
\[%
\begin{array}
[c]{ll}%
b_{0}^{2n}=b_{1}^{2n-1}, & k=0,\\
b_{2k}^{2n}=-b_{2k-1}^{2n-1}+\left(  2k+1\right)  b_{2k+1}^{2n-1}, & 1\leq
k\leq n-1,\\
b_{2n}^{2n}=-b_{2n-1}^{2n-1}. & k=n.
\end{array}
\]

and from \ref{Ap028},%
\[%
\begin{array}
[c]{ll}%
b_{2k-1}^{2n-1}=-b_{2k-2}^{2n-2}+\left(  2k\right)  b_{2k}^{2n-2}, & 1\leq
k\leq n-2,\\
b_{2n-1}^{2n-1}=-b_{2n-2}^{2n-2}, & k=n-1.
\end{array}
\]

so we have%
\begin{align*}
b_{2k}^{2n}  & =-b_{2k-1}^{2n-1}+\left(  2k+1\right)  b_{2k+1}^{2n-1}\\
& =b_{2k-2}^{2n-2}-\left(  2k\right)  b_{2k}^{2n-2}+\left(  2k+1\right)
\left(  -b_{2k}^{2n-2}+\left(  2k+2\right)  b_{2k+2}^{2n-2}\right) \\
& =b_{2k-2}^{2n-2}-\left(  4k+1\right)  b_{2k}^{2n-2}+\left(  2k+1\right)
\left(  2k+2\right)  b_{2k+2}^{2n-2}.
\end{align*}

which for $n\geq2$ can be expressed as%
\begin{equation}%
\begin{array}
[c]{ll}%
b_{2k}^{2n}=b_{2k-2}^{2n-2}-\left(  4k+1\right)  b_{2k}^{2n-2}+\left(
2k+1\right)  \left(  2k+2\right)  b_{2k+2}^{2n-2}, & 1\leq k\leq n-1,\\
b_{0}^{2n}=-b_{0}^{2n-2}, & \\
b_{2n}^{2n}=b_{2n-2}^{2n-2}. &
\end{array}
\label{Ap093}%
\end{equation}

In particular%
\begin{align*}
b_{0}^{4}  & =b_{1}^{3},\\
b_{2}^{4}  & =-b_{1}^{3}+3b_{3}^{3},\\
b_{4}^{4}  & =-b_{3}^{3}.
\end{align*}

and%
\begin{align*}
b_{1}^{3}  & =-b_{0}^{2}+2b_{2}^{2},\\
b_{3}^{3}  & =-b_{2}^{2}.
\end{align*}

so%
\begin{align*}
b_{0}^{4}  & =-b_{0}^{2}+2b_{2}^{2},\\
b_{2}^{4}  & =b_{0}^{2}-2b_{2}^{2}-3b_{2}^{2}\\
& =b_{0}^{2}-5b_{2}^{2},\\
b_{4}^{4}  & =b_{2}^{2},
\end{align*}

which can be written%
\[
\left(  b_{2k}^{4}\right)  =%
\begin{pmatrix}
-1 & 2\\
1 & -5\\
0 & 1
\end{pmatrix}
\left(  b_{2k}^{2}\right)  .
\]

In general, for $n\geq1$: $b_{0}^{0}=1$ and%
\begin{align}
&  \left(  b_{2k}^{2n}\right)  =\nonumber\\
&
\begin{pmatrix}
-1 & 2 & 0 &  &  &  &  &  & \\
1 & -5 & 12 &  &  &  &  &  & \\
0 & 1 & -9 &  &  &  &  &  & \\
&  & \ddots & \ddots &  &  &  &  & \\
&  &  & 1 & -\left(  4k+1\right)  & \left(  2k+1\right)  \left(  2k+2\right)
&  &  & \\
&  &  &  &  & \ddots & \ddots &  & \\
&  &  &  &  &  & -\left(  4n-7\right)  & \left(  2n-4\right)  \left(
2n-3\right)  & 0\\
&  &  &  &  &  & 1 & -\left(  4n-3\right)  & \left(  2n-2\right)  \left(
2n-1\right) \\
&  &  &  &  &  & 0 & 1 & -\left(  4n+1\right)
\end{pmatrix}
\nonumber\\
&  \times%
\begin{pmatrix}
b_{2k}^{2n-2}\\
0
\end{pmatrix}
.\label{Ap036}%
\end{align}

Numerically we create the square matrix \ref{Ap036} and then remove the last column.

\subsection{Calculating $p_{s}^{\left(  2n+1\right)  }$ - for numerical use}

This subsection is needed to calculate iterative formulas \ref{Ap036} for the
constants $b_{2k+1}^{2n+1}$ of the polynomials $p_{s}^{\left(  2n+1\right)  }$
and so obtain Conjecture \ref{Conj_argmax_p_odd} above. These results can be
also used to prove the representations \ref{Ap045} and \ref{Ap052} by induction.

From \ref{Ap028},%
\[%
\begin{array}
[c]{ll}%
b_{2k+1}^{2n+1}=-b_{2k}^{2n}+\left(  2k+2\right)  b_{2k+2}^{2n}, &
k=0,1,2,\ldots,n-1.\\
b_{2n+1}^{2n+1}=-b_{2n}^{2n}. &
\end{array}
\]

and from \ref{Ap033},%
\begin{align*}
b_{0}^{2n}  & =b_{1}^{2n-1},\\
b_{2k}^{2n}  & =-b_{2k-1}^{2n-1}+\left(  2k+1\right)  b_{2k+1}^{2n-1},\\
b_{2n}^{2n}  & =-b_{2n-1}^{2n-1},
\end{align*}

so that%
\begin{align*}
b_{2k+1}^{2n+1}  & =-b_{2k}^{2n}+\left(  2k+2\right)  b_{2k+2}^{2n}\\
& =-\left(  -b_{2k-1}^{2n-1}+\left(  2k+1\right)  b_{2k+1}^{2n-1}\right)
+\left(  2k+2\right)  \left(  -b_{2k+1}^{2n-1}+\left(  2k+3\right)
b_{2k+3}^{2n-1}\right) \\
& =b_{2k-1}^{2n-1}-\left(  2k+1\right)  b_{2k+1}^{2n-1}-\left(  2k+2\right)
b_{2k+1}^{2n-1}+\left(  2k+2\right)  \left(  2k+3\right)  b_{2k+3}^{2n-1}\\
& =b_{2k-1}^{2n-1}-\left(  4k+3\right)  b_{2k+1}^{2n-1}+\left(  2k+2\right)
\left(  2k+3\right)  b_{2k+3}^{2n-1},\\
b_{1}^{2n+1}  & =-3b_{1}^{2n-1}+6b_{3}^{2n-1},\\
b_{2n-1}^{2n+1}  & =b_{2n-3}^{2n-1}-\left(  4n-1\right)  b_{2n-1}^{2n-1},\\
b_{2n+1}^{2n+1}  & =b_{2n-1}^{2n-1}.
\end{align*}

i.e. for $n\geq3$,%
\begin{align}
b_{1}^{2n+1}  & =-3b_{1}^{2n-1}+6b_{3}^{2n-1},\label{Ap035}\\
b_{3}^{2n+1}  & =b_{1}^{2n-1}-7b_{3}^{2n-1}+20b_{5}^{2n-1},\nonumber\\
b_{5}^{2n+1}  & =b_{3}^{2n-1}-11b_{5}^{2n-1}+42b_{7}^{2n-1},\nonumber\\
& \vdots\nonumber\\
b_{2k+1}^{2n+1}  & =b_{2k-1}^{2n-1}-\left(  4k+3\right)  b_{2k+1}%
^{2n-1}+\left(  2k+2\right)  \left(  2k+3\right)  b_{2k+3}^{2n-1},\quad
k=1,2,\ldots,n-2,\nonumber\\
& \vdots\nonumber\\
b_{2n-1}^{2n+1}  & =b_{2n-3}^{2n-1}-\left(  4n-1\right)  b_{2n-1}%
^{2n-1},\nonumber\\
b_{2n+1}^{2n+1}  & =b_{2n-1}^{2n-1},\nonumber
\end{align}

i.e.%
\begin{equation}%
\begin{array}
[c]{lll}%
b_{1}^{2n+1} & =-3b_{1}^{2n-1}+6b_{3}^{2n-1}, & k=1,\\
b_{2k+1}^{2n+1} & =b_{2k-1}^{2n-1}-\left(  4k+3\right)  b_{2k+1}%
^{2n-1}+\left(  2k+2\right)  \left(  2k+3\right)  b_{2k+3}^{2n-1}, &
k=1,2,\ldots,n-2,\\
b_{2n+1}^{2n+1} & =b_{2n-1}^{2n-1}, & k=2n+1,
\end{array}
\label{Ap098}%
\end{equation}

so in matrix form the iterations \ref{Ap035} become for $n\geq2$,%
\begin{align}
&
\begin{pmatrix}
b_{1}^{2n+1}\\
b_{3}^{2n+1}\\
b_{5}^{2n+1}\\
\vdots\\
b_{2k+1}^{2n+1}\\
\vdots\\
b_{2n-3}^{2n+1}\\
b_{2n-1}^{2n+1}\\
b_{2n+1}^{2n+1}%
\end{pmatrix}
\nonumber\\
&  =%
\begin{pmatrix}
-3 & 6 & 0 &  &  &  &  &  & \\
1 & -7 & 20 &  &  &  &  &  & \\
0 & 1 & -11 & \ddots &  &  &  &  & \\
&  & \ddots & \ddots &  &  &  &  & \\
&  &  & 1 & -\left(  4k+3\right)  & \left(  2k+2\right)  \left(  2k+3\right)
&  &  & \\
&  &  &  &  & \ddots & \ddots &  & \\
&  &  &  &  & \ddots & -\left(  4n-5\right)  & \left(  2n-2\right)  \left(
2n-1\right)  & 0\\
&  &  &  &  &  & 1 & -\left(  4n-1\right)  & 2n\left(  2n+1\right) \\
&  &  &  &  &  & 0 & 1 & -\left(  4n+3\right)
\end{pmatrix}
\nonumber\\
&  \qquad\times%
\begin{pmatrix}
b_{1}^{2n-1}\\
b_{3}^{2n-1}\\
b_{5}^{2n-1}\\
\vdots\\
b_{2k+1}^{2n-1}\\
\vdots\\
b_{2n-3}^{2n-1}\\
b_{2n-1}^{2n-1}\\
0
\end{pmatrix}
,\label{Ap037}\\
b_{1}^{1} &  =1.\nonumber
\end{align}

\textbf{Check}: from \ref{Ap040} we know that%
\begin{align*}
b_{1}^{1}  & =1\\
b_{odd}^{3}  & =\left(  3,-1\right)  ^{\prime},\\
b_{odd}^{5}  & =\left(  -15,10,-1\right)  ^{\prime},\\
b_{odd}^{7}  & =\left(  105,-105,21,-1\right)  ^{\prime}.
\end{align*}

and we check that%
\[
b_{odd}^{3}=%
\begin{pmatrix}
-3\\
1
\end{pmatrix}
b_{1}^{1},\quad b_{odd}^{5}=%
\begin{pmatrix}
-3 & 6\\
1 & -7\\
0 & 1
\end{pmatrix}
b_{odd}^{3},\quad b_{odd}^{7}=%
\begin{pmatrix}
-3 & 6 & 0\\
1 & -7 & 20\\
0 & 1 & -11\\
0 & 0 & 1
\end{pmatrix}
b_{odd}^{5}.
\]

Numerically, we create the square matrix \ref{Ap037} and then remove the last column.

\subsection{Calculating $b_{2k+1}^{2n+1}/b_{1}^{2n+1}$ and $b_{2k}^{2n}%
/b_{0}^{2n}$ \label{SbSect_calc_C^n_k}}

The values of $b_{2k+1}^{2n+1}/b_{1}^{2n+1}$ and $b_{2k}^{2n}/b_{0}^{2n}$ are
much, much smaller than $b_{2k+1}^{2n+1}$ and $b_{2k}^{2n}$. Accordingly
define%
\begin{equation}
g_{2k+1}^{2n+1}:=b_{2k+1}^{2n+1}/b_{1}^{2n+1},\quad g_{2k}^{2n}:=b_{2k}%
^{2n}/b_{0}^{2n}.\label{Ap094}%
\end{equation}

From \ref{Ap083},%
\begin{equation}
b_{0}^{2n}=b_{1}^{2n-1}=\left(  -1\right)  ^{n}2^{n}\left(  \frac{1}%
{2}\right)  _{n},\quad n\geq1,\label{Ap025}%
\end{equation}

so that%
\begin{align*}
b_{0}^{2n}  & =\left(  -1\right)  ^{n}2^{n}\left(  \frac{1}{2}\right)
_{n}=\left(  -2\right)  \left(  -1\right)  ^{n-1}2^{n-1}\left(  \frac{1}%
{2}\right)  _{n-1}\left(  \frac{1}{2}+n-1\right)  =\\
& =\left(  -2\right)  \left(  -1\right)  ^{n-1}2^{n-1}\left(  \frac{1}%
{2}\right)  _{n-1}\left(  n-\frac{1}{2}\right)  =-\left(  2n-1\right)  \left(
-1\right)  ^{n-1}2^{n-1}\left(  \frac{1}{2}\right)  _{n-1}=\\
& =-\left(  2n-1\right)  b_{0}^{2n-2},
\end{align*}

and%
\[
b_{1}^{2n+1}=\left(  -1\right)  ^{n+1}2^{n+1}\left(  \frac{1}{2}\right)
_{n+1}=\left(  -2\right)  \left(  \frac{1}{2}\right)  _{n}\left(  \frac{1}%
{2}+n\right)  =-\left(  2n+1\right)  b_{1}^{2n-1},
\]

and we have shown that%
\begin{equation}
b_{0}^{2n}=-\left(  2n-1\right)  b_{0}^{2n-2},\quad b_{1}^{2n+1}=-\left(
2n+1\right)  b_{1}^{2n-1}.\label{Ap095}%
\end{equation}

From \ref{Ap093},%
\[%
\begin{array}
[c]{ll}%
b_{2k}^{2n}=b_{2k-2}^{2n-2}-\left(  4k+1\right)  b_{2k}^{2n-2}+\left(
2k+1\right)  \left(  2k+2\right)  b_{2k+2}^{2n-2}, & 1\leq k\leq n-1,\\
b_{0}^{2n}=-b_{0}^{2n-2}, & \\
b_{2n}^{2n}=b_{2n-2}^{2n-2}, &
\end{array}
\]

so%
\begin{align*}
g_{2k}^{2n}=b_{2k}^{2n}/b_{0}^{2n}  & =b_{2k-2}^{2n-2}/b_{0}^{2n}-\left(
4k+1\right)  b_{2k}^{2n-2}/b_{0}^{2n}+\left(  2k+1\right)  \left(
2k+2\right)  b_{2k+2}^{2n-2}/b_{0}^{2n}\\
& =-\frac{1}{2n-1}\left(  b_{2k-2}^{2n-2}/b_{0}^{2n-2}-\left(  4k+1\right)
b_{2k}^{2n-2}/b_{0}^{2n-2}+\left(  2k+1\right)  \left(  2k+2\right)
b_{2k+2}^{2n-2}/b_{0}^{2n-2}\right) \\
& =-\frac{1}{2n-1}\left(  g_{2k-2}^{2n-2}-\left(  4k+1\right)  g_{2k}%
^{2n-2}+\left(  2k+1\right)  \left(  2k+2\right)  g_{2k+2}^{2n-2}\right)  ,
\end{align*}

and%
\begin{equation}%
\begin{array}
[c]{ll}%
g_{2k}^{2n}=\frac{1}{2n-1}g_{2k-2}^{2n-2}-\frac{4k+1}{2n-1}b_{2k}^{2n-2}%
+\frac{\left(  2k+1\right)  \left(  2k+2\right)  }{2n-1}g_{2k+2}^{2n-2}, &
1\leq k\leq n-1,\\
g_{0}^{2n}=\frac{1}{2n-1}g_{0}^{2n-2}, & \\
g_{2n}^{2n}=-\frac{1}{2n-1}g_{2n-2}^{2n-2}. &
\end{array}
\label{Ap096}%
\end{equation}

Hence we can use the same iterative matrix that was used to calculate the
$b_{2k}^{2n}$ in \ref{Ap036}, except that we now have the scalar multiplier
$-\frac{1}{2n-1}$ i.e.%
\begin{align}
&  \left(  g_{2k}^{2n}\right) \nonumber\\
&  =-\frac{1}{2n-1}\times\nonumber\\
&
\begin{pmatrix}
-1 & 2 & 0 &  &  &  &  &  & \\
1 & -5 & 12 &  &  &  &  &  & \\
0 & 1 & -9 &  &  &  &  &  & \\
&  & \ddots & \ddots &  &  &  &  & \\
&  &  & 1 & -\left(  4k+1\right)  & \left(  2k+1\right)  \left(  2k+2\right)
&  &  & \\
&  &  &  &  & \ddots & \ddots &  & \\
&  &  &  &  &  & -\left(  4n-7\right)  & \left(  2n-4\right)  \left(
2n-3\right)  & 0\\
&  &  &  &  &  & 1 & -\left(  4n-3\right)  & \left(  2n-2\right)  \left(
2n-1\right) \\
&  &  &  &  &  & 0 & 1 & -\left(  4n+1\right)
\end{pmatrix}
\nonumber\\
&  \times%
\begin{pmatrix}
g_{2k}^{2n-2}\\
0
\end{pmatrix}
,\label{Ap097}\\
g_{0}^{0} &  =-1.\nonumber
\end{align}

For numerical purposes we first multiply each element of the matrix by
$-\frac{1}{2n-1}$ and \textbf{then} multiply by $%
\begin{pmatrix}
g_{2k}^{2n-2}\\
0
\end{pmatrix}
$: this involves much smaller numbers.

Similarly, from \ref{Ap095} $b_{1}^{2n+1}=-\left(  2n+1\right)  b_{1}^{2n-1}$
so the iteration \ref{Ap098} and its matrix form \ref{Ap037} become,%
\begin{equation}%
\begin{array}
[c]{lll}%
g_{1}^{2n+1} & =-\frac{1}{2n+1}\left(  -3c_{1}^{2n-1}+6c_{3}^{2n-1}\right)
, & k=1,\\
g_{2k+1}^{2n+1} & =-\frac{1}{2n+1}\left(  g_{2k-1}^{2n-1}-\left(  4k+3\right)
g_{2k+1}^{2n-1}+\left(  2k+2\right)  \left(  2k+3\right)  g_{2k+3}%
^{2n-1}\right)  , & k=1,2,\ldots,n-2,\\
g_{2n+1}^{2n+1} & =-\frac{1}{2n+1}g_{2n-1}^{2n-1}, & k=2n+1,
\end{array}
\label{Ap099}%
\end{equation}

and%
\begin{align}
&  \left(  g_{2k+1}^{2n+1}\right)  =-\frac{1}{2n+1}\times\label{Ap090}\\
&
\begin{pmatrix}
-3 & 6 & 0 &  &  &  &  &  & \\
1 & -7 & 20 &  &  &  &  &  & \\
0 & 1 & -11 & \ddots &  &  &  &  & \\
&  & \ddots & \ddots &  &  &  &  & \\
&  &  & 1 & -\left(  4k+3\right)  & \left(  2k+2\right)  \left(  2k+3\right)
&  &  & \\
&  &  &  &  & \ddots & \ddots &  & \\
&  &  &  &  & \ddots & -\left(  4n-5\right)  & \left(  2n-2\right)  \left(
2n-1\right)  & 0\\
&  &  &  &  &  & 1 & -\left(  4n-1\right)  & \left(  2n\right)  \left(
2n+1\right) \\
&  &  &  &  &  & 0 & 1 & -\left(  4n+3\right)
\end{pmatrix}
\nonumber\\
&  \times%
\begin{pmatrix}
g_{2k-1}^{2n-1}\\
0
\end{pmatrix}
,\nonumber\\
g_{1}^{1} &  =-1.\nonumber
\end{align}

\section{Bounds for $\left(  \widehat{a}D\right)  ^{m}\left\vert x\right\vert
^{k},$ $m,k=0,1,2,\ldots$}

The conjectured bounds are stated in Theorem \ref{Thm_(aD)^m_|x|^n}. Several
conjectures concerning the polynomials $p_{s}^{\left(  k\right)  }$, inspired
by Matlab experiments, are required to obtain this theorem. We have two
conjectures about $p_{s}^{\left(  2n\right)  }$ when $s=-1,-2,-3,\ldots$:

\begin{conjecture}
\label{Conj_argmax_p^2m_neg2k}\textbf{When} $s=-2k$ \textbf{for}
$k=1,2,3,\ldots$ From Matlab experiments we conjecture that%
\[
\operatorname*{argmax}_{t\in\left[  -1,1\right]  }\left\vert p_{s}^{\left(
2n\right)  }\left(  t\right)  \right\vert =\left\{
\begin{array}
[c]{ll}%
\pm1, & 2\leq2n\leq-s,\\
0, & 2n>-s,
\end{array}
\right.
\]

and%
\[
\operatorname*{argmax}_{t\in\left[  -1,1\right]  }\left\vert p_{s}^{\left(
2n+1\right)  }\left(  t\right)  \right\vert =\left\{
\begin{array}
[c]{ll}%
\pm1, & 1\leq2n+1<-s-1.\\
\left[  -1,1\right]  , & 2n+1=-s-1,\\
0, & 2n+1>-s,
\end{array}
\right.
\]

\end{conjecture}

and:

\begin{conjecture}
\label{Conj_argmax_p^(2m+1)_neg2k}\textbf{When} $s=-\left(  2k+1\right)  $
\textbf{for} $k=1,2,3,\ldots$ From Matlab experiments we conjecture that%
\[
\operatorname*{argmax}_{t\in\left[  -1,1\right]  }\left\vert p_{s}^{\left(
2n+1\right)  }\left(  t\right)  \right\vert =\left\{
\begin{array}
[c]{ll}%
\pm1, & 3\leq2n+1<-s,\\
0, & 2n+1\geq-s,
\end{array}
\right.
\]

and%
\[
\operatorname*{argmax}_{t\in\left[  -1,1\right]  }\left\vert p_{s}^{\left(
2n\right)  }\left(  t\right)  \right\vert =\left\{
\begin{array}
[c]{ll}%
\pm1, & 3\leq2n<-s,\\
0, & 2n>-s,
\end{array}
\right.
\]

\end{conjecture}

We consider four cases:\medskip

\fbox{Bound for $\left(  \widehat{a}D\right)  ^{2n}\left\vert x\right\vert
^{2k}$} From \ref{Ap059},%
\[
\left(  \widehat{a}D\right)  ^{2n}\left\vert x\right\vert ^{-s}=2^{n}\left(
\frac{s}{2}\right)  _{n}p_{s}^{\left(  2n\right)  }\left(  \widehat{a}%
\widehat{x}\right)  \left\vert x\right\vert ^{-\left(  s+2n\right)  },
\]

so that%
\[
\left(  \widehat{a}D\right)  ^{2n}\left\vert x\right\vert ^{2k}=2^{n}\left(
-k\right)  _{n}p_{-2k}^{\left(  2n\right)  }\left(  \widehat{a}\widehat
{x}\right)  \left\vert x\right\vert ^{2k-2n},
\]

and from Conjecture \ref{Conj_argmax_p^2m_neg2k},%
\begin{align*}
\left\vert \left(  \widehat{a}D\right)  ^{2n}\left\vert x\right\vert
^{2k}\right\vert  & \leq\left\{
\begin{array}
[c]{ll}%
2^{n}\left\vert \left(  -k\right)  _{n}\right\vert p_{-2k}^{\left(  2n\right)
}\left(  1\right)  \left\vert x\right\vert ^{2k-2n}, & 2\leq2n\leq2k,\\
2^{n}\left\vert \left(  -k\right)  _{n}\right\vert p_{-2k}^{\left(  2n\right)
}\left(  0\right)  \left\vert x\right\vert ^{2k-2n}, & 2n>2k,
\end{array}
\right. \\
& =\left\{
\begin{array}
[c]{ll}%
2^{n}\left\vert \left(  -k\right)  _{n}\right\vert p_{-2k}^{\left(  2n\right)
}\left(  1\right)  \left\vert x\right\vert ^{2k-2n}, & 1\leq n\leq k,\\
2^{n}\left\vert \left(  -k\right)  _{n}\right\vert p_{-2k}^{\left(  2n\right)
}\left(  0\right)  \left\vert x\right\vert ^{2k-2n}, & n>k.
\end{array}
\right.
\end{align*}

But $\left(  \widehat{a}D\right)  ^{2n}\left\vert \cdot\right\vert ^{2k}=0$
when $n>k$ so%
\[
\left\vert \left(  \widehat{a}D\right)  ^{2n}\left\vert x\right\vert
^{2k}\right\vert \leq\left\{
\begin{array}
[c]{ll}%
2^{n}\left\vert \left(  -k\right)  _{n}\right\vert p_{-2k}^{\left(  2n\right)
}\left(  1\right)  \left\vert x\right\vert ^{2k-2n}, & 1\leq n\leq k,\\
0, & n>k.
\end{array}
\right.
\]

From \ref{Ap032},%
\[
p_{-2k}^{\left(  2n\right)  }\left(  1\right)  =\frac{\left(  -2k\right)
_{2n}}{2^{n}\left(  -k\right)  _{n}}=2^{n}\left(  \frac{1}{2}-k\right)
_{n},\quad n=0,1,2,\ldots,
\]

so%
\begin{equation}%
\begin{array}
[c]{lll}%
\left(  \widehat{a}D\right)  ^{2n}\left\vert x\right\vert ^{2k} & =0, & n>k,\\
\left\vert \left(  \widehat{a}D\right)  ^{2n}\left\vert x\right\vert
^{2k}\right\vert  & \leq\left\vert \left(  -2k\right)  _{2n}\right\vert
\left\vert x\right\vert ^{2k-2n}, & n\leq k.
\end{array}
\label{Ap068}%
\end{equation}

\fbox{Bound for $\left(  \widehat{a}D\right)  ^{2n+1}\left\vert x\right\vert
^{2k}$} From \ref{Ap022},%
\[
\left(  \widehat{a}D\right)  ^{2n+1}\left\vert x\right\vert ^{-s}%
=2^{n+1}\left(  \frac{s}{2}\right)  _{n+1}\left(  \widehat{a}\widehat
{x}\right)  p_{s}^{\left(  2n+1\right)  }\left(  \widehat{a}\widehat
{x}\right)  \left\vert x\right\vert ^{-s-\left(  2n+1\right)  },\quad n\geq0,
\]

so that when $s=-2k$,%
\[
\left(  \widehat{a}D\right)  ^{2n+1}\left\vert x\right\vert ^{2k}%
=2^{n+1}\left(  -k\right)  _{n+1}\left(  \widehat{a}\widehat{x}\right)
p_{-2k}^{\left(  2n+1\right)  }\left(  \widehat{a}\widehat{x}\right)
\left\vert x\right\vert ^{2k-\left(  2n+1\right)  },\quad n,k\geq0.
\]

From Conjecture \ref{Conj_argmax_p^2m_neg2k},%
\begin{align*}
&  \left\vert \left(  \widehat{a}D\right)  ^{2n+1}\left\vert x\right\vert
^{2k}\right\vert \\
&  \leq2^{n+1}\left\vert \left(  -k\right)  _{n+1}\right\vert \left(
\widehat{a}\widehat{x}\right)  \left\vert x\right\vert ^{2k-\left(
2n+1\right)  }\left\{
\begin{array}
[c]{ll}%
p_{-2k}^{\left(  2n+1\right)  }\left(  0\right)  , & 2n+1>2k,\\
p_{-2k}^{\left(  2n+1\right)  }\left(  t\right)  ,t\in\left[  -1,1\right]  , &
2n+1=2k-1,\\
p_{-2k}^{\left(  2n+1\right)  }\left(  1\right)  , & 1\leq2n+1<2k-1,
\end{array}
\right. \\
&  =2^{n+1}\left\vert \left(  -k\right)  _{n+1}\right\vert \left(  \widehat
{a}\widehat{x}\right)  \left\vert x\right\vert ^{2k-\left(  2n+1\right)
}\left\{
\begin{array}
[c]{ll}%
p_{-2k}^{\left(  2n+1\right)  }\left(  0\right)  , & n\geq k,\\
p_{-2k}^{\left(  2n+1\right)  }\left(  t\right)  ,t\in\left[  -1,1\right]  , &
n=k-1,\\
p_{-2k}^{\left(  2n+1\right)  }\left(  1\right)  , & 0\leq n\leq k-2,
\end{array}
\right. \\
&  =2^{n+1}\left\vert \left(  -k\right)  _{n+1}\right\vert \left(  \widehat
{a}\widehat{x}\right)  \left\vert x\right\vert ^{2k-\left(  2n+1\right)
}\left\{
\begin{array}
[c]{ll}%
p_{-2k}^{\left(  2n+1\right)  }\left(  0\right)  , & n\geq k,\\
p_{-2k}^{\left(  2n+1\right)  }\left(  1\right)  , & 0\leq n\leq k-1,
\end{array}
\right. \\
&  =2^{n+1}\left\vert \left(  -k\right)  _{n+1}\right\vert \left(  \widehat
{a}\widehat{x}\right)  \left\vert x\right\vert ^{2k-\left(  2n+1\right)
}\left\{
\begin{array}
[c]{ll}%
0, & n\geq k,\\
p_{-2k}^{\left(  2n+1\right)  }\left(  1\right)  , & 0\leq n\leq k-1.
\end{array}
\right.
\end{align*}

From \ref{Ap026}, $p_{s}^{\left(  2n+1\right)  }\left(  1\right)
=-\frac{\left(  s\right)  _{2n+1}}{2^{n+1}\left(  \frac{s}{2}\right)  _{n+1}}%
$, so $p_{-2k}^{\left(  2n+1\right)  }\left(  1\right)  =-\frac{\left(
-2k\right)  _{2n+1}}{2^{n+1}\left(  -k\right)  _{n+1}}$ and%
\begin{equation}%
\begin{array}
[c]{lll}%
\left\vert \left(  \widehat{a}D\right)  ^{2n+1}\left\vert x\right\vert
^{2k}\right\vert  & \leq\left\vert \left(  -2k\right)  _{2n+1}\right\vert
\left(  \left\vert \left(  \widehat{a}\widehat{x}\right)  \right\vert
\left\vert x\right\vert ^{2k-\left(  2n+1\right)  }\right)  , & n\leq k,\\
\left(  \widehat{a}D\right)  ^{2n+1}\left\vert x\right\vert ^{2k} & =0, & n>k.
\end{array}
\label{Ap069}%
\end{equation}
\medskip

\fbox{Bound for $\left(  \widehat{a}D\right)  ^{2n}\left\vert x\right\vert
^{2k+1}$} From \ref{Ap059},%
\[
\left(  \widehat{a}D\right)  ^{2n}\left\vert x\right\vert ^{-s}=2^{n}\left(
\frac{s}{2}\right)  _{n}p_{s}^{\left(  2n\right)  }\left(  \widehat{a}%
\widehat{x}\right)  \left\vert x\right\vert ^{-\left(  s+2n\right)  },
\]

so that when $s=-\left(  2k+1\right)  $,%
\[
\left(  \widehat{a}D\right)  ^{2n}\left\vert x\right\vert ^{2k+1}=2^{n}\left(
-k-\frac{1}{2}\right)  _{n}p_{-\left(  2k+1\right)  }^{\left(  2n\right)
}\left(  \widehat{a}\widehat{x}\right)  \left\vert x\right\vert ^{2k+1-2n},
\]

From Conjecture \ref{Conj_argmax_p^(2m+1)_neg2k},%
\begin{align*}
\left\vert \left(  \widehat{a}D\right)  ^{2n}\left\vert x\right\vert
^{2k+1}\right\vert  & \leq2^{n}\left\vert \left(  -k-\frac{1}{2}\right)
_{n}\right\vert \left\vert x\right\vert ^{2k+1-2n}\times\left\{
\begin{array}
[c]{ll}%
\left\vert p_{-\left(  2k+1\right)  }^{\left(  2n\right)  }\left(  1\right)
\right\vert , & 3\leq2n<2k+1,\\
\left\vert p_{-\left(  2k+1\right)  }^{\left(  2n\right)  }\left(  0\right)
\right\vert , & 2n>2k+1,
\end{array}
\right. \\
& =2^{n}\left\vert \left(  -k-\frac{1}{2}\right)  _{n}\right\vert \left\vert
x\right\vert ^{2k+1-2n}\times\left\{
\begin{array}
[c]{ll}%
\left\vert p_{-\left(  2k+1\right)  }^{\left(  2n\right)  }\left(  1\right)
\right\vert , & 3/2\leq n<k+1/2,\\
\left\vert p_{-\left(  2k+1\right)  }^{\left(  2n\right)  }\left(  0\right)
\right\vert , & n>k+1/2,
\end{array}
\right. \\
& =2^{n}\left\vert \left(  -k-\frac{1}{2}\right)  _{n}\right\vert \left\vert
x\right\vert ^{2k+1-2n}\times\left\{
\begin{array}
[c]{ll}%
\left\vert p_{-\left(  2k+1\right)  }^{\left(  2n\right)  }\left(  1\right)
\right\vert , & 2\leq n\leq k,\\
\left\vert p_{-\left(  2k+1\right)  }^{\left(  2n\right)  }\left(  0\right)
\right\vert , & n\geq k+1.
\end{array}
\right.
\end{align*}

From \ref{Ap032},%
\[
p_{s}^{\left(  2n\right)  }\left(  1\right)  =2^{n}\left(  \frac{s+1}%
{2}\right)  _{n},\quad n=0,1,2,\ldots,
\]

so that%
\[
p_{-\left(  2k+1\right)  }^{\left(  2n\right)  }\left(  1\right)
=2^{n}\left(  -k\right)  _{n},\quad n=0,1,2,\ldots,
\]

and from \ref{Ap031},%
\[
p_{s}^{\left(  2n\right)  }\left(  0\right)  =\left(  -1\right)  ^{n}%
2^{n}\left(  \frac{1}{2}\right)  _{n},\quad n=1,2,3,\ldots;\text{ }%
s\in\mathbb{R}^{1},
\]

so that%
\[
p_{-\left(  2k+1\right)  }^{\left(  2n\right)  }\left(  0\right)  =\left(
-1\right)  ^{n}2^{n}\left(  \frac{1}{2}\right)  _{n}.
\]

Thus%
\begin{align}
\left\vert \left(  \widehat{a}D\right)  ^{2n}\left\vert x\right\vert
^{2k+1}\right\vert  & \leq2^{n}\left\vert \left(  -k-\frac{1}{2}\right)
_{n}\right\vert \left\{
\begin{array}
[c]{ll}%
\left\vert 2^{n}\left(  -k\right)  _{n}\right\vert , & 2\leq n\leq k,\\
2^{n}\left(  \frac{1}{2}\right)  _{n}, & n\geq k+1,
\end{array}
\right\}  \left\vert x\right\vert ^{2k+1-2n}\nonumber\\
& =\left\{
\begin{array}
[c]{ll}%
2^{2n}\left\vert \left(  -k-\frac{1}{2}\right)  _{n}\right\vert \left\vert
\left(  -k\right)  _{n}\right\vert , & 2\leq n\leq k,\\
2^{2n}\left\vert \left(  -k-\frac{1}{2}\right)  _{n}\right\vert \left(
\frac{1}{2}\right)  _{n}, & n\geq k+1.
\end{array}
\right\}  \left\vert x\right\vert ^{2k+1-2n}\label{Ap071}%
\end{align}
\medskip

\fbox{Bound for $\left(  \widehat{a}D\right)  ^{2n+1}\left\vert x\right\vert
^{2k+1}$} From \ref{Ap022},%
\[
\left(  \widehat{a}D\right)  ^{2n+1}\left\vert x\right\vert ^{-s}%
=2^{n+1}\left(  \frac{s}{2}\right)  _{n+1}\left(  \widehat{a}\widehat
{x}\right)  p_{s}^{\left(  2n+1\right)  }\left(  \widehat{a}\widehat
{x}\right)  \left\vert x\right\vert ^{-s-\left(  2n+1\right)  },\quad n\geq0,
\]

so that when $s=-\left(  2k+1\right)  $,%
\[
\left(  \widehat{a}D\right)  ^{2n+1}\left\vert x\right\vert ^{2k+1}%
=2^{n+1}\left(  -k-\frac{1}{2}\right)  _{n+1}\left(  \widehat{a}\widehat
{x}\right)  p_{-\left(  2k+1\right)  }^{\left(  2n+1\right)  }\left(
\widehat{a}\widehat{x}\right)  \left\vert x\right\vert ^{2k-2n},\quad
n,k\geq0.
\]
From Conjecture \ref{Conj_argmax_p^(2m+1)_neg2k},%
\begin{align*}
\left\vert \left(  \widehat{a}D\right)  ^{2n+1}\left\vert x\right\vert
^{2k+1}\right\vert  & \leq2^{n+1}\left\vert \left(  -k-\frac{1}{2}\right)
_{n+1}\right\vert \left\{
\begin{array}
[c]{ll}%
p_{-\left(  2k+1\right)  }^{\left(  2n+1\right)  }\left(  1\right)  , &
3\leq2n+1<2k+1,\\
p_{-\left(  2k+1\right)  }^{\left(  2n+1\right)  }\left(  0\right)  , &
2n+1\geq2k+1,
\end{array}
\right\}  \left\vert \widehat{a}\widehat{x}\right\vert \left\vert x\right\vert
^{2k-2n}\\
& \leq2^{n+1}\left\vert \left(  -k-\frac{1}{2}\right)  _{n+1}\right\vert
\left\{
\begin{array}
[c]{ll}%
p_{-\left(  2k+1\right)  }^{\left(  2n+1\right)  }\left(  1\right)  , & 1\leq
n<k,\\
p_{-\left(  2k+1\right)  }^{\left(  2n+1\right)  }\left(  0\right)  , & n\geq
k.
\end{array}
\right\}  \left\vert \widehat{a}\widehat{x}\right\vert \left\vert x\right\vert
^{2k-2n}%
\end{align*}
But from \ref{Ap031},
\[
p_{s}^{\left(  2n+1\right)  }\left(  0\right)  =\left(  -1\right)
^{n+1}2^{n+1}\left(  \frac{1}{2}\right)  _{n+1},\quad n=0,1,2,\ldots;\text{
}s\in\mathbb{R}^{1},
\]
and from \ref{Ap026}, $p_{s}^{\left(  2n+1\right)  }\left(  1\right)
=2^{n}\left(  \frac{s+1}{2}\right)  _{n}$ so%
\[
p_{-\left(  2k+1\right)  }^{\left(  2n+1\right)  }\left(  1\right)
=2^{n}\left(  -k\right)  _{n},
\]

and these two equations imply%
\begin{align}
\left\vert \left(  \widehat{a}D\right)  ^{2n+1}\left\vert x\right\vert
^{2k+1}\right\vert  & \leq2^{n+1}\left\vert \left(  -k-\frac{1}{2}\right)
_{n+1}\right\vert \left\{
\begin{array}
[c]{ll}%
2^{n}\left\vert \left(  -k\right)  _{n}\right\vert , & 1\leq n<k,\\
2^{n+1}\left(  \frac{1}{2}\right)  _{n+1}, & n\geq k.
\end{array}
\right\}  \left\vert \widehat{a}\widehat{x}\right\vert \left\vert x\right\vert
^{2k-2n}\nonumber\\
& =\left\{
\begin{array}
[c]{ll}%
2^{2n+1}\left\vert \left(  -k-\frac{1}{2}\right)  _{n+1}\right\vert \left\vert
\left(  -k\right)  _{n}\right\vert , & 1\leq n<k,\\
2^{2n+2}\left\vert \left(  -k-\frac{1}{2}\right)  _{n+1}\right\vert \left(
\frac{1}{2}\right)  _{n+1}, & n\geq k.
\end{array}
\right\}  \left\vert \widehat{a}\widehat{x}\right\vert \left\vert x\right\vert
^{2k-2n}\label{Ap070}%
\end{align}

In summary, if the Conjectures \ref{Conj_argmax_p^(2m+1)_neg2k} and
\ref{Conj_argmax_p^(2m+1)_neg2k} are true then from \ref{Ap068}, \ref{Ap069},
\ref{Ap070} and \ref{Ap071} we have:

\begin{theorem}
\label{Thm_(aD)^m_|x|^n}\textbf{CONJECTURES for upper bounds for }$\left(
\widehat{a}D\right)  ^{m}\left\vert x\right\vert ^{j}$ \textbf{when}
$m,j=0,1,2,\ldots$

Assume $k,n,m\geq0$ are integers and suppose $a,x\in\mathbb{R}^{d}$ and
$\widehat{a}=a/\left\vert a\right\vert $.

\textbf{For even powers of }$\left\vert x\right\vert $:%
\begin{equation}%
\begin{array}
[c]{lll}%
\left\vert \left(  \widehat{a}D\right)  ^{2n}\left\vert x\right\vert
^{2k}\right\vert  & \leq\left\vert \left(  -2k\right)  _{2n}\right\vert
\left\vert x\right\vert ^{2k-2n}, & n\leq k,\\
\left\vert \left(  \widehat{a}D\right)  ^{2n+1}\left\vert x\right\vert
^{2k}\right\vert  & \leq\left\vert \left(  -2k\right)  _{2n+1}\right\vert
\left\vert \widehat{a}\widehat{x}\right\vert \left\vert x\right\vert
^{2k-\left(  2n+1\right)  }, & n\leq k,
\end{array}
\label{Ap226}%
\end{equation}

and%
\begin{equation}
\left(  \widehat{a}D\right)  ^{m}\left\vert x\right\vert ^{2k}=0,\quad
m>2k.\label{Ap227}%
\end{equation}

\textbf{For odd powers of }$\left\vert x\right\vert $:%
\begin{equation}%
\begin{array}
[c]{ll}%
\left\vert \left(  \widehat{a}D\right)  ^{2n}\left\vert x\right\vert
^{2k+1}\right\vert  & \leq\left\{
\begin{array}
[c]{ll}%
2^{2n}\left\vert \left(  -k-\frac{1}{2}\right)  _{n}\right\vert \left\vert
\left(  -k\right)  _{n}\right\vert \left\vert x\right\vert ^{2k+1-2n}, & 0\leq
n\leq k,\\
2^{2n}\left\vert \left(  -k-\frac{1}{2}\right)  _{n}\right\vert \left(
\frac{1}{2}\right)  _{n}\left\vert x\right\vert ^{2k+1-2n}, & n\geq k+1,
\end{array}
\right. \\
\left\vert \left(  \widehat{a}D\right)  ^{2n+1}\left\vert x\right\vert
^{2k+1}\right\vert  & \leq\left\{
\begin{array}
[c]{ll}%
\left(  2^{2n+1}\left\vert \left(  -k-\frac{1}{2}\right)  _{n+1}\right\vert
\left\vert \left(  -k\right)  _{n}\right\vert \right)  \left\vert \widehat
{a}\widehat{x}\right\vert \left\vert x\right\vert ^{2k-2n}, & 0\leq n<k,\\
\left(  2^{2n+2}\left\vert \left(  -k-\frac{1}{2}\right)  _{n+1}\right\vert
\left(  \frac{1}{2}\right)  _{n+1}\right)  \left\vert \widehat{a}\widehat
{x}\right\vert \left\vert x\right\vert ^{2k-2n}, & n\geq k.
\end{array}
\right.
\end{array}
\label{Ap228}%
\end{equation}

\end{theorem}

By applying the approximation $\left\vert \widehat{a}\widehat{x}\right\vert
\leq1$,\ these inequalities are weakened to:

\begin{corollary}
\label{Cor_Thm_(aD)^m_|x|^n}\textbf{CONJECTURES for upper bounds for }$\left(
\widehat{a}D\right)  ^{m}\left\vert x\right\vert ^{j}$ \textbf{when}
$m,j=0,1,2,\ldots$.

Assume $k,n,m\geq0$ are integers and suppose $a,x\in\mathbb{R}^{d}$ and
$\widehat{a}=a/\left\vert a\right\vert $. Then%
\[
\left\vert \left(  \widehat{a}D\right)  ^{m}\left\vert x\right\vert
^{n}\right\vert \leq k_{m,n}\left\vert x\right\vert ^{n-m},
\]

where \textbf{for even powers of }$\left\vert x\right\vert $:%
\begin{equation}
k_{m,2k}=\left\{
\begin{array}
[c]{ll}%
\left\vert \left(  -2k\right)  _{m}\right\vert \left\vert x\right\vert
^{2k-m}, & m\leq2k,\\
0, & m>2k,
\end{array}
\right. \label{p14}%
\end{equation}

and \textbf{for odd powers of }$\left\vert x\right\vert $:%
\begin{equation}%
\begin{array}
[c]{ll}%
k_{2n,2k+1} & =2^{2n}\left\vert \left(  -k-\frac{1}{2}\right)  _{n}\right\vert
\times\left\{
\begin{array}
[c]{ll}%
\left\vert \left(  -k\right)  _{n}\right\vert , & 0\leq n\leq k,\\
\left(  \frac{1}{2}\right)  _{n}, & n\geq k+1,
\end{array}
\right. \\
k_{2n+1,2k+1} & =2^{2n+1}\left\vert \left(  -k-\frac{1}{2}\right)
_{n+1}\right\vert \times\left\{
\begin{array}
[c]{ll}%
\left\vert \left(  -k\right)  _{n}\right\vert , & 0\leq n<k,\\
2\left(  \frac{1}{2}\right)  _{n+1}, & n\geq k.
\end{array}
\right.
\end{array}
\label{p12}%
\end{equation}

\end{corollary}

\section{Estimating $D_{k}^{m}\left\vert x\right\vert ^{-s}$}

We have the following simple consequences of \ref{Ap057} and \ref{Ap055}. From
\ref{Ap057},%
\begin{align*}
\left\vert \left(  \widehat{a}D\right)  ^{2n+1}\left\vert x\right\vert
^{-s}\right\vert  & =\left\{
\begin{array}
[c]{ll}%
2^{2n+2}\left\vert \left(  \frac{s}{2}\right)  _{n+1}\right\vert \left(
\frac{1}{2}\right)  _{n+1}, & -4\leq s\leq2,\\
2^{2n+1}\left(  \frac{s}{2}\right)  _{n+1}\left(  \frac{s+1}{2}\right)
_{n}, & s\geq2,
\end{array}
\right\}  \left\vert \widehat{a}\widehat{x}\right\vert \left\vert x\right\vert
^{-s-\left(  2n+1\right)  }\\
for\text{ }n  & \geq0,
\end{align*}

so letting $\widehat{a}=\mathbf{e}_{k}$ gives%
\begin{align}
\left\vert D_{k}^{2n+1}\left\vert x\right\vert ^{-s}\right\vert  & =\left\{
\begin{array}
[c]{ll}%
2^{2n+2}\left\vert \left(  \frac{s}{2}\right)  _{n+1}\right\vert \left(
\frac{1}{2}\right)  _{n+1}, & -4\leq s\leq2,\\
2^{2n+1}\left(  \frac{s}{2}\right)  _{n+1}\left(  \frac{s+1}{2}\right)
_{n}, & s\geq2,
\end{array}
\right\}  \left\vert \widehat{x}_{k}\right\vert \left\vert x\right\vert
^{-s-\left(  2n+1\right)  }\label{Ap063}\\
for\text{ }n  & \geq0.\nonumber
\end{align}

From \ref{Ap055},%
\begin{align*}
\left\vert \left(  \widehat{a}D\right)  ^{2n}\left\vert x\right\vert
^{-s}\right\vert  & \leq\left\{
\begin{array}
[c]{ll}%
2^{2n}\left\vert \left(  \frac{s}{2}\right)  _{n}\right\vert \left(  \frac
{1}{2}\right)  _{n}\left\vert x\right\vert ^{-s-2n}, & -2\leq s\leq0,\\
2^{2n}\left(  \frac{s}{2}\right)  _{n}\left(  \frac{s+1}{2}\right)
_{n}\left\vert x\right\vert ^{-s-2n}, & s\geq0,
\end{array}
\right. \\
for\text{ }n  & \geq1.
\end{align*}

so letting $\widehat{a}=\mathbf{e}_{k}$ gives%
\begin{align}
\left\vert D_{k}^{2n}\left\vert x\right\vert ^{-s}\right\vert  & \leq\left\{
\begin{array}
[c]{ll}%
2^{2n}\left\vert \left(  \frac{s}{2}\right)  _{n}\right\vert \left(  \frac
{1}{2}\right)  _{n}, & -2\leq s\leq0,\\
2^{2n}\left(  \frac{s}{2}\right)  _{n}\left(  \frac{s+1}{2}\right)  _{n}, &
s\geq0,
\end{array}
\right\}  \left\vert x\right\vert ^{-s-2n}\label{Ap064}\\
for\text{ }n  & \geq1.\nonumber
\end{align}

\section{Upper bounds for $D^{\alpha}\left\vert x\right\vert ^{-s}$}

We want to obtain meaningful upper bounds for $D^{\alpha}\left\vert
x\right\vert ^{-s}$. When $\max\alpha\leq1$ we prove the formula \ref{Ap021}
for $D^{\alpha}\left\vert x\right\vert $. When $\left\vert \alpha\right\vert
+s>0$ we conjecture the bounds \ref{Ap115} and \ref{Ap113} for $D^{\alpha
}\left\vert x\right\vert ^{-s}$ when $\left\vert \alpha\right\vert +s>0$; this
is achieved using the integral formula \ref{Ap020} for the gamma function.

\subsection{A bound for $D^{\alpha}\left\vert x\right\vert $ when $\max
\alpha\leq1$.}

Permute the multi-index $\alpha$, if necessary, to write:%
\[
\alpha=\beta_{1,d^{\prime}}\mathbf{0}_{d^{\prime}+1,d},
\]

where $\beta_{1,d^{\prime}}=1$, $\beta_{d^{\prime}+1,d}=0$.\medskip

\fbox{Concerning $D^{\beta}\left\vert x\right\vert $}%
\begin{align*}
D_{1}\left\vert x\right\vert  & =x_{1}\left\vert x\right\vert ^{-1},\\
D_{2}D_{1}\left\vert x\right\vert  & =x_{1}D_{2}\left\vert x\right\vert
^{-1}=\left(  1\right)  \left(  -1\right)  x_{1}x_{2}\left\vert x\right\vert
^{-3},\\
D_{3}D_{2}D_{1}\left\vert x\right\vert  & =-x_{1}x_{2}D_{3}\left\vert
x\right\vert ^{-3}=\left(  1\right)  \left(  -1\right)  \left(  -3\right)
x_{1}x_{2}x_{3}\left\vert x\right\vert ^{-5},\\
D_{4}D_{3}D_{2}D_{1}\left\vert x\right\vert  & =\left(  1\right)  \left(
1-2\right)  \left(  1-4\right)  \left(  1-6\right)  x_{1}x_{2}x_{3}%
x_{4}\left\vert x\right\vert ^{-7},\\
& \vdots\\
D_{k}\ldots D_{3}D_{2}D_{1}\left\vert x\right\vert  & =\left(  1\right)
\left(  1-2\right)  \left(  1-4\right)  \left(  1-6\right)  \ldots\left(
1-\left(  2k-2\right)  \right)  x_{1}x_{2}\ldots x_{k}\left\vert x\right\vert
^{-\left(  2k-1\right)  }\\
& =\left(  1\right)  \left(  -1\right)  \left(  -3\right)  \left(  -5\right)
\ldots\left(  3-2k\right)  x_{1}x_{2}\ldots x_{k}\left\vert x\right\vert
^{-\left(  2k-1\right)  }\\
& =\left(  -1\right)  ^{k-1}\left(  1\right)  \left(  1\right)  \left(
3\right)  \left(  5\right)  \ldots\left(  2k-3\right)  x_{1}x_{2}\ldots
x_{k}\left\vert x\right\vert ^{-\left(  2k-1\right)  }\\
& =\left(  -1\right)  ^{k-1}\frac{1.1.2.3.4.5.\ldots\left(  2k-4\right)
\left(  2k-3\right)  \left(  2k-2\right)  }{2.4.6.\ldots\left(  2k-4\right)
\left(  2k-2\right)  }x_{1}x_{2}\ldots x_{k}\left\vert x\right\vert ^{-\left(
2k-1\right)  }\\
& =\left(  -1\right)  ^{k-1}\frac{\left(  2k-2\right)  !}{2^{k-1}\left(
k-1\right)  !}x_{1}x_{2}\ldots x_{k}\left\vert x\right\vert ^{-\left(
2k-1\right)  },
\end{align*}

if we define $0!:=1$. Now a permutation argument implies%
\begin{equation}
D^{\alpha}\left\vert x\right\vert =\left(  -1\right)  ^{\left\vert
\alpha\right\vert -1}\frac{\left(  2\left\vert \alpha\right\vert -2\right)
!}{2^{\left\vert \alpha\right\vert -1}\left(  \left\vert \alpha\right\vert
-1\right)  !}\widehat{x}^{\mathbf{\alpha}}\left\vert x\right\vert
^{1-\left\vert \alpha\right\vert },\quad0\leq\alpha\leq\mathbf{1},\alpha
\neq\mathbf{0,}\label{Ap021}%
\end{equation}

and so%
\begin{equation}
\left\vert D^{\alpha}\left\vert x\right\vert \right\vert \leq\frac{\left(
2\left\vert \alpha\right\vert -2\right)  !}{2^{\left\vert \alpha\right\vert
-1}\left(  \left\vert \alpha\right\vert -1\right)  !}\left\vert x\right\vert
^{1-\left\vert \alpha\right\vert },\quad0\leq\alpha\leq\mathbf{1,}\text{
}\alpha\neq\mathbf{0.}\label{a023}%
\end{equation}

\subsection{Conjectured bounds for $D^{\alpha}\left\vert x\right\vert ^{-s}$
when $\max\alpha_{i}\geq2$ and $\left\vert \alpha\right\vert +s>0$.}

Using a permutation argument we can assume without loss of generality that ??
NEEDED? ??
\begin{equation}
\alpha_{1}\geq2,\alpha_{i}\geq1\text{ }when\text{ }i\leq d^{\prime}\text{
}and\text{ }\alpha_{i}=0\text{ }when\text{ }i>d^{\prime}.\label{a024}%
\end{equation}

From \ref{Ap045}, for $n=0,1,2,\ldots$,
\begin{align*}
\left(  \widehat{a}D\right)  ^{2n+1}\left\vert x\right\vert ^{-s}  & =\left(
\sum_{k=0}^{n}b_{2k+1}^{2n+1}\left\{  s\left(  s+2\right)  \ldots\left(
s+2n+2k\right)  \right\}  \left(  \widehat{a}\widehat{x}\right)
^{2k+1}\right)  \left\vert x\right\vert ^{-s-\left(  2n+1\right)  }\\
& =\left(  \sum_{k=0}^{n}b_{2k+1}^{2n+1}2^{n+k+1}\left\{  \frac{s}{2}\left(
\frac{s}{2}+1\right)  \ldots\left(  \frac{s}{2}+n+k\right)  \right\}  \left(
\widehat{a}\widehat{x}\right)  ^{2k+1}\right)  \left\vert x\right\vert
^{-s-\left(  2n+1\right)  }\\
& =\left(  \sum_{k=0}^{n}b_{2k+1}^{2n+1}2^{n+k+1}\left(  \frac{s}{2}\right)
_{n+k+1}\left(  \widehat{a}\widehat{x}\right)  ^{2k+1}\right)  \left\vert
x\right\vert ^{-s-\left(  2n+1\right)  }\\
& =\sum_{k=0}^{n}b_{2k+1}^{2n+1}2^{n+k+1}\left(  \frac{s}{2}\right)
_{n+k+1}\left(  \widehat{a}\frac{x}{\left\vert x\right\vert }\right)
^{2k+1}\left\vert x\right\vert ^{-s-\left(  2n+1\right)  }\\
& =\sum_{k=0}^{n}b_{2k+1}^{2n+1}2^{n+k+1}\left(  \frac{s}{2}\right)
_{n+k+1}\left(  \widehat{a}x\right)  ^{2k+1}\left\vert x\right\vert
^{-2k-1}\left\vert x\right\vert ^{-s-\left(  2n+1\right)  }\\
& =\sum_{k=0}^{n}b_{2k+1}^{2n+1}2^{n+\frac{1}{2}+k+\frac{1}{2}}\left(
\frac{s}{2}\right)  _{n+\frac{1}{2}+k+\frac{1}{2}}\left(  \widehat{a}x\right)
^{2k+1}\left\vert x\right\vert ^{-\left(  s+2n+1+2k+1\right)  },
\end{align*}

which we can now write as%
\begin{equation}
\left(  \widehat{a}D\right)  ^{m}\left\vert x\right\vert ^{-s}=\sum_{j=0}%
^{m}c_{j}^{m}\left(  s\right)  \left(  \widehat{a}x\right)  ^{j}\left\vert
x\right\vert ^{-\left(  s+m+j\right)  },\quad m=1,3,5,\ldots,\label{Ap073}%
\end{equation}

where%
\begin{equation}
c_{j}^{m}\left(  s\right)  =\left\{
\begin{array}
[c]{ll}%
b_{j}^{m}2^{\frac{m+j}{2}}\left(  \frac{s}{2}\right)  _{\frac{m+j}{2}}, &
j\text{ }odd,\\
0, & j\text{ }even.
\end{array}
\right. \label{Ap074}%
\end{equation}

From \ref{Ap052}, for $n=1,2,3,\ldots$,
\begin{align*}
\left(  \widehat{a}D\right)  ^{2n}\left\vert x\right\vert ^{-s}  & =\left(
\sum_{k=0}^{n}b_{2k}^{2n}s\left(  s+2\right)  \ldots\left(  s+2n+2k-2\right)
\left(  \widehat{a}\widehat{x}\right)  ^{2k}\right)  \left\vert x\right\vert
^{-s-2n}\\
& =\left(  \sum_{k=0}^{n}b_{2k}^{2n}2^{n+k}\left\{  \frac{s}{2}\left(
\frac{s}{2}+1\right)  \ldots\left(  \frac{s}{2}+n+k-1\right)  \right\}
\left(  \widehat{a}\widehat{x}\right)  ^{2k}\right)  \left\vert x\right\vert
^{-s-2n}\\
& =\left(  \sum_{k=0}^{n}b_{2k}^{2n}2^{n+k}\left(  \frac{s}{2}\right)
_{n+k}\left(  \widehat{a}\widehat{x}\right)  ^{2k}\right)  \left\vert
x\right\vert ^{-s-2n}\\
& =\sum_{k=0}^{n}b_{2k}^{2n}2^{n+k}\left(  \frac{s}{2}\right)  _{n+k}\left(
\widehat{a}\frac{x}{\left\vert x\right\vert }\right)  ^{2k}\left\vert
x\right\vert ^{-s-2n}\\
& =\sum_{k=0}^{n}b_{2k}^{2n}2^{n+k}\left(  \frac{s}{2}\right)  _{n+k}\left(
\widehat{a}\frac{x}{\left\vert x\right\vert }\right)  ^{2k}\left\vert
x\right\vert ^{-\left(  s+2n\right)  }\\
& =\sum_{k=0}^{n}b_{2k}^{2n}2^{\frac{2n+2k}{2}}\left(  \frac{s}{2}\right)
_{\frac{2n+2k}{2}}\left(  \widehat{a}x\right)  ^{2k}\left\vert x\right\vert
^{-\left(  s+2n+2k\right)  },
\end{align*}

which we can now write in the form%
\begin{equation}
\left(  \widehat{a}D\right)  ^{m}\left\vert x\right\vert ^{-s}=\sum_{j=0}%
^{m}c_{j}^{m}\left(  s\right)  \left(  \widehat{a}x\right)  ^{j}\left\vert
x\right\vert ^{-\left(  s+m+j\right)  },\quad m=0,2,4,\ldots,\label{Ap075}%
\end{equation}

where%
\begin{equation}
c_{j}^{m}\left(  s\right)  =\left\{
\begin{array}
[c]{ll}%
b_{j}^{m}2^{\frac{m+j}{2}}\left(  \frac{s}{2}\right)  _{\frac{m+j}{2}}, &
j\text{ }even,\\
0, & j\text{ }odd.
\end{array}
\right. \label{Ap076}%
\end{equation}

Equations \ref{Ap073} to \ref{Ap076} can be written more concisely as:%
\begin{equation}
\left(  \widehat{a}D\right)  ^{m}\left\vert x\right\vert ^{-s}=\sum_{j=0}%
^{m}c_{j}^{m}\left(  s\right)  \left(  \widehat{a}x\right)  ^{j}\left\vert
x\right\vert ^{-\left(  s+m+j\right)  },\quad m\geq0,\label{Ap077}%
\end{equation}

where%
\begin{equation}
c_{j}^{m}\left(  s\right)  =\left\{
\begin{array}
[c]{ll}%
b_{j}^{m}2^{\frac{m+j}{2}}\left(  \frac{s}{2}\right)  _{\frac{m+j}{2}}, &
j,m\text{ }both\text{ }odd\text{ }or\text{ }even,\\
0, & \left\{
\begin{array}
[c]{c}%
j\text{ }even\text{ }and\text{ }m\text{ }odd\\
or\\
j\text{ }odd\text{ }and\text{ }m\text{ }even,
\end{array}
\right.  .
\end{array}
\right. \label{Ap078}%
\end{equation}

See estimates of Subsection \ref{SbSect_estim_(aD)^m|x|^-s_integ_gamma}.

Using the gamma function form \ref{Ap065} i.e. $\left(  s\right)  _{k}%
=\frac{\Gamma\left(  s+k\right)  }{\Gamma\left(  s\right)  }$, we can simplify
\ref{Ap078} by generalizing $\left(  s\right)  _{m}$\ to:%
\begin{equation}
\left(  s\right)  _{r}=\frac{\Gamma\left(  s+r\right)  }{\Gamma\left(
s\right)  },\quad r,s\in\mathbb{R}^{1},\label{Ap079}%
\end{equation}

so that%
\begin{equation}
c_{j}^{m}\left(  s\right)  =b_{j}^{m}2^{\frac{m+j}{2}}\left(  \frac{s}%
{2}\right)  _{\frac{m+j}{2}}=b_{j}^{m}2^{\frac{m+j}{2}}\frac{\Gamma\left(
\frac{m+s+j}{2}\right)  }{\Gamma\left(  \frac{s}{2}\right)  },\quad
j=0,1,\ldots,m.\label{Ap080}%
\end{equation}

When $a=\mathbf{e}_{k}$, \ref{Ap077} becomes
\begin{align*}
D_{k}^{m}\left\vert x\right\vert ^{-s}  & =\sum\limits_{j=0}^{m}c_{j}%
^{m}\left(  s\right)  x_{k}^{j}\left\vert x\right\vert ^{-\left(
s+m+j\right)  },\quad k=1,\ldots,d,\\
& implies\text{ }that\\
D_{k}^{\alpha_{k}}\left\vert x\right\vert ^{-s_{k}}  & =\sum\limits_{\beta
_{k}\leq\alpha_{k}}c_{\beta_{k}}^{\alpha_{k}}\left(  s_{k}\right)
x_{k}^{\beta_{k}}\left\vert x\right\vert ^{-\left(  s_{k}+\alpha_{k}+\beta
_{k}\right)  }=\sum\limits_{\beta_{k}\leq\alpha_{k}}c_{\beta_{k}}^{\alpha_{k}%
}\left(  s_{k}\right)  x_{k}^{\beta_{k}}\left\vert x\right\vert ^{-s_{k+1}},
\end{align*}

where we have defined%
\begin{align*}
s_{1}  & =s,\\
s_{2}  & =s_{1}+\alpha_{1}+\beta_{1}\\
& =s+\alpha_{1}+\beta_{1},\\
s_{3}  & =s_{2}+\alpha_{2}+\beta_{2}\\
& =s+\left(  \alpha_{1}+\alpha_{2}\right)  +\left(  \beta_{1}+\beta
_{2}\right)  ,\\
& \vdots\\
s_{k+1}  & =s+\left(  \alpha_{1}+\alpha_{2}+\ldots+\alpha_{k}\right)  +\left(
\beta_{1}+\beta_{2}+\ldots+\beta_{k}\right)  <s+2\left(  \alpha_{1}+\alpha
_{2}+\ldots+\alpha_{k}\right)  ,\\
& \vdots\\
s_{d^{\prime}+1}  & =s+\left(  \alpha_{1}+\alpha_{2}+\ldots+\alpha_{d^{\prime
}}\right)  +\left(  \beta_{1}+\beta_{2}+\ldots+\beta_{d^{\prime}}\right)
=s+\left\vert \alpha\right\vert +\left\vert \beta\right\vert .
\end{align*}

so that%
\begin{align*}
D_{1}^{\alpha_{1}}\left\vert x\right\vert ^{-s}  & =\sum\limits_{\beta_{1}%
\leq\alpha_{1}}c_{\beta_{1}}^{\alpha_{1}}\left(  s_{1}\right)  x_{1}%
^{\beta_{1}}\left\vert x\right\vert ^{-\left(  s_{1}+\alpha_{1}+\beta
_{1}\right)  }=\sum\limits_{\beta_{1}\leq\alpha_{1}}c_{\beta_{1}}^{\alpha_{1}%
}\left(  s_{1}\right)  x_{1}^{\beta_{1}}\left\vert x\right\vert ^{-s_{2}},\\
D_{2}^{\alpha_{2}}\left\vert x\right\vert ^{-s_{2}}  & =\sum\limits_{\beta
_{2}\leq\alpha_{2}}c_{\beta_{2}}^{\alpha_{2}}\left(  s_{2}\right)
x_{2}^{\beta_{2}}\left\vert x\right\vert ^{-\left(  s_{2}+\alpha_{2}+\beta
_{2}\right)  }=\sum\limits_{\beta_{2}\leq\alpha_{2}}c_{\beta_{2}}^{\alpha_{2}%
}\left(  s_{2}\right)  x_{2}^{\beta_{2}}\left\vert x\right\vert ^{-s_{3}%
}\Rightarrow\\
D_{2}^{\alpha_{2}}D_{1}^{\alpha_{1}}\left\vert x\right\vert ^{-s}  &
=\sum\limits_{\beta_{1}\leq\alpha_{1}}c_{\beta_{1}}^{\alpha_{1}}\left(
s_{1}\right)  x_{1}^{\beta_{1}}D_{2}^{\alpha_{2}}\left\vert x\right\vert
^{-s_{2}}\\
& =\sum\limits_{\beta_{1}\leq\alpha_{1}}c_{\beta_{1}}^{\alpha_{1}}\left(
s_{1}\right)  x_{1}^{\beta_{1}}\sum\limits_{\beta_{2}\leq\alpha_{2}}%
c_{\beta_{2}}^{\alpha_{2}}\left(  s_{2}\right)  x_{2}^{\beta_{2}}\left\vert
x\right\vert ^{-s_{3}}\\
& =\sum\limits_{\beta_{1}\leq\alpha_{1}}\sum\limits_{\beta_{2}\leq\alpha_{2}%
}c_{\beta_{1}}^{\alpha_{1}}\left(  s_{1}\right)  c_{\beta_{2}}^{\alpha_{2}%
}\left(  s_{2}\right)  x_{1}^{\beta_{1}}x_{2}^{\beta_{2}}\left\vert
x\right\vert ^{-s_{3}},\\
& and\text{ }in\text{ }general\\
D^{\alpha}\left\vert x\right\vert ^{-s}  & =\sum\limits_{\beta\leq\alpha
}c_{\beta_{1}}^{\alpha_{1}}\left(  s_{1}\right)  c_{\beta_{2}}^{\alpha_{2}%
}\left(  s_{2}\right)  \ldots c_{\beta_{d^{\prime}}}^{\alpha_{d^{\prime}}%
}\left(  s_{d^{\prime}}\right)  x^{\beta}\left\vert x\right\vert
^{-s_{d^{\prime}+1}}\\
& =\sum\limits_{\beta\leq\alpha}c_{\beta_{1}}^{\alpha_{1}}\left(
s_{1}\right)  c_{\beta_{2}}^{\alpha_{2}}\left(  s_{2}\right)  \ldots
c_{\beta_{d^{\prime}}}^{\alpha_{d^{\prime}}}\left(  s_{d^{\prime}}\right)
x^{\beta}\left\vert x\right\vert ^{-\left(  s+\left\vert \alpha\right\vert
+\left\vert \beta\right\vert \right)  }\\
& =\sum\limits_{\beta\leq\alpha}c_{\beta_{1}}^{\alpha_{1}}\left(
s_{1}\right)  c_{\beta_{2}}^{\alpha_{2}}\left(  s_{2}\right)  \ldots
c_{\beta_{d^{\prime}}}^{\alpha_{d^{\prime}}}\left(  s_{d^{\prime}}\right)
\widehat{x}^{\beta}\left\vert x\right\vert ^{-s-\left\vert \alpha\right\vert
}\\
& =\left(  \sum\limits_{\beta\leq\alpha}c_{\beta_{1}}^{\alpha_{1}}\left(
s_{1}\right)  c_{\beta_{2}}^{\alpha_{2}}\left(  s_{2}\right)  \ldots
c_{\beta_{d^{\prime}}}^{\alpha_{d^{\prime}}}\left(  s_{d^{\prime}}\right)
\widehat{x}^{\beta}\right)  \left\vert x\right\vert ^{-s-\left\vert
\alpha\right\vert }.
\end{align*}

Using \ref{Ap080} i.e. $c_{j}^{m}\left(  s\right)  =b_{j}^{m}2^{\frac{m+j}{2}%
}\frac{\Gamma\left(  \frac{m+s+j}{2}\right)  }{\Gamma\left(  \frac{s}%
{2}\right)  }$, we have%
\[
c_{\beta_{k}}^{\alpha_{k}}\left(  s_{k}\right)  =b_{\beta_{k}}^{\alpha_{k}%
}2^{\frac{\alpha_{k}+\beta_{k}}{2}}\frac{\Gamma\left(  \frac{\alpha_{k}%
+\beta_{k}+s_{k}}{2}\right)  }{\Gamma\left(  \frac{s_{k}}{2}\right)
}=b_{\beta_{k}}^{\alpha_{k}}2^{\frac{s_{k+1}-s_{k}}{2}}\frac{\Gamma\left(
\frac{s_{k+1}}{2}\right)  }{\Gamma\left(  \frac{s_{k}}{2}\right)  },
\]

so that%
\begin{align*}
D^{\alpha}\left\vert x\right\vert ^{-s}  & =\left(  \sum\limits_{\beta
^{\prime}\leq\alpha^{\prime}}b_{\beta_{1}}^{\alpha_{1}}2^{\frac{s_{2}-s_{1}%
}{2}}\frac{\Gamma\left(  \frac{s_{2}}{2}\right)  }{\Gamma\left(  \frac{s_{1}%
}{2}\right)  }b_{\beta_{2}}^{\alpha_{2}}2^{\frac{s_{3}-s_{2}}{2}}\frac
{\Gamma\left(  \frac{s_{3}}{2}\right)  }{\Gamma\left(  \frac{s_{2}}{2}\right)
}\times\ldots\times b_{\beta_{d^{\prime}}}^{\alpha_{d^{\prime}}}%
2^{\frac{s_{d^{\prime}+1}-s_{d^{\prime}}}{2}}\frac{\Gamma\left(
\frac{s_{d^{\prime}+1}}{2}\right)  }{\Gamma\left(  \frac{s_{d^{\prime}}}%
{2}\right)  }\widehat{x}^{\beta}\right)  \left\vert x\right\vert
^{-s-\left\vert \alpha\right\vert }\\
& =\left(  \sum\limits_{\beta\leq\alpha}2^{\frac{s_{d^{\prime}+1}-s_{1}}{2}%
}\frac{\Gamma\left(  \frac{s_{d^{\prime}+1}}{2}\right)  }{\Gamma\left(
\frac{s_{1}}{2}\right)  }b_{\beta_{1}}^{\alpha_{1}}b_{\beta_{2}}^{\alpha_{2}%
}\ldots b_{\beta_{d^{\prime}}}^{\alpha_{d^{\prime}}}\widehat{x}^{\beta
}\right)  \left\vert x\right\vert ^{-s-\left\vert \alpha\right\vert }\\
& =\left(  \sum\limits_{\beta\leq\alpha}2^{\frac{s+\left\vert \alpha
\right\vert +\left\vert \beta\right\vert -s}{2}}\frac{\Gamma\left(
\frac{s+\left\vert \alpha\right\vert +\left\vert \beta\right\vert }{2}\right)
}{\Gamma\left(  \frac{s}{2}\right)  }b_{\beta_{1}}^{\alpha_{1}}b_{\beta_{2}%
}^{\alpha_{2}}\ldots b_{\beta_{d^{\prime}}}^{\alpha_{d^{\prime}}}\widehat
{x}^{\beta}\right)  \left\vert x\right\vert ^{-s-\left\vert \alpha\right\vert
}\\
& =2^{\frac{\left\vert \alpha\right\vert }{2}}\left(  \sum\limits_{\beta
\leq\alpha}2^{\frac{\left\vert \beta\right\vert }{2}}\frac{\Gamma\left(
\frac{s+\left\vert \alpha\right\vert +\left\vert \beta\right\vert }{2}\right)
}{\Gamma\left(  \frac{s}{2}\right)  }b_{\beta_{1}}^{\alpha_{1}}b_{\beta_{2}%
}^{\alpha_{2}}\ldots b_{\beta_{d^{\prime}}}^{\alpha_{d^{\prime}}}\widehat
{x}^{\beta}\right)  \left\vert x\right\vert ^{-s-\left\vert \alpha\right\vert
}\\
& =2^{\frac{\left\vert \alpha\right\vert }{2}}\left(  \sum\limits_{\beta
\leq\alpha}2^{\frac{\left\vert \beta\right\vert }{2}}\frac{\Gamma\left(
\frac{s+\left\vert \alpha\right\vert +\left\vert \beta\right\vert }{2}\right)
}{\Gamma\left(  \frac{s}{2}\right)  }b_{\beta_{1}}^{\alpha_{1}}b_{\beta_{2}%
}^{\alpha_{2}}\ldots b_{\beta_{d}}^{\alpha_{d}}\widehat{x}^{\beta}\right)
\left\vert x\right\vert ^{-s-\left\vert \alpha\right\vert }\\
& =2^{\frac{\left\vert \alpha\right\vert }{2}}\left(  \sum\limits_{\beta
\leq\alpha}2^{\frac{\left\vert \beta\right\vert }{2}}\left(  \frac{s}%
{2}\right)  _{\frac{\left\vert \alpha\right\vert +\left\vert \beta\right\vert
}{2}}b_{\beta_{1}}^{\alpha_{1}}b_{\beta_{2}}^{\alpha_{2}}\ldots b_{\beta_{d}%
}^{\alpha_{d}}\widehat{x}^{\beta}\right)  \left\vert x\right\vert
^{-s-\left\vert \alpha\right\vert },
\end{align*}

Note that if $\left\vert \alpha\right\vert $ is even then $\left\vert
\beta\right\vert $ is even and if $\left\vert \alpha\right\vert $ is odd then
$\left\vert \beta\right\vert $ is odd. Next define%
\begin{equation}
b_{\beta}^{\alpha}:=b_{\beta_{1}}^{\alpha_{1}}b_{\beta_{2}}^{\alpha_{2}}\ldots
b_{\beta_{d}}^{\alpha_{d}},\label{Ap081}%
\end{equation}

so that%
\begin{equation}
D^{\alpha}\left\vert x\right\vert ^{-s}=2^{\frac{\left\vert \alpha\right\vert
}{2}}\left(  \sum\limits_{\beta\leq\alpha}2^{\frac{\left\vert \beta\right\vert
}{2}}\left(  \frac{s}{2}\right)  _{\frac{\left\vert \alpha\right\vert
+\left\vert \beta\right\vert }{2}}b_{\beta}^{\alpha}\widehat{x}^{\beta
}\right)  \left\vert x\right\vert ^{-s-\left\vert \alpha\right\vert
}.\label{Ap082}%
\end{equation}

Now when $\left\vert \alpha\right\vert +s>0$,%
\begin{equation}
\left(  \frac{s}{2}\right)  _{\frac{\left\vert \alpha\right\vert +\left\vert
\beta\right\vert }{2}}=\frac{1}{\Gamma\left(  \frac{s}{2}\right)  }%
\Gamma\left(  \frac{\left\vert \alpha\right\vert +\left\vert \beta\right\vert
+s}{2}\right)  =\frac{1}{\Gamma\left(  \frac{s}{2}\right)  }\int_{0}^{\infty
}e^{-t}t^{\frac{\left\vert \alpha\right\vert +\left\vert \beta\right\vert
+s}{2}-1}dt,\label{Ap020}%
\end{equation}

and hence%
\begin{align*}
\sum\limits_{\beta\leq\alpha}2^{\frac{\left\vert \beta\right\vert }{2}}\left(
\frac{s}{2}\right)  _{\frac{\left\vert \alpha\right\vert +\left\vert
\beta\right\vert }{2}}b_{\beta}^{\alpha}\widehat{x}^{\beta}  & =\sum
\limits_{\beta\leq\alpha}2^{\frac{\left\vert \beta\right\vert }{2}}\frac
{1}{\Gamma\left(  \frac{s}{2}\right)  }\int_{0}^{\infty}e^{-t}t^{\frac
{\left\vert \alpha\right\vert +\left\vert \beta\right\vert +s}{2}-1}dt\text{
}b_{\beta}^{\alpha}\widehat{x}^{\beta}\\
& =\frac{1}{\Gamma\left(  \frac{s}{2}\right)  }\int_{0}^{\infty}e^{-t}%
\sum\limits_{\beta\leq\alpha}2^{\frac{\left\vert \beta\right\vert }{2}%
}t^{\frac{\left\vert \alpha\right\vert +\left\vert \beta\right\vert +s-2}{2}%
}b_{\beta}^{\alpha}\widehat{x}^{\beta}dt\\
& =\frac{1}{\Gamma\left(  \frac{s}{2}\right)  }\int_{0}^{\infty}e^{-t}%
t^{\frac{\left\vert \alpha\right\vert +s-2}{2}}\left(  \sum\limits_{\beta
\leq\alpha}2^{\frac{\left\vert \beta\right\vert }{2}}t^{\frac{\left\vert
\beta\right\vert }{2}}b_{\beta}^{\alpha}\widehat{x}^{\beta}\right)  dt\\
& =\frac{1}{\Gamma\left(  \frac{s}{2}\right)  }\int_{0}^{\infty}e^{-t}%
t^{\frac{\left\vert \alpha\right\vert +s-2}{2}}\left(  \sum\limits_{\beta
\leq\alpha}\left(  2t\right)  ^{\frac{\left\vert \beta\right\vert }{2}%
}b_{\beta}^{\alpha}\widehat{x}^{\beta}\right)  dt\\
& =\frac{1}{\Gamma\left(  \frac{s}{2}\right)  }\int_{0}^{\infty}e^{-t}%
t^{\frac{\left\vert \alpha\right\vert +s}{2}-1}\prod\limits_{k=1}^{d}\left(
\sum\limits_{\beta_{k}=0}^{\alpha_{k}}b_{\beta_{k}}^{\alpha_{k}}\left(
2t\right)  ^{\frac{\beta_{k}}{2}}\widehat{x}_{k}^{\beta_{k}}\right)  dt\\
& =\frac{1}{\Gamma\left(  \frac{s}{2}\right)  }\int_{0}^{\infty}e^{-t}%
t^{\frac{\left\vert \alpha\right\vert +s}{2}-1}\prod\limits_{k=1}^{d}\left(
\sum\limits_{j=0}^{\alpha_{k}}b_{j}^{\alpha_{k}}\left(  2t\right)  ^{\frac
{j}{2}}\widehat{x}_{k}^{j}\right)  dt.
\end{align*}

\textbf{We can assume that }$x>0$\textbf{\ so }$\widehat{x}>0$ and%
\begin{equation}
\sum\limits_{\beta\leq\alpha}2^{\frac{\left\vert \beta\right\vert }{2}}\left(
\frac{s}{2}\right)  _{\frac{\left\vert \alpha\right\vert +\left\vert
\beta\right\vert }{2}}b_{\beta}^{\alpha}\widehat{x}^{\beta}=\frac{1}%
{\Gamma\left(  \frac{s}{2}\right)  }\int_{0}^{\infty}e^{-t}t^{\frac{\left\vert
\alpha\right\vert +s}{2}-1}\prod\limits_{k=1}^{d}\sum\limits_{j=0}^{\alpha
_{k}}b_{j}^{\alpha_{k}}\left(  2t\left(  \widehat{x}_{k}\right)  ^{2}\right)
^{j/2}dt.\label{Ap091}%
\end{equation}

Define%
\[
f_{n}\left(  u\right)  :=\sum\limits_{j=0}^{n}b_{j}^{n}u^{j/2},\quad u\geq0,
\]

i.e.%
\begin{equation}
\left.
\begin{array}
[c]{l}%
f_{2m}\left(  u\right)  =\sum\limits_{j=0}^{m}b_{2j}^{2m}u^{j},\\
f_{2m+1}\left(  u\right)  =u^{1/2}\sum\limits_{j=0}^{m}b_{2j+1}^{2m+1}u^{j},
\end{array}
\right\}  \quad u\geq0.\label{Ap087}%
\end{equation}

Thus%
\begin{equation}
\sum\limits_{\beta\leq\alpha}2^{\frac{\left\vert \beta\right\vert }{2}}\left(
\frac{s}{2}\right)  _{\frac{\left\vert \alpha\right\vert +\left\vert
\beta\right\vert }{2}}b_{\beta}^{\alpha}\widehat{x}^{\beta}=\frac{1}%
{\Gamma\left(  \frac{s}{2}\right)  }\int_{0}^{\infty}e^{-t}t^{\frac{\left\vert
\alpha\right\vert +s}{2}-1}\prod\limits_{k=1}^{d}f_{\alpha_{k}}\left(
2t\left(  \widehat{x}_{k}\right)  ^{2}\right)  dt.\label{Ap089}%
\end{equation}

The next step is to estimate the $f_{\alpha_{k}}$.

\subsubsection{Bounding $f_{n}$\label{SbSbSect_bound_fn}}

There are two cases:

\fbox{\textbf{Case }$f_{2m}$} Using Matlab we calculate $\ln\left(  \left\vert
f_{2m}\left(  u\right)  \right\vert /\left\vert f_{2m}\left(  0\right)
\right\vert \right)  $ and \textbf{conjecture} that for each $m\geq1 $ there
is a unique straight line $l_{2m}$ which is tangent to all the peaks and such
that $\ln\left(  \left\vert f_{2m}\left(  u\right)  \right\vert /\left\vert
f_{2m}\left(  0\right)  \right\vert \right)  \leq l_{2m}\left(  u\right)  $
when $u\geq0$. We \textbf{also} \textbf{conjecture} that this line also passes
through $\left(  0,\left\vert f_{2m}\left(  0\right)  \right\vert \right)  $.
Denote the slope of the line corresponding to $2m$ by $v_{2m}$. We
\textbf{next} \textbf{conjecture} that $v_{2m+2}>v_{2m}>\ldots>v_{2}$.
Experiments show that%
\begin{equation}
0.2784<v_{2}<0.2785.\label{Ap109}%
\end{equation}

Observe that from \ref{Ap087} and \ref{Ap083},
\begin{equation}
f_{2m}\left(  0\right)  =b_{0}^{2m}=b_{1}^{2m-1}=\left(  -1\right)  ^{m}%
2^{m}\left(  \frac{1}{2}\right)  _{m}.\label{Ap088}%
\end{equation}

So we have%
\[
\ln\left(  \left\vert f_{2m}\left(  u\right)  \right\vert /\left\vert
b_{0}^{2m}\right\vert \right)  \leq v_{2}u,\quad u\geq0.
\]

i.e.%
\begin{equation}
\left\vert f_{2m}\left(  u\right)  \right\vert \leq2^{m}\left(  \frac{1}%
{2}\right)  _{m}e^{v_{2}u},\quad u\geq0,\label{Ap092}%
\end{equation}

or equivalently%
\begin{equation}
\left\vert f_{n}\left(  u\right)  \right\vert \leq2^{\frac{n}{2}}\left(
\frac{1}{2}\right)  _{\frac{n}{2}}e^{v_{2}u},\quad n\text{ }even,\text{ }%
u\geq0\label{Ap105}%
\end{equation}

\fbox{\textbf{Case }$f_{2m+1}$} We start by considering the situation when
$n\geq3$.

Using Matlab we calculate $\ln\left(  \left\vert f_{2m+1}\left(  u\right)
\right\vert /\left\vert b_{1}^{2m+1}\right\vert \right)  $ and
\textbf{conjecture} that for each $m\geq1$ there is a unique straight line
$l_{2m+1} $ which is tangent to all the peaks and such that $\ln\left(
\left\vert f_{2m+1}\left(  u\right)  \right\vert /\left\vert b_{1}%
^{2m+1}\right\vert \right)  \leq l_{2m+1}\left(  u\right)  $ when $u\geq0$
with equality only at the tangents. Suppose now that $l_{2m+1}$ has slope
$v_{2m+1}$ and \textit{y}-intercept $b_{2m+1}$. We \textbf{also conjecture}
that $v_{2m+1}<v_{2m-1}<\ldots<v_{3}$ and $b_{2m+1}<b_{2m-1}<\ldots<b_{3}$.
Experiments show that%
\begin{equation}
v_{3}\simeq0.273358,\text{ }b_{3}\simeq-0.6351,\label{Ap102}%
\end{equation}

which gives a vertical distance between the curves of less than $10^{-3}$. It
is now clear that%
\begin{equation}
\ln\left(  \left\vert f_{2m+1}\left(  u\right)  \right\vert /\left\vert
b_{1}^{2m+1}\right\vert \right)  \leq l_{3}\left(  u\right)  ,\quad
u\geq0,\text{ }m\geq1.\label{Ap100}%
\end{equation}

Now when $n=1$ we have $m=0$ and $f_{1}\left(  u\right)  =\sqrt{u}\left\vert
b_{1}^{1}\right\vert =\sqrt{u}$. Graphically we see that $l_{3}$ can be moved
upwards so that it is tangential to $f_{1}$ and that in this case \ref{Ap100}
also holds when $m=0$. However, using Matlab we can conjecture that there is a
line $l_{1}$, tangential $f_{1}$, that has smaller slope than $l_{3}$ and
which satisfies%
\begin{equation}
\ln\left(  \left\vert f_{2m+1}\left(  u\right)  \right\vert /\left\vert
b_{1}^{2m+1}\right\vert \right)  \leq l_{1}\left(  u\right)  ,\quad
u\geq0,\text{ }m\geq0,\label{Ap101}%
\end{equation}

where%
\begin{equation}
v_{1}\simeq0.247,\text{ }b_{1}\simeq-0.144.\label{Ap103}%
\end{equation}

Thus from \ref{Ap088} and \ref{Ap101},%
\[
\left\vert f_{2m+1}\left(  u\right)  \right\vert \leq\left\vert b_{1}%
^{2m+1}\right\vert e^{l_{1}\left(  u\right)  },
\]

i.e.%
\begin{equation}
\left\vert f_{2m+1}\left(  u\right)  \right\vert \leq2^{m+1}\left(  \frac
{1}{2}\right)  _{m+1}e^{b_{1}+v_{1}u},\quad u\geq0,\text{ }m\geq
0,\label{Ap104}%
\end{equation}

or equivalently%
\begin{equation}
\left\vert f_{n}\left(  u\right)  \right\vert \leq2^{\frac{n+1}{2}}\left(
\frac{1}{2}\right)  _{\frac{n+1}{2}}e^{b_{1}+v_{1}u},\quad u\geq0,\text{
}n\text{ }odd.\label{Ap106}%
\end{equation}

The key property of the estimates \ref{Ap092} and \ref{Ap104} is that the
exponent of $e$ is a linear function of $u$.

We now bound the right side of \ref{Ap089}. First some definitions: suppose
$u\in\mathbb{R}^{d}$ and $u\geq0$; $n_{o}\left(  \alpha\right)  $ is the
number of odd-valued $\alpha$ and $n_{e}\left(  \alpha\right)  $ is the number
of even-valued $\alpha$; $\alpha_{o}^{!}=\prod\limits_{\alpha_{k}\text{ }%
odd}\alpha_{k}$; $\left\vert u\right\vert _{1}=\sum\limits_{k}u_{k}$,
$\left\vert u\right\vert _{1,e}=\sum\limits_{\alpha_{k}\text{ }odd}u_{k}$ and
$\left\vert u\right\vert _{1,e}=\sum\limits_{\alpha_{k}\text{ }even}u_{k}$.
Then using the estimates \ref{Ap105} and \ref{Ap106},%
\begin{align*}
\left\vert \prod\limits_{k=1}^{d}f_{\alpha_{k}}\left(  u_{k}\right)
\right\vert  & =\prod\limits_{\alpha_{j}\text{ }even}\left\vert f_{\alpha_{j}%
}\left(  u_{j}\right)  \right\vert \prod\limits_{\alpha_{k}\text{ }%
odd}\left\vert f_{\alpha_{k}}\left(  u_{k}\right)  \right\vert \\
& \leq\prod\limits_{\alpha_{j}\text{ }even}2^{\frac{\alpha_{j}}{2}}\left(
\frac{1}{2}\right)  _{\frac{\alpha_{j}}{2}}e^{v_{2}u_{j}}\text{ }%
\prod\limits_{\alpha_{k}\text{ }odd}2^{\frac{a_{k}+1}{2}}\left(  \frac{1}%
{2}\right)  _{\frac{a_{k}+1}{2}}e^{b_{1}+v_{1}u_{k}}\\
& <\left(  \prod\limits_{\substack{\alpha_{j} \\even}}2^{\frac{\alpha_{j}}{2}%
}\prod\limits_{\substack{\alpha_{k} \\odd}}2^{\frac{a_{k}+1}{2}}\right)
\left(  \prod\limits_{\substack{\alpha_{j} \\even}}\left(  \frac{1}{2}\right)
_{\frac{\alpha_{j}}{2}}\prod\limits_{\substack{\alpha_{k} \\odd}}\left(
\frac{1}{2}\right)  _{\frac{\alpha_{k}+1}{2}}\right)  \left(  \prod
\limits_{\substack{\alpha_{j} \\even}}e^{v_{2}u_{j}}\prod
\limits_{\substack{\alpha_{k} \\odd}}e^{b_{1}+v_{1}u_{k}}\right) \\
& =2^{\frac{\left\vert \alpha\right\vert }{2}}2^{\frac{n_{o}\left(
\alpha\right)  }{2}}e^{n_{o}\left(  \alpha\right)  b_{1}}\left(
\prod\limits_{\alpha_{j}\text{ }even}\left(  \frac{1}{2}\right)
_{\frac{\alpha_{j}}{2}}\prod\limits_{\alpha_{k}\text{ }odd}\left(  \frac{1}%
{2}\right)  _{\frac{\alpha_{k}+1}{2}}\right)  e^{v_{2}\left\vert u\right\vert
_{1}}.
\end{align*}

But%
\[
\left(  \frac{1}{2}\right)  _{\frac{\alpha_{k}+1}{2}}=\frac{\Gamma\left(
\frac{\alpha_{k}+2}{2}\right)  }{\Gamma\left(  1/2\right)  }=\frac
{\Gamma\left(  \frac{\alpha_{k}}{2}+1\right)  }{\Gamma\left(  1/2\right)
}=\frac{\alpha_{k}}{2}\frac{\Gamma\left(  \frac{\alpha_{k}}{2}\right)
}{\Gamma\left(  1/2\right)  }=\frac{\alpha_{k}}{2}\left(  \frac{1}{2}\right)
_{\frac{\alpha_{k}}{2}},
\]

so we can write%
\begin{align*}
\left\vert \prod\limits_{k=1}^{d}f_{\alpha_{k}}\left(  u_{k}\right)
\right\vert  & <2^{\frac{\left\vert \alpha\right\vert }{2}}2^{\frac
{n_{o}\left(  \alpha\right)  }{2}}e^{n_{o}\left(  \alpha\right)  b_{1}}\left(
\prod\limits_{\alpha_{j}\text{ }even}\left(  \frac{1}{2}\right)
_{\frac{\alpha_{j}}{2}}\prod\limits_{\alpha_{k}\text{ }odd}\frac{\alpha_{k}%
}{2}\left(  \frac{1}{2}\right)  _{\frac{\alpha_{k}}{2}}\right)  e^{v_{2}%
\left\vert u\right\vert _{1}}\\
& =2^{\frac{\left\vert \alpha\right\vert }{2}}2^{\frac{n_{o}\left(
\alpha\right)  }{2}}e^{n_{o}\left(  \alpha\right)  b_{1}}\left(
\prod\limits_{\alpha_{k}\text{ }odd}\frac{\alpha_{k}}{2}\right)  \left(
\frac{1}{2}\right)  _{\frac{\alpha}{2}}e^{v_{2}\left\vert u\right\vert _{1}}\\
& =2^{\frac{\left\vert \alpha\right\vert }{2}}2^{-\frac{n_{o}\left(
\alpha\right)  }{2}}e^{n_{o}\left(  \alpha\right)  b_{1}}\alpha_{o}^{!}\left(
\frac{1}{2}\right)  _{\frac{\alpha}{2}}e^{v_{2}\left\vert u\right\vert _{1}},
\end{align*}

where%
\begin{equation}
\left(  s\right)  _{\beta}:=\left(  s\right)  _{\beta_{1}}\left(  s\right)
_{\beta_{2}}\ldots\left(  s\right)  _{\beta_{d}},\quad\beta=\left(  \beta
_{i}\right)  .\label{Ap112}%
\end{equation}

Further%
\begin{align*}
2^{-\frac{n_{o}\left(  \alpha\right)  }{2}}e^{n_{o}\left(  \alpha\right)
b_{1}}  & =\exp\left(  -\frac{n_{o}\left(  \alpha\right)  }{2}\ln2\right)
\exp\left(  n_{o}\left(  \alpha\right)  b_{1}\right) \\
& =\exp\left(  \left(  b_{1}-\frac{1}{2}\ln2\right)  n_{o}\left(
\alpha\right)  \right) \\
& \leq1,
\end{align*}

since%
\[
b_{1}-\frac{1}{2}\ln2\simeq-0.144-0.3466\simeq-0.490,
\]

so now we have%
\begin{equation}
\left\vert \prod\limits_{k=1}^{d}f_{\alpha_{k}}\left(  u_{k}\right)
\right\vert \leq2^{\frac{\left\vert \alpha\right\vert }{2}}e^{-\left(
-b_{1}+\frac{1}{2}\ln2\right)  n_{o}\left(  \alpha\right)  }\alpha_{o}%
^{!}\left(  \frac{1}{2}\right)  _{\alpha/2}e^{v_{2}\left\vert u\right\vert
_{1}}.\label{Ap107}%
\end{equation}

Thus \ref{Ap089} can be estimated by%
\begin{align}
&  \left\vert \sum\limits_{\beta\leq\alpha}2^{\frac{\left\vert \beta
\right\vert }{2}}\left(  \frac{s}{2}\right)  _{\frac{\left\vert \alpha
\right\vert +\left\vert \beta\right\vert }{2}}b_{\beta}^{\alpha}\widehat
{x}^{\beta}\right\vert \nonumber\\
&  \leq\frac{1}{\Gamma\left(  \frac{s}{2}\right)  }\int_{0}^{\infty}%
e^{-t}t^{\frac{\left\vert \alpha\right\vert +s}{2}-1}\left\vert \prod
\limits_{k=1}^{d}f_{\alpha_{k}}\left(  2t\left(  \widehat{x}_{k}\right)
^{2}\right)  \right\vert dt\nonumber\\
&  \leq\frac{1}{\Gamma\left(  \frac{s}{2}\right)  }\int_{0}^{\infty}%
e^{-t}t^{\frac{\left\vert \alpha\right\vert +s}{2}-1}2^{\frac{\left\vert
\alpha\right\vert }{2}}e^{-\left(  -b_{1}+\frac{1}{2}\ln2\right)  n_{o}\left(
\alpha\right)  }\alpha_{o}^{!}\left(  \frac{1}{2}\right)  _{\alpha/2}%
e^{v_{2}\left\vert 2t\left(  \widehat{x}_{k}\right)  ^{2}\right\vert _{1}%
}dt\nonumber\\
&  =2^{\frac{\left\vert \alpha\right\vert }{2}}e^{-\left(  -b_{1}+\frac{1}%
{2}\ln2\right)  n_{o}\left(  \alpha\right)  }\alpha_{o}^{!}\frac{\left(
\frac{1}{2}\right)  _{\alpha/2}}{\Gamma\left(  \frac{s}{2}\right)  }\int%
_{0}^{\infty}e^{-t}t^{\frac{\left\vert \alpha\right\vert +s}{2}-1}%
e^{v_{2}\left\vert 2t\left(  \widehat{x}_{k}\right)  ^{2}\right\vert _{1}%
}dt\nonumber\\
&  =2^{\frac{\left\vert \alpha\right\vert }{2}}e^{-\left(  -b_{1}+\frac{1}%
{2}\ln2\right)  n_{o}\left(  \alpha\right)  }\alpha_{o}^{!}\frac{\left(
\frac{1}{2}\right)  _{\alpha/2}}{\Gamma\left(  \frac{s}{2}\right)  }\int%
_{0}^{\infty}e^{-t}t^{\frac{\left\vert \alpha\right\vert +s}{2}-1}e^{2v_{2}%
t}dt\nonumber\\
&  =2^{\frac{\left\vert \alpha\right\vert }{2}}e^{-\left(  -b_{1}+\frac{1}%
{2}\ln2\right)  n_{o}\left(  \alpha\right)  }\alpha_{o}^{!}\frac{\left(
\frac{1}{2}\right)  _{\alpha/2}}{\Gamma\left(  \frac{s}{2}\right)  }\int%
_{0}^{\infty}e^{-\left(  1-2v_{2}\right)  t}t^{\frac{\left\vert \alpha
\right\vert +s}{2}-1}dt.\label{Ap108}%
\end{align}

Noting that \ref{Ap109} implies $0<2v_{2}<1$, the next step is to express the
integral in \ref{Ap108} as a gamma function using the change of variables
$r=\left(  1-2v_{2}\right)  t$, $dr=\left(  1-2v_{2}\right)  dt$:%
\begin{align}
\int_{0}^{\infty}e^{-\left(  1-2v_{2}\right)  t}t^{\frac{\left\vert
\alpha\right\vert +s}{2}-1}dt  & =\int_{0}^{\infty}e^{-r}\left(  \frac
{r}{1-2v_{2}}\right)  ^{\frac{\left\vert \alpha\right\vert +s}{2}-1}\frac
{dr}{1-2v_{2}}\nonumber\\
& =\left(  1-2v_{2}\right)  ^{-\frac{\left\vert \alpha\right\vert +s}{2}}%
\int_{0}^{\infty}e^{-r}r^{\frac{\left\vert \alpha\right\vert +s}{2}%
-1}dr\nonumber\\
& =\left(  1-2v_{2}\right)  ^{-\frac{\left\vert \alpha\right\vert +s}{2}%
}\Gamma\left(  \frac{\left\vert \alpha\right\vert +s}{2}\right)
,\label{Ap116}%
\end{align}

so that%
\begin{align}
& \left\vert \sum\limits_{\beta\leq\alpha}2^{\frac{\left\vert \beta\right\vert
}{2}}\left(  \frac{s}{2}\right)  _{\frac{\left\vert \alpha\right\vert
+\left\vert \beta\right\vert }{2}}b_{\beta}^{\alpha}\widehat{x}^{\beta
}\right\vert \nonumber\\
& \leq2^{\frac{\left\vert \alpha\right\vert }{2}}e^{-\left(  -b_{1}+\frac
{1}{2}\ln2\right)  n_{o}\left(  \alpha\right)  }\alpha_{o}^{!}\left(  \frac
{1}{2}\right)  _{\alpha/2}\frac{\Gamma\left(  \frac{\left\vert \alpha
\right\vert +s}{2}\right)  }{\Gamma\left(  \frac{s}{2}\right)  }\left(
1-2v_{2}\right)  ^{-\frac{\left\vert \alpha\right\vert +s}{2}}\nonumber\\
& =2^{\frac{\left\vert \alpha\right\vert }{2}}e^{-\left(  -b_{1}+\frac{1}%
{2}\ln2\right)  n_{o}\left(  \alpha\right)  }\alpha_{o}^{!}\left(  \frac{1}%
{2}\right)  _{\alpha/2}\left(  \frac{s}{2}\right)  _{\left\vert \alpha
\right\vert /2}\left(  1-2v_{2}\right)  ^{-\frac{\left\vert \alpha\right\vert
+s}{2}},\label{Ap110}%
\end{align}

where from \ref{Ap109},%
\begin{equation}
0.4430<1-2v_{2}<0.4432.\label{Ap117}%
\end{equation}

Thus we can now estimate \ref{Ap082} by%
\begin{align}
\left\vert D^{\alpha}\left\vert x\right\vert ^{-s}\right\vert  & \leq
2^{\frac{\left\vert \alpha\right\vert }{2}}\left\vert \sum\limits_{\beta
\leq\alpha}2^{\frac{\left\vert \beta\right\vert }{2}}\left(  \frac{s}%
{2}\right)  _{\frac{\left\vert \alpha\right\vert +\left\vert \beta\right\vert
}{2}}b_{\beta}^{\alpha}\widehat{x}^{\beta}\right\vert \left\vert x\right\vert
^{-s-\left\vert \alpha\right\vert }\nonumber\\
& \leq2^{\left\vert \alpha\right\vert }e^{-\left(  -b_{1}+\frac{1}{2}%
\ln2\right)  n_{o}\left(  \alpha\right)  }\alpha_{o}^{!}\left(  \frac{1}%
{2}\right)  _{\alpha/2}\left(  \frac{s}{2}\right)  _{\left\vert \alpha
\right\vert /2}\left(  1-2v_{2}\right)  ^{-\frac{\left\vert \alpha\right\vert
+s}{2}}\left\vert x\right\vert ^{-s-\left\vert \alpha\right\vert
}\label{Ap113}\\
when\text{ }\left\vert \alpha\right\vert +s  & >0.\nonumber
\end{align}

\fbox{A bound for $\left(  \frac{1}{2}\right)  _{\alpha/2}$} Write%
\begin{equation}
\left(  \frac{1}{2}\right)  _{\alpha/2}=\prod\limits_{k=1}^{d}\frac
{\Gamma\left(  \frac{\alpha_{k}+1}{2}\right)  }{\Gamma\left(  \frac{1}%
{2}\right)  }=\frac{1}{\Gamma\left(  \frac{1}{2}\right)  ^{d}}\prod
\limits_{k=1}^{d}\Gamma\left(  \frac{\alpha_{k}+1}{2}\right)  =\frac{1}%
{\pi^{d/2}}\prod\limits_{k=1}^{d}\Gamma\left(  \frac{\alpha_{k}+1}{2}\right)
.\label{Ap114}%
\end{equation}

Now%
\[
\prod\limits_{k=1}^{d}\Gamma\left(  \frac{\alpha_{k}+1}{2}\right)
=\prod\limits_{k=1}^{d}\int_{0}^{\infty}e^{-t_{k}^{2}}t_{k}^{2\left(
\frac{\alpha_{k}+1}{2}\right)  -1}dt_{k}=\prod\limits_{k=1}^{d}\int%
_{0}^{\infty}e^{-t_{k}^{2}}t_{k}^{\alpha_{k}}dt_{k}=\int_{\mathbb{R}_{+}^{d}%
}e^{-\left\vert \tau\right\vert ^{2}}\tau^{\alpha}d\tau,
\]

so that clearly%
\[
\left\vert \prod\limits_{k=1}^{d}\Gamma\left(  \frac{\alpha_{k}+1}{2}\right)
\right\vert \leq\int_{\mathbb{R}_{+}^{d}}e^{-\left\vert \tau\right\vert ^{2}%
}\left\vert \tau\right\vert ^{\left\vert \alpha\right\vert }d\tau,
\]

and hence%
\begin{equation}
\left(  \frac{1}{2}\right)  _{\alpha/2}\leq\frac{1}{\pi^{d/2}}\int%
_{\mathbb{R}_{+}^{d}}e^{-\left\vert \tau\right\vert ^{2}}\left\vert
\tau\right\vert ^{\left\vert \alpha\right\vert }d\tau,\label{Ap111}%
\end{equation}

so we use the spherical polar coordinates change of variables e.g. Section
5.43 of Adams \cite{Adams75}:%
\[
d\tau=d\tau_{1}d\tau_{2}\ldots d\tau_{d-1}=r^{d-1}\prod\limits_{j=1}^{d}%
\sin^{j-1}\phi_{j}drd\phi_{1}\ldots d\phi_{d-1},
\]

to obtain for $d\geq2$,%
\begin{align*}
\int\limits_{\mathbb{R}_{+}^{d}}e^{-\left\vert \tau\right\vert ^{2}}\left\vert
\tau\right\vert ^{\left\vert \alpha\right\vert }d\tau & =\underset{\left(
d-2\right)  \times}{\underbrace{\int\limits_{0}^{\pi}\ldots\int\limits_{0}%
^{\pi}}}\int\limits_{-\pi}^{\pi}\int\limits_{0}^{\infty}e^{-\rho^{2}%
}r^{\left\vert \alpha\right\vert }r^{d-1}\prod\limits_{j=1}^{d-1}\sin
^{j-1}\phi_{j}drd\phi_{1}\ldots d\phi_{d-1}\\
& =2\pi\left(  \prod\limits_{j=2}^{d-1}\int\limits_{0}^{\pi}\sin^{j-1}\theta
d\theta\right)  \int\limits_{0}^{\infty}e^{-r^{2}}r^{\left\vert \alpha
\right\vert +d-1}dr\\
& =2\pi\left(  \prod\limits_{j=1}^{d-2}\int\limits_{0}^{\pi}\sin^{j}\theta
d\theta\right)  \int\limits_{0}^{\infty}e^{-r^{2}}r^{2\left(  \frac{\left\vert
\alpha\right\vert +d}{2}\right)  -1}dr\\
& =2\pi\left(  \prod\limits_{j=1}^{d-2}\int\limits_{0}^{\pi}\sin^{j}\theta
d\theta\right)  \int\limits_{0}^{\infty}e^{-r^{2}}r^{2\left(  \frac{\left\vert
\alpha\right\vert +d}{2}\right)  -1}dr\\
& =2\pi\left(  2^{d-2}\prod\limits_{j=1}^{d-2}\int\limits_{0}^{\pi/2}\sin
^{j}\theta d\theta\right)  \Gamma\left(  \frac{\left\vert \alpha\right\vert
+d}{2}\right)  .
\end{align*}

From exercise 10.4.9 (b) of Arfken \cite{Arfken70},%
\[
\int\limits_{0}^{\pi/2}\sin^{j}\theta d\theta=\frac{\sqrt{\pi}}{2}%
\frac{\left(  \left(  j-1\right)  /2\right)  !}{\left(  j/2\right)  !},
\]

so that%
\begin{align*}
\prod\limits_{j=1}^{d-2}\int\limits_{0}^{\pi/2}\sin^{j}\theta d\theta=\left(
\frac{\sqrt{\pi}}{2}\right)  ^{d-2}\frac{0!}{\left(  \left(  d-2\right)
/2\right)  !} &  =\left(  \frac{\sqrt{\pi}}{2}\right)  ^{d-2}\frac{1}{\left(
\frac{d}{2}-1\right)  !}\\
&  =\left(  \frac{\sqrt{\pi}}{2}\right)  ^{d-2}\frac{1}{\Gamma\left(  \frac
{d}{2}\right)  },
\end{align*}

and hence%
\begin{align*}
\int\limits_{\mathbb{R}_{+}^{d}}e^{-\left\vert \tau\right\vert ^{2}}\left\vert
\tau\right\vert ^{\left\vert \alpha\right\vert }d\tau\leq2\pi2^{d-2}\left(
\frac{\sqrt{\pi}}{2}\right)  ^{d-2}\frac{1}{\Gamma\left(  \frac{d}{2}\right)
}\Gamma\left(  \frac{\left\vert \alpha\right\vert +d}{2}\right)   & =2\pi
\frac{\left(  \sqrt{\pi}\right)  ^{d-2}}{\Gamma\left(  \frac{d}{2}\right)
}\Gamma\left(  \frac{\left\vert \alpha\right\vert +d}{2}\right) \\
& =\frac{2\pi^{d/2}}{\Gamma\left(  \frac{d}{2}\right)  }\Gamma\left(
\frac{\left\vert \alpha\right\vert +d}{2}\right) \\
& =2\pi^{d/2}\left(  \frac{d}{2}\right)  _{\left\vert \alpha\right\vert /2}.
\end{align*}

From \ref{Ap111},%
\[
\left(  \frac{1}{2}\right)  _{\alpha/2}\leq\frac{1}{\pi^{d/2}}\int%
_{\mathbb{R}_{+}^{d}}e^{-\left\vert \tau\right\vert ^{2}}\left\vert
\tau\right\vert ^{\left\vert \alpha\right\vert }d\tau\leq\left(  \frac{d}%
{2}\right)  _{\left\vert \alpha\right\vert /2},
\]

and the right side depends on $\left\vert \alpha\right\vert $. As a
consequence \ref{Ap113} now becomes%
\begin{align}
\left\vert D^{\alpha}\left\vert x\right\vert ^{-s}\right\vert  &
\leq2^{\left\vert \alpha\right\vert }e^{-\left(  -b_{1}+\frac{1}{2}%
\ln2\right)  n_{o}\left(  \alpha\right)  }\alpha_{o}^{!}\left(  \frac{1}%
{2}\right)  _{\alpha/2}\left(  \frac{s}{2}\right)  _{\left\vert \alpha
\right\vert /2}\left(  1-2v_{2}\right)  ^{-\frac{\left\vert \alpha\right\vert
+s}{2}}\left\vert x\right\vert ^{-s-\left\vert \alpha\right\vert }\nonumber\\
& \leq2^{\left\vert \alpha\right\vert }e^{-\left(  -b_{1}+\frac{1}{2}%
\ln2\right)  n_{o}\left(  \alpha\right)  }\alpha_{o}^{!}2\left(  \frac{d}%
{2}\right)  _{\left\vert \alpha\right\vert /2}\left(  \frac{s}{2}\right)
_{\left\vert \alpha\right\vert /2}\left(  1-2v_{2}\right)  ^{-\frac{\left\vert
\alpha\right\vert +s}{2}}\left\vert x\right\vert ^{-s-\left\vert
\alpha\right\vert }\nonumber\\
& =2^{\left\vert \alpha\right\vert +1}e^{-\left(  -b_{1}+\frac{1}{2}%
\ln2\right)  n_{o}\left(  \alpha\right)  }\alpha_{o}^{!}\left(  \frac{d}%
{2}\right)  _{\left\vert \alpha\right\vert /2}\left(  \frac{s}{2}\right)
_{\left\vert \alpha\right\vert /2}\left(  1-2v_{2}\right)  ^{-\frac{\left\vert
\alpha\right\vert +s}{2}}\left\vert x\right\vert ^{-s-\left\vert
\alpha\right\vert }\label{Ap115}\\
when\text{ }\left\vert \alpha\right\vert +s  & >0.\nonumber
\end{align}

Further%
\begin{align*}
\left\vert D^{\alpha}\left\vert x\right\vert ^{-s}\right\vert  &
\leq2^{\left\vert \alpha\right\vert +1}\alpha_{o}^{!}\left(  \frac{d}%
{2}\right)  _{\left\vert \alpha\right\vert /2}\left(  \frac{s}{2}\right)
_{\left\vert \alpha\right\vert /2}\left(  1-2v_{2}\right)  ^{-\frac{\left\vert
\alpha\right\vert +s}{2}}\left\vert x\right\vert ^{-s-\left\vert
\alpha\right\vert }\\
& \leq2^{\left\vert \alpha\right\vert +1}\alpha!\left(  \frac{d}{2}\right)
_{\left\vert \alpha\right\vert /2}\left(  \frac{s}{2}\right)  _{\left\vert
\alpha\right\vert /2}\left(  1-2v_{2}\right)  ^{-\frac{\left\vert
\alpha\right\vert +s}{2}}\left\vert x\right\vert ^{-s-\left\vert
\alpha\right\vert }.
\end{align*}

But from \ref{Ap029}%
\[
\sqrt{2\pi}e^{-u+\frac{1}{1+12u}}u^{u+\frac{1}{2}}<\Gamma\left(  u+1\right)
<\sqrt{2\pi}e^{-u+\frac{1}{12u}}u^{u+\frac{1}{2}},\quad u>0,
\]

so that%
\begin{align*}
\alpha!  & =\alpha_{1}!\alpha_{2}!\ldots\alpha_{d}!\\
& <\left(  \sqrt{2\pi}e^{-\alpha_{1}+\frac{1}{12\alpha_{1}}}\alpha_{1}%
^{\alpha_{1}+\frac{1}{2}}\right)  \left(  \sqrt{2\pi}e^{-\alpha_{2}+\frac
{1}{12\alpha_{2}}}\alpha_{2}^{\alpha_{2}+\frac{1}{2}}\right)  \ldots\left(
\sqrt{2\pi}e^{-\alpha_{d}+\frac{1}{12\alpha_{d}}}\alpha_{d}^{\alpha_{d}%
+\frac{1}{2}}\right) \\
& \leq\left(  \sqrt{2\pi}\right)  ^{d}e^{\frac{d}{12}}e^{-\left\vert
\alpha\right\vert }\alpha_{1}^{\alpha_{1}+\frac{1}{2}}\alpha_{2}^{\alpha
_{2}+\frac{1}{2}}\ldots\alpha_{d}^{\alpha_{d}+\frac{1}{2}}\\
& \leq\left(  e^{\frac{1}{12}}\sqrt{2\pi}\right)  ^{d}e^{-\left\vert
\alpha\right\vert }\left\vert \alpha\right\vert ^{\alpha_{1}+\frac{1}{2}%
}\left\vert \alpha\right\vert ^{\alpha_{2}+\frac{1}{2}}\ldots\left\vert
\alpha\right\vert ^{\alpha_{d}+\frac{1}{2}}\\
& =\left(  e^{\frac{1}{12}}\sqrt{2\pi}\right)  ^{d}e^{-\left\vert
\alpha\right\vert }\left\vert \alpha\right\vert ^{\left\vert \alpha\right\vert
+\frac{d}{2}},
\end{align*}

and hence we \textbf{conjecture}: \underline{when $\left\vert \alpha
\right\vert +s>0$},%
\begin{equation}
\left\vert D^{\alpha}\left\vert x\right\vert ^{-s}\right\vert \leq\left(
e^{\frac{1}{12}}\sqrt{2\pi}\right)  ^{d}2^{\left\vert \alpha\right\vert
+1}e^{-\left\vert \alpha\right\vert }\left\vert \alpha\right\vert ^{\left\vert
\alpha\right\vert +\frac{d}{2}}\left(  \frac{d}{2}\right)  _{\left\vert
\alpha\right\vert /2}\left(  \frac{s}{2}\right)  _{\left\vert \alpha
\right\vert /2}\left(  1-2v_{2}\right)  ^{-\frac{\left\vert \alpha\right\vert
+s}{2}}\left\vert x\right\vert ^{-s-\left\vert \alpha\right\vert
}.\label{Ap024}%
\end{equation}

\subsection{Estimating $\left(  \widehat{a}D\right)  ^{m}\left\vert
x\right\vert ^{-s}$ using integral form of the $\Gamma$
function\label{SbSect_estim_(aD)^m|x|^-s_integ_gamma}}

?? Condition on $s$? ?? \textbf{Negative result}

From \ref{Ap077} and \ref{Ap080},%
\[
\left(  \widehat{a}D\right)  ^{m}\left\vert x\right\vert ^{-s}=\sum_{j=0}%
^{m}c_{j}^{m}\left(  s\right)  \left(  \widehat{a}x\right)  ^{j}\left\vert
x\right\vert ^{-\left(  s+m+j\right)  },\quad m\geq0,
\]

and%
\[
c_{j}^{m}\left(  s\right)  =b_{j}^{m}2^{\frac{m+j}{2}}\left(  \frac{s}%
{2}\right)  _{\frac{m+j}{2}}.
\]

But from \ref{Ap029},%
\[
\left(  \frac{s}{2}\right)  _{\frac{m+j}{2}}=\frac{1}{\Gamma\left(  \frac
{s}{2}\right)  }\Gamma\left(  \frac{m+j+s}{2}\right)  =\frac{1}{\Gamma\left(
\frac{s}{2}\right)  }\int_{0}^{\infty}e^{-t}t^{\frac{m+j+s}{2}-1}dt.
\]

Now write%
\[
\left(  \widehat{a}D\right)  ^{m}\left\vert x\right\vert ^{-s}=\sum_{j=0}%
^{m}c_{j}^{m}\left(  s\right)  \left(  \widehat{a}x\right)  ^{j}\left\vert
x\right\vert ^{-\left(  s+m+j\right)  }=\left(  \sum_{j=0}^{m}c_{j}^{m}\left(
s\right)  \left(  \widehat{a}\widehat{x}\right)  ^{j}\right)  \left\vert
x\right\vert ^{-\left(  s+m\right)  },
\]

and we consider the polynomial $\sum_{j=0}^{m}c_{j}^{m}\left(  s\right)
u^{j}$:%
\begin{align*}
\sum_{j=0}^{m}c_{j}^{m}\left(  s\right)  u^{j}  & =\sum_{j=0}^{m}b_{j}%
^{m}2^{\frac{m+j}{2}}\left(  \frac{s}{2}\right)  _{\frac{m+j}{2}}u^{j}\\
& =\sum_{j=0}^{m}b_{j}^{m}2^{\frac{m+j}{2}}\frac{u^{j}}{\Gamma\left(  \frac
{s}{2}\right)  }\int_{0}^{\infty}e^{-t}t^{\frac{m+j+s}{2}-1}dt\\
& =\frac{2^{m/2}}{\Gamma\left(  \frac{s}{2}\right)  }\int_{0}^{\infty}%
e^{-t}\sum_{j=0}^{m}b_{j}^{m}2^{\frac{j}{2}}t^{\frac{m+j+s}{2}-1}u^{j}dt\\
& =\frac{2^{m/2}}{\Gamma\left(  \frac{s}{2}\right)  }\int_{0}^{\infty}%
e^{-t}\sum_{j=0}^{m}b_{j}^{m}2^{\frac{j}{2}}t^{\frac{m+j+s}{2}-1}u^{j}dt\\
& =\frac{2^{m/2}}{\Gamma\left(  \frac{s}{2}\right)  }\int_{0}^{\infty}%
e^{-t}t^{\frac{m+s}{2}-1}\sum_{j=0}^{m}b_{j}^{m}2^{\frac{j}{2}}t^{\frac{j}{2}%
}u^{j}dt\\
& =\frac{2^{m/2}}{\Gamma\left(  \frac{s}{2}\right)  }\int_{0}^{\infty}%
e^{-t}t^{\frac{m+s}{2}-1}\sum_{j=0}^{m}b_{j}^{m}2^{\frac{j}{2}}t^{\frac{j}{2}%
}u^{j}dt\\
& =\frac{2^{m/2}}{\Gamma\left(  \frac{s}{2}\right)  }\int_{0}^{\infty}%
e^{-t}t^{\frac{m+s}{2}-1}\sum_{j=0}^{m}b_{j}^{m}\left(  2tu^{2}\right)
^{j/2}dt\\
& =\frac{2^{m/2}}{\Gamma\left(  \frac{s}{2}\right)  }\int_{0}^{\infty}%
e^{-t}t^{\frac{m+s}{2}-1}f_{m}\left(  2tu^{2}\right)  dt.
\end{align*}

using definition \ref{Ap027} of $f_{m}$.\medskip

\fbox{If $m=2n$} then this equation becomes%
\[
\sum_{j=0}^{2n}c_{j}^{2n}\left(  s\right)  u^{j}=\frac{2^{n}}{\Gamma\left(
\frac{s}{2}\right)  }\int_{0}^{\infty}e^{-t}t^{\left(  n+\frac{s}{2}\right)
-1}f_{2n}\left(  2tu^{2}\right)  dt,\quad2n+s>0,
\]

and from \ref{Ap092},%
\[
\left\vert f_{2n}\left(  \tau\right)  \right\vert \leq2^{n}\left(  \frac{1}%
{2}\right)  _{n}e^{v_{2}\tau},\quad\tau\geq0,
\]

so that%
\begin{align*}
\left\vert \sum_{j=0}^{2n}c_{j}^{2n}\left(  s\right)  u^{j}\right\vert  &
\leq\frac{2^{n}}{\left\vert \Gamma\left(  \frac{s}{2}\right)  \right\vert
}\int_{0}^{\infty}e^{-t}t^{\left(  n+\frac{s}{2}\right)  -1}2^{n}\left(
\frac{1}{2}\right)  _{n}e^{v_{2}2tu^{2}}dt\\
& =2^{2n}\left(  \frac{1}{2}\right)  _{n}\frac{1}{\left\vert \Gamma\left(
\frac{s}{2}\right)  \right\vert }\int_{0}^{\infty}e^{-\left(  1-2v_{2}%
u^{2}\right)  t}t^{\left(  n+\frac{s}{2}\right)  -1}dt.
\end{align*}

From \ref{Ap116},%
\[
\int_{0}^{\infty}e^{-\left(  1-2v_{2}\right)  t}t^{\frac{\left\vert
\alpha\right\vert +s}{2}-1}dt=\left(  1-2v_{2}\right)  ^{-\frac{\left\vert
\alpha\right\vert +s}{2}}\Gamma\left(  \frac{\left\vert \alpha\right\vert
+s}{2}\right)  ,
\]

so now%
\begin{align*}
\left\vert \sum_{j=0}^{2n}c_{j}^{2n}\left(  s\right)  u^{j}\right\vert  &
\leq2^{2n}\left(  \frac{1}{2}\right)  _{n}\frac{1}{\left\vert \Gamma\left(
\frac{s}{2}\right)  \right\vert }\left(  1-2v_{2}u^{2}\right)  ^{-\left(
n+\frac{s}{2}\right)  }\Gamma\left(  n+\frac{s}{2}\right) \\
& =2^{2n}\left\vert \left(  \frac{s}{2}\right)  _{n}\right\vert \left(
\frac{1}{2}\right)  _{n}\left(  1-2v_{2}u^{2}\right)  ^{-\left(  n+\frac{s}%
{2}\right)  }\\
& \leq2^{2n}\left\vert \left(  \frac{s}{2}\right)  _{n}\right\vert \left(
\frac{1}{2}\right)  _{n}\left(  1-2v_{2}\right)  ^{-\left(  n+\frac{s}%
{2}\right)  },\\
when\text{ }s  & >-2n.
\end{align*}

From \ref{Ap117}, $0.4430<1-2v_{2}<0.4432$, so%
\[
\left\vert \sum_{j=0}^{2n}c_{j}^{2n}\left(  s\right)  u^{j}\right\vert
\leq2^{2n}\left\vert \left(  \frac{s}{2}\right)  _{n}\right\vert \left(
\frac{1}{2}\right)  _{n}\frac{1}{\left(  0.4430\right)  ^{n+\frac{s}{2}}}.
\]

Compare this result with: from \ref{Ap061} and \ref{Ap055},%
\begin{align*}
\left\vert \left(  \widehat{a}D\right)  ^{2n}\left\vert x\right\vert
^{-s}\right\vert  & \leq\left\{  \left\{
\begin{array}
[c]{ll}%
2^{2n}\left\vert \left(  \frac{s}{2}\right)  _{n}\right\vert \left(  \frac
{1}{2}\right)  _{n}, & -2\leq s\leq0,\\
\left(  s\right)  _{2n}, & s\geq0,
\end{array}
\right.  \right\}  \left\vert x\right\vert ^{-s-2n}\\
& =\left\{
\begin{array}
[c]{ll}%
2^{2n}\left\vert \left(  \frac{s}{2}\right)  _{n}\right\vert \left(  \frac
{1}{2}\right)  _{n}, & -2\leq s\leq0,\\
2^{2n}\left(  \frac{s}{2}\right)  _{n}\left(  \frac{s+1}{2}\right)  _{n}, &
s\geq0.
\end{array}
\right\}  \left\vert x\right\vert ^{-s-2n}.
\end{align*}

\textbf{We conclude that} using the integral form of the gamma function does
not provide a sufficiently good estimate.

Also from \ref{Ap025} and \ref{Ap094},%
\[
b_{2k}^{2n}=b_{0}^{2n}g_{2k}^{2n}=\left(  -1\right)  ^{n}2^{n}\left(  \frac
{1}{2}\right)  _{n}g_{2k}^{2n},\quad n\geq1,
\]

so%
\begin{align*}
\sum_{j=0}^{2n}c_{j}^{2n}\left(  s\right)  u^{j}  & =\frac{2^{n}}%
{\Gamma\left(  \frac{s}{2}\right)  }\int_{0}^{\infty}e^{-t}t^{\left(
n+\frac{s}{2}\right)  -1}\sum_{k=0}^{n}\left(  -1\right)  ^{n}2^{n}\left(
\frac{1}{2}\right)  _{n}g_{2k}^{2n}\left(  2tu^{2}\right)  ^{k}dt\\
& =\left(  -1\right)  ^{n}2^{2n}\left(  \frac{1}{2}\right)  _{n}\frac
{1}{\Gamma\left(  \frac{s}{2}\right)  }\int_{0}^{\infty}e^{-t}t^{\left(
n+\frac{s}{2}\right)  -1}\sum_{k=0}^{n}g_{2k}^{2n}\left(  2tu^{2}\right)
^{k}dt
\end{align*}

\section{More properties of the $b_{k}^{n}$}

For values of $b_{k}^{n}$ see \ref{Ap040}.

Plotting $\log\left\vert b_{k}^{n}\right\vert $ gives a very nice plot which suggests:

\begin{conjecture}
For $n\geq2$,%
\begin{align*}
\left\vert b_{2k+1}^{2n-1}\right\vert  & <\left\vert b_{2k}^{2n}\right\vert
<\left\vert b_{2k+1}^{2n+1}\right\vert ,\quad1\leq k\leq n-1,\\
\left\vert b_{2n}^{2n}\right\vert  & =\left\vert b_{2n+1}^{2n+1}\right\vert
,\\
\left\vert b_{1}^{2n-1}\right\vert  & =\left\vert b_{2}^{2n}\right\vert .
\end{align*}

\end{conjecture}

?? \textbf{What to do next? Corresponding formula for the constants }%
$g_{k}^{j}$? ??

\section{Orthonormal transform}

From \ref{Ap022},%
\[
\left(  \widehat{a}D\right)  ^{2n+1}\left\vert x\right\vert ^{-s}%
=2^{n+1}\left(  \frac{s}{2}\right)  _{n+1}\left(  \widehat{a}\widehat
{x}\right)  p_{s}^{\left(  2n+1\right)  }\left(  \widehat{a}\widehat
{x}\right)  \left\vert x\right\vert ^{-s-\left(  2n+1\right)  }.
\]

If $\mathcal{O}$ is an orthonormal transformation then $\left\vert
\mathcal{O}y\right\vert =\left\vert y\right\vert $ and $x\mathcal{O}%
^{T}y:=\left(  x,\mathcal{O}^{T}y\right)  =\left(  \mathcal{O}x,y\right)
=\mathcal{O}y\,x$.

Note that an orthogonal transformation $\mathcal{O}$ satisfies $\mathcal{O}%
^{T}=\mathcal{O}^{-1}$ where $\left(  \mathcal{O}x,y\right)  =\left(
x,\mathcal{O}^{T}y\right)  $ for the Euclidean inner product.

Note that $\mathcal{O}:\mathbb{R}^{d}\rightarrow\mathbb{R}^{d}$ is an
isometric isomorphism.

Note that $\widehat{\mathcal{O}a}=\mathcal{O}\widehat{a}$.%
\begin{align*}
\left(  \left(  \widehat{\mathcal{O}a}D\right)  ^{2n+1}\left\vert
\cdot\right\vert ^{-s}\right)  \left(  x\right)   & =\left(  \widehat
{\mathcal{O}a}D\right)  ^{2n+1}\left\vert x\right\vert ^{-s}\\
& =2^{n+1}\left(  \frac{s}{2}\right)  _{n+1}\left(  \widehat{\mathcal{O}%
a}\widehat{x}\right)  p_{s}^{2n+1}\left(  \widehat{\mathcal{O}a}\widehat
{x}\right)  \left\vert x\right\vert ^{-s-\left(  2n+1\right)  }\\
& =2^{n+1}\left(  \frac{s}{2}\right)  _{n+1}\left(  \mathcal{O}\widehat
{a},\widehat{x}\right)  p_{s}^{2n+1}\left(  \left(  \mathcal{O}\widehat
{a},\widehat{x}\right)  \right)  \left\vert x\right\vert ^{-s-\left(
2n+1\right)  }\\
& =2^{n+1}\left(  \frac{s}{2}\right)  _{n+1}\left(  \widehat{a},\mathcal{O}%
^{T}\widehat{x}\right)  p_{s}^{2n+1}\left(  \left(  \widehat{a},\mathcal{O}%
^{T}\widehat{x}\right)  \right)  \left\vert x\right\vert ^{-s-\left(
2n+1\right)  }\\
& =2^{n+1}\left(  \frac{s}{2}\right)  _{n+1}\left(  \widehat{a},\widehat
{\mathcal{O}^{T}x}\right)  p_{s}^{2n+1}\left(  \left(  \widehat{a}%
,\widehat{\mathcal{O}^{T}x}\right)  \right)  \left\vert \mathcal{O}%
^{T}x\right\vert ^{-s-\left(  2n+1\right)  }\\
& =\left(  \left(  \widehat{a}D\right)  ^{2n+1}\left\vert \cdot\right\vert
^{-s}\right)  \left(  \mathcal{O}^{T}x\right)  .
\end{align*}

\[
\left(  \widehat{a}D\right)  ^{2n+1}\left\vert x\right\vert ^{-s}=\left(
\left(  \widehat{\mathcal{O}a}D\right)  ^{2n+1}\left\vert \cdot\right\vert
^{-s}\right)  \left(  \mathcal{O}x\right)  .
\]

\section{?? Estimating $D^{\beta}\left\vert x\right\vert ^{-s}$ - matrix
approach}%

\[
\sum_{\left\vert \beta\right\vert =n}\frac{a^{\beta}}{\beta!}D^{\beta}%
f=\frac{1}{n!}\left(  aD\right)  ^{n}f,\quad a\in\mathbb{R}^{d}.
\]

\[
I_{n}:=\left\{  \beta:\left\vert \beta\right\vert =n\right\}  ,\quad
\#I_{n}=\binom{n+d-1}{n}.
\]

\textbf{Approach} Choose $\binom{n+d-1}{n}$ values of $a$ such that we have
$\binom{n+d-1}{n}$ linear equations for $\left\{  D^{\beta}f\right\}
_{\left\vert \beta\right\vert =n}$ and the matrix is regular.

What about choosing%
\[
\sum_{\left\vert \beta\right\vert =n}\frac{\alpha^{\beta}}{\beta!}D^{\beta
}f=\frac{1}{n!}\left(  \alpha D\right)  ^{n}f,\quad\alpha\in I_{n}?
\]

Is the matrix $M_{n}=\left(  \frac{\alpha^{\beta}}{\beta!}\right)
_{\left\vert \alpha\right\vert ,\left\vert \beta\right\vert =n}$ regular?

\section{Another idea - unfinished ??}%

\[
\frac{1}{m!}\left(  \left(  a.D\right)  \cdot\left(  a.D\right)  \right)
^{m}\left\vert x\right\vert ^{-s}=\sum_{\left\vert \beta\right\vert =m}%
\frac{a^{\beta}D^{\beta}a^{\beta}D^{\beta}}{\beta!}\left\vert x\right\vert
^{-s}=\sum_{\left\vert \beta\right\vert =m}\frac{a^{2\beta}}{\beta!}D^{2\beta
}\left\vert x\right\vert ^{-s},
\]

implies%
\[
D_{a}^{2\gamma}\frac{1}{m!}\left\vert a.D\right\vert ^{2m}\left\vert
x\right\vert ^{-s}=\frac{\left(  2\gamma\right)  !}{\gamma!}D^{2\gamma
}\left\vert x\right\vert ^{-s},\quad\left\vert \gamma\right\vert =m.
\]

\[
\frac{1}{m!}D_{a}^{2\gamma}\left\vert D_{x}\right\vert ^{2m}\left\vert
a.x\right\vert ^{-s}=\frac{\left(  2\gamma\right)  !}{\gamma!}D^{2\gamma
}\left\vert x\right\vert ^{-s},\quad\left\vert \gamma\right\vert =m.
\]

\[
D_{i}\left(  u\left(  \left\vert x\right\vert \right)  \right)  =u^{\prime
}\left(  \left\vert x\right\vert \right)  D_{i}\left\vert x\right\vert
=\left\vert x\right\vert ^{-1}u^{\prime}\left(  \left\vert x\right\vert
\right)  x_{i}%
\]
and so
\[
a_{i}^{2}D_{i}^{2}\left(  u\left(  \left\vert x\right\vert \right)  \right)
=a_{i}^{2}\left\vert x\right\vert ^{-1}u^{\prime}\left(  \left\vert
x\right\vert \right)  +\left\vert x\right\vert ^{-2}u^{\prime\prime}\left(
\left\vert x\right\vert \right)  a_{i}^{2}x_{i}^{2}-\left\vert x\right\vert
^{-3}u^{\prime}\left(  \left\vert x\right\vert \right)  a_{i}^{2}x_{i}^{2}.
\]

Hence
\begin{align*}
\left\vert \widehat{a}.D\right\vert ^{2}\left(  u\left(  \left\vert
x\right\vert \right)  \right)   & =\left\vert \widehat{a}_{i}D_{i}\right\vert
^{2}\left(  u\left(  \left\vert x\right\vert \right)  \right) \\
& =\left(  \left\vert x\right\vert ^{-1}-\left\vert \widehat{a}.x\right\vert
^{2}\left\vert x\right\vert ^{-3}\right)  u^{\prime}\left(  \left\vert
x\right\vert \right)  +\left\vert \widehat{a}.x\right\vert ^{2}\left\vert
x\right\vert ^{-2}u^{\prime\prime}\left(  \left\vert x\right\vert \right) \\
& =\left(  \left\vert x\right\vert ^{2}-\left\vert \widehat{a}.x\right\vert
^{2}\right)  \left\vert x\right\vert ^{-3}u^{\prime}\left(  \left\vert
x\right\vert \right)  +\left\vert \widehat{a}.x\right\vert ^{2}\left\vert
x\right\vert ^{-2}u^{\prime\prime}\left(  \left\vert x\right\vert \right) \\
& =\mathcal{L}_{\widehat{a}}u.
\end{align*}

When $u\left(  t\right)  =t^{-s}$:%
\begin{align*}
\left\vert \widehat{a}.D\right\vert ^{2}\left(  u\left(  \left\vert
x\right\vert \right)  \right)   & =\left(  \left\vert x\right\vert
^{2}-\left\vert \widehat{a}.x\right\vert ^{2}\right)  \left\vert x\right\vert
^{-3}\left(  -s\left\vert x\right\vert ^{-s-1}\right)  +\left\vert \widehat
{a}.x\right\vert ^{2}\left\vert x\right\vert ^{-2}s\left(  s+1\right)
\left\vert x\right\vert ^{-s-2}\\
& =\left(  -s\left(  \left\vert x\right\vert ^{2}-\left\vert \widehat
{a}.x\right\vert ^{2}\right)  +s\left(  s+1\right)  \left\vert \widehat
{a}.x\right\vert ^{2}\right)  \left\vert x\right\vert ^{-s-4}\\
& =s\left(  -\left\vert x\right\vert ^{2}+\left\vert \widehat{a}.x\right\vert
^{2}+\left(  s+1\right)  \left\vert \widehat{a}.x\right\vert ^{2}\right)
\left\vert x\right\vert ^{-s-4}\\
& =s\left(  -\left\vert x\right\vert ^{2}+\left(  s+2\right)  \left\vert
\widehat{a}.x\right\vert ^{2}\right)  \left\vert x\right\vert ^{-s-4},
\end{align*}

and so%
\[
\left\vert a.D\right\vert ^{2}\left(  u\left(  \left\vert x\right\vert
\right)  \right)  =s\left(  -\left\vert a\right\vert ^{2}\left\vert
x\right\vert ^{2}+\left(  s+2\right)  \left\vert a.x\right\vert ^{2}\right)
\left\vert x\right\vert ^{-s-4}.
\]

\begin{align*}
\left\vert b.D\right\vert ^{2}\left(  \left\vert b.x\right\vert ^{2}\right)
& =\sum\limits_{k}b_{k}^{2}D_{k}^{2}\left(  \ldots+b_{k}^{2}x_{k}^{2}%
+\ldots\right)  =2\sum\limits_{k}b_{k}^{4}.\\
\left\vert b.D\right\vert ^{2}\left(  \left\vert b.x\right\vert ^{2}\left\vert
x\right\vert ^{-s}\right)   & =??
\end{align*}

\section{Bounds for $D^{\alpha}\left\vert x\right\vert ^{k}$ and $D^{\alpha
}\left(  \left\vert x\right\vert ^{k}\log\left\vert x\right\vert \right)  $}

??

\begin{lemma}
\label{Lem_bnd_deriv_rad_func}For $\alpha\geq0,$ $t\in\mathbb{R}^{1}$ there
exist constants $c_{\alpha,t}$ and $c_{\alpha,t}^{\prime}$ such that:

\begin{enumerate}
\item $\left\vert D^{\alpha}\left\vert x\right\vert ^{t}\right\vert \leq
c_{\alpha,t}\left\vert x\right\vert ^{t-\left\vert \alpha\right\vert }$.

\item $\left\vert D^{\alpha}\left\vert x\right\vert ^{t}\log\left\vert
x\right\vert \right\vert \leq c_{\alpha,t}^{\prime}\left\vert x\right\vert
^{t-\left\vert \alpha\right\vert }$.
\end{enumerate}
\end{lemma}

\begin{proof}
Since $\left\vert x\right\vert ^{t}$ and $\left\vert x\right\vert ^{t}%
\log\left\vert x\right\vert $ are both homogeneous of order $t$ it follows
that the derivative of order $\alpha$ is homogeneous of order $t-\left\vert
\alpha\right\vert $.
\end{proof}

Put bounds using constants depending on $\left\vert \alpha\right\vert $ in
this lemma.

\begin{lemma}
\label{Lem_deriv_rad_func}:

\begin{enumerate}
\item \textbf{CONJECTURE} There exist positive constants $\left\{
c_{m,t}:t\in\mathbb{R}^{1},\text{ }m=0,1,2,\ldots\right\}  $ such that%
\[
\left\vert D^{\alpha}\left\vert x\right\vert ^{t}\right\vert \leq
c_{\left\vert \alpha\right\vert ,t}\left\vert x\right\vert ^{t-\left\vert
\alpha\right\vert },\quad\alpha\geq0,\text{ }t\in\mathbb{R}^{1}.
\]

Also $c_{\alpha,t}=0$ when $t>0$ is an even integer and $\left\vert
\alpha\right\vert >k$.

\item For all $\alpha$,%
\[
\left\vert D^{\alpha}\log\left\vert x\right\vert \right\vert \leq
c_{\left\vert \alpha\right\vert }^{\prime}\left\vert x\right\vert
^{-\left\vert \alpha\right\vert },
\]

where%
\[
c_{\left\vert \alpha\right\vert }^{\prime}:=\left\{
\begin{array}
[c]{ll}%
2^{\left\vert \alpha\right\vert }\left\vert \alpha\right\vert !, & \alpha
\leq\mathbf{1},\\
c_{\left\vert \alpha\right\vert -1,-2}+\left(  \left\vert \alpha\right\vert
-1\right)  c_{\left\vert \alpha\right\vert -2,-2}, & \alpha\nleqslant
\mathbf{1}.
\end{array}
\right.
\]

Now define%
\begin{equation}
b_{n,k}:=\max_{m\leq n}c_{m,k},\quad b_{n}^{\prime}:=\max\limits_{m\leq
n}c_{m}^{\prime}.\label{Ap229}%
\end{equation}

\item For all $\alpha$ and integer $k>0$,%
\[
\left\vert D^{\alpha}\left(  \left\vert x\right\vert ^{k}\log\left\vert
x\right\vert \right)  \right\vert \leq b_{\left\vert \alpha\right\vert
,k}\left(  \left(  2^{\left\vert \alpha\right\vert }-1\right)  b_{\left\vert
\alpha\right\vert }^{\prime}+\left\vert \log\left\vert x\right\vert
\right\vert \right)  \left\vert x\right\vert ^{k-\left\vert \alpha\right\vert
}.
\]

\item If $k>0$ is even then
\[
\left\vert D^{\alpha}\left(  \left\vert x\right\vert ^{k}\log\left\vert
x\right\vert \right)  \right\vert \leq\left\{
\begin{array}
[c]{ll}%
b_{\left\vert \alpha\right\vert ,k}\left(  \left(  2^{\left\vert
\alpha\right\vert }-1\right)  b_{\left\vert \alpha\right\vert }^{\prime
}+\left\vert \log\left\vert x\right\vert \right\vert \right)  \left\vert
x\right\vert ^{k-\left\vert \alpha\right\vert }, & \left\vert \alpha
\right\vert \leq k,\\
\left(  2^{\left\vert \alpha\right\vert }-1\right)  b_{\left\vert
\alpha\right\vert ,k}b_{\left\vert \alpha\right\vert }^{\prime}\left\vert
x\right\vert ^{k-\left\vert \alpha\right\vert }, & \left\vert \alpha
\right\vert >k.
\end{array}
\right.
\]

\item For $k=0,1,2,\ldots$,%
\[
\left\vert D^{\beta}\left\vert x\right\vert ^{2k}\right\vert \leq2^{2k}%
\pi^{-d/2}d^{2k-\left\vert \beta\right\vert }\frac{\left(  k!\right)  ^{2}%
}{\left(  2k-\left\vert \beta\right\vert \right)  !}\left\vert x\right\vert
^{2k-\left\vert \beta\right\vert },\quad\left\vert \beta\right\vert \leq2k.
\]

\item For $k=0,1,2,\ldots$,%
\[
\left\vert D^{\beta}\left\vert x\right\vert ^{2k+1}\right\vert \leq
2^{2k+\left\vert \beta\right\vert }\pi^{-d/2}d^{2k}\left(  k!\right)
^{2}b_{\left\vert \beta\right\vert ,1}\left\vert x\right\vert
^{2k+1-\left\vert \beta\right\vert },\quad\left\vert \beta\right\vert
\leq2k+1.
\]

\end{enumerate}
\end{lemma}

\begin{proof}
\textbf{Part 1 }From the conjecture \ref{Ap024} there \textbf{MAY} exist such
constants $c_{\left\vert \alpha\right\vert ,t}$.\medskip

\textbf{Part 2} \fbox{Case $\alpha\leq\mathbf{1}$} Suppose $\alpha_{k_{j}}=1$
for $j=1,2,\ldots,d^{\prime}\leq d$.

Then
\begin{align*}
D_{k_{1}}\log\left\vert x\right\vert  & =\left\vert x\right\vert ^{-1}%
D_{k_{1}}\left\vert x\right\vert =\left\vert x\right\vert ^{-2}x_{k_{1}},\\
D_{k_{2}}D_{k_{1}}\log\left\vert x\right\vert  & =-2\left\vert x\right\vert
^{-4}x_{k_{1}}x_{k_{2}},
\end{align*}

so that
\begin{align*}
D^{\alpha}\log\left\vert x\right\vert  & =\left(  -1\right)  ^{\left\vert
\alpha\right\vert +1}2\times4\times\ldots\times\left(  2\left\vert
\alpha\right\vert \right)  \left\vert x\right\vert ^{-2\left\vert
\alpha\right\vert }x^{\alpha}\\
& =\left(  -1\right)  ^{\left\vert \alpha\right\vert +1}2^{\left\vert
\alpha\right\vert }\left\vert \alpha\right\vert !\left\vert x\right\vert
^{-2\left\vert \alpha\right\vert }x^{\alpha}.
\end{align*}

Hence%
\[
\left\vert D^{\alpha}\log\left\vert x\right\vert \right\vert \leq2^{\left\vert
\alpha\right\vert }\left\vert \alpha\right\vert !\left\vert x\right\vert
^{-2\left\vert \alpha\right\vert }\left\vert x\right\vert ^{\left\vert
\alpha\right\vert }=2^{\left\vert \alpha\right\vert }\left\vert \alpha
\right\vert !\left\vert x\right\vert ^{-\left\vert \alpha\right\vert },
\]

so choose $c_{\left\vert \alpha\right\vert }^{\prime}=2^{\left\vert
\alpha\right\vert }\left\vert \alpha\right\vert !$.\medskip

\fbox{Case $\alpha\nleqslant\mathbf{1}$} Since $D_{1}\log\left\vert
x\right\vert =\left\vert x\right\vert ^{-1}D_{1}\left\vert x\right\vert
=\left\vert x\right\vert ^{-2}x_{1}$ we have%
\begin{align*}
D_{1}^{\alpha_{1}}\log\left\vert x\right\vert  & =D_{1}^{\alpha_{1}-1}\left(
\left\vert x\right\vert ^{-2}x_{1}\right)  =\sum\limits_{j\leq\alpha_{1}%
-1}\tbinom{\alpha_{1}-\mathbf{1}}{j}D_{1}^{\alpha_{1}-1-j}\left(  \left\vert
x\right\vert ^{-2}\right)  D_{1}^{j}x_{1}\\
& =\sum\limits_{j=0}^{1}\tbinom{\alpha_{1}-\mathbf{1}}{j}D_{1}^{\alpha
_{1}-1-j}\left(  \left\vert x\right\vert ^{-2}\right)  D_{1}^{j}x_{1}\\
& =\tbinom{\alpha_{1}-1}{0}D_{1}^{\alpha_{1}-1-0}\left(  \left\vert
x\right\vert ^{-2}\right)  D_{1}^{0}x_{1}+\tbinom{\alpha_{1}-1}{1}%
D_{1}^{\alpha_{1}-1-1}\left(  \left\vert x\right\vert ^{-2}\right)  D_{1}%
^{1}x_{1}\\
& =x_{1}D_{1}^{\alpha_{1}-1}\left(  \left\vert x\right\vert ^{-2}\right)
+\left(  \alpha_{1}-1\right)  D_{1}^{\alpha_{1}-2}\left(  \left\vert
x\right\vert ^{-2}\right)  ,
\end{align*}

so%
\[
D^{\alpha}\log\left\vert x\right\vert =x_{1}D^{\alpha-\mathbf{e}_{1}}\left(
\left\vert x\right\vert ^{-2}\right)  +\left(  \alpha_{1}-1\right)
D^{\alpha-2\mathbf{e}_{1}}\left(  \left\vert x\right\vert ^{-2}\right)  .
\]

Now from part 1,%
\begin{align*}
\left\vert D^{\alpha}\log\left\vert x\right\vert \right\vert  & \leq\left\vert
x_{1}\right\vert \left\vert D^{\alpha-\mathbf{e}_{1}}\left(  \left\vert
x\right\vert ^{-2}\right)  \right\vert +\left(  \left\vert \alpha\right\vert
-1\right)  \left\vert D^{\alpha-2\mathbf{e}_{1}}\left(  \left\vert
x\right\vert ^{-2}\right)  \right\vert \\
& \leq c_{\left\vert \alpha-\mathbf{e}_{1}\right\vert ,-2}\left\vert
x_{1}\right\vert \left\vert x\right\vert ^{-\left\vert \alpha\right\vert
-1}+\left(  \alpha_{1}-1\right)  c_{\left\vert \alpha-2\mathbf{e}%
_{1}\right\vert ,-2}\left\vert x\right\vert ^{-\left\vert \alpha\right\vert
}\\
& =c_{\left\vert \alpha\right\vert -1,-2}\left\vert x_{1}\right\vert
\left\vert x\right\vert ^{-\left\vert \alpha\right\vert -1}+\left(  \alpha
_{1}-1\right)  c_{\left\vert \alpha\right\vert -2,-2}\left\vert x\right\vert
^{-\left\vert \alpha\right\vert }\\
& \leq c_{\left\vert \alpha\right\vert -1,-2}\left\vert x\right\vert
\left\vert x\right\vert ^{-\left\vert \alpha\right\vert -1}+\left(  \left\vert
\alpha\right\vert -1\right)  c_{\left\vert \alpha\right\vert -2,-2}\left\vert
x\right\vert ^{-\left\vert \alpha\right\vert }\\
& =\left(  c_{\left\vert \alpha\right\vert -1,-2}+\left(  \left\vert
\alpha\right\vert -1\right)  c_{\left\vert \alpha\right\vert -2,-2}\right)
\left\vert x\right\vert ^{-\left\vert \alpha\right\vert },
\end{align*}

so choose $c_{\left\vert \alpha\right\vert }^{\prime}:=c_{\left\vert
\alpha\right\vert -1,-2}+\left(  \left\vert \alpha\right\vert -1\right)
c_{\left\vert \alpha\right\vert -2,-2}$.\medskip

\textbf{Part 3 }If $k>0$ then
\begin{align*}
D^{\alpha}\left(  \left\vert x\right\vert ^{k}\log\left\vert x\right\vert
\right)   & =\sum\limits_{\gamma\leq\alpha}\binom{\alpha}{\gamma}\left(
D^{\alpha-\gamma}\left\vert x\right\vert ^{k}\right)  D^{\gamma}\log\left\vert
x\right\vert \\
& =\sum\limits_{\mathbf{0}<\gamma\leq\alpha}\binom{\alpha}{\gamma}\left(
D^{\alpha-\gamma}\left\vert x\right\vert ^{k}\right)  D^{\gamma}\log\left\vert
x\right\vert +\left(  D^{\alpha}\left\vert x\right\vert ^{k}\right)
\log\left\vert x\right\vert .
\end{align*}

Thus from parts 1 and 2,%
\begin{align*}
\left\vert D^{\alpha}\left(  \left\vert x\right\vert ^{k}\log\left\vert
x\right\vert \right)  \right\vert  & \leq\sum\limits_{\mathbf{0}<\gamma
\leq\alpha}\binom{\alpha}{\gamma}\left\vert D^{\alpha-\gamma}\left\vert
x\right\vert ^{k}\right\vert \left\vert D^{\gamma}\log\left\vert x\right\vert
\right\vert +\left\vert \left(  D^{\alpha}\left\vert x\right\vert ^{k}\right)
\log\left\vert x\right\vert \right\vert \\
& \leq\sum\limits_{\mathbf{0}<\gamma\leq\alpha}\binom{\alpha}{\gamma
}c_{\left\vert \alpha-\gamma\right\vert ,k}\left\vert x\right\vert
^{k-\left\vert \alpha-\gamma\right\vert }c_{\left\vert \gamma\right\vert
}^{\prime}\left\vert x\right\vert ^{-\left\vert \gamma\right\vert
}+c_{\left\vert \alpha\right\vert ,k}\left\vert x\right\vert ^{k-\left\vert
\alpha\right\vert }\left\vert \log\left\vert x\right\vert \right\vert \\
& =\left(  \sum\limits_{\mathbf{0}<\gamma\leq\alpha}\binom{\alpha}{\gamma
}c_{\left\vert \alpha\right\vert -\left\vert \gamma\right\vert ,k}%
c_{\left\vert \gamma\right\vert }^{\prime}+c_{\left\vert \alpha\right\vert
,k}\left\vert \log\left\vert x\right\vert \right\vert \right)  \left\vert
x\right\vert ^{k-\left\vert \alpha\right\vert }\\
& \leq\left(  \sum\limits_{\mathbf{0}<\gamma\leq\alpha}\binom{\alpha}{\gamma
}\left(  \max_{m<\left\vert \alpha\right\vert }c_{m,k}\right)  \left(
\max_{0<m\leq\left\vert \alpha\right\vert }c_{m}^{\prime}\right)
+c_{\left\vert \alpha\right\vert ,k}\left\vert \log\left\vert x\right\vert
\right\vert \right)  \left\vert x\right\vert ^{k-\left\vert \alpha\right\vert
}\\
& \leq\left(  \max_{m\leq\left\vert \alpha\right\vert }c_{m,k}\right)  \left(
\sum\limits_{\mathbf{0}<\gamma\leq\alpha}\binom{\alpha}{\gamma}\left(
\max_{0<m\leq\left\vert \alpha\right\vert }c_{m}^{\prime}\right)  +\left\vert
\log\left\vert x\right\vert \right\vert \right)  \left\vert x\right\vert
^{k-\left\vert \alpha\right\vert }\\
& =\left(  \max_{m\leq\left\vert \alpha\right\vert }c_{m,k}\right)  \left(
\sum\limits_{\mathbf{0}<\gamma\leq\alpha}\binom{\alpha}{\gamma}\mathbf{1}%
^{\gamma}\mathbf{1}^{\alpha-\gamma}\left(  \max_{0<m\leq\left\vert
\alpha\right\vert }c_{m}^{\prime}\right)  +\left\vert \log\left\vert
x\right\vert \right\vert \right)  \left\vert x\right\vert ^{k-\left\vert
\alpha\right\vert }\\
& =\left(  \max_{m\leq\left\vert \alpha\right\vert }c_{m,k}\right)  \left(
\left(  \left(  \mathbf{1}+\mathbf{1}\right)  ^{\alpha}-1\right)  \left(
\max_{0<m\leq\left\vert \alpha\right\vert }c_{m}^{\prime}\right)  +\left\vert
\log\left\vert x\right\vert \right\vert \right)  \left\vert x\right\vert
^{k-\left\vert \alpha\right\vert }\\
& =\left(  \max_{m\leq\left\vert \alpha\right\vert }c_{m,k}\right)  \left(
\left(  2^{\left\vert \alpha\right\vert }-1\right)  \max_{0<m\leq\left\vert
\alpha\right\vert }c_{m}^{\prime}+\left\vert \log\left\vert x\right\vert
\right\vert \right)  \left\vert x\right\vert ^{k-\left\vert \alpha\right\vert
}\\
& =b_{\left\vert \alpha\right\vert ,k}\left(  \left(  2^{\left\vert
\alpha\right\vert }-1\right)  b_{\left\vert \alpha\right\vert }^{\prime
}+\left\vert \log\left\vert x\right\vert \right\vert \right)  \left\vert
x\right\vert ^{k-\left\vert \alpha\right\vert }.
\end{align*}
\medskip

\textbf{Part 4} Apply part 1 to part 3.\smallskip

\textbf{Part 5} If $\left\vert \beta\right\vert \leq2k$ then $D^{\beta
}\left\vert x\right\vert ^{2k}=k!\sum\limits_{\left\vert \alpha\right\vert
=k}\frac{D^{\beta}x^{2\alpha}}{\alpha!}$. But $\frac{1}{\beta!}D^{\beta
}x^{2\alpha}=\binom{2\alpha}{\beta}x^{2\alpha-\beta}$ if $\beta\leq2\alpha$ so%
\begin{align*}
D^{\beta}\left\vert x\right\vert ^{2k}  & =k!\sum\limits_{\substack{\left\vert
\alpha\right\vert =k \\2\alpha\geq\beta}}\frac{D^{\beta}x^{2\alpha}}{\alpha
!}=k!\sum\limits_{\substack{\left\vert \alpha\right\vert =k \\2\alpha\geq
\beta}}\frac{\beta!}{\alpha!}\tbinom{2\alpha}{\beta}x^{2\alpha-\beta}\\
& =k!\sum\limits_{\substack{\left\vert \alpha\right\vert =k \\2\alpha\geq
\beta}}\frac{\beta!}{\alpha!}\frac{\left(  2\alpha\right)  !}{\left(
2\alpha-\beta\right)  !\beta!}x^{2\alpha-\beta}\\
& =k!\sum\limits_{\substack{\left\vert \alpha\right\vert =k \\2\alpha\geq
\beta}}\frac{1}{\alpha!}\frac{\left(  2\alpha\right)  !}{\left(  2\alpha
-\beta\right)  !}x^{2\alpha-\beta}.
\end{align*}

Thus%
\[
\left\vert D^{\beta}\left\vert x\right\vert ^{2k}\right\vert \leq\left(
k!\sum\limits_{\substack{\left\vert \alpha\right\vert =k \\2\alpha\geq\beta
}}\frac{1}{\alpha!}\frac{\left(  2\alpha\right)  !}{\left(  2\alpha
-\beta\right)  !}\right)  \left\vert x\right\vert ^{2k-\beta}.
\]

Now%
\begin{align*}
\sum\limits_{\substack{\left\vert \alpha\right\vert =k \\2\alpha\geq\beta
}}\frac{1}{\alpha!}\frac{\left(  2\alpha\right)  !}{\left(  2\alpha
-\beta\right)  !}  & \leq\left(  \max_{_{\substack{\left\vert \alpha
\right\vert =k \\2\alpha\geq\beta}}}\frac{\left(  2\alpha\right)  !}{\alpha
!}\right)  \sum\limits_{\substack{\left\vert \alpha\right\vert =k
\\2\alpha\geq\beta}}\frac{1}{\left(  2\alpha-\beta\right)  !}\\
& =\left(  \max_{_{\substack{\left\vert \alpha\right\vert =k \\2\alpha
\geq\beta}}}\alpha!\frac{\left(  2\alpha\right)  !}{\alpha!\alpha!}\right)
\sum\limits_{\substack{\left\vert \alpha\right\vert =k \\2\alpha\geq\beta
}}\frac{1}{\left(  2\alpha-\beta\right)  !}\\
& \leq\left(  \max_{\left\vert \alpha\right\vert =k}\alpha!\right)  \left(
\max_{_{\substack{\left\vert \alpha\right\vert =k \\2\alpha\geq\beta}}}%
\frac{\left(  2\alpha\right)  !}{\alpha!\alpha!}\right)  \sum
\limits_{\substack{\left\vert \alpha\right\vert =k \\2\alpha\geq\beta}%
}\frac{1}{\left(  2\alpha-\beta\right)  !}\\
& =\left(  \max_{\left\vert \alpha\right\vert =k}\alpha!\right)  \left(
\max_{_{\substack{\left\vert \alpha\right\vert =k \\2\alpha\geq\beta}}}%
\tbinom{2n}{n}\right)  \sum\limits_{\substack{\left\vert \alpha\right\vert =k
\\2\alpha\geq\beta}}\frac{1}{\left(  2\alpha-\beta\right)  !}.
\end{align*}

But $\max\limits_{\left\vert \alpha\right\vert =k}\alpha!\leq k!$ and%
\[
\tbinom{2n}{n}<2^{2n}\pi^{-1/2}\Longrightarrow\tbinom{2\alpha}{\alpha
}<2^{2\left\vert \alpha\right\vert }\pi^{-d/2}=2^{2k}\pi^{-d/2},
\]

so $\left(  \max\limits_{\left\vert \alpha\right\vert =k}\alpha!\right)
\max\limits_{_{\substack{\left\vert \alpha\right\vert =k \\2\alpha\geq\beta}%
}}\tbinom{2n}{n}\leq2^{2k}\pi^{-d/2}k!$ and hence%
\[
k!\sum\limits_{\substack{\left\vert \alpha\right\vert =k \\2\alpha\geq\beta
}}\frac{1}{\alpha!}\frac{\left(  2\alpha\right)  !}{\left(  2\alpha
-\beta\right)  !}\leq2^{2k}\pi^{-d/2}\left(  k!\right)  ^{2}.
\]

Further%
\[
\sum\limits_{\substack{\left\vert \alpha\right\vert =k \\2\alpha\geq\beta
}}\frac{1}{\left(  2\alpha-\beta\right)  !}\leq\sum\limits_{\left\vert
\gamma\right\vert =2k-\left\vert \beta\right\vert }\frac{1}{\gamma!}%
=\sum\limits_{\left\vert \gamma\right\vert =2k-\left\vert \beta\right\vert
}\frac{\mathbf{1}^{2\gamma}}{\gamma!}=\frac{\left\vert \mathbf{1}\right\vert
^{2\left(  2k-\left\vert \beta\right\vert \right)  }}{\left(  2k-\left\vert
\beta\right\vert \right)  !}=\frac{d^{2k-\left\vert \beta\right\vert }%
}{\left(  2k-\left\vert \beta\right\vert \right)  !},
\]

so%
\[
k!\sum\limits_{\substack{\left\vert \alpha\right\vert =k \\2\alpha\geq\beta
}}\frac{1}{\alpha!}\frac{\left(  2\alpha\right)  !}{\left(  2\alpha
-\beta\right)  !}\leq2^{2k}\pi^{-d/2}d^{2k-\left\vert \beta\right\vert }%
\frac{\left(  k!\right)  ^{2}}{\left(  2k-\left\vert \beta\right\vert \right)
!}.
\]

\[
\left\vert D^{\beta}\left\vert x\right\vert ^{2k}\right\vert \leq\left(
2^{2k}\pi^{-d/2}d^{2k-\left\vert \beta\right\vert }\frac{\left(  k!\right)
^{2}}{\left(  2k-\left\vert \beta\right\vert \right)  !}\right)  \left\vert
x\right\vert ^{2k-\left\vert \beta\right\vert }.
\]

\textbf{Part 6 }Using the estimates of part 5 and part 1:
\begin{align*}
\left\vert D^{\beta}\left\vert x\right\vert ^{2k+1}\right\vert  & =\left\vert
D^{\beta}\left\vert x\right\vert ^{2k}\left\vert x\right\vert \right\vert
=\left\vert \sum\limits_{\alpha\leq\beta}\tbinom{\beta}{\alpha}D^{\alpha
}\left(  \left\vert x\right\vert ^{2k}\right)  D^{\beta-\alpha}\left(
\left\vert x\right\vert \right)  \right\vert \leq\\
& \leq\sum\limits_{\substack{\alpha\leq\beta\\\left\vert \alpha\right\vert
\leq2k}}\tbinom{\beta}{\alpha}\left\vert D^{\alpha}\left(  \left\vert
x\right\vert ^{2k}\right)  \right\vert \left\vert D^{\beta-\alpha}\left(
\left\vert x\right\vert \right)  \right\vert \\
& \leq\sum\limits_{\substack{\alpha\leq\beta\\\left\vert \alpha\right\vert
\leq2k}}\tbinom{\beta}{\alpha}\left(  2^{2k}\pi^{-d/2}d^{2k-\left\vert
\alpha\right\vert }\frac{\left(  k!\right)  ^{2}}{\left(  2k-\left\vert
\alpha\right\vert \right)  !}\right)  \left\vert x\right\vert ^{2k-\left\vert
\alpha\right\vert }c_{\left\vert \beta-\alpha\right\vert ,1}\left\vert
x\right\vert ^{1-\left\vert \beta-\alpha\right\vert }\\
& =2^{2k}\pi^{-d/2}d^{2k}\left(  k!\right)  ^{2}\left(  \sum
\limits_{\substack{\alpha\leq\beta\\\left\vert \alpha\right\vert \leq
2k}}\tbinom{\beta}{\alpha}\frac{d^{-\left\vert \alpha\right\vert }}{\left(
2k-\left\vert \alpha\right\vert \right)  !}c_{\left\vert \beta\right\vert
-\left\vert \alpha\right\vert ,1}\right)  \left\vert x\right\vert
^{2k+1-\left\vert \beta\right\vert }\\
& \leq2^{2k}\pi^{-d/2}d^{2k}\left(  k!\right)  ^{2}\left(  \max
_{\substack{\left\vert \alpha\right\vert \leq\left\vert \beta\right\vert
\\\left\vert \alpha\right\vert \leq2k}}c_{\left\vert \alpha\right\vert
,1}\right)  \left(  \sum\limits_{\substack{\alpha\leq\beta\\\left\vert
\alpha\right\vert \leq2k}}\tbinom{\beta}{\alpha}\frac{d^{-\left\vert
\alpha\right\vert }}{\left(  \left\vert \beta\right\vert -\left\vert
\alpha\right\vert \right)  !}\right)  \left\vert x\right\vert
^{2k+1-\left\vert \beta\right\vert }.
\end{align*}

But if $\left\vert \beta\right\vert \leq2k$ then using the identity: if
$\left\vert xy\right\vert <1$%
\[
\left(  1+xy\right)  ^{m}=\sum_{\left\vert \alpha\right\vert \leq m}\binom
{m}{\alpha}x^{\alpha}y^{\alpha},\text{\quad}\binom{m}{\alpha}=\frac{m!}%
{\alpha!\left(  m-\left\vert \alpha\right\vert \right)  !}.
\]

we obtain%
\begin{align*}
\sum\limits_{\alpha\leq\beta}\tbinom{\beta}{\alpha}\frac{d^{-\left\vert
\alpha\right\vert }}{\left(  \left\vert \beta\right\vert -\left\vert
\alpha\right\vert \right)  !}  & =\sum\limits_{\alpha\leq\beta}\frac
{1}{\left(  \beta-\alpha\right)  !}\frac{\beta!}{\alpha!\left(  \left\vert
\beta\right\vert -\left\vert \alpha\right\vert \right)  !}\left(  \frac{1}%
{d}\mathbf{1}\right)  ^{\alpha}\mathbf{1}^{\alpha}\\
& \leq\sum\limits_{\left\vert \alpha\right\vert \leq\left\vert \beta
\right\vert }\binom{\left\vert \beta\right\vert }{\alpha}\left(  \frac{1}%
{d}\mathbf{1}\right)  ^{\alpha}\mathbf{1}^{\alpha}\\
& =\left(  1+\left(  \frac{1}{d}\mathbf{1},\mathbf{1}\right)  \right)
^{\left\vert \beta\right\vert }\\
& =2^{\left\vert \beta\right\vert },
\end{align*}

we means that when $\left\vert \beta\right\vert \leq2k+1$,
\[
\left\vert D^{\beta}\left\vert x\right\vert ^{2k+1}\right\vert \leq
2^{2k+\left\vert \beta\right\vert }\pi^{-d/2}d^{2k}\left(  k!\right)
^{2}b_{\left\vert \beta\right\vert ,1}\left\vert x\right\vert
^{2k+1-\left\vert \beta\right\vert },\quad.
\]

\end{proof}

\section{Bounds for $\left(  aD\right)  ^{n}\left(  \left\vert x\right\vert
^{t}\right)  $ and $\left(  aD\right)  ^{n}\left(  \left\vert x\right\vert
^{t}\log\left\vert x\right\vert \right)  $}

??

\begin{lemma}
\label{Lem_op_aD_estim}Define the operator $aD=a_{1}D_{1}+\ldots+a_{d}D_{d}$
where $D=\left(  D_{i}\right)  $. Hence, noting identity \ref{p08},%
\[
\frac{1}{n!}\left(  \left(  aD\right)  ^{n}u\right)  \left(  x\right)
=\sum\limits_{\left\vert \alpha\right\vert =n}\frac{a^{\alpha}}{\alpha
!}D^{\alpha}u\left(  x\right)  ,
\]

and:

\begin{enumerate}
\item $\left(  aD\right)  \left\vert x\right\vert =ax\left\vert x\right\vert
^{-1}$,

\item $\left(  aD\right)  f\left(  \left\vert x\right\vert \right)
=f^{\prime}\left(  \left\vert x\right\vert \right)  \left(  aD\right)
\left\vert x\right\vert =a\widehat{x}f^{\prime}\left(  \left\vert x\right\vert
\right)  $,

\item $\left(  aD\right)  \left\vert x\right\vert ^{-k}=-k\left(  ax\right)
\left\vert x\right\vert ^{-k-2}$,

\item $\left(  aD\right)  ax=\left\vert a\right\vert ^{2}$,

\item $\left(  aD\right)  ^{2}\left\vert x\right\vert ^{-k}=k\left(
k+2\right)  \left(  ax\right)  ^{2}\left\vert x\right\vert ^{-\left(
k+4\right)  }-k\left\vert x\right\vert ^{-\left(  k+2\right)  }$.
\end{enumerate}

\textbf{CONJECTURES} for upper bounds:

\begin{enumerate}
\item For $m,n=0,1,2,\ldots$,
\[
\left\vert \left(  \widehat{a}D\right)  ^{m}\left\vert x\right\vert
^{n}\right\vert \leq k_{m,n}\left\vert x\right\vert ^{n-m},
\]
where the constants $k_{m,n}$ is given in Corollary \ref{Cor_Thm_(aD)^m_|x|^n}%
.\medskip

\item For $n\geq1$,%
\begin{equation}
\left\vert \left(  \widehat{a}D\right)  ^{2n}\left\vert x\right\vert
^{-s}\right\vert \leq c_{2n}\left(  -s\right)  \left\vert x\right\vert
^{-s-2n},\quad s\geq-2,\nonumber
\end{equation}

where%
\begin{equation}
c_{2n}\left(  -s\right)  :=\left\{
\begin{array}
[c]{ll}%
2^{2n}\left\vert \left(  \frac{s}{2}\right)  _{n}\right\vert \left(  \frac
{1}{2}\right)  _{n}, & -2\leq s\leq0,\\
2^{2n}\left(  \frac{s}{2}\right)  _{n}\left(  \frac{s+1}{2}\right)  _{n}, &
s\geq0.
\end{array}
\right. \label{Ap000}%
\end{equation}

\item For $n\geq0$,%
\[
\left\vert \left(  \widehat{a}D\right)  ^{2n+1}\left\vert x\right\vert
^{-s}\right\vert \leq c_{2n+1}\left(  -s\right)  \left\vert x\right\vert
^{-s-\left(  2n+1\right)  },\quad s\geq-4,
\]

where%
\begin{equation}
c_{2n+1}\left(  -s\right)  =\left\{
\begin{array}
[c]{ll}%
2^{2n+2}\left\vert \left(  \frac{s}{2}\right)  _{n+1}\right\vert \left(
\frac{1}{2}\right)  _{n+1}, & -4\leq s\leq2,\\
2^{2n+1}\left(  \frac{s}{2}\right)  _{n+1}\left(  \frac{s+1}{2}\right)
_{n}, & s\geq2.
\end{array}
\right. \label{Ap018}%
\end{equation}

\item For all $a\in\mathbb{R}^{d}$,%
\begin{align*}
\widehat{a}D\log\left\vert x\right\vert  & =\left(  \widehat{a}\widehat
{x}\right)  \left\vert x\right\vert ^{-1},\\
\left(  \widehat{a}D\right)  ^{m}\log\left\vert x\right\vert  & =\left(
\widehat{a}\widehat{x}\right)  \left\vert x\right\vert \left(  \widehat
{a}D\right)  ^{m-1}\left\vert x\right\vert ^{-2}+\left(  m-1\right)  \left(
\widehat{a}D\right)  ^{m-2}\left\vert x\right\vert ^{-2},\quad m\geq2,
\end{align*}

and we have the estimates%
\[
\frac{1}{m!}\left\vert \left(  \widehat{a}D\right)  ^{m}\log\left\vert
x\right\vert \right\vert \leq\left(  2-\frac{1}{m}\right)  \left\vert
x\right\vert ^{-m},\quad m\geq1.
\]

\item For all $t\in\mathbb{R}^{1}$ and $n\geq1$,%
\[
\frac{1}{n!}\left\vert \left(  \widehat{a}D\right)  ^{n}\left(  \left\vert
x\right\vert ^{t}\log\left\vert x\right\vert \right)  \right\vert \leq\frac
{1}{n!}\left\vert \left(  \widehat{a}D\right)  ^{n}\left\vert x\right\vert
^{t}\right\vert \left\vert \log\left\vert x\right\vert \right\vert
+\sum\limits_{j=0}^{n-1}\left(  2-\frac{1}{n-j}\right)  \frac{1}{j!}\left\vert
\left(  \widehat{a}D\right)  ^{j}\left\vert x\right\vert ^{t}\right\vert
\left\vert x\right\vert ^{-m}.
\]

\item For $n\geq1$ and $s\geq-2$,%
\[
\frac{1}{\left(  2n\right)  !}\left\vert \left(  \widehat{a}D\right)
^{2n}\left(  \left\vert x\right\vert ^{-s}\log\left\vert x\right\vert \right)
\right\vert \leq\left(  \frac{c_{2n}\left(  -s\right)  }{\left(  2n\right)
!}\left\vert \log\left\vert x\right\vert \right\vert +\sum\limits_{j=0}%
^{2n-1}\left(  2-\frac{1}{2n-j}\right)  \frac{c_{j}\left(  -s\right)  }%
{j!}\right)  \left\vert x\right\vert ^{-s-2n}.
\]

\item For $i=0,1,2,\ldots$ and $n\geq0$,%
\[
\frac{1}{n!}\left\vert \left(  \widehat{a}D\right)  ^{n}\left(  \left\vert
x\right\vert ^{i}\log\left\vert x\right\vert \right)  \right\vert \leq\left(
\frac{k_{n,i}}{n!}\left\vert \log\left\vert x\right\vert \right\vert
+\sum\limits_{j=0}^{n-1}\left(  2-\frac{1}{n-j}\right)  \frac{k_{j,i}}%
{j!}\right)  \left\vert x\right\vert ^{i-n}.
\]

where the constants $k_{j,i}$ are defined in Corollary
\ref{Cor_Thm_(aD)^m_|x|^n}. Further%
\[
\frac{1}{n!}\left\vert \left(  \widehat{a}D\right)  ^{n}\left(  \left\vert
x\right\vert ^{2i}\log\left\vert x\right\vert \right)  \right\vert \leq\left(
\sum\limits_{j=0}^{2i}\left(  2-\frac{1}{n-j}\right)  \frac{k_{j,2i}}%
{j!}\right)  \left\vert x\right\vert ^{2i-n},\quad n>2i.
\]

\end{enumerate}
\end{lemma}

\begin{proof}
\textbf{Parts 1 to 5}%
\begin{align*}
\left(  aD\right)  f\left(  \left\vert x\right\vert \right)  =\sum
\limits_{i=1}^{d}a_{i}D_{i}f\left(  \left\vert x\right\vert \right)
=\sum\limits_{i=1}^{d}a_{i}f^{\prime}\left(  \left\vert x\right\vert \right)
D_{i}\left\vert x\right\vert  &  =f^{\prime}\left(  \left\vert x\right\vert
\right)  \sum\limits_{i=1}^{d}a_{i}D_{i}\left\vert x\right\vert \\
&  =f^{\prime}\left(  \left\vert x\right\vert \right)  \left(  aD\right)
\left\vert x\right\vert .
\end{align*}

\begin{align*}
\left(  aD\right)  \left\vert x\right\vert  & =\sum\limits_{i=1}^{d}a_{i}%
D_{i}\left\vert x\right\vert =\left\vert x\right\vert ^{-1}\sum\limits_{i=1}%
^{d}a_{i}x_{i}=ax\left\vert x\right\vert ^{-1},\\
\left(  aD\right)  \left\vert x\right\vert ^{-1}  & =-\left\vert x\right\vert
^{-2}\left(  aD\right)  \left\vert x\right\vert =-ax\left\vert x\right\vert
^{-3},\\
\left(  aD\right)  \left\vert x\right\vert ^{-k}  & =-k\left\vert x\right\vert
^{-k-1}\left(  aD\right)  \left\vert x\right\vert =-kax\left\vert x\right\vert
^{-k-2},\\
\left(  aD\right)  ax  & =\sum\limits_{i=1}^{d}a_{i}D_{i}\left(  ax\right)
=\sum\limits_{i=1}^{d}a_{i}^{2}=\left\vert a\right\vert ^{2}.\\
\left(  aD\right)  ^{2}\left\vert x\right\vert  & =\left\vert x\right\vert
^{-1}-\left(  ax\right)  ^{2}\left\vert x\right\vert ^{-3}.\\
\left(  aD\right)  ^{3}\left\vert x\right\vert  & =-2ax\left\vert x\right\vert
^{-3}-ax\left\vert x\right\vert ^{-3}+3\left(  ax\right)  ^{3}\left\vert
x\right\vert ^{-5}\\
& =-3ax\left\vert x\right\vert ^{-3}+3\left(  ax\right)  ^{3}\left\vert
x\right\vert ^{-5}.\\
\left(  aD\right)  ^{4}\left\vert x\right\vert  & =-3\left\vert x\right\vert
^{-3}+9\left(  ax\right)  ^{2}\left\vert x\right\vert ^{-5}+9\left(
ax\right)  ^{2}\left\vert x\right\vert ^{-5}-15\left(  ax\right)
^{4}\left\vert x\right\vert ^{-7}\\
& =-3\left\vert x\right\vert ^{-3}+18\left(  ax\right)  ^{2}\left\vert
x\right\vert ^{-5}-15\left(  ax\right)  ^{4}\left\vert x\right\vert ^{-7}\\
\left(  aD\right)  ^{5}\left\vert x\right\vert  & =36ax\left\vert x\right\vert
^{-5}+9ax\left\vert x\right\vert ^{-5}-60\left(  ax\right)  ^{3}\left\vert
x\right\vert ^{-7}-90\left(  ax\right)  ^{3}\left\vert x\right\vert
^{-7}+105\left(  ax\right)  ^{5}\left\vert x\right\vert ^{-9}\\
& =45ax\left\vert x\right\vert ^{-5}-150\left(  ax\right)  ^{3}\left\vert
x\right\vert ^{-7}+105\left(  ax\right)  ^{5}\left\vert x\right\vert ^{-9},
\end{align*}

which strongly suggest the general equations of part 6.\medskip

\textbf{Part 6} Restatement of Corollary \ref{Cor_Thm_(aD)^m_|x|^n}.\medskip

\textbf{Part 7} This is inequality \ref{Ap055} which is a consequence of
Conjecture \ref{Conj_argmax_p_even}.\medskip

\textbf{Part 8} This is inequality \ref{Ap057} which is a consequence of
Conjecture \ref{Conj_argmax_p_odd}.\medskip

\textbf{Part 9 }When $k=1$,%
\begin{equation}
\widehat{a}D\log\left\vert x\right\vert =\left\vert x\right\vert ^{-1}%
\widehat{a}D\left\vert x\right\vert =\left(  \widehat{a}x\right)  \left\vert
x\right\vert ^{-2}=\left(  \widehat{a}\widehat{x}\right)  \left\vert
x\right\vert ^{-1},\label{p29}%
\end{equation}

and when $k\geq2$, since $\left(  \widehat{a}D\right)  \left(  \widehat
{a}x\right)  =1$,%
\begin{align*}
\left(  \widehat{a}D\right)  ^{k}\log\left\vert x\right\vert  & =\left(
\widehat{a}D\right)  ^{k-1}\widehat{a}D\log\left\vert x\right\vert =\left(
\widehat{a}D\right)  ^{k-1}\left(  \left(  \widehat{a}x\right)  \left\vert
x\right\vert ^{-2}\right)  =\\
& =\sum\limits_{j=0}^{k-1}\tbinom{k-1}{j}\left(  \widehat{a}D\right)
^{j}\left(  \widehat{a}x\right)  \left(  \widehat{a}D\right)  ^{k-1-j}%
\left\vert x\right\vert ^{-2}\\
& =\left(  \widehat{a}x\right)  \left(  \widehat{a}D\right)  ^{k-1}\left\vert
x\right\vert ^{-2}+\sum\limits_{j=1}\tbinom{k-1}{j}\left(  \widehat
{a}D\right)  ^{j}\left(  \widehat{a}x\right)  \left(  \widehat{a}D\right)
^{k-1-j}\left\vert x\right\vert ^{-2}\\
& =\left(  \widehat{a}x\right)  \left(  \widehat{a}D\right)  ^{km-1}\left\vert
x\right\vert ^{-2}+\tbinom{k-1}{1}\left(  \widehat{a}D\left(  \widehat
{a}x\right)  \right)  \left(  \widehat{a}D\right)  ^{k-2}\left\vert
x\right\vert ^{-2}\\
& =\left(  \widehat{a}\widehat{x}\right)  \left\vert x\right\vert \left(
\widehat{a}D\right)  ^{k-1}\left\vert x\right\vert ^{-2}+\left(  k-1\right)
\left(  \widehat{a}D\right)  ^{k-2}\left\vert x\right\vert ^{-2},
\end{align*}

i.e.%
\[
\left(  \widehat{a}D\right)  ^{k}\log\left\vert x\right\vert =\left(
\widehat{a}\widehat{x}\right)  \left\vert x\right\vert \left(  \widehat
{a}D\right)  ^{k-1}\left\vert x\right\vert ^{-2}+\left(  m-1\right)
\left\vert \left(  \widehat{a}D\right)  ^{k-2}\left\vert x\right\vert
^{-2}\right\vert ,\quad k\geq2.
\]

From part 7
\begin{equation}
\left\vert \left(  \widehat{a}D\right)  ^{2k}\left\vert x\right\vert
^{-2}\right\vert \leq c_{2k}\left(  -2\right)  \left\vert x\right\vert
^{-2k-2},\nonumber
\end{equation}

and so%
\[
\left\vert \left(  \widehat{a}D\right)  ^{2k-2}\left\vert x\right\vert
^{-2}\right\vert \leq c_{2k-2}\left(  -2\right)  \left\vert x\right\vert
^{-2-\left(  2k-2\right)  }=c_{2k-2}\left(  -2\right)  \left\vert x\right\vert
^{-2k},
\]

and from part 8%
\[
\left\vert \left(  \widehat{a}D\right)  ^{2k+1}\left\vert x\right\vert
^{-2}\right\vert \leq c_{2k+1}\left(  -2\right)  \left\vert x\right\vert
^{-2k-3},
\]

so that%
\[
\left\vert \left(  \widehat{a}D\right)  ^{2k-1}\left\vert x\right\vert
^{-2}\right\vert \leq c_{2k-1}\left(  -2\right)  \left\vert x\right\vert
^{-2k-1}.
\]

Hence%
\begin{align}
\left\vert \left(  \widehat{a}D\right)  ^{2k}\log\left\vert x\right\vert
\right\vert  & \leq\left\vert \left(  \widehat{a}D\right)  ^{2k-1}\left\vert
x\right\vert ^{-2}\right\vert +\left(  2k-1\right)  \left\vert \left(
\widehat{a}D\right)  ^{2k-2}\left\vert x\right\vert ^{-2}\right\vert
\nonumber\\
& \leq c_{2k-1}\left(  -2\right)  \left\vert \widehat{a}\widehat{x}\right\vert
\left\vert x\right\vert ^{-2k-1}+\left(  2k-1\right)  c_{2\left(  k-1\right)
}\left(  -2\right)  \left\vert x\right\vert ^{-2k}\nonumber\\
& =\left(  c_{2k-1}\left(  -2\right)  +\left(  2k-1\right)  c_{2k-2}\left(
-2\right)  \right)  \left\vert x\right\vert ^{-2k},\label{a004}%
\end{align}

and%
\begin{align}
\left\vert \left(  \widehat{a}D\right)  ^{2k+1}\log\left\vert x\right\vert
\right\vert  & \leq\left\vert \left(  \widehat{a}D\right)  ^{2k}\left\vert
x\right\vert ^{-2}\right\vert +2k\left\vert \left(  \widehat{a}D\right)
^{2k-1}\left\vert x\right\vert ^{-2}\right\vert \nonumber\\
& \leq c_{2k}\left(  -2\right)  \left\vert x\right\vert ^{-2k-1}%
+2kc_{2k-1}\left(  -2\right)  \left\vert x\right\vert ^{-2k-1}\nonumber\\
& =\left(  c_{2k}\left(  -2\right)  +2kc_{2k-1}\left(  -2\right)  \right)
\left\vert x\right\vert ^{-2k-1},\label{a003}%
\end{align}

so that%
\[
\left\vert \left(  \widehat{a}D\right)  ^{m}\log\left\vert x\right\vert
\right\vert \leq\left\{
\begin{array}
[c]{ll}%
\left\vert x\right\vert ^{-1}, & m=1,\\
\left(  c_{m-1}\left(  -2\right)  +\left(  m-1\right)  c_{m-2}\left(
-2\right)  \right)  \left\vert x\right\vert ^{-m}, & m\geq2,
\end{array}
\right.
\]

But
\[%
\begin{array}
[c]{ll}%
c_{2k}\left(  -2\right)  =2^{2k}\left(  1\right)  _{k}\left(  \frac{3}%
{2}\right)  _{k}=k!\frac{\left(  2k+1\right)  !}{k!}=\left(  2k+1\right)  !, &
k\geq1,\\
c_{2k+1}\left(  -2\right)  =2^{2k+1}\left(  1\right)  _{k+1}\left(  \frac
{3}{2}\right)  _{k}=\left(  2k+2\right)  !, & k\geq0,
\end{array}
\]

so that%
\begin{equation}
c_{m}\left(  -2\right)  =\left(  m+1\right)  !,\quad m\geq0,\label{a002}%
\end{equation}

and%
\begin{align*}
\left\vert \left(  \widehat{a}D\right)  ^{m}\log\left\vert x\right\vert
\right\vert  & \leq\left\{
\begin{array}
[c]{ll}%
\left\vert x\right\vert ^{-1}, & m=1,\\
\left(  m!+\left(  m-1\right)  \left(  m-1\right)  !\right)  \left\vert
x\right\vert ^{-m}, & m\geq2,
\end{array}
\right. \\
& =\left\{
\begin{array}
[c]{ll}%
\left\vert x\right\vert ^{-1}, & m=1,\\
\left(  2m-1\right)  \left(  m-1\right)  !\left\vert x\right\vert ^{-m}, &
m\geq2,
\end{array}
\right. \\
& =\left(  2-\frac{1}{m}\right)  m!\left\vert x\right\vert ^{-m},\text{ }%
m\geq1.
\end{align*}
\medskip

\textbf{Part 10} For $n\geq1$, all $t\in\mathbb{R}^{1}$ and $x\neq\mathbf{0}$,%
\begin{align*}
\left(  \widehat{a}D\right)  ^{n}\left(  \left\vert x\right\vert ^{t}%
\log\left\vert x\right\vert \right)   & =\sum\limits_{m=0}^{n}\tbinom{n}%
{m}\left(  \left(  \widehat{a}D\right)  ^{n-m}\left\vert x\right\vert
^{t}\right)  \left(  \widehat{a}D\right)  ^{m}\log\left\vert x\right\vert \\
& =\left(  \left(  \widehat{a}D\right)  ^{n}\left\vert x\right\vert
^{t}\right)  \log\left\vert x\right\vert +\sum\limits_{m=1}^{n}\tbinom{n}%
{m}\left(  \left(  \widehat{a}D\right)  ^{n-m}\left\vert x\right\vert
^{t}\right)  \left(  \widehat{a}D\right)  ^{m}\log\left\vert x\right\vert ,
\end{align*}

and from part 9%
\begin{align*}
\left\vert \left(  \widehat{a}D\right)  ^{n}\left(  \left\vert x\right\vert
^{t}\log\left\vert x\right\vert \right)  \right\vert  & \leq\left\vert \left(
\widehat{a}D\right)  ^{n}\left\vert x\right\vert ^{t}\right\vert \left\vert
\log\left\vert x\right\vert \right\vert +\sum\limits_{m=1}^{n}\tbinom{n}%
{m}\left\vert \left(  \widehat{a}D\right)  ^{n-m}\left\vert x\right\vert
^{t}\right\vert \left\vert \left(  \widehat{a}D\right)  ^{m}\log\left\vert
x\right\vert \right\vert \\
& \leq\left\vert \left(  \widehat{a}D\right)  ^{n}\left\vert x\right\vert
^{t}\right\vert \left\vert \log\left\vert x\right\vert \right\vert
+\sum\limits_{m=1}^{n}\tbinom{n}{m}\left\vert \left(  \widehat{a}D\right)
^{n-m}\left\vert x\right\vert ^{t}\right\vert \left(  2-\frac{1}{m}\right)
m!\left\vert x\right\vert ^{-m}\\
& =\left\vert \left(  \widehat{a}D\right)  ^{n}\left\vert x\right\vert
^{t}\right\vert \left\vert \log\left\vert x\right\vert \right\vert
+n!\sum\limits_{m=1}^{n}\frac{2-\frac{1}{m}}{\left(  n-m\right)  !}\left\vert
\left(  \widehat{a}D\right)  ^{n-m}\left\vert x\right\vert ^{t}\right\vert
\left\vert x\right\vert ^{-m}\\
& =\left\vert \left(  \widehat{a}D\right)  ^{n}\left\vert x\right\vert
^{t}\right\vert \left\vert \log\left\vert x\right\vert \right\vert
+n!\sum\limits_{j=0}^{n-1}\left(  2-\frac{1}{n-j}\right)  \frac{1}%
{j!}\left\vert \left(  \widehat{a}D\right)  ^{j}\left\vert x\right\vert
^{t}\right\vert \left\vert x\right\vert ^{-\left(  n-j\right)  },
\end{align*}

so that%
\[
\frac{1}{n!}\left\vert \left(  \widehat{a}D\right)  ^{n}\left(  \left\vert
x\right\vert ^{t}\log\left\vert x\right\vert \right)  \right\vert \leq\frac
{1}{n!}\left\vert \left(  \widehat{a}D\right)  ^{n}\left\vert x\right\vert
^{t}\right\vert \left\vert \log\left\vert x\right\vert \right\vert
+\sum\limits_{j=0}^{n-1}\left(  2-\frac{1}{n-j}\right)  \frac{1}{j!}\left\vert
\left(  \widehat{a}D\right)  ^{j}\left\vert x\right\vert ^{t}\right\vert
\left\vert x\right\vert ^{-\left(  n-j\right)  }.
\]
\smallskip

\textbf{Part 11} From part 10 and then parts 7 and 8, when $s\geq-2$,
\begin{align*}
&  \left\vert \left(  \widehat{a}D\right)  ^{2n}\left(  \left\vert
x\right\vert ^{-s}\log\left\vert x\right\vert \right)  \right\vert \\
&  \leq\left\vert \left(  \widehat{a}D\right)  ^{2n}\left\vert x\right\vert
^{-s}\right\vert \left\vert \log\left\vert x\right\vert \right\vert +\left(
2n\right)  !\sum\limits_{j=0}^{2n-1}\frac{2-\frac{1}{2n-j}}{j!}\left\vert
\left(  \widehat{a}D\right)  ^{j}\left\vert x\right\vert ^{-s}\right\vert
\left\vert x\right\vert ^{-\left(  2n-j\right)  }\\
&  \leq c_{2n}\left(  -s\right)  \left\vert x\right\vert ^{-s-2n}\left\vert
\log\left\vert x\right\vert \right\vert +\left(  2n\right)  !\sum
\limits_{j=0}^{2n-1}\frac{2-\frac{1}{2n-j}}{j!}c_{j}\left(  -s\right)
\left\vert x\right\vert ^{-s-j}\left\vert x\right\vert ^{-\left(  2n-j\right)
}\\
&  =c_{2n}\left(  -s\right)  \left\vert x\right\vert ^{-s-2n}\left\vert
\log\left\vert x\right\vert \right\vert +\left(  2n\right)  !\sum
\limits_{j=0}^{2n-1}\frac{2-\frac{1}{2n-j}}{j!}c_{j}\left(  -s\right)
\left\vert x\right\vert ^{-s-2n}\\
&  =\left(  c_{2n}\left(  -s\right)  \left\vert \log\left\vert x\right\vert
\right\vert +\left(  2n\right)  !\sum\limits_{j=0}^{2n-1}\frac{2-\frac
{1}{2n-j}}{j!}c_{j}\left(  -s\right)  \right)  \left\vert x\right\vert
^{-s-2n},
\end{align*}

the $c_{j}$ are given by \ref{Ap000} and \ref{Ap018}. Finally%
\[
\frac{1}{\left(  2n\right)  !}\left\vert \left(  \widehat{a}D\right)
^{2n}\left(  \left\vert x\right\vert ^{-s}\log\left\vert x\right\vert \right)
\right\vert \leq\left(  \frac{c_{2n}\left(  -s\right)  }{\left(  2n\right)
!}\left\vert \log\left\vert x\right\vert \right\vert +\sum\limits_{j=0}%
^{2n-1}\left(  2-\frac{1}{2n-j}\right)  \frac{c_{j}\left(  -s\right)  }%
{j!}\right)  \left\vert x\right\vert ^{-s-2n}%
\]
\medskip

\textbf{Part 12} Part 10 with $t=i$ and then using part 6 gives%
\begin{align*}
\frac{1}{n!}\left\vert \left(  \widehat{a}D\right)  ^{n}\left(  \left\vert
x\right\vert ^{i}\log\left\vert x\right\vert \right)  \right\vert  & \leq
\frac{1}{n!}\left\vert \left(  \widehat{a}D\right)  ^{n}\left\vert
x\right\vert ^{i}\right\vert \left\vert \log\left\vert x\right\vert
\right\vert +\sum\limits_{j=0}^{n-1}\left(  2-\frac{1}{n-j}\right)  \frac
{1}{j!}\left\vert \left(  \widehat{a}D\right)  ^{j}\left\vert x\right\vert
^{i}\right\vert \left\vert x\right\vert ^{-\left(  n-j\right)  }\\
& \leq\frac{1}{n!}k_{n,i}\left\vert x\right\vert ^{i-n}\left\vert
\log\left\vert x\right\vert \right\vert +\sum\limits_{j=0}^{n-1}\left(
2-\frac{1}{n-j}\right)  \frac{1}{j!}k_{j,i}\left\vert x\right\vert
^{i-j}\left\vert x\right\vert ^{-\left(  n-j\right)  }\\
& =\left(  \frac{k_{n,i}}{n!}\left\vert \log\left\vert x\right\vert
\right\vert +\sum\limits_{j=0}^{n-1}\left(  2-\frac{1}{n-j}\right)
\frac{k_{j,i}}{j!}\right)  \left\vert x\right\vert ^{i-n}.
\end{align*}

From part 6, $k_{n,2i}=0$ when $n>2i$ so%
\[
\frac{1}{n!}\left\vert \left(  \widehat{a}D\right)  ^{n}\left(  \left\vert
x\right\vert ^{2i}\log\left\vert x\right\vert \right)  \right\vert \leq\left(
\sum\limits_{j=0}^{n-1}\left(  2-\frac{1}{n-j}\right)  \frac{k_{j,2i}}%
{j!}\right)  \left\vert x\right\vert ^{2i-n}.
\]

\end{proof}

\section{Bounds for $\left\vert D\right\vert ^{2n}\left(  \left\vert
x\right\vert ^{t}\right)  $ and $\left\vert D\right\vert ^{2n}\left(
\left\vert x\right\vert ^{t}\log\left\vert x\right\vert \right)  $}

\begin{lemma}
\label{Lem_iterLaplacian_rad}The action of \textbf{iterated Laplacian}
$\left\vert D\right\vert ^{2n}$ on $\left\vert x\right\vert ^{k}$, $\left\vert
x\right\vert ^{-k}$ and $\left\vert x\right\vert ^{2}\log\left\vert
x\right\vert $:

\begin{enumerate}
\item
\begin{align*}
\left\vert D\right\vert ^{2}\left(  u\left(  \left\vert x\right\vert \right)
\right)   & =\left(  d-1\right)  \left\vert x\right\vert ^{-1}u^{\prime
}\left(  \left\vert x\right\vert \right)  +u^{\prime\prime}\left(  \left\vert
x\right\vert \right) \\
& =\left(  \left(  d-1\right)  s^{-1}u^{\prime}\left(  s\right)
+u^{\prime\prime}\left(  s\right)  \right)  \left(  s=\left\vert x\right\vert
\right)  .
\end{align*}

If%
\[
\mathcal{L}u:=\left(  d-1\right)  s^{-1}u^{\prime}+u^{\prime\prime},
\]

then
\[
\left\vert D\right\vert ^{2n}\left(  u\left(  \left\vert x\right\vert \right)
\right)  =\left(  \mathcal{L}^{n}u\right)  \left(  \left\vert x\right\vert
\right)  .
\]

These formulas will be used to prove:\smallskip

\item $\left\vert D\right\vert ^{2}\left(  \left\vert x\right\vert
^{2}\right)  =2d$.

\item $\left\vert D\right\vert ^{2}\log\left\vert x\right\vert =\left(
d-2\right)  \left\vert \cdot\right\vert ^{-2}$.

\item $\left\vert D\right\vert ^{2}\left(  \left\vert x\right\vert ^{2}%
\log\left\vert x\right\vert \right)  =d+2+2d\log\left\vert x\right\vert $.

\item For all $n$ and real $s$,%
\[
\left\vert D\right\vert ^{2n}\left(  \left\vert x\right\vert ^{-s}\right)
=2^{2n}\left(  \frac{s}{2}\right)  _{n}\left(  \frac{s}{2}-\frac{d}%
{2}+1\right)  _{n}\left\vert x\right\vert ^{-s-2n},
\]

and $\left\vert D\right\vert ^{2n}\left(  \left\vert x\right\vert ^{s}\right)
=0$ iff $s\in\left\{  0,2,4,\ldots,2n-2\right\}  \cup\left\{  0,2,4,\ldots
,2n-2\right\}  +2-d$.\medskip

\item $\left\vert D\right\vert ^{2n}\left(  \left\vert \cdot\right\vert
^{-k}\right)  =0$ iff $k\in\left\{  d-2,d-4,\ldots,d-2n\right\}  $.

\item If $0!:=0$ then%
\[
\left\vert D\right\vert ^{2n}\left(  \left\vert x\right\vert ^{2}%
\log\left\vert x\right\vert \right)  =\left\{
\begin{array}
[c]{ll}%
d+2+2d\log\left\vert x\right\vert , & n=1,\\
2d\left(  d-2\right)  \left\vert x\right\vert ^{-2}, & n=2,\\
2^{2n-1}\left(  n-2\right)  !\left(  -\frac{d}{2}\right)  _{n}\left\vert
x\right\vert ^{2-2n}, & n\geq2.
\end{array}
\right.
\]

\item For integer $n,k\geq0$, using the generalized binomial coefficient
notation%
\[
\left\vert D\right\vert ^{2n}\left(  \left\vert x\right\vert ^{2k}\right)
=\left\{
\begin{array}
[c]{ll}%
2^{2n}\left(  n!\right)  ^{2}\tbinom{k}{n}\tbinom{k+d/2-1}{n}\left\vert
x\right\vert ^{2k-2n} & n\leq k,\\
0, & n>k.
\end{array}
\right.
\]

\item For all integer $n,k\geq1$,%
\[
\left\vert D\right\vert ^{2n}\left(  \left\vert x\right\vert ^{2k-1}\right)
=??2^{2n}\left(  -k+\frac{1}{2}\right)  _{n}\left(  -k-\frac{d}{2}+\frac{3}%
{2}\right)  _{n}\left\vert x\right\vert ^{2k-1-2n}.
\]

\end{enumerate}
\end{lemma}

\begin{proof}
\textbf{Part 1} $D_{i}\left(  u\left(  \left\vert x\right\vert \right)
\right)  =u^{\prime}\left(  \left\vert x\right\vert \right)  D_{i}\left\vert
x\right\vert =\left\vert x\right\vert ^{-1}u^{\prime}\left(  \left\vert
x\right\vert \right)  x_{i}$ and so
\[
D_{i}^{2}\left(  u\left(  \left\vert x\right\vert \right)  \right)
=\left\vert x\right\vert ^{-1}u^{\prime}\left(  \left\vert x\right\vert
\right)  +\left\vert x\right\vert ^{-2}u^{\prime\prime}\left(  \left\vert
x\right\vert \right)  x_{i}^{2}-\left\vert x\right\vert ^{-3}u^{\prime}\left(
\left\vert x\right\vert \right)  x_{i}^{2}.
\]
\smallskip

Hence $\left\vert D_{i}\right\vert ^{2}\left(  u\left(  \left\vert
x\right\vert \right)  \right)  =d\left\vert x\right\vert ^{-1}u^{\prime
}\left(  \left\vert x\right\vert \right)  +u^{\prime\prime}\left(  \left\vert
x\right\vert \right)  -\left\vert x\right\vert ^{-1}u^{\prime}\left(
\left\vert x\right\vert \right)  =\left(  d-1\right)  \left\vert x\right\vert
^{-1}u^{\prime}\left(  \left\vert x\right\vert \right)  +u^{\prime\prime
}\left(  \left\vert x\right\vert \right)  $.\smallskip

\textbf{Part 2} $u=s^{2}$ so $\mathcal{L}u=\left(  d-1\right)  s^{-1}%
u^{\prime}\left(  s\right)  +u^{\prime\prime}\left(  s\right)  =2\left(
d-1\right)  +2=2d $.\smallskip

\textbf{Part 3} $u=\log s$ so $\mathcal{L}u=\left(  d-1\right)  s^{-1}%
u^{\prime}\left(  s\right)  +u^{\prime\prime}\left(  s\right)  =\left(
d-1\right)  s^{-2}-s^{-2}=\left(  d-2\right)  s^{-2}$.\smallskip

\textbf{Part 4 }$u=s^{2}\log s$ so\textbf{\ }$u^{\prime}=2s\log s+s$ and
$u^{\prime\prime}=2\log s+2+1=3+2\log s$ so
\begin{align*}
\mathcal{L}u  & =\left(  d-1\right)  s^{-1}u^{\prime}\left(  s\right)
+u^{\prime\prime}\left(  s\right)  =\left(  d-1\right)  s^{-1}\left(  s+2s\log
s\right)  +\left(  3+2\log s\right) \\
& =\left(  d-1\right)  \left(  1+2\log s\right)  +\left(  3+2\log s\right)
=d+2+2d\log s.
\end{align*}
\smallskip

\textbf{Part 5 }For convenience we start by assuming that $u\left(  t\right)
=t^{s}$ then change $s$ to $-s$.%
\begin{align*}
\mathcal{L}u=\left(  d-1\right)  t^{-1}u^{\prime}+u^{\prime\prime}  & =\left(
d-1\right)  t^{-1}\left(  st^{s-1}\right)  +s\left(  s-1\right)  t^{s-2}\\
& =\left(  \left(  d-1\right)  s+s^{2}-s\right)  t^{s-2}\\
& =s\left(  s+d-2\right)  t^{s-2},
\end{align*}

so that when $x\neq0$,
\begin{align}
&  \left\vert D\right\vert ^{2n}\left(  \left\vert x\right\vert ^{s}\right)
\nonumber\\
&  =\left(  \mathcal{L}^{n}u\right)  \left(  \left\vert x\right\vert \right)
\nonumber\\
&  =s\left(  s+d-2\right)  \left(  s-2\right)  \left(  s+d-4\right)
\ldots\left(  s-2n+2\right)  \left(  s+d-2n\right)  \left\vert x\right\vert
^{s-2n}\nonumber\\
&  =\left\{  s\left(  s-2\right)  \ldots\left(  s-2n+2\right)  \right\}
\left\{  \left(  s+d-2\right)  \left(  s+d-4\right)  \ldots\left(
s+d-2n\right)  \right\}  \left\vert x\right\vert ^{s-2n},\label{Ap136}%
\end{align}

for all $s$.

Thus if $s\in\left\{  0,2,4,\ldots,2n-2\right\}  \cup\left\{  2-d,4-d,\ldots
,2n-d\right\}  $ then $\left\vert D\right\vert ^{2n}\left(  \left\vert
x\right\vert ^{s}\right)  =0$.

Further%
\begin{align*}
\left\vert D\right\vert ^{2n}\left(  \left\vert x\right\vert ^{-s}\right)   &
=\left\{  -s\left(  -s-2\right)  \ldots\left(  -s-2n+2\right)  \right\}
\times\\
& \quad\times\left\{  \left(  -s+d-2\right)  \left(  -s+d-4\right)
\ldots\left(  -s+d-2n\right)  \right\}  \left\vert x\right\vert ^{-s-2n}\\
& =\left\{  s\left(  s+2\right)  \ldots\left(  s+2\left(  n-1\right)  \right)
\right\}  \left\{  \left(  s-d+2\right)  \left(  s-d+4\right)  \ldots\left(
s-d+2n\right)  \right\}  \left\vert x\right\vert ^{-s-2n}\\
& =2^{2n}\left\{  \frac{s}{2}\left(  \frac{s}{2}+1\right)  \ldots\left(
\frac{s}{2}+n-1\right)  \right\}  \times\\
& \quad\times\left\{  \left(  \frac{s-d+2}{2}\right)  \left(  \frac{s-d+2}%
{2}+1\right)  \ldots\left(  \frac{s-d+2}{2}+n-1\right)  \right\}  \left\vert
x\right\vert ^{-s-2n}\\
& =2^{2n}\left(  \frac{s}{2}\right)  _{n}\left(  \frac{s-d+2}{2}\right)
_{n}\left\vert x\right\vert ^{-s-2n}.
\end{align*}
\medskip

\textbf{Part 6} From part 5
\[%
\begin{array}
[c]{lll}%
\left\vert D\right\vert ^{2n}\left(  \left\vert \cdot\right\vert ^{-k}\right)
=0 & iff & k-d+2j=0\text{ }for\text{ }some\text{ }integer\text{ }1\leq j\leq
n\\
& iff & j=\frac{d-k}{2}\text{ }for\text{ }some\text{ }integer\text{ }1\leq
j\leq n\\
& iff & \frac{d-k}{2}\in\left\{  1,2,\ldots,n\right\} \\
& iff & k\in\left\{  d-2,d-4,\ldots,d-2n\right\}
\end{array}
\]
\medskip

\textbf{Part 7} Using part 4 and then part 3 we get
\begin{align*}
\left\vert D\right\vert ^{2n}\left(  \left\vert x\right\vert ^{2}%
\log\left\vert x\right\vert \right)   & =\left\vert D\right\vert
^{2n-2}\left\vert D\right\vert ^{2}\left(  \left\vert x\right\vert ^{2}%
\log\left\vert x\right\vert \right)  =\left\vert D\right\vert ^{2n-2}\left(
d+2+2d\log\left\vert x\right\vert \right)  =\\
& =2d\left\vert D\right\vert ^{2n-2}\log\left\vert x\right\vert \\
& =2d\left\vert D\right\vert ^{2n-4}\left\vert D\right\vert ^{2}\log\left\vert
x\right\vert \\
& =2d\left\vert D\right\vert ^{2\left(  n-2\right)  }\left(  \left(
d-2\right)  \left\vert x\right\vert ^{-2}\right) \\
& =2d\left(  d-2\right)  \left\vert D\right\vert ^{2\left(  n-2\right)
}\left(  \left\vert x\right\vert ^{-2}\right)  ,
\end{align*}

and from part 5, when $n\geq3$,%
\begin{align*}
\left\vert D\right\vert ^{2m}\left(  \left\vert x\right\vert ^{-k}\right)   &
=2^{2m}\left(  k/2\right)  _{m}\left(  k/2-d/2+1\right)  _{m}\left\vert
x\right\vert ^{-\left(  k+2m\right)  },\\
\left\vert D\right\vert ^{2m}\left(  \left\vert x\right\vert ^{-2}\right)   &
=2^{2m}m!\left(  2-d/2\right)  _{m}\left\vert x\right\vert ^{-2-2m},\\
\left\vert D\right\vert ^{2\left(  n-2\right)  }\left(  \left\vert
x\right\vert ^{-2}\right)   & =2^{2\left(  n-2\right)  }\left(  n-2\right)
!\left(  2-d/2\right)  _{n-2}\left\vert x\right\vert ^{2-2n},
\end{align*}

so that when $n\geq3$,%
\begin{align*}
\left\vert D\right\vert ^{2n}\left(  \left\vert x\right\vert ^{2}%
\log\left\vert x\right\vert \right)   & =d\left(  d-2\right)  2^{2n-3}\left(
n-2\right)  !\left(  2-\frac{d}{2}\right)  _{n-2}\left\vert x\right\vert
^{2-2n}\\
& =2^{2}\frac{d}{2}\left(  1-\frac{d}{2}\right)  2^{2n-3}\left(  n-2\right)
!\left(  2-\frac{d}{2}\right)  \left(  3-\frac{d}{2}\right)  \ldots\left(
n-1-\frac{d}{2}\right)  \left\vert x\right\vert ^{2-2n}\\
& =2^{2n-1}\left(  n-2\right)  !\left(  -\frac{d}{2}\right)  _{n}\left\vert
x\right\vert ^{2-2n}.
\end{align*}
\medskip

\textbf{Parts 8} From \ref{Ap136},%
\begin{align*}
\left\vert D\right\vert ^{2n}\left(  \left\vert x\right\vert ^{2m}\right)   &
=\left\{  2m\left(  2m-2\right)  \left(  2m-2n+2\right)  \right\}  \times\\
& \quad\times\left\{  \left(  2m+d-2\right)  \left(  2m+d-4\right)
\ldots\left(  2m+d-2n\right)  \right\}  \left\vert x\right\vert ^{2m-2n}\\
& =2^{2n}\left\{  m\left(  m-1\right)  \ldots\left(  m-n+1\right)  \right\}
\times\\
& \quad\times\left\{  \left(  m+\frac{d}{2}-1\right)  \left(  m+\frac{d}%
{2}-2\right)  \ldots\left(  m+\frac{d}{2}-n\right)  \right\}  \left\vert
x\right\vert ^{2m-2n}\\
& =2^{2n}\left(  n!\right)  ^{2}\tbinom{m}{n}\tbinom{m+d/2-1}{n}\left\vert
x\right\vert ^{2m-2n}.
\end{align*}
\medskip

\textbf{Parts 9} ??
\end{proof}

\section{Bounds for $\left(  \widehat{\cdot}D\right)  ^{n}\left(  \left\vert
x\right\vert ^{t}\log\left\vert x\right\vert \right)  $ and $\left(
\widehat{\cdot}D\right)  ^{n}\left(  \left\vert x\right\vert ^{t}\right)  $}

\begin{lemma}
\label{Lem_deriv_rad_funcs}The action of the operator $\widehat{\cdot}D$ on
$\left\vert \cdot\right\vert ^{t}$ and $\left\vert \cdot\right\vert ^{t}%
\log\left\vert \cdot\right\vert $:

The operators $\cdot D$ and $\widehat{\cdot}D$ will be defined by%
\begin{align*}
\left(  \left(  \cdot D\right)  u\right)  \left(  x\right)   & =\sum
\limits_{k=1}^{d}x_{k}D_{k}u\left(  x\right)  ,\\
\left(  \left(  \widehat{\cdot}D\right)  u\right)  \left(  x\right)   &
=\sum\limits_{k=1}^{d}\widehat{x}_{k}D_{k}u\left(  x\right)  ,
\end{align*}

where $\widehat{x}=x/\left\vert x\right\vert $. Hence%
\begin{align*}
\left(  \left(  \cdot D\right)  ^{n}u\right)  \left(  x\right)   &
=\sum\limits_{k=1}^{d}\frac{x^{\alpha}}{\alpha!}D^{\alpha}u\left(  x\right)
,\\
\left(  \left(  \widehat{\cdot}D\right)  ^{n}u\right)  \left(  x\right)   &
=\sum\limits_{k=1}^{d}\frac{\left(  \widehat{x}\right)  ^{\alpha}}{\alpha
!}D^{\alpha}u\left(  x\right)  .
\end{align*}

\begin{enumerate}
\item If $\psi\left(  x\right)  =\psi_{\circ}\left(  \left\vert x\right\vert
\right)  $ then $\left(  \left(  \widehat{\cdot}D\right)  \psi\right)  \left(
x\right)  =\left(  D\psi_{\circ}\right)  \left(  \left\vert x\right\vert
\right)  $ and
\begin{equation}
\left(  \left(  \widehat{\cdot}D\right)  ^{n}\psi\right)  \left(  x\right)
=\left(  D^{n}\psi_{\circ}\right)  \left(  \left\vert x\right\vert \right)
.\label{Ap131}%
\end{equation}

This result allows us to use the following 1-dimensional results:\medskip

\item For $k=1,2,3,\ldots$,%
\[
D^{n}s^{k}=\left\{
\begin{array}
[c]{ll}%
\left(  -1\right)  ^{n}\left(  -k\right)  _{n}s^{k-n}, & n\leq k,\\
0, & n>k.
\end{array}
\right.
\]

\item For any $t\in\mathbb{R}^{1}$,%
\[
D^{n}s^{-t}=\left(  -1\right)  ^{n}\left(  t\right)  _{n}s^{-t-n},\quad
n\geq0.
\]

\item If $n$ and $k$ are non-negative integers:%
\begin{align*}
D^{n}\left(  s^{k}\phi\left(  s\right)  \right)   & =\left\{
\begin{array}
[c]{ll}%
k!\sum\limits_{j=0}^{n}\tbinom{n}{j}\frac{s^{k-n+j}}{\left(  k-\left(
n-j\right)  \right)  !}D^{j}\phi\left(  s\right)  , & 1\leq n\leq k,\\
k!\sum\limits_{j=n-k}^{n}\tbinom{n}{j}\frac{s^{k-n+j}}{\left(  k-n+j\right)
!}D^{j}\phi\left(  s\right)  , & n>k,
\end{array}
\right. \\
& =k!\sum\limits_{j=\max\left\{  0,n-k\right\}  }^{n}\tbinom{n}{j}%
\frac{s^{k-n+j}}{\left(  k-n+j\right)  !}D^{j}\phi\left(  s\right)  .
\end{align*}

\item
\[
D^{n}\left(  s^{t}\phi\left(  s\right)  \right)  =\sum\limits_{j=0}^{n}\left(
-1\right)  ^{j}\tbinom{n}{j}\left(  -t\right)  _{j}s^{t-j}D^{n-j}\phi\left(
s\right)  .
\]

\item
\[
D^{j}\log s=\left(  -1\right)  ^{j-1}\left(  j-1\right)  !s^{-j},\quad j\geq1,
\]

and%
\[
D^{n}\left(  s^{t}\log s\right)  =\left(  -1\right)  ^{n}\left(  \left(
-t\right)  _{n}\log s-\sum\limits_{l=1}^{n}\frac{\left(  -1\right)  ^{l}}%
{l}\left(  -n\right)  _{l}\left(  -t\right)  _{n-l}\right)  s^{t-n},\quad
t\in\mathbb{R}^{1}.
\]

\item
\[
D^{n}\left(  s^{k}\log s\right)  =\left\{
\begin{array}
[c]{ll}%
n!\left(  \tbinom{k}{n}\log s+\sum\limits_{j=1}^{n}\frac{\left(  -1\right)
^{j+1}}{j}\tbinom{k}{n-j}\right)  s^{k-n}, & 1\leq n\leq k,\\
\left(  -1\right)  ^{n-k+1}k!\left(  n-k-1\right)  !s^{k-n}, & n>k.
\end{array}
\right.
\]

\end{enumerate}
\end{lemma}

\begin{proof}
\textbf{Part 1}%
\begin{align*}
\left(  \widehat{\cdot}D\right)  \psi_{\circ}\left(  \left\vert x\right\vert
\right)   & =\sum\limits_{i=1}^{d}\widehat{x}_{i}D_{i}\psi_{\circ}\left(
\left\vert x\right\vert \right)  =\sum\limits_{i=1}^{d}\widehat{x}_{i}\left(
D\psi_{\circ}\right)  \left(  \left\vert x\right\vert \right)  D_{i}\left\vert
x\right\vert =\\
& =\sum\limits_{i=1}^{d}\widehat{x}_{i}^{\prime}\psi_{\circ}\left(  \left\vert
x\right\vert \right)  \widehat{x}_{i}=\psi_{\circ}^{\prime}\left(  \left\vert
x\right\vert \right)  \sum\limits_{i=1}^{d}\left(  \widehat{x}_{i}\right)
^{2}=\psi_{\circ}^{\prime}\left(  \left\vert x\right\vert \right)  .
\end{align*}
\medskip

\textbf{Part 2} Expand $D^{n}\left(  s^{k}\psi_{\circ}\left(  s\right)
\right)  $ by differentiating using Leibniz's theorem.\medskip

\textbf{Part 3}.%
\begin{align*}
Ds^{-t}  & =-ks^{-\left(  t+1\right)  },\\
D^{2}s^{-t}  & =t\left(  t+1\right)  s^{-\left(  t+2\right)  },\\
& and\text{ }in\text{ }general\\
D^{n}s^{-t}  & =\left(  -1\right)  ^{n}t\left(  t+1\right)  \ldots\left(
t+n-1\right)  s^{-\left(  t+2\right)  }\\
& =\left(  -1\right)  ^{n}\left(  t\right)  _{n}s^{-t-n}\\
& =\left(  -1\right)  ^{n}n!\tbinom{t+n-1}{n}s^{-t-n}.
\end{align*}
\medskip

\textbf{Part 4} ??\medskip

\textbf{Part 5} ??\medskip

\textbf{Part 6} Examples: $D\log s=s^{-1}$, $D^{2}\log s=-s^{-2}$, $D^{3}\log
s=2!s^{-3}$ etc. so
\[
D^{j}\log s=\left(  -1\right)  ^{j+1}\left(  j-1\right)  !s^{-j}.
\]

Thus%
\begin{align*}
D^{n}\left(  s^{t}\log s\right)   & =\sum\limits_{j=0}^{n}\left(  -1\right)
^{j}\tbinom{n}{j}\left(  -t\right)  _{j}s^{t-j}D^{n-j}\log s\\
& =\left(  -1\right)  ^{n}\left(  -t\right)  _{n}s^{t-n}\log s+\\
& \quad+\sum\limits_{j=0}^{n-1}\left(  -1\right)  ^{j}\tbinom{n}{j}\left(
-t\right)  _{j}s^{t-j}\left(  -1\right)  ^{n-j-1}\left(  n-j-1\right)
!s^{j-n}\\
& =\left(  -1\right)  ^{n}\left(  -t\right)  _{n}s^{t-n}\log s-\left(
-1\right)  ^{n}\sum\limits_{j=0}^{n-1}\tbinom{n}{j}\left(  -t\right)
_{j}\left(  n-j-1\right)  !s^{t-n}\\
& =\left(  -1\right)  ^{n}\left(  \left(  -t\right)  _{n}\log s-\sum
\limits_{j=0}^{n-1}\tbinom{n}{j}\left(  -t\right)  _{j}\left(  n-j-1\right)
!\right)  s^{t-n}\\
& =\left(  -1\right)  ^{n}\left(  \left(  -t\right)  _{n}\log s-\sum
\limits_{j=0}^{n-1}\frac{1}{n-j}\frac{n!}{j!}\left(  -t\right)  _{j}\right)
s^{t-n}\\
& =\left(  -1\right)  ^{n}\left(  \left(  -t\right)  _{n}\log s-\sum
\limits_{l=0}^{n-1}\frac{1}{l+1}\frac{n!}{\left(  n-l-1\right)  !}\left(
-t\right)  _{n-l-1}\right)  s^{t-n}\\
& =\left(  -1\right)  ^{n}\left(  \left(  -t\right)  _{n}\log s-\sum
\limits_{l=0}^{n-1}\frac{\left(  -1\right)  ^{l+1}}{l+1}\left(  -n\right)
_{l+1}\left(  -t\right)  _{n-l-1}\right)  s^{t-n}\\
& =\left(  -1\right)  ^{n}\left(  \left(  -t\right)  _{n}\log s-\sum
\limits_{l=1}^{n}\frac{\left(  -1\right)  ^{l}}{l}\left(  -n\right)
_{l}\left(  -t\right)  _{n-l}\right)  s^{t-n}.
\end{align*}
\medskip

\textbf{Part 7} ?? But from part 4,%
\[
D^{n}\left(  s^{k}\log s\right)  =k!\sum\limits_{j=0}^{n}\tbinom{n}{j}%
\frac{s^{k-n+j}}{\left(  k-n+j\right)  !}D^{j}\log s,\quad1\leq n\leq k,
\]

so%
\begin{align*}
D^{n}\left(  s^{k}\log s\right)   & =k!\sum\limits_{j=0}^{n}\tbinom{n}{j}%
\frac{s^{k-n+j}}{\left(  k-n+j\right)  !}D^{j}\log s\\
& =k!\left(  \frac{s^{k-n}}{\left(  k-n\right)  !}\log s+\sum\limits_{j=1}%
^{n}\tbinom{n}{j}\frac{s^{k-n+j}}{\left(  k-n+j\right)  !}\left(  -1\right)
^{j+1}\left(  j-1\right)  !s^{-j}\right) \\
& =k!s^{k-n}\left(  \frac{1}{\left(  k-n\right)  !}\log s+\sum\limits_{j=1}%
^{n}\left(  -1\right)  ^{j+1}\tbinom{n}{j}\frac{\left(  j-1\right)  !}{\left(
k-\left(  n-j\right)  \right)  !}\right) \\
& =s^{k-n}\left(  \frac{k!}{\left(  k-n\right)  !}\log s+\sum\limits_{j=1}%
^{n}\left(  -1\right)  ^{j+1}\frac{n!}{\left(  n-j\right)  !j!}\frac{k!\left(
j-1\right)  !}{\left(  k-\left(  n-j\right)  \right)  !}\right)  ,
\end{align*}

so that%
\begin{align*}
\frac{1}{n!}D^{n}\left(  s^{k}\log s\right)   & =s^{k-n}\left(  \tbinom{k}%
{n}\log s+\sum\limits_{j=1}^{n}\frac{\left(  -1\right)  ^{j+1}}{j}\frac
{1}{\left(  n-j\right)  !}\frac{k!}{\left(  k-\left(  n-j\right)  \right)
!}\right) \\
& =s^{k-n}\left(  \tbinom{k}{n}\log s+\sum\limits_{j=1}^{n}\frac{\left(
-1\right)  ^{j+1}}{j}\tbinom{k}{n-j}\right)  .
\end{align*}

On the other hand, if $n>k$,%
\begin{align*}
\frac{1}{n!}D^{n}\left(  s^{k}\log s\right)   & =\frac{1}{n!}D^{n-k}%
D^{k}\left(  s^{k}\log s\right) \\
& =\frac{1}{n!}D^{n-k}s^{k-k}\left(  \tbinom{k}{k}\log s+\sum\limits_{j=1}%
^{k}\frac{\left(  -1\right)  ^{j+1}}{j}\tbinom{k}{k-j}\right) \\
& =\frac{1}{n!}D^{n-k}\log s\\
& =\left(  -1\right)  ^{n-k+1}\frac{\left(  n-k-1\right)  !}{n!}s^{-\left(
n-k\right)  }.
\end{align*}

\end{proof}

\bibliographystyle{amsplain}
\bibliography{acompat,BasisFunc}

\newif\ifabfull\abfulltrue
\providecommand{\bysame}{\leavevmode\hbox to3em{\hrulefill}\thinspace}
\providecommand{\MR}{\relax\ifhmode\unskip\space\fi MR }
\providecommand{\MRhref}[2]{%
  \href{http://www.ams.org/mathscinet-getitem?mr=#1}{#2}
}
\providecommand{\href}[2]{#2}
\begin{thebibliography}{10}

\bibitem{PrudBryMar86}
Yu.~Brychkov A.~Prudnikov and O.~Marichevl, \emph{Integrals and series}, vol.
  1: Elementary Functions, Gordon And Breach Science Publishers, New York,
  1986.

\bibitem{AbramowStegun70}
M.~Abramowitz and I.~Stegun, \emph{Handbook of mathematical functions}, Dover
  Publications, 1970.

\bibitem{Adams75}
R.~A. Adams, \emph{Sobolev spaces}, 1st ed., Pure and applied mathematics
  series, vol.~65, Academic Press, London, 1975.

\bibitem{Archbold64}
J.~W. Archbold, \emph{Algebra}, 3rd ed., Pitman, Melbourne, 1964.

\bibitem{Arfken70}
G.~Arfken, \emph{Mathematical methods for physicists}, 2nd ed., Academic Press,
  New York, London, 1970.

\bibitem{DenyLions54}
J.~Deny and J.~L. Lions, \emph{Les espaces du type de beppo levi}, Ann. Inst.
  Fourier \textbf{5} (1954), 305--370 (French).

\bibitem{Dwight61}
H.~B. Dwight, \emph{Tables of integrals and other mathematical data}, 4th ed.,
  Macmillan, 1961.

\bibitem{Dyn89}
N.~Dyn, \emph{Interpolation and approximation by radial and related functions},
  Approximation Theory VI (C.~Chui., L.~Schumaker., and J.~Ward, eds.), Acad.
  Press, 1989, pp.~211--234.

\bibitem{GelfandVilenkin64}
I.~M. Gelfand. and N.~Ya. Vilenkin, \emph{Generalised functions}, Acad. Press,
  New York, 1964.

\bibitem{Hardy71}
R.~L. Hardy, \emph{Multiquadric equations of topography and other irregular
  surfaces}, J. Geophysical Res. \textbf{76} (1971), 1905--1915.

\bibitem{Hardy90}
\bysame, \emph{Theory and applications of the multiquadric-biharmonic method},
  Comput. Math. Appl. \textbf{19} (1990), 163--208.

\bibitem{LightWayne95ErrEst}
W.~Light and H.~Wayne, \emph{Error estimates for approximation by radial basis
  functions.}, Approximation theory, wavelets and applications (Maratea, 1994),
  NATO Adv. Sci. Inst. Ser. C Math. Phys. Sci., vol. 454, Kluwer Acad. Publ.,
  Dordrecht, 1995, pp.~215--246.

\bibitem{LightWayne98PowFunc}
\bysame, \emph{On power functions and error estimates for radial basis function
  interpolation}, J. Approx. Th. \textbf{92} (1998), no.~2, 245--266.

\bibitem{LightWayneX98Weight}
\bysame, \emph{Spaces of distributions, interpolation by translates of a basis
  function and error estimates}, Numer. Math. \textbf{81} (1999), no.~3,
  415--450.

\bibitem{MadychNelson88}
W.~Madych and S.~Nelson, \emph{Multivariate interpolation and conditionally
  positive definite functions}, J. Approx. Th. Appl. \textbf{4} (1988), 77--89.

\bibitem{MadychNelson90}
\bysame, \emph{Multivariate interpolation and conditionally positive definite
  functions ii}, Math. Comp. \textbf{54} (1990), 211--230.

\bibitem{Malliavin95}
P.~Malliavin, H.~Airault, L.~Kay, and G.~Letac, \emph{Integration and
  probability}, Springer-Verlag, Berlin, 1995.

\bibitem{Micchelli86}
C.~Micchelli, \emph{Interpolation of scattered data: Distance matrices and
  conditionally positive definite functions}, Const. Approx. \textbf{2} (1986),
  11--22.

\bibitem{NarcWardWend2004}
F.~Narcowich, J.~Ward, and H.~Wendland, \emph{Sobolev bounds on functions with
  scattered zeros, with applications to radial basis function surface fitting},
  Math. Comput. \textbf{74} (2004), no.~250, 743--763.

\bibitem{Petersen83}
B.~Petersen, \emph{Introduction to the fourier transform and
  pseudo-differential operators}, Monographs and studies in mathematics,
  vol.~19, Pitman Advanced Pub. Program, Boston, 1983.

\bibitem{ReedSimon72}
M.~Reed and B.~Simon, \emph{Methods of modern mathematical physics i:
  Functional analysis}, Academic Press, New York, 1972.

\bibitem{Schaback95}
R.~Schaback, \emph{Multivariate interpolation and approximation by translates
  of a basis function}, Approximation Theory VIII, Vol I; Approximation and
  interpolation (C.K. Chui and Larry Schumaker, eds.), World Scientific
  Publishing Co., 1995, pp.~491--514.

\bibitem{Schaback99}
\bysame, \emph{Native {H}ilbert spaces for radial basis functions i}, Int. Ser.
  Num. Math., vol. 132, pp.~255--282, Birkhauser-Verlag, 1999.

\bibitem{Schoenberg38}
I.~J. Schoenberg, \emph{Metric spaces and completely monotonic functions}, Ann.
  of Math. \textbf{39} (1938), 811--841.

\bibitem{Schwartz66}
L.~Schwartz, \emph{Theorie des distributions}, Hermann, Paris, 1966.

\bibitem{SteinWeiss71}
E.~M. Stein and G.~Weiss, \emph{Fourier analysis and {E}uclidean spaces},
  Princeton University Press, Princeton, 1971.

\bibitem{Vladimirov}
V.~Vladimirov, \emph{Equations of mathematical physics}, Pure and Applied
  Mathematics series, vol.~3, Marcel Dekker, New York, 1971.

\bibitem{Wendland05}
H.~Wendland, \emph{Scattered data approximation}, Cambridge Monographs on
  Computational and Applied Mathematics, vol.~17, Cambridge Univ. Press,
  Cambridge, 2005.

\bibitem{WilliamsZeroOrdSmthV4}
P.~Y. Williams, \emph{A weight function theory of zero order basis function
  interpolants and smoothers}, unpublished arXiv:0708.0780v4 [math.NA], 2014.

\bibitem{Yosida58}
K.~Yosida, \emph{Functional analysis}, Springer-Verlag, Berlin, 1958.

\bibitem{Yosida68}
\bysame, \emph{Functional analysis}, Springer, Berlin, Heidelberg, New York,
  1968.

\bibitem{WuZongSchab93}
W.~Zong-min and R.~Schaback, \emph{Local error estimates for radial basis
  function interpolation to scattered data}, J. Numer. Anal. \textbf{13}
  (1993), 13--27.

\end{thebibliography}

\end{document}